\PassOptionsToPackage{bookmarks=true,colorlinks=true,pdfstartview=FitV,linkcolor=blue, citecolor=red,urlcolor=green}{hyperref}
\documentclass[12pt, A4]{memo-l}
\usepackage{geometry}

\usepackage{amsmath}
\usepackage{amssymb}
\usepackage{mathptmx}
\usepackage{helvet}
\usepackage{courier}
\usepackage{mathrsfs} 
\usepackage{dsfont} 
\usepackage{mathtools}
\usepackage{bbm}
\usepackage{type1cm}     
\usepackage{graphicx}        
\usepackage{multicol}        
\usepackage[bottom]{footmisc}

\usepackage[numbered]{bookmark}

\usepackage{epstopdf}
\usepackage{bm}

\usepackage{lpic}

\usepackage{subfigure}

\usepackage{syntonly}
\usepackage{xr}
\usepackage{setspace}

\usepackage{framed}

 \usepackage{comment}

\newtheorem{theorem}{Theorem}[section]
\newtheorem{lemma}[theorem]{Lemma}
\newtheorem{proposition}[theorem]{Proposition} 
\newtheorem{corollary}[theorem]{Corollary} 

\theoremstyle{definition}
\newtheorem{definition}[theorem]{Definition}
\newtheorem{example}[theorem]{Example}

\newtheorem{hyp}[theorem]{Hypothesis}
\newtheorem{conjecture}[theorem]{Conjecture}

\theoremstyle{remark}
\newtheorem{remark}[theorem]{Remark}

\numberwithin{section}{chapter}
\numberwithin{equation}{chapter}

\numberwithin{section}{chapter}
\numberwithin{equation}{section}

\newcommand\bR{{\mathbb{R}}}
\newcommand\bC{{\mathbb C}}
\newcommand\bZ{{\mathbb Z}}
\newcommand\bQ{{\mathbb Q}}

\newcommand\bV{{\mathbf{V}}}
\newcommand\bW{{\mathbf{W}}}

\newcommand\bA{{\mathbf{A}}}
\newcommand\bU{{\mathbb{U}}}
\newcommand\bP{{\mathbb{P}}}
\newcommand\Uu{{\mathbb{U}}}
\newcommand\bB{{\mathbb{B}}}
\newcommand\bF{{\mathbb{F}}}
\newcommand\bG{{\mathbb{G}}}
\newcommand\bX{{\mathbb{X}}}
\newcommand\boldx{{\mathbf{x}}}

\newcommand\Hom{{\rm Hom}}
\newcommand\dev{{\bf dev}}
\newcommand\SI{{\mathbb{S}}}

\newcommand\Bd{{\rm bd}}
\newcommand\clo{{\rm Cl}}
\newcommand\bdd{{\mathbf{d}}}

\newcommand\ra{\rightarrow}

\newcommand\che{\check}
\newcommand\emp{\emptyset}
\newcommand\eps{\epsilon}

\newcommand\Aff{{\mathbf{Aff}}}
\newcommand\ovl{\overline}
\newcommand\Aut{{\mathbf{Aut}}}
\newcommand\Idd{{\rm I}}

\newcommand\CP{{\mathbb{CP}}}
\newcommand\tr{{\mathrm{tr}}}

\newcommand\CN{{\mathcal{N}}}

\newcommand\bv{{\vec{v}}}
\newcommand\bu{{\vec{u}}}
\newcommand{\bx}{\vec{x}}
\newcommand{\bw}{\vec{w}}

\newcommand\tri{\triangle}

\newcommand\rp{\mathbb{RP}}

\newcommand\rpn{\mathbb{RP}^n}

\newcommand\rpno{\mathbb{RP}^{n-1}}
\newcommand\RP{\mathbb{RP}}

\newcommand\PO{{\mathsf{PO}}}

\newcommand\rep{\mathrm{rep}}

\newcommand\SL{{\mathsf{SL}}}

\newcommand\PGL{{\mathsf{PGL}}}
\newcommand\PSL{{\mathsf{PSL}}}
\newcommand\SO{{\mathsf{SO}}}
\newcommand\PSO{\mathsf{PSO}}
\newcommand\Ort{{\mathsf{O}}}
\newcommand\SLnp{{\mathsf{SL}}_\pm(n+1, \bR)}
\newcommand\SLn{{\mathsf{SL}}_\pm(n, \bR)}
\newcommand\SLf{{\mathsf{SL}}_\pm(4, \bR)}

\newcommand\GL{{\mathsf{GL}}}
\newcommand\GLnp{{\mathsf{GL}}(n+1, \bR)}
\newcommand\PGLnp{{\mathsf{PGL}}(n+1, \bR)}

\newcommand\Def{\mathrm{Def}}
\newcommand\CDef{\mathrm{CDef}}

\newcommand\SDef{\mathrm{SDef}}
\newcommand\hol{\mathrm{hol}}
\newcommand\orb{\mathcal{O}} 
\newcommand\torb{\tilde{\mathcal{O}}}
\newcommand\bGamma{{\bm\Gamma}}
\newcommand\leng{{\mathrm{length}}}

\newcommand\mx{{\mathrm{max}}}
\newcommand\cwl{{\mathrm{cwl}}}

\newcommand\supp{{\mathrm{supp}}}

\newcommand\cR{{\mathcal{R}}}
\newcommand\cT{{\mathcal{T}}}
\newcommand\SLpm{{\mathrm{SL}}_{\pm}(n+1, \bR)}
\newcommand\Ag{{\mathrm{Ag}}}
\newcommand\rlp{\,\rm(}
\newcommand\rrp{\,\rm)}

\newcommand{\llrrbracket}[1]{
	\left[\mkern-2mu\left[#1\right]\mkern-2mu\right]}
\newcommand{\llrrparen}[1]{
	\left(\mkern-4mu\left(#1\right)\mkern-4mu\right)}

\newcommand{\llrrV}[1]{
	\left|\mkern-2mu\left|#1\right|\mkern-2mu\right|}

\DeclareMathOperator{\rank}{rank}

\newcommand*{\Scale}[2][4]{\scalebox{#1}{$#2$}}%

\DeclareMathOperator{\Ima}{Im}

\newcommand{\orar}{\overrightarrow}

\newcommand{\CH}{\mathcal{CH}}

\newcommand{\romu}{\mathrm{u}}
\newcommand{\mbv}{\mathrm{v}}
\newcommand{\lh}{\mathrm{lh}}

\newcommand{\SSn}{[$\SI^n$S]}
\newcommand{\SnT}{[$\SI^n$T]}
\newcommand{\SnP}{[$\SI^n$P]}
\newcommand{\cpr}{[$\mathbf{cpr}$]}

\makeindex             

\begin{document}
\pagenumbering{arabic}

\frontmatter
\title{Real projective orbifolds with ends and their deformation theory: 	\protect\\  The deformation theory for nicest ones}

\author{Suhyoung Choi}
\address{Department of Mathematical Sciences, KAIST, Daejeon, South Korea}
\curraddr{}
\email{schoi@math.kaist.ac.kr} 
\thanks{This work was supported by the National Research Foundation	of Korea (NRF) grant funded by the Korea government (MEST) (No.2010-0027001).} 



\date{\today}

\maketitle

	

	
	\tableofcontents

\chapter*{Preface}




 


Let $G$ be a Lie group acting transitively and effectively on a manifold $X$. 
An $(X, G)$-geometry is given by this pair. Furthermore, 
an $(X, G)$-structure on an orbifold or a manifold is an atlas of charts to $X$ 
with transition maps in $G$. Here, we are concerned with 
$G = \PGL(n+1, \bR)$ and $X = \RP^n$. 

Cartan, Ehresmann, and others started the field of $(X, G)$-structures. 
Subjects of $(X, G)$-structures were popularized by Thurston and Goldman among many 
other people. 
These structures provide alternative viewpoints for understanding representations and their deformations, beyond just algebraic ones.
Our deformation spaces often parameterize significant parts of the spaces of representations. 

Since the examples are easier to construct, even now, we will be studying orbifolds, 
a natural generalization of manifolds. Also, computations can be done 
fairly well for simple examples. 
We began our study with Coxeter orbifolds where the computations 
are probably the simplest possible. 

Thurston did use the theory of orbifolds in a deep way. 
The hyperbolization of Haken $3$-manifold requires 
the use of the deformation theory of orbifolds where we build from hyperbolic structures 
from handlebodies with ``scalloped'' orbifold structures.
(See Morgan \cite{Morgan79}.) We do not yet know how to bypass this step, 
which was a very subtle point that some experts misunderstood.  
Also, orbifolds are natural objects obtained when we take quotients of 
manifolds by fibrations and so on. 
These are some of the reasons we study orbifolds instead of only manifolds. 

Classically, conformally flat structures were studied much by differential geometers.  Projectively flat structures were also studied from Cartan's time. However, our techniques differ 
significantly from their approaches. 

Convex real projective orbifolds are quotient spaces of convex domains on 
a projective space $\RP^n$ by a discrete group of projective automorphisms.
Hyperbolic manifolds and many symmetric manifolds are natural examples. 
These can be deformed into ones not coming from simple constructions. 
The study was initiated by Kuiper \cite{Kuiper53}, Koszul \cite{Koszul68}, Benz\'ecri \cite{Benzecri60}, Vey \cite{Vey}, and Vinberg \cite{Vinberg71},
accumulating some class of results.
Closed manifolds or orbifolds admit many such structures as shown first by Kac-Vinberg \cite{Vinberg67}, 
Goldman \cite{Goldman90}, and Cooper-Long-Thistlethwaites \cite{CLT06}, \cite{CLT07}. Some parts of the theory for closed orbifolds 
were completed by Benoist \cite{Benoist01} in the 1990s. 

The topics of convex real projective structures on manifolds and orbifolds are currently developing. We present some parts. 
This book is mainly written for researchers in this field. 

There are surprisingly many such structures coming from hyperbolic ones and deforming
as shown by Vinberg for Coxeter orbifolds, Goldman for surfaces, and later by Cooper-Long-Thistlethwaite for $3$-manifolds. 

Most of these theories focuses on closed manifolds, and we concentrate on 
open manifolds. The Smale-Hirsch immersion theory \cite{Smale59}, \cite{Hirsch59}, 
and the h-principle of Gromov \cite{Gromov86} more generally tells us that each open manifold admits 
$(X, G)$-structure for any $(X, G)$-geometry provided obvious topological obstructions on
the manifold vanish. So we need some type of boundary 
conditions. 

We compare these theories 
to the Mostow or Margulis-type rigidity for symmetric spaces. 
The rigidity can be replaced by what is called 
the {\em Ehresmann-Thurston-Weil principle} that 
\begin{itemize} 
\item a subspace of the $G$-character space (variety) of the 
fundamental group of a manifold or orbifold $M$ classifies the $(X,G)$-structures
on $M$ under the map 
\[ \hol: \Def_c(M) \ra \Hom_c(\pi_1(M),G)/G \]
where 
\begin{itemize} 
\item we define the deformation space
\[\Def_c(M):= \{ \hbox{$(X, G)$-structures on $M$ satisfying some boundary conditions denoted by $c$} \}/\sim\] where 
$\sim$ is the isotopy equivalence relation, and 
\item $\Hom_c(\pi_1(M),G)/G$ is the subspace of 
the character space $\Hom(\pi_1(M),G)/G$ satisfying
the corresponding conditions to $c$. 
\end{itemize} 
\end{itemize} 

For closed real projective orbifolds, Benoist's work \cite{Benoist01}
is widely regarded as quite an encompassing one. 
Hence, we won't prove much about this topic. (See Choi-Lee-Marquis for a survey \cite{CLM18}.)

We focus on tame convex real projective orbifolds with ends, where some topologists have now 
accumulated some classes of examples. 
Basically, we will prove an Ehresmann-Thurston-Weil principle: 
We will show that the deformation space of properly convex real projective 
structures on an orbifold with some end conditions is identified under a map with 
the union of components of the subset of character spaces of the orbifold satisfying 
the corresponding conditions on end holonomy groups. 
Our conditions on the ends are probably very generic ones, and 
we have many examples of such deformations. 

In fact, we are focusing on generic cases of lens-type or horospherical ends. 
To complete the picture, we need to consider all types of ends. Even with end vertex conditions, 
we are still to complete the picture leaving out the NPNC-ends. We hope to allow these types for 
our deformation spaces in the near future. 

The book is divided into three parts: 
\begin{description} 
	\item[(Part I)] We give some introduction and survey some recent results.  
	\item[(Part II)] We classify the types of ends we will work with. 
	We use the uniform middle-eigenvalue conditions. The condition is used to prove the 
the existence of convex domains called lens. 
	\item[(Part III)] We aim to prove the Ehresmann-Thurston-Weil principle for the deformation spaces for our type of orbifolds. We show the local homeomorphism property and 
	the closedness of the images for the maps from the deformations spaces to 
	the character spaces restricted by the end conditions. 
Finally, we give some examples where our theory is applicable. 
	\end{description} 
We adhere to the strictly logical progression of the material. 
However, for each chapter, we will introduce the main results first. 

As an application, 
we will use the results of the whole of the monograph at 
Chapter \ref{ex-sec-nicecase}, which are the nice cases.   
One can consider these as the conclusions of the monograph.

We describe what are still to develop in these theories so far to 
understand the convexity of the orbifolds still: 
\begin{itemize} 
\item One should be able to generalize the radial or totally geodesic boundary conditions to 
more broader classes, perhaps by adding ideal boundary stratas of various codimensions. 
\item Lens type or weak middle-eigenvalue conditions might be generalized or even 
dropped for all these types of ends. 
\end{itemize}
Many parts of these were carried out under the condition that the end fundamental groups are amenable by Ballas, Cooper, Leitner, Long, and so on as describe in the monograph. 
General philosophy here is that the convexity of some end neighborhoods will 
be equivalent to the convexity of the whole orbifold. 
Tillmann and I talked about this in Melbourne near 2009 when I was in the sabbatical in 
Melbourne. 
However, the conditions are sometimes not compatible. 

As a further motivation for our study, we say about some long-term goal: 
Deforming a real projective structure on an orbifold to an unbounded situation results in the actions of the fundamental group on affine buildings. This will lead us to some understanding of orbifolds and manifolds in particular of dimension three
as indicated by Cooper. 

The logical dependence of the monograph 
is as ordered by the order of appearance. 
Appendix \ref{ch-abelian} depends only on Chapter \ref{ch-prelim},
and the results are used in the monograph except for Chapter 
\ref{ch-prelim}.  

We give an outline at the beginning of each part.

We need to lift the objects to $\SI^n$ using Section \ref{prelim-sub-lifting}. 
We give proofs in the book by considering objects to be 
in $\SI^n$ and using the projective automorphism group 
$\SLpm$. 
We will use proof symbols: 
\begin{description}
	\item[[$\SI^n$S]] at the end of the proof 
	indicates that it is sufficient to prove for $\SI^n$
	since the conclusion does not involve $\RP^n$ nor $\SI^n$.
	\item[[$\SI^n$T]] indicates that the version of the theorem, proposition, or lemma 
	for $\SI^n$ implies one for $\RP^n$
	often with the help of Theorem \ref{prelim-thm-lifting}
and Proposition \ref{prelim-prop-closureind}. 
	\item[[$\SI^n$P]]   indicates that the proof of the theorem for 
	$\SI^n$ implies one for $\RP^n$
		often with the help of Theorem \ref{prelim-thm-lifting}
and Proposition \ref{prelim-prop-closureind}. 
\item[\cpr] indicates that the proof did not end here but will continue.
\end{description}

If we do not need to go to $\SI^n$ to prove the result, we leave 
no mark except for the end of the proof.


There is a concurrent work by the group consisting of 
Cooper, Long, and Tillmann with Ballas and Leitner 
on the same subjects but with different conditions on ends. They impose
the condition that the end fundamental groups to be amenable. 
However, we do not require the same conditions in this paper but 
instead we will require geometric conditions and 
use some type of norms of eigenvalue conditions
to guarantee the convexity during the deformations. 
We note that their deformation spaces are somewhat differently defined. 
Of course, we benefited much from their work and insights 
in this book and are very grateful for their generous help and guidance.

This book generalizes and simplifies the earlier preprints of the author. 
We were able to drop many conditions in the earlier versions of the theorems
overcoming many limitations. 
Some of the results were announced in some survey articles \cite{End1} and 
\cite{convsurv}. 

Also, there are some related later developments generalizing the convex compact real projective orbifolds, termed convex cocompact actions of real projective groups as developed 
by Danciger, Gu\'eritaud, Kassel, Lee, and Marquis \cite {DGK21},\cite{DGKLM}. These are generalizations of some of the work here with possibly similar techniques. 
Of course, there are growing numbers of other later related works by other authors, which we mentioned in various parts of this monograph. 

The project itself was very old starting from Melbourne in 2008. However, there were too many difficulties with many ideas and organizing these ideas well. 
Finally, to better communicate the ideas, 
the author made some effort to make the material clearer and more precise, 
entailing the trade-off of the writing being long, somewhat technical, and sometimes redundant.
Finally, we compiled an extensive indices, which the readers can use to find definitions and terminologies.

\vspace{\baselineskip}
\begin{flushright}\noindent
\hfill {\it Daejeon,}\\
September, 2025\hfill {\it  Suhyoung Choi}\\
\end{flushright}




\chapter*{Acknowledgments}

We also appreciate much help from Crampon and Marquis working also 
independently of the above group and us. We also thank Cooper for reading earlier 
 versions and pointing out a few mistakes, which were very helpful. 

Of course we benefit much from interactions and communications with Samuel Ballas, Yves Benoist, Jeffrey Danciger, David Fried, William Goldman, Gye-Seon Lee, Ludovic Marquis, and so on.  
We partly worked with many people including Gye-Seon Lee, Stephen Tillmann, Yves Carri\`ere, David Fried, William Goldman, and Daryl Cooper in a chronological order. However, 
these people all wished for me to be the sole author of this work.  We would have never come up with some of the ideas in this monograph without their generous help and support. 

We also thank the graduate students  who were in my course during Spring of 2024 going over 
the materials in the book and checking the text and so on. They were very helpful in catching many errors and mistakes. We thank in particular Juseop Lee,  Seungyeol Park, and Yongho Seo
for reading the monograph rather thoroughly for mistakes and so on. 

This work was supported by the National Research Foundation
of Korea (NRF) grant funded by the Korea government (MEST) (No.2010-0027001).



\mainmatter

%
%
%
\part{Introduction to orbifolds and real projective structures.}

Part I aims to survey some preliminary definitions and elementary facts used. These are standard materials, and there are no new result.

In Chapter \ref{ch-prelim}, we go over basic preliminary materials. 
We begin with defining geometric structures and real projective structures, 
in particular convex ones. 
Affine orbifolds and affine suspensions of real projective orbifolds are defined.  
We discuss the linear algebra and estimations derived from it,
orthopotent actions of Lie groups, 
higher-convergence group actions, attracting and repelling sets, 
 convexity, 
the Benoist theory on convex dividing actions, 
and so on. 
We discuss the dual orbifolds of  given convex real projective orbifolds as given by Vinberg. Finally, we extend the concept of duality to all convex compact 
sets and discuss the geometric limits of the dual convex sets. 
Here we employ a slightly more generalized version of convexity.

In Chapter \ref{ch-ex}, we provide some examples, where our full theory applies.
We will explain this in Chapter \ref{ex-sec-nicecase}.
Coxeter orbifolds and the orderability theory for Coxeter orbifolds using 
the Vinberg theory will be explained. 
We discuss the work jointly done with Gye-Seon Lee, Hodgson, and Greene. 
We state the work of Heusner-Porti on projective deformations of 
hyperbolic link complements. 
Also, we present results by Cooper-Long-Tillmann and Crampon-Marquis on finite-volume convex real projective structures and their admission of thick and thin decompositions."



\chapter{Preliminaries} \label{ch-prelim}


We survey the underlying theory.
	In Section \ref{prelim-sec-prelim}, we discuss the Hausdorff convergences of sequences of compact sets, Hilbert metrics, some orbifold topology, geometric structures on orbifolds, real projective structures on orbifolds, spherical real projective structures and liftings. 
	We also classify compact convex subsets of $\SI^n$ in
	Proposition \ref{prelim-prop-classconv}. 
	In Section \ref{prelim-sec-affinorb},  we discuss  affine structures and affine suspensions of real projective orbifolds. 
In Section \ref{prelim-sec-linearalg},  
	we discuss the linear algebra and estimation to find convergences, orthopotent groups, proximal and semiproximal actions, 
	 semisimplicity, and the higher convergence groups. 
	 Higher convergence groups are generalizations of convergence groups. 
	 In Section \ref{prelim-sec-conv},
we explain the comprehensive Benoist theory 
on properly convex real projective orbifolds, where he 
completed theories of Kuiper, 
Koszul, 
Vey, 
Vinberg, 
and so on,  
on dividing actions on convex linear cones as he terms them.  
In particular, the strict-join decomposition of properly convex real projective orbifolds will be explained. 
Lemma \ref{prelim-lem-domainIn} shows 
that the properly convex real projective structures are uniquely determined by holonomy groups, which is a somewhat commonly overlooked fact. 
In Section \ref{prelim-sec-duality}, we explain the duality theory of Vinberg. We introduce the augmented boundary of properly convex domains 
as the set of boundary points and the sharply supporting hyperspaces
associated with these points. The duality map is extended to the 
augmented boundary. Duality is extended to sweeping actions also. 
 The duality is extended to every compact convex set in $\SI^n$, 
and we discuss the relationship between the duality and the geometric convergences of the sequences of properly convex sets.

\section{Preliminary definitions} \label{prelim-sec-prelim}

Given a vector space $\bV$, we define $\bP(\bV)$ as the quotient space of \index{p@$\bP(\cdot)$|textbf}
$\bV -\{O\}$ under the equivalence relation
\[\bv\sim \bw \hbox{ for } \bv, \bw \in \bV  -\{O\}  \hbox{ iff } \bv = s \bw, \hbox{ for } s \in \bR -\{0\}.\] 
We denote as $[\bv]$ the equivalence class of $\bv \in \bV  -\{O\}$. 
For a subspace $\bW$ of $\bV$, we denote by $\bP(\bW)$ the image of $\bW-\{O\}$ under the quotient map, also said to be a {\em subspace}. 
\index{subspace} 
\label{Pparanthesis} 

As usual, we denote by $\rpn$ the projectivization of $\bR^{n+1}$. 
There is a group $\PGL(n+1, \bR)$ acting effectively and transitively on $\rpn$.  \label{npr}

Given a vector space $V$, we denote by $\SI(\bV)$ the quotient space
\[ (\bV - \{O\})/\sim \hbox{ where } v \sim w \hbox{ iff } v = s w \hbox{ for } s > 0. \]
We denote by $\SI^n := \SI(\bR^{n+1})$. 
We will represent each element of $\PGL(n+1, \bR)$ by a matrix of determinant $\pm 1$; 
i.e., $\PGL(n+1, \bR) = \SLpm/\langle\pm \Idd \rangle$. 
Recall the covering map $p_{\SI^n}: \SI^{n}= \SI(\bR^{n+1}) \ra \rpn$. 
\label{sphere} 
 \index{sn@$\SI^n$}
\index{sphericalization|textbf} 
 \index{S@$\SI$|textbf}


We let $\llrrparen{\vec{v}}$ denote the equivalence class of $\vec{v} \in \bR^{n+1} -\{O\}$. 
Given a vector subspace  $V \in \bR^{n+1}$, we denote by $\SI(V)$ the image of $V -\{O\}$ under the quotient map. \index{s@$\SI(\cdot)$}
The image is called a {\em subspace}. A set of antipodal points forms a subspace of dimension $0$. 
\label{Sparanthesis} 
There is a double covering map $p_{\SI^n}:\SI^n \ra \RP^n$ with the deck transformation group generated by
$\mathcal{A}$. 
This mapping induces a projective structure on $\SI^n$. 
The group of projective automorphisms is identified with $\SL_\pm(n+1, \bR)$. \index{SLplusm@$\SL_\pm(n+1, \bR)$}
The notion of geodesics is defined as in the projective geometry: they correspond 
to arcs in great circles in $\SI^n$. 
\index{geodesic} 

A collection of subspaces $\SI(V_1), \dots, \SI(V_m)$
(resp. $\bP(V_1), \dots, \bP(V_m)$) is {\em independent} 
if the subspaces $V_1, \dots, V_m$ are independent. 
\index{independent subspaces|textbf} 




The group $\SL_{\pm}(n+1, \bR)$ of linear transformations of determinant $\pm 1$ 
maps to the projective group $\PGL(n+1, \bR)$ by a double covering homomorphism $\hat q$, 
and $\SL_{\pm}(n+1, \bR)$ acts on $\SI^n$ lifting the projective transformations. 
The elements are also {\em projective transformations}. 

\begin{remark}\label{intro-rem-SL}
For each $g \in \PGL(n+1, \bR)$ acting on a convex open domain $\Omega$, 
there is a unique lift in $\SLpm$ preserving each component of the inverse image of $\Omega$
under $\SI^{n} \ra \rpn$.  We will use this representative. 
\end{remark}

A cone $C$ in $\bR^{n+1} -\{O\}$ is a subspace such that given a vector $x \in C$, 
$s x \in C$ for every $s \in \bR_+$. 
A {\em convex cone} is a cone that is a convex subset of $\bR^{n+1}$ in the usual sense. 
A {\em properly convex cone} is a convex cone not containing a complete affine line. 
\index{cone|textbf} 
\index{cone!convex|textbf} 

\begin{itemize}
	\item Given a vector $\vec{v} \in \bR^{n+1} -\{O\} $, we denote by $[\vec{v}] \in \RP^n$ the equivalence class. 
	Let $\Pi: \bR^{n+1} -\{O\} \ra \RP^n$ denote the projection. 
\item Given a connected subset $A$ of an affine subspace of 
	$\RP^n$, a cone $C(A) \subset \bR^{n+1} $ of $A$ is given 
	as a connected cone in $\bR^{n+1}$ mapped onto $A$ under the projection $\Pi: \bR^{n+1} -\{O\} \ra \RP^n$.
	\item $C(A)$ is unique up to the antipodal map $\mathcal{A}: \bR^{n+1} \ra \bR^{n+1}$ given by $\vec{v} \ra -\vec{v}$. 
\end{itemize}
\index{rpn@$\RP^n$} 
\index{C@$C(\cdot)$}
\index{A@$\mathcal{A}$}

\label{pi} 
The general linear group $\GLnp$ acts on $\bR^{n+1}$, and $\PGLnp$ acts on $\rpn$ effectively. 
Denote by $\bR_+ =\{ r \in \bR| r > 0\}$.
The {\em real projective sphere} $\SI^n$ is defined as the quotient of $\bR^{n+1} -\{O\}$ under the quotient relation 
$\vec{v} \sim \vec{w}$ iff $\vec{v} = s\vec{w}$ for $s \in \bR_+$. 
We will also use $\SI^n$ as the double cover of $\rpn$. 
Then $\Aut(\SI^n)$, isomorphic to the subgroup $\SLnp$ of $\GLnp$ of 
determinant $\pm 1$, double-covers $\PGLnp$. 
$\Aut(\SI^{n})$ acts as a group of projective automorphisms of $\SI^n$. 
A {\em projective map} of a real projective orbifold to another is a map that is projective by charts to $\rpn$. 
Let $\Pi: \bR^{n+1}-\{O\} \ra \RP^n$ be a projection and let $\Pi':  \bR^{n+1}-\{O\} \ra \SI^n$ denote one for $\SI^n$. 
\index{pi@$\Pi$|textbf} 
\index{pi'@$\Pi'$|textbf}
An infinite subgroup $\Gamma$ of $\PGLnp$ (resp. $\SLnp$) is {\em strongly irreducible} if every finite-index subgroup is irreducible. 
A {\em subspace} $S$ of $\RP^n$ (resp. $\SI^n$) is the image of a subspace with the origin removed under the projection $\Pi$ (resp. $\Pi'$).
\index{subspace|textbf}
\index{projective map|textbf}
\index{projective automorphism|textbf}
\index{irreducible!strongly|textbf}
\index{strongly irreducible|textbf}
\label{piprime} 


A line in $\rpn$ or $\SI^n$ is an embedded arc in a $1$-dimensional subspace. 
A {\em projective geodesic} is an arc in a projective orbifold developing into a line in $\rpn$
or to a one-dimensional subspace of $\SI^n$. 
A {\em great segment} is an embedded geodesic connecting a pair of 
antipodal points in $\SI^n$ or the complement of a point in a $1$-dimensional 
subspace in $\RP^n$. Sometimes an open great segment 
is called a {\em complete affine line}. 
An affine space $\mathds{A}^n$ can be identified with the complement of a codimension-one subspace 
$\rpno$ such that the geodesic structures are same up to parametrizations. 
A {\em convex subset} of $\rpn$ is a convex subset of an affine subspace in this paper. 
A {\em properly convex subset} of  $\rpn$ is a precompact convex subset of an affine subspace. 
$\bR^n$ is identified with an open half-space in $\SI^n$ defined by a linear function on $\bR^{n+1}$. 
(In this paper an affine subspace is either embedded in $\rpn$ or $\SI^n$.)
\index{projective geodesic|textbf} 
\index{great segment|textbf} 
\index{complete affine line|textbf} 
\index{convex subset|textbf}
\index{convex subset!properly|textbf}
\index{affine space|textbf} 
\index{An@$\mathds{A}^n$|textbf} 
\index{affine subspace|textbf}
\label{affine space}

An {\em $i$-dimensional complete affine subspace} is 
a subspace of a projective $n$-orbifold projectively diffeomorphic to 
an $i$-dimensional affine subspace in some affine subspace 
$\mathds{A}^n$ of $\rpn$ or $\SI^n$. 
\index{idim@$i$-dimensional complete affine subspace|textbf}

Again an affine subspace in $\SI^n$ is a lift of an affine subspace in $\rpn$, 
which is the interior of an $n$-hemisphere. 
Convexity and proper convexity in $\SI^n$ are defined in the same way as in $\rpn$. 
 \index{hemisphere} 
 \index{complete affine space|textbf} 
 


The complement of a codimension-one subspace $W$ in $\rpn$ can be considered as an affine 
space $\mathds{A}^n$ by correspondence 
\[[1, x_1, \dots, x_n] \ra (x_1, \dots, x_n)\] for a coordinate system where $W$ is given by $x_0=0$. 
The group $\Aff(\mathds{A}^n)$ of projective automorphisms acting on $\mathds{A}^n$ is identical with 
the group of affine transformations of form 
\[ \vec{x} \mapsto A \vec{x} + \vec{b} \] 
for a linear map $A: \bR^n \ra \bR^n$ and $\vec{b} \in \bR^n$. 
The projective geodesics and the affine geodesics agree up to parametrizations.
\index{group of affine transformations|textbf}
\index{affine coordinate system|textbf}
\index{affA@$\Aff(\mathds{A}^n)$|textbf}

A subset $A$ of $\RP^n$ or $\SI^n$ {\em spans} a subspace $S$ if $S$ is the smallest subspace containing $A$. We write $S = \langle A \rangle$. 
Of course, we use the same term for affine and vector spaces as well. 
\index{span} 
\index{lA@$\langle A \rangle$}

The following notation is used in the monograph. 
For a subset $A$ of a space $X$, we denote by $\clo_X(A)$
the closure of $A$ in $X$ and $\Bd_X A$ the boundary of $A$ in $X$. 
We omit the subscript $X$ if $X$ is clear from the context. 
If $A$ is a domain of a subspace of $\rpn$ or $\SI^n$, we denote by $\Bd A$ the topological boundary in the subspace. 
\index{Bd@$\Bd{ }$|textbf}
The closure $\clo(A)$ of a subset $A$ of $\rpn$ or $\SI^n$ is the topological closure in 
the respective spaces. 
\index{cl@$\clo(\cdot)$|textbf}
We will also denote by $K^o$ the manifold or orbifold interior for an $n$-manifold or $n$-orbifold $K$. Also, we may use $K^o$ as the interior relative to the topology of $P$ 
when $K$ is a domain $K$ in a totally geodesic subspace $P$ in $\SI^n$ or $\rp^n$.  
Define $\partial A$ for a manifold or orbifold $A$ to be the {\em manifold or orbifold boundary}. (See Section \ref{prelim-sub-ends}.)
\index{partial@$\partial$|textbf}
\index{boundary!manifold or orbifold} 
\index{boundary!topological} 
\label{manifoldboundary}

Let $p, q \in \SI^n$.
We also denote by $\overline{pq}$ a minor arc connecting $p$ and $q$ in a great circle in $\SI^n$. 
If $q \ne p_-$, this is unique. Otherwise, we need to specify a point in $\SI^n$ not antipodal to both. 
We denote by $\overline{pzq}$ the unique minor arc connecting $p$ and $q$ passing $z$. 
This is called a {\em great segment} 
\index{great segment|textbf} 

If $p, q \in \RP^n$, then $\overline{pq}$ denotes one of the closures of a component of $l - \{p, q\}$ for 
a one-dimensional projective line containing $p, q$. 
\index{Overline@$\overline{\ast\ast}$|textbf} 
\label{overline}



\subsection{Convex sets in $\rpn$ and $\SI^n$. }

Recall that an affine path in $\rpn$ is a complement of a codimension-one subspace. 
It has a canonical geodesic structure where each projective geodesic corresponds to 
affine geodesics up to parametrizations and conversely.
A {\em convex set } in $\rpn$ is a convex subset of an affine patch of $\rpn$.

We use a slightly different definition of convexity for $\SI^n$. 
\begin{definition} \label{prelim-defn-convexitySn} 
A {\em convex segment} is an arc contained in \index{great segment} a great segment. 
\index{convex segment|textbf} 
A {\em convex subset} of $\SI^n$ is a subset $A$ where every pair of points of $A$ is connected by a convex segment in $A$. 
\end{definition} 
\index{convex set for $\SI^n$|textbf}



\begin{definition}\label{prelim-defn-join}
	Given a convex set $D$ in $\rpn$, we obtain a connected cone $C(D)$ in $\bR^{n+1}-\{O\}$
mapping onto $D$,
	determined up to the antipodal map. For a convex domain $D \subset \SI^n$, we have a unique domain $C(D) \subset \bR^{n+1}-\{O\}$. 
For a nonzero vector $\vec x$, we denote by $\llrrparen{\vec{x}}$ the equivalence class of $\vec{x}$
where $\vec{x} \sim \vec{y}$ if and only if $\vec{x} = s\vec{y}, s> 0$.

	A {\em join} of two properly convex subsets $A$ and $B$ in a convex domain $D$ of $\rpn$ (resp. $\SI^n$) is defined as
	\[A \ast B := \{[ t \vec{x} + (1-t) \vec{y}]| \vec{x}, \vec{y} \in C(D),  [\vec{x}] \in A, [\vec{y}] \in B, t \in [0, 1] \} \]
	\[\hbox{(resp. }
	A \ast B := \{\llrrparen{t \vec{x} + (1-t) \vec{y}}| \vec{x}, \vec{y} \in C(D),  \llrrparen{\vec{x}} \in A, \llrrparen{\vec{y}} \in B, t \in [0, 1] \} 
	\hbox{)}\]
	where $C(D)$ is a cone corresponding to $D$ in $\bR^{n+1}$. 
The definition of the join is independent of the choice of $C(D)$ in $\SI^n$. 
	In $\RP^n$, the join may depend on the choice of $C(D)$. 
	Note we use $p \ast B = \{p\} \ast B$ interchangeably for a point $p$.  
\end{definition} 
\index{join|textbf} 


\begin{definition}
	Let $C_1, \dots, C_m$ respectively be cones in a set of independent vector subspaces $V_1, \dots, V_m$ of $\bR^{n+1}$. 
	In general, a {\em sum} of convex sets $C_1, \dots, C_m$ in $\bR^{n+1}$ 
	in independent subspaces $V_i$ is defined by
	\[ C_1+ \dots + C_m := \{v | v = \vec{c}_1+ \cdots + \vec{c}_m, \vec{c}_i \in C_i \}.\]
	A {\em strict join} of convex sets $\Omega_i$ in $\SI^n$ (resp. in $\RP^n$) is given as 
	\[\Omega_1 \ast \cdots \ast \Omega_m := \Pi'(C_1 + \cdots + C_m -\{O\}) \hbox{ (resp. } \Pi(C_1 + \cdots + C_m -\{O\})  ) \]
	where each $C_i-\{O\}$ is a convex cone with image $\Omega_i$ for each $i$
	for the projection $\Pi'$ (resp. $\Pi$). 
\end{definition}
\index{join!strict|textbf} 
\index{sum} 
\index{projection}

Note that each set $A$ in $\SI^n$ (resp. $\rpn$) 
has a minimal great sphere containing it, called 
{\em span} of $A$ and denoted by $\langle A \rangle$. 
\index{span} 
\index{langle@$\langle \cdot \rangle$} 

By the following, it is easy to see that either
a convex subset of $\SI^n$ is contained in an affine subspace, it is in a closed hemisphere, or it is a great sphere of dimension $\geq 1$. 
In the first case, the set embeds to a convex set in $\RP^n$ under the covering map.

Since an affine patch of $\rpn$ always lifts to an open hemisphere in $\SI^n$,
a convex subset of $\rpn$ always lifts to a convex subset of $\SI^n$ which maps to it homeomorphically
under the projection $\SI^n \ra \rpn$.  Hence, each convex subset of $\rpn$ corresponds to 
a convex subset of $\SI^n$ contained in an open hemisphere. 

\begin{proposition} \label{prelim-prop-classconv} 
	A closed convex subset $K$ of $\SI^n$ is either 
	a great sphere $\SI^{i_0}$ of dimension $i_0\geq 1$, 
	or is contained in a closed hemisphere $H^{i_0}$ in $\SI^{i_0}$ 
	and the following hold: 
	\begin{itemize} 
\item A closure of a convex set in $\SI^n$ is a convex compact subset of $\SI^n$. 
		\item There exists a great sphere $\SI^{j_0}$ of dimension $j_0 \geq 0$ 
		in the boundary $\Bd K$ and a compact properly convex domain $K_K$ 
		in an independent subspace of $\SI^{j_0}$ and 
		$K = \SI^{j_0} \ast K_K$, a strict join. 
		Moreover, $\SI^{j_0}$ is a unique maximal great sphere in $K$. 
		\item If $K$ is not a great sphere, then $K^o$ is a convex domain in the interior of an $(i_0+1)$-hemisphere for some $i_0$. 
\item  Unless $K$ is a great sphere, $\partial K = \Bd _{\langle K \rangle} K$, 
and $K$ is homeomorphic to a cell. 
\item For a properly convex compact domain $K$ in $\rpn$,  we have
$\partial K = \Bd_{\rpn} K$. 
		\end{itemize} 

	\end{proposition}
\begin{proof} 
The first item follows since we can use sequence arguments on pairs of points and 
the fact that a sequence of convex segments converges to a convex segment up to a choice of 
subsequences. 

	Let $\SI^{i_0}$ be the span of $K$. Then $K^o$ is not an empty domain
	in $\SI^{i_0}$. The map 
	$x \mapsto \bdd(x, K)$ is continuous on $\SI^{i_0}$. 
	Choose a maximum point $x_0$. 
	If the maximum $R$ is $< \pi/2$, then the spherical geometry tells us that 
	there are at least two points $y, z$ of $K$ closest to $x_0$ of 
	a same distance from $x_0$: 
	Suppose not. 
We have $\{y\} = K \cap B^\bdd_R(x_0)$ for a unique point $y$. 
The $A$ be a small connected open neighborhood of $y$ in $\partial B^\bdd_R(x_0)$. 
Then the arc $B := \partial B^\bdd_R(x_0) -A$ has a $\delta$-neighborhood $U$ 
disjoint from $K$. 
We let $s$ be the great circle through $y$ and $x_0$. 
Then there is $\eps', \eps'> 0$, such that for $x'_0$ in 
$\{x'_0| \bdd(x'_0, R) < \eps', x'_0 \in s - \ovl{x_0y}\}$, 
we have $B^\bdd(x'_0) \subset B^\bdd_R(x_0) \cup U$, 
and $B^\bdd(x'_0, R)$ is disjoint from $K$. Now, $R$ can be increased
which is a contradiction. 
%
%
%

	Given two points $y, z$, there is a closer point on $\overline{yz}^o$ in $K$ to $x_0$.
	This is a contradiction. 
	Hence, 
	$K = \SI^{i_0}$. Otherwise, $K$ is a subset of an $i_0$-hemisphere in
	$\SI^{i_0}$. (See \cite{cdcr1} also.)  
	
	The second part follows 
	from Section 1.4 of \cite{ChCh93}. 
	(See also \cite{GV58}.) 
	Hence, we obtain a unique maximal great sphere $\SI^{j_0}$ in $K$ 
	which is contained in $\Bd K$, and 
	$K$ is a union of $j_0+1$-hemispheres with common boundary 
	$\SI^{j_0}$.
	
	By choosing an independent subspace $\SI^{n-j_0-1}$ 
	to $\SI^{j_0}$, each $j_0+1$-hemisphere in $K$ is 
	transverse to  $\SI^{n-j_0-1}$ and hence meets it in 
	a unique point. We let $K_K$ denote the set of intersection 
	points. Therefore, $K =  \SI^{j_0} \ast K_K$. 
	
	There is a map $K \ra K_K$ given by sending a $j_0+1$-hemisphere to 
	its intersection point. Obviously, this is a restriction of projective 
	diffeomorphism from the space of $j_0+1$-hemispheres with boundary $\SI^{j_0}$ 
	to $\SI^{n-j_0 -1}$. 
	Since $K$ cannot contain a higher-dimensional great sphere, it follows that
	 $K_K$ is also properly convex. 

For the third, foruth, and fifth item, the fact that $K$ is a join of a great sphere with a properly convex domain implies this since the items are true for properly convex domains. 
	\end{proof}

\subsection{The Hausdorff distances used} \label{prelim-sub-Haus}

We will be using the standard elliptic metric $\bdd$ on $\RP^n$ (resp. in $\SI^n$) where
the set of geodesics coincides with the set of projective geodesics 
up to parametrizations. Sometimes, these are called Fubini-Study metrics. 
 \index{d@$\bdd$} \index{elliptic metric} 
 \label{bdd} 


\begin{definition}\label{prelim-defn-Haus} 
Given a subset $A$ of $\SI^n$ (resp. $\RP^n$), 
we define
\[N_\eps(A):= \{x\in \SI^n \mid \bdd(x, A) < \eps\}\, \big(
\hbox{ resp. }N_\eps(A):= \{x\in \RP^n \mid \bdd(x, A) < \eps\}.\big)
\] 
Given two subsets $K_1$ and $K_2$ of $\SI^n$ {\rm (}resp. $\RP^n${\rm ),}
we define the {\em Hausdorff distance} $\bdd_H(K_1, K_2)$ between $K_1$ and $K_2$ to be 
\[ \inf\{\eps > 0| K_2 \subset N_\eps(K_1), K_1 \subset N_\eps(K_2)  \}.\] 
The simple distance $\bdd(K_1, K_2)$ is defined as 
\[ \inf\{ \bdd(x, y)| x \in K_1, y\in K_2 \}.\]
\end{definition}
\index{dh@$\bdd_H$|textbf} \index{Hausdorff distance|textbf} 


A sequence $\{A_{i}\}$ of compact sets 
is said to {\em converges} to a compact subset $A$ if $\{\bdd_{H}(A_{i}, A)\} \ra 0$. 
In this case, the limit is unique. 
Recall that every sequence of compact sets $\{A_i\}$ in  $\SI^n$
(resp. $\RP^n$) has
a convergent subsequence. 
The limit $A$ is characterized as follows if it exists: 
\begin{equation} \label{prelim-eqn-characterization} 
A := \{a\in H|\, a \hbox{ is a limit point of some sequence } \{a_i| a_i \in A_i\} \}.  
\end{equation}
(See Proposition E.12 of \cite{BP92} for a proof since the Chabauty topology for a compact space is the Hausdorff topology. See also Munkres \cite{Munkres75}.)
\index{convergence!geometrically} 

We will use the same notation even when $A_i$ and $A$ are closed subsets of a fixed open domain
using $\bdd_H$ and $\bdd$.

\begin{proposition}[Benedetti-Petronio] \label{prelim-prop-BP} 
	A sequence $\{A_i\}$ of compact sets in  $\RP^n$ {\rm (}resp.  $\SI^n${\rm )} 
	converges to $A$ in the Hausdorff topology if and only if 
	both of the following hold\/{\rm :} 
	\begin{itemize} 
		\item If $x_{i_j} \in A_{i_j}$ and $\{x_{i_j}\} \ra x$, where 
		$i_j\ra \infty$, then $x \in A$. 
		\item If $x \in A$, then there exists $x_i \in A_i$ for each $i$ 
		such that $\{x_i\} \ra x$. 
	\end{itemize} 
\end{proposition}
\begin{proof} 
	Since $\SI^n$ and $\RP^n$ are compact, the Chabauty topology is the same as the Hausdorff topology. Hence, this follows from Proposition E.12 of Benedetti-Petronio \cite{BP92}.
\end{proof}

\begin{lemma} \label{prelim-lem-geoconv} 
	Let $\{g_i\}$ be a sequence of elements of $\PGL(n+1, \bR)$
	{\rm (}resp. $\SLpm$\/{\rm )}
	converging to $g_\infty$ in $\PGL(n+1, \bR)$ {\rm (}resp. $\SLpm$\/{\rm ).} 
	Let $\{K_i\}$ be a sequence of compact sets, and let $K$ be another one.
	Then $\{K_i\} \ra K$ if and only if $\{g_i(K_i)\} \ra g_\infty(K)$.
	\end{lemma} 
\begin{proof} 
	We use the above point description of the geometric limit.
	\end{proof} 

An $n$-hemisphere $H$ in $\SI^n$ {\em supports} a set $D$ if 
$H$ contains $D$. $H$ is called a {\em supporting } hemisphere. 
An oriented hyperspace $S$ in $\SI^n$ {\em supports} a set $D$ if 
the closed hemisphere bounded in an inner-direction by $S$ contains $D$. 
$S$ is called a {\em supporting} hyperspace. 
If a supporting hyperspace contains at least one boundary point $x$ of $D$,
then it is called a {\em sharply supporting} hyperspace at $x$.
If the boundary of a supporting $n$-hemisphere is sharply supporting at $x$, 
then the hemisphere is called a {\em sharply-supporting hemisphere} at $x$. 
\index{supporting} 
\index{hyperspace} 
\index{hyperspace!supporting}


\index{supporting hyperspace|textbf} 
\index{supporting hemisphere|textbf} 
\index{supporting hyperspace!sharply|textbf} 
\index{supporting hemisphere!sharply|textbf}


\begin{proposition} \label{prelim-prop-convC} 
	Let $K_i$ be a sequence of compact convex sets {\rm (}resp. cells{\rm )} of $\SI^n$. 
	Then up to choosing a subsequence, we have
	$K_i \ra K$ for a compact convex set {\rm (}resp. cell{\rm )}  $K$ of $\SI^n$. 
	Also, a geometric limit must be a compact convex cell of possibly lower dimension when $K_i$ are compact convex cells. 
	If $K_i$ is in a fixed $n$-hemisphere, then 
	so is $K$. 
	\end{proposition} 
\begin{proof} 
By Proposition \ref{prelim-prop-BP}, we can show this when $K_i$ is a hemisphere. 
 For other cases, consider sequences of segments and Proposition \ref{prelim-prop-BP}. 
	\end{proof}

The following is probably well-known. 
\begin{lemma} \label{prelim-lem-bdconv}
Suppose that one of the following holds{\rm :} 
	\begin{itemize} 
		\item $K$ and $K_i$ for each $i$, $i=1, 2,\dots$, 
		is a compact convex domain in $\SI^n$. 
		\item $K$ and $K_i$ is a convex open domain in $\SI^n$
		\item $K$ and $K_i$ is a properly convex domain in $\RP^n$. 
		\end{itemize} 
	Suppose that a sequence $\{\clo(K_i)\}$ 
	geometrically converges to 
	$\clo(K)$ with a nonempty interior.
	Then $\{\Bd K_i\} \ra \Bd K$. 
\end{lemma} 
\begin{proof}
	We prove for $\SI^n$. 
	Suppose that a point $p$ is in $\Bd K$. 
	Let $B_\eps(p)$ be an open $\epsilon$-ball of $p$. 
	Suppose $B_\eps(p)\cap K_i = \emp$ for infinitely many $i$. 
	Then $p$ cannot be a limit point of $K$ by Proposition \ref{prelim-prop-BP}. 
	This is a contradiction. 
	Thus, $B_\eps(p) \cap K_i \ne \emp$ for $i > N$ for some $N$. 
	Suppose that $B_\epsilon(p) \subset K_i$ for infinitely many $i$.
	Then each point in $B_\epsilon(p) - K$ is 
	a limit point of some sequence $p_i, p_i \in K_i$. 
	and hence we have $B_\epsilon(p) \subset K, p \in K^o$, a contradiction. 
	Hence, given $\epsilon> 0$,  $B_\epsilon(p) \cap \Bd K_i \ne \emp$ for $i > M$ for some $M$. 
	Then $p$ is a limit of a sequence $p_i, p_i\in \Bd K_i$. 
	
	Conversely, suppose that a sequence $\{p_{i_j}\}, p_{i_j}\in \Bd K_{i_j}$ where 
	$i_j \ra \infty$ as $j \ra \infty$, converges to $p$. 
	Then $p\in K$ clearly.  
	Suppose that $p \in K^o$. 
	Then there is $\epsilon, \epsilon> 0,$ with $B_\epsilon(p) \subset K$.
	Now, $K_{i_j}$ has a sharply supporting closed hemisphere $H_{i_j}$ at $p_{i_j}$ with 
	$ K_{i_j} \subset H_{i_j}$. 
	Since $\{p_{i_j}\} \ra p$, we may choose a subsequence $k_j$ such that 
	$\{H_{k_j}\} \ra H_\infty$ and 
	$\bdd_H(H_{k_j}, H_\infty) < \epsilon/4$ for a hemisphere $H_\infty$. 
Clearly, $H_\infty$ is a closed hemisphere with a boundary point $p$.
	Let $q \in B_{3\epsilon/4}(p) - H_\infty$ such that $\bdd_{H}(q, H_\infty) > \epsilon/4$. 
	Hence, $B_{\epsilon/4}(q) \in B_\epsilon(p)  - H_{k_j}$ for all $j$. 
	Since $K_{k_j} \subset H_{k_j}$, no sequence $\{q_{k_j}\}, q_{k_j} \in K_{k_j}$ converges 
	to $q$. However, since $\{K_{k_j}\} \ra K$ and $q \in K$, this is a contradiction to Proposition \ref{prelim-prop-BP}.
	Hence, $p \in \Bd K$. Now, Proposition \ref{prelim-prop-BP} proves 
	$\{\Bd K_i\} \ra \Bd K$. 
	
	When $K_i$ is an open domain in $\SI^n$, we just need to take its closure and use the first part. 
	
	For the $\RP^n$-version, we lift $K_i$ to $\SI^n$ to 
	properly convex domains $K_i'$. Now, we may also choose a subsequence 
	such that $\{K'_i\}$ geometrically converges to a choice of a lift $K'$ of $K$
	by Proposition \ref{prelim-prop-BP}. 
	Since $K'$ is properly convex, 
	$K'$ is in a bounded subset of an affine subspace of $\SI^n$. 
	Then the result follows from the $\SI^n$-version. 
\end{proof} 
We note that the last statement is false if $\{K_i\}$ geometrically converges to a hemisphere when lifted to 
$\SI^n$. 

\begin{theorem}\label{prelim-thm-partialKi}
	Suppose that $K_i$ and $K$ are compact convex balls of the same dimension in $\SI^n$. 
	Suppose that $\{K_i\} \ra K$. 
	It follows that 
	\begin{equation} \label{prelim-eqn-partialKi}
	\{\partial K_i\} \ra \partial K.
	\end{equation} 
	This holds also provided $K_i$ and $K$ are properly compact convex in $\RP^n$
	with $\{K_i\} \ra K$. 
\end{theorem}
\begin{proof} 
	Since $K_i$ and $K$ are of the same dimension, 
	we find $g_i \in \SLpm$ such that 
	$g_i(\langle K_i\rangle) = \langle K\rangle$ and $\{g_i\} \ra g_\infty$ for 
	$g_\infty \in  \SLpm$. Then 
	$\{g_i(K_i)\} \ra g_\infty(K)$. 
	Thus we obtain $\{\partial g_i(K_i)\} \ra \partial g_\infty(K)$ 
	by Lemma \ref{prelim-lem-bdconv}. 
	Hence, $\{\partial K_i\} \ra \partial K$ by Lemma \ref{prelim-lem-geoconv}. 
	
	The last statement follows from  Proposition \ref{prelim-prop-closureind}. 
	\end{proof}

\subsection{The Hilbert metric} \label{prelim-sub-Hilbert}   
Let $\Omega$ be a convex open domain. A line or a one dimensional subspace
in $\RP^n$ has a $2$-dimensional homogeneous coordinate system. 
Let $[o, s, q, p]$ denote the cross ratio of four points on a line as defined by 
\[ \frac{\bar o - \bar q}{\bar s - \bar q} \frac{\bar s - \bar p}{\bar o - \bar p} \] 
where $\bar o, \bar p, \bar q, \bar s$ denote respectively the first coordinates of the homogeneous coordinates \index{Hilbert metric} 
of $o, p, q , s$ provided that the second coordinates are normalized to be $1$. 
Define a pseudo-metric for $p, q \in \Omega$, 
$d_\Omega(p, q)= \log|[o,s,q,p]|$ where $o$ and $s$ are 
endpoints of the maximal segment in $\Omega$ containing $p, q$
where $o, q$ separates $p, s$.  
If $\Omega$ is properly convex,  then it is a metric and  a Finsler metric 
(See \cite{Kobpaper}.) If $\Omega$ is complete affine, $d_\Omega$ is identically zero. 

\begin{lemma} \label{prelim-lem-convmetric}
Assume that $\{K_i\} \ra K$, $K_i, K \subset \SI^n$ {\rm (}resp. $\subset \rpn$\/{\rm )} geometrically for a sequence of properly convex compact domains $K_i$ and a convex compact domain $K$ which is not necessarily properly convex.   
Suppose that two sequences of points $\{x_i| x_i \in K_i^o\}$ and $\{y_i| y_i \in K_i^o\}$ 
converge to $x, y \in K^o$ respectively. 
Then 
\begin{equation} \label{prelim-eqn-HiHa}
\{d_{K_i^o}(x_i, y_i)\} \ra d_{K^o}(x, y),
\end{equation}   
where $d_{K^o}$ can be a pseudo-metric. 
\end{lemma} 
\begin{proof}
Let $z_i$ and $t_i$ denote the endpoints of the maximal line containing $x_i$ and $y_i$ in $K_i$. 
Let $z$ and $t$ denote the endpoint of one containing $x$ and $y$ in $K$. 
It is easy to see $z_i \ra z$ and $t_i \ra t$. 
Let $l_i$ and $m_i$ denote the sharply 
supporting great hyperspheres at $z_i$ and $t_i$ for $K_i$. 
Then $l_i$ and $m_i$ must converge up to subsequences to a supporting great hypersphere 
at $z$ and $t$ respectively since the closures of components of the complements of $l_m$ and $m_i$ are disjoint from $K_i^o$. 
This implies that $\overline{z_it_i}$ converges to a subsegment of $\overline{zt}$
up to a choice of subsequences. 
However, a limit cannot be a proper segment since otherwise a boundary point of $K_i$ converges to 
an interior point of $K$ contradicting Theorem \ref{prelim-thm-partialKi}.
This implies the result.  \hfill \SSn {\parfillskip0pt\par}
\end{proof} 


\begin{lemma}[Cooper-Long-Tillmann \cite{CLT15}] \label{prelim-lem-nhbd} 
Let $U$ be a convex subset of a properly convex domain $V$ in $\SI^n$ {\rm (}resp. $\RP^n${\rm ).} 
Let \[U_\eps := \{x \in V | d_{V}(x, U) \leq \eps\}\] for $\eps > 0$. 
Then $U_\eps$ is properly convex. 
\end{lemma}
\begin{proof} 
Given $u, v \in U_\eps$, we find 
\[w, t \in V  \hbox{ such that } d_V(u, w) < \eps, d_V(v, t) < \eps.\]
Then each point of $\ovl{uv}$ is within $\eps$ of $\ovl{wt} \subset U$ in the $d_V$-metric.
By Lemma 1.8 of \cite{CLT15}, this follows. \hfill 
\SSn {\parfillskip0pt\par}
\end{proof} 

\begin{proposition} \label{prelim-prop-AutK} 
	Let $\Omega$ be a properly convex domain in $\SI^n$ {\rm (}resp. $\rpn$\/{\rm ).}
Then 	the group $\Aut(\Omega)$ of projective automorphisms of 
	$\Omega$ is closed in  $\SLnp$ {\rm (}resp. $\PGL(n+1, \bR)$\/{\rm )}. 
	Also, the set of elements of $g$ of $\Aut(\Omega)$ such that 
	$g(x)\in K$ for a compact subset $K$ of $\Omega$ is compact. 
\end{proposition}
\begin{proof} 
	We prove for $\SI^n$. 
	Clearly, the limit of a sequence of elements in $\Aut(\Omega)$ 
	is an isometry of the Hilbert metric of $\Omega$. 
	Hence, it acts on $\Omega$ since the Hilbert metric is complete according to 
\cite{Kobpaper}. 
	
	For the second part, we take an  $n$-simplex $\sigma$ with 
	a point $x$ in the interior as a base point.
	
	The space $\mathcal{S}^n$ of nondegenerate 
	convex $n$-simplices with base points in their interiors with 
	the Hausdorff topology is homeomorphic to 
	$\SL_\pm(n+1, \bR)$ since the action of $\SL_\pm(n+1, \bR)$ is 
	simply transitive on $\mathcal{S}^n$.
	
	The subspace of simplices of the form 
	$g(\sigma)$ for $g$ with $g(x) \in K$, $g\in \Aut(\Omega)$ is compact 
	by the existence of the Hilbert metric: We can show this by 
	using the invariants. The edge lengths are invariants. 
	The distance from each vertex to the hyperspace 
	containing the remaining vertices is an invariant of the action. 

We can see that the set of such $g$ is bounded: We can find the bounded set  
\[\{h_g\in \Aut(\SI^n)|h_g\circ g(\sigma, x) = (\sigma, x)\}\]
since any sequence of sets of the form 
$g(\sigma)$ does not degenerate and $g(x)$ is uniformly bounded away from the boundary of 
$g(\sigma)$.  Since the simplex $\sigma$ and the basepoint $x$ are fixed, 
we have $h_g\circ g = \Idd$ for $g(x) \in K$. Hence the set $\{g| g(x)\in K\}$ is uniformly bounded. 

Since $\mathcal{S}^n$ is diffeomorphic to $\SL_\pm(n+1, \bR)$, 
the closedness of $\Aut(\Omega)$ proves the result. 
\hfill \SnP {\parfillskip0pt\par}
	\end{proof} 



\subsection{Topology of orbifolds} \label{prelim-sub-ends}
We summarize Chapter 4 of \cite{Cbook}. We will only briefly outline it. 
An $n$-dimensional {\em orbifold structure} on a Hausdorff space $X$ is given by a maximal collection of charts $(U, \phi, G)$ satisfying the following conditions: 
\begin{itemize}
	\item $U$ is an open subset of $\bR^n$ and $\phi:U \ra X$ is 
	a map and $G$ is a finite group acting on $U$,
	\item the chart $\phi: U\ra X$ induces a homeomorphism 
	$U/G$ to an open subset of $X$, 
	\item the collection of the sets of the form $\phi(U)$ from a triple 
$(U, \phi, G)$ in the collection  covers $X$, 
	\item for any pair of models $(U, \phi, G)$ and $(V, \psi, H)$ 
	with an inclusion map $\iota: \phi(U) \ra \psi(V)$
    lifts to an embedding $U\ra V$ that is equivariant with respect to 
    an injective homomorphism $G \ra H$. (compatibility condition)
    \end{itemize}
An {\em orbifold} $\orb$ is a topological space with an orbifold structure. 
The boundary $\partial \orb$ of an orbifold is defined as the set of points with only half-open sets as models. 
(This is often distinct from topological boundary.)
An  {\em $i$-dimensional suborbifold} $N$ of $\orb$ is a subspace of $X$ equipped with a maximal collection of 
charts containing the orbifold charts of form 
\[(U\cap \phi^{-1}(N), \phi| U\cap \phi^{-1}(N), G| U \cap \phi^{-1}(N)) \hbox{ where } 
U\cap \phi^{-1}(N) \subset \bR^i \subset \bR^n\]
from charts of form $(U, \phi, G)$ of $\orb$. (See Definition 4.4.2 of \cite{Cbook}.)
(Note this is more general than other definitions. See \cite{Borzellino12}, \cite{Borzellino15} 
for various approaches.)
Every boundary component of $\orb$ is a suborbifold.    
\index{orbifold} \index{orbifold!boundary|textbf} 

Orbifolds are stratified by manifolds. 
Let $\orb$ denote an $n$-dimensional orbifold with finitely many ends. 
(See \cite{Thurston80}, \cite{ALR07}, \cite{Kapovich09}
and \cite{Cbook} for details.)

An {\em orbifold covering map} $p: \orb_1\ra \orb$ for orbifolds $\orb$ and $\orb_1$ is a map such that
for any point on $\orb$, 
\begin{itemize} 
\item there is a connected open set $U \subset X$ with 
a model $(\tilde U, \phi, G)$ as above 
\item whose inverse 
image $p^{-1}(U)$ is a union of connected open sets $U_i$ for $i \in I$ for some index set $I$ 
of $\orb_1$ with models $(\tilde U, \phi_i, G_i)$ for a subgroup $G_i \subset G$
and the induced chart $\phi_i: \tilde U \ra \tilde U_i$. 
\end{itemize} 
\index{orbifold!covering} 

We say that an orbifold is a manifold if it has a subatlas of charts with trivial local groups. 
We will consider good orbifolds only, i.e., covered by simply connected manifolds. 
In this case, the universal covering orbifold $\torb$ is a manifold
with a projective orbifold-covering map $p_\orb: \torb \ra \orb$. \index{po@$p_{\orb}$}
The group of deck transformations is denoted by $\pi_1(\orb)$ or $\bGamma$, 
and is said to be the {\em fundamental group} of $\orb$. 
It acts properly discontinuously on $\torb$ but not necessarily freely. 

We will follow  Section 4.4.2 of \cite{Cbook}. (See Chapter 4 of \cite{Hirsch} for manifolds.) 
A {\em neat suborbifold} $N \subset \orb$ is a suborbifold such that $\partial N \subset \partial \orb$
and the tangent space to $N$ at $\partial N$ is transversal to the tangent space of $\partial \orb$.  
Of course, if $\partial N=\emp$, a suborbifold $N$ is still considered neat. 
Let $N(N)$ denote the subspace of the tangent bundle of $\orb$ over $N$ consisting of vectors perpedicular to $N$. Let $\eps: N \ra [0, \infty)$ denote a real valued function. 
We denote by $N_\eps(N)$ the subspace of normal vectors to $N$ of length $\leq \eps(x)$ at each 
$T_x \orb$, $x\in N$. The exponential map is an embedding of $N_\eps(N)$ to $\orb$
for sufficiently small $\eps$. 
We call the image a {\em tubular neighborhood} of $N$. 
\index{suborbifold!neat} 
\index{neighborhood!tubular}

\begin{proposition}\label{prelim-prop-neat} $ $
\begin{itemize} 
\item We can give a Riemannian metric on an orbifold 
$\orb$ such that $\partial \orb$ is totally geodesic
and a neat suborbifold $N$ to be totally geodesic perpendicular to $\partial \orb$. 
\item  Suppose that $N$ is a connected neat suborbifold of $\orb$. 
A tubular neighborhood of $N$ is always diffeomorphic to the orbifold product $N \times J$ 
where $J$ is an orbifold covered by an open ball of dimension $\dim \orb - \dim N$. 
\begin{itemize} 
\item If $J$ is $[0, 1)$, then $N$ is a boundary component. 
\item If $J$ is doubly covered by $(-1, 1)$, then $N$ is a silvered submanifold. 
\end{itemize} 
\end{itemize} 
\end{proposition} 
\begin{proof}
See Section 4.4.2 and Lemma 4.4.1 of \cite{Cbook}.
\end{proof} 

Recall that a quotient space of a manifold by a properly discontinuous action of a group
has an orbifold structure (see Theorem 4.2.3 of \cite{Cbook}.)  
In these cases, we will consider the quotient space
to have an orbifold structure automatically. 

\begin{proposition} \label{prelim-prop-prdisc} 
	Let $\Omega$ be a properly convex domain in $\SI^n$ {\rm (}resp. $\rpn$\/{\rm ).}
	Suppose that a discrete subgroup $\Gamma$ of $\SLnp$ {\rm (}resp. $\PGL(n+1, \bR)$\/{\rm )} acts on $\Omega$. 
	Then $\Omega/\Gamma$ has a structure of an $n$-orbifold. 
\end{proposition}
\begin{proof}
	The second part of Proposition \ref{prelim-prop-AutK} implies that 
	$\Gamma$ acts properly discontinuously. 
(Or Corollary 3.5.11 of \cite{Thurston97} shows this.)
We obtain that $\Omega/\Gamma$ has a structure of an orbifold. 
\hfill \SnT {\parfillskip0pt\par}
\end{proof}

\subsection{Geometric structures on orbifolds}\label{prelim-sub-geostr}

Let $G$ be a Lie group acting on a space $X$ transitively and effectively. 
An {\em $(X,G)$-structure} on an orbifold $\mathcal{O}$ is a maximal 
atlas of charts from open subsets of 
\index{xgs@$(X, G)$-structure|textbf} 
\index{geometric structure|textbf} 
$X$ with finite subgroups of $G$ acting on them, and the inclusions always lift to restrictions of 
elements of $G$ in open subsets of $X$. This is equivalent to saying that the orbifold $\mathcal{O}$ has 
a simply connected manifold cover $\tilde{\mathcal{O}}$ with an immersion 
$\dev:\tilde{\mathcal{O}} \ra X$ and the fundamental group $\pi_1(\mathcal{O})$ acts on $\tilde{\mathcal{O}}$ properly 
discontinuously such that $h: \pi_1(\orb) \ra G$ is a homomorphism satisfying 
$\dev\circ \gamma = h(\gamma)\circ \dev$ for 
each $\gamma \in \pi_1(\orb)$.  
Here, $\pi_1(\orb)$ is allowed to have fixed points with finite stabilizers. 
(We shall use this second more convenient definition here.)
The map $\dev$ is said to be a {\em developing map} and 
$h$ is said to be a {\em holonomy homomorphism}
and the pair $(\dev, h)$ is called a {\em development pair}. 
Tor a given $(X,G)$-structure, it is determined only up to an action 
\[(\dev, h(\cdot)) \mapsto (k\circ \dev, kh(\cdot)k^{-1}) \hbox{ for } k\in G.\] 
Conversely, a development pair completely 
determines the $(X,G)$-structure. 
(See Chapter 6 of \cite{Cbook} or 
Thurston \cite{Thurston97} for the general theory of geometric structures.)

Thurston showed that an orbifold with an $(X,G)$-structure is always good, i.e., covered \index{orbifold!good} 
by a manifold with an $(X,G)$-structure. (See Proposition 13.2.1 
of Chapter 13 of Thurston 
\cite{Thurston80}.) 
Hence, every geometric orbifold is of form $\tilde M/\Gamma$ for a discrete 
group $\Gamma$ acting on a simply connected manifold $\tilde M$.
Here, we have to understand $\tilde M/\Gamma$ as having an orbifold structure 
coming from an atlas 
where each model set is based on a precompact open cell of $\tilde M$ on which 
a finite subgroup of $\Gamma$ acts. (See Theorem 4.2.3 
of \cite{Cbook} for details.)

\index{geometric structure|textbf} 

Let $(X_i, G_i)$ be a pair where $G_i$ is a Lie group acting effectively and transitively on $X_i$ 
for $i=1, 2$. 
Suppose that $X_1 \subset X_2$ and $G_1 \subset G_2$. 
Given an orbifold $\orb$ an $(X_1, G_1)$-structure is {\em compatible} with 
an $(X_2, G_2)$-structure if charts in the $(X_1, G_1)$-structure is in
the charts of $(X_2, G_2)$-structure. 

\index{geometric structure!compatible|textbf}

\subsection{Real projective structures on orbifolds.} \label{prelim-sec-rps}





We consider an orbifold $\orb$ with a real projective structure, i.e., 
($\rpn, \PGLnp)$-structure: 
This can be expressed as 
\begin{itemize}
\item having a pair $(\dev, h)$ where 
$\dev:\torb \ra \rpn$ is an immersion equivariant with respect to 
\item the homomorphism $h: \pi_1(\orb) \ra \PGLnp$ where 
$\torb$ is the universal cover, and $\pi_1(\orb)$ is the group of deck transformations acting on $\torb$. 
\end{itemize}

We will regard $\torb$ to have an induced real projective structure, where the covering map
is a projective. Converly, any convex open domain projectively covering $\orb$ is a universal cover. 

We will usually use only one pair where $\dev$ is an embedding for this paper and hence 
identify $\torb$ with its image. 
A {\em holonomy} is an image of an element under $h$. 
The {\em holonomy group} is the image group $h(\pi_1(\orb))$.  
\index{holonomy|textbf} 
\index{holonomy group|textbf}
\index{developing map|textbf} 
\index{development pair|textbf} 
\index{holonomy homomorphism|textbf}
\index{real projective structure|textbf} 

We denote by $\Aut(K)$ the group of projective automorphisms of a set $K$ in 
$\RP^n$ (or $\SI^n$). 
\label{Aut} 
The Klein model of the hyperbolic geometry is given as follows: 
Let $x_0, x_1, \dots, x_n$ be the standard coordinates of $\bR^{n+1}$. 
Let $\bB$ be the interiorof a standard ball in $\rpn$  or $\SI^n$ that is the image of the positive cone $x_0^2 > x_1^2 + \dots + x_n^2$ in $\bR^{n+1}$.
Then $\bB$ can be identified with a hyperbolic $n$-space. The isometry group of
the hyperbolic space 
equals the group $\Aut(\bB)$ of projective automorphisms acting on $\bB$. 
Thus, a complete hyperbolic manifold carries a unique real projective structure and is 
projectively diffeomorphic to $\bB/\Gamma$ for $\Gamma \subset \Aut(\bB)$.  
In fact, $g(\bB)$ for any $g \in \PGLnp$ will serve as a Klein model of the hyperbolic space, and 
$\Aut(g\bB) = g\Aut(\bB)g^{-1}$ is the isometry group.
(See \cite{Cbook} for details.) \index{b@$\bB$|textbf}
\index{Klein model} 

A real projective orbifold is {\em convex} (resp. {\em properly convex})  if it is projectively diffeomorphic to a  projective quotient of a convex domain (resp. {\em prolery convex domain}). 
in an affine subspace in $\RP^n$.  (Note that this definition is stricter than ones in \cite{psconv} but conforms to definitions in current literature.) 
\index{real projective structure!convex|textbf} 
\index{real projective structure!properly convex|textbf} 

Above hyperbolic manifolds form examples of properly convex real projective manifolds.

A {\em totally geodesic hypersurface} (resp. {\em suborbifold}) $A$ in $\torb$ is 
a suborbifold of codimension-one (resp. codimension $\geq 1$) 
where each point $p$ in $A$ has a neighborhood $U$ in $\torb$ 
such that $D|A\cap U$ has the image in a hyperspace (resp. subspace). 
A suborbifold $A$ pf $\orb$ is a {\em totally geodesic hypersurface} (resp. {\em suborbifold})  
if it is covered by  a one in $\torb$. 
\index{hypersurface!totally geodesic|textbf}

An {\em affine $n$-orbifold} is an $n$-orbifold with an affine structure, i.e., a geometric structure modeled on $(\mathds{A}^n, \Aff(\mathds{A}^n))$. 
\index{affine orbifold|textbf} 
A Euclidean manifold admits an affine structure. The Euclidean structure is compatible with 
the affine structure as the charts in the Euclidean geometric structures are also in
the charts for the affine structure. 
Clearly, an affine structure is compatible with a real projective structure.

\subsection{Spherical real projective structures} \label{prelim-sub-orp}

A {\em spherical real projective structure} is an $(\SI^n, \SL_\pm(n+1, \bR)$-structure. 
\index{spherical real projective structure|textbf} 
We now discuss the standard lifting: 
A real projective structure on $\mathcal{O}$ provides us 
with a development pair $(\dev, h)$ where 
$\dev: \tilde{\mathcal{O}} \ra \RP^n$ is an immersion and $h: \pi_1(\mathcal{O}) \ra \PGL(n+1, \bR)$ is a homomorphism. \index{developing map!lifting} 
Since $p_{\SI^n}$ is a covering map and $\tilde{\mathcal{O}}$ is a simply connected manifold,  \index{holonomy homomorphism!lifting} 
$\mathcal{O}$ being a good orbifold, 
there exists a lift $\dev': \tilde{\mathcal{O}} \ra \SI^n$ unique up to the action of $\{\Idd, \mathcal{A}\}$.
This induces a spherical real projective structure on $\tilde{\mathcal{O}}$ and $\dev'$ is a developing map for this real projective structure.
Given a deck transformation $\gamma:\tilde{\mathcal{O}} \ra \tilde{\mathcal{O}}$, the composition
$\dev' \circ \gamma$ is again a developing map for the real projective structure
and hence equals $h'(\gamma) \circ \dev'$ for $h'(\gamma) \in \SL_\pm(n+1, \bR)$. 
We verify that $h':\pi_1(\mathcal{O}) \ra \SL_\pm(n+1, \bR)$ is a homomorphism. 
Hence, $(\dev', h')$ gives us a spherical real projective structure, which 
induces the original real projective structure. 
\index{developing map!lifting} 
\index{holonomy homomorphism!lifting}

Given a real projective structure where $\dev: \torb \ra \rpn$ is an embedding onto a properly convex 
open subset, the developing map  
$\dev$ lifts to an embedding $\dev': \torb \ra \SI^n$ onto 
an open domain $D$ with no pair of antipodal  points. $D$ is determined up to $\mathcal{A}$. 

We will identify $\torb$ with $D$ or $\mathcal{A}(D)$ and 
$\pi_1(\orb)$ with $\bGamma$. 
Then the holonomy group $\bGamma$ lifts to a subgroup $\bGamma'$ of $\SLnp$ acting effectively and discretely on $\torb$.
There is a unique way to lift such that $D/\bGamma'$ is projectively diffeomorphic to $\torb/\bGamma$. 

%

\begin{theorem}\label{intro-thm-doubledef} 
	There is a one-to-one correspondence between the space of real projective structures on an orbifold $\orb$ 
	and the space of spherical real projective structures on $\orb$. 
	Moreover, a real projective diffeomorphism
	of real projective orbifolds is an $(\SI^n, \SL_\pm(n+1, \bR))$-diffeomorphism of 
corresponding spherical real projective orbifolds, 
	and vice versa. 
\end{theorem} 
\begin{proof} 
	Straightforward. See p. 143 of Thurston \cite{Thurston97} (see also Section \ref{prelim-sub-lifting}). 
\end{proof}



\begin{proposition}[Selberg-Malcev]\label{prlim-prop-Malcev} 
	The holonomy group of a real projective $n$-orbifold is residually finite and virtually
torsion-free.
The fundamental group 
of a convex real projective orbifold is residually finite and virtually torsion-free. 
\end{proposition} 
\begin{proof} 
The holonomy group $\Gamma$ also injectively lifts to a group of 
	projective automorphisms of the domain in 
	$\SLpm$. The lifted group is residually finite by 
	by Malcev \cite{Malcev}. 
	Hence, $\Gamma$ is thus always residually finite. 

Now, Selberg's Lemma \cite{Selberg} shows that the holonomy group is virtually torsion-free.

Suppose that $\orb$ is convex. Then 
$\dev: \torb \ra \RP^n$ a diffeomorphism to a convex domain $\Omega$, which is
simply connected.  Hence, the holonomy homomorphism must be injective. 

\end{proof}
\index{residually finite}


\begin{theorem}[Selberg]\label{prelim-thm-vgood} 
A real projective orbifold $S$ is covered finitely by a real projective manifold $M$
and $S$ is real projectively diffeomorphic to $M/G_1$ for a finite group $G_1$ of real projective automorphisms of $M$. 
An affine orbifold $S$ is covered finitely by an affine manifold $N$, and $S$ is affinely 
diffeomorphic to $N/G_2$ for a finite group $G_2$ of affine automorphisms of $N$. 
Finally, given a two convex real projective or affine orbifold $S_1$ and $S_2$ with isomorphic fundamental groups, 
one is a closed orbifold if and only if the other is closed as well.
\end{theorem} 
\begin{proof} 
Since $\Aff(\mathds{A}^n)$ is a subgroup of the general linear group $\GL(n+1, \bR)$, Selberg's Lemma \cite{Selberg} guarantees the existence of 
a torsion-free subgroup of the holonomy group of finite index,  
and we can choose the group to be a normal subgroup. 
We take the inverse image in the fundamental group. 
 and the second item follows. We obtain the orbifold regular covering by a manifold.


The first item follows similarly by Proposition \ref{prlim-prop-Malcev}. 

For the final item, we can take a torsion-free subgroups and the finite covers of 
$S_1$ and $S_2$ are manifolds serving as $K(G, 1)$-spaces for an identical group $G$. 
Hence, the conclusion follows. \hfill
\SSn {\parfillskip0pt\par}
\end{proof}
\index{Selberg's lemma}

\subsection{A comment on lifting real projective structures and conventions} 
\label{prelim-sub-lifting}

We sharpen Theorem \ref{intro-thm-doubledef}.  
Let $\SL_{-}(n+1, \bR)$ denote the component of $\SL_{\pm}(n+1, \bR)$ that does not contain $\Idd$. 
A projective automorphism $g$ of $\SI^n$ is orientation preserving if and only if $g$ has a matrix in $\SL(n+1, \bR)$. 
For even $n$, the quotient map $ \SL(n+1, \bR) \ra \PGL(n+1, \bR)$ is an isomorphism and 
the map $\SL_{-}(n+1, \bR) \ra \PGL(n+1, \bR)$ is a diffeomorphism for the component 
of $\SL_\pm(n+1, \bR)$ with determinants equal to $-1$. For odd $n$, the quotient map 
$ \SL(n+1, \bR) \ra \PGL(n+1, \bR)$ 
is a $2$ to $1$ covering map onto its image component with deck transformations given by $A \ra \pm A$. 
\index{lifting projective structure}


\begin{theorem}\label{prelim-thm-lifting} 
	Let $M$ be a strongly tame $n$-orbifold.
	Suppose that $h: \pi_1(M) \ra \PGL(n+1, \bR)$ is a holonomy homomorphism of a real projective structure on $M$. 
 \index{real projective structure!lifting}
Then the following statements hold\,{\rm :} 
	\begin{itemize}
		\item Suppose that $M$ is orientable. 
		Then we can lift $h$ to a homomorphism $h': \pi_1(M) \ra \SL(n+1, \bR)$, which is a holonomy homomorphism of 
		the $(\SI^n, \SL_\pm(n+1, \bR))$-structure lifting the real projective structure. 
		\item Suppose that $M$ is not orientable. Then we can lift $h$ to a homomorphism $h': \pi_1(M) \ra \SL_\pm(n+1, \bR)$
		that is the holonomy homomorphism of the $(\SI^n, \SL_\pm(n+1, \bR))$-structure lifting the real projective structure
		such that the condition {\rm ($\ast$)} is satisfied. 
		\begin{itemize}
			\item[({$\ast$})] a deck transformation maps to a negative determinant matrix if and only if it reverses orientations. 
		\end{itemize}
	\end{itemize}
In general a lift $h'$ is unique if we require it to be the holonomy homomorphism 
of the lifted structure. For even $n$, the lifting is unique if we require the condition ($\ast$).
\end{theorem} 
\begin{proof}
For the first part,	recall $\SL(n+1, \bR)$ is the group of orientation-preserving linear automorphisms of $\bR^{n+1}$ and hence is precisely 
	the  group of orientation-preserving projective automorphisms of $\SI^n$. 
	Since the deck transformations of the universal cover $\tilde M$ of the lifted $(\SI^n, \SL_\pm(n+1, \bR))$-orbifold 
	are orientation-preserving, the holonomy of the lift are in $\SL(n+1, \bR)$. 
	We use $h' $ as the holonomy homomorphism of the lifted structure. 
	
	For the second part, we can construct an orientable double cover $M'$ of $M$ by an orientable with an orientation-reversing 
	$\bZ_2$-action of the projective automorphism group generated by an automorphism 
$\phi:M'\ra M'$.
	$\phi$ lifts to $\tilde \phi: \tilde M'\ra \tilde M'$ for the universal covering manifold $\tilde M'= \tilde M$  and hence 
	$h(\tilde \phi) \circ \dev = \dev \circ \tilde \phi$ for  the developing map $\dev$ 
	and the holonomy  \[h(\tilde \phi) \in \SL_-(n+1, \bR).\] 
	Then it follows from the first item since $\dev$ preserves orientations for a given orientation of $\tilde M$. 
	(See p. 143 of Thurston \cite{Thurston97}.)
	
	The proof of the uniqueness is straightforward. 
\end{proof}



\begin{remark}[Convention on using spherical real projective structures] \label{prelim-rem-convention}
	Suppose we are given 
	a convex real projective orbifold of form $\Omega/\Gamma$ for 
	$\Omega$ a convex domain in $\RP^n$ and $\Gamma$ a subgroup of $\PGL(n+1, \bR)$.
	We can also think of $\Omega$ as a domain in $\SI^n$ and 
	$\Gamma \subset \SL_{\pm}(n+1, \bR)$. 
	We will use a convenient one for the purposes ahead.
	\end{remark}


\subsubsection{Convex hulls}\label{prelim-subsub-convh} 

\begin{definition}\label{prelim-defn-convhull} 
	Let $K$ be a subset of a properly convex domain $\Omega$ 
in an affine subspace $\mathds{A}^n$ in $\SI^n$ 
	{\rm (}resp. $\RP^n${\rm ),} the {\em convex hull} $\CH(K)$ of \index{ch@$\CH(\cdot)$}
	$K$ is defined as the smallest convex set containing $K$ in 
$\clo(\Omega) \subset \mathds{A}^n$. 
\end{definition}
The convex hull is well-defined on $\SI^{n}$ as long as $\Omega$ is properly convex
independent of the choice of $\Omega$ or $\mathds{A}^n$. 
For $\RP^{n}$, the convex hull may not be well-defined. 
Once the affine suspace or a properly convex domain containing $K$ is specified, it is
well-defined. 
(Usually it will be clear what $\Omega$ is by context, but we will mention these.) 

Also, it is well-known that each point of the convex hull of a set $K$
has a direction vector equal to a linear sum of at most $n+1$ vectors in
the direction of $K$. 
Hence, the convex hull is a union of $n$-simplices with vertices in 
$K$. 
Also, if $K$ is compact, then its convex hull is also compact
(See Berger \cite{Berger}.) 

\index{convex hull|textbf} 

\begin{lemma} \label{prelim-lem-fch} 
	Let $\Omega$ be a compact convex domain.
	Let $\{K_i\}$ be  a sequence of compact sets in $\Omega$.
	For each $i$, $K_i $ in a properly convex domain.
	Suppose that $\{K_i\}$ geometrically converges to a compact set 
	$K \subset \Omega$, 
	Then $\{\CH(K_i)\}$ geometrically converges to $\CH(K)$.  
	\end{lemma} 
\begin{proof} 
	It suffices to prove for $\SI^n$. 
	We write each element of $\CH(K_i)$ as a finite sum  
	$\llrrparen{\sum_{j=1}^{n+1}\lambda_{i,j} \vec{v}_{i,j}}$ for $\vec{v}_{i,j}$ in 
	the direction of $K_i$ and $\lambda_{i, j} \geq 0$. 
	Proposition \ref{prelim-prop-BP} implies the result.
\hfill	\SnT {\parfillskip0pt\par}
	\end{proof}

\subsubsection{Needed convexity facts} 
We will use the following lemma many times in the monograph. 

\begin{proposition} \label{prelim-prop-closureind} 
	Let $\Omega$ be a properly convex domain in $\SI^n$, and 
	let $\Omega'$ be the image of $\Omega$ under 
	the double covering map $p_{\SI^n}$. 
	Then the restriction $\clo(\Omega) \ra \clo(\Omega')$ is 
	one-to-one and onto. 
	\end{proposition} 
\begin{proof} 
	This follows since we can find 
	an affine subspace $\mathds{A}^n$ containing $\clo(\Omega)$. 
	Since the covering map restricts to a homeomorphism on 
	$\mathds{A}^n$, the result follows. 
	\end{proof}

\begin{lemma}[Chapter 11 of \cite{Thurston80}] \label{prelim-lem-locconv}
	Let $K$ be a closed subset of a convex domain $\Omega$ in 
	$\RP^n$ {\rm (}resp. $\SI^n$\/{\rm )}
	such that each point of $\Bd K$ has a convex neighborhood. 
	Then $K$ is a convex domain. 
\end{lemma} 
\begin{proof} 
	Assume $\Omega \subset \SI^n$. 
	We can connect each pair of points by a broken projective geodesic. 
	Then the local convexity shows that we can reduce the number of geodesic segments 
	by one using triangles. Finally, we obtain a geodesic 
	segment connecting the pair of points. \hfill \SSn  {\parfillskip0pt\par}
\end{proof}

\begin{lemma} \label{prelim-lem-simplexbd}
	Let $\Omega$ be 
	a convex domain in $\RP^n$ {\rm (}resp. in $\SI^n$\/{\rm ).} 
	Let $\sigma$ be a convex domain in $\clo(\Omega) \cap P$ for a subspace $P$. 
	Then either $\sigma \subset  \Bd \Omega$ or $\sigma^o$ is in $\Omega$. 
\end{lemma} 
\begin{proof} 
	Assume $\Omega \subset \SI^n$. 
	Suppose that $\sigma^o$ meets 
	$\Bd \Omega$ and is not entirely contained in it.  
Since $\Omega^o$ is convex, $\sigma \cap \Omega^o$ is a convex open set, possibly empty. 
Suppose that $\sigma \cap \Omega^o$ is not empty. 
	Then we can find a segment $s \subset \sigma^o$ with a point $z$ such that a component $s_1$ of $s-\{z\}$ 
	is in $ \Bd \Omega$ and the other component $s_2$ is disjoint from it.
	 
	We may perturb $s$ in a $2$-dimensional totally geodesic space 
	containing $s$ and 
	such that the new segment $s' \subset \clo(\Omega)$ 
	meets $\Bd \Omega$ only at its interior point. 
	This contradicts the fact that $\Omega$ is convex by Theorem A.2 of \cite{psconv}. \hfill \SnP {\parfillskip0pt\par}
\end{proof}

\section{Affine orbifolds}\label{prelim-sec-affinorb} 

Recall that an orbifold is a topological space stratified by open manifolds
(See Chapter 4 of \cite{Cbook}). 
An affine or projective $n$-orbifold is {\em triangulated} if 
\begin{itemize} 
\item there is a smoothly embedded $n$-cycle composed of geodesic $n$-simplices 
on the compactified orbifold relative to ends by adding an ideal point to a radial end and an ideal boundary component to each totally geodesic ends
\item 
where the interiors of $i$-simplices in the cycle are mutually disjoint and are embedded in
strata of the same or higher dimension. 
\end{itemize} 
\index{orbifold!triangulated|textbf} 

\subsection{Affine suspension constructions} \label{prelim-sub-asusp}
The affine subspace $\bR^{n+1}$ is a dense open subset of $\RP^{n+1}$ 
which is the complement of an $n$-dimensional projective subspace.
Thus, an affine transformation is a restriction of a unique projective automorphism acting on $\bR^{n+1}$. 
The group of affine transformations $\Aff(\mathds{A}^{n+1})$ is isomorphic to the group of projective automorphisms 
acting on $\bR^{n+1}$ by the restriction homomorphism. 

A {\em dilatation} $\gamma$ in an affine subspace 
$\bR^{n+1}$ is a linear transformation with respect 
to an affine coordinate system 
such that all its eigenvalues have norm $> 1$ or $< 1$.
In this setting, $\gamma$ acts as an expanding map in the dynamical sense. 
A {\em scalar dilatation} is a diagonalizable dilatation with a single eigenvalue.
\index{dilatation|textbf} 
\index{dilatation!scalar|textbf} 

An affine $(n+1)$-orbifold $\orb$ is {\em radiant} if $h(\pi_1(\orb))$ fixes a point in $\bR^{n+1}$ for the holonomy homomorphism $h: \pi_1(\orb) \ra \Aff(\mathds{A}^{n+1})$. \index{radiant affine manifold}
A real projective orbifold $\mathcal{O}$ of dimension $n$ 
has a developing map $\dev': \torb \ra \SI^n$ and the holonomy homomorphism 
$h': \pi_1(\orb) \ra \SL_\pm(n+1, \bR)$. 
We regard that $\SI^n$ is embedded as a unit sphere in $\bR^{n+1}$
temporarily. 
We obtain 
a radiant affine $(n+1)$-orbifold by taking $\tilde{\mathcal{O}}$ and $\dev'$ and $h'$: 
Define $D'':\tilde{\mathcal{O}} \times \bR_+ \ra \bR^{n+1}$ by sending $(x, t)$ to $t\dev'(x)$. 
For each element of $\gamma \in \pi_1(\mathcal{O})$, 
we define the transformation $\gamma'$ on $\tilde{\mathcal{O}} \times \bR_+$ from 
\begin{align} 
\gamma'(x,t) = & ( \gamma(x), \theta(\gamma)||h'(\gamma)(t\dev'(x))||) \nonumber \\ 
                         &\hbox{ for a homomorphism } \theta: \pi_1(\orb) \ra \bR_+.
\end{align}
 Also, consider the transformation $S_s:\tilde{\mathcal{O}} \times \bR_+ \ra \tilde{\mathcal{O}} \times \bR_+$ 
sending $(x, t)$ to $(x, st)$ for $s \in \bR_+$. 
Thus, \[\tilde{\mathcal{O}}\times \bR_+/\langle S_\rho,\pi_1(\mathcal{O}) \rangle, \rho\in \bR_+, \rho > 1\] 
is an affine $(n+1)$-orbifold with the fundamental group 
isomorphic to $\pi_1(\mathcal{O}) \times \bZ$ where 
\begin{itemize}
\item the developing map is given by $D''$, and 
\item the holonomy homomorphism is given by $h'$ 
which sends the generator of $\bZ$ to $S_\rho$. 
\end{itemize} 
We call the result the {\em affine suspension} of $\mathcal{O}$, which necessarily is radiant. 
The representation of $\pi_1(\orb) \times \bZ$ with the center $\bZ$ mapped to a scalar dilatation is called an {\em affine suspension} 
of $h$. 
A {\em special affine suspension} is an affine suspension in which $\theta \equiv 1$ holds identically. \index{affine suspension|textbf}  \index{affine suspension!special|textbf}

We deduce that
 \[\tilde{\mathcal{O}}\times \bR_+/\langle \gamma,\pi_1(\mathcal{O}) \rangle\] 
 is an affine orbifold with the fundamental group isomorphic to 
 $\langle \pi_1(\orb), \bZ \rangle$. 

There is a variation called {\em generalized affine suspension}. 
Here we use any $\gamma$ that is a dilatation and normalizes 
$h'(\pi_1(\mathcal{O}))$. We will not use this concept. 
(See Sullivan-Thurston \cite{ST83}, Barbot \cite{Barbot00} and Choi \cite{rdsv} also.)  
\index{affine suspension!generalized|textbf}

\begin{definition}\label{prelim-defn-affs}
We denote by $C(\torb)$ the manifold $\torb \times \bR$ with the structure 
given by $D''$. We say that $C(\torb)$ is the {\em affine suspension} of $\torb$. \index{C@$C(\cdot)$|textbf}
\end{definition} 

Let $S_t:\bR^{n+1} \ra \bR^{n+1}$, given by $\vec{v} \mapsto t\vec{v}$, 
$t \in \bR_+$, be a one-parameter family of dilations fixing a common point. 
A family of self-diffeomorphisms $\Psi_t$, $in \bR$, 
on an affine orbifold $M $ \index{radiant flow diffeomorphism}
lifting to $\hat \Psi_t: \widetilde M  \ra \widetilde M $
such that $D \circ \hat \Psi_t = S_{e^{ct}} \circ D$ for every $t \in \bR$ and some $c \ne 0$
is called a group of {\em radiant flow diffeomorphisms}.  \index{radiant flow diffeomorphism} 


\begin{lemma}\label{prelim-lem-affsus}
Let $\orb$ be a strongly tame real projective $n$-orbifold. 
\begin{itemize}
\item An affine suspension $\mathcal{O}'$ of $\orb$ always 
admits a group $\{\Phi_t|t\in \bR\}$ of radiant flow diffeomorphisms.
Here, $\{\Phi_t\}$ forms a circle, and all flow lines are closed. 
\item Conversely, if there exists a group of radiant flow diffeomorphisms where all orbits are closed and have the homology class $\llrrbracket{\ast\times \SI^1}\in H_1(\orb\times \SI^1)$ with an affine structure, then  $\orb \times \SI^1$ is 
affinely diffeomorphic to one obtained by an affine suspension construction
from a real projective structure on $\orb$. 
\end{itemize} 
\end{lemma} 
\begin{proof} 
	The first item directly fallows by the above construction.
	
The generator of the $\pi_1(\SI^1)$-factor goes to a scalar dilatation since it induces 
the identity map on the space of directions of radial segments from the global fixed point. 
Thus, each closed curve along $\ast \times \SI^1$ 
represents a nontrivial homology. 
The homology direction of the flow equals $\llrrbracket{\ast \times \SI^1} \in \SI(H_1(\orb \times \SI^1; \bR))$. 
By Theorem D of \cite{Fried82}, there exists a connected cross-section homologous 
to \[[\orb \times \ast] \in H_n(\orb\times \SI^1, V \times \SI^1; \bR)\cong H^1(\orb\times \SI^1; \bR)\]
where $V$ is the union of disjoint end neighborhoods of product forms in $\orb$. 
By Theorem C of \cite{Fried82}, any cross-section is isotopic to $\orb \times \ast$.  
The radial flow is transverse to the cross-section isotopic to 
$\orb\times \ast$, and hence $\orb$ admits a real projective 
structure. It follows easily now that $\orb \times \SI^1$ is an affine suspension. (See \cite{Barbot00} for examples.)
\end{proof}

\section{The needed linear algebra} \label{prelim-sec-linearalg}


Here, we collect the linear algebra we need in this monograph. 
A source is a comprehensive book by Hoffman and Kunze \cite{HK71}. 

\begin{definition}\label{prelim-defn-jordan}
	Given an eigenvalue $\lambda$ of an element $g \in \SLnp$, 
	a {\em $\bC$-eigenvector} $\vec v$ is a nonzero vector in 
	\[\bR E_{\lambda}(g) := \bR^{n+1} \cap \big(\ker (g - \lambda I) 
	+ \ker (g - \bar \lambda I)\big), \lambda \ne 0, \Im \lambda \geq 0\]
	A {\em $\bC$-fixed point} is the direction of a $\bC$-eigenvector in
	$\RP^n$ (resp. $\SI^n$ or $\CP^n$). \index{Ce@$\bC$-eigenvector|textbf}
	
	Any element $g$ has a primary decomposition. (See Section 6.8 of \cite{HK71}.)
	Write the minimal polynomial of $g$ as $\prod_{i=1}^{m} (x- \lambda_{i})^{r_{i}}$ for $r_{i} \geq 1$ and 
	mutually distinct complex numbers $\lambda_{1}, \dots, \lambda_{m}$. 
	Define
	\[C_{\lambda_{i}}(g) := \ker (g- \lambda_{i}\Idd)^{r_{i}} \subset \bC^{n+1}\]
	where $r_{i} = r_{j}$ if $\lambda_{i} = \bar \lambda_{j}$. 
	Then the primary decomposition theorem states that
	\[\bC^{n+1} = \bigoplus_{i=1}^{m} C_{\lambda_{i}}(g),\]
which is a canonical decomposition. 
	\index{Clambda@$C_{\lambda}(\cdot)$} 
	\index{primary decomposition theorem} 
	
	A {\em real primary subspace} is the sum 
	$\bR^{n+1} \cap (C_{\lambda}(g) + C_{\bar \lambda}(g))$ for 
	$\lambda$ an eigenvalue of $g$. 
\index{real primary subspace}

	A point $[\vec v], \vec v \in \bR^{n+1},$ is {\em affiliated} with a norm $\mu$ of an eigenvalue if 
	\begin{equation} \label{prelim-eqn-affil} 
	\vec v \in {\mathcal{R}}_{\mu}(g):= \left(\bigoplus_{i \in \{j| |\lambda_{j}| = \mu\}} C_{\lambda_{i}}(g)\right) \cap \bR^{n+1}.
	\end{equation} 
	Let $\mu_{1}, \dots, \mu_{l}$ denote the set of distinct norms of eigenvalues of $g$. 
	We also have $\bR^{n+1} = \bigoplus_{i=1}^{l} {\mathcal{R}}_{\mu_{i}}(g)$. 
	Here, ${\mathcal{R}}_{\mu}(g) \ne \{0\}$ if $\mu$ equals $|\lambda_{i}|$ for some $i$. 
\end{definition}  \index{affiliated|textbf} \index{real primary space|textbf} 
\index{Rmug@${\mathcal{R}}_{\mu}(\cdot)$}

Given a subset $K$ in $\RP^n$ {\rm (}resp. $\SI^n${\rm )} 
and a set of projective automorphims $g_i$,
we say that $g_i(K)$ {\em accumulates to a set $A$ only} 
if accumulations points of any sequence of the form 
$\{g_i(k_i)| k_i \in K\}$ lies in $A$ only.  
\index{accumulate} 

\begin{proposition}\label{prelim-prop-attract} 
	Let $g$ be an element of $\PGL(n+1, \bR)$ {\rm (}resp. $\SLpm${\rm )} acting on
	$\RP^n$ {\rm (}resp. $\SI^n${\rm ).} 
	Let $V$ and $W$ be independent complementary subspaces on which $g$ acts. 
	Suppose that the norm of every eigenvalue of any eigenvector 
	in the direction of 
	$V$ is strictly larger than any norms of the eigenvalues of
	the vectors in the direction of $W$.
	Let $V^S$ be the subspace that is the sum of the $\bC$-eigenspaces of $V$. 
	Then
	\begin{itemize}
		\item for any compact set $K \in \RP^n - \Pi(V) - \Pi(W)$
		{\rm (}resp. $\SI^n - \Pi'(V) -\Pi'(W)${\rm ),} 
		$\{g^n(K)\}$ accumulates only at points in $\Pi(V^S)$
		{\rm (}resp. $\Pi'(V^S)${\rm )} as $n \ra \infty$. 
		\item Let $U$ be a neighborhood of $x$ in $\RP^n - \Pi(V) - \Pi(W)$
		{\rm (}resp. $\SI^n - \Pi'(V) -\Pi'(W)${\rm ).}  There exists an open
		subset $U'$ of $\Pi(V^S)$ {\rm (}resp. $\Pi'(V^S)${\rm )} where each point of $U'$ 
is realized as a limit point of 
		$\{g^n(y)\}$ as $n \ra \infty$ for some $y$ in $U$. 
	\end{itemize} 
\end{proposition} 
\begin{proof} 
	It is sufficient to prove the statements for $\bC^{n+1}$ and $\CP^n$. 
	We express 
the minimal polynomial of $g$ as $\prod_{i=1}^{m} (x- \lambda_{i})^{r_{i}}$ where 
$r_{i} \geq 1$ and 
$\lambda_{1}, \dots, \lambda_{m}$ are	mutually distinct complex numbers. 
	Let $W^{\bC}$ be the complexification of the subspace corresponding 
	to $W$ and $V^{\bC}$ be the one for $V$. 
	Then $C_{\lambda_{i}}(g)$ is a subspace of $W^{\bC}$ or $V^{\bC}$
	by elementary linear algebra.
	
	Now, we write the matrix of $g$ determined only up to $\pm \Idd$ 
	in terms of above primary decomposition
	spaces. Then we write the matrix of $g$ in the Jordan form as an upper 
	triangular matrix. 
	The diagonal terms of the matrix of $g^n$ 
	dominate nondiagonal terms in terms of ratios of the absolute values. 
	Thus the lemma easily follows. 
	
	The last part follows from writing $x$ and $y$ in terms of 
	vectors in directions of $V$ and $W$ and other $g$-invariant subspaces. 
\end{proof}

\subsection{Nilpotent and orthopotent groups}

Let $\bU$ denote a maximal unipotent subgroup of $\SL_\pm(n+1, \bR)$ 
given by upper triangular matrices with diagonal entries equal to $1$.
We let $\bU_{\bC}$ denote the group of 
upper triangular matrices with diagonal entries equal to $1$
in $\SL_{\pm}(n+1, \bC)$. 


Let $\Ort(n+1)$ denote the orthogonal group of $\bR^{n+1}$ with the standard hermitian
inner product.    
\begin{lemma}[Kostant Decomposition]\label{prelim-lem-orth} 
	The matrix of $g \in \Aut(\SI^n)$ can be written under a coordinate system as $k(g) a(g) n(g)$, where 
	$k(g)$ is an element of $\Ort(n+1)$, $a(g)$ is a positive diagonal element, and $n(g)$ is a real unipotent matrix. 
	Also, diagonal elements of $a(g)$ are the norms of eigenvalues of $g$ as elements of $\Aut(\SI^n)$.  The coordinate system does depend on $g$,
and the coordinate vectors are orthonomal ones with respect to the standard 
coordinate system. 
\end{lemma} 
\begin{proof}  
	Recall the Jordan decomposition of the matrix of $g$. 
	Then $g = k'(g) a'(g) n'(g)$ where $k'(g)$ is conjugate to an orthogonal element
	$a'(g)$ is conjugate to a diagonal one and $n'(g)$ is unipotent
	and $k'(g), a'(g), $ and $n'(g)$ mutually commute,
	and each acts on the cyclic decomposition subspaces of $g$.  
	
	Let $\vec{v}_1, \dots, \vec{v}_{n+1}$ denote the basis vectors of $\bC^{n+1}$ 
	that are chosen from the Jordan decomposition of $g$ with the same norms of eigenvalues.
	where $\vec{v}_j \in V^{i_0+1}_\infty$ for $j =1, \dots, i_0+1$. 
	Let $\lambda_j$ denote the norms of eigenvalues associated 
	with $\vec{v_j}$ in the nonincreasing order for $j = 1, \dots, l+1$. 
	Assume that $\lambda_1, \dots, \lambda_l$ are real, and  
	while $\lambda_{l+2j-1} = \bar \lambda_{l+2j}$ for each $j = 1, \dots, (n+1-l)/2$
	and are nonreal.  
	Assume that $\vec{v}_1, \dots, \vec{v}_{l}$ are real vectors. 
	
	Now we fix the standard hermitian inner product on $\bC^{n+1}$. 
	We obtain vectors 
	\[\vec{v}_1' , \dots, \vec{v}_{n+1}'\]
	by applying the Gram-Schmidt orthogonalization process 
	starting from $v_1$ and continuing to $v_{2}, \dots$. 
	
	Then $g$ can be written as an upper triangular matrix with 
	diagonals unchanged. Hence, we factor $g$ as
	$k'(g) a'(g) n'(g)$ where $k'(g)$ is a diagonal with unit complex 
	eigenvalues, $a'(g)$ is a positive diagonal matrix whose 
diagonal entries are the norms of the eigenvalues of $g$, and $n'(g) \in \bU_{\bC}$. 
	We permute the coordinates if necessary such that 
	$k'(g)$ is blocked by $l$ $1\times 1$-matrices and 
	$(n+1 - l)/2$ number of $2\times 2$-matrices with conjugate eigenvalues. 
	For each pair of coordinates corresponding to these $2\times 2$-matrices, we 
	introduce a real coordinate system such that $k'(g)$ is a real rotation 
	on each block with 
	respect to the standard metric. 
	Under the new coordinate system, $a'(g)$ is unchanged and 
	$n'(g)$ is still unipotent. 
	Finally, we can verify that 
	the new coordinate system is real and orthogonal. 
	$n'(g)$ is real since $k'(g)$ and $a'(g)$ are real
	and $g$ has a real matrix under the real coordinate system. 
	(See also Proposition 2.1 of Kostant \cite{Kostant73}.) 
\end{proof}

Recall that all maximal unipotent subgroups are conjugate 
in $\SL_{\pm}(n+1, \bR)$. (See Section 21.3 of Humphreys \cite{Humphreys}.)
We define 
\begin{equation} \label{prelim-eqn-U}
\bU' := \bigcup_{k \in \Ort(n+1)} k \bU k^{-1} 
= \bigcup_{ k \in \SL_{\pm}(n+1, \bR)} k\bU k^{-1}.
\end{equation} 
 The second equality is explained: 
A unipotent group conjugated 
by any element is still unipotent.
  Each maximal unipotent subgroup is characterized by a maximal flag. 
 Each maximal unipotent subgroup is conjugate to 
 a standard lower triangular unipotent group by 
 an orthogonal element in $\Ort(n+1)$ since $\Ort(n+1)$ acts 
 transitively on the maximal flag space. 
Another way to prove this equality is to write each 
$k\in \SL_{\pm}(n+1, \bR)$ by the fixed KAN decomposition.

\begin{remark}\label{prelim-rem-invmetric} 
We use the Riemannian metric that is left-invariant and invariant 
under the conjugations by the compact group $\Ort(n+1)$.
For example, $\tr((g^{-1}dg)^T g^{-1}dg)$ gives us a positive definite Riemannian 
metric left-invariant by $\SL_\pm(n+1, \bR)$ and invariant by the right action of  $\Ort(n+1)$ as Robert Bryant states \cite{Bryant}. 
Any left-invariant metric gives the global 
distance metric that are uniformly bounded above and below by this metric 
since such metrics are determined at the identity.
\end{remark} 
 
\begin{corollary}\label{prelim-cor-bdunip} 
	Suppose that there exists a positive constant $C_1$, and an element 
	$g \in \SL_\pm(n+1, \bR)$, 
	\[ \frac{1}{C_1} \leq \lambda_{n+1}(g) \leq \lambda_1(g) \leq C_1\]
	for the minimal norm $\lambda_{n+1}(g)$ of the eigenvalues of $g$
	and the maximal norm $\lambda_1(g)$ of the eigenvalues of $g$.
	Then $g$ lies within a bounded distance from $\bU'$ with the bound depending only on $C_1$.
Furthermore, if $C_1 = 1$, then $g\in \bU'$. 
\end{corollary} 
\begin{proof} 
	Let us fix an Iwasawa decomposition 
	$\SL_\pm(n+1, \bR) = \Ort(n+1) D_{n+1} \bU$ for a positive diagonal group 
	$D_{n+1}$.
	By Lemma \ref{prelim-lem-orth}, we can find an element $k\in \Ort(n+1)$ such that 
	\[ g  = k k(g) k^{-1} k a(g) k^{-1} k n(g) k^{-1},\] 
	where $ k(g) \in \Ort(n+1), a(g)  \in D^+_n$,  and $k n(g) k^{-1} \in \bU$
where $a(g)$ is a diagonal matrix with
 norms of eigenvalues as the diagonal.
	Then $k k(g) k^{-1} \in \Ort(n+1),$ and $k a(g) k^{-1}$ is uniformly bounded from $\Idd$ by 
	a constant depending only on $C_1$ by assumption.  

The final part follows from \eqref{prelim-eqn-U}. 
\end{proof}

A subset $A$ of a Lie group is of {\em polynomial growth} if the volume 
$A \cap B_R(\Idd)$ of the ball $B_R(\Idd)$ radius $R$ is less than or
equal to a polynomial of $R$. 
\index{polynomial growth}

\begin{lemma} \label{prelim-lem-bU} 
	$\bU'$ is of polynomial growth in terms of the distance from $\Idd$
in terms of its path-metric induced from restricting the left-invariant metric
in $\Aut(\SI^n)$.
The same is also true for any $C$-neighborhood of $\bU'$ for $C> 0$. 
\end{lemma}
\begin{proof}
	Let $\Aut(\SI^n)$ have a left-invariant Riemannian metric that is also invariant under conjugation by $\Ort(n+1)$.  By Remark \ref{prelim-rem-invmetric}, the metric is compatible with any left-invariant metric. 

	Clearly $\bU$ is of polynomial growth by Gromov \cite{Gromov81} since $\bU$ is nilpotent.
	Choose a unipotent element $u \in \bU'$.  We can write $u(s) = \exp(s \vec{u}), s \geq 0$ where 
	$\vec{u}$ is a nilpotent matrix of unit norm. $g(t) := \exp(t \vec{x}), t \geq 0$ for $\vec{x}$ in the Lie 
	algebra of $\Ort(n+1)$ of unit norm. 
	For a family of $g(t) \in \Ort(n+1)$, we define 
	\begin{equation}\label{prelim-eqn-conj}
	u(t, s) = g(t) u(s) g(t)^{-1} =  \exp(s \mathrm{Ad}_{g(t)} \vec{u}).
	\end{equation}
	We compute 
	\[ u(t, s)^{-1} \frac{d u(t, s)}{d t} =u(t, s)^{-1} (\vec{x} u(t, s) - u(t, s) \vec{x}) = (\mathrm{Ad}_{u(t, s)^{-1}} - \Idd)( \vec{x} ).\]
	Since $\vec{u}$ is nilpotent, $\mathrm{Ad}_{u(t, s)^{-1}} - \Idd$ is a polynomial of variables $t, s$. 
	The norm of $d u(t, s)/ dt$ is bounded above by a polynomial in $s$ and $t$. 
	The conjugation orbits of $\Ort(n+1)$ in $\Aut(\SI^n)$ are compact.
	Also, the conjugation by $\Ort(n+1)$ preserves the distances of elements from $\Idd$
	since the left-invariant metric $\mu$ is preserved by conjugation at $\Idd$ by 
$\Ort(n+1)$ by our choice of the metric,  
	and geodesics from $\Idd$ go to geodesics from $\Idd$ of same $\mu$-lengths under 
	the conjugations by \eqref{prelim-eqn-conj}.
	Hence, we obtain a parametrization of $\bU'$ by $\bU$ and $\Ort(n+1)$ where 
	the volume of each orbit of $\Ort(n+1)$ grows polynomially. 
	Since $\bU$ is of polynomial growth, $\bU'$ is of polynomial growth
	in terms of the distance from $\Idd$. 

For a $C$-neighborhood, we can deduce the same results hold for a $C$-neighborhood $N_C$ of $\bU'$: Instead of a ball $B_R$ of radius $R$ in $N_C'$, we take a ball of radius $R$ in $\bU'$ and 
take balls of radius $R$ in the normal bundle of $\bU'$. The union $B'_R(C)$ of the exponential images 
of the balls contain $B_R$. We can find the upper bound to 
the volume of $B'_R(C)$ by equation (2.0.1) of \cite{HK78}. 
since the set of balls in the normal bundles are compact up to left 
translations and so are their exponential images. 


Then the same argument as above will apply. 
\end{proof}

\begin{theorem}[Zassenhaus \cite{Zassenhaus}]\label{prelim-thm-Zassenhaus}
	For every discrete subgroup $G$ of $\GL(n+1, \bR)$, whose elements have the shape
	in a complex basis in $\bC^n$ 
	\[\left(\begin{array}{cccc}
	e^{i\theta_1} & \ast         & \cdots & \ast \\ 
	0     & e^{i\theta_2}    & \cdots & \ast \\ 
	0     & 0    & \ddots & \vdots \\ 
	0     &   0 &  \cdots   & e^{i\theta_{n+1}} 
	\end{array}\right), \]
	there exists a positive number $\eps$, such that all matrices $A\in G$ 
satisfying the inequalities $|e^{i\theta_j} -1 | < \epsilon$ for every $j=1, \dots, n+1$ 
are contained in the radical of the group, i.e., the subgroup $G_u$ of 
	elements of $G$ with only unit eigenvalues. 
\end{theorem}

An element $g$ of $\GL(n+1, \bR)$ (resp. $\PGL(n+1, \bR)$) is said to be {\em unit-norm-eigenvalued} if 
it (resp. its representative) has only eigenvalues of norm $1$. 
A group is unit-norm-eigenvalued if all of its elements are {\em unit-norm-eigenvalued}. 
\index{unit-norm-eigenvalued} 

A subgroup $G$ of $\SL_\pm(n+1, \bR)$ is {\em orthopotent} if 
there is a flag of subspaces $0=Y_0 \subset Y_1 \subset \cdots \subset Y_m = \bR^{n+1}$
preserved by $G$ such that $G$ acts as an orthogonal group 
on $Y_{j+1}/Y_j$ for each $j=0, \dots, m-1$
for some choices of inner-products. 
(See D. Fried \cite{Fried86}.) 
\index{orthopotent|textbf} 

Recall that a {\em distal group} is a linear group whose elements  do not decrease norms of vectors. 
A distal group is orthopotent by  \cite{CG74} or \cite{Moore68}. 

\begin{theorem}\label{prelim-thm-orthopotent} 
	Let $G$ be a unit-norm-eigenvalued subgroup of $\SL_\pm(n+1, \bR)$.
	Then $G$ is orthopotent and distal, and the following hold{\rm :} 
	\begin{itemize}
		\item If $G$ is discrete, then $G$ is virtually unipotent.
		\item If $G$ is a connected Lie group, 
		then $G$ is an extension of a solvable group by a compact group; i.e., 
		$G/S$ is a compact group for a normal solvable group $S$ in $G$.  
		\item If $G$ is contractible, then $G$ is a simply connected 
		solvable Lie group. 
	\end{itemize}
	\end{theorem} 
\begin{proof} 
	By Corollary \ref{prelim-cor-bdunip}, 
	$G$ is in $\bU' \subset \Ort(n+1) \bU \Ort(n+1)$.  
Hence, the group is a distal group and orthopotent. 
	
	Suppose that $G$ is discrete. Then 
	$G$ is of polynomial growth by Lemma \ref{prelim-lem-bU}. 
	By Gromov \cite{Gromov81}, 
	$G$ is virtually nilpotent. 
	
	Choose a finite-index normal nilpotent subgroup $G'$ of $G$.
	Since $G'$ is solvable, Theorem 3.7.3 of \cite{Varadarajan84}
	shows that $G'$ can be put into an upper triangular form for 
	a complex basis. 
    Let $G'_u$ denote the subset of $G'$ with only elements with all eigenvalues equal to $1$. 
    $G'_u$ is a normal subgroup since it is in an upper triangular form. 
	The map $G' \ra G'/G'_u$ factors
	into a map $G' \ra (\SI^1)^n$ by taking the complex eigenvalues. 
	By Theorem \ref{prelim-thm-Zassenhaus}, the image is 
	a discrete subgroup of $(\SI^1)^n$. Hence, $G'/G'_u$ is finite 
	where $G'_u$ is unipotent. 
	(An alternative proof is given in the Remark of page 124 of Jenkins \cite{Jenkins}.) 
	 
	Suppose that $G$ is a connected Lie group.
		Since $|\lambda(g)| =1$  for every eigenvalue $\lambda(g)$ 
	for all $g\in G$, 
	Corollary \ref{prelim-cor-bdunip} shows that
	$G \subset K \bU K$ for a compact Lie group $K$.

	Since $\bU$ is a distal group, 
	$G$ is a distal group, and hence $G$ is orthopotent by 
	\cite{CG74} or \cite{Moore68}. 
	
	

Since $G \subset K \bU K$, it is of polynomial growth by Lemma \ref{prelim-lem-bU}, 
By Corollary 2.1 of Jenkins \cite{Jenkins}, 
$G$ is an extension of a solvable Lie group
by a compact Lie group. 

If $G$ is contractible, then by the second item 
$G$ can only be an extension of a solvable Lie group by 
a finite group. Since $G$ is determined by its Lie algebra, 
$G$ must be solvable by the second item. 
	\end{proof}



%


\subsection{Elements of dividing groups} \label{prelim-sub-semi} 

Suppose that $\Omega$ is an open domain in $\SI^n$ (resp. $\subset \RP^n$) that is 
properly convex but not necessarily strictly convex. 
Let $\Gamma$ be a discrete group of $\SLpm$ (resp. $\PGL(n+1, \bR)$)
 acting on $\Omega$ such that 
$\Omega/\Gamma$ is compact. Then $\Omega/\Gamma$ has the induced orbifold structure
by Proposition \ref{prelim-prop-prdisc}. 


An element of $\Gamma$ is said to be 
{\em elliptic} if it is conjugate to an element 
of a compact subgroup of $\PGL(n+1, \bR)$ or $\SLpm$. 
\index{elliptic|textbf} 


\begin{lemma}\label{prelim-lem-elliptic} 
	Suppose that $\Omega$ is a properly convex domain in 
	$\RP^n$ {\em (}resp. in $\SI^n${\em ),} and $\Gamma$ 
	is a group of projective automorphisms of $\Omega$.
	Suppose that $\Gamma$ is discrete. 
Then an element $g$ of $\Gamma$ is elliptic if and only if $g$ fixes a point of 
$\Omega$ if and only if $g$ is of finite order. 
\end{lemma} 
\begin{proof}
Assume that $\Omega \subset \SI^n$. 
Let $g$ be an elliptic element of $\Gamma$. 
Take a point $x \in \Omega$. Let $\vec{x}$ denote a vector in 
the cone $C(\Omega)\subset \bR^{n+1}$ corresponding to $x$. 
Then the orbit $\{g^n(\vec{x})| n \in \bZ \}$ has a compact closure. 
There is a fixed vector in $C(\Omega)$, which corresponds to 
a fixed point of $\Omega$ by averaging the orbit of $x$. 
 
If $x$ is a point of $\Omega$ fixed by $g$, then $g$ is in 
the stabilizer group of $x$. Since $\Omega/\Gamma$ is an orbifold, 
$g$ is of finite order. 

If $g$ is of finite order, $g$ certainly is elliptic. 
\hfill \SnT {\parfillskip0pt\par}
\end{proof}



%

The {\em multiplicity} of the norm of an eigenvalue of a matrix $g$ 
is the sum of multiplicities of all eigenvalues of that norm. 
We recall the definitions of Benoist \cite{Benoist05}:
For an element $g$ of $\SL_\pm(n+1, \bR)$, we denote by 
$\lambda_1(g), \dots, \lambda_{n+1}(g)$ the nonincreasing 
sequence of the norms of  its eigenvalues with repetitions by their respective 
multiplicities. 
The first one $\lambda_1(g)$ is called the {\em spectral radius} of $g$. 

\begin{itemize}
	\item $g$  is {\em proximal} if $\lambda_1(g)$ has multiplicity one.
	\item  $g$ is {\em positive proximal } if $g$ is proximal and $\lambda_1(g)$ is an eigenvalue of $g$. 
	\item  $g$ is {\em semiproximal} if either
	$\lambda_1(g)$ or $-\lambda_1(g)$ is an eigenvalue of $g$. 
	\item $g$ is {\em positive semiproximal} if 
	$\lambda_1(g)$ is an eigenvalue of $g$. (Definition 3.1 of \cite{Benoist05}.)
\index{proximal}
\index{proximal!positive} 
\index{semiproximal} 
\index{semiproximal!positive} 
\item $g$ is called {\em positive bi-proximal} if 
$g$ and $g^{-1}$ are both positive proximal. 
\item $g$ is called {\em positive bi-semiproximal} if 
$g$ and $g^{-1}$ are both positive semiproximal. 
\end{itemize} 
\index{bi-semiproximal!positive} 
\index{bi-proximal!positive} 
For an element $g$ of $\PGL(n+1, \bR)$, 
we use the same terminology if a lift of $g$ satisfies one of the above. 


Let $g$ be a positive bi-semiproximal element in $\SL_\pm(n+1, \bR)$. 
\begin{itemize}
\item we denote by 
$V_{g}^A:=\left(\oplus_{i=1}^{k(g)} \ker \left(g - \lambda^{(i)}_1(g)\Idd\right)^{m^{(i)}_1}\right) \cap \bR^{n+1}$
where $m^{(i)}_1$ is the multiplicity of the eigenvalue 
$\lambda^{(i)}_1(g)$, $i=1, \dots, k(g)$, of 
the norm $\lambda_1(g)$ 
in the characteristic polynomial of $g$, and  
\item by $V_{g}^R = \left(\oplus_{i=1}^{k'(g)} \ker \left(g - \lambda^{(i)}_n(g)\Idd\right)^{m^{(i)}_{n+1}}\right)\cap  \bR^{n+1} $
where $m^{(i)}_{n+1}$ is the multiplicity of the eigenvalue $\lambda^{(i)}_{n+1}(g)$,
$i=1, \dots, k'(g)$ of the norm $\lambda_{n+1}(g)$.
\end{itemize} 
We denote by $\hat A_g = \llrrparen{V_g^A}  \cap \clo(\Omega)$ and 
$\hat R_g = \llrrparen{V_g^R} \cap \clo(\Omega)$. 

We define
\[A_g:=A \cap \clo(\Omega) \hbox{ and } R_g:=R \cap \clo(\Omega)\] 
for the eigenspace $A$ associated with positive eigenvalue of $g$ having the largest norm, 
and the eigenspace $R$ associated with positive eigenvalue of $g$ having the smallest norm. 

Suppose that $\Omega \subset \RP^n$. 
The sets $A_g, R_g, \hat A_g, \hat R_g$ are defined as the images of the corresponding sets 
under $p_{\SI^n}$. 
By Proposition \ref{prelim-prop-closureind}, the map restricts to homeomorphisms.  

\begin{lemma}[Lemma 3.2 of \cite{Benoist05}] \label{prelim-lem-prox} 
	Suppose that a nonidentity 
	projective automorphism $g$ acts on a properly convex domain. 
Assume that $\lambda_1(g) \ne \lambda_{n+1}(g)$. 
Then 	$g$ is positive bi-semiproximal.
Moreover, $A_g$ and $R_g$ are not empty. 
	\end{lemma} 
\begin{proof} 
Since $g$ acts on a proper cone in $\bR^{n+1}$ corresponding to $\Omega$, 
an eigenvalue corresponding to $A_g$ is positive by Lemma 5.3 of \cite{Benoist05}
and there must be a fixed point associated with it by taking averages of vectors in $V^A_g$, 
and so are one corresponding to $R_g = A_{g^{-1}}$. 
This proves that $g$ is bi-semiproximal. 
\hfill \SnT  {\parfillskip0pt\par}
%
%
\end{proof} 





The following propositions 
are related to Section 2, 3 of \cite{CLT15} where they use somewhat different methods. 
We denote by $\llrrV{\cdot}$ a standard Euclidean norm of a vector space over $\bR$. 
\label{twovert2} 

\begin{lemma} \label{prelim-lem-unitnolower}
	Suppose that $\Omega$ is an open properly convex domain in 
$\SI^n$ {\rm (}resp. $\rpn${\rm).}
	Suppose that an infinite-order 
	element $g$ acts on $\Omega$ with only unit norm of 
	eigenvalues. Then 
		\[\inf_{y\in \Omega}\{d_{\Omega}(y, g(y))| y\in\Omega \} = 0.\]
	\end{lemma} 
\begin{proof} 
Assume that $\Omega \subset \SI^n$. Now, 
$g$ fixes a point $x$ in $\clo(\Omega)$ by the Brouwer fixed-point theorem. 
If $g$ fixes a point in $\Omega$, the result follows.
Assume $x\in\Bd \Omega$ is a fixed point of $g$. 


We prove this by induction. When $\dim \Omega = 0, 1$, this is obvious. 
Suppose that we proved the conclusion when $\dim \Omega = n-1$, $n \geq 2$.
We now assume that $\dim \Omega = n$. 

	
	We may assume that 
	$\Omega\subset \mathds{A}^n$ for some affine space 
	$\mathds{A}^n$ since $\Omega$ is properly convex. 
	We choose a coordinate system where $x$ is the origin of $\mathds{A}^n$. 
	Then $g$ takes a form of a rational map. 
	We denote by $Dg_x$ the linear map that is the differential of $g$ 
	at $x$. 
	
	Let $S_r$ denote the similarity transformation of $\mathds{A}^n$ 
	fixing $x$ with a factor $r> 0$. 
Then we obtain 
	\[S_r \circ g \circ S_{1/r}:S_r(\Omega) \ra S_r(g(\Omega)).\]  
We introduce a coordinate system on $\mathds{A}^n$ where $x$ is the origin. 
	
	Recall that the differential linear map $Dg_x: \bR^n\ra \bR^n$ 
	satisfies
	\[ \lim_{y, u\ra \vec{0}} \frac{\llrrV{ g(y) -g(u) - Dg_x(y-u)}}{\llrrV{y-u}} \ra 0. \]
Hence, 
	\[ \lim_{r \ra \infty}  r\llrrV{ g(S_{1/r}(y)) -g(S_{1/r}(u)) - S_{1/r}Dg_x(y-u)}\ra 0. \]
	Setting $u=\vec{0}$, 
	we obtain 
		\[ \lim_{r \ra \infty} \llrrV{S_r\circ g\circ S_{1/r}(y)  - Dg_x(y)}\ra 0. \]
		We obtain that 
	as $r \ra \infty$, $\{S_r \circ g \circ S_{1/r}\}$ converges to $Dg_x$
uniformy on any open ball around $x$. 
	
	Also, it is easy to show that as $r \ra \infty$, 
	$\{S_r(\Omega)\}$ geometrically converges to a cone 
	$\Omega_{x, \infty}$ with the vertex at $x$ on which
	$Dg_x$ acts. 
	
	
	Let $x_n$ be an affine coordinate function for a sharply-supporting 
	hyperspace of $\Omega$ that takes the value $0$ at $x$. 
For now, any such function will suffice.
	Let ${\mathbf{x}}(t)$ be a projective geodesic with ${\mathbf{x}}(0) =x$
	at $t = 0$ and ${\mathbf{x}}(t) \in \Omega, x_n({\mathbf{x}}(t))= t$ for $t > 0$ 
	and let $\vec{u} = d {\mathbf{x}}(t)/dt \ne 0$ at $t = 0$.
	We assumed in the premise that $g$ is unit-norm-eigenvalued. 
	Then 
	\begin{equation} \label{prelim-eqn-coneD} 
	 \lim_{t \ra 0} d_\Omega(g({\mathbf{x}}(t)), {\mathbf{x}}(t)) = 
	d_{\Omega_{x, \infty}}(Dg_x(\vec{u}), \vec{u})        
	\end{equation} 
	considering $\vec{u}$ as an element of the cone 
$	\Omega_{x, \infty}$: This follows from
	\[  d_{\Omega}(g \circ S_{1/r} \circ S_r({\mathbf{x}}(t)), {\mathbf{x}}(t)) = 
       d_{S_r(\Omega)}(S_r\circ g \circ S_{1/r} ( S_r({\mathbf{x}}(t))), S_r({\mathbf{x}}(t)))  \]
since $S_r:(\Omega, d_{\Omega}) \ra (S_r(\Omega), d_{S_r(\Omega)})$ is an isometry. We set ${\mathbf{x}}(t) = {\mathbf{x}}(1/r)$ and obtain $S_r({\mathbf{x}}(1/r)) \ra \vec{u}$
as $r \ra \infty$. 
Since $S_{r}(\Omega) \ra \Omega_{x, \infty}$ as $r \ra \infty$, 
Lemma \ref{prelim-lem-convmetric} shows that 
 \begin{equation} \label{prelim-eqn-Sr} 
 \{ d_{S_r(\Omega)}(S_r\circ g \circ S_{1/r} ( S_r({\mathbf{x}}(1/r))), S_r({\mathbf{x}}(1/r))) \}
 \ra d_{\Omega_{x, \infty}}(Dg_x(\vec{u}), \vec{u}) \hbox{ as } r \ra \infty.
\end{equation} 

Now, $\Omega_{x, \infty}^o$ is a convex cone
of form $C(U)$ for a convex open domain $U$ in $\Bd \mathds{A}^n$ with the 
origin as the vertex. 
The space $U^\ast$ of 
sharply-supporting hyperspaces of $\Omega_{x, \infty}^o$ at $x$ 
is a compact convex set. 


Since $g$ acts on a ball $U^\ast$, $g$ fixes a point by the Brouwer-fixed-point theorem, which corresponds to a hyperspace. 
Let $P$ be a hyperspace in $\mathds{A}^n$ passing $x$  and 
sharply-supporting $\Omega_{x, \infty}^o$ invariant under $Dg_x$. 

	We now choose the affine coordinate $x_n$ for $P$
	such that $P$ is the zero set of $x_n$. 
There are three possibilities for $\Omega_{x, \infty}$ as a cone of an open convex domain
by Proposition \ref{prelim-prop-classconv}: 
\begin{itemize}
\item  a cone over complete affine space, 
\item a cone over a properly convex domain, or 
\item a cone cover a convex domain that is neither properly convex nor complete affine.
\end{itemize}

First, suppose that $\Omega_{x, \infty}$ is a properly convex cone.
	Projecting  $\Omega_{x, \infty}^o$  to the space 
	$\SI^{n-1}_x$ of rays starting from $x$, 
	we obtain a properly convex open domain 
	$\Omega_1 = R_x(\Omega)$. Here, $\dim \Omega_1 \leq n-1$.

	By the induction hypothesis on
	dimension $\dim \Omega_1$, since the $Dg_x$-action has only one-norm 
	of the eigenvalues, 
	we can find a sequence $\{z_i\}$ in $\Omega_1 $ such that 
	$\{d_{\Omega_1}(Dg_x(z_i), z_i)\} \ra 0$. 
	Since $\Omega_{x, \infty}$ is a properly convex cone in $\mathds{A}^n$, 
	we choose a sequence $u_i \in \Omega_{x, \infty} $ with 
	$x_n(u_i) = 1$ and $u_i$ has the direction of $z_i$ from $x$.  
	Let $\vec{u_i}$ denote the vector in the direction of 
	$\overrightarrow{xz_i}$ on $\mathds{A}^n$ where $x_n(\vec{u}_i)=1$. 
	Since $g$ is unit-norm-eigenvalued, 
	$x_n(Dg_x(\vec{u}_i)) =1$ also. 
	Hence, the geodesic to measure the Hilbert metric 
	from $\vec{u}_i$ to $Dg_x(\vec{u}_i)$ is on the hyperspace given by $x_n=1$ in $\mathds{A}^n$. 
	Let $P_1$ denote the affine subspace given by $x_n = 1$.
	The projection $\Omega_{x, \infty} \cap P_1 \ra U$ from $x$ is a projective diffeomorphism and hence is an isometry. 
Therefore, 
$d_{\Omega_{x, \infty}}(Dg_x(\vec{u}_i), \vec{u}_i) \ra 0$. 

We can find arcs ${\mathbf{x}}_i(t)$ with 
\[ x_n({\mathbf{x}}_i(t)) = t \hbox{ and }  d {\mathbf{x}}_i(t)/dt = \vec{u}_i \hbox{ at } t =0.\] 
Also, we find a sequence of points $\{{\mathbf{x}}_i(t_i)\}$, $t_i \ra 0$, such that 
\[ d_{\Omega}(g({\mathbf{x}}_i(t_i)), {\mathbf{x}}_i(t_i))= 
d_{S_{1/t_i}(\Omega)}(S_{1/t_i}\circ g \circ S_{t_i} ( S_{1/t_i}({\mathbf{x}}_i(t_i))), S_{1/t_i}({\mathbf{x}}_i(t_i))).\]
Since $\{S_{1/t_j}(\clo(\Omega))\} \ra \clo(\Omega_{x, \infty})$, 
and $\{S_{1/t_j}({\mathbf{x}}_i(t_j))\} \ra \vec{u}_i$ as $j \ra \infty$,  
we obtain 
\[ \{d_{\Omega}(g({\mathbf{x}}_i(t_j)), {\mathbf{x}}_i(t_j))\} \ra 
d_{\Omega_{x, \infty}}(Dg_x(\vec{u}_i), \vec{u}_i)\] 
by Lemma \ref{prelim-lem-convmetric}. 

By choosing $j_i$ sufficiently large for each $i$, we obtain
\[ \{d_{\Omega}(g({\mathbf{x}}_i(t_{j_i})), {\mathbf{x}}_i(t_{j_i}))\} \ra 0. \]

Suppose that $\Omega_{x, \infty}$ is not a properly convex cone.
If $\Omega_1$ is a complete affine space, then the pseudometric $d_{\Omega_1}$ is $0$.
Hence, Lemma \ref{prelim-lem-convmetric} shows the result as above. 


Now, $\Omega_1= R_x(\Omega)$ is a convex domain neither properly convex nor complete affine. 
By Proposition \ref{prelim-prop-classconv}, 
such a set is foliated by complete affine spaces of dimension $j$, $0< j< n-1$
or it  is a complete affine space of dimension $n-1$. 
The quotient space $O_x:= \Omega_1/\sim$ with equivalence 
relationship given by complete affine subspace is a properly convex open domain
of dimension $< n$. 
Recall that there is a pseudo-metric 
$d_{\Omega_1}$ on $\Omega_1$. 
Since $\Omega_1$ is not properly convex, we have  $\dim O_x \geq 1$. 
Note that the projection $\pi: \Omega_1 \ra O_x$ is projective 
and $d_{\Omega_1}(y, z) = d_{O_x}(\pi(y), \pi(z))$ for
all $y, z \in \Omega_1$, which is fairly easy to show. 
The differential $Dg_x$ induces a projective map 
$D'g_x: O_x \ra O_x$.
Since $\dim O_x < \dim \Omega$, we have  
by induction a sequence $y_i \in O_x$  
such that $d_{O_x}(y_i, D'g_x(y_i)) \ra 0$ as $i \ra \infty$. 
We take an inverse image $z_i$ in $\Omega_1$ of $y_i$. 
Then $d_{\Omega_1}(z_i, Dg_x(z_i)) = d_{O_x}(y_i, D'g_x(y_i)) \ra 0$
as $i \ra \infty$.  
As above, we obtain the desired result. 
%
%
%
%
\hfill \SnT
\end{proof}


We believe that these were already well known by Benoist and Cooper-Long-Tillmann \cite{CLT15}. 
\begin{proposition} \label{prelim-prop-nonu} 
	Let $\Omega$ be a properly convex domain in $\SI^n$ {\rm (} resp. $\rpn$ {\rm).}
	Suppose that $\Gamma \subset \SL_\pm(n+1, \bR)$ {\rm (} resp. $\PGL(n+1, \bR)$ {\rm)}
is a discrete group acting on $\Omega$ such that $\Omega/\Gamma$ is compact. 
	Let $g$ be a non-torsion and non-identity element
acting on a subspace $Q$ with $Q \cap \Omega \ne \emp$. 
	Then the following hold\/{\rm :}  
\begin{itemize} 
\item  $g$ has at least two distinct positive eigenvalues associated with $Q$
(resp. up to mutliplying by $\pm 1$). 
\item The largest and the smallest norms of $g$ are realized by positive eigenvalues 
greater than $1$ and less than $1$, and $g$ is positive semiproximal. 
\item 	In particular $g|Q$ is not orthopotent, and, hence, $g|Q$ cannot be unipotent.
\end{itemize} 
	\end{proposition} 
\begin{proof} 
Again by Proposition \ref{prelim-prop-prdisc} $\Omega/\Gamma$ has 
a structure of a compact orbifold.
Suppose that $g$ acts with a single norm of eigenvalues on 
a subspace $Q$ with $Q \cap \Omega \ne \emp$. 
Applying Lemma \ref{prelim-lem-unitnolower} where $n$ is replaced by  
the dimension of $Q$, we obtain $0$ as the infimum of the Hilbert 
lengths of closed curves in a compact orbifold $\Omega/\Gamma$. 
Since $\Omega/\Gamma$ is a compact orbifold, there should be a positive 
lower bound. This is a contradiction. 

Lemma  \ref{prelim-lem-prox} and the fact that the product of the norms of eigenvalues is
$1$ prove this result by taking $g$ and $g^{-1}$. 
\hfill \SnP {\parfillskip0pt\par}
\end{proof}

 We generalize Proposition 5.1 of 
Benoist \cite{Benoist04}.
By Theorem \ref{prelim-thm-vgood} following from Selberg's Lemma \cite{Selberg}, 
there is a finite-index subgroup $\Gamma' \subset \Gamma$
whose elements of $\Gamma$ are not elliptic. 
(In fact, a finite manifold cover is enough.)

\begin{theorem}[Benoist \cite{Benoist05}] \label{prelim-thm-semi} 
	Suppose that $\Omega$ is properly convex but not necessarily strictly convex in $\SI^n$
	{\rm (}resp. $\rpn$ {\rm).}
	Let $\Gamma$ be a discrete group acting on $\Omega$ so 
	that $\Omega/\Gamma$ is compact and Hausdorff. 
	Then each nonidentity nontorsion element $g\in\Gamma$ is positive bi-semiproximal
	with the following properties\/{\rm :} 
	\begin{itemize} 
		\item $\lambda_1(g) > 1, \lambda_{n+1}(g) < 1$,
and eigenvalues of norm $\lambda_i(g)$ are positive for $i=1$ or $i=n+1$. 
\item 	$A_g, \hat A_g \subset \Bd \Omega$, 
$R_g, \hat R_g \subset \Bd \Omega$ are nonempty properly convex subsets in
the boundary. 
\item 
$\dim \hat A_g = \dim V_g^A - 1$, $\dim \hat R_g = \dim V_g^R-1$, 
$\dim A_g = \dim \ker (g-\lambda_1 \Idd) -1$, 
and  $\dim R_g = \dim \ker (g-\lambda_{n+1} \Idd) -1$.
\item Let $K$ be a compact set in $\Omega$. 
Then $\{g^i(K)|i \geq 0\}$ has the limit set in $A_g$, 
and  $\{g^i(K)|i < 0\}$ has the limit set in $R_g$.
\end{itemize} 
Furthermore, if $\Omega$ is strictly convex, then 
$A_g = \hat A_g$ is a singleton  in $\Bd \Omega$ 
and $R_g = \hat R_g$ is a singleton in $\Bd \Omega$. In this cse, 
$g$ is positive bi-proximal. 
	\end{theorem} 
\begin{proof} 
By 	Proposition \ref{prelim-prop-nonu}, 
every nonidentity nontorsion element $g$ of $\Gamma$ has an eigenvalue of norm $> 1$. 
By Lemma \ref{prelim-lem-prox}, $g$ is positive bi-semiproximal. 
Also, the eigenvalues of norm $\lambda_i$ for $i=1, n+1$ cannot be negative:
We can consider orbits of generic points, and the convex hull of their limit set cannot 
be properly convex. Since $\Omega$ is properly convex, this cannot happen.

Proposition \ref{prelim-prop-attract} shows that 
$\hat A_g$ contains a limit point of $\{g^i(K)|i > 0 \}$ for any compact set $K \subset \Omega$
and also confirms the dimension of $\hat A_g$.  
Hence, $\hat A_g$ is not empty and $\hat A_g \subset \Bd \Omega$. 
Similarly, $\hat R_g$ is not empty as well. 
Also, the same proposition shows that 
$\hat A_g$ has the same dimension as $V_{g}^A-1$. 
Similarly, taking $g^{-1}$ instead of $g$, we obtain that  
$\hat R_g$ has the same dimension as $V_{g}^R-1$. 
Also, $\hat A_g$ and $\hat R_g$ are convex compact sets since 
they are intersection of $\clo(\Omega)$ with the subspaces corresponding 
to $V_g^A$ and $V_g^R$ respectively. 


We can consider the orbits in $\hat A_g$. The convex hull of each orbit must contain 
a point of $A_g$. Hence $A_g$ is not empty. Similarly $R_g$ is not empty.

Then $A_g$ equals the intersection $\langle V_1\rangle  \cap \clo(\Omega)$ 
for the eigenspace $V_1$ of $g$ corresponding to $\lambda_1(g)$. 
Since $g$ fixes each point of $\langle V_1\rangle$, it follows that 
$A_g$ is a compact convex subset of $\Bd \Omega$. 
Similarly, $R_g$ is also a compact convex subset of $\Bd \Omega$. 


Suppose that $\hat A_g \cap \Omega \ne \emp$. 
Then $g$ acts on the open 
properly convex domain $\hat A_g \cap \Omega $
as a unit-norm-eigenvalued element. 
By Lemma \ref{prelim-lem-unitnolower} applied to $\hat A_g \cap \Omega$,
we obtain a contradiction again by obtaining a sequence of closed 
curves of $d_{\Omega}$-lengths in $\Omega/\Gamma$ converging to $0$ 
which is impossible for a closed orbifold. 
Thus, $\hat A_g \subset \Bd \Omega$. 
As above, it is a compact convex subset. 
Similarly, $\hat R_g$ is a compact convex subset of 
$\Bd \Omega$. 

The third item follows from the dimension of $\hat A_g$ since 
$A_g$ is obtained by taking 
generic orbits in $\hat A_g$ and finding the averages in $\bR^{n+1}$. 
The dimension of $R_g$ is also clear from that of $\hat R_g$. 


Suppose that $\Omega$ is strictly convex. 
Then \[\dim A_g =0, \dim \hat A_g = 0, \dim R_g=0, 
\dim \hat R_g =0\] by the strict convexity. 
Proposition 5.1 of \cite{Benoist04} proves that $g$ is proximal. 
$g^{-1}$ is also proximal by the same proposition.
These are positive proximal since $g$ acts on a proper cone. 
Hence, $g$ is positive bi-proximal. 
\hfill \SSn {\parfillskip0pt\par}
	\end{proof} 

Note here that $\hat A_g$ may properly contain $A_g$ and 
$\hat R_g$ may also properly contain $R_g$.


\subsection{The higher-convergence-group} \label{prelim-sub-convG}

In this section, we only work with $\SI^n$ as this is the only version needed.
Considering $\GL(n+1, \bR)$ as a dense open subspace 
of $M_{n+1}(\bR) -\{O\}$, we can project $\SL_{\pm}(n+1, \bR)$ 
as a dense open subset of $\SI(M_{n+1}(\bR))$. 
We can compactify $\SL_{\pm}(n+1, \bR)$ as 
$\SI(M_{n+1}(\bR))$. Denote by $\llrrparen{g}$ the equivalence class of 
$g \in \SL_{\pm}(n+1, \bR)$. 
\index{higher convergence group} 

\index{compactification of $\SL_{\pm}(n+1, \bR)$} 

\begin{theorem}[The higher-convergence-group property] \label{prelim-thm-converg}
	Let $g_i$ be an unbounded sequence of projective automorphisms of 
	a properly convex domain $\Omega$ in $\SI^n$. 
	We consider $g_i \in \SL_\pm(n+1, \bR)$ according 
	to convention \ref{intro-rem-SL}. 
	Then we can choose a subsequence of $\{\llrrparen{g_i}\}$, 
	$g_i \in \SL_\pm(n+1, \bR)$, 
	converging to $\llrrparen{g_\infty}$ in $\SI(M_{n+1}(\bR))$ for 
	$g_\infty \in M_{n+1}(\bR)$ where the following hold\/ {\rm :} 
	\begin{itemize} 
		\item	$\llrrparen{g_\infty}$ is undefined on 
		$\SI(\ker g_\infty)$ and the range is 
		$\SI(\Ima g_\infty)$. 
		\item $\dim \SI(\ker g_\infty) + \dim \SI(\Ima g_\infty) = n-1.$
		\item For every compact subset $K$ of 
		$\SI^n -  \SI(\ker g_\infty)$, $\{g_i(K)\}\ra K_\infty$ 
		for a subset $K_\infty$ of $\SI(\Ima g_\infty)$. 
		
		\item Given a convergent subsequence $\{g_i\}$ as above, 
		$\{\llrrparen{g g_i}\}$ also converges to $\llrrparen{g g_\infty}$ and 
		$\SI(\ker g g_\infty) = \SI(\ker g_\infty)$ 
		and $\SI(\Ima g g_\infty) = g\SI(\Ima g_\infty)$
		\item Given a convergent subsequence $\{g_i\}$ as above, 
$\{\llrrparen{g_i g}\}$ also converges to $\llrrparen{g_\infty g}$ and 
		\[\SI(\ker g_\infty g) = g^{-1}(\SI(\ker g_\infty))
		\hbox{ and } \SI(\Ima g_\infty g) = \SI(\Ima g_\infty).\] 
	\end{itemize}
\end{theorem} 
\begin{proof} 
	Since $\SI(M_{n+1}(\bR))$ is compact, we can find a subsequence 
	of $g_i$ converging to an element  $\llrrparen{g_\infty}$. 
	The second item is the consequence of the rank and nullity 
	of $g_\infty$. 
	The third item follows from considering the compact-open topology 
	of maps and $g_i$ divided by its maximal norm of the matrix entries.
	
	The final two items are straightforward.  	
\end{proof} 


\begin{lemma} \label{prelim-lem-convcomp} 
	$\llrrparen{g_\infty}$ can be obtained by 
	taking the limit of $g_i/m(g_i)$ in $M_{n+1}(\bR)$ first 
	and then taking the direction 
	where $m(g_i)$ is the maximal norm of 
	elements of $g_i$ in the matrix form of $g_i$. 
	\end{lemma} 
\begin{proof} 
	This follows since $g_i/m(g_i)$ does not go to zero. 
	\end{proof} 


\index{higher convergence group} 

This definition was suggested by Goldman.
\begin{definition} \label{prelim-defn-convlimit} 
	An unbounded sequence $\{g_i\}$, $g_i \in \SL_{\pm}(n+1, \bR)$, 
	such that $\{\llrrparen{g_i}\}$ is convergent in $\SI(M_{n+1}(\bR))$ 
	is called a {\em convergence sequence}.
	In this case, $g_\infty\in M_{n+1}(\bR)$ is called a {\em convergence limit},
	determined only up to a positive constant. 
	The element $\llrrparen{g_\infty} \in \SI(M_{n+1}(\bR))$ where 
	$\{\llrrparen{g_i}\} \ra \llrrparen{g_\infty}$ is called a {\em convergence limit}.
	
	We may also apply this definition for $\PGL(n+1, \bR)$. 
	An unbounded sequence $\{g_i\}$, $g_i \in \PGL(n+1, \bR)$
	such that $\{[g_i]\}$ is convergent in $\bP(M_{n+1}(\bR))$
	is called a {\em convergence sequence}. 
	Also, the element $[g_\infty] \in \bP(M_{n+1}(\bR))$ is 
	$\{[g_i]\} \ra [g_\infty]$ is called a {\em convergence limit}. 
	
	We can provide additional interpretations: 
	We use the KTK-decomposition (or polar decomposition) 
	of Cartan for $\SL_{\pm}(n+1, \bR)$. 
(See Theorem 7.39 of \cite{Knapp} and \cite{Helgason}.) 
	We may express $g_i = k_i d_i \hat k_i^{-1}$ where 
	$k_i, \hat k_i\in \Ort(n+1, \bR)$ and $d_i$ is a positive diagonal matrix
	with a nonincreasing set of elements 
	\[a_{1,i} \geq a_{2,i} \geq \cdots \geq a_{n+1,i}. \]
\index{KTK-decomposition|textbf}

	Let $\SI([1,m])$ denote the subspace spanned by the points corresponding 
to $\vec{e}_1, \dots, \vec{e}_m$, 
	and let $\SI([m+1, n+1])$ denote the subspace spanned by
the points corresponding to $\vec{e}_{m+1}, \dots, \vec{e}_{n+1}$. 
	We assume that $\{k_i\}$ converges to $k_\infty$,  
	$\{\hat k_i\}$ converges to $\hat k_\infty$, 
	and \hfil \break 
 $\{[a_{1,i}, a_{2,i}, \cdots, a_{n+1,i}]\}$ is convergent in $\RP^n$. 
	We will further require this for convergence sequences.  
	\index{so@$\SI([1,m])$|textbf }
	\index{se@$\SI([m+1,n+1])$|textbf }
	
	For a sequence $\{g_i\}$ in $\PGL(n+1, \bR)$, we may also write $g_i = k_i d_i \hat k_i^{-1}$
	where $d_i$ is represented by positive diagonal matrices as above. Then we require as above. 
\end{definition} 
This of course generalizes the concept of the convergence sequence, without the second set of 
requirements above, for 
$\PSL(2, \bR)$ as given by Tukia (see \cite{AC96}). 

Given a convergence sequence $\{g_i\}, g_i \in \Aut(\SI^n)$, 
we define 
\begin{align} 
\hat A_\ast(\{g_i\}) := \SI(\Ima g_\infty)  \\
\hat N_\ast(\{g_i\}) := \SI(\ker g_\infty) \\ 
A_\ast(\{g_i\}) := \SI(\Ima g_\infty) \cap \clo(\Omega) \\
N_\ast(\{g_i\}) := \SI(\ker g_\infty) \cap \clo(\Omega)
\end{align} 
\index{aa@$A_\ast(\cdot)$|textbf}
\index{na@$N_\ast(\cdot)$|textbf}
\index{aa@$\hat A_\ast(\cdot)$|textbf}
\index{na@$\hat N_\ast(\cdot)$|textbf}
\index{convergence sequence|textbf} 




%

For a matrix $A$, we denote by $|A|$ the maximum of the norms of entries of $A$. 
\label{twoverts} 
Let $U$ be an orthogonal matrix in $\Ort(n+1, \bR)$. 
Then we obtain 
\begin{equation} \label{prelim-eqn-comporth} 
\frac{1}{n+1}|A| \leq |AU| \leq (n+1)|A|.
\end{equation} 
The second inequality follows since the entries of $AU$ are dot products of rows of $A$ with 
elements of $U$ whose entries are bounded above by $1$ and below by  $-1$. 
Also, we can multiply $U^{-1}$ to $AU$ to obtain the first inequality. 
Hence, we obtain for $g = k D \hat k^{-1}$ for 
$k, \hat k \in \Ort(n+1, \bR)$ and $D$ positive diagonal as above. 
\begin{equation} \label{prelim-eqn-gD}
\frac{1}{(n+1)^2} |D| \leq |g| \leq (n+1)^2 |D|
\end{equation}


Recall Definition \ref{prelim-defn-jordan}, we then obtain
\begin{theorem} \label{prelim-thm-convGU} 
	Let  $\{g_i\}, g_i \in \Aut(\SI^n)$, be a convergence sequence.
	We consider $g_i \in \SL_\pm(n+1, \bR)$ according 
	to convention \ref{intro-rem-SL}. 
We write $g_i = k_i d_i \hat k_i^{-1}$ for $k_i, \hat k_i \in \Ort(n+1, \bR)$ and 
positive diaongal matrix $d_i$ with nonincreasing sequence of positive diagonal 
entries $a_{1, i}, \dots, a_{n+1, i}$. 
	Then we may assume that the following hold
	up to a choice of subsequence of $\{g_i\}$\/{\rm :}
	\begin{itemize} 
	\item there exists $m_a$, $1 \leq m_a < n+1$, where 
	$\{a_{j,i}/a_{1, i}\} \ra 0$ for $j > m_a$ 
	and $a_{j,i}/a_{1,i} > \eps$ for $j \leq m_a$ for a uniform $\eps > 0$. 
	\item there exists $m_r$, $1\leq m_a < m_r \leq n+1$, where 
	$a_{j,i}/a_{n+1,i} < C$ for $j \geq m_r$ for a uniform $C> 1$,  
	and $\{a_{j,i}/a_{n+1,i}\} \ra \infty$ for $j < m_r$.
	\item 
	$\hat N_\ast(\{g_i\})$ is the geometric limit of $\hat k_i(\SI([m_a+1,n+1]))$.
	\item 
	$\hat A_\ast(\{g_i\})$ is the geometric limit of $A^p(g_i)=k_i(\SI([1,m_a]))$.
	\item $\{g_i^{-1}\}$ is also a convergence sequence up to a choice of 
	subsequences, and 
	$\hat A_\ast(\{g_i^{-1}\}) \subset \hat N_\ast(\{g_i\})$. 
	\end{itemize} 
	\end{theorem} 
\begin{proof} 
	We choose a subsequence such that $m_a$ and $m_r$ are defined respectively, 
	and $\{k_i\}, \{\hat k_i\}$ form convergent sequences. 
	We denote $D_\infty $ as the limit of $\{D_i/|D_i|\}$ and 
	$k_\infty$ and $\hat k_\infty$ as the respective limits of $\{k_i\}$ and $\{\hat k_i\}$. 
Then we obtain by \eqref{prelim-eqn-gD}, 
\[ \frac{1}{(n+1)^3} |k_\infty\circ  D_\infty \circ \hat k_\infty^{-1} (\vec{v})| \leq |g_\infty (\vec{v})| \leq 
(n+1)^3  |k_\infty\circ D_\infty \circ \hat k_\infty^{-1}(\vec{v})|    \]
for every $\vec{v} \in \bR^{n+1}$. 
Thus, 
$k_\infty\circ  D_\infty \circ \hat k_\infty^{-1} (\vec{v}) =0$
if and only if $g_\infty(\vec{v}) = 0$
Thus, null spaces of 
$k_\infty\circ  D_\infty \circ \hat k_\infty^{-1}$ and $g_\infty$  coincide. 
The image of these maps also coincide since we are taking the matrix limits. 
We obtain that 
\begin{multline*} 
\SI(\Ima g_\infty) = 
k_\infty \SI(\Ima D_\infty) = k_\infty \SI([1, m_a])
\hbox{ and } \\ 
\SI(\ker g_\infty) = \hat k_\infty(\SI(\ker D_\infty)
= \hat k_\infty((\SI([m_{a}+1, n+1]))).
\end{multline*}
Hence, the first four items follow. 

%

	The last item follows from considering the third and fourth items
	and using the fact that $m_r \geq m_a$. 
	\end{proof} 


When $m_a$ and $m_r$ exist for $\{g_i\}$, and 
$\{k_i\}$ and $\{\hat k_i\}$ for some choices of $k_i$ and $\hat k_i$ for each $g_i$ 
are convergent for a convergence sequence, 
we say that $g_i$ are {\em set-convergent}.
\index{set-convergent}  
\index{convergent!set-convergent} 

We define for each $i$, 
\begin{equation*}
F^p(g_i) := k_i \SI([1,m_r -1]), 
\hbox{ and } R^p(g_i):= \hat k_i\SI([m_r, n+1]).
\end{equation*}
%
We define 
$\hat R_\ast(\{g_i\})$ as the geometric limit of 
$\{R^p(g_i)\}$, and 
$\hat F_\ast(\{g_i\})$ as the geometric limit of 
$\{F^p(g_i)\}$. 
We also define
\[R_\ast(\{g_i\}):=\hat R_\ast(\{g_i\}) \cap \clo(\Omega),
F_\ast(\{g_i\}): = \hat F_\ast(\{g_i\}) \cap \clo(\Omega).\]
\index{ap@$A^p(g_i)$|textbf } 
\index{rp@$R^p(g_i)$|textbf }
\index{ra@$R_\ast(\{g_i\})$|textbf } 
\index{ra@$\hat R_\ast(\{g_i\})$|textbf } 
\index{fp@$F^p(g_i)$|textbf } 
\index{fa@$F_\ast(\{g_i\})$|textbf } 
\index{fa@$\hat F_\ast(\{g_i\})$|textbf } 

\begin{lemma} \label{prelim-lem-Aginverse} 
	Suppose that $\{g_i\}$ and $\{g_i^{-1}\}$ are set-convergent sequences. Then 
	$\hat R_\ast (\{g_i\}) = \hat A_\ast(\{g_i^{-1}\})$ and 
	$\hat F_\ast(\{g_i\}) = \hat N_{\ast}(\{g_i^{-1}\})$. 
\end{lemma} 
\begin{proof} 
	For	$g_i = k_i d_i \hat k_i^{-1}$, we have
	$g_i^{-1}= \hat k_i^{-1} d_i^{-1} k_i$. 
	Hence, $A^p(g_i^{-1}) = \hat k_i \SI([m_r, n+1]) = R^p(g_i)$. 
	We also have $F^p(g_i^{-1}) = \hat k_i \SI([m_a+1, n+1])$. 
	The result follows.    
\end{proof} 

\begin{lemma}\label{prelim-lem-Ag}
	We also have 
	$\hat A_\ast(\{g^i\}) \subset \hat A_g$ and 
	$\hat R_\ast(\{g^i\}) \subset \hat R_g$ for positive bi-semiproximal element $g$
	with $\lambda_1(g)> 1$.  
\end{lemma} 
\begin{proof} 
	Let $V^A_g$ denote $\mathcal{R}_{\lambda_1(g)}$ in $\bR^{n+1}$
	which is a $g$-invariant subspace from \eqref{prelim-eqn-affil}. 
There is a complementary $g$ 
which is a direct sum of $\mathcal{R}_{\mu}(g)$ for $\mu < \lambda_1(g)$. 
	Then we use Proposition \ref{prelim-prop-attract}
	applied to $\Pi'(V^A_g)$ and $\Pi'(N^A_g)$.
	
	For the second part, we use $g^{-1}$ and argue using obvious facts that 
	$\hat R_g = \hat A_{g^{-1}}$ and 
	$\hat A_\ast(\{ g^{-i}\})
	= \hat R_\ast(\{g^i\})$. 
\end{proof}

\begin{proposition} \label{prelim-prop-Aint} 
We obtain that $A_\ast(\{g_i\})$ contains an open subset of $\hat A_\ast(\{g_i\})$
	and hence $\langle A_\ast(\{g_i\})\rangle = \hat A_\ast(\{g_i\})$. 
	Also, $R_\ast(\{g_i\})$ contains an open subset of $\hat R_\ast(\{g_i\})$
	and hence $\langle R_\ast(\{g_i\})\rangle = \hat R_\ast(\{g_i\})$.
\end{proposition}
\begin{proof}
	We write $g_i = k_i D_i \hat k_i^{-1}$. 
	By Theorem \ref{prelim-thm-convGU},  
	$\hat A_\ast(\{g_i\})$ is the geometric limit of 
	$k_i(\SI[1, m_a])$ for some $m_a$ as above. 
We have 
	$g_i(U) = k_i D_i(V_i)$ for an open set $U \subset \torb$ where
	$V_i = \hat k_i^{-1}(U)$. 
	Since $\hat k_i^{-1}$ is a $\bdd$-isometry, 
	$V_i$ is an open set containing a closed ball $B_i$ of fixed radius
	$\epsilon$. 
	$\{D_i\}$ converges to a diagonal matrix $D_\infty$. 
	Without loss of generality, we may assume that 
	\[\{B_i\} \ra B_\infty, B_i\cap B_\infty \supset B \hbox{ for sufficiently large } i\]
where $B_\infty$ is a ball of a fixed radius $\epsilon$ and $B$ is a fixed ball
	of radius $\eps/2$. 
	Then $\{D_i(B)\} \ra D_\infty(B) \subset \SI([1, m_a])$.  
	Here, $D_\infty(B)$ is a subset of $\SI([1, m_a])$
	containing an open set. 
	Since $\{D_i(B_i)\}$ geometrically converges to a subset containing
	$D_\infty(B)$, up to a choice of subsequence, 
	$\{k_i D_i(V_i)\}$ geometrically converges to a subset 
	containing $k_\infty D_\infty(B)$
	by Lemma \ref{prelim-lem-geoconv}. 
	
	For the second part, we use the sequence $g_i^{-1} = \hat k_i D_i^{-1} k_i^{-1}$
	and Lemma \ref{prelim-lem-Aginverse}. 
\end{proof}

\begin{lemma} \label{prelim-lem-low}
	Suppose that $\Gamma$ acts properly discontinuously on a properly convex 
	open domain $\Omega$, and $\{g_i\}$ is a set-convergent sequence in $\Gamma$. 
	Suppose that $\{g_i\}$ is not bounded in $\SL_{\pm}(n+1, \bR)$, and 
	is a set-convergent sequence. Then
	the following hold{\rm :} 
	\begin{enumerate} 
	\item[(i)]  $\hat R_\ast(\{g_i\}) \cap \Omega = \emp$, 
	\item[(ii)] $\hat A_\ast(\{g_i\}) \cap \Omega = \emp$,
	\item[(iii)] $\hat F_\ast(\{g_i\}) \cap \Omega = \emp$, and
		\item[(iv)] $\hat N_\ast(\{g_i\}) \cap \Omega = \emp$.
	\end{enumerate} 
\end{lemma} 
\begin{proof} 
	(ii) 
	Suppose not. 
    Since $\hat A_\ast(\{g_i\}) \cap \Omega \ne \emp$,
    $A_{\ast}(\{g_i\})$ meets $\Omega$. 
    Since $A_{\ast}(\{g_i\})$ 
    is contained in the set of limit points of $g_i(x)$ for $x \in \Omega$, 
    the proper discontinuity of the action of $\Gamma$ shows 
    that $A_{\ast}(\{g_i\})$ does not meet $\Omega$. 
	
	(iv) For each $x$ in $\Omega$, a fixed ball $B$ in $\Omega$ 
	centered at $x$ does not meet
$\hat k_i(\SI([m_a+1, n+1]))$ for infinitely many $i$. 
	Otherwise $\{g_i^m(B)\}$ converges to a nonproperly convex set in 
	$\clo(\Omega)$ as $m\ra \infty$,  a contradiction.
	Hence, the second item follows.

	The remaining items follow by replacing $g_i$ to $g_i^{-1}$
	and using Lemma \ref{prelim-lem-Aginverse}.  
\end{proof} 

\begin{theorem} \label{prelim-thm-AR} 
	Let $\{g_i\}$ be a set-convergence sequence in $\Gamma$ acting 
	properly discontinuously on a properly convex domain $\Omega$. 
	Then 
	\begin{alignat}{3}
	& A_\ast(\{g_i\}) && = \hat A_\ast(\{g_i\}) \cap \clo(\Omega) && = \hat A_\ast(\{g_i\}) \cap \Bd\Omega , \\
	&  N_\ast(\{g_i\}) && = \hat N_\ast(\{g_i\}) \cap \clo(\Omega) && = \hat N_\ast(\{g_i\}) \cap \Bd \Omega, \\
	&  R_\ast(\{g_i\}) && = \hat R_\ast(\{g_i\}) \cap \clo(\Omega) && = \hat R_\ast(\{g_i\}) \cap \Bd \Omega, \\
	& F_\ast(\{g_i\}) && = \hat F_\ast(\{g_i\}) \cap \clo(\Omega)  && = \hat F_\ast(\{g_i\}) \cap \Bd \Omega 
	 \end{alignat} 
	are subsets of $\Bd \Omega$ and 
	they are nonempty sets. 
	Also, we have 
	\begin{align} 
	\hat A_\ast(\{g_i\}) & \subset \hat F_\ast(\{g_i\}), &
	A_\ast(\{g_i\}) & \subset F_\ast(\{g_i\}), \\
	\hat R_\ast(\{g_i\}) & \subset \hat N_\ast(\{g_i\}), &
	R_\ast(\{g_i\}) & \subset N_\ast(\{g_i\}).
	\end{align} 
	
	\end{theorem} 
\begin{proof} 
By Lemma \ref{prelim-lem-low}, we only need to show  
the respective sets are not empty. 
By the third item of Theorem \ref{prelim-thm-converg}, 
a point $x$ in $\hat A_\ast(\{g_i\}) \cap \clo(\Omega)$ is a limit of 
$g_i(y)$ for some $y \in \Omega$. Since $\Gamma$ acts 
properly discontinuously, $x \not\in \Omega$ 
and $x \in \Bd \Omega$. 
By taking $\{g_i^{-1}\}$, 
we obtain $\hat R_\ast(\{g_i\}) \cap \clo(\Omega)\ne \emp$.
Since $\hat R_\ast(\{g_i\}) \subset \hat N_\ast(\{g_i\})$ 
and $\hat A_\ast(\{g_i\}) \subset \hat F_\ast(\{g_i\})$, 
the remaining part follows. 
The last collections follow directly from definitions. 
\end{proof}

\begin{proposition} \label{prelim-prop-Aast} 
	For an automorphism $g$ of $\Omega$, and 
	a set-convergence sequence $\{g_i\}$, 
	the following hold\/{\rm :} 
	\begin{align} \label{prelim-eqn-gAgR} 
	\hat A_\ast(\{g g_i\}) = g(\hat A_\ast(\{g_i\})), 
	\hat A_\ast(\{g_i g\}) = \hat A_\ast(\{g_i\}), \\
	\hat N_\ast(\{g g_i\}) = \hat N_\ast(\{g_i\}),
	\hat N_\ast(\{g_i g\}) = g^{-1}(\hat N_\ast(\{g_i\})),\\
		\hat F_\ast(\{g g_i\}) = g(\hat F_\ast(\{g_i\})), 
	\hat F_\ast(\{g_i g\}) = \hat F_\ast(\{g_i\}), \\
	\hat R_\ast(\{g g_i\}) = \hat R_\ast(\{g_i\}),
	\hat R_\ast(\{g_i g\}) = g^{-1}(\hat R_\ast(\{g_i\})). 
	\end{align}
\end{proposition} 
\begin{proof} 
	%
	The fourth and fifth items of Theorem \ref{prelim-thm-converg} imply
	the first and second equations here. 
	The third line follows from the second line by Lemma \ref{prelim-lem-Aginverse}. 
	Also, the fourth line follows from the first line by Lemma \ref{prelim-lem-Aginverse}. 
\end{proof}

Of course, there are $\rpn$-versions of the results here. However, 
we do not state these.

\section{Convexity of real projective orbifolds}\label{prelim-sec-conv}

\subsection{Convexity}


In the following, a zero-dimensional sphere 
$\SI^0_\infty$ denotes a pair of antipodal points. 
\begin{proposition}\label{prelim-prop-projconv} $ $ 
	\begin{itemize}
		\item A real projective $n$-orbifold is convex if and only if the developing map sends 
		the universal cover to a convex domain in $\rpn$ {\rm (}resp. $\SI^n${\rm ).} 
		\item A real projective $n$-orbifold is properly convex if and only if the developing map sends 
		the universal cover to a precompact properly convex open domain in an affine patch of $\rpn$
		{\rm (}resp. $\SI^{n}${\rm ).}
		\item If a convex real projective $n$-orbifold is neither properly convex nor complete affine, then
		its holonomy is reducible in $\PGL(n+1, \bR)$ {\rm (}resp. $\SLnp${\rm ).} In this case, $\torb$ is foliated 
		by affine subspaces $l$ of dimension $i$ with the common boundary $\clo(l) - l$ equal to 
		a fixed subspace  $\RP^{i-1}_{\infty}$ 
		{\rm (}resp. $\SI^{i-1}_{\infty}${\rm )} in $\Bd \torb$.  
		Furthermore, this holds for any convex domain in $\rpn$ {\rm (}resp. $\SI^n${\rm )} 
		and the projective group action on it. 
	\end{itemize}
\end{proposition}
\begin{proof} 
	We prove for $\SI^n$ first. 
Since the universal cover is projectively diffeomorphic to a convex open domain, 
the developing map must be an embedding.
	The converse is also trivial. 
	(See Proposition A.2 of \cite{psconv}. ) 
	
	The second follows immediately. 

	For the final item, 
	a convex subset of $\SI^n$ is a convex subset of an affine subspace $\mathds{A}^n$, isomorphic to an affine space, 
	which is the interior of a hemisphere $H$. 
	We may assume that $D^o \ne \emp$ by restricting to a spanning subspace of $D$ 
	in $\SI^n$. 
	Let $D$ be a convex subset of $H^o$. 
	If $D$ is not properly convex, the closure $\clo(D')$ must be of 
the form $\SI^{i_0} \ast K$ for a properly convex domain by Proposition \ref{prelim-prop-classconv}. 
	
	Since $\SI^{i_0}$ must be holonomy invariant, 
	the holonomy group is reducible. 
	
 For the $\RP^n$-version, 
 we use the double covering map $p_{\SI^n}$ mapping an open
	hemisphere to an affine subspace. \hfill	\SnT  {\parfillskip0pt\par}
\end{proof} 


We define $\SI^{n\ast}:= \SI(\bR^{n+1\ast}-\{O\})$, i.e, the sphericalization of
the dual space $\bR^{n+1\ast}$. 
\label{Snast} 
\index{Snast@$\SI^{n\ast}$|textbf} 

\begin{theorem}\label{prelim-thm-sweepU}
	Suppose that $G$ acts on a convex domain 
	$\Omega$ in $\RP^n$ {\em (}resp. in $\SI^n${\em ),} 
	such that $\Omega/G$ is a compact orbifold. 
	Then if $G$ has only elements with unit-norm eigenvalues, 
	$\Omega$ is complete affine. 
\end{theorem} 
\begin{proof} 
%
%
%
Theorem \ref{prelim-thm-orthopotent} shows that $G$ is orthopotent
and has a unipotent subgroup $U$ of finite index. 
	A unipotent group has a global fixed point in $\SI^{n\ast}$, and so does 
	$U^\ast$. 
	Thus, there exists a hyperspace $P$ in $\SI^n$
	where $U$ acts on. $P \cap \clo(\Omega) \ne \emp$ and $P \cap \Omega =\emp$ 
	by Proposition \ref{prelim-prop-sweep}. 
	Thus, $\Omega$ lies in an affine subspace $\mathds{A}$ bounded by $P$. 
	Also, $G$ acts on $\mathds{A}$ 
	as a group of affine transformations since
	every projective action on a complete affine space is 
	affine. (See Berger \cite{Berger}.)
	Discrete orthopotent affine groups are distal by Theorem \ref{prelim-thm-orthopotent}.
	Lemma 2 of Fried \cite{Fried86} implies the conclusion.
\hfill	\SnT
	
%
\end{proof}

\subsection{The flexibility of boundary} 

The following lemma gives us some flexibility of boundary.  
A smooth hypersurface embedded in a real projective manifold is called {\em strictly convex} if 
there is no open line segment in the hypersurface containing 
a point of the hypersurface. \index{convex!strictly!hypersurface|textbf}

\begin{lemma} \label{prelim-lem-pushing} 
	Let $M$ be a strongly tame properly convex real projective orbifold with strictly convex $\partial M$.
	We can modify $\partial M$ inward $M$ and the result bound 
	a strongly tame or compact properly convex real projective orbifold $M'$ with strictly convex $\partial M'$
\end{lemma} 
\begin{proof}
	Let $\Omega$ be a properly convex domain covering $M$. 
	We may assume that $\Omega \subset \SI^n$. 
	We may modify $M$ by pushing $\partial M$ inward. 
	We take an arbitrary inward vector field defined on a tubular neighborhood of $\partial M$.
	(See Section 4.4 of \cite{Cbook} for the definition of tubular neighborhoods.)
	We use the flow defined by this vector field
to modify $\partial M$. By the $C^{2}$-convexity condition, 
	for sufficiently small change, $\partial M$ still remains strictly convex and smooth. 
	Let the resulting compact $n$-orbifold be denoted by $M'$.
	$M'$ is covered by a subdomain $\Omega'$ in $\Omega$. 
	
	Since $M'$ is a compact suborbifold of $M$, 
	$\Omega'$ is a properly embedded domain in $\Omega$ 
	and thus, $\Bd \Omega' \cap \Omega = \partial \Omega'$. 
	$\partial \Omega'$ is a strictly convex hypersurface since so is $\partial M'$.
	This means that $\Omega'$ is locally convex. 
	A locally convex closed subset of a convex domain is convex 
	by Lemma \ref{prelim-lem-locconv}. 
	
	Hence, $\Omega'$ is convex and hence is properly convex since it is 
a subset of a properly convex domain. 
	So is $M'$. 
\hfill	\SSn 
\end{proof}

\begin{remark}\label{prelim-rem-bdconv}
	Thus, by choosing a point in the interior, we may assume without loss of generality that a strictly convex 
	boundary component can be pushed out 
	to a strictly convex boundary component. 
\end{remark}

\subsection{The Benoist theory.} \label{prelim-sub-ben}

In the late 1990s, Benoist more or less completed the theory 
of the dividing action as initiated by Benz\'ecri, Vinberg, Koszul, Vey, and so on
in the series of papers \cite{Benoist04}, \cite{Benoist03}, \cite{Benoist05}, \cite{Benoist062}, \cite{Benoist00}, \cite{Benoist97}.
We will generalize these to sweeping actions in Lemma \ref{prelim-lem-domainIn}. 
The comprehensive theory will aid us much in this paper. 
(We will use ``dividing'' instead of ``divisible'' by using the active form.)

\begin{proposition}[Corollary 2.13 \cite{Benoist05}]\label{prelim-prop-Benoist}  
Suppose that $\Gamma$ is a discrete subgroup  of $\SLn$ {\rm (}resp. $\PGL(n, \bR)${\rm ),} 
$n \geq 2$,
acts on a properly convex $(n-1)$-dimensional open domain $\Omega$ in $\SI^{n-1}$ {\rm (}resp. $\RP^{n-1}${\rm )} 
such that $\Omega/\Gamma$ is a compact orbifold. Then the following statements are equivalent. 
\begin{itemize} 
\item Every finite-index subgroup of $\Gamma$ has a finite center. 
 \item Every finite-index subgroup of $\Gamma$ has a trivial center. 
\item Every finite-index subgroup of $\Gamma$ is irreducible in $\SLn$ {\rm (}resp. in $\PGL(n, \bR)${\rm ).} 
That is, $\Gamma$ is strongly irreducible. 
\item The Zariski closure of $\Gamma$ is semisimple as a linear group. 
\item $\Gamma$ does not contain an infinite nilpotent normal subgroup. 
\item $\Gamma$ does not contain an infinite abelian normal subgroup.
\end{itemize}
\end{proposition}
\begin{proof}
Corollary 2.13 of \cite{Benoist05} considers $\PGL(n, \bR)$ and $\RP^{n-1}$. 
However, the corresponding results for $\SI^{n-1}$ follows from this since we can always lift 
a properly convex domain in $\RP^{n-1}$ to one $\Omega$ in $\SI^{n-1}$ and 
the group to one in $\SLn$ acting on $\Omega$
by Theorem \ref{prelim-thm-lifting}. 
\end{proof}

\label{zariski} 

The center of a group $G$ is denoted by $\bZ(G)$.
\index{z@$\bZ(\cdot)$ }
A {\em virtual center} of a group $G$ 
is a subgroup of $G$ whose elements commute with 
a finite-index subgroup of $G$.
A group with properties above is said to be a group with 
a {\em trivial virtual center}. 
\index{virtual center|textbf} 
\label{center}

Given a group $G$ acting on a manifold $N$, the action is said to be {\em cocompact} if 
there is a compact subset $K$ of $N$ such that $\bigcup_{g\in G} gK= N$.
\index{action!cocompact} 

\begin{theorem}[Theorem 1.1 of \cite{Benoist05}] \label{prelim-thm-Benoist} 
Suppose that a discrete subgroup $\Gamma$ of $\SLn$ {\rm (}resp. $\PGL(n, \bR)${\rm )}, $n \geq 2$,  with trivial virtual center acts on 
a properly convex $(n-1)$-dimensional open domain $\Omega \subset \SI^{n-1}$ so 
that $\Omega/\Gamma$ is a compact orbifold. 
Then every representation of the component of $\Hom(\Gamma, \SLn)$ {\rm (}resp. $\Hom(\Gamma, \PGL(n, \bR))${\rm ) }
containing the inclusion representation also acts on a properly convex $(n-1)$-dimensional open domain whose quotient under the representation is a compact orbifold. 
\end{theorem}
For $n=3$, Choi-Goldman  \cite{CG93} proved this result. 
When $\Omega/\Gamma$ admits a hyperbolic structure and $n=4$, 
Inkang Kim \cite{Kimi01} proved this result simultaneously for a union of components.





A group $G$ acts on a space $X$ {\em cocompactly} if there is a compact subset $Y$ of $X$ 
such that $X = \bigcup_{g\in G} gY$. 
\index{action!cocompact} 
Provided $X$ is locally compact, this condition is equivalent to the condition
that $X/G$ is compact. 

\begin{proposition}[Theorem 1.1. of Benoist \cite{Benoist03}] \label{prelim-prop-Ben2}  
Suppose that a discrete subgroup $\Gamma$ of $\SLn$ 
{\rm (}resp. $\PGL(n, \bR)${\rm ),} $n \geq 2$,
acts on a properly convex $(n-1)$-dimensional open domain $\Omega$ in $\SI^{n-1}$ 
{\rm (}resp, $\RP^{n-1}${\rm )} 
such that $\Omega/\Gamma$ is a compact orbifold.  
Then 
\begin{itemize}
\item $\Omega$ is projectively diffeomorphic to the interior of 
a strict join $K:=K_1 * \cdots * K_{l_0}$ where $K_i$ is a properly convex open domain of dimension $n_i \geq 0$ in 
the subspace $\SI^{n_i}$ in $\SI^n$ {\rm (}resp. $\RP^{n_i}$ in $\RP^n${\rm )}.
$K_{i}$ corresponds to a convex cone $C_i \subset \bR^{n_i+1}$ for each $i$. 
\item $\Omega$ is the interior of the image of $C_1 + \cdots + C_{l_{0}}$. 
\item Let $\Gamma_{i}'$ be the image of $\Gamma'$ to $K_i$ for
the restriction map of the subgroup $\Gamma'$ of $\pi_{1}(\Sigma)$ acting on each $K_{j}$, 
$j=1, \dots, l_0$.  
We denote by $\Gamma_i$ an arbitrary extension of $\Gamma_i'$ by
requiring it to act trivially on $K_j$ for $j \ne i$
and to have $1$ as the eigenvalue associated with vectors in
their directions. 
\item The subgroup corresponding to $\bR^{l_0-1}$ acts trivially on each $K_j$
and forms a positive diagonalizable matrix group. 
\item The fundamental group 
$\pi_1(\Sigma)$ is virtually isomorphic to a subgroup of 
$\bR^{l_0-1} \times \Gamma_1 \times \cdots \times \Gamma_{l_0}$ for 
$(l_0 -1) + \sum_{i=1}^{l_{0}} n_i = n$. 
\item $\pi_{1}(\Sigma)$ acts on $K^o$ cocompactly and discretely
and in a semisimple manner {\em (}Theorem 3 of Vey \cite{Vey}{\em )}.
\item 
The Zariski closure of $\Gamma'$ equals 
$\bR^{l_0-1}\times G_1\times \cdots \times G_{l_0}$. 
Each $\Gamma_j$ acts on $K_j^o$ cocompactly, and $G_{j}$ is an simple 
Lie group {\em (}Remark after Theorem 1.1 of \cite{Benoist03}{\em ),} 
and $G_{j}$ acts trivially on $K_m$ for $m \ne j$. 
\item A virtual center of $\pi_1(\Sigma)$ of maximal rank is isomorphic to $\bZ^{l_0-1}$ corresponding to the subgroup of $\bR^{l_0-1}$.
 {\em  (}Proposition 4.4 of \cite{Benoist03}\/{\em .)}
\end{itemize} 
\end{proposition}

%

We will often indicate by $\bZ^{l_0-1}$  the virtual center
of $\pi_1(\Sigma)$. 
See Example 5.5.3 of Morris \cite{Morris15} for a group that acts  properly on a product of two hyperbolic spaces but restricts to a nondiscrete group for each factor space.

\begin{corollary} \label{prelim-cor-Benoist}
	Assume as in Proposition \ref{prelim-prop-Ben2}. 
	Then every normal solvable 
	subgroup of a finite-index subgroup $\Gamma'$ of 
	$\Gamma$ is virtually central in $\Gamma$.
	\end{corollary} 
\begin{proof} 
Assume that $\Gamma$ is in $\Aut(\SI^n)$. 
	If $\Gamma$ is virtually abelian, this is obvious. 
	
	Suppose that $\Omega$ is properly convex. 
	Let $G$ be a normal solvable subgroup of a finite-index subgroup $\Gamma'$ of 
	$\Gamma$. We may assume without loss of 
	generality that $\Gamma'$ acts on each $K_i$ by taking 
	a further finite-index subgroup and replacing $G$ by a finite-index subgroup of $G$. 
Now, $G$ is a normal solvable subgroup of 
the Zariski closure $\mathcal{Z}(\Gamma')$. 
By Theorem 1.1 of \cite{Benoist03}, 
$\mathcal{Z}(\Gamma')$ equals $G_1\times \cdots \times G_l \times \bR_+^{l-1}$
and $K = K_1\ast \cdots \ast K_l$ 
where $G_i$ is reductive and 
the following holds:  
\begin{itemize} 
	\item if $K_i$ is homogeneous, then $G_i$ is simple and $G_i$ 
	is commensurable with $\Aut(K_i)$. 
	\item Otherwise, $K_i^o$ is divisible and $G_i$ is a union of 
	components of $\SL_{\pm}(V_i)$
	\end{itemize} 
The image of $G$ into $G_i$ by the restriction homomorphism to $K_i$ 
is a normal solvable subgroup of $G_i$. 
Since $G_i$ is virtually simple, the image is a finite group. 
Hence, $G$ must be virtually a subgroup of the diagonalizable group 
$\bR_+^{l-1}$ and is virtually central in $\Gamma'$. 
\hfill\SSn {\parfillskip0pt\par}
\end{proof}


If $l_0 = 1$, $\Gamma$ is strongly irreducible as shown by Benoist. \index{irreducible!strongly}
However, the images of these groups will be subgroups of $\PGL(m, \bR)$ and 
$\SL_\pm(m, \bR)$ for $m \leq n$. 
If $l_0 > 1$, we say that such an image in $\Gamma$ is {\em virtually factorizable}. 
Otherwise, such an image  is a {\em non-virtually-factorizable} group. 

\index{factorizable!virtually|textbf} 
\index{factorizable!non-virtually|textbf} 

An action of a projective group $G$ on a properly convex domain 
$\Omega$ is {\em sweeping} if the $G$-action is cocompact but $\Omega/G$ is not required 
to be Hausdorff. A dividing action is sweeping.
\index{action!sweeping} 

Recall that a {\em commutant} $H$ of a group acting on a properly convex 
domain is a maximal diagonalizable group commuting with the group.
(See Vey \cite{Vey}.)
\index{commutant|textbf}

We have a useful theorem: 
\begin{theorem}[Proposition 3 of Vey \cite{Vey}] \label{prelim-thm-Vey} 
	Suppose that a projective group $\Gamma$ acts on a properly convex open domain $\Omega$ in $\RP^n$ {\em (}resp. in $\SI^n${\em ),}  with a sweeping action.
	Then $\Omega$ equals a convex hull of the orbit 
	$\Gamma(x)$ for any $x \in \Omega$. 
	\end{theorem}

We generalize Proposition \ref{prelim-prop-Ben2}. 
\begin{proposition} \label{prelim-prop-sweep}
Suppose that a projective group $G$ acts on an $n$-dimensional properly convex open domain $\Omega$ in $\SI^n$ {\em (}resp. $\RP^n${\em )} as 
a sweeping action. Then the following hold\/{\em :}
	\begin{itemize} 
	\item Let $L$ be any proper subspace where $G$ acts on. 
Then $L \cap \clo(\Omega) \ne \emp$
but $L \cap \Omega = \emp$. 
\item If $G$ acts on a compact convex set $K$, 
then $K$ must meet $\clo(\Omega) \cup \mathcal{A}(\clo(\Omega))$. 
\item Suppose that $G$ is semisimple. 
Then 
all the items up to the last one 
in the conclusion of Proposition \ref{prelim-prop-Ben2}
with $G$ replacing $\pi_1(\tilde E)$ 
without the discreteness assumption hold.
In particular, $\clo(\Omega) = K_1 \ast \cdots \ast K_{l_0}$ for properly convex sets 
$K_1, \dots, K_{l_0}$. 
\item Suppose that $G$ is semisimple. 
Then the closure of $G$  has a virtual center containing a group of 
diagonalizable projective automorphisms that is isomorphic to $\bZ^{l_0-1}$ acting 
trivially on each $K_i$. 
\end{itemize} 
\end{proposition} 
\begin{proof} 
	Assume $\Omega \subset \SI^n$. 
	Suppose that $L \cap \clo(\Omega)= \emp$. 
	Then there is a lower bound for the $\bdd$-distance from 
	$\Bd \Omega$ to $L$. 
	Let $x \in \Omega$. 
	We denote the space of 
	oriented maximal open segments containing $x$ and ending 
	at a point of $L$ by $L_{\Omega, x}$. This is a set homeomorphic to $\SI^{\dim L}$. 
	
	Let $l_+$ denote the endpoint of $L \cap l$ ahead of $x$. 
	Let $l_{\Omega, 0}$ denote the endpoint of $l \cap \Omega$ ahead of 
	$x$, and $l_{\Omega, 1}$ denote the endpoint of $l \cap \Omega$ after  
	$x$.
	We define a function $f: \Omega \ra (0, \infty)$ 
	given by 
	\[ f(x) = \inf\{\log(l_+,l_{\Omega, 0}, x,l_{\Omega, 1})| l \in L_{\Omega, x}  \} \]
	where the logarithm measures the Hilbert distance between 
	$l_{\Omega, 0}$ and $x$ on the properly convex segment with endpoints $l_+$ 
	and $l_{\Omega, 1}$. 
Since the map from the space of oriented maximal segments in $\Omega$ to the set of their 
endpoints is continuous, and $\bigcup_{x\in K}L(\Omega, x)$ is compact for each compact 
subset $K$ of $\Omega$, this is a continuous positive function. 

Also, we can choose a subsequence $x_i \ra \Bd \Omega$ such that $f(x_i) \ra 0$: 
Choose a maximal line in $\Omega$ and choose $x_i$ on it converging to one of 
its endpoints.

	Since $f(g(x)) = f(x)$ for all $x \in \Omega$ and $g \in G$, $f$ induces 
	a continuous map $\bar f: \Omega/G \ra (0, \infty)$. 
	Here, $\bar f$ can take values arbitrarily to $0$ as one wishes. 
	This contradicts the compactness of $\Omega/G$. 
	
	Suppose that $L \cap \Omega \ne \emp$. Then $G$ acts on
	the convex domain $L\cap \Omega$ open in $L$. 
	We define a function $f: \Omega \ra [0, \infty)$ given 
	by measuring the Hilbert distance from $L \cap \Omega$. 
	Then $f(x) \ra \infty$ as $x \ra \Bd \Omega - L$. 
	Again $f(g(x)) = f(x)$ for all $x\in \Omega, g\in G$. 
	This induces $\bar f: \Omega/G \ra [0, \infty)$. 
	Since $\bar f$ can take arbitrarily small positive values, 
	this contradicts the compactness of $\Omega/G$. 

For the second item, suppose that such a set $K$ exists. 
$K$ and $\mathcal{A}(K)$ are disjoint from 
$\clo(\Omega)$. For $x \in \Omega$, 
we define $K_{\Omega, x}$ as the space of 
oriented open segments containing $x$ and ending in $K$ and $\mathcal{A}(K)$.  
We define $l_+$ as the first point of $K \cap l$ ahead of $x$. 
Then a similar argument as in the above proof applies and we obtain a contradiction.

Now, we go to the third item.	
Let $G$ have a $G$-invariant decomposition $\bR^n=V_1\oplus \cdots \oplus V_{l_0}$ 
where $G$ acts irreducibly. 
This item follows from Lemma 2 of \cite{Vey} since
any decomposition of $\bR^n$ gives rise to a diagonalizable 
commutant of rank $l_0$. 

For the fourth item, we prove the case where
$G$ has a $G$-invariant decomposition $\bR^n= V_1 \oplus V_2$.
Then by the second item, $G$ acts on $K = K_1 \ast K_2$ where $K_i$ is a
properly convex domain in $\SI(V_i)$ for $i = 1, 2$. 
$G$ acts cocompactly on $K^o$
and $G$ is a subgroup of $G_1 \times G_2 \times \bR_+$ where 
$G_i$ is isomorphic to $G|K_i$ extended to act trivially on $K_{i+1}$
with $G_i| V_{i+1} = \Idd$ with the indices modulo $2$. 

The closure $\bar G$ of $G$ in $\Aut(K)$ is a subgroup of 
$\bar G_1\times \bar G_2 \times \bR_+$ 
for the closure $\bar G_i$ of $G_i$ in $\Aut(K_i)$ for $i=1,2$.  
\[\bar G \subset \{ (g_1, g_2, r)| g_i \in \bar G_i, i=1, 2, r \in \bR_+\}.\]
For a fixed pair $(g_1, g_2)$, if there are more than one associated 
$r$, then we obtain by taking differences 
that $(\Idd, \Idd, r)$ is in the group $\bar G$ for 
$r \ne 1$. This implies that $\bar G$ contains a nontrivial subgroup of $\bR_+$. 


Otherwise, $\bar G$ is in the graph of homomorphism 
$\lambda:\bar G_1\times \bar G_2 \ra \bR_+$. 
An orbit of this action on the manifold $\bF K_1^o\times \bF K_2^o\times \bR_+$ 
lies in the orbit of the image of $\lambda$. 
Hence, each orbit of a compact set meets $(y_1, y_2) \times (0, 1)$ for $y_i \in \bF K_i^o$, $i=1,2$, at  a compact set. 
Thus, we do not have a cocompact action. 

Furthermore, if we have a $G$-invariant decomposition 
$K_1\ast \cdots \ast K_m$, we can use the decomposition 
$K_1 \ast (K_2 \ast \cdots \ast K_m)$. 
Now, we use induction to obtain the result. 
\hfill \SnT {\parfillskip0pt\par}
%
\end{proof}

\begin{proposition}\label{prelim-prop-convatt}
Assume as in Proposition \ref{prelim-prop-Ben2}. Then $K$ is the closure of 
the convex hull of $\bigcup_{g \in \bZ^{l_0-1}} A_g$ for the attracting 
limit set $A_g$ of $g$. 
Also, for any partial join $\hat K:= K_{i_1} \ast \dots \ast K_{i_j}$
for a subcollection $\{i_1, \dots, i_j\}$, the closure of 
the convex hull of $\bigcup_{g \in \bZ^{l_0-1}} A_g \cap \hat K$ equals $\hat K$. 
	\end{proposition} 
\begin{proof} 
	We take a finite-index normal subgroup $\Gamma'$ of $\Gamma$ such that $\bZ^{l_0-1}$ is 
	the center of $\Gamma'$. Using Theorem \ref{prelim-thm-vgood}, we may assume that 
	$\Gamma'$ is torsion-free. 
	Note that $k A_g$ for any $k \in \Gamma'$ equals $A_{kgk^{-1}} = A_g$ 
	since $k g k^{-1} = g$. Thus, $\Gamma'$ acts on 
	$\bigcup_{g \in \bZ^{l_0-1}} A_g$ since it is a $\Gamma'$-invariant set. 
	The interior $C$ of the convex hull of $\bigcup_{g \in \bZ^{l_0-1}} A_g$ is a subdomain in 
	$K^o$. Since $C/\Gamma' \ra K^o/\Gamma'$ is a homotopy equivalence of closed manifolds,
	we obtain $C = K^o$ and $\clo(C) = K$ by Lemma 
	\ref{prelim-lem-domainIn}.
	
	For the second part, if the closure of the convex hull of $\bigcup_{g \in \bZ^{l_0-1}} A_g \cap \hat K$
	is a proper subset of $\hat K$, then the closure of the convex hull of $\bigcup_{g \in \bZ^{l_0-1}} A_g$
	is a proper subset of $K$. This is a contradiction.  
%
\hfill	\SnP {\parfillskip0pt\par}
	\end{proof}

%

\begin{theorem}[Kobayashi \cite{Kobpaper}]\label{prelim-thm-Kobayashi} 
	Suppose that a closed real projective $n$-orbifold has a developing map into 
	a properly convex domain $D$ in $\RP^n$ {\em (}resp. in $\SI^n${\em )}. Then the orbifold
	is projectively diffeomorphic to $\Omega/\Gamma$ for 
	the holonomy group $\Gamma$ and for the unique nonempty $\Gamma$-invariant open set 
$\Omega$ in $D^o$, and $\Omega$ must be a properly convex domain. 
	\end{theorem} 
\begin{proof} 
	This follows since all maximal segments in $D$ are 
	of $\bdd$-length $\leq \pi-\eps_0$ for a uniform $\eps_0>0$. 
	Hence, the Kobayashi metric is well-defined proving that the orbifold 
	is properly convex. Hence, the developing image $\Omega$ is a convex open domain 
by \cite{Kobpaper}. Clearly, $\Omega$ is $\Gamma$-invariant. 

Let $U$ be any $\Gamma$-invariant open set in $D$. Then the convex 
hull $U'$ of $U\cup \Omega$ is a convex open domain containing $\Omega$. 
Again, $\Gamma$ acts properly discontinuously on $U'$ since there is a Hilbert metric on $U'$. 
There is a homotopy equivalence $\Omega/\Gamma \ra U'/\Gamma$.
This implies $U'=\Omega = U$. 
\hfill \SnP {\parfillskip0pt\par}
	\end{proof} 

\begin{lemma}[Domains for holonomy]\label{prelim-lem-domainIn} 
	Suppose that $\Omega$ is an open 
	domain  in an open hemisphere in $\SI^{n-1}$ 
	{\rm (}resp. in $\RP^{n-1}${\rm )} where a projective group 
	$\Gamma$ acts on such that $\Omega/\Gamma$ is a closed orbifold. 
	Suppose that $\Gamma$ acts on a compact properly convex domain $K$ where $K^o\ne \emp$.
	Then 
	\begin{itemize} 
		\item $K^o = R(\Omega)$ where $R$ is a diagonalizable projective automorphism 
		commuting with a finite-index subgroup of $\Gamma$ with eigenvalues $\pm 1$ only
		and is a composition of reflections commuting with one-another.  
		\item In fact $K = K_1 \ast \cdots \ast K_k$ where 
		$K_j = K\cap P_j$, $J=1, \dots, k$, for 
		a virtually invariant subspace $P_j$ of $\Gamma$ where $R$ equals 
		$\Idd$ or $\mathcal{A}$. 
		\item  $\Omega$ must be properly convex. 
		\item If $K^o$ meets with $\Omega$, then $K^o=\Omega$. 
	\end{itemize} 
\end{lemma} 
\begin{proof} 
	Suppose that $\Omega \subset \SI^{n-1}$. 
	We prove by induction on dimensions. 
	For $\SI^1$, it is clear. 
	
	By Theorem \ref{prelim-thm-vgood}, there is a torsion-free finite-index subgroup 
	$\Gamma'$ in $\Gamma$.  
	Suppose that $\Omega \cap K^o \ne \emp$. 
	Then $(\Omega\cap K^o)/\Gamma'$ is homotopy equivalent to 
	$\Omega/\Gamma'$, a closed manifold. Hence, $\Omega\cap K^o = \Omega$ and 
	$\Omega \subset K^o$. Similarly, $K^o\subset \Omega$. 
	We obtain $K^o = \Omega$. 
	Also, if $\Omega \cap \mathcal{A}(K^o) \ne \emp$, 
	then $\mathcal{A}(K^o) = \Omega^o$.  
	The lemma has been proved for this case. 
	
	Proposition \ref{prelim-prop-prdisc} and the second part of 
Theorem \ref{prelim-thm-vgood} show that 
	$K^o/\Gamma$ is again a closed orbifold. 
	Suppose that $\Gamma$ is not virtually factorizable with respect to $K$. 
	Then $\Gamma$ is strongly irreducible by Benoist \cite{Benoist03}. 
	$K$ contains the attracting fixed points of bi-semiproximal element $g$ of 
	$\Gamma$. This implies that $\clo(\Omega)\cap \clo(K) \ne \emp$ or 
	$\clo(\Omega) \cap \mathcal{A}(\clo(K))\ne \emp$ 
	since $\Omega$ contains a generic point of $\SI^n$. 
	As stated in the above paragraph, we may suppose that 
	the intersection is in $\Bd \Omega \cap \Bd K$ or 
	$\Bd \Omega \cap \mathcal{A}(\Bd K)$. 
	Then this is a compact convex 
	set invariant under $\Gamma$. Hence, $\Gamma$ is reducible, a contradiction.

	Now, suppose that $\Gamma$ is virtually factorizable. Then 
	there exists a diagonalizable free abelian group $D$ of rank $k-1$ for some $k \geq 2$ 
	in the virtual center of $\Gamma$ by Proposition \ref{prelim-prop-Ben2}. 
	$D$ acts trivially on a finite set of subspaces
	$P_1, \dots, P_k$  by Proposition	\ref{prelim-prop-Ben2}
that are minimal invariant subspaces.
	Since $\Gamma$ permutes these subspaces, a torsion-free 
	finite-index subgroup $\Gamma'$ of $\Gamma$ acts on each of 
	$P_1, \dots, P_k$. 
	We denote $\hat P_j := P_1\ast \cdots \ast P_{j-1} \ast P_{j+2} \ast \cdots \ast P_k$.  
	Then $\Omega$ is disjoint from $\hat P_j$ for each $j$ since otherwise
	$(\hat P_j\cap \Omega)/\Gamma' \ra \Omega/\Gamma'$ is a homotopy equivalence
	of different dimensional manifolds. 
	
	However, $P_j \cap \clo(\Omega)\ne \emp$ since we can choose a 
	sequence $\{g_i\}$ of elements $g_i \in D$ such that the associated eigenvalues for 
	$P_j$ go to zero and the other eigenvalues go to infinity while their 
	ratios are uniformly bounded as $i \ra \infty$. 
	
	Again, define $K_j :=K \cap P_j$. 
	Since $K$ is properly convex, $K_1\ast \cdots \ast K_k \subset K$. 
	Since the action of $\Gamma$ on $K^o$ is cocompact and proper, 
	Proposition \ref{prelim-prop-Ben2} shows that $K = K_1\ast \cdots \ast K_k$. 
	We have a projection for $K^o \ra K_j^o$ for each $j$ obtained from the join
	structure. 
	Then the action of $\Gamma$ on $K_j$ is cocompact 
	since otherwise $K^o/\Gamma$ cannot be compact.
	Also, the action of $\Gamma'$ on $P_j$ is irreducible by Benoist \cite{Benoist03}. 
The second item has been proved.

	
		We can find a sequence in $D$ converging to a 
	projection $\Pi_j$ to each $P_j$ with the undefined space $\hat P^j$ in the sense of 
Section \ref{prelim-sub-convG}. 
		We define the domains 
	$\Omega_j := \clo(\Pi_j(\Omega))$ in $P_j$.  
	Since $\Omega$ is in an open subset in a hemisphere, 
	there exists a convex hull of $\clo(\Omega)$, and 
	hence so has $\Omega_j$ for each $j=1, \dots m$. 
%
%
Then $\Omega_j$ is properly convex
by the third item of Proposition \ref{prelim-prop-projconv} and the irreducibilty of the action in each factor $K_j$ 
in Proposition \ref{prelim-prop-Ben2}. 	
%
	Hence, $\Omega$ is in a properly convex domain 
	$\Omega_1 \ast \cdots \ast \Omega_m$. 
	
	By Theorem \ref{prelim-thm-Kobayashi} of Kobayashi, 
	$\Omega$ equals the interior of $\Omega_1 \ast \cdots \ast \Omega_m$.
Since $\Gamma$ acts cocompactly on $\Omega$, it acts on its projection $\Omega_j$ in a sweeping manner. 
	Suppose that 
	$\Omega_j^o\cap K_j^o \ne \emp$ or $\Omega^j \cap \mathcal{A}(K_j^o) \ne \emp$.  
	Theorem \ref{prelim-thm-Vey} shows that 
	$K_j^o = \Omega_j^o$ or $\mathcal{A}(K_j^o) = \Omega_j^o$
	since $\Gamma$ acts on the convex domain $\Omega_j^o$ as a sweeping action. 
	Suppose that $\Omega_j^o \cap K_j^o$ or $\Omega_j^o \cap \mathcal{A}(K_j^o)$
	are empty for all $j$. 
	Then 
	$\clo(\Omega_j)\cap K_j \ne \emp$ or 
	$\clo(\Omega_j)\cap \mathcal{A}(K_j)\ne \emp$ 
	by Proposition \ref{prelim-prop-sweep}. 
	Since such intersection has a unique minimal subspace containing it, 
	this contradicts the irreducibility of $\Gamma'$-action on $P_j$.

	
	
	Hence, it follows that 
	$K':= (K_1^{\prime}\ast \cdots \ast K_k^{\prime})^o$ 
	is a subset of $\Omega$ for $K_j' = K_j$ or $K'_j = \mathcal{A}(K_j)$. 
	Again $K'/\Gamma' \ra \Omega/\Gamma'$ is a homotopy equivalence, 
	and hence $K' = \Omega$.   Hence, the first item is proved. 
	
	The action of projective automorphisms which restricts to $\Idd$ or $\mathcal{A}$ on
	each $P_j$ 	gives us the final part. 
	(See Theorem 4.1 of \cite{CLM18} also.) 
\hfill	\SnT {\parfillskip0pt\par}
\end{proof}

We have the following useful result. 
\begin{corollary} \label{prelim-cor-uniqueness}
Let $\{h_i\}$ be a sequence of faithful
discrete representations from $\Gamma$ 
to $\SL_\pm(n+1, \bR)$ {\rm (}resp. $\PGL(n+1, \bR)${\rm )} such that 
$\orb_i:= \Omega_i/h_i(\Gamma)$ is a closed real projective orbifold
for a properly convex domain $\Omega_i$ in $\SI^n$
for each $i$. 
Suppose that $\{h_i\} \ra h_\infty$ algebraically, 
$h_\infty$ is faithful with a discrete image, 
and $\{\clo(\Omega_i)\}$ geometrically converges to 
a properly convex domain $\clo(\Omega_{\infty})$ 
with a nonempty interior $\Omega_{\infty}$. 
Then $h_\infty(\Gamma)$ acts on
the interior $\Omega_{\infty}$
such that the following statements hold\/{\rm :}
\begin{itemize} 
\item  $\Omega_{\infty}/h_\infty(\Gamma)$ is a closed real projective
orbifold. 
\item $\Omega_{\infty}/h_\infty(\Gamma)$ is diffeomorphic to 
$\orb_i$ for sufficiently large $i$. 
\item If $U$ is a properly convex domain where $h_\infty(\Gamma)$ acts 
such that $U/h_\infty(\Gamma)$ is an orbifold, then 
$U = \Omega_{\infty}$ or ${J}(\Omega_{\infty})$
where $J$ is a projective automorphism commuting with a finite-index subgroup of 
$h_\infty(\Gamma)$. In particular, if $\Gamma$ is non-virtually-factorizable, 
then $J = \Idd$ or $\mathcal{A}$. 
\end{itemize}
\end{corollary} 
\begin{proof} 
By Proposition \ref{prelim-prop-prdisc}, the quotient 
$\Omega_\infty/h_\infty(\Gamma)$ is an orbifold. 
For the second item, 
see the proof of Theorem 4.1 of \cite{CLM18}. 
The third item follows from Lemma \ref{prelim-lem-domainIn}.  
\hfill \SnT {\parfillskip0pt\par}
\end{proof}

\subsection{Technical propositions.} \label{prelim-sub-technical} 


\begin{proposition} \label{prelim-prop-joinred} 
	If a group $G$ of projective automorphisms 
	acts on a strict join $A= A_1 $ with $A_1 \ast A_2$ for two compact convex sets $A_1$ and $A_2$ in $\SI^n$ {\rm (}resp. in $\RP^n${\rm )} with $\dim A_1 + \dim A_2 = n-1$, 
then $G$ is virtually reducible. 
\end{proposition} 
\begin{proof} 
	We prove for $\SI^n$.
	Let $x_1, \dots, x_{n+1}$ denote the homogeneous coordinates. 
	There is at least one pair of strict join sets $A_{1}, A_{2}$. 
	We choose a maximal collection of compact convex sets $A'_1, \dots, A'_m$ such that 
	$A$ is a strict join $A'_1 \ast \cdots \ast A'_m$. 
	Here, we have $A'_i \subset S_i$ for a subspace $S_i$ corresponding to 
	a subspace $V_i \subset \bR^{n+1}$ that forms an independent set of subspaces. 
	
	We claim that $g \in G$ permutes the collection $\{A'_1, \dots, A'_m\}$:
	Suppose not. We give coordinates such that 
	for each $i$, there exists some index set $I_i$ 
	such that elements of $A'_i$ satisfy $x_j = 0$ 
	for $j \in I_i$ 
	and elements of $A$ satisfy $x_i \geq 0$.
	Then we form a new collection of nonempty sets 
	\[J':= \{ A'_i \cap g(A'_j)| 0 \leq i, j \leq n, g \in G\}\]
	with more elements. 
	Since \[A = g(A) = g(A'_1) \ast \cdots \ast g(A'_l),\] 
it follows 
	that each $A'_i$ is a strict join of nonempty sets in \[J'_i := \{ A'_i \cap g(A'_j)| 0 \leq j \leq l, g \in G\}\]
	using coordinates.
	$A$ is a strict join of the collection of the sets in $J'$, a contraction to the maximal property. 
	
	Hence, by taking a finite-index subgroup $G'$ of $G$ acting trivially on the collection, 
	$G'$ is reducible. \hfill \SnT {\parfillskip0pt\par}
\end{proof}

\begin{proposition} \label{prelim-prop-decjoin} 
	Suppose that a set $G$ of projective automorphisms  in $\SI^n$ {\rm (}resp. in $\RP^n${\rm )}
	acts on subspaces $S_1, \dots, S_{l_0}$ and a properly convex domain $\Omega \subset \SI^n$ 
	{\rm (}resp.  $\subset \RP^n${\rm )}, corresponding 
	to independent 
	subspaces $V_1, \dots, V_{l_0}$ such that $V_i \cap V_j =\{0\}$ for $i \ne j$ 
	and $V_1 \oplus \cdots \oplus V_{l_0} = \bR^{n+1}$. Let $\Omega_i : = \clo(\Omega) \cap S_i$ for each $i$, $i=1, \dots, l_0$. 
	Let $\lambda_i(g)$ denote the largest norm 
	of the eigenvalues of $g$ restricted to $V_i$. 
	We assume that  
	\begin{itemize}
		\item for each $S_i$, $G_i := \{g|S_i| g \in G\}$ forms a bounded set of automorphisms, and 
		\item for each $S_i$, there exists a sequence $\{g_{i, j} \in G\}$ such that
		\[\left\{\frac{\lambda_i(g_{i, j})}{\lambda_k(g_{i, j})}\right\} \ra \infty
		\hbox{ for each } k, k\ne i \hbox{ as } j \ra \infty.\]  
	\end{itemize}
	Then $\clo(\Omega) = \Omega_1 * \cdots * \Omega_{l_0}$ for $\Omega_j \ne \emp, j=1, \dots, l_0$. 
\end{proposition} 
\begin{proof}
	First, we have $\Omega_i \subset \clo(\Omega)$ by definition. 
	Each element of a strict join has a vector that is a linear combination of 
	vectors in the directions of $\Omega_1, \dots, \Omega_{l_0}$, 
	Hence, we obtain \[\Omega_1 \ast \cdots \ast \Omega_{l_0} \subset \clo(\Omega)\] since
	$\clo(\Omega)$ is convex. 
	
	Let $z = [ \vec{v}_z]$ for a vector $\vec{v}_z$ in $\bR^{n+1}$.
	We express as a unique sum
\[\vec{v}_z= \vec{v}_1 + \cdots + \vec{v}_{l_0}, \vec{v}_j \in V_j \hbox{ for each } j, j=1, \dots, l_0.\] 
	Then $z$ uniquely determines $z_i = [\bv_i]$. 
	
	Let $z $ be any point.  
	We choose a subsequence of $\{g_{i, j}\}$ such that $\{g_{i, j}|S_i\}$ converges to a projective automorphism 
	$g_{i, \infty}: S_i \ra S_i$ and $\lambda_{i, j} \ra \infty$ as $j \ra \infty$. 
	Then $g_{i, \infty}$ also acts on $\Omega_i$. 
By Proposition \ref{prelim-prop-attract},  
	$\{g_{i, j}(z_i)\} \ra g_{i, \infty}(z_i) = z_{i, \infty}$  
	for some point $z_{i, \infty} \in S_i$. 
	We also have
	\begin{equation} \label{pr-eqn-lim}
	z_i = g_{i, \infty}^{-1}(g_{i, \infty}(z_i)) = g_{i, \infty}^{-1}(\lim_j g_{i, j}(z_i)) = g_{i, \infty}^{-1}(z_{i, \infty}).
	\end{equation}
	
	Now suppose that $z \in \clo(\Omega)$. 
	We have   $\{g_{i, j}(z)\} \ra z_{i, \infty}$ by the eigenvalue condition. Thus, 
	we obtain $z_{i, \infty} \in \Omega_i$ as $z_{i, \infty}$ is the limit of a sequence of orbit points of $z$. 
	Hence we also obtain $z_i \in \Omega_i$ by \eqref{pr-eqn-lim}. 
	We obtain $\Omega_i \ne \emp$.
	Therefore, we conclude that $\clo(\Omega) = \Omega_1 * \cdots * \Omega_{l_0}$
	since $z \in \{z_1\} \ast \cdots \ast \{z_{l_0}\}$. 
	
	

\hfill	\SnT {\parfillskip0pt\par}
\end{proof} 


\section{The Vinberg duality of real projective orbifolds}\label{prelim-sec-duality} 

The duality is a natural concept in real projective geometry, and it continues to play an essential role 
in this theory as well. \index{duality!Vinberg}

\subsection{The duality} \label{prelim-sub-duality} 
We begin with linear duality. Let $\Gamma$ be a group of linear transformations $\GL(n+1, \bR)$. 
Let $\Gamma^\ast$ be the {\em dual group} defined by $\{g^{\ast -1}| g \in \Gamma \}$. 

Suppose that $\Gamma$ acts on a properly convex cone $C$ in $\bR^{n+1}$ with the vertex $O$.
\begin{itemize} 
\item An open convex cone  $C^{\ast}$ in $\bR^{n+1 \ast}$  is {\em dual} to an open convex cone $C $ in $\bR^{n+1}$  if 
$C^{\ast} \subset \bR^{n+1 \ast}$ is the subset of all linear functionals that takes positive values on $\clo(C)-\{O\}$.
The cone $C^{\ast}$ has the origin as the vertex again. 
Note $(C^{\ast})^{\ast} = C$, and $C$ must be properly convex
since otherwise $C^o$ cannot be open. 
We generalize the notion in Section \ref{prelim-sub-Eduality}. 
\index{dual cones} 

\item Now,  $\Gamma^{\ast}$ acts on $C^{\ast}$.
A {\em central dilatational extension} $\Gamma'$ of $\Gamma$ by $\bZ$ is given by adding a scalar dilatation by a factor $s > 1$ for the set $\bR_+$ of positive real numbers. 
\item The dual $\Gamma^{\prime \ast}$ of $\Gamma'$ is a central dilatational extension of $\Gamma^{\ast}$. 
 Also, $\Gamma'$ acts cocompactly on $C$ if and only if $\Gamma^{\prime \ast}$ acts so 
on $C^{\ast}$. 
 (See Chapter 4 of  \cite{goldmanbook} for details.)
 \index{central dilatational extension} 

\item  Given a subgroup $\Gamma$ in $\PGLnp$, the dual group $\Gamma^{\ast}$ is the image in $\PGLnp$ of the dual of 
the inverse image of $\Gamma$ in $\SLpm$. 
\index{group!dual|textbf}
\index{dual group|textbf}  
\item We define $\RP^{n\ast}$ as $\bP(\bR^{n+1 \ast})$. \label{dualrpn}
\item A properly convex open domain $\Omega$ in $\bP(\bR^{n+1})$ is {\em dual} to a properly convex open domain
$\Omega^{\ast}$ in $P(\bR^{n+1 \ast})$ if $\Omega$ corresponds to an open convex cone $C$ 
and $\Omega^{\ast}$ corresponds to its dual $C^{\ast}$. 
We say that $\Omega^{\ast}$ is dual to $\Omega$. 
We also have $(\Omega^{\ast})^{\ast} = \Omega$ and $\Omega$ is properly convex if and only if so is $\Omega^{\ast}$. 
\index{dual domains} 
\label{duality}

\item We call $\Gamma$ a {\em dividing group} of $\Omega$ if a central dilatational extension of $\Gamma$ acts cocompactly on $C$ with a Hausdorff quotient. 
$\Gamma$ is dividing if and only if so is $\Gamma^{\ast}$. 
Also, we say that $\Omega$ is divisible by $\Gamma$.
\index{dividing group|textbf} 
\index{divisible|textbf} 

\item 
Define $\SI^{n\ast} := \SI(\bR^{n+1\ast})$. 
For an open properly convex subset $\Omega$ in $\SI^{n}$, the dual domain is defined as the quotient  in $\SI^{n\ast}$
of the dual cone of the cone $C_\Omega$ corresponding to $\Omega$. 
The dual set $\Omega^{\ast}$ is also open and properly convex
but the dimension may not change.  
Again, we have $(\Omega^{\ast})^{\ast} =\Omega$. 

\item If $\Omega$ is a properly convex compact domain, not necessarily open, 
then we define $\Omega^\ast$ as the closure of the dual domain of $\Omega^o$.  
This definition is consistent with the definition given in Section \ref{prelim-sub-Eduality}
for any compact convex domains since a sharply supporting hyperspace 
can be perturbed to a supporting hyperspace that is not sharply supporting. 
(See Berger \cite{Berger}.)

\item Given a properly convex domain $\Omega$ in $\SI^n$ (resp. $\RP^n$), 
we define the {\em augmented boundary} of $\Omega$ as follows:
\begin{multline} 
\Bd^{\Ag} \Omega  := \{ (x, H)| x \in \Bd \Omega, x \in H,  \\ 
H \hbox{ is an oriented sharply supporting hyperspace of } \Omega \}  
\subset \SI^{n}\times \SI^{n\ast}.
\end{multline}  \index{bd@$\Bd^{\Ag}\Omega$|textbf}
Obviously this forms a compact space. 
Define the projection \[\Pi^{\Ag}_\Omega: \Bd^{\Ag} \Omega  \ra \Bd \Omega\] 
by $(x,H) \mapsto x$.  \index{pi@$\Pi_{\Ag}$|textbf}
Each $x \in \Bd \Omega$ has at least one oriented sharply supporting hyperspace.
An oriented hyperspace is an element of $\SI^{n \ast}$ since it is represented by a linear functional. 
Conversely, an element of $\SI^n$ represents an oriented hyperspace in $\SI^{n \ast}$. 
(Clearly, we can do this for $\RP^{n}$ and the dual space $\RP^{n \ast}$
but we consider only nonoriented supporting hyperspaces here.) 
\index{boundary!augmented|textbf} 
\end{itemize} 
\label{PiAg} 

\begin{theorem} \label{prelim-thm-Serre} 
Let $A$ be a subset of $\Bd \Omega$. Let $A':= \Pi_\Omega^{\Ag,-1}(A)$ be a subset of $\Bd^{\Ag}(A)$. 
Then $\Pi^\Ag_\Omega| A': A' \ra A$ is a quasi-fibration. 
\end{theorem} 
\begin{proof}
We take a Euclidean metric on an affine subspace containing $\clo(\Omega)$. 
Each sharply supporting hyperspaces at $x$ correspond to a unit normal vector at $x$. 
Each fiber $\Pi_\Omega^{\Ag, -1}(x)$ is 
a properly convex compact domain in a sphere of unit vectors through $x$. 
We find a continuous section defined on $\Bd \Omega$ by taking the center of mass of each fiber with respect 
to the Euclidean metric. This gives a local coordinate system on each fiber by assigning it 
the origin, and each fiber is a compact convex domain 
containing the origin. Then the quasi-fibration property follows. 
\hfill \SnT  {\parfillskip0pt\par}
\end{proof}

\begin{remark} 
We notice that for 
properly convex open or compact domains $\Omega_1$ and $\Omega_2$ in $\SI^n$ 
(resp. in $\RP^n$) we have 
\begin{equation}\label{prelim-eqn-dualinc}
\Omega_1 \subset \Omega_2 \hbox{ if and only if } \Omega_2^{\ast} \subset \Omega_1^{\ast}  
\end{equation}

\end{remark}

Consider decomposition of $\bR^{n+1}$ as a direct sum of two 
subspaces $\bR^{k+1}$ of dimension $k+1$ and  $\bR^{n-k}$ of dimension $n-k$. 
They correspond to $\RP^k$ and $\RP^{n-k-1}$ in $\RP^n$. 
$\RP^{k\ast}$ embeds into $\RP^{n\ast}$ as $\bP(V_1)$ for the subspace 
$V_1$ of linear functionals nullifying 
vectors in directions of $\RP^{n-k-1}$, and $\RP^{n-k-1\ast}$ embeds into  
$\RP^{n\ast}$ as $\bP(V_2)$ for the subspace $V_2$ of linear functionals 
nullifying the vectors in directions of $\RP^k$. 
We denote these by $\RP^{k\dagger}$ and $\RP^{n-k-1\dagger}$
respectively and refer them as {\em intrinsic duals}. 
Similarly, we can define the intrinsic duals $\SI^{k\dagger}$ similarly. 

\begin{remark} \label{prelim-rem-intrinsic_dual}
A {\em intrinsic-subspace dual $K^{\dagger}_X$ with respect to $X=\RP^k$ or $\SI^{k}$} of a convex domain $K$ in $X= \RP^k$ or $\SI^{k}$ is the dual domain as obtained from considering $X$ and  its corresponding vector subspace only. 
\label{dagger} 
\index{dagger@$(\cdot)^\dagger$|textbf}
\index{dual!instrinsic} 
\index{intrinsic dual}

Consider a strict join $A\ast B$ for a properly convex compact $k$-dimensional domain $A$ in $\RP^k\subset \RP^n$ and a properly convex compact $(n-k-1)$-dimensional domain $B$ in the complementary $\RP^{n-k-1} \subset \RP^n$.
Let $A^\dagger_{\RP^k}$ denote the
subspace dual in $\RP^{k\ast}$ of $A$ in $\RP^k$ defined intrinsically 
and $B^\dagger_{\RP^{n-k-1}}$ the intrinsic subspace dual domain in $\RP^{n-k-1\ast}$ 
of $B$ in $\RP^{n-k-1}$ defined intrinsically. 

 Then we have 
\begin{equation}\label{prelim-eqn-dualjoin} 
(A \ast B)^{\ast} = A^{\dagger}_{\RP^k} \ast B^{\dagger}_{\RP^{n-k-1}}.
\end{equation} 
This follows from the definition as 
every linear functional as 
a sum of linear functionals in the direct-sum subspaces.

 Suppose that $A \subset \SI^k$ and $B \subset \SI^{n-k-1}$ are 
 $k$-dimensional and $(n-k-1)$-dimensional domains respectively,
 where $\SI^k$ and $\SI^{n-k-1}$ are complementary subspaces in $\SI^n$. 
 We realize $\SI^{k\ast}$ and $\SI^{n-k-1\ast}$ as subspaces of $\SI^{n\ast}$ 
 as above using the linear functionals. 
We denote these subspaces by $\SI^{k\dagger}$ and $\SI^{n-k-1\dagger}$ 
 respectively. 
 Suppose that $A$ and $B$ have respective dual sets 
 $A^\dagger_{\SI^k} \subset \SI^{k\ast}, B^\dagger_{\SI^{n-k-1}} 
 \subset \SI^{n-k-1\ast}$. 
 Then the above equation also holds with the subscripts exchanged appropriately.

\end{remark}

An element $(x, H)$ is in $\Bd^{\Ag} \Omega$ if and only if $x \in \Bd \Omega$ and $H$ is represented 
by a linear functional $\alpha_H$ such that $\alpha_H(\vec{y}) > 0$ for all $\vec{y}$ in the open cone $C(\Omega)$ corresponding to $\Omega$ and 
$\alpha_H(\bv_x) =0$ for a vector $\bv_x$ representing $x$. 

Let $(x, H) \in \Bd^{\Ag} \Omega$. 
The dual cone $C(\Omega)^{\ast}$ consists of every nonzero $1$-form $\alpha$ such that $\alpha(\vec{y}) > 0$ for all $\vec{y} \in \clo(C(\Omega)) - \{O\}$. 
Thus, $\alpha(\bv_x) > 0$ for all $\alpha \in C^{\ast}$ and $\alpha_H(\bv_x) = 0$, 
and  
$\alpha_H \not\in C(\Omega)^{\ast}$ since $\bv_x \in \clo(C(\Omega))-\{O\}$. 
However, $H \in \Bd \Omega^{\ast}$ as we can perturb $\alpha_H$ such that it is in $C^{\ast}$. 
Thus, $x$ is a sharply supporting hyperspace at $H \in \Bd \Omega^{\ast}$.
We define a {\em duality map}
\[ {\mathcal{D}}^\Ag_{\Omega}: \Bd^{\Ag} \Omega \ra \Bd^{\Ag} \Omega^{\ast} \] 
given by sending $(x, H)$ to $(H, x)$ for each $(x, H) \in \Bd^{\Ag} \Omega$. 
\index{duality map|textbf}
\index{do@$\mathcal{D}^\Ag_{\Omega}(\cdot)$|textbf} 
\index{D@$\mathcal{D}(\cdot)$|textbf} 



\begin{proposition} \label{prelim-prop-duality}
Let $\Omega$ and $\Omega^{\ast}$ be dual open domains in $\SI^{n}$ and $\SI^{n \ast}$ {\rm (}resp. $\RP^{n}$ 
and $\RP^{n \ast}${\rm ).} 
\begin{enumerate}  
\item[{\rm (i)}]\, There is a proper map $\Pi^{\Ag}: \Bd^{\Ag} \Omega \ra \Bd \Omega$
given by sending $(x, H)$ to $x$. 
\item[{\em (ii)}]\,\, A projective automorphism group 
$\Gamma$ acts properly on a properly convex open domain $\Omega$ if and only if 
so $\Gamma^{\ast}$ acts on $\Omega^{\ast}$\/  {\rm (}see Vinberg's Theorem \ref{prelim-thm-dualdiff}\/ {\rm ).}
\item[{\rm (iii)}]\,\, There exists a duality map 
\[ {\mathcal{D}}^\Ag_{\Omega}: \Bd^{\Ag} \Omega \ra \Bd^{\Ag} \Omega^{\ast} \] 
which is a homeomorphism. 
\item[{\rm (iv)}]\,\, Let $A \subset \Bd^{\Ag} \Omega$ be a subspace, and let $A^{\ast}\subset \Bd^{\Ag} \Omega^{\ast}$
be the corresponding dual subspace $\mathcal{D}^\Ag_{\Omega}(A)$. A group $\Gamma$ acts on $A$ such that $A/\Gamma$ is compact 
if and only if $\Gamma^{\ast}$ acts on $A^{\ast}$ and $A^{\ast}/\Gamma^{\ast}$ is compact. 
\end{enumerate} 
\end{proposition} 
\begin{proof} 
We will prove for $\SI^n$ first. 
(i) Each fiber is a closed set of hyperspaces at a point and forms a compact set.
The set of sharply supporting hyperspaces at a compact subset of $\Bd \Omega$ is closed. 
The closed set of hyperspaces having a point in a compact subset of $\SI^{n+1}$ is compact. 
Thus, $\Pi^{\Ag}$ is proper.  Clearly, $\Pi^{\Ag}$ is continuous and is an open map
since $\Bd^{\Ag} \Omega$ is given the subspace topology from 
$\SI^n \times \SI^{n \ast}$
with a product topology where $\Pi^{\Ag}$ extends to a projection.

(ii) See Theorem 4.4.10  \cite{goldmanbook} or Vinberg \cite{Vinberg63}.

(iii) ${\mathcal{D}}^\Ag_{\Omega}$ has the inverse map ${\mathcal{D}}^\Ag_{\Omega^{\ast}}$. 

(iv) This follows directly from (iii). 
\hfill \SnT {\parfillskip0pt\par}
\end{proof} 

\index{PiAg@$\Pi^{\Ag}$|textbf} 

\begin{definition}
The two subgroups $G_1$ of $\Gamma$ and $G_2$ of $\Gamma^{\ast}$ are {\em dual} if 
sending $g \mapsto g^{-1, T}$ gives us an isomorphism $G_1 \ra G_2$. 
A set in $A \subset \Bd \Omega$ is {\em dual} to a subset $B \subset \Bd \Omega^{\ast}$ if 
 $\mathcal{D}^\Ag_{\Omega}: \Pi_\Ag^{-1}(A) \ra \Pi_{\Ag}^{-1}(B)$ is a one-to-one and onto map. 
\end{definition} 
\index{dual groups}
\index{dual sets} 

\begin{remark}\label{prelim-rem-duallens}
	For an open subspace $A \subset \Bd \Omega$ that is $C^1$ and strictly convex, 
	${\mathcal{D}}^\Ag_{\Omega}$ induces a well-defined map 
	\[A \subset \Bd \Omega \ra A' \subset \Bd \Omega^{\ast}\]
	since each point has a unique sharply supporting hyperspace 
	for an open subspace $A'$.
	The image $A'$ of the map is also smooth and strictly convex by Lemma \ref{pr-lem-predual}. 
	We will simply say that $A'$ is the {\em image}  of $\mathcal D$. 
\end{remark}




Let $\RP^{n-1}_{x}$ denote the space of concurrent lines to a point $x$
where two lines are equivalent if they coincide in a neighborhood of $x$. 
Here, $\RP^{n-1}_{x}$ is projectively diffeomorphic to $\RP^{n-1}$. 
The real projective transformations fixing $x$ induce real projective transformations of $\RP^{n-1}_{x}$. 
\index{RPnminusonex@$\RP^{n-1}_{x}$|textbf} 
Let $x\in \SI^n$. 
The space $\SI^{n-1}_x$ denotes the space of equivalence classes of concurent lines ending 
at $x$ with orientation away from $x$ where two such lines are considered 
equivalent if they agree on an open subset with a common boundary point $x$. 
Each equivalence class is called a {\em direction from} $x$. 
Note that $\SI^{n-1}_x$ is well-defined on $\RP^n$ as well for $x \in \RP^n$. 
\index{direction}
\index{Snminusonex@$\SI^{n-1}_x$|textbf}


\label{page:Rx}

For a subset $K$ in a convex domain $\Omega$ in $\RP^n$ or 
$\SI^n$, let $x$ be a boundary point. 
For a subset $K$ of $\Omega$,
we define $R_x(K)$ as the space of 
directions of open rays from $x$ in $\Omega$ ending at $K$. 
We defined $R_x(K) \subset \SI^{n-1}_x$. 
\index{rx@$R_x(\cdot)$|textbf}
Any projective group fixing $x$ induces an action on 
$\SI^{n-1}_x$.

We say that a two-sided open hypersurface is {\em convex polyhedral} if it is a union of a locally 
finite collection of 
compact polytopes in hyperspaces meeting one another in strictly convex 
angles where the convexity is towards a fixed side. 
\index{hypersurface!convex polyhedral}

\begin{lemma}\label{pr-lem-predual}
	Let $\Omega^*$ be the dual of a properly convex open domain $\Omega$ 
	in $\RP^n$ {\em (}resp. in $\SI^n${\em ).} 
	Then 
	\begin{enumerate}
		\item[{\rm (i)}] $\Bd \Omega$ is $C^1$ and strictly convex at a point $p \in \Bd \Omega$ if and only if 
		$\Bd \Omega^*$  is $C^1$ and strictly convex at the unique corresponding point $p^{*}$. 
		\item[{\rm (ii)}] $\Omega$ is an ellipsoid if and only if so is $\Omega^*$. 
		\item[{\rm (iii)}] $\Bd \Omega^*$ contains a properly convex domain $D = P \cap \Bd \Omega^*$ open in a totally geodesic hyperspace $P$ 
		if and only if $\Bd \Omega$ contains 
		a vertex $p$ with $R_p(\Omega)$ a properly convex domain. 
		In this case, ${\mathcal{D}}^\Ag_{\Omega}$ sends the pair of $p$ and the associated sharply supporting hyperspace of $\Omega$
		to the pairs of the totally geodesic hyperspace containing $D$ and points of $D$.  
		Moreover, $D$ and $R_p(\Omega)$ are properly
		convex, and the projective dual of $D$ is projectively diffeomorphic to 
		$R_p(\Omega)$. 
		\item[{\rm (iv)}] Let $S$ be a convex polyhedral 
		open subspace of $\Bd \Omega$. Then $S$  
		is dual to a convex polyhedral open subspace of $\Bd \Omega^\ast$. 
	\end{enumerate}
\end{lemma}
\begin{proof} 
	We suppose that $\Omega \subset \SI^n$. 
	(i) 
	The one-to-one map ${\mathcal{D}}^\Ag_{\Omega}$ sends each pair $(x, H)$ of a 
	point of $\Bd \Omega$ and the sharply supporting hyperspace to 
	the pair $(H, x)$ where $H$ is a point of $\Bd \Omega^\ast$ and 
	$x$ is a sharply supporting hyperspace at $H$ of $\Omega^\ast$. 
	
\begin{itemize} 
\item The fact that $\Bd \Omega$ is $C^1$ imlpies that
for $x \in \Bd \Omega$, $H$ is unique, and 
\item the fact that $\Bd \Omega$ strictly convex implies that  for $H$, 
	there is only one point of $\Bd \Omega$ where $H$ meets $\Bd \Omega$. 
\end{itemize} 
	Also, this is equivalent to the fact that 
	for each $H \in \Bd \Omega^\ast$, the supporting hyperspace $x$ is unique 
	and for each $x$, there is one point of $\Bd \Omega^\ast$ where $x$ 
	meets $\Bd \Omega^\ast$. This shows that $\Bd \Omega^\ast$ is 
	strictly convex and $C^1$. 
	
	
	(ii) Let $\bR^{n+1}$ have the standard Lorentz inner product $B$. 
	Let $C$ be the open positive cone. Then the space of linear functionals
	positive on $C$ is in one-to-one correspondence with 
	vectors in $C$ using the isomorphism $C^\ast \ra C$ given
	by $\phi \mapsto \bv_\phi$ such that $\phi = B(\bv_\phi, \cdot)$.
	(See \cite{goldmanbook}.)
	
	(iii) Suppose that $R_{p}(\Omega)$ is properly convex. 
	We consider the set of hyperspaces sharply supporting $\Omega$ at $p$. 
	This forms a properly convex domain: 
	Let $\vec{v}$ be the vector in $\bR^{n+1}$ in a direction of $p$. 
	Then we find the set of  linear functionals positive on $C(\Omega)$ but vanishing on $\vec{v}$. 
	Let $\bV$ be a complementary space of $\vec{v}$ in $\bR^{n+1}$. 
	Let $\mathds{A}$ be given as the affine subspace 
	$\bV + \{\vec{v}\}$ of $\bR^{n+1}$. 
	Choose $\bV$ such that $C_{\vec{v}}:= C(\clo(\Omega)) \cap \mathds{A}$ is a bounded convex domain in $\mathds{A}$. 
	We give $\mathds{A}$ a linear structure such that $\vec{v}$ corresponds to the origin.
	We identify this space with $\bV$. 
	The set of linear functionals positive on $C(\Omega)$ and $0$ at $\vec{v}$ is identical with 
	the set of linear functionals on $\bV$ positive on $C_{\vec{v}}$:
	we define 
	\begin{multline} 
	C(D) := \big\{ f\in \bR^{n+1\ast}\big| \, f| C(\clo(\Omega)) -\{t\vec{v}| t\geq 0\} > 0, f(\vec{v}) = 0  \big\} \nonumber \\ 
	\cong \widehat C_{\vec{v}}^{\ast} :=\big\{ g \in \bV^{\ast} \big| \, g| C_{\vec{v}}-\{O\} > 0 \big\} \subset \bR^{n+1\ast}. 
	\end{multline}
	Here $\cong$ indicates a linear isomorphism, 
	which follows from the decomposition $\bR^{n+1} = \{t\vec{v}| t \in \bR\} \oplus \bV$.
	Define $R'_{\vec{v}}(C_{\vec{v}})$ as the equivalence classes of properly convex segments in $C_{\vec{v}}$ ending at $\vec{v}$ 
	where two segments are equivalent if they coincide in an open neighborhood of $\vec{v}$. 
	$R_{p}(\Omega)$ is identical to $R'_{\vec{v}}(C_{\vec{v}})$ by 
the sphericalization map $\SI: \bR^{n+1}-\{O\} \ra \SI^n$.  
	Hence, $R'_{\vec{v}}(C_{\vec{v}})$ is a properly convex open domain in $\SI(\bV)$. 
	Since $R'_{\vec{v}}(C_{\vec{v}})$ is properly convex, 
	the interior of its spherical projectivization 
	$\SI(\widehat C_{\vec{v}}^{\ast}) \subset \SI(\bV^{\ast})$ 
	is dual to the properly convex domain $R'_{\vec{v}}(C_{\vec{v}}) \subset \SI(\bV)$.  
	
	Again, we have a projection 
	$\SI: \bR^{n+1\ast}-\{O\} \ra \SI^{n\ast}$. 
	Define $D:= \SI(C(D)) \subset \SI^{n\ast}$. 
	Since $R'_{\vec{v}}(C_{\vec{v}})$ corresponds to $R_{p}(\Omega)$,  and 
	$\SI(\widehat C_{\vec{v}}^{\ast})$ corresponds to $D$, the duality follows. 
	Also, $D \subset \Bd \Omega^\ast$ since points of $D$ are oriented 
	sharply supporting hyperspaces to $\Omega$ by Proposition \ref{prelim-prop-duality} (iii).
	(Here, we can also use 
	Proposition \ref{pr-prop-dualtube}.) 

The converse follows from reversing the argument above. 
	
	(iv) From (iii) each vertex of a convex polyhedral subspace of $S$ corresponds to 
	a compact convex polytope in the dual subspace. Also, we can check that each 
	side of dimension $i$ corresponds to a side of dimension $n-i-1$.  
\hfill	\SnT {\parfillskip0pt\par}
\end{proof}

\subsection{The duality of convex real projective orbifolds with strictly convex boundary} 

Since $\orb = \Omega/\Gamma$ for an open properly convex domain $\Omega$
in $\RP^n$ (resp. in $\SI^n$), 
the dual orbifold $\orb^{\ast} = \Omega^{\ast}/\Gamma^{\ast}$ is a properly convex real projective orbifold. 
The dual orbifold is well-defined up to projective diffeomorphisms. 




\begin{theorem}[Vinberg, see Theorem 4.4.10 in Chapter 4 of \cite{goldmanbook}] \label{prelim-thm-dualdiff} 
	Let $\orb$ be a strongly tame properly convex real projective orbifold that is either open or closed. 
	The dual orbifold $\orb^{\ast}$ is diffeomorphic to $\orb$.
\end{theorem}

The map in \cite{Vinberg63} is {\em called}
the {\em Vinberg duality diffeomorphism}.
For an orbifold $\orb$ with boundary, 
the map is a diffeomorphism in the interiors $\orb^{o} \ra \orb^{\ast o}$.
Let $\torb$ denote the properly convex projective domain covering $\orb$ projectively.
Also,  ${\mathcal{D}}^\Ag_{\torb^o}$ gives us the diffeomorphism $\partial \orb \ra \partial \orb^{\ast}$. 
(We conjecture that they form a diffeomorphism $\orb \ra \orb^{\ast}$ up to isotopies. We also remark that 
${\mathcal{D}}_{\orb^\ast}\circ {\mathcal{D}}_{\orb}$ 
may not be identity as shown by Vinberg.) 
\index{Vinberg duality diffeomorphism|textbf} 

For each $p \in \Omega$, let $\vec{p}_{V, \Omega}$ denote the vector in 
$C(\Omega)$ with $f_V^{-1}(\vec{p}_{V, \Omega})= 1$ for the Koszul-Vinberg function 
$f_{V, \Omega}$ for $C(\Omega)$. 
(See \eqref{cl-eqn-kv}.)
Define $\vec{p}_{V, \Omega}^\ast$ as the $1$-form 
$Df_{V, \Omega}(\vec{p}_V)$, and  also define $p^\ast$ as
$\llrrparen{Df_{V, \Omega}(\vec{p}_{V, \Omega})}$.
We obtain a compactification of $\Omega$ by defining 
$\clo^\Ag(\Omega) := \Omega \cup \Bd^\Ag \Omega$.
For any sequence $p_i \in \Omega$, we form a pair 
$(\vec{p}_{i}, \vec{p}_{i, V}^\ast)$ where $\vec{p}_{i, V}^\ast$ is the $1$-form in $\bR^{n\ast}$ given by 
\[D f_V(\vec{p}_{i, V}).\] 
	Clearly, a limit point of $\{\vec{p}_{i, V, \Omega}^\ast\}$ is a supporting $1$-form of 
	$C(\Omega)$ since it supports a properly convex domain 
	$f_V^{-1}(1, \infty)\subset C(\Omega)$ and contains a point of $\Bd \Omega$. 
	We say that $p_i$ converges to an element of $\Bd^\Ag \Omega$ if 
	this augmented sequence converges to it. 
\index{sequence!augmented}
\index{closure!augmented}
\index{Vinberg duality diffeomorphism!augmented|textbf} 

\begin{theorem}\label{prelim-thm-AugVinberg} 
	Let $\Omega$ be a properly convex domain
	in $\RP^n$ {\em (}resp. in $\SI^n${\em ),}
	and let $\Omega^\ast$ be its dual 
	in $\RP^{n\ast}$ {\em (}resp. in $\SI^{n\ast}${\em )}. 
	Then the Vinberg duality diffeomorphism ${\mathcal{D}}_{\Omega}: \Omega \ra \Omega^\ast$ 
	extends to a homeomorphism $\bar{\mathcal{D}}_\Omega^\Ag: \clo^\Ag(\Omega) \ra \clo^\Ag(\Omega)$. 
	Moreover, for any projective group $\Gamma$ acting on it, 
	$\bar{\mathcal{D}}_\Omega^\Ag$ is equivariant with respect to
	the duality map $\Gamma \ra \Gamma^\ast$ given by $g \mapsto g^{\ast-1}$. 
\end{theorem} 
\begin{proof} 
	We assume $\Omega \subset \SI^n$. 
	The continuity follows from 
	the paragraph above the theorem 
	since $\mathcal{D}_{\torb}$ is induced 
	by $(p, p_{V, \Omega}^\ast) \ra (p_{V, \Omega}^\ast, p)$,
	and $\mathcal{D}^\Ag_{\torb}$ also switches the orders of the pairs.
	
	Proposition \ref{prelim-prop-duality} shows the injectivity of 
	$\bar{\mathcal{D}}_\Omega^\Ag$. 
	The map is surjective since $\mathcal{D}_\Omega$ and 
	$\mathcal{D}_\Omega^\Ag$ are surjective. 
	
	The equivariance follows since so are $\mathcal{D}_{\Omega}$ and 
	$\mathcal{D}_{\Omega}^\Ag$.  \hfill \SnT {\parfillskip0pt\par}
	\end{proof}

\subsection{Sweeping actions} \label{prelim-sub-sweeping}

The properly convex open set $D$ in $\RP^n$ 
(resp. $\SI^n$) has a Hilbert metric. Also the group $\Aut(K)$ of projective automorphisms of $K$ in $\SL_\pm(n+1, \bR)$ is a locally compact closed group. 

\begin{lemma}\label{np-lem-invmet} 
	Let $D$ be a properly convex open 
	domain in $\RP^n$ {\em (}resp. $\SI^n${\em )}
	with $\Aut(D)$ of smooth projective automorphisms of $D$.
	Let a group $G$ act isometrically and effectively on an open domain $D$ such that 
the  map $G \rightarrow \Aut(D)$ is an embedding. 
	Suppose that the action of $G$ on $D$ is cocompact. 
	Then the closure $\bar G$ of $G$ is a Lie subgroup 
	acting on $D$ properly, and there exists a smooth Riemannian metric on $D$ that is 
	$\bar G$-invariant. 
\end{lemma} 
\begin{proof} 
	Assume $D \subset \SI^n$. 
	Since $\bar G$ is contained in $\SL_\pm(n+1, \bR)$, the closure 
	$\bar G$ is a Lie subgroup acting on $D$ properly. 
	
	We can construct a Riemannian metric $\mu$ with bounded entries. 
	Let $\phi$ be a function supported on a compact set $F$ such that 
	$G(F) \supset D$ where $\phi|F > 0$.
	Given a bounded subset of $\bar G$, the elements are in a bounded subset of the projective
	automorphism group $\SL_{\pm}(n+1, \bR)$. A bounded subset of projective automorphisms 
	has uniformly bounded set of derivatives on $\SI^{n}$ up to the $m$-th order for any $m$. 
	We can assume  that the derivatives of the entries of $\phi \mu$ up to the $m$-th order are uniformly bounded above.
	Let $d\eta$ be the left-invariant measure on $\bar G$.
	
	Then $\{g^*\phi\mu| g \in \bar G\}$ is an equicontinuous family on any compact 
	subset of $D^{o}$ up to any order. 
	For any compact $J \subset D^o$, $\supp (g^\ast \phi \mu) \cap J \ne \emp$ for $g$ 
	in a compact set of $\bar G$. 
	Thus, the integral 
	\[ \int_{g \in \bar G} g^* \phi\mu d \eta \]
	of $g^* \phi\mu$ for $g \in \bar G$ is a $C^\infty$-Riemannian metric and that is positive definite. 
	This bestows us a $C^\infty$-Riemannian metric $\mu_D$ on $D$ invariant under $\bar G$-action.
	\hfill \SnT {\parfillskip0pt\par}
\end{proof}


By Lemma \ref{np-lem-invmet}, there exists a Riemannian metric 
on a properly convex domain $\Omega$ invariant under $\Aut(\Omega)$. 
Hence, we can define a frame bundle $\bF \Omega$
on which $\Aut(\Omega)$ acts freely.

\begin{proposition}[Lemma 1 of Vey \cite{Vey}] \label{prelim-prop-sweepduality}
	Suppose that a projective group $G$ acts on an $(n-1)$-dimensional  
	properly convex open domain $\Omega$ as 
	a sweeping action. Then the dual group $G^*$ acts on $\Omega^*$ as 
	a sweeping action also. 
\end{proposition} 
\begin{proof} 	
	The Vinberg duality map in Theorem \ref{prelim-thm-dualdiff} 
	is a diffeomorphism $\Omega \ra \Omega^\ast$. 
	This map is equivariant under the duality homomorphism 
	$g \mapsto g^{\ast -1}$ for each $g \in G$. 
	Here, $G$ does not need to be a dividing action. 
\hfill \SnT {\parfillskip0pt\par}
	\end{proof}

\begin{theorem}[Generalized Theorem \ref{prelim-thm-Kobayashi}]\label{prelim-thm-Kobayashi2}
	Suppose that a projective group $G$  acts sweepingly on
	a properly convex open domain $D$ in $\RP^n$ {\em (}resp. in $\SI^n${\em )}. 
Then for any properly convex open domain 
$\Omega$ where the $G$-action of $\Omega$ is cocompact and $\Omega \cap D \ne \emp$, 
we must have $\Omega = D$. 
	\end{theorem} 
\begin{proof} 
Suppose not. Then $G$ acts on $D\cap \Omega$ as a sweeping manner and $D\cap \Omega$ is a proper subset of $D$. Let $x \in D\cap \Omega$. 
By Theorem \ref{prelim-thm-Vey}, the convex hull of $Gx$ must equal both $D\cap \Omega$ and $D$. Hence, $D\subset \Omega$. The converse also holds by the same reason.  \hfill 
\SnT {\parfillskip0pt\par}
\end{proof} 

\subsection{Extended duality} \label{prelim-sub-Eduality}

We can generalize the duality for convex domains
as in the beginning of Section \ref{prelim-sec-duality};
however, we can generalize for $\SI^n$ but not for $\RP^n$. 
Given a closed convex cone $C_1$ in $\bR^{n+1}$, 
consider the set of linear functionals on $\bR^{n+1\ast}$ taking 
nonnegative values in $C_1$. 
This forms a  closed convex cone. We call this the {\em dual cone} of $C_1$
and denote it by $C_1^\ast$. 

A closed cone $C_2$ in $\bR^{n+1 \ast}$ is {\em dual} to a closed 
convex cone $C_1$ in $\bR^{n+1}$ if $C_2$ is the set of 
linear functionals taking nonnegative values on $C_1$. 

For a convex compact set $U$ in $\SI^n$, we form a corresponding convex cone 
$C(U)$. Then we form $C(U)^\ast$ and the image of its projection 
a convex compact set 
$U^\ast$ in $\SI^{n\ast}$. 
Clearly, $(U^\ast)^\ast = U$ for a compact convex set $U$ by 
definition since this operation just interchange the roles of vectors and dual vectors. 

 
 Also, the definition agrees with the previous definition
 for properly convex domains. This is straightforward: 
Suppose that $U$ is properly convex. 
Functions in $C(U)^\ast$ can be approximated arbitrarily by 
functions strictly positive on $C(U)$. 

 
Recall the classification of compact convex 
sets in Proposition \ref{prelim-prop-classconv}. 
%
%
%
%
Recall the intrinsic duality in Remark \ref{prelim-rem-intrinsic_dual}. 
 
 \begin{proposition} \label{prelim-prop-NPduality}$ $
\begin{itemize}  	
\item  	Let $\SI^{i_0}$ be a great sphere of dimension $i_0$. 
\item Let $\SI^{j_0}$ be one of dimension $j_0$
 	with $i_0 + j_0+1 \leq n$ independent of $\SI^{i_0}$. 
\item We also have 
 	the join $\SI^{i_0+ j_0+1}$ of $\SI^{i_0}$ and $\SI^{j_0}$ 
 	and its complementary subspace $\SI^{n-i_0-j_0 -2}$. 	
\item Let $\SI^{n-i_0-1}$ be one of dimension $n-i_0-1$
 complementary to $\SI^{i_0}$
 where $\SI^{n-i_0-j_0-2}\subset \SI^{n-i_0-1}$. 
\item  Let us identify the intrinsic dual $\SI^{i_0\dagger}$ of 
$\SI^{i_0}$ as $\SI^{n-i_0 -1 \ast}$ by taking restrictions of 
linear maps and the intrinsic dual
$\SI^{n- i_0 - j_0 -2\dagger}$ of $\SI^{n-i_0-j_0-2}$ 
as $\SI^{i_0+j_0+1\ast}$. 
\end{itemize} 
 	Let $U$ be a convex compact proper set in $\SI^n$. 
 	Then the following hold\/{\em :}
 	\begin{description} 
 		\item[(i)] $U$ is a great $i_0$-sphere if and only if $U^\ast$ is a great $n-i_0-1$-sphere. 
 		$U^\ast$ is not convex if and only if $i_0= n-1$. 
 		\item[(ii)] If $U$ is a strict join of a properly convex domain $K$ of dimension 
 		$i_0$ in the great sphere $\SI^{i_0}$ and its complementary great sphere 
 		$\SI^{j_0}$ for $i_0, j_0\geq  0$, then 
 		\begin{itemize} 
 		\item $U^\ast$ is a strict join of $\SI^{n- i_0 - j_0 -2\dagger}= \SI^{i_0+j_0+1\ast}$ and 
 		a properly convex domain $K^{\dagger}$ in $\SI^{i_0\dagger}= \SI^{n-i_0-1\ast}$ of dimension $i_0$ intrinsically dual to $K$ in $\SI^{i_0}$ if $i_0 + j_0 +1 < n$. 
 		\item $U^\ast$ is $K^{\dagger} \subset  \SI^{i_0\dagger}$ if $i_0+j_0+1 = n$. 
 		\end{itemize} 
 	\item[(iii)] If $U$ is a properly convex domain $K$ in the great sphere $\SI^{i_0}$ of 
 	dimension $i_0, i_0 < n$, then $U^\ast$ is a strict join of the intrinsic dual $K^{\dagger}$ of $K$ in $\SI^{i_0\dagger}= \SI^{n-i_0-1\ast}$ and a great sphere $\SI^{i_0\ast}$ of dimension $n-i_0-1$
 	for any choice of the complement $\SI^{n-i_0-1}$ of 
 	$\SI^{i_0}$.  
 	\item[(iv)] $U$ is a properly convex compact $n$-dimensional domain  	if and only if so is $U^\ast$.
\item[(v)]  If $U$ is not properly convex and has a 
 nonempty interior, then $U^\ast$ has an empty interior. 
\item[(vi)]  If $U$ has an empty interior, then $U^\ast$ is not properly convex 
 and has a nonempty interior provided $U$ is properly convex. 
\item[(vii)]  In particular, if $U$ is an $n$-hemisphere, then $U^\ast$ is a point 
 and vice versa. 
\item[(viii)] 	If $U^o\ne \emp$ and $U^{\ast o} \ne \emp$, then 
 both $U$ and $U^\ast$ are properly convex domains in $\SI^n$. 
 \item[(xi)] $U$ is contained in a hemisphere if and only if $U^\ast$ is contained in a hemisphere. 
 \end{description} 
 \end{proposition} 
\begin{proof} 
(i) Suppose that $U$ is a great $i_0$-sphere. 
Then $C(U)$ is a subspace of dimension $i_0+1$. 
The set of linear functionals that takes value $0$ on $C(U)$ form
a subspace of dimension $n-i_0$. Hence, $U^\ast= \SI(C(U)^\ast)$ is 
a great sphere of dimension $n-i_0-1$. 
The converse also holds. 

(ii) Suppose that $U$ is not a great sphere. 
Proposition \ref{prelim-prop-classconv} shows us that 
$U$ is contained in an $n$-hemisphere. 

Let $\SI^{m_0}$ be the span of $U$. Here, $m_0= i_0 + j_0+1$. 
Then $U = \SI^{j_0} \ast K^{i_0}$ for a great sphere 
$\SI^{j_0}$ and a properly convex domain $K^{i_0}$ 
in a great sphere of dimension $i_0$ in $\SI^{m_0}$, independent of 
the first one by Proposition \ref{prelim-prop-classconv}. 

$C(U)$ is a closed cone in the vector subspace $\bR^{m_0+1}$. 
Then $C(U) = \bR^{j_0+1} + C(K^{i_0})$
where $C(K^{i_0}) \subset \bR^{i_0+1}$
for independent subspaces $\bR^{j_0+1}$ and $\bR^{i_0+1}$. 
Let $C(U)'$ denote the dual of $C(U)$ in $ \bR^{m_0+1\ast}$. 
For $f\in C(U)'$, $f = 0$ on $\bR^{j_0+1}$, 
and $f| \bR^{i_0+1}$ takes a nonnegative value on $C(K^{i_0})$. 
Hence, 
\[ f: \bR^{m_0+1} = \bR^{j_0+1} \oplus \bR^{i_0+1} \ra \bR
\hbox{ is in } \{0\}\oplus C(K^{i_0})^\ast.\] 
Denote the projection of $C(U)'$ in $\SI^{m_0}$ by $U'$. 

Suppose $m_0 = n$. Then we have shown the second case of (ii). 

Suppose $m_0 < n$. 
Then decompose $\bR^{n+1} $ as $\bR^{n-m_0} \oplus \bR^{m_0+1}$. 
We obtain that 
$f \in C(U)^\ast$ is a sum $f_1 + f_2$ where 
$f_1$ is an element of  $C(U)'$ extended by setting 
$f_1| \bR^{n-m_0}\oplus \{O\}=0$ and $f_2$ is 
any linear functional satisfying $f_2|\{O\}\oplus\bR^{m_0+1}=0$
where we indicate by $\{O\}$ the trivial subspaces of the complements. 
Hence, $\llrrparen{f_2} \in \SI^{m_0\ast}= \SI^{n-i_0-j_0-2\dagger}$. 
Hence, $U^\ast$ is a strict join of $U'$ and 
a great sphere $\SI^{n-i_0-j_0-2\dagger}$.

(iii) This follows from taking the dual of the second case of (ii).

(iv) 
Since the definition is consistent with the classical one for properly convex domains, 
this follows. Also, one can derive this as contrapositive of (ii) and (iii)
since domains and their duals not covered by (ii) and (iii) are the properly convex domains.
 

(v) $U$ corresponds to the second case of (ii).

(vi) If $U$ has an empty interior, then $U$ is covered by (ii)
and (iii), or is a great sphere of dimension $< n$. 
(iii) corresponds to the case 
when the dual of $K$ has nonempty interior. 
 
(vii) The forward part is given by (iii) where $K$ is a singleton in $\SI^0$
and $i_0=0$. 
The converse part is the second case of (ii) where $K$ is of 
dimension zero and $j_0= n-1$.  

(viii) 	This follows from the item (v) since $(U^\ast)^\ast = U$.

(xi)  Proposition \ref{prelim-prop-classconv} shows that 
(ii), (iii), and (iv) cover all compact convex sets that are not great spheres. 
%
\end{proof} 


We also note that for any properly convex 
domain $K$, $K \subset H^k \subset \SI^k$ for 
an open hemisphere $H_k$, and 
a great sphere $\SI^j$ in an independent space, the interior of 
$K \ast \SI^j \subset H^k \ast \SI^j$ is in an affine space 
$H^{k+j+1} = H^k \times H^{j}\times \bR \subset \SI^n$. 
Hence, a join can be viewed as a product of a certain form.
We call this an {\em affine form} of a strict join.

\index{join!affine form} 

\subsection{Duality and geometric limits}

Define the {\em thickness} of a properly convex domain $\Delta$ is given as 
\[\min\{ \max\{\bdd(x, \Bd \Delta)|x \in \Delta\}, \max\{\bdd(y, \Bd \Delta^*)|y \in \Delta^*\}\} \]
for the dual $\Delta^*$ of $\Delta$. 

\begin{figure}[h]
	\centering 
	\includegraphics[height=5cm]{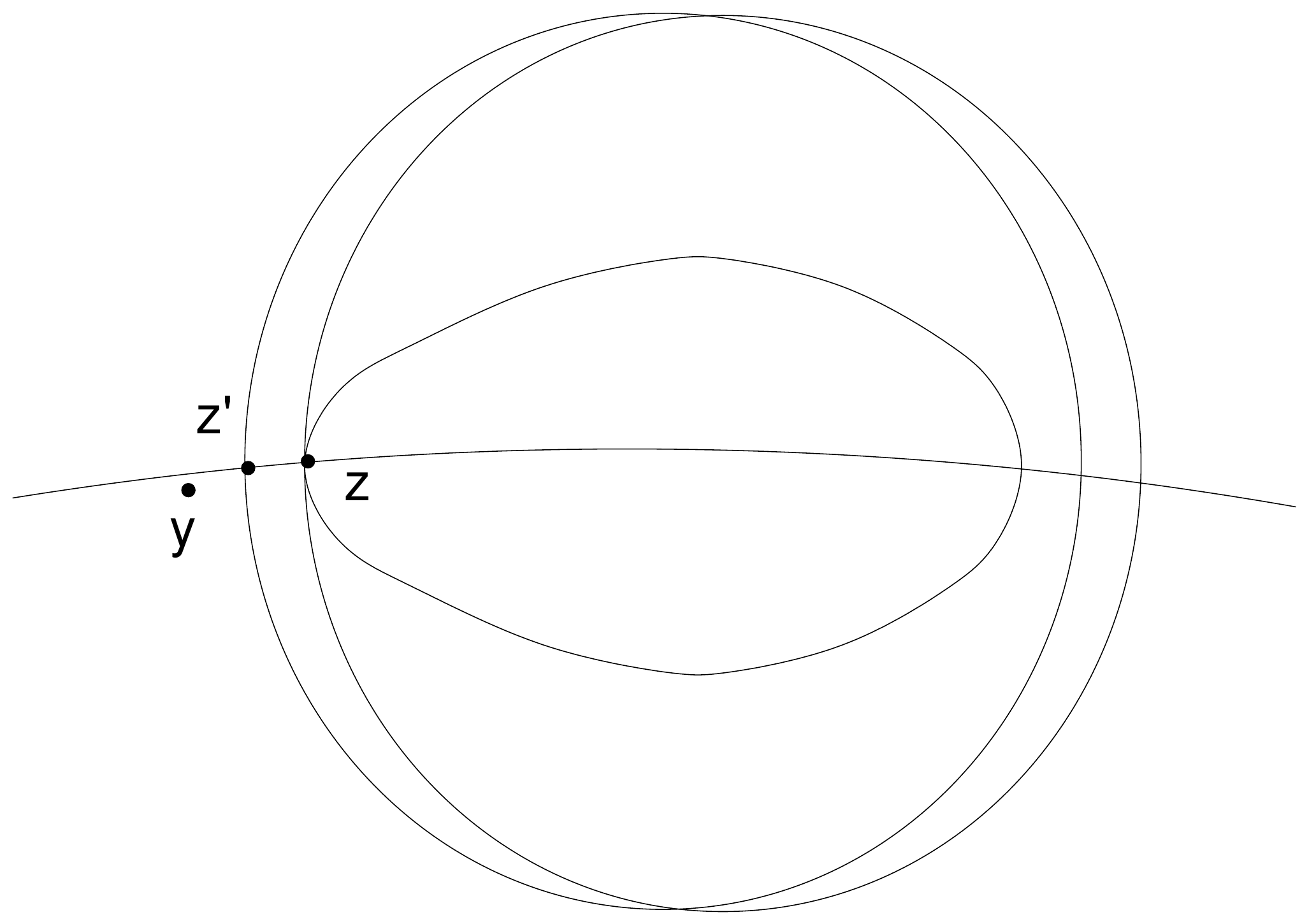}
	
	\caption{The diagram for Lemma \ref{cl-lem-convHaus}. }
	\label{cl-fig-figdual}
\end{figure}


\begin{lemma}\label{cl-lem-convHaus} 
	Let $\Delta$ be a properly convex open {\rm (}resp. compact{\rm)} domain in $\RP^n$ {\rm (}resp. $\SI^n$ {\rm )} and its dual $\Delta^*$  in $\RP^{n *}$ {\rm (}resp. $\SI^{n\ast}$ {\rm ).}
	Let $\eps$ be a positive number less than the thickness of $\Delta$
	and less than $\frac{1}{2}\bdd(\Delta', \mathcal{A}(\Delta'))$ and 
	$\frac{1}{2}\bdd(\Delta^{\prime \ast}, \mathcal{A}(\Delta^{\prime \ast}))$
	for a lift $\Delta'$ of $\Delta$ to $\SI^n$
	{\rm (}resp. $\Delta' = \Delta${\rm )}.  

	Then the following statements hold\/{\rm :} 
	\begin{itemize}
		\item $N_\eps(\Delta) \subset (\Delta^* - \clo(N_\eps(\Bd \Delta^*)))^*$.
		\item 
 Let $\Delta_1$ and $\Delta_2$ be two properly convex open domains. 
If the Hausdorff distance of $\Delta_1$ and $\Delta_2$ is less than  
		$\eps$ less than the thickness of $\Delta_1$ and $\Delta_2$, 
		then $\Delta^*_1$ and $\Delta_2^*$ are of Hausdorff distance $< \eps$.
		\item Furthermore, if 
		$\Delta_2 \subset N_{\eps'}(\Delta_1)$ and $\Delta_1 \subset N_{\eps'}(\Delta_2)$ for $0< \eps'< \eps$, 
		then we have $\Delta_2^* \supset \Delta_1^* - \clo(N_{\eps'}(\Bd \Delta_1^*))$
		and $\Delta_1^* \supset \Delta_2^* - \clo(N_{\eps'}(\Bd \Delta_2^*))$.
	\end{itemize} 
\end{lemma} 
\begin{proof} 
	Using the double covering map $p_{\SI^n}$ 
	and $p_{\SI^{n\ast}}:\SI^{n\ast} \ra \RP^{n\ast}$ of  the unit spheres in 
	$\bR^{n+1}$ and $\bR^{n+1\ast}$,
	we take components of $\Delta$ and $\Delta^*$. 
We can verify the result for properly convex open 
	domains in $\SI^n$ and $\SI^{n*}$ is sufficient. 
	
	For elements $\phi \in \SI^{n*}$ and $x \in \SI^{n}$, 
	we say $\phi(x) < 0$ if $f(v) < 0$ for $\phi =[f], x =[\vec{v}]$ 
	for $f \in \bR^{n+1 *}, \vec{v} \in \bR^{n+1}$. 
	Similarly, we say $\phi(x) > 0$ if $f(v) > 0$ for $\phi =[f], x =[\vec{v}]$ 
	for $f \in \bR^{n+1 *}, \vec{v} \in \bR^{n+1}$. 
	
	For the first item, let $y \in N_\eps(\Delta)$. Suppose that $\phi(y) < 0$ for 
	\[\phi \in \clo((\Delta^* - \clo(N_\eps(\Bd \Delta^*))) \ne \emp.\] 
	Since $\phi \in \Delta^*$, the set of points of $\SI^n$ where $\phi$ is positive 
is an open hemisphere $H$
	containing $\Delta$ but not $y$. 
	The boundary $\Bd H$ of $H$ has a closest point $z \in \Bd \Delta$ 
	of distance $< \eps$. The closest point $z'$ to $z$  on $\Bd H$ is in $N_\eps(\Delta)$ 
	since $y$ is in $N_\eps(\Delta) - H$ and $z'$ is closest to $\Bd \Delta$.
	The great circle $\SI^1$ containing $z$ and $z'$ is perpendicular to $\Bd H$
	since $\ovl{zz'}$ minimizes length.  
	Hence  $\SI^1$ passes through the center of the hemisphere.
	One can push the center of the hemisphere on $\SI^1$
	until it becomes a sharply supporting hemisphere to $\Delta$. The corresponding $\phi'$ is in $\Bd \Delta^*$ and 
	the distance between $\phi$ and $\phi'$ is less than $\eps$. This is a contradiction. 
	Thus, the first item holds (See Figure \ref{cl-fig-figdual}.)
	
	

	For the final item, we have that 
	\[\Delta_2 \subset N_{\eps'}(\Delta_1),  \Delta_1 \subset N_{\eps'}(\Delta_2) \hbox{ for } 0< \eps' <  \eps.\] 
	Hence, $\Delta_2 \subset (\Delta_1^* - \clo(N_{\eps'}(\Bd \Delta_1^*)))^*$, and 
	$\Delta_2^* \supset \Delta_1^* - \clo(N_{\eps'}(\Bd \Delta_1^*))$
	by \eqref{prelim-eqn-dualinc},
	which proves the third item where we need to switch $1$ and $2$ also. 
By a distance argument,	we can show that 
$\Delta_2^* \supset \Delta_1^* - \clo(N_{\eps'}(\Bd \Delta_1^*))$ also shows 
 $N_\eps(\Delta_2^*) \supset \Delta_1^*$ and vice versa. 
	The second item follows. \hfill \SnT   {\parfillskip0pt\par}
\end{proof}

The following may not hold for $\RP^n$: 
\begin{proposition} \label{prelim-prop-dualHausdorff}
	Suppose that $\{K_i\}$ is a sequence of properly convex domains in $\SI^n$ 
	geometrically converging to a compact convex set $K$. 
	Then $\{K_i^\ast\}$ geometrically converges to the compact convex set $K^\ast$ that is
	dual to $K$. 
%
	\end{proposition} 
\begin{proof} 
%
	Recall the compact metric space of all compact subsets of $\SI^n$ with 
	the Hausdorff metric $\bdd_H$. (See p.280-281 of Munkres \cite{Munkres75}.)
	$K_i$ is a Cauchy sequence under the Hausdorff metric $\bdd_H$.
Since the space is a compact metric space, we may assume without loss of generality that 
$K_i^\ast$ is also a convergent sequence 
	under the Hausdorff metric $\bdd_H$ on $\SI^{n\ast}$ 
by choosing further subsequences if necessary. 
To further explain, if we show that any subsequence of a subsequence converges to $K^*$, 
this will prove that $K_i^\ast$ converges to $K^\ast$. 

The space of all compact subsets of
	$\SI^{n\ast}$ with the  Hausdorff metric 
is a compact metric space. (See Munkres \cite{Munkres75}.)
	
	Since each $K_i$ is contained in an $n$-hemisphere corresponding 
	to the linear functional $\phi_i$ with $\phi_i|C(K_i) \geq 0$, 
	we deduce that $K$ is contained in an $n$-hemisphere. 

Let $K^\infty$ denote the limit of the Cauchy sequence $\{K_i^\ast\}$. 
	We will show $K^\infty = K^\ast$. 
	
	First, we show $K^\infty \subset K^\ast$: 
	Let $\phi_\infty$ be the limit of a sequence $\phi_i$ where
	$\phi_i \in C(K_i^\ast)$ for each $i$. 
	By Proposition \ref{prelim-prop-BP}, it is sufficient to show 
	$\llrrparen{\phi_\infty} \in K^\ast$ for every such $\phi_\infty$. 
	We may assume that their Euclidean norms are always $1$
	under the standard Euclidean metric on $\bR^{n+1}$.
	We will show that 
	$\phi_\infty| C(K) \geq 0$. 
	
	Let $\SI^n_1$ denote the unit sphere in $\bR^{n+1}$ with a Fubini-Study path-metric 
	$\bdd_1$. 
We identify $\SI^n$ with $\SI^n_1$ for convenience. 
	The projection $\SI^n \ra \SI^n_1$ is an isometry from $\bdd$ to $\bdd_1$. 
	Then $C(K_i)\cap \SI^n_1 \ra C(K) \cap \SI^n_1$ geometrically
	under the Hausdorff metric $\bdd_{H, 1}$ associated with $\bdd_1$. 
	Let $N_\eps(U)$ denote the $\eps$-neighborhood of a subset $U$ of $\SI^n_1$ 
	under $\bdd_1$. 
	Since $K_i \ra K$, 
	we find a sequence $\eps_i$ such that 
	$N_{\eps_i}(C(K_i))\cap \SI^n_1 \supset C(K) \cap \SI^n_1$
	and $\eps_i \ra 0$ as $i \ra \infty$. 
	
	For any $\phi$ in $\bR^{n+1\ast}$ of unit norm and 
	every pair of points $x, y \in \SI^n_1$ 
	with $\phi(x) \geq 0$,   we obtain that
	\begin{equation} \label{prelim-eqn-phidelta}
	\bdd_1(x, y)\leq \delta \hbox{ implies } \phi(y) \geq -\delta. 
	\end{equation}
	We can prove this by integrating along the geodesic from $x$ to $y$ 
	by considering $\phi$ as a $1$-form of norm $1$. 
	
	Since $\min\{\phi_i|N_{\eps_i}(C(K_i)) \cap \SI^n_1\} \geq -\eps_i$ by 
	  \eqref{prelim-eqn-phidelta}, 
	  we obtain $\phi_i| C(K) \cap \SI^n_1 \geq - \eps_i$ for sufficiently large $i$. 
	  Since $\eps_i \ra 0$, 
	  we obtain $\phi_\infty| C(K) \cap \SI^n_1 \geq 0$,  
	and $\phi_\infty \in C(K)^\ast$.

	Conversely, we show $K^\ast \subset K^\infty$. 
	Let $\phi \in C(K)^\ast$. Then $\phi| C(K) \cap \SI^n_1 \geq 0$. 
	Define $\eps_i = \min\{\phi(C(K_i) \cap \SI^n_1)\} $.
	If $\eps_i \geq 0$ for sufficiently large $i$, 
	then $\phi \in C(K_i)^\ast$ 
	and $\llrrparen{\phi} \in K^\infty$ by Proposition \ref{prelim-prop-BP}, which
completes the proof in this case. 
	
	Suppose $\eps_i < 0$ for infinitely many $i$. 
	By taking a subsequence if necessary, we assume that $\eps_i < 0$ for all $i$. 
	Let $H_\phi$, $H_\phi \subset \SI^n_1$, 
	be the hemisphere determined by the nonnegative condition of $\phi$. 
	Then $K_i - H_\phi \ne \emp$ for every $i$. 
	Choose a point $y_i$ in $K_i$ that maximizes the distance from $H_\phi$.
	Then $\bdd_1(y_i, H_\phi) = \delta_i$ for $0 < \delta_i \leq \pi/2$. 
	Since $\{K_i\} \ra K$, we deduce $\{\delta_i\} \ra 0$ as $i\ra \infty$ obviously. 
	Assume $\delta_i < \pi/4$ without loss of generality.
	
	Define a distance function $f_1(\cdot):=\bdd_1(\cdot, H_\phi): \SI^n_1  \ra \bR_+$. 
	Then $y_i$ is contained in a smooth sphere $S_{\delta_i}$ at the level $\delta_i$
	with the center $x_\phi$ of the complementary hemisphere of $H_\phi$. 
	Also, $K_i$ is contained in the complement of the convex open ball 
	$B_{\delta_i}$  bounded by $S_{\delta_i}$. 
	
	Now, we apply the convex geometry. 
	Let $H_i$ denote the hemisphere whose boundary contains $y_i$ and 
	tangent to $S_{\delta_i}$ and disjoint from $B_{\delta_i}$.
	Then $y_i$ is a unique maximum point of $ f_1|K_i$ since otherwise 
	there would exist a point with a larger $f_1$-value by the convexity of $K_i$ and $B_{\delta_i}$.  
	Moreover, $K_i$ is disjoint from $B_{\delta_i}$ as $y_i$ is the unique maximum point. 
	Since $K_i$ and $\clo(B_{\delta_i})$ are both convex and meets only at $y_i$, 
	$\partial H_i$ supports $K_i$ and $\clo(B_{\delta_i})$ by the 
	hyperspace separation theorem applied to 
	$C(K_i)$ and $C(\clo(B_{\delta_i}))$. 
	
	
	We obtained $K_i \subset H_i$. 
	Let $\phi_i$ be the linear functional of unit norm corresponding to $H_i$. 
	Then $\phi_i| C(K_i) \geq 0$.
	Let $s_i$ be the shortest segment from $y_i$ to $\partial H_\phi$
	whose other endpoint is $x_i \in \partial H_\phi$. 
	The center $\mathcal{A}(x_{\phi})$ 
	of $H_\phi$ is on the great circle $\hat s_i$ containing $s_i$. 
	The center of $H_i$ is on $\hat s_i$ and of distance $\delta_i$ 
	from $\mathcal{A}(x_{\phi})$ 
	since $\bdd_1(y_i, x_i) = \delta_i$. 
	
	This implies that $\bdd(\phi, \phi_i) = \delta_i$. 
	Since $\delta_i \ra 0$, we obtain $\{\phi_i\} \ra \phi$ and
	 $K^\ast \subset K^\infty$
	by Proposition \ref{prelim-prop-BP}. 
	\end{proof}

\chapter{Examples of properly convex real projective orbifolds with ends: cusp openings}
\chaptermark{Examples of properly convex real projective orbifolds}
\label{ch-ex} 


We present examples where our theory applies. 
We discuss the theory of convex projective structures on 
	Coxeter orbifolds  
	as well as the orderability theory for Coxeter orbifolds.
	We explain our work jointly done with Gye-Seon Lee and Craig Hodgson generalizing the work of 	Benoist and Vinberg.
	We also introduce the vertex orderable Coxeter orbifolds.
	We state the work of Heusner-Porti on projective deformations of 
	the hyperbolic link complement and the subsequent work by Ballas. 
	Also, we state some nice results on 
	finite volume convex real projective structures by Cooper-Long-Tillmann 
	and Crampon-Marquis on horospherical ends
	and thick and thin decomposition. 

How, these examples well-fitting into this monograph are explained in Chapter \ref{ex-sec-nicecase}. 


\section{History of examples} 

Originally, Vinberg \cite{Vinberg71} investigated convex real projective 
Coxeter orbifolds as 
linear groups acting on convex cones. The groundbreaking work also produced
many examples of real projective orbifolds and manifolds $\orb$ relevant to our study. 
For example, see Kac and Vinberg \cite{Vinberg67} for the deformation of triangle groups. 
However, the work was reduced to studying some Cartan forms with rank $n+1$ for $n=\dim \orb$. 
The method proves to be challenging in computing explicit examples. 

Later, Benoist \cite{Benoist062} worked out some examples on prisms.  
Generalizing this, 
Choi \cite{Choi06} studied the orderability of Coxeter orbifolds
after conversing with Kapovich about the deformability. 
This produced many examples of noncompact orbifolds with 
properly convex projective structures 
based on the work of Vinberg. 
Later, Marquis \cite{Marquis10} generalized the technique 
to study the convex real projective structures based on Coxeter 
orbifolds with truncation polytopes as base spaces. 
Since these are compact orbifolds, we will not discuss these in this section. 

For compact hyperbolic $3$-manifolds, Cooper-Long-Thistlethwaite \cite{CLT06}
and \cite{CLT07} produced many examples with deformations using numerical
methods. Some of these are exact computations. 

We now discuss the noncompact strongly tame orbifolds with
convex real projective structures. 

Also, Choi, Hodgson, and Lee \cite{CHL12} computed the restricted 
deformation spaces of convex 
real projective structures of some complete 
hyperbolic Coxeter orbifolds with ideal vertices but no finite vertices, and 
Choi and Lee \cite{CL15} showed that all compact hyperbolic 
weakly orderable Coxeter orbifolds have local deformation spaces 
of dimension $e_+ - 3$ where $e_+$ is the number of ridges with order $\geq 3$.
These Coxeter orbifolds form a large class of Coxeter orbifolds.  

Our future plan is to generalize these results to complete hyperbolic 
Coxeter orbifolds that are weakly orderable with respect to ideal vertices. 

For noncompact hyperbolic $3$-manifolds, Porti and Tillmann \cite{PTp}, 
Cooper-Long-Tillmann \cite{CLT15}, and Crampon-Marquis \cite{CM14}
made theories where the ends were restricted to be horospherical. 
Ballas \cite{Ballas14} and \cite{Ballas15} made initial studies
of deformations of complete hyperbolic $3$-manifolds to 
convex real projective ones. 
Cooper, Long, and Tillmann \cite{CLT18} have produced a deformation theory
for convex real projective manifolds  
parallel to ours with different types of restrictions on
ends, such as requiring the end holonomy group to be abelian. 
They focus on the openness of the deformation spaces. 
We will provide our theory in Part 3.

\section{Examples and computations}
\label{ex-sec-examples}

We will give some series of examples due to the author and many other 
people. 
Here, we will not give compact examples as we already gave some survey 
in Choi-Lee-Marquis \cite{CLM18}.


Given a polytope $P$, a {\em face} is a codimension-one side of $P$.
A {\em ridge} is the codimension-two side of $P$. 
When $P$ is $3$-dimensional, a ridge is called an {\em edge}. 

We will concentrate on $n$-dimensional orbifolds whose base spaces 
are homeomorphic to convex Euclidean polyhedrons and whose faces are 
silvered and each ridge is given an order. 
For example, a hyperbolic polyhedron with edge angles of form 
$\pi/m$ for positive integers $m$ will have a natural orbifold 
structure. 


\begin{definition}\label{ex-defn-Coxeter} 
	A {\em Coxeter group} $\Gamma$ is an abstract group
	defined by a group presentation of the form 
	\[(R_i;(R_iR_j)^{n_{ij}}), i, j \in I\]
	where $I$ is a countable index set, 
	$n_{ij} \in \mathbf{N}$ is symmetric for $i, j$ and $n_{ii}=1$, 
	$n_{ij} \geq 2$ for $i\ne j$. 
\end{definition} \index{Coxeter group}

The fundamental group of the orbifold is a Coxeter group with 
a presentation 
\[ R_i, i=1,2,\dots, f, (R_iR_j)^{n_{ij}} = 1\] 
where each $R_i$ is associated with a silvered sidesand 
$R_iR_j$ has order $n_{ij}$ associated with the edge 
formed by the intersection of the $i$-th and $j$-th sides.

Let us consider only the $3$-dimensional orbifolds for now. 
Let $P$ be a fixed convex $3$-polyhedron. 
Let us assign an order at each edge. 
We let $e$ be the number of edges 
and $e_2$ be the number of edges of order-two. 
Let $f$ be the number of sides. 

For any vertex of $P$, we can remove it unless
the link in $P$ forms a spherical Coxeter $2$-orbifold of codimension one. 
This makes $P$ into a $3$-dimensional orbifold. 

Let $\hat P$ denote the differentiable orbifold with 
silvered sides, and the edge orders realized as assigned
from $P$ with the above vertices removed. 
We say that $\hat P$ has a {\em Coxeter orbifold structure}.  
\index{Coxeter orbifold structure}

In this chapter, we will exclude a {\em cone-type Coxeter orbifold}
whose polyhedron has a side $F$ and a vertex $v$ where all other sides 
are adjacent triangles 
to $F$ and contain $v$ and all ridge orders of $F$ are $2$.
Another type we will {\em not} study is a {\em product-type Coxeter orbifold} 
whose polyhedron is topologically a product of a polygon and an interval 
and edge orders of top and the bottom sides are all $2$. 
These orbifolds are essentially of lower dimensions.
Finally, we will {\em not} study Coxeter orbifolds with finite fundamental groups.
If $\hat P$ does not belong to any of the above types and does not admit 
an affine structure, then $\hat P$ is said to be a {\em normal-type} 
Coxeter orbifold. 


\begin{remark} \label{rem:geomded} 
Actually, there was an error in the author's paper \cite{Choi06}. We need to exclude Coxeter orbifolds admitting affine structures. 
We will try to mend the proof here. 
This correction is necessary for the proof of Proposition 2 in \cite{Choi06}:
First, note that if a reflection with a sphere of fixed points contains an antipodal fixed point of another reflection,
then they must commute, their associated sides must meet, and their edge order 
must be $2$.

Let $P$ be a properly convex fundamental polytope, and let $F_i$ denote the sides of $P$
for $i$ in an index set $I$. 
Let $R_i$ be the reflections on the sides $F_i$ of $P$. 

We need to show that there does not exist a holonomy invariant union of one or two $1$-dimensional subspace  or 
holonomy invariant $2$-dimensional subspace. 

Suppose that $l$ is a holonomy invariant union of one or two $1$-dimensional subspaces.
The case of
 two $1$-dimensional subspaces reduces to the first one because a reflection must act on 
each $1$-dimensional subspace if it acts on a union of two $1$-dimensional subspaces. 
Also, if two faces are adjacent at an edge, then the edge order must be finite. 

Let $l$ be a holonomy invariant $1$-dimensional subspace. 
If a face is contained in a $2$-dimensional subspace containing $l$, we call it a {\em parallel} face.
The associated reflection is also called {\em parallel}. 
Otherwises, 
it is called a {\em transverse} face. The associated reflection is called {\em transverse}. 
Then the antipodal fixed point must lie on $l$ and the fixed point subspace must intersect $l$. 

Suppose that $\Omega \cap l = \emptyset$. Then we can take a union of open hemispheres with boundary $l$ 
meeting $\Omega$. Then $P$ must have at least two parallel sides $F_i$ and $F_j$. 
In this case, $R_i$ and $R_j$ commute with all other 
reflections at the sides. Hence, we violated the normality. 
Otherwise, we cannot have a compact $P$: If there is one parallel face, then
 the group is simply an extension of 
$2$-dimensional Coxter group by a reflection of the face,
and we have a cone-type orbifold, contradicting the normality assumption. 
This result follows from the second paragraph above. 
If there is no parallel face, then each sphere containing $l$ 
is invariant by the nature of transversal reflections, 
and again we have a $2$-dimensional Coxeter group. 
These do not admit properly convex fundamental polytope with some vertices removed.

Suppose that $\Omega \cap l \ne \emptyset$. 
If there exists a pair of adjacent transversal faces $F_i$ and $F_j$,
then $R_i$ and $R_j$ generate a finite group, and by their action and convexity, we must have $\Omega \supset l$.
This implies $\Omega = \SI^3$ by convexity, and we must have a finite Coxeter group, contradiction the normality. 
Also, there can be at most two parallel faces since $P$ is a properly convex polytope. 
From these two facts, it follows that there must be at most one transverse face.
We cannot constuct properly convex $P$ in this situation. 

Suppose that $S$ is a holonomy invariant $2$-dimensional subspace. 
If $S \cap \Omega = \emptyset$, then there exists an 
invariant affine subspace containing $\Omega$, 
implying that our orbifold admits an affine structure, contradicting the normality assumption. 

The old proof correctly rules out  $S \cap \Omega \ne \emptyset$. We explain a bit more. 
In the old proof, we chose the fundamental polytope $P$ such that $P\cap S \ne \emptyset$. Also, to deduce 
the infinite edge orders, we tacitly used the fact that 
our holonomy group acts as a $2$-dimensional Coxeter group on $S\cap\Omega$. 

\end{remark}


A huge class of examples are obtainable from complete hyperbolic $3$-polytopes with dihedral angles that are submultiples of $\pi$. (See Andreev \cite{Andreev70} and Roeder \cite{Roeder07}.)


\begin{definition}\label{ex-defn-deformation}
	The {\em deformation space} $\mathfrak{D}(\hat P)$ of projective structures on a Coxeter orbifold $\hat P$
	is the space of all projective structures on $\hat P$ quotient 
	by isotopy group actions of $\hat P$.
\end{definition}
\index{deformation space} \index{d@$\mathfrak{D}(\cdot)$} 

This definition has also used in a number of papers \cite{Choi06}, \cite{Choi12}, 
and \cite{CHL12}. 
The topology on $\mathfrak{D}(\hat P)$ is defined by as follows: 
$\mathfrak{D}(\hat P)$ is the quotient space of 
the space of the development pairs $(\dev, h)$ with 
the compact open $C^r$-topology, 
$r \geq 2$, for the maps $\dev: \tilde P \ra \RP^n$.
\index{crtopology@$C^r$-topology} 

In Proposition \ref{op-prop-DPCDef}, 
we will explain that the space is identical with 
$\CDef_{{\mathcal E}}(\hat P)$. 

%


%
	%


A point $p$ of $\mathfrak{D}(\hat{P})$ gives a fundamental polyhedron $P$ in $\RP^3$,
well-defined up to projective automorphisms. 
By Proposition \ref{op-prop-DPCDef}, $\mathfrak{D}(\hat{P})$ can be identified 
with $\CDef_{{\mathcal E}}(\hat P)$. 
We concentrate on the space
of $p \in \mathfrak{D}(\hat{P})$ giving
a fundamental polyhedron $P$ fixed up to projective automorphisms. 
This space is called the \emph{restricted deformation space} of $\hat{P}$ and denoted by $\mathfrak{D}_P(\hat{P)}$. A point $t$ in $\mathfrak{D}_P(\hat{P)}$ is said to be \emph{hyperbolic} if a hyperbolic structure on $\hat{P}$
induces the projective structure; 
that is, it is projectively diffeomorphic to 
$\bB/\Gamma$ for a standard unit ball $\bB$ and a discrete group 
$\Gamma \subset \Aut(\bB)$. 
A point $p$ of $\mathfrak{D}(\hat P)$ always 
determines a fundamental polyhedron $P$ up to projective 
automorphisms  because $p$ also determines 
reflections corresponding to sides up to conjugations.  
\index{restricted deformation space} \index{dp@$\mathfrak{D}_P(\hat P)$}


The work of Vinberg \cite{Vinberg71} implies that each element of $\mathfrak{D}(\hat P)$ gives a convex projective structure (see Theorem 2 of \cite{Choi06}).
That is, the image of the developing map of the orbifold universal cover
of $\hat P$
is projectively diffeomorphic to a convex domain in $\RP^3$, 
and the holonomy homomorphism  is a discrete faithful representation.

Now, we state the key property in this chapter:
\begin{definition}\label{ex-def-orderable}
	We say that $P$ is {\em orderable} if we can order the sides of $P$ such that 
	each side meets sides of higher index  in at most $3$ edges. 
	We say that a Coxeter orbifold $\hat P$ is {\em orderable} if we can order the sides of $\hat P$ 
such that the set of order-two edges and an edge from higher order sides
on a side of $P$ has the cardinality $\leq 3$. 
\end{definition}
\index{Coxeter orbifold!orderable} 


A pyramid with a complete hyperbolic structure and 
dihedral angles that are submultiples of $\pi$ 
is an obvious example. See Proposition 4 of \cite{Choi06} worked out with 
J. R. Kim. 

Another example is a drum-shaped convex polyhedron which has top and bottom sides 
of the same polygonal type, and each vertex of the bottom side is connected to two vertices on
the top side and vice versa. Another example is  a convex polyhedron where 
the union of triangles separates each pair of the interiors of nontriangular sides. 
In these examples, since all nontriangular sides are separated by the union of triangular sides, 
the sides are either level $0$ or level $1$, and hence they satisfy the orderability condition.


If $P$ is compact, Marquis \cite{Marquis10} showed that $P$ 
is a truncation polytope. That is, one starts from a tetrahedron, 
one cuts a neighborhood of a vertex so as to change the combinatorial 
type near that vertex only, and one iterates this process finitely many times. 
Many of these can be realized as compact
hyperbolic polytopes with dihedral angles equal to submultiples of $\pi$
by the Andreev theorem \cite{Roeder07}. 
If $P$ is not compact, we do not have a complete classification.
Also, infinitely many of these can be realized as complete 
hyperbolic polytopes with dihedral angles that are submultiples of $\pi$. 
(D. Choudhury was first to show this.)

\begin{definition}
	We denote by $k(P)$ the dimension of the projective group acting on a convex polyhedron $P$.
\end{definition}

The dimension $k(P)$ of the subgroup of $G$ acting on 
$P$ equals $3$ if $P$ is a tetrahedron and equals $1$ if 
$P$ is a cone with base a convex polyhedron which is not a triangle. Otherwise, $k(P) = 0$. 


\begin{definition}\label{ex-defn-deformable}
	Let $P$ be a 3-dimensional hyperbolic Coxeter polyhedron, and let $\hat{P}$ denote its Coxeter orbifold structure.
	Suppose that $t$ is the corresponding hyperbolic point of $\mathfrak{D}_P(\hat{P})$.
	We call a neighborhood of $t$ in $\mathfrak{D}_P(\hat{P})$ the \emph{local restricted deformation space} of $P$. We say that $\hat{P}$ is \emph{projectively deformable relative to the mirrors}, or
	simply \emph{deforms rel mirrors},
	if the dimension of its local restricted deformation space is positive.
	Conversely, we say that $\hat{P}$ is \emph{projectively rigid relative to the mirrors}, or  \emph{rigid rel  mirrors}, if the dimension of
	its local restricted deformation space is $0$.
\end{definition} 
\index{projectively rigid relative to the mirrors} 

\begin{theorem}\label{ex-thm-main} 
	Let $P$ be a convex polyhedron, and let $\hat P$ be given a normal-type Coxeter orbifold 
	structure. 
	Suppose that $\hat P$ is orderable.
	Then the restricted deformation space of projective 
	structures on the orbifold $\hat P$ is a smooth manifold
	of dimension $3f - e - e_2 - k(P)$ if it is not empty.
\end{theorem}

If we start from a complete hyperbolic polytope whose dihedral angles 
are submultiples of $\pi$, we know that the restricted deformation space 
is not empty. 

If we assume that $P$ is compact, then we refer to Marquis \cite{Marquis10} for the complete theory. However, the topic lies beyond the scope of this monograph.


The following theorem describes the local restricted deformation space for a class of
Coxeter orbifolds arising from  {\em ideal hyperbolic polytope}, i.e. a polytope with all vertices on the sphere at infinity.

\begin{theorem}[Choi-Hodgson-Lee \cite{CHL12}]\label{ex-thm-CHL}
	Let $P$ be an ideal $3$-dimensional hyperbolic polyhedron
	whose dihedral angles are all equal to $\pi/3$, and suppose that $\hat{P}$ is given its Coxeter orbifold structure.
	If $P$ is not a tetrahedron,
	then a neighborhood of the hyperbolic point in
	$\mathfrak{D}_P(\hat{P)}$ is a smooth $6$-dimensional manifold.
\end{theorem}


The main ideas in the proof of Theorem \ref{ex-thm-CHL} are as follows.
We first show that
$\mathfrak{D}_P(\hat{P)}$
is isomorphic to the solution set of  a system of polynomial equations, following ideas of Vinberg \cite{Vinberg71} and Choi \cite{Choi06}. Since the faces of $P$ are fixed,
each projective reflection
in a face of the polyhedron is determined by a reflection vector $b_i$.
We then compute the Jacobian matrix of the equations for the $b_i$ at the hyperbolic point. This
reveals that the matrix has exactly the same rank as the Jacobian matrix of the equations for
the Lorentzian unit normals of a hyperbolic polyhedron with the given dihedral angles. By the infinitesimal rigidity of the hyperbolic structure on $\hat P$, this matrix is of full rank and has the kernel of dimension six; the result then follows from the  implicit function theorem.
In fact, we can interpret the infinitesimal projective deformations as applying infinitesimal hyperbolic isometries to the reflection vectors

As Hodgson later pointed out,
we can generalize the above theorem slightly: 

\begin{definition} \label{ex-dfn-hyperideal}
	Let $X$ be a hyperbolic $n$-orbifold with totally geodesic boundary component 
	diffeomorphic to an $(n-1)$-orbifold $\Sigma$.
	Let $\tilde X$ denote the universal cover in the Klein model $\bB$ in 
	$\SI^n$. Let $\Gamma$ be the group of deck transformations considered 
	as projective automorphisms of $\SI^n$. 
	Then a complete hyperbolic hyperspace $\tilde \Sigma$ covers $\Sigma$. 
	Every component of the inverse image of $\Sigma$ is of the form $g(\tilde \Sigma)$
	for $g \in \pi_1(X)$.  
	A point $\mbv_{\tilde \Sigma} \in \SI^n - \bB -\mathcal{A}(\bB)$ 
	is projectively dual to the hyperspace containing $\tilde \Sigma$
	with respect to the bilinear form $B$. 
	(See Section \ref{intro-sub-horo}.)
	Then we form the join 
	$C:=\{\mbv_{\tilde \Sigma}\} \ast \tilde \Sigma -\{\mbv_{\tilde \Sigma}\}$. 
	Then we form $\hat C:= X \cup \bigcup_{g \in \Gamma} g(C)$. 
	$\hat C/\Gamma$ is an $n$-orbifold with radial ends. 
	We call the ends the {\em hyperideal ends}. 	
	(See also  \cite{BaoBonahon}, \cite{Rousset}, and \cite{Schlenker}.) 
	\end{definition}
\index{end!hyperideal}

A point of $\mathfrak{D}_P(\hat{P)}$ corresponding to a hyperbolic $n$-orbifold 
with hyperideal ends added will again be called a {\em hyperbolic point}. 
An {\em  $3$-dimensional hyperbolic polyhedron with possibly 
hyperideal vertices} is a compact convex polyhedron with 
vertices outside $\bB$ removed where no $1$-dimensional open edge is outside $\bB$.



\begin{corollary}[Choi-Hodgson-Lee]\label{ex-thm-CHL2}
	Let $P$ be an ideal $3$-dimensional hyperbolic polyhedron
	with possibly hyperideal vertices whose dihedral angles are of the form $\pi/p$ for integers $p \geq 3$, and suppose that $\hat{P}$ is given its Coxeter orbifold structure.
	If $P$ is not a tetrahedron,
	then a neighborhood of the hyperbolic point in
	$\mathfrak{D}_P(\hat{P)}$ is a smooth $6$-dimensional manifold.
\end{corollary}
We did not give a proof for the case where some edges orders are greater than equal to $4$
in the article \cite{CHL12}. 
We can allow any of our end orbifold to be 
a $(p, q, r)$-triangle reflection orbifold for $p, q, r \geq 3$.
The same proof will apply as first observed by Hodgson:
We modify the proof of Theorem 1 of the article in Section 3.3 of \cite{CHL12}. 
Let $\partial_\infty \hat P$ denote the union of end orbifolds of $\hat P$ which 
are orbifolds based on $2$-sphere with singularities that admit either a Euclidean
or hyperbolic structures. Let $h:\pi_1(\hat P) \ra \PO(3, 1)\subset \PGL(4, \bR)$ denote 
the holonomy homomorphism associated with the convex real projective structure 
induced from the hyperbolic structure. 
We only need to show 
\[H^1(\hat P, so(3, 1)_{Ad_h}) = 0, 
H^1(\partial_\infty P, so(3, 1)_{Ad_h}) = 0.\] 
Recall that a $(p, q, r)$-triangle reflection orbifold with $1/p + 1/q + 1/r < 1$ 
has a rigid hyperbolic and conformal structure. 
By Corollary 2 of \cite{Sullivan82}, the representation into $\PO(3, 1)$ is rigid. 
The first part of the equation follows. 
The second part also follows from Corollary 2 of \cite{Sullivan82}.
These examples are convex as shown by the work of Vinberg \cite{Vinberg71}.  
Corollary \ref{ex-cor-closed2} implies the proper convexity. 

We comment that 
we are using Theorem 7 (Sullivan rigidity) of \cite{Sullivan82}
as the generalization of the Garland-Raghunathan-Weil rigidity \cite{GR70} \cite{Weil62}. 



\subsection{Vertex orderable Coxeter orbifolds} 

\subsubsection{Vinberg theory} 

Let $\hat P$ be a Coxeter orbifold of dimension $n$. 
Let $P$ be the fundamental convex polytope of  $\hat P$. 
The reflection is determined by a point, called a reflection point, and a hyperspace. 
Let $R_i$ be a projective
 reflection on a hyperspace $S_i$ containing a side of $P$.
 Then we can write
\[ R_i:= \Idd - \alpha_i \otimes \vec{v}_i       \]
where $\alpha_i$ is zero on $S_i$ and $\vec{v}_i$ is 
the reflection vector and $\alpha_i(\vec{v}_i) = 2$. 

Given a reflection group $\Gamma$, 
we form a Cartan matrix $A(\Gamma)$ given by $a_{ij} := \alpha_i(\vec{v}_j)$. 
Vinberg \cite{Vinberg71} 
proved that the following conditions are necessary and sufficient for 
$\Gamma$ to be a linear Coxeter group:
\begin{enumerate}
	\item[(C1)] $a_{ij} \leq 0$ for $i \ne j$, and $a_{ij} =0$ if and only if 
	$a_{ji} = 0$.
	\item[(C2)] $a_{ii} = 2$, and 
	\item[(C3)] for $i \ne j$, $a_{ij}a_{ji} \geq 4$ or
	$a_{ij}a_{ji} = 4\cos^2(\frac{\pi}{n_{ij}})$ for an integer $n_{ij}$.
\end{enumerate}
The Cartan matrix is an $f\times f$-matrix when $P$ has $f$ sides. 
Also, $a_{ij} = a_{ji}$ for all $i, j$ if $\Gamma$ is conjugate to a reflection group
in $O^+(1, n)$. This condition ensures that $\hat P$ is a hyperbolic 
Coxeter orbifold. 

The Cartan matrix is determined only up to an action of  
the group $D_{f, f}$ of nonsingular diagonal matrices: 
\[A(\Gamma) \ra DA(\Gamma)D^{-1} \hbox{ for } D \in D_{f,f}.\] 
This is due to the ambiguity of choices
\[ \alpha_i \mapsto c_i \alpha_i, \vec{v}_i \mapsto \frac{1}{c_i} \vec{v}_i, c_i > 0.\]
Vinberg showed that the set of all {\em cyclic invariants} of the form 
$a_{i_1i_2}a_{i_2i_3} \cdots a_{i_ri_1}$ 
classifies isomorphic semisimple linear Coxeter groups generated by reflections up to 
conjugation. (See Fact 3.28, Fact 3.29 of \cite{DGKLM} and Proposition 16 of \cite{Vinberg71}.) 

\subsubsection{The classification of convex real projective structures on 
	triangular reflection orbifolds} 
We will follow Kac-Vinberg \cite{Vinberg67}. 
Let $\hat T$ be a $2$-dimensional Coxeter orbifold based on a triangle $T$.
Let the edges of $T$ be silvered. Let the vertices be given orders 
$p, q, r$ where $1/p + 1/q + 1/r \leq 1$. 
If $1/p + 1/q + 1/r \leq  1$, then 
the universal cover $\tilde T$ of $\hat T$ is a properly convex domain or a complete affine plane by Vinberg \cite{Vinberg71}. 
We can find the topology of $\mathfrak{D}(\hat T)$ as Goldman did in his senior thesis \cite{Gthesis}. 
We may put $T$ as a standard triangle with vertices $\vec{e}_1:=[1,0,0], \vec{e}_2:=[0,1,0], \vec{e}_3:= [0,0,1]$. 

Let $R_i$ be the reflection on a line containing $[\vec{e}_{i-1}], [\vec{e}_{i+1}]$ and with 
a reflection vertex $[\vec{v}_i]$. 
Let $\alpha_i$ denote the linear function on $\bR^3$ that takes zero values on 
$\vec{e}_{i-1}$ and $\vec{e}_{i+1}$. We choose $\vec{v}_i$ to satisfy $\alpha_i(\vec{v}_i)= 2$. 

When $1/p + 1/q + 1/r = 1$, 
the triangular orbifold admits a compatible Euclidean structure. 
When $1/p + 1/q + 1/r < 1$, 
the triangular orbifold admits a hyperbolic structure not
necessarily compatible with the real projective structure. 

A linear Coxeter group $\Gamma$ is hyperbolic if and only if the Cartan matrix A of $\Gamma$ is indecomposable, of negative type, and equivalent to a symmetric matrix of signature $(1,n)$.

Assume that none of $p, q, r$ is $2$ and $1/p + 1/q + 1/r < 1$. 
Let $a_{ij}$ denote the entries of the Cartan matrix. 
	It satisfies \[a_{12}a_{21} = 4 \cos^2 \pi/p, 
	a_{23}a_{32} = 4 \cos^2 \pi/q, a_{13}a_{31} = 4 \cos^2 \pi/r.\]
There are only two cyclic invariants 
$a_{12}a_{23}a_{31}$ and $a_{13}a_{32}a_{21}$
satisfying 
\[a_{12}a_{23}a_{31}a_{13}a_{32}a_{21}= 64 \cos^2 \pi/p \cos^2 \pi/q \cos^2 \pi/r.\] 
Then 
the {\em triple invariant} $a_{12}a_{23}a_{31} \in \bR_+$ classifies 
the conjugacy classes of $\Gamma$. 
A single point of $\bR_+$ corresponds to a hyperbolic structure.
For different points, they are properly convex by \cite{cgorb}. 
\index{triple invariant} 

Since $a_{ij} = a_{ji}$ for geometric cases, 
we obtain that 
\[a_{12}a_{23}a_{31} = 2^3  \cos(\pi/p) \cos(\pi/q)\cos(\pi/r)\]
gives the unique hyperbolic points. 

We conclude that for this orbifold 
$\mathfrak{D}(\hat T) \cong \bR_+$ the space of the triple invariants. 

\begin{example}[Lee's example] \label{ce-exmp-Lee}
Consider the Coxeter orbifold $\hat P$ with the underlying space on a polyhedron $P$ 
which has the combinatorics of a cube with all sides mirrored and
all edges given order $3$ but with vertices removed. 
By the Mostow-Prasad rigidity and the Andreev theorem, 
the orbifold has a unique complete hyperbolic structure. 
By Theorem \ref{ex-thm-CHL}, 
there exists a six-dimensional space of real projective structures on it 
where one has a projectively fixed fundamental domain 
in the universal cover.

There are eight ideal vertices of $P$ corresponding to eight ends of $\hat P$. 
Each end orbifold is a $2$-orbifold based on a triangle with edges mirrored, 
and vertex orders are all $3$. 
Each end orbifold has a real projective structure, and hence it is characterized 
by the triple invariant. 
Thus, each end has a neighborhood diffeomorphic 
to the $2$-orbifold multiplied by $(0, 1)$.
The eight triple invariants are related when we are working on the restricted deformation space 
since the deformation space is only six-dimensional. 


\end{example}


\subsubsection{The end mappings} 
We will give some explicit conjectural class of examples where we can control 
the end structures. We worked this out with Greene, Gye-Seon Lee, and Marquis starting from
the workshop at the ICERM in 2014. 

Let $\mathcal{F}$ be the set of faces of $C$ and give a total order $\leqslant$ on $\mathcal{F}$. A face $F'$ is $E_2$-greater than $F$ if $F < F'$ and $F \cap F'$ is an edge of label $2$.

A {\em flexible} vertex of a Coxeter orbifold is a vertex of the base polytope where 
no edge of order $2$ ends. 
Let $\mathcal{V}_f$ be a set of flexible vertices in $C$, and let $\mathcal{V}$ be a subset of $\mathcal{V}_f$. A face $F'$ is $\mathcal{V}$-greater than $F$ if $F < F'$ and there exists a face $F''$ such that $F < F''$ and $F \cap F' \cap F''$ is a vertex in $\mathcal{V}$. 

A Coxeter $3$-orbifold $\hat P$ is $\mathcal{V}$-orderable if there is no triangular face all of whose vertices are in $\mathcal{V}$ and the faces of $C$ can be ordered such that for each face $F$ of $C$, the number of the set of order $2$ edges of $F$ and edges in a face $\mathcal{V}$-greater than $F$ is less than or equal to $3$.



Let $\partial_{\mathcal{V}} \orb$ denote the disjoint union of end orbifolds corresponding to 
the set of ideal vertices $\mathcal{V}$. 

\begin{conjecture}[Choi-Greene-Lee-Marquis \cite{CGLMp}] \label{ex-thm-endinvmap}
We conjecture as follows: 
	Suppose that $P$ with a set of vertices $\mathcal{V}$ is $\mathcal{V}$-orderable, 
	and $P$ admits a Coxeter orbifold structure with a convex real projective 
	structure. Then the following map is onto
\[{\mathfrak{D}}(\orb) \ra {\mathfrak{D}}(\partial_{\mathcal{V}} \orb).\]
	\end{conjecture}


A Coxeter $3$-orbifold $\hat P$ is  {\em weakly $\mathcal{V}$-orderable} if there is no triangular face all vertices of which are in $\mathcal{V}$ and the faces of $C$ can be ordered such that for each face $F$ of $C$, the number of the set of edges in faces $E_2$-greater and $\mathcal{V}$-greater than $F$ is less than or equal to $3$. 

\begin{conjecture}[Choi-Greene-Lee-Marquis \cite{CGLMp}] \label{ex-thm-endinvmap2}
	Suppose that $P$ with a set of vertices $V$ is weakly $\mathcal{V}$-orderable. 
	Suppose $P$ admits a Coxeter orbifold structure with an ideal or hyperideal end
	structure. 
	Then the function $\mathfrak{D}(\orb) \ra \mathfrak{D}(\partial_{\mathcal{V}} \orb)$ is locally surjective 
	at the hyperbolic point.  
\end{conjecture}

\section{Some relevant results} 

The deformation spaces of 
convex structures on closed hyperbolic manifolds were extensively studied by 
Cooper-Long-Thistlethwaite \cite{CLT06} and \cite{CLT07}.

\subsection{The work of Heusener-Porti}

\begin{definition} 
Let $N$ be a closed hyperbolic $3$-manifold.  
We consider the holonomy representation of $N$
 \[\rho:  \pi_1(N) \ra \PSO(3, 1) \hookrightarrow   \PGL(4, \bR).\]
A closed hyperbolic three manifold $N$ is called {\em infinitesimally projectively rigid} if 
\[H^1(\pi_1(N), \mathfrak{sl}(4, \bR)_{{\mathrm{Ad}} \rho}) = 0. \]
\end{definition} 
\index{infinitesimally projectively rigid} 

\begin{definition} 
Let $M$ denote a compact $3$-manifold with boundary a union of tori and whose interior is hyperbolic with finite volume.
$M$ is called {\em infinitesimally projectively rigid relative to the cusps} if the inclusion $\partial M \ra   M$ induces an injective homomorphism
\[H^1(\pi_1(M), \mathfrak{sl}(4, \bR)_{{\mathrm{Ad}} \rho}) \ra H^1(\pi_1(\partial M), \mathfrak{sl}(4, \bR)_{{\mathrm{Ad}} \rho}).\]
\end{definition}

\begin{theorem}[Heusener-Porti \cite{HP11}] 
Let $M$ be an orientable $3$-manifold whose interior has a complete hyperbolic metric with finite volume. If $M$ is infinitesimally projectively rigid relative to the cusps, then infinitely many Dehn fillings on M are infinitesimally projectively rigid.
\end{theorem} 

\begin{theorem}[Heusener-Porti \cite{HP11}] 
Let $M$ be an orientable $3$-manifold whose interior has a complete hyperbolic metric of finite volume. If a hyperbolic Dehn filling $N$ on $M$ satisfies\/{\rm :}
\begin{enumerate}
\item[(i)] $N$ is infinitesimally projectively rigid,
\item[(ii)] the Dehn filling slope of N is contained in the {\rm (}connected\/{\rm )} hyperbolic Dehn filling space of M ,
\end{enumerate}
then infinitely many Dehn fillings on $M$ are infinitesimally projectively rigid.
\end{theorem}

The complete hyperbolic manifold $M$ that is the complement of a figure-eight knot in $\SI^3$ is infinitesimally projectively rigid. 
Then infinitely many Dehn fillings on $M$ are infinitesimally projectively
rigid. 

They showed the following: 
\begin{itemize} 
\item For a sufficiently large positive integer $k$, the homology sphere obtained by {$\frac{1}{k}$}–Dehn filling on the figure-eight knot is infinitesimally not projectively rigid.
Since the Fibonacci manifold $M_k$ is a branched cover 
of $\SI^3$ over the figure eight knot complements, 
for any $k \in N$, the Fibonacci manifold $M_k$ is not projectively rigid. 
\item All but finitely many punctured torus bundles with tunnel number one are infinitesimally projectively rigid relative to the cusps.
All but finitely many twist-knot complements are infinitesimally projectively rigid relative to the cusps.
\end{itemize}

\subsection{Ballas's work on ends.} 
The following are from Ballas \cite{Ballas14} and \cite{Ballas15}. 
\begin{itemize}
\item Let $M$ be the complement in $\SI^3$ of $4_1$ (the figure-eight knot), $5_2$, $6_1$, or $5^2_1$ (the Whitehead link). Then $M$ does not admit strictly convex deformations of its complete hyperbolic structure.
\item Let $M$ be the complement of a hyperbolic amphichiral knot. sSppose that $M$ is infinitesimally projectively rigid relative to the boundary, and the longitude is a rigid slope. Then for sufficiently large $n$, there is a one-dimensional family of strictly convex deformations of the complete hyperbolic structure on $M(m/0)$ for $m \in \bZ$. 
\item Let $M$ be the complement in $\SI^3$ of the figure-eight knot. There exists $\epsilon$ such that for each $s \in (-\epsilon, \epsilon)$, $\rho_s$ is the holonomy of a finite volume properly convex projective structure on $M$
for a parameter $\rho_s$ of representations $\pi_1(M) \ra \PGL(4, \bR)$.  
Furthermore, when $s \ne 0$, this structure is not strictly convex.
\end{itemize}
We also note the excellent 
work of Ballas, Danciger, and Lee \cite{BDL} experimenting 
with more of these and finding a method to glue along tori for 
deformed hyperbolic $3$-manifolds to produce convex real projective $3$-manifolds that 
does not admit hyperbolic structures.

\subsection{Finite volume strictly convex real projective orbifolds with ends} 
 We summarize the main results from two independent groups. 
The Hilbert metric is a complete Finsler metric on a properly convex set $\Omega$. This is the hyperbolic metric in the Klein model when $\Omega$ is projectively diffeomorphic to a standard ball. A simplex with its Hilbert metric is isometric to a normed vector space, and appears in a natural geometry on the Lie algebra $\mathfrak{sl}(n, \bR)$. A singular version of this metric arises in the study of certain limits of projective structures. The Hilbert metric has a Hausdorff measure and hence a notion of finite volume.
(See \cite{CLT15}.)

\begin{theorem}[Choi \cite{marg}, Cooper-Long-Tillmann \cite{CLT15}, Crampon-Marquis \cite{CM14}]
For each dimension $n \geq 2$, there is a Margulis constant $\mu_n > 0$ with the following 
property. If $M$ is a properly convex projective $n$-manifold and $x$ is a point in $M$, then the subgroup of $\pi_1(M,x)$ generated by loops based at $x$ of length less than $\mu_n$ is virtually nilpotent.
In fact, there is a nilpotent subgroup of index bounded above by $m = m(n)$. Furthermore, if $M$ is strictly convex and has finite volume, this nilpotent subgroup is abelian. If $M$ is strictly convex and closed, this nilpotent subgroup is either trivial or infinite cyclic.
\end{theorem}

\begin{theorem}[Cooper-Long-Tillmann \cite{CLT15}, Crampon-Marquis \cite{CM14}]
Each end of a strictly convex projective manifold or orbifold of finite volume is 
horospherical.
\end{theorem}

\begin{theorem}[(Relatively hyperbolic). Cooper-Long-Tillmann \cite{CLT15}, Crampon-Marquis \cite{CM14}]
Suppose that  $M = \Omega/\Gamma$ is a properly convex manifold of finite volume which is the interior of a compact manifold $N$, and the holonomy of each component of $\partial N$ is topologically parabolic. Then the following are equivalent\/{\rm :}
\begin{itemize} 
\item[1] $\Omega$ is strictly convex,
\item[2] $\partial \Omega$ is $C^1$,
\item[3] $\pi_1(N)$ is hyperbolic relative to the subgroups of the boundary components.
\end{itemize}
\end{theorem}
Here, the definition of the term ``topologically parabolic'' is according to \cite{CLT15}. 
This is not a Lie group definition but rather a topological term. 
We have found a generalization Theorem \ref{rh-thm-relhyp} and its converse 
Theorem \ref{rh-thm-converse} in Chapter \ref{ch-rh}.

%
%


%
%



\part{The classification of radial and totally geodesic ends.}




The purpose of this part is to understand 
the structures of ends of real projective $n$-dimensional orbifolds for $n \geq 2$.
In particular, we consider 
the radial or totally geodesic ends. Hyperbolic manifolds with cusps 
and hyperideal ends are examples. 
For this, we will study the natural conditions on eigenvalues of holonomy 
representations of ends when these ends are comprehensive. 
This is the most technical part of the monograph, containing a large number of results useful in the other two parts. 

We begin our study of radial ends in Chapter \ref{ch-ce}.
We categorize radial ends into three classes: complete affine radial ends, properly convex ends, and convex ends that are neither properly convex nor
complete affine. We define the lens and horospherical conditions for these ends. 
We give some examples of these radial ends.


In Chapter \ref{ch-du}, we study the theory of affine actions.   
This is the major technical section of this part. 
We consider the case when there is a discrete affine action of a group $\Gamma$
 acting cocompactly on
a properly convex domain $\Omega$ in the boundary of the affine subspace
$\mathds{A}^n$ in $\RP^n$ or $\SI^n$.
We study the convex domain $U$ in an affine space 
$\mathds{A}^n$ whose closure meets with $\Bd \mathds{A}^n$ in $\Omega$. 
We can find a domain $U$ having an asymptotic hyperspace at each point of 
$\Bd \Omega$ if and only if $\Gamma$ satisfies the uniform middle-eigenvalue condition
with respect to $\Bd \mathds{A}^n$. 
To prov this, we study the flow on the affine bundle over the unit tangent 
space over $\Omega$ generalizing parts of the work by Goldman-Labourie-Margulis on complete flat Lorentz $3$-manifolds \cite{GLM09}. 
We end with showing that a T-end has a lens neighborhood if it satisfies the uniform middle-eigenvalue condition. 

In Chapter \ref{ch-pr}, we study the theory of properly convex R-ends. 
Tubular actions and the dual theory of affine actions are discussed. We show that distanced actions
and asymptotically nice actions are dual. We explain that the uniform middle-eigenvalue condition 
implies the existence of a distanced action. The main result here is the characterization of 
R-ends whose end holonomy groups
satisfy 
uniform middle-eigenvalue conditions. 
In other words, they are generalized lens-shaped R-ends. 
We also discuss some important properties of lens-shaped R-ends. 
Finally, we show that lens-shaped T-ends and lens-shaped R-ends are dual. 
We conclude by discussing the properties of T-ends 
derived from this duality. 

In Chapter \ref{ch-app}, we investigate the applications of the radial end theory
such as the stability condition. We discuss the expansion and shrinking of the end neighborhoods.
We demonstrate the openness of the lens condition in Theorem \ref{app-thm-qFuch}, which is one of the central results required in Part III.
 We also prove Theorem \ref{intro-thm-sSPC},
the strong irreducibility of strongly tame properly convex real projective orbifolds with
generalized lens-shaped $\cR$-ends or horospherical $\cR$- or $\cT$-ends.

In major technical Chapter \ref{ch-np}, we discuss the R-ends that are neither properly convex 
nor complete affine (NPNC). 
First, we show that the end holonomy group of an NPNC-end $E$ will have an exact sequence 
\[ 1 \ra N \ra h(\pi_1(\tilde E)) \longrightarrow N_K \ra 1\] 
where $N_K$ velong 
to the projective automorphism group $\Aut(K)$ of a properly convex compact set $K$, 
$N$ is the normal subgroup of elements mapping to the trivial automorphism of $K$,
and $N_K$ acts cocompactly on $K^o$. 
We show that $\Sigma_{\tilde E}$ is foliated by 
complete affine subspaces of dimension $\geq 1$. 
We explain that an NPNC-end satisfying the transverse weak middle-eigenvalue condition
(see Definition \ref{np-defn-weakmec})
for NPNC-ends is a quasi-joined R-end under some natural conditions.
A quasi-joined end is an end with an end neighborhood 
covered by the join of a properly convex action and a horoball action 
twisted by translations (see Definition \ref{np-defn-quasijoin}.) 
For virtually abelian groups, Ballas-Cooper-Leitner \cite{BCLp}, \cite{BCL22} had covered much of these materials but not in our generality.

We will also classify the complete affine ends in the final chapter \ref{ce-sub-horo} of this part.

\chapter[Theory of convex radial ends]{Introduction to the theory of convex radial ends}  \label{ch-ce}


In Section \ref{intro-sec-ends},
 we will discuss the convex radial ends of orbifolds, covering the most elementary aspects of the theory. 
For a properly convex real projective orbifold, 
the space of rays of  each R-end gives us a closed real projective orbifold 
of dimension $n-1$. 
The orbifold is convex. The universal cover can be a complete affine subspace (CA)
or a properly convex domain (PC) or a convex domain that is neither (NPNC) 
in Section \ref{ce-sec-ends}. 
We discuss objects associated with R-ends, 
and examples of ends; horospherical ones, totally geodesic ones, 
and bendings of ends to obtain more general examples of ends. 

In particular, in Section \ref{np-sub-unitnorm}, we concentrated on 
the ends whose holonomies have only unit norm eigenvalues. We show
that they have to be cuspidal. 
In Section \ref{ce-sub-affineorb}, we discuss the ends of affine orbifolds.  

In Section \ref{intro-exmp}, we discuss some examples.

\section{End structures} \label{intro-sec-ends} 


%
%

\subsection{End fundamental groups} \label{intro-sub-endf} 
%
%
%
Let $\orb$ be a strongly
 tame real projective orbifold with universal cover $\torb$ and projective covering map $p_{\torb}$. 
A compact smooth orbifold $\bar \orb$ whose interior is $\orb$ 
is called a compactification of $\orb$. 
We don't assume that $\bar \orb$ has any compatible projective structure. 
There might be more than one compactification. 
(There might be more than one diffeomorphism classes of compactifications 
as smooth orbifolds.)

Throughout this book, we assume: 
\begin{hyp}\label{intro-hyp-compactification} 
A strongly tame orbifold $\orb$ is equipped 
with a chosen {\em compactification} $\bar \orb$
which is a smooth orbifold with boundary. 
\end{hyp} 

When we say $\orb$, we actually mean $\orb$ with 
$\bar \orb$. Each boundary component of $\bar \orb$ 
is the {\em ideal boundary component} of $\orb$ and 
is an {\em end} of $\orb$. 
\index{end|textbf} 
\index{ideal boundary component|textbf} 
An {\em end neighborhood} $U$ of $\orb$ is an open set $U$ where 
$\Sigma_E \cup U$ forms a neighborhood of 
an ideal boundary component $\Sigma_E$ corresponding 
to an end $E$. 
\index{end neighborhood|textbf} 

Let $\hat \orb$ denote the universal cover of $\bar \orb$ with 
the covering map $\hat p_{\torb}$ where we assume that $\torb \subset \bar \orb$ and 
$\hat p_{\torb}$ restricts to $p_{\torb}$ on $\torb$. 
Let $\Gamma$ be the deck transformation group of $\hat \orb \ra \bar \orb$
which also restricts to the deck transformation group of 
$\torb \ra \orb$.

An end $E$ is {\em proper} if the end has the following property:  
\begin{itemize}
\item The end $E$ has an end neighborhood homeomorphic to a closed connected $(n-1)$-dimensional orbifold $B$ times a half-open interval $(0, 1)$.
\item  An end neighborhood completes to an orbifold $U'$ that contains
a suborbifold diffeomorphic to $B \times (0, 1]$ in the compactification 
orbifold $\bar \orb$ with $B\times \{1\}$ a boundary component of $\bar \orb$.  
We say that $U'$ is an end neighborhood of $E$ {\em compatible} with $\bar \orb$. 
This is the {\em compatibility condition} with
the compactification $\bar \orb$. 
 \index{end neighborhood!compatible|textbf} 
 \index{end!proper|textbf} 
\item We say that the subset of $U'$ corresponding to $B \times \{1\}$ is an 
{\em ideal boundary component}. We will denote it by $S_E$. 
\end{itemize} 
The completion is called a {\em compactified end neighborhood} of the end $E$. 
The boundary component $S_E$ is called the {\em ideal boundary}
{\em component} of \index{boundary!ideal|textbf}
the end. Such ideal boundary components may not be uniquely determined
as there are two projectively nonequivalent ways to add boundary components of 
elementary annuli (see Section 1.4 of \cite{cdcr2}). \index{elementary annulus} 
Two compactified end neighborhoods of an end are {\em equivalent} if the
end neighborhood contains a common end neighborhood
whose compactification embeds into the compactified 
end neighborhoods. 
(See Section \ref{op-sub-endstr} for more detail.)

In this paper, we will not study the ends which are not proper.
An orbifold is strongly tame if and only if every end is proper and there are finitely many ends. 
However, there might be some orbifolds that have some proper ends and some ends are not proer. 

Let $E$ be a proper end. 
Each end neighborhood $U$, diffeomorphic to $S_{E} \times (0, 1)$ for an $(n-1)$-orbifold $S_E$, 
of an end $E$ lifts to a connected open set 
$\tilde U$ in $\torb$. 
We choose $U$ and the diffeomorphism $f_U: U \ra S_E \times (0, 1)$
such that $S_{E}\times (0, 1]$ is also diffeomorphic to 
a tubular neighborhood of a boundary component of 
$\bar \orb$ corresponding to $U$. 
A subgroup $\bGamma_{\tilde U}$ of $\Gamma$ acts on $\tilde U$ where 
\[p_{\torb}^{-1}(U) = \bigcup_{g\in \pi_1(\orb)} g(\tilde U).\] 
Each component
$\tilde U$ is said to be a {\em proper pseudo-end neighborhood}. 
\begin{itemize} 
	\item A {\em super-exiting sequence} of sets $U_{1}, U_{2}, \cdots $ in $\torb$ is 
	a sequence such that for each compact subset $K$ of $\orb$
	there exists an integer $N$ satisfying $p_\orb^{-1}(K) \cap U_i = \emp$ for $i > N$.  
	\item A {\em pseudo-end neighborhood sequence} is a super-exiting sequence of proper pseudo-end neighborhoods 
	\[\{U_{i}|i=1, 2, 3, \dots\}, \hbox{ where } U_{i+1} \subset U_{i} \hbox{ for every } i.\]
	\item Two pseudo-end sequences $\{ U_{i}\}$ and $\{V_{j}\}$ are {\em compatible} if 
	for each $i$, there exists $J$ such that $ V_{j} \subset U_{i}$ for all $j > J$
	and conversely for each $j$, there exists $I$ such that $U_{i} \subset V_{j}$ for all $i > I$. 
	\item A compatibility class of a proper pseudo-end sequence is called 
	a {\em p-end} of $\torb$.
	Each of these corresponds to an end of $\orb$ under the universal covering map $p_{\orb}$.
	\index{pseudo-end|textbf} 
	\index{p-end|textbf} 
	\item For a pseudo-end $\tilde E$ of $\torb$, we denote by $\bGamma_{\tilde E}$ the subgroup $\bGamma_{\tilde U}$ where 
	$U$ and $\tilde U$ are as above. We call $\bGamma_{\tilde E}$ a {\em pseudo-end fundamental group}.
	We will also denote it by $\pi_{1}(\tilde E)$.  They are independent of the choice of $U$ 
up to natural canonical inclusion homomorphisms by following 
Proposition \ref{intro-prop-endf}.
	\index{pseudo-end fundamental group|textbf} 
	\index{p-end fundamental group|textbf} 
	\item A {\em pseudo-end neighborhood } $U$ of a pseudo-end $\tilde E$ is 
a $\bGamma_{\tilde E}$-invariant open set that contains
	a proper pseudo-end neighborhood of $\tilde E$. 
	\index{end!p-end neighborhood|textbf} 
	A proper pseudo-end neighborhood is an example. 
\end{itemize}
(From now on, we will replace ``pseudo-end'' with the abbreviation ``p-end''.)
\index{end fundamental group|textbf} 
\index{pseudo-end|textbf} 
\index{p-end|textbf} 

As a summary, the set of boundary components of $\bar \orb$ 
has a one-to-one correspondence with the set of 
p-ends of $\orb$. 

\begin{proposition} \label{intro-prop-endf} 
	Let $\tilde E$ be a p-end of a strongly tame orbifold $\orb$. 
	The p-end fundamental group $\bGamma_{\tilde E}$ of $\tilde E$ is independent of the choice of $U$. 
\end{proposition}
\begin{proof}
	Given end neighborhoods $U$ and $U'$ for an end $E$, 
	let $\tilde U$ and $\tilde U'$ be p-end neighborhoods of a p-end $\tilde E$
	that are components of $p^{-1}(U)$ and $p^{-1}(U')$ respectively. 
	Let $\tilde U''$ be the component of $p^{-1}(U'')$ that is a p-end neighborhood of $\tilde E$. 
	Then $\bGamma_{\tilde U''}$ injects into $\bGamma_{\tilde U}$
	since both are subgroups of $\Gamma$. 
	Any $\mathcal{G}$-path in $U$ in the sense of Bridson-Haefliger \cite{BH99} is homotopic 
	to a $\mathcal{G}$-path in $U''$ by a translation in the $I$-factor. 
	Thus, $\pi_{1}(U'') \ra \pi_{1}(U)$ is surjective.  
	Since $\tilde U$ is connected, every element $\gamma$ of 
	$\bGamma_{\tilde U}$ can be represented by a $\mathcal{G}$-path connecting 
	$x_{0}$ to $\gamma(x_{0})$.  (See Example 3.7 in Chapter III.$\mathcal{G}$ of \cite{BH99}.) 
	Thus, $\bGamma_{\tilde U}$ 
	is isomorphic to the image of $\pi_{1}(U) \ra \pi_{1}(\orb)$. 
	Since $\bGamma_{\tilde U''}$ is surjective to the image of $\pi_{1}(U'') \ra \pi_{1}(\orb)$, 
	it follows that $\bGamma_{\tilde U''}$ is isomorphic to both $\bGamma_{\tilde U}$
and $\Gamma_{\tilde U'}$. 
\end{proof}

\subsection{Totally geodesic ends} \label{intro-sub-totgeo} 
Let $E$ be  an end  of a real projective orbifold $\orb$ of dimension $n \geq 2$ 

Suppose that in addition, we have the following: 
\begin{itemize} 
\item  $\orb \cup B$ is an orbifold with boundary $B$. 
We assume that it has a compatible real projective structure extending that of $\orb$, and  that each point of the added boundary component has a neighborhood projectively diffeomorphic to 
the quotient orbifold of an open set $V$ in an affine half-space $P$ 
such that $V \cap \partial P \ne \emp$
by a projective action of a finite group. This implies that the developing map extends to 
the universal cover of the orbifold with $U'$ attached. 
\end{itemize} 
The equivalence class of compactified end neighborhoods is called a {\em totally geodesic end structure} ({\em T-end structure}) for an end $E$. 
\index{end!totally geodesic!completion} 
\index{end!structure!totally geodesic}
\index{end!totally geodesic|textbf} \index{end!type-T|textbf}

\index{end orbifold!totally geodesic} 
\index{end!ideal boundary component}

Now, $\rpn$ admits a Riemannian metric of constant curvature called the 
Fubini-Study metric. 
Recall that the universal cover $\torb$ of $\orb$ has 
a path-metric induced by $\dev: \torb \ra \rpn$,
where  we assume that $\dev$ extends smoothly to 
$\hat \orb$ as well. 
We can Cauchy complete $\torb$ of this path-metric.
This Cauchy completion is called the {\em Kuiper completion} of $\torb$. 
(See \cite{psconv}.)
\index{Kuiper completion}
Note that we may sometimes use a lift $\dev: \torb \ra \SI^n$ 
of the developing map and use the same notation. 

A {\em T-end}
is a proper end equipped with a T-end structure. 
A {\em T-p-end} is a p-end $\tilde E$ 
corresponding to a T-end $E$. 
Using the definition of the T-end structure, we can easily show that 
there is a totally geodesic $(n-1)$-dimensional domain $\tilde S_{\tilde E}$ 
in the Cauchy completion of $\torb$ in the closure of
a p-end neighborhood of $\tilde E$. Of course, 
$\tilde S_{\tilde E}$ covers $S_E$. We call $\tilde S_{\tilde E}$
the p-end ideal boundary component. 
We will identify it with a domain in a hyperspace in $\RP^n$ (resp. $\SI^n$)
when $\dev$ is a fixed map to $\RP^n$ (resp. $\SI^n$). 
\index{p-end!ideal boundary component} 

Note that $\orb$ does not need to be strongly tame for this definition. 

\begin{definition}\label{intro-defn-lens}
A {\em topological lens} is a properly convex domain $L$ in 
$\RP^n$ such that $\partial L$ is a union of two strictly convex open $(n-1)$-dimensional cells. \index{lens!topolocial|textbf}
$L$ is a {\em generalized lens} if $\partial L$ is a union of two open disks, one of which is strictly convex and smooth
and the other is allowed to be just a topological $(n-1)$-cell.
\index{lens!generalized|textbf}
$L$ is a {\em PL-lens} if $\partial L$ is a union of two open $(n-1)$-cells, 
which is a union of compact convex 
$(n-1)$-dimensional polyhedrons meeting one another in strictly convex dihedral angles. 
$L$ is a {\em lens} if $\partial L$ is a union of two smooth strictly convex open $(n-1)$-cell. 
\index{lens|textbf}
A {\em lens-orbifold} (or {\em lens}) 
is a compact quotient orbifold of a lens by a properly discontinuous action of a projective group $\Gamma$ acting on
each boundary component as well. 
\index{lens!orbifold|textbf} 
\end{definition}


\begin{description} 
\item[(Lens condition for T-ends)] {\em The ideal boundary component} is identified 
as a totally geodesic suborbifold in the interior of a lens-orbifold 
in the ambient real projective orbifold containing $\bar \orb$. 
\end{description}
If the lens condition is satisfied for a T-end, 
we call it the {\em lens-shaped T-end}. 
The intersection of a lens with $\orb$ is called a {\em lens end neighborhood} of the T-end. A corresponding T-p-end is said to be a {\em lens-shaped T-p-end}. 
\index{T-end!lens-shaped} 
\index{T-p-end!lens-shaped} 
\index{T-end!neighborhood!lens} 

In these cases, $\tilde S_{\tilde E}$ is a properly convex $(n-1)$-dimensional domain, 
and $S_E$ is a $(n-1)$-dimensional properly convex 
real projective orbifold. 
We furthermore assume that 
$\pi_1(\tilde E)$ acts properly and cocompactly on the lens. 
Then $L \cap \torb$ is said to be {\em lens p-end neighborhood} of $\tilde E$
or $\tilde S_{\tilde E}$. 
\index{lens!cocompactly acted} 
\index{T-p-end!neighborhood!lens} 

We remark that for each component $\partial_i L$ for $i=1,2$ of $L$, 
$\partial_i L/\Gamma$ is compact and both are homotopy equivalent 
up to a virtual manifold cover $L/\Gamma'$ of $L/\Gamma$
for some finite-index subgroup $\Gamma'$. 
Also, the ideal boundary component of $L/\Gamma'$ has the same homotopy type as 
$L/\Gamma'$ and is a compact manifold. (See Selberg's Theorem \ref{prelim-thm-vgood}.)

\subsubsection{p-end ideal boundary components}
We recall Section \ref{intro-sub-totgeo}.
Let $E$ be a T-end of a real projective orbifold $\orb$. 
Given a totally geodesic end of $\orb$ and an end neighborhood $U$ diffeomorphic to $S_{E} \times [0, 1)$
with an end-completion by a totally geodesic orbifold $S_E$, 
we take a component $U_1$ of $p^{-1}(U)$ and a convex domain $\tilde S_{\tilde E}$, the ideal boundary component, 
developing into a totally geodesic hypersurface under $\dev$.
Here, $\tilde E$ is the p-end corresponding to $E$ and $U_1$.
There exists a subgroup $\bGamma_{\tilde E}$ acting on $\tilde S_{\tilde E}$. 
Again, $S_{\tilde E}:= {\tilde S_{\tilde E}}/\bGamma_{\tilde E}$ is projectively diffeomorphic to the {\em end orbifold} to be denoted by $S_E$ or 
$S_{\tilde E}$.

\begin{itemize} 
	\item We call $\tilde S_{\tilde E}$ a  p-end ideal boundary component of $\torb$. \index{end!p-ideal boundary component} 
	\item We call $S_{\tilde E}$ an ideal boundary component of $\orb$. \index{end!ideal boundary component} 
\end{itemize} 
We may regard $\tilde S_{\tilde E}$ as a domain in a hyperspace in 
$\RP^n$ or $\SI^n$. 
	\index{stildeE@$\tilde S_{\tilde E}$}
\index{stildeE@$S_{\tilde E}$}
\index{sE@$S_{E}$}




\subsection{Radial ends} \label{intro-sub-R-ends}
A {\em segment} is a convex arc in a $1$-dimensional subspace of $\RP^{n}$ or $\SI^{n}$. 
We denote the closed segment 
by $\ovl{xy}$ if $x $ and $y$ are endpoints. It is uniquely determined by $x$ and $y$ if $x$ and $y$ are not antipodal in $\SI^n$. 
In the following, all the sets are required to be inside an affine subspace $\mathds{A}^n$ and 
their closures to be either in $\RP^n$ or $\SI^n$. 

Let $\torb$ denote the universal cover of $\orb$ with the developing map $\dev$. 
Suppose that 
an end  $E$ of a real projective orbifold satisfies the following: 
\begin{itemize}
\item The end has an end neighborhood $U$ foliated by properly embedded projective geodesics. \index{end!radial|textbf} \index{end!type-R|textbf}
\item Choose any map $f: \bR \times [0, 1] \ra \orb$ such that $f|\bR\times \{t\}$
 is a geodesic leaf of such a foliation of $U$ for each $t$.
Then $f$ lifts to $\tilde f:\bR \times [0,1] \ra \torb$ where $\dev \circ \tilde f| \bR \times \{t\}$ for each $t \in [0, 1]$, 
maps to a geodesic in $\rpn$ ending at a point of concurrency common for every $t$. 
\end{itemize} 
The foliation  is called a {\em radial foliation} and leaves {\em radial lines} of $E$. 
Two such radial foliations $\mathcal F_{1}$ and $\mathcal F_{2}$ 
of radial end neighborhoods of an end are {\em equivalent} if the restrictions of $\mathcal{F}_{1}$ and $\mathcal{F}_{2}$ 
in an end neighborhood agree with each other. 
A {\em radial end structure} is an equivalence class of radial foliations.

Remember that $\orb$ always comes with a smooth compact orbifold $\bar \orb$ with boundary.
We will fix a radial end structure for each end of $\orb$
coming from a smooth foliation whose leaves 
end transversely to the boundary component of $\bar \orb$. 
This is the {\em compatibility condition} of the R-end structure to 
$\bar \orb$. \index{compatibility condition for the R-end structure|textbf} 

 \index{end neighborhood!compatible|textbf} 

To explain further, 
 an end neighborhood $U$ is {\em compatible} with $\bar \orb$ 
 if $U$ has a diffeomorphism $f:U \ra \Sigma_E \times (0, 1)$ 
 where the image of  $\Sigma_E \times \{t\}$ is transverse 
 to the radial foliation for $t \in (0, 1)$, and $f$ 
 extends to diffeomorphism $U\cup \Sigma_E \ra \Sigma_E \times (0, 1]$. 
 Here, $\Sigma_E$ is called an end orbifold of the radial end $E$. 
 \index{end orbifold!radial} 
 \index{R-end orbifold}



An {\em R-end} is a proper end with a radial end structure. 
An {\em R-p-end} is a p-end with a p-end neighborhood covering a radial end neighborhood with 
induced foliation. Of course, we do not forget the foliation. 
Each lift of the radial foliation has a finite path-length induced by 
$\dev$.   
A {\em pseudo-end }({\em p-end}\/) {vertex} of a radial p-end neighborhood or a radial p-end is 
the common endpoint of concurrent lifts of leaves of the radial foliation,
which we obtain by a Cauchy completion along the leaves.
The point will be denoted by $\mbv_{\tilde E}$ for an R-p-end $\tilde E$. 
Note that $\dev$ always extends to the pseudo-end vertex.
The p-end vertex is defined independently of the choice of $\dev$. 
$\mbv_{\tilde E}$ 
also denotes the image point of $\RP^n$ (resp. $\SI^n$) under $\dev$ if $\dev$
is a fixed developing map to $\RP^n$ (resp. $\SI^n$). We will often fix $\dev$. 
 \index{end!structure!radial|textbf}
 \index{end!pseudo!vertex|textbf} 
 \index{p-end!vertex|textbf}
 \index{end!p-end vertex|textbf} 

(See Section \ref{op-sub-endstr} for more detail.)
Note that $\orb$ does not need to be strongly tame for this definition. 
\index{end!radial!completion} 

\begin{remark}{(End-compactification structures)}\label{intro-rem-fixcomp} 
	If we have a compactification $\bar \orb'$ of $\orb$ not 
	diffeomorphic to $\bar \orb$, and choose $\bar \orb'$ instead of 
	$\bar \orb$, all these discussions have to take 
	place with respect to $\bar \orb'$. 
	By the s-cobordism theorem of Mazur \cite{Mazur}, Barden \cite{Barden} and 
	Stallings as well as the existence theorem 11.1 of 
	Milnor \cite{Milnor}, there are tame manifolds with more than 
	one compactifications. 
	(This is due to Benoit Kloeckner in the mathematics overflow.
	See also Section \ref{op-sub-endstr}.) 
However, notice that the radial structure determines the diffeomorphism type of $\bar \orb$ 
since each flow line determines the unique boundary points, and the set of flow lines passing 
a codimension-one transversal ball determines the diffeomorphism type. 
However,  given a fixed $\bar \orb$, the radial structure is determined only up to isotopies by
the isotopy uniqueness part of the tubular neighborhood theorem, which also holds for orbifolds (See Section 4.5 of \cite{Hirsch} and Section 4.4 of \cite{Cbook}.)

	\end{remark}

Two rays $l$ and $m$ with some arclength-parametrizations in $\Omega$ are
{\em asymptotic} if 
\[d_{\Omega}(l(t), m(t)) < C \hbox{ for some constant } C \hbox{ for all } t \geq 0.\]
However, note that 
the property does not require that the two geodesics share an endpoint.  
(See Section 3.11.3 of \cite{DK18}).


\begin{lemma}[Benoist \cite{Benoist04}] 
	\label{intro-lem-decrease} 
Let $l$ be a line in a properly convex open domain $\Omega$ 
in $\RP^n$ {\em (}resp. $\SI^n${\em ),} $n \geq 2$, ending 
at $x \in \Bd \Omega$. Let $m$ be a line ending at $x$ also. 
Then for a parametrization $l(t)$ of $l$, there is a parametrization 
$m(t)$ of $m$ such that $d_{\Omega}(m(t), l(t)) < C$ for some constant $C$ 
independent of $t$. Furthermore, $m$ and $l$ are asymptotic rays. 
\end{lemma}  
\begin{proof} 
	We will prove this for $\SI^n$. 
	We choose a supporting hyperspace $P$ at $x$. 
	Then $P \cap \clo(\Omega)$ is a properly convex domain. 
	Let $Q$ be a codimension-one subspace of $P$ disjoint from $P\cap \clo(\Omega)$. 
	There is a parameter of hyperspaces $P_t$ passing $l(t)$ and 
	containing $Q$. We denote $m(t) = m \cap P_t$. 
	For convenience, we may suppose our interval is $[0, 1)$ and
	that $l(0)$ and $m(0)$ are the beginning points 
	of $l$ and $m$ respectively. 
	Let $J$ denote the $2$-dimensional subspace containing both $l$ and $m$. 
	Now, $m(t), l(t)$ lie on the line $P_t \cap J$. 

By the convexity of the $2$-dimensional domain $\Omega \cap J$, 
	the function $t \mapsto d_{\Omega}(l(t), m(t))$ is eventually 
	decreasing: 
	We can draw four segments from $x$ to $l(t), m(t)$ and 
	the endpoints of $\overline{m(t)l(t)} \cap \Omega$.  
	The segments to the endpoints always move outward, and 
	the the segments to $l(t), m(t)$ are constant.  
	Hence, $d_{\Omega}(l(t), m(t)) \leq C' d_{\Omega}(l(0), m(0))$
	for a constant $C'\geq 1$. 	
	(See also Section 3.2.6 and 3.2.7 of \cite{Benoist04}.
	For eventual decreasing property, we don't need the 
	$C^1$-boundary property of $\Omega$.) 

Finally, we may use the arclength parameterizations by taking 
discrete equidistantly placed elements in $l$ and $m$ respectively, and 
a triangle inequality argument: We can show that the parameterizations are related by 
constants. Then we increase the intervals. 
	\hfill \SSn  {\parfillskip0pt\par}
	\end{proof}




Following Lemma \ref{intro-lem-actionRend} gives us 
another characterization of R-end and the condition  
$R_x(\tilde U)= R_x(\Omega)$.

 \begin{lemma}\label{intro-lem-actionRend} 
 	Let $\Omega$ be a properly convex open domain in $\RP^n$ {\em (}resp. $\SI^n${\em ),} $n \geq 2$. 
 	Suppose that $\orb= \Omega/\Gamma$ is a noncompact strongly tame orbifold.
 	Let $U$ be a proper end neighborhood, and let $\tilde U$ be a connected 
 	open set in $\Omega$ 
 	covering $U$. Let $\Gamma_{\tilde U}$ denote the subgroup of 
 	$\Gamma$ acting on $\tilde U$. 
 	Suppose that $\tilde U$ is foliated by segments with 
 	a common endpoint $x$ in $\Bd \Omega$. 
 	Suppose that $\Gamma_{\tilde U}$ fixes $x$.  
 	Then the following hold{\rm :} 
 	\begin{itemize}
 	\item 
 	 $\Gamma_{\tilde U}$ acts properly on 
 	$R_x(\tilde U)$ if and only if 
 	every radial ray in $\tilde U$ ending at $x$ 
 	maps to a properly embedded arc in 
 	$U$. 
 	
 	\item 
 	If the above item holds, then $R_x(\Omega) = R_x(\tilde U)$, 
 	and $x$ is an R-p-end vertex of $\torb$. 
 	\end{itemize} 
 \end{lemma} 
\begin{proof}
	It is sufficient to prove  the statements for $\SI^n$. 
	The forward direction of the first item is clear:  If a leaf $l$ is not properly embedded, 
then there exists a sequence $g_i \in \Gamma_{\tilde U}$ such that the direction of 
$g_i(l)$ accumulates to a point of $R_x(\tilde U)$. 
The properness of the action of $\Gamma_{\tilde U}$ contradicts this. 
	
	For the converse, suppose that $g_i(K) \cap K \ne \emp$ for infinitely 
	many mutually distinct $g_i \in \Gamma_{\tilde U}$ for a compact set $K \subset R_x(\tilde U)$.
	Then there exists a sequence $p_i \in K$ such that $g_i(p_i) \ra p_\infty$,
	$p_\infty \in K$. We can choose a compact set $\hat K \subset \Omega$ 
	such that the ray $l_i$ ending at $x$ in the direction of $p_i$ has an endpoint 
	$\hat p_i \in \hat K_i$ for each $i$.

	
	We can choose $q_i$ on $l_i$ such that $g_i(q_i) \ra q_\infty$ for 
	$q_\infty$ in $\Omega$ since the direction of $g_i(l_i)$ is in $K$ 
	and their endpoints are uniformly bounded away from $x$ for all $i$. 
	
	By Lemma \ref{intro-lem-decrease}, we can choose a point $\bar q_i$ on $l_1$
	such that $d_{\Omega}(q_i, \bar q_i) < C$ for some uniform constant $C$. 
	Thus, $d_{\Omega}(g_i(q_i), g_i(\bar q_i)) < C$. 
	We can choose a subsequence such that $g_i(\bar q_i) \ra q'_\infty$ 
	for a point $q'_\infty$ in $\Omega$ with $d_{\Omega}(q_\infty, q'_\infty)
	\leq C$. This implies that $l_1$ is not properly embedded in $U$. 
	This is a contradiction. 
	
	For the second item, we have $R_x(\tilde U)\subset R_x(\Omega)$ clearly. 
	Let $l$ be a line from $x$ in $U$, and
	let $m$ be any line from $x$ in $\Omega$. 
	For a parametrization of $l$ by $[0, 1)$, 
	we obtain $d_{\Omega}(l(t), m(t)) < C$, $t \in [0, 1)$,
where $C$ is a uniform constant $C> 0$ and $m(t)$ is a parameterization of $m$
	by Lemma \ref{intro-lem-decrease}. 
	Since $l$ maps to a  properly embedded arc in $U$, and $\Bd_\orb U$ is compact, 
it follows that	$d_{\Omega}(l(t), \Bd \tilde U \cap \Omega) \ra \infty$ as 
	$t \ra 1$. This implies that $m(t)\in \tilde U$ for sufficiently large $t$. 
	Therefore, $m$ has a direction in $R_x(\tilde U)$.  
	Hence, we showed $R_x(\tilde U) = R_x(\Omega)$.  
	\hfill  \SnT {\parfillskip0pt\par}
	\end{proof}  


 Let $\Omega$ be a properly convex domain in $\RP^n$ (resp. in $\SI^n$)  such that 
 $\orb = \Omega/\Gamma$ for a discrete subgroup $\Gamma$ of automorphisms 
 of $\Omega$. 
 The space of radial lines in an R-end lifts to a space $R_x(\Omega)$ consisting
 of lines in $\Omega$ ending at a point $x$ of $\Bd \Omega$.
 By Lemma \ref{intro-lem-actionRend}, $\Gamma_x$ acts properly on $R_x(\Omega)$
since we assume that we have radial ends only.  
 The quotient space  $R_x(\Omega)/\Gamma_x$ has an $(n-1)$-orbifold structure
 by Lemma \ref{intro-lem-actionRend}. 
\index{rx@$R_x(\cdot)$}
The {\em end orbifold} $\Sigma_E$ associated with an R-end is defined as the space of radial lines in $\orb$. It is clear that $\Sigma_E$ can be identified with 
 $R_x(\Omega)/\Gamma_x$. 
 By the compatibility condition stated at the beginning of Section \ref{intro-sub-R-ends},
$\Sigma_E$ is diffeomorphic to 
 the component of $\bar \orb - \orb$ corresponding to $E$. 
The space of radial lines in an R-end has the local structure of $\RP^{n-1}$ since we can lift a local 
neighborhood to $\torb$, and these radial lines lift to lines developing into concurrent lines.  
The end orbifold has an induced real projective structure of one lower dimension.

For the following, we may assume that all subsets here are bounded subsets of 
an affine subspace $\mathds{A}^n$ in $\RP^n$ (resp. in $\SI^n$).
\begin{itemize}
	\item An $n$-dimensional submanifold $L$ of $\mathds{A}^n$ is said to be a {\em pre-horoball} if it is strictly convex, 
	and the boundary $\partial L$ is diffeomorphic to $\bR^{n-1}$ 
	and $\Bd L - \partial L$ is a single point. 
	The boundary $\partial L$ is said to be a {\em pre-horosphere}. 
	\index{pre-horoball|textbf} \index{pre-horosphere|textbf}
	\item Recall that an $n$-dimensional subdomain $L$ of $\mathds{A}^{n}$ is called a  lens if $L$ is a convex domain 
	and $\partial L$ is a disjoint union of two smooth strictly convex embedded open
	$(n-1)$-cells $\partial_+ L$ and $\partial_{-} L$. 
It is a generalized lens if $\partial_+ L$ is allowed to be topological cell. 
	\index{lens} 
\index{lens-shaped!strict|textbf}
	A topological lens $L$ is said to be a {\em strict} one if the following hold: 
\begin{multline} 
	\partial \clo(\partial_+ L) = \clo(\partial_+ L) - \partial_+ L  = \partial \clo(\partial_- L)= \clo(\partial L_-) - \partial_- L ,\, \\
	\partial_+ L \cup \partial L_- = \partial L, \hbox{ and } \clo(\partial_+ L) 
\cup \clo(\partial_- L) = \partial \clo(L).
\end{multline} 
		\index{lens!strict|textbf}
			\index{lens!generalized!strict|textbf}
Any subdomain in $\torb$ mapping diffeomorphic to these under the developing map
will be named the same. 
			
	\item A {\em cone} is a bounded domain $D$ in an affine patch with a point in the boundary, called an {\em end vertex} $v$
	such that every other point $x \in D$ has an open segment $\ovl{vx}^{o} \subset D$.
A {\em trivial one-dimensional cone} is an open half-space in $\bR^1$ given by $x > 0$ or $x < 0$. 
\index{trival one-dimensional cone} 
	\index{cone|textbf} \index{end vertex} 
	A cone $D$ is a {\em join} $\{v\} \ast A$ for a subset $A$ of $D$ 
	if $D$ is a union of segments 
	starting from $v$ and ending at $A$. (See Definition \ref{prelim-defn-join}.)
	\item The cone $\{p\} \ast L$ over a lens-shaped domain $L$ in $\mathds{A}^n$, 
	$p\not\in \clo(L)$ 
	is called a {\em lens cone} if it is a convex domain and satisfies
	\begin{itemize}
		\item $\{p\}\ast L = \{p\} \ast \partial_+ L$ for one boundary component $\partial_+ L$ of $L$,  and 
		\item every segment from $p$ to $\partial_{+} L$ meets the other boundary component
		$\partial_{-} L$ of $L$ at a unique point.
	\end{itemize} 
	\item As a consequence, each line segment from $p$ to $\partial_{+} L$ is transverse to $\partial_{+} L$.
	$L$ is called the {\em lens} of the lens cone. 
	(Here, different lenses may give the identical lens cone.)
	Also, $\{p\}\ast L - \{p\}$ is a manifold with boundary $\partial_{+} L$.  
	\index{lens!cone|textbf} 
	\item Each of  two boundary components of $L$ is called a {\em top} or {\em bottom} hypersurface, 
	depending on whether it is further away from $p$. 
The top component is denoted by $\partial_+ L$
	and the bottom one by $\partial_{-}L$. 
	\index{top component|textbf} \index{bottom component|textbf} 
	\item An $n$-dimensional subdomain $L$ of $\mathds{A}^{n}$ is a {\em generalized lens} if $L$ is a convex domain 
	and $\partial L$ is a disjoint union of a strictly convex smoothly embedded open
	$(n-1)$-cell $\partial_- L$ and an embedded open $(n-1)$-cell $\partial_{+} L$, which is not necessarily smooth. 
	\index{lens!generalized|textbf} 
	\item 
	A cone  $\{p\} \ast L$ is said to be a {\em generalized lens cone}  if 
	\begin{itemize} 
		\item $\{p\} \ast L  = \{p\} \ast \partial_+ L, p \not\in \clo(L)$ 
		is a convex domain for a generalized lens $L$, and 
		\item every segment from $p$ to $\partial_{+} L$ meets $\partial_{-} L$ at a unique point.
	\end{itemize} 
	A lens cone will, of course, be considered a generalized lens cone. 
	\index{lens!cone!generalized|textbf}
	\item  We again define the top hypersurface  and the bottom one as above.
	They are denoted by $\partial_{+} L$ and $\partial_{-} L$ respectively. 
	$\partial_+ L$ can be non-smooth; however, $\partial_-L$ is required to be smooth. 
	\item A {\em totally-geodesic submanifold} is a convex domain in a subspace. 
	A {\em cone-over} a totally-geodesic submanifold $D$ is a union of all segments 
	with  one endpoint $x$ not in the subspace spanned by $D$ 
	and the other endpoint in $D$. We denote it by $\{x\} \ast D$. 
	\index{cone over a totally-geodesic submanifold} 
\end{itemize}


Recall that we name the image of the ideal vertex $\mbv_{\tilde E}$ 
by the same name as we fix the developing map $\dev$ from the beginning of 
Section \ref{intro-sub-R-ends}. 
We apply these to ends: 

\begin{definition} \label{intro-defn-R-ends} $ $
\begin{description}
	\item [Pre-horospherical R-end] An R-p-end $\tilde E$ of $\torb$ is {\em pre-horospherical} if it has a pre-horoball in $\torb$ as a p-end neighborhood, or equivalently, an open p-end neighborhood $U$ in $\torb$ 
	such that $\Bd U \cap \torb = \Bd U - \{\mbv_{\tilde E}\}$. 
	$\tilde E$ is {\em pre-horospherical} if it has a pre-horoball in $\torb$ as a p-end neighborhood. 
	We require that the radial foliation of $\tilde E$ with leaves ending 
at the $p$-end vertex $\mbv_{\tilde E}$ is the one where 
	each leaf ends at $\mbv_{\tilde E}$ under the developing map.  
	\item[Lens-shaped R-end]  An R-p-end $\tilde E$ is {\em lens-shaped} (resp.
	{\em generalized-lens-shaped}\,),  if it has a p-end neighborhood that is 
	projectively diffeomorphic to $L \ast \{\mbv_{\tilde E}\} -\{\mbv_{\tilde E}\}$ 
(resp. the interior of the same set\,), under $\dev$ where 
	\begin{itemize} 
	\item $L$ is a strict lens (resp. strict generalized lens)
	and 
	\item $h(\pi_1(\tilde E))$ acts properly and cocompactly on $L$, 
	\end{itemize} 
and every leaf of the radial foliation of the p-end neighborhood with ending 
at the p-end vertex 
maps to a radial segment ending at $v$ under the developing map.
A p-end end neighborhood of $\tilde E$ is ({\em generalized lens-shaped})
 {\em lens-shaped} if it has 
the interior projectively diffeomorphic to the interior of $L \ast \{\mbv_{\tilde E}\} -\{\mbv_{\tilde E}\}$.
\end{description}
\index{lens} 
\index{lens!cocompactly acted} 
\index{lens!generalized!cocompactly acted} 

An R-end of $\orb$ is {\em lens-shaped}  (resp. {\em totally geodesic cone-shaped},
{\em generalized lens-shaped}\,) if 
the corresponding R-p-end is lens-shaped  (resp. {\em totally geodesic cone-shaped},
{\em generalized lens-shaped}\,).  
An end neighborhood of an end $E$ is  ({\em generalized}) {\em lens-shaped} 
if so is a corresponding p-end neighborhood $\tilde E$. 
\index{R-p-end!lens-shaped|textbf} 
\index{R-p-end!generalized lens-shaped|textbf}  
\index{R-p-end!totally geodesic cone-shaped|textbf} 
\end{definition} 
\index{R-end!pre-horospherical|textbf} 
\index{R-end!neighborhood!generalized lens-shaped|textbf} 
\index{R-end!neighborhood!lens-shaped|textbf} 
\index{R-end!generalized lens-shaped|textbf} 
\index{R-end!lens-shaped|textbf} 
\index{end!generalized lens-shaped|textbf} 
An end neighborhood is {\em lens-shaped} if 
it is a lens-shaped R-end neighborhood or T-end neighborhood. 
A p-end neighborhood is {\em lens-shaped} if 
it is a lens-shaped R-p-end neighborhood or T-p-end neighborhood. 
\index{end!lens-shaped} 
\index{end!neighborhood!lens-shaped} 
\index{p-end!neighborhood!lens-shaped} 
\index{p-end!lens-shaped|textbf} 
Of course, it is redundant to say that an R-end or a T-end satisfies the 
lens condition depending on its radial or totally geodesic end structure. 

Also, the domain on $\torb$ with the $p$-end vertex 
added corresponding to $L \ast \{v\} $ will be called 
lens cone. The domains in the lens cone corresponding to (resp. generalized) 
lens will be called ({\em generalized}) {\em lens}.  

\begin{definition}\label{intro-defn:endtype} 
	A {\em real projective orbifold with radial or totally geodesic ends} is a strongly tame 
	orbifold with a real projective structure, where each end is an R-end or a T-end with a given end structure.
	An end of a real projective orbifold is (resp. {\em generalized} ) {\em lens-shaped} or {\em pre-horospherical} if it 
	is a (resp. {\em generalized}) {\em lens-shaped} or {\em pre-horospherical} R-end, 
	or if it is a lens-shaped T-end. 
\end{definition}


\subsubsection{p-end vertices} \label{prelim-subsub-pendv}
Let $\orb$ be a real projective orbifold with the universal cover $\torb$. 
We fix a developing map $\dev$ in this subsection and identify with its image. 
Given a radial end of $\orb$ and an end neighborhood $U$ of a product form $E \times [0, 1)$ with a radial foliation, \index{p-end} 
we take a component $U_1$ of $p^{-1}(U)$ and the lift of the radial foliation. 
The developing images of leaves of the foliation end at a common point $x$ in $\RP^n$. 
\begin{itemize}
	\item Recall that a p-end vertex of $\torb$ is the ideal point of 
	leaves of $U_1$. (See Section \ref{intro-sub-R-ends}.)
	\item Let $\SI^{n-1}_{\mbv_{\tilde E}}$ denote the space of equivalence classes of rays from $\mbv_{\tilde E}$ diffeomorphic to an $(n-1)$-sphere \index{snminusonev@$\SI^{n-1}_{\mbv_{\tilde E}}$}
	where $\pi_1({\tilde E})$ acts as a group of projective automorphisms.  
	Here, $\pi_1({\tilde E})$ acts on $\mbv_{\tilde E}$ and sends the leaves to leaves in $U_1$. 
	\item Given a p-end $\tilde E$ corresponding to $\mbv_{\tilde E}$, 
	we define $ \tilde \Sigma_{\tilde E} := R_{\mbv_{\tilde E}}(\torb)$ the space of directions of developed leaves under $\dev$  oriented away from $\mbv_{\tilde E}$ into a p-end neighborhood of $\tilde{\mathcal{O}}$
	corresponding to $\tilde E$. 
	The space develops to $\SI^{n-1}_{x}$ by $\dev$ as an embedding 
	into a convex open domain. 
	\index{rvO@$R_{\mbv_{\tilde E}}(\torb)$} 
	\index{sigmaE@$\Sigma_{E}$}
	\index{sigmatildeE@$\tilde \Sigma_{\tilde E}$}
	 \index{sigmatildeE@$\Sigma_{\tilde E}$}
	
	\item Recall that $\tilde \Sigma_{\tilde E}/\bGamma_{\tilde E}$ is projectively diffeomorphic to the end orbifold  $\Sigma_E$ or by $\Sigma_{\tilde E}$. (See Lemma \ref{intro-lem-actionRend}.)
	\item We may use the lift of $\dev$ to $\SI^n$. 
	The endpoint $x'$ of the lift of radial lines will be identified with 
	the {\em p-end vertex} when the lift of $\dev$ is fixed. 
	Here, we can canonically identify
	$\SI^{n-1}_{x'}$ and $\SI^{n-1}_x$ and 
	the group actions of $\bGamma_{\tilde E}$ on them. 
\end{itemize} 

\subsection{R-Ends} \label{ce-sec-ends}

We classify $R$-ends into three classes: complete affine ends, properly convex ends, and
nonproperly convex and not complete affine ends. 
We also introduce T-ends.


Recall that an R-p-end $\tilde E$ is convex 
if $\tilde \Sigma_{\tilde E}$ is convex. 
Since $\tilde \Sigma_{\tilde E}$ is a convex open domain,
it is contractible by Proposition \ref{prelim-prop-classconv}, and
it always lifts to $\SI^n$ as an embedding. 
By Proposition \ref{prelim-prop-classconv}, 
a convex R-end is either
\begin{enumerate}
\item[(i):]  complete affine (CA), 
\item[(ii):]  properly convex (PC), or 
\item[(iii):] convex but neither properly convex nor complete affine (NPNC). 
We may drop ``convex but''. 
\end{enumerate} 
 \index{end!complete affine|textbf} 
\index{end!properly convex|textbf}
 \index{end!convex but neither properly convex nor complete affine|textbf} 
\index{end!NPNC|textbf} 
\index{end!neither properly convex nor complete affine|textbf} 

%
%
%

\subsection{Cusp ends}\label{intro-sub-horo} 

A {\em parabolic algebra} $\mathfrak{p}$ is an algebra in a semisimple Lie algebra $\mathfrak{g}$
whose complexification contains a maximal solvable subalgebra of $\mathfrak{g}$  (p. 279--288 of \cite{Varadarajan84}).
A {\em parabolic group} $P$ of a semisimple Lie group $G$ is the full normalizer of a parabolic subalgebra. \index{parabolic group|textbf}


An {\em ellipsoid} in $\RP^n=\bP(\bR^{n+1})$ (resp. in $\SI^n=\SI(\bR^{n+1})$)  \index{ellipsoid|textbf}
is the projection $C -\{O\}$ of the null cone 
\[C:=\{\bx \in \bR^{n+1}| B(\bx, \bx)=0\}\] 
for a nondegenerate symmetric
bilinear form $B: \bR^{n+1} \times \bR^{n+1} \ra \bR$ of 
signature $(1, n)$. 
Ellipsoids are always equivalent by projective automorphisms of $\rpn$.
An {\em ellipsoid ball} is the closed contractible domain in an affine subspace $\mathds{A}^n$ of 
$\RP^n$ (resp. $\SI^n$)  bounded by an ellipsoid contained in $\mathds{A}^n$. 
A {\em horoball} is an ellipsoid ball with a point $p$ of the boundary removed. \index{horoball|textbf}
An ellipsoid with a point $p$ on it removed is called a {\em horosphere}. The {\em vertex} of \index{horosphere|textbf}
the horosphere or the horoball is defined as $p$.

Let $U$ be a horoball with a vertex $p$ in the boundary of $B$. 
A real projective orbifold that is projectively diffeomorphic to an orbifold
$U/\Gamma_p$ for a discrete subgroup $\Gamma_p \subset \PO(1, n)$ fixing 
a point $p \in \Bd B$ is called a {\em horoball orbifold}. \index{horoball!orbifold|textbf}
A {\em cusp or horospherical end} is an end with an end neighborhood that is such an orbifold. 
A {\em cusp group} is a subgroup of a parabolic subgroup of an isomorphic copy of $\PO(1, n)$ in $\PGL(n+1, \bR)$
or in $\SO^+(1, n)$ in $\SL_{\pm}(n+1, \bR)$.  
A cusp group is a {{\em unipotent cusp-group }} if it is unipotent as well. 
\index{cusp group|textbf} 
\index{cusp group!unipotent|textbf} 
\index{end!cusp|textbf} 
\index{end!horospherical|textbf} 

By Corollary \ref{ce-cor-cusphor}, 
an end is pre-horospherical if and only if it is a cusp end.
We will use the term interchangeably
except in Chapter \ref{ce-sub-horo} where we will prove this fact.

\subsubsection{Lie group invariant p-end neighborhoods}
We will need the following lemma later.

A {\em p-end holonomy group} is the image of a p-end fundamental group
under the holonomy homomorphism. If its universal cover 
$\torb$ embeds to $\SI^n$ or $\RP^n$, then $h$ is injective, and 
hence p-end holonomy group is isomorphic to 
the p-end fundamental group. 
An {\em end holonomy group} is the image of an end fundamental group. 
\index{p-end!holonomy group|textbf} 
\index{end!holonomy group|textbf} 

We say that a subgroup $\Gamma$ of a group $G$ is {\em cocompact} if 
$G = \Gamma \cdot K$ for a compact subset $K$ of $G$. 
\index{subgroup!cocompact|textbf} 

\begin{lemma} \label{prelim-lem-endnhbd}
	Let $\orb$ be a convex real projective $n$-orbifold,
	and let $\torb$ be its universal cover in $\RP^n$ {\em (}resp. in $\SI^n${\em ).} 
	Let $U$ be an R-p-end neighborhood of a p-end $\tilde E$ in 
	$\torb$
	where a p-end holonomy group $\bGamma_{\tilde E}$ acts on. 
	Let $Q$ be a discrete subgroup of $\bGamma_{\tilde E}$. 
	Suppose that $G$ is a connected 
	Lie subgroup of $\PGL(n+1, \bR)$  
 {\em (}resp. $\SLpm${\em )} that 
virtually contains $Q$ such that 
	$G \cap Q$ is cocompact in $G$. 
	Assume that $G$ fixes the R-p-end vertex $\mbv_{\tilde E}$ 
	and acts on $\tilde \Sigma_{\tilde E}$. 
	Then $\bigcap_{g \in G} g(U)$ is a non-empty $G$-invariant 
	p-end neighborhood of $\tilde E$. 	
\end{lemma} 
\begin{proof} 
	We first assume that $\torb \subset \SI^n$
	and $Q\subset G$. 
	Let $F$ be a compact fundamental domain of $G$ under 
	$G\cap Q$. 
	It is sufficient to prove the case when 
	$U$ is a proper p-end neighborhood
	since for any open set $V$ containing $U$, 
	$\bigcap_{g\in G} g(V)$ contains a p-end neighborhood 
	 $\bigcap_{g \in G} g(U)$.
	Hence, we assume that $\Bd_{\torb} U/\bGamma_{\tilde E}$ is 
	a smooth compact hypersurface. Let $F_U$ denote the fundamental 
	domain of $\Bd_{\torb} U$. 

	Let $F$ be a compact fundamental domain of $G$ with respect to 
	$Q$. 
	Let $L$ be a compact subset of $\tilde \Sigma_{\tilde E}$
	and let $\hat L$ denote the union of all maximal open segments 
	with endpoints $\mbv_{\tilde E}$ and $\mbv_{\tilde E-}$ 
	in the direction of $L$. 

	We claim that $\bigcap_{g\in F} g(U) = \bigcap_{g\in G} g(U)$ 
	contains an open set in $\hat L$. 	
	We show this by proving that for any maximal $l$ in $\hat L$, 
	$\bigcap_{g\in F} g(U) \cap l$ has a lower bound on its 
	$\bdd$-length. This lower bound is uniform over $L$. 
	
	Suppose not. Then there exists a sequence $g_i\in F$ and 
	maximal segment $l_i$ in $\hat L$ such that 
	the sequence of $\bdd$-length of $g_i(U) \cap l_i$ 
	from $\mbv_{\tilde E}$ goes to $0$ as $i \ra \infty$.
	The endpoint of $g_i(U)\cap l_i$ equals $g_i(y_i)$ for $y_i \in \Bd_{\torb} U$.
	This implies that $\{g_i(y_i)\} \ra \mbv_{\tilde E}$. 
	
	Now, $y_i$ corresponds to a direction $u_i \in \tilde \Sigma_{\tilde E}$. 
	Since $F$ is compact, $u_i$ corresponds to 
	a point of a compact set $F^{-1}(L)$, which corresponds 
	to a compact set $\hat F_U$ of $\Bd_{\torb} U$ with directions in $F^{-1}(L)$.   
	Hence, $y_i \in \hat F_U$, a compact set. 
	Since $\mbv_{\tilde E}$ is a fixed point of $G$, and 
	$y_i \subset \hat F_U$ for a compact subset $\hat F_U$ of $\SI^n$ not containing 
	$\mbv_{\tilde E}$, this shows that $g_i$ form an unbounded 
	sequence in $\SLnp$. This is a contradiction to $g_i \in F$. 
	
	



	We have a nonempty set
	\[\hat U:=\bigcap_{g\in G} g(U) = \bigcap_{g\in F} g(U)\]
	containing an open subset of $U$. 
	$G$ acts on $\hat U$ clearly. We take the interior of $\hat U$. 
	If $G$ only virtually contains $\bGamma_{\tilde E}$, 
	we just need to add finitely many elements to the above arguments. 
\hfill	\SSn {\parfillskip0pt\par}
	\end{proof} 

\begin{lemma}\label{prelim-lem-cuspsegment}  
Suppose that $\orb$ is properly convex. Then 
 a p-end vertex of a horospherical p-end cannot be an endpoint of a segment in $\Bd \tilde{\mathcal{O}}$. 
\end{lemma}
\begin{proof} 
Suppose that $\Bd \torb$ contains a segment $s$ ending at the p-end vertex $\mbv_{\tilde E}$. 
Then $s$ is on an invariant hyperspace of $\bGamma_{\tilde E}$.
Now, conjugate $\bGamma_{\tilde E}$ into a cusp subgroup $P$ of $\SO(n,1)$ fixing 
$(1,-1,0,\dots, 0) \in \bR^{n+1}$ by using an element $\eta$ of $\SL_{\pm}(n+1,\bR)$. 
By simple computations using the matrix forms of $\bGamma_{\tilde E }$, 
we can find a sequence 
$\{g_i\}, g_i \in \eta\bGamma_{\tilde E}h^{-1} \subset P$ such that $\{g_i(\eta(s))\}$ geometrically converges to 
a great segment. Thus, for
 $\eta^{-1} g_{i} h\eta\in \bGamma_{\tilde E}$,
 the sequence $\{\eta^{-1}g_{i}\eta(s)\}$ geometrically 
converges to a great segment in $\clo(\torb)$. 
This contradicts the proper convexity of $\torb$. 
\hfill  \SnT {\parfillskip0pt\par}
\end{proof}

\subsection{Unit-norm eigenvalued actions on ends}\label{np-sub-unitnorm}

Here, we collect useful results on unit-normed actions resulting in
Proposition \ref{ce-prop-orthouni2} and Lemma \ref{ce-lem-unithoro}.

 \begin{lemma} \label{ce-lem-transitive} $ $
 	\begin{itemize} 
 		\item Suppose that a closed connected projective group $G$ acts properly and cocompactly on 	a convex domain $\Omega$ in $\SI^n$\/ {\rm (}resp. $\RP^n$\/{\rm ).} Then $G$ acts transitively on $\Omega$. 
 		\item Suppose that $\Gamma$ is a uniform lattice in a closed connected group $G$	acting on a convex domain $\Omega$ in $\SI^n$\/ {\rm (} resp. $\RP^n$ {\rm ).} Also suppose that $\Gamma$ acts properly and cocompactly on $\Omega$. Then $G$ acts transitively on $\Omega$. 
 	\end{itemize} 
 \end{lemma} 
 \begin{proof} 
The first part is Lemma 2.5 of \cite{Benoist03}.  	

 For the second item, we also claim that $G$ acts properly. 
 Let $\hat F$ be the fundamental domain of $G$ with $\Gamma$ action. 
 Let $x \in \Omega$. 
 Let $F'$ be the image $F(x):= \{g (x)| g \in F\}$ in $\Omega$. 
 Then $F'$ is a compact set. 
 Define
 \[\Gamma_{F'} := \{g \in \Gamma|  g(F(x)) \cap F(x) \ne \emp\}.\] 
 Then $\Gamma_{F'}$ is finite by the properness of the 
 action of $\Gamma$. 
 Sinc eeach element of $G$ is a product of an element $g'$ of $\Gamma$ 
 and $f \in F$, and $g' f(x) = x$, it follows that 
 \[g'F(x) \cap F(x) \ne \emp \hbox{ and } g' \in \Gamma_{F'}.\] 
 Hence, the stabilizer $G_{x} $ is a subset of $\Gamma_{F'}F$, and $G_{x}$ is compact.
 $G$ becomes a Riemannian isometry group with respect to a metric on $\Omega$
by Lemma 3.4.11 of \cite{Thurston80}. Hence, $G$ acts properly. 
 	The second part follows from the first part
 		since $G$ must act both properly and cocompactly.
\hfill 	\SnT {\parfillskip0pt\par}
 \end{proof} 

\begin{lemma} \label{ce-lem-transitive2} 
Suppose that a simply connected Lie group $G$ acts smoothly and isometrically 
on a simply connected manifold $M$ with a metric such that each stabilizer is trivial. 
Suppose that $\dim G = \dim M$ and $G$ is a closed subgroup of the isometry group of $M$. 
Then $G$ acts transitively on $M$.  
\end{lemma}
\begin{proof} 
The orbit should be open by the Invariance of domain theorem. 
Since $G$ is a closed subgroup of the isometry group, 
an orbit is a closed set also. Since $M$ is connected, the orbit should be $M$. 
\end{proof}

 \begin{proposition} \label{ce-prop-orthouni} 
 	Let $N$ be either a discrete group or 
an $(n-1)$-dimensional connected Lie group where all the elements 
 	have only eigenvalues of unit norms acting on 
 	a convex $(n-1)$-dimensional domain $\Omega$
 	in $\SI^n$ {\em (}resp. $\RP^n${\em )}
 	 projectively, properly, and cocompactly. 
 	Then  $\Omega$ is a complete affine space. 
Moreover, 
 	if $N$ is a connected Lie group, then 
$N$ is a simply connected orthopotent solvable group. 
 \end{proposition}
 \begin{proof} 
 	Again, we first prove for the case when $\Omega \subset \SI^n$. 
 	First assume that $N$ is discrete. 
 	By Theorem \ref{prelim-thm-orthopotent}, $N$ is an orthopotent Lie group. 
 	Theorem \ref{prelim-thm-sweepU} implies that $\Omega$ is a complete affine space. 

Now, consider the case when $N$ is a connected Lie group. 
	By Lemma \ref{ce-lem-transitive}, $N$ acts transitively on $\Omega$. 
$N$ has an $N$-invariant metric on $\Omega$ by the properness of the action.
Consider an orbit map $N \ra N(x)$ for $x \in \Omega$. 
If a stabilizer of a point $x$ of $\Omega$ contains a group of dimension 
$\geq 1$, then $\dim N > \dim \Omega$. 
The stabilizer is a finite group. 
Hence, $N$ covers $\Omega$ finitely. Since $\Omega$ is contractible, 
the orbit map is a diffeomorphism. Hence, $N$ is contractible. 

Since the stabilizer group must be a trivial group as we see here,
the action of $N$ on $\Omega$ is simply transitive. 

By Theorem \ref{prelim-thm-orthopotent}, $N$ is a solvable Lie group. 

 	Now, $\Omega$ cannot be properly convex: 
 	Otherwise, 
 	by Fait 1.5 of \cite{Benoist03}, 
 	$N$ either acts irreducibly on $\Omega$ or $\Omega$ is the interior of a join of 
 	domain $\Omega_1, \dots, \Omega_n$ 
 	where $N$ acts irreducibly on each $\Omega_i$. 
 	Since a solvable group never acts irreducibly unless the domain is 
 	$0$-dimensional by the Lie-Kolchin theorem, 
 	$\Omega$ is a simplex or a point. (See Theorem 17.6 of \cite{Humphreys}.)
 	Then $N$ must be diagonalizable, which is a contradiction 
 	to the unit-norm-eigenvalued property. 
 	If $n-1=0$, the conclusion is true.

 	Finally, suppose that $\Omega$ is neither properly convex nor complete affine. 
 	Then $\Omega$ is foliated by $i_0$-dimensional complete affine spaces
 	for $i_0 < n$. 
 	The space of affine leaves is a properly convex domain $K$ 
 	by discussions on R-ends in Section \ref{ce-sec-ends}.
 	Hence, $N$ acts on $K$. The subgroup $N_l$ of $N$ acting on $l$ 
is $i_0$-dimensional since  the $N$-action is simply transitive. 
 	Hence, $N/N_l$ acts on a properly convex set $K^o$ satisfying 
 	the premises. Again, this is a contradiction as in the above paragraph. 
 	
 	Hence, $\Omega$ is complete affine. 
 \hfill	\SSn {\parfillskip0pt\par}
 \end{proof}

Given a virtually solvable subgroup $G$ of an algebraic Lie group, 
a {\em syndetic hull} ${\mathcal{S}}(G)$ of $G$ is a virtually solvable Lie group 
with finitely many components such that $G$ is cocompact in it. 
(See Theorem 1.6 of 
Fried and Goldman \cite{FG83} and D. Witte \cite{Witte95}.)
\index{syndetic hull} \index{s@$\mathcal{S}(\cdot)$}

 
 \begin{lemma} \label{ce-lem-orthnil}
 	Let $N$ be a closed 
 	orthopotent Lie group in $\SL_{\pm}(n, \bR)$ 
 	acting on $\bR^n$ inducing a proper action on
 	 an $(n-1)$-dimensional affine space $\mathds{A}^{n-1}$ that is 
 	the upper half-space of $\bR^n$ quotient by 
 	 the scalar multiplications.
 	Suppose that $N$ acts cocompactly on $\mathds{A}^{n-1}$. 
 	Then  there is a connected group $N_u$
 	with the following properties{\em :} 
 	\begin{itemize} 
 		\item $N \cap N_u$ in cocompact in both $N$ and $N_u$.  
 	\item $N_u$ is homeomorphic to a cell, 
 	\item $N_u$ acts simply transitively on $\mathds{A}^{n-1}$,
 	\item $N_u$ is the unipotent subgroup of dimension $n-1$ in $\SL_{\pm}(n, \bR)$ 
 	of $N$ normalized by $N$. 
 	\end{itemize} 
 	\end{lemma} 
 \begin{proof} 
 Since $N$ is orthopotent, 
 there is a flag of 
 vector subspaces $\{0\} = V_0 \subset V_1 \subset \cdots \subset V_m = \bR^n$
 preserved by $N$ 
 where $N$ acts as an orthogonal group on $V_{i+1}/V_i$
 for each $i= 0, \dots, m-1$.
 Here, $\mathds{A}^{n-1}$ is parallel to 
 the vector subspace $V_{m-1}$ of dimension $n-1$. 
 (See Chapter 2 of Berger \cite{Berger}.)
 
 Hence, there is a homomorphism 
 $N \ra \bigoplus_{i=0}^{m-1} \Ort(V_{i+1}/V_i)$. 
 Let $N_u'$ denote its kernel. 
 Then $N_u'$ is a unipotent group that is cocompact in $N$. 
 
 We define $N_u$ to be the Zariski closure of $N_u'$ in 
 $\SL_{\pm}(n, \bR)$.  Now, $N_u$ is a unipotent Lie group, 
 and $N_u'$ is cocompact in $N_u$ by Malcev \cite{Malcev49}. 
 
 Since $N_u$ also acts on $\mathds{A}^{n-1}$, $N'_u$ is cocompact in $N_u$, 
 and $N'_u$ acts properly, it follows that 
 $N_u$ acts properly on $\mathds{A}^{n-1}$. 
 By Lemma \ref{ce-lem-transitive}, 
 $N_u$ acts transitively on $\mathds{A}^{n-1}$. 
 The action has the trivial stabilizer since $N_u$ is unipotent. 
 This implies $N_u$ is homeomorphic to $\mathds{A}^{n-1}$.
\end{proof}

The following is related to Section 9 of \cite{CLT15}, and we made uses of the work of 
\cite{CM14}. However, they often assume the finiteness of certain types of 
volumes. 

\begin{proposition} \label{ce-prop-orthouni2} 
	Let $U$ be a properly convex open domain in $\SI^n$ {\rm (}resp. $\RP^n$\/{\rm )}
	radially foliated by lines with a common end point $p \in \Bd U$
	with smooth $\Bd U -\{p\}$. 
	Let $N$ be an $(n-1)$-dimensional connected Lie group 
	with only unit norm eigenvalues 
	acting on $U$ and fixing $p$.
	Suppose that $N$ acts on 
	$R_p(U)$ properly and cocompactly. 
	Then $U$ is the interior of an ellipsoid, 
	and $N$ is a unipotent cusp group acting transitively and freely on
	$\Bd U -\{p\}$ 
	\end{proposition} 
\begin{proof}  
We first assume $U \subset \SI^n$. 
By Proposition \ref{ce-prop-orthouni},  $R_p(U)= \mathds{A}^{n-1}$ is complete affine. 
Thus, $N$ is a simply-connected solvable Lie group
where each element is unit-norm eigenvalued 
by Proposition \ref{ce-prop-orthouni}.  

By Lemma \ref{ce-lem-orthnil}, there is a unipotent group $N_u$ of dimension $n-1$ 
where both $N/N\cap N_u$ and $N_u/N \cap N_u$ are compact. 
Since $N_u$ is isomorphic to a unipotent subgroup, 
$N \cap N_u$ is a lattice both in $N$ and $N_u$.  

The compactness of $N_u/N\cap N_u$  implies that 
the tangent space to $N\cap N_u$ must be same as the tangent space of 
the identity at $N_u$ as we can see from the central series 
extensions by free abelian groups. 
It follows that $N\cap N_u = N_u$ and so $N_u \subset N$.
Since they have the identical dimension and are both connected, it follows that 
$N_u = N$. Hence, $N$ acts simply transitively on $\mathds{A}^{n-1}= R_p(U)$.

We now show that $U$ is the interior of an ellipsoid. 
We identify $p$ with $[1, 0, \dots, 0]$.
Let $W$ denote the hyperspace in $\SI^n$ 
containing $p$ sharply supporting $U$. Here, 
$W$ corresponds to a supporting hyperspace in $\SI^{n-1}_p$ 
of the set of directions of an open hemisphere 
 $R_p(U)$, and hence it is the unique supporting hyperspace at $p$ 
 and, thus, $N$-invariant. Also, 
$U_1:= W \cap \clo(U)$ is a properly convex subset of $W$. 


Suppose that $U_1$ contains more than one point $p$. 
Let $V$ be the smallest subspace containing $p$ and $U_{1}$. 
Again, the unipotent group $N$ acts on $V$. 
Now, $V$ is divided into disjoint open hemispheres of various dimensions where $N$ acts on:
By Theorem 3.5.3 of \cite{Varadarajan84}, $N$ preserves a full flag structure 
$V_0 \subset V_1 \subset \dots \subset V_k = V$.
We take components of the complement $V_i - V_{i-1}$. 
Let $H_V:=V - V_{k-1}$. 

Suppose that $\dim V = n-1$ for contradiction. 
Then $H_V \cap U_1$ is not empty; otherwise, $V$ would have smaller dimension. 
Let $h_V$ be the component of $H_V$ meeting $U_1$.
Since $N$ is unipotent, $h_V$ has an $N$-invariant metric by Theorem 3 of Fried \cite{Fried86}.

We claim that the orbit of the action of $N$ is of dimension $n-1$ and hence is locally transitive on $H_V$: 
If not, then a one-parameter subgroup $N'$ fixes a point of $h_V$.  
This group acts trivially on $h_V$ since the unipotent group contains a trivial orthogonal subgroup. 
Since $N'$ is not trivial, it acts as a group of nontrivial translations on the affine subspace $H^o$.
It follows that $N'(U)$ is not properly convex, which is absurd. 
Hence, an orbit of $N$ is open in $h_{V}$, and 
$N$ acts locally simply-transitively without fixed points.

Since $N$ is unipotent, $N$ has trivial stabilizer at each point of $h_V$. 
There is an $N$-invariant Riemannian metric on $h_V$. 
The orbit of $N$ in $h_V$ is closed  since $h_{V}$ has an $N$-invariant metric, 
and $N$ is closed in the isometry group of $h_V$. 
Thus, $N$ acts transitively on $h_V$ since $\dim N= \dim h_V$. 

Hence, the orbit $N(y)$ of $N$ for $y \in H_V \cap U_1$ contains a component of $H_V$. 
This contradicts the assumption that $\clo(U)$ is properly convex
(compare this with arguments in \cite{CM14}.)


Suppose that $1 \leq \dim V \leq n-2$. 
Let $J$ be a subspace of dimension $1$ greater than $\dim V$ and containing $V$ and  meeting $U$. Let $J_{\mathds{A}}$ denote the subspace of $\mathds{A}^{n-1}$ corresponding 
to the directions in $J$. 
Then $J_{\mathds{A}}$ is mapped to itself or 
disjoint subspaces under the action of $N$. 
Since $N$ acts on $\mathds{A}$ transitively, 
a nilpotent subgroup $N_J$ of $N$ acts on $J_{\mathds{A}}$ transitively. 
Hence, \[\dim N_J =  \dim J_{\mathds{A}} = \dim V,\] 
and we are in the situation immediately 
above. 
The orbit $N_J(y)$ for a limit point $y \in H_V$ contains a component of $V - V_{k-1}$
as above. Thus, $N_J(y)$ contains the same component, which is an affine subspace. 
As above, we have a contradiction to the proper convexity since the above
argument applies to $N_J$. 

Therefore, points of form $x_0 \in W \cap \Bd(\torb)  - \{p\}$ do not exist. 
Hence, for any sequence of elements $g_i \in \bGamma_{\tilde E}$, we have
$\{g_i(y)\} \ra p$. 
Hence, \[\Bd U = (\Bd U \cap \torb) \cup \{ p\}.\]  
Clearly, $\Bd U$ is homeomorphic to an $(n-1)$-sphere.

Since $U$ is radial, this means that $U$ is a pre-horospherical p-end neighborhood. 
(See Definition \ref{intro-defn-R-ends}.)
Since $N$ acts transitively on a complete affine space $R_p(U)$, 
and there is a $1$ to $1$ radial correspondence of $R_p(U)$ and $\Bd U -\{p\}$,
it acts so on $\Bd U - p$. 
Since $N$ is unipotent and acts transitively on 
$\Bd U -\{p\}$, 
Lemma 7.12 of \cite{CM14} shows that $U$ is bounded by an ellipsoid. 
Choose $x \in U$, then $N(x)\subset U$ is a horospherical p-end neighborhood 
also. Since $\Aut(U)$ is the group of hyperbolic isometries of $U$ 
with the Hilbert metric, it follows that $N$ is a unipotent cusp group. 
\hfill \SnT {\parfillskip0pt\par}
\end{proof}


The following answers a question that we discussed with J. Porti at the UAB, 
Barcelona, in 2013
whether there is a noncuspidal unipotent group acting as an end holonomy group
of an R-end. 
%
Note that this was also proved by Theorem 5.7 in Cooper-Long-Tillmann \cite{CLT15}
using the duality theory of ends. Here we do not need to use duality. 
The following is related to the result of D. Fried \cite{Fried86}.

 \begin{lemma}\label{ce-lem-unithoro}
 	Assume that $\orb$ is a properly 
 	convex real projective orbifold with an end $E$ and
 	the universal cover $\torb$ in $\SI^n$ {\em (}resp. $\RP^n${\em )}.
 	Suppose that $E$ is a convex end with a corresponding 
 	p-end $\tilde E$. 
 	Suppose that the eigenvalues of elements of  $\bGamma_{\tilde E}$ have unit norm only. 
 	Then $\bGamma_{\tilde E}$ is conjugate to a subgroup of a parabolic subgroup in $\SO(n, 1)$ {\em (}resp. $\PO(n, 1)${\em ),} and
 	a finite-index subgroup of $\bGamma_{\tilde E}$ is unipotent
 	and $\tilde E$ is horospherical, i.e., cuspidal. 
 \end{lemma}
 \begin{proof}
 	We will assume first $\torb \subset \SI^n$. 
  	By Theorem \ref{prelim-thm-orthopotent}, 
 	$\bGamma_{\tilde E}$ is virtually orthopotent. 
 	By Proposition \ref{ce-prop-orthouni}, $\tilde \Sigma_{\tilde E}$ is complete 
 	affine, and $\bGamma_{\tilde E}$ acts on it as an affine transformation group.  
 	By Theorem 3 in Fried \cite{Fried86}, $\bGamma_{\tilde E}$ is virtually unipotent. 
 	Since $\tilde \Sigma_{\tilde E}/\bGamma_{\tilde E}$ is a compact complete-affine manifold, 
 	a finite-index subgroup $F$ of $\bGamma_{\tilde E}$ is contained in a unipotent Lie subgroup
 	acting on $\tilde \Sigma_{\tilde E}$. 
 	Now, by Malcev \cite{Malcev49}, it follows that the same group is contained in 
 	a simply connected unipotent group $N$ acting on $\SI^n$ since $F$ is unipotent. 
 	The dimension of $N$ is $n-1 = \dim \tilde \Sigma_{\tilde E}$ by Theorem 3 of \cite{Fried86}. 
 	
 	Let $U$ be a component of the inverse image of a p-end neighborhood 
 	such that ${\mbv_{\tilde E}} \in \Bd U$. Assume that $U$ is a radial p-end neighborhood of $\mbv_{\tilde E}$. 
 	The group $N$ acts on a smaller open set covering a p-end neighborhood
 	by Lemma \ref{prelim-lem-endnhbd}. 
 	From now on, let $U$ be this open set. Consequently, $\Bd U \cap \torb$ is smooth.
 	
 	Now $N$ acts transitively and properly on $\tilde \Sigma_{\tilde E}$ 
 	by Lemma \ref{ce-lem-transitive} since $F$ acts properly on it and $N/F$ is compact. 
 	$N$ acts cocompactly on $\tilde \Sigma_{\tilde E}$ since so does $F$, $F \subset N$. 
 	
 	By Proposition \ref{ce-prop-orthouni2}, 
 	$N$ is a cusp group, and $U$ is a p-end neighborhood bounded by an ellipsoid. 
 	
 	Since $\bGamma_{\tilde E}$ has a finite extension of $N$ 
 	as the Zariski closure, the connected identity component $N$
 	is normalized by $\bGamma_{\tilde E}$. 
 	
 	Also, for an element $g \in \bGamma_{\tilde E} - F$, 
 	suppose $g(x)\in U$
 	for $x \in \Bd U- \{  \mbv_{\tilde E}\})$. 
 	Now,  $g(\Bd U - \{  \mbv_{\tilde E}\})$ 
 	is an $N$-orbit of $g(x)$. 
 	Hence, $g(\Bd U - \{  \mbv_{\tilde E}\})) \subset U$ since these are a $N$-orbit ellipsoid and the interior of one. 
 	Hence $g^n$ is not in $F$ for all $n$, which leads to a contradiction. 
 	Also, $g(x)$ cannot be outside $\clo(U)$ similarly.
 	Hence, $\bGamma_{\tilde E}$ acts on $U$. 
 	Now, it is clear that $\bGamma_{\tilde E}$ is a cusp group. 
 \hfill \SnT {\parfillskip0pt\par}
 \end{proof}

\subsection{The affine orbifolds and ends} \label{ce-sub-affineorb} 

Recall that a geodesic is {\em complete} in a direction
if the affine geodesic parameter is infinite in the direction. 

\begin{itemize} 
\item An affine orbifold has a {\em parallel end} if the corresponding end has an end neighborhood foliated by properly embedded affine geodesics \index{end!parallel}
parallel to one another in charts and if each leaf is complete in one direction. 
We assume that the affine geodesics are leaves assigned as above. 
\begin{itemize}
\item We obtain a smooth complete vector field $X_E$ in a neighborhood of $E$ for each end following the affine geodesics, \index{end!end vector field} 
which is affinely parallel in the flow; i.e., leaves have parallel tangent vectors. We call this 
an {\em end vector field}. 
\item We denote by $X_\mathcal{O}$ the vector field partially defined on $\mathcal{O}$ by taking the union of 
vector fields defined on some mutually disjoint neighborhoods of the ends using the partition of unity. 
\item The oriented direction of 
the parallel end is uniquely determined in the developing image of each p-end neighborhood of the 
universal cover of $\mathcal{O}$.
\item Finally, we put a fixed complete Riemannian metric on $\mathcal{O}$ such that for each end there is an open 
neighborhood where the metric is invariant under the flow generated by $X_\mathcal{O}$.
Note that such a Riemannian metric always exists. 
\end{itemize} 
\item An affine orbifold has a {\em totally geodesic end} $E$ if each end can be completed by a totally geodesic affine hypersurface. \index{end!totally geodesic} 
That is, there exists a neighborhood of the end $E$ diffeomorphic to $\Sigma_E \times [0, 1)$ for an $(n-1)$-orbifold $\Sigma_E$ that compactifies to an orbifold
diffeomorphic to $\Sigma_E \times [0, 1]$, and each point of 
$\Sigma_E \times \{1\}$ has a neighborhood affinely diffeomorphic 
to a neighborhood of a point $p$ of $\partial H$ for a half-space $H$ of an affine space.  
This implies that the corresponding p-end holonomy group $h(\pi_1(\tilde E))$ 
for a p-end $\tilde E$ going to $E$
acts on a hyperspace $P$ corresponding to $E \times \{1\}$. 
\end{itemize} 

Note that an affine structure has a unique compatible real projective structure
since the atlas of charts in the affine structure can be considered as one for
the real projective structures. 
These ends of affine orbifolds 
have radial end-structure or totally geodesic end-structures
when they are considered as real projective $(n+1)$-dimensional orbifolds.

An affine suspension of a horospherical orbifold is called a {\em suspended horoball orbifold}. \index{suspended horoball orbifold} 
An end of an affine orbifold with an end neighborhood affinely diffeomorphic to this is said 
to be of {\em suspended horoball type}. It also has a parallel end since the fixed point on the boundary of $\bR^n$ gives a unique direction. 

\index{affine suspension!horoball|textbf}

\begin{proposition} \label{prelim-prop-susp} 
Under the affine suspension construction, a real projective $n$-orbifold has radial, totally geodesic, or horospherical ends 
if and only if the affine $(n+1)$-orbifold affinely suspended from it 
has parallel,  totally geodesic, or suspended horospherical ends. 
\end{proposition} 

Again, an suspended affine $(n+1)$-orbifold has {\em type $\cR$- or $\cT$-ends} if 
the corresponding real projective $n$-orbifold has $\cR$- or $\cT$-ends in
correspondingly.

\section{Examples} \label{intro-exmp} 

We begin with two examples. We will elaborate these in this section, which we will fully justify later.

\begin{example}\label{intro-exmp-hyp}
	The interior of a finite-volume hyperbolic $n$-orbifold with rank $n-1$ horospherical ends and
	totally geodesic boundary forms an example of a noncompact strongly tame properly convex real projective orbifold 
	with radial or totally geodesic ends. 
	For horospherical ends, the corresponding end orbifolds have Euclidean structures. 
	(Also, we could allow hyperideal ends by attaching radial ends. See below)
\end{example} 


\begin{example} \label{intro-exmp-endv}
Suppose that the end orbifold of an R-end $E$ is a $2$-orbifold based on a sphere with
	three singularities of order $3$. Then a line of singularity is a leaf of a radial foliation. 
 End orbifolds of the Porti-Tillmann orbifold \cite{PTp} and the double of a tetrahedral reflection orbifold
	are examples. 
	A double orbifold of a cube with edges of orders $3$ only has eight such end orbifolds.
	(See Proposition 4.6 of \cite{End1} and their deformations are computed in \cite{CHL12}.
	Also, see Ryan Greene \cite{greene} for the theory.  These are explained again in
Section \ref{ex-sec-examplesII}.)

\end{example}


%
%
%




Recall the Klein model of hyperbolic geometry: It is a pair $(\bB, \Aut(\bB))$ where 
$\bB$ is the interior of an ellipsoid in $\rpn$ or $\SI^{n}$ and 
$\Aut(\bB)$ is the group of projective automorphisms of $\bB$. 
Now,  
$\bB$ has a Hilbert metric which in this case is two times  
the hyperbolic metric. 
Then $\Aut(\bB)$ is the group of isometries of $\bB$. (See Section \ref{prelim-sec-rps}.)
\index{ellipsoid} 

From hyperbolic manifolds, we obtain some examples of ends. 
Let $M$ be a complete hyperbolic manifold with cusps. 
$M$ is a quotient space of the interior $\bB$ of an ellipsoid in $\rpn$ or $\SI^n$
under the action of a discrete subgroup $\bGamma$ of $\Aut(\bB)$. 
Then some horoballs are p-end neighborhoods of the horospherical R-ends. 

We generalize Definition \ref{ex-dfn-hyperideal}. 
Suppose that a noncompact strongly tame convex real projective orbifold $M$ has totally geodesic embedded surfaces $S_1,.., S_m$ that are homotopic to the ends.
Let $M$ be covered by a properly convex domain
$\tilde M$ in an affine subspace of $\SI^n$. 
\begin{itemize}
\item We remove the outside of $S_j$s to obtain a properly convex 
real projective orbifold $M'$ with totally geodesic boundary.
Suppose that each $S_j$ can be considered a lens-shaped T-end. 
\item Each $S_i$ corresponds to a disjoint union of totally geodesic domains $\bigcup_{j \in J} \tilde S_{i, j}$ in $\tilde M$ 
for a collection $J$. For each $\tilde S_{i, j} \subset \tilde M$, a group 
$\Gamma_{i,j}$ acts on it where $\tilde S_{i, j}/\Gamma_{i, j}$ is 
a closed orbifold projectively diffeomorphic to $S_i$. 
\item Suppose that $\Gamma_{i, j}$ fixes a point $p_{i,j}$ outside $\tilde M$.  
\item Hence, we form the cone 
$M_{i, j} := \{p_{i, j}\} \ast \tilde S_{i, j}$. 
\item We obtain the quotient 
$M_{i, j}/\Gamma_{i, j} -\{p_{i, j}\}$ and identify $\tilde S_{i, j}/\Gamma_{i, j}$ with $S_{i,j}$ in 
$M'$ to obtain the examples of real projective manifolds with R-ends. 
\item $(\{p_{i, j}\}\ast \tilde S_{i, j})^{o}$ is an R-p-end neighborhood and the end is a totally geodesic R-end. 
\end{itemize}
The result is convex by Lemma \ref{rh-lem-commsupp}
since we can regard $S_{j}$ as an ideal boundary component of $M'$ 
and that of $M_{i, j}/\Gamma_{i, j} -\{p_{i, j}\}$. 
This orbifold is called the {\em hyperideal extension} of the convex real projective orbifold. 
When $M$ is hyperbolic, 
each $S_j$ is lens-shaped by Proposition \ref{ce-prop-lensend}. 
Hence, the hyperideal extensions of hyperbolic orbifolds are 
properly convex.  (See also  \cite{BaoBonahon}, \cite{Rousset}, and \cite{Schlenker}.) 
\index{end!hyperideal}
\index{hyperideal extension} 





We will fully generalize the following in Chapter \ref{ch-pr}. 
We remark that Proposition \ref{ce-prop-lensend} also follows from Lemma \ref{ex-lem-niceend}. However, we used more elementary 
results to prove it here. 

\begin{proposition}\label{ce-prop-lensend}
Suppose that $M$ is a convex real projective orbifold.
Let $\tilde E$ be an R-p-end of $M$. 
Suppose that 
\begin{itemize}
\item the p-end holonomy group of $\pi_{1}(\tilde E)$ 
\begin{itemize} 
\item is generated by the homotopy classes of finite order, 
\item is simple, or 
\item satisfies the unit middle-eigenvalue condition
\end{itemize} 
and 
\item $\tilde E$ has a $h(\pi_{1}(\tilde E))$-invariant $(n-1)$-dimensional totally geodesic 
properly convex open domain $D$ in a p-end neighborhood, and 
\item The closure of $D$ does not contain the p-end vertex. 
\end{itemize} 
Then the R-p-end $\tilde E$ is lens-shaped.
\end{proposition}
\begin{proof} 
Let $\tilde M$ be the universal cover of $M$ in $\SI^n$. 
$\tilde E$ is an R-p-end of $M$, and $\tilde E$ has a $\pi_{1}(\tilde E)$-invariant $(n-1)$-dimensional totally geodesic 
properly convex domain $D$. 
Since $D$ is properly convex, $h(\pi_1(\tilde E))$ acts properly discontinuously on it. 
Up to a finite cover, the quotient $D/h(\pi_1(\tilde E))$ is a closed manifold
by Theorem \ref{prelim-thm-vgood}. 
Since $D/\pi_1(\tilde E)$ is homotopy equivalent to 
$\tilde \Sigma_{\tilde E}/\pi_1(\tilde E)$ for the end vertex $\mbv_{\tilde E}$ 
up to a finite cover, 
$D$ cannot just project to a subspace of codimension greater than $1$. 
Hence, $D$ projects to an open domain.
By Theorem \ref{prelim-thm-Kobayashi}, 
$D$ projects onto $\tilde \Sigma_{\tilde E}$, and hence 
$D$ is transverse to radial lines from $\mbv_{\tilde E}$.

Under the first assumption, 
since the end holonomy group $\bGamma_{\tilde E}$
is generated by elements of finite order, 
the eigenvalues of the generators  corresponding to the p-end vertex  $\mbv_{\tilde E}$ equal $1$,
and hence every element of the end holonomy group has $1$ as the eigenvalue at
$\mbv_{\tilde E}$.

Now, assume that the end holonomy groups
fix the p-end vertices with eigenvalues equal to $1$. 

Then the p-end neighborhood $U$ can be chosen to be 
the open cone over the totally geodesic domain with vertex $\mbv_{\tilde E}$. 
Now, $U$ is projectively diffeomorphic  to the interior of a properly convex cone in 
an affine subspace $\mathds{A}^n$. 
The end holonomy group acts on 
$U$ as a discrete linear group of unit determinant. 
The theory of convex cones applies, and using the level sets of 
the Koszul-Vinberg function, we obtain 
a one-sided convex neighborhood $N$ in $U$ with smooth boundary 
(see Section 4.4 of Goldman \cite{goldmanbook}).
Let $F$ be a fundamental domain of $N$ with a compact closure in $\torb$. 

We obtain a one-sided neighborhood on the other side as follows: 
We take $R(N)$ by applying a reflection $R$ that fixes each point of the hyperspace containing $\tilde \Sigma$ and the p-end vertex.
Then we choose a diagonalizable transformation $\mathcal{D}$ fixing the p-end vertex, and 
every point of $\tilde \Sigma$ such that the image $\mathcal{D}\circ R(F)$ is in $\torb$. It follows that $\mathcal{D}\circ R(N) \subset \torb$ as well. 
Thus, $N \cup \mathcal{D}\circ R(N)$ is the lens we needed. 
The interior of the cone $ \{\mbv_{\tilde E}\}  \ast (N\cup \mathcal{D}\circ R(N))$ is the lens cone neighborhood for $\tilde E$. 
\hfill \SSn {\parfillskip0pt\par}
\end{proof}

A more specific example is below. 
Let $S_{3,3,3}$ denote the $2$-orbifold with base space homeomorphic to a $2$-sphere and 
three cone-points of order $3$. 
The $3$-orbifolds satisfying the following properties are 
the examples of Porti-Tillmann \cite{PTp} 
and the hyperbolic Coxeter $3$-orbifolds based on 
 ideal hyperbolic $3$-polytopes of dihedral angles $\pi/3$. 
(See Choi-Hodgson-Lee \cite{CHL12}.)

The following is a more specific version of Lemma \ref{ex-lem-niceend}.  
We provide a much more elementary proof that does not depend on the full theory of 
this monograph. 
\begin{proposition} \label{ce-prop-lensauto}
Let $\mathcal{O}$ be a convex real projective $3$-orbifold with R-ends where each end orbifold 
is diffeomorphic  to a sphere $S_{3,3,3}$ or a disk with three silvered edges and three corner-reflectors of orders $3, 3, 3$. 
Assume that the holonomy group of $\pi_{1}(\orb)$ is strongly irreducible. 
Then the orbifold has only lens-shaped R-ends or horospherical R-ends. 
\end{proposition}
\begin{proof}
Again, it is sufficient to prove this for the case $\torb \subset \SI^3$. 
Let $\tilde E$ be an R-p-end corresponding to an R-end whose end orbifold is diffeomorphic to $S_{3, 3, 3}$. 
It is sufficient to consider only $S_{3,3,3}$ since it double-covers the disk orbifold. 
Since $\bGamma_{\tilde E}$ is generated by finite order elements fixing a p-end vertex $\mbv_{\tilde E}$,
every holonomy element has the eigenvalue equal to $1$ at $\mbv_{\tilde E}$.
Take a finite-index free abelian group $A$ of rank two in $\bGamma_{\tilde E}$. 
Since $\Sigma_E$ is convex, 
a convex projective torus $T^2$ covers $\Sigma_E$ finitely.
Therefore, $\tilde \Sigma_{\tilde E}$ is projectively diffeomorphic either to 
\begin{itemize}
\item a complete affine subspace or 
\item the interior of a properly convex triangle or
\item a half-space 

\end{itemize} 
by the classification of 
convex tori by Nagano-Yagi \cite{NY} 
found in many places including \cite{goldmanbook}, \cite{Benoist94}, and 
Proposition \ref{prelim-prop-projconv}. 
Since there exists a holonomy automorphism of order $3$ fixing a point of 
$\tilde \Sigma_{\tilde E}$, 
it cannot be a quotient of a half-space with a distinguished foliation by lines.
Thus, the end orbifold admits a complete affine structure or is a quotient of a properly convex triangle. 

 Suppose that $\Sigma_{\tilde E}$ has a complete affine structure.
 Since $\lambda_{\mbv_{\tilde E}}(g) = 1$ for all $g \in \bGamma_{\tilde E}$, 
 the only possibility according to Theorem \ref{ce-thm-comphoro} is when $\bGamma_{\tilde E}$ is virtually nilpotent  and we have a horospherical p-end for $\tilde E$. 
 
Suppose that $\Sigma_{\tilde E}$ has a properly convex  open triangle $T'$ as its universal cover. 
Since the diagonalizable abelian group $A$ acts on $T'$ cocompactly, 
we can deduce that 
$A$ contains an element $g'$ with the largest norm eigenvalue $>1$, and the smallest norm eigenvalue $<1$ as a transformation 
in $\SL_\pm(3, \bR)$ the group of projective automorphisms at $\SI^2_{\mbv_{\tilde E}}$.
As an element of $\SL_{\pm}(4, \bR)$, we have $\lambda_{\mbv_{\tilde E}}(g') = 1$, and the product 
of the remaining eigenvalues is $1$, the corresponding largest and smallest eigenvalues are $> 1$ and $< 1$ respectively. 
Thus, an element of $\SL_{\pm}(4, \bR)$, $g'$ fixes $v_1$ and $v_2$ other than $\mbv_{\tilde E}$
in directions of vertices of $T'$.
Since $\bGamma_{\tilde E}$ contains an order-three element exchanging the vertices of $T'$, 
there are three fixed points of an element of $A$ different from $\mbv_{\tilde E}, \mbv_{\tilde E-}$. 
By commutativity, there is a properly convex compact triangle $T\subset \SI^{3}$ 
that has these three fixed points as vertices 
where $A$ acts on. Hence, $A$ is diagonalizable over the reals.

We can make any vertex of $T$ to be an attracting fixed point of an element of $A$. 
Each element $g \in \bGamma_{\tilde E}$ conjugates elements of $A$ to $A$. 
Therefore $g$ sends the attracting fixed points of elements of $A$ to those of elements of $A$. 
Hence $g(T) = T$ for all $g \in \bGamma_{\tilde E}$. 

Each point on the edge $E$ of $\clo(T)$ is an accumulation point of an orbit of $A$
by taking a sequence $g_{i}$ such that the sequence of the largest eigenvalue norms $\lambda_{1}(g_{i})$ 
and the sequence of second largest eigenvalue norms $\lambda_{2}(g_{i})$ are 
going to $+\infty$ while the sequence $\log|\lambda_{1}(g_i)/\lambda_{2}(g_i)|$ 
is bounded. Since $\lambda_{\mbv_{\tilde E}} = 1$, writing every vector 
as a linear combination of vectors in the direction of the four vectors, 
this follows. 
Hence, $\partial T \subset \Bd \torb$ and $T \subset \clo(\orb)$. 

If $T^{o} \cap \Bd \orb \ne \emp$, then $T \subset \Bd \orb$ by Lemma \ref{prelim-lem-simplexbd}. 
Then each segment from $\mbv_{\tilde E}$ ending in $\Bd \orb$ has 
the direction in $\clo(\Sigma_{\tilde E})= T'$. It must end at a point of $T$. 
Hence, $\torb = (T \ast \mbv_{\tilde E})^{o}$, an open tetrahedron $\sigma$. 
Since the holonomy group acts on it, we can take a finite-index group fixing each vertex of $\sigma$.
Thus, the holonomy group is virtually reducible.  This is a contradiction. 

Therefore, $T \subset \orb$ as $T \cap \Bd \orb = \emp$. 
We have a totally geodesic R-end, and by Proposition \ref{ce-prop-lensend}, the end is lens-shaped. 
(See also \cite{End2}.) 
\hfill \SSn {\parfillskip0pt\par}
\end{proof}

The following construction is called ``bending'' and was first investigated by Johnson and Millson \cite{JM87}. This construction gives us examples of R-ends that are not 
totally geodesic R-ends. 
See Ballas and Marquis \cite{BM20} for additional examples. 
\index{R-end!bending|textbf} 


\begin{example}[Bending] \label{ce-exmp-bending} 
Let $\orb$ have the usual assumptions. 
We will concentrate on an end and not take into consideration the rest of the orbifold. 
Certainly, the deformation given here may not extend to the rest.
(If the totally geodesic hypersurface exists on the orbifold, the bending does extend to the rest.)
\index{end!bending} 

Suppose that $\orb$ is an oriented hyperbolic manifold with a hyperideal end $E$. 
Then $E$ is a totally geodesic R-end with an R-p-end $\tilde E$.
Let  the associated orbifold 
$\Sigma_{E}$ for $E$ of $\orb$ be a closed $2$-orbifold and 
let $c$ be a  two-sided simple closed geodesic in $\Sigma_{E}$. 
Suppose that $E$ has an open end neighborhood $U$ in $\orb$ diffeomorphic to $\Sigma_{E} \times (0,1)$
with totally geodesic boundary $\Bd U \cap \orb$ diffeomorphic to $\Sigma_E$.
Let $\tilde U$ be a p-end neighborhood in $\torb$ corresponding to $\tilde E$
bounded by $\tilde \Sigma_{\tilde E}$ covering $\Sigma_E$.
Then $U$ has a radial foliation whose leaves lift to radial lines in $\tilde U$ from $\mbv_{\tilde E}$. 

Let $A$ be an annulus in $U$ diffeomorphic to $c \times (0, 1)$, foliated by leaves of the radial foliation of $U$.
Now a lift $\tilde c$ of $c$ is in an embedded disk $A'$ that covers $A$.
Let $g_c$ be the deck transformation corresponding to $\tilde c$ and $c$. 
Suppose that $g_{c}$ is orientation-preserving. 
Since $g_{c}$ is a hyperbolic isometry of the Klein model, 
the holonomy $g_c$ is conjugate 
to a diagonal matrix with entries $\lambda,  \lambda^{-1}, 1, 1, $ where $\lambda > 1$,
and the last $1$ corresponds to the vertex $\mbv_{\tilde E}$.  
We take an element $k_b$ of $\SLf$ of the following form in this coordinate system: 
\begin{equation}\label{ce-eqn-bendingm} 
\left(
\begin{array}{cccc}
1           &       0              & 0   & 0  \\
 0          &       1              & 0  & 0 \\ 
 0           &      0              & 1  & 0 \\ 
 0           &      0               & b & 1   
\end{array}
\right)
\end{equation}
where $b \in \bR$. 
$k_b$ commutes with $g_c$. 
Let us just work on the end $E$. 
We can ``bend'' $E$ by $k_b$: 

Since $k_{b}$ commutes with $g_{c}$, 
$k_{b}$ induces a diffeomorphism $\hat k_{b}$ of an open neighborhood of 
$A$ in $U$ to another one  of $A$ 
We find tubular neighborhoods $N_1$ of $A$ in $U$
and $N_{2}$ of $A$. 
We choose $N_{1}$ and $N_{2}$ such that they are
diffeomorphic by a projective map $\hat k_b$. 
Then we obtain two copies $A_1$ and $A_2$ of $A$ 
by completing $U - A$. 

Give orientations on $A$ and $U$. 
Let $N_{1,-}$ denote the inner component of $N_{1} - A$ and 
let $N_{2, +}$ denote the outer component of $N_{2} -A$. 


We take a disjoint union $(U - A) \sqcup N_1 \sqcup N_2$ and 
\begin{itemize} 
\item identify the projectively diffeomorphic copy of $N_{1,-}$ in $N_{1}$ with $N_{1, -}$ in $U- A $ by the identity map
and 
\item identify the projectively diffeomorphic  copy of $N_{2,+}$ in $N_2$ with $N_{2,+}$ in $U - A$ by the identity also.
\end{itemize} 
We glue back $N_1$ and $N_2$ by 
the real projective diffeomorphism $\hat k_b$ of $N_1$ to $N_2$. 
Then $N_{1} - (N_{1,-} \cup A)$ is identified with $N_{2,+}$, 
and $N_{2} - (N_{2, +} \cup A)$ is identified with $N_{1, -}$ by  $\hat k_b$. 
We obtain a new manifold. 

For sufficiently small $b$, we see that the end is still lens-shaped, 
and it is not a totally geodesic R-end. (This follows since the condition of being 
a lens-shaped R-end is an open condition. 
See  Section \ref{cl-sec-openness}.)




For the same $c$, 
let $k_s$ be given by 
\begin{equation}\label{ce-eqn-bendingm2} 
\left(
\begin{array}{cccc}
s           &       0              & 0   & 0  \\
 0          &       s              & 0  & 0 \\ 
 0           &      0              & s  & 0 \\ 
 0           &       0            & 0 & 1/s^3   
\end{array}
\right)
\end{equation}
where $s \in \bR_+$. 
These give us bendings of the second type. 
For $s$ sufficiently close to $1$, the property of being lens-shaped is preserved 
and as is the property of being g a totally geodesic R-end. 
(However, these will be understood by cohomology.)

If $s \lambda < 1$ for the maximal eigenvalue $\lambda$ of a closed curve $c_1$
meeting $c$ odd number of times, we have that the holonomy along $c_1$ has the attracting 
fixed point at $\mbv_{\tilde E}$. This implies that we no longer have lens-shaped R-ends even if we have 
started with a lens-shaped R-end. 

\end{example} 

There are some recent developments: 
See Martin Bobb \cite{Bobb} for examples of these bendings carried out on the whole orbifolds. He and James Farr \cite{BFp2024} have recently developed bending coordinates for convex real projective surfaces,  for our first type of bending here. 









\chapter[Properly convex affine actions]{The affine action on a properly convex domain 
	whose boundary is in the ideal boundary} \label{ch-du}


In this chapter, we discuss the asymptotic niceness of the affine actions
when the affine group $\Gamma$ 
acts on a convex domain $\Omega$ in $\mathds{A}^n$ and 
a properly convex domain in the ideal boundary of 
$\mathds{A}^n$. We find a properly convex domain 
in $\mathds{A}^n$ with boundary in $\Omega$ under some conditions.  
	The main tools are generalized Anosov flows on the affine bundles over 
	the unit tangent bundles as in Goldman-Labourie-Margulis \cite{GLM09}. 
	We introduce a flat bundle and decompose it in an Anosov-type manner. 
	Then we find an invariant section.  
	We  prove the asymptotic niceness using 
	the sections.   
In Section \ref{du-sec-affine}, 
we  define asymptotic niceness and flow decomposition of 
the vector bundles over $\Uu \Omega/\Gamma$. 
In Section \ref{du-sub-anosov},
we begin with a strictly convex domain $\Omega$
with a hyperbolic $\Gamma$ 
and the main result Theorem \ref{du-thm-asymnice}. 
We define 
proximal flows and decompose the vector bundle flows into contracting, repelling, and neutral subbundles. 
In Section \ref{du-sub-Anosov}, we show the contracting 
and expansion properties of 
the contracting and repelling subbundles, 
with a somewhat technical argument involving pulling-back. 
However, the neutral subbundles here are more of
a generalized type than what they had. 
We obtain the neutralized sections
as Goldman-Labourie-Margulis did. 
We will prove the main result for strictly convex $\Omega$ at the end of 
this section using the neutralized sections to obtain asymptotic hyperspaces. 
In Section \ref{du-sec-gen}, we will generalize these results to the case when 
$\Omega$ is not necessarily strictly convex. They are 
Theorems \ref{du-thm-asymniceII} and \ref{du-thm-ASunique}. 
A basic technique here is to enlarge the unit tangent bundle to 
an augmented unit tangent bundle by blowing 
up using the compact sets of hyperspaces at the endpoints of geodesics. 
The strategy to prove the second main result
is analogous to the strictly convex case for $\Omega$. 
In Section \ref{du-sec-lensT-end}, we discuss the lens condition for T-ends
obtained by the uniform middle-eigenvalue condition.
We will end by finding strictly convex smooth hypersurfaces approximating 
any convex boundary components for these types of domains. 
Except for Section \ref{du-sec-lensT-end}, we will work only in 
$\SI^n$ for simplicity. 

Note that there recently appeared related works by Nie, Seppi \cite{NieSeppi23}, and 
Ablondi \cite{Ablondi25} which lead to and confirm partly the main results in this chapter using the techinques of Cheng and Yau \cite{ChengYau} and Barbot \cite{Barbot2005}. 


\section{Affine actions} \label{du-sec-affine} 
Let $\Gamma$ be an affine group acting on the affine subspace $\mathds{A}^n$, which 
is an open hemisphere, with boundary $\Bd \mathds{A}^n = \SI^{n-1}_\infty $ in $\SI^n$. 
Let $U$ be a properly convex $\Gamma$-invariant domain 
in $\mathds{A}^n$ with the following property: 
\[\clo(U) \cap \Bd \mathds{A}^n = \clo(\Omega) \subset \Bd \mathds{A}^n\]
for a properly convex open domain $\Omega$.
We also assume that $\Omega/\Gamma$ is a closed orbifold. 
The action of $\Gamma$ on $\SI^n$ or $\RP^n$ is said to be a {\em properly convex affine action}.
Also, $(\Gamma, U, \Omega)$ is said to be a {\em properly convex affine triple}.
\index{properly convex affine action} 

A sharply supporting hyperspace $P$ at $x\in \partial \clo(\Omega)$ is 
{\em asymptotic} to $U$ if there is no other sharply supporting hyperspace $P'$ 
at $x$ such that $P'\cap \mathds{A}^n$ separates $U$ from $P \cap \mathds{A}^n$. 
In this case, we say that hyperspace $P$ is {\em asymptotic} to $U$. 
We will use the abbreviation {\em AS-hyperspace} to indicate 
for asymptotic sharply supporting hyperspace.
\index{supporting hyperspace!asymptotic|textbf} 
\index{AS-hyperspace|textbf} 

Let $(\Gamma, U, \Omega)$ be a properly convex affine triple. 
A properly convex affine action of $\Gamma$ is said to be {\em asymptotically nice}  \index{asymptotically nice} with respect to $U$ if 
$\Gamma$ acts on a compact subset
\[J:=\{ H| H \hbox{ is an AS-hyperspace in $\SI^n$ at } x \in \partial \clo(\Omega), H \not\subset {\SI^{n-1}_\infty}\}\] 
where we require that every sharply supporting $(n-2)$-dimensional space 
of $\Omega$ in ${\SI^{n-1}_\infty}$ is 
contained in at least one element of $J$. 
As a consequence, for any sharply supporting $(n-2)$-dimensional space $Q$
of $\Omega$, the set 
\[H_Q :=\{H \in J| H \supset Q\} \]
is compact and bounded away from $\Bd \mathds{A}^n$
in the Hausdorff metric $\bdd_H$.


\begin{definition} \label{du-defn-excon} 
A subspace $U$ of $\bR^n$ is {\em expanding} under a linear map $L$ 
if $\llrrV{L(u)} \geq C \llrrV{u}$ for every $u \in \bR^n$ for 
a fixed norm $\llrrV{\cdot}$ of $\bR^n$ and $C> 1$.  

A subspace $U$ of $\bR^n$ is {\em contracting} under a linear map $L$ 
if $\llrrV{L(u)} \leq C \llrrV{u}$ for every $u \in U$ for 
a fixed norm $\llrrV{\cdot}$ of $\bR^n$ and $0 < C < 1$. 
\end{definition} \index{map!expanding} \index{map!contracting} 
The expanding condition 
is equivalent to the condition that 
all the norms of eigenvalues of $L|U$ are
strictly larger than $1$. (See Corollary 1.2.3 of Katok and Hasselblatt
\cite{KH95}.) 

In this section, we will work with $\SI^n$ only, while
the $\RP^n$ versions of the results follows from the results here
in an obvious manner by Results in Section \ref{prelim-sub-lifting}
and then projecting back to $\RP^n$. 

For each element of $g \in \Gamma$, 
\begin{equation}\label{du-eqn-bendingm4} 
h(g)= 
\left(
\begin{array}{cc}
\frac{1}{\lambda_{{\tilde E}}(g)^{1/n}} \hat h(g)          &       \vec{b}_g     \\
\vec{0}          &     \lambda_{{\tilde E}}(g)                  
\end{array}
\right)
\end{equation}
where $\vec{b}_g$ is $n\times 1$-vector and $\hat h(g)$ is an $n\times n$-matrix of determinant $\pm 1$
and $\lambda_{{\tilde E}}(g) > 0$.  
\index{lambdaE@$\lambda_{{\tilde E}}(\cdot)$|textbf } 
In the affine coordinates, it is of the form 
\begin{equation}\label{du-eqn-affact} 
x \mapsto \frac{1}{\lambda_{{\tilde E}}(g)^{1+ \frac{1}{n}}} \hat h(g) x + \frac{1}{\lambda_{{\tilde E}}(g)} \vec{b}_g. 
\end{equation}
Let $\lambda_1(g), \dots, \lambda_n(g), \lambda_{\tilde E}(g)$ denote the norms of eigenvalues of $h(g)$ where 
$\lambda_1(g)$ denotes the maximal norm of the eigenvalues 
and $\lambda_n(g)$ the least norm of the eigenvalues
for $h(g)$, $g\in \Gamma$. 
If there exists a uniform constant $C > 1$ such that 
\begin{equation} \label{du-eqn-umecDii}
C^{-1} \leng_\Omega(g) \leq \log \frac{\lambda_1(g)}{\lambda_{{\tilde E}}(g)} \leq C \leng_\Omega(g), \quad
g \in \bGamma_{\tilde E} -\{\Idd\},
\end{equation} 
then we say that $\Gamma$  satisfy 
the {\em uniform middle-eigenvalue condition} with respect to 
the boundary hyperspace. 
\index{middle-eigenvalue condition!uniform|textbf} 
\index{uniform middle-eigenvalue condition|textbf} 

By taking $g^{-1}$ instead, we obtain the equivalent condition: 
\begin{equation} \label{du-eqn-umecDii-inverse}
C^{-1} \leng_\Omega(g) \leq \left|\log \frac{\lambda_n(g)}{\lambda_{{\tilde E}}(g)}\right| \leq C \leng_\Omega(g), \quad
g \in \bGamma_{\tilde E} -\{\Idd\},
\end{equation} 

This implies that 
\begin{equation} \label{du-eqn-umec}
\lambda_1(g)/\lambda_{\tilde E}(g) > 1 
\hbox{ and } \lambda_n(g)/\lambda_{\tilde E}(g) < 1
\end{equation} 
as we can see by taking the inverse of $g$.

Daryl Cooper also mentioned that 
\begin{equation}
\log\left(\frac{\lambda_1(g)}{\lambda_{\mbv_{\tilde E}}(g)}\right) \leq  \leng_K(g) 
\end{equation}
is true since 
$\leng_K(g) = \log\left(\frac{\lambda_1(g)}{\lambda_{n+1}(g)}\right)$. 
(See Prop. 2.1 of \cite{CLT15}.) 

We can replace $\leng_K(g)$ with $\cwl(g)$ for properly convex ends
by  Svarc-Milnor \cite{Svarc} and \cite{Milnor2} 
 (A. Swartz is currently a professor at UC. Davis) or Milnor-Svarc theorem. 
(See Theorem 8.1 of Farb-Margalit \cite{FM} also mentioning 
Efremovi\v{c} \cite{Efremovic}.)
\index{length@$\leng_K(g)$|textbf} 

We denote by $\mathcal{L}: \Aff(\mathds{A}^n) \ra \GL(n, \bR)$ 
the homomorphism $g \mapsto M_g$ taking the {\em linear part} of 
an affine transformation $g: x \mapsto M_gx + \vec{b}_g$ to 
$M_g \in \GL(n, \bR)$. 
\index{linear-part homomorphism|textbf} \index{l@$\mathcal{L}(\cdot)$|textbf}

\begin{lemma} \label{du-lem-umec} 
Let $\Gamma$ be a real projective 
group acting on the affine subspace $\mathds{A}^n$ with boundary $\Bd \mathds{A}^n$ in $\SI^n$ satisfying the uniform middle-eigenvalue 
condition with respect to $\Bd \mathds{A}^n$. Then  the linear part 
of $g \in \Gamma$ given by $\frac{1}{\lambda_{{\tilde E}}(g)^{1+\frac{1}{n}}} \hat h(g)$ 
has a nonzero expanding subspace and a contracting subspace
in $\bR^n$ provided $\leng_\Omega(g) > 0$. 
\hfill $\square$
\end{lemma} 

\begin{theorem}\label{du-thm-asymnice}
	We assume that $\Gamma$ is a word-hyperbolic group with a properly convex 
	affine action. 
	Let $\Gamma$ have an affine action on the affine subspace 
	$\mathds{A}^n \subset \SI^n$,
	acting on a properly convex domain $\Omega$ in 
	$ \Bd \mathds{A}^n$. 
	Suppose that $\Omega/\Gamma$ is a closed $(n-1)$-dimensional orbifold, and 
	that $\Gamma$ satisfies the uniform middle-eigenvalue condition. 
	Then $\Gamma$ acts on a properly convex 
	open domain $U$ with the following properties: 
	\begin{itemize}
		\item $(\Gamma, U, \Omega)$ is a properly convex triple, and 
 $\Gamma$ is asymptotically nice with the properly convex open domain $U$,
		and 
		\item For any open set $U'$ such that $(\Gamma, U', \Omega)$ is a properly convex triple,
		the  AS-hyperspace at each point of $\partial \clo(\Omega)$ exists and 
		coincides with that of $U$.  That is, $\Gamma$ is also asymptotically nice with respect to $U'$. 
	\end{itemize} 
\end{theorem}

The word hyperbolicity condition will be removed in Theorem \ref{du-thm-asymniceII}. 
However, there might be multiple sharpely supporting hyperspace at each point 
of $\partial \Omega$. 

When the linear parts of the affine maps are unimodular, 
then the uniform middle-eigenvalue condition automatically is true since 
$\lambda_E(\gamma) = 1$ for all $\gamma \in \Gamma$. 
This follows from the dual of Remark \ref{pr-rem-eigenlem}, 
which we will explain later. 

Theorem 8.2.1 of Labourie \cite{Labourie07} shows that such a domain $U$ exists but without 
establishing the asymptotic niceness of the group. 
Also, when the linear part of $\Gamma$ is a geometrically finite 
Kleinian group in $\SO(n, 1)$, 
Barbot showed this result in Theorem 4.25 of 
\cite{Barbot2005} in the context of 
globally hyperbolic Lorentzian spacetimes.
We believe our theory also generalizes to the case when 
$\mathcal{L}(\Gamma)$ is convex cocompact. 
Fried also found a solution when the linear   part of $\Gamma$ is in $\SO(2, 1)$ using 
cocycles \cite{Friedlens} with informal notes
in the same context but in the dual picture of R-ends as in Chapter \ref{ch-pr}.

As mentioned above Nie and Seppi \cite{NieSeppi23} and Ablondi \cite{Ablondi25} generalize
these results when the linear part is in $\SL_\pm(n+1, \bR)$ using the techniques of 
Cheng-Yau \cite{ChengYau} and the techniques of Barbot. They of course do not assume the word hyperbolicity here. 
Hence, the results of 
our Theorems \ref{du-thm-asymnice} and \ref{du-thm-asymniceII} overlap with their result.

\subsection{Flow setup} \label{du-sub-flow} 
The following flow setup will be applicable in the following. A slight modification 
is required later in Section \ref{du-sec-gen}.

We generalize the work of Goldman-Labourie-Margulis \cite{GLM09}
using Anosov flows:
We assume that $\Gamma$ has a properly convex 
affine action with the triple $(\Gamma, U, \Omega)$ for $U \subset \mathds{A}^n$. 
Since $\Omega$ is properly convex, 
$\Omega$ has a Hilbert metric. 
Let $T\Omega$ denote the tangent space of $\Omega$. 
Let $\Uu\Omega$ denote the {\em unit tangent bundle over} $\Omega$.
It has a smooth structure as a quotient space of $T\Omega - O/\sim$ where 
\begin{itemize}
	\item $O$ is the image of the zero-section, and 
	\item $\vec{v} \sim \vec{w}$ if $\vec{v}$ and $\vec{w}$ are over the same point of $\Omega$
	and $\vec{v} = s \vec{w}$ for a real number $s > 0$.
\end{itemize} 
\index{unit tangent bundle} 
\index{tangent bundle}

Let $\mathds{A}^n$ be the $n$-dimensional affine subspace. 
Let $h: \Gamma \ra \Aff(\mathds{A}^n)$ denote the representation 
as described in \eqref{du-eqn-affact}.  
We form the product $\Uu\Omega \times \mathds{A}^n$ that is an affine bundle over $\Uu\Omega$. 
We take the quotient of $\tilde \bA := \Uu\Omega \times \mathds{A}^n$ by the diagonal action 
\[g(x, \vec u)= (g(x), h(g) \vec u) \hbox{ for } g \in \Gamma, x \in \Uu\Omega, \vec u \in \mathds{A}^n.\] 
We denote the quotient by $\bA$, which fibers over
the smooth orbifold $\Uu\Omega/\Gamma$ with fiber $\mathds{A}^n$. 
\index{a@$\bA$|textbf} \index{a@$\tilde \bA$|textbf}

Let $V^n$ be the vector space associated with $\mathds{A}^n$.
Then we can form $\tilde \bV:= \Uu\Omega \times V^n$ and take the quotient under 
the diagonal action:
\[g(x, \vec u)= (g(x), {\mathcal L}(h(g)) \vec u) \hbox{ for } g \in \Gamma, x\in \Uu\Omega, \vec u \in V^n\]
where $\mathcal L$ is the homomorphism taking the linear part of $g$.  
We denote by $\bV$ the fiber bundle over $\Uu\Omega/\Gamma$ 
with fiber $V^n$. 
\index{v@$\bV$|textbf} \index{v@$\tilde \bV$|textbf}


There exists a flow $\hat \Phi_t: \Uu\Omega/\Gamma \ra \Uu\Omega/\Gamma$ for $t \in \bR$ given 
by sending $\vec v$ at $\alpha(0)$ to the unit tangent vector at $\alpha(t)$ where 
$\alpha$ is a geodesic tangent to $\vec v$. 
This flow is induced from the geodesic flow $\widetilde{\hat \Phi}_t: \Uu \Omega \ra \Uu \Omega$.

We define a flow on $\tilde \Phi_t: \tilde \bA \ra \tilde \bA$ by considering 
a unit-speed-geodesic flow-line $\vec{l}$ in $\Uu\Omega$ 
and $\vec{l} \times \mathds{A}^n$ and acting trivially on the second factor as we go from $\vec v$ to $\hat \Phi_t(\vec v)$
(See remarks in the beginning of Section 3.3 and equations in Section 4.1 of \cite{GLM09}.)
Each flow line in $\Uu\Sigma$ lifts to a flow line on $\bA$ from every point along it. This induces a flow $\Phi_t: \bA \ra \bA$.  

We define a flow on $\mathcal{L}(\tilde{\Phi}_t): \tilde \bV \ra \tilde \bV$ by considering a unit-speed geodesic-flow line $\vec{l}$ in $\Uu\Omega$ and 
considering $\vec{l} \times V^n$ and acting trivially on the second factor as we go from $\vec v$ to 
$\widetilde{\hat \Phi}_t(\vec v)$
for each $t$. 
This induces a flow $\mathcal{L}(\Phi_t): \bV \ra \bV$. 
(This generalizes the flow on \cite{GLM09}.)

We let $\llrrV{\cdot}_{\mathrm{fiber}}$ 
denote a metric on these bundles over $\Uu\Sigma/\Gamma$ defined as a fiberwise inner product:
We choose a cover of $\Omega/\Gamma$ by compact sets $K_i$, 
 choose a metric over $K_i \times \mathds{A}^n$, and use the partition of unity. 
This induces a fiberwise metric on $\bV$ as well. 
By pulling the metric back to $\tilde \bA$ and $\tilde \bV$, we obtain a fiberwise metric
to be denoted by $\llrrV{\cdot}_{\mathrm{fiber}}$. 

We recall the trivial product structure. 
$\Uu\Omega \times \mathds{A}^n$ is a flat $\mathds{A}^n$-bundle over $\Uu\Omega$ with a flat affine connection $\nabla^{\tilde \bA}$, 
and $\Uu\Omega \times V^n$ has a flat linear connection $\nabla^{\tilde \bV}$.
The action described above preserves the connections. 
We have a flat affine connection $\nabla^{\bA}$ on the bundle $\bA$ over $\Uu\Sigma$, 
and a flat linear connection $\nabla^{\bV}$ on the bundle $\bV$ over $\Uu\Sigma$. 

\begin{remark}\label{du-rem-GLM}
	In \cite{GLM09}, the authors use the term "recurrent geodesic.''  A geodesic is ``recurrent'' in their sense	if it accumulates to compact subsets in both directions. 
	They work in a compact subsurface where geodesics are recurrent in both directions. 
	In our work, since $\Omega/\Gamma$ is a closed orbifold, every geodesic is recurrent in their sense. 
	Hence, their theory generalizes to our case. 
\end{remark}

\section{The proximal flow.}\label{du-sub-anosov}

We will start with the case when $\Gamma$ is word-hyperbolic and hence
 $\Omega$ must be strictly convex 
with $\partial \clo(\Omega)$ being $C^1$ 
by Theorem 1.1 of \cite{Benoist04}.

When the linear parts of the affine maps are unimodular, 
Theorem 8.2.1 of Labourie \cite{Labourie07} shows that such a domain $U$ exists but without 
establishing the asymptotic niceness of the group. 
Also, when the linear part of $\Gamma$ is a geometrically finite 
Kleinian group in $\SO(n, 1)$, 
Barbot showed this result in Theorem 4.25 of 
\cite{Barbot2005} in the context of 
globally hyperbolic Lorentzian spacetimes.
We believe our theory also generalizes to the case when 
$\mathcal{L}(\Gamma)$ is convex cocompact. 
Fried also found a solution when the linear   part of $\Gamma$ is in $\SO(2, 1)$ using 
cocycles \cite{Friedlens} with informal notes
in the same context but in the dual picture of R-ends as in Chapter \ref{ch-pr}.

As mentioned above Nie and Seppi \cite{NieSeppi23} and Ablondi \cite{Ablondi25} generalizes 
these results when the linear part is in $\SL_\pm(n+1, \bR)$ using the techniques of 
Cheng-Yau \cite{ChengYau} and the techniques of Barbot. They of course do not assume the word hyperbolicity here.



The word-hyperbolicity of $\Gamma$ shows that $\Omega$ is strictly 
convex by Benoist \cite{Benoist04}. 
We will generalize the theorem to Theorem \ref{du-thm-asymniceII} 
without the word-hyperbolicity condition of $\Gamma$. 
Furthermore, we will 
show that the middle-eigenvalue condition actually implies the existence 
of the properly convex domain $U$ in Theorem \ref{du-thm-asymniceII}. 
Also, the uniqueness of the set of asymptotic hyperspaces is given by 
Theorem \ref{du-thm-ASunique}.

The reason for presenting weaker Theorem \ref{du-thm-asymnice}
is to convey the basic idea of the proof of the generalized theorem.  





\subsection{The decomposition of the flow}
%
We assume that $\Gamma$ is hyperbolic.  
Since $\Sigma:= \Omega/\Gamma$ is a closed strictly convex real projective orbifold, 
$\Uu\Sigma := \Uu\Omega/\Gamma$ is also a compact smooth orbifold. 
A geodesic flow on $\Uu\Omega/\Gamma$ is Anosov 
and hence topologically mixing. 
Therefore, the flow is nonwandering everywhere. (See \cite{Benoist01}.)
$\Gamma$ acts irreducibly on $\Omega$, and $\partial \clo(\Omega)$ is of class $C^1$. 
Denote by $\Pi_{\Omega}: \Uu \Omega \ra \Omega$ the projection to the base points. 
\index{us@$\Uu \Sigma$|textbf}
\index{uo@$\Uu \Omega$|textbf} \index{pi@$\Pi_\Omega$|textbf}
We identify $\Bd \mathds{A}^n = \SI(V^{n})= \SI^{n-1}$, where $g$ acts by  ${\mathcal{L}}(g) \in \GL(n, \bR)$. 
We give a decomposition of $\tilde \bV$ into three parts $\tilde \bV_+, \tilde \bV_0, \tilde \bV_-$: 
\begin{itemize} 
\item For each vector $\vec u \in \Uu\Omega$, we find the maximal oriented geodesic 
$l$ ending at the backward endpoint $\partial_+ l$ 
and the forward endpoint $\partial_- l \in \partial\clo(\Omega)$. They correspond to 
the $1$-dimensional vector subspaces $\tilde \bV_+(\vec u)$ and $\tilde \bV_-(\vec u) \subset V$. 
\item Recall that $\partial \clo( \Omega)$ is of class $C^{1}$ since $\Omega$ is strictly convex (see \cite{Benoist04}).
There exists a unique pair of sharply supporting hyperspheres $H_+$ and $H_-$ in $\Bd \mathds{A}^n$ at each of $\partial_+ l$ and $\partial_- l$. 
We denote by $H_0 = H_+ \cap H_-$. It is a codimension-two great sphere in $\Bd \mathds{A}^n$
and corresponds to a vector subspace $\tilde \bV_0(\vec u)$ of codimension-two in $\tilde \bV$. 
\item 
For each vector $\vec u$, we find the decomposition of $V$ as $\tilde \bV_+(\vec u) \oplus \tilde \bV_0(\vec u) \oplus \tilde \bV_-(\vec u)$
and hence we can form the subbundles $\tilde \bV_+, \tilde \bV_0, \tilde \bV_-$ over $\Uu\Omega$
where \[\tilde \bV = \tilde \bV_+ \oplus \tilde \bV_0 \oplus \tilde \bV_-.\] 
\end{itemize} 
The map $\Uu\Omega \ra \partial \clo(\Omega)$ sending a vector to the endpoint of the geodesic tangent to it is $C^{1}$.
The map $\partial \clo(\Omega) \ra \mathcal{H}$ sending a boundary point to its sharply supporting hyperspace 
in the space $\mathcal{H}$ of hyperspaces in $\Bd \mathds{A}^n$ is continuous.   
Hence $\tilde \bV_{+}, \tilde \bV_{0},$ and $\tilde \bV_{-}$ are 
continuous bundles. 
Since the action preserves the decomposition of $\tilde \bV$, 
$\bV$ also decomposes as 
\begin{equation} \label{du-eqn-Vdec} 
\bV = \bV_{+} \oplus \bV_0 \oplus \bV_-.
\end{equation} 
\index{v@$\bV$} \index{vp@$\bV_+$|textbf} \index{vz@$\bV_0$|textbf} \index{vm@$\bV_-$|textbf}

For each complete  geodesic $l$ in $\Omega$, let $\vec l$ denote the set of unit vectors on $l$ in one of the two directions.  
On $\vec{l}$, we have a decomposition
\begin{align*}
&\tilde \bV|\vec l = \tilde \bV_{+}|\vec l \oplus \tilde \bV_{0}|\vec l \oplus \tilde \bV_{-}|\vec{l}  \hbox{ of form }  \\
&\vec l \times \tilde \bV_{+}(\vec u), \vec l \times \tilde \bV_{0}(\vec u), \vec l\times \tilde \bV_{-}(\vec u)
\hbox{ for a vector $\vec u$ tangent to $l$}
\end{align*} 
where  we recall:
\begin{itemize}
\item $\tilde \bV_{+}(\vec u)$ is the space of vectors in the direction of the backward endpoint of $\vec{l}$. 
\item $\tilde \bV_{-}(\vec u)$ is the space of vectors in the direction of the forward endpoint of $\vec{l}$. 
\item $\tilde \bV_{0}(\vec u)$ is the space vectors in directions of $H_{0}= H_{+}\cap H_{-}$ for $\partial l$. 
\end{itemize}
That is, these bundles have constant fibers along $l$. 

Suppose that $g \in \Gamma$ acts on a complete geodesic $l$ with a unit vector 
$\vec u$ in the direction of the action of $g$. 
Then $\tilde \bV_-(\vec u)$ and $\tilde \bV_+(\vec u)$ corresponding to endpoints of $l$ 
are respectively eigenspaces of the eigenvalue of 
the largest norm  $\lambda_1(g)$ 
and the eigenvalue of the smallest norm $\lambda_n(g)$
of the linear part ${\mathcal{L}}(g)$ of $g$. 
Hence 
\begin{itemize}
\item on $\tilde \bV_-(\vec u)$, $g$ acts by expanding by 
\begin{equation} \label{du-eqn-exp} 
\frac{\lambda_1(g)}{\lambda_{{\tilde E}}(g)} > 1,
\end{equation} 
and 
\item on $\tilde \bV_+(\vec u)$, $g$ acts by contracting by 
\begin{equation} \label{du-eqn-cont}
\frac{\lambda_n(g)}{\lambda_{{\tilde E}}(g)} < 1.
\end{equation} 
\end{itemize}

There exists a flow $\hat \Phi_t: \Uu\Omega \ra \Uu\Omega$ for $t \in \bR$ given 
by sending $\vec v$ at $\alpha(0)$ to the unit tangent vector at $\alpha(t)$ where 
$\alpha$ is a geodesic tangent to $\vec v$ at $\alpha(0)$.

We define a flow $\tilde \Phi_t: \tilde \bA \ra \tilde \bA$ by considering 
a unit-speed geodesic-flow line $\vec{l}$ in $\Uu\Omega$ 
and considering $\vec{l} \times \mathds{A}^n$ and acting trivially on the second factor as we go from $\vec v$ to $\hat \Phi_t(\vec v)$.
(See remarks in the beginning of Section 3.3 and equations in Section 4.1 of \cite{GLM09}.)
Each flow line in $\Uu\Sigma$ lifts to a flow line on $\bA$ from every point on it. This induces a flow
 $\Phi_t: \bA \ra \bA$.  
 
We defined a flow $\tilde \Phi_t: \tilde \bV \ra \tilde \bV$ by considering 
a unit-speed geodesic-flow line $\vec{l}$ in $\Uu\Omega$ and 
 considering $\vec{l} \times V$, and acting trivially on the second factor as we go from $\vec v$ to $\tilde \Phi_t(\vec v)$
for each $t$. (This generalizes the flow on \cite{GLM09}.)
Furthermore, $\mathcal{L}(\tilde{\Phi}_t)$ preserves $\tilde \bV_+$, $\tilde \bV_0$, and 
$\tilde \bV_-$ since the endpoints $\partial_\pm l$ on the line $l$ remain unchanged during the flow. 
Again, this induces a flow 
\[\mathcal{L}(\Phi_{t}):  \bV \ra \bV, \bV_{+} \ra \bV_{+}, \bV_{0} \ra \bV_{0}, \bV_{-} \ra \bV_{-}.\]
\index{phit@$\Phi_t$} 
\index{phit@$\tilde\Phi_t$} 


We let $\llrrV{\cdot}_S$ denote a metric on these bundles over $\Uu\Sigma/\Gamma$ defined as a fiberwise inner product:
We choose a cover of $\Omega/\Gamma$ by compact sets $K_i$, 
choose a metric over $K_i \times \mathds{A}^n$ and use the partition of unity. 
This induces a fiberwise metric on $\bV$ as well. 
Pulling back the metric to $\tilde \bA$ and $\tilde \bV$, we obtain fiberwise metrics
to be denoted by $\llrrV{\cdot}_{S}$.

We will prove the following: 
By the uniform middle-eigenvalue condition,  
$\bV$ satisfies the following properties for $\vec u \in \Uu \Omega/\Gamma$: 
\begin{itemize} 
\item the flat linear connection $\nabla^{\bV}$ on $\bV$ 
is bounded with respect to $\llrrV{\cdot}_{\mathrm{fiber}}$. \label{vertfiber} 
\item hyperbolicity: There exist constants $C, k > 0$ such that 
\begin{align} \label{du-eqn-Anosov1} 
\llrrV{\mathcal{L}(\Phi_t)(\vec{v})}_{\mathrm{fiber}} \geq \frac{1}{C} \exp(kt) \llrrV{\vec{v}}_{\mathrm{fiber}} \hbox{ as } t \ra \infty
\end{align}
for $\vec{v} \in \bV_+$ and 
\begin{align} \label{du-eqn-Anosov2} 
 \llrrV{\mathcal{L}(\Phi_t)(\vec{v})}_{\mathrm{fiber}} \leq C \exp(-kt) \llrrV{\vec{v}}_{\mathrm{fiber}} \hbox{ as } t \ra \infty
\end{align}
for $\vec{v} \in \bV_-$
\end{itemize} 
Using 
Proposition \ref{du-prop-contr}, 
we prove these properties by taking $C$ sufficiently large according to $t_1$, 
which is a standard technique. 
\index{fiberwise metric|textbf}


\subsection{The proof of the proximal property.} \label{du-sub-Anosov}


We may assume that $\Gamma$ has no element of finite order  by 
taking a finite-index subgroup using Theorem \ref{prelim-thm-vgood}. 
Also, by Benoist \cite{Benoist04}, the elements of 
$\Gamma -\{\Idd\}$ are positive bi-proximal. (See Theorem \ref{prelim-thm-semi}.)

We can apply this property 
to $\bV_-$ and $\bV_+$ by possibly reversing the direction of 
the flow. 
The Anosov property follows from the following proposition. 

Let $\bV_{-,1}$ denote the subset of $\bV_-$ of the unit length under $\llrrV{\cdot}_{\mathrm{fiber}}$.

\begin{proposition} \label{du-prop-contr} 
Let $\Omega/\Gamma$ be a closed strictly convex real projective $(n-1)$-orbifold with hyperbolic fundamental group $\Gamma$. 
Then there exists a constant $t_1$ such that 
\[ \llrrV{\mathcal{L}(\Phi_t)(\bv)}_{\mathrm{fiber}} \leq \tilde C \llrrV{\bv}_{\mathrm{fiber}}, \bv \in \bV_- \hbox{ and } \llrrV{\mathcal{L}(\Phi_{-t})(\bv)}_{\mathrm{fiber}} \leq \tilde C \llrrV{\bv}_{\mathrm{fiber}}, \bv \in \bV_+\] 
for $t \geq t_1$ and a uniform $\tilde C$, $0 < \tilde C < 1$.
\end{proposition}
\begin{proof} 
It sufficies to prove the first part of the inequalities since we can substitute $t \ra -t$
and switch $\bV_{+}$ with $\bV_{-}$ as the direction of the vector changed to the opposite one. 

By following Lemma \ref{du-lem-S_0}, the uniform convergence implies that for given $0< \eps < 1$, 
for every vector $\bv$ in $\bV_{-, 1}$, there exists a uniform $T$ such that for $t > T$, 
$\mathcal{L}(\Phi_t)(\bv)$ is in 
the $\eps$-neighborhood $U_\eps(S_0)$ of the image $S_0$ of the zero section. 
Hence, we obtain that $\mathcal{L}(\Phi_{t})$ is uniformly contracting near $S_{0}$, which implies the result. 
\end{proof} 

Now, we will prove Lemma  \ref{du-lem-S_0}; but we need some preliminary material: 
\begin{remark} \label{pr-rem-virtual}
	We only need to prove the following for a finite-index subgroup of $\bGamma$ since 
	the contracting properties are invariant under finite regular covering maps. 
	Hence, we may assume that each element is positive proximal or positive semiproximal 
	by Theorem \ref{prelim-thm-semi} depending on whether 
$\Omega$ is strictly convex or just properly convex 
	respectively. 
\end{remark} 

\begin{lemma} \label{du-lem-length} 
	Let $\Gamma$ act properly discontinuously on a strictly convex domain 
	$\Omega$. 
	Assume that $g_i$ is a sequence of distinct positive bi-proximal elements of $\Gamma$. 
	Suppose that the sequence of attracting fixed points $a_i$ and that of repelling fixed point $\{r_i\}$ of $g_i$ converges to distinct pair of points respectively. 
	Then 
\begin{equation}\label{du-eqn-length}
\{\leng_\Omega(g_i)\} \ra \infty
\end{equation}
\end{lemma} 
\begin{proof} 
	Since $\{a_i\} \ra a_\ast$ and 
	$\{r_i\} \ra r_\ast$, the segments 
	$\overline{a_ir_i}$ pass through a fixed compact domain $U$ in $\Omega$ for 
	sufficiently large $i$. 
Suppose that $\leng_\Omega(g_i) < C$ for a constant $C$. 
Then $g_i(U)$ passes through $\overline{a_ir_i}$ for each $i$. 
Hence, $g_i(U)$ is a subset of a ball of radius $2L + C$. 
Since $\{g_i\}$ is a sequence of mutually distinct elements,
this contradicts the proper discontinuity of the action of $\Gamma$. 
\end{proof} 

\begin{lemma} \label{du-lem-convergence} 
	Let $\Omega$ be strictly convex and $C^1$.  
We choose a subsequence $\{g_i\}$  of positive-bi-proximal elements of 
$\Gamma$ such that the sequences $\{a_i\}$ and $\{r_i\}$
are convergent for the attracting fixed point $a_i \in \clo(\Omega)$ 
and the repelling fixed point $r_i \in \clo(\Omega)$ of each $g_i$. 
Suppose that 
\[\{a_i\} \ra a_\ast \hbox{ and } \{r_i\} \ra r_\ast \hbox{ for } a_\ast, r_\ast \in \partial \clo(\Omega), a_\ast \ne r_\ast.\] 
Suppose that $g_i$ is an unbounded sequence. 
Then for every compact $K \subset \clo(\Omega) - \{r_\ast\}$,
\begin{equation}\label{du-eqn-giK2}
\{g_i(K)\} \ra \{a_\ast\} 
\end{equation} 
uniformly. 
\end{lemma} 
\begin{proof} 
Each $g_i$ acts on an $(n-3)$-dimensional subspace $W_{g_i}$ in 
$\SI^{n-1}_\infty$ disjoint from $\Omega$. 
Here, $W_{g_i}$ is the intersection of two sharply supporting hyperspaces of 
$\Omega$ at $a_i$ and $r_i$. 
The set $\{W_{g_i}\}$ is precompact by our condition. 
By the $C^1$-property, 
we may assume that $\{W_{g_i}\} \ra W_\ast$ for an $n-3$-dimensional 
subspace $W_\ast$ that is the intersection of 
two hyperspaces supported at $a_\ast$ and $r_\ast$. 
Also, $W_{g_i} \cap \clo(\Omega) = \emptyset$ by this property. 

Let $\eta_i$ denote the complete geodesic connecting $a_i$ and $r_i$, and 
let $\eta_\infty$ denote the one connecting 
$a_\ast$ and $r_\ast$. 
Since $W_\ast$ is the intersection of two sharply supporting hyperspaces of 
$\Omega$ at $a_\ast$ and $r_\ast$, $\eta_\infty$ has endpoints 
$a_\ast, r_\ast$, and $\Omega$ is strictly convex,  it follows 
$\langle \eta_\infty \rangle  \cap W_\ast = \emp$. 

For each $g_i$, we call $P_i \cap \Omega$  a {\em slice } of $g_i$
for the $(n-2)$-dimensional subspace $P_i$ containing $W_{g_i}$. 
The closure of a component of $\Omega$ with a slice of $g_i$ removed 
is called a {\em half-space} of $g_i$. 

Let $H_i$ denote the half-space of $g_i$ containing $K$. 
Since $\{\clo(\eta_i)\}$ and $\{W_{g_i}\}$ are 
geometrically convergent respectively, 
and $\langle \eta_\infty \rangle  \cap W_\infty = \emp$,
it follows that  
$\{g_i(P_i)\}$ geometrically converges to a hyperspace containing 
$W_\infty$ passing $a_\ast$. 
Therefore, one deduces easily that 
$\{g_i(H_i)\} \ra \{a_\ast\}$ geometrically.
%
Since $K \subset H_i$, the lemma follows. 
\end{proof}

The line bundle $\bV_-$ lifts to $\tilde \bV_-$ where for each unit vector $\bu$ on $\Omega$, 
we associate the line $\bV_{-,\bu}$ corresponding to the starting point in $\partial \clo(\Omega)$ of the oriented geodesic $l$ tangent to it.
$\tilde \bV_{-}| \vec l$ equals $\vec l \times \bV_{-, \bu}$. 
$\mathcal{L}(\Phi_{t})$ lifts to a parallel translation or constant flow $\mathcal{L}(\tilde{\Phi}_{t})$
of the form 
\[ (\bu, \vec v) \ra (\hat \Phi_{t}(\bu), \vec v).\]

\begin{figure}
\centering
\includegraphics[height=8cm]{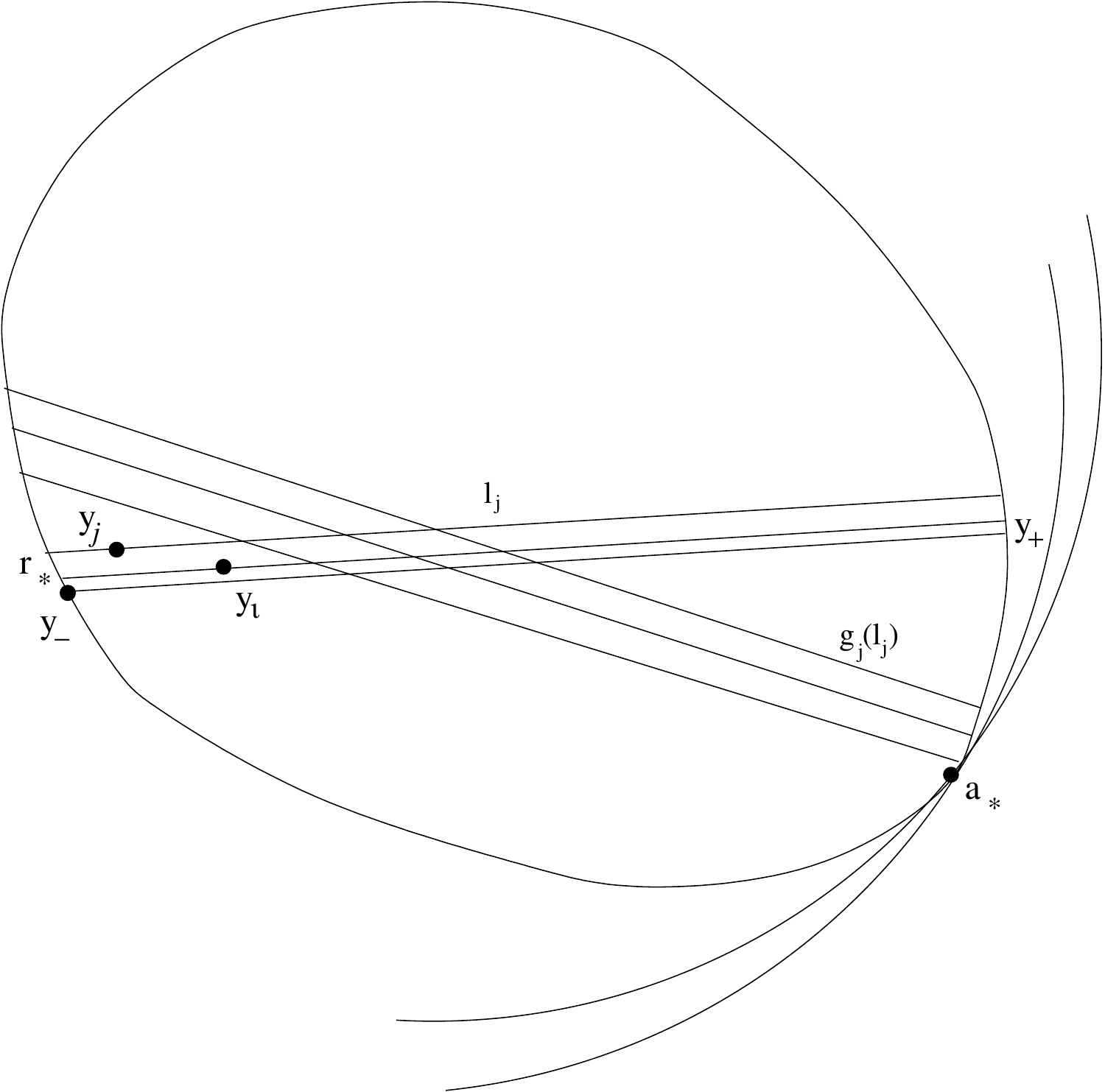}
\caption{The figure for Lemma \ref{du-lem-S_0}. 
Here $y_i$ and $y_j$ denote the images under $\Pi_{\Omega}$ of 
 the named points in the proof of Lemma \ref{du-lem-S_0}.}

\label{du-fig-contr}

\end{figure}




We will generalize the following to Lemma \ref{du-lem-S_0II}.

\begin{lemma} \label{du-lem-S_0} 
Suppose that $\Omega$ is strictly convex with $\partial \clo(\Omega)$ being $C^1$, 
and $\Gamma$ acts properly discontinuously and cocompactly on $\Omega$ 
satisfying the uniform middle-eigenvalue condition. Then 
$\{\llrrV{\mathcal{L}(\Phi_{t})| \bV_-}_{\mathrm{fiber}}\}  \ra 0 $ uniformly as $t \ra \infty$. 
\end{lemma}
\begin{proof} 
%
Let $F$ be a fundamental domain of $\Uu\Omega$ under $\Gamma$.
It sufficies to prove this for $\mathcal{L}(\tilde{\Phi}_{t})$ on the fibers over 
$F$ of $\Uu\Omega$ with a fiberwise metric $\llrrV{\cdot}_{\mathrm{fiber}}$.

We choose an arbitrary sequence $\{x_i\}$, $\{x_i\} \ra x$ in $F$.
For each $i$, let $\bv_{-, i}$ be a Euclidean unit vector in  $V_{-,i}:= \tilde \bV_{-}(x_{i})$ for the unit vector $x_i \in \Uu\Omega$.
That is, $\bv_{-, i}$ is in the $1$-dimensional subspace in $\bR^{n}$,  
corresponding to the backward endpoint of the geodesic $l_i$ in $\Omega$ 
determined by $x_i$ in $\partial\clo(\Omega)$ and in a direction of $\clo(\Omega)$.

We will show that 
\[\left\{\llrrV{\mathcal{L}(\tilde{\Phi}_{t_i})(x_{i}, \bv_{-,i})}_{\mathrm{fiber}}\right\} \ra 0 \hbox{ for any sequence } t_i \ra \infty,\] 
which is sufficient to prove the uniform convergence to $0$ by the compactness of  $\bV_{-,1}$. 

It sufficies 
to show that any sequence of $\{t_i\} \ra \infty$ has a subsequence $\{t_j\}$ such that 
\[\{\llrrV{\mathcal{L}(\tilde{\Phi}_{t_j})(x_{i}, \bv_{-,  j})}_{\mathrm{fiber}}\} \ra 0.\]
This follows since if the uniform convergence did not hold, then we can easily find a sequence without such subsequences. 

Denote $y_i := \hat{\Phi}_{t_i}(x_i)$ for the lift of the flow $\hat \Phi$. 
By construction, we recall that each $\Pi_{\Omega}(y_i)$ is in the geodesic $l_i$. 
Since we have the sequence of vectors $\{x_i\} \ra x$, $x_{i}, x \in F$, 
we obtain that $\{l_i\}$ geometrically converges to a line $l_\infty$ passing $\Pi_\Omega(x)$ in $\Omega$. 
Let $y_+$ and $y_-$ be the endpoints of $l_\infty$ where $\{\Pi_{\Omega}(y_i)\} \ra y_-$.
Hence, \[\{\llrrparen{\bv_{+, i}}\} \ra y_{+}, \{\llrrparen{\bv_{-, i}}\} \ra y_{-}.\] 

Find a deck transformation $g_i$ such that 
$g_i(y_i) \in F$ and $g_i$ acts on the line bundle $\tilde \bV_-$ by the linearization of the matrix of the form in \eqref{du-eqn-bendingm4}: 
\begin{align} \label{du-eqn-Lgi} 
\mathcal{L}(g_{i}) &: \tilde \bV_{-} \ra \tilde \bV_{-} \hbox{ with } \nonumber \\
 (y_{i}, \bv) & \mapsto (g_{i}(y_{i}), {\mathcal L}(g_{i})(\bv)) \hbox{ where } \nonumber \\
{\mathcal L}(g_i) & = \frac{1}{\lambda_{{\tilde E}}(g_i)^{1+\frac{1}{n}}} \hat h(g_i) : \tilde \bV_{-}(y_i) = \tilde \bV_{-}(x_{i}) \ra \tilde \bV_{-}(g_i(y_i)).
\end{align}
\begin{description}
\item[{\rm (Goal)}] We will show $\{(g_{i}(y_{i}), {\mathcal L}(g_i)(\bv_{-, i}))\} \ra 0$ under $\llrrV{\cdot}_{\mathrm{fiber}}$. 
This will complete the proof since $g_{i}$ acts as isometries on $\tilde \bV_{-}$ with $\llrrV{\cdot}_{\mathrm{fiber}}$. 
\end{description} 
Also, we may assume that $\{g_i\}$ is a sequence of mutually distinct elements up to a choice of 
subsequences since 
$g_i^{-1}(F)$ contains $y_i$, and $\{y_i\}$ is an unbounded sequence. 

Since $g_i(l_i) \cap F \ne \emp$, we choose a subsequence of $g_i$ and relabel it as $g_i$ such that
$\{g_i(l_i)\}$ converges to a nontrivial line $\hat l_\infty$ 
in $\Omega$.



By our choice of $l_i$, $y_i$, $g_i$ as above, 
and by Remark \ref{pr-rem-virtual}, we may assume without loss of generality 
that each $g_i$ is positive bi-proximal 
since $\Omega$ is strictly convex. 

We let $\hat g_i$ to denote $\mathcal{L}(g_i)$. 
We choose a subsequence of $\{\hat g_i\}$ such that the sequences $\{a_i\}$ and $\{r_i\}$
are convergent for the attracting fixed point $a_i \in \clo(\Omega)$ 
and the repelling fixed point $r_i \in \clo(\Omega)$ of each $\hat g_i$. 
Then 
\[\{a_i\} \ra a_\ast \hbox{ and } \{r_i\} \ra r_\ast \hbox{ for } a_\ast, r_\ast \in \partial \clo(\Omega).\] 
(See Figure \ref{du-fig-contr}.) 

Suppose that $a_{\ast} = r_{\ast}$. Then we choose an element $g \in \Gamma$ such that $\hat g(a_{\ast}) \ne r_{\ast}$
and replace the sequence by $\{g \hat g_i\}$ while replacing $F$ by $F \cup g(F)$. 
The above uniform convergence condition still holds. 
Thus for the new attracting fixed points $a'_i$ of $\hat g\hat g_i$, 
we have $\{a'_i\} \ra \hat g(a_\ast)$.
Furthermore, the sequence $\{r'_i\}$ of repelling fixed point $r'_i$ of $\hat g\hat g_i$ converges to $r_\ast$ as shown by Lemma \ref{pr-lem-gatt}.
Hence, we may assume, without loss of generality, that \[a_\ast \ne r_\ast\]
by replacing our sequence $\hat g_i$.

Now, Lemma \ref{du-lem-convergence} shows that 
for every compact $K \subset \clo(\Omega) - \{r_\ast\}$,
\begin{equation}\label{du-eqn-giK} 
\{\hat g_i(K)\} \ra \{a_\ast\} 
\end{equation} 
uniformly. 

Suppose that both $y_+, y_- \ne r_\ast$. Then $\{\hat g_i(l_i)\}$ converges to a singleton $\{a_\ast\}$ by \eqref{du-eqn-giK} and this cannot happen 
since $\hat l_\infty \subset \Omega$. 
If \[r_\ast = y_+ \hbox{ and } y_- \in \partial\clo(\Omega) - \{r_\ast\},\] then 
$\{\hat g_i(y_i)\} \ra a_\ast$ by \eqref{du-eqn-giK} again. 
Since $\hat g_i(y_i) \in F$, this is a contradiction. 
Therefore \[r_\ast = y_- \hbox{ and } y_+ \in \partial \clo(\Omega)-\{r_\ast\}.\] 

Let $d_i= \llrrparen{\bv_{+, i}}$ 
denote the other endpoint of $l_i$ from $\llrrparen{\bv_{-, i}}$. 
\begin{itemize} 
\item Since $\{\llrrparen{\bv_{-, i}}\} \ra y_-$ and $\{l_i\}$ converges to a nontrivial line $l_\infty$, it follows that 
$\{d_i = \llrrparen{\bv_{+, i}}\}$ is in a compact set in $\partial \clo(\Omega) -\{r_\ast\}$,
and $\{d_i\} \ra y_+$. 
\item Then $\{\hat g_i(d_i)\} \ra a_{\ast}$ as $\{d_i\}$ is in a compact set in $\partial \clo(\Omega) -\{r_\ast\}$.
\item Thus, $\{\hat g_i(\llrrparen{\bv_{-, i}})\} \ra y' \in \partial \clo(\Omega)$ where $a_\ast \ne y'$ holds since 
$\{\hat g_i(l_i)\}$ converges to a nontrivial line $\hat l_{\infty}$ in $\Omega$
as we mentioned above.
\end{itemize} 


Also, $\hat g_i$ has an invariant great sphere $\SI^{n-2}_i \subset \Bd \mathds{A}^n$ 
containing the attracting fixed point $a_i$ and sharply supporting $\Omega$ at $a_{i}$.
Thus, $r_i$ is uniformly bounded at a distance from $\SI^{n-2}_i$ since $\{r_i\} \ra y_-= r_\ast$
and $\{a_{i}\} \ra a_{\ast}$ where $\{\SI^{n-2}_{i}\}$ geometrically converging to a sharply supporting sphere $\SI^{n-2}_{\ast}$
at $a_{\ast}$. 

Let $\llrrV{\cdot}_E$ denote the standard Euclidean metric of $\bR^{n}$. 
\begin{itemize}
\item Since $\{\Pi_\Omega(y_i)\} \ra y_-$, $\Pi_\Omega(y_i)$ 
is also uniformly bounded away from $a_i$ and the tangent sphere $\SI^{n-1}_i$ at $a_i$. 
\item Since $\{\llrrparen{\bv_{-, i}}\} \ra y_-$,
the vector $\bv_{-, i}$ has the component $\bv_i^p$ parallel to $r_i$ and the component $\bv_i^S$ 
in the direction of $\SI^{n-2}_i$ where $\bv_{-, i} = \bv_i^p + \bv_i^S$. 
\item Since $\{r_i\} \ra r_\ast= y_-$ and $\{\llrrparen{\bv_{-, i}}\} \ra y_-$, we obtain $\{\llrrV{\bv_i^S}_E\} \ra 0$ and that 
\[ \frac{1}{C} < \llrrV{\bv_i^p}_E < C  \]
for some constant $C> 1$.
\item $\hat g_i$ acts by preserving the directions of $\SI^{n-2}_i$ and $r_i$. 
\end{itemize} 
Since $\{g_i(\llrrparen{\bv_{-, i}})\}$ converges to $y'$, $y'\in \partial \clo(\Omega)$, and is uniformly bounded away from $\SI^{n-2}_i$, we obtain that the sequence of
the Euclidean norms of 
\[\frac{{\mathcal L}(g_i)(\bv_i^S)}{\llrrV{{\mathcal L}(g_i)(\bv_i^p)}_E}\] is bounded above uniformly
by considering the homogeneous coordinates
\[\llrrparen{{\mathcal L}(g_i)(\bv_i^S): {\mathcal L}(g_i)(\bv_i^p) }.\] 
Since $r_i$ is a repelling fixed point of $g_i$ and $\{\llrrV{\bv_i^p}_E\}$ is uniformly bounded above, 
$\{{\mathcal L}(g_i)(\bv_i^p)\} \ra 0$ by \eqref{du-eqn-umecDii} and
\eqref{du-eqn-length}.
\[\{{\mathcal L}(g_i)(\bv_i^p)\} \ra 0 \hbox{ implies } \{{\mathcal L}(g_i)(\bv_i^S)\} \ra 0\]
for $\llrrV{\cdot}_{E}$.  
Hence, we obtain 
 \begin{equation}\label{du-eqn-bv}
\left\{\llrrV{\mathcal{L}(g_{i})(\bv_{-,i}))}_{E}\right\} \ra 0.
\end{equation}

Recall that $\mathcal{L}(\tilde{\Phi}_{t})$ is the identity map on the second factor of $\Uu\Omega \times V^n$. 
\[g_{i}(\mathcal{L}(\tilde{\Phi}_{t_{i}})(x_{i}, \bv_{-,i})) = (g_{i}(y_{i}), {\mathcal L}(g_{i})(\bv_{-, i})) \in F \times \bV_{-} \]
 is a vector over the compact fundamental domain $F$ of $\Uu\Omega$. 
\[\left\{\llrrV{\mathcal{L}(g_i)(\tilde \Phi_{t_i})(x_i, \bv_{-, i})}_E\right\} \ra 0
\hbox{ implies } 
\left\{\llrrV{\mathcal{L}(g_i)(\tilde \Phi_{t_i})(x_i, \bv_{-, i})}_{\mathrm{fiber}}\right\} \ra 0\]
since $g_i(x_i) \in F$ and $\llrrV{\cdot}_{\mathrm{fiber}}$ and 
$\llrrV{\cdot}_E$ are compatible over points of $F$. 
Since $g_i$ is an isometry of $\llrrV{\cdot}_{\mathrm{fiber}}$,   
we conclude that $\{\llrrV{\mathcal{L}(\tilde{\Phi}_{t_i})(x,\bv_{-,i})}_{\mathrm{fiber}} \}\ra 0$:
\end{proof}






\subsection{The neutralized section.} \label{du-sub-neutral}

We denote by $\Gamma(\bV)$ the space of sections 
$\Uu\Sigma \ra \bV$, and by $\Gamma(\bA)$ the space of
sections $\Uu\Sigma \ra \bA$. 
As introduced in \cite{GLM09}, the one parameter-group of 
operators $D\Phi_{t, *}$ on $\Gamma(\bV)$ and 
$\Phi_{t, *}$ on $\Gamma(\bA)$. 
In our terminology $D\Phi_{t, *} = \mathcal{L}(\Phi_t)$. 
Recall Lemma 8.3 of \cite{GLM09} also. 
We denote by $\phi$ the vector field generated by 
this flow on $\Uu\Sigma$.

A section $s:\Uu\Sigma \ra \bA$ is {\em neutralized} if 
\begin{equation}\label{du-eqn-neu}
\nabla^{\bA}_{\phi} s \in \bV_0 . 
\end{equation}
\index{neutralized section} 

\begin{lemma}\label{du-lem-lem83}
 If $\psi \in \Gamma(\bA)$ and 
\[t \mapsto D\Phi_{t, *}(\psi) \] 
is a path in $\Gamma(\bV)$ that is differentiable at $t =0$, then 
\[ \left.\frac{d}{dt}\right|_{t=0}  (D\Phi_t)_* (\psi)   = \nabla^{\bA}_\phi(\psi). \]
\end{lemma} 

Recall that $\Uu\Sigma$ is a recurrent set under the geodesic flow. 
\begin{lemma}\label{du-lem-exist} 
There exists a neutralized section $s_0:\Uu\Sigma \ra \bA$. 
This lifts to a map $\tilde s_0: \Uu\Omega \ra \tilde \bA$  
such that $\tilde s_0 \circ \gamma = \gamma \circ \tilde s_0$
for each $\gamma$ in $\Gamma$ acting on 
$\tilde \bA = \Uu \Omega \times \mathds{A}^n$. 
\end{lemma}
\begin{proof} 
Let $s$ be a continuous section $\Uu\Sigma \ra \bA$. 
We construct $\nabla^{\bA_\pm}$ by projecting the values of 
$\nabla$ to $\bV_{\pm}$ and $\nabla^{\bA_0}$ by projecting the values 
of $\nabla$ to $\bV_{0}$. 
We decompose 
\[\nabla^{\bA}_\phi(s) = \nabla^{\bA_+}_\phi(s) + \nabla^{\bA_0}_\phi(s) + \nabla^{\bA_-}_\phi(s) \in \bV \]
where 
$\nabla^{\bA_\pm}_\phi(s) \in \bV_\pm$ and $\nabla^{\bA_0}_\phi(s)  \in \bV_0$ hold. 
This can be done since along the vector field $\phi$,
$\bV_{\pm}$ and $\bV_0$ are constant bundles.  
By the uniform convergence property of \eqref{du-eqn-Anosov1} and \eqref{du-eqn-Anosov2}, 
the following integrals converge to smooth functions over $\Uu\Sigma$. 
Again
\[ s_0 = s + \int_0^\infty (D\Phi_t)_*(\nabla^{\bA_-}_\phi(s)) dt - \int_0^\infty (D\Phi_{-t})_*(\nabla^{\bA_+}_\phi(s)) dt\] 
is a continuous section, and 
$\nabla^{\bA}_\phi(s_0) = \nabla^{\bA_0}_\phi(s_0) \in \bV_0$ as shown in
Lemma 8.4 of \cite{GLM09}. 

Since $\Uu\Sigma$ is connected, there exists a fundamental domain $F$ 
such that we can lift $s_0$ to $\tilde s_0'$ defined on $F$ mapping to $\bA$. 
We can extend $\tilde s_0'$ to $\Uu\Omega \ra \Uu \Omega \times \mathds{A}^n$. 
\end{proof} 

Let $N_2(\mathds{A}^n)$ denote the space of codimension-two affine subspaces of $\mathds{A}^n$. 
We denote by $G(\Omega)$ the space of maximal oriented geodesics in $\Omega$.
We use the quotient topology on both spaces. 
There exists a natural action of $\Gamma$ on both spaces. 
\index{N2@$N_2(\cdot)$} 

For each element $g\in \Gamma -\{\Idd\}$, we define $N_2(g)$: 
Now, $g$ acts on $\Bd \mathds{A}^n$ with invariant subspaces corresponding to 
those of the linear part ${\mathcal L}(g)$ of $g$. 
Since $g$ and $g^{-1}$ are positive proximal,  
\begin{itemize}
\item a unique fixed point of $g$ in $\Bd \mathds{A}^{n}$ corresponds
to an eigenvector of an eigenvalue of the largest norm and an attracting fixed point of $g$ in $\Bd \mathds{A}^n$, 
and 
\item a unique fixed point of $g$ in $\Bd \mathds{A}^{n}$ corresponds to 
an eigenvector of an eigenvalue of the smallest norm and a repelling fixed point of $g$. 
\end{itemize} 
by \cite{Benoist01} or \cite{Benoist97}.
There exists an ${\mathcal L}(g)$-invariant 
vector subspace $\bV_g^0$ complementary to the sum of the subspace generated 
by these eigenvectors. (This space equals $\bV_0(\vec u)$ for the unit tangent vector $\vec u$
to the unique maximal geodesic $l_g$ in $\Omega$ on which $g$ acts.)
It corresponds to a $g$-invariant 
subspace $M(g)$ of codimension two in $\Bd \mathds{A}^n$.

Let $\tilde c$ be the geodesic in $\Uu\Sigma$ that is $g$-invariant for $g \in \Gamma$. 
$\tilde s_0(\tilde c)$ lies on a fixed affine subspace parallel to $V_g^{0}$ by neutrality, i.e., Lemma \ref{du-lem-exist}. 
There exists a unique affine subspace $N_2(g)$ of codimension two in $\mathds{A}^n$ 
containing $\tilde s_0(\tilde c)$. 
It immediately follows that 
$N_2(g) = N_2(g^{m}), m \in \bZ -\{0\}$, 
and that $g$ acts  on $N_{2}(g)$. 

\begin{definition}\label{du-defn-tau} 
We define $S'(\partial \clo(\Omega))$ as the space of hyperspaces $P$ meeting 
$\mathds{A}^n$ 
where $P \cap \Bd \mathds{A}^n$ is a sharply supporting hyperspace in 
$\Bd \mathds{A}^n$ to $\Omega$. 
We denote by $S(\partial \clo(\Omega))$ the space of pairs 
$(x, H)$ where $H \in S'(\partial \clo(\Omega))$, and $x$ is in the boundary of 
$H$ and in $\partial \clo(\Omega)$. 
\end{definition}

Define $\Delta$ to be the diagonal set of $\partial\clo(\Omega) \times \partial \clo(\Omega)$. 
We denote by $\Lambda^* = \partial\clo(\Omega) \times \partial \clo(\Omega)- \Delta$. 
Let $G(\Omega)$ denote the space of maximal oriented geodesics in $\Omega$. 
Then $G(\Omega)$ is in a one-to-one correspondence with $\Lambda^*$ by 
the map taking the maximal oriented geodesic to the ordered pair of its endpoints. 

\begin{lemma} \label{du-lem-Gconnect} 
	$G(\Omega)$ is a connected subspace.
	\end{lemma} 
\begin{proof} 
The proof follows from generalizing Lemma 1.3 of \cite{GLM09}
using a bi-proximal element of $\Gamma$.
\cite{Benoist04}. Or alternatively, one can use the fact that 
$\Uu \Omega$ is connected. 
	\end{proof}

\begin{proposition}\label{du-prop-mapgh}  $ $ 
\begin{itemize}
\item There exists  a continuous function $\hat s: \Uu\Omega \ra N_2(\mathds{A}^n)$ 
equivariant with respect to $\Gamma$-actions.
\item For each $g \in \Gamma$, let $\vec{l}_g$ denote the unique unit-speed geodesic  in 
$\Uu\Omega$ lying over a geodesic $l_g$ where $g$ acts on. 
Then $\hat s(\vec{l}_g) = N_2(g)$. 
\item This gives a continuous map 
\[\bar s: G(\Omega)=\partial\clo(\Omega) \times \partial \clo(\Omega) - \Delta \ra N_2(\mathds{A}^n)\] 
again equivariant with respect to the $\Gamma$-actions. 
There exists a continuous function 
\[\tau:\Lambda^* \ra S(\partial\clo(\Omega).\]
\end{itemize}
\end{proposition}
\begin{proof} 
Given a vector ${\vec u} \in \Uu\Omega$, we consider $\tilde s_0(\vec{u})$.
There exists a lift $\tilde \phi_t: \Uu\Omega \ra \Uu\Omega$ of 
the geodesic flow $\phi_t$.
Then $\tilde s_0(\tilde \phi_t({\vec u}))$ is in 
an affine subspace $H^{n-2}$ parallel to $V_0$ for $\vec u$ by the neutrality 
condition \eqref{du-eqn-neu}.
We define $\hat s(\vec u)$ to be this $H^{n-2}$. 

For any unit vector $\vec u'$ on the maximal (oriented) geodesic 
in $\Omega$ determined by $\vec u$, we obtain
$\hat s(\vec u') = H^{n-2}$. 
Hence, this determines the continuous map $\bar s: G(\Omega) \ra N_2(\mathds{A}^n)$. 
The $\Gamma$-equivariance follows from that of $\tilde s_0$. 

Let $g \in \Gamma$. 
Suppose that $\vec u$ is tangent to ${\vec l}_g$. 
Then $\vec u$ and $g(\vec u)$ lie on the lift ${\vec l}'_g$. 
Since $g(\tilde s_0(\vec u)) = \tilde s_0(g(\vec u))$ by equivariance, 
$g(\tilde s_0(\vec u))$ lies on
$\hat s(\vec u) = \hat  s(g(\vec u))$ in $\Uu \Omega$ 
by the paragraph before Definition \ref{du-defn-tau}. 
We conclude $g(\bar s({\vec l}'_g)) = \bar s({\vec l}'_g)$, 
which shows $N_2(g) = \hat s({\vec l}_g)$.

The map $\bar s$ is defined since
$\partial\clo(\Omega) \times\partial\clo(\Omega) - \Delta$ is in one-to-one correspondence with 
the space $G(\Omega)$. 
The map $\tau$ is defined by taking for each pair $(x, y) \in \Lambda^*$
\begin{itemize} 
\item the geodesic $l$ with endpoints $x$ and $y$,
and 
\item the hyperspace containing $\bar s(l)$ and its boundary containing $x$.  
\end{itemize} 
\end{proof}


\subsection{The asymptotic niceness.} \label{du-sub-asymnice}


We denote by $h(x, y)$ the hyperspace part in 
$\tau(x, y) = (x, h(x, y))$. 

\begin{lemma}\label{du-lem-hdisj} 
Let $U$ be a $\bGamma_{\tilde E}$-invariant 
properly convex open domain in $\bR^n$ such that 
$\Bd U \cap \Bd \mathds{A}^n = \clo(\Omega)$. 
Suppose that $x$ and $y$ are attracting and repelling fixed points of an element $g$ of $\Gamma$ in $\partial\clo(\Omega)$. 
Then $h(x, y)$ is disjoint from $U$. 
\end{lemma} 
\begin{proof}
Suppose not. 
$h'(x, y):= h(x, y) \cap \mathds{A}^n$ is a $g$-invariant open hemisphere, and $x$ is an attracting fixed point of $g$ in it.
If necessary, we can replace $g$ with $g^{-1}$.
Then $U \cap h(x, y)$ is a $g$-invariant properly convex open domain containing $x$ in its boundary. 

Suppose first that $h'(x, y)$ has a limit point $z$ of 
$g^{-n}(u)$ for some point $u \in h'(x, y) \cap U$. 
 
Here, $y$ is the fixed point of the smallest eigenvalue of the linear part of $g$.
Then $y$ is the attracting fixed point of 
a component of $\mathds{A}^n - h'(x, y)$ containing $U$ for $g^{-1}$. 
The antipodal point $y_-$ is the attracting fixed 
point of the other component $\mathds{A}^n - h'(x, y)$. 
Take a ball $B$ in $U$ centered at $u$ in the convex set $U \cap h(x, y)$. 
Then $\{g^{-n}(u)\}$ converges to $z$ as $n\ra \infty$. 
Let $u_1, u_2$ be two nearby points in $B$ such that 
$\ovl{u_1u_2}$ is separated by $h'(x, y)$ and $\ovl{u_1u_2}\cap h'(x,y) = u$.  
Then $\{g^{-n_i}(\ovl{u_1u_2})\}$ geometrically coverges to 
$\ovl{yz}\cup \ovl{y_-z}$ for some sequence $n_i$. 
Hence $\clo(U)$ cannot be properly convex. 

If the above assumption does not hold, then 
an orbit $g^{-n}(u)$ for $u \in U\cap h'(x, y)$ 
has a limit point only in the boundary of $h'(x, y)$. 
Since $g$ is bi-proximal, $x$ is the repelling fixed point 
of $h'(x, y)$ under $g^{-1}$. Hence, a limit point $y'$ 
cannot be either $x$ or $x_-$. 

Since $y'$ is a limit point, 
$y' \in \clo(U)$. It follows that $y' \in \clo(\Omega)$. 
Now, $x, y' \in \clo(\Omega)$ implies
$\ovl{x y'} \subset \Bd \mathds{A}^n \subset \clo(\Omega)$.
Finally, $\ovl{x y'} \subset \partial h'(x, y)$ for 
the sharply supporting subspace $\partial h'(x, y)$ of $\clo(\Omega)$
violates the strict convexity of $\Omega$.
(See Definition \ref{intro-defn-strict} and Benoist \cite{Benoist01}.)


\end{proof} 


Lemma \ref{du-lem-inde} will be generalized to Lemma \ref{du-lem-inde2} with a proof generalized word-for-word.  

\begin{lemma} \label{du-lem-inde}
Let $U$ be a $\bGamma_{\tilde E}$-invariant 
properly convex open domain in $\bR^n$ such that 
$\Bd U \cap \Bd \mathds{A}^n = \clo(\Omega)$. 
	Let $\Gamma$ act on strictly convex domain $\Omega$ 
	with $\Bd \Omega$ being $C^1$ 
	in a cocompact 
	manner. 
	Let $(x, y) \in \partial\clo(\Omega)\times \partial(\clo(\Omega) - \Delta$. Then 
\begin{itemize}
\item $\tau(x, y)$ is independent of $y$ and is uniquely determined by each $x$. 
\item $h'(x, y):= h(x,y)\cap \mathds{A}^n$ 
contains $\bar s(\overline{xy})$, but is independent of $y$
and $h(x, y) = h(x)$. 
\item The map $\tau': \partial \clo(\Omega) \ra S(\partial \clo(\Omega))$ induced from $\tau$ is continuous. 
	\item For an open set $U'$ such that $(\Gamma, U', \Omega)$ is a properly convex triple,
		the  AS-hyperspace at each point of $\partial \clo(\Omega)$ exists and 
		coincides that of $U$. 
\end{itemize}
\end{lemma}
\begin{proof} 
Let $l_1$ be an augmented geodesic in $\Uu\Omega$ 
	with endpoints $x$ and $z$, oriented towards $x$. 
	Consider a connected subspace $\mathcal{L}_x$ of 
	$\Uu \Omega$ consisting of points of maximal 
	augmented geodesics in $\Omega$ ending at $x$. 
The space of geodesic leaves in 
	$\mathcal{L}_x$ is in one-to-one correspondence with 
	$\partial \clo(\Omega) -\{x\}$.
	We will show that $\tau$ is locally constant on $\mathcal{L}_x$,  
	which shows that it is constant. 
	
	Let $\tilde l_1$ denote the lift of $l_1$ in $\Uu \Omega$. 
	Let $S$ be a compact neighborhood 
	in $\mathcal{L}_x$ of 
	a point $y$ on $\tilde l_1$ transverse to $\tilde l_1$. 
	Any two rays of geodesic flow 
	$\Phi: S \times \bR \ra \Uu \Omega$ are asymptotic 
	on $\mathcal{L}_{(x, h_1)}$ by Lemma \ref{intro-lem-decrease}. 
	

Let $y \in \tilde l_1$. 
Consider another point $y' \in S \subset \Uu \Omega$ with 
endpoints $x$ and $z'$ 
where  \[(x, z), (x, z') \in \Lambda^{\ast}.\] 
	
	
	Choose a fixed fundamental domain $F$ of $\hat\Uu \Omega$. 
	Let $\{y_i = \Phi_{t_i}(y)\}, y_i \in \tilde l_1,$ be a sequence whose projection under 
	$\Pi_\Omega$ converges to $x$.
	We use a deck transformation $g_i$ such that 
	$g_i(y_i) \in F$. 
	Then $g_i(\tau(l_1)) = \tau(g_i(l_1))$ is a hyperspace
	containing $g_i(x)$ and $\hat s(g_i(\tilde l_1))$. 
	
	
	Let $\bv_+$ denote a vector in the direction of the end of 
	$l_1$ other than $x$. 
	Equation \eqref{du-eqn-Anosov1} shows that 
	$\{\llrrV{\mathcal{L}(\Phi_t)| \bV_{+}}_{\mathrm{fiber}}\} \ra \infty$ as $t \ra \infty$. 
	Since $g_i$ is an isometry under $\llrrV{\cdot}_{\mathrm{fiber}}$, 
	and $\hat \Phi_{t_i}(y) = y_i$ and $g_i(y_i) \in F$, 
	it follows that 
	the $\tilde \bV_+$-component of $g_i(y_i, \bv_+)$ satisfies 
	\begin{equation}\label{du-eqn-LinearP} 	
	\{\llrrV{\mathcal{L}(g_i)(\bv_+)}_{\mathrm{fiber}}\} \ra \infty.
	\end{equation} 
		Since $g_i(y_i) \in F$ and 
	the metrics are compatible with a uniform constant
	under the Euclidean norm over the compact set $F$, 
	we obtain
	$\{\llrrV{\mathcal{L}(g_i) (\bv_{+} )}_E\} \ra \infty$. 
	
	Since the affine hyperspace in
	$\tau(x,z)$ and $\tau(x, z')$ 
	contain $x$ in their boundary sets, 
	they restrict to parallel affine hyperspaces in $\mathds{A}$. 
	Suppose that the affine hyperspace part of $\tau(x,z)$ 
	differs from one of  $\tau(x, z')$ 
	by a translation by $\bv_{+}$ multiplied by a constant. 
	This implies that the sequence of the Euclidean distances between 
	the respective affine hyperspaces corresponding to 
	\[g_i(\tau(x, z)) \hbox{ and } g_i(\tau(x, z'))\]
	goes to infinity as $i \ra \infty$. 
	
	Now we consider that $\Phi(S\times [t_i, t_i+1]) \subset \Uu\Omega$, 
	and we have obtained $g_i$ such that 
	$g_i(\Phi(S\times [t_i, t_i+1]))$ is in a fixed compact subset 
	$\hat P$ of $\Uu\Omega$ by the uniform boundedness of 
	$\Phi(S\times [t_i, t_i+1]))$ shown in the second paragraph of this proof. 
	There is a map $E: \Uu \Omega \ra \Lambda^{\ast}$
	given by sending the vector in $\Uu \Omega$ to the ordered pair of endpoints and supporting hyperspaces of the geodesic passing the vector.
	Since $\hat s$ is continuous, $\tau \circ E| \hat P$ is uniformly bounded. 
	The previous paragraph shows that the sequence of the diameters of 
	$\tau\circ E|g_i(\Phi(S\times [t_i, t_i+1]))$ can become arbitrarily large.
	This is a contradiction. Hence, $\tau$ is constant on 
	$\mathcal{L}_{x}$. 
	
	This proves the first two items. 
	The third item follows since $\tau^{\prime}$ is an induced map. 

	Define $H(x)$ to be the open $n$-dimensional hemisphere in $\SI^n$ bounded by 
	the great sphere containing the affine hyperspace $\tau^{\prime}(x)$
	and containing $\Omega$. Let 
	\[ \hat U:= \bigcap_{x \in \partial \clo(\Omega)} H(x) \cap \mathds{A}^n. \]

We claim that for any $x, y $ in $\partial \clo(\Omega)$, $h'(x, y)$ is disjoint from $U$: 
By Theorem 1.1 of Benoist \cite{Benoist01}, the geodesic flow on $\Omega/\Gamma$ is Anosov, and hence the set of 
closed geodesics in $\Omega/\Gamma$ is dense in the space of geodesics
by the basic property of the Anosov flow. 
Since the fixed points lie in $\partial \clo(\Omega)$, 
we can find sequences  $\{x_i\} \ra x$ and 
$\{y_i\} \ra y$ where 
$x_i$ and $y_i$ are fixed points of an element $g_i \in \Gamma$ for each $i$. 
If $h'(x, y) \cap U \ne \emp$, then $h'(x_i, y_i) \cap U \ne \emp$ for 
sufficiently large $i$ by the continuity of the map $\tau$
from Proposition \ref{du-prop-mapgh}. 
This is a contradiction by Lemma \ref{du-lem-hdisj} 

In particular, this also holds for $\hat U$ and $U'$ in the premise.

	
	Now, we show that the affine hyperspace part of 
	$\tau(x)$ is an AS-hyperspace for $U$:  
	Suppose that for $x\in \partial \clo(\Omega)$, 
	the AS-hyperspace $Q$ satisfies  
	$Q \ne \tau(x)$.  
	Then the hemisphere $H_Q$ bounded by $Q$ contains $U$.
	By the disjointness of $h(x, y)$ to $U$ stated above,  
we have $H_Q \subset H(x)$. 
	Again, we choose a segment $l_1$ ending at $x$. 
	Then we choose sequences $g_i$ as above in the proof 
	before \eqref{du-eqn-LinearP}. This shows, as above, 
	the sequence of the Euclidean distances between the respective affine 
	hyperspace parts of 
	\[ g_i(\tau(x)) \hbox{ and } g_i(Q)    \]
	goes to infinity. 
		Proposition \ref{du-prop-mapgh2} shows that 
	$\{g_i(\tau(x))\}$  
	lies in the image of $\tau$, which is a compact set. 
	The set of supporting hyperspace of $U$ is 
	bounded away from $\Bd \mathds{A}^n$ since 
	they must lie between those of the image of 
	$\tau$ and $U$.
	Since $g_i(Q)$ is still a supporting 
	hyperspace of $U$, 
	this leads to a contradiction. 

Now, we can replace $U$ by $U'$. 
The existence of an AS-hyperspace at each point follows  
from the argument above that $h'(x, y)$ is disjoint from $U'$. 
The fourth item follows from the paragraph above also.

\end{proof}

\begin{proof}[{\sl Proof of Theorem \ref{du-thm-asymnice}}] 
Let	\[ \hat U:= \bigcap_{x \in \partial \clo(\Omega)} H(x) \cap \mathds{A}^n. \]
Then this follows from Lemma \ref{du-lem-inde}. 
\end{proof} 


\section{Generalization to nonstrictly convex domains} \label{du-sec-gen}

\subsection{Main argument}
Now, we drop the condition of hyperbolicity on $\Gamma$. 
Hence, $\Omega \subset \Bd \mathds{A}^n$
is not necessarily strictly convex.  
Also, $\Omega$ is allowed to be the interior of 
a strict join.
Here, we do not assume that $\Gamma$ 
is not necessarily hyperbolic, and 
%
hence, it is more general. 
Also, we obtain an asymptotically nice properly convex domain $U$ in $\mathds{A}^n$ 
where $\Gamma$ acts properly on.

\begin{theorem}\label{du-thm-asymniceII}
	Let $\Gamma$ have an affine action on the affine subspace 
	$\mathds{A}^n$, $\mathds{A}^n \subset \SI^n$,
	acting on a properly convex domain $\Omega$ in 
	$ \Bd \mathds{A}^n$. 
	Suppose that $\Omega/\Gamma$ is a closed $(n-1)$-dimensional orbifold, and 
that $\Gamma$ satisfies the uniform middle-eigenvalue condition. 
	Then $\Gamma$ acts on a properly convex 
	open domain $U$ with the following properties\/{\rm :} 
	\begin{itemize}
		\item $(\Gamma, U, \Omega)$ is a properly convex triple, and 
 $\Gamma$ is asymptotically nice with the properly convex open domain $U$,
		and 
		\item if any open set $U'$ such that $(\Gamma, U', \Omega)$ is a properly convex triple, then 
		the  AS-hyperspace at each pair consisting of a point $x$ of $\partial \clo(\Omega)$ and a strictly supporting hyperspace of $\Omega$ in $\Bd \mathds{A}^n$ at $x$ exists and
		is the same as that of $U$.  That is, $U'$ is also asymptotically nice. 
	\end{itemize} 
\end{theorem}

	The proof is analogous to that of Theorem \ref{du-thm-asymnice}.
	Now, $\Omega$ is not strictly convex and hence for each 
	point of $\partial \clo(\Omega)$, there might be more than one 
	sharply supporting hyperspace in $\Bd \mathds{A}^n$.
	We generalize $\Uu \Omega$ to the {\em augmented unit tangent bundle} 
	\begin{multline} 
	\Uu^\Ag \Omega := \{(\vec x, H_a, H_b)|   
     \hbox{ $\vec x \in \Uu \Omega$ is a direction vector at a point }\\ 
     \hbox{ of a maximal oriented geodesic $l_{\vec x}$ in $\Omega$,} \\
	\hbox{ $H_a$ is a sharply supporting hyperspace in 
		 $\Bd \mathds{A}^n$ at the starting point 
		of $l_{\vec x}$,} \\
	\hbox{ $H_b$ 
		is a sharply supporting hyperspace in $\Bd \mathds{A}^n$ at the ending point of 
		$l_{\vec x}$}
	 \}.
	 \end{multline} 
	 Here, we regard $\vec x$ as a based vector, meaning that it has the information on 
	 where its base point is on $l$,
	 and $H_a$ and $H_b$ are given orientations such that $\Omega$ is in the 
	 interior direction to them. 
	 The space is not a manifold but a locally compact Hausdorff space and is metrizable. 
	 Since the set of sharply supporting hyperspaces of $\Omega$ at a point 
	 of $\partial \clo(\Omega)$ is compact, 
	 $\Uu^\Ag \Omega/\Gamma$ is a compact Hausdorff space
	 fibering over $\Omega/\Gamma$ with compact fibers. 	
	 The obvious metric is induced from $\Omega$ and 
	 the space $\SI^{n\ast}$ of oriented hyperspaces in $\SI^n$. 
	 We also write $\Pi^\Ag: \Uu^\Ag \Omega \ra \Omega$ the obvious projection
	 $(\vec{x}, H_a, H_b) = \Pi_{\Omega}(\vec{x})$.
	 \index{augmented unit tangent bundle|textbf} 
	 	
	
	From Section \ref{prelim-sec-duality}, we recall the augmented boundary 
	$\partial^\Ag \clo(\Omega)$. 
	We define 
	\[\Lambda^{\ast \Ag}  = \partial^\Ag \clo(\Omega) \times \partial^\Ag \clo(\Omega) 
	- (\Pi_{\Ag}\times\Pi_{\Ag})^{-1}(\Delta^\Ag)\]
	where $\Delta^\Ag$ is defined as the closed subset
	\[\{(x,y)| x,y \in \partial \clo(\Omega), x =y \hbox{ or } \overline{xy} \subset \partial \clo(\Omega) \}.\] 
	Let $G^\Ag(\Omega)$ denote the set of oriented 
	maximal geodesics in $\Omega$ where 
	endpoint is augmented with the sharply supporting hyperspace. 
	The elements are called {\em augmented geodesics}.
	There is a one-to-one and onto correspondence between 
	$\Lambda^{\ast \Ag}$ and $G^\Ag(\Omega)$. 
	We denote by
	$\overline{(x, h_1)(y, h_2)}$ the complete geodesic in 
	$\Omega$ with endpoints $x, y$ and sharply supporting hyperspaces $h_1$ at $x$ 
	and $h_2$ at $y$. 
	\index{lambda@$\Lambda^{\ast \Ag}$|textbf}
	\index{augmented geodesics|textbf} 
	
	\begin{lemma} \label{du-lem-hatGconnect} 
		$G^\Ag(\Omega)$ is a connected subspace.
	\end{lemma} 
	\begin{proof} 
		We can use the fact that 
		$\Uu^\Ag \Omega$ is connected. 
	\end{proof}


Now, we follow Section \ref{du-sub-flow} and define 
$\tilde \bA = \Uu^\Ag \Omega \times \mathds{A}^n$, 
$\tilde \bV =  \Uu^\Ag \Omega \times V^n$, 
$\bA$ by $\Uu^\Ag \Omega \times \mathds{A}^n/\Gamma $ 
and $\bV := \Uu^\Ag \Omega \times V^n/\Gamma$
and corresponding subbundles
 $\tilde \bV_+, \tilde \bV_-, \tilde \bV_0$, $\bV_+$, $\bV_-$, and $\bV_0$. 
We define 
the flows $\hat \Phi_t, \Phi_t, \tilde \Phi_t, \mathcal{L}(\Phi_t), \mathcal{L}(\tilde \Phi_t)$ by replacing $\Uu \Omega$ with $\Uu^\Ag \Omega$ and geodesics by augmented geodesics 
and so on in a natural way. 
	
	For each point $\boldx = (\vec{x}, H_a, H_b)$ of $\Uu^\Ag \Omega$,
	\begin{itemize}  
	\item we define $\tilde \bV_+(\boldx)$ to be the space of vectors in the direction of the backward endpoint of $l_{\boldx}$, 
	\item $\tilde \bV_-(\boldx)$ that for the forward endpoint of $l_{\boldx}$, 
\item $\tilde \bV_0(\boldx)$ to be the space of vectors in 
directions of $H_a\cap H_b$.
\end{itemize}

	For each $\boldx \in \Uu^\Ag \Omega$, 
	\[V^n = \tilde \bV_+(\boldx)\oplus \tilde \bV_0(\boldx) \oplus 
	\tilde \bV_-(\boldx).\] 
	This gives us a decomposition. 
	$\tilde \bV = \tilde \bV_+ \oplus \tilde \bV_0 \oplus \tilde \bV_-$, 
	and 	$\bV =  \bV_+ \oplus  \bV_0 \oplus \bV_-$.
	Clearly, $\bV_+$ and $\bV_-$ are topological line bundles
	since the beginning and the endpoints depend continuously on 
	points of $\Uu^\Ag \Omega$. 
	Also, $\tilde \bV_0$ is a subspace of $\bR^n$ whose 
	directions of nonzero vectors lies in $H_a \cap H_b$. 
	Since $(H_a, H_b)$ depends continuously on points of 
	$\Uu^\Ag \Omega$, we obtain that 
	$\bV_0$ is a continuous bundle over $\Uu^\Ag \Omega$. 
	
	Obviously, the geodesic flow exists on 
	$\Uu^\Ag \Omega$ using the ordinary geodesic flow with respect to 
	the geodesics and not considering the augmented boundary.
	
	There exist constants $C, k > 0$ such that 
	\begin{align} \label{du-eqn-Anosov3} 
	\llrrV{\mathcal{L}(\Phi_t) ({\vec{v}})}_{\mathrm{fiber}} \geq \frac{1}{C} \exp(kt) \llrrV{{\vec{v}}}_{\mathrm{fiber}} \hbox{ as } t \ra \infty
	\end{align}
	for $\vec{v} \in \bV_+$ and 
	\begin{align} \label{du-eqn-Anosov3ii} 
	\llrrV{\mathcal{L}(\Phi_t) (\vec{v})}_{\mathrm{fiber}} \leq C \exp(-kt) \llrrV{\vec{v}}_{\mathrm{fiber}} \hbox{ as } t \ra \infty
	\end{align}
	for $\vec{v} \in \bV_-$. 
	
	We prove this by showing
	$\{\llrrV{\mathcal{L}(\Phi_{t})| \bV_-}_{\mathrm{fiber}}\} \ra 0 $ uniformly as $t \ra \infty$
	i.e., Proposition \ref{du-prop-contr} 
	under the more general conditions that $\Omega$ is properly convex 
	but not necessarily strictly convex.  
	We generalize Lemma \ref{du-lem-S_0}.
	We will repeat the strategy of the proof since it is important to check. 
	However, the proof follows the same philosophy 
	with some technical differences.

	\begin{lemma} \label{du-lem-S_0II}
		Assume that $\Omega$ is properly convex 
		and $\Gamma$ acts properly discontinuously satisfying the uniform 
		middle-eigenvalue condition with respect to $\Bd \mathds{A}^n$. Then 
	$\{\llrrV{\Phi_t| \bV_-}_{\mathrm{fiber}}\} \ra 0$ uniformly as $t \ra \infty$. 
	\end{lemma} 
	\begin{proof} 
			
		We proceed as in the proof of Lemma \ref{du-lem-S_0}.
	It sufficies to prove the uniform convergence to $0$ by the compactness of  $\bV_{-,1}$. 
		Let $F$ be a fundamental domain of $\Uu^\Ag\Omega$ under $\Gamma$.
	It is sufficient to prove this for $\mathcal{L}(\tilde{\Phi}_{t})$ on the fibers of over 
	$F$ of $\Uu^\Ag\Omega$ with a fiberwise metric $\llrrV{\cdot}_{\mathrm{fiber}}$.
	
	We choose an arbitrary sequence $\{(\vec{x}_i, H_{a_i}, H_{b_i}) \in \Uu^\Ag \Omega\}$, 
	$\{\boldx_i\} \ra \boldx$ in $F$ where $a_i, b_i$ are the backward and forward point of 
	the maximal oriented geodesic passing through $\vec{x}_i$ in $\Omega$. 
	For each $i$, let $\bv_{-, i}$ be a Euclidean unit vector in  $\tilde \bV_{-}(\boldx_{i})$ for the unit vector $\vec{x}_i \in \Uu^\Ag\Omega$.
	That is, $\bv_{-, i}$ is in the $1$-dimensional subspace in $\bR^{n}$,  
	corresponding to the forward
	endpoint of the geodesic determined by $\vec{x}_i$ in $\partial \clo(\Omega)$. 
	
	Here, $\llrrparen{\bv_{-, i}}$ is an endpoint of $l_i$ in the direction given by $\boldx_i$.
	For this, we just need to show that any sequence of $\{t_i\} \ra \infty$ has a subsequence $\{t_j\}$ such that 
	$\left\{\llrrV{\mathcal{L}(\tilde{\Phi}_{t_j})
		(\boldx_{i}, \bv_{-,  j})}_{\mathrm{fiber}}\right\} \ra 0$.
	This follows since if the uniform convergence did not hold, then we can easily find a sequence without such subsequences. 
	
	Define $y_i := \widehat{\tilde \Phi}_{t_i}(\boldx_i)$ for the lift $\widehat{\tilde \Phi}$ 
of the flow $\hat \Phi$. 
	By construction, each $\Pi^\Ag_\Omega(y_i)$ lie on the geodesic $l_i$. 
	Since we have the sequence $\{\boldx_i\} \ra \boldx$, $\boldx_{i}, \boldx \in F$, 
	we obtain that $\{l_i$\} geometrically converges to a line $l_\infty$ passing $\Pi^\Ag_\Omega(\boldx)$ in $\Omega$. 
	Let $y_+$ and $y_-$ be the endpoints of $l_\infty$ where 
	$\{\Pi^\Ag_\Omega(y_i)\} \ra y_-$.
	Hence, \[\{\llrrparen{\bv_{+, i}}\} \ra y_{+}, \{\llrrparen{\bv_{-, i}}\}\ra y_{-}.\] 
	(See Figure \ref{du-fig-contr} for a similar situation.)
	
	Find a deck transformation $g_i$ such that 
	$g_i(y_i) \in F$, and $g_i$ acts on the line bundle $\tilde \bV_-$ by the linearization of the matrix of the form \eqref{du-eqn-bendingm4}: 
	\begin{align} 
	g_{i} &: \tilde \bV_{-} \ra \tilde \bV_{-} \hbox{ given by } \nonumber \\
	(y_{i}, \bv) & \ra (g_{i}(y_{i}), {\mathcal L}(g_{i})(\bv)) 
	\hbox{ where } \nonumber \\
	{\mathcal L}(g_i) &:= 
	\frac{1}{\lambda_{{\tilde E}}(g_i)^{1+\frac{1}{n}}} \hat h(g_i) : \tilde \bV_{-}(y_i) = \tilde \bV_{-}(x_{i}) \ra \tilde \bV_{-}(g_i(y_i)). \label{du-eqn-LgiII} 
	\end{align}

	 We will show $\{(g_{i}(y_{i}), {\mathcal L}(g_i)(\bv_{-, i}))\} \ra 0$ under $\llrrV{\cdot}_{\mathrm{fiber}}$. 
		This will complete the proof since $g_{i}$ acts as isometries on $\tilde \bV_{-}$ with $\llrrV{\cdot}_{\mathrm{fiber}}$. 

	
	\begin{figure}
		\centering
		\includegraphics[height=8cm]{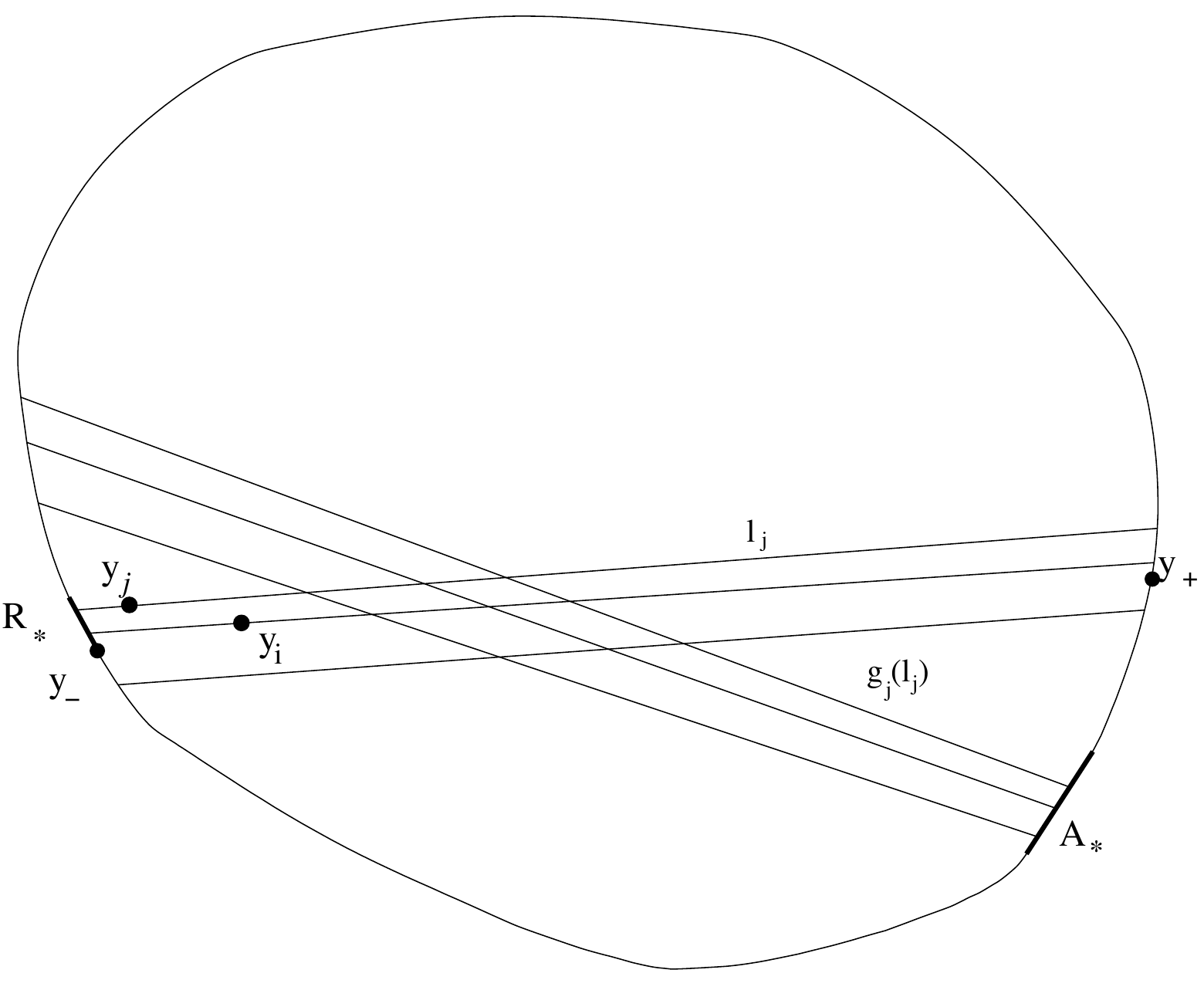}
		\caption{The figure for Lemma \ref{du-lem-S_0II}.}
		
		\label{du-fig-contrii}
		
	\end{figure}

		Since $g_i(y_i) \in F$ for every $i$, we obtain 
	\[g_i(l_i)\cap F\ne \emp.\]
	Since $g_i(l_i) \cap F \ne \emp$, we may pass to 
a subsequence of $g_i$ and still denoted as $g_i$ such that
	$\{\Pi^\Ag_\Omega(g_i(l_i))\}$ converges to a nontrivial line 
	$\hat l_\infty$ in $\Omega$.
	
	Remark \ref{pr-rem-virtual} 
	shows that we may assume without loss of generality that each element of $\Gamma$ is positive bi-semiproximal.
We recall facts from Section \ref{prelim-sub-convG}. 
Given a generalized 
convergence sequence $\hat g_i:= \mathcal{L}(g_i)$, we obtain an endomorphism $\hat g_\infty$ in 
$M_{n}(\bR)$ such that $\{\llrrparen{\mathcal{L}(g_i)}\} \ra \llrrparen{\hat g_\infty} \in \SI(M_{n}(\bR))$. 
Recall
\[A_\ast(\{\hat g_i\}, N_\ast(\{\hat g_i\}), R_\ast(\{\hat g_i\}), F_\ast(\{\hat g_i\}),\] 
which are subsets of $\partial \clo(\Omega)$ by Theorem \ref{prelim-thm-AR}. 

	Up to passing to a subsequence, 
ratio arguments as in the proof of 
Theorem  \ref{prelim-thm-converg}
	imply
	that for any compact subset $K$ of $\clo(\Omega) - R_\ast(\{\hat g_i\})$,
	there is a convex compact subset $K_\ast$ in $F_\ast(\{\hat g_i\})$,
	\begin{equation}\label{du-eqn-giK3} 
	\{\hat g_i(K)\} \ra K_\ast \subset F_\ast(\{g_i\}). 
	\end{equation} 

	Suppose that $y_- \in \clo(\Omega) - R_\ast(\{\hat g_i\})$. 
	Then $\{\hat g_i(y_{-,i})\}\ra \hat y \in F_\ast(\{\hat g_i\})$ since $y_{-.i}$ are in a compact subset of 	$\clo(\Omega) - R_\ast(\{\hat g_i\})$ and \eqref{du-eqn-giK3}. 
This means that $\overline{y_i y_{-, i}} \subset l_i$ is eventually of some $\bdd$-distance from $F_i(\hat g_i)$, and hence
$\hat g_i(\overline{y_i y_{-, i}})$ geometrically converges to a segment or a singleton in 
$F_\ast(\{\hat g_i\})$ disjoint from $\Omega$, which contradicts $\hat g_i(y_i) \in F$.
	Hence, $y_- \in R_\ast(\{\hat g_i\})$. 
Furthermore,  \[a_{n}(\hat g_i) \ra 0\]
since otherwise $\hat g_i$ would be bounded. (Recall Section \ref{prelim-sub-convG}.)
	
	Let $d_i=\llrrparen{\bv_{+,i}}$ denote the endpoint of $l_i$ distinct from
	$\llrrparen{\bv_{-,i}}$ as above. Let $d_\infty$ denote the limit of $d_i$ 
	in $\partial \clo(\Omega)$. 
	We deduce as above up to a choice of a subsequence: 
	\begin{itemize} 
		\item Since $\{\llrrparen{\bv_{-, i}}\}\ra y_-$, $y_- \in R_\ast(\{\hat g_i\})$ 
		and since $\{l_i\}$ converges to a nontrivial line $l_\infty \subset \Omega$
		and $R_\ast(\{\hat g_i\})$ is a convex compact subset of $\partial \clo(\Omega)$, 
		it follows that 
		$\{d_i\}$ lies  in a compact set in $\partial \clo(\Omega) -R_\ast(\{\hat g_i\})$.
		\item Then $\{\hat g_i(d_i)\} \ra a_{\ast} \in F_\ast(\{\hat g_i\})$ by \eqref{du-eqn-giK3}
		since $\{d_i\}$ is in a compact set in $\partial \clo(\Omega) -R_\ast(\{\hat g_i\})$.
		\item Thus, $\{\hat g_i(\llrrparen{\bv_{-, i}})\} \ra y'_- \in \partial \clo(\Omega) - F_\ast({\hat g_i})$ holds since 
		$\{\hat g_i(l_i)\}$ converges to a nontrivial line in $\Omega$.
	\end{itemize} 

Let $m_r$ be obtained from $\{\hat g_i\}$ as in Theorem \ref{prelim-thm-convGU}.
	Recall that $\llrrV{\cdot}_E$ denotes the standard Euclidean metric of $\bR^{n}$. 
	We write $\hat g_i = k_i D_i \hat k_i^{-1}$ for 
	$k_i, \hat k_i \in \Ort(n, \bR)$, and $D_i$ is a positive diagonal matrix
	of determinant $\pm 1$ with nonincreasing diagonal entries. 
	\begin{itemize}
		\item Since $\{\llrrparen{\bv_{-, i}}\} \ra y_-$,
		the vector $\bv_{-, i}$ has the component $\bv_i^p$ parallel to 
		$R^p(g_i) = \hat k_i(\SI([m_r, n]))$ and the component $\bv_i^S$ in the orthogonal complement $\left(R^p(g_i)\right)^\perp = \hat k_i(\SI([1, m_r-1])) $ 
		where $\bv_{-, i} = \bv_{-, i}^p + \bv_{-.i}^S$. 
		We may require $\llrrV{\bv_{-,i}}_E = 1$.
		(We remark 
		\[\{\llrrV{\bv_i^p}_E\} \ra 1 \hbox{ and } \left\{\llrrV{\bv_i^S}_E\right\} \ra 0\]
		since $\{\bv_{-, i}\}$ converges to a point of 
		$R_\ast(\{g_i\})$. 
	\item $\left\{\llrrV{{\mathcal L}(g_i)(\bv_{-, i}^p)}_E\right\} \ra 0$ 
	by Theorem \ref{prelim-thm-convGU} 
	since $a_{n}(\mathcal{L}(g_i))  \ra 0$. 
	\item Since $\{g_i(\llrrparen{\bv_{-, i}})\}$ converges to $y'_-$, 
	$y'_- \in \partial \clo(\Omega) - F_\ast({g_i})$,  
	$\{g_i(\llrrparen{\bv_{-, i}})\}$ is uniformly bounded away from $F_\ast(\{g_i\})$. 
	\end{itemize} 
	
	Because of the orthogonal decomposition $\hat k_i(\SI([m_r, n])$ and 
	$\hat k_i(\SI([1, m_r-1]))$, and the fact that 
	$\hat g_i = k_i D_i \hat k_i^{-1}$, 
	and $\{\llrrparen{\mathcal{L}(g_i)}\} \ra \llrrparen{\hat g_\infty}$ in $\SI(M_n(\bR))$, 
	it follows that 
	either $\{\mathcal{L}(g_i)(\bv^S_i)\}$  converges to zero, 
	or $\{\llrrparen{\mathcal{L}(g_i)(\bv^S_i)}\}$ converges to $F_\ast(\{g_i\})$ by 
using ratios of elements of $D_i$ 
as in the proof of Theorem \ref{prelim-thm-convGU}.  
	We have
	\[\{{\mathcal L}(g_i)(\bv_i^p)\} \ra 0 \hbox{ implies }
	 \{{\mathcal L}(g_i)(\bv_i^S)\} \ra 0\]
	for $\llrrV{\cdot}_{E}$ since otherwise 
	$\{\llrrparen{{\mathcal L}(g_i)(\bv_{-, i})}\}$ converges to a point of 
	$F_\ast(\{g_i\})$ and hence 
	$\{g_i(\llrrparen{\bv_{-, i}})\}$ cannot be converging to $y'$.
	
	Hence, we obtain $\{{\mathcal L}(g_i)(\bv_{-, i})\} \ra 0$ under $\llrrV{\cdot}_E$. Now, we can deduce the result as in 
	the final part of the proof of Lemma \ref{du-lem-S_0}

	\end{proof}


As in Section \ref{du-sub-neutral}, we find  the neutralized section 
	$s:\Uu^\Ag\Sigma \ra \bA$ such that $\nabla^{\bA}_\phi s \in \bV_0$. 
		
		Since we are considering $\Uu^\Ag \Omega$,
		the section $s: \Uu^\Ag \Omega \ra N_2(\mathds{A}^n)$, 
		we need to look at the boundary point and 
		a sharply supporting hyperspace at the point and find the affine 
		subspace of dimension $n-2$ in $\bR^n$,
		generalizing Proposition \ref{du-prop-mapgh}. 	
	We now generalize Definition \ref{du-defn-tau}: 
	\begin{definition}\label{du-defn-tauAg} 
		We denote by $S^\Ag(\partial^\Ag \clo(\Omega))$ the space of pairs 
		$((x, H \cap \Bd \mathds{A}^n), H)$ 
		where $H \in S'(\partial \clo(\Omega))$, and $x$ is in the boundary of 
		$H$ and $(x, H \cap \Bd \mathds{A}^n) \in \partial^\Ag \clo(\Omega)$. 
	\end{definition}

\begin{proposition}\label{du-prop-mapgh2}  $ $
	\begin{itemize}
		\item There exists  a continuous function 
$\hat s: \Uu^\Ag\Omega \ra N_2(\mathds{A}^n)$ that is 
		equivariant with respect to $\Gamma$-actions.
		\item For each $g \in \Gamma$, let $\vec{l}_g$ denote 
		the unique unit-speed geodesic in 
		$\Uu^\Ag\Omega$ lying over an augmented geodesic $l_g$ where $g$ acts on, 
		$\hat s(\vec{l}_g) = \{N_2(g)\}$. 
		\item This gives a continuous map 
		\[\bar s^\Ag: \partial^\Ag \clo(\Omega) \times \partial^\Ag \clo(\Omega) - 
		(\Pi_{\Ag}\times \Pi_{\Ag})^{-1}(\Delta^\Ag) \ra N_2(\mathds{A}^n)\] 
		again equivariant under the $\Gamma$-actions. 
		There exists a continuous function 
		\[\tau^\Ag:\Lambda^{\ast \Ag} \ra S^\Ag(\partial \clo(\Omega)).\]
	\end{itemize}
\end{proposition}
\begin{proof}
	The proof is entirely similar to that of Proposition \ref{du-prop-mapgh}
	but using a straightforward generalization of Lemma \ref{du-lem-exist}. 
\end{proof}

We generalize Proposition \ref{du-lem-inde}. 
We define $\tau': \Uu^\Ag \Omega \ra  S^\Ag(\partial \clo(\Omega))$ as 
a composition of $\tau^\Ag$ with the map from $\Uu^\Ag \Omega$ 
to $\Lambda^{\ast \Ag}$. This map is continuous. 
Here, we don't assume that $\Gamma$ acts on a properly convex domain 
in $\mathds{A}^n$ with boundary $\Omega$. Hence, the situation is more general,
and we need a different proof. We just need that the orbit closures are compact. 

\begin{lemma} \label{du-lem-inde2}
	Let an affine group $\Gamma$ act on an affine subspace $\mathds{A}^n$
	on a properly convex domain $\Omega$ in the boundary 
	of an affine subspace $\mathds{A}^n$. 
	Let $\Gamma$ act on a properly convex domain $\Omega$ 
	a cocompactly and properly 
	 and satisfies the uniform middle-eigenvalue 
	condition with respect to $\Bd \mathds{A}^n$. 
	Let $((x, h_1), (y, h_2)) \in \Lambda^{\ast \Ag}$. Then 
	\begin{itemize}
		\item $\tau^\Ag((x, h_1), (y, h_2))$ does not depend on $(y, h_2)$ and is unique for each $(x, h_1)$. 
		\item $h( (x, h_1), (y, h_2))$ contains $\bar s^\Ag((x, h_1), (y, h_2))$ but is independent of $(y, h_2)$. 
		\item $h((x, h_1), (y, h_2))$ is never a hemisphere in $\Bd \mathds{A}^n$ for every 
		$((x, h_1), (y, h_2)) \in \Lambda^{\ast \Ag}$.
		\item $\tau^\Ag$ induces a map 
		$\tau^{\Ag \prime}: \partial^\Ag \clo(\Omega) \ra S^\Ag(\partial \clo(\Omega))$ that is continuous.  
			\item There exists an asymptotically nice convex $\Gamma$-invariant 
		open domain $U$ in $\mathds{A}^n$ with 
		$\Bd U \cap \partial \mathds{A}^n = \clo(\Omega)$. 
		For every $(x, h_1) \in \partial^\Ag \clo(\Omega)$, 
		 $\tau(x, h_1)$ is an AS-hyperspace of $U$.
	\end{itemize}
\end{lemma}
\begin{proof} 
	Let $l_1$ be an augmented geodesic in $\Uu^\Ag\Omega$ 
	with endpoints $(x, h_1)$ and $(z, h_2)$
	oriented towards $x$. 
	Consider a connected subspace $\mathcal{L}_{(x, h_1)}$ of 
	$\Uu^\Ag \Omega$ of points of maximal 
	augmented geodesics in $\Omega$ ending at $(x, h_1)$. 
The space of geodesic leaves in 
	$\mathcal{L}_{(x, h_1)}$ is in one-to-one correspondence with 
	$\Bd^{\Ag}\Omega - \Pi_{\Ag}^{-1}(K_x)$
	for the maximal flat $K_x$ in $\partial \clo(\Omega)$ containing $x$.  
	We will show that $\tau^\Ag$ is locally constant on $\mathcal{L}_{(x, h_1)}$ 
	which implies that it is constant. 
	
	Let $\tilde l_1$ denote the lift of $l_1$ in $\Uu^\Ag \Omega$. 
	Let $S$ be a compact neighborhood 
	in $\mathcal{L}_{(x, h_1)}$ of 
	a point $y$ of $\tilde l_1$ transverse to $\tilde l_1$. 
	Any two rays of geodesic flow 
	$\Phi: S \times \bR \ra \Uu^\Ag \Omega$ are asymptotic 
	on $\mathcal{L}_{(x, h_1)}$ by Lemma \ref{intro-lem-decrease}. 
	

Let $y \in \tilde l_1$. 
Consider another point $y' \in S \subset \Uu^\Ag \Omega$
with endpoints $x$ and $z'$ in 
a sharply supporting hyperspace $h'_2$
where  \[((x, h_1), (z, h_2)), ((x, h_1), (z', h'_2)) \in \Lambda^{\ast \Ag}.\] 
	
	
	Choose a fixed fundamental domain $F$ of $\hat\Uu \Omega$. 
	Let $\{y_i = \Phi_{t_i}(y)\}, y_i \in \tilde l_1,$ be a sequence whose projection under 
	$\Pi_\Omega$ converges to $x$.
	We use a deck transformation $g_i$ such that 
	$g_i(y_i) \in F$. 
	Then $g_i(\tau^\Ag(l_1)) = \tau^\Ag(g_i(l_1))$ is a hyperspace
	containing $g_i(x)$, $g_i(h_1)$, and $\hat s(g_i(\tilde l_1))$. 
	
	
	Let $\bv_+$ denote a vector in the direction of the end of 
	$l_1$ other than $x$. 
	Equation \eqref{du-eqn-Anosov3} shows that 
	$\{\llrrV{\mathcal{L}(\Phi_t)| \bV_{+}}_{\mathrm{fiber}}\} \ra \infty$ as $t \ra \infty$. 
	Since $g_i$ is an isometry under $\llrrV{\cdot}_{\mathrm{fiber}}$, 
	and $\hat \Phi_{t_i}(y) = y_i$ and $g_i(y_i) \in F$, 
	it follows that 
	the $\tilde \bV_+$-component of $g_i(y_i, \bv_+)$ satisfies 
	\begin{equation}\label{du-eqn-Linear} 	
	\{\llrrV{\mathcal{L}(g_i)(\bv_+)}_{\mathrm{fiber}}\} \ra \infty.
	\end{equation} 
		Since $g_i(y_i) \in F$ and 
	under the Euclidean norm since over a compact set $F$ 
	the metrics are compatible by a uniform constant, 
	we obtain that
	$\{\llrrV{\mathcal{L}(g_i) (\bv_{+} )}_E\} \ra \infty$. 
	
	Since the affine hyperspaces in
	$\tau^\Ag((x, h_1),(z, h_2))$ and $\tau^\Ag((x, h_1), (z', h'_2))$ 
	contain $x$ and $h_1$ in their boundary spheres, 
	they restrict to parallel affine hyperspaces in $\mathds{A}$. 
	Suppose that the affine hyperspace part of $\tau^\Ag((x, h_1),(z, h_2))$ 
	differs from one of  $\tau^\Ag((x, h_1), (z', h'_2))$ 
	by a translation by $\bv_{+}$ scaled by a constant. 
	This implies that the sequence of the Euclidean distances between 
	the respective affine hyperspaces corresponding to 
	\[g_i(\tau^\Ag((x, h_1),(z, h_2))) \hbox{ and } g_i(\tau^\Ag((x, h_1), 
	(z', h'_2)))\]
	goes to infinity as $i \ra \infty$. 
	
	Now consider $\Phi(S\times [t_i, t_i+1]) \subset \Uu^\Ag \Omega$. 
	We have obtained $g_i$ such that 
	$g_i(\Phi(S\times [t_i, t_i+1]))$ is in a fixed compact subset 
	$\hat P$ of $\Uu^\Ag \Omega$ by the uniform boundedness of 
	$\Phi(S\times [t_i, t_i+1])$ shown in the second paragraph of this proof. 
	Define a map $E: \Uu^\Ag \Omega \ra \Lambda^{\ast \Ag}$
	given by sending the vector in $\Uu^\Ag \Omega$ to the ordered pair of endpoints and supporting hyperspaces of the geodesic it lies on. 
	Since $\hat s$ is continuous, $\tau^\Ag\circ E| \hat P$ is uniformly bounded. 
	The above paragraph shows that the sequence of the diameters of 
	$\tau^\Ag\circ E|g_i(\Phi(S\times [t_i, t_i+1]))$ can become arbitrarily large.
	This is a contradiction. Hence, $\tau^\Ag$ is constant on 
	$\mathcal{L}_{(x, h_1)}$. 
	
	This proves the first two items. 
	The fourth item follows since $\tau^{\Ag \prime}$ is an induced map. 
	The image of $\tau^{\Ag \prime}$ is compact 
	since $\partial \clo(\Omega)$ is compact. This implies the third item. 

	Define $H(x, h_1)$ to be the open $n$-dimensional hemisphere in $\SI^n$ bounded by 
	the great sphere containing the affine hyperspace $\tau^{\Ag \prime}(x, h_1)$
	and also containing $\Omega$. 
	We define 
	\[ U:= \bigcap_{(x, h_1) \in \partial^\Ag \clo(\Omega)} H(x, h_1) \cap \mathds{A}^n. \]
	
	Now, we show that the affine hyperspace part of 
	$\tau^\Ag(x, h_1)$ is an AS-hyperspace for $U$:  
	Suppose that for $(x, h_1)\in \partial^\Ag \clo(\Omega)$, 
there exists 
	the AS-hyperspace $Q$ with $Q \cap \Bd \mathds{A}^n = h_1$
	and $Q \ne \tau^\Ag(x, h_1)$.  
	Then the hemisphere $H_Q$ bounded by $Q$ contains $U$.
	By definition, $H_Q \subset H(x, h_1)$. 
	Then again we choose a segment $l_1$ ending at $x$. 
	Then we choose sequences $g_i$ as above in the proof 
	before \eqref{du-eqn-Linear}. This shows that 
	the sequence of the Euclidean distances between the respective affine 
	hyperspace parts of 
	\[ g_i(\tau^\Ag(x, h_1)) \hbox{ and } g_i(Q)    \]
	goes to infinity. 
		Proposition \ref{du-prop-mapgh2} shows that 
	$\{g_i(\tau^\Ag(x, h_1))\}$  
	is in the image of $\tau^\Ag$, a compact set. 
	The set of suppporting hyperspaces of $U$ is 
	bounded away from $\Bd \mathds{A}^n$ since 
	they have to be between those of the image of 
	$\tau^\Ag$ and $U$.
	Since $g_i(Q)$ is still a supporting 
	hyperspace of $U$, 
	this leads to to a contradiction. 
	\end{proof}


\begin{proof}[Proof of Theorem \ref{du-thm-asymniceII}]
	First, we obtain a properly convex domain where $\Gamma$ acts.
	Since $\partial^\Ag \clo(\Omega)$ is compact, the image under $\tau^{\Ag}$ is compact. 
	Then $U:=\bigcap_{(x, h)\in \partial^\Ag \clo(\Omega)} H^o_{(x, h)}\cap \mathds{A}^n$ 
	contains $\Omega$.  This is an open set since the compact set of 
	$H_{(x, h)}$, $(x, h)\in \partial^\Ag \clo(\Omega)$ has 
	a lower bound on angles with $\Bd \mathds{A}^n$. 
	Thus, $U$ is asymptotically nice.  
Now, the proof is identical with that of Theorem  \ref{du-thm-asymnice} with 
Lemma \ref{du-lem-inde2} replacing Lemma \ref{du-lem-inde}. 

The uniqueness part is established  in Theorem \ref{du-thm-ASunique} immediately below. 
	\end{proof}

\subsection{Uniqueness of AS-hyperspaces}
Finally, we end with some uniqueness properties. 
The following generalizes Lemma \ref{du-lem-inde}.
\begin{theorem}\label{du-thm-ASunique}
	Let $(\Gamma, U, D)$ be a properly convex affine triple. 
	Suppose that $\Gamma$ satisfies the uniform middle-eigenvalue condition.
Then for any properly convex triple $(\Gamma, U', D)$ for an open set $U'$, 
 the set of AS-hyperspaces for $U'$ containing all 
	sharply supporting hyperspaces of $\Omega$ in $\Bd \mathds{A}^n$ exists and is 
	independent of the choice of $U'$.  That is, $U'$ is also asymptotically nice. 
\end{theorem}
\begin{proof} 
	For each $x\in \partial \clo(\Omega)$, 
	let $h$ be a supporting hyperspace of $\Omega$ in $\Bd \mathds{A}^n$, and 
	let $S_{(x, h)}$ be the AS-hyperspace in $\SI^n$ for $U$
	such that $S_{(x, h)} \cap \Bd \mathds{A}^n = h$. 
	Again for each $(x, h)$, we 
	let $S'_{(x, h)}$ is the AS-hyperspace in $\SI^n$ for $U'$ if it exists
	such that $S'_{(x, h)} \cap \Bd \mathds{A}^n = h$. 

If $(x, h)$ is a fixed point of a nontrivial $g\in \Gamma$, 
then the same argument as in the proof of Lemma \ref{du-lem-hdisj}
shows that $U'$ and $S_{(x, h)}$ are disjoint.  

Suppose that $\Gamma$ is not virtually factorizable when it 
restricted to $\Bd \mathds{A}^n$. 
(See Proposition \ref{prelim-prop-Ben2} .) 
Every extreme point of the dual domain $\Omega^\ast$ in $\SI^{n-1\ast}_\infty$ is the limit of 
a sequence of fixed points of $\Gamma^\ast$ by Proposition 5.1 of \cite{Zimmer23}. 
Since $\tau^\Ag$ is continuous, $U'$ and $S_{(x', h')}$ are disjoint if $h'$ is an extreme point of 
$\Omega^\ast$ by the density of fixed points. 
Since every boundary point of $\Omega^\ast$ is a nonnegative linear sum of extreme points of
$\Omega^\ast$, it follows that $U'$ and $S_{(x', h')}$ is disjoint for any pair $(x', h')$. 
This is because $S_{(x', h')}$ is also supporting $\Omega$ at $x'' \in \partial \clo(\Omega)$ 
that $h'$ contains. 
This proves that $U'$ is also asymptotically nice. 

Suppose that $\Gamma$ is virtually factorizable restricted to $\Bd \mathds{A}^n$. 
Then $\Omega$ is the interior of the strict join $K_1\ast \dots \ast K_m$. 
Then there are supporting hyperspaces $H_j$ to 
$\Omega$ in $\SI^{n-1}_\infty$ that contain all the factors except for $K_j$ 
for each $j =1, \dots, m$, 
and a supporting hyperspace to $K_j$ in the span of $K_j$. 
Then the virtual center, isomorphic to $\bZ^{m-1}$, acts on $H_j$. 
If $U'$ intersects with $H_j$, then the second paragraph above will again show the contraction. 
Since any supporting hyperspaces are positive linear combination of some $H_i'$, we have the 
disjointness of $U'$ with any $S_{(x', h')}$.


	Since $U$ and $U'$ are both asymptotically nice,
	the sets of AS-hyperspaces are compact. 
	For each $(x, h)$, $S_{(x, h)}$ and $S'_{(x, h)}$ differ by 
	a uniformly bounded distance in $\SI^{n\ast}$. 
	
	Suppose that $S_{(x, h_1)}$ is different from $S'_{(x, h_1)}$ 
	for some $x, h_1 \in \partial \clo(\Omega)$. 
	Now, we follow the argument from the proof of Lemma \ref{du-lem-inde2}. 
	We again obtain a sequence $g_i \in \Gamma$ such that 
	\[g_i(S_{(x, h_1)}\cap \mathds{A}^n) \hbox{ and }
	g_i(S'_{(x, h_1)}\cap \mathds{A}^n)\] 
	are parallel affine planes, 
	and the sequence of their Euclidean distances goes to $\infty$ as $i\ra \infty$. 
	By compactness, we know both sequences $\{g_i(S_{(x, h_1)}\cap \mathds{A}^n)\}$
    and $\{g_i(S'_{(x, h_1)}\cap \mathds{A}^n)\}$ respectively 
	converge to two hyperspaces up to a choice of a subsequence. 
	This means that their Euclidean distances are uniformly bounded. 
	Again \eqref{du-eqn-Linear} contradicts this. 
\end{proof} 

\begin{remark} 
	Theorems \ref{du-thm-asymnice} and \ref{du-thm-asymniceII} 
	also generalize to the case when 
	$\Gamma$ acts on $\Omega$ as a convex cocompact group: 
	i.e., there is a convex domain $C \subset \Omega$ 
	such that $C/\Gamma$ is compact.
	We work on the set of geodesics in $C$ only and the set $\Lambda$ 
	of endpoints of these.
	In this case, $\Lambda$ may be disconnected. 
	Definitions such as asymptotic niceness should be 
	applied to points of $\Lambda$ only. 
	Here we do need the connectedness of 
	Lemmas \ref{du-lem-Gconnect} and \ref{du-lem-hatGconnect} 
	to be generalized to this case. 
	However, the proofs indicated there will work. 

\end{remark} 

\section{Lens type T-ends} \label{du-sec-lensT-end}

\subsection{Existence of lens-neighborhood}

\begin{proposition}\label{du-prop-lensn} 
Let $\Gamma$ be a discrete group in $\SLnp$ {\rm (}resp. $\PGL(n+1, \bR)$\/{\rm )} 
acting on a properly convex domain $\Omega$
cocompactly and properly, 
$\Omega \subset \Bd \mathds{A}^n \subset \SI^n$ {\rm (}resp. $\subset \RP^n$\/{\rm ),}
such that $\Omega/\Gamma$ is 
a closed $n$-orbifold.
\begin{itemize}
\item Suppose that $\Gamma$ satisfies the uniform middle-eigenvalue condition
with respect to $\Bd \mathds{A}^n$.
\item Let $P$ be the hyperspace containing $\Omega$. 
\end{itemize}
Let $U$ be any one-sided $\Gamma$-invariant open neighborhood of $\Omega$.
Then $\Gamma$ acts on a properly convex manifold $L$ in $\SI^n$\/
{\rm (}resp. in $\RP^n$\/{\rm )} with sboundary $\partial L$  and 
 following properties: 
\begin{itemize} 
\item $\partial L$ is either strictly convex or is a union of polyhedrons meeting 
one another in angles $< \pi$, 
\item $L^o$ is a neighborhood of $\Omega$, 
\item $L -P$ is closed in $\SI^n -P$ {\rm (}resp. in $\RP^n _P$\/{\rm )}, 
\item 
\[ \Omega \subset L \subset U, \partial L \subset \SI^{n}-P \,
\mathrm{(} \hbox{resp. } \subset \RP^n - P \mathrm{).} \]
\item  $L$ satisfies $\Bd_{\partial \clo(L)} \partial L \subset P$
\end{itemize} 
\end{proposition} 
\begin{proof} 
	We prove for $\SI^n$ first. 
We will prove only for the general case since the special case when 
$\Omega$ is strictly convex and $C^1$, the augmented boundary is 
given as the set of all $(x, h)$ where $x$ is in $\partial \clo(\Omega)$ 
and $h$ is the unique supporting hyperspace of $\Omega$ at $x$. 
Assume without loss of generality 
that $U$ is an asymptotically nice open domain for $\Gamma$. 

For each $(x, h)$ in the augmented boundary of $\Omega$, 
define a half-space $H(x, h) \subset \mathds{A}^n$ bounded by $\tau^{\Ag}(x, h)$, 
and containing $\Omega$ in the boundary. 
In the proofs of Theorems \ref{du-thm-asymnice} and 
\ref{du-thm-asymniceII},
for each $H(x, h)$, $(x, h) \in \partial^\Ag \clo(\Omega)$, 
there exists
an open $n$-hemisphere $H'(x, h) \subset \SI^{n}$ such that 
$H'(x, h) \cap \mathds{A}^{n} = H(x, h)$. 
Then we define
\[V:=\bigcap_{(x, h) \in \partial^\Ag \clo(\Omega)} H'(x, h) \subset \SI^{n}\] 
is a convex open domain containing $\Omega$
as in the proof of Lemma \ref{du-lem-inde}.

$\Gamma$ acts on a compact set 
\[\mathcal{H}:=
\{ h'| h' \hbox{ is an AS-hyperspace to $V$ at } (x, h) \in \partial^\Ag \clo(\Omega), h'\cap \SI^{n-1}_\infty = h\}.\] 
Let $\mathcal{H}'$ denote the set of hemispheres bounded by an element of $\mathcal{H}$ and 
containing $\Omega$. 
 Then we define
\[V:=\bigcap_{H \in \mathcal{H}'} H \subset \SI^{n}\] 
is a convex open domain containing $\Omega$.
Again the set of AS-hyperspaces to $V$ is closed and bounded away from $\SI^{n-1}_{\infty}$. 

First suppose that $V$ is properly convex. 
Then $V$ has a $\Gamma$-invariant 
Hilbert metric $d_V$ that is also Finsler. (See \cite{goldmanbook} and \cite{Kobpaper}.)
Then
\[N_\eps =\{ x\in V| d_V(x, \Omega) < \eps\}\]
 is a convex subset of $V$ by Lemma \ref{prelim-lem-nhbd}. 



A compact tubular neighborhood $M$ of $\Omega/\Gamma$ in $V/\Gamma$ is
diffeomorphic to $\Omega/\Gamma \times [-1,1]$. (See Section 4.4.2 of \cite{Cbook}.)
We choose $M$ in $(U \cap V)/\Gamma$. 
Since 
$\Omega/\Gamma$ is compact, the regular neighborhood of $M$ has a compact closure. 
Thus, $d_V(\Omega/\Gamma, \Bd M) > \eps_0$ for some $\eps_0 > 0$. 
If $\eps < \eps_0$, then $N_\eps \subset M$. We obtain that 
$\Bd_V N_\eps/\Gamma$ is compact.


Clearly, $\Bd_V N_\eps/\Gamma$ has two components, 
each lying in two distinct components of $(V - \Omega)/\Gamma$.
Let $F_1$ and $F_2$ be 
compact fundamental domains of respective components of 
$\Bd_V N_\eps$ with respect to $\Gamma$. 
For each $j=1,2$, 
we choose the set ${\mathcal  H}_{j}$ of
finitely many open hemispheres $H_i$, $H_i \supset \Omega$, 
such that open sets $(\SI^n - \clo(H_i) )\cap N_\eps$ cover $F_j$. 
By Lemma \ref{du-lem-locfin}, 
the following is an open set containing $\Omega$
\[W := \bigcap_{g \in \Gamma} \bigcap_{H_{i} \in \mathcal{H}_{1}\cup \mathcal{H}_{2}} g(H_i) \cap V.\]
Since any path in $V$ from $\Omega$ to $\Bd_V N_\eps$ must meet $\Bd_V W -P$ first, $N_\eps$ contains $W$ and $\Bd_V W$. 
Now, $\Bd_V W/\Gamma$ is a union of 
compact totally geodesic polyhedrons meeting one another only at
angles $< \pi$ or in the empty set. 
Let $L$ be $\clo(W) \cap \mathds{A}^n$. Then by Lemma \ref{du-lem-attracting3}, 
the hypersurface $\partial L$ has boundary only in $\Bd \mathds{A}^{n}$
since $\Gamma$ satisfies the uniform middle-eigenvalue condition with respect to 
$\Bd \mathds{A}^n$.


Suppose that $V$ is not properly convex. Then $\Bd V$ contains $v, v_-$.
$V$ is a convex domain of the form $\{v, v_-\} \ast \Omega$. 
This follows from Proposition \ref{prelim-prop-classconv} 
where we take the closure of $V$ and then the interior. 
We take any two open hemispheres $S_1$ and $S_2$ containing $\clo(\Omega)$ such that 
$\{v, v_-\} \cap S_1 \cap S_2 = \emp$. 
Then $\bigcap_{g \in \Gamma} g(S_1 \cap S_2) \cap V$ is a properly convex open domain containing $\Omega$, 
and we can apply the same argument as above. 

To prove for $\RP^n$, we need to find $L$ in a sufficiently 
thin neighborhood of $\Omega$, which the theorem for $\SI^n$ 
provides. Then we can project to obtain the desired set. 
\end{proof}

\begin{lemma}\label{du-lem-attracting3} 
	Let $\Gamma$ be a discrete group in $\SLnp$ acting on a properly convex domain
	$\Omega \subset \Bd \mathds{A}^n$ such that $\Omega/\Gamma$ is 
	a closed $(n-1)$-orbifold.  Suppose that $\Gamma$ satisfies the uniform middle-eigenvalue condition with respect to $\Bd \mathds{A}^n$ and acts on a properly convex domain $V$ in $\SI^n$ such that $\clo(V) \cap \Bd \mathds{A}^n = \clo(\Omega)$ holds. 
		Suppose that $g_i$ is a sequence of mutually distinct elements of $\Gamma$ acting on $\Omega$. 
	Let $J$ be a compact subset of $V$. Then 
	$\{g_i(J)\}$ has accumulation points only at  $\clo(\Omega)$. 
\end{lemma} 
\begin{proof} 
	Since the $\Gamma$-action on $\Omega$ is cocompact, $\hat h(g_i)$ is an unbounded 
	sequence of elements in $\SL_{\pm}(n, \bR)$. 
	Recalling \eqref{du-eqn-bendingm4}, we can write  
	each $g_i$ in the form
	\begin{equation}\label{du-eqn-bendingm5} 
	\left(
	\begin{array}{cc}
	\frac{1}{\lambda_{{\tilde E}}(g_i)^{1/n}} \hat h(g_i)          &       \vec{b}_{g_i}     \\
	\vec{0}          &     \lambda_{{\tilde E}}(g_i)                  
	\end{array}
	\right)
	\end{equation}
	where $\vec{b}_{g_i}$ is an $n\times 1$-vector and $\hat h(g_i)$ is 
	an $n\times n$-matrix of determinant $\pm 1$
	and $\lambda_{{\tilde E}}(g_i) > 0$.  
%
	Let $m(\hat h(g_i))$ denote the maximal modulus of 
	the entries of $\hat h(g_i)$ in $\SL_\pm(n, \bR)$. 
	We may assume without loss of generality 
	$\{\hat h(g_i)/m(\hat h(g_i))\} \ra g_{n-1,\infty}$ in $M_n(\bR)$.
	A matrix analysis easily tells us 
	$\lambda_1(\hat h(g_i)) \leq n m(\hat h(g_i))$
	since the later term bounds the amount of stretching of norms of vectors
	under the action of $\hat h(g_i)$.   
	By Lemma \ref{prelim-lem-convcomp}, dividing by $\lambda_1(\hat h(g_i))$, 
	we have 
	\[ \hat h(g_i)/\lambda_1(\hat h(g_i)) \ra A g_{n-1,\infty}  \]
	for 
	 $A \geq 1/n$, or the sequence is unbounded in $M_n(\bR)$. 

Now, we let 
	$\{\llrrparen{g_i}\} \ra \llrrparen{g_\infty}$ where 
	$g_\infty$ is obtained as a limit of 
	{\large 
		\begin{equation}\label{du-eqn-bendingm7} 
	\frac{\lambda_1({\hat h}(g_i))}{m(\hat h(g_i)}  	
	\left(
		\begin{array}{cc}
		\frac{1}{\lambda_1({\hat h}(g_i))} \hat h(g_i)          &   
		\frac{\lambda_{{\tilde E}}(g_i)^{1/n}}{\lambda_1({\hat h}(g_i))}   \vec{b}_{g_i}    \\
		\vec{0}          &    \frac{\lambda_{{\tilde E}}(g_i)^{1/n}\lambda_{{\tilde E}}(g_i)}{\lambda_1(\hat h(g_i))}                
		\end{array}
		\right).
		\end{equation}
	}	
	By the uniform middle-eigenvalue condition \eqref{du-eqn-umecDii-inverse}, we obtain 
	that $\{\llrrparen{g_i}\} \ra \llrrparen{g_\infty}$ where
	$g_\infty$ is of the form 
	\begin{equation}\label{du-eqn-bendingm8} 
	\left(
	\begin{array}{cc}
	g_{n-1,\infty}       &       \vec{b}   \\
	\vec{0}          &    0                 
	\end{array}
	\right).
	\end{equation}
	by rescaling if necessary. 
	Now, $R_\ast(\{g_i\})$ is a subset of $\SI_{\infty}^{n-1}$ 
	since the lower row is zero. 
	%
	Hence, $\{g_i(J)\}$ geometrically converges to $g_\infty(J) \subset \SI_{\infty}^{n-1}$. 
	Since $\Gamma$ acts on $\clo(V)$ and $\clo(V) \cap \Bd \mathds{A}^n = \clo(\Omega)$, 
	we obtain $g_\infty(J) \subset \clo(\Omega)$. 
\end{proof}

\begin{lemma}\label{du-lem-locfin}
	Let $\Gamma$ be a discrete group of projective automorphisms of 
	a properly convex domain $V$ and 
	a domain $\Omega \subset V$ of dimension $n-1$. 
	Assume that $\Omega/\Gamma$ is a closed $(n-1)$-orbifold. 
	Suppose that $\Gamma$ satisfies the uniform middle-eigenvalue condition
	with respect to the hyperspace containing $\Omega$. 
	Let $\hat P$ be a finite collection of hyperspaces of 
$\SI^{n}$ such that $P \cap \clo(\Omega) = \emp$ for $P \in \hat P$
	and $P \cap V \ne \emp$. 
	Then the collection 
\[\{g(P)\cap V|g \in \Gamma, P\in \mathcal{P}\}\] is a locally finite in $V$, 
and consists of closed subsets of $V$. 
\end{lemma}
\begin{proof} 
	Suppose not. Then there exists a sequence $\{x_i\}$ where $x_{i}\in g_i(P_i) \cap V$, 
$P_i \in \hat P$, mutually distinct collection $g_{i}\in \Gamma$ such that $\{x_i\} \ra x_\infty \in V$. 
	
	We have 
	$x_{i} \in F$ for a compact set $F \subset V$. 
	Then Lemma \ref{du-lem-attracting3} applies. 
	$\{g_{i}^{-1}(F)\}$ accumulates to $\partial \clo(\Omega)$.
	In particular, $g_i^{-1}(x_i)$ accumulates to points of $\partial \clo(\Omega)$.
	Since $g_i^{-1}(x_{i}) \in P_i \cap V$, an accumulating hyperspace 
$P_\infty\in \hat P$ of $P_i$ 
meets $\clo(\Omega)$.  This is a contradiction. 
\end{proof}

\subsection{Approximating a convex hypersurface by strictly convex hypersurfaces}


See also Chapter 9 of \cite{DGK21} where they obtain the  $C^1$-property only. 

\begin{theorem}\label{du-thm-lensnpre} 
		We assume that $\Gamma$ is a projective group with a properly convex 
	affine action with the triple $(\Gamma, U, D)$ for $U \subset \mathds{A}^n
	\subset \SI^n$.
	Assume the following{\em :} 
	\begin{itemize}  
\item  $\Gamma$ satisfies the uniform middle-eigenvalue condition
with respect to $\Bd \mathds{A}^n$.
	\item $U$ is an asymptotically nice properly convex 
	domain closed in $\mathds{A}^n$,
	\item the boundary $\Bd_{\mathds{A}^n} U
=\Bd U \cap \mathds{A}^n$, which is an $(n-1)$-manifold, 
	is in an asymptotically nice properly convex open domain $V'$ where $\Gamma$ acts on,
	and $\clo(U)\cap \mathds{A}^n \subset V'$, 
	\item $\clo(V')\cap \Bd \mathds{A}^{n} = D$, and  
	\item $\Bd_{\mathds{A}^n} U/\Gamma$ is a compact convex 
	hypersurface.
	\end{itemize}  
	Then there exists an  
	asymptotically nice properly convex domain $V$ closed in $V'$ containing $U$
	such that $\partial V/\Gamma$ is a compact hypersurface with strictly convex 
	smooth boundary. 
	Furthermore, 
for any complete Riemannian metric on $V'/\Gamma$, 
$\partial V/\Gamma$ can be chosen to be arbitrarily close to 
	$\Bd_{\mathds{A}^n} U/\Gamma$ in $V/\Gamma$.
\end{theorem} 
\begin{proof} 
%
We  prove by showing that the conditions of Lemma \ref{du-lem-approx} hold. 

Let $V''$ be a properly convex domain such that $\partial V''' \subset U$, 
and suppose that $(\bGamma, V'', D)$ is a properly convex affine triple.
We can construct such a domain by the proof of Proposition \ref{du-prop-lensn}.

Let $\mathcal{P}$ denote the set of hyperspaces 
sharply supporting $V'''$ at $\Bd_{\mathds{A}^n} U$.
For each $P\in \mathcal{P}$, 
let $H_P$ denote the open half-space bounded by $P$ 
disjoint from the interior of $U$; that is,  the exterior of $P$. 
The dual $P^\ast$ of $P\in \mathcal{P}$ in $\SI^{n\ast}$ is a point of 
the properly convex domain $V^{\prime \prime \ast}$ 
by \eqref{prelim-eqn-dualinc}. 
Hence, 
by Proposition \ref{prelim-prop-duality} applied to 
$\clo(U)$ and the hyperspace $\Bd_{\mathds{A}^n} U$, 
the set $\mathcal{P}^\ast$ of dual points corresponding to elements of
$\mathcal{P}$ is a properly embedded hypersurface in the interior of 
$\clo(V^{\prime \prime\ast})$.  

Also, $\Bd_{\mathds{A}^n}  U/\bGamma$ is a compact orbifold by 
Proposition \ref{prelim-prop-duality}(iv).
By the duality described above, $\mathcal{P}^\ast/\bGamma^\ast$ is a compact orbifold,
and so is $\mathcal{P}/\bGamma$. 
There is a fundamental domain $F$ of $\mathcal{P}$ under $\bGamma$. 


For each point of $\Bd_{\mathds{A}^n} V'''$, there is a neighborhood $N$ 
whose closure lies entirely within $U$. 

%
%

For any $\eps> 0$, there exists a compact set $K'$ 
such that elements of $\mathcal{P} - K'$ are $\eps$-$\bdd_H$-close 
to the hyperspaces asymptotic to $V'$ by the second paragraph above. 
Therefore, for each $x\in \Bd_{\mathds{A}^n} V'''$, 
there exists a neighborhood $N$ such that 
the set of sharply supporting hyperspaces $P$ 
in $\mathcal{P}$ with $N \subset H_P$ is compact.  
Hence, Lemma \ref{du-lem-approx} implies the result. 
\end{proof}

We have a generalization:

\begin{theorem}\label{du-thm-horonpre} 
	Assume the following\/{\em :} 
	\begin{itemize} 
		\item  $\Gamma$ is a projective discrete group 
		acting properly on a properly convex domain $U$ closed in $\mathds{A}^n
		\subset \SI^n$, 
		\item  $\{p\}=\clo(U)- U$ is a singleton, 
		\item  $\Gamma$ is a cusp group, 
		\item the manifold boundary $\Bd_{\mathds{A}^n} U $ 
		is in a properly convex open domain $V'$ where $\Gamma$ acts on, 
		\item $U$ is closed in $V'$, and
		\item  $\Bd_{\mathds{A}^n}  U/\Gamma$ is a compact convex 
		hypersurface. 
	\end{itemize} 
	Then there exists an  
	properly convex domain $V$ closed in $V'$ containing $U$
	such that $\partial V/\Gamma$ is a compact hypersurface with strictly convex
	smooth boundary. 
	Furthermore, $\partial V/\Gamma$ can be chosen to be arbitrarily close to 
	$\Bd_{\mathds{A}^n}  U/\Gamma$ in $V/\Gamma$ with any complete Riemannian metric on 
	$V'/\Gamma$. 
\end{theorem}
\begin{proof} 
	Let $\mathds{A}^n$ be the affine space bounded by a hyperspace 
	tangent to $\clo(U)$ at $p$.
	Let $\mathcal{P}$ denote the set of hyperspaces 
	sharply supporting $U$. 
The union $\mathcal{P}$ can only accumulate to itself or 
to points of a unique supporting hyperspace $\partial \mathds{A}^n$
	to $U$ at $p$. 
	This follows since the dual of 
$\mathcal{P}$ is a horosphere where $\Gamma$ acts as a cusp group. 
Then the analogous proof as 
	that of Theorem \ref{du-thm-lensnpre}, we can show that 
for any $\eps> 0$, there exists a compact set $K'$ 
such that elements of $\mathcal{P} - K'$ are $\eps$-$\bdd_H$-close 
to $\partial \mathds{A}^n$. Now, Lemma \ref{du-lem-approx} is applicable exactly. 
\end{proof} 

\begin{lemma}\label{du-lem-approx} 
		We assume that $\Gamma$ is a projective group
acting on $\mathds{A}^n$ and an open domain $U \subset \mathds{A}^n	\subset \SI^n$.
	Suppose that $U$ is a properly convex domain in an affine space $\mathds{A}^n$
	in $\SI^n$
	with $\Bd_{\mathds{A}^n} U$ a properly embedded hypersurface in it.
	Let $\mathcal{P}$ denote the set of 
	sharply supporting hyperspaces of $U$ meeting $\Bd_{\mathds{A}^n} U$.
	Assume the following\/{\rm :}
	\begin{itemize}
\item  $\Bd_{\mathds{A}^n} U/\Gamma$ is a compact orbifold.
\item $\clo(U) \cap \mathds{A}^n$ is in a convex open domain $V'$ where $\Gamma$ acts on, 
\item Each point $x$ of 
$V'$ has a neighborhood $N$ such that 
the set of sharpley supporting hyperspaces $P$ 
in $\mathcal{P}$ with $N$ meeting $H_P$ for the exterior $H_P$ of $P$ is compact.  
\end{itemize} 
Then there exists an  
properly convex domain $V$ closed in $V'$ containing $U$
such that $\partial V/\Gamma$ is a compact hypersurface with strictly convex
smooth boundary. 
Furthermore, $\partial V/\Gamma$ can be chosen to be arbitrarily close to 
$\Bd_{\mathds{A}^n} U/\Gamma$ in $V/\Gamma$ with any complete Riemannian metric on 
$V'/\Gamma$. 
\end{lemma} 
\begin{proof} 
%
For each $P\in \mathcal{P}$, 
we define an affine function $f_P$ on $\mathds{A}^n$ 
such that $f_P^{-1}(0) = P$ and 
$f_P> 0 $ on $H_P$.
It follows that if $P' = g(P)$ for $P \in \mathcal{P}$ and $g \in \Gamma$, 
then $f_{P'}\circ g^{-1} = f_P$.


We define a smooth function 
\[g(t) = t^2 \exp(1/t^2) \hbox{ for } t > 0, \hbox{ and } g(t)=0 \hbox{ for } t \leq 0.\] 
We let $g_P = g \circ f_P$.
Then by the premise, $g_P(x)$  for each $x \in V'$ 
is nonzero for only compact subset of $\mathcal{P}$.

The $\Gamma$-action on $V'$ preserves the Hilbert metric of $V'$ and hence the action 
is properly discontinuous on $V'$. 
Since the action is properly discontinuous, we can put a $\Gamma$-invariant Riemannian metric
on $V'$. 
The dual $\mathcal{P}^\ast \subset V^{\prime \ast}$ is a dual set of $\mathcal{P}$ considering 
each hyperspace as a linear function in $\bR^{n+1}$. 

Since $\mathcal{P}$ is the boundary of the supporting hyperspaces of $U$, it is 
the boundary of $U^*$ by duality.   
Hence, $\mathcal{P}^\ast$ can be considered a topological manifold in $V^{\prime\ast}$.  
Since $\Gamma^\ast$-action on $\mathcal{P}^\ast \subset V^{\prime\ast}$ is 
properly discontinuous, $\mathcal{P}^\ast/\Gamma^\ast$ is a compact topological orbifold. 
Also, we may assume that $f_P$ for $P\in \mathcal{P}$ is chosen continuously with respect to $P$ 
by taking the fundamental domain of $\mathcal{P}^\ast$ under the $\Gamma^\ast$-action.
There is a $\Gamma$-invariant measure $d\mu$ on $\mathcal{P}^\ast$ 
compatible with a positive continuous function 
times a volume on each chart of 
$\mathcal{P}$.  
We define a smooth function 
\[\chi_U: V' \ra \bR \hbox{ by } \int_{P\in \mathcal{P}} g_P d\mu.\]
Hence, $\chi_U$ is well-defined in $V'$ by the above paragraph. 
Moreover, 
\[\chi_U^{-1}(0)= \bigcap_{P\in \mathcal{P}} (\mathds{A}^n - H_P)= \clo(U)\cap \mathds{A}^n.\]  
The third item of the premise and 
the proof of Proposition 2.1 of Ghomi \cite{Ghomi} imply that 
$\chi_U$ is strictly convex on $V' - U$ since only compact subset of $\mathcal{P}$ is involved in
the computations for each neighborhood of the third item of the premise. 
By our definition, $\chi_U$ is $\Gamma$-invariant. 

We give an arbitrary Riemannian metric $\mu'$ on $V'/\Gamma$.
There exists a neighborhood $N$ of $\clo(U) \cap V'/\Gamma$ in $V'/\Gamma$
where 
$\chi_U$ has a nonzero differential in $N - \clo(U)/\Gamma$
as we can see from the integral $\int_{P\in \mathcal{P}} Dg_P d\mu$
where the elements denoted by $Dg_P$ are in a properly convex cone $C_{\mathcal{P}}$ 
in $\bR^{n\ast}$ spanned by $\{u_P| P \in \mathcal{P}\}$
for each point of $V' - \clo(U)$. 
Then as $\eps \ra 0$,
\[\Sigma_\eps := \{\chi^{-1}_U(\eps)\}/\Gamma \ra \Bd_{\mathds{A}^n} U/\Gamma\]
geometrically
since $N$ has a compact closure and the gradient vectors are
uniformly bounded with respect to $\mu'$ and are zero only at points of
$U/\Gamma$ in the closure and hence we can isotopy $\Sigma_\eps$ 
along the gradient vector field to as close to $\Bd_{\mathds{A}^n}  U/\Gamma$ as we wish. 
(See Batyrev \cite{Batyrev} and Ben-Tal \cite{Ben-Tal} also.)
Furthermore, since $\chi_U$ is strictly convex, 
$\chi^{-1}_U(\eps)$ is a  strictly convex smooth hypersurface on 
which $\Gamma$ acts.

\end{proof}

The following is a consequence of Theorems \ref{du-thm-asymnice} and \ref{du-thm-asymniceII}
and Proposition \ref{du-prop-lensn}: 

\begin{theorem}\label{du-thm-lensn} 
 Suppose that $\Gamma$ satisfies the uniform middle-eigenvalue condition
with respect to $\Bd \mathds{A}^n$.
We assume that $\Gamma$ is a projective group with a properly convex 
	affine action acting on a properly convex domain $\Omega$ in $\Bd \mathds{A}^n$. 
Then we can choose a lens $L$ containing $\Omega$ 
such that $L/\Gamma$ is a lens neighborhood with strictly convex smooth boundary. 
\end{theorem} 
\begin{proof}
We can smooth $\Bd_VW$ to a strictly convex hypersurface 
to obtain a lens-neighborhood $W' \subset W$ of $\Omega$ in $N_\eps$ 
where $\Gamma$ acts cocompactly by Theorem \ref{du-thm-lensnpre}. 
\end{proof}


%
%
%




\chapter[Properly convex R-ends and T-ends]{Properly convex radial ends and totally geodesic ends: lens properties} \label{ch-pr}


	We will consider properly convex ends in this chapter.
In Section \ref{pr-sec-mainresult}, we define the uniform middle-eigenvalue conditions for R-ends and T-ends. We state the main results of this chapter
Theorem \ref{pr-thm-equiv}: the equivalence of
these conditions with the generalized lens conditions for R-ends or T-ends.
The generalized lens conditions often improve to
lens conditions, as shown in Theorem \ref{pr-thm-secondmain}. 
In Section \ref{pr-sec-endth}, we start to study the R-end theory. First, we discuss the holonomy representation spaces.
Tubular actions and the dual theory of affine actions are discussed. 
We show that distanced actions
and asymptotically nice actions are dual. Hence, the uniform middle
eigenvalue condition implies the distanced action 
deduced from the dual theory in Chapter \ref{ch-du}.
In Section \ref{pr-sec-chlens}, we prove the main results. 
In Section \ref{pr-sub-eigen}, 
we estimate the largest norm $\lambda_{1}(g)$ of eigenvalues 
in terms of word length.  In Section \ref{pr-sub-umecorbit},  we study orbits 
under the action satisfying the uniform middle-eigenvalue conditions. In Section \ref{app-sec-Koszul}, 
we prove a minor extension of Koszul's openness for bounded manifolds, well-known 
to many people.
In Section \ref{pr-subsec-redlens}, 
we show how to prove the strictness of the boundary of lenses 
and prove our main result Theorem \ref{pr-thm-equiv} using the orbit results 
and the Koszul's openness. 
In Section \ref{pr-sub-umecl}, 
we prove one of the main results: Theorem \ref{pr-thm-secondmain}. 
In Section \ref{pr-sec-lens}, 
we show that 
the lens-shaped ends have concave end neighborhoods, 
and we discuss the properties of lens-shaped ends in 
Theorems \ref{pr-thm-lensclass} and \ref{pr-thm-redtot}.
If the generalized lens-shaped end is virtually factorizable, 
it can be made into a lens-shaped totally-geodesic R-end, which is a surprising result. 
In Section \ref{pr-sec-dualT},  
we obtain the duality between the lens-shaped T-ends and generalized lens-shaped R-ends. 




 The main reason that we are studying the lens-shaped ends is to understand deformations preserving the convexity properties. 
 These objects are useful in trying to understand this phenomenon.

We also remark that 
sometimes a lens-shaped p-end neighborhood may not exist for an R-p-end
within a given convex real projective orbifold. 
However, a generalized lens-shaped p-end neighborhood  may exist for the R-p-end.

%

\section{Main results} \label{pr-sec-mainresult}
Let $\orb$ be a convex real projective $n$-orbifold, and 
let $\torb$ be a convex domain in $\SI^n$ projectively covering $\orb$. 
Let $h: \pi_1(\orb) \ra \SLpm$ denote the holonomy homomorphism 
with its image $\bGamma$. 
All statements will take place in $\SI^n$ by default 
However, readers can easily modify these to $\RP^n$-versions
by Proposition \ref{prelim-prop-closureind}  results in Section \ref{prelim-sub-lifting}.

\begin{definition}\label{pr-defn-conc} 
	Suppose that $\tilde E$ is an R-p-end of generalized lens-type. 
 Then $\tilde E$ has a p-end neighborhood that is projectively diffeomorphic 
 to the interior of $\{p\} \ast L^o -\{p\}$ under $\dev$ where 
 $\{p\} \ast L$ is a generalized lens cone over a generalized lens $L$
 where $\partial (\{p\}\ast L -\{p\}) = \partial_+L$ for a boundary component 
 $\partial_+ L$ of $L$, and let 
 $h(\pi_1(\tilde E))$ acts on $L$ properly and cocompactly. 
 A {\em concave pseudo-end neighborhood} of $\tilde E$ is the open pseudo-end neighborhood in $\torb$ projectively diffeomorphic to 
 $\{p\}\ast L -\{p\} - L$ for some choice of a lens $L$. 
A {\em concave end neighborhood} of an end $E$ is an end neighborhood covered by 
a concave pseudo-end neighborhood.
 \index{end!neighborhood!concave|textbf} 
 \index{p-end!neighborhood!concave|textbf} 
 \index{pseudo-end!neighborhood!concave|textbf} 
 \end{definition}


\subsection{Uniform middle-eigenvalue conditions} 
The following applies to both R-ends and T-ends. 
Let $\tilde E$ be a p-end and $\bGamma_{\tilde E}$ the associated p-end holonomy group.
We say that 
$\tilde E$ is {\em non-virtually-factorizable} if any finite-index subgroup has a finite center or 
$\bGamma_{\tilde E}$ is virtually center-free;
otherwise, $\tilde E$ is virtually factorizable by Theorem 1.1 of \cite{Benoist03}. 
(See Section \ref{prelim-sub-ben}.)

\index{factorizable!virtually} 
\index{factorizable!non-virtually}

Let $\tilde \Sigma_{\tilde E}$ denote the universal cover of the end orbifold $\Sigma_{\tilde E}$ associated with $\tilde E$. 
We recall Proposition \ref{prelim-prop-Ben2}
(Theorem 1.1 of Benoist \cite{Benoist05}). 
If $\bGamma_{\tilde E}$ is virtually factorizable, then 
$\bGamma_{\tilde E}$ satisfies the following conditions: 
%
%
\begin{itemize}
\item  $\clo(\tilde \Sigma_{\tilde E}) = K_{1}\ast \cdots \ast K_{k}$ 
where each $K_{i}$ is properly convex or is a singleton. 
\item Let $G_{i}$ be the restriction of the $K_{i}$-stabilizing subgroup of $\bGamma_{\tilde E}$ to $K_{i}$. 
Then $G_{i}$ acts on $K_{i}^{o}$ cocompactly. 
(Here $K_{i}$ can be a singleton, and $\Gamma_{i}$ a trivial group.)
\item A finite-index subgroup $G'$ of $\bGamma_{\tilde E}$ is isomorphic to a cocompact subgroup of 
$\bZ^{k-1}\times G_{1}\times \cdots \times G_{k}$.
\item  The center $\bZ^{k-1}$ of $G'$ is a subgroup acting trivially on each $K_{i}$. 
\end{itemize} 

Note that there are examples of discrete groups of form $\bGamma_{\tilde E}$  where 
$G_{i}$ are nondiscrete. (See also Example 5.5.3 of \cite{Morris15} as pointed out by M. Kapovich.)

We will use simply $\bZ^{k-1}$ to represent the corresponding group on $\bGamma_{\tilde E}$.
Here, $\bZ^{k-1}$ is called a virtual center of $\bGamma_{\tilde E}$. 
\index{virtual center} 


Let $\Gamma$ be generated by finitely many elements $g_1, \ldots, g_m$. 
Let $w(g)$ denote the minimum word length of $g \in G$ written as words of $g_{1}, \dots, g_{m}$. 
The {\em conjugate word length} $\cwl(g)$ of $g \in \pi_1(\tilde E)$ is
\[\min \{ w(cgc^{-1})| c \in \pi_{1}(\tilde E) \}. \]
 \index{word length!conjugate} \index{cwl@$\cwl(\cdot)$}

Let $d_K$ denote the Hilbert metric of the interior $K^o$ of a properly convex domain $K$ in $\RP^n$ or $\SI^n$. 
Suppose that a projective automorphism group $\Gamma$ acts on $K$ properly. 
Let $\leng_K(g)$ denote the infimum of $\{ d_K(x, g(x))| x \in K^o\}$, compatible with $\cwl(g)$. 

For the following, we do assume that $\Sigma_{\tilde E}$ is properly convex. 

\begin{definition}\label{pr-defn-umec}
Let $\mbv_{\tilde E}$ be a p-end vertex of an R-p-end $\tilde E$. 
Let $K := \clo(\tilde \Sigma_{\tilde E})$. 
The p-end holonomy group $\bGamma_{\tilde E}$ satisfies the {\em uniform middle-eigenvalue condition  {\rm (}umec{\rm )} with respect to $\mbv_{\tilde E}$, to the R-p-end $\tilde E$, 
or the corresponding R-end $E$} if the following hold: 
\begin{itemize}
\item each $g\in \bGamma_{\tilde E}$ satisfies for a uniform  $C> 1$ independent of $g$
\begin{equation}\label{pr-eqn-umec}
C^{-1} \leng_K(g) \leq \log\left(\frac{\lambda_1(g)}{\lambda_{\mbv_{\tilde E}}(g)}\right) 
\leq C \leng_K(g) , 
\end{equation}
for  the largest norm $\lambda_1(g)$  of the eigenvalues of $g$
and the eigenvalue $\lambda_{\mbv_{\tilde E}}(g)$ of $g$ at $\mbv_{\tilde E}$.
\index{lambda@$\lambda_{\mbv_{\tilde E}}(\cdot)$} 
\index{middle-eigenvalue condition!uniform|textbf} 
\index{uniform middle-eigenvalue condition|textbf}

\end{itemize} 
We required here that $\tilde \Sigma_{\tilde E}$ is properly convex and cocompact under
$\bGamma_{\tilde E}$ as $\Sigma_{\tilde E}$ is properly convex and compact.
Of course, we choose the matrix of $g$ such that $\lambda_{\mbv_{\tilde E}}(g) > 0$.
See Remark \ref{intro-rem-SL} as we are looking for the lifting of $g$ that acts on 
p-end neighborhood. 
We remark that the condition does depend on the choice of $\mbv_{\tilde E}$;
however, the radial end structures will determine the end vertices.



\end{definition}
The definition of course applies to the case when $\bGamma_{\tilde E}$ has the finite-index subgroup with the above properties. 


We recall a dual definition identical with the definition in 
Section \ref{du-sec-affine} but adopted to the cases of T-p-ends. 
\begin{definition} \label{pr-defn-umecD}
Suppose that $\tilde E$ is a properly convex T-p-end. 
 Suppose that the ideal boundary component $\tilde \Sigma_{\tilde E}$ of 
 the T-p-end is properly convex. Let $K = \clo(\tilde \Sigma_{\tilde E})$. 
Let $g^*:\bR^{n+1 \ast} \ra \bR^{n+1 \ast}$ be the dual transformation of $g: \bR^{n+1} \ra \bR^{n+1}$. 
The p-end holonomy group $\bGamma_{\tilde E}$ satisfies the {\em uniform middle-eigenvalue condition {\rm (}umec{\rm )} with respect to $\tilde \Sigma_{\tilde E}$, 
the T-p-end $\tilde E$, or the corresponding T-end $E$}
\begin{itemize}
\item  if each $g\in \bGamma_{\tilde E}$ has a uniform  constant $C> 1$ independent of $g$
such that 
\begin{equation}\label{pr-eqn-umecD}
C^{-1} \leng_K(g) \leq \log\left(\frac{\lambda_1(g)}{\lambda_{K^{*}}(g^{*})}\right) \leq C \leng_K(g)
\end{equation}
 for the largest norm $\lambda_1(g)$ 
of the eigenvalues of $g$ and  
the eigenvalue $\lambda_{K^{*}}(g^{*})$ of $g^*$ in the vector in the direction of $K^*$, the point dual 
to the hyperspace containing $K$. 

\end{itemize} 

\end{definition} 

Again, the condition depends on the choice of the hyperspace containing 
$\tilde \Sigma_{\tilde E}$, i.e., the T-p-end structure. 
(We again lift $g$ such that $\lambda_{K^\ast}(g) > 0$.)

We mention that this condition is stronger than one where we replace $\lambda_1$ with 
the largest singular value of $g$ by obvious matrix inequalities. 
The reverse is not true of course. 
This condition is similar to the Anosov condition studied by Guichard and Wienhard \cite{GW12}, and the results also seem similar. We do not use their theories.
They also use word length instead. 
One may look at the paper of Kassel-Potrie \cite{KP2002} to understand 
the relationship between eigenvalues and singular values when the group is word-hyperbolic. 
We use the eigenvalues to obtain conjugacy invariant conditions
which are needed in proving the converse part of Theorem \ref{pr-thm-equiv}. 
Our main tools to understand these questions are in Chapter \ref{ch-du}
which we use here.


We will see that the condition is an open condition; and hence a ``structurally stable one."
(See Corollary \ref{app-cor-mideigen}.)  

\subsection{Lens and the uniform middle-eigenvalue condition}

As holonomy groups, the condition for being a generalized lens R-p-end and one for 
being a lens R-p-end are equivalent. 
For the following, we are not concerned with a lens cone being in $\torb$. 




\begin{theorem}[Lens holonomy]\label{pr-thm-equiv}
Let $\tilde E$ be an R-p-end of a  convex real projective $n$-orbifold. 
Then 
 the holonomy group 
 $h(\pi_{1}(\tilde E))$ satisfies the uniform middle-eigenvalue condition
 for the R-p-end vertex $\mbv_{\tilde E}$ 
 if and only if it acts on a topological lens cone with vertex $\mbv_{\tilde E}$ 
 and its topological lens properly and cocompactly. 
 Moreover, in this case, 
 the lens cone exists in the union of great segments with 
 the vertex $\mbv_{\tilde E}$
 in the direction of a properly convex 
 domain $\Omega \subset \SI^{n-1}_{\mbv_{\tilde E}}$, where 
  $h(\pi_{1}(\tilde E))$ acts properly discontinuously. 
\end{theorem} 

For the following, we are concerned with a lens cone being in $\torb$. 

\begin{theorem}[Theorem \ref{pr-thm-equ} (Actual lens cone )]\label{pr-thm-secondmain} 
Let $\mathcal{O}$ be a convex real projective $n$-orbifold. 
\begin{itemize} 
\item Let $\tilde E$ be a properly convex R-p-end. 
\begin{itemize} 
\item The p-end holonomy group satisfies the uniform middle-eigenvalue condition with respect to 
R-p-end $\tilde E$ 
if and only if $\tilde E$ is a generalized lens-shaped R-p-end.
\end{itemize} 
\item Assume that the holonomy group of $\mathcal{O}$ is strongly irreducible, strongly tame, 
and $\orb$ is properly convex. 
If $\orb$ satisfies the triangle condition 
{\rm (}see Definition \ref{pr-defn-tri}{\rm )} 
or $\tilde E$ is virtually factorizable or is a totally geodesic R-end, 
then we can replace the term ``generalized lens-shaped'' with ``lens-shaped''
in the above statement. 
\end{itemize} 
\end{theorem}
We will prove the analogous result for totally geodesic ends in Theorem \ref{pr-thm-equ2}. 
\index{triangle condition} 

Notice that there is no condition on $\orb$ to be properly convex.

Another main result concerns the duality of lens-shaped ends: 
Recalling from Section \ref{prelim-sub-duality}, 
we have $\RP^{n \ast}={\bP}(\bR^{n+1 \ast})$ the dual real projective space of 
$\RP^n$. Recall also $\SI^{n \ast} = \SI(\bR^{n+1 \ast})$ 
as the dual spherical projective space of $\SI^n$. 
In Section \ref{pr-sec-endth}, we define the projective dual domain $\Omega^*$ in $\RP^{n \ast}$ 
to a properly convex domain $\Omega$ in $\RP^n$ where 
the dual group $\Gamma^*$ to $\Gamma$ acts on. 
Vinberg showed that there is a duality diffeomorphism between $\Omega/\Gamma$ and $\Omega^{\ast}/\Gamma^{\ast}$. 
The ends of $\orb$ and $\orb^*$ are in a one-to-one correspondence. 
Horospherical ends are dual to themselves, i.e., ``self-dual types'', 
and properly convex R-ends and T-ends are dual to one another. (See Proposition \ref{pr-prop-dualend}.)
We will see that generalized lens-shaped 
properly convex R-ends 
are always dual to lens-shaped T-ends by Corollary \ref{pr-cor-duallens2}.

We mention that Fried also solved this question when the linear holonomy is 
in $\SO(2, 1)$ \cite{Friedlens}. 

\section{The end theory}
\label{pr-sec-endth}
In this section, we discuss the properties of lens-shaped radial and totally geodesic ends and their duality also.


\subsection{The holonomy homomorphisms of the end fundamental groups: the tubes.} \label{pr-sub-holfib}

Although we restrict our discussion to $\SI^n$ only here but the obvious corresponding $\RP^n$-version of the theory exists.
Let $\tilde E$ be an R-p-end of $\torb$. 
Let $\SLnp_{\mbv_{\tilde E}}$ be the subgroup of $\SLnp$ fixing a point $\mbv_{\tilde E} \in \SI^n$.
This group can be understood as follows from letting $\mbv_{\tilde E} = [0, \ldots, 0, 1]$ 
as a group of matrices: For $g \in \SLnp_{\mbv_{\tilde E}}$, we have 
\begin{equation} \label{pr-eqn-mForm}
\left( \begin{array}{cc} 
        \frac{1}{\lambda_{\mbv_{\tilde E}}(g)^{1/n}} \hat h(g) & \vec{0} \\ 
        \vec{v}_g                & \lambda_{\mbv_{\tilde E}}(g)
        \end{array} \right) 
        \end{equation}
where $\hat h(g) \in \SLn, \vec{v} \in \bR^{n \ast}, \lambda_{\mbv_{\tilde E}}(g) \in \bR_+ $. 
Here, \[\lambda_{\mbv_{\tilde E}}: g \mapsto \lambda_{\mbv_{\tilde E}}(g) 
\hbox{ for } g \in \SLnp_{\mbv_{\tilde E}}\] is a homomorphism.
In particular, it is trivial in the commutator group $[\bGamma_{\tilde E}, \bGamma_{\tilde E}]$. 
There is a group homomorphism 
\begin{align} 
{\mathcal L}': \SLnp_{\mbv_{\tilde E}} & \ra \SLn \times \bR_+ \nonumber \\
g &\mapsto (\hat h(g), \lambda_{\mbv_{\tilde E}}(g)) 
\end{align} 
with the kernel equal to $\bR^{n \ast}$, a dual space to $\bR^n$. 
Thus, we obtain a diffeomorphism \[\SLnp_{\mbv_{\tilde E}} \ra \SLn \times \bR^{n \ast} \times \bR_+.\]
We note the multiplication rules
\begin{equation}\label{pr-eqn-multrule}  
 (A, \vec{v}, \lambda) (B, \vec{w}, \mu) = (AB, \frac{1}{ \mu^{1/n} } \vec{v}B + \lambda \vec{w}, \lambda \mu).
 \end{equation}  
We denote by ${\mathcal L}_{1}: \SLnp_{\mbv_{\tilde E}} \ra \SLn$ the further projection to $\SLn$.
\index{Lone@${\mathcal L}_1$} 

Let $\Sigma_{\tilde E}$ be the end $(n-1)$-orbifold. 
Given a representation 
\[\hat h: \pi_1(\Sigma_{\tilde E}) \ra \SLn \hbox{ and a homomorphism } 
\lambda_{v_{\tilde E}}: \pi_1(\Sigma_{\tilde E}) \ra \bR_+,\] 
we denote by $\bR^{n}_{\hat h, \lambda_{v_{\tilde E}}}$
the $\bR$-module with the $\pi_1(\Sigma_{\tilde E})$-action given 
by \[g\cdot \vec v = \frac{1}{\lambda_{v_{\tilde E}}(g)^{1/n}}\hat h(g)(\vec v).\] 
And we denote by $\bR^{n \ast}_{\hat h, \lambda_{v_{\tilde E}}}$ the dual vector space
with the right dual action given by 
\[g\cdot \vec v = \frac{1}{{\lambda_{v_{\tilde E}}(g)^{1/n}}}\hat h(g)^{\ast}(\vec v).\] 
Let $H^1(\pi_1(\tilde E), \bR^{n \ast}_{\hat h, \lambda_{v_{\tilde E}}})$ denote the cohomology 
space of $1$-cocycles 
\[ \Gamma \ni g \mapsto \vec v(g) \in  \bR^{n \ast}_{\hat h, \lambda_{v_{\tilde E}}}.\]  


As $\Hom(\pi_1(\Sigma_{\tilde E}), \bR_+)$ equals $H^1(\pi_1(\Sigma_{\tilde E}), \bR)$, we obtain: 

\begin{theorem} \label{pr-thm-defspace}
Let $\orb$ be a  convex  real projective $n$-orbifold, and 
let $\torb$ be its universal cover. 
Let $\Sigma_{\tilde E}$ be the end orbifold associated with an R-p-end $\tilde E$ of $\torb$. 
Then the space of representations 
\[\Hom(\pi_1(\Sigma_{\tilde E}), \SLnp_{\mbv_{\tilde E}})/\SLnp_{\mbv_{\tilde E}}\] 
has a natural surjection to 
\[\left(\Hom(\pi_1(\Sigma_{\tilde E}), \SLn)/\SLn\right) \times H^1(\pi_1(\Sigma_{\tilde E}), \bR)\]
where each fiber isomorphic to 
$H^1(\pi_1(\Sigma_{\tilde E}), \bR^{n \ast}_{\hat h, \lambda_{v_{\tilde E}}}) $ 
for each $([\hat h], \lambda)$. 
\end{theorem} 



On a Zariski open subset of $\Hom(\pi_1(\Sigma_{\tilde E}), \SLn)/\SLn $, 
the dimensions of the fibers are constant (see Johnson-Millson \cite{JM87}).
A similar idea is given by Mess  \cite{Mess07}. 
In fact, dualizing these matrices gives us
a representation to $\Aff(\mathds{A}^n)$. (See Chapter \ref{ch-du}.)
In particular, when 
the linear parts are in $\SO(n, 1)$, then we are exactly in the cases studied by Mess. 
(The concept of the duality is explained in Section \ref{pr-sub-affdualtub}) 

%

If $\Sigma$ is a closed $2$-orbifold with negative Euler characteristic, one can compute 
the dimension of $H^1(\pi_1(\Sigma_{\tilde E}), \bR^{n \ast}_{\hat h, \lambda_{v_{\tilde E}}}) $ 
using the twisted orbifold Euler characteristic of Porti \cite{Porti20}. 

\subsection{Tubular actions.} \label{pr-sub-tubular}

Consider a pair of antipodal points $\mbv$ and $\mbv_-$. 
If a group $\Gamma$ of projective automorphisms fixes a pair of fixed points $\mbv$ and $\mbv_-$, 
then $\Gamma$ is said to be {\em tubular}.
There is a projection $\Pi_{\mbv}: \SI^n -\{\mbv, \mbv_-\} \ra \SI^{n-1}_{\mbv}$ given 
by sending every great segment with endpoints $\mbv$ and $\mbv_-$
to a point of the sphere of directions at $\mbv$. 
\index{piv@$\Pi_{\mbv}$|textbf} 
\index{snminusone@$\SI^{n-1}_\mbv$}

A {\em tube} in $\SI^n$ (resp. in $\RP^n$) is the closure of the preimage 
$\Pi^{-1}_{\mbv}(\Omega)$ of a domain $\Omega$
in $\SI^{n-1}_{\mbv}$ (resp. in $\RP^{n-1}_{\mbv}$).
We often 
denote the closure in $\SI^{n}$ by ${\mathcal T}_{\mbv}(\Omega)$, and 
we call it a {\em tube domain}. 
Given an R-p-end $\tilde E$ of $\torb$, let $\mbv := \mbv_{\tilde E}$. 
The corresponding {\em end domain} is denoted $R_{\mbv}(\torb)$.
 \index{tv@${\mathcal T}_{\mbv}(\cdot)$|textbf } 
If an R-p-end $\tilde E$ has the end domain $\tilde \Sigma_{\tilde E} = R_{\mbv}(\torb)$, 
the group $h(\pi_1(\tilde E))$ acts on $B:={\mathcal T}_{\mbv}(\Omega)$. 
We also denote by $R_{\mbv}(B)$ the space of directions of great segments in $B$ from 
$\mbv$ to its antipode. 
\index{rv@$R_{\mbv}(\torb)$} \index{sigmatildeE@$\tilde \Sigma_{\tilde E}$} 
\index{tubular action} 
\index{tube domain}

The image of the tube domain ${\mathcal T}_{\mbv}(\Omega)$ in $\RP^n$ is still called a {\em tube domain} and 
denoted by ${\mathcal T}_{[v]}(\Omega)$ where $[v]$ is the image of $\mbv$.

We now discuss for the $\SI^n$-version but the $\RP^n$ version is obviously obtained from this 
by a minor modification. 

Letting $\mbv$ have the coordinates $[0, \dots, 0, 1]$, we obtain 
the matrix of $g$ of $\pi_1(\tilde E)$ of form 
\begin{equation}\label{pr-eqn-bendingm3} 
\left(
\begin{array}{cc}
\frac{1}{\lambda_{\mbv}(g)^{\frac{1}{n}}} \hat h(g)          &       0                \\
\vec{b}_g           &      \lambda_{\mbv}(g)                  
\end{array}
\right)
\end{equation}
where $\vec{b}_g$ is an $n\times 1$-vector and $\hat h(g)$ is an $n\times n$-matrix of determinant $\pm 1$
and $\lambda_{\mbv}(g) $ is a positive constant. 

Note that the representation $\hat h: \pi_1(\tilde E) \ra \SLn$ is given by 
$g \mapsto \hat h(g)$. Here we have $\lambda_{\mbv}(g) > 0$.  
If $\tilde \Sigma_{\tilde E}$  is properly convex, 
then the convex tubular domain and the action is said to be {\em properly tubular}.
\index{tubular action}
\index{tubular action!properly}

\subsection{Affine actions  dual to tubular actions.}\label{pr-sub-affdualtub}






Let ${\SI^{n-1}}$ in $\SI^{n} = \SI(\bR^{n+1})$ be a great sphere of dimension $n-1$.
A component of the complement of ${\SI^{n-1}}$
can be identified with an affine space $\mathds{A}^{n}$. 
The subgroup of projective automorphisms preserving ${\SI^{n-1}}$ and the components of the complements is equal to
the affine group $\Aff(\mathds{A}^n)$.

By duality, a great $(n-1)$-sphere ${\SI^{n-1}}$ corresponds to a point $\mbv_{\SI^{n-1}}$. 
Thus, for a group $\Gamma$ in $\Aff(\mathds{A}^n)$, 
its dual group $\Gamma^*$ acts on $\SI^{n\ast}:=\SI(\bR^{n+1, *})$ fixing $\mbv_{\SI^{n-1}}$.
(See Proposition \ref{prelim-prop-duality} also.)

Let $\SI^{n-1}_\infty$ denote a hyperspace in $\SI^n$. 
Suppose that $\Gamma$ acts on a properly convex open domain $U$ where $\Omega := \Bd U \cap {\SI^{n-1}_\infty}$
is a properly convex domain. 
We recall that $\Gamma$ has a properly convex affine action. 
Let us recall some facts from Section \ref{prelim-sub-Eduality}
\begin{itemize}
\item A great $(n-2)$-sphere $P \subset \SI^{n}$ is dual to a great circle $P^{\ast}$ in $\SI^{n\ast}$
given as the set of hyperspheres containing $P$. 
\item The great sphere $\SI^{n-1}_{\infty} \subset \SI^{n}$ with an orientation is dual to a point $\mbv \in \SI^{n\ast}$ 
and it with an opposite orientation is dual to $\mbv_{-}\in \SI^{n\ast}$. 
\item An oriented hyperspace $P \subset \SI^{n-1}_{\infty}$ of dimension $n-2$ is dual to an oriented great circle 
passing $\mbv$ and $\mbv_{-}$, giving us an element $P^{\dagger}$ of the linking sphere 
$\SI^{n-1\ast}_{\mbv}$ of rays from $\mbv$ in $\SI^{n\ast}$.  
\item The space $S$ of oriented hyperspaces in $\SI^{n-1}_{\infty}$ equals 
$\SI^{n-1 \dagger}_{\infty}$. 
Thus, there is a projective isomorphism 
\[\mathcal{I}_{2}: S= \SI^{n-1\ast}_{\infty} \ni P \leftrightarrow P^{\dagger} \in \SI^{n-1\ast}_{\mbv}.\]
\end{itemize} 
\index{i2@$\mathcal{I}_{2}$|textbf}


For the following, we say that an oriented hyperspace $V$ in $\SI^{i}$ {\em supports} 
an open submanifold $A$ if  it bounds an open $i$-hemisphere $H$ in the right orientation containing $A$. 
In projective geometry, a {\em pencil} is a one-parameter family of hyperspaces containing 
a codimension-two subspace. 

\begin{proposition}\label{pr-prop-dualtube}
Suppose that $\Gamma \subset \SLnp$ acts on a properly convex open domain $\Omega \subset {\SI^{n-1}_\infty}$
cocompactly.
Then the dual group $\Gamma^*$ acts on a properly tubular
domain $B$ with vertices $\mbv:= \mbv_{\SI^{n-1}_\infty}$ and $\mbv_- := \mbv_{{\SI^{n-1}_\infty}, -}$ dual to 
$\SI^{n-1}_{\infty}$. Moreover, 
the domain $\Omega$ and domain $R_{\mbv}(B)$ in the linking sphere $\SI^{n-1}_{\mbv}$ from $\mbv$ in the directions of $B^o$
are projectively diffeomorphic to a pair of dual domains in $\SI_{\infty}^{n-1}$ respectively. 
\end{proposition} 
 \begin{proof} 
Given $\Omega \subset \SI^{n-1}_\infty$, we obtain the properly convex open dual domain 
$\Omega^{\ast}$ in $\SI^{n-1 \dagger}_\infty$. 
An oriented great $(n-2)$-sphere sharply supporting $\Omega$  
in $\SI^{n-1}_\infty$  corresponds to 
a point of $\Bd \Omega^{\ast}$ and vice versa. (See Section \ref{prelim-sec-duality}.) 
An oriented great $(n-1)$-sphere in $\SI^{n}$ supporting $\Omega$
but not containing $\Omega$ 
meets a great $(n-2)$-sphere $P$ in $\SI^{n-1}_{\infty}$
supporting $\Omega$. 
The dual $P^{*}$ of $P$ is the set of hyperspaces containing $P$, 
a great circle in $\SI^{n\ast}$. 
The set of oriented great $(n-1)$-spheres containing $P$ supporting $\Omega$
but not containing $\Omega$ 
forms a pencil; in this case, a great open segment $I_{P^{\ast}}$ in $\SI^{n \ast}$ with endpoints $\mbv$ and $\mbv_-$.  
(See Section \ref{prelim-sub-Haus} for the definition of supporting hyperspaces.)
Let $P^{\ddagger} \in \SI^{n-1\dagger}_{\infty}$ denote the dual of $P$ in  $\SI^{n-1}_{\infty}$.
Then $P^{\dagger}:= \mathcal{I}_{2}(P^{\ddagger})$ is the direction of $P^{\ast}$ at $\mbv$ as we can see from the projective isomorphism 
$\mathcal{I}_{2}$. 
Recall from the beginning of Section \ref{prelim-sub-duality} 
$P$ supports $\Omega$ if and only if $P^{\ddagger} \in \Omega^\ast$.
Hence, there is a homeomorphism
\begin{multline} \label{pr-eqn-dualsupp}
I_{P} := \{ Q | Q \hbox{ is an oriented great $(n-1)$-sphere supporting but not meeting } \clo(\Omega), \\ Q \cap \SI^{n-1}_{\infty }= P\}  \rightarrow \nonumber \\
S_{P^{\ast}} = \{ p| p \hbox{  is a point of a great open segment in } P^\ast 
\hbox{ with endpoints } \mbv, \mbv_{-}  \nonumber \\ 
 \hbox{ where  the direction $P^{\dagger}= \mathcal{I}_{2}(P^{\ddagger}), P^{\ddagger} \in \Omega^\ast$}\}.
\end{multline}

The set $B$ of oriented hyperspaces supporting $\Omega$ possibly containing
$\Omega$ meets an oriented $(n-2)$-hyperspace in $\SI^{n-1}_{\infty}$
supporting $\Omega$. Denote by 
$\alpha_{x}$ the great segment with vertices $\mbv$ and 
$\mbv_-$ in the direction of $x \in \SI^{n-1}_\mbv$. 
Thus, we obtain
 \[ B^{\ast} = \bigcup_{P \in \Omega^\ast} S_{P^{\ast}} 
=\bigcup_{x \in \mathcal{I}_2(\Omega^\ast)} \alpha_{x}
\subset \SI^{n\ast}.\]
Let ${\mathcal T}(\Omega^{\ast})$ denote the union of open great segments with endpoints $\mbv$ and $\mbv_{-}$ 
in direction of $\Omega^{\ast}$. Thus, $B^{\ast} = {\mathcal T}(\Omega^{\ast})$.
Thus, there is a homeomorphism
\begin{align} \label{pr-eqn-supp}
I&:= \{ Q | Q \hbox{ is an oriented great $(n-1)$-sphere sharply  supporting } \Omega\} \rightarrow  \nonumber \\
S &= \{ p| p \in S_{P^{\ast}}, P^{\ddagger} \in \Bd \Omega^{\ast}\} = \Bd B^{\ast} -\{\mbv, \mbv_{-}\}.
\end{align} 
Also, $R_{\mbv}(B^{\ast}) = \mathcal{I}_2(\Omega^{\ast})$ by the above equation. 
Thus, the following are equivalent: 
\begin{itemize} 
\item  $\Gamma$ acts on $\Omega$. 
\item $\Gamma$ acts on $I$.
\item $\Gamma^{\ast}$ acts on $S$. 
\item $\Gamma^{\ast}$ acts on $B^{\ast}$ and on $\Omega^{\ast}$.
\end{itemize} 
\end{proof} 



\subsection{Distanced tubular actions and asymptotically nice affine actions.} 


The approach is similar to what we did in Chapter \ref{ch-du}
but is in the dual setting. 

\begin{definition}\label{pr-defn-tubular} $ $ 
\begin{description}
\item[Radial action:] A properly tubular action of $\Gamma$ is said to be {\em distanced} if a $\Gamma$-invariant tubular domain contains 
a properly convex compact $\Gamma$-invariant subset disjoint from the vertices of 
the tubes. 
\index{distanced action|textbf} 
\item[Affine action:] We recall from Chapter \ref{ch-du} the following fact: 
A properly convex affine action of $\Gamma$ is said to be {\em asymptotically nice} if 
$\Gamma$ acts on a properly convex open domain $U'$ in $\mathds{A}^n$ with boundary in 
$\Omega \subset {\SI^{n-1}_\infty}$, 
and $\Gamma$ acts on a compact subset
\[J:= \{ H| H \hbox{ is an AS-hyperspace passing through } x \in \Bd \Omega, H \not\subset {\SI^{n-1}_\infty}\}\] 
 where we require that every sharply supporting $(n-2)$-dimensional space of $\Omega$ in ${\SI^{n-1}_\infty}$ is 
 contained in at least one of the element of $J$. 
 \index{asymptotically nice action} 
\end{description}
\end{definition}


The following is a simple consequence of the homeomorphism 
given by equation \eqref{pr-eqn-supp}.
\begin{proposition}\label{pr-prop-dualDA} 
Let $\Gamma$ and $\Gamma^*$ be a pair of groups dual to each other 
where $\Gamma$ has an affine action on $\mathds{A}^n$ and $\Gamma^*$ is tubular with 
the vertex $\mbv = \mbv_{\SI^{n-1}_\infty}$ dual to the boundary $\SI^{n-1}_\infty$ of $\mathds{A}^n$.
Let $\Gamma= (\Gamma^*)^*$ acts on a convex open domain $\Omega$ with 
a closed $n$-orbifold $\Omega/\Gamma$.
Then $\Gamma$ act asymptotically nicely if and only if 
$\Gamma^*$ acts on a properly tubular domain $B$ and this action is distanced. 
\end{proposition}
\begin{proof} 
From the definition of asymptotic niceness, we can obtain the following: 
for each point $x$ and a sharply supporting hyperspace $P$ of $\Bd \Omega$ passing $x$ in $\SI^{n-1}$, we choose 
a great $(n-1)$-sphere in $\SI^n$ sharply supporting $\Omega$ at $x$ containing $P$ and 
uniformly bounded at a distance in the $\bdd_H$-sense from ${\SI^{n-1}_\infty}$. 
This forms a compact $\Gamma$-invariant set $J$ of hyperspaces. 

Let $B$ be a tube domain with the vertex $\mbv$ dual to ${\SI^{n-1}_\infty}$.
Then by  \eqref{pr-eqn-supp} (see the proof of Proposition \ref{pr-prop-dualtube})
the dual of the hyperspaces corresponds to points of $B$. 
The hyperspaces in $J$ sharply supporting $U$ at $x \in \Bd \Omega$, 
are bounded at a distance from ${\SI^{n-1}_\infty}$  in the $\bdd_H$-sense. 
The dual points are uniformly bounded away
at a distance from the vertices $\mbv$ and $\mbv_-$. We take the closure of the set of 
hyperspaces in the dual space of $\SI^{n\ast}$. Let us call this compact set $K$. 
Let $\Omega^* \subset \SI^{n-1}_{\mbv}$ be the dual domain of $\Omega$. 
Then for every point of $\Bd \Omega^*$,
we have a point of $K$ in the corresponding great segment from $\mbv$ to $\mbv_-$. Then
$K$ is uniformly bounded at a distance from $\mbv$ and $\mbv_-$  in the $\bdd$-sense.
The convex hull of $K$ in $\clo(\torb)$  is a compact convex set bounded at a uniform distance from $\mbv$ and $\mbv_-$
since the tube domain is properly tubular and convex:
To explain more, we can just consider an affine patch $A$ containing $K$ with origin $v$ 
and then the tube domain meets in $A$ in a strictly convex cone, which implies the result. 
Since $K$ is $\Gamma^{\ast}$-invariant, so is the convex hull in $\clo(\torb)$. 
Therefore, $\Gamma^{\ast}$ acts on $B$ as a distanced action.

The converse also follows directly by \eqref{pr-eqn-supp} in the proof of Proposition \ref{pr-prop-dualtube}. We can take the intersection $U$ of the inner components of hyperspaces involved here and obtain an open set. Also, $U$ is not empty by an elementary geometric argument  since the angles between $\SI^{n-1}_\infty$ and the strictly supporting hyperspaces are uniformly bounded below. 
Therefore, we conclude that the action is the asymptotically nice.
\end{proof}


\begin{theorem}\label{pr-thm-distanced}
Let $\Gamma$ have a nontrivial properly convex tubular action at vertex $\mbv = \mbv_{\SI^{n-1}_\infty}$ on 
$\SI^n$ {\rm (}resp. in $\RP^n${\rm )} 
and acts on a properly convex tube $B$
and satisfies the uniform middle-eigenvalue condition with respect to $\mbv_{{\SI^{n-1}_\infty}}$. 
We assume that $\Gamma$ acts on a convex open domain $\Omega \subset \SI^{n-1}_{\mbv}$ 
where $B = \mathcal{T}_{\mbv}(\Omega)$ and $\Omega/\Gamma$ is a closed $n$-orbifold. 
Then $\Gamma$ is distanced inside the tube $B$, 
and $B$ contains a distanced $\Gamma$-invariant compact set $K$. 
Finally, if $\Gamma$ is virtually factorizable, we can choose the distanced set 
$K$ to be in a hypersphere disjoint from $\mbv, \mbv_{-}$ 
\end{theorem} 
\begin{proof}
We will again prove for $\SI^n$. 
Let $\Omega$ denote the convex domain in $\SI^{n-1}_{\mbv}$ corresponding to $B^{o}$. 
By Theorems \ref{du-thm-asymnice} and \ref{du-thm-asymniceII}, 
$\Gamma^*$ is asymptotically nice. 
Proposition \ref{pr-prop-dualDA} implies the result.

Now, we prove the final part to show the total geodesic property of virtually factorizable ends: 
Suppose that $\Gamma$ acts virtually reducibly on $\SI^{n-1}_{\mbv}$ on a properly convex domain $\Omega$. 
Then $\Gamma$ is virtually isomorphic to a cocompact subgroup of 
 \[\bZ^{l_0-1} \times \bGamma_1 \times \dots \times \bGamma_{l_0}\]
where $\bGamma_i$ is irreducible
by Proposition \ref{prelim-prop-Ben2}.  Also, 
$\Gamma$ acts on \[K:= K_{1}\ast \cdots \ast K_{l_{0}}= \clo(\Omega) \subset \SI^{n-1}_{\mbv} \]
where $K_{i}$ denotes the properly convex compact set in $\SI^{n-1}_{\mbv}$ where $\bGamma_{i}$ acts on for each $i$. 
Here, $K_{i}$ is $0$-dimensional for $i=s+1, \dots, l_0$ for $s+1 \leq l_0$.  
Let $B_{i}$ be the convex tube with vertices $\mbv$ and $\mbv_{-}$ corresponding to $K_{i}$.
Each $\bGamma_i$ for $i = 1, \dots, s$ acts on a nontrivial tube $B_i$ with vertices $\mbv$ and $\mbv_-$ in a subspace. 


For each $i$, $s+1 \leq i \leq r$,  $B_i$ is a great segment with endpoints $\mbv$ and $\mbv_-$. 
A point $p_i$ corresponds to $B_i$ in $\SI^{n-1}_{\mbv}$. 




The virtual center isomorphic to $\bZ^{l_0-1}$ is in the group
$\Gamma$ by Proposition \ref{prelim-prop-Ben2}. 
Recall that a nontrivial element $g$ of the virtual center acts 
trivially on the subspace $K_{i}$ of $\SI^{n-1}_{\mbv}$; that is, 
$g$ has only one associated eigenvalue in points of $K_{i}$. 
There exists a nontrivial element $g$ of 
the virtual center whose the largest norm eigenvalue is associated with the direction 
$K_{i}$ for the induced 
$g$-action on $\SI^{n-1}_{\mbv}$ 
since the action of $\Gamma$ on $\Omega$ is cocompact.  
By the middle-eigenvalue condition, 
for each $i$, we can find $g$ in the center such that $g$ has a hyperspace $K'_{i} \subset B_{i}$ corresponding to 
the  largest norm eigenvalues. 
%
%
The convex hull of
\[K'_1 \cup \cdots \cup K'_{l_0}\] 
in $\clo(B)$ is a distanced $\Gamma$-invariant compact convex set. 
For $(\zeta_1, \cdots, \zeta_{l_0}) \in \bR_+^{l_0}$, we define
\begin{equation} \label{pr-eqn-diag} 
\zeta(\zeta_1, \cdots, \zeta_{l_0}):= 
\left(\begin{array}{ccccc}
\zeta_{1} \Idd_{n_1+1} & 0 & \cdots & 0 \\
0   & \zeta_{2} \Idd_{n_2+1} & \cdots & 0 \\ 
\vdots & \vdots  &\ddots    & \vdots \\ 
0 &   0 & \cdots & \zeta_{l_0} \Idd_{n_{l_0}+1} 
\end{array}\right), 
\zeta_{1}^{n_1+1}\zeta_{2}^{n_2+1}\cdots\zeta_{l_0}^{n_{l_0}+1} =1, 
\end{equation}
using the coordinates where each $K_i$ corresponds to a block. 

Now, we consider the general case. 
The element $x$ of $K^o \subset \SI^{n-1}_{\mbv}$ has coordinates 
\[\llrrparen{\lambda_1, \dots, \lambda_{l_0}, \vec{x}_1, \dots, \vec{x}_{l_0}}, 
\hbox{ where } 
\sum_{i=1}^{l_0} \lambda_i = 1,
x = \llrrparen{\sum_{i=1}^{l_0} \lambda_i \vec{x}_i}\] 
for $\vec{x}_i$ is a unit vector in the direction of $K_i^o$
for $i=1, \dots, l_0$.

For any subgroup $G$ of $\SL_{\pm}(n+1, \bR)$, 
let $\mathcal{Z}(G)$ denote the Zariski closure in 
$\SL_{\pm}(n+1, \bR)$. 
\index{Zariski closure} 
\index{z@$\mathcal{Z}(\cdot)$|textbf}

Let $\Gamma'$ denote the finite-index normal subgroup acting on 
each of $K_i$ in $\Gamma$. 
Then the Zariski closure of $\Gamma'$ is isomorphic to 
\[\bR^{l_0-1} \times {\mathcal{Z}}(\bGamma_1) \times \cdots {\mathcal{Z}}(\bGamma_{l_0})\]
where ${\mathcal{Z}}(\bGamma_i)$ is the Zariski closure of $\bGamma_i$
easily derivable from Theorem 1.1 of Benoist \cite{Benoist03}
for our setting. 
The elements of $\bR^{l_0}$ commute with elements of 
$\bGamma_i$ and hence with $\Gamma'$. 

Let $\bZ^{l_0}$ be the virtual center of $\bGamma$, 
as obtained by Proposition \ref{prelim-prop-Ben2}.  
There is a linear map $Z: \bZ^{l_0-1} \ra \bR^{l_0}$ such that 
an isomorphism $\bZ^{l_0-1} \ra \Gamma'$ is 
represented by $\zeta \circ \exp \circ Z$. 

Let $\log \lambda_1: \bZ^{l_0-1} \ra \bR$ denote 
a map given by taking the log of the largest norm of the elements,  
$\log \lambda_n: \bZ^{l_0-1} \ra \bR$
given by taking the log of the smallest norm of the elements, 
and $\log \lambda_{\mbv}: \bZ^{l_0-1} \ra \bR$ 
the log of the eigenvalue at $\mbv$ of the elements. 
Now, $\log \lambda_1$ and $\log \lambda_n$ extends to
piecewise linear functions on $\bR^{l_0-1}$ 
that are linear over cones with origin as the vertex.

$\log \lambda_1$ has only nonnegative values, while 
$\log \lambda_n$ has nonpositive values.
The uniform middle-eigenvalue condition for $\bZ^{l_0-1}$
is equivalent to the condition
that $\log \lambda_1 > \log \lambda_{\mbv} > \log \lambda_n$ 
holds over $\bR^n-\{O\}$.

Let $B_i$ denote the tube $\mathcal{T}_{\mbv}(K_i)$. 
We choose an element $g$ of the virtual center having the largest norm of 
the eigenvalue at points of $K_i$ as an automorphism of 
$\SI^{n-1}_{\mbv}$. 
$g$ acts on $B_i$. 
By the uniform middle-eigenvalue condition, 
$g$ fixes a subspace $\hat K_i$ equal to $B_i \cap P_i$ 
where $P_i$ is a hyperspace  in the span of $B_i$
corresponding to the largest norm eigenvalue of $g$
as an element of $\SLnp$. 
By commutativity, the virtual center also acts on $\hat K_i$. 
Define $P$ as the join of $P_1, \dots, P_{l_0}$.
Hence, the center acts on 
the join $\hat K_1\ast \cdots \ast \hat K_{l_0}$, which 
equals $\mathcal{T} \cap P$ for a hyperspace $P$.
By commutativity, $\Gamma'$ also acts on $B_i$. 

Suppose that for some $g \in \Gamma$, $g(P) \ne P$. 
Then $g(P) \cap B_j$ contains a point $x$ closer to $\mbv$ or $\mbv_-$ 
than any points on $P \cap B_j$ for some $j$. 
Assume that it is closer to $\mbv$ without loss of generality. 
We find a sequence $\{\vec{k}_i\}$ such that 
$g_i = \zeta \circ \exp \circ Z(\vec{k}_i)$ has the largest
eigenvalue at points of $B_i$ and $\lambda_1(g_i)/\lambda_{\mbv}(g_i) \ra \infty$.
Since $\lambda_n(g_i^{-1}) = \lambda_1(g_i)^{-1}$ and 
$\lambda_{\mbv}(g_i^{-1}) = \lambda_{\mbv}(g_i)^{-1}$, 
we obtain that 
$\{g_i^{-1}(x)\} \ra \mbv$ as $i \ra \infty$.   
Then we obtain that $g_i(g(P)) \cap \mathcal{T}$ is not distanced. 
This contradicts the first paragraph of the proof. 

Hence, $\Gamma$ acts on the hyperspace $P$. Hence, letting $K = B \cap P$ completes the proof. 
\end{proof}






\section{The characterization of lens-shaped representations} \label{pr-sec-chlens}

The main purpose of this section is to characterize the lens-shaped representations
in terms of eigenvalues, a major result of this monograph.

First, we prove the eigenvalue estimation in terms of lengths for non-virtually-factorizable and hyperbolic ends. 
We show that the uniform middle-eigenvalue condition implies the existence of limit sets. 
This proves Theorem \ref{pr-thm-equiv}.
Finally, we prove  in Theorem \ref{pr-thm-equ} that the lens condition 
and the uniform middle-eigenvalue condition are equivalent 
for both R-ends and T-ends under very general conditions. In particular, this establishes
Theorem \ref{pr-thm-secondmain}.


Techniques here are somewhat related to the works of Guichard-Wienhard \cite{GW12}
and Benoist \cite{Benoist00}.


\subsection{The eigenvalue estimations}\label{pr-sub-eigen} 

Let $\orb$ be a convex real projective $n$-orbifold 
and $\torb$ be the universal cover in $\SI^n$. 
Let $\tilde E$ be a properly convex R-p-end of $\torb$, and let $\mbv_{\tilde E}$ be the p-end vertex. 
Let \[h: \pi_1(\tilde E) \ra \SLnp_{\mbv_{\tilde E}}\] be a homomorphism, and suppose that $\pi_1(\tilde E)$ is hyperbolic. 


In this article, we assume that $h$ satisfies the middle-eigenvalue condition.  
We denote the norms of eigenvalues of $g$ by
\begin{equation}
\lambda_1(g), \ldots , \lambda_n(g), \lambda_{\mbv_{\tilde E}}(g), \hbox{where } \lambda_1(g) \cdots \lambda_n(g) \lambda_{\mbv_{\tilde E}}(g)= \pm 1, \hbox{ and }
\lambda_1(g)\geq \ldots \geq \lambda_n(g),
\end{equation} 
where we allow repetitions. 

Recall the linear part homomorphism 
${\mathcal L}_{1}$ from the beginning of Section \ref{pr-sec-endth}. 
We denote by $\hat h: \pi_1(\tilde E) \ra \SLn$ the homomorphism 
${\mathcal L}_1 \circ h$. Since $\hat h$ is a holonomy of a closed convex real projective $(n-1)$-orbifold,
and $\Sigma_{\tilde E}$ is assumed to be properly convex, 
$\hat h(\pi_1(\tilde E))$ divides a properly convex domain $\tilde \Sigma_{\tilde E}$ in $\SI^{n-1}_{\mbv_{\tilde E}}$.

We denote by $\tilde \lambda_1(g), \ldots, \tilde \lambda_n(g)$ the norms of eigenvalues of 
$\hat h(g)$ such that 
\begin{equation} 
\tilde \lambda_1(g) \geq  \ldots \geq \tilde \lambda_n(g), \tilde \lambda_1(g)   \ldots \tilde \lambda_n(g) = \pm 1 
\end{equation} 
hold.
These are called the {\em relative norms of eigenvalues} of $g$.
Note that \[\lambda_i(g) =\frac{\tilde \lambda_i(g)}{\lambda_{\mbv_{\tilde E}}(g)^{\frac{1}{n}}} 
\hbox{ for } i=1, \dots, n.\]  

For each nontorsion element $g$, eigenvalues corresponding to 
\[\lambda_1(g), \tilde \lambda_1(g), \lambda_n(g), \tilde \lambda_n(g), \lambda_{\mbv_{\tilde E}}(g)\]
are all positive, and $g$ is positive bi-semiproximal  by Proposition \ref{prelim-prop-nonu}.
(See also Theorem \ref{prelim-thm-semi}.) 
We define 
\[\leng_{\Sigma_{\tilde E}}(g):= \log\left(\frac{\tilde\lambda_1(g)}{\tilde \lambda_n(g)}\right) = \log\left(\frac{\lambda_1(g)}{\lambda_n(g)}\right).\]
This equals the infimum of the Hilbert metric lengths of the associated closed curves in $\tilde \Sigma_{\tilde E}/\hat h(\pi_1(\tilde E))$
as first shown by Kuiper. (See Proposition 2.1  of \cite{CLT15} or  \cite{Benoist97} for example.)

We recall the notions in 
Section \ref{prelim-sub-semi}. (See 
\cite{Benoist97} and \cite{Benoist00} also.)



When $\pi_1(\tilde E)$ is hyperbolic, 
all infinite order elements of $\hat h(\pi_1(\tilde E))$ are positive bi-proximal
and a finite-index subgroup has only positive bi-proximal elements and the identity.


Assume that $\bGamma_{\tilde E}$ is hyperbolic. 
Suppose that $g \in \bGamma_{\tilde E}$ is proximal. 
We define 
\begin{equation}\label{pr-eqn-alphabetag}
\alpha_g := \frac{\log \tilde \lambda_1(g)- \log \tilde \lambda_n(g)}{\log \tilde \lambda_1(g) - \log \tilde \lambda_{n-1}(g)}, 
\beta_g :=   \frac{\log \tilde \lambda_1(g)- \log \tilde \lambda_n(g)}{\log \tilde \lambda_1(g) - \log \tilde \lambda_{2}(g)},
\end{equation} 
and denote by $\bGamma_{\tilde E}^p$ the set of proximal elements. We define
\[\beta_{\bGamma_{\tilde E}} := \sup_{g \in \bGamma_{\tilde E}^p} \beta_g, 
\alpha_{\bGamma_{\tilde E}} := \inf_{g\in \bGamma_{\tilde E}^p} \alpha_g. \]
Proposition 20 of Guichard \cite{Guichard05} shows that  
we have 
\begin{equation}\label{pr-eqn-betabound}
1 < \alpha_{\tilde \Sigma_{\tilde E}} \leq \alpha_\Gamma \leq 2 \leq \beta_\Gamma \leq \beta_{\tilde \Sigma_{\tilde E}} < \infty 
\end{equation}
for constants $\alpha_{\tilde \Sigma_{\tilde E}}$ and $\beta_{\tilde \Sigma_{\tilde E}}$ depending only on $\tilde \Sigma_{\tilde E}$
since $\tilde \Sigma_{\tilde E}$ is properly and strictly convex.
\index{alpha@$\alpha_{\bGamma_{\tilde E}}$|textbf}  \index{bg@$\beta_{\bGamma_{\tilde E}} $|textbf } 

Here, it follows that $\alpha_{\bGamma_{\tilde E}}$ and $\beta_{\bGamma_{\tilde E}}$
depend on $\hat h$, and they form positive-valued functions on the union of components of  
\[\Hom(\pi_1(\tilde E), \SL_{\pm}(n, \bR))/\SL_{\pm}(n, \bR)\] 
consisting of convex dividing representations
with the algebraic convergence topology as given by Benoist \cite{Benoist05}. 




\begin{theorem}\label{pr-thm-eignlem} 
Let $\orb$ be a convex real projective $n$-orbifold. 
Let $\tilde E$ be a properly convex R-p-end of the universal cover $\torb$, 
where $\torb \subset \SI^n$, $n \geq 2$.
Let $\bGamma_{\tilde E}$ be a hyperbolic group. 
Then for every nonelliptic element $g \in \hat h(\pi_1(\tilde E))$, we have
\[ \frac{1}{n}\left(1+ \frac{n-2}{\beta_{\bGamma_{\tilde E}}}\right) \leng_{\Sigma_{\tilde E}}(g)
 \leq \log \tilde \lambda_1(g)   \leq  \frac{1}{n}\left(1+ \frac{n-2}{\alpha_{\bGamma_{\tilde E}}}\right) \leng_{\Sigma_{\tilde E}}(g).\]
\end{theorem}
\begin{proof} 
Since there is a positive bi-proximal subgroup of finite index, we concentrate on positive bi-proximal elements only.
We obtain from above that 
\[ \frac{\log \frac{\tilde \lambda_1(g)}{\tilde \lambda_n(g)}}{\log \frac{\tilde \lambda_1(g)}{\tilde \lambda_2(g)}} 
\leq \beta_{\bGamma_{\tilde E}}.\] 
We deduce that 
\begin{equation}\label{pr-eqn-eigratio} 
\frac{\tilde \lambda_1(g)}{\tilde \lambda_2(g)} \geq \left( \frac{\lambda_1(g)}{\lambda_n(g)} \right)^{\frac{1}{\beta_{\bGamma_{\tilde E}}}}
=  \left( \frac{\tilde \lambda_1(g)}{\tilde \lambda_n(g)} \right)^{\frac{1}{\beta_{\bGamma_{\tilde E}}}} = \exp\left(\frac{\leng(g)}{\beta_{\bGamma_{\tilde E}}}\right).
\end{equation}
Since we have $\tilde \lambda_i \leq \tilde \lambda_2 $ for $i\geq 2$, we obtain
\begin{equation}\label{pr-eqn-betab} 
\frac{\tilde \lambda_1(g)}{\tilde \lambda_i(g)} \geq \left( \frac{\lambda_1(g)}{\lambda_n(g)} \right)^{\frac{1}{\beta_{\bGamma_{\tilde E}}}}
\end{equation}
and since $\tilde \lambda_1 \cdots \tilde \lambda_n = 1$, 
we have 
\[ \tilde \lambda_1(g)^n = \frac{\tilde \lambda_1(g)}{\tilde \lambda_2(g)} \cdots  \frac{\tilde \lambda_1(g)}{\tilde \lambda_{n-1}(g)}
 \frac{\tilde \lambda_1(g)}{\tilde \lambda_n(g)} \geq \left(  \frac{\tilde \lambda_1(g)}{\tilde \lambda_n(g)} \right)^{\frac{n-2}{\beta_{\bGamma_{\tilde E}}} + 1}.\] 
 We obtain 
 \begin{equation}\label{pr-eqn-betabd}
  \log \tilde \lambda_1(g) \geq \frac{1}{n}\left(1+ \frac{n-2}{\beta_{\bGamma_{\tilde E}}}\right) \leng_{\Sigma_{\tilde E}}(g).
  \end{equation}
By similar reasoning, we also obtain 
\begin{equation}\label{pr-eqn-alphabd}
\log \tilde \lambda_1(g) \leq \frac{1}{n}\left(1+ \frac{n-2}{\alpha_{\bGamma_{\tilde E}}}\right) \leng_{\Sigma_{\tilde E}}(g). 
\end{equation}
\end{proof}

\begin{remark} \label{pr-rem-eigenlem}
Under the assumption of Theorem \ref{pr-thm-eignlem}, if we do not assume that $\pi_1(\tilde E)$ is hyperbolic, then 
we obtain 
\begin{equation} \label{pr-eqn-eigenlem}
 \frac{1}{n} \leng_{\Sigma_{\tilde E}}(g) \leq \log \tilde \lambda_1(g)   \leq  
\frac{n-1}{n} \leng_{\Sigma_{\tilde E}}(g)
 \end{equation}
for every semiproximal element $g \in \hat h(\pi_1(\tilde E))$.
\end{remark} 
\begin{proof} 
Let $\tilde \lambda_i(g)$ denote the norms of $\hat h(g)$ for $i=1, 2, \dots, n$. 
\[\log \tilde \lambda_1(g) \geq  \ldots \geq \log \tilde \lambda_n(g), 
\log \tilde \lambda_1(g)  + \cdots + \log \tilde \lambda_n(g) = 0\] 
hold.
We deduce 
\begin{alignat}{3} 
\log \tilde \lambda_n(g) &=& -\log \lambda_1 - \cdots - \log \tilde \lambda_{n-1}(g) \nonumber \\
& \geq & -(n-1) \log \tilde \lambda_1 \nonumber\\ 
\log \tilde \lambda_1(g) & \geq & -\frac{1}{n-1} \log \tilde \lambda_n(g) \nonumber\\ 
\left(1+ \frac{1}{n-1}\right) \log \tilde \lambda_1(g) & \geq & \frac{1}{n-1} \log \frac{\tilde \lambda_1(g)}{\tilde \lambda_n(g)}\nonumber \\ 
\log \tilde \lambda_1(g) & \geq & \frac{1}{n} \leng_{\Sigma_{\tilde E}}(g).
\end{alignat}
We also deduce 
\begin{alignat}{3} 
-\log \tilde \lambda_1(g) & = & \log \tilde \lambda_2(g) + \cdots + \log \tilde \lambda_{n}(g) \nonumber \\
 & \geq & (n-1) \log \tilde \lambda_{n}(g) \nonumber \\ 
-(n-1) \log \tilde \lambda_{n}(g) & \geq & \log \tilde \lambda_1(g) \nonumber \\ 
(n-1) \log \frac{\tilde \lambda_1(g)}{\tilde \lambda_{n}(g)} & \geq & n \log \tilde \lambda_1(g) \nonumber \\ 
\frac{n-1}{n} \leng_{\Sigma_{\tilde E}}(g) & \geq & \log \tilde \lambda_1(g). \tag*{$\square$} 
\end{alignat} 
\end{proof}

\begin{remark}
We cannot show that the middle-eigenvalue condition implies 
the uniform middle-eigenvalue condition. This could be false.
For example, we  could obtain a sequence of elements $g_i \in \Gamma$ such that 
$\{\lambda_1(g_i)/ \lambda_{\mbv_{\tilde E}}(g_i)\} \ra 1$ while $\Gamma$ satisfies the middle-eigenvalue 
condition. Certainly, we could have an element $g$ where 
$\lambda_1(g) = \lambda_{\mbv_{\tilde E}}(g)$. 
However, even if there is no such element, we might still have 
a counter-example. 
For example, suppose that we might have 
\[\left\{\frac{\log \left(\frac{\lambda_1(g_i)}{\lambda_{\mbv_{\tilde E}}(g_i)}\right)}{\leng_{\Sigma_{\tilde E}}(g)} \right\}\ra 0.\] 
This could happen
by changing $\lambda_{\mbv_{\tilde E}}$ considered as a homomorphism 
$\pi_{1}(\Sigma_{\tilde E}) \ra \bR_{+}$.   
Such assignments are not well understood globally,
but see Benoist \cite{Benoist97}. An analogous phenomenon also appears
in the study of the Margulis space-time and diffused Margulis invariants as investigated by 
Charette, Drumm, Goldman, Labourie, and Margulis 
recently.  See \cite{GLM09})
\end{remark}


\subsection{The uniform middle-eigenvalue conditions and the orbits.} \label{pr-sub-umecorbit}

Let $\tilde E$ be a properly convex R-p-end of the universal cover $\torb$ of 
 a convex real projective $n$-orbifold $\orb$. 
Assume that $\bGamma_{\tilde E}$ satisfies the uniform middle-eigenvalue condition
with respect to $\tilde E$. 
There exists a $\bGamma_{\tilde E}$-invariant compact set 
denoted by 
$L_{\tilde E}$ distanced from $\{\mbv_{\tilde E}, \mbv_{\tilde E-}\}$
by Theorem \ref{pr-thm-distanced}. 
For the corresponding tube 
${\mathcal T}_{\mbv_{\tilde E}}(\tilde \Sigma_{\tilde E})$, 
	$L_{\tilde E} \cap \Bd {\mathcal T}_{\mbv_{\tilde E}}(\tilde \Sigma_{\tilde E})$ is a compact 
subset distanced from $\{\mbv_{\tilde E}, \mbv_{\tilde E-}\}$.
\index{LambdaE@$\Lambda_{\tilde E}$} 
Let $\CH(L)$ be the convex hull of $L$ in the tube ${\mathcal T}_{\mbv_{\tilde E}}(\tilde \Sigma_{\tilde E})$ obtained by 
Theorem \ref{pr-thm-distanced}.
Then $\CH(L)$ is a $\bGamma_{\tilde E}$-invariant 
distanced subset of ${\mathcal T}_{\mbv_{\tilde E}}(\tilde \Sigma_{\tilde E})$. 

Throughout this section, we work with $\SI^n$; however, the corresponding results for 
$\RP^n$ are fairly obvious to be obtained. 

\begin{definition} \label{pr-defn-lambdaE} 
A {\em transversal set} is a compact subset of 
${\mathcal T}_{\mbv_{\tilde E}}(\tilde \Sigma_{\tilde E}) -\{\mbv_{\tilde E}, \mbv_{\tilde E -}\}$
that meets the interior of every  great segment from $\mbv_{\tilde E}$ to $\mbv_{\tilde E -}$ in it. 
A {\em transversal  boundary set} is a compact subset $C$ of 
	$\Bd {\mathcal T}_{\mbv_{\tilde E}}(\tilde \Sigma_{\tilde E}) -\{\mbv_{\tilde E}, \mbv_{\tilde E -}\}$
that meets the interior of every  great segment from $\mbv_{\tilde E}$ to $\mbv_{\tilde E -}$ that is the intersection of a properly convex compact set with 
${\mathcal T}_{\mbv_{\tilde E}}(\tilde \Sigma_{\tilde E}) -\{\mbv_{\tilde E}, \mbv_{\tilde E -}\}$.
	We define the {\em limit set} $\Lambda_{\tilde E}$ of 
	a properly convex R-p-end $\tilde E$ as the smallest nonempty
	compact $\bGamma_{\tilde E}$-invariant transversal boundary set.
	\end{definition}

Corollary \ref{pr-cor-LambdaW}  shows that the limit set is 
well-defined. Compare it also to Definition \ref{app-defn-limitset}. 
\index{end!limit set} 
\index{p-end!limit set}


The following main result of this subsection shows that 
$\Lambda_{\tilde E}$ is characterized. 
\begin{corollary}\label{pr-cor-LambdaW} 
A compact transversal boundary $\bGamma_{\tilde E}$-invariant set $C$ in $\Bd {\mathcal T}_{\mbv_{\tilde E}}(\tilde \Sigma_{\tilde E}) -\{\mbv_{\tilde E}, \mbv_{\tilde E -}\}$  exists and is unique. 
Also, it satisfies $\CH(C)\cap (\Bd {\mathcal T}_{\mbv_{\tilde E}}(\tilde \Sigma_{\tilde E}) -\{\mbv_{\tilde E}, \mbv_{\tilde E -}\}) = C$. 
\end{corollary}
\begin{proof} 
Proposition \ref{pr-prop-orbit} shows that $C$ is independent of the choice, 
and $C$ meets each great segment for any distanced compact convex set $L$
in ${\mathcal T}_{\mbv_{\tilde E}}(\tilde \Sigma_{\tilde E})$ at a unique point. 

Since $\CH(C)\cap (\Bd {\mathcal T}_{\mbv_{\tilde E}}(\tilde \Sigma_{\tilde E})
-\{\mbv_{\tilde E}, \mbv_{\tilde E -}\})$ contains $C$ 
and is also a $\bGamma_{\tilde E}$-invariant distanced compact set,
the uniqueness part of  
Proposition \ref{pr-prop-orbit} shows that it must equal $C$.  
%
\end{proof} 

Moreover, by the 
$\bGamma_{\tilde E}$-invariance and the middle-eigenvalue condition,  
$\Lambda_{\tilde E} \cap \Bd {\mathcal T}_{\mbv_{\tilde E}}(\tilde \Sigma_{\tilde E})$ contains all attracting and repelling fixed points of $\gamma \in \bGamma_{\tilde E}$ 




\subsubsection{Hyperbolic groups}\label{pr-subsub-hyperbolic} 

We first consider the case where $\bGamma_{\tilde E}$ is hyperbolic. 
In this case, $\Omega$ strictly convex, and for any nontorsion $g\in \bGamma_{\tilde E}$, 
$g$ has exactly two fixed points in $\Bd \Omega$, an attracting and a repelling one.
This follows from  
Proposition \ref{prelim-prop-nonu}, Theorem \ref{prelim-thm-semi}  and 
the strict convexity. (See \cite{Benoist04}.) 

\begin{lemma}\label{pr-lem-attracting} 
Let $\orb$ be a convex real projective $n$-orbifold. 
Let $\tilde E$ be a properly convex R-p-end. 
Assume that $\bGamma_{\tilde E}$ is hyperbolic
 and satisfies the uniform middle-eigenvalue conditions with respect to $\tilde E$. 
\begin{itemize}
\item Suppose that $\gamma_i$ is a sequence of elements of $\bGamma_{\tilde E}$ acting on ${\mathcal T}_{\mbv_{\tilde E}}(\tilde \Sigma_{\tilde E})$. 
\item The sequence of attracting fixed points $a_i$ and the sequence of  repelling fixed points $b_i$ are such that 
$\{a_i\} \ra a_\infty$ and $\{b_i\} \ra b_\infty$ where $a_\infty, b_\infty$ are not in $\{ \mbv_{\tilde E}, \mbv_{\tilde E-}\}$.
\item Suppose that the sequence $\{\lambda_i\}$ of eigenvalues where 
$\lambda_i$ corresponds to $a_i$ converges to $+\infty$. 
\end{itemize} 
Then the point $a_\infty$ 
in $\Bd {\mathcal T}_{\mbv_{\tilde E}}(\tilde \Sigma_{\tilde E}) 
-\{ \mbv_{\tilde E}, \mbv_{\tilde E -}\}$ is
the geometric limit of $\{\gamma_i(K)\}$ for any compact subset $K \subset M$. 
\end{lemma} 
\begin{proof} 
	We may assume without loss of 
	generality that $a_\infty \ne b_\infty$
	since otherwise we replace $\{g_i\}$ with $\{gg_i\}$ where 
	$g(a_\infty) \ne b_\infty$. Proving for this case implies the general cases. 
	
Let $k_{i}$ be the inverse of the factor 
\[\min \left\{\frac{\tilde \lambda_1(\gamma_i)}{\tilde \lambda_2(\gamma_i)}, 
\frac{\tilde \lambda_1(\gamma_i)}{\lambda_{\mbv_{\tilde E}}(\gamma_i)^{\frac{n+1}{n}}}
= \frac{\lambda_{1}(\gamma_{i})}{\lambda_{\mbv_{\tilde E}}(\gamma_i)}
\right\}.\]
Then $\{k_i\} \ra 0$ by the uniform middle-eigenvalue condition and  \eqref{pr-eqn-eigratio}. 

There exists a totally geodesic sphere $\SI^{n-1}_i$ sharply supporting ${\mathcal T}_{\mbv_{\tilde E}}(\tilde \Sigma_{\tilde E})$ at $b_i$. 
$a_i$ is uniformly bounded away from $\SI^{n-1}_i$ for $i$ sufficiently large.
$\SI^{n-1}_i$ bounds an open hemisphere $H_i$ containing $a_i$ where $a_i$ 
is the attracting fixed point 
by Corollary 1.2.3 of \cite{KH95} or by Proposition \ref{prelim-prop-attract}.
For a Euclidean metric $d_{E, i}$, we have for
$\gamma_i| H_i: H_i \ra H_i$ 
\begin{equation}\label{pr-eqn-kcont}
d_{E, i}(\gamma_{i}(x), \gamma_{i}(y)) \leq k_{i} d_{E, i}(x, y), x, y \in H_{i}.
\end{equation}
Note that $\{\clo(H_i)\}$ converges geometrically to $\clo(H_\infty)$ for an open hemisphere containing $a$ in its interior. 

Actually, we can choose a Euclidean metric $d_{E, i}$ on $H_i^o$ 
such that $\{d_{E, i}| J \times J \}$ is uniformly convergent for any compact subset $J$ of 
$H_\infty$.
Hence there exists a uniform positive constant $C'$ such that 
\begin{equation}\label{pr-eqn-Cp} 
\bdd(a_i, K)  < C' d_{E_i}(a_i, K). 
\end{equation} 
for $a_{i} \in J, K \subset J$ and sufficiently large $i$. 

Since $\bGamma_{\tilde E}$ is hyperbolic, the domain $\Omega$ corresponding to ${\mathcal T}_{\mbv_{\tilde E}}(\tilde \Sigma_{\tilde E})$
in $\SI^{n-1}_{\mbv_{\tilde E}}$ is strictly convex. 
For any compact subset $K$ of $M$, the condition $K \subset M$ is equivalent to 
\[K \cap \clo(\bigcup_{i=1}^\infty \ovl{b_i\mbv_{\tilde E}} \cup \ovl{b_i\mbv_{\tilde E-}}) = \emp.\]
Since the boundary sphere $\Bd H_{\infty}$ meets ${\mathcal T}_{\mbv_{\tilde E}}(\tilde \Sigma_{\tilde E})$ in this set only
by the strict convexity of $\Omega$, we obtain $K \cap \Bd H_{\infty} = \emp$. Moreover,
 $K \subset H_\infty$ since ${\mathcal T}_{\mbv_{\tilde E}}(\tilde \Sigma_{\tilde E}) \subset \clo(H_{\infty})$. 

We have $\bdd(K, \Bd H_{\infty}) > \eps_{0}$ for $\eps_{0}> 0$. 
Thus, the distance $\bdd(K, \Bd H_i)$ is uniformly bounded by a constant $\delta$. 
$\bdd(K, \Bd H_i) > \delta$ implies that 
$d_{E_i}(a_i, K) \leq C/\delta$ for a positive constant $C> 0$
Acting by $g_i$, we obtain 
$d_{E_i}(g_i(K), a_i) \leq k_i C/\delta$ by \eqref{pr-eqn-kcont}, which implies  
$\bdd(g_i(K), a_i) \leq C' k_i C/\delta$ by \eqref{pr-eqn-Cp}.
Since $\{k_i\} \ra 0$, and $\{a_i\} \ra a$, it follows 
that $\{g_i(K)\}$ geometrically converges to $a$. 
%
 \end{proof}

\begin{lemma}\label{pr-lem-attractingII} 
	Let $\orb$ be a convex real projective $n$-orbifold. 
	Let $\tilde E$ be a properly convex R-p-end. 
	Assume that $\bGamma_{\tilde E}$ is hyperbolic, 
	and satisfies the uniform middle-eigenvalue conditions with respect to $\tilde E$. 
	Suppose that $\{\gamma_i\}$ is 
	sequence of elements of $\bGamma_{\tilde E}$ acting on 
	${\mathcal T}_{\mbv_{\tilde E}}(\tilde \Sigma_{\tilde E})$
 and forms a convergent sequence acting on $\SI^{n-1}_{\mbv_{\tilde E}}$. 
	Then any transversal boundary $\Gamma_{\tilde E}$-invariant set $L_{\tilde E}$ for $\tilde E$
	contains the geometric limit of
	any subsequence of $\{\gamma_i(K)\}$ for any compact subset $K \subset {\mathcal{T}}_{\mbv_{\tilde E}}^o$. 
	Furthermore, 
	\[A_\ast(\{\gamma_i\}) \cap {\mathcal T}_{\mbv_{\tilde E}}(\tilde \Sigma_{\tilde E}), R_\ast(\{\gamma_i\})\cap {\mathcal T}_{\mbv_{\tilde E}}(\tilde \Sigma_{\tilde E}) 
	\subset L_{\tilde E}.\]
\end{lemma} 
\begin{proof} 
	Let $z \in {\mathcal T}^{o}_{\mbv_{\tilde E}}$.
	Let $\llrrparen{z'}$ denote the element in $\Sigma_{\tilde E}$
	corresponding to the ray from $\mbv_{\tilde E}$ to $z$. 
	Let $\{\gamma_i\}$ be any sequence in $\bGamma_{\tilde E}$ 
	such that the corresponding sequence $\{\gamma_i(\llrrparen{z})\}$
	in $\Sigma_{\tilde E} \subset \SI^{n-1}_{\mbv_{\tilde E}}$ converges to a point $z'$ in 
	$\Bd \Sigma_{\tilde E} \subset \SI^{n-1}_{\mbv_{\tilde E}}$. 
	
	Clearly, a fixed point of $g \in \bGamma_{\tilde E} -\{\Idd\}$ 
	in $\Bd {\mathcal T}_{\mbv_{\tilde E}}(\tilde \Sigma_{\tilde E})  - \{\mbv_{\tilde E}, \mbv_{\tilde E -}\}$ is in $L_{\tilde E}$
	since $g$ has at most one fixed point on each open segment in the boundary. 
	Let $a_i$ and $r_i$ denote the attracting fixed points  of $\gamma_i$ respectively, 
	we may assume that  
	\[\{a_i\} \ra a, \{r_i\} \ra r \hbox{ for } a_i, r_i \in  L_{\tilde E} \]
We have $a, r \in L_{\tilde E}$ by the $\Gamma_{\tilde E}$-invariance and 
the closedness of $L_{\tilde E}$. 
	Assume $a \ne r$ first.
	By Lemma \ref{pr-lem-attracting}, we have $\{\gamma_i(z)\} \ra a$, and hence the limit 
	$z_\infty = a$. 
	
	However, it could be that $a = r$. In this case, we choose $\gamma_0 \in \bGamma_{\tilde E}$
	such that $\gamma_0(a) \ne r$. Then $\gamma_0\gamma_i$ has 
the attracting fixed point $a'_i$ 
	satisfying $\{a'_i\} \ra \gamma_0(a)$ 
	and repelling fixed points $r'_i$ satisfying $\{r'_i \}\ra r$ 
	by Lemma \ref{pr-lem-gatt}.
	This implies the first part. 
	
	As above since $\{\gamma_0 \gamma_i(z) \} \ra \gamma_0(a)$,
	we apply $\gamma_0^{-1}$ now
	to conclude $\{\gamma_i(z) \} \ra a$. 
	Thus, the limit set is contained in $L_{\tilde E}$. 
\end{proof}

\begin{lemma} \label{pr-lem-gatt}
	Let $\{g_i\}$ be a sequence of projective automorphisms 
	acting on a strictly convex domain $\Omega$ in $\SI^n$. 
	Suppose that we have a sequence of attracting fixed points 
	$\{a_i \in \Bd \Omega\} \ra a$ and a sequence of 
	repelling fixed points $\{r_i \in \Bd \Omega\} \ra r$. 
	Assume that the corresponding sequence of eigenvalues of 
	$a_i$ limits to $+\infty$ and that of $r_i$ limits to $0$. 
	Let $g$ be any projective automorphism of $\Omega$. 
	Then $\{gg_i\}$ has the sequence of attracting fixed points 
	$\{a'_i\}$ converging to $g(a)$ and the sequence of repelling 
	fixed points converging to $r$. 
\end{lemma}
\begin{proof}
	Recall that $g$ is a quasi-isometry. 
	Given $\eps >0$ and 
	a compact ball $B$ disjoint from a ball around $r$, 
	we obtain that $g g_{i}(B)$ lies within an $\eps$-ball of of $g(a)$
	for sufficiently large $i$. 
	For a choice of $B$ and sufficiently large $i$, 
	we obtain $g g_{i}(B) \subset B^{o}$. 
	Since $g g_{i}(B) \subset B^{o}$, we obtain 
	\[(g g_{i})^{n}(B) \subset (g g_{i})^{m}(B)^{o} \hbox{ for } n > m\] by induction, 
	There exists an attracting fixed point $a'_{i}$ of $g g_{i}$ in $g g_{i}(B)$. 
	Since the sequence of the diameters of $g g_{i}(B)$ converges to $0$, we obtain that $\{a'_{i}\} \ra g(a)$. 
	
	Also, given $\eps > 0$ and a compact ball $B$ disjoint from a ball around $g(a)$, 
	$g_{i}^{-1}g^{-1}(B)$ is in the ball of radius $\eps$ of $r$. 
	Similarly to the above, we obtain the needed conclusion. 
	%
\end{proof}

\subsubsection{Non-hyperbolic groups}\label{pr-subsub-nonhyperbolic} 

Now, we generalize to not necessarily hyperbolic  
$\bGamma_{\tilde E}$. 
A $\bGamma_{\tilde E}$-invariant 
distanced set $L_{\tilde E}$ contains the attracting fixed set $A_i$ and 
the repelling fixed set $R_i$ of any $g\in \bGamma_{\tilde E}$
by invariance and sequence arguments. 



\begin{lemma}\label{pr-lem-attracting2} 
	Let $\orb$ be a convex real projective $n$-orbifold. 
	Let $\tilde E$ be a properly convex R-p-end. 
	Assume that $\bGamma_{\tilde E}$ is non-hyperbolic or virtually-factorizable
	and satisfies the uniform middle-eigenvalue conditions with respect to 
	$\mbv_{\tilde E}$. 
		Suppose that $\{\gamma_i\}$ is a generalized convergence 
		sequence in  $\bGamma_{\tilde E}$ acting on 
		${\mathcal T}_{\mbv_{\tilde E}}(\tilde \Sigma_{\tilde E})$. 
	Let $L_{\tilde E}$ be a transverse boundary set for $\tilde E$. 
	Then $L_{\tilde E}$ contains the geometric limit of
	any subsequence of $\{\gamma_i(K)\}$ for any compact subset $K \subset  {\mathcal T}_{\mbv_{\tilde E}}(\tilde \Sigma_{\tilde E})^o$. 
	Furthermore, 
	\[A_\ast(\{\gamma_i\}) \cap {\mathcal T}_{\mbv_{\tilde E}}(\tilde \Sigma_{\tilde E}), 
	R_\ast(\{\gamma_i\})\cap {\mathcal T}_{\mbv_{\tilde E}}(\tilde \Sigma_{\tilde E}) 
	\subset L_{\tilde E}.\]
\end{lemma} 
\begin{proof} 
	Let $L = \CH(L_{\tilde E}) \cap {\mathcal T}_{\mbv_{\tilde E}}(\tilde \Sigma_{\tilde E})$. 
	Then $L$ is a convex set 
	uniformly bounded away from $\mbv_{\tilde E}$ and its antipode 
	by a geometric consideration. 

Given any sequence $g_i$, we can extract a convergence sequence 
$\{g_i\}$ with a convergence limit $g_\infty$. 

Suppose that $L^o=\emp$. Then $L$ is a convex domain on a hyperspace $P$ 
disjoint from $\mbv_{\tilde E}$. 
We use a coordinate system where each $\gamma \in \Gamma$ is of the form \eqref{pr-eqn-bendingm3} where $\vec{b}_g = 0$. 
By dividing $g_i$ by $\lambda_1(g_i)$ and taking a limit, we obtain that 
$g_\infty$ equals 
\begin{equation}\label{pr-eqn-gired}
\left(
\begin{array}{cc}
\hat g_\infty         &       0                \\
0        &      0                 
\end{array}
\right)
\end{equation}
by the uniform middle-eigenvalue condition
and Lemma \ref{prelim-lem-convcomp}. 
Hence $A_\ast(\{g_i\}) \subset P$. 
By Theorem \ref{prelim-thm-AR},  
\[A_\ast(\{g_i\}) \subset 
P \cap \Bd {\mathcal T}_{\mbv_{\tilde E}}(\tilde \Sigma_{\tilde E}))= L_{\tilde E}.\] 
The remaining parts are straighforward to show. 

Suppose that $L^o$ is not empty. Then $L^o \cap N_\ast(\{g_i\})= \emp$
by Lemma \ref{prelim-lem-low}.
Given any convergence sequence $\{g_i\}, g_i\in \Gamma$ converging to $g_\infty$, 
the sequence $g_i(x)$ for $x \in L$ converges to a point of $A_\ast(\{g_i\})$.

By Lemma \ref{prelim-lem-convcomp}, $\mbv_{\tilde E} \in N_\ast(\{g_i\})$ 
since $\{\lambda_{\mbv_{\tilde E}}(g_i)/\lambda(g_i)\} \ra 0$
by the uniform middle-eigenvalue condition. 
Dividing $g_i$ by $\lambda_1(g_i)$ and taking a limit, we obtain that 
$g_\infty$ equals 
\begin{equation}\label{pr-eqn-ginftyb} 
\left(
\begin{array}{cc}
\hat g_\infty         &       0                \\
\hat b        &      0                 
\end{array}
\right)
\end{equation}
by the uniform middle-eigenvalue condition
and Lemma \ref{prelim-lem-convcomp}
dualizing the proof of Lemma \ref{du-lem-attracting3}.
Here $\hat g_\infty$ is not zero; otherwise 
we have uniform convergence to $\mbv_{\tilde E}$ or $\mbv_{\tilde E-}$ for 
any compact set disjoint from $\{\mbv_{\tilde E}, \mbv_{\tilde E-}\}$
which contradicts the invariance of the set $L$. 
Since $\hat g_\infty \ne 0$, 
the image of $g_\infty$ is now a subspace 
of the same dimension as $A_\ast(\{\llrrparen{\hat g_i}\})$. 
Actually, it is a graph over $A_\ast(\{\llrrparen{\hat g_i}\})$
where the vertical direction is given by the direction to 
$\mbv_{\tilde E}$ for a linear function given by $\hat b$.

Since $\bGamma_{\tilde E}$ acts on $L$, 
$g_\infty(x) \in \clo(L) \cap \Bd {\mathcal T}_{\mbv_{\tilde E}}(\tilde \Sigma_{\tilde E})$. 
Hence, $g_\infty(L) = A_\ast(\{g_i\}) \subset L_{\tilde E}$. 
Using $\{g^{-1}_i\}$, we obtain 
$R_\ast(\{g_i\}) \subset L_{\tilde E}$. 

Any element $x \in {\mathcal T}_{\mbv_{\tilde E}}(\tilde \Sigma_{\tilde E})$ 
can be written $x = \llrrparen{\vec{v}_x}, \vec{v}_x = \vec{v}_L + c \vec{v}_{\tilde E}$
for a constant $c> 0$, a vector $\vec{v}_L$ in the direction of a point of $L$, 
and a vector $\vec{v}_{\tilde E}$ in the direction of $\mbv_{\tilde E}$.
Then \[ g_\infty(\llrrparen{\vec{v}_x}) 
= \llrrparen{g_\infty(\vec{v}_L) + c g_\infty (\vec{v}_{\tilde E})}.\]
Since $c g_\infty  (\vec{v}_{\tilde E})=0$ from \eqref{pr-eqn-ginftyb}, we obtained that 
$g_\infty({\mathcal T}_{\mbv_{\tilde E}}(\tilde \Sigma_{\tilde E})) = g_\infty(L)$. 
Since \[g_\infty({\mathcal T}_{\mbv_{\tilde E}}(\tilde \Sigma_{\tilde E}))
\subset A_\ast(\{g_i\}) \cap {\mathcal T}_{\mbv_{\tilde E}}(\tilde \Sigma_{\tilde E}),\]
and $g_\infty(L) = A_\ast(\{g_i\})$, 
we obtain the result. The final statement is also proved by taking 
the sequence $g_i^{-1}$.
	\end{proof}

In the following, $\bGamma_{\tilde E}$ can be virtually factorizable.
By Proposition \ref{pr-prop-orbit}, 
$\Lambda_{\tilde E}$ is well-defined and independent of 
the choice of $K$.
\begin{proposition}\label{pr-prop-orbit}
Let $\orb$ be a convex real projective $n$-orbifold. 
Let $\tilde E$ be a properly convex R-p-end. 
Assume that $\bGamma_{\tilde E}$ satisfies the uniform middle-eigenvalue condition
with respect to the R-p-end structure.  
Let $\mbv_{\tilde E}$ be the R-end vertex
and  $z \in {\mathcal T}_{\mbv_{\tilde E}}(\tilde \Sigma_{\tilde E})^o$. 
Let $L_{\tilde E}$ be a transversal $\bGamma_{\tilde E}$-invariant 
boundary set for $\tilde E$,  
and let $L$ be the closure of $\CH(L_{\tilde E})$. 
Then the following properties are satisfied\,{\rm :} 
\begin{enumerate} 
\item[{\rm (i)}] 
$L_{\tilde E}$ contains all the limit points of orbits of 
each compact subset of ${\mathcal T}_{\mbv_{\tilde E}}(\tilde \Sigma_{\tilde E})^o$.
$L_{\tilde E}$ contains all attracting fixed points of elements of 
$\bGamma_{\tilde E}$. 
If $\bGamma_{\tilde E}$ is hyperbolic, then 
the set of attracting fixed points is dense in the set. 
\item[{\rm (ii)}] For each segment $s$ in $\Bd {\mathcal T}_{\mbv_{\tilde E}}(\tilde \Sigma_{\tilde E})$ with 
an endpoint $\mbv_{\tilde E}$, the great segment containing $s$ meets 
$L_{\tilde E}$ 
at a unique point other than $\mbv_{\tilde E}, \mbv_{\tilde E-}$. 
That is, there is a one-to-one correspondence between $\partial \clo(\Sigma_{\tilde E})$
and $L_{\tilde E}$. 
\item[{\rm (iii)}] $L_{\tilde E}$ is homeomorphic to $\SI^{n-2}$. 
\item[{\rm (iv)}] For any $\bGamma_{\tilde E}$-distanced compact set $L'$ in 
$\Bd \mathcal{T}_{\bv_{\tilde E}}(\tilde \Sigma_{\tilde E}) - \{ \mbv_{\tilde E}, \mbv_{\tilde E-}\}$ 
meeting every great segment in ${\mathcal T}_{\mbv_{\tilde E}}(\tilde \Sigma_{\tilde E})$,  
we have $L_{\tilde E} = L'$. {\rm (}uniqueness{\rm )} 
\end{enumerate} 
\end{proposition} 
\begin{proof} 
%

We will first prove (i),(ii), and (iii) for various cases, and then prove (iv) 
all together: 

(A) Consider first when $\bGamma_{\tilde E}$ is hyperbolic.
Proposition \ref{pr-lem-attractingII} proves (i) here. 
Let $L$ be the closure of $\CH(L_{\tilde E})$, which 
is $\bGamma_{\tilde E}$-invariant. 
Let $K' = L \cap \Bd {\mathcal T}_{\mbv_{\tilde E}}(\tilde \Sigma_{\tilde E}) - \{\mbv_{\tilde E}, \mbv_{\tilde E -}\}$. 
Clearly $L_{\tilde E}\subset K'$.


Since $\bGamma_{\tilde E}$ is hyperbolic,  
each point $y$ of $\Bd \tilde \Sigma_{\tilde E} \subset \SI^{n-1}_{\mbv_{\tilde E}}$ 
is a limit point of some sequence $\{g_{i}(x)\}$ for 
$x \in \tilde \Sigma_{\tilde E}$ by \cite{Benoist04}. 
Let $l_y$ denote the segment with endpoints $\mbv_{\tilde E}$ and $\mbv_{\tilde E -}$ 
in the direction of $y$. 
At least one point of $l_y$ 
is a limit point of some subsequence of $\{g_{i}(x)\}$ by Lemma \ref{pr-lem-attracting}. 
Thus, $l_{y} \cap L_{\tilde E} \ne \emp$. 
and $l_y \cap K' \ne \emp$.

Let us choose a standard Euclidean metric $\llrrV{\cdot}_E$ for $\bR^{n+1}$. 
We identify $\SI^{n-1}_{\mbv_{\tilde E}}$ with a hyperspace $V$ not passing $\mbv_{\tilde E}$
for convenience during this proof.  We consider $V$ to correspond to $\bR^n$ 
and $\mbv_{\tilde E}$ to be the $(n+1)$-st unit vector.

We claim that $l_{y} \cap K'$ is unique: Suppose not. 
Let $z$ and $z'$ be the two points of $l_y \cap K'$.  
We choose a line $l$ in $\tilde \Sigma_{\tilde E}$ 
ending at $y$. 
Let $y_i$ be the sequence of points on $l$ converging to $y$. 
We choose $g_i$ as in the proof of 
Lemma \ref{du-lem-S_0} 
such that $g_i(y_i) \in F$ for a compact 
fundamental domain $F$ of $\tilde \Sigma_{\tilde E}$. 
We may assume that $\{g_i\}$ is a convergence sequence by passing to 
a subsequence.
(Here, $\llrrparen{\vec v_{-,i}}$ is fixed to be a single point $y=\llrrparen{\vec v_-}$
modifying the proof of Lemma \ref{du-lem-S_0}.)
Given the other endpoint $t$ of $l$, we have 
\[\{g_i(t)\} \ra a_\ast(\{g_i\})\] for 
an attractor of $a_\ast(\{g_i\})$ of $\{g_i\}$. 
This follows from the same reasoning as in the proof of 
Lemma \ref{du-lem-S_0}. 
This means that $\{g_i(y)\}$ is uniformly bounded away from $a_\ast(\{g_i\})$
since $g_i(l)$ passes $F$ with $\{g_i(t)\}$ converging to 
$a_\ast(\{g_i\})$. 
Since $\{g_i(y)\}$ is bounded away from $a_\ast(\{g_i\})$ uniformly, 
as in the proof of Lemma \ref{du-lem-S_0}
using \eqref{du-eqn-Lgi} and similarly to proving 
the conclusion of the lemma, 
we obtain 
 \begin{equation}\label{pr-eqn-lhg}
 \left\{ 
  \frac{1}{\lambda_{\mbv_{\tilde E}}(g_i)^{1+\frac{1}{n}}}
 \hat h(g_i)(\vec v_{-}) \right\}\ra 0 
 \end{equation}
 in the Euclidean metric.
 To explain more, we write $\vec v_-$ as a sum of $\vec v_i^p + \vec v_i^S$ as 
in the referenced proof.
The rest is analogous.
 
 
 Let $\mbv_{\tilde E}$ denote the unit vector in the direction of 
 $\vec v_{\tilde E}$. We consider $\bR^{n}$ to be a complementary subspace
 to this vector under the norm $\llrrV{\cdot}$. 
 We write the vector for $z$ as 
 $\vec{v}_z = \lambda \vec v_- +  \vec v_{\tilde E}$ and similarly the vector for $z'$ as 
 $\vec{v}_{z'} = \lambda' \vec v_- +  \vec v_{\tilde E}$. 
 Then 
 \[g_i(\vec{v}_z) = 
  \lambda\frac{1}{\lambda_{\mbv_{\tilde E}}^{\frac{1}{n}}(g_i)} \hat h(g_i)(\vec v_-) + (\lambda \vec{b}_{g_i}\cdot \vec v_- + 
 \lambda_{\mbv_{\tilde E}}(g_i)) \vec v_{\tilde E}\] 
 by \eqref{pr-eqn-mForm}.
 Let us denote 
 \[c_i := \llrrV{\frac{1}{\lambda_{\mbv_{\tilde E}}^{\frac{1}{n}}(g_i)} \hat h(g_i)(\vec v_-)}_E.\]
 Since the direction of $g_i(\vec{v}_z)$ is bounded away from 
 $\vec v_{\tilde E}$, 
\[\left| \lambda \frac{\vec{b}_{g_i}\cdot \vec v_-}{c_i} 
+ \frac{\lambda_{\mbv_{\tilde E}}(g_i)}{c_i} \right| \]
is uniformly bounded. 
By \eqref{pr-eqn-lhg}, we obtain
\[\left\{ \left|\frac{\lambda_{\mbv_{\tilde E}}(g_i)}{c_i}\right| \right\} \ra \infty.  \]
 Hence, 
 \[\left\{\left| \frac{\vec{b}_{g_i}\cdot \vec v_-}{c_i}    \right| \right\}\ra \infty 
 \hbox{ as }  i \ra \infty.\]  
 
 We also have 
 \[g_i(\vec{v}_{z'}) = 
 \lambda'\frac{1}{\lambda_{\mbv_{\tilde E}}^{\frac{1}{n}}(g_i)} \hat h(g_i)(\vec v_-) + (\lambda' \vec{b}_{g_i}\cdot \vec v_- + 
 \lambda_{\mbv_{\tilde E}}(g_i)) \vec v_{\tilde E}.\] 
 Since $\lambda' \ne \lambda$, 
 and $\left\{\left| \frac{\vec{b}_{g_i}\cdot \vec v_-}{c_i}    \right|\right\} \ra \infty$, 
it follows that 
 \[ \left\{\left|\lambda' \frac{\vec{b}_{g_i}\cdot \vec v_-}{c_i} + 
 \frac{\lambda_{\mbv_{\tilde E}}(g_i) \vec v_{\tilde E}}{c_i} \right| 
 = \left|(\lambda' - \lambda) \frac{\vec{b}_{g_i}\cdot \vec v_-}{c_i} + 
 \lambda \frac{\vec{b}_{g_i}\cdot \vec v_-}{c_i} + 
 \frac{\lambda_{\mbv_{\tilde E}}(g_i) \vec v_{\tilde E}}{c_i} \right| \right\}
 \]
cannot be uniformly bounded.
This implies that $g_i(z')$ converges to 
$\mbv_{\tilde E}$ or $\mbv_{\tilde E-}$.
Since $z' \in K'$ and $K'$ is 
$\bGamma_{\tilde E}$-invariant, this is a contradiction. 

By Lemma \ref{pr-lem-attracting}, $L_{\tilde E}$ meets every
great segment in $\mathcal{T}$.  
Thus, $K' \cap l_y  = L_{\tilde E}\cap l_y$ for 
every $y$ in $\partial \clo(\tilde \Sigma_{\tilde E})$. 
Thus, $K'= L_{\tilde E}$, and 
 (i) and (ii) hold for $L_{\tilde E}$.

 (iii) Since $L_{\tilde E}$ is closed, compact, and bounded away 
 from $\mbv_{\tilde E}, \mbv_{\tilde E-}$, 
 the section $s: \partial \clo(\Omega_{\tilde E}) \ra \Bd \mathcal{T}_{\mbv_{\tilde E}}(\tilde \Sigma_{\tilde E})$ 
 is continuous. Otherwise, we can contradict (ii) by taking two sequences converging 
 to distinct points on a great segment from $\bv_{\tilde E}$ to its antipode.




(B) Now suppose that $\bGamma_{\tilde E}$ is neither virtually factorizable 
nor hyperbolic. 
Lemma \ref{pr-lem-attracting2} proves that the orbits have limits points in  
$L_{\tilde E}$ only. 
Every attracting attracting 
fixed points in $\mathcal{T}_{\mbv_{\tilde E}}(\tilde \Sigma_{\tilde E})$ is in 
$L_{\tilde E}$ as in case (A). 

First suppose that a great segment $\eta$ in $\mathcal{T}_{\mbv_{\tilde E}}(\tilde \Sigma_{\tilde E})$ 
with endpoints $\mbv_{\tilde E}$ and $\mbv_{\tilde E -}$ 
corresponds to an element $y$ of $\partial \clo(\tilde \Sigma_{\tilde E})$. 
Now we take a line $l$ in $\tilde \Sigma_{\tilde E}$ as in the hyperbolic case. Then 
\eqref{pr-eqn-lhg} holds as above using Lemma \ref{du-lem-S_0II} 
instead of Lemma \ref{du-lem-S_0}. 
The identical argument 
will show that $\eta^o$ meets with $L_{\tilde E}$ at a unique point. 
This proves (i) and (ii). 
(iii) follows as above.

(C) Suppose that $\bGamma_E$ is virtually factorizable. 
We follow the proof of Theorem \ref{pr-thm-distanced}.
Now, the space of open great segments with an endpoint $\mbv_{\tilde E}$ in
 ${\mathcal T}_{\mbv_{\tilde E}}(\tilde \Sigma_{\tilde E})^o$ corresponds to 
a properly convex domain $\Omega$ that is the interior of 
the strict join $K_1\ast \cdots \ast K_l$. 
Then a totally geodesic $\bGamma_{\tilde E}$-invariant 
hyperspace $H$ is disjoint from $\{\mbv_{\tilde E}, \mbv_{\tilde E-}\}$ 
as shown in the proof of Theorem \ref{pr-thm-distanced}. 
Here, we may regard $K_i \subset H$ for each $i=1, \dots, l$. 
Then consider any sequence $\{g_i\}$ such that 
$\{g_i(x)\} \ra x_0$ for a point 
$x \in  {\mathcal T}_{\mbv_{\tilde E}}(\tilde \Sigma_{\tilde E})^o$ and $x_{0} \in {\mathcal T}_{\mbv_{\tilde E}}(\tilde \Sigma_{\tilde E})$.
Let $x'$ denote the corresponding point of $\tilde \Sigma_{\tilde E}$ for $x$. 
Then $\{g_i(x')\}$ converges to a point $y \in \SI^{n-1}_{\mbv_{\tilde E}}$. 
Let $\vec x\in \bR^{n+1}$ be the vector in the direction of $x$. 
We write \[\vec x = \vec x_E + \vec x_H\] where $\vec x_H$ is in the direction of $H$ and $\vec x_E$ is in the direction of $\mbv_{\tilde E}$. 
By the uniform middle-eigenvalue condition and estimating the size of 
vectors, we obtain the same situation as in \eqref{pr-eqn-gired} and hence
$\{g_i(x)\} \ra x_{0}$ for $x \in L_{\tilde E}$ and $x_0\in H$. 
Hence, $x_{0} \in H \cap L_{\tilde E}$. 
Thus, every limit point of an orbit of $x$ is in $H$. 

If there is a point $y$ in $(L_{\tilde E} -H) \cap 
\partial{\mathcal T}_{\mbv_{\tilde E}}(\tilde \Sigma_{\tilde E})^o - 
\{\mbv_{\tilde E}, \mbv_{\tilde E-}\}$, 
then there is a strict join $K_{i_1}\ast \cdots \ast K_{i_l} \ast \{\mbv_{\tilde E}\}$ 
over a proper collection containing $y$. 
As in the proof of Theorem \ref{pr-thm-distanced}, 
by Proposition \ref{prelim-prop-Ben2}, 
we can find a sequence $g_i$, virtually 
central $g_i\in \bGamma_{\tilde E}$, such that 
$\{g_i| K_{i_1}\ast \cdots \ast K_{i_l}\} $ converges to the identity, and 
the maximal norm of $g_i$ is supported in the complementary domains. 
The norms of eigenvalues associated with 
$K_{i_1}\ast \cdots \ast K_{i_l}$ have uniformly  bounded ratios with 
the minimal one $\lambda_{n}(g_i)$ in this case. 
The maximal norm $\bar \lambda_1(g_i)$ of the eigenvalue associated with 
$K_{i_1}\ast \cdots \ast K_{i_l}$ and $\lambda_{\tilde E}(g_i)$
satisfy $\{\bar \lambda_1(g_i)/\lambda_{\tilde E}(g_i)\} \ra 0$
by the uniform middle-eigenvalue condition. 
Hence $\{g_i(y)\} \ra \mbv_{\tilde E}$. Again this is a contradiction. 
Hence, we obtain $L_{\tilde E} \subset H$. 
(i), (ii), and (iii) follow easily now. 

(iv)  Suppose that we have another distanced set $L'$. We take 
a convex hull of $L_E\cup L'$ and apply the same reasoning as above. 
\end{proof}

\subsection{An extension of Koszul's openness} \label{app-sec-Koszul}

Here, we state and prove a well-known 
minor modification of Koszul's openness result. 

For a real projective orbifolds, two immersed oriented domains in $(n-1)$-dimnsional orbifolds
{\em meet in a dihedral angle} $< \pi$ at their common boundary point if 
their lifts in the universal cover develops to two oriented domains in hyperspaces meeting each 
other in a dihedral angle $< \pi$. 

A {\em radial affine connection} is an affine connection on $\bR^{n+1} -\{O\}$ invariant under 
a scalar dilatation $S_t: \vec{v} \ra t\vec{v}$ for every $t > 0$. 

\begin{proposition}[Koszul] \label{app-prop-koszul} 
	Let $M$ be a compact properly convex real projective $n$-orbifold with strictly convex
smooth (resp. polyhedral) boundary. 
	Let 
	\[h:\pi_{1}(M) \ra \PGL(n+1, \bR)\: (\hbox{resp. } h: \pi_1(M) \ra \SLnp)\] denote the holonomy homomorphism
	acting on a properly convex domain $\Omega_{h}$ in $\RP^{n}$ {\rm (}resp.  in $\SI^{n}${\rm )}.  
	Assume that $M$ is projectively diffeomorphic to $\Omega_{h}/h(\pi_{1}(M))$. 
	Then there exists a neighborhood $U$ of $h$ in 
	\[\Hom(\pi_1(M), \PGL(n+1, \bR))\:
	(\hbox{resp. } \Hom(\pi_1(M), \SLnp) )\] 
	such that every $h' \in U$ acts on a properly convex domain $\Omega_{h'}$
	and $\Omega_{h'}/h'(\pi_{1}(M))$ is 
	a compact properly convex real projective $n$-orbifold 
	with strictly smooth (resp. polyhedral) convex boundary. 
	Also,  $\Omega_{h'}/h'(\pi_{1}(M))$ is diffeomorphic to $M$. 
\end{proposition}
\begin{proof} We prove for $\SI^{n}$. 
First, we assume that $\partial M$ is smooth. 
	Let $\Omega_{h}$ be a properly convex domain projectively covering $M$. 
	We may modify $M$ by pushing $\partial M$ inward:
	Let $\Omega'_{h}$ be the inverse image of $M'$ in $M$. 
	Then $M'$ and  $\Omega'_{h}$ are properly convex by Lemma \ref{prelim-lem-pushing}. 
	
	The linear cone $C(\Omega^{o}_{h})\subset \bR^{n+1} = \Pi^{-1}(\Omega^{o}_{h})$ over $\Omega^{o}_{h}$ 
	has a smooth, strictly convex Hessian function $V$ 
	by Vey \cite{Vey} or Vinberg \cite{Vinberg67}. 
	Let $C(\Omega'_{h})$ denote the linear cone over $\Omega'_{h}$.
	We extend the group $\mu(\pi_1(M))$ by adding
	a transformation $\gamma: \vec{v} \mapsto 2\vec{v}$ to $C(\Omega^{o}_{h})$. 
	For the fundamental domain $F'$ of $C(\Omega'_{h})$ under this group, 
	the Hessian matrix of $V$ restricted to $F \cap C(\Omega'_{h})$ 
has a constant positive definite matrix as a lower bound. 
	Also, the boundary $\Bd C(\Omega'_{h})$ is strictly convex 
	in its intersection to any transverse affine subspace to the radial directions at any point.
	
	Let $N'$ be a compact orbifold 
$C(\Omega'_{h})/\langle \mu(\pi_1(\tilde E)), \gamma \rangle$ equipped with 
a torsion-free flat connection $\nabla_h$. 
	Note that $S_t$, $t \in \bR_+$, becomes an action of a circle on $M$.
	The change of representation from $h$ to $h': \pi_1(M) \ra \SLnp$ 
	is realized by a change of a holonomy representation of $M$ and hence by
	a change of an affine connection on $C(\Omega'_{h})$. Since $S_t$ commutes with the images of $h$ and $h'$, 
	$S_t$ still gives us a circle action on $N'$ with a different affine connection. 
This can be proved by covering $N$ into a circle times $(n-1)$-balls 
and using the patching as in Lok \cite{Lok} (see \cite{dgorb}
for the orbifold version).
This can be accomplished for $h'\in U$ for 
a sufficiently small neighborhood $U$ of $h$ in  $\Hom(\pi_1(M), \SLnp)$.
	Without loss of generality, we may assume 
	that the circle action is fixed and $N'$ equipped with new holonomy 
$h'\in U$ is invariant under this action.
	
	We may assume that $N'$ and $\partial N'$ are foliated by circles that are flow curves of the circle action. 
	The change corresponds to a sufficiently small $C^{r}$-change in the affine connection for $r \geq 2$ as shown in \cite{dgorb}. 
We find a transversal hyperspace $O_h:= P\cap C(\Omega'_h)$ in a convex domain 
within an affine patch in $P$.
	Now, the strict positivity of 
	the Hessian of $V$ on $O_h$ and the boundary convexity
in the transversal hyperspaces are preserved
where we are changing $V$ and the connection simultaneously. 
	Let $C(\Omega''_{h'})$ denote the universal cover of $N'$ equipped 
with the new affine connection
$\nabla_{h'}, h'\in U$ that is invariant under the circle action.
Hence, $C$ under a new affine connection has locally convex boundary. 
This implies that $C$ is still convex with a function with a positive definite Hessian. 
	Thus, $C(\Omega''_{h'})$ is also a properly convex affine cone by \cite{Koszul68}. 
	Also, it is a cone over a properly convex domain $\Omega''_{h'}$ in $\SI^{n}$.
The strict convexity of $\partial O_h$ implies that $\partial \Omega''_{h'}$ is strictly convex. 

Now, assume that $\partial M$ is convex polyhedral. 
First, we perturb the vertices of $\partial M$ generically so that $\partial M$ is a union of 
$(n-1)$-simplices meeting one another in dihedral angles $< \pi$:
We consider $M$ to be in an open real projective manifold $N$. 
We triangulate $M$ by simplices as in Chapter 4 of \cite{Cbook}. 
We may assum that each simplex is in the interior of a compact convex subset of $N$. 
Let $N$ have some Riemannian distance metric $\mu_N$. 
We cover $M$ by a finite set $\mathcal{U}_N$ of the interiors of compact convex balls of $N$
so that each simplex is inside one where the minimal $\mu_N$ distance from the simplex to
the boundary of a compact convex set is at least $\eps$ for some $\eps >0$. 
We perturb vertices of each boundary polyhedron within these open sets. 
Then we take a convex hull of vertices in each former polyhedron.  
The dihedral angles change very small amount
in the universal cover, and we need to look at only finitely many of these. 
Taking sufficiently small pertubations of the vertices, 
we obtain a triangulated submanifold $M'$ of $N$ with convex polyhedral boundary. 
By choosing the vertices inside $M$, we obtain that $M'$ is properly convex
by Lemma \ref{prelim-lem-locconv}. 

Now we change the holonomy: We interpret this as a chang e of projective connection of $N$. 
The compact convex balls in $\mathcal{U}_N$ also can be changed to ones whose interior still covers $M$ following the arguments in \cite{dgorb}. (See \cite{Lok} also.)
We still have the property that minimal $\mu_N$-distance from each simplex to the boundary 
of compact convex balls are $> \eps/2$
if the connection by a sufficiently $C^r$-close, $r \geq 2$, to the original one.
As above, we obtain a triangulated submanifold $M''$ of $N$ with the new connection
with convex polyhedral boundary. 

The universal covers 
of special affine suspection of $M''$ and $N$ have 
a function with positive definite Hessian since $N$ has one. 
Again \cite{Koszul68} applies. 
\hfill	\SnT {\parfillskip0pt\par}
\end{proof} 

We denote by $\PGL(n+1, \bR)_{v}$ the subgroup of $\PGL(n+1, \bR)$ fixing a point $v \in \RP^n$, 
and denote by $\SLnp_{v}$ the subgroup of $\SLnp$ fixing a point $v \in \SI^n$. 
Let	$\Hom_C(\Gamma, \PGL(n+1, \bR)_{v})$ denote the space of 
representations dividing on a properly convex domain in 
$\SI^{n-1}_{v}$. We define $\Hom_C(\pi_1(M), \SLnp_{v})$
similarly. 

\begin{proposition}\label{app-prop-lensP}
	Let $T$ be a tube domain over a properly convex domain 
$\Omega \subset \RP^{n-1}$ {\em (}resp. $\subset \SI^{n-1}${\em ).} 
	Let $B$ be a strictly convex smooth (resp. polyhedral) 
hypersurface bounding a properly convex domain in 
	a tube domain $T$ where a projective group $\Gamma \subset \PGL(n+1, \bR)$ 
 {\em (}$\subset \SLnp${\em )} acts on.  
Let $v$ be a vertex of $T$ fixed by $\Gamma$. 
	Assume that $B$ meets each radial ray in $T$ from $v$ transversely. 
	Assume that the induced action of
$\Gamma$ on $\Omega$ is properly discontinuous and cocompact.
	Then there exists a neighborhood of the inclusion map in 
	\[\Hom_C(\Gamma, \PGL(n+1, \bR)_{v})\:
	(\hbox{resp. }\Hom_C(\pi_1(M), \SLnp_{v}))\]
	where every element $h$ acts on 
	a strictly convex hypersurface $B_{h}$ in a tube domain ${\mathcal T}_{h}$ meeting 
	each radial ray at a unique point, and bounding a properly convex domain in 
${\mathcal T}_{h}$. 
\end{proposition}
\begin{proof}
Suppose that $B$ is smooth. 
	We first assume that $B, T \subset \SI^n$. 
	For sufficiently small neighborhood $V$ of $h$ in $\Hom_C(\Gamma, \SLnp_{v})$,
	$h(\Gamma)$, $h'\in V$ acts on a properly convex domain $\Omega_{h'}$ properly discontinuously and cocompactly by Theorem 1 of \cite{dgorb}
	(see Koszul \cite{Koszul68}). 
	A large compact subset $K$ of $\Omega$ flows to a compact 
	subset $K_h$ by a diffeomorphism using a method of Section 5 of  
	\cite{dgorb}. 
	Let ${\mathcal T}_{h}$ denote the tube over $\Omega_{h}$. 
	Since $B/\Gamma$ is a compact orbifold, 
	we choose $V' \subset V$ such that for the projective connections on a compact neighborhood of $B/\Gamma$ 
	corresponding to elements of $V'$, $B/\Gamma$ is still strictly convex
	and transverse to radial lines. 
	For each $h\in V'$, we obtain an immersion $\iota_{h}: B \ra {\mathcal T}_{h}$ to a strictly convex domain, 
	transverse to radial lines, since we can view the change of holonomy as
	small $C^1$-change of connections. 
	(Or we can use the method described in Section 5 of \cite{dgorb}.)
	Let $p_{{\mathcal T}_h}: {\mathcal T}_{h} \ra \Omega_{h}$ denote the projection with fibers equal to the radial lines. 
	Also, in this way of viewing as the connection change, 
	$p_{{\mathcal T}_h}\circ \iota_{h}$ is a proper immersion to $\Omega_{h}$, it is a diffeomorphism $B\ra \Omega_h$. 
(Here again, we can use Section 5 of \cite{dgorb}.)
	Each point of $B$ is transverse to a radial segment from $v$.
	By considering the compact fundamental domains of $B$,  
	we see that the same holds for $B_{h}$ for $h$ sufficiently near $\Idd$.  
	Also, $B_{h}$ is strictly convex and smooth. 
By Proposition \ref{app-prop-koszul}, the conclusion follows.  

If $B$ is convex polyhedral hypersurface, entirely similar argument will apply. 
\hfill	\SnT {\parfillskip0pt\par}
\end{proof}

\subsection{Convex cocompact actions of the p-end holonomy groups.} \label{pr-subsec-redlens}


In $\SI^n$,	recall from Definition \ref{intro-defn-R-ends} that 
a (resp. generalized) lens-shaped R-p-end with the p-end vertex $\mbv_{\tilde E}$ 
identified with a point $\SI^n$ is ({\em resp. generalized}) {\rm lens-shaped} 
	if we can choose a (resp. generalized) lens domain $D$ 
	in $\SI^n$ such that 
\begin{itemize} 
\item the interior of a lens cone $D \ast \mbv_{\tilde E}$ 
	is a developing image of a p-end neighborhood 
		
\item where the top hypersurfaces $A$ and the bottom one $B$ of $L$
satisfy 
\[\clo(A) - A = \clo(B) - B, A\cup B = \partial L \hbox{ and }
\clo(A) \cup \clo(B) = \Bd L.\]  
\end{itemize} 
We deduce that 

\begin{lemma} \label{pr-lem-strictlens} 
Let $\tilde E$ be a lens-shaped R-end with $D \ast \mbv_{\tilde E}$ as a developing image 
of a (generalized) 
lens cone p-end neighborhood. Then each great open segment in $\SI^n$ from 
$p$ in the direction of 
		$\partial \clo(\tilde \Sigma_{\tilde E})$ meets 
		$\clo(D) - A - B$ at a unique point. 
Conversely, this property implies the strictness of $D$. 
\end{lemma} 
\begin{proof}
$\dev$ restricted to $A$ and $B$ is a homeomorphism to its image.
Hence, each boundary of $A$ and $B$ maps to the boundary of this domain.
 Definition \ref{intro-defn-R-ends} requires $L$ to be strict lens or strict generalized lens, 
which implies the result. 
The converse is also straightforward. 
\end{proof} 


	Recall that in order that $L$ is to be a lens, we assume that $\pi_1(\tilde E)$ acts cocompactly on $L$. 
	Moreover,  by 
	Corollary \ref{pr-cor-LambdaW}, $\clo(A)-A$ must equal the limit set $\Lambda_{\tilde E}$ of 
	$\tilde E$ 
		\index{lens!strict}
			\index{lens!generalized!strict}
			
			Also, the images of these under $p_{\SI^n}$ are called by 
			the same names respectively. 

In this section, we will prove Proposition \ref{pr-prop-convhull2} obtaining a lens. 




We will use the Hilbert metric $d_{V}$ for any properly convex $V$. (See Section \ref{rh-sub-Hilbert}.)
Let $U$ be a properly convex set of 
${\mathcal T}_{\mbv_{\tilde E}}(\tilde \Sigma_{\tilde E})$ containing 
a lens cone with the vertex at $\mbv_{\tilde E}$. 
Let $\Sigma_1$ and $\Sigma_2$ be two hypersurfaces in $U$ 
meeting every radial line at exactly one point
in ${\mathcal T}_{\mbv_{\tilde E}}(\tilde \Sigma_{\tilde E})$.
Then $\Sigma_1$ and $\Sigma_2$ are said to be
{\em $\eps$-approximate with respect to $d_U$} (resp. $\bdd$), for $\eps > 0$,  
if for each maximal radial segment $l$ in $U$, the shortest $d_U$-path 
(resp. $\bdd$-path) length 
between $l\cap \Sigma_1$ and $l\cap \Sigma_2$ is less than $\eps$. 
\index{epsilon-approximate@$\eps$-approximate|textbf}

In the following, $\orb$ needs not be properly convex but merely convex. 
\begin{proposition}\label{pr-prop-convhull2} 
Let $\orb$ be a convex
real projective $n$-orbifold where 
$\torb \subset \SI^n$ {\em (}resp. $\RP^n${\em )}.
\begin{itemize} 
\item Let $\bGamma_{\tilde E}$ be the holonomy group of a properly convex R-p-end $\tilde E$.
\item Let ${\mathcal T}_{\mbv_{\tilde E}}(\tilde \Sigma_{\tilde E})$ be an open tube corresponding to $R(\mbv_{\tilde E})$. {\rm (}See page \ref{page:Rx} for the definition.{\rm )}
\item Suppose that $\bGamma_{\tilde E}$ satisfies the uniform middle-eigenvalue condition with respect to the R-p-end structure, 
and acts on a distanced compact convex set 
$K$ in ${\mathcal T}_{\mbv_{\tilde E}}(\tilde \Sigma_{\tilde E})$  
where $K \cap {\mathcal T}_{\mbv_{\tilde E}}(\tilde \Sigma_{\tilde E}) \subset \torb$. 
\end{itemize} 
Then the following  statements hold\/{\rm :}
\begin{itemize}  
\item Any connected open  p-end neighborhood $U$ containing a lift to $\torb$ of 
$K \cap {\mathcal T}_{\mbv_{\tilde E}}(\tilde \Sigma_{\tilde E})$ also contains 
a lens $L'$ and a lens cone p-end neighborhood 
$L' \ast \{\mbv_{\tilde E}\} - \{\mbv_{\tilde E}\}$ of 
the R-p-end $\tilde E$. 
{\rm (}\/Note that $U$ need not be in the lift of $\torb$.{\rm )}
\item We can choose the lens $L'$ in $U$  such that $\Bd L' \cap \mathcal{T}^o = A \cup B$
for strictly convex smooth connected hypersurfaces $A$ and $B$. 
\item  We can choose $A$ such that $A$ radially $\eps$-approximate any outer boundary component of a generalized lens cone with respect to $d_U$ provided $U$ is properly convex.
\item The lens $L'$ of the cone is a strict lens
and $L' \ast \mbv_{\tilde E}$ is a strict lens-cone. 
\end{itemize}  
\end{proposition}
\begin{proof}\renewcommand{\qedsymbol}{}
	First suppose that $\torb \subset \SI^n$. 
We may assume that $U$ embeds into a neighborhood of $L$ under 
a developing map
by taking $U$ sufficiently small. 
We continue to denote the image by $U$. 
The projection of $K$ to $\tilde \Sigma_{\tilde E}$ must be onto since
it must be a $\Gamma_{\tilde E}$-invariant convex set. 
So, $K$ must meet each great segment with endpoints $\mbv_{\tilde E}$ in the directions of 
$\tilde \Sigma_{\tilde E}$.  
Hence, $K \cap {\mathcal T}_{\mbv_{\tilde E}}(\tilde \Sigma_{\tilde E})$ 
is a separating set in $\torb$, and 
 $U - K$ has two components since the boundary of $K$ has two components in $\torb$. 

%
Let $\Lambda_{\tilde E}$ denote $\Bd {\mathcal T}_{\mbv_{\tilde E}}(\tilde \Sigma_{\tilde E}) \cap K$. 
Let us choose finitely many points $z_1, \dots, z_m \in U - K$ in the two components of $U - K$.

Proposition \ref{pr-prop-orbit} shows that the orbits of $z_i$ for each $i$ accumulate to points of $\Lambda_{\tilde E}$ only. 
Hence, a totally geodesic hypersphere separates $\mbv_{\tilde E}$ with these orbit points
and another one separates $\mbv_{\tilde E-}$ and the orbit points. 
Define the convex hull $C_1:= \CH(\bGamma_{\tilde E}(\{z_1, \dots, z_m\})\cup K)$. 
Thus, $C_1$ is a compact convex set disjoint from 
$\mbv_{\tilde E}$ and $\mbv_{\tilde E-}$ 
and satisfies 
$C_1 \cap \Bd {\mathcal T}_{\mbv_{\tilde E}}(\tilde \Sigma_{\tilde E}) = \Lambda_{\tilde E}$. 
(See Definition \ref{prelim-defn-convhull}.) 

We need the following lemma:

\hfill \cpr {\parfillskip0pt\par}
\end{proof} 

\begin{lemma}\label{pr-lem-push} 
We continue to assume as in Proposition \ref{pr-prop-convhull2}
where $\torb \subset \SI^n$. 
Then we can choose 
$z_1, \dots, z_m$ in $U$ such that 
for $C_1:= \CH(\bGamma_{\tilde E}(\{z_1, \dots, z_m\})\cup K)$, 
$\Bd C_1 \cap \torb$ is disjoint from $K$, 
and $C_1 \subset U$.
\end{lemma}
\begin{proof} 
First, suppose that $K$ is not in a hyperspace. Then 
$(\Bd K \cap {\mathcal T}_{\mbv_{\tilde E}}(\tilde \Sigma_{\tilde E}))/\bGamma_{\tilde E}$ is diffeomorphic to 
a disjoint union of two copies of $\Sigma_{\tilde E}$. 
We can cover a compact fundamental domain of 
$\Bd K \cap {\mathcal T}_{\mbv_{\tilde E}}(\tilde \Sigma_{\tilde E})$ by the interior of $n$-balls in $\torb$ that are convex hulls of finite sets
of points in $U$. Since the $\bGamma_{\tilde E}$-action on $L$ is cocompact for a lens $L$ 
containing $K\cap \torb$, it follows that 
so is $(K\cap \torb)/\bGamma_{\tilde E}$,
and 
there exists a positive lower bound of $\{d_{\torb}(x, \Bd U\cap \torb)| x \in K\}$.
Let $F$ denote the union of these finite sets. 
We can choose $\eps > 0$ such that 
the $\eps$-$d_{\torb}$-neighborhood $U'$ of $K$ in $\torb$ is a subset of $U$. 
Moreover, $U'$ is convex by Lemma \ref{prelim-lem-nhbd} following \cite{CLT15}. 
Let $z_1, \dots, z_m$ denote the points of $F$. 
If we choose $F$ to be in $U'$, then $C_1$ is in $U'$ since $U'$ is convex. 

The disjointness of $\Bd C_1$ from 
$K \cap  {\mathcal T}_{\mbv_{\tilde E}}(\tilde \Sigma_{\tilde E})$ follows since 
the $\bGamma_{\tilde E}$-orbits of the above balls cover $\Bd K \cap {\mathcal T}_{\mbv_{\tilde E}}(\tilde \Sigma_{\tilde E})$. 

If $K$ lies in a hyperspace, the reasoning is similar to the above. 
\end{proof}



We continue:
\begin{lemma} \label{pr-lem-infiniteline} 
Let $C$ be a compact $\bGamma_{\tilde E}$-invariant 
distanced convex set with boundary in 
where the $\bGamma_{\tilde E}$-action on $(C \cap {\mathcal{T}^o_{\tilde E}})$ is cocompact. 
There are two connected hypersurfaces $A$ and $B$ of
$\Bd C \cap {\mathcal T}^o_{\tilde E}$ meeting every great segment in ${\mathcal T}^o_{\tilde E}$. 
Suppose that $A$ and $B$ are disjoint from another $C'$ $\bGamma_{\tilde E}$-invariant 
distanced compact convex set with boundary in 
where $\bGamma_{\tilde E}$-action on $(C' \cap {\mathcal{T}^o_{\tilde E}})$ is cocompact. 
Then neither $A$ nor $B$ contains no line ending in $\Bd \torb$. 
\end{lemma} 
 \begin{proof} 
 It suffices to prove for $A$. 
Suppose that there exists 
a line $l$ in  $A$ ending at a point of 
$\Bd {\mathcal T}_{\mbv_{\tilde E}}(\tilde \Sigma_{\tilde E})$.
 Assume $l \subset A$. 
The line $l$ projects to a line $l'$ in ${\tilde E}$. 

Let $C_1 = C \cap {\mathcal T}_{\mbv_{\tilde E}}(\tilde \Sigma_{\tilde E})$. 
Since $A/\bGamma_{\tilde E}$ and $B/\bGamma_{\tilde E}$ are both compact, 
there exists a fibration $C_1/\bGamma_{\tilde E} \ra A/\bGamma_{\tilde E}$ 
induced from $C_1 \ra A$ using the foliation by great segments with endpoints $\mbv_{\tilde E}, \mbv_{\tilde E-}$. 

Since the $\bGamma_{\tilde E}$-action on $A$ is cocompact,  
we choose a compact fundamental domain $F$ in $A$, and 
choose a sequence $\{x_i \in l\}$ 
whose image sequence in $l'$ converges to the endpoint of $l'$ in 
$\partial \clo(\tilde \Sigma_{\tilde E})$. 
We choose $\gamma_i \in \bGamma_{\mbv_{\tilde E}}$ such that $\gamma_i(x_i) \in F$ 
where $\{\gamma_i(\clo(l'))\}$ geometrically converges to a segment $l'_\infty$ with both endpoints in $\partial \clo(\tilde \Sigma_{\tilde E})$. 
Hence, $\{\gamma_i(\clo(l))\}$ geometrically converges to a segment $l_\infty$ in $A$.
We may assume that the endpoint $z$ of $l$ in $A$ satisfies $\{\gamma_i(z)\} \ra p_1$. 
Proposition \ref{pr-prop-orbit} implies that the endpoint $p_1$ of $l_\infty$ is in $L_{\tilde E}:= L \cap \Bd {\mathcal T}_{\tilde E}$. 
Let $t$ be the endpoint of $l$ distinct from $z$. Then $t \in A$. 
Since $\gamma_{i}$ is not a bounded sequence, $\gamma_i(t)$ converges to a point of $\Lambda_{\tilde E}$. 
Thus, both endpoints of $l_\infty$ are in $\Lambda_{\tilde E}$
and hence $l_\infty^{o} \subset C'$ by the convexity of $C'$.
However, $l \subset A$ implies that $l_\infty^{o} \subset A$. As $A$ is disjoint from $C'$, 
this is a contradiction.  The similar arguments apply to $B$. 
\end{proof} 

\begin{proof}[Proof of Proposition \ref{pr-prop-convhull2} continued]  
We denote by $C_1$ the compact convex subset 
$C_1 = C \cap {\mathcal T}_{\mbv_{\tilde E}}(\tilde \Sigma_{\tilde E})$ for $C$ obtained by Lemma
\ref{pr-lem-push}. 
Since $C_1$ meets in a compact segment any great segment in 
$ {\mathcal T}_{\mbv_{\tilde E}}(\tilde \Sigma_{\tilde E})$, it follows that 
$\Bd C_1 \cap {\mathcal T}_{\mbv_{\tilde E}}(\tilde \Sigma_{\tilde E})^o$ is the union of 
two hypersurfaces $A$ and $B$. 
Since  $C$ is defined as the convex hull of 
$(\bGamma_{\tilde E}(\{z_1, \dots, z_m\})\cup K)$, and balls that are convex hulls of
some points of $z_1, \dots, z_m$ and their images cover 
$\Bd K \cap {\mathcal T}_{\mbv_{\tilde E}}(\tilde \Sigma_{\tilde E})^o$, 
it follows that 
the extreme points of $A$ or $B$ must be vertices of the images of 
$z_1, \dots, z_m$ and points of $K \cap \Bd {\mathcal T}_{\mbv_{\tilde E}}
(\tilde \Sigma_{\tilde E})$. 
Since the ball cover $\Bd K \cap  {\mathcal T}_{\mbv_{\tilde E}}(\tilde \Sigma_{\tilde E})^o$, 
$A$ and $B$ are disjoint from $K$. 

Since $A$ and analogously $B$ do not contain any geodesic ending at $\Bd \torb$, 
by Lemma \ref{pr-lem-infiniteline},  
$A \cup B = \Bd C_1 \cap {\mathcal T}_{\mbv_{\tilde E}}(\tilde \Sigma_{\tilde E})^o$ is a union of compact $(n-1)$-dimensional polytopes which meet one another in dihedral angles $<pi$

Following Proposition \ref{pr-prop-lenssmooth} completes the proof of Proposition \ref{pr-prop-convhull2}. The approximation property of $A$ 
is the consequence of being inside $L_1^o$. and outside $L_2$ 
and the cocompact action of $\bGamma_{\tilde E}$ on $L_1 - L_2^o$.
\end{proof} 

\begin{proposition}\label{pr-prop-lenssmooth}
Assume the premise of Proposition \ref{pr-prop-convhull2}.
Suppose that a lens cone $L_1\ast \{\bv_{\tilde E}\} - \{\mbv_{\tilde E}\}$ 
is in a convex p-end neighborhood $U$ of a p-end $\tilde E$
a p-end neighborhood of $\tilde E$. 
Suppose that $L_1$ contains a lens $L$ in its interior
where $L\ast \{\bv_{\tilde E}\} - \{\mbv_{\tilde E}\}$ is again 
 a R-p-end neighborhood of $\tilde E$. 
Suppose that $L_1$ is bounded by two connected convex polyhedral hypersurfaces. 
Then there exists a lens $L_2$ bounded by two connected strictly convex smooth
hypersurfaces such that $L_2 \subset U$ and $L \subset L_2^o, L_2 \subset L_1^o$.	
\end{proposition}
\begin{proof} 
	First, assume $\torb \subset \SI^n$. 
Let us take the dual domain $U_L$ of $(L \ast \{\mbv_{\tilde E}\})^o$.
The dual $U_1$ of $(L_1 \ast \{\mbv_{\tilde E}\})^o$ is an open subset 
of $U_L$ by \eqref{prelim-eqn-dualinc}. 
By Proposition \ref{pr-prop-dualDA},  the dual action is asymptotically nice. 
By the uniqueness parts of Theorems \ref{du-thm-asymnice} and \ref{du-thm-asymniceII},
$U_L$ and $U_1$ are asymptotically-nice properly convex domains. 
By Lemma \ref{pr-lem-predual} (iii), the hyperspaces sharply supporting 
$(L_1 \ast \{\mbv_{\tilde E}\})^o$ at $\mbv_{\tilde E}$ 
correspond to points of a totally geodesic domain $D$, 
\[ D = \clo(U_1)\cap P = \clo(U_L)\cap P\]
for a $\bGamma_{\tilde E}^\ast$-invariant hyperspace $P$. 
Hence, $U_L$ and $U_1$ are asymptotically nice domains with respect to $D^o$. 
(See Section \ref{du-sec-affine}.) 

By the premise, we have a connected convex polyhedral open subspace 
\[S_1:=\Bd (L_1 \ast \{\mbv_{\tilde E}\}) \cap {\mathcal T}_{\mbv_{\tilde E}}(\tilde \Sigma_{\tilde E})^o  
\subset \Bd (L_1 \ast \{\mbv_{\tilde E}\}).\] 
By Lemma \ref{pr-lem-predual} (iv),
$S_1$ corresponds to a connected convex polyhedral hypersurface 
\[S_1^{\ast} \subset \Bd U_1\] by 
 $\mathcal{D}^\Ag_{(L_1 \ast \{\mbv_{\tilde E}\})^o}$. 
 Since $S_1$ is disjoint from $L$ by the premise, 
 it follows $S_1^{\ast} \subset \Bd U_1 \cap U_L$ by 
 \eqref{prelim-eqn-dualinc}.
Since $S_1/\bGamma_{\tilde E}$ is compact, 
so is $S_1^{\ast}/\bGamma_{\tilde E}^\ast$ by Proposition \ref{prelim-prop-duality}. 
Theorem \ref{du-thm-lensn} shows that $D \cup S_1' = \Bd U_2$. 
Hence, $\Bd U_1 \cap U_L = S_1^{\ast}$. 


By Theorem \ref{du-thm-lensnpre}, we obtain an asymptotically nice closed
domain $U_2$ with connected strictly convex smooth hypersurface boundary $S_2$ 
in $U_L$ with $U_1 \cup S_1 \subset U_2^o$. 
The dual $U_2^{\ast}$ of $U_2$ has a connected strictly convex smooth hypersurface boundary 
$S_2^\ast$ in ${\mathcal T}_{\mbv_{\tilde E}}(\tilde \Sigma_{\tilde E})^o$ disjoint from $(L\ast \{\mbv_{\tilde E}\})^o$ and inside $(L_1\ast \{\mbv_{\tilde E}\})^o$ 
by \eqref{prelim-eqn-dualinc}. 
This is what we wanted.  

Also, considering $(L \ast \{\mbv_{\tilde E -}\})^o$ and 
$(L_1 \ast  \{\mbv_{\tilde E -}\})^o$, we obtain a connected strictly convex 
smooth hypersurface in the other component of ${\mathcal T}_{\mbv_{\tilde E}}(\tilde \Sigma_{\tilde E})^o - L$ in $U$.  
The union of the two hypersurfaces bounds a lens $L_2$ in $U$. 
(See Section \ref{prelim-sub-duality}.)

Let $F$ denote the compact fundamental domain 
of the boundary of the lens.  
The strictness of the lens follows 
from Proposition \ref{pr-prop-orbit} since the boundary of the lens 
is a union of orbits of $F$ and the limit points are only
in $\Lambda_{\tilde E}$. 

%
%
%
%
%
Again Proposition \ref{prelim-prop-closureind}  completes the proof for 
$\RP^n$. 
\end{proof} 

\begin{proof}[{\sl Proof of Theorem \ref{pr-thm-equiv}}] 
The forward direction follows from 
Proposition \ref{pr-prop-convhull2} by taking $\torb := \mathcal{T}_{\mbv_{\tilde E}}(\tilde \Sigma_{\tilde E})$.

Now, we show the converse. 
Let $L$ be a topological lens of the 
topological lens cone where $\bGamma_{\tilde E}$ acts cocompactly on. 
We may assume without loss of generality that 
the boundary is a union of two connected convex polyhedral hypersurfaces. 
Let ${\mathcal T}_{\tilde E}$ denote the tube ${\mathcal T}_{\mbv_{\tilde E}}(\tilde \Sigma_{\tilde E})$.

We denote by $h: \pi_{1}(\tilde E) \ra \SL_{\pm}(n+1, \bR)$ the holonomy
homomorphism of the end fundamental group with image $\bGamma_{\tilde E}$. 
We assume that the image of $h$ consists of matrices of the form \eqref{pr-eqn-mForm}.

By Theorem \ref{pr-thm-defspace}, $h$ is determined by $\hat h(\pi_1(\tilde E)) \ra \SL_\pm(n, \bR)$
and $\lambda_{\mbv_{\tilde E}}: \pi_1(\tilde E) \ra \bR_+$ evaluating the eigenvalue at 
$\mbv_{\tilde E}$ and a cocycle 
$\vec{b} \in H^1(\pi_1(\Sigma_{\tilde E}), \bR^{n \ast}_{\hat h, \lambda_{v_{\tilde E}}})$ 
for each $([\hat h], \lambda)$. 

(I) Suppose that we set the cocycle to be zero. 
Then we claim that there still exists 
a PL-lens for new homorphism $h': \pi_1(\tilde E) \ra \SL_\pm(n, \bR)$
where the boundary is a union of two convex polyhedral hypersurfaces. 
Here, we have an invariant hyperspace $P$ not meeting $\mbv_{\tilde E}$: 

Let $K$ be a compact set in the topologcal lens $L$. 
We know that 
$\{h(g_i)(K')\}$ converges to a compact set in the limit set 
$\Lambda_{\tilde E} :=\clo(L) \cap \Bd \mathcal{T}_{\tilde E}$ disjoint from $\mbv_{\tilde E}$ and 
$\mbv_{\tilde E -}$ by passing to a subsequence. 
By the definition of a topological lens, we know that $\Lambda_{\tilde E}$ satisfies the property of 
Lemma \ref{pr-lem-strictlens}. 

Also, we may deduce that any compact subset $K''$ in $\mathcal{T}_{\tilde E}$ is in a 
$h(\pi_1(\tilde E)$-invariant strict topological lens.
This can be obtained by taking a component of $\mathcal{T}_{\tilde E} - \partial_- L$
which is properly convex. We take the Hilbert metric,
form a neighborhood of $\partial_+ L$ and  the union with $L$. 
Similar construction applies by switching $\partial_- L$ and $\partial_+ L$.
The strictness of the newly obtained genralized lens follows from considering
 the Hilbert metric on the lines from $\mbv_{\tilde E}$ to its antipode.


Now, if we make $\vec b$ be zero, then the same happens with $\Lambda$ replaced with $P \cap \Bd \mathcal{T}_{\tilde E}$: 
The projected images of $h(g_i)(K')$ to $\Omega$ is the same as that of 
$h'(g_i)(K')$.


We use the coordinates on $\mathcal{T}_{\tilde E}$ given by $\mbv_{\tilde E}$ as the last coordinate axis direction and $P$ giving us the complementary subspace with $n$ coordinates. 
Let $\llrrparen{u, x} \in K$ where  $||u||=1, |x|$ is bounded above, 
$u$ is a vector in the direction of $P$, and 
$x$ is the coordinate in the direction of $\mbv_{\tilde E}$. 
Denote by $A_{g_i}$ the upper left $n\times n$-submatrix of $h(g_i)$
and $b_{g_i}$ the $n$-vector given by the $(n+1)$-row vector with the last element removed. 
Then 
\[ h(g_i)\llrrparen{u, x}  = \llrrparen{A_{g_i} u, b_{g_i} u + \lambda_{\mbv_{\tilde E}}(g_i) x} 
\ra \llrrparen{y, z}, ||y||=1\] 
where $|z|$ is bounded above since the limit is in $\Lambda$. 
Therefore, $||A_{g_i} u||/|(b_{g_i} u + \lambda_{\mbv_{\tilde E}}(g_i)x)|$
is uniformly bounded below. 

We can choose $x$ such that $b_g u$ and $\lambda_{\mbv_{\tilde E}}(g_i)x$ have the same sign
because we can find the strict topological lens to contain any compact set as we showed above.
So we can drop $b_{g_i} u$ and 
the above quantity is still uniformly bounded below. 


Suppose that some subsequence of $||A_{g_i} u||/|\lambda_{\mbv_{\tilde E}}(g_i) x|$ does not tend to infinity. Then we have $h'(g_i)(\llrrparen{u, x})$ does not accumulate only to $\Lambda'$
and $h'(g_i)(\llrrparen{u, -x})$ converges to distinct points from that of 
$h'(g_i)(\llrrparen{u, x})$. 
Even after restoring $b_{g_i} u$ back in
and we still have two distinct limits for the sequences 
$h(g_i)(\llrrparen{u, -x})$ and
$h(g_i)(\llrrparen{u, x})$. 
This contradicts the strict topological lens property of Lemma \ref{pr-lem-strictlens}.

So we proved that each orbit under $h'(\pi_1(\tilde E))$ only accumulates to the standard limit set. 
Then we can form a convex hull $C$ of the orbit of two points in each component of 
$\mathcal{T}_{\tilde E}- P$. 
The boundary $\partial C$ contains no infinite line. 
Otherwise, pull back using the elements of $\bGamma_{\tilde E}$,
 the sequence of points along the line to obtain
a biinfinite line in the boundary of $C$ as a geometric limit. 
This lines has to lie on $P$ since its two endpoints are on $\Lambda'$. 
However, since $C$ is $\bGamma_{\tilde E}$-invariant, this is a contradiction. 

Choose locally finite discrete sets in $\partial C$ invariant under $h'(\pi_1(E))$. 
Then take the convex hull of this to get a locally finite triangulated convex hull.

(II) Now, we prove the uniform middle eigenvalue condition: 
There is an abelianization map 
\[A: \pi_1(\tilde E) \ra H_1(\pi_1(\tilde E), \bR)\]
obtained by taking a homology class. 
The above map $g \ra \log \lambda_{\mbv_{\tilde E}}(h(g))$ 
induces homomorphism
\[\Lambda^{\prime h}: H_1(\pi_1(\tilde E), \bR) \ra \bR \]
that depends on the holonomy homomorphism $h$. 

Also, we may assume that $L$ is a PL-lens by taking a finite set of 
points in both components of $\partial L$ and taking the convex hull. 



Let us give an arbitrary Riemannian metric $\mu$ on $\Sigma_{\tilde E}$. 
Recall that a {\em current} is a transverse measure on a partial foliation by 
$1$-dimensional subspaces in the compact 
space $\Uu \Sigma_{\tilde E}$ on the transverse measure. (See \cite{RS}.)
These are not necessarily geodesic currents as in 
Bonahon \cite{Bonahon}. 
We denote the space of currents with the weak topology by $\mathcal{C}(\Uu \Sigma_{\tilde E})$

Hence, we assume that $\vec b$ is trivial, and $\Lambda$ lies on a hyperspace $P$ disjoint from
$\mbv_{\tilde E}$. 

The abelianization map 
$\pi_1(\Sigma_{\tilde E}) \ra H_1(\pi_1(\tilde E), \bR)$
can be understood as sending a closed curve to a current 
the corresponding homology class. 
This map extends to 
$\mathcal{C}(\Uu \Sigma_{\tilde E}) \ra H_1(\pi_1(\tilde E), \bR)$. 
(See Proposition 1 of \cite{RS} and Theorem 14 of \cite{deRham}.)
Also, $\Lambda^h:\pi_1(\tilde E)\ra \bR$ induces the continuous map 
$\hat \Lambda^h:\mathcal{C}(\Uu \Sigma_{\tilde E})\ra \bR$ 
which restricts to $\Lambda^h$ on the image currents of $\pi_1(\Sigma_{\tilde E})$: 
$\Lambda^h$ can be realized 
by integrating a $1$-form on $\Sigma_{\tilde E}$ along
the closed curve representing $\pi_1(\Sigma_{\tilde E})$ 
since $\Hom(H_1(\pi_1(\tilde E), \bR), \bR) = 
H^1(\pi_1(\tilde E), \bR)$. 
Since the integrations along currents are well-defined, the construction is complete.


Let $\lambda^{\mathrm{ul}}_{1}(h(g))$ denote the maximal norm of 
the eigenvalues $h(g)$
of the upper-left corner of $h$ in \eqref{pr-eqn-mForm}.
Obviously, $\lambda^{\mathrm{ul}}_{1}(h(g)) \geq \lambda_{\mbv_{\tilde E}}(h(g))$ 
for each $g\in \pi_1(\tilde E)$:  
If the eigenvalue of the upper-left corner matrix of $h(g)$ is strictly smaller than $\lambda^{\mathrm{ul}}_1(h(g))$, 
Proposition \ref{prelim-prop-attract}
shows that the closure of $L$ contains $\mbv_{\tilde E}$ 
or $\mbv_{\tilde E-}$
considering the orbit of $\{g^n\}$, which is a contradiction. 

Let 
$g \in \bGamma_{\tilde E}$, and 
let $[g]$ denote the current supported on a closed curve $c_{g}$ on $\Sigma_{\tilde E}$ corresponding to $g$ 
lifted to $\Uu \Sigma_{\tilde E}$.
Define $\leng_\mu(g)$ as the infimum of the $\mu$-length of 
such closed curves corresponding to $g$. 
Suppose that $\bGamma_{\tilde E}$ does not satisfy 
\[  C^{-1}\leng_\mu(g) \leq \log\frac{\lambda^{\mathrm{ul}}_{1}(h(g))}{\lambda_{\mbv_{\tilde E}}(h(g))} \leq C \leng_\mu(g)  \]
for a uniform constant $C> 1$. 
Then there exists a sequence $g_i$ of elements of $\bGamma_{\tilde E}$ such that 
\begin{equation}\label{pr-eqn-lambdaul}  
\left\{\frac{\log\left(\frac{\lambda^{\mathrm{ul}}_{1}(h(g_i))}{\lambda_{\mbv_{\tilde E}}(h(g_i))}\right)}{ \leng_\mu(g_i)} \right\} \ra 0 \hbox{ as } i \ra \infty.
\end{equation} 

Let $[g_{\infty}]$ denote a limit point of $\{[g_i]/\leng_\mu(g_i)\}$ in the space of currents on 
$\Uu\Sigma_{\tilde E}$. Since $\Uu \Sigma_{\tilde E}$ is compact, a limit point exists. 
We may modify $h$ by changing the homomorphism 
$g \mapsto \lambda_{\mbv_{\tilde E}}(h(g))$ only; 
that is, we only modify the $(n+1)\times (n+1)$-entry of the matrices form \eqref{pr-eqn-mForm} with corresponding changes. 
By Proposition \ref{app-prop-lensP}, 
the perturbed cocompactly-acted PL-lens $L'$,
constructed with a modified holonomy representation, 
is still a properly convex domain with the same tube. 
Since its boundary polyhedrons are still tranversal to the great segments from 
$\mbv_{\tilde E}$,  its closure does not contain $\mbv_{\tilde E}$.

Suppose that $\hat \Lambda^h([g_\infty])=0$. Then 
we have 
\[ \frac{\log \lambda_{\mbv_{\tilde E}}(g_i)}{\leng_K(g_i)} \ra 0.\]
This implies 
\[  
 {\lambda_{\mbv_{\tilde E}}(g_i)}^{\frac{1}{\leng_K(g_i)}} \ra 1.
\] 
Equation \ref{pr-eqn-lambdaul} and this imply 
 \[  
 {\lambda_{1}^{ul}(g_i)}^{\frac{1}{\leng_K(g_i)}} \ra 1.
\] 
Denote by $\tilde \lambda_1(g_i)$ the largest eigenvalue of $\hat h(g_i)$. 
Since
$ {\lambda_{1}}^{ul}(g_i) = \frac{\tilde \lambda_1(g_i)}{ \lambda_{\mbv_{\tilde E}}^{\frac{1}{n}}(g_i)}$
By equation \eqref{pr-eqn-eigenlem}, 
we have 
\[ 
 {\tilde \lambda_{1}(g_i)}^{\frac{1}{\leng_K(g_i)}} \geq e^{1/n}.
\]
This implies that 
\[ 
\frac{{\tilde \lambda_{1}(g_i)}^{\frac{1}{\leng_K(g_i)}}}{  {\lambda_{\mbv_{\tilde E}}(g_i)}^{\frac{1}{n\leng_K(g_i)}} } \ra  C \geq e^{1/n}
\]
for some $C> 0$.  This is a contradiction. 
Hence, $[g_\infty] \ne 0$ in $H_1(\Sigma_{\tilde E}, \bR)$, and 
$\Lambda^h$ is not zero on $[g_\infty]$. 

By considering the image of $[g_\infty]$ in $H_1(\Sigma_{\tilde E}, \bR)$, 
we can make a sufficiently 
small change of $h$ to $h'$ by changing $\lambda_{\mbv_{\tilde E}}$ and normalizing $h$ 
such that $\Lambda^{h'}([g_{\infty}]) > \Lambda^{h}([g_\infty])$.

From this, we obtain that
\begin{equation} \label{pr-eqn-negative} 
\log\left(\frac{\lambda^{\mathrm{ul}}_{1}(h'(g_{i}))}{\lambda_{\mbv_{\tilde E}}(h'(g_i))}\right)< 0  \hbox{ for sufficiently large } i. 
\end{equation} 


By \eqref{pr-eqn-negative},  
we obtain that $\lambda^{\mathrm{ul}}_{1}(h'(g)) < \lambda_{\mbv_{\tilde E}}(h'(g))$ for some $g$
and that $\lambda_{\mbv_{\tilde E}}(h'(g))$ at $\mbv_{\tilde E}$.
Hence, we can decompose $\SI^n$ into a hyperspace $S'$ and the complementary
$\{\mbv_{\tilde E}, \mbv_{\tilde E-}\}$. 
The norms of eigenvalues associated with $S'$ are strictly less than 
that of $\mbv_{\tilde E}$. 
Proposition \ref{prelim-prop-attract} shows 
that the closure of $L$ contains 
$\mbv_{\tilde E}$ or $\mbv_{\tilde E-}$ by 
considering the orbits under $\{g^i\}$.
\end{proof}

We thank Daryl Cooper for pointing out a gap in an earlier version. 





\subsection{The uniform middle-eigenvalue conditions and the lens-shaped ends.} \label{pr-sub-umecl} 



Now, we aim to prove Theorem \ref{pr-thm-secondmain} restated 
as Theorem \ref{pr-thm-equ}.  
A {\em radially foliated end neighborhood system} of a real projective $n$-orbifold $\mathcal{O}$
is a collection of end neighborhoods of $\mathcal{O}$ that is radially foliated 
and outside a compact suborbifold of $\mathcal{O}$ 
whose interior is isotopic to $\mathcal{O}$.  


\begin{definition} \label{pr-defn-tri}
We say that a convex 
$\mathcal{O}$ with $\torb \subset \SI^n$ (resp. $\subset \RP^n$) 
satisfies the {\em triangle condition} if 
for any fixed end neighborhood system of $\mathcal{O}$ and 
every triangle $T \subset \clo(\torb)$, 
if \[\partial T \subset  \Bd \torb, T^{o} \subset \torb, \hbox{ and  }
\partial T \cap \clo(U)\ne \emp\] 
for a radial p-end neighborhood $U$,  then 
$\partial T$ is a subset of $\clo(U) \cap \Bd \torb$. 
\end{definition}  
\index{triangle condition|textbf} 

For example, 
by Corollary \ref{app-cor-stLens}, strongly tame strict SPC-orbifolds with 
generalized lens-shaped or horospherical ends satisfy this condition. 
The converse does not necessarily hold. 


A minimal $\bGamma_{\tilde E}$-invariant distanced compact set is 
the smallest compact $\bGamma_{\tilde E}$-invariant distanced set in 
${\mathcal T}_{\tilde E}$: This follows since we can take intersections of such sets to obtain such a set.  

\begin{lemma} \label{pr-lem-coneseq} 
Suppose that $\mathcal{O}$ is a convex real projective $n$-orbifold 
and satisfies the triangle condition. 
Then no triangle $T$ with $T^o \subset \torb, \partial T \subset \Bd \torb$ 
has a vertex equal to an R-p-end vertex.
\end{lemma} 
\begin{proof} 
	Assume $\torb \subset \SI^n$. 
Let $\mbv_{\tilde E}$ be a p-end vertex. Choose a fixed radially foliated p-end neighborhood system. 
Suppose that a triangle $T$ with $\partial T \subset \Bd \torb$ contains a vertex equal to a p-end vertex. 
Let $U$ be a radial p-end neighborhood of a p-end $\tilde E$ with 
a p-end vertex $\mbv_{\tilde E}$. 
By the triangle condition, $\partial T \subset \clo(U) \cap \Bd \torb$. 

Since $U$ is foliated by radial lines from $\mbv_{\tilde E}$, 
we choose $U$ such that $\Bd U \cap \torb$ covers a compact hypersurface
in $\orb$. 
Let $\mathcal U$ denote the set of segments in $\clo(U)$ from $\mbv_{\tilde E}$. 
Every segment in $\mathcal U$ in the direction of 
$\tilde \Sigma_{\tilde E}$ ends in $\Bd U \cap \torb$. 
Also, the segments $\mathcal U$ in directions of 
$\Bd \tilde \Sigma_{\tilde E}$ lie in $\Bd U \cap \Bd \torb$
by the definition of $\tilde \Sigma_{\tilde E}$. 
Also, $\clo(U)$ is a union of segments in 
$\mathcal U$. 
Thus, $\clo(U) \cap \Bd \torb$ is a union of segments 
in directions of $\Bd \tilde \Sigma_{\tilde E}$. 

Since $T^o \subset \torb$, each segment in $\mathcal U$ with 
interior in $T^o$ lies in directions not contained in  $\partial \clo(\tilde \Sigma_{\tilde E})$. 
Let $w$ be the endpoint of the maximal extension in $\torb$ of such 
a segment. Then $w$ is not in $\clo(U) \cap \Bd \torb$
by the conclusion of the above paragraph. 
This contradicts
$\partial T \subset \clo(U) \cap \Bd \torb$.

The proof for $\RP^n$ case follows from Proposition \ref{prelim-prop-closureind}. 
%
\hfill \SnP {\parfillskip0pt\par}
\end{proof}

\begin{lemma} \label{pr-lem-genlens} 
Suppose that $\mathcal{O}$ is a convex real projective $n$-orbifold 
and satisfies the triangle condition, or alternatively, an R-p-end $\tilde E$ is virtually factorizable.  
Then the R-p-end $\tilde E$ is generalized lens-shaped if and only if it is  lens-shaped. 
\end{lemma}
\begin{proof} 
	Again, we prove the statements for $\SI^n$. 
If $\tilde E$ is virtually factorizable, this follows from Theorem \ref{pr-thm-redtot}. 

Suppose that $\tilde E$ is not virtually factorizable. 
Now assume the triangle condition. 
Given a generalized cocompactly-acted lens $L$, let $L^{b}$ denote $\clo(L) \cap {\mathcal T}_{\mbv_{\tilde E}}(\tilde \Sigma_{\tilde E})$. 
We obtain the convex hull $M$ of $L^{b}$. $M$ is a subset of $\clo(L)$. 
The lower boundary component of $L$ is a smooth, strictly convex hypersurface. 

Let $M_{1}$ be the outer component of $\Bd M \cap {\mathcal T}_{\mbv_{\tilde E}}(\tilde \Sigma_{\tilde E})$. 
Suppose that $M_{1}$ meets $\Bd \torb$. 
$M_{1}$ is a union of interiors of simplices. 
By Lemma \ref{prelim-lem-simplexbd}, either a simplex $\sigma$ in $\clo(\torb)$ 
is in $\Bd \torb$ or 
its interior $\sigma^o$ is disjoint from it. 
Hence, there is a simplex $\sigma$ in $M_{1}\cap \Bd \torb$. 
Taking the convex hull of $\mbv_{\tilde E}$ and an edge in $\sigma$, 
we obtain a triangle $T$ with $\partial T \subset \Bd \torb$ and $T^{o}\subset \torb$. 
This contradicts the triangle condition by Lemma \ref{pr-lem-coneseq}. 
Thus, $M_{1}\subset \torb$. 
By Theorem \ref{pr-thm-equ}, the end satisfies the uniform middle-eigenvalue condition. 
By Proposition \ref{pr-prop-convhull2}, we obtain a lens cone in $\torb$. 
\hfill \SSn {\parfillskip0pt\par}
\end{proof}

\begin{theorem}\label{pr-thm-equ} 
Let $\orb$ be a convex real projective $n$-orbifold. 
Let $\bGamma_{\tilde E}$ be the holonomy group of a properly convex R-end $\tilde E$
and the end vertex $\mbv_{\tilde E}$. 
Then the following are equivalent\/{\rm :}
\begin{enumerate}
\item[{\rm (i)}] $\tilde E$ is a generalized lens-shaped R-end.
\item[{\rm (ii)}] $\bGamma_{\tilde E}$ 
satisfies the uniform middle-eigenvalue condition
with respect to $\mbv_{\tilde E}$. 
\end{enumerate}
If ${\mathcal{O}}$ furthermore 
satisfies the triangle condition or, alternatively, assume that $\tilde E$ is virtually factorizable, then the following holds\/{\rm :} 
\begin{itemize} 
\item $\bGamma_{\tilde E}$ is lens-shaped if and only if 
 $\bGamma_{\tilde E}$ satisfies the uniform middle-eigenvalue condition with respect to $\tilde E$. 
\end{itemize}
\end{theorem}
\begin{proof} 
	Assume $\torb \subset \SI^n$. 
(ii) $\Rightarrow$ (i): This follows from Theorem \ref{pr-thm-equiv} since 
we can intersect the lens with $\torb$
and add the top hypersurface
 to obtain a generalized lens and  a generalized lens cone from it. (In particular, $\pi_1(\tilde E)$ acts cocompactly on
the generalized lens.) 
(i) $\Rightarrow$ (ii):   
Let $L$ be a generalized cocompactly-acted lens in the generalized 
lens cone $L\ast \mbv_{\tilde E}$. 
Let $B$ be the lower boundary component of $L$ in
the tube $\mathcal{T}_{\mbv_{\tilde E}}(\tilde \Sigma_{\tilde E})$. Since $B$ is strictly convex, 
the upper component of 
$\mathcal{T}_{\mbv_{\tilde E}}(\tilde \Sigma_{\tilde E}) -B$ is a properly convex domain, which we denote by 
$U$. 
Let $l_x$ denote the maximal segment from $\mbv_{\tilde E}$ passing $x$ for 
$x \in U - L$. 

We define a function $f: U - L \ra \bR$ given 
by $f(x)$ as the Hilbert distance on line 
$l_x$ from $x$ to $L \cap l_x$. 
Then a level set of $f$ is always strictly convex: 
This follows from considering a $2$-plane $P$ containing $\mbv_{\tilde E}$ and 
passing through $L$. Let $x,y $ be points of $f^{-1}(c)$ for a constant $c> 0$. 
Let $x'$ be the point of $\clo(L) \cap l_x$ closest to $x$ 
and $y'$ be one of $\clo(L)\cap l_y$ closest to $y$. 
Let $x''$ be one of $\clo(L)\cap l_x$ furthest from $x$.
Let $y''$ be one of $\clo(L)\cap l_y$ furthest from $x$.
Since $f(x) = f(y)$, a cross-ratio argument shows 
that the lines extending $\overline{xy}, \overline{x'y'}$ and 
$\overline{x''y''}$ are concurrent outside $U\cap P$. 
The strict convexity of $B$ and Lemma 1.8 of \cite{CLT15} show that 
$f(z) < \eps$ for $z \in \overline{xy}^o$. 



We can approximate each level set by a convex polyhedral hypersurface in $U- L$ 
by selecting a finite set in the level set and taking the convex hull of the orbit
taking advantage of the strict convexity of the level sets. 
We do the same for the lower hypersurface of $L$. 
We take the convex hull and take a boundary component which is 
a convex polyhedral hypersurface.
Let $V$ denote the domain bounded by this hypersurface and $B$.   
Theorem \ref{pr-thm-equiv} implies (ii).

The final part follows from Lemma \ref{pr-lem-genlens}.
\hfill \SSn {\parfillskip0pt\par}
\end{proof}

\section{The properties of lens-shaped ends.} \label{pr-sec-lens}

\begin{lemma}\label{pr-lem-concavebd} 
	Let $\orb$ be a properly convex $n$-orbifold. Then
	given any end neighborhood, there exists a concave end neighborhood in it. 
Furthermore, the $d_\orb$-diameter of the boundary of a concave end neighborhood of
an R-end $E$ is 
bounded by the Hilbert diameter of the end orbifold $\Sigma_E$ of $E$. 	
	\end{lemma} 
\begin{proof} 
	It is sufficient to prove for the case $\torb \subset \SI^n$. 
		Suppose that we have a generalized lens cone $V$
	that is a p-end neighborhood equal to the interior of $L \ast v_{\tilde E}$ 
	where $L$ is a generalized cocompactly-acted lens bounded away from $v_{\tilde E}$. 
	
	Now take a p-end neighborhood $U'$. We assume without loss of generality that 
	$U'$ covers a product end neighborhood with compact boundary. 
	
	By taking smaller $U'$ if necessary, we may assume that $U'$ and $L$ are disjoint. 
	Since both $(\Bd U'\cap \torb)/h(\pi_1(\tilde E))$ and $L/h(\pi_1(\tilde E))$ are compact, 
	$\eps > 0$. 
	Let
	\[L' := \{x \in V | d_{V}(x, L) \leq \eps\}.\] 
	Since a lower component of $\partial L$ is strictly convex, 
	we can show that $L'$ can be polyhedrally approximated and 
	smoothed into be a cocompactly-acted lens 
	by Proposition \ref{pr-prop-lenssmooth}.
	
	Clearly, $h(\pi_1(\tilde E))$ acts on $L'$. 
	
	We choose sufficiently large $\eps'$ such that
	$\Bd U \cap \torb \subset L'$. 
Hence $V-L'\subset U$ form a concave p-end neighborhood as above.


	Let $\tilde E$ be a p-end corresponding to $E$. 
	Let $U$ be a concave p-end neighborhood of $\tilde E$ that is a cone: 
	$U $ is the interior of $\{\mbv\}\ast L - L$ where $L$ is a generalized cocompactly-acted lens, and $v$ is the p-end vertex corresponding to $U$. 
	Let $T$ denote the tube with vertex $\mbv$ in the direction of $L$. 
	Then $B:= \Bd U \cap \torb$ is a smooth lower boundary component of $L$. 
	
	Any tangent hyperspace $P$ at a point of $B$ meets $\Bd T$ in a sphere of 
	dimension $n-2$. By the convexity of $L$ and 
	the strict convexity of $B$, it follows that  
	$P \cap L$ is a point.
	We claim that $P\cap \Bd \torb \subset P\cap \Bd T$: 
	We put $T$ into an affine space $\mathds{A}^n$ with 
	vertices in $\Bd \mathds{A}^n$. Then $T$ is foliated 
	by parallel complete affine lines. Consider these as vertical lines.
	$B$ is a strictly convex hypersurface 
	meeting these vertical lines transversely. 
	Hence, the property of $P$ is evident.
	
	Thus any maximal segment in $\orb$ tangent to $B$ at $x$ must end in 
	$\Bd T\cap \Bd \torb$.
	There is a projection 
	$\Pi_{\mbv_{\tilde E}}: B \ra \tilde \Sigma_{\tilde E}$
	that is a diffeomorphism. 
	Hence, under $\Pi_{\mbv}$, the maximal segment is sent to a maximal segment of 
	$\tilde \Sigma_{\tilde E}$ , which forms an isometry 
	on the segments with the Hilbert metrics. 
	Moreover this is a Finsler isometry by comparing the Finsler metric 
	restricted to the tangent space to $B$ at $x$ to 
	that of the tangent space to $\tilde \Sigma_{\tilde E}$ 
	at $\Pi_{\mbv_{\tilde E}}(x)$.  
	The conclusion follows. \hfill \SSn {\parfillskip0pt\par}
	\end{proof}

\subsection{The properties for a lens cone in non-virtually-factorizable cases} 
Recall that each infinite-order $g \in \bGamma_{\tilde E}$ 
is positive bi-semiproximal 
by Proposition \ref{prelim-prop-nonu}. 

\begin{theorem}\label{pr-thm-lensclass}
Let $\mathcal{O}$ be a properly convex real projective $n$-orbifold. 
Let $\tilde E$ be an R-p-end of $\tilde{\mathcal{O}} \subset \SI^n$
{\rm (}resp. in $\RP^n${\rm )} 
with a generalized lens p-end neighborhood. 
Let $\mbv_{\tilde E}$ be the p-end vertex. 
Assume that $\pi_1(\tilde E)$ is non-virtually-factorizable. 
Then 
$\bGamma_{\tilde E}$ satisfies the uniform middle-eigenvalue condition with 
respect to $\mbv_{\tilde E}$, and
there exists a generalized cocompactly-acted lens $D$ disjoint from $\mbv_{\tilde E}$
with the following properties\/{\em :}
\begin{enumerate} 
\item[{\rm (i)}] 
\begin{itemize}
\item  The set $\Bd D -\partial D = \Lambda_{\tilde E}$ is independent of the choice of $D$
where $\Lambda_{\tilde E}$ is from Proposition \ref{pr-prop-orbit}. 
\item $D$ is strictly generalized lens-shaped. 
\item Each nontorsion nonidentity element $g \in \bGamma_{\tilde E}$ has an attracting 
fixed set in $\Bd D$ intersected with the union of some great segments from $\mbv_{\tilde E}$ in the directions in $\Bd \tilde \Sigma_{\tilde E}$. 
\item 
The closure of the union of all attracting fixed sets is a subset of $\Bd D - A - B$ for the top and the bottom hypersurfaces $A$ and $B$. 
The closure equals $\Bd D - A - B$ if $\bGamma_{\tilde E}$ is 
hyperbolic. 
\end{itemize}
\item[{\rm (ii)}] 
\begin{itemize} 
\item Let $l$ be a segment $l \subset \Bd \tilde{\mathcal{O}}$ 
with $l^o \cap \clo(U) \neq \emp$ for  any concave p-end neighborhood $U$ of $\mbv_{\tilde E}$. Then $l$ lies in $\bigcup S(\mbv_{\tilde E})$ and also
in the closure in $\clo(V)$ of any concave or proper p-end neighborhood $V$ of $\mbv_{\tilde E}$. 
\item The set $S(\mbv_{\tilde E})$ of maximal segments from $\mbv_{\tilde E}$ in $\clo(V)$ is independent of 
a concave or proper p-end neighborhood $V$
{\rm (}In fact, it is the set of maximal segments from $\mbv_{\tilde E}$ ending 
in $\Bd D - A - B$\/{\rm ).}
\item \[\bigcup S(\mbv_{\tilde E}) = \clo(V) \cap \Bd \tilde{\mathcal{O}}.\]
\end{itemize} 
\item[{\rm (iii)}] $S(g(\mbv_{\tilde E}))=g(S(\mbv_{\tilde E}))$ for $g \in \pi_1(\tilde E)$. 
\item[{\rm (iv)}] 
Given $g \in \pi_1(\orb)$, we have 
\begin{equation} 
\left(\bigcup S(g(\mbv_{\tilde E}))\right)^o
\cap \bigcup S(\mbv_{\tilde E}) = \emp \hbox{ or else } 
\bigcup S(g(\mbv_{\tilde E})) =
\bigcup S(\mbv_{\tilde E}) \hbox{ with } g \in \bGamma_{\tilde E}. 
\end{equation} 
\item[{\rm (v)}] A concave p-end neighborhood is a proper p-end neighborhood. 
\item[{\rm (vi)}] Assume that $\bv_{\tilde E'}$ is 
the p-end vertex of an R-p-end $\tilde E'$. 
We can choose mutually disjoint concave p-end neighborhoods for 
every R-p-ends.
Then 
\[\left(\bigcup S(\mbv_{\tilde E})\right)^o \cap \bigcup S(\mbv_{\tilde E'}) = \emp 
	\hbox{ or } \bigcup S(\mbv_{\tilde E})= \bigcup S(\mbv_{\tilde E'}) \hbox{ with } \mbv_{\tilde E} = \mbv_{\tilde E'}, \tilde E= \tilde E'\] 
for an R-p-end vertice $\mbv_{\tilde E}$. 
{\rm (}This is a sharping of (iv){\rm .)}
\end{enumerate}
\end{theorem}   
\begin{proof}  
	Suppose first $\torb \subset \SI^n$. 
%
Theorem \ref{pr-thm-equ} implies the uniform middle-eigenvalue condition. 

(i) 
%
Let $U_1$ be a concave end neighborhood. 
Since $\bGamma_{\tilde E}$ acts on $U_1$, 
$U_1$ is a component of the complement of a generalized lens $D$ in 
a generalized R-end of form $D \ast \{\mbv_{\tilde E}\}$ by definition. 
The action on $D$ is cocompact and proper since we can use 
a foliation by great segments in a tube corresponding to 
$\tilde E$. 

Proposition \ref{pr-prop-orbit} 
implies that the lens is strict.
This implies (i).

(ii) Consider any segment $l$ in $\Bd \tilde{\mathcal O}$ with $l^o$ 
meeting $\clo(U_1)$ for a concave p-end neighborhood $U_1$ of $\mbv_{\tilde E}$. 
Here, the generalized lens $D$ has boundary components $A$ and $B$ 
where $B$ is also a boundary component of $U_1$ in $\torb$. 
Let $T$ be the open tube corresponding to $\tilde \Sigma_{\tilde E}$. 
Then $\torb \subset T$ since $\tilde \Sigma_{\tilde E}$ is the direction of 
all segments in $\torb$ starting from $\mbv_{\tilde E}$. 
Let $T_{1}$ be a component of $\Bd T - \partial_{1} B $ containing $\mbv_{\tilde E}$.
Then $T_{1} \subset \clo(U_{1}) \cap \Bd \torb$ by the definition of concave p-end neighborhoods. 
In the closure of $U_{1}$, an endpoint of $l$ is in $T_{1}$. 
Then $l^{o} \subset \Bd T$ since  $l^{o}$ is tangent to $\Bd T - \{\mbv_{\tilde E}, \mbv_{\tilde E-}\}$. 
Any convex segment $s$ from $\mbv_{\tilde E}$ to any point of $l$
must be in $\Bd T$. By the convexity of $\clo(\torb)$, 
we have $s \subset \clo(\torb)$. 
Thus, $s$ is in $\Bd \torb$ since $\Bd T \cap \clo(\torb) \subset \Bd \torb$. 
Therefore, the segment $l$ is contained in the union of segments in $\Bd \torb$ from $\mbv_{\tilde E}$. 

We now suppose that $l$ is a segment from 
$\mbv_{\tilde E}$ containing a segment $l_0$ in $\clo(U_1)\cap \Bd \tilde{\mathcal O}$ from $\mbv_{\tilde E}$, 
and we will show that $l$ is in $\clo(U_1) \cap \Bd \torb$, which 
will be sufficient to prove (ii). 
$l^{o}$ contains a point $p$ of $\Bd D - A - B$, 
which is a subset of $\Bd \mathcal{T}_{\mbv_{\tilde E}}(\tilde \Sigma_{\tilde E}) \cap D$. 
Since $l \subset \clo(\torb)$, we obtain
$\bigcup_{g \in \bGamma_{\tilde E}} g(l) \subset \clo(\torb)$, 
a properly convex subset. 
Hence, $\bigcup_{g \in \bGamma_{\tilde E}}  
g(l) - U_1$ is a distanced set, and 
has a distanced compact closure. Then the convex hull of the closure 
meets $\Bd \mathcal{T}_{\mbv_{\tilde E}}(\tilde \Sigma_{\tilde E})$ in a way contradicting  
Proposition \ref{pr-prop-orbit} (ii) 
where $D$ is $\Lambda_{\tilde E}$ in the proposition. 
Thus, $l^o$ does not meet  $\Bd D - A - B$. 
 Thus, \[l \subset \clo(U_1)\cap \Bd \tilde{\mathcal O}.\]

We define $S(\mbv_{\tilde E})$ as the set of maximal segments
in $\clo(U_1) \cap \Bd \torb$. 
Such a maximal segment is also maximal in $\clo(U) \cap \Bd \torb$ by 
the above paragraph. 
Hence, we can characterize $S(\mbv_{\tilde E})$ as the set of 
maximal segments in $\Bd \tilde{\mathcal O}$ from $\mbv_{\tilde E}$ 
ending at points of  $\Bd D - A - B$.
Also, $\bigcup S(\mbv_{\tilde E})= \clo(U_1)\cap \Bd \torb$. 

For any other concave affine neighborhood $U_2$ of $U_1$, 
we have 
\[U_2 = \{\mbv_{\tilde E}\}\ast D_2 -D_2 -\{\mbv_{\tilde E}\}\] for a generalized cocompactly-acted lens $D_2$. 
Since $\clo(D_2) - \partial D_2= \clo(D) - \partial D$, it follows
that $\clo(U_2) \cap \Bd \torb = \clo(U_1)\cap \Bd \torb = \bigcup S(\mbv_{\tilde E})$.

%

Let $U'$ be any proper p-end neighborhood associated with 
$\mbv_{\tilde E}$. 
By Lemma \ref{pr-lem-concavebd}, there exists a concave p-end 
neighborhood $U_1 \subset U'$. 
Again, $U_1 = \{\mbv_{\tilde E}\} \ast D - D - \{\mbv_{\tilde E}\}$ where 
$D$ is a generalized cocompactly-acted lens and $\mbv_{\tilde E} \not\in \clo(D)$.  
Hence, $\clo(U_1)\cap \Bd \torb \subset \clo(U')\cap \Bd \torb$. 
Moreover, every maximal segment in $S(\mbv_{\tilde E})$ is in 
$\clo(U')$. 
 
%

We can form $S'(\mbv_{\tilde E})$ as
the set of maximal segments from $\mbv_{\tilde E}$ in $\clo(U') \cap \Bd \tilde{\mathcal O}$.
Then no segment $l$ in $S'(\mbv_{\tilde E})$ has interior points in $\Bd D - A - B$ as above. 
 Thus, \[S(\mbv_{\tilde E}) = S'(\mbv_{\tilde E}).\] 
 Also, since every points of $\clo(U') \cap \Bd \torb$ has a segment in the direction of $\Bd \tilde \Sigma_{\tilde E}$, 
 we obtain 
 \[\bigcup S(\mbv_{\tilde E}) = \clo(U') \cap \Bd \torb.\]

 (iii) Since $g(D)$ is the generalized cocompactly-acted lens for the the generalized lens neighborhood $g(U)$ of $g(\mbv_{\tilde E})$, 
we obtain $g(S(\mbv_{\tilde E}))= S(g(\mbv_{\tilde E}))$ for any p-end vertex $\mbv_{\tilde E}$.

(iv)  Let $U$ be a proper p-end neighborhood of $\tilde E$ covering an end neighborhood 
of a product form and having compact boundary.
We choose a generalized cocompactly-acted lens $L$ of a generalized lens cone 
such that $C_{\tilde E}:= \{\mbv_{\tilde E}\}\ast L - L-  \{\mbv_{\tilde E}\}$ lies in $U$ 
by Lemma \ref{pr-lem-concavebd}.  
We can choose $C_{\tilde E}$ to be a proper concave p-end neighborhood
since such a choice is always possible in a proper p-end neighborhood.  
The properness of $U$ shows that
 \begin{equation} \label{pr-eqn-CE} 
 g(C_{\tilde E})=C_{\tilde E} \hbox{ for } g \in \bGamma_{\tilde E}, 
 \hbox{ or else } g(C_{\tilde E})\cap C_{\tilde E} = \emp
 \hbox{ for every } g \in \bGamma_{\tilde E}.
 \end{equation} 
 



Let $B_L$ denote the boundary component of $L$ meeting the closure of 
$C_{\tilde E}$.
Now, $\bigcup S(\mbv_{\tilde E})^o$ has an open neighborhood 
$C_{\tilde E} \cup \bigcup S(\mbv_{\tilde E})^o$ in $\orb$ 
since $B_L$ is a separating hypersurface in $\torb$.  
We obtain the conclusion since  
the intersection of the two sets implies the intersections of the neighborhoods of
the sets. 

(v) Let $C_{\tilde E}$ be a concave p-end neighborhood
$\{\mbv_{\tilde E}\}\ast L - L-  \{\mbv_{\tilde E}\}$ for a lens $L$. 
We will now show that $C_{\tilde E}$ is a proper p-end neighborhood.  
Suppose for contradiction that 
\[g (C_{\tilde E}) \cap C_{\tilde E} \ne \emp \hbox{ and } g(C_{\tilde E}) \ne C_{\tilde E} \hbox{ for } 
g \in \bGamma_{\tilde E}.\]



Since $C_{\tilde E}$ is concave, 
each point $x$ of $\Bd C_{\tilde E} \cap \torb$ is contained in a sharply supporting  hyperspace 
$H$ such that 
a component $C$ of $\torb - H$ is in $C_{\tilde E}$ 
where $\clo(C) \ni \mbv_{C_{\tilde E}}$ for the p-end vertex $\mbv_{C_{\tilde E}}$ of $C_{\tilde E}$. 
Similar statements hold for $g(C_{\tilde E})$. 

Since $g(C_{\tilde E}) \cap C_{\tilde E} \ne \emp$ 
and $g(C_{\tilde E}) \ne C_{\tilde E}$, 
it follows that 
\[\Bd g(C_{\tilde E}) \cap C_{\tilde E} \ne \emp \hbox{ or } g(C_{\tilde E}) \cap \Bd C_{\tilde E} \ne \emp.\]
Assume the second case without loss of generality. 
Let $x \in \Bd C_E$ in $g(C_E)$, and choose $H, C$ as above. 
Let $\clo(C)$ be the closure 
of a component $C$ of $\clo(\torb) - H$ containing $\mbv_{\tilde E}$ 
where $H$ is a separating hyperspace.
$C\cap \Bd \torb$ is a union of lines in $S(\mbv_{\tilde E})$.
Now, $H \cap g(C_{\tilde E})$ contains an open neighborhood in $H$ 
of $x$. 

Since $H$ contains a point of a concave p-end neighborhood 
$g(C_{\tilde E})$ of $g(\tilde E)$, 
it meets a points of $\{g(\mbv_{\tilde E})\}\ast g(D) - g(D)- \{g(\mbv_{\tilde E})\}$  for a lens $D$ of $\tilde E$ 
and a ray from $g(\mbv_{\tilde E})$ in $g(C_{\tilde E})$. 
We deduce that 
$H \cap g(C_{\tilde E})$ separates $g(C_{\tilde E})$ 
into two open sets $C_1$ and $C_2$ in the direction of one of the sides of $H$  
where $\clo(C_1)-H$ and $\clo(C_2) - H$ meet 
$g\left( \bigcup S(\mbv_{\tilde E}) \right)^o$
at nonempty sets. 
One of $C_1$ and $C_2$ is in $C$ since $C$ is a component of 
$\torb - H$. 
Also, $\clo(C)-H$ meets the set at 
$(\clo(C)-H) \cap \Bd \orb \subset \bigcup S(\mbv_{\tilde E})$. 
Hence, this implies 
\[g\left( \bigcup S(\mbv_{\tilde E}) \right)^o \cap 
\bigcup S(\mbv_{\tilde E}) \ne \emp.\] 
By (iv), this means $g \in \pi_1(\tilde E)$. 
Hence, $g(C_{\tilde E}) = C_{\tilde E}$ and this is absurd.
We have 
\[ g(C_{\tilde E}) \cap C_{\tilde E} = \emp \hbox{ or }
g(C_{\tilde E}) = C_{\tilde E} \hbox{ for } g \in \pi_1(\orb).
\] 
Since $g$ acts on $C_{\tilde E}$ and the maximal segments in 
$S(\mbv_{\tilde E})$ must go to maximal segments, and the
interior points of maximal segments cannot be an image of $\mbv_{\tilde E}$,
we must have $g(\mbv_{\tilde E}) = \mbv_{\tilde E}$. 
Hence, $g(U) \cap U \ne \emp$ for any proper p-end neighborhood of 
$\tilde E$, and $g \in \pi_1(\tilde E)$.  






(vi) Suppose that 
$S(\mbv_{\tilde E})^o\cap S(\mbv_{\tilde E'}) \ne \emp$. Then
$S(\mbv_{\tilde E})^o \cup C_{\tilde E}$ is 
a neighborhood of $S(\mbv_{\tilde E})^o$ for 
a proper concave p-end neighborhood $C_{\tilde E}$ of $\tilde E$. 
Also, $S(\mbv_{\tilde E'})^o \cup C_{\tilde E}$ is 
a neighborhood of $S(\mbv_{\tilde E'})^o$ for 
a proper concave p-end neighborhood $C_{\tilde E'}$ of $\tilde E'$. 

The above argument in (iv) applies within this situation 
to show that $\tilde E = \tilde E'$ and 
$\mbv_{\tilde E} = \mbv_{\tilde E'}$. 
%

Proposition \ref{prelim-prop-closureind} implies the version for 
$\RP^n$. 
\hfill \SSn {\parfillskip0pt\par}
\end{proof}

\subsection{The properties of lens cones for factorizable case}

\begin{theorem}\label{pr-thm-redtot}
Let $\mathcal{O}$ be a properly convex real projective $n$-orbifold. 
Suppose that 
\begin{itemize}
\item $\clo(\orb)$ is not of form $\mbv_{\tilde E} \ast D$ for a totally geodesic 
properly convex domain $D$, or 
\item the holonomy group $\bGamma$ is strongly irreducible. 
\end{itemize} 
Let $\tilde E$ be an R-p-end of the universal cover $\torb$, $\torb \subset \SI^n$ 
{\rm (}resp. $\subset \RP^n$\/{\rm )} with a generalized lens p-end neighborhood. 
Let $\mbv_{\tilde E}$ be the p-end vertex,  and  $\Sigma_{\tilde E}$  
the p-end orbifold of $\tilde E$. 
Suppose that the p-end holonomy group $\bGamma_{\tilde E}$ is 
virtually factorizable.
Then $\bGamma_{\tilde E}$ satisfies the uniform middle-eigenvalue condition
with respect to $\mbv_{\tilde E}$, 
and the following statements hold\,{\rm :}
\begin{enumerate} 
\item[{\rm (i)}] The R-p-end is  totally geodesic. 
$D_i \subset \SI^{n-1}_{\mbv_{\tilde E}}$ 
is projectively diffeomorphic by the projection $\Pi_{\mbv_{\tilde E}}$ 
to totally geodesic convex domain $D'_i$ in a hyperspace in $\SI^n$ {\rm (} resp. in $\RP^n$\/{\rm )}
disjoint from $\mbv_{\tilde E}$. Moreover, 
$\bGamma_{\tilde E}$ is virtually a subgroup of 
$\bR^{l_0-1} \times \prod_{i=1}^{l_0} \bGamma_i$ where 
$\bGamma_i$ acts on $D'_i$ irreducibly, and acts trivially on $D'_j$ for $j \ne i$, 
and $\bR^{l_0-1}$ acts trivially on $D'_j$ for every $j=1, \dots, l_0$. 
\item[{\rm (ii)}] The R-p-end is strictly lens-shaped, and
each $C'_i$ corresponds to a cone $C^*_i = \mbv_{\tilde E}\ast D'_{i}$. 
The R-p-end has a p-end neighborhood equal to the interior of 
\[\{\mbv_{\tilde E}\} \ast D \hbox{ for } D:=  \clo(D'_1) \ast \cdots \ast \clo(D'_{l_0})\]
where the interior of $D$ 
forms the boundary of the p-end neighborhood in $\torb$. 
\item[{\rm (iii)}] The set $S(\mbv_{\tilde E})$ of maximal segments in $\Bd \torb$ from $\mbv_{\tilde E}$ in the closure of a p-end neighborhood of 
$\mbv_{\tilde E}$ is independent of the p-end neighborhood. 
\[\bigcup S(\mbv_{\tilde E}) = \bigcup_{i=1}^{l_{0}} \{\mbv_{\tilde E}\} *\clo(D'_1)*\cdots *\clo(D'_{i-1})* \partial \clo(D'_i) * \clo(D'_{i+1}) * \cdots * \clo(D'_{l_0}).\]
Finally, the statements {\rm (i)-(vi)} of Theorem \ref{pr-thm-lensclass} also hold. 
\end{enumerate}
\end{theorem} 
\begin{proof} 
Again the $\SI^n$-version is enough by Proposition \ref{prelim-prop-closureind}. 
The final statement on {\rm (i)-(vi)} of Theorem \ref{pr-thm-lensclass}
can be proved similiarly to ones given in the proof there. 

(i)  This follows from 
Proposition \ref{prelim-prop-Ben2} 
(see Benoist \cite{Benoist04}). 

As in the proof of Theorem \ref{pr-thm-lensclass}, 
Theorem \ref{pr-thm-equ} implies 
that $\bGamma_{\tilde E}$ satisfies 
the uniform middle-eigenvalue condition.
Proposition \ref{pr-prop-convhull2} implies that the cocompactly-acted lens 
is a strict one. 
Theorem \ref{pr-thm-distanced} implies that the distanced 
$\bGamma_{\tilde E}$-invariant set is contained in a hyperspace $P$
disjoint from $\mbv_{\tilde E}$. 


(ii) By the uniform middle-eigenvalue condition, 
the largest norm of the eigenvalue $\lambda_1(g)$ is strictly 
bigger than $\lambda_{\mbv_{\tilde E}}(g)$. 


By Proposition \ref{prelim-prop-Ben2}, 
$\bGamma_{\tilde E}$ is virtually a subgroup of 
$\bR^{l_0-1} \times \Gamma_1 \times \cdots \times \Gamma_{l_0}$ 
with $\bR^{l_0-1}$ acting as a diagonalizable group, 
and  there are subspaces $\hat \SI_j$, $j=1, \dots, l_0$, 
in $\SI^{n-1}_{\mbv_{\tilde E}}$ where 
the factor groups $\bGamma_1, ...,\bGamma_{l_0}$ act irreducibly
by Benoist \cite{Benoist04}.
Let $\SI_j$, $j=1, \dots, l_0$, be the projective subspaces containing  $\mbv_{\tilde E}$ 
which goes to $\hat S_j$ under $\Pi_{\mbv_{\tilde E}}$. 
Now, 
$\clo(\tilde \Sigma_{\tilde E}) \cap \SI_j$ is a properly convex domain $K_i$
by Benoist \cite{Benoist04}. 
Let $C_i$ denote the union of great segments from $\mbv_{\tilde E}$ with directions in $K_i$ in $\SI_i$ for each $i$.  
The abelian center isomorphic to $\bZ^{l_0-1}$ acts as the identity on the subspace corresponding to $C_{i}$
in the projective space $\SI^{n-1}_{\mbv_{\tilde E}}$. 

We denote by $D'_i := C_i\cap P$. We denote by 
$D = D'_1\ast \cdots \ast D'_{l_0} \subset P$. 
Since the set of directions from $\mbv_{\tilde E}$ to points of $D^o$ is 
exacty $\tilde \Sigma_{\tilde E}$, 
the interior of $\mbv_{\tilde E}\ast D$ is a p-end neighborhood of 
$\tilde E$. This proves (i) and (ii).

From (ii), we see that $\mbv_{\tilde E} \ast \partial D =\Bd \torb \cap \clo(U)$.
Since $\partial D$ equals 
\[\bigcup_{i=1}^{l_0}
 \clo(D'_1)*\cdots *\clo(D'_{i-1})*\partial \clo(D'_i) \ast \clo(D'_{i+1})* \cdots * \clo(D'_{l_0}),\] 
we obtain that 
\[\bigcup_{i=1}^{l_0}
\mbv_{\tilde E} * \clo(D'_1)*\cdots *\clo(D'_{i-1})*\partial \clo(D'_i) \ast \clo(D'_{i+1})* \cdots * \clo(D'_{l_0}) = \Bd \torb \cap \clo(U).\] 
(iii) follows since $S(\mbv_{\tilde E})$ is independent of choice of 
R-p-end neighborhood. 
\hfill \SSn {\parfillskip0pt\par}
\end{proof}

\subsection{Uniqueness of vertices outside the lens}

We will need this later in Chapter \ref{ch-cl}. 

\begin{proposition} \label{cl-prop-section} 
	Suppose that $h:\pi_1(E) \ra \SL_{\pm}(n+1, \bR)$ is a holonomy representation of 
	end fundamental group $\pi_1(E)$ of a properly convex real projective $n$-orbifold.
	Let $h(\pi_1(E))$ act on either a generalized 
	lens cone $\{v\} \ast L$ for vertex $v$ and  a generalized lens $L$
where the action is proper and cocompact, 
	or 	act on a horosphere as a lattice in a cusp group. 
	Then the following hold\/{\rm :}
	\begin{itemize} 
		\item the vertex of the lens cone is determined up to the antipodal map. 
		\item If the lens cone 
		is given an outward direction, then the vertex of any lens cone where 
		$h(\pi_1(E))$ acts as a p-end neighborhood equals the vertex of the 
		lens cone and is uniquely determined. 
		\item The vertex of the horospherical end is uniquely determined. 
	\end{itemize} 
\end{proposition}
\begin{proof} 
	The horospherical case can be understood from 
	the horospherical action on a ball of a Klein model
	where fixed points form a pair of antipodal points. 
	
	Suppose that $h(\pi_1(E))$ acts on a generalized lens cone $\{v\}\ast L$ for a generalized lens $L$ as in the premise and	a vertex $v$. 
By Theorem \ref{pr-thm-equiv}, $h(\pi_1(E))$ satisfies 
	the uniform middle-eigenvalue condition with respect to $v$. 
	Suppose that there exists another point $w$ fixed by 
	$h(\pi_1(E))$ such that	$\{w\} \ast L'$ is a generalized lens cone for another generalized lens $L'$
on which the action is proper and cocompact.
	Let $\overrightarrow{vw}$ denote a vector tangent to $\overline{vw}$ oriented 
	away from $v$. 
	Then  $\overline{vw}$ goes to a point $\llrrparen{\overrightarrow{vw}}$ on $\SI^{n-1}_v$ fixed by
	$h(\pi_1(E))$. 
	Hence, $\pi_1(E)$ acts reducibly on 
	$\SI^{n-1}_v$, and $h(\pi_1(E))$ is virtually factorizable. 
	Thus, $h(\pi_1(E))$ acts on a hyperspace $S$ disjoint from $v$
	by Theorem \ref{pr-thm-redtot}. There is a properly convex domain 
	$D$ in $S$ where $h(\pi_1(E))$ acts properly discontinuously. 
	Also, $\clo(D) = K_1 \ast \cdots \ast K_m$ for properly convex 
	domain $K_j$ where $h(\pi_1(E))$ acts irreducibly
	by Proposition \ref{prelim-prop-Ben2}.
	
	Suppose that $w \in S$. Then $\{w\} $ must be $K_j$ or $K_{j-}$ for some $j$
by the irreducibility of the action on each factor. 
	The uniform middle-eigenvalue condition 
	with respect to $v$ implies that $\pi_1(E)$ does not 
	satisfy the same condition with respect to $w$. 
	Hence, $w \not\in S$. 
	
	Hence $\pi_1(E)$ acts on a domain $\Omega$ equal to 
	the interior of $K:= K_1\ast \cdots \ast K_m$ where $K_j$ is compact and 
	convex and a finite-index subgroup $\Gamma'_E$ of $\pi_1(E)$ acts on 
	each $K_j$ irreducibly. 
	The great segment from $v$ containing $w$ meets $S$ in a point $w'$. 
	There exists a virtual-center diagonalizable group $D$ 
	acting on each $K_j$ as the identity by 
	Proposition \ref{prelim-prop-Ben2}
	(more precisely, Proposition 4.4 of \cite{Benoist03}). 
	Hence $\{w'\}$ must be one of $K_k$ or its anitpode $K_{k-}$ 
	since otherwise we can find an element of $D$ not fixing $w'$. 
	
	Since the action is cocompact on $\Omega$, 
	there must be an element of $D$ acting with the largest norm eigenvalue 
	on $K_k$.
	
	Since $h(\pi_1(E))$ acts on $\{w\} \ast L'$ with a compact set $L'$ 
	disjoint from $w$, we construct a tube domain $T$ and 
	$L'\cap T$ gives us a cocompactly-acted lens in the tube.  
	Hence, by Theorem \ref{pr-thm-equiv}
	$\lambda_w(g)$ satisfies the uniform middle-eigenvalue condition. 
	We choose $g \in D$ with a largest norm eigenvalue at $K_k$. 
	Since $v, w, w'$ are distinct points in a properly convex segment 
	and are fixed points of $g$, it follows that 
	$\lambda_1(g) = \lambda_w(g) =\lambda_v(g)$. 
	This contradicts the uniform middle-eigenvalue condition for $v$ under 
	$\pi_1(E)$. 
	Thus, we obtain $v=w$.
%
\end{proof} 

By duality, we obtain the following: 

\begin{proposition} \label{cl-prop-section2} 
	Suppose that $h:\pi_1(E) \ra \SLnp$ is a representation 
	\begin{itemize} 
	\item acting properly, cocompactly, 
and effectively
on a lens neighborhood of a totally geodesic $(n-1)$-dimensional domain $\Omega$, 
	and $\Omega/h(\pi_1(E))$ is 
	a compact orbifold or 	
	\item acting on a horosphere
	as a cocompact cusp group. 
	\end{itemize} 
	Then $h(\pi_1(E))$ uniquely determines the hyperspace $P$ 
	with one of the following properties\/{\rm :} 
	\begin{itemize} 
		\item  $P$ meets a lens domain $L'$ with the property that 
		$(L' \cap P)/h(\pi_1(E))$ is a compact orbifold
		with $L'\cap P = L^{\prime o}\cap P$. 
		\item $P$ is tangent to 
		the $h(\pi_1(E))$-invariant horosphere at the cusp point of the horosphere. 
	\end{itemize} 
\end{proposition}
\begin{proof} 
	The duality will prove this by Proposition \ref{pr-prop-dualend}
	and Corollary \ref{pr-cor-duallens} in the next section. 
	The vertex and the hyperspace exchange roles. 
\end{proof}

\section{Duality and lens-shaped T-ends}\label{pr-sec-dualT}

We first discuss the duality map. We show a lens cone p-end neighborhood  of an R-p-end 
is dual to a lens p-end neighborhood of a T-p-end. 
Using this duality, we prove Theorem \ref{pr-thm-equ2} dual to Theorem \ref{pr-thm-equ}, i.e., 
Theorem \ref{pr-thm-secondmain}. 

\subsection{Duality map.} 

The Vinberg duality diffeomorphism induces a one-to-one correspondence between p-ends of $\torb$ and $\torb^*$ 
by considering the dual relationship $\bGamma_{\tilde E}$ and $\bGamma^*_{\tilde E'}$ for each pair of 
p-ends $\tilde E$ and $\tilde E'$ with dual p-end holonomy groups. (See Section \ref{prelim-sec-duality}.)

Given a properly convex domain $\Omega$ in $\SI^n$ (resp. $\RP^n$), 
we recall the augmented boundary of $\Omega$
\begin{align} 
\Bd^{\Ag} \Omega  &:= \{ (x, H)| x \in \Bd \Omega, x\in H, \nonumber \\ 
 & H \hbox{ is an oriented sharply supporting hyperspace of } \Omega \}  \subset \SI^{n}\times \SI^{n \ast}.
\end{align}
This is a closed subspace. 
Each $x \in \Bd \Omega$ is contained in at least one sharply supporting hyperspace
oriented towards $\Omega$, 
and an element of $\SI^n$ represents an oriented hyperspace in $\SI^{n \ast}$. 
\index{boundary!augmented} 

We recall a duality map.
\begin{equation} \label{pr-eqn-dualmap}  
{\mathcal{D}}^\Ag_{\Omega}: \Bd^{\Ag} \Omega \rightarrow \Bd^{\Ag} \Omega^* 
\end{equation}
given by sending $(x, H)$ to $(H, x)$ for each $(x, H) \in \Bd^{\Ag} \Omega$. 
This is a diffeomorphism since $\mathcal{D}^\Ag_{\Omega}$ has an inverse given by 
switching factors by Proposition \ref{prelim-prop-duality} (iii).
\index{duality map} 
\index{do@${\mathcal{D}}^\Ag_{\Omega(\cdot)}$}

\begin{figure}
\centering
\includegraphics[trim = 10mm 80mm 1mm 10mm, clip, width=10cm, height=5cm]{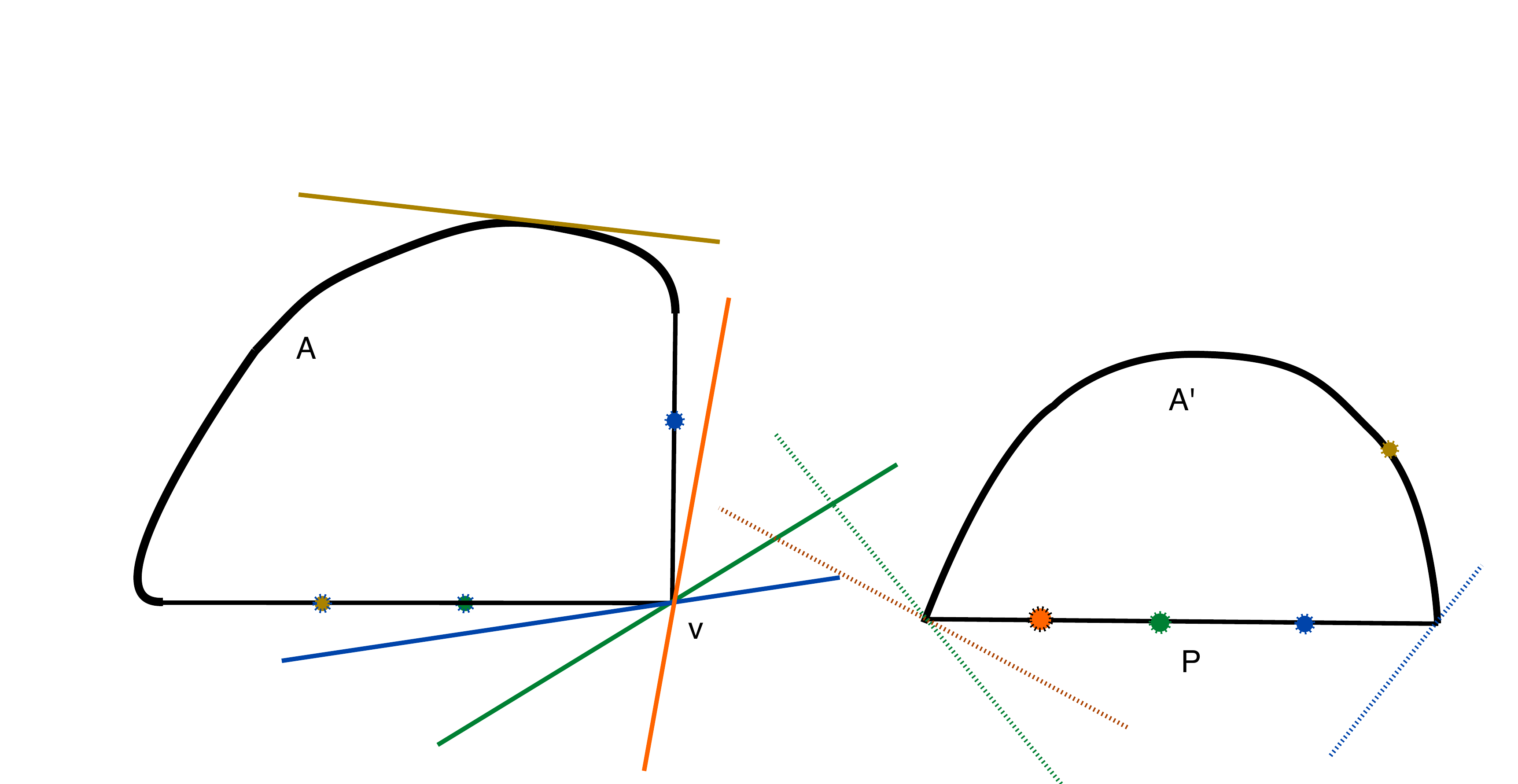}

\caption{The figure for Corollary \ref{pr-cor-duallens}. Lines passing $v$ in the left figure correspond to points on $P$ in the left. The line passing a point of $A$ corresponds 
to a point on $A'$. Lines in the right figure correspond to points in the right figure. }

\label{pr-fig-duallens}
\end{figure}

%

We need the corollary about the duality of lens cone and lens neighborhoods. 
Recall that given a properly convex domain $D$ in $\SI^{n}$ or $\RP^{n}$, 
the dual domain is the closure of the open set given by the collection of (oriented) hyperspaces in $\SI^{n}$ or $\RP^{n}$ not meeting $\clo(D)$. 
Let $\Omega$ be a properly convex domain. 
We recall the duality from Section \ref{prelim-sub-duality}
with the projection map 
\[  \Pi^{\Ag}_{\Omega}: \Bd^\Ag \Omega \ra \Bd \Omega \]
sending each pair $(x, h)$ of a point $x$, $x\in \Bd \Omega$, 
and sharply supporting hyperspace $h$ at $x$. 

\begin{corollary} \label{pr-cor-duallens} 
The following hold in $\SI^n$\/{\rm :}  
\begin{itemize} 
\item Let $L$ be a lens and $\{\mbv\} \not\in \clo(L)$ such that $\mbv \ast L$ is a properly convex lens cone. 
Suppose that the smooth, strictly convex 
boundary component $A$ of $L$ is tangent to a segment from $\mbv$ at each point of 
$\partial \clo(A)$
and $\{\mbv\}\ast L = \{\mbv\} \ast A$. 
Then the following hold\/{\em :}
\begin{itemize} 
\item the dual domain of $\clo(\{\mbv\}\ast L)$ is the closure of a component $L_{1}$ of $L' - P$ where $L'$ is a lens and $P$ is a hyperspace meeting $L^{\prime o}$. 
\item $A$ corresponds to a hypersurface $A' \subset \Bd L'$ under the duality \eqref{pr-eqn-dualmap}. 
\item $A' \cup D$ is the boundary of $L_1$ for a totally geodesic 
properly convex $(n-1)$-dimensional compact domain $D$ dual to $R_{\mbv}(\{\mbv\} \ast L)$ where $D$ is given by
 $\Pi_{\{\mbv\}\ast L}^\Ag\circ \mathcal{D}_{\{\mbv\}\ast L}((\Pi^{\Ag})^{-1}_{\{\mbv\} \ast L}(\{\mbv\}))$.
\end{itemize}  
\item Conversely, we are given a lens $L'$ and a hyperspace $P$ meeting $L^{\prime o}$ 
but not meeting the boundary of $L'$. 
Let $L_{1}$ be a component of $L'-P$ with smooth, strictly convex boundary 
$\partial L_{1}$ such that $\partial \clo(\partial L_{1}) \subset P$. 
Here, we assume $\partial L_1$ is an open hypersurface. 
Then the following hold\/{\em :}
\begin{itemize} 
\item The dual of the closure of a component $L_{1}$ of $L' -P$ is 
the closure of $\mbv \ast L$ for a lens $L$ and $\mbv \not\in L$ such that 
$\mbv\ast L$ is a properly convex lens cone. 
Here, $\{\mbv\} = \Pi^\Ag_{\clo(L_1)}\circ \mathcal{D}_{\clo(L_1)}(P)$. 
\item The outer boundary component $A$ of $L$ is tangent to a segment from $\mbv$ at each point of $\partial \clo(A)$. 
\end{itemize} 
\item In the above, the vertex denoted by $v$ corresponds to a hyperspace denoted by $P$ uniquely. 
\end{itemize}
\end{corollary} 
\begin{proof} 
	In the proof, all hyperspaces are oriented such that $L^o$ is in its interior direction. 
Let $A$ denote the boundary component of $L$ such that $\{\mbv\} \ast L = \{\mbv\} \ast A$. 
We determine the dual domain $(\clo(\{\mbv\} \ast L))^\ast$ of
 $\{\mbv\} \ast L$ by finding the boundary of $D$ using the duality map $\mathcal{D}_{\{\mbv\} \ast L}$. 
The set of hyperspaces sharply supporting $\clo(\{\mbv\}\ast L)$ at $\mbv$ forms a properly totally geodesic domain $D$ in $\SI^{n\ast}$
contained in a hyperspace $P$ dual to $\mbv$ by Lemma \ref{pr-lem-predual}. 
Also, the set of hyperspaces sharply supporting $\clo(\{\mbv\} \ast L)$ at points of $A$ goes to the strictly convex hypersurface $A'$ in $\Bd (\mbv \ast L)^\ast$
by Lemma \ref{pr-lem-predual} since $\mathcal{D}_{\{\mbv\} \ast L}$ is a diffeomorphism. 
(See Remark \ref{prelim-rem-duallens} and Figure \ref{pr-fig-duallens}.)
The subspace
$S:=\Bd (\{\mbv\} \ast A) - A$ is a union of segments from $\mbv$. 
The sharply supporting hyperspaces containing these segments 
go to points in $\partial D$. 
Each point of $\clo(A') - A'$ is a limit of a sequence $\{p_{i}\}$ of points of $A'$, corresponding to a sequence of sharply
supporting hyperspheres $\{h_{i}\}$ to $A$. 
The tangency condition of $A$ at $\partial \clo(A)$ implies that the limit hypersphere contains the segment in $S$ from $\mbv$. 
We obtain that $\clo(A') - A'$ equals the set of hyperspheres containing the segments in $S$ from $\mbv$, and 
they go to points of $\partial D$ with $\partial \clo(A') = \partial D$. 
We conclude 
\[\Pi_{\{\mbv\}\ast L}^\Ag\circ\mathcal{D}_{\{\mbv\} \ast L}(\Bd  (\{\mbv\} \ast L)) = A' \cup D.\] 

Let $P$ be the unique hyperspace containing $D$. 
Then each point of $\partial \clo(A)$ is mapped to a sharply supporting hyperspace 
at a point of $\partial \clo(A')$
distinct from $P$. 
Let $L^{\ast}$ denote the dual domain of $\clo(L)$. 
Since $\clo(L) \subset \clo(\{\mbv\}\ast L)$, we obtain 
$(\clo(\{\mbv\}\ast L))^\ast \subset (\clo(L))^{\ast}$ by  \eqref{prelim-eqn-dualinc}.
Since $A \subset \Bd L$, we obtain 
\[\Pi_{\{\mbv\}\ast L}^\Ag\circ\mathcal{D}_{\{\mbv\}\ast L}(\Bd (\{\mbv\}\ast L))  \subset A' \cup P, \hbox{ and } A' \subset \Bd L^{\ast}.\]
Proposition \ref{prelim-prop-duality} implies that 
$((\{\mbv\}\ast L)^o)^\ast $ is a component $L_1$ 
of $(L^o)^{\ast} - P$ since the first 
domain can have boundary points in $A' \cup P$ only and cannot have points outside 
the component. Hence, the dual of $\clo(\{\mbv\} \ast L)$ is 
$\clo(L_1)$. 
Moreover,  $A' \subset \Bd L_{1}$ since $A'$ is a strictly convex hypersurface 
with boundary in $P$. 

The second item is proved similarly to the first. Now hyperspaces are oriented 
such that $L_1^o$ is in its interior. 
Then $\partial L_{1}$ is mapped to a hypersurface $A$ in the
boundary of the dual domain
$L_1^\ast$ of $L_{1}$ under $\mathcal{D}_{L_1}$. 
Again $A$ is a smooth, strictly convex boundary component. 
Since $\partial \clo(L_{1}) \subset P$ and $L_{1}$ is a component of $L'-P$, we obtain
$\Bd L_{1} - \partial L_{1} = \clo(L_{1}) \cap P$. 
This is a totally geodesic properly convex domain $D$ in $P$.


Suppose that $l \subset P$ is an $n-2$-dimensional space disjoint from $L_1^o$. 
Then a space of oriented hyperspaces containing $l$, and bounding an open hemisphere 
containing $L_1^o$
forms a parameter dual to a convex projective geodesic in $\SI^{n\ast}$. 
An {\em $L_{1}$-pencil} $P_{t}$ with ends $P_{0}, P_{1}$ is a parameter satisfying 
\[P_{t} \cap P = P_{0 }\cap P,  P_{t}\cap L_{1}^{o} = \emp \hbox{ for  all } t \in [0, 1]\]
where $P_t$ is oriented such that it bounds an open hemisphere containing $L_1^o$. 

There is a one-to-one correspondence 
\[ \{P'| P' \hbox{ is an oriented hyperspace that supports } L_{1} \hbox{ at points of } \partial D \} 
\leftrightarrow \{\mbv\} \ast \partial \clo(A): \]

Every supporting hyperspace $P'$ to $L_{1}$ at points of $\partial D$ is contained in  an $L_{1}$-parameter $P_{t}$ with $P_{0} = P', P_{1} = P$.
Now, $\mbv$ is the dual to $P$ in $\SI^{n\ast}$. 
Each path $P_{t}$ is a segment in $\SI^{n\ast}$ with an endpoint $\mbv$. 

Under the duality map $\Pi^\Ag_{\clo(L_1)} \circ \mathcal{D}_{\clo(L_1)}$, the image of $\Bd L_1$ is a union of $A$ and these segments. 
Given any hyperspace $P'$ disjoint from $L_{1}^{o}$, we find a one-parameter family of hyperspaces containing $P' \cap P$. 
Thus, we find an $L_1$-pencil family $P_{t}$ 
with $P_{0} = P', P_{1} = P$.
We can extend the $L_1$-pencil such that the ending hyperspace $P''$ of the $L_1$-pencil 
meets  $\partial L_1$ tangentially, and $P'' \cap P$ is a 
supporting hyperspace of $D$ in $P$. 
Since the hyperspaces are disjoint from $L_{1}$, the segment is in $L_1^\ast$. 
Since $L_1$ is a properly convex domain,  
we can deduce that $(\clo(L_1))^\ast$ is the closure of the cone  $\{\mbv\}\ast A$. 

Let $L''$ be the dual domain of $\clo(L')$.  
Since $L_1 \subset L'$, we obtain $L'' \subset (\clo(L_1))^\ast$ by \eqref{prelim-eqn-dualinc}. 
Since $\partial L_{1}\subset \Bd L'$, we obtain $A \subset \Bd L''$ by the duality map $\mathcal{D}_{\clo(L_1)}$.  
We obtain that $L^{\prime \prime} \cup A \subset \clo(\{\mbv\}\ast A)$.  

Let $B$ be the image of the other boundary component $B'$ of $L'$ under $\mathcal{D}_{L'}$. 
We take a sharply supporting hyperspace $P_{y}$ at $y\in B'$. 
Then $P_y \cap P$ is disjoint from $\clo(D)$ by the strict convexity of $B'$.
We find an $L_1$-pencil $P_{t}$ of hyperspaces containing $P_{y} \cap P$
with $ P_{0}=P_{y}, P_{1}=P$.
This $L_1$-pencil goes into the segment from $\mbv$ to a point of $B$ under the duality.
We can extend the $L_1$-pencil such that the ending hyperspace meets $\partial L_1$ tangentially.  The dual pencil is a segment from $\mbv$ to a point of $A$. 
Thus, each segment from $\mbv$ to a point of $A$ meets $B$. 
Thus, $L^{\prime \prime o} \cup A \cup B$ is a lens of the lens cone $ \{\mbv\}\ast A$.
This completes the proof. 
\end{proof}

\subsection{The duality of T-ends and properly convex R-ends.}  \label{pr-sub-dualend}

%



%
Let $\Omega$ be the properly convex domain projectively covering $\orb$. 
For a T-end $E$, the totally geodesic ideal boundary component $S_{E}$ of $E$ 
is covered by a properly convex open domain in $\Bd \Omega$ corresponding 
to a T-p-end $\tilde E$.  We denote it by $\tilde S_{\tilde E}$. 


Recall R-end structure from Section \ref{intro-sub-R-ends}.

\begin{lemma} \label{pr-lem-endnhbd}
	Let $\orb$ be a properly convex real projective orbifold with ends
	where 
	$\orb = \torb/\Gamma$ for a properly convex domain $\torb \subset \SI^n$ and
	a discrete projective group $\Gamma$. Let $p_{\orb}: \torb \ra \orb$ 
	denote the projective covering map.
	Suppose that a p-end fundamental group $\Gamma_{\tilde E}$ for 
	a p-end $\tilde E$ acts on a simply connected hypersurface $\tilde \Sigma$ in $\torb^o$. 
	Then a component $U'$ of $\torb - \tilde \Sigma$ is a p-end neighborhood of 
	$\tilde E$. Furthermore, the following hold\/{\rm :}
	\begin{itemize}
	\item Suppose that $U'$ is a horoball such that $\Bd U' \cap \torb 
	= \Bd U' -\{p\}= \Sigma$ for the common fixed point $p$ of $\Gamma_{\tilde E}$. 
	Then $p_{\orb}(U')$ can be given the structure of 
	a horospherical R-end neighborhood of $\tilde E$ 
	and $p=\mbv_{\tilde E}$. 
	\item Suppose that
	 $U'$ equals $U_L:= L \ast \{p\} - L$ for a lens $L$ where
	\begin{itemize}
		\item 
	$h(\pi_1(\tilde E))$ acts properly and cocompactly on $U'$, 
	\item we have 
	a lens cone $L\ast \{p\}$ for a common fixed point $p$ of 
	$\Gamma_{\tilde E}$, and	
	\item $\Bd U_L \cap \torb = \tilde \Sigma$.  
	\end{itemize} 
	Then $U_L$ can be given the structure of
	a concave R-p-end neighborhood of $\tilde E$, 
and $p = \mbv_{\tilde E}$. 
\item Suppose that $U'$ equals a component $L_1$ of $L - P$ for a lens $L$ 
where 
\begin{itemize}
\item $h(\pi_1(\tilde E))$ acts properly and cocompactly on $L$, 
\item $h(\pi_1(\tilde E))$ acts on the hyperspace $P$,
\item  $L_1 \subset \torb$, and 
\item $\Bd L_1 \cap \torb =\tilde  \Sigma$. 
\end{itemize} 
 Then $L_1^o$ can be given the structure of a T-p-end neighborhood of 
$\tilde E$ and, moreover, $\tilde \Sigma_{\tilde E}$ can be identified with $L \cap P$.  
\end{itemize} 	
Moreover, the corresponding end completions give us a compact smooth 
orbifold $\bar \orb$ whose interior is $\orb$.  
	
\end{lemma}
\begin{proof} 
	It is sufficient to prove for the case when the orbifold is an orientable
manifold since we can take a finite quotient by 
	Theorem \ref{prelim-thm-vgood}. 
	Let $U$ be a proper p-end neighborhood of $\tilde E$.
Recall that $\Gamma_{\tilde E}$ acts on $\tilde \Sigma_{\tilde E}$ with 
a closed $(n-1)$-manifold as a quotient. 
	We take $\torb/\Gamma_{\tilde E}$, which is homotopy eqiuivalent to  
	$A:= B \times \bR$, 
	where $B:= (\Bd U \cap \torb)/\Gamma_{\tilde E}$ is a closed submanifold of codimension-one.
	We can take an exiting sequence $U_i$ of p-end neighborhoods in $U$.  
	Then $C:=\tilde \Sigma/\Gamma_{\tilde E}$ is a closed submanifold homotopy equivalent to $B \times \bR$.
Hence, one $U_B$ of the two components $\torb/\Gamma_{\tilde E} - B$ 
	contains $U_i$ for sufficiently large $i$.  
	Hence, a component of the inverse image of $U_B$ which is also a component of 
	$\torb - \Sigma$ is a p-end neighborhood, perhaps not a proper one. 
	
		In the first case, we have a horoball $U'$ inside $\torb^o$ and in $H$ 
since the sharply supporting hyperspaces at the vertex of $H$
	must coincide by the invariance under $h(\pi_1(\tilde E))$
	by a limiting argument. 
	By above, one of the two component of $\torb^o - \Bd U'$
	is a p-end neighborhood. 
	One cannot put the outside component into $U'$ by an element of 
	$h(\pi_1(\orb))$. 
	Hence, $U$ is a horospherical p-end neighborhood of $\tilde E$. 
	
	For the second item, Lemma \ref{intro-lem-actionRend} implies that $U_L$ is 
	foliated by radial lines mapping to properly embedded 
	lines in $\orb$. 
	Thus, $U_L$ is the p-end neighborhood of $\tilde E$, and   			 
	$\tilde E$ has a radially foliated p-end neighborhood
	with $\mbv_{\tilde E}$ as the p-end vertex.  
	
	In the third case,	let $D$ denote $\clo(L_1)\cap P$, a properly convex domain. 
By premise, 
	$D^o/h(\pi_1(\tilde E))$ is a closed orbifold of codimension-one.
	Then one of the two components of $\torb - \Bd L_1$ is a p-end neighborhood of $\tilde E$ 
as shown above. 
	$D^o$ is totally geodesic and $\Bd L_1 \cap \torb$ is not. 
	Hence, $L_1^o$ is a p-end neighborhood of $\tilde E$. 
	By premise, $h(\pi_1(\tilde E))$ acts properly on 
	$L_1 \cup D^o$. 
	The $T$-end structure is given by $(D^o \cup L_1)/h(\pi_1(\tilde E))$, 
	which is the completion of the end neighborhood $p(L_1^o)$, 
	projectively diffeomorphic to $L_1^o/h(\pi_1(\tilde E))$. 
	(See Section \ref{intro-sub-totgeo}.)
\end{proof} 

%
%
%


\begin{proposition}\label{pr-prop-dualend2}
	Let $\orb$ be a properly convex real projective $n$-orbifold. 
	The following conditions are equivalent\,{\rm :} 
	\begin{enumerate} 
		\item[{\rm (i)}] A properly convex R-end $E$ of $\orb$ satisfies the uniform middle-eigenvalue condition.
		\item[{\rm (ii)}] The corresponding T-end $E^*$ of $\orb^*$ satisfies this condition where the correspondence of 	the vertex of the p-end $\tilde E$ of $E$ to the hyperspace of the p-end 
		$\tilde E^\ast$, which is given the unique hyperspace containing
		$\Pi^\Ag_{\clo(\torb)}\circ\mathcal{D}_{\torb}((\Pi^{\Ag})^{-1}_{\clo(\torb)}(\mbv_{\tilde E}))$. 
	\end{enumerate} 
\end{proposition}
\begin{proof}
	The items (i) and (ii) are equivalent by considering \eqref{pr-eqn-umec} and 
	\eqref{pr-eqn-umecD}. 
	Proposition \ref{prelim-prop-closureind} implies the $\RP^n$-version.
	\hfill \SnT  {\parfillskip0pt\par}
\end{proof}


We now prove the dual to Theorem \ref{pr-thm-equ}. For this, we do not need the triangle condition or the
reducibility of the end. 


\begin{theorem}\label{pr-thm-equ2}
	Let $\orb$ be a properly convex real projective $n$-orbifold. 
	Let $\tilde S_{\tilde E}$ be a totally geodesic ideal boundary component of a T-p-end $\tilde E$ of $\torb$. 
	Then the following conditions are equivalent\,{\rm :} 
	\begin{enumerate} 
		\item[{\rm (i)}] The end holonomy group of 
		$\tilde E$ satisfies the uniform middle-eigenvalue condition
		with respect to the T-p-end structure of $\tilde E$. 
		\item[{\rm (ii)}] $\tilde S_{\tilde E}$ has a lens neighborhood in an ambient open manifold containing $\torb$ with cocompact action of $\pi_1(\tilde E)$, 
		and hence $\tilde E$ has a lens-shaped p-end neighborhood in $\torb$. 
	\end{enumerate} 
\end{theorem}
\begin{proof}
	We prove the result for the $\SI^n$-version.  
	Assuming (i), we can deduce the existence of a lens neighborhood from Theorem \ref{du-thm-lensn} and Lemma \ref{pr-lem-endnhbd}. 
	
	Assuming (ii), we obtain a totally geodesic $(n-1)$-dimensional properly convex domain $\tilde S_{\tilde E}$ in
	a subspace $\SI^{n-1}$ on which $\bGamma_{\tilde E}$ acts. 
	Let $U$ be a lens-neighborhood of it on which $\bGamma_{\tilde E}$ acts. 
	Since $U$ is a neighborhood, the sharply supporting hemisphere at each point of 
	$\clo(\tilde S_{\tilde E})- \tilde S_{\tilde E}$ is now transverse to $\SI^{n-1}$. 
	Let $P$ be the hyperspace containing $\tilde S_{\tilde E}$, and 
	let $U_{1}$ be the component of $U - P$. Then the dual $U_{1}^{\ast}$ is 
	a lens cone by the second part of Corollary \ref{pr-cor-duallens} where 
	$P$ corresponds to a vertex of the lens cone.
	The dual $U^*$ of $U$  is aa lens contained in a lens cone $U_{1}^{\ast}$ where 
	$\bGamma_E$ acts on $U^{\ast}$. 
	We apply the part (i) $\Rightarrow$ (ii) of Theorem \ref{pr-thm-equ}.
By Proposition \ref{pr-prop-dualend2}, the result follows.
	Proposition \ref{prelim-prop-closureind} implies the result for $\RP^n$. 
\hfill	\SnT {\parfillskip0pt\par}
\end{proof}

\begin{proposition}\label{pr-prop-dualend} 
Let $\orb$ be a properly convex real projective $n$-orbifold with R-ends or T-ends with universal covering domain $\Omega$. Let $\orb^\ast$ be the dual orbifolds with a universal covering domain $\torb^\ast$. 
Then $\orb^*$ can be given R-end and T-end structures to the ends 
such that the following hold\/{\em :}
\begin{itemize} 
\item there exists a one-to-one correspondence $\mathcal{C}$ between the set of p-ends of $\torb$ and the set of p-ends of $\torb^*$ which sends a p-end neighborhood to 
a p-end neighborhood using the Vinberg diffeomorphism of Theorem \ref{prelim-thm-dualdiff}. 
\item $\mathcal{C}$ restricts to such a one  between the subset of horospherical p-ends of $\torb$ and the subset of horospherical ones of $\torb^*$.
Also, the augmented Vinberg duality homeomorphism $\bar{\mathcal{D}}^\Ag$ sends the p-end vertex to the p-end vertex of the dual p-end. 
\item $\mathcal{C}$ restricts to a one-to-one correspondence   between 
the subset of all generalized lens-shaped R-ends of $\orb$ and the subset of all
lens-shaped T-ends of $\orb^*$. 
Moreover, 
$\tilde \Sigma_{\tilde E}$ of an R-p-end is projectively dual to the ideal boundary component $\tilde S_{\tilde E^*}$ 
for the corresponding dual T-p-end $\tilde E^*$ of $\tilde E$. 
Moreover, $\mathcal{D}^\Ag_{\torb}$ gives a one to one correspondence 
between $\Pi^{\Ag -1}_{\torb}(\bigcup S(\mbv_{\tilde E}))$ in $\Bd^\Ag \torb$ 
to $\Pi^{\Ag -1}_{\torb^\ast}(\clo(\tilde S_{\tilde E^\ast}))$  of $\Bd^\Ag \torb^\ast$. 
\item $\mathcal{C}$ restricts to a one-to-one correspondence   between the set of lens-shaped
T-p-ends of $\torb$ with the set of p-ends of generalized lens-shaped R-p-ends of $\torb^*$.
The ideal boundary component $\tilde S_{\tilde E}$ for a T-p-end $\tilde E$ 
is projectively diffeomorphic to the properly convex open domain dual to the domain
$\tilde \Sigma_{\tilde E^*}$ for the corresponding R-p-end $\tilde E^*$ of $\tilde E$. 
Also, $\mathcal{D}^\Ag_{\torb}$ gives one to one correspondence 
between $\Pi^{\Ag -1}_{\torb}(\clo(\tilde S_{\tilde E}))$ in $\Bd^\Ag \torb$ 
to $\Pi^{\Ag-1}_{\torb^\ast}(\bigcup S(\tilde E^\ast))$ of $\Bd^\Ag \torb^\ast$. 
\end{itemize} 
\end{proposition}
\begin{proof}
We prove for the result the $\SI^n$-version first. 
Let $\torb$ be the universal cover of $\orb$. Let $\torb^*$ be the dual domain. 
The first item follows from the fact that 
this diffeomorphism sends p-end neighborhoods to p-end neighborhoods. 

Let $\tilde E$ be a horospherical R-p-end with $x$ as the end vertex. 
Since there is a subgroup $\Gamma_{\tilde E}$ 
of a cusp group acting on $\clo(\torb)$ with
a unique fixed point, 
the intersection of the unique sharply supporting hyperspace $P$ with $\clo(\torb)$ at $x$ is a singleton $\{x\}$. 
(See Theorem \ref{ce-thm-affinehoro}.)
The dual subgroup $\Gamma_{\tilde E}^\ast$
is also a cusp group and acts on $\clo(\torb^*)$ with the fixed point $P^\ast$ dual to $P$.
So the corresponding $\torb^*$ has the tangent hyperspace $x^*$ dual to $x$ 
meeting $\clo(\torb^\ast)$ at the unique point $P^\ast$. 
There is a horosphere $S$ where the end fundamental group $\Gamma_{\tilde E}^\ast$ acts on. 
By Lemma \ref{pr-lem-endnhbd}, $S$ bounds a horospherical $p$-end neighborhood of $\tilde E$.
Hence, $x^*$ is the vertex of a horospherical end. 
$\mathcal{D}_{\torb}^\Ag(x) = x^\ast$ since $\Gamma_{\tilde E} \ra \Gamma_{\tilde E}^\ast$ 
and these are unique fixed points. 

An R-p-end $\tilde E$ of $\torb$ has a p-end vertex $\mbv_{\tilde E}$. 
$\tilde \Sigma_{\tilde E}$ is a properly convex domain in $\SI^{n-1}_{\mbv_{\tilde E}}$. 
The space of sharply supporting hyperspaces of $\torb$ at 
$\mbv_{\tilde E}$ forms a properly convex domain of dimension $n-1$ 
since they correspond to hyperspaces in $\SI^{n-1}_{\mbv_{\tilde E}}$
not intersecting $\tilde \Sigma_{\tilde E}$. 
Under the duality map ${\mathcal{D}}^\Ag_{\torb}$ in Proposition \ref{prelim-prop-duality},  
$(\mbv_{\tilde E}, P)$ for a sharply supporting hyperspace $P$ is sent 
to $(P^{\ast}, \mbv_{\tilde E}^{\ast})$ for a point $P^\ast$ and 
a hyperspace $\mbv_{\tilde E}^\ast$. 
Lemma \ref{pr-lem-predual} 
shows that 
$P^{\ast}$ is a point in a properly convex $(n-1)$-dimensional domain 
$D:= \Bd \torb^{\ast} \cap P$ for $P = \mbv_{\tilde E}^{\ast}$, a hyperspace. 

Corollary \ref{pr-cor-duallens} implies the fact about
$\mathcal{D}_{\torb}^\Ag$. 


Since $D$ is a properly convex domain with a Hilbert metric, 
$\pi_1(\tilde E)$ acts properly on $D^o$. 
The $n$-orbifold $(\torb\cup D^o)/\pi_1(\tilde E)$ 
has closed-orbifold boundary $D^o/\pi_1(\tilde E)$. 
There is a Riemannian metric on the $n$-orbifold such that 
$D^o/\pi_1(\tilde E)$ is totally geodesic. 
Using the exponential map, we obtain a tubular neighborhood of 
$D^o/\pi_1(\tilde E)$. Hence, $\torb$ has a p-end neighborhood that corresponds 
to $\pi_1(\tilde E)$ and contains $D^o$ in the boundary. 
The dual group $\bGamma_{\tilde E}^\ast$ satisfies the uniform middle 
eigenvalue condition since $\bGamma_{\tilde E}$ satisfies the condition. 
By Theorem \ref{du-thm-lensn} and Lemma \ref{pr-lem-endnhbd}, 
we can find a p-end neighborhood $U$
in $\torb^\ast$ bounded by a strictly convex 
hypersurface  $\Bd U \cap \torb$, 
where $\clo(\Bd U\cap \torb)-(\Bd U \cap \torb) \subset \partial D$.

By Lemma \ref{pr-lem-endnhbd},  
$\tilde S_{\tilde E^\ast} \subset \Bd \Omega^\ast$, and 
$\tilde E^{\ast}$ is a totally geodesic end with 
$\tilde S_{\tilde E^{\ast}}$ dual to $\tilde \Sigma_{\tilde E}$. 
This proves the third item.

The fourth item follows similarly. 
Take a T-p-end $\tilde E$. We take the ideal p-end boundary
$\tilde \Sigma_{\tilde E}$. The map $\mathcal{D}^\Ag_{\torb}$ sends
$P$ to a singleton $P^\ast$ in $\Bd \torb^\ast$, and 
points of $\clo(\tilde S_{\tilde E}))$ to 
hyperspaces supporting $\torb^\ast$ at $P^\ast$. 
Since $\bGamma_{\tilde E}^\ast$ satisfies the uniform middle-eigenvalue 
condition with respect to $P^\ast$, Theorem \ref{pr-thm-equiv} shows that 
there exists a lens cone on which $\bGamma_{\tilde E}^\ast$ acts. 
Also, 
$\bGamma_{\tilde E}$ acts on a tube domain $\mathcal{T}_{P^\ast}(D^o)$ for
a properly convex domain $D^o$. Then $\clo(\torb^\ast) \cap 
\Bd \mathcal{T}_{P^\ast}(D^o)$ is a $\bGamma_{\tilde E}^\ast$-invariant 
closed set. Also, this set is the image under  $\mathcal{D}^\Ag_{\torb}$
of all hyperspaces supporting $\torb$ at points of $\clo(\tilde \Sigma_{\tilde E})$
by Corollary \ref{pr-cor-duallens}. 
Hence, $R_{P^\ast}(\torb^\ast) = D^o$ by convexity. 
Since $D^o$ is properly convex, $\bGamma_{\tilde E}^\ast$ acts properly on it. 
By Lemma \ref{intro-lem-actionRend}, $P^\ast$ is a p-end vertex of 
a p-end neighborhood. There is a lens $L$ such that $U_L:= P^\ast \ast L - L$ is 
$\bGamma$-invariant. There is a boundary component $\partial_- L$ of 
this $U_L$ in $\torb^\ast$. 
By Lemma \ref{pr-lem-endnhbd}, this implies that $U_L$ is a p-end neighborhood 
corresponding to $\bGamma_{\tilde E}$. 
Corollary \ref{pr-cor-duallens} implies the fact about
$\mathcal{D}_{\torb}^\Ag$. 
%

The proof for $\RP^n$-version follows from
Proposition \ref{prelim-prop-closureind}. 
\hfill  \SnT {\parfillskip0pt\par}
\end{proof} 





\begin{remark}\label{pr-rem-comp} 
	We also remark that  the map induced from the limit points 
	of p-end neighborhoods of $\Omega$ to those of $\Omega^{*}$ by
	$\bar{\mathcal{D}}_{\Omega}^\Ag$ is compatible with the Vinberg diffeomorphism
	by the continuity part of Theorem \ref{prelim-thm-AugVinberg}. 
	That is, the limit points of $\Bd^\Ag \Omega$ of a p-end neighborhood of 
	a p-end $\tilde E$ 
	go to the limit points of $\Bd^\Ag \Omega^\ast$ 
	a p-end neighborhood of a dual p-end $\tilde E^\ast$ of $\tilde E$
	by $\bar{\mathcal{D}}_{\Omega}^\Ag$. 
\end{remark}

$\mathcal{C}$ restricts to a correspondence between the lens-shaped R-ends and
lens-shaped T-ends. See Corollary \ref{pr-cor-duallens2}  for details.


Theorems \ref{pr-thm-equ} and \ref{pr-thm-equ2} and Propositions \ref{pr-prop-dualend} and \ref{pr-prop-dualend2} imply 
\begin{corollary}\label{pr-cor-duallens2} 
Let $\orb$ be a properly convex real projective $n$-orbifold, and 
let $\orb^{\ast}$ be its dual orbifold. 
Then we can give the structure of R-ends and T-ends to ends of 
$\orb$ and $\orb^\ast$ such that 
dual end correspondence $\mathcal{C}$ restricts to a correspondence between the generalized lens-shaped R-ends and
lens-shaped T-ends and horospherical ends to themselves.  
If $\orb$ satisfies the triangle condition, or if every end is virtually factorizable, 
then $\mathcal{C}$ restricts to 
a correspondence between the lens-shaped R-ends and
lens-shaped T-ends and horospherical ends to themselves. 
\end{corollary} 

\begin{corollary} \label{pr-cor-tangency} 
	Let $\orb$ be a properly convex real projective $n$-orbifold.
	Let $\tilde E$ be a lens-shaped p-end. Then for a lens cone p-end neighborhood 
	$U$ of form $\{\mbv_{\tilde E}\} \ast L -\{\mbv_{\tilde E}\}$ 
	for lens $L$, we have 
	the upper boundary component $A$ is tangent to radial rays 
	from $\mbv_{\tilde E}$ at $\Bd A$.
	\end{corollary} 
\begin{proof} 
By Corollary \ref{pr-cor-duallens2}, 
$\tilde E$ corresponds to a lens-shaped T-p-end $\tilde E^\ast$ of $\orb^\ast$. 
Now, we take the dual domain $L_1'$ of $\{\mbv_{\tilde E}\} \ast L$ of 
$\tilde E^\ast$. The second part of Corollary \ref{pr-cor-duallens} applied 
to $L_1'$ gives us
the result for $\{\mbv_{\tilde E}\} \ast L$. 

Proposition \ref{prelim-prop-closureind} completes the proof. 
\hfill \SnT {\parfillskip0pt\par}
\end{proof}

\chapter[Applications]{Application: The openness of the lens properties, and expansion and shrinking of end neighborhoods} 
\label{ch-app}



This chapter lists applications of the main theory of Part 2, except for Chapter \ref{ch-np}, 
which contains the results we need in Part 3.
In Section \ref{app-sec-results}, 
We show that the lens-shaped property is stable under changes in holonomy representations. 
In Section \ref{app-sec-limitset},
we define limit sets of ends and discuss the properties. 
We obtain the exhaustion of $\torb$ by a sequence of p-end neighborhoods 
of $\torb$.
We show that any end neighborhood contains either a horospherical or a concave end neighborhood. 
We discuss maximal concave end neighborhoods.
In Section \ref{app-sec-strirr}, 
Corollary \ref{app-cor-disjclosure} shows that
the closures of the p-end neighborhoods are disjoint
in the closures of the universal cover in $\SI^n$ (resp. in $\RP^n$). 
From this, we prove the strong irreducibility of $\orb$, 
Theorem \ref{intro-thm-sSPC}, under
the conditions (IE) and (NA).


For results in this chapter, we don't necessarily assume that the holonomy group of $\pi_{1}(\orb)$ is strongly irreducible. 
In addition, we will not explicitly mention
Proposition \ref{prelim-prop-closureind} since its usage 
is well-established.

\subsection{SPC-structures and its properties} \label{intro-sub-spc}

\begin{definition} \label{intro-defn-IE}
For an orbifold $\mathcal{O}$, 
\begin{enumerate}
\item[(IE)] $\orb$ or $\pi_1(\orb)$ satisfies 
{\em infinite end index condition IE}
\index{end!condition!IE|textbf}
\index{end!condition!infinite end index|textbf}
if $[\pi_1(\mathcal O): \pi_1(E)] = \infty$ for the end fundamental group $\pi_1(E)$ of each end $E$. 
\item[(NA)] $\orb$ or $\pi_1(\orb)$ satisfies the {\em nonparallel end condition NA}  \index{end!condition!NA|textbf} 
\index{end!condition!nonparallel end|textbf} 
if $\pi_{1}(\tilde E_{1}) \cap \pi_{1}(\tilde E_{2})$ is finite for two distinct p-ends $\tilde E_{1}, \tilde E_{2}$ of $\orb$. 
\end{enumerate} 
\end{definition} 
These conditions are satisfied by complete hyperbolic manifolds with cusps.
These are group-theoretical properties with respect to the end groups. 
\index{IE|textbf} \index{NA|textbf} 

(IE) and (NA) imply that 
for each end $E$, 
$\pi_1(E)$ contains every element $g \in \pi_1(\orb)$ normalizing $\langle k \rangle$ for an infinite order $k \in \pi_1(E)$. 
We can choose a mutually disjoint set of end neighborhoods on $\orb$. 
Hence, the set of p-end neighborhoods that are components of the inverse images 
is mutually disjoint. 
If $k \not \in \pi_1(E)$, then $k$ sends any p-end corresponding to 
$E$ to a different p-end. Hence, $g = kgk^{-1}$ belongs to 
the intersection of two distinct p-end-fundamental groups. This is a contradiction. 


A {\em character} of a group $\Gamma$ to a Lie group $G$ 
is a conjugacy class of a representation $G$.
A representation or a character is said to be {\em stable} if the orbit of its representative is closed and the stabilizer is finite under the conjugation action in  $\Hom(\Gamma, G)$.
By Theorem 1.1 of Johnson-Millson \cite{JM87}, a representation $\rho$ is stable if and only if it is irreducible and no proper parabolic subgroup 
\index{character!stable|textbf} \index{representation!stable|textbf}
\index{stable representation|textbf} 
contains the image of $\rho$. 
Both stability and irreducibility are open conditions in the Zariski topology. 
Also, if the image of $\rho$ is Zariski dense, then $\rho$ is stable. 
$\PGL(n+1, \bR)$ acts properly on the open set of stable representations
in $\Hom(\pi_1(\mathcal{O}), \PGL(n+1,\bR))$. Similarly, 
$\SLpm$ acts properly on $\Hom(\pi_{1}(\orb), \SLpm)$.
(See \cite{JM87} for more details.)
(For $\SLpm$ and $\PGL(n+1,\bR)$, irreducibility is equivalent to 
the stability. However, we may use these terms simultaneously for emphasis. )

\begin{definition} \label{intro-defn-SPC} 
A {\em stable properly convex  real projective} ({\em SPC}-)structure on an $n$-orbifold 
is a structure of a properly convex real projective 
$n$-orbifold with a stable and irreducible holonomy group.
\end{definition}
\index{SPC-structure|textbf} 
\index{stable properly-convex real projective structure|textbf} 

\begin{definition}\label{intro-defn-strict} 
Suppose that an $n$-orbifold $\mathcal{O}$ has an SPC-structure. Let $\tilde U$ be 
the inverse image in $\tilde{\mathcal{O}} \subset \RP^n$ of the union $U$ of some choice of a collection of disjoint end neighborhoods of $\orb$. 
\index{SPC-structure|textbf}
If every straight arc and every non-$C^1$-point   in $\Bd \tilde{\mathcal{O}}$ 
are contained in the closure of a component of $\tilde U$, 
then $\mathcal{O}$ is said to be {\em strictly convex} with respect to the collection of the ends.  \index{convex!strictly|textbf}
In this case, $\mathcal{O}$ is also said to have a 
{{\em strict stable properly convex {\rm (}strict SPC-\/{\rm )}}}structure with respect to the collection of ends. \index{strict SPC-structure!strict|textbf}
\end{definition}
Currently, we can show that $U$ can be a union of any collection of end neighborhoods
when $\orb$ has only lens-type or horospherical ends
since $\clo(U_i) \cap \Bd \torb$ is independent of the choice of the end neighborhood
$U_i$. 
This is proved for the R ends by Theorems \ref{pr-thm-lensclass} and \ref{pr-thm-redtot}, and for the T ends by Theorem \ref{pr-thm-equ2}, using the lens property.
For horospherical ends, we can show that the closure of a p-end neighborhood
meeting the boundary of $\torb$ at a singleton since it has a cusp group 
as an end holonomy group by Section \ref{intro-sub-horo}.  We omit the obvious details. 


By a real projective orbifold {\em with generalized lens-shaped or horospherical $\cR$- or $\cT$-ends}, 
we mean one with a real projective structure that has a $\cR$-type or $\cT$-type assigned for
each end, where each $\cR$-end is either generalized lens-shaped or horospherical
and each $\cT$-end is lens-shaped or horospherical. 

Notice that the definition depends on the choice of $U$. However, we will show that
if each end is required to be a lens-shaped or horospherical
$\cR$- or $\cT$-end, then
we show that the definition is independent of $U$ in
Corollary \ref{app-cor-independence}. 


We will prove in Section \ref{app-sec-strirr} 
that topological conditions imply stability: 
\begin{theorem}\label{intro-thm-sSPC} 
Let $\orb$ be a noncompact strongly tame properly convex real projective
$n$-orbifold, $n \geq 2$, with 
generalized lens-shaped or horospherical $\cR$- or $\cT$-ends. 
Assume that $\orb$ satisfies  {\rm (}IE\/{\rm )} and {\rm (}NA\/{\rm )}.   
Then the holonomy group is strongly irreducible and not 
contained in any proper parabolic subgroup of $\PGL(n+1, \bR)$  {\rm (}resp.  $\SLpm${\rm ).} That is, the holonomy is stable. 
\end{theorem} 


\section{The openness of lens properties.} \label{app-sec-results}





As conditions in the representations of $\pi_1(\tilde E)$, the conditions for
generalized lens-shaped ends and one for lens-shaped ends are the same. 
Given a holonomy group of $\pi_1(\tilde E)$ acting on a generalized lens-shaped cone p-end neighborhood,  
it satisfies the uniform middle-eigenvalue condition by Theorem \ref{pr-thm-equ}. 
We can find a lens cone by choosing our orbifold to be ${\mathcal T}_{\mbv_{\tilde E}}(\tilde \Sigma_{\tilde E})^{o}/\pi_1(\tilde E)$ 
by Proposition \ref{pr-prop-convhull2}.

Let 
\[\Hom_{\mathcal{E}}(\pi_1(\tilde E), \SLnp)\: (\hbox{resp. } \Hom_{\mathcal{E}}(\pi_1(\tilde E), \PGL(n+1, \bR)))\]
denote the space of representations of the 
fundamental group of an $(n-1)$-orbifold $\Sigma_{\tilde E}$
fixing a common point if $\tilde E$ is a R-p-end and acting on a common hyperspace 
if $\tilde E$ is a T-p-end.
\index{homE@$\Hom_{\mathcal{E}}(\pi_1(\tilde E), \SLnp)$|textbf} 
\index{homE@$\Hom_{\mathcal{E}}(\pi_1(\tilde E), \PGL(n+1, \bR))$|textbf} 


Recall Definition \ref{intro-defn-R-ends}  
 for generalized lens-shaped R-ends. 
A {\em lens-shaped} representation for an R-end fundamental group is a representation 
acting on a (resp. generalized) lens cone as a p-end neighborhood. 
Hence, it acts properly, effectively and cocompactly on a lens with a compact Hausdorff quotient. 

We remark that the condition of having a lens cone or 
a generalized lens cone are the same by Theorems \ref{pr-thm-secondmain}
and \ref{pr-thm-equiv}. 
They are both equivalent to the uniform middle eigenvalue condition.

\index{representation!lens-shaped|textbf} 
\index{representation!generalized lens-shaped|textbf}

\begin{theorem}\label{app-thm-qFuch}
Let $\orb$ be a properly convex real projective $n$-orbifold. 
Assume that the universal cover $\torb$ is a subset of $\SI^n$ {\rm (}resp.\, $\RP^n${\rm ).} 
Let $\tilde E$ be a properly convex R-p-end of the universal cover $\torb$. 
Then 
\begin{enumerate}
\item[{\rm (i)}] $\tilde E$ is a generalized lens-shaped R-end if and only if $\tilde E$ is a strictly generalized lens-shaped R-end.
\item[{\rm (ii)}] The subspace of generalized lens-shaped  representations of an R-end is open in \[\Hom_{\mathcal{E}}(\pi_1(\tilde E), \SLnp) \:
(\hbox{resp. } \Hom_{\mathcal{E}}(\pi_1(\tilde E), \PGL(n+1, \bR))).\]
\end{enumerate}
Finally, if $\orb$ is properly convex and 
satisfies the triangle condition or $\tilde E$ is virtually factorizable, then we can replace the term  ``generalized lens-shaped''
with  ``lens-shaped'' in each of the above statements. 
\end{theorem}
\begin{proof} 
	We will assume $\torb \subset \SI^n$ first. 
(i) If $\pi_1(\tilde E)$ is non-virtually-factorizable, then the equivalence follows from Theorem \ref{pr-thm-lensclass} (i), and 
if $\pi_1(\tilde E)$ is virtually factorizable, 
then it follows from Theorem \ref{pr-thm-redtot} (ii). The converse is obvious. 

(ii) Let $\mu$ be a representation $\pi_1(\tilde E) \ra \SLnp$ associated with a generalized lens cone.  
By Theorem \ref{pr-thm-secondmain}, $\pi_1(\tilde E)$ satisfies the uniform middle eigenvalue condition with respect to $\mbv_{\tilde E}$.  
By Theorem \ref{pr-thm-equiv}, 
we obtain a lens $K$ in ${\mathcal T}_{\mbv_{\tilde E}}(\tilde \Sigma_{\tilde E})$ with smooth convex boundary components
$A \cup B$ since ${\mathcal T}_{\mbv_{\tilde E}}(\tilde \Sigma_{\tilde E})$ itself satisfies the triangle condition, although it is not properly convex. 
(Note, we don't need $K$ to be in $\torb$ for the proof.)

$K/\mu(\pi_1(\tilde E))$ is a compact orbifold whose boundary is the union of two closed $n$-orbifold components
$A/\mu(\pi_1(\tilde E)) \cup B/\mu(\pi_1(\tilde E))$.
Suppose that $\mu'$ is sufficiently near $\mu$. 
We may assume that $\mbv_{\tilde E}$ is fixed by conjugating $\mu'$ by a bounded projective transformation. Hence, $K$ deforms to a lens under $\mu'$. 
Considering the radial segments in $K$, we also obtain a foliation by the radial lines in
$\{\mbv_{\tilde E}\} \ast K$. 
By Proposition \ref{app-prop-koszul}, applying Proposition \ref{app-prop-lensP} to both boundary components of 
the lens, we obtain a lens cone in a tube domain 
${\mathcal T}'_{\mbv_{\tilde E}}$ 
deformed from the original one. 
This implies that a sufficiently small change of holonomy keeps $\tilde E$ having a concave p-end neighborhood. 
Proposition \ref{pr-prop-orbit} proves the strictness of the new lens deformed from $K$. 
This completes the proof of (ii). 


The final statement follows from Lemma \ref{pr-lem-genlens}. 
\hfill \SnT {\parfillskip0pt\par}
\end{proof} 






\begin{theorem}\label{app-thm-qFuch2}
Let $\orb$ be a strongly tame properly convex real projective $n$-orbifold. 
Assume that the universal cover $\torb$ is a subset of $\SI^n$ {\rm (}resp. of $\RP^n${\rm ).}
Let $\tilde E$ be a T-p-end of the universal cover $\torb$. 
Then the subspace of  lens-shaped representations for a T-p-end is open in 
\[\Hom_{\mathcal{E}}(\pi_1(\tilde E), \SLnp) \:
(\hbox{resp. } \Hom_{\mathcal{E}}(\pi_1(\tilde E), \PGL(n+1, \bR))).\]
\end{theorem}
\begin{proof} 
By Theorem \ref{pr-thm-equ2}, the condition of the lens T-p-end
is equivalent to the uniform middle-eigenvalue condition for the end. 
Proposition \ref{pr-prop-dualend2}, Theorem \ref{pr-thm-secondmain},
Proposition \ref{pr-cor-duallens2}, and  
Theorem \ref{app-thm-qFuch} complete the proof. 
\hfill \SnT{\parfillskip0pt\par}
\end{proof} 




\begin{corollary}\label{app-cor-mideigen}
We are given a properly convex R- or T-end $\tilde E$ of 
a strongly tame convex $n$-orbifold $\orb$.
Assume that $\torb \subset \SI^n$ {\rm (}resp.\, $\torb \subset \RP^n${\rm ).}  
Then the subset of 
\[\Hom_{\mathcal{E}}(\pi_1(\tilde E), \SL_\pm(n+1, \bR)) \hbox{ {\rm (}resp.}\, \Hom_{\mathcal{E}}(\pi_1(\tilde E), \PGL(n+1, \bR)) ) \]
consisting of  representations satisfying the uniform middle-eigenvalue condition
with respect to certain choices of fixed points or fixed hyperspaces of
the holonomy group is open.
\end{corollary} 
\begin{proof} 
For R-p-ends, this follows from Theorems \ref{pr-thm-equ} and \ref{app-thm-qFuch}.
For T-p-ends, this follows from dual results: Theorem \ref{pr-thm-equ2} and Theorems \ref{app-thm-qFuch2}. 
\hfill \SnT {\parfillskip0pt\par}
\end{proof}

\section{The end and the limit sets.} \label{app-sec-limitset}

\begin{definition} \label{app-defn-limitset}  $  $ 
\begin{itemize}
\item Define the {\em limit set} $\Lambda(\tilde E)$ of an R-p-end $\tilde E$ with a generalized p-end neighborhood to be 
$\Bd D - \partial D$ for a generalized lens $D$ of $\tilde E$ in $\SI^n$ {\rm (}resp. $\RP^n${\rm ).}
This is identical with the set $\Lambda_{\tilde E}$ in Definition \ref{pr-defn-lambdaE} by Corollary \ref{pr-cor-LambdaW}. 
\item The {\em limit set} $\Lambda(\tilde E)$ of a lens-shaped T-p-end $\tilde E$ to be 
$\clo(\tilde S_{\tilde E})- \tilde S_{\tilde E}$ for the ideal boundary component $\tilde S_{\tilde E}$ of $\tilde E$. 
\item The limit set of a horospherical end is the set of the end vertex. 
\end{itemize} 
\end{definition} 
\index{end!limit set} 

The definition does depend on whether we work on $\SI^n$ or $\RP^n$. 
However, by Proposition \ref{prelim-prop-closureind}, 
there are always straightforward one-to-one correspondences. 
We remark that this may not equal the closure of the union
of the attracting fixed sets in some cases. 


\begin{corollary}\label{app-cor-independence} 
Let $\mathcal{O}$ be a convex real projective $n$-orbifold
where $\torb \subset \SI^n$ {\rm (resp. }$\subset \RP^n$\/{\rm ).} 
Let $U$ be a p-end neighborhood of $\tilde E$ where $\tilde E$ is a lens-shaped T-p-end 
or a generalized lens-shaped or lens-shaped or horospherical R-p-end. 
Then $\clo(U) \cap \Bd \torb$
equals $\clo(\tilde S_{\tilde E})$ or $\clo(\bigcup S(\mbv_{\tilde E}))$ or $\{\mbv_{\tilde E}\}$ depending on 
whether $\tilde E$ is a lens-shaped T-p-end or a generalized lens-shaped or horospherical R-p-end. Furthermore, 
this set is independent of the choice of $U$ 
and so is the limit set $\Lambda(\tilde E)$ of $\tilde E$.
\end{corollary}
\begin{proof} 
	We first assume $\torb \subset \SI^n$. 
Let $\tilde E$ be a generalized lens-shaped R-p-end. Then, by Theorem \ref{pr-thm-equ}, 
$\bGamma_{\tilde E}$ satisfies the uniform middle-eigenvalue condition with respect to $\tilde E$. Suppose that $\pi_1(\tilde E)$ is not virtually factorizable. 
Let $L^b$ denote $\partial {\mathcal T}_{\mbv_{\tilde E}}(\tilde \Sigma_{\tilde E}) \cap L$ for a distanced compact convex set $L$ on which $\bGamma_{\tilde E}$ acts.
We have $L^b = \Lambda(\tilde E)$ by Proposition \ref{pr-prop-orbit}. 
Since $S(\mbv_{\tilde E})$ is 
an $h(\pi_1(\tilde E))$-invariant set, and the convex hull of 
$\Bd \bigcup S(\mbv_{\tilde E})$ is a distanced compact convex set
by the proper convexity of $\tilde \Sigma_{\tilde E}$,
Theorems \ref{pr-thm-lensclass} and \ref{pr-thm-redtot}
show that the 
limit set is determined by the set $L^b$ in $\bigcup S(\mbv_{\tilde E})$, 
and $\clo(U) \cap \Bd \torb = \bigcup S(\mbv_{\tilde E})$. 

Suppose now that $\pi_1(\tilde E)$ is virtually factorizable. Then, by Theorem \ref{pr-thm-redtot}, $\tilde E$ is a totally geodesic R-p-end. 
Proposition \ref{pr-prop-orbit} and Theorem 
\ref{pr-thm-redtot} again imply the result. 

Let $\tilde E$ be a T-p-end. Theorems \ref{pr-thm-equ2} and \ref{du-thm-lensn} imply
\[\clo(A) - A \subset \clo(\tilde S_{\tilde E}) \hbox{ for } A= \Bd L \cap \torb\]
for a cocompactly-acted lens neighborhood $L$ by the strictness of the lens. 
Thus, $\clo(U) \cap \Bd \torb$  equals $\clo(\tilde S_{\tilde E})$. 

For horospherical ones, 
we simply use the definition to obtain the result. 
\hfill \SnT  {\parfillskip0pt\par}
\end{proof}

Proposition \ref{prelim-prop-closureind} shows that this definition is 
equivalent to Definition \ref{intro-defn-strict}.  
Corollary \ref{app-cor-strictconv} shows the independence of the definition
with respect to the choice of the end neighborhoods
when the ends are generalized lens-type $\cR$-end or lens-shaped $\cT$-ends. 
We conjecture that this also holds  for the ends of four types given by 
Ballas-Cooper-Leitner \cite{BCLp}.

\begin{corollary} \label{app-cor-strictconv} 
Suppose that $\mathcal O$ is an $n$-dimensional strongly tame strictly SPC-orbifold
with generalized lens-shaped R-ends or lens-shaped T-ends or horospherical ends. 
Let $\torb$ be a properly convex domain in $\RP^n$ {\rm (} resp. in $\SI^n$\,
{\rm )} 
projectively covering $\orb$.
Choose any disjoint collection of end neighborhoods in $\orb$. Let $U$ denote their union. \index{po@$p_{\orb}$}
Let $p_{\orb}: \torb \ra \orb$ denote the universal cover. 
Then any segment or a non-$C^1$-point of $\Bd \torb$ is contained in the closure of a component of $p_{\orb}^{-1}(U)$ for 
any choice of $U$. 
\end{corollary}
\begin{proof}
%
We first assume $\torb \subset \SI^n$. 
By the definition of a strict SPC-orbifold, any segment or a non-$C^1$-point has to be in the closure of 
a p-end neighborhood.  Corollary \ref{app-cor-independence} proves the claim. 
\hfill \SnT {\parfillskip0pt\par}
\end{proof}

%
%

\subsection{Convex hulls of ends.} \label{app-sub-convh}

We will sharpen Corollary \ref{app-cor-independence} and the convex hull part of Lemma \ref{app-lem-expand}.   
Again, these sets are all defined in $\SI^n$, and
we define the corresponding objects for $\RP^n$ by their images in $\RP^n$. 
by Proposition \ref{prelim-prop-closureind}.

One can associate a {\em convex hull $I(\tilde E)$ of a p-end} $\tilde E$ 
of $\tilde{\mathcal{O}}$ as follows: 
\index{p-end! convex hull|textbf} 
\begin{itemize}
	\item For horospherical p-ends, the convex hull of each is defined to be the set of the end vertex. 
	\item The convex hull of a lens-shaped totally geodesic p-end $\tilde E$ is the closure $\clo(\tilde S_{\tilde E})$
	the totally geodesic ideal boundary component $\tilde S_{\tilde E}$ corresponding to $\tilde E$. 
	\item For a generalized lens-shaped p-end $\tilde E$,  the convex hull 
	of $\tilde E$ is the convex hull of $\bigcup S(\mbv_{\tilde E})$ in $\clo(\torb)$; that is, 
	\[ I(\tilde E):= \CH(\bigcup S(\mbv_{\tilde E})).\]  
\end{itemize} 
The first two equal $\clo(U) \cap \Bd \torb$ for any p-end neighborhood $U$ of $\tilde E$ by Corollary \ref{app-cor-independence}.
Corollary \ref{app-cor-independence}  and Proposition \ref{app-prop-I} imply that the convex hull of an end is well defined.
We can also characterize it as the intersection 
\[I(\tilde E)=\bigcap_{U_{1}\in \mathcal U}\CH(\clo(U_1))\] for 
the collection $\mathcal U$ of p-end neighborhoods $U_1$ of $\mbv_{\tilde E}$
by Proposition \ref{app-prop-I}.
\index{end!convex hull|textbf} 
\index{itildeE@$I(\tilde E)$|textbf}

We define $\partial_S I(\tilde E)$ as the set of endpoints of 
maximal rays from $\mbv_{\tilde E}$ ending at $\Bd I(\tilde E)$ 
and in the directions of $\tilde \Sigma_{\tilde E}$. 
It is homeomorphic to $\tilde \Sigma_{\tilde E}$ by rays
and has a compact quotient under $\Gamma_{\tilde E}$.
\index{partialSitildeE@$\partial_S I(\tilde E)$|textbf} 
	Since the convex hull of $\bigcup S(\mbv_{\tilde E})$ is 
a subset of the tube with a vertex $\mbv_{\tilde E}$ 
in the directions of elements of 
$\clo(\Pi_{\mbv_{\tilde E}}(R_{\mbv_{\tilde E}}(\Sigma_{\tilde E})))$, 
we obtain 
\begin{equation}\label{app-eqn-pI}  
\Bd I(\tilde E) = \partial_S I(\tilde E) \cup 
\bigcup S(\mbv_{\tilde E}).
\end{equation}  

	\begin{lemma}\label{app-lem-sigma}
		If $\sigma \in S_i$ meets $\partial_S I(\tilde E)$, 
		then $\sigma^o \subset \partial_S I(\tilde E)$
		and the vertices of $\sigma$ are endpoints of maximal segments in 
		$S(\mbv_{\tilde E})$. 
		\end{lemma} 
	\begin{proof} 
	Suppose that $x \in \partial_S I(\tilde E)$ and $x \in \sigma^o$ for 
	a simplex $\sigma\in S_i$ for minimal $i$. 
	The vertices of $\sigma$ lie in $\bigcup S(\mbv_{\tilde E})$. 
	If at least one vertex $v_1$ is in the interior of a segment in 
	$S(\mbv_{\tilde E})$, then by taking points in the neighborhood 
	of $x$ in $\bigcup S(\mbv_{\tilde E})$, we can deduce that $\sigma^o$ is not
	in the boundary of the convex hull. Moreover, 
	$\sigma^o$ is a subset of $\partial I_{\tilde E}$ by 
	Lemma \ref{prelim-lem-simplexbd}. 
	Hence, the vertices are the endpoints of the maximal segments in 
	$S(\mbv_{\tilde E})$. 
	
	Suppose that a point of $\sigma^o$ lies in a segment in 
	$S({\mbv_{\tilde E}})$. Then an interior point of 
	$\Pi_{\mbv_{\tilde E}}(\sigma)$ meets 
	the boundary of $\clo(R_{\mbv_{\tilde E}}(\torb))$. 
	By  Lemma \ref{prelim-lem-simplexbd}, 
	$\Pi_{\mbv_{\tilde E}}(\sigma) \subset 
	\Bd (R_{\mbv_{\tilde E}}(\torb))$. 
	Thus, $\sigma$ is in a union of segments from 
	$\mbv_{\tilde E}$ in the directions of $\Bd (R_{\mbv_{\tilde E}}(\torb))$. 
	By Theorems \ref{pr-thm-lensclass} and \ref{pr-thm-redtot}, these segments are contained in $\bigcup S(\mbv_{\tilde E})$.
	We obtain $\sigma \subset \bigcup S(\mbv_{\tilde E})$.  
	This is a contradiction. 
	By \eqref{app-eqn-pI},  
	$\sigma^o \subset \partial_S I_{\tilde E}$. 
\hfill \SSn {\parfillskip0pt\par}
	\end{proof}

A {\em topological orbifold} is one in which we are allowed to use
continuous maps as charts.
We say that two very good topological orbifolds $\orb_1$ and $\orb_2$ are {\em homeomorphic} if 
there is a homeomorphism of the base spaces and 
that is induced from a homeomorphism $f:M_1 \ra M_2$ of respective 
very good manifold covers $M_1$ and $M_2$ 
where $f$ induces the isomorphism of the deck transformation group of  
$M_1 \ra \orb_1$ to that of $M_2 \ra \orb_2$.
\index{orbifold!homeomorphism}

\begin{proposition}\label{app-prop-I} 
	Let $\orb$ be a properly convex real projective $n$-orbifold with radial ends or lens-shaped totally geodesic ends. 
	Let $\tilde E$ be a generalized 
	lens-shaped R-p-end, and let $\mbv_{\tilde E}$ be an associated 
	p-end vertex. Let $I(\tilde E)$ be the convex hull of $\tilde E$. 
	\begin{enumerate}
		\item[{\rm (i)}] Suppose that $\tilde E$ is a lens-shaped radial p-end. 
		Then $\partial_S I(\tilde E)= \Bd I(\tilde E) \cap \tilde{\mathcal{O}}$, and 
		$\partial_S I(\tilde E)$ is contained in the lens in a lens-shaped p-end neighborhood. 
		\item[{\rm (ii)}] $I(\tilde E)$ contains any concave p-end neighborhood of $\tilde E$ and 
		\begin{gather*} 
		I(\tilde E) = \CH(\clo(U)) \\
		I(\tilde E) \cap \torb = \CH(\clo(U)) \cap \torb
		\end{gather*}
		for a concave p-end neighborhood $U$ of $\tilde E$. 
		Thus, $I(\tilde E)$ has a nonempty interior. 
		\item[{\rm (iii)}] Each segment from $\mbv_{\tilde E}$ maximal in $\tilde{\mathcal{O}}$ 
		meets the set $\partial_S I(\tilde E)$ exactly once and
		$\partial_S I(\tilde E)/\bGamma_{\tilde E}$ is homeomorphic to
		 $\Sigma_E$. 
		\item[{\rm (iv)}] There exists a nonempty interior of the convex hull $I(\tilde E)$ of $\tilde E$ 
		where $\bGamma_{\tilde E}$ acts such that $I(\tilde E)^o/\bGamma_{\tilde E}$ is homeomorphic to the end orbifold times an interval. 
	\end{enumerate}
\end{proposition}
\begin{proof}
	Assume first that $\torb \subset \SI^n$. 
	(i) Suppose that $\tilde E$ is lens-shaped. 
	We define $S_1$ as the set of $1$-simplices with endpoints in segments in $\bigcup S({\mbv_{\tilde E}})$, and we inductively define
	$S_i$ to be the set of $i$-simplices with faces in $S_{i-1}$. 
	Then 
	\[I(\tilde E) = \bigcup_{\sigma \in S_1 \cup S_2 \cup \cdots \cup S_m} \sigma.\] 
	

	
	Since any point of $\partial_S I(\tilde E)$ is in some simplex 
	$\sigma$, $\sigma \in S_i$, 
	we obtain that $\partial_S I(\tilde E)$ is the union 
		\[\bigcup_{\sigma \in S_1 \cup S_2 \cup \cdots \cup S_m, \, \sigma^o \subset \partial_S I(\tilde E)} \sigma^o\]
	by Lemma \ref{app-lem-sigma}.
	
	Suppose that $\sigma \in S_i$ with $\sigma^o \subset \partial_S I(\tilde E)$. Then 
	each of its vertices must be in an endpoint of a segment in $S({\mbv_{\tilde E}})$
	by Lemma \ref{app-lem-sigma}. 
	By Theorems \ref{pr-thm-lensclass} and \ref{pr-thm-redtot},
	the endpoints of the segments in $S(\mbv_{\tilde E})$ 
	are in $\Lambda(\tilde E)$.
	Hence,	by the convexity of $L$, 
 $\sigma^o$ is contained in a cocompactly-acted lens-shaped domain $L$ 
since the vertices of $\sigma$ all lie in $\Bd L - \partial L = \Lambda(\tilde E)$

	
	Thus, each point of $\partial_S I(\tilde E)$ is in $L^o \subset \tilde \Sigma$. 
	Hence, $\partial_S I(\tilde E) \subset \Bd I(\tilde E) \cap \torb$. 
	Conversely, a point of $ \partial I(\tilde E) \cap \torb$ is an endpoint of 
	a maximal segment in a direction of $\tilde \Sigma_{\tilde E}$.
	By \eqref{app-eqn-pI}, we obtain
	$\partial_S I(\tilde E) = \partial I(\tilde E) \cap \torb$.
	
	(ii) Since $I(\tilde E)$ contains the segments in $S({\mbv_{\tilde E}})$ and is convex, and so does a concave p-end neighborhood $U$, 
	we obtain $\Bd U \subset I(\tilde E)$: 
	Otherwise, let $x$ be a point of $\Bd U \cap \Bd I(\tilde E) \cap \torb$ where some neighborhood 
	in $\Bd U$ is not in $I(\tilde E)$. Then, since $\Bd U$ is a union of a
	strictly convex hypersurface 
	$\Bd U \cap \torb$ and $S({\mbv_{\tilde E}})$, 
	each sharply supporting hyperspace at $x$ of the convex set $\Bd U \cap \torb$
	meets a segment in $S({\mbv_{\tilde E}})$ in its interior: 
	consider the lens $L$ such that one of the boundary components is 
	$\Bd U$. The supporting hyperspace at the boundary component 
	cannot meet the closure of $L$ at any other point by the strict
	convexity. 
	
	This is a contradiction since $x$ must then be in $I(\tilde E)^{o}$. 
	Thus, $U \subset I(\tilde E)$. Thus, \[\CH(\clo(U)) \subset I(\tilde E).\] 
	Conversely, since $\clo(U) \supset \bigcup S({\mbv_{\tilde E}})$ by Theorems \ref{pr-thm-lensclass} and \ref{pr-thm-redtot}, we obtain that 
	\[\CH(\clo(U)) \supset I(\tilde E).\] 
	
	
	(iii) 
	We again use Proposition \ref{prelim-thm-vgood}. It is sufficient
	to prove the result by taking a very good cover of 
	$\orb$ and considering the corresponding p-end to $\tilde E$. 
	Each point of it meets a maximal segment from ${\mbv_{\tilde E}}$ in the end but not in $S({\mbv_{\tilde E}})$ at exactly one point since a maximal segment must eventually leave the lens cone.
	Thus, $\partial_S I(\tilde E)$ is homeomorphic to an $(n-1)$-cell and the result follows. 
	
	(iv) This follows from (iii) since we can use rays from $x$ meeting $\Bd I(\tilde E) \cap  \tilde{\mathcal{O}}$ at unique points 
	and use them as leaves of a fibration. 
\hfill	\SnT  {\parfillskip0pt\par}
\end{proof}


\begin{figure}
	\centering
	\includegraphics[height=7cm]{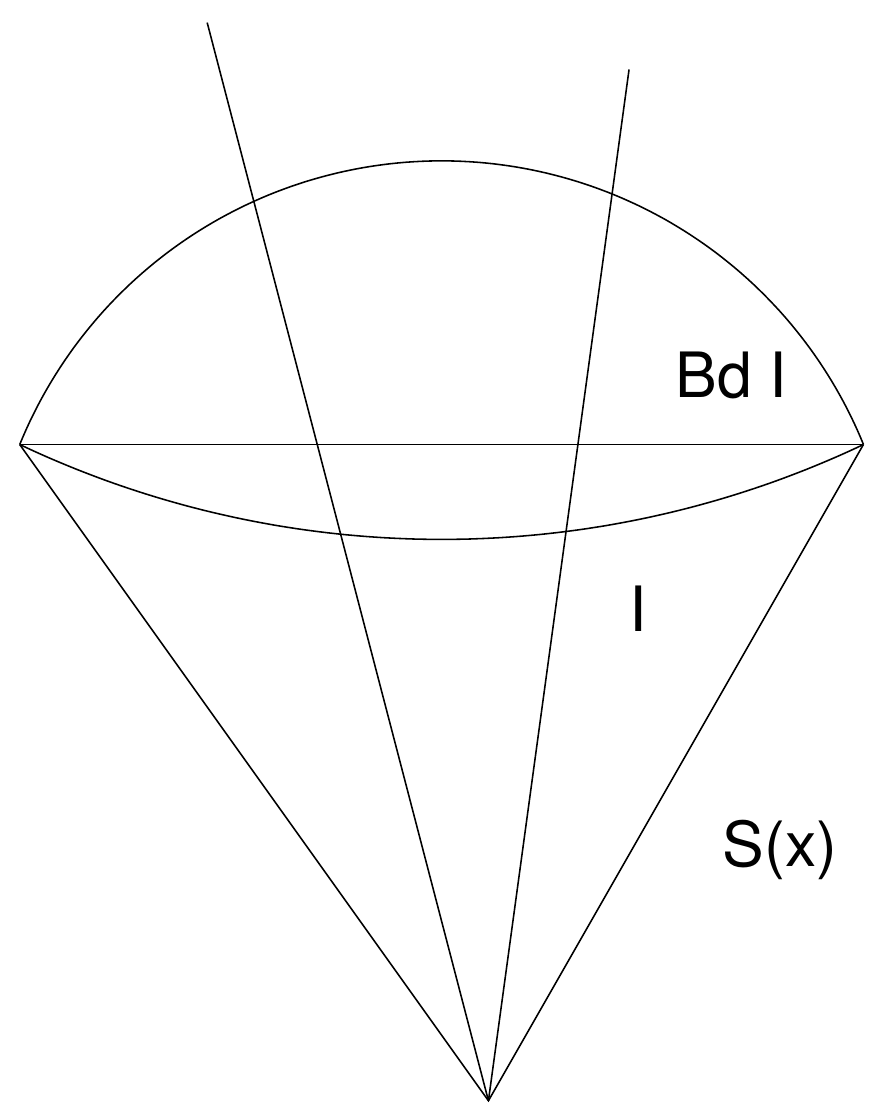}
	
	\caption{The structure of a lens-shaped p-end.}
	\label{app-fig-lenslem2}
\end{figure}


\subsection{Expansion of lens or horospherical p-end neighborhoods.} 



\begin{lemma}\label{app-lem-expand}  
Let $\mathcal{O}$ have a properly convex real projective structure $\mu$. 
\begin{itemize} 
\item Let $U_1$ be a p-end neighborhood of a horospherical or a lens-shaped R-p-end $\tilde E$ with the p-end vertex $\mbv_{\tilde E}${\rm ;} 
or 
\item Let $U_1$ be a lens-shaped p-end neighborhood of a T-p-end $\tilde E$.
\end{itemize} 
Let $\bGamma_{\tilde E}$ denote the p-end holonomy group corresponding to $\tilde E$. 
Then we can construct a sequence of lens cone or lens p-end neighborhoods $U_{i}$, $i=1, 2, \dots, $ satisfying 
$U_{i} \subset U_{j} \subset  \torb$ for every pair $i, j$, $i > j$, where 
the following hold\,{\rm :} 
\begin{itemize} 
\item Given a compact subset of $\tilde{\mathcal{O}}$, there exists an integer $i_0$ such that 
$U_i$ for $i > i_0$ contains it. 
\item The Hausdorff distance between $U_i$ and $\tilde{\mathcal{O}}$ can be made as small as possible, i.e., 
\[ \forall \eps > 0, \exists\, J > 0, \hbox{ such that }  \bdd_H (U_i, \torb) < \epsilon \hbox{ for } i > J. \]
\item There exists a sequence of convex open p-end neighborhoods $U_i$ of $\tilde E$ in $\tilde{\mathcal{O}}$ 
such that $(U_i - U_j)/\bGamma_{\tilde E}$ for a fixed $j$ and $i> j$ is diffeomorphic to a product of an open interval with 
the end orbifold. 
\item We can choose $U_i$ such that $\Bd U_i \cap \torb$ is smoothly embedded and strictly convex with 
$\clo(\Bd U_i) - \torb \subset \Lambda(\tilde E)$. 
\end{itemize}
\end{lemma}
\begin{proof} 
%
Suppose that $\torb \subset \SI^n$ first. 
Suppose that $\tilde E$ is a lens-shaped R-end first. 
Let $U_1$ be a  lens cone.
Take a union of finitely many geodesic leaves $L$ from $\mbv_{\tilde E}$ in $\torb$ 
of $d_{\torb}$-length $t$ outside the lens cone $U_1$ and
take the convex hull of $U_1$ and $\bGamma_{\tilde E}(L)$ in $\tilde{\mathcal{O}}$. 
Denote the result by $\Omega_t$. Thus, the endpoints of $L$ not equal to $\mbv_{\tilde E}$ are in $\torb$.

We claim that 
\begin{itemize}
\item $\Bd \Omega_t \cap \tilde{\mathcal{O}}$ is a connected $(n-1)$-cell, 
\item $(\Bd \Omega_t \cap \tilde{\mathcal{O}})/\bGamma_{\tilde E}$ is 
a compact $(n-1)$-orbifold diffeomorphic to $\Sigma_{\tilde E}$, and  
\item $\Bd U_1 \cap \torb$ bounds 
a compact orbifold diffeomorphic to the product of a closed interval with  
$(\Bd \Omega_t \cap \tilde{\mathcal{O}})/\bGamma_{\tilde E}$: 

\end{itemize} 
First, for $l$ in $L$, 
 any converging subsequence of
$\{g_i(l)\}, g_i\in \bGamma_{\tilde E}$, converges to a segment in
$S(\mbv_{\tilde E})$ for a sequence $\{g_i\}$ of mutually distinct elements. 
This follows since the limit is a segment in $\Bd \tilde{\mathcal{O}}$ with an endpoint $\mbv_{\tilde E}$ 
and must belong to $S(\mbv_{\tilde E})$ 
by Theorems \ref{pr-thm-lensclass} and \ref{pr-thm-redtot}.  

Let $S_1$ be the set of segments with endpoints in $\bGamma_{\tilde E}(L) \cup \bigcup S(\mbv_{\tilde E})$.
Inductively, we define $S_i$ as the set of simplices all of whose sides are in $S_{i-1}$. 
Then the convex hull of $\bGamma_{\tilde E}(L)$ in $\clo(\torb)$ is a union of $S_1 \cup \cdots \cup S_m$. 

We claim that for each maximal segment $s$ in $\clo(\torb)$ from $\mbv_{\tilde E}$ 
not in $S(\mbv_{\tilde E})$, $s^o$ meets $\Bd \Omega_t \cap \torb$ at a unique point: 
Suppose not. 
Then let $v'$ be the other endpoint of $s$ in  $\Bd \tilde{\mathcal{O}}$ 
with $s^{o}\cap \Bd \Omega_t \cap \torb=\emp$. 
Thus, $v' \in \Bd \Omega_{t}$. 

Now, $v'$ is contained in the interior of a simplex $\sigma$ in $S_i$ for some $i$.
Since $\sigma^o \cap \Bd \torb \ne \emp$, it follows that $\sigma\subset \Bd \torb$ according to Lemma \ref{prelim-lem-simplexbd}.
Since the endpoints $\bGamma_{\tilde E}(L)$ are in $\torb$, the only possibility is that
the vertices of $\sigma$ are in $\bigcup S(\mbv_{\tilde E})$. 
Also, $\sigma^{o}$ is transverse to radial rays since, otherwise, $v'$ is not in $\Bd \torb$. 
Thus, $\sigma^{o}$ projects to an open simplex of the same dimension in $\tilde \Sigma_{\tilde E}$.
Since $U_1$ is convex and contains $\bigcup S(\mbv_{\tilde E})$ in its boundary, 
$\sigma$ is in the lens cone $\clo(U_1)$. 
Since a lens cone has a boundary that is a union of a strictly convex open hypersurface $A$ and $\bigcup S(\mbv_{\tilde E})$, 
and $\sigma^{o}$ cannot tangentially meet $A$, it follows that
$\sigma^{o}$ is in the interior of the lens cone, 
and no interior point of $\sigma$ is in $\Bd \tilde{\mathcal{O}}$, a contradiction. 
Therefore, each maximal segment $s$ from $\mbv_{\tilde E}$ exactly reaches the boundary $\Bd \Omega_t \cap \tilde{\mathcal{O}}$ once. 

As in Lemma \ref{pr-lem-infiniteline}, $\Bd \Omega_t \cap \torb$ contains no segment ending in $\Bd \torb$.  
The strict convexity of $\Bd \Omega_t \cap \torb$ follows 
 as in the proof of Proposition \ref{pr-prop-convhull2} by smoothing. 
Taking sufficiently many leaves for $L$ with sufficiently large $d_{\torb}$-lengths $t_i$, we can show that any compact subset is contained in $\Omega_t$. 
Choose some sequence $\{t_i\}$ such that $\{t_i\} \ra \infty$ as $i \ra \infty$. 
Now, let $U_i := \Omega_{t_i}$. 
From this, the final item follows. 
The first three items now follow if $\tilde E$ is an R-end. 

%

Suppose now that $\tilde E$ is horospherical and
$U_1$ is a horospherical p-end neighborhood. 
We can smooth its boundary to be strictly convex. 
$\bGamma_{\tilde E}$ is in a cusp subgroup by the premise. 
Taking $L$ sufficiently dense, we can choose similarly as above a sequence $\Omega_i$ of polyhedral convex horospherical open sets at $\mbv_{\tilde E}$ 
such that eventually any compact subset of $\tilde{\mathcal{O}}$ is in it for sufficiently large $i$.
Theorem \ref{du-thm-horonpre} provides a smooth, strictly convex horospherical p-end neighborhood 
$U_i$. 

Suppose now that $\tilde E$ is totally geodesic. We then use the dual domain $\torb^*$ and the group $\bGamma_{\tilde E}^*$. 
Let $\mbv_{\tilde E^*}$ denote the dual vertex of the hyperspace
containing $\tilde S_{\tilde E}$. 
By the diffeomorphism induced by great segments with the common endpoint $\mbv_{\tilde E}^*$, we obtain an orbifold homeomorphism 
 \[(\Bd \torb^* - \bigcup S(\mbv_{\tilde E^*}))/\bGamma_{\tilde E}^* \cong \Sigma_{\tilde E^\ast}/\bGamma_{\tilde E}^*,\] 
between compact orbifolds.
Then we obtain $U_i$ that contains $\torb^*$ in ${\mathcal T}_{\mbv_{\tilde E}}(\tilde \Sigma_{\tilde E})$
by taking finitely many hyperspace $F_{i}$ disjoint from $\torb^*$ but meeting
${\mathcal T}_{\mbv_{\tilde E}}(\tilde \Sigma_{\tilde E})^o$. Let $H_{i}$ be the open hemisphere that contains $\torb^*$ bounded by $F_{i}$. 
Then we form $U_1 := \bigcap_{g\in \bGamma_{\tilde E}} g(H_i)$. 
By taking more hyperspaces, we obtain a sequence 
\[U_1 \supset U_2 \supset \cdots \supset U_i \supset U_{i+1} \supset \cdots \supset \torb^* \]
such that $\clo(U_{i+1}) \subset U_i$ and 
\[\bigcap_i \clo(U_i) =\clo(\torb^*). \] 
That is, by using sufficiently many hyperspaces, we can make 
$U_i$ disjoint from any compact subset disjoint from $\clo(\torb^*)$.
Now taking the dual $U_i^*$ of $U_i$ and using equation \eqref{prelim-eqn-dualinc} we obtain
\[ U_1^* \subset U_2^* \subset \cdots \subset U_i^* \subset U_{i+1}^* \subset \cdots \subset \torb.\]
Then $U_{i}^* \subset \torb$ is an increasing sequence that eventually contains all compact subsets of $\torb$ by duality from the above disjointness.
This completes the proof for the first three items.

The fourth item follows from Corollary \ref{app-cor-independence}.
\hfill \SnP {\parfillskip0pt\par}
\end{proof}

\subsection{Shrinking of lens and horospherical p-end neighborhoods.} 

We now discuss the ``shrinking'' of p-end neighborhoods. These repeat some results. 
In this subsection, we do not necessarily assume the strong tameness for later use. 

\begin{corollary} \label{app-cor-shrink} 
Suppose that $\orb$ is a properly convex real projective $n$-orbifold 
and let $\torb$ be a properly convex domain in $\SI^n$ {\rm (}resp. $\RP^n${\rm )}
projectively covering $\orb$. 
Then the following statements hold\,{\rm :} 
\begin{enumerate} 
\item[{\rm (i)}] If $\tilde E$ is a horospherical R-p-end, then 
every p-end neighborhood of $\tilde E$ contains a horospherical p-end neighborhood.
\item[{\rm (ii)}] Suppose that $\tilde E$ is a generalized lens-shaped or lens-shaped R-p-end. 
Let $I(\tilde E)$ be the convex hull of $\bigcup S(\mbv_{\tilde E})$, 
and let $V$ be a p-end neighborhood where $(\Bd V \cap \torb)/\pi_1(\tilde E)$ is a compact orbifold. 
If $V^o \supset I(\tilde E) \cap \torb$, then 
$V$ contains a lens cone p-end neighborhood of $\tilde E$. 
\item[{\rm (iii)}] If $\tilde E$ is a generalized lens-shaped R-p-end or satisfies the uniform middle-eigenvalue condition, then 
every p-end neighborhood of $\tilde E$ contains a concave p-end neighborhood.
\item[{\rm (iv)}] Suppose that $\tilde E$ is a lens-shaped T-p-end or satisfies the uniform middle-eigenvalue condition.
Then every p-end neighborhood contains a lens p-end neighborhood $L$ with strictly convex boundary in $\torb$. 
\item[{\rm (v)}] Suppose further that $\orb$ is strongly tame. 
We can choose a collection of mutually disjoint end neighborhoods for 
all ends that are lens-shaped T-end neighborhood, concave R-end neighborhood or a horospherical ones. 
\end{enumerate} 
\end{corollary} 
\begin{proof}
Suppose that $\torb \subset \SI^n$ first. 

(i) Let $\mbv_{\tilde E}$ denote the R-p-end vertex corresponding to $\tilde E$. 
By premise, we obtain a conjugate $G$ of a subgroup of 
a cusp subgroup of $\SO(n, 1)$ as the finite-index subgroup of $h(\pi_1(\tilde E))$
acting on $U$, a p-end neighborhood of $\tilde E$. 
For any $\eps$, 
we can choose a $G$-invariant ellipsoid of $\bdd$-diameter $\leq \eps$ 
in $U$ containing $\mbv_{\tilde E}$.    

(ii) This follows from Proposition \ref{pr-prop-convhull2} since 
the convex hull of $\bigcup S(\mbv_{\tilde E})$  contains 
a generalized lens with the right properties.



(iii) This was proved in Proposition \ref{pr-lem-concavebd}.  





(iv) The existence of a lens-shaped p-end neighborhood of $\tilde S_{\tilde E}$ follows from Theorem \ref{du-thm-lensn}. 

(v) We choose mutually disjoint end neighborhoods for all ends. 
Then we choose the desired ones as above. 
\hfill \SnT {\parfillskip0pt\par}
\end{proof}

\subsection{The mc-p-end neighborhoods}

The mc-p-end neighborhoods will be useful in other papers. 

\index{end!mc-p-end neighborhood}
\begin{definition}\label{app-defn-lambda}
 Let $\tilde E$ be a lens-shaped R-end of 
a convex projective $n$-orbifold $\orb$ with a projective universal cover 
 $\torb \subset \SI^n$ (resp. $\RP^n$). 
 Let $\CH(\Lambda(\tilde E))$ denote the convex hull of $\Lambda(\tilde E)$. 
Let $U'$ be any p-end neighborhood of $\tilde E$ containing $\CH(\Lambda(\tilde E)) \cap \torb$. 
We define a {\em maximal concave p-end neighborhood} (or {\em mc-p-end neighborhood}) $U$ 
 to be one of the components of $U' - \CH(\Lambda(\tilde E))$ containing 
a p-end neighborhood of $\tilde E$. 
The {\em closed maximal concave p-end neighborhood} is 
$\clo(U) \cap \torb$. 
An $\eps$-$d_{\torb}$-neighborhood $U''$ of a maximal concave p-end neighborhood is called 
an {\em $\eps$-mc-p-end neighborhood}.  \index{end!emc@$\eps$-mc-p-end neighborhood}
\end{definition} 
In fact, these are independent of the choice of $U'$. 
Note that a maximal concave p-end neighborhood $U$ is uniquely determined since so is $\Lambda(\tilde E)$. 


Each radial segment $s$ in $\torb$ from $\mbv_{\tilde E}$ meets $\Bd U \cap \torb$ at a unique point
since the point $s \cap \Bd U$ is in an $(n-1)$-dimensional ball 
$D = P \cap U$ for a hyperspace $P$
sharply supporting $\CH(\Lambda(\tilde E))$ 
with $\partial \clo(D) \subset \bigcup S(\mbv_{\tilde E})$. 

\begin{lemma} \label{app-lem-mcc}
Let $D$ be an $i$-dimensional totally geodesic compact convex domain with $i \geq 1$. 
Let $\tilde E$ be a generalized lens-shaped R-p-end with a p-end vertex ${\mbv_{\tilde E}}$. 
Suppose $\partial D \subset \bigcup S({\mbv_{\tilde E}})$. 
Then $D^o \subset V$ for 
a maximal concave p-end neighborhood $V$. Moreover, 
for sufficiently small $\eps>0$, an $\eps$-$d_{\torb}$-neighborhood of 
$D^o $ is contained in $V'$ for any $\eps$-mc-p-end neighborhood $V'$.
\end{lemma}
\begin{proof} 
	Suppose that $\torb \subset \SI^n$ first. 
Assume that $U$ is a generalized lens cone of ${\mbv_{\tilde E}}$. 
Then $\Lambda(\tilde E)$ is the set of endpoints of segments in 
$S({{\mbv_{\tilde E}}})$ 
which are not ${\mbv_{\tilde E}}$ by Theorems \ref{pr-thm-lensclass} and
\ref{pr-thm-redtot}. 
Let $P$ be the subspace spanned by $D \cup \{{\mbv_{\tilde E}}\}$. 
Since \[\partial D, \Lambda(\tilde E) \cap P  \subset \bigcup S({\mbv_{\tilde E}}) \cap P,\] 
and $\partial D \cap P$ is closer than $\Lambda(\tilde E) \cap P$ from ${\mbv_{\tilde E}}$, 
it follows that 
$P\cap \clo(U) -D$ has a component $C_1$ containing ${\mbv_{\tilde E}}$ and  
$\clo(P\cap \clo(U) - C_1)$ contains $\Lambda(\tilde E) \cap P$.
Hence, \[\clo(P\cap \clo(U) - C_1) \supset \CH(\Lambda(\tilde E)) \cap P\] by the convexity of $\clo(P \cap \clo(U) - C_1)$. 
Since $\CH(\Lambda(\tilde E))\cap P$ is a convex set in $P$, we have 
only two possibilities\/{\rm :}
\begin{itemize}
\item $D$ is disjoint from $\CH(\Lambda(\tilde E))^o$ or 
\item $D$ contains $\CH(\Lambda(\tilde E)) \cap P$.

\end{itemize} 
Let $V$ be an mc-p-end neighborhood of $U$. 
Since $\clo(V)$ includes the closure of the component of $U - \CH(\Lambda(\tilde E))$ with a limit point ${\mbv_{\tilde E}}$, 
it follows that $\clo(V)$ includes $D$. 

Since $D$ is in $\clo(V)$, the boundary $\Bd V' \cap \torb$ of the 
$\eps$-mc-p-end neighborhood $V'$ does not meet $D$. Hence $D^o \subset V'$. 
\hfill \SnT {\parfillskip0pt\par}
\end{proof} 


The following provides a characterization of $\eps$-mc-p-end neighborhoods of 
$\tilde E$. 

\begin{corollary} \label{app-cor-mcn} 
Let $\orb$ be a properly convex real projective $n$-orbifold with 
generalized lens-shaped or horospherical $\cR$- or $\cT$-ends 
satisfying  {\rm (}IE\/{\rm )}. 
Let $\tilde E$ be a generalized lens-shaped R-end. 
Then the following statements hold\/{\em :}  
\begin{enumerate}
\item[{\rm (i)}] A concave p-end neighborhood of $\tilde E$ is always a subset of an mc-p-end neighborhood of the same R-p-end. 
\item[{\rm (ii)}] The closed mc-p-end neighborhood of $\tilde E$ is 
the closure in $\torb$ of the union of all concave end neighborhoods of $\tilde E$.
\item[{\rm (iii)}] The mc-p-end neighborhood $V$ of $\tilde E$ is a 
proper p-end neighborhood, and covers an end neighborhood with compact boundary in $\orb$.
\item[{\rm (iv)}]  For sufficiently small $\eps > 0$, an $\eps$-mc-p-end neighborhood of $\tilde E$ is a proper p-end neighborhood. 
\item[{\rm (v)}] For sufficiently small $\eps> 0$, 
the image end neighborhoods in $\orb$ of $\eps$-mc-p-end neighborhoods of R-p-ends
are either mutually disjoint or identical. 
\end{enumerate}
\end{corollary} 
\begin{proof}
Suppose first that $\torb \subset \SI^n$. 
(i) Since the limit set $\Lambda(\tilde E)$ lies in every generalized cocompactly-acted lens $L$ by Corollary \ref{app-cor-independence}, 
a generalized lens cone p-end neighborhood $U$ of $\tilde E$ contains $\CH(\Lambda) \cap \torb$. 
Hence, a concave end neighborhood is contained in an mc-p-end neighborhood. 



(ii) 
Let $V$ be an mc-p-end neighborhood of $\tilde E$.
Then define $S$ to be the subspace of endpoints in $\clo(\torb)$ 
of maximal segments in $V$ from ${\mbv_{\tilde E}}$ in directions of $\tilde \Sigma_{\tilde E}$. 
Then $S$ is homeomorphic to $\tilde \Sigma_{\tilde E}$ by the map induced by radial segments
as shown in the previous paragraph. 
Thus, $S/\pi_1(\tilde E)$ is a compact set since $S$ is contractible and $\tilde \Sigma_{\tilde E}/\pi_1(\tilde E)$ is a $K(\pi_1(\tilde E))$-space
up to taking a torsion-free finite-index subgroup by Theorem \ref{prelim-thm-vgood} (Selberg's lemma). 
We can approximate $S$ in the $d_{\torb}$-sense 
by smooth convex boundary component $S_{\eps}$ outward
of a generalized cocompactly-acted lens
by Theorem \ref{du-thm-lensnpre}
since $\tilde E$ satisfies the uniform middle-eigenvalue condition. 
A component $U-S_{\eps}$ is a concave p-end neighborhood. 
(ii) follows from this. 

(iii) Since a concave p-end neighborhood is a proper p-end neighborhood by Theorems \ref{pr-thm-lensclass}(iv) and \ref{pr-thm-redtot}, 
and $V$ is a union of concave p-end neighborhoods, 
we obtain either
\[g(V) \cap V = \emp \hbox{ or } g(V) = V \hbox{ for } g \in \pi_1(\orb) \hbox{ by (ii).} \] 

Suppose that $g(\clo(V) \cap \torb) \cap \clo(V) \ne \emp$. Then $g(V) = V$. 
In that case, $g \in \pi_1(\tilde E)$:
Otherwise, $g(V) \cap V =\emp$, and $g(\clo(V) \cap \torb)$ meets $\clo(V)$ 
in a totally geodesic hypersurface $S$ equal to $\CH(\Lambda)^o$
by the concavity of $V$. 
Furthermore, for every $g \in \pi_1(\orb)$, $g(S) = S$, since $S$ is a maximal 
totally geodesic hypersurface in $\torb$.
Hence, $g(V) \cup S \cup V = \torb$ since these are subsets of a properly convex domain $\torb$, the boundary of $V$ and $g(V)$ are in $S$, 
and $S$ itself is now in the interior of $\torb$. 
Then $\pi_1(\orb)$ acts on $S$, and $S/G$ is homotopy 
equivalent to $\torb/G$ for a torsion-free  finite-index subgroup $G$ of $\pi_1(\orb)$ by Theorem \ref{prelim-thm-vgood} (Selberg's lemma). 
This contradicts condition {\rm (IE)}.

Hence, the only possibility is that
$\clo(V) \cap \torb = V \cup S$ for a hypersurface $S$ and  
\[g(V\cup S) \cap V \cup S = \emp \hbox{ or } g(V\cup S) = V \cup S \hbox{ for } g \in \pi_{1}(\orb).\]



Now suppose that $S \cap \Bd \torb \ne \emp$. 
Let $S'$ be a maximal totally geodesic domain in $\clo(V)$ that contains $S$. 
Then $S' \subset \Bd \torb$ by convexity and Lemma \ref{prelim-lem-simplexbd}, 
i.e., $S'=S \subset \Bd \torb$.  
In this case, $\torb$ is a cone over $S$ with the end vertex ${\mbv_{\tilde E}}$ of $\tilde E$.
For each $g \in \pi_1(\orb)$, 
\[g(V) \cap V \ne \emp \hbox{ implies } g(V)=V\] since $g({\mbv_{\tilde E}})$ is on $\clo(S)$. 
Thus, $\pi_1(\orb) = \pi_1(\tilde E)$. 
This contradicts the condition (IE) for $\pi_1(\tilde E)$. 

We showed that $\clo(V) \cap \torb = V \cup S$ for a hypersurface $S$ 
and covers a submanifold in $\orb$ which is a closure of 
an end neighborhood covered by $V$. 
Thus, an mc-p-end neighborhood $\clo(V) \cap \torb$ is a proper end neighborhood of $\tilde E$
with compact embedded boundary $S/\pi_1(\tilde E)$. 

(iv) Obviously,
we can choose positive $\eps$ such that an $\eps$-mc-p-end neighborhood is also a proper p-end neighborhood.

(v) For two mc-p-end neighborhoods $U$ and $V$ for different R-p-ends, 
we have $U \cap V =\emp$
by reasoning as in (iii) replacing $g(V)$ with $U$: 
We showed that $\clo(V) \cap \torb$ for an mc-p-end neighborhood $V$ covers an end neighborhood in $\orb$. 


Suppose that $U$ is an mc-p-end neighborhood different from $V$.
We claim that $\clo(U) \cap \clo(V) \cap \torb = \emp$:
Suppose not.  For $g \not\in \bGamma_{\tilde E}$, $g(\clo(V))$
must be a subset of $U$ since otherwise we have 
a situation of (iii) for $V$ and $g(V)$. 
Since the preimages of the end neighborhoods are disjoint, 
$g(V)$ is a p-end neighborhood of the same end as $U$.
Since both are $\eps$-mc-p-end neighborhoods which are 
canonically defined, we obtain $U = g(V)$. This was ruled out in (iii). 
\hfill \SnT {\parfillskip0pt\par}
\end{proof}


\section{The strong irreducibility of the real projective orbifolds.}\label{app-sec-strirr}

The main purpose of this section is to prove 
Theorem \ref{intro-thm-sSPC}, the strong irreducibility result. 
But we shall first discuss convex hull of the ends. 
We show that the closures of convex hulls of p-end neighborhoods are disjoint in $\Bd \torb$. 
The fact that there are infinitely many of these will show the strong irreducibility.


For the following, we need a stronger condition of lens-shaped ends,
and not just the generalized lens-shaped condition, 
to obtain the disjointedness of the closures of p-end neighborhoods. 
For now, we cannot do the following for the generalized lens. 

\begin{corollary} \label{app-cor-disjclosure} 
Let $\orb$ be a strongly tame  properly convex real projective $n$-orbifold with
lens-shaped or horospherical $\cR$- or $\cT$-ends 
and satisfy  {\rm (}NA\/{\rm )}. 
Let $\mathcal U$ be the collection of the components of the inverse image in $\torb$ 
of the union of the disjoint collection of the end neighborhoods of $\orb$. 
Now for each concave p-end neighborhood, replace it by 
R-p-end neighborhood of the collection $\mathcal U$ 
using Corollary \ref{app-cor-shrink} {\rm (iii).} 
Then the following statements hold\,{\em :} 
\begin{enumerate} 
\item[{\rm (i)}] Given horospherical, concave, or one-sided lens p-end neighborhoods $U_1$ and $U_2$ contained in $\bigcup \mathcal U$, 
we have $U_1 \cap U_2 =\emp$ or $U_1= U_2$. 
\item[{\rm (ii)}] Let $U_1$ and $U_2$ be in $\mathcal U$. Then either 
$\clo(U_1) \cap \clo(U_2) \cap  \Bd \tilde{\mathcal{O}} = \emp$ 
or $U_1 = U_2$ holds. 
\end{enumerate}
\end{corollary}
\begin{proof} 
	Suppose first that $\torb \subset \SI^n$. Suppose that $\orb$ has ends $E_1, \dots, E_m$.
Since the neighborhoods in $\mathcal U$ are mutually disjoint, 
\begin{itemize}
\item $\clo(U''_1) \cap \clo(U''_2) \cap \torb = \emp$ or
\item For each pair $U''_1, U''_2 \in \mathcal U$, we have $U''_1 = U''_2$. 
\end{itemize} 
%
%
(i) Assume, without loss of generality, that $E_1$ and $E_2$ are R-ends. 
Suppose that $U_1$ and $U_2$ are concave p-end neighborhoods of R-p-ends
$\tilde E_1$ and $\tilde E_2$ respectively. 
 Let $U'_1$ be the interior of the associated generalized lens cone of $U_1$ in $\clo(\torb)$, 
and let $U'_2$ be that of $U_2$. 
Let $U''_i$ be the concave p-end neighborhood of $U'_i$ for $i=1,2$ 
by Corollary \ref{app-cor-shrink} (iii)
that covers a respective end neighborhood in $\orb$.


Assume that $U''_i \in {\mathcal{U}}$, $i=1, 2$, and $U''_1 \ne U''_2$. 
Suppose that the closures of $U''_1$ and $U''_2$ intersect in $\Bd \torb$.
Suppose that they are R-p-end neighborhoods. 
Then 
the respective closures of convex hulls $I_1$ and $I_2$ as obtained by Proposition \ref{app-prop-I} intersect as well. 
Take a point $z \in \clo(U''_1) \cap \clo(U''_2) \cap \Bd \torb$. 
Suppose that we choose $p_i = \mbv_{\tilde E_i}$, $i=1,2$, for each p-end $\tilde E_i$. 

Suppose that $\ovl{p_1p_2}^o \subset \Bd \torb$. 
Then there exists a segment in $\Bd \torb$ from 
$\mbv_{\tilde E_1}$ to $\mbv_{\tilde E_2}$. 
By Theorems \ref{pr-thm-lensclass} and \ref{pr-thm-redtot}, 
these vertices are in the closures of p-end neighborhoods 
of one another. 
Since $\clo(U''_1)$ and $\clo(U''_2)$ contain some open metric ball neighborhoods of $p_1$ and $p_2$ respectively for the metric $\bdd$ restricted to $\clo(\torb)$, and $U''_j$ is dense in $\clo(U''_j)$, $j=1,2$, 
we obtain $U''_1\cap U''_2 \ne \emp$. 
This is a contradiction. 

Hence, $\ovl{p_1p_2}^o \subset \torb$ holds. 
Then $\ovl{p_1z} \subset S(\mbv_{\tilde E_2})$, and 
$\ovl{p_2z} \subset S(\mbv_{\tilde E_1})$ and these segments
are maximal since otherwise $U''_1 \cap U''_2 \ne \emp$.
The segments intersect transversely at $z$ 
since otherwise we violate the maximality in Theorems \ref{pr-thm-lensclass} and 
\ref{pr-thm-redtot}.
We obtain a triangle $\tri(p_1p_2z)$ in $\clo(\torb)$ with vertices $p_1, p_2, z$. 

There is a lens $L_i$ for each $\tilde E_i$ such that $L_i \ast p_i$ is a p-end neighborhood of 
$\tilde E_i$. 
Proposition \ref{pr-cor-tangency} contradicts the existence of the above triangle dues to 
the uniqueness of the supporting hyperspace at $z$, which is in both $\clo(L_i) - L_i$ for $i=1,2$.

Now, consider when $U_1$ is a one-sided lens neighborhood of a T-p-end and 
$U_2$ is a concave R-p-end neighborhood of an R-p-end of $\torb$.
Let $z$ be the intersection point in $\clo(U_1) \cap \clo(U_2)$. 
Suppose that the hyperspace containing the ideal boundary component of $U_1$ is transversal to the line 
$\ovl{z\mbv_{\tilde E_2}}$. 
We can use the same reasoning as above by choosing any $p_1$ in $\tilde \Sigma_{\tilde E_1}$ 
such that $\ovl{p_1z}$ passes the interior of $\tilde E_1$. Let $p_2$ be the R-p-end vertex of $U_2$. 
We again obtain the triangle with vertices $p_1, p_2$, and $z$ as above. 
Proposition \ref{pr-cor-tangency} again gives us a contradiction. 

Suppose that the hyperspace containing the ideal boundary component of $U_1$ is tangent to the line
$\ovl{z\bv_{\tilde E_2}}$. 
We take the convex hull of $z$ with the convex domain $\tilde \Sigma_{\tilde E_2}$ in
the hyperspace $S$ containing both. The convex hull is in $\clo(\torb)\cap S$ where $\pi_1(\tilde E_2)$ acts on.  The uniqueness part of Theorem \ref{prelim-thm-Kobayashi} shows that the interior of the convex hull is still $\tilde \Sigma_{\tilde E_2}$. This is a contradiction.  


Next, consider the case where $U_1$ and $U_2$ both are one-sided  lens neighborhoods of T-p-ends.  If the hyperspaces $S_i$ containing the ideal boundary component $\tilde \Sigma_{\tilde E_i}$ for $i=1, 2$, are identical, then 
we can again take the convex hull as in the above paragraph and obtain a contradiction. 

Suppose that these hyperspaces $S_1$ and $S_2$ are transversal. 
Let $z$ be a  point of intersection of $\clo(\Omega_1)$ and $\clo(\Omega_2)$. 
Then $S_i$ must be an asymptotic hyperspace of $\Omega_j$ at $z$ for $i\ne j, i, j=1, 2$.

We can see that $\clo(\torb)$ meets $S_i$ at $\clo(\tilde \Sigma_{\tilde E_i})$ again by
the uniqueness part of Theorem \ref{prelim-thm-Kobayashi}. Hence, 
$(\pi_1(\tilde E_i), \tilde \Sigma_{\tilde E_i}, \torb^o)$ is a properly convex triple. 

Let $H_i$ denote the open half-space containing $\torb^o$ bounded by $S_i$. 
Let $B =H_1 \cap H_2$ be the convex domain containing $\torb^o$.

Let $L = S_1 \cap S_2$. Then $L$ is sharply supporting $\Omega_i$ at $z$ in $S_i$. 
Let $A_i$ denote the asymptotic hyperspace to $U_i$ at $z$ containing $L$. 
By using Lemma \ref{du-lem-hdisj} as in the proof of Theorem \ref{du-thm-ASunique}, 
we conclude that $\torb^o$ cannot intersect the supporting hyperspaces $A_1$ and $A_2$. 

Since $A_i$ is asymptotic to $U_i$, $S_j$, $j\ne i$, 
it follows that $S_i \cap B$ cannot separate $B$. 
Since $S_i$ must meet the boundary of $\torb$, 
the only possibility is $S_1 = A_2$ and $S_2 = A_1$. 
This proves our claim. 

This means that $L_1^o \cap L_2^o \ne \emp$ by the asymptoticity of $A_1$ and $A_2$. 
Moreover, for any choice of lenses $L'_1$ and $L'_2$ their interiors meet.
This contradicts the notion that there are some p-end neighborhoods of lens type in any proper end neighborhoods according to Theorem \ref{du-thm-lensn}. 


We finally consider when $U$ is a horospherical R-p-end. Since $\clo(U) \cap \Bd \torb$ is a unique point, Lemma \ref{prelim-lem-cuspsegment}.  
implies the result.
\hfill \SnT {\parfillskip0pt\par}
\end{proof}





We modify Theorem \ref{pr-thm-redtot} 
by replacing some conditions. 
In particular, we do not assume that $h(\pi_{1}(\orb))$ is strongly irreducible. 
\begin{lemma}\label{app-lem-redtot2}
Let $\orb$ be a strongly tame 
properly convex real projective $n$-orbifold  
satisfying {\rm (}IE\/{\rm )}. 
Let $\tilde E$ be a virtually factorizable R-p-end of the universal cover 
$\torb$ of generalized lens-shaped.
Then the following hold\/{\rm :}
\begin{itemize}
\item there exists a totally geodesic hyperspace $P$ on which $h(\pi_1(\tilde E))$ acts, 
\item $D:=P \cap \torb$ is a properly convex domain, 
\item $D^{o} \subset \torb$, and
\item $D^{o}/\pi_1(\tilde E)$ is a compact orbifold. 
\end{itemize} 
\end{lemma}
\begin{proof}
	Assume first that $\torb \subset \SI^n$. 
The proof of Theorem \ref{pr-thm-redtot}   shows that  
\begin{itemize}
\item either $\clo(\torb)$ is a strict join $\{\mbv_{\tilde E}\} \ast D'$ for 
a properly convex domain $D$ in a hyperspace, or 
\item the conclusion of Theorem \ref{pr-thm-redtot}   holds.

\end{itemize} 
In both cases, $\pi_1(\tilde E)$ acts on a compact, totally geodesic convex domain $D$ of codimension $1$.
$D$ is the intersection $P_{\tilde E} \cap \clo(\torb)$ for a $\pi_1(\tilde E)$-invariant subspace $P_{\tilde E}$. 
Suppose that $D^{o}$ is not a subset of $\torb$. Then by Lemma \ref{prelim-lem-simplexbd}, 
$D\subset \Bd \torb$. 

In the former case, $\clo(\torb)$ is the join ${\mbv_{\tilde E}} \ast D$. 
For each $g \in \pi_{1}(\tilde E)$ satisfying $g(\{\mbv_{\tilde E}\}) \ne {\mbv_{\tilde E}}$, we have $g(D) \ne D$
since $g({\mbv_{\tilde E}})\ast g(D) = \{\mbv_{\tilde E}\}\ast D$. 
The intersection 
$g(D) \cap D$ is a proper compact convex subset
of $D$ and $g(D)$. 
Moreover, 
\[\clo(\torb) =\{\mbv_{\tilde E}\}\ast g(\{\mbv_{\tilde E}\}) \ast (D\cap g(D)).\]
We can continue as many times as there is a mutually distinct collection of vertices of the form 
$g({\mbv_{\tilde E}})$. This process cannot continue indefinitely, 
We have a contradiction since 
by Condition (IE), there are infinitely many distinct end vertices of the form $g({\mbv_{\tilde E}})$ for $g \in \pi_{1}(\orb)$. 

Now, we only have the alternative $D^{o} \subset \torb$
where $D^o/\bGamma_{\tilde E}$ is projectively diffeomorphic to
$\tilde \Sigma_{\tilde E}/\bGamma_{\tilde E}$. 
\hfill \SSn {\parfillskip0pt\par}
\end{proof}

\begin{proof}[{\sl Proof of Theorem \ref{intro-thm-sSPC}}]
	It is sufficient to prove for the assumption $\torb \subset \SI^n$. 
Let $h:\pi_1(\orb) \ra \SLnp$ be the holonomy homomorphism. 
Suppose that $h(\pi_1(\orb))$ is virtually reducible. Then we can choose a finite cover 
$\orb_1$ such that $h(\pi_1(\orb_1))$ is reducible
since it is sufficient to prove for the finite-index groups.  

Each end holonomy group $h(\pi_1(\tilde E)$ for R-p-end $\tilde E$ is reducible. 
Suppose that it is not factorizable. Then it acts on a hyperplane $P$ disjoint from 
$\mbv_{\tilde E}$. 
The interior of $\mbv_{\tilde E} \ast \Omega$ for a stricly convex domain 
$\Omega$ in $P$ must be in $\torb$. We also know that $\partial \Omega$ is $C^1$ 
by \cite{Benzecri60}. 
This interior cannot be $\torb$ since this will imply that $\mbv_{\tilde E}$ is fixed 
by $h(\pi_1(\orb))$ violating (IE). 
Therefore, $\Omega \subset \torb$, and $\tilde E$ is a totally geodesic R-p-end. 

Theroem \ref{pr-thm-equ} shows that $h(\pi_1(\tilde E))$ satisfies the uniform middle eigenvalue
condition with respect to $\mbv_{\tilde E}$. 
$\clo(\Omega)$ is a distanced compact convex set. 
Proposition \ref{pr-prop-convhull2} 
 implies that $\tilde E$ is a lens-shaped R-p-end.

Since $\tilde E$ can only be factorizable otherwise, all generalized lens $\cR$-ends are lens $\cR$-ends
by Theorem \ref{pr-thm-redtot}. 

We may assume that $\pi_1(\orb)$ is torsion-free 
by taking a finite cover by Theorem \ref{prelim-thm-vgood}.

We denote $\orb_1$ by $\orb$ for simplicity. 
Let $S$ denote a proper subspace where $\pi_1(\orb)$ acts on. 
Suppose that $S$ meets $\torb$. 
Then $\pi_1(\tilde E)$ acts on a properly convex open domain $S\cap \torb$ for each p-end 
$\tilde E$. 
Then $S \cap \torb$ for any p-end neighborhood gives a suborbifold of 
a closed end orbifold homotopy equivalent to it. 
Thus, $(S\cap \torb)/\pi_1(\tilde E)$ is a compact orbifold homotopy equivalent to one of the end orbifold, and $S$ must be of codimension one. 
However, $S \cap \torb$ is $\pi_1(\tilde E')$-invariant and cocompact 
for any other p-end $\tilde E'$. 
Hence, each p-end fundamental group $\pi_{1}(\tilde E)$ is virtually identical to
any other p-end fundamental group. 
This contradicts (NA). 
Therefore, 
\begin{equation}\label{app-eqn-K}
K:=S \cap \clo(\torb) \subset \Bd \torb,
\end{equation} 
where $g(K)=K$ for every $g \in h(\pi_1(\orb))$.

We divide into steps:
\begin{enumerate}
	\item[(A)] First, we show $K\ne \emp$. 
	\item[(B)] We show $K= D_j$ or $K = \{\mbv_{\tilde E}\}\ast D_j$ 
where $D_j$ is a properly convex domain in $\Bd \torb \cap \clo(U)$ for a p-end neighborhood $U$ of $\tilde E$. 
	\item[(C)] Finally, we show that 
there exists a finite-index subgroup $\bGamma'$ of 
$\bGamma$ such that $g(D_j) = D_j$ for every $g \in \bGamma'$, 
	and we use Corollary \ref{app-cor-disjclosure} to obtain a contradiction.
	 
\end{enumerate}

(A) We show that $K:= \clo(\torb) \cap S \ne \emp$: 
Let $\tilde E$ be a p-end. If $\tilde E$ is horospherical, $\pi(\tilde E)$ acts on a great sphere $\hat S$ tangent 
to an end vertex. 
Since $S$ is $\bGamma$-invariant, 
$S$ must be a subspace in $\hat S$ containing the end vertex 
by Theorem \ref{ce-thm-affinehoro}(iii). 
This implies that every horospherical p-end vertex is in $S$. 
Let $p$ be one. 
Since there is no nontrivial segment in $\Bd \torb$ containing
$p$ by Theorem \ref{ce-thm-affinehoro}(v), 
$p$ equals $S \cap \clo(\torb)$. 
Hence, $p$ is $\bGamma$-invariant and $\bGamma = \bGamma_{\tilde E}$. 
This contradicts Condition (IE). 

Suppose that $\tilde E$ is a generalized lens-shaped R-p-end. Then, by the existence of 
attracting subspaces of some elements of $\bGamma_{\tilde E}$, we have
\begin{itemize}
\item either $S$ passes the end vertex $\mbv_{\tilde E}$ or 
\item there exists a subspace $S'$ containing $S$ and $\mbv_{\tilde E}$ that is $\bGamma_{\tilde E}$-invariant. 

\end{itemize} 
In the first case, we have $S \cap \clo(\torb) \ne \emp$, and we completed Step (A).

In the second case, 
$S'$ corresponds to a subspace of $\SI^{n-1}_{\mbv_{\tilde E}}$
and $S$ is a hyperspace of dimension $\leq n-1$ disjoint from $\mbv_{\tilde E}$. 
Thus, $\tilde E$ is a virtually factorizable R-p-end. 
By Theorem \ref{prelim-thm-semi},
we obtain some attracting fixed points 
in the limit sets of $\pi_1(\tilde E)$. 
If $S'$ is a proper subspace, then $\tilde E$ is factorizable, 
and $S'$ contains the 
attracting fixed set of some positive bi-semiproximal $g$, $g \in \bGamma_{E}$. 
The uniform middle-eigenvalue condition shows that 
positive bi-semiproximal $g$ has attracting fixed points in $\clo(L)$.  
Since $g$ acts on $S$, 
we obtain $S\cap \clo(L)\ne \emp$
by the uniform middle-eigenvalue condition. 

If $S'$ is not a proper subspace, then $g$ acts on $S$. 
Moreover, 
$S$ contains attracting fixed points of $g$ by the uniform middle-eigenvalue condition. 
Thus, $S \cap \clo(L)\ne \emp$. 


If $\tilde E$ is a lens-shaped T-p-end, we can apply a similar argument using the attracting fixed sets. 
Therefore, $S \cap \clo(\torb)$ is a subset $K$ of $\Bd \torb$ 
with $\dim K \geq 0$ and $K$ is not empty. 
In fact, we have shown that the closure of each p-end neighborhood meets $K$. 




(B) Taking a dual orbifold if necessary, 
we assume without loss of generality that there exists a generalized lens-shaped R-p-end $\tilde E$ 
with a radial p-end vertex $\mbv_{\tilde E}$. 

As in (A), suppose that $\mbv_{\tilde E} \in K$. 
There exists $g \in \pi_{1}(\orb)$,
such that 
\[g(\mbv_{\tilde E}) \ne \mbv_{\tilde E}, 
\hbox{ and } g(\mbv_{\tilde E}) \in K \subset \Bd \torb\]
since $g$ acts on $K$. The point
$g(\mbv_{\tilde E})$ is outside the closure of
the concave p-end neighborhood of $\tilde E$ by Corollary
\ref{app-cor-disjclosure}. Since $K$ is connected, $K$ meets $\clo(L)$ for the CA-lens or generalized
cocompactly-acted lens $L$ of $\tilde E$. 

If $\mbv_{\tilde E} \not\in K$, then again $K \cap \clo(L) \ne \emp$ as in (A) using attracting fixed sets of 
some elements of $\pi_{1}(\tilde E)$. 
Hence, we conclude $K \cap \clo(L) \ne \emp$ for a 
generalized cocompactly-acted lens $L$ of $\tilde E$. 

Let $\Sigma'_{\tilde E}$ denote $D^{o}$ from Lemma \ref{app-lem-redtot2}. 
Since $K \subset \Bd \orb$, it follows that
$K$ cannot contain $\Sigma'_{\tilde E}$. 
Thus, $K \cap \clo(\Sigma'_{\tilde E})$ is a proper subspace of $\clo(\Sigma'_{\tilde E})$, and 
$\tilde E$ must be a virtually factorizable end. 

By Lemma \ref{app-lem-redtot2}, 
there exists a totally geodesic domain $\Sigma'_{\tilde E}$ in the cocompactly-acted lens. 
A p-end neighborhood of $\mbv_{\tilde E}$ equals 
$U_{\mbv_{\tilde E}}:=(\{\mbv_{\tilde E}\} \ast \Sigma'_{\tilde E})^{o}$. 
Since $\pi_{1}(\tilde E)$ acts reducibly, 
\[\clo(\Sigma'_{\tilde E}) = D_{1}\ast \cdots \ast D_{m},\] 
where 
$K \cap \clo(U_{\mbv_{\tilde E}})$ contains a join $D_{J}:= \ast_{i\in J}D_{i}$ for a proper subcollection 
$J$ of $\{1, \dots, m\}$. Moreover, $K \cap \clo(\Sigma'_{\tilde E}) = D_{J}$. 

 Since $g(U_{\mbv_{\tilde E}})$ is a p-end neighborhood of $g(\mbv_{\tilde E})$, we obtain $g(U_{\mbv_{\tilde E}}) = U_{g(\mbv_{\tilde E})}$.
Since $g(K) = K$ for all $g \in \Gamma$, we obtain 
\[K \cap g(\clo(\Sigma'_{\tilde E}))  = g(D_{J}).\] 

Lemma \ref{app-lem-redtot2} implies that 
\begin{align} \label{app-eqn-Uv}
U_{g(\mbv_{\tilde E})} \cap U_{\mbv_{\tilde E}} = \emp \hbox{ for } g \not\in \pi_{1}(\tilde E) \hbox{ or } \nonumber \\ 
U_{g(\mbv_{\tilde E})} = U_{\mbv_{\tilde E}} \hbox{ for } g \in \pi_{1}(\tilde E) 
\end{align} 
by the similar properties of $S(g(\mbv_{\tilde E}))$ and $S(\mbv_{\tilde E})$ and the fact that
$\Bd U_{\mbv_{\tilde E}} \cap \torb$ and $\Bd U_{g(\mbv_{\tilde E})}\cap \torb$ are 
totally geodesic domains. 

Let $g$ be normalized to have determinant $\pm 1$. 
Let $\lambda_{J}(g)$ denote the $(\dim D_{J}+1)$-th root of the norm of the determinant of the submatrix of $g$ associated with $D_{J}$.
By Proposition 4.4 of \cite{Benoist03}, there exists 
a sequence of virtually central diagonalizable 
elements $\gamma_{i}\in \pi_{1}(\tilde E)$ such that 
\[\{\gamma_{i}| D_{J}\} \ra \Idd, \{\gamma_{i}| D_{J^{c}}\} \ra \Idd \hbox{ satisfying }
\left\{\frac{\lambda_{J}(\gamma_{i})}{\lambda_{J^{c}}(\gamma_{i})}\right\} \ra \infty \]
for the complement $J^{c}:= \{1, 2, \dots, m\} - J$. 
Since the lens-shaped ends satisfy the uniform middle-eigenvalue condition by Theorem \ref{pr-thm-redtot}, we obtain
\begin{align} 
\{\gamma_{i}| D_{J}\} \ra \Idd, \{\gamma_{i}| D_{J^{c}}\} \ra \Idd \hbox{ for the complement } J^{c}:= \{1, 2, \dots, m\} - J, \nonumber \\
\left\{\frac{\lambda_{J}(\gamma_{i})}{\lambda_{\mbv_{\tilde E}}(\gamma_{i})}\right\} \ra \infty, \left\{\frac{\lambda_{J^{c}}(\gamma_{i})}{\lambda_{\mbv_{\tilde E}}(\gamma_{i})}\right\} \ra 0, 
\left\{\frac{\lambda_{J}(\gamma_{i})}{\lambda_{J^{c}}(\gamma_{i})}\right\} \ra \infty. 
\end{align}  
(See Theorem \ref{prelim-prop-Ben2}.)

Since $D_{J}\subset K$, the uniform middle-eigenvalue condition 
for $\bGamma_{\tilde E}$ implies that one of the following holds: 
\[K = D_{J}, K = \{\mbv_{\tilde E}\} \ast D_{J} \hbox{ or } K = \{\mbv_{\tilde E}\}\ast D_{J} \cup \{\mbv_{\tilde E-}\} \ast D_{J}\]
by the invariance of $K$ under $\gamma_{i}^{-1}$,
and the fact that $K \cap \clo(\Sigma'_{\tilde E}) = D_{J}$. 
Since $K \subset \clo(\torb)$, the third case is not possible
by the proper convexity of $\clo(\torb)$. 
We obtain 
\begin{equation} \label{app-eqn-twoDJ}
K = D_{J} \hbox{ or } K= \{\mbv_{\tilde E}\} \ast D_{J}.
\end{equation}

(C) We now explore the two cases of \eqref{app-eqn-twoDJ}. 
Assume that the second case holds. 
Let $g$ be an arbitrary element of $\pi_{1}(\orb) - \pi_{1}(\tilde E)$. 
Since $D_{J} \subset K$, we obtain $g(D_{J}) \subset K$. 
Recall that $U_{\mbv_{\tilde E}} \cup S(\mbv_{\tilde E})^{o}$ is a neighborhood of the points of $S(\mbv_{\tilde E})^{o}$ in $\clo(\torb)$.
Thus, $g(U_{\mbv_{\tilde E}} \cup S(\mbv_{\tilde E})^{o})$ is a neighborhood of the points of $g(S(\mbv_{\tilde E})^{o})$. 

Recall that $D_{J}^{o}$ lies in the closure of $U_{\mbv_{\tilde E}}$. 
If $D_{J}^{o}$ meets 
\[g(\{\mbv_{\tilde E}\} \ast D_{J} - D_{J}) \subset g(U_{\mbv_{\tilde E}}\cup S(\mbv_{\tilde E})^{o}) \supset g(S(\mbv_{\tilde E})^{o}),\] 
then 
\[U_{\mbv_{\tilde E}} \cap g(U_{\mbv_{\tilde E}}) \ne \emp, 
\hbox{ and } S(\mbv_{\tilde E})^{o} \cap g(S(\mbv_{\tilde E})^{o}) \ne \emp\] 
since 
these are components of $\torb$ with certain totally geodesic hyperspaces removed.
Hence, $\mbv_{\tilde E} = g(\mbv_{\tilde E})$ by Theorems \ref{pr-thm-lensclass} 
and \ref{pr-thm-redtot}. Finally, we obtain $D_{J} = g(D_{J})$ since
\[K = \{\mbv_{\tilde E}\} \ast D_{J} = g(\{\mbv_{\tilde E}\}) \ast g(D_{J}).\] 

If $D_{J}^{o}$ is disjoint from $g(\{\mbv_{\tilde E}\} \ast D_{J} - D_{J})$, then $g(D_{J}) \subset D_{J}$
since $K = \{\mbv_{\tilde E}\} \ast D_{J}$ and $g(K) = K$. 
Since both $D_{J}$ and $g(D_{J})$ arise as intersections of a hyperspace with $\Bd \torb$,  
we obtain $g(D_{J}) = D_{J}$. 

Both cases of \eqref{app-eqn-twoDJ} imply that
 $g(D_{J}) = D_{J}$ for $g \in \pi_{1}(\orb)$.
Since $\mbv_{\tilde E}$ and $g(\mbv_{\tilde E})$ are not equal for $g \in \pi_{1}(\orb) - \pi_{1}(\tilde E)$, we obtain 
\[ \clo(U_1) \cap g(\clo(U_1)) \ne \emp.\]
%
Since all virtually reducible R-ends are totally geodesic, 
all R-ends are of lens-type ones according to Theorem \ref{pr-thm-redtot}.  
Corollary \ref{app-cor-disjclosure} gives us a contradiction. 
Therefore, we deduced that the $h(\pi_{1}(\orb))$-invariant subspace $S$ does not exist. 

Since parabolic subgroups of $\PGL(n+1, \bR)$ or 
$\SL_\pm(n+1, \bR)$ are reducible, the proof is complete. 
\hfill \SSn  {\parfillskip0pt\par}

\end{proof}

\subsection{Equivalence of lens-ends and generalized lens-ends for strict SPC-orbifolds}

\begin{corollary} \label{app-cor-stLens}
	Suppose that $\mathcal O$ is an $n$-dimensional strongly tame strictly SPC-orbifold
	with generalized lens-shaped or horospherical $\cR$- or $\cT$-ends,
	and satisfying the conditions  {\rm (}IE\/{\rm )} and {\rm (}NA\/{\rm )}. 
	Then $\orb$ satisfies the triangle condition, and 
	every generalized lens-shaped R-ends are lens-shaped R-ends. 
\end{corollary}
\begin{proof} 
	Assume first $\torb \subset \SI^n$. 
	Let $\tilde E$ be a generalized lens-shaped p-end neighborhood of $\torb$.
	Let $L$ be the generalized cocompactly-acted lens such that the interior $U$ of 
	$ \{\mbv_{\tilde E}\}\ast L$ is 
	a lens p-end neighborhood. Then $U- L$ is a concave p-end neighborhood. 
	Recall the triangle condition of Definition \ref{pr-defn-tri}. 
	Let $T$ be a triangle with 
	\[\partial T \subset \Bd \torb, T^o \subset \torb
	\hbox{ and }\partial T \cap \clo(U)\ne \emp\] for an R-p-end neighborhood $U$. 
	By the strict convexity $\torb$, each edge of $T$ 
	must be inside a set of the form 
	$\clo(V) \cap \Bd \torb$ for a p-end neighborhood $V$. 
	Corollary \ref{app-cor-disjclosure} implies that all edges lie in 
	$\clo(U) \cap \Bd \torb$ for a single R-p-end neighborhood $U$. 
	Hence, the triangle condition is satisfied. 
	By Theorem \ref{pr-thm-equ}, $\tilde E$ is a lens-shaped p-end.  

	\hfill \SnT {\parfillskip0pt\par}
\end{proof} 











\newcommand{\Sio}{\SI^{i_0}_\infty}

\hyphenation{neighbor-hood}










\setcounter{tocdepth}{3} 






\chapter[The NPNC-ends]{The convex but nonproperly convex and non-complete-affine radial ends} \label{ch-np}

In previous chapters, we classified 
properly convex or complete radial ends under suitable conditions. 
In this chapter, we will study radial ends that are convex but neither properly convex nor complete affine. 
The main techniques are the theory of Fried and Goldman on affine manifolds, and a generalization of 
the work on Riemannian foliations by Molino, Carri\`ere, and so on. 
We will show that these are quasi-joins of horospheres and totally geodesic radial ends under the transverse weak middle-eigenvalue conditions (see Definition \ref{np-defn-weakmec}).
These are suitable deformations of joins of horospheres and totally geodesic radial ends. 
Since this is the most technical chapter, we will provide outlines in places
in addition to the main outline in Section \ref{np-sub-outline}. 
 







\section{Introduction} \label{np-sec-intro} 



\subsection{General setting} \label{np-sub-general} 

In this chapter, we work with $\SI^n$ and $\SL_{\pm}(n+1, \bR)$
with only a few exceptions
since our purpose is to classify some objects modulo projective automorphisms. 
However, the corresponding results for $\RP^n$ can be obtained 
easily from results in Section \ref{prelim-sub-lifting} 
and then projecting back to $\RP^n$. 
Let $\tilde E$ be an R-p-end of a convex real projective $n$-orbifold $\orb$ 
with end orbifold $\Sigma_{\tilde E}$ and its universal cover 
$\tilde \Sigma_{\tilde E}$ and the p-end vertex $\mbv_{\tilde E}$. 
 We recall Proposition \ref{prelim-prop-classconv}.
If $\tilde \Sigma_{\tilde E}$ is convex but neither properly convex 
and not complete affine, then we call $E$ a  nonproperly convex and noncomplete end (NPNC-end.)
\index{end!nonproperly convex and noncomplete end}
\index{end!NPNC}
The closure $\clo(\tilde \Sigma_{\tilde E})$ contains 
a great $(i_0-1)$-dimensional sphere 
$\SI^{i_0-1}_\infty$, for $0< i_0 <n-1$, 
and the convex open domain $\tilde \Sigma_{\tilde E}$ is foliated by $i_0$-dimensional hemispheres 
with this boundary $\SI^{i_0-1}_\infty$. 
(These follow from Section 1.4 of \cite{ChCh93}. See also \cite{GV58}.) 

The space of $i_0$-dimensional hemispheres in $\SI^{n-1}_{\mbv_{\tilde E}}$ with boundary $\SI^{i_0-1}_\infty$ forms 
a projective sphere $\SI^{n-i_0-1}$: 
This follows since a complementary subspace $S$ isomorphic to 
$\SI^{n-i_0-1}$ parameterizes the space by the intersection points
with $S$.
The fibration with fibers open hemispheres of dimension $i_0$ with boundary 
$\SI^{i_0-1}_\infty$ 
\begin{alignat}{2} \label{np-eqn-pik}
	\hat \Pi_{K}:\SI^{n-1}_{\mbv_{\tilde E}} - \SI^{i_0-1}_\infty &  \, \longrightarrow  \,\, && \SI^{n-i_0-1} \quad  \quad\\
	\quad \uparrow \quad  &                     && \uparrow \quad\quad \nonumber \\
	\tilde \Sigma_{\tilde E}  \quad  & \, \longrightarrow  \,\, &&  K^{o}  \quad \quad\nonumber
\end{alignat}
gives us an image of $\tilde \Sigma_{\tilde E}$ 
that is the interior $K^{o}$ of a  properly convex compact set $K$. 

Let $\SI^{i_0}_\infty$ be a great $i_0$-dimensional sphere in $\SI^{n}$ containing $\mbv_{\tilde E}$ corresponding to the directions
of $\SI^{i_0-1}_\infty$ from $\mbv_{\tilde E}$. 
The space of  $(i_0+1)$-dimensional hemispheres in $\SI^n$
with boundary $\SI^{i_0}_\infty$ again has the structure of the projective sphere $\SI^{n-i_0-1}$, 
which is identifiable with the above one. 



Each $i_0$-dimensional 
hemisphere $H^{i_0}$ in $\SI^{n-1}_{\mbv_{\tilde E}}$ with 
$\Bd H^{i_0} = \SI^{i_0-1}_\infty$ corresponds to an $(i_0+1)$-dimensional hemisphere 
$H^{i_0+1}$ in $\SI^n$ with common boundary $\SI^{i_0}_\infty$ that contains $\mbv_{\tilde E}$. 

There is also a fibration whose fibers are open hemispheres of dimension $i_0+1$ with
common boundary $\SI^{i_0}_\infty$: 
\begin{alignat}{2} \label{np-eqn-pik2}
	\Pi_{K}:\SI^n - \SI_\infty^{i_0} &  \, \longrightarrow  \,\, && \SI^{n-i_0-1} \quad  \quad\\
	\quad \uparrow \quad  &                     && \uparrow \quad\quad \nonumber \\
	U \quad  & \, \longrightarrow  \,\, &&  K^{o}  \quad \quad\nonumber
\end{alignat}
since $\SI_\infty^{i_0-1}$ corresponds to $\SI_\infty^{i_0}$ in 
the projection $\SI^{n} - \{\mbv_{\tilde E}, \mbv_{\tilde E-}\} \ra 
\SI^{n-1}$. 
\index{PiK@$\Pi_K$} 

Let $\SL_\pm(n+1, \bR)_{\SI^{i_0}_\infty, \mbv_{\tilde E}}$ 
denote the subgroup of $\Aut(\SI^n)$ acting on 
$\SI^{i_0}_\infty$ and $\mbv_{\tilde E}$.  
The projection $\Pi_K$ induces a homomorphism 
\[\Pi_K^*: \SL_\pm(n+1, \bR)_{\SI^{i_0}_\infty, \mbv_{\tilde E}} 
\ra \SL_\pm( n-i_{0}, \bR).\]
\index{PiKast@$\Pi_K^\ast$} 

Suppose that $\SI^{i_0}_\infty$ is $h(\pi_1(\tilde E))$-invariant. 
We let $N$ be the subgroup of $h(\pi_1(\tilde E))$ of elements inducing trivial actions on $\SI^{n-i_0-1}$. 
The exact sequence 
\begin{equation}\label{np-eqn-exact}
	1 \ra N \ra h(\pi_1(\tilde E)) \stackrel{\Pi^*_K}{\longrightarrow} N_K \ra 1
\end{equation}
satisfies 
that the kernel normal subgroup $N$ acts trivially on $\SI^{n-i_0-1}$ but acts on each hemisphere with 
boundary equal to $\SI^{i_0}_\infty$,
and $N_K$ acts effectively by the action induced from $\Pi^*_K$.

Since $K$ is a properly convex domain, $K^{o}$ admits a Hilbert metric $d_{K}$ and 
$\Aut(K)$ is a subgroup of isometries of $K^{o}$. 
Here $N_K$ is a subgroup of the group $\Aut(K)$ of the group of projective automorphisms of $K$, and 
 is called the {\em properly convex quotient } of $h(\pi_1(\tilde E))$ or $\bGamma_{\tilde E}$. 

We showed: 
\begin{theorem}\label{np-thm-folaff}
	Let $\Sigma_{\tilde E}$ be the end orbifold of an NPNC R-p-end $\tilde E$ of a	properly convex $n$-orbifold $\orb$ with radial or totally geodesic ends. Let $\torb$ be the universal cover in $\SI^n$. 
	We consider the induced action of $h(\pi_1(\tilde E))$ 
	on $\Aut(\SI^{n-1}_{\mbv_{\tilde E}})$ 
for the corresponding end vertex $\mbv_{\tilde E}$. 
	Then the following hold\,{\em :} 
	\begin{itemize} 
		\item $\tilde \Sigma_{\tilde E}$ is foliated by complete affine subspaces of dimension
 $i_0 > 0$. We denote by $K$ the properly convex compact domain of dimension $n-i_0-1$ in $\SI^{n-i_0-1}$ whose interior is the space of complete affine subspaces of dimension $i_0$. 
		\item $h(\pi_1(\tilde E))$ acts on the great sphere $\SI^{i_0-1}_\infty$ of dimension $i_0-1$ in $\SI^{n-1}_{\mbv_{\tilde E}}$. 
		\item There exists an exact sequence 
		\[ 1 \ra N \ra \pi_1(\tilde E) \stackrel{\Pi^*_K}{\longrightarrow} N_K \ra 1 \] 
		where $N$ acts trivially on the quotient great sphere $\SI^{n-i_0-1}$, and 
		$N_K$ acts effectively on a properly convex domain $K^o$ in $\SI^{n-i_0-1}$ isometrically 
		with respect to the Hilbert metric $d_K$. 
	\end{itemize} 
\hfill $\square$ 
\end{theorem}

\subsection{Main results.}

\begin{figure}
	\centering
	\includegraphics[height=6cm]{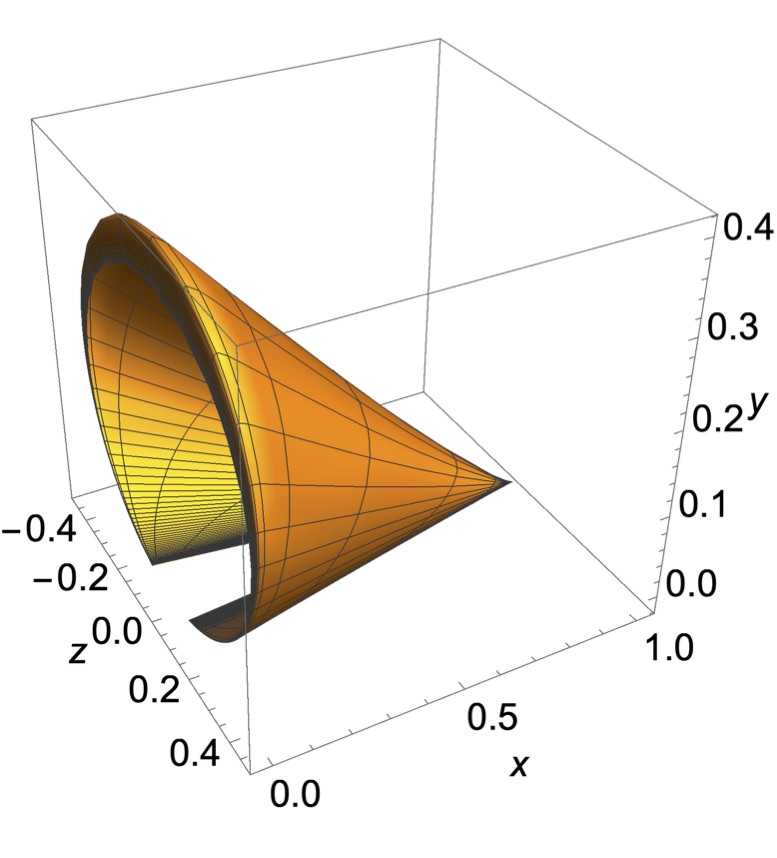}
	\caption{A partial development of 
		the boundary of a quasi-joined R-p-end neighborhood
		in an affine patch with Euclidean coordinates 
		where $\llrrparen{0,0,0,1}$ in $\SI^3$ corresponds to $(0,0,0)$. 
		In the notation of Definition \ref{np-defn-quasijoin}, 
	$\mbv$ is the point $(0,0,0)$, $\hat K$ is the singleton of $(1, 0, 0)$, 
	$H_x$ contains the upper-half of the $yz$-plane, and $\SI^{1}_\infty$ contains the $y$-axis. 
	A cusp group acts on each hyperspace containing the $y$-axis.  
	$H_x$ and $\SI^1_\infty$. These hyperspaces are not 
	in the affine space of the projection. 
	(See the mathematica file \cite{quasijoin})
}
	\label{np-fig:quasi-j}
\end{figure}
We begin with a definition of quasi-joined R-ends. 

\begin{definition}\label{np-defn-quasijoin} 
	We have the following: 
	\begin{itemize} 
	\item Let $\hat K$ be a compact properly convex subset of dimension 
$n-i-2$ in $\SI^n$ for $i \geq 1$. (Here, $\hat K$ will correspond to a codimension-one subspace of $K$ under the projection $\Pi_K$. )
\item 	Let $\SI_{\infty}^i$ be a subspace of dimension $i$ complementary to
	the subspace containing $\hat K$.
	\item Let $\mbv$ be a point in $\SI_{\infty}^i$. 
	\item A group $G$ acts on $\hat K$, $\SI_{\infty}^i$, and $\mbv$ and on  
	an open set $U$ with $\mbv\in \Bd U$. 
		\item $U/G$ is required to be diffeomorphic to a compact orbifold times an interval.
	\item There is a fibration $\Pi_K: \SI^n - \SI_{\infty}^i \ra \SI^{n-i-1}$ 
	with fibers equal to open $(i+1)$-hemispheres with boundary $\SI_{\infty}^i$. 
	\begin{itemize} 
	\item The set of fibering open $(i+1)$-hemispheres $H$ with $H\cap U \ne \emp$ 
is projected onto the interior of  $\{x\} \ast \Pi_K(\hat K)$
	in $\SI^{n-i-1}$ for a certain point $x$ not in $\Pi_K(\hat K)$. 
	\item  For each fibering open hemisphere $H$, 
	$H\cap U$ is an open set bounded by an ellipsoid
	containing $\mbv$ unless $H \cap U$ is empty.
	\item $\clo(H_x) \cap \clo(U) = \{\mbv\}$ for $H_x$ the fibering open $(i+1)$-hemisphere 
	over $x$. 
	\end{itemize} 

	\end{itemize}  
	If an R-end $E$ of a real projective orbifold has an end neighborhood 
	projectively diffeomorphic to $U/G$ with the induced radial foliation
	corresponding to $\mbv$, 
	then $E$ is called a {\em quasi-joined end}\/ 
	({\em of a totally geodesic R-end and a horospherical end
	with respect to $\mbv$}\/), 
	and a corresponding R-p-end is said to be a {\em quasi-joined R-p-end} also. 
Also, any R-end with an end neighborhood 
covered by an end neighborhood of a quasi-joined R-end is called by 
the same name. In these cases, the end holonomy group 
is a {\em quasi-joined end group}
({\em of a totally geodesic R-end and a horospherical end
	with respect to $\mbv$}\/)
	\end{definition} 
In this chapter, we will characterize 
the NPNC-ends and see some examples. See Proposition \ref{np-prop-qjoin} and Remark \ref{np-defn-qjoin} 
for detailed understanding of quasi-joined ends. 
\index{end!quasi-joined|textbf}

Let $\tilde E$ be an NPNC-end. 
 Recall from Chapter \ref{ch-ce} that the universal cover $\tilde \Sigma_{\tilde E}$ of 
 the end orbifold $\Sigma_{\tilde E}$ is foliated by complete affine $i_{0}$-dimensional totally geodesic leaves 
 for $i_{0} > 1$. 
 The end fundamental group $\pi_{1}(\tilde E)$ acts on a properly convex domain $K$
 that is the space of $i_{0}$-dimensional totally geodesic hemispheres foliating $\tilde \Sigma_{\tilde E}$.

\begin{definition} \label{np-defn-NS} 
	A countable group $G$ satisfies the property {\em (}NS\/{\em )} if 
	every normal solvable subgroup $N$ of a finite-index subgroup $G'$ is virtually central in $G'$;  
	that is, $N\cap G''$ is central in $G''$ for a finite-index subgroup $G''$ of $G'$. 
\end{definition} 
By Corollary \ref{prelim-cor-Benoist}, the fundamental group of 
a closed orbifold admitting a properly convex structure has 
the property (NS). 
Clearly, a virtually abelian group satisfies (NS). 
Obviously, the groups of Benoist are somewhat related to this condition. 
(see Proposition \ref{prelim-prop-Benoist}.)
\index{NS|textbf} 
\index{end!condition!NS|textbf} 
 
The main result of this chapter is the following. We need the proper convexity of $\torb$. The following theorem shows that given that 
$\tilde \Sigma_{\tilde E}$ is convex but not properly convex, 
the transverse weak middle-eigenvalue condition
implies that the end is quasi-joined type. For example, 
this implies that the holonomy group 
decomposes into a semisimple part and a horospherical part. 
(See \eqref{np-eqn-gformII}). 
This type is much easier to 
understand. 
A later Proposition \ref{np-prop-qjoin} gives some detailed discussions. 

\begin{theorem}\label{np-thm-thirdmain} 
Let $\mathcal{O}$ be a strongly tame properly convex real projective $n$-orbifold with radial or totally geodesic ends.
Assume that the holonomy group of $\mathcal{O}$ is strongly irreducible.
Let $\tilde E$ be an NPNC R-p-end with the end orbifold 
$\Sigma_{\tilde E}$. The universal cover $\tilde \Sigma_{\tilde E}$ is
foliated by $i_0$-dimensional totally geodesic hemispheres for $0< i_0 < n-1$. 
The leaf space naturally identifies with 
the interior of a compact convex set $K$. 
Suppose that the following hold\/{\rm :}
\begin{itemize} 
\item the end fundamental group satisfies the property {\rm (}NS\/{\rm )}
{\rm (}for example $\dim K^o =0, 1${\rm )}
for the leaf space $K^o$ of $\tilde E$. 
\item the p-end holonomy group $h(\pi_{1}(\tilde E))$ virtually satisfies the transverse weak middle-eigenvalue condition with respect to a p-end vertex $\mbv_{\tilde E}$
\end{itemize}
Then $\tilde E$ is a quasi-joined type R-p-end for $\mbv_{\tilde E}$. 
\end{theorem}
See Definition \ref{np-defn-weakmec} for the transverse weak middle-eigenvalue condition
for NPNC-ends. Without this condition, we doubt that we can obtain this type of results. 
However, it is open to investigations. 
In this case, $\tilde E$ does not satisfy the uniform middle-eigenvalue condition
as stated in Chapter \ref{ch-ce} for properly convex ends. 

For amenable groups, Ballas-Cooper-Leitner \cite{BCLp}, \cite{BCL22} had covered much of the material.



Recall that the dual orbifold $\orb^{\ast}$ is obtained from 
a properly convex real projective $n$-orbifold $\orb$. 
(See Section \ref{pr-sub-dualend}.) 
The set of ends of $\orb$ is in a one-to-one correspondence with the set of ends of $\orb^{\ast}$. 
We show that the dual of a quasi-joined NPNC R-p-end is a quasi-joined NPNC R-p-end. 

\begin{corollary}\label{np-cor-dualNPNC} 
Let $\mathcal{O}$ be a properly convex real projective $n$-orbifold with radial or totally geodesic ends.
Assume that the holonomy group of $\mathcal{O}$ is strongly irreducible.
Let $\tilde E$ be a quasi-joined NPNC R-p-end for an end $E$ of $\orb$ virtually satisfying 
the transverse weak middle-eigenvalue condition with respect to the p-end vertex. 
Suppose that 
the end fundamental group satisfies the property {\rm (}NS\/{\rm )}
{\rm (} in particular if $\dim K^o =1$\/{\rm )}
for the leaf space $K^o$ of $\tilde E$. 

Similarly, let $\orb^{\ast}$ denote the dual real projective orbifold of $\orb$. 
Let $\tilde E^{\ast}$ be a p-end corresponding to a dual end of $E$. 
Then $\tilde E^{\ast}$ has a p-end neighborhood of a quasi-joined type R-p-end for 
the universal cover of $\orb^\ast$ for a unique choice of a p-end vertex.  
\end{corollary} 
In short, we are saying that $\tilde E^{\ast}$ can be considered a quasi-joined type R-p-end 
by choosing its p-end vertex well. 
However, this does involve artificially introducing a radial foliation structure
in an end neighborhood.
We mention that the choice of the p-end vertex is uniquely 
determined for $\tilde E^{\ast}$ to be quasi-joined.

\subsection{Outline} \label{np-sub-outline}

In Section \ref{np-sec-notprop}, we discuss the R-ends that are NPNC. 
We introduce the transverse weak middle-eigenvalue condition. 
We explain the main eigenvalue estimates following from the transverse weak middle-eigenvalue condition for 
NPNC-ends. Then we explain our plan to prove Theorem \ref{np-thm-thirdmain}.

In Section \ref{np-sec-gendiscrete}, 
we introduce the example of the joining of horospherical and totally geodesic R-ends. 
We now study a bit more general situation introducing Hypothesis \ref{np-h-norm}. 
We try to obtain the splitting 
under some hypothesis. We will outline the subsection there.
By computations involving the normalization conditions, 
we show that the above exact sequence is virtually split under the condition (NS).
Surprisingly, we can demonstrate that 
the R-p-ends are of strictly joined or quasi-joined types. 
Then we show using the irreducibility of the holonomy group
of $\pi_{1}(\orb)$ that they can only be of quasi-joined type. 
We divide the tasks:  
\begin{itemize} 
\item In Section \ref{np-sub-hyp}, we introduce
Hypothesis \ref{np-h-norm} under which we work.  
We show that $K$ has to be a cone, and the conjugation action on 
$\CN$ has to be scalar orthogonal type changes.
Finally, we show the splitting of the NPNC-ends. 
We outline more completely there. 
\item In Section \ref{np-subsub-qjoin}, 
we introduce Hypothesis \ref{np-h-qjoin}, requiring 
more than Hypothesis \ref{np-h-norm}. 
We prove Proposition \ref{np-prop-qjoin} that the ends are 
quasi-joins under the hypothesis. 
\end{itemize} 
As a final part of this section in Section \ref{np-sub-discretecase}, 
we discuss the case when $N_K$ is discrete. We prove Theorem \ref{np-thm-thirdmain} for this case by showing that the above 
two hypotheses are satisfied.

In Section \ref{np-sec-nondiscrete}, we discuss when $N_K$ is not discrete. There is a foliation 
by complete affine subspaces as above. 
We use some estimates on eigenvalues to show that each leaf 
is of polynomial growth.
The leaf closures are suborbifolds $V_l$ 
by the theory of Carri\`ere \cite{Carriere88} and Molino \cite{Molino88} on Riemannian foliations. 
They form the fibration with compact fibers. 
$\pi_1(V_l)$ is solvable using the work of Carri\`ere \cite{Carriere88}. 
One can then take the syndetic closure to 
obtain a larger group that acts transitively on each leaf following Fried and Goldman \cite{FG83}. 
We find a unipotent cusp group acting on each leaf transitively normalized by $\bGamma_{\tilde E}$. 
Then we show that 
the end also splits virtually using the theory of Section \ref{np-sec-gendiscrete}. This proves Theorem \ref{np-thm-thirdmain}
for this case. 

In Section \ref{np-sec-dualNPNC}, we prove Corollary \ref{np-cor-dualNPNC}
showing that the duals of NPNC-ends are NPNS-ends, and 
This will be needed in the proof of Theorem \ref{ce-thm-mainaffine}. 

 In Section \ref{np-sub-counter}, we discuss a counterexample 
 to Theorem \ref{np-thm-NPNCcase} when the condition (NS) is dropped.

\subsubsection{The plan of the proof of Theorem \ref{np-thm-thirdmain}} \label{np-subsub-plan}
We show that our NPNC-ends are quasi-joined type ones; i.e., we prove
Theorem \ref{np-thm-thirdmain} 
For discrete $N_K$, we do this in 
Section \ref{np-sub-discretecase}, 
For nondiscrete $N_K$, we do this in Section \ref{np-sub-nondproof}. 

\begin{itemize} 
	\item Assume Hypotheses \ref{np-h-norm}.
	\begin{itemize} 
	\item We show that $\bGamma_{\tilde E}$ acts as a subgroup of
$\bR_+$ times an orthogonal group 
	on a Lie group $\CN$, which is 
realized as an real unipotent abelian group $\bR^{i_{0}}$. 
See Lemmas \ref{np-lem-similarity}. 
	This is done by computations and coordinate change arguments and 
	the distal group theory of Fried \cite{Fried86}. 
	\item We show that $K$ is a cone in Lemma \ref{np-lem-conedecomp1}. 
	\item We refine the matrix forms in Lemma \ref{np-lem-matrix} when $\mu_{g}=1$. 
	Here the matrices are in almost desired forms. 
	\item Proposition \ref{np-prop-decomposition} shows the splitting of the representation of $\bGamma_{\tilde E}$. 
	One uses the transverse weak middle-eigenvalue condition to realize the compact $(n-i_{0}-1)$-dimensional totally geodesic domain independent of $\SI_\infty^{i_0}$ where $\bGamma_{\tilde E}$ acts on. 
	\end{itemize}
\item Now we can assume Hypothesis  \ref{np-h-qjoin} additionally.
 In Section \ref{np-subsub-qjoin}, we discuss joins and quasi-joins. The idea is to show that 
	the join cannot occur by Propositions  \ref{prelim-prop-joinred}  and \ref{prelim-prop-decjoin}. 
\item This settles the cases of discrete $N_{K}$ in Theorem \ref{np-thm-NPNCcase} in Section \ref{np-sub-discretecase}. 
\item In Section \ref{np-sec-nondiscrete}, we settle the cases of nondiscrete $N_{K}$. See above for the outline of this section.
	
\end{itemize}

We remark that we can often take a finite-index subgroup of 
$\bGamma_{\tilde E}$ during our proofs 
since Definition \ref{np-defn-qjoin} is 
a definition given up to finite-index subgroups.

\begin{remark}
The main result of this chapter, Theorem \ref{np-thm-thirdmain} and Corollary \ref{np-cor-dualNPNC}, are stated without references to $\SI^n$ or 
$\RP^n$. We will work in the space $\SI^n$ only. 
Often the result for $\SI^n$ implies 
the result for $\RP^n$. In this case, we only prove for $\SI^n$. In other cases, we can easily modify 
the $\SI^n$-version proof to the one for the $\RP^n$-version proof. 
\end{remark}



\section{The transverse weak middle-eigenvalue conditions for NPNC ends} \label{np-sec-notprop}

We now study ends $E$ where corresponding $\Sigma_E$ are neither properly convex 
nor projectively diffeomorphic to a complete affine subspace.
Let $\tilde E$ be an R-p-end of $\orb$, and let $U$ be the corresponding p-end neighborhood in $\torb$
with the p-end vertex $\mbv_{\tilde E}$. 
Let $\tilde \Sigma_{\tilde E}$ denote the universal cover of the p-end orbifold $\Sigma_{\tilde E}$ as 
a domain in $\SI^{n-1}_{\mbv_{\tilde E}}$. 





We denote by $\bGamma_{\tilde E}$ the p-end holonomy group acting on $U$ fixing $\mbv_{\tilde E}$. 
We denote by $\widetilde{\mathcal{F}}_{\tilde E}$ the foliation on $\Sigma_{\tilde E}$, and 
 $\mathcal{F}_{\tilde E}$
the corresponding one in $\Sigma_{\tilde E}$. 
We recall 
\[\leng_K(g):= \inf \{ d_K(x, g(x))| x \in K^o \}, 
g \in \bGamma_{\tilde E}.\]
\index{length@$\leng_K(g)$} 

We recall Definition \ref{prelim-defn-jordan}. 
Let $V^{i_{0}+1}_\infty$ denote the subspace of $\bR^{n+1}$ corresponding 
to $\SI^{i_{0}}_\infty$. 
By invariance of $\SI^{i_{0}}_\infty$, 
if \[{\mathcal{R}}_{\mu}(g)\cap V^{i_{0}+1}_\infty \ne \{0\}, \mu > 0,\] 
then ${\mathcal{R}}_{\mu}(g) \cap V^{i_{0}+1}_\infty$ always contains a $\bC$-eigenvector of $g$. 

\begin{definition}\label{np-defn-eig} 
	Let $\Sigma_{\tilde E}$ be the end orbifold of a nonproperly convex R-p-end $\tilde E$ of a 
	properly convex $n$-orbifold $\orb$.
	Let $\bGamma_{\tilde E}$ 
	be the p-end holonomy group. 
	\begin{itemize}
		\item Let $\lambda^{Tr}_{\max}(g)$ denote the largest norm of the eigenvalue of $g \in \bGamma_{\tilde E}$ affiliated
		with $\vec{v} \ne 0$, $\llrrparen{\vec{v}} \in \SI^n - \SI^{i_0}_\infty$,
		i.e., \[\lambda^{Tr}_{\max}(g):= \max\{\mu\,| \exists\, \vec v \in  {\mathcal{R}}_{\mu}(g) -  V^{i_0+1}_\infty\},\]
		which is the maximal norm of transverse eigenvalues.  \index{eigenvalue!transverse|textbf} 
		\index{lambdaTrmax@$\lambda^{Tr}_{\max}(\cdot)$|textbf}
		\item Also, let $\lambda^{Tr}_{\min}(g)$ denote the smallest one affiliated with a nonzero vector $\vec v$, $\llrrparen{\vec{v}} \in \SI^n - \SI^{i_0}_\infty$,
		i.e., \[\lambda^{Tr}_{\min}(g):= \min \{\mu\,| \exists \vec v \in {\mathcal{R}}_{\mu(g)} - V^{i_0+1}_\infty\}, \]
		which is the minimal norm of transverse eigenvalues.
			\index{lambdaTrmin@$\lambda^{Tr}_{\min}(\cdot)$|textbf}
		\item Let $\lambda^{\Sio}_{\max}(g)$ be the largest of the norms of the eigenvalues of $g$ with $\bC$-eigenvectors 
		of form $\vec v$, $\llrrparen{\vec v} \in \Sio$
		and $\lambda^{\Sio}_{\min}(g)$ the smallest such one. 
		\index{lambdaSiomax@$\lambda^{\Sio}_{\max}(\cdot)$|textbf}
		\index{lambdaSiomin@$\lambda^{\Sio}_{\min}(\cdot)$|textbf}
	\end{itemize} 
\end{definition}

\begin{definition}\label{np-defn-weakmec}
Let $\lambda_{\mbv_{\tilde E}}(g)$ denote the eigenvalue of $g$ at $\mbv_{\tilde E}$. 
The {\em transverse weak middle-eigenvalue condition} with 
respect to $\mbv_{\tilde E}$ or the R-p-end structure 
is that each element $g$ of $h(\pi_{1}(\tilde E))$ satisfies
\begin{equation} \label{np-eqn-mec}
\lambda^{Tr}_{\max}(g) \geq \lambda_{\mbv_{\tilde E}}(g).
\end{equation} 
\end{definition} 
\index{middle-eigenvalue condition!weak!transverse|textbf} 
\index{transverse weak middle eigenvalue condition|textbf} 

Theorem \ref{abelian-thm-justify} provides some justification to our approach. 
We do believe that the weak middle-eigenvalue conditions
imply  the transverse ones at least for relevant cases.

The following proposition is very important in this chapter and shows that 
$\lambda^{Tr}_{\max}(g)$ and $\lambda^{Tr}_{\min}(g)$ are true largest and smallest norms of the eigenvalues of $g$. 

\begin{proposition}\label{np-prop-eigSI}  
Let $\Sigma_{\tilde E}$ be the end orbifold of an NPNC R-p-end $\tilde E$ of 
a properly convex real projective $n$-orbifold $\orb$ with radial or totally geodesic ends. 
Suppose that $\torb$ in $\SI^n$\/ {\rm (}resp. $\RP^n$\/{\rm )} 
covers $\orb$ as a universal cover. 
Let $\bGamma_{\tilde E}$ be the p-end holonomy group satisfying the transverse weak middle-eigenvalue condition for the R-p-end structure. 
Let $g \in \bGamma_{\tilde E}$. 
Then the following hold\/{\em :}
\begin{gather} \label{np-eqn-eigen} 
\lambda^{Tr}_{\max}(g) \geq \lambda^{\Sio}_{\max}(g) \geq  \lambda_{\mbv_{\tilde E}}(g) \geq \lambda^{\Sio}_{\min}(g) \geq \lambda^{Tr}_{\min}(g), \nonumber \\
\lambda_1(g) = \lambda^{Tr}_{\max}(g) \hbox{, and }
\lambda_{n+1}(g) = \lambda^{Tr}_{\min}(g).
\end{gather} 
\end{proposition} 
\begin{proof}
We may assume that $g$ has an infinite order. 
Suppose that $\lambda^{\Sio}_{\max}(g) > \lambda^{Tr}_{\max}(g)$. 
We have $\lambda^{\Sio}_{\max}(g) > \lambda_{\mbv_{\tilde E}}(g)$ by 
the transverse weak uniform middle-eigenvalue condition. 

Now, $\lambda^{Tr}_{\max}(g) < \lambda^{\Sio}_{\max}(g)$ implies that
\[R_{\lambda^{\Sio}_{\max}(g)}:= \bigoplus_{|\mu|=\lambda^{\Sio}_{\max}(g)}  \mathcal{R}_{\mu}(g)\]
is a subspace of $V^{i_0+1}_\infty$ and corresponds to a great sphere $\SI^j$, 
$\SI^j \subset \SI^{i_0}_\infty$.  
Hence, a great sphere $\SI^j$, $j \geq 0$, in $\SI^{i_0}_\infty$ is disjoint from 
$\{\mbv_{\tilde E}, \mbv_{\tilde E-}\}.$ 
Since $\mbv_{\tilde E} \in \SI^{i_0}_\infty$ is not contained in $\SI^j$, we obtain $j+1 \leq i_0$.

A vector space $V_1$ corresponds $\bigoplus_{|\mu|<\lambda^{\Sio}_{\max}(g)}  \mathcal{R}_{\mu}(g) $ where 
$g$ has strictly smaller norm eigenvalues and is complementary to $R_{\lambda^{\Sio}_{\max}(g)}$. 
Let $C_{1} =\SI(V_{1})$. 
The great sphere $C_1$ is disjoint from $\SI^j$; however, $C_1$ contains $\mbv_{\tilde E}$. 
Moreover, $C_1$ is of complementary dimension to $\SI^j$, i.e., $\dim C_1 = n - j-1$. 


Since $C_1$ is complementary to $\SI^j \subset \SI^{i_0}_\infty$, $C_1$ contains 
a complementary subspace 
$C'_1$ to $\SI^{i_0}_\infty$ of dimension $n-i_0-1$ in $\SI^n$. 
Considering the sphere $\SI^{n-1}_{\mbv_{\tilde E}}$ at $\mbv_{\tilde E}$, 
it follows that $C'_1$ is mapped to an $(n-i_0-1)$-dimensional subspace $C''_1$ in 
$\SI^{n-1}_{\mbv_{\tilde E}}$ disjoint from $\partial l$ for any complete affine leaf $l$. 
Each complete affine leaf $l$ of $\tilde \Sigma_{\tilde E}$ has the dimension $i_0$ 
and meets $C''_1$ in $\SI^{n-1}_{\mbv_{\tilde E}}$ by the dimension consideration. 


\begin{figure}
\centering
\includegraphics[trim = 15mm 30mm 10mm 10mm, clip,  height=5cm]{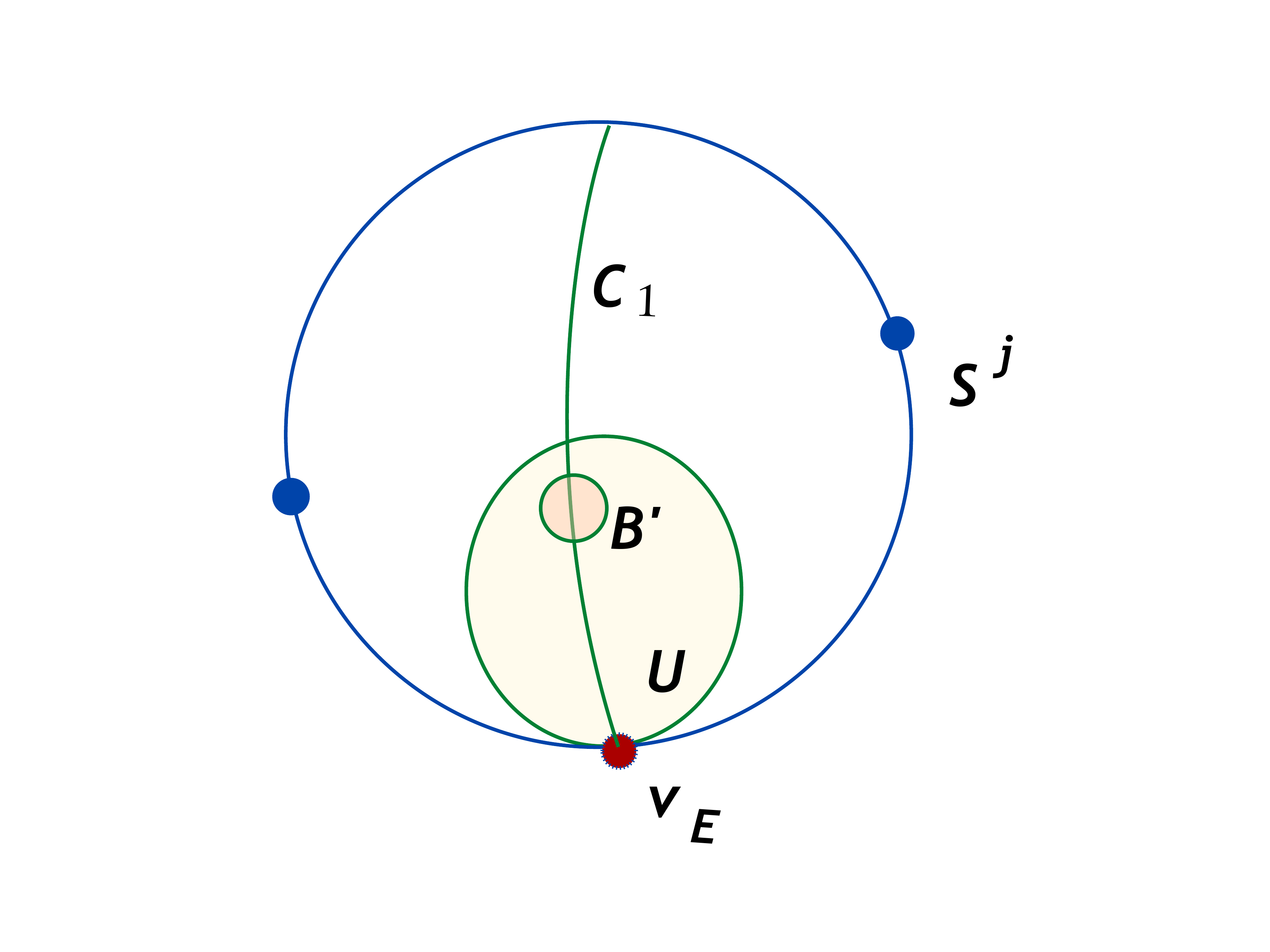}
\caption{The figure for the proof of Proposition \ref{np-prop-eigSI}. }
\label{np-fig:compl}
\end{figure}

Hence, a small ball $B'$ in $U$ meets $C_1$ in its interior.
\begin{align}
& \hbox{For any } \llrrparen{\vec{v}} \in B', \vec{v} \in \bR^{n+1}, \vec{v} = \vec{v}_1 + \vec{v}_2 \hbox{ where } \llrrparen{\vec{v}_1} \in C_1 \hbox{ and } \llrrparen{\vec{v}_2} \in \SI^j. \nonumber \\
& \hbox{We obtain } g^k(\llrrparen{\vec{v}}) = \llrrparen{g^k(\vec{v}_1) + g^k(\vec{v}_2)},
\end{align}
where we used $g$ to represent the linear transformation of determinant $\pm 1$ as well (See Remark \ref{intro-rem-SL}.) 
By Proposition \ref{prelim-prop-attract}, 
the action of $g^k$ as $k \ra \infty$ makes the former vectors very small compared to the latter ones, i.e., 
\[ \left\{\llrrV{g^k(\vec{v}_1)}/\llrrV{g^k(\vec{v}_2)}\right\} \ra 0 \hbox{ as } k \ra \infty.\]
Hence, $\{g^k(\llrrparen{\vec{v}})\}$ converges to the same limit as $\{g^k(\llrrparen{\vec{v}_2})\}$ if it exists. 

Now choose $\llrrparen{\vec{w}}$ in $C_1 \cap B'$ and $\vec{v}, \llrrparen{\vec{v}}\in \SI^j$. 
For small $\eps>0$, 
set $\vec{w}_1 = \llrrparen{\vec{w} + \eps \vec{v}}$ and $\vec{w}_2 =\llrrparen{\vec{w} -\eps \vec{v}}$ in $B'$. 
Choose a subsequence $\{k_i\}$ such that $\{g^{k_i}(\vec{w}_1)\}$ converges to a point of $\SI^n$. 
The above estimation shows that $\{g^{k_i}(\vec{w}_1)\}$ and $\{g^{k_i}(\vec{w}_2)\}$ converge to 
an antipodal pair of points in $\clo(U)$ respectively. 
This contradicts the proper convexity of $U$ 
as $g^k(B') \subset U$ and its geometric limit is in $\clo(U)$. 

Also the consideration of $g^{-1}$ yields the reverse inequality, and 
the second equation follows from the first one.
\hfill \SnT  {\parfillskip0pt\par}
\end{proof}



\section{The general theory} \label{np-sec-gendiscrete}



\subsection{Examples} 
First, we give some examples. 

\subsubsection{The standard quadric in $\bR^{i_0+1}$ and the group acting on it.} \label{np-subsub-quadric}

Let us consider an affine subspace $\mathbb{A}^{i_0+1}$ of $\SI^{i_0+1}$ with 
coordinates $x_0, x_1, \dots, x_{i_0+1}$ given by $x_0 > 0$.
The standard quadric in $\mathbb{A}^{i_0+1}$ is given by 
\[ x_{i_0+1}x_0 = \frac{1}{2}(x_1^2 + \cdots + x_{i_0}^2). \] 
Clearly, the group of orthogonal maps $\Ort(i_0)$ acting on the planes given by $x_0=C_1, x_{i_0+1} = C_2$ for constants $C_1, C_2$ acts on
the quadric also. Also, the group of matrices of the form 
\[
\left(
\begin{array}{ccc}
 1                                 & 0                    & 0  \\
 \vec{v}^T                    & \Idd_{i_0}        & 0   \\
\frac{\llrrV{\vec{v}}^2}{2}  & \vec{v}   & 1  
\end{array}
\right)
\]
acts on the quadric. 

The group of affine transformations that acts on the quadric is exactly the Lie group
generated by the above cusp group and $\Ort(i_0)$. The action is transitive 
and each stabilizer is a conjugate of
$\Ort(i_0)$ by elements of the above group. 

The proof of this fact is simply that such  a group of  affine transformations is conjugate to a parabolic isometry group in
the $(i_0+1)$-dimensional complete hyperbolic space $H$ where the ideal fixed point is identified with 
$\llrrparen{0, \dots, 0, 1} \in \SI^{i_0+1}$ and with $\Bd H$ tangent to $\Bd A^{i_0}$.

The group 
of projective automorphisms of the following forms is called a unipotent cusp group. 
\begin{equation} \label{np-eqn-generator}
\CN'(\vec{v}):= \left( \begin{array}{ccc}                
                                                         1 &  \vec{0}  & 0 \\
                                                        \vec{v}^T  & \Idd_{i_0-1} & \vec{0}^T \\ 
                                                \frac{\llrrV{\vec{v}}^2}{2}  &  \vec{v}       & 1 
                                    \end{array} \right) \hbox{ for } \vec{v} \in \bR^{i_0}. 
 \end{equation}
 (see \cite{CM14} for details.) \index{cusp group!unipotent|textbf} 
 We can make each translation direction of a generator of $\CN$ in $\tilde \Sigma_{\tilde E}$ to be one of the standard vectors. 
Therefore, we can find 
a coordinate system of $V^{i_0+2}$ such that the generators are the 
$(i_0+2) \times (i_0+2)$-matrices of forms 
\begin{equation} \label{np-eqn-generator2}
\CN'_j(s):=  \left( \begin{array}{cccc}                1 &   \vec{0} & 0 \\ 
                                                        s\vec{e}^T_j  & \Idd_{i_0}    & \vec{0}^T \\
                                                       \frac{1}{2} &  s\vec{e}_j     & 1 
                                    \end{array} \right)
 \end{equation}
 where $s \in \bR$ and 
 $(\vec{e}_{j})_{k} = \delta_{jk}$, the row $i$-vector for $j=1, \dots, i_{0}$.  
 That is, 
 \[\CN'(\vec{v}) = \CN'_1(v_{1})\cdots \CN'_{i_0}(v_{i_{0}}).\]


\subsubsection{Examples of generalized joined ends}

We begin with some  examples. 

\begin{example} \label{np-exmp-joined} 
Let us consider two ends $E_1$, a totally geodesic R-end, 
with the p-end neighborhood $U_1$ in the universal cover of a real projective orbifold $\mathcal{O}_1$ in $\SI^{n-i_0-1}$
of dimension $n-i_0-1$ with the p-end vertex $\mbv_1$, 
and $E_2$  with the p-end neighborhood $U_2$, a horospherical type one, in the universal cover
of a real projective orbifold $\mathcal{O}_2$ of dimension $i_0+1$
with the p-end vertex $\mbv_2$. 
\begin{itemize} 
\item Let $\bGamma_1$ denote the projective automorphism group in $\Aut(\SI^{n-i_0-1})$ acting on $U_1$ 
corresponding to $E_1$. 
We assume that $\bGamma_1$ acts on a great sphere $\SI^{n-i_0-2} \subset \SI^{n-i_0-1}$ disjoint from $\mbv_1$.
There exists a properly convex open domain $K'$ in $\SI^{n-i_0-2}$ where $\bGamma_1$ acts cocompactly
but not necessarily freely. 
We change $U_1$ to be the interior of the join of $K'$ and $\mbv_1$. 
\item Let $\bGamma_2$ denote the one in $\Aut(\SI^{i_0+1})$ acting on $U_2$ 
as a subgroup of a cusp group. 
\item We embed $\SI^{n-i_0-1}$ and $\SI^{i_0+1}$ in $\SI^n$ such that 
they meet transversely at $\mbv = \mbv_1 = \mbv_2$. 
\item We embed $U_2$ in $\SI^{i_0+1}$ and $\bGamma_2$ in $\Aut(\SI^n)$ fixing each point of 
$\SI^{n-i_0-1}$. 
\item We can embed $U_1$ in $\SI^{n-i_0-1}$ and $\bGamma_1$ in $\Aut(\SI^n)$ acting 
on the embedded $U_1$ such that 
$\bGamma_1$ acts on $\SI^{i_0-1}$ normalizing $\bGamma_2$ and acting on $U_1$. 
One can find some of such embeddings by choosing an arbitrary 
homomorphism $\rho: \bGamma_1 \ra N(\bGamma_2)$ for a normalizer $N(\bGamma_2)$ of 
$\bGamma_2$ in $\Aut(\SI^n)$. 
\end{itemize}

We denote by $\Sigma_1$ the quotient orbifold $K^{\prime o}/\bGamma_1$ and
$\Sigma_2$ the quotient orbifold $\partial U_2/\bGamma_2$. 
We find an element $\zeta \in \Aut(\SI^n)$ fixing each point of $\SI^{n-i_0-2}$ and 
acting on $\SI^{i_0+1}$ as a unipotent element normalizing $\bGamma_2$ 
such that the corresponding matrix has only two norms of eigenvalues. 
Then $\zeta$ centralizes $\bGamma_1| \SI^{n-i_0-2}$ and normalizes $\bGamma_2$. 
Let $U$ be the strict join of $U_1$ and $U_2$, which is a properly convex domain. 
$U/\langle \bGamma_1, \bGamma_2, \zeta \rangle$ gives us an R-p-end of dimension $n$
diffeomorphic to $\Sigma_{E_1} \times \Sigma_{E_2} \times \SI^1 \times \bR$, 
and the transverse real projective orbifold is diffeomorphic to $\Sigma_{E_1} \times \Sigma_{E_2} \times \SI^1$. 
We call the results the {\em generalized joined} end and the generalized joined end neighborhoods. 
Those ends with end neighborhoods finitely covered by these are also called 
{\em  generalized joined} ends. The generated group $\langle \bGamma_1, \bGamma_2, \zeta \rangle$ is 
called a {\em generalized joined group}. 

Now we generalize this construction slightly: 
As above, assume $U_1 \subset \SI^{n-i_0-1}$ and $U_2 \subset \SI^{i_0+1}$. 
Suppose now that $U_i$ for each $i=1, 2$
 is a domain (not necessarily open) radially foliated by segments from 
a common point $\mbv \in \Bd U_i$. 
There is a manifold boundary $\partial U_i$ transversal
to these radial segments. We do not assume $\partial U_i$ is dense in $\Bd U_i$. 

Suppose that $\bGamma_1$ 
and $\bGamma_2$ are discrete subgroups of $\Aut(\SI^n)$, and they 
act on $\partial U_1$ and $\partial U_2$  respectively, 
and we have a parameter of $\zeta^t$ for $ t\in \bR$ 
normalizing both $\bGamma_1$ and $\bGamma_2$. 
Let $U = U_1 \ast U_2$. 
Let $G$ denote the group generated by 
$\bGamma_1$ and $\bGamma_2$ and $\{\zeta_t\}$. 
Suppose there is a subgroup of $G$ acting on 
$\partial U_1 \ast \partial U_2$ cocompactly and properly discontinuously. 
Then $G$ acts as a {\em generalized joined action} on $U$, and 
$U / G$ is also called a {\em generalized joined end}. 


\index{end!joined!generalized} 
\end{example} 

\begin{remark}
We will deform this construction to a construction to obtain quasi-joined ends
(see Definition \ref{np-defn-quasijoin}).
This will be done by adding some translations appropriately. 
\end{remark}

We continue the above example to a more specific situation. 

\begin{example}\label{np-exmp-nonexmp}
Let $N$ be as in \eqref{np-eqn-nilmatstd}. 
Let $N$ be a cusp group 
conjugate to one in $\SO(i_0+1, 1)$ acting on an $i_0$-dimensional ellipsoid in $\SI^{i_0+1}$. 

We find a closed properly convex real projective 
orbifold $\Sigma$ of dimension $n-i_0-2$ and find a homomorphism from $\pi_1(\Sigma)$
to a subgroup of $\Aut(\SI^{i_0+1})$ normalizing $N$. 
(We begin by using  a trivial homomorphism.)
Using this, we obtain a group $\Gamma$ such that we have an exact sequence
\[ 1 \ra N \ra \Gamma \ra \pi_1(\Sigma) \ra 1. \] 
We assume that this sequence splits; 
i.e., $\pi_1(\Sigma)$ acts trivially on $N$.

We now consider an example where $i_0 = 1$. 
Let $N$ be $1$-dimensional, generated by $N_1$ as in \eqref{np-eqn-nilmat3}. 
\renewcommand{\arraystretch}{1.5}
\begin{equation} \label{np-eqn-nilmat3}
\newcommand*{\tempV}{\multicolumn{1}{r|}{}}
N_1:= \left( \begin{array}{ccccc}         \Idd_{n-i_0-1} & 0 & \tempV &  0 & 0 \\ 
                                                       \vec{0}           & 1  & \tempV &  0  & 0\\
                                                       \cline{1-5}
                                                       \vec{0}         & 1 & \tempV & 1   & 0 \\
                                                       \vec{0}       &\frac{1}{2} & \tempV & 1  & 1 
                                    \end{array} \right)
 \end{equation}
 where $i_0 = 1$. 

We take a discrete faithful proximal semisimple 
representation 
\[\tilde h: \pi_1(\Sigma) \ra \GL(n-i_0, \bR)\]
acting on a convex cone $C_\Sigma$ in $\bR^{n-i_0}$.
We define 
\[h: \pi_1(\Sigma) \ra \GL(n+1, \bR)\] by matrices
\begin{equation} \label{np-eqn-gammaJp}
h(g):= \left( \begin{array}{ccc}  \tilde h(g)         & 0                          & 0\\ 
                                                 \vec{d}_1(g)    & a_1(g)           & 0 \\
                                                 \vec{d}_2(g)    & c(g)  & \lambda_{\mbv_{\tilde E}}(g)    
                                                               \end{array} \right)
\end{equation}
where $\vec{d}_1(g)$ and $\vec{d}_2(g)$ are $(n-i_0)$-vectors, $a_1(g)$ is a constant 
depending on $g$, 
$g \mapsto \lambda_{\mbv_{\tilde E}}(g)$ is a homomorphism as defined above
for the p-end vertex, where $(\det \tilde h(g)) a_1(g) \lambda_{\mbv_{\tilde E}}(g) = 1$ holds.
 \begin{equation} \label{np-eqn-gammaJpi}
h(g^{-1}):= \left( \begin{array}{cc}  \tilde h(g)^{-1}         & \begin{array}{cc} \,0  & \quad \quad \quad 0 \end{array}   \\ 
                             -   \left(   \begin{array}{cc} \frac{\vec{d}_1(g)}{a_1(g)} \\ \frac{-c(g)\vec{d}_1(g)}{a_1(g) \lambda_{\mbv_{\tilde E}}(g)} 
                                + \frac{\vec{d}_2(g)}{\lambda_{\mbv_{\tilde E}}(g)}
                                                                 \end{array}   \right) \tilde h(g)^{-1} 

                                          & \begin{array}{cc} \frac{1}{a_1(g)} & 0 \\ 
                                          \frac{-c(g)}{a_1(g) \lambda_{\mbv_{\tilde E}}(g)} &  
                                          \frac{1}{\lambda_{\mbv_{\tilde E}}(g)} \end{array}\\ 
                         
                                                               \end{array} \right).
\end{equation}
 
Then the conjugation of $N_1$ by $h(g)$ gives us
 \begin{equation} \label{np-eqn-nilmat4}
 \left( \begin{array}{cc}         \Idd_{n-i_0}  & \begin{array}{cc} \quad 0  & \quad \,  0 \end{array} \\ 
                                                      \left(\begin{array}{cc} \vec{0} & a_1(g) \\ \vec{\ast} & \ast \end{array} \right) \tilde h(g)^{-1}&
                                                      \begin{array}{cc} 1   & 0 \\ \frac{\lambda_{\mbv_{\tilde E}}(g)}{a_1(g)}  & 1 \end{array} 
                                                      \end{array}
                                                      \right).
 \end{equation}
 Our condition on the form of $N_1$ implies that 
 $\underbrace{(0,0,\dots,0,1)}_{n-i_0}$ 
 has to be a common eigenvector by $\tilde h(\pi_1(\tilde E))$.  
 We also assume that $a_1(g) = \lambda_{\mbv_{\tilde E}}(g)$ for simplicity. 
 The last row of $\tilde h(g)$ equals 
 $(\vec{0}, \lambda_{\mbv_{\tilde E}}(g))$. Thus, the upper left 
 $(n-i_0)\times (n-i_0)$-block of $h(\pi_1(\tilde E))$ is reducible. 

Suppose that $\pi_1(\tilde E)$ is generated by $N_1$ and another group $G'$. 
Some further computations show that 
we can take any 
\[h: \pi_1(\tilde E) \ra \SL(n-i_0, \bR)\] with matrices of the form 
\renewcommand{\arraystretch}{1.2}
\begin{equation} \label{np-eqn-nilmat5}
\newcommand*{\tempV}{\multicolumn{1}{r|}{}}
h(g):= \left( \begin{array}{ccccc}   S_{n-i_0-1}(g) & 0 & \tempV &  0 & 0 \\ 
                                                       \vec{0}           & \lambda_{\mbv_{\tilde E}}(g)  &   \tempV &  0  & 0\\
                                                        \cline{1-5}
                                                       \vec{0}         & 0 & \tempV & \lambda_{\mbv_{\tilde E}}(g)   & 0 \\
                                                       \vec{0}          & 0 & \tempV & 0  & \lambda_{\mbv_{\tilde E}}(g) 
                                    \end{array} \right)
 \end{equation}
 for $g \in G'$ by a choice of coordinates by the semisimple property of the $(n-i_0) \times (n-i_0)$-upper left part of $h(g)$. 
 
 Since $\tilde h(\pi_1(\tilde E))$ has a common eigenvector, 
 Theorem 1.1 of Benoist \cite{Benoist03} shows that the open convex domain $K$, 
 that is the image of $\Pi_K$, is reducible. 
 We assume that 
 $N_K = N'_K \times \bZ$ for another subgroup $N'_K$,
and the image of the homomorphism $g \in N'_K \ra S_{n-i_0-1}(g)$ gives 
a discrete projective automorphism group acting properly discontinuously on a properly convex subset 
$K'$ in $\SI^{n-i_0-2}$ with a compact quotient. 

Let $\mathcal{E}$ be the one-dimensional ellipsoid 
where lower-right $3\times 3$-blocks of $N$ acts on. 
From this, the end is of the join form 
$K^{\prime o}/N'_K \times \SI^1 \times {\mathcal{E}}/\bZ$ 
by taking a double cover if necessary. 
Moreover, $\pi_1(\tilde E)$ is isomorphic to $N'_K \times \bZ \times \bZ$
up to taking an index two subgroups. 

We can think of this as the join of $K^{\prime o}/N'_K$ with $\mathcal{E}/\bZ$ as
$K'$, and $\mathcal{E}$ lies on disjoint complementary projective spaces
of respective dimensions $n-3$ and $2$ to be denoted $\SI(K')$ and $\SI(\mathcal{E})$ respectively. 


\end{example}

\subsection{Hypotheses to derive the splitting result} \label{np-sub-hyp} 
These hypotheses will help us to obtain the splitting. 
Afterward, we will show the NPNC-ends with transverse weak middle-eigenvalue conditions will satisfy these. 

We outline this subsection. 
In Section \ref{np-subsub-matrixform},  we ntroduce a standard coordinate system to work on, where we introduce the 
unipotent cusp group $\CN \cong \bR^{i_{0}}$. Our group
$\bGamma_{\tilde E}$ normalizes $\CN$ by the hypothesis. 
{\em Similarity} Lemma \ref{np-lem-similarity} shows that 
a conjugation in $\CN$ by an element of $\bGamma_{\tilde E}$ acts as a similarity, 
a simple consequence of the normalization property. 
We use this similarity and the Benoist theory \cite{Benoist03} to prove
{\em $K$-is-a-cone} Lemma \ref{np-lem-conedecomp1} that $K$ decomposes into a cone $\{k\} \ast K'' $ where 
$\CN$ has a nice expression for the adopted coordinates.
(If an orthogonal group acts cocompactly on an open manifold, then the manifold is zero-dimensional.) 
In Section \ref{np-subsub-split}, 
{\em splitting} Proposition \ref{np-prop-decomposition} shows that 
the end holonomy group splits. To do that, we find a sequence of elements of the virtual center 
expanding neighborhoods of a copy  of $K''$.   
Here, we explicitly find a part corresponding 
to $K'' \subset \Bd \torb$, and we realize $k$ as an $(i_{0}+1)$-dimensional hemisphere 
where $\CN$ acts on. 


\subsubsection{The matrix form of the end holonomy group.} \label{np-subsub-matrixform}

Let $\bGamma_{\tilde E}$ be an R-p-end holonomy group,
and let $l \subset \SI^{n-1}_{\mbv_{\tilde E}}$ be a complete $i_0$-dimensional leaf in $\tilde \Sigma_{\tilde E}$ to be specified a bit later. 
Then a unique great sphere $\SI^{i_0+1}_l$ contains the great segments from 
$\mbv_{\tilde E}$ in the direction of $l$. 
Let $V^{i_0+1}$ denote the subspace corresponding to $\SI^{i_0}_\infty$ containing $\mbv_{\tilde E}$, 
and $V^{i_0+2}_l$ the subspace corresponding to $\SI^{i_0+1}_l$.
We choose a coordinate system such that 
\[\mbv_{\tilde E} = \underbrace{\llrrparen{0, \cdots, 0, 1}}_{n+1}\]
and points of $V^{i_0+1}$, and those of $V^{i_0+2}$ respectively correspond to 
\[ \llrrparen{\overbrace{0, \dots, 0}^{n-i_0}, \ast, \cdots, \ast}, \quad \llrrparen{\overbrace{0, \dots, 0}^{n-i_0-1}, \ast, \cdots, \ast}.\]

Since $\SI^{i_0}_\infty$ and $\mbv_{\tilde E}$ are $g$-invariant, 
$g\in \bGamma_{\tilde E}$, is of {\em standard} form 
\renewcommand{\arraystretch}{1.2}
\begin{equation}\label{np-eqn-matstd}
\newcommand*{\tempV}{\multicolumn{1}{r|}{}}
\left( \begin{array}{ccccccc} 
S(g) & \tempV & s_1(g) & \tempV & 0 & \tempV & 0 \\ 
 \cline{1-7}
s_2(g) &\tempV & a_1(g) &\tempV & 0 &\tempV & 0 \\ 
 \cline{1-7}
C_1(g) &\tempV & a_4(g) &\tempV & A_5(g) &\tempV & 0 \\ 
 \cline{1-7}
c_2(g) &\tempV & a_7(g) &\tempV & a_8(g) &\tempV & a_9(g) 
\end{array} 
\right)
\end{equation}
where $S(g)$ is an $(n-i_0-1)\times (n-i_0-1)$-matrix,
$s_1(g)$ is an $(n-i_0-1)$-column vector, 
$s_2(g)$ and $c_2(g)$ are $(n-i_0-1)$-row vectors, 
$C_1(g)$ is an $i_0\times (n-i_0-1)$-matrix, 
$a_4(g)$ is an  $i_0$-column vector, 
$A_5(g)$ is an $i_0\times i_0$-matrix, 
$a_8(g)$ is an $i_0$-row vector, 
and $a_1(g), a_7(g)$, and $ a_9(g)$ are scalars. 

Denote 
\[ \hat S(g) =
\left( \begin{array}{cc} 
S(g) & s_1(g)\\ 
s_2(g) & a_1(g)
\end{array} \right),\]
and is called an {\em upper-left part} of $g$.  
\index{upper-left part} 

Let $\CN$ be a unipotent group acting on $\SI^{i_{0}}_{\infty}$, 
inducing $\Idd$ on $\SI^{n-i_{0}-1}$, and restricting to a cusp group for at least one great  $(i_{0}+1)$-dimensional sphere
$\SI^{i_{0} + 1}_l$ containing $\SI^{i_{0}}_{\infty}$.

We can write each element  $g \in \CN$ as an $(n+1)\times (n+1)$-matrix
\begin{equation} \label{np-eqn-nilmat}
 \left( \begin{array}{ccc}                \Idd_{n-i_0-1} &  0 &   0 \\ 
                                                       \vec{0}      & 1  &  0 \\
                                                     C_g          & *   & U_g 
                                    \end{array} \right)
 \end{equation}
where $C_g > 0$ is an $(i_0+1)\times (n-i_0-1)$-matrix, 
$U_g$ is a unipotent $(i_0+1) \times (i_0+1)$-matrix, 
$0$ indicates various zero row or column vectors, 
$\vec{0}$ denotes the zero row-vector of dimension $n-i_0-1$, and 
$\Idd_{n-i_0-1}$ is  the $(n-i_0-1) \times (n-i_0-1)$-identity matrix. 
This follows since $g$ acts trivially on $\bR^{n+1}/V^{i_0+1}_l$ 
and $g$ acts as a cusp group element on the subspace $V^{i_0+2}_l$.  



For $\vec{v} \in \bR^{i_0}$, we define 
\renewcommand{\arraystretch}{1.5}
\begin{equation} \label{np-eqn-nilmatstd}
\newcommand*{\tempV}{\multicolumn{1}{r|}{}}
\CN(\vec{v}):= \left( \begin{array}{ccccccc}         
	\Idd_{n-i_0-1} & 0 &\tempV &  0 & 0& \dots & 0 \\ 
    \vec{0}           & 1  &\tempV &  0  & 0&  \dots & 0\\
                    \cline{1-7}
        \vec{c}_1(\vec{v})    & {v}_1 &\tempV & 1   & 0 &  \dots & 0 \\
       \vec{c}_2(\vec{v})   & {v}_2 &\tempV & 0   & 1 & \dots & 0\\
      \vdots   & \vdots &\tempV & \vdots & \vdots & \ddots & \vdots \\
      \vec{c}_{i_0+1}(\vec{v})  & \frac{1}{2}\llrrV{\vec{v}}^2& \tempV & {v}_1 & v_2 & \dots  & 1
          \end{array} \right)
 \end{equation}
where $\llrrV{v}$ is the norm of $\vec{v} = (v_1, \cdots, v_{i_0}) \in \bR^{i_0}$. 
We assume that \[\CN := \{ \CN(\vec{v})| \vec v \in \bR^{i_{0}}\}\] is a group, which must be nilpotent. 
The elements of our nilpotent group 
$\CN$ are of this form since $\CN(\vec{v})$ is the product  
$\prod_{j=1}^{i_0} \CN(e_j)^{v_j}$.  
By the way we defined this, 
$\vec{c}_k:\bR^{i_0} \ra \bR^{n-i_0-1}$ for each $k=1, \dots, i_0$ is a  linear function 
of $\vec{v}$, defined as 
\[\vec{c}_k(\vec{v}) = \sum_{j=1}^{i_0} \vec c_{kj} v_j \hbox{ for } \vec{v} = (v_1, v_2, \dots, v_{i_0})\]
such that we form a group. 
(We do not need the property of $\vec{c}_{i_0+1}$ at the moment.)

From now on, 
\begin{itemize} 
\item we denote by $C_1(\vec{v})$ the $(n-i_0-1) \times i_0$-matrix given by 
the $j$-th row $\vec{c}_j(\vec{v})$ for $j= 1, \dots, i_0$
and 
\item by $c_2(\vec{v})$ the row $(n-i_0-1)$-vector $\vec{c}_{i_0+1}(\vec{v})$. 
\item The lower-right $(i_0+2)\times (i_0+2)$-matrix is called the {\em standard cusp matrix form}.
\end{itemize}

We denote by $\hat A$ the matrix
\begin{equation}
\newcommand*{\tempV}{\multicolumn{1}{r|}{}}
\left( \begin{array}{ccccccc} 
\Idd_{n-i_0-1} & \tempV & 0 & \tempV & 0 & \tempV & 0 \\ 
\cline{1-7}
0 &\tempV & 1 &\tempV & 0 &\tempV & 0 \\ 
\cline{1-7}
0 &\tempV &0 &\tempV &  A &\tempV & 0 \\ 
\cline{1-7}
0 &\tempV & 0 &\tempV &  0 &\tempV & 1
\end{array} 
\right)
\end{equation}
for $A$ an $i_0\times i_0$-matrix.
Denote the group of the form 
\[ \{\hat O_5| O_5 \in \Ort(i_0)\}\] 
by $\hat \Ort(i_0)$.

If $\mathcal{N}$ can be put in the form \eqref{np-eqn-nilmatstd} with 
 $C_1(\vec{v})=0, c_2(\vec{v})=0$ for 
all $\vec{v}$, we call $\CN$ the {\em standard cusp group} 
({\em of type $(n+1, i_0)$}) in the standard form. 
The {\em standard parabolic group} (of type $(n+1, i_0)$) 
is a  group conjugate to $\CN \hat \Ort(i_0)$, 
where $\CN$ is in the standard cusp group in the standard form. 
Groups conjugate to these are called {\em standard cusp groups} 
and {\em standard extended cusp groups} respectively. 
\index{cusp group!standard|textbf} 
\index{cusp group!standard!extended|textbf} 
\index{Ort@$\hat \Ort(i_0)$|textbf}

The assumptions for this subsection are as follows: 
We assume that the group satisfies the condition virtually only
since this is sufficient for our purposes. 
\begin{hyp}\label{np-h-norm} $ $
\begin{itemize} 
\item  Let $K$ be defined as above for an R-p-end $\tilde E$.
Assume that $N_K$ acts on $K^{o}$ cocompactly. 
\item $\bGamma_{\tilde E}$ satisfies the transverse weak middle-eigenvalue condition
for the R-p-end structure. 
\item  Under a common coordinate system, 
elements take the matrix form given in  \eqref{np-eqn-matstd} 
\item A group $\CN$ of the form \eqref{np-eqn-nilmatstd} in the same coordinate 
system 
acts on each hemisphere with boundary $\SI^i_{\infty}$, and fixes
$\mbv_{\tilde E} \in \SI^i_\infty$ with coordinates 
$\llrrparen{0, \cdots, 0, 1}$. 
\item  
$N \subset \CN$ in the same coordinate system as above. 
\item The p-end holonomy group $\bGamma_{\tilde E}$ normalizes 
$\CN$. 
\item $\CN$ acts on a p-end neighborhood $U$ of $\tilde E$, 
and acts on $U \cap \SI^{i_{0}+1}$
for each great sphere $\SI^{i_{0}+1}$ containing $\SI^{i_{0}}_{\infty}$
whenever  $U \cap \SI^{i_{0}+1} \ne \emp$. 
\item $\CN$ acts freely and transitively 
on the space of $i_{0}$-dimensional leaves of $\tilde \Sigma_{\tilde E}$
by an induced action. 
\end{itemize}
\end{hyp}

Let $U$ be a p-end neighborhood of $\tilde E$. 
Let $l'$ be an $i_{0}$-dimensional leaf of $\tilde  \Sigma_{\tilde E}$. 
The consideration of the projection $\Pi_{K}$ shows us that 
the leaf $l'$ corresponds to a hemisphere $H^{i_0+1}_{l'}$ where
\begin{equation} \label{np-eqn-HlU} 
U_{l'}:= (H^{i_0+1}_{l'} - \SI^{i_0}_\infty) \cap U \ne \emp
\end{equation} 
 holds. 

\begin{lemma}[Cusp]\label{np-lem-bdhoro} 
Assume Hypothesis \ref{np-h-norm}. 
Let $l'$ be an $i_{0}$-dimensional leaf of $\tilde  \Sigma_{\tilde E}$. 
Let $H^{i_{0}+1}_{l'}$ denote the $(i_{0}+1)$-dimensional hemisphere with boundary $\SI^{i_{0}}_{\infty}$ 
corresponding to $l'$.  
Then 
the action of $\CN$ on $H^{i_0+1}_{l'}$ has the form of the lower-right $(i_0+2)\times(i_0+2)$-submatrices
in \eqref{np-eqn-nilmatstd} under a coordinate system.
\end{lemma}
\begin{proof} 


Since $l'$ is an $(i_{0}+1)$-dimensional leaf of $\tilde \Sigma_{\tilde E}$, we obtain
$H^{i_{0}}_{l'} \cap U \ne \emp$. 
Let $J_{l'}:= H^{i_{0}+1}_{l'} \cap U \ne \emp$ where $\mathcal{N}$ acts on.

Now, $l'$ corresponds to an interior point of $K$. 
We need to change the coordinates of $\SI^{n-i_{0}-1}$ such that $l'$ goes to 
\[\underbrace{\llrrparen{0, \cdots, 0, 1}}_{n-i_0} \hbox{ under } \Pi_{K}.\]
This involves the coordinate changes of its first $n-i_{0}$ coordinates. 
Now, we can restrict $g$ to $H^{i_{0}+1}_{l'}$ such that the matrix form is 
with respect to $U_{l'}$. 
Now give a coordinate system on the open hemisphere $H^{i_0+1, o}_{l'}$ 
given by 
\[\underbrace{\llrrparen{1, \vec{x}, x_{i_0+1}}}_{i_0+2} 
\hbox{ for } \vec{x} \in \bR^{i_0} \hbox{ where }
\underbrace{\llrrparen{1, 0, \dots, 0}}_{i_0+2} \hbox{ is the origin,}\] and 
$(\vec{x}, x_{i_0+1})$ 
gives the affine coordinate system on $H^{i_0+1, o}_{l'}$. 

Using \eqref{np-eqn-nilmatstd} restricted to $\SI^{i_0}$, 
the lowest row of the lower-right $(i_0+1)\times (i_0+1)$-block
has to be of the form $(\ast, \vec{v}, 1)$. 
We obtain that 
each $g \in \CN$ then has the form in $H_{l'}^{i_0+1}$ as 
\[
\left(
\begin{array}{ccc}
1              & 0         &  0 \\
L(\vec{v}^T)  & \Idd_{i_0}  & 0  \\
\kappa(v) & \vec{v}  & 1  
\end{array}
\right)
\]
where $L: \bR^{i_0} \ra \bR^{i_0}$ is a linear map. The linearity of $L$ is a consequence of the group property. 
$\kappa: \bR^{i_0} \ra \bR$ is some function. We consider $L$ as an $i_0\times i_0$-matrix. 

Suppose that $L$ has a nontrivial kernel $K_1$. We use $t\vec{v} \in K_1-\{O\}$. 
As $t \ra \infty$, consider each orbit of the subgroup 
$\mathcal{N}(\bR\vec{v}) \subset \mathcal{N}$ given by 
\[\llrrparen{1, \vec{x}, x_{i_0+1}} \ra \llrrparen{1, \vec{x}, 
	\kappa(t\vec{v}) + t \vec{v}\cdot \vec{x} + x_{i_0+1}}.\]
This action fixes coordinates from the second to the $(i_0+1)$-st ones to be
$\vec{x}$. Hence, each orbit lies on an affine line from $\llrrparen{0, 0, \dots, 1}$. 
Since all eigenvalues of elements of $\CN$ equal $1$, 
the action is unipotent. 
Since the action is unipotent, either the action is trivial or 
the orbit is the entire complete affine line on $H^{i_0}_{l'}$. 
Since the action on each leaf $l'$ is free, the action cannot be trivial. 
Thus, the orbit is the complete affine line, and 
this contradicts the proper convexity of $\torb$. 



Also, since $\CN$ is abelian, the computations of 
\[\CN(\vec{v})\CN(\vec{w}) = \CN(\vec{w})\CN(\vec{v})\] 
show that $\vec v L \vec w^T = \vec w L \vec v^T$
for every pair of vectors $\vec v$ and $\vec w$ in $\bR^{i_0}$. 
Thus, $L$ is a symmetric matrix. 

We may obtain new coordinates $x_{n-i_0+1}, \dots, x_n$ by taking 
linear combinations of these. 
Since $L$ hence is nonsingular, we can find 
new coordinates $x_{n-i_0+1}, \dots, x_n$ such that $\CN$ now take the standard form: 
We conjugate $\CN$ by 
\[
\left(
\begin{array}{ccc}
1              & 0         &  0 \\
0 & A  & 0  \\
0 &  0 & 1  
\end{array}
\right)
\]
for nonsingular $A$. 
We obtain
\[
\left(
\begin{array}{ccc}
1              & 0         &  0 \\
AL\vec{v}^T  & \Idd_{i_0}  & 0  \\
\kappa(\vec{v})  & \vec{v} A^{-1}  & 1  
\end{array}
\right).
\]
We thus need to solve for $A^{-1} A^{-1 T} = L$.  
This is possible since $L$ is nonsingular and symmetric as shown above. 

Now, we conjugate as we wished to. 
We can factorize each element of $\CN$ into forms 
\[
\left(
\begin{array}{ccc}
1              & 0         &  0 \\
0  & \Idd_{i_0}  & 0  \\
\kappa(\vec{v}) - \frac{\llrrV{\vec{v}}^2}{2}  & 0  & 1  
\end{array}
\right)
\left(
\begin{array}{ccc}
1              & 0         &  0 \\
\vec{v}^T  & \Idd_{i_0}  & 0  \\
 \frac{\llrrV{\vec{v}}^2}{2}   & \vec{v}  & 1  
\end{array}
\right).
\]
Again, by the group property, $\aleph_7(\vec{v}) : = \kappa(\vec{v}) - \frac{\llrrV{\vec{v}}^2}{2}$ 
gives us a linear function $\aleph_7: \bR^{i_0} \ra \bR$. 
Hence $\aleph_7(\vec{v}) = \kappa_\alpha \cdot \vec{v}$
for $\kappa_\alpha \in \bR^{i_0}$. 
Now, we conjugate $\CN$ by 
the matrix 
\[
\left(
\begin{array}{ccc}
 1 & 0  & 0  \\
 0 & \Idd_{i_0}  & 0  \\
0  & -\kappa_\alpha  & 1  
\end{array}
\right), 
\]
and this puts $\CN$ into the standard form. 

Now it follows that the orbit of $\CN(x_0)$ for any point $x_0$ of $J_{l'}$ is an ellipsoid with a point removed:  
$\CN$ in the standard form can be recognized as 
that of the parabolic group in the hyperbolic space 
in the Klein model in an appropriate coordinates.
\end{proof}

Recall the standard cusp group from the end of 
Section \ref{np-subsub-matrixform}. 
For later purposes, we need: 

\begin{lemma} \label{np-lem-cuspgroup} 
	Let $C$ be a standard cusp group acting on a hemisphere $H$ of dimension
	$i_0+1$ with boundary $\SI^{i_0}$ fixing a point $p$ in $\SI^{i_0}$. 
	Then the following hold: 
	\begin{itemize} 
	\item There exists an affine space $\mathds{A}^n_C$ 
	with $\SI^{i_0}\subset \Bd \mathds{A}^n_C$ and $H^o \subset \mathds{A}^n_C$
	where orbits of points have three types\/{\em :}  
	\begin{itemize} 
	\item The orbit of each point in $\mathds{A}^n_C$ is an ellipsoid 
	in an affine subspace of dimension $i_0$ parallel to 
	the affine subspace $H^o$. 
	\item The orbit of each point of 
	a great sphere $\SI^{n-i_0-1} \subset \Bd \mathds{A}^n$ of dimension $(n-i_0-1)$ containing $p$ is transverse to $\SI^{i_0}$
	where the orbits are singletons. 
	\item The orbit of every point of $\Bd \mathds{A}^n_C - \SI^{n-i_0-1}$ is not 
	contained in a properly convex domain. 
\end{itemize} 
	\item The affine space $\mathds{A}^n_C$ with these orbit properties is 
	uniquely determined. 
	\end{itemize} 
	\end{lemma} 
\begin{proof} 
	We choose the affine space $\mathds{A}^n_C$ given by $x_{n-i_0}> 0$ for
	the coordinate system where $C$ is written
	as in \eqref{np-eqn-nilmatstd} with $c_i$, $i=1, \dots, i_0+1$, being zero. 
	Since $C$ is standard, 
	There exists a sphere of dimension $\SI^{n-i_0-2}$ complementary 
	to $\SI^{i_0}$ in $\Bd \mathds{A}^n_C$ where $C$ acts trivially. 
	Let $\SI^{n-i_0-1}$ be the join $\{p, p_-\} \ast \SI^{n-i_0-2}$. 
	On $\SI^{n-i_0-1}$ the orbits are singletons. 
	Each orbit of a point of $\Bd \mathds{A}^n_C - \SI^{n-i_0-1}$ can be understood 
	by the matrix form. Each orbit always contain 
	a pair of antipodal points in the closures by considering 
	$\CN$ with the $(n-i_0)$-th row and the $(n-i_0)$-th column removed.
	The affine translations commute with each element of $C$. 
	This shows that each orbit in $\mathds{A}^n_C$ is as claimed. 
	
	Also, since the orbit types are characterized, 
	$\mathds{A}^n_C$ is uniquely determined. 
\end{proof}

Let $a_5(g)$ denote $\left| \det(A^5_g)\right|^{\frac{1}{i_0}}$.
Define $\mu_g:= \frac{a_5(g)}{a_1(g)} = \frac{a_9(g)}{a_5(g)}$ for $g \in \bGamma_{\tilde E}$
from Lemma \ref{np-lem-similarity}. 

\begin{lemma}[Similarity] \label{np-lem-similarity}
Assume Hypothesis \ref{np-h-norm} and $g \in \bGamma_{\tilde E}$ is in 
the form of \eqref{np-eqn-matstd}. 
Then 
any element $g \in \bGamma_{\tilde E}$ 
induces an $(i_0\times i_0)$-matrix $M_g$ given by
\[g \CN(\vec{v}) g^{-1} = \CN(\vec{v}M_g) \hbox{ where } \] 
\[M_g = \frac{1}{a_1(g)} (A_5(g))^{-1} = \mu_g O_5(g)^{-1} \]
for $O_5(g)$ in a compact Lie group $G_{\tilde E}$, and 
the following hold. 
\begin{itemize} 
\item $(a_5(g))^2 = a_1(g) a_9(g)$ or equivalently $\frac{a_5(g)}{a_1(g)}= \frac{a_9(g)}{a_5(g)}$.
\item Finally, $a_1(g), a_5(g),$ and $a_9(g)$ are all positive. 
\end{itemize}  
\end{lemma} 
\begin{proof}
Since the conjugation by $g$ sends elements of $\CN$ to themselves in a one-to-one manner, 
the correspondence between the set of $\vec{v}$ for $\CN$ and $\vec{v'}$ is one-to-one.

Since we have $g \CN({\vec{v}}) = \CN({\vec{v}'}) g$ for vectors $\vec{v}$ and $\vec{v'}$ 
in $\bR^{i_0}$ by Hypothesis \ref{np-h-norm},
we consider
\begin{equation} \label{np-eqn-firstm}
\newcommand*{\tempV}{\multicolumn{1}{r|}{}}
\left( \begin{array}{ccccccc} 
S(g) & \tempV & s_1(g) & \tempV & 0 & \tempV & 0 \\ 
 \cline{1-7}
s_2(g) &\tempV & a_1(g) &\tempV & 0 &\tempV & 0 \\ 
 \cline{1-7}
C_1(g) &\tempV & a_4(g) &\tempV & A_5(g) &\tempV & 0 \\ 
 \cline{1-7}
c_2(g) &\tempV & a_7(g) &\tempV & a_8(g) &\tempV & a_9(g) 
\end{array} 
\right)
\left( \begin{array}{ccccccc} 
\Idd_{n-i_0-1} & \tempV & 0 & \tempV & 0 & \tempV & 0 \\ 
 \cline{1-7}
0 &\tempV & 1 &\tempV & 0 &\tempV & 0 \\ 
 \cline{1-7}
C_1({\vec{v}}) &\tempV & \vec{v}^T &\tempV & \Idd_{i_0} &\tempV & 0 \\ 
 \cline{1-7}
c_2({\vec{v}}) &\tempV & \frac{\llrrV{\vec{v}}^2}{2} &\tempV & \vec{v} &\tempV & 1 
\end{array} 
\right)
\end{equation}
where $C_1({\vec{v}})$ is an $(n-i_0-1)\times i_0$-matrix where
each row is a linear function of $\vec{v}$, 
$c_2({\vec{v}})$ is a $(n-i_0-1)$-row vector, and
$\vec{v}$ is an $i_0$-row vector. 
This must equal the following matrix
for some $\vec{v'}\in \bR$
\begin{equation} \label{np-eqn-secondm}
\newcommand*{\tempV}{\multicolumn{1}{r|}{}}
\left( \begin{array}{ccccccc} 
\Idd_{n-i_0-1} & \tempV & 0 & \tempV & 0 & \tempV & 0 \\ 
 \cline{1-7}
0 &\tempV & 1 &\tempV & 0 &\tempV & 0 \\ 
 \cline{1-7}
C_1({\vec{v'}}) &\tempV & \vec{v'}^T &\tempV & \Idd_{i_0} &\tempV & 0 \\ 
 \cline{1-7}
c_2({\vec{v'}}) &\tempV & \frac{\llrrV{\vec{v'}}^2 }{2} &\tempV & \vec{v'} &\tempV & 1 
\end{array} 
\right)
\left( \begin{array}{ccccccc} 
S(g) & \tempV & s_1(g) & \tempV & 0 & \tempV & 0 \\ 
 \cline{1-7}
s_2(g) &\tempV & a_1(g) &\tempV & 0 &\tempV & 0 \\ 
 \cline{1-7}
C_1(g) &\tempV & a_4(g) &\tempV & A_5(g) &\tempV & 0 \\ 
 \cline{1-7}
c_2(g) &\tempV & a_7(g) &\tempV & a_8(g) &\tempV & a_9(g) 
\end{array} 
\right).
\end{equation}
From \eqref{np-eqn-firstm}, 
we compute the $(4, 3)$-block of the result 
to be $a_8(g) + a_9(g) \vec{v}$. 
From \eqref{np-eqn-secondm}, 
the $(4, 3)$-block is
$\vec{v'} A_5(g) + a_8(g)$. We obtain the relation
$a_9(g) \vec{v} = \vec{v'} A_5(g)$ for every $\vec{v}$. 
Since the correspondence between $\vec{v}$ and $\vec{v'}$ is one-to-one, 
we obtain 
\begin{equation}\label{np-eqn-vp1}
\vec{v'} = a_9(g) \vec{v} (A_5(g))^{-1}
\end{equation}
for the $i_0\times i_0$-matrix $A_5(g)$ and we also infer $a_9(g) \ne 0$
and $\det(A_5(g)) \ne 0$. 
The $(3, 2)$-block of the result of \eqref{np-eqn-firstm} 
equals 
\[a_4(g) + A_5(g) \vec{v}^T. \]  
The $(3, 2)$-block of the result of \eqref{np-eqn-secondm} 
equals
\begin{equation}
C_1({\vec{v}'}) s_1(g) + a_1(g) \vec{v}^{\prime T} + a_4(g). 
\end{equation}
Thus,
\begin{equation} \label{np-eqn-sim0}
A_5(g) \vec{v}^T =  C_1({\vec{v}'}) s_1(g)  + a_1(g) \vec{v}^{\prime T}.
\end{equation} 


For each $g$, we can choose a coordinate system such that $s_1(g) = 0$ 
since by the Brouwer fixed point theorem, there is a fixed point in the compact 
convex set $K \subset \SI^{n-i_0-1}$.
This involves the coordinate changes of the first $n-i_0$ coordinate functions only. 
This may change the matrix form of $\CN$, which we do below. 



Let $l'$ denote one of the $i_0$-dimensional leaves. 
Since $\CN$ acts on $\SI^{i_0+1}_{l'}$ for $l'$ as a cusp group, 
Lemma \ref{np-lem-bdhoro} shows that 
there exists a coordinate change involving the last $(i_0+2)$-coordinates 
\[x_{n-i_0+1}, \dots, x_n, x_{n+1}\]
such that the matrix form of the lower-right $(i_0+2)\times(i_0+2)$-matrix of each element $\CN$ is of the standard cusp form.
This does not affect $s_1(g) = 0$ as we can check from conjugating matrices used in 
the proof of Lemma \ref{np-lem-bdhoro}
as the change involves the above coordinates only. 
Denote this coordinate system by $\Phi_{g, l'}$. 

We may assume that the transition to this coordinate system from the original one 
is uniformly bounded: First, we change for $\SI^{n-i_0-1}$ with a bounded elliptic 
coordinate change since we are only picking out a single point to be a coordinate 
axis. 

Also, for the later coordinate changes for $l'$, we need to use
 $L$ and $\kappa$ in the proof of Lemma \ref{np-lem-bdhoro} applied here 
to be uniformly bounded functions since we are working with only one $l'$ here. 
Hence, $A$ and $\kappa_\alpha$ in the proof applied here
are also uniformly bounded. 


Let us use $\Phi_{g, l'}$ for a while
using primes for new set of coordinates functions. 
Now $A'_5(g)$ is conjugate to $A_5(g)$ as we can check in the proof of Lemma \ref{np-lem-bdhoro}. 
Under this coordinate system for given $g$, 
we obtain $a'_1(g) \ne 0$ and we can recompute to show that 
$a'_9(g) \vec{v} = \vec{v'} A'_5(g)$ for every $\vec{v}$ as in \eqref{np-eqn-vp1}. 
By \eqref{np-eqn-sim0} recomputed for this case, we obtain 
\begin{equation}\label{np-eqn-vp2}
\vec{v'} = \frac{1}{a'_1(g)} \vec{v} (A'_5(g))^T
 \end{equation}
 as $s'_1(g) = 0$ here since we are using the coordinate system $\Phi_{g, l'}$.
 Since this holds for every $\vec{v} \in \bR^{i_0}$, 
 we obtain 
 \[a'_9(g) (A'_5(g))^{-1} = \frac{1}{a'_1(g)} (A'_5(g))^T.\] 
 Hence $\frac{1}{ |\det(A'_5(g))|^{1/i_0}} A'_5(g) \in \Ort(i_0)$. 
 Also, \[\frac{a'_9(g)}{a'_5(g)} = \frac{a'_5(g)}{a'_1(g)}.\] 

 Here, $A'_5(g)$ is a conjugate of the original matrix $A_5(g)$ by linear coordinate changes 
 as we can see from the above processes to obtain the new coordinate system.

 
 This implies that the original matrix $A_5(g)$ is conjugate to an orthogonal matrix multiplied by a positive scalar for every $g$. 
The set of matrices $ \{ A_5(g)| g \in \bGamma_{\tilde E}\}$ forms a group since every $g$ is of a standard matrix form 
(see \eqref{np-eqn-matstd}). 
  Given such a group of matrices normalized to have determinant $\pm 1$, we obtain 
 a compact group 
 \[G_{\tilde E}:= \Bigg\{ \frac{1}{|\det A_5(g)|^{\frac{1}{i_0}}} A_5(g) \Bigg| g \in \Gamma_{\tilde E} \Bigg\}\]
  by Lemma \ref{np-lem-cpt}. 
This group has a coordinate system where every element is orthogonal by a coordinate change
of $x_{n-i_0+1}, \dots, x_n$. 


Also, $a_1(g), a_5(g), a_9(g)$ are conjugation invariant. Hence, 
we proved 
\begin{equation}\label{np-eqn-a9a5a1}  
\frac{a_9(g)}{a_5(g)} = \frac{a_5(g)}{a_1(g)} \left(= \mu_g\right)
\end{equation} 

We have $a_9(g) = \lambda_{\mbv_{\tilde E}}(g) > 0$. 
Since $a_5(g)^2= a_1(g) a_9(g)$, we obtain $a_1(g) > 0$. Finally, 
$a_5(g) > 0$ by definition.
\end{proof} 


\begin{lemma}\label{np-lem-cpt} 
Suppose that $G$ is a subgroup of a linear group $\GL(i_0, \bR)$ where each 
element is conjugate to an orthogonal element by a uniformly bounded set of 
conjugating matrices. Then $G$ is in a compact Lie group. 
\end{lemma}
\begin{proof} 
Clearly, the norms of eigenvalues of $g \in G$ are all $1$. Hence, 
$G$ is virtually an orthopotent group by Theorem \ref{prelim-thm-orthopotent}. 
(See \cite{Moore68} and \cite{CG74}).) 
Hence, $\bR^{i_0}$ has subspaces 
$\{0\} =V_0\subset V_1 \subset \cdots \subset V_m = \bR^{i_0}$ 
where $G$ acts as orthogonal on $V_{i+1}/V_{i}$ up to a choice of coordinates. 
Hence, the Zariski closure $\mathcal{Z}(G)$ of $G$ is also orthopotent. 


If $G$ is discrete, Theorem \ref{prelim-thm-orthopotent} shows that 
$G$ is virtually unipotent. The unipotent subgroup of $G$ is trivial 
since the elements must be conjugate to orthogonal elements. 
Thus, $G$ is a finite group, and we finished the proof.


Suppose now that the closure $\bar G$ of $G$ in $\GL(i_0, \bR)$ 
is a Lie group of dimension $\geq 1$. 
Each element of the identity component $\bar G$ is again elliptic or is the identity 
since each accumulation point is the limit of a sequence of 
elliptic elements conjugated by a uniformly bounded collection of elements. 
%
%
Let $\Ort(\oplus_{i=1}^m V_i/V_{i-1})$ denote the group of linear transformations 
acting on each $V_i/V_{i-1}$ orthogonally for each $i=1, \dots, m$. 
By Theorem \ref{prelim-thm-orthopotent}, 
there is a homomorphism $\mathcal{Z}(G) \ra \Ort(\oplus_{i=1}^m V_i/V_{i-1})$
whose kernel $U_{i_0}$ is the a group of unipotent matrices. 
Let $\hat G$ denote the image of $\bar G$ in the second group.
Then $\bar G \cap U_{i_0}$ is trivial since every element of $\bar G$ is elliptic or is the identity.
Thus, $\bar G$ is isomorphic to a compact group $\hat G$. 
\end{proof}


From now on,
we denote by $(C_1({\vec{v}}), \vec{v}^T)$ the matrix obtained 
from $C_1(\vec{v})$ by adding a column vector $\vec{v}^T$
and denote 
$O_5(g):= \frac{1}{|\det A_5(g)|^{\frac{1}{i_0}}} A_5(g)$. 
We also let 
\[ \hat S(g) := \left(\begin{array}{cc} 
S(g), & s_1(g) \\ s_2(g), & a_1(g) \end{array} \right). 
\]


\begin{lemma}[$K$ is a cone] \label{np-lem-conedecomp1}
Assume Hypothesis \ref{np-h-norm}. 
Suppose that $\bGamma_{\tilde E}$ acts semisimply on $K^o$. 
Then the following hold{\rm :} 
\begin{itemize} 
\item $K$ is a cone over a totally geodesic 
$(n-i_0-2)$-dimensional domain $K''$. 
\item The rows of $(C_1({\vec{v}}), \vec{v}^T)$ are proportional to a single vector, and 
we can find a coordinate system where $C_1({\vec v}) = 0$
not changing any entries of the lower-right $(i_0+2)\times (i_0+2)$-submatrices for all $\vec v \in \bR^{i_0}$. 
\item We can find a common coordinate system where 
\begin{align}\label{np-eqn-O5coor}
& O_5(g)^{-1} = O_5(g)^T, O_5(g) \in \Ort(i_0), \nonumber \\ 
& s_1(g) = s_2(g) = 0 \hbox{ for all } g \in \bGamma_{\tilde E}
\end{align} 
where $O_5(g) = \left| \det(A^5_g)\right|^{\frac{1}{i_0}} A_5(g)$. 
\item In this coordinate system, we have
\begin{multline} \label{np-eqn-conedecomp1}
s_1(g)=0, s_2(g)=0,  \\ 
 a_9(g) c_2({\vec{v}})  
 = c_2({\mu_g\vec{v}O_5(g)^{-1}}) S(g)  + \mu_g \vec{v} O_5(g)^{-1} C_1(g)
\hbox{ for all } g\in \bGamma_{\tilde E}.
\end{multline} 


\end{itemize}
\end{lemma}
\begin{proof} 

The assumption implies that $M_g =  \mu_g O_5(g)^{-1}$ by Lemma \ref{np-lem-similarity}.
We consider the equation 
\begin{equation} \label{np-eqn-conj0}
g \CN(\vec{v}) g^{-1} = \CN(\mu_g \vec{v} O_5(g)^{-1}).
\end{equation}



We change to 
\begin{equation}\label{np-eqn-gN}
g \CN(\vec v) = \CN(\mu_g \vec v O_5(g)^{-1}) g.
\end{equation}
Considering the lower left $(n-i_0) \times (i_0+1)$-matrix of the left side of  \eqref{np-eqn-gN},  
we obtain 
\begin{equation} \label{np-eqn-mm1}
\left(\begin{array}{cc} 
C_1(g), & a_4(g) \\ c_2(g), & a_7(g) \end{array} \right) 
+ 
\left( \begin{array}{cc} 
a_5(g) O_5(g) C_1({\vec{v}}), & a_5(g) O_5(g) \vec{v}^T\\ 
a_8(g) C_1({\vec{v}}) + a_9 c_2({\vec{v}}),  
& a_8(g) {\cdot} \vec{v}^T + a_9(g) \vec{v}\cdot\vec{v}/2
\end{array} \right) \end{equation}
where the entry sizes are clear. 
From the right side of \eqref{np-eqn-gN}, we obtain
\begin{equation} \label{np-eqn-mm2} 
\begin{split} 
\left( \begin{array}{cc} 
C_1({\mu_g \vec{v} O_5(g)^{-1}}), & \mu_g O_5(g)^{-1, T}\vec{v}^T \\ 
c_2({\mu_g\vec{v}}O_5(g)^{-1}),     & \vec{v}\cdot\vec{v}/2
\end{array} \right) \hat S(g) + \\
 \left(\begin{array}{cc} 
C_1(g), & a_4(g) \\ \mu_g\vec{v} O_5(g)^{-1} \cdot C_1(g) + c_2(g), & a_7(g) + \mu_g \vec{v} O_5(g)^{-1} \cdot a_4(g) 
\end{array} \right).
\end{split}
\end{equation}
From the top rows of \eqref{np-eqn-mm1} and \eqref{np-eqn-mm2},  we obtain that 
\begin{equation}\label{np-eqn-C1} 
\begin{split}
 \bigg(a_5(g) O_5(g) C_1({\vec{v}}) , \, & a_5(g) O_5(g) \vec{v}^T \,\bigg) = \\
& \bigg( \mu_g C_1\left({\vec{v}}  O_5(g)^{-1}\right),  \mu_g O_5(g)^{-1, T} \vec{v}^T  \bigg) \hat S(g). 
 \end{split}
 \end{equation} 
 We multiply the both sides by $O_5(g)^{-1}$ from the left 
and by $\hat S(g)^{-1}$ from the right to obtain
 \begin{equation}\label{np-eqn-C2} 
 \begin{split} 
\bigg(a_5(g) C_1({\vec{v}}) , \, & a_5(g) \vec{v}^T\bigg) \hat S(g^{-1})
=   \\
& \bigg(\mu_g O_5(g)^{-1} C_1({\vec{v}}  O_5(g)^{-1} ),  \mu_g O_5(g)^{-1}  O_5(g)^{-1, T} \vec{v}^T\,\bigg).
\end{split}
\end{equation}
Let us form the subspace $V_C$ in the dual sphere $\bR^{n-i_0 \ast}$ spanned by 
row vectors of 
\[\{(C_1({\vec{v}}), \vec{v}^T)| \vec{v}\in \bR^{i_0}\}.\] 
Let $\SI_C^\ast$ denote the corresponding subspace 
in $\SI^{n-i_0-1 \ast}$. Then 
\[\left\{\frac{1}{\det \hat S(g)^{\frac{1}{n-i_0-1}}}\hat S(g) \vert g \in \bGamma_{\tilde E}\right\}\] 
acts on $V_C$ as a group of bounded 
linear automorphisms by \eqref{np-eqn-C2} 
since $O_5(g) \in G$ for a compact group $G$.
Therefore, $\{\hat S(g)| g\in \bGamma_{\tilde E}\}$ on $\SI_C^\ast$ is in a compact group of 
projective automorphisms.  

We recall that the dual group $N_K^*$ of $N_K$ acts on the properly convex dual domain $K^*$ of $K$ cocompactly 
by Proposition \ref{prelim-prop-sweepduality}.
Notice that a compact irreducible group 
cannot act cocompactly on a convex open set unless the set is a singleton. 
Since $N_K^\ast$ acts as a compact group on $\SI_C^\ast$, 
it must be that $N_K^*$ is reducible. 





Now, we apply the theory of Vey \cite{Vey} and Benoist \cite{Benoist03}: 
Since $N_K^*$ is semisimple by above premises, 
$N_K^*$  acts on 
a complementary subspace of $\SI^{\ast}_N$. 
By Proposition \ref{prelim-prop-sweep}, 
$K^*$ has an invariant subspace $K^*_1$ and $K^*_2$ 
such that we have the strict join 
 \[K^* = K^*_1 \ast K^*_2 \hbox{ where }  \dim K^*_1 = \dim \SI_C^{\ast}, \dim K^*_2 = \dim \SI_N^{\ast}\] 
where
\[K^*_1 = K^* \cap \SI_C^{\ast}, K^*_2 = K^{\ast}\cap \SI_N^{\ast}.\] 
Also, $N_K^*$ is isomorphic to a subgroup of 
\[N_{K, 1}^*\times N_{K, 2}^* \times A\] 
where 
\begin{itemize} 
\item $A$ is a diagonalizable subgroup with positive eigenvalues only isomorphic 
to a subgroup of $\bR_+$, 
\item $N_{k, i}^*$ is the restriction image of $N_K^\ast$ to $K_i^\ast$ for $i=1,2$, and
\item $N_{K, i}^*$ acts on 
the interior of $K^*_i$ properly and cocompactly. 

\end{itemize} 
But since $N_{K, 1}^*$ acts orthogonally on $\SI_C^*$, 
the only possibility is that $\dim \SI_C^* = 0$: 
Otherwise $K^{\ast o}_1/N_K$ is not compact contradicting Proposition 
\ref{prelim-prop-sweep}. 
Hence, $\dim \SI_{C}^* =0$
and $K_1^\ast$ is a singleton and $K_2^\ast$ is $n-i_0-2$-dimensional properly
convex domain. 

Rows of $(C_{1}(\vec{v}), \vec{v}^{T})$ are elements of the $1$-dimensional subspace in $\bR^{n-i_{0}-1\ast}$ 
corresponding to $\SI_{C}^{\ast}$. 
Therefore this shows that the rows of $(C_1({\vec{v}}), \vec{v}^T)$ are proportional to a single row vector.

Since $(C_1({\vec{e}_j}), \vec{e}_j^T)$ has $0$ as the last column element except 
for the $j$th one, only the $j$th row of $C_{1}(\vec{e}_{j})$ is nonzero. 
Let $C_1({1, \vec{e}_1})$ be the first row of $C_1({\vec{e}_1})$. 
Thus, each row of $(C_{1}(\vec{e}_{j}), \vec{e}_{j}^{T})$ equals to a scalar multiple of 
$(C_1({1, \vec{e}_1}), 1)$ for every $j$. 
Now we can choose coordinates of $\bR^{n-i_0 \ast}$ such that 
this $(n-i_0)$-row-vector now has coordinates $(0, \dots, 0, 1)$. 
We can also choose the coordinates such that $K^*_2$ is 
in the zero set of the last coordinate.  
With this change, we need to do conjugation by matrices 
with the top left $(n-i_0-1)\times (n-i_0-1)$-submatrix being different from $\Idd$ and 
the rest of the entries staying the same. 
This does not affect the expressions of matrices of 
 lower right $(i_0+2)\times (i_0+2)$-matrices involved here. 
Thus, $C_1({\vec{v}}) =0$ in this coordinate for all $\vec{v} \in \bR^{i_0}$.
Also, $\underbrace{ \llrrparen{0, \dots, 0, 1}}_{n-i_0}$ is an eigenvector of every element of $N_{K}^{\ast}$. 

The hyperspace containing $K^{\ast}_{2}$ is also $N_{K}^{\ast}$-invariant. 
Thus, the $(n-i_0)$-vector $(0, \dots, 0, 1)$ corresponds to an eigenvector 
of every element of $N_{K}$. 
In this coordinate system, $K$ is a strict join of a point 
for an $(n-i_0)$-vector 
\[ k = \underbrace{ \llrrparen{0, \dots, 0, 1}}_{n-i_0}\]
and a domain $K''$ given by
setting $x_{n-i_0} = 0$ in a totally geodesic 
sphere of dimension $n-i_0-2$ by duality. 
We also obtain
\[s_1(g) =0, s_2(g) =0.\]

For the final item we have under our coordinate system. 
\begin{equation} \label{np-eqn-form1}
g = \left( \begin{array}{cccc} 
 S(g) & 0 & 0 & 0 \\ 
 0 & a_1(g) & 0 & 0 \\ 
C_1(g) & a_4(g) & a_5(g) O_5(g) & 0  \\
c_2(g) & a_7(g) & a_8(g) & a_9(g) 
\end{array} \right), 
\end{equation} 
\begin{equation} \label{np-eqn-form2}
 \CN(\vec{v}) = \left( \begin{array}{cccc} 
 \Idd_{n-i_0-1} & 0 & 0 & 0 \\ 
 0 &    1 & 0 & 0 \\ 
0 & \vec{v}^T & \Idd & 0 \\ 
 c_2({\vec{v}}) & \frac{1}{2} \llrrV{\vec{v}}^2 & \vec{v} & 1 
 \end{array} \right). 
 \end{equation} 
  
 The normalization of $\CN$ shows as in the proof of Lemma \ref{np-lem-similarity} that $O_5(g)$ is orthogonal now.
 (See  \eqref{np-eqn-vp1} and  \eqref{np-eqn-sim0}.)
By \eqref{np-eqn-conj0}, we have 
 \[g\CN(\vec{v}) = \CN(\vec{v'})g, \vec{v}'= \mu_g \vec{v} O_5(g)^{-1}.\]
 We consider the lower-right $(i_0+1)\times (n-i_0)$-submatrices of 
 $g\CN(\vec{v})$ and $\CN(\vec{v'})g$. 
 For the first one, we obtain 
 \[
 \left( \begin{array}{cc} C_1(g), & a_4(g) \\ c_2(g), & a_7(g) \end{array} \right ) 
 + \left( \begin{array}{cc} a_5(g) O_5(g), & 0 \\ a_8(g), & a_9(g) \end{array} \right) 
 \left( \begin{array}{cc} 0, & \vec{v}^T \\ c_2({\vec{v}}), &  \frac{1}{2} \llrrV{\vec{v}}^2 \end{array} \right) 
 \]
 For $\CN(\vec{v'})g$, we obtain 
  \[
  \left( \begin{array}{cc} 0, & \vec{v'}^T \\ c_2({\vec{v'}}), & \frac{1}{2} \llrrV{\vec{v}'}^2 \end{array} \right) 
 \left( \begin{array}{cc} S(g), & 0 \\ 0, & a_1(g) \end{array} \right ) 
 + \left( \begin{array}{cc}  \Idd, & 0 \\ \vec{v'}, & 1 \end{array} \right) 
 \left( \begin{array}{cc} C_1(g), & a_4(g) \\ c_2(g), & a_9(g) \end{array} \right).
 \]
Considering $(2, 1)$-blocks, we obtain 
\[ c_2(g) + a_9(g) c_2({\vec{v}}) = c_2({\vec{v'}}) S(g) 
+ \vec{v'} C_1(g) 
+ c_2(g). \]

\end{proof}



\begin{lemma} \label{np-lem-matrix}
Assume Hypothesis \ref{np-h-norm} and that $N_K$ acts semisimply. 
Then we can use the coordinates found in Lemma \ref{np-lem-conedecomp1}
such that the following holds for all $g$\,{\rm :} 
\begin{align}
\frac{a_{9}(g)}{a_{5}(g)} O_5(g)^{-1} a_4(g) &= a_8(g)^{T} \hbox{ or }\frac{a_{9}(g)}{a_{5}(g)}   a_4(g)^T O_5(g) = a_8(g).   \label{np-eqn-a4a8} \\
\hbox{If } \mu_{g} =1, \hbox{ then } &
a_1(g) = a_9(g)  = \lambda_{\mbv_{\tilde E}}(g) \hbox{ and } 
A_5(g) = \lambda_{\mbv_{\tilde E}}(g) O_5(g). \label{np-eqn-mug} 
\end{align}
\end{lemma} 
\begin{proof} 

Again, we use \eqref{np-eqn-firstm} and \eqref{np-eqn-secondm}. 
We only need to consider lower right $(i_0+2)\times (i_0+2)$-matrices. 
\begin{align} 
&\left( \begin{array}{ccc} 
a_1(g) & 0 & 0 \\ a_4(g) & a_5(g) O_5(g)  & 0 \\ a_7(g) & a_8(g) & a_9(g) \end{array} \right) 
\left(\begin{array}{ccc} 
1 & 0 & 0 \\ \vec{v}^T & \Idd & 0 \\ \frac{1}{2}\llrrV{\vec{v}}^2 & \vec{v} & 1 \end{array} \right) \\
& = \left( \begin{array}{ccc} 
a_1(g) & 0 & 0 \\ a_4(g) + a_5(g)  O_5(g) \vec{v}^T & a_5(g) O_5(g) & 0 
\\ a_7(g) + a_8(g) \vec{v}^T + \frac{a_9(g)}{2} \llrrV{\vec{v}}^2 & a_8(g) + a_9(g) \vec{v} & a_9(g) 
\end{array} \right). 
\end{align}
This equals 
\begin{align}
& \left( \begin{array}{ccc} 
1 & 0 & 0 \\ \vec{v'}^T & \Idd & 0 \\ \frac{1}{2}\llrrV{\vec{v'}}^2 & \vec{v'} & 1 \end{array} \right) 
\left( \begin{array}{ccc} 
a_1(g) & 0 & 0 \\ a_4(g) & a_5(g) O^5_g & 0 \\ a_7(g) & a_8(g) & a_9(g) \end{array} \right) \\
& = \left(\begin{array}{ccc} 
a_1(g) & 0 & 0 \\ a_1(g) \vec{v'}^T + a_4(g) & a_5(g) O_5(g) & 0 \\ 
\frac{a_1(g)}{2} \llrrV{\vec{v'}}^2 + \vec{v'} a_4(g) + a_7(g) & 
a_5(g) \vec{v'} O_5(g) + a_8(g) & a_9(g) \end{array} \right).
\end{align}
Then by comparing the $(3, 2)$-blocks, 
we obtain 
\[a_8(g) + a_9(g) \vec{v} = a_8(g) + a_5(g) \vec{v'} O_5(g) .\]
Thus, $\vec{v} =  \frac{a_{5}(g)}{a_{9}(g)}\vec{v'} O_5(g).$

From the $(3, 1)$-blocks, we obtain 
\[ a_1(g) \vec{v'}\cdot\vec{v'}/2 + \vec{v'}a_4(g) = a_8(g)\vec{v}^T + a_9(g) \vec{v}\cdot\vec{v}/2. \]
Since the quadratic forms have to equal each other, we obtain 
\[\frac{a_{9}(g)}{a_{5}(g)}\vec{v} O_5(g)^{-1} \cdot a_4(g) = \vec{v} \cdot a_8(g) \hbox{ for all } \vec{v} \in \bR^{i_{0}}.\] 
Thus, $\frac{a_{9}(g)}{a_{5}(g)}(O_5(g)^T a_4(g))^T = a_8(g)^T$.

Since we have $\mu_g = 1$, we obtain 
\[a_1(g) = a_9(g) = a_5(g) = \lambda_{\mbv_{\tilde E}}(g)
\hbox{ and }
A_5(g) = \lambda_{\mbv_{\tilde E}}(g) O_5(g)\] 
by Lemma \ref{np-lem-similarity}. 
Also, $a_1(g) = a_9(g) = a_5(g) = \lambda_{\mbv_{\tilde E}}(g)$.
\end{proof}

Under Hypothesis \ref{np-h-norm} and assuming that 
$N_K$ acts semisimply, 
we conclude 
by \eqref{np-eqn-a4a8}
and \eqref{np-eqn-a9a5a1}  
that each $g \in \bGamma_{\tilde E} $ has the form 
\begin{equation} \label{np-eqn-formgi}
\newcommand*{\tempV}{\multicolumn{1}{r|}{}}
\left( \begin{array}{ccccccc} 
S(g) & \tempV & 0 & \tempV & 0 & \tempV & 0 \\ 
 \cline{1-7}
0 &\tempV & a_{1}(g) &\tempV & 0 &\tempV & 0 \\ 
 \cline{1-7}
C_1(g) &\tempV & a_{1}(g) \vec{v}^T_g &\tempV & a_{5}(g) O_5(g) &\tempV & 0 \\ 
 \cline{1-7}
c_2(g) &\tempV & a_7(g) &\tempV &  a_{5}(g) \vec{v}_g O_5(g) &\tempV & a_{9}(g)
\end{array} 
\right)
\end{equation}
defining $\vec{v}_g:= \frac{a_4(g)}{a_1(g)}$.  

\begin{remark} \label{np-rem-mug} 
Since the matrices are of form \eqref{np-eqn-formgi}, 
$g \mapsto \mu_g$ is a homomorphism. 
\end{remark} 


\begin{corollary}\label{np-cor-formg2} 
	If $g$ of form \eqref{np-eqn-formgi} centralizes a Zariski dense subset $A'$ of $\CN$, 
	then $\mu_{g}=1$ and $O_5(g) = \Idd_{i_0}$. 
\end{corollary} 
\begin{proof} 
	$\CN$ is isomorphic to $\bR^{i_{0}}$.
	The subset $A''$ of $\bR^{i_{0}}$ corresponding to $A'$ is also Zariski dense in $\bR^{i_{0}}$. 
	$g \CN(\vec{v}) = \CN(\vec{v}) g$ shows that 
	$\vec{v} = \vec{v} O_5(g)$ for all $\vec{v} \in A''$. 
	Hence $O_5(g) =\Idd$. 
\end{proof}

\subsubsection{Invariant $\aleph_7$.} \label{np-sub-alpha7} 
We assume $\mu_g=1$, $g \in \bGamma_{\tilde E}$, 
identically in this subsubsection. 
When $\mu_{g}=1$ for all $g \in \bGamma_{\tilde E}$, by 
taking a finite-index subgroup of $\bGamma_{\tilde E}$, 
we conclude that each $g \in \bGamma_{\tilde E} $ has the following form 
by Lemma \ref{np-lem-matrix}
\begin{equation} \label{np-eqn-formgii}
\newcommand*{\tempV}{\multicolumn{1}{r|}{}}
M(g):= 
\left( \begin{array}{ccccccc} 
S(g) & \tempV & 0 & \tempV & 0 & \tempV & 0 \\ 
 \cline{1-7}
0 &\tempV & \lambda_{\mbv_{\tilde E}}(g) &\tempV & 0 &\tempV & 0 \\ 
 \cline{1-7}
C_1(g) &\tempV & \lambda_{\mbv_{\tilde E}}(g)\vec{v}^T_g &\tempV & \lambda_{\mbv_{\tilde E}}(g) O_5(g) &\tempV & 0 \\ 
 \cline{1-7}
c_2(g) &\tempV & a_7(g) &\tempV & \lambda_{\mbv_{\tilde E}}(g) \vec{v}_g O_5(g) &\tempV & \lambda_{\mbv_{\tilde E}}(g)
\end{array} 
\right).
\end{equation}

We define an invariant: 
\[\aleph_7(g):= \frac{a_7(g)}{\lambda_{\mbv_{\tilde E}}(g)} - \frac{\llrrV{\vec{v}_g}^2}{2}
\hbox{ for } g\in \bGamma_{\tilde E}. \]
We denote by $\hat M(g)$ the lower right $(i_0+1)\times(i_0+1)$-submatrix of $M(g)$
for $g \in \bGamma_{\tilde E}$. 
An easy computation shows that 
$\hat M(g)\hat M(h) = \hat M(gh)$ where 
$\vec{v}_{gh} = \vec{v} + O_5(g)\vec{v}_h$ holds 
for $g, h \in \bGamma_{\tilde E}$.
Then it is easy to show that 
\begin{equation} \label{eqn:alpha_add}
\aleph_7(gh) = \aleph_7(g) + \aleph_7(h), \hbox{ whenever }g, h
\in \bGamma_{\tilde E}.
\end{equation} 
Hence, we obtain a homomorphism to the additive group $\bR$ 
\begin{equation}
\aleph_7: \bGamma_{\tilde E} \ra \bR.
\end{equation} 
(See \eqref{np-eqn-kernel}.)

Here $\aleph_7(g)$ is also determined by factoring  
the matrix of $g$ into commuting matrices of form
\renewcommand{\arraystretch}{1.5}
\begin{equation}\label{np-eqn-kernel}
\newcommand*{\tempV}{\multicolumn{1}{r|}{}}
\left( \begin{array}{ccccccc}
\Idd_{n-i_0-1} & \tempV & 0 & \tempV & \, 0 & \tempV & \, 0 \\   
\cline{1-7}
0                    &\tempV & 1 & \tempV & \, 0 & \tempV & \, 0  \\ 
\cline{1-7}
0                    & \tempV &   0 &\tempV & \, \Idd_{i_0} &\tempV & \, 0 \\ 
\cline{1-7}
0                    &\tempV &  \aleph_7(g) &\tempV & \, \vec{0} &\tempV & \, 1 \\ 
\end{array} 
\right) 
\left( \begin{array}{ccccccc}
S_g & \tempV & 0 & \tempV & 0 &\tempV & 0 \\ 
\cline{1-7}
0 & \tempV & \lambda_{\mbv_{\tilde E}}(g) & \tempV & 0 & \tempV & 0  \\ 
\cline{1-7}
C_1(g) & \tempV & \lambda_{\mbv_{\tilde E}}(g) \vec{v}_g &\tempV & \lambda_{\mbv_{\tilde E}}(g) O_5(g) &\tempV & 0\\
\cline{1-7} 
c_2(g) & \tempV & \lambda_{\mbv_{\tilde E}}(g)  \frac{\llrrV{\vec{v_g}}^2}{2}  
&\tempV & \lambda_{\mbv_{\tilde E}}(g) \vec{v}_g O_5(g) &\tempV & \lambda_{\mbv_{\tilde E}}(g) \\ 
\end{array} 
\right).
\end{equation}

\begin{remark} \label{np-rem-alpha7}
	We give a bit more explanations. 
	Recall that the space of segments in a hemisphere $H^{i_0+1}$ with the vertices $\mbv_{\tilde E}, \mbv_{\tilde E-}$ 
	forms an affine subspace $\mathds{A}^i$ one-dimension lower, and the group $\Aut(H^{i_0+1})_{\mbv_{\tilde E}}$ of projective automorphisms 
	of the hemisphere fixing $\mbv_{\tilde E}$ maps to $\Aff(\mathds{A}^{i_0})$ with kernel $K$ equals the group of 
	transformations of the $(i_0+2)\times (i_0+2)$-matrix forms
	\renewcommand{\arraystretch}{1.5}
	\begin{equation}\label{np-eqn-kernel2}
	\newcommand*{\tempV}{\multicolumn{1}{r|}{}}
	\left( \begin{array}{ccccc} 
	1 & \tempV & \vec{0}^T & \tempV & 0  \\ 
	\cline{1-5}
	O &\tempV & \Idd_{i_0} &\tempV & \vec{0} \\ 
	\cline{1-5}
	b &\tempV & \vec{0}^T &\tempV & 1 \\ 
	\end{array} 
	\right)
	\end{equation}
	where $\mbv_{\tilde E}$ is given coordinates $\llrrparen{0, 0, \dots, 1}$, 
	$\vec{0}$ denote the $0$-vector in $\bR^{i_0}$ and 
	a center point of $H^{i_0+1}_l$ the coordinates $\llrrparen{1, 0, \dots, 0}$. 
	In other words, the transformations are of form 
	\begin{align}\label{np-eqn-temp}
	\left[
	\begin{array}{c}
	1\\
	x_1  \\
	\vdots \\ 
	x_{i_0} \\ 
	x_{i_0+1}   
	\end{array}
	\right]
	\mapsto 
	\left[
	\begin{array}{c}
	1 \\
	x_1\\
	\vdots \\ 
	x_{i_0} \\ 
	x_{i_0+1}+b
	\end{array}
	\right]
	\end{align}
	and hence $b$ determines the kernel element. 
	Hence $\aleph_7(g)$ indicates the translation towards $\mbv_{\tilde E}=\llrrparen{0,\dots, 1}$. (Actually the vertex corresponds to 
	$(1, 0, \dots, +\infty)$-point in this view. )
\end{remark}

We denote by $T(n+1, n-i_0)$ the 
group of matrices restricting to \eqref{np-eqn-kernel2} in the lower-right
$(i_0+2)\times (i_0+2)$-submatrices and equal to $\Idd$ on upper-left $(n-i_0-1)\times(n-i_0-1)$-submatrices 
and zero elsewhere. 
\index{Tnplusone@$T(n+1, n-i_0)$|textbf}


\subsubsection{Splitting the NPNC end} \label{np-subsub-split}

%






\begin{proposition}[Splitting] \label{np-prop-decomposition}
Assume Hypothesis \ref{np-h-norm} for 
$\bGamma_{\tilde E }$. 
Suppose additionally the following{\rm :} 
\begin{itemize} 
	\item Suppose that $N_K$ acts on $K$ in a semisimple manner. 
\item $K = \{k\} \ast K''$ a strict join, and $K^{o}/N_{K}$ is compact
with $k$ a common fixed point of $N_K$. 
\item Let $H$ be a commutant of a finite-index subgroup of $N_K$ that is positive 
diagonalizable. Assume that $N_K \cap H$ contains a free abelian group of 
rank $l_0-1$ provided $K''$ is a strict join of compact convex subsets $K_1, \dots, K_{l_0}$
where $H$ acts trivially on each $K_j$, $j=1, \dots, l_0$. 
\end{itemize} 
Then the following hold\/{\rm :}
\begin{itemize} 
\item there exists an exact sequence 
\[ 1 \ra \mathcal{N} \ra \langle \bGamma_{\tilde E}, \mathcal{N} \rangle \stackrel{\Pi_{K}^\ast}{\longrightarrow} N_K \ra 1. \]
\item $K''$ embeds projectively in the closure of $\Bd \torb$ whose image is 
$\bGamma_{\tilde E}$-invariant, and 
\item one can find a coordinate system such that every $\CN(\vec v)$ for 
$\vec v \in \bR^{i_0}$ is in the standard form and each element $g$ of 
$\bGamma_{\tilde E}$ is written such that
\begin{itemize}
\item $C_1(\vec v)=0, c_2({\vec v})=0$, and 
\item $C_1(g)=0$ and $c_2(g) = 0$.
\end{itemize}
\end{itemize} 
\end{proposition}
\begin{proof}
(A) Let $Z$ denote $\langle \bGamma_{\tilde E}, \mathcal{N} \rangle$.
Since $N \subset \mathcal{N}$, 
we have homomorphism
\[ Z \, \stackrel{\Pi^\ast_K}{\longrightarrow} \, N_K \ra 1 \]
extending $\Pi_K^\ast$ of \eqref{np-eqn-exact}.

Since the kernel of $\Pi_K^\ast| Z$ is generated by $\mathcal{N}$ and $N$, 
we proved the first item. 

Now, we set up some matrix notations: 
The function $\lambda_{\mbv_{\tilde E}}: \bGamma_{\tilde E} \ra \bR_+$ extends to 
$\lambda_{\mbv_{\tilde E}}: Z \ra \bR_+$. 
By \eqref{np-eqn-formgi}, we deduce that 
every element $g$ of $Z$ is of form: 
\begin{equation} \label{np-eqn-formgiiZ}
\newcommand*{\tempV}{\multicolumn{1}{r|}{}}
\left( \begin{array}{ccccccc} 
S(g) & \tempV & 0 & \tempV & 0 & \tempV & 0 \\ 
\cline{1-7}
0 &\tempV & a_{1}(g) &\tempV & 0 &\tempV & 0 \\ 
\cline{1-7}
C_1(g) &\tempV & a_{1}(g) \vec{v}^T_g &\tempV & a_{5}(g) O_5(g) &\tempV & 0 \\ 
\cline{1-7}
c_2(g) &\tempV & a_7(g) &\tempV &  a_{5}(g) \vec{v}_g O_5(g) &\tempV & a_{9}(g)
\end{array} 
\right)
\end{equation}
for some functions $C_1, c_2, a_i:Z \ra \bR$ for $i=1, 5, 9$.
Notice that 
$\aleph_7$ is identical zero on $\mathcal{N}$. 
Since $N \subset \mathcal{N}$ by Hypothesis \ref{np-h-norm}, 
$\aleph_7$ is zero on the kernel $N$.  
For $g\in \mathcal{N}$, there is an element $\vec{v}_g \in \bR^{i_0}$ such that 
$\exp(\vec{v}_g) = g$. We define $C_1, c_2: \bR^{i_0}\ra \bR$ by 
setting $C_1(\vec{v}_g) = C_1(g), c_2(\vec{v}_g) = c_2(g)$ for each $g\in \mathcal{N}$. 
Hence, $g \in \mathcal{N}$ is of form 
\begin{equation} \label{np-eqn-formgiiZii}
\newcommand*{\tempV}{\multicolumn{1}{r|}{}}
g=
\left( \begin{array}{ccccccc} 
\Idd_{n-i_0-1} & \tempV & 0 & \tempV & 0 & \tempV & 0 \\ 
\cline{1-7}
0 &\tempV & 1 &\tempV & 0 &\tempV & 0 \\ 
\cline{1-7}
C_1(\vec{v}_g) &\tempV &\vec{v}^T_g &\tempV &  O_5(g) &\tempV & 0 \\ 
\cline{1-7}
c_2(\vec{v}_g) &\tempV & \frac{\llrrV{\vec{v}_g}^2}{2}  &\tempV &  \vec{v}_g O_5(g) &\tempV & 1
\end{array} 
\right)
\end{equation}
since $\aleph_7(g)=0$, and $S(g) =  \lambda_g \Idd_{n-i_0-1}$ 
and the $(n-i_0, n-i_0)$-term must be $\lambda_g$ for some $\lambda_g> 0$
such that it goes to $\Idd$ in $K$.

(B)	
Lemma \ref{np-lem-conedecomp1} shows that $C_1(\vec{v}) =0$ for all $\vec{v} \in \bR^{i_0}$ for a coordinate system where $k$ has 
the form \[\llrrparen{0,\dots, 0, 1} \in \SI^{n-i_0-1}.\]

Let $\lambda_{S}$ denote the minimal norm of the eigenvalues of 
the upper-left part $S_g$ of $g$. 
We define 
\[\bGamma_{\tilde E, S}:= \left\{g \in \bGamma_{\tilde E}\big| \lambda_{S}(g) > a_1(g)\right\}.\]
There is always an element like this because $\bGamma_{\tilde E}$ 
acts on a subspace of dimension $n-i_0-1$ containing a compact set projectively 
diffeomorphic to $K$. 
In particular, we take the inverse image of 
suitable diagonalizable elements of the center  $H \cap N_K$  by the premise. 
We take the diagonalizable element in $N_K$ with $K''$ having a largest norm
eigenvalue. Let $g$ be such an element. 
Since $O_5(g)$ is orthogonal, the transversal weak middle-eigenvalue condition tells us 
\[\max\{\lambda_S(g), a_5(g), a_1(g)\} \geq a_9(g) =  \lambda_{\mbv_{\tilde E}}(g).\]
We have either $a_9(g) \geq a_5(g) \geq a_1(g)$ or 
$a_1(g) \geq a_5(g) \geq a_9(g)$ by \eqref{np-eqn-a9a5a1}  
depending on $\mu_g \geq 1$ or $\leq 1$. 
The second case can be ignored
since $a_9(g)$ is the smallest eigenvalue in that case and we can consider $g^{-1}$ 
to obtain that $a_1(g)= a_5(g) = a_9(g)$  by Proposition \ref{np-prop-eigSI} 
reducing to the first case. 
Hence,for $g \in \bGamma_{\tilde E, S}$, we have $\mu_g \geq 1$. 
Hence, we have
\begin{equation} \label{np-eqn-a1a5a9} 
a_{1}(g) \leq a_5(g) \leq a_{9}(g) \leq \lambda_{S}(g) \hbox{ and } 
\mu_{g}  \geq 1 \hbox{ for } g \in \bGamma_{\tilde E, S}.
\end{equation} 

(C) 
Applying Lemma \ref{np-lem-conedecomp1},
we modify the coordinates such that 
the elements of $\mathcal{N}$ 
are of form: 
\begin{equation}\label{np-eqn-formn-5}
\newcommand*{\tempV}{\multicolumn{1}{r|}{}}
\eta= 
\left( \begin{array}{ccccccc} 
\Idd_{n-i_0-1} & \tempV & 0 & \tempV & 0 & \tempV & 0 \\ 
\cline{1-7}
0 &\tempV & 1 &\tempV & 0 &\tempV & 0 \\ 
\cline{1-7}
C_1(\vec{v}_\eta) &\tempV & \vec{v}_\eta^T &\tempV & \Idd_{i_0} &\tempV & 0 \\ 
\cline{1-7}
c_2(\vec{v}_\eta) &\tempV & \frac{\llrrV{\vec{v}_\eta}^2}{2} &\tempV & \vec{v}_\eta &\tempV & 1 
\end{array} 
\right) \hbox{ where } C_1(\vec{v}_\eta) = 0.
\end{equation} 
By the group property, $\vec{v} \mapsto c_2(\vec{v})$ is a linear map. 

We have coordinates such that $K''\subset \SI^{n-i_0-2}$. 
There exists a sequence of elements $z_{i}$ of 
$N_K\cap H$ in the virtual center  $H$
such that a largest norm eigenvalue has a direction 
in $K''$ and $z_i|K'' \ra \Idd_{K''}$. 

Since $\clo(U)$ is in properly convex $\clo(\torb)$, 
it is in an affine patch where 
$\mbv_{\tilde E}$ is the origin. That is, 
$\mbv_{\tilde E} = \llrrparen{0, 0, 0, 1}$. 
Let $g_i \in \bGamma_{\tilde E+}$ be the element going to $z_i$ under $\Pi^\ast_K$. 
Then $\{g_{i}(x)\}$ for a point $x$ of $U$ 
converges to $\llrrparen{\lambda\vec{a}, 0, \vec{w}, 1}$ 
for $\llrrparen{a} \in K''$ and some 
$\vec{w} \in \bR^{i_0}$, $\lambda \in \bR$. 
Here, $\lambda> 0$ since our point must project to the limit of $z_i(\Pi_K(x))$ as $i\ra \infty$.
Hence by \eqref{np-eqn-a1a5a9} and 
\begin{equation}\label{np-eqn-conva} 
\llrrparen{\lambda\vec{a}, 0, \vec{w}, 1} \in \clo(U)
\end{equation} 
Since $z_i|K'' \ra \Idd_{K''}$, we may 
assume that an open subset of $K''$ can be realized 
as $\llrrparen{\vec{a}}$ in the above. 

By \eqref{np-eqn-formn-5}
\begin{equation} \label{np-eqn-Nc2zero} 
\CN(\vec{v}_1)^k \llrrparen{\lambda\vec{a}, 0, \vec{w}, 1}  = 
\llrrparen{\lambda\vec{a}, 0, \vec{w}, k \lambda \vec{a} \cdot 
	c_2(\vec{v}_1) + k \vec{v}_1 \cdot \vec{w} + 1}
\end{equation} 
as $C_1(\vec{v}_1) =0$ for every $\vec{v}_1 \in \bR^{i_0}$. 
Suppose that $\lambda \vec{a} \cdot c_2(\vec{v}_1)+ \vec{v}_1 \cdot \vec{w} \ne 0$. 
Then
as $k \ra \infty$, $\{\CN(k\vec{v}_1)(\llrrparen{\vec{a}})\}$ converges to a point 
and as $k \ra -\infty$, it converges to its antipode. 
The limits form an antipodal pair of points in $\clo(U)$. 
This contradicts the proper convexity of $\torb$. 

Hence, $\lambda \vec{a} \cdot c_2(\vec{v}_1) + \vec{w}\cdot\vec{v}_1 = 0$ holds 
for every $\vec{v}_1 \in \bR^{i_0}$.
We obtain that $c_2(\vec{v}_1)$ is a linear function of $\vec{v}_1$ to be 
written as $\hat c_2 \cdot \vec{v}_1$ where  $\hat c_2$ is an $(n-i_0)\times i_0$-matrix. 
Then $\vec{w}^T := -\lambda\hat c_2^T \vec{a}^T. $
Let $\hat K$ denote the image of 
\[\llrrparen{\vec{a}} \mapsto \llrrparen{\lambda\vec{a}, 0, \vec{w}, 0}, 
\vec{w}^T := -\lambda\hat c_2^T \vec{a}^T, \llrrparen{\vec{a}} \in K''.\]
Under $\Pi_K$, the compact convex set $\hat K$ embeds onto a compact convex set $K''$ of the same dimension. 

Since every point of $K^{''o}$ 
is a limit point of the orbit of a point of $K^o$ under $z_i$, 
every point of $\hat K^o$ is a limit point of a point of $U$ under $g_i$. 
Hence, we obtain $\hat K\subset \Bd \torb$ by convexity
by \eqref{np-eqn-conva}.  

Also, $\CN$ acts on $\hat K$ by our discussion and hence on $\hat K$. 
We choose the coordinates such that $\hat K$ corresponds to 
$x_1 = x_2 = \cdots = x_{n-i_0-1} = 0$. 
Under this coordinate system, 
\begin{equation} \label{np-eqn-C1c2} 
C_1(\vec{v})=0, c_2(\vec{v})=0
\hbox{ for every } \vec{v} \in \bR^{i_0}. 
\end{equation} 




(D)
Consider a sequence $\{g_i\}$ of elements $g_i \in \bGamma_{\tilde E, S}$
with $\{\Pi_K^\ast(g_i)(y)\}$ converging to $x\in K''$. 
We claim that every limit point $x'$ of $g_i(u)$ for $u\in U$ 
is in $\hat K$: In our coordinates as above \eqref{np-eqn-C1c2}, 
we have $x' =\llrrparen{\lambda \vec{a}, 0, \vec{w}, 1}$. 
\eqref{np-eqn-Nc2zero} still holds, and since $c_2(\vec{v}_1)= 0$, 
$\vec{v}_1\cdot \vec{w} = 0$ for every $\vec{v}_1\in \bR^{i_0}$. 
Thus, $\vec{w}=0$, and $x'\in \hat K$. 

Since the set of such sequences is invariant under the conjugation 
by $\bGamma_{\tilde E}$, it follows that the set of accumulation points of 
such sequence of elements
in $\hat K$ is  $\bGamma_{\tilde E}$-invariant. 
Since each point of $K''$ can be an accumulation point of some sequence of elements 
in $N_K$, it follows that 
$\hat K$ is $\bGamma_{\tilde E}$-invariant. 
This implies that $\hat K$ is  $\langle \bGamma_{\tilde E}, \CN \rangle$-invariant.

We may assume in our chosen coordinates that 
\begin{equation}\label{np-eqn-CCoor} 
C_1(g)=c_2(g) = C_1(\vec{v}) = c_2(\vec{v}) = 0
\hbox{ for every } g \in \bGamma_{\tilde E}, \vec{v} \in \bR^{i_0}.
\end{equation}

\end{proof}


%
%





\subsection{Strictly joined and quasi-joined ends for $\mu\equiv 1$} \label{np-subsub-qjoin}

We now discuss joins and their generalizations in-depth in this subsection. 
That is, we only consider when $\mu_{g}= 1$ for all $g\in \bGamma_{\tilde E}$.
We use a hypothesis and later show that the hypothesis is true in our cases
to prove the main results. 
Again, we assume the hypothesis virtually since it is sufficient. 


\begin{hyp}[$\mu \equiv1$]\label{np-h-qjoin} 
Let $\bGamma_{\tilde E}$ be a p-end holonomy group. 
We continue to assume Hypothesis \ref{np-h-norm} for $\bGamma_{\tilde E}$.  
\begin{itemize} 
\item Every $g \in \bGamma_{\tilde E} \ra M_g$ is 
such that $M_g$ is in a fixed orthogonal group $\Ort(i_0)$. Thus, $\mu_g = 1$ identically. 
\item $\bGamma_{\tilde E}$ acts on a subspace $\SI^{i_0}_\infty$ containing 
$\mbv_{\tilde E}$ and the properly convex 
domain $K'''$ in the subspace $\SI^{n-i_0-2}$ 
forming an independent pair with $\SI^{i_0}_\infty$ 
mapping homeomorphic to the factor $K''$ of 
$K = \{k\} \ast K'' $ under $\Pi_{K}$.
\item $\CN$ acts on these two subspaces, fixing every point 
of $\SI^{n-i_0-2}$.
\end{itemize} 
\end{hyp} 




By the convexity of $\tilde \Sigma_{\tilde E}$, 
we can choose $H$ such that $U \subset H^o$,
$K''' \subset H$ and $\SI^{i_0}_\infty \subset \Bd H$. 
We may assume that
$H$ is the closed $n$-hemisphere defined by $x_{n-i_0} \geq 0$. 
We identify $H^o$ with an affine space $\mathds{A}^n$. 
(See Section \ref{prelim-sec-rps}.)

By Hypothesis \ref{np-h-qjoin}, 
the elements of $\CN$ have the form of \eqref{np-eqn-nilmatstd} with 
\[C_1(\vec{v})=0, c_2(\vec{v})=0 \hbox{ for all } \vec{v} \in \bR^{i_0},\] 
and the elements of $\bGamma_{\tilde E}$ have the form of \eqref{np-eqn-formgii} with 
\[s_1(g) =0, s_2(g) = 0, C_1(g) = 0, \hbox{ and } c_2(g) = 0.\] 





We recall the invariants from the form of \eqref{np-eqn-kernel}. 
\[\aleph_7(g):= \frac{a_7(g)}{\lambda_{\mbv_{\tilde E}}(g)} - \frac{\llrrV{\vec{v}_g}^2}{2} \]
for every $g \in \bGamma_{\tilde E}$. Recall 
$\aleph_7(gh) = \aleph_7(g) + \aleph_7(h), \hbox{ whenever }g, h \in \bGamma_{\tilde E}.$



Under Hypothesis \ref{np-h-qjoin}, Lemma \ref{np-lem-matrix} shows that
every $g \in \bGamma_{\tilde E}$ is of the following form: 
\renewcommand{\arraystretch}{1.5}
\begin{equation}\label{np-eqn-gformII}
\newcommand*{\tempV}{\multicolumn{1}{r|}{}}
\left( \begin{array}{ccccccc} 
S_{g} & \tempV & 0 & \tempV & 0 & \tempV & 0 \\ 
\cline{1-7}
0 &\tempV & \lambda_g &\tempV & 0 &\tempV & 0 \\ 
\cline{1-7}
0 &\tempV & \lambda_g \vec{v}^T_g &\tempV & \lambda_g O_5(g) &\tempV & 0 \\ 
\cline{1-7}
0 &\tempV & \lambda_g\left(\aleph_7(g) + \frac{\llrrV{\vec{v}_g}^2}{2}\right) &\tempV &
\lambda_g \vec{v}_g O_5(g) &\tempV & \lambda_g 
\end{array} 
\right),
\end{equation}
and every element of $\mathcal{N}$ is of form 
\renewcommand{\arraystretch}{1.5}
\begin{equation}\label{np-eqn-nformII}
\newcommand*{\tempV}{\multicolumn{1}{r|}{}}
\mathcal{N}(\vec{v}) = 
\left( \begin{array}{ccccccc} 
	\Idd_{n-i_0-1} & \tempV & 0 & \tempV & 0 & \tempV & 0 \\ 
	\cline{1-7}
	0 &\tempV & 1 &\tempV & 0 &\tempV & 0 \\ 
	\cline{1-7}
	0 &\tempV & \vec{v}^T &\tempV & \Idd_{i_0} &\tempV & 0 \\ 
	\cline{1-7}
	0 &\tempV & \frac{\llrrV{\vec{v}}^2}{2} &\tempV & \vec{v} &\tempV & 1 
\end{array} 
\right).
\end{equation}

%

We assumed $\mu\equiv 1$. 
We define 
\begin{equation}\label{np-eqn-lambdak} 
\lambda_{k}(g):=\lambda_{\mbv_{\tilde E}}(g) \hbox{ for } g\in \bGamma_{\tilde E}.
\end{equation}
We define $\lambda_{K''}(g)$ as the maximal norm of the eigenvalue occurring for $S(g)$. 
We define $\bGamma_{\tilde E, +}$ as a subset of $\bGamma_{\tilde E}$ consisting of elements $g$ such that the following hold:  
\begin{itemize} 
\item the largest norm $\lambda^{Tr}_{\max}(g)$ of the eigenvalues occurs at the vertex $k$, i.e., $\lambda_{1}(g) = \lambda_{k}(g)$,  and 
\item all other norms of the eigenvalues occurring at $K'''$ are strictly less than $\lambda_{\mbv_{\tilde E}}(g)$.  
\end{itemize} 
Then, since $\mu_g=1$, we necessarily have
$\lambda_k(g)= a_1(g) = a_5(g)=\lambda_{\mbv_{\tilde E}}(g)$ and
hence $\lambda^{Tr}_{\max}(g) = \lambda_{\mbv_{\tilde E}}(g)$
for $g \in \bGamma_{\tilde E, +}$. 

The second largest norm $\lambda_2(g)$ is equal to $\lambda_{K''}(g)$. 
Thus, $\bGamma_{\tilde E, +}$ is a semigroup.     
The condition that $\aleph_7(g) \geq 0$ for $g \in \bGamma_{\tilde E, +}$ is said to be the 
{\em nonnegative translation condition}. 
\index{nonnegative translation condition } 

Again, we define 
\begin{equation}\label{np-eqn-mu7} 
\mu_7(g) : = \frac{\aleph_7(g)}{\log\frac{\lambda_{\mbv_{\tilde E}}(g)}{\lambda_2(g)}}
\hbox{ for } g \in \bGamma_{\tilde E, +}.
\end{equation} 
The condition 
\begin{equation}\label{np-eqn-uptc}
\mu_7(g) > C_0, g \in  \bGamma_{\tilde E, +} \hbox{ for a uniform constant } C_0, C_0> 0 
\end{equation}
is called the {\em uniform positive translation condition}. 
(Heuristically, the condition means that we do not translate in the negative direction by too much
for bounded $\frac{\lambda_{\mbv_{\tilde E}}(g)}{\lambda_2(g)}$.)
\index{uniform positive translation condition}

\begin{lemma} \label{np-lem-i}	
The condition $\aleph_7(g) \geq 0$ for $g \in \bGamma_{\tilde E,+}$ 
is a necessary condition 
such that $\bGamma_{\tilde E}$ acts on a properly convex domain.
\end{lemma} 
\begin{proof}
Suppose that $\aleph_7(g) < 0$ for some $g \in \bGamma_{\tilde E, +}$.
Let $k'\in K^o$. Then $g^n(k') \ra k$ as $n \ra \infty$ in $K$ by the property of 
being in $\bGamma_{\tilde E, +}$. 
Now, we use \eqref{np-eqn-HlU} and see that $\{g^n(U_{k'})\}$ 
converges geometrically to an $(i_0+1)$-dimensional hemisphere
since $\{\aleph_7(g^n)\} \ra -\infty$ as $n \ra \infty$ implies that
$g$ sends the affine subspace $H^o_{k'}$  to $H^o_{g^n(k')}$, 
and points in it towards $\llrrparen{-1,0,\dots, 0}$ in the above coordinate system by \eqref{np-eqn-gformII}.
Thus, $\bGamma_{\tilde E}$ cannot act on a properly convex domain. 
(See also Remark \ref{np-rem-alpha7}.)
\end{proof}

From the matrix equation \eqref{np-eqn-gformII},
we define $\vec{v}_{g}$ for every $g \in \langle \bGamma_{\tilde E}, 
\CN \rangle$. 
(We just need to do this under a single coordinate system. )

\begin{lemma} \label{np-lem-vbound}
	Given $\bGamma_{\tilde E}$ satisfying Hypotheses \ref{np-h-norm} and \ref{np-h-qjoin}, 
	let $\gamma_{m}$ be any sequence of elements of $\bGamma_{\tilde E, +}$
	such that $\{\lambda_{k}(\gamma_{m})/\lambda_{K''}(\gamma_{m})\} \ra \infty$. 
	Then we can replace it with another sequence $\{g_{m}^{-1}\gamma_m\}$ for$g_m \in \bGamma_{\tilde E}$ such that 
	\[\llrrV{\vec{v}_{g_{m}^{-1}\gamma_m}} \hbox{ and }
\Pi^{\ast}_{K}(g_{m}) \in \Aut(K)  \] 
	are uniformly bounded, and
                          \[ \left\{\frac{\lambda_k(g_m^{-1}\gamma_m)}{\lambda_{K''}(g_m^{-1}\gamma_m)}\right\} \ra \infty. \]
\end{lemma} 
\begin{proof}
    Denote $\vec{v}_m:= \vec{v}_{\gamma_m}$. 	
	Suppose that $N_{K}$ is discrete. 
Then, since its action on the interior of $K$ is properly discontinuous, 
we have an orbifold bundle $\Bd U/\bGamma_{\tilde E} \ra K^o/N_K$. 
	This means that the subgroup $\bGamma_{\tilde E, l}$of $\bGamma_{\tilde E}$
	acts cocompactly on a complete affine leaf $l$ giving us the fibers. 
	Since the stabilizer of $N_K$ on each point of $K^o$ is finite, 
	$\bGamma_{\tilde E}\cap \ker \Pi_K^\ast$ acts on $l$ cocompactly. 
	The action of $ \mathcal{N}\hat O(i_0)$ is proper on each leaf. 
	Hence,$\bGamma_{\tilde E}\cap \mathcal{N}\hat O(i_0)$ is a lattice in $\mathcal{N} \hat O(i_0)$. 
	By the cocompact action of $\bGamma_{\tilde E}\cap 
	\mathcal{N} \hat O(i_0)$
	in $\mathcal{N}\hat O(i_0)$, we can multiply $\gamma_{m}$ by $g_{m}^{-1}$ for 
	an element $g_{m}$ of $\bGamma_{\tilde E}\cap \mathcal{N}\hat O(i_0)$
	nearest to $\mathcal{N}(\vec{v}_{m})$.  The result follows since
	the action on $\SI^{n-1}_{\mbv_{\tilde E}}$ is given by only $\vec{v}_g$ and 
	$S_g$ for $g \in \bGamma_{\tilde E}$ as we can see from 
	\eqref{np-eqn-gformII}. 
	The last convergence follows since $S(g_m) = \Idd_{n-i_0}$ and 
	the matrix multiplication form of $g_m^{-1}\gamma_m$
	considering the top left $(n-i_0)\times (n-i_0)$-block 
	and the bottom right $(i_0+2)\times (i_0+2)$-block. 
	
	We now assume that $N_K$ is nondiscrete. 
	$\tilde \Sigma_{\tilde E}$ has a compact fundamental domain $F$.
	Thus, given $\vec{v}_m$,  for $x \in F$, 
	\[\CN(\vec{v}_m)(x) \in g_m(F) \hbox{ for some } g_m \in \bGamma_{\tilde E}.\] 
	Then $g_m^{-1}\CN(\vec{v}_m)(x) \in F$. Since 
	\[g_m(y) = \CN({\vec{v}_m})(x) \in g_m(F)  \hbox{ for } y \in F \hbox{ and } x \in F,\] 
	it follows that
	\begin{equation}\label{np-eqn-dK}
	d_K\bigg(\Pi_K(y), \Pi^{\ast}_K(g_m)(\Pi_{K}(y))(= \Pi_{K}(x))\bigg) < C_{F}
	\end{equation}
	for a constant $C_F$ depending on $F$. 
	\begin{enumerate}
		\item [(i)] $g_m$ is of the form of matrix \eqref{np-eqn-gformII}
		\item[(ii)] $S_{g_m}$ is in a bounded neighborhood of $\Idd$ by above \eqref{np-eqn-dK}
		since $\hat S_{g_m}$ moves a point of a compact set $F$ to 
		a uniformly bounded set. (This follows from considering the Hilbert metric.)
	
	\end{enumerate} 
	From the linear block form of $g_m^{-1}\CN(\vec{v}_m)$ and the fact that  $g_m^{-1}\CN(\vec{v}_m)(x) \in F$, 
	we obtain that the corresponding $\vec{v}_{g_m^{-1}\CN(\vec{v}_m)}$ can be made uniformly bounded independent of $\vec{v}_m$. 
	
	For the element $\gamma_m$ above, we take its vector $\vec{v}_{\gamma_m}$ and find our $g_{m}$ for $\mathcal{N}(\vec{v}_{\gamma_m})$.
	We obtain $\gamma'_{m}:= g_{m}^{-1}\gamma_m$. Then the corresponding 
	$\vec{v}_{g_{m}^{-1}\gamma_m}$ is uniformly bounded as
	we can see from the block multiplications and the action on 
	$\tilde \Sigma_{\tilde E}$ in $\SI^{n-1}_{\tilde E}$. 
The final part follows from (ii) and the fact that
$\{\lambda_{k}(\gamma_{m})/\lambda_{K''}(\gamma_{m})\} \ra \infty$
and the matrix multiplication form of $g_m^{-1} \gamma_m$
if we consider the top left $(n-i_0)\times (n-i_0)$-submatrix 
and the bottom right $(i_0+2)\times (i_0+2)$-submatrix.
\end{proof} 


\begin{lemma}\label{np-lem-noball} 
	Suppose that the holonomy group of $\orb$ is strongly irreducible. 
	Given $\bGamma_{\tilde E}$ satisfying Hypotheses \ref{np-h-norm} and \ref{np-h-qjoin}. 
	Let $U$ be the properly convex p-end neighborhood of $\tilde E$.
	Let $H_k$ be the $(i_0+1)$-hemisphere mapping to 
	the vertex $k$ of Hypothesis \ref{np-h-qjoin} under $\Pi_K$. 
	Then the interior of 
$\clo(U) \cap H_{k}$ is not a domain $B$ with a nonempty interior in $H^k$ with $\Bd B \ni \mbv_{\tilde E}$.
\end{lemma} 
\begin{proof} 
	Since $\mathcal{N}$ acts on $H_k$, it also acts on $B$. 
	The matrix form of $\mathcal{N}$ is given by the coordinates 
	where $H_k$ is the 
	projectivization of the span of $e_{n-k}, \dots, e_{n+1}$
	as we can see from \eqref{np-eqn-form1} in the proof of 
	Lemma \ref{np-lem-conedecomp1}.
	Hence, $B$ is an ellipsoid as we can see from the form of $\mathcal{N}$
	in \eqref{np-eqn-nformII}. 
	Firstly, 
	\begin{equation}\label{np-eqn-alpha7} 
	\aleph_7(h) =0 \hbox{ for all } h \in \bGamma_{\tilde E} 
	\end{equation}
	by Lemma \ref{np-lem-i} since otherwise by \eqref{np-eqn-gformII}
	\[\{\gamma^i(B)\} \ra H_{k} \hbox{ as } i \ra \infty \hbox{ or } -\infty \hbox{ for } 
	\gamma \hbox{ with } \aleph_7(\gamma) \ne 0.\]
	Since $\tilde \Sigma_{\tilde E}/\bGamma_{\tilde E}$ is compact,
we have the sequence $\gamma_i \in \bGamma_{\tilde E, +}$ where 
	\[\left\{\frac{\lambda_{\mbv_{\tilde E}}(\gamma_i)}{\lambda_2(\gamma_i) } \right\} \ra \infty, \aleph_7(\gamma_{i})=0,
	\hbox{ and $\{\gamma_i| K''\}$ are uniformly bounded.} \]
	Now modify $\gamma_{i}$ to $g_i^{-1} \gamma_i$ 
	by Lemma \ref{np-lem-vbound} for $g_i$ obtained there. 
	Hence, rewriting $\gamma_i$ as the modified one,
	we have 
	\[\left\{\frac{\lambda_{\mbv_{\tilde E}}(\gamma_i)}{\lambda_2(\gamma_i) }\right\} \ra \infty, \aleph_7(\gamma_{i})=0, \bv_{\gamma_i} 
	\hbox{ and $\gamma_i| K''$ is uniformly bounded.} \]
	
	Recall that $K$ is a strict join $K'' \ast \{k\}$ for a properly convex domain $K''\subset \Bd \torb$ of dimension 
	$n-i_{0}-2$ and a vertex $k$.
	Denote by $\SI(K'')$ and $\SI(H)$ the subspaces spanned by $K'$ and $H_{k}$
	forming a pair of complementary subspaces in $\SI^n$. 
	
	From the form of the lower-right $(i_0+2)\times (i_0+2)$-matrix of the above matrix, 
	$h_{i}$ must act on the horosphere $H \subset \SI(H_{k})$. 
The group $\CN$ also acts transitively on $H_{k}$. 
	Hence, for any such matrix, we can find an element of $\CN$ such that 
	the product is in the orthogonal group acting on $H_{k}$. 
	
	Now, this is the final part of the proof: 
	Let $H_{\mx}$ denote $\SI(H_{k}) \cap \clo(\torb)$ and $K^{\prime \prime}_{\mx}$ the set $\SI(K'') \cap \clo(\torb)$.
	Since $\{\vec{v}_{\gamma_m}\}$ is bounded and 
	$\aleph_7(\gamma_m)= 0$, 
	we have the sequence $\{\gamma_m\}$
	\begin{itemize}
		\item acting on $K''_{\mx}$ is uniformly bounded and 
		\item $\gamma_m$ acting on $H_{\mx}$ in a uniformly bounded manner 
		as $m \ra \infty$. 
		
	\end{itemize}
	By Proposition \ref{prelim-prop-decjoin} for the $l = 2$ case,
    $\clo(\torb)$ equals the join of $H_{\mx}$ and $K''_{\mx}$. 
	This implies that $\bGamma$ is virtually reducible by Proposition \ref{prelim-prop-joinred}, which contradicts  the premise of the lemma.
	%
\end{proof}





For this proposition, we do not assume that $N_K$ is discrete. 
Also, we do not need the assumption of
the proper convexity of $\orb$. 

\begin{proposition}[Quasi-joins] \label{np-prop-qjoin}
Let $\Sigma_{\tilde E}$ be the end orbifold of an NPNC R-end $\tilde E$ of a properly convex 
real projective $n$-orbifold $\orb$. 
Let $\bGamma_{\tilde E}$ be the p-end holonomy group. 
Let $\tilde E$ be an NPNC R-p-end
and $\bGamma_{\tilde E}$ and $\mathcal N$ act on a p-end neighborhood 
$U$ fixing $\mbv_{\tilde E}$. Let $K, K'', K''', \SI^{i_0}_\infty,$ and $\SI^{i_0+1}$ be as in Hypotheses \eqref{np-h-norm} and \eqref{np-h-qjoin}. 
Assume that 
\begin{itemize} 
\item $\bGamma_{\tilde E}$ satisfies the transverse weak middle-eigenvalue condition
with respect to $\mbv_{\tilde E}$. 
\item $\bGamma_{\tilde E}$ acts on $K'''$ and $k$. 
\item $\mu_{g} = 1$ for all $g \in \bGamma_{\tilde E}$. 
\item Elements of $\bGamma_{\tilde E}$ and $\CN$ respectively are of the form \eqref{np-eqn-form1} and \eqref{np-eqn-form2}
with \[C_1(\vec{v}) = 0, c_2(\vec{v})=0, C_1(g) = 0, c_2(g) =0\]
for every $\vec{v} \in \bR^{i_0}$ and $g \in \bGamma_{\tilde E}$.
\item $\bGamma_{\tilde E}$ normalizes $\CN$, and $\CN$ acts on $U$ 
and each leaf of $\widetilde{\mathcal{F}}_{\tilde E}$ of $\tilde \Sigma_{\tilde E}$. 
\end{itemize} 
Then the following hold\/{\rm :} 
\begin{enumerate} 
\item[(i)] 
The uniform positive translation condition is equivalent to the existence of
a properly convex p-end neighborhood $U'$ whose closure meets $\SI^{i_0+1}_k$ at $\mbv_{\tilde E}$ only.
\item[(ii)] In these cases, $U'$ is radially foliated by 
line segments from $\mbv_{\tilde E}$, 
\item[(iii)] $\clo(U') \cap \Bd \mathds{A}^n = K''' \ast \{\mbv_{\tilde E}\}$.
\item[(iv)] $\aleph_7$ is identically zero if and only if $\CH(U)$ is the interior of the join $K''' \ast B$ for an open ball $B$ in $\SI^{i_0+1}_k$, 
and $\CH(U)$ is properly convex. 

\end{enumerate}
\end{proposition}
\begin{proof} 
(A) First, we find some coordinate system and describe
the action of $\bGamma_{\tilde E}$ on it. 

Let $H$ be the unique $n$-hemisphere containing segments in 
directions of $\tilde \Sigma_{\tilde E}$ from $\mbv_{\tilde E}$ 
 where $\partial H$ contains $\SI^{i_{0}}_{\infty}$ and $K'''$ 
 in a general position by Hypothesis
\ref{np-h-qjoin}. 
Then $H^{o}$ is an affine subspace that is also denoted by
$\mathds{A}^n$ containing $U$. 
Since $\bGamma_{\tilde E }$ and $\mathcal{N}$ act on $K'''$ and $\SI^{i_{0}}_{\infty}$ spanning $\partial H$, we deduce that
$\bGamma_{\tilde E }$ and $\mathcal{N}$ 
act as an affine transformation group on $\mathds{A}^n$. 

Let $H_{l}$ denote the hemisphere with boundary $\SI^{i_{0}}_{\infty}$ and corresponding to a leaf $l$
of the foliation on $\tilde \Sigma_{\tilde E}$.  
Recall that $K^o$ is the leaf space. 
Let $\mathds{A}^n$ be the affine space whose boundary
contains $\SI^{i_0}_{\infty}$ and a leaf $l$ is the $(i_0+1)$-dimensional 
affine subspace with a transverse affine space of dimension
$n-(i_0+1)$ meeting it at a point, considered the origin. 
We chose $\mathds{A}^n$ such that 
 $\mbv_{\tilde E}$ in $\Bd \mathds{A}^n$.
We may further require that from $\mbv_{\tilde E}$
the space of directions to $\mathds{A}^n$ contains
those to $U$. 
Furthermore, by projection to the affine space of
dimension $n-i_0-1$ with kernel vector space parallel to $l$, 
we obtain an affine space $\mathds{A}^{n-i_0-1}$ with the origin
the image of $\mbv_{\tilde E}$. 
The image
of $K^o$ is a cone from the origin in $\mathds{A}^{n-i_0-1}$.
The space of the directions of the cone is isomorphic to $K^{''o}$.
We may regard $K'''$ as a subset of the ideal boundary of
$\Bd \mathds{A}^{n-i_0-1}$ and in $\Bd \mathds{A}^n$. 
Hence, we projectively identify 
\[\bigcup_{l \, \in \, K^o} H_{l}^{o} = K^o \times \bR^{i_{0}+1} \subset \mathds{A}^{n}.\] 
(It might be helpful to see Figure \ref{np-fig:quasi-j} where
we convert the coordinates such that $(0,\vec{0},0)$, $\SI^1_\infty$, and $\hat K$ are now in
$\Bd \mathds{A}^n$.) 

There exists a family of affine subspaces in $\mathds{A}^n$ parallel to
$\mathds{A}^{n-i_0-1}$.
In addition, there is a transverse family of affine spaces forming a foliation
$\mathcal{V}^{i_0+1}$ where each leaf is
the complete affine space parallel to $\mathds{A}^{i_0+1}$. 
We denote it by $\{(\vec{x}_1, 1)\}\times \bR^{i_0+1}$ for
$\vec{x}_1 \in \bR^{n-i_0-1}$. 
The affine coordinates are given by $(\vec{x}_1, 1, \vec{x}_2, x_{n+1})$ 
where $1$ is at the $(n-i_0)$-th position, $\vec{x}_1$ is 
an $(n-i_0-1)$-vector and $\vec{x}_2$ is an $i_0$-vector. 


Now we describe $H^o_l \cap U$ for each $l \in K^o$. 
We use the affine coordinate system on $\mathds{A}^n$ such that 
$H_l^o$ are parallel affine $(i_0+1)$-dimensional spaces 
with origins in $\mathds{A}^{n-i_0-1}$. We use the parallel affine 
coordinates. 
According to the matrix form \eqref{np-eqn-nformII}, $\CN$ acts 
on each $x \times \bR^{i_0+1}$, $x \in \hat K^o$. 

We denote each point in $H_l^o$ by $(\tilde x, x_{n+1})$ where 
$\tilde x$ is a point of $\mathds{A}^{i_0}$.  
Each of $E_l := H_l \cap U$ is given by 
\begin{equation} \label{np-eqn-paraboloid} 
x_{n+1} > \llrrV{\vec{x}_2}^2/2 + C_l, C_l \in \bR 
\end{equation}
since $\CN$ acts on each
where $C_l$ is a constant depending on $l$ and $U$
by \eqref{np-eqn-kernel}.
(The vertex $\mbv_{\tilde E}$ corresponds to the ideal 
point in the positive infinity in terms of the $x_{n+1}$-coordinate.
See Section \ref{np-subsub-quadric}.)

There is a family of quadrics of the form 
$Q_{l, C}$ defined by $x_{n+1} =  \llrrV{\vec{x}_2}^2/2 + C$ for each $C\in \bR$
on each leaf $l$ of $\mathds{V}^{i_0+1}$ using the affine coordinate system. 
The family forms a foliation $\mathcal{Q}_l$ for each $l$.

Now we describe the $\bGamma_{\tilde E}$-action. 
Since $\mu_g = 1$ for all $g \in \bGamma_{\tilde E}$ by our hypothesis \ref{np-h-qjoin},  
it follows that 
$\lambda_{\mbv_{\tilde E}} = \lambda_k$ by definition \eqref{np-eqn-lambdak}. 
Given a point $x = \llrrparen{\vec v} \in U'\subset \SI^n$ where 
$\vec v = \vec v_s + \vec v_h$ 
where $\vec v_s$ is in the direction of $\SI(K''')$ and
$\vec v_h$ is in one of $\SI_k^{i_0+1}$. 
If $g \in \bGamma_{\tilde E, +}$, then we obtain
\begin{equation} \label{np-eqn-gplus}
g\llrrparen{\vec v} = \llrrparen{g \vec v_s + g \vec v_h} \hbox{ where } 
\llrrparen{g \vec v_s} \in K'''
\hbox{ and } \llrrparen{g \vec v_h} \in H_k.
\end{equation} 


Let $\Pi_{i_0}: U \ra \bR^{i_0+1}$ be the projection to the last 
$i_0+1$ coordinates $x_{n-i_0}, \dots, x_n$. 
We obtain a commutative diagram and an affine map $L_g$ induced from $g$
\begin{alignat}{3} \label{np-eqn-affine}
H_l^o  \, & \, \stackrel{g}{\longrightarrow} & \, \, g(H_l^o) \nonumber \\
\, \Pi_{i_0}\downarrow \, &  & \, \Pi_{i_0} \downarrow\, \, \nonumber\\ 
\bR^{i_0+1} \, & \, \stackrel{L_g}{\longrightarrow}  & \, \, \bR^{i_0+1}. 
\end{alignat}
By \eqref{np-eqn-kernel}, 
$L_g$ preserves the family of quadrics $\mathcal{Q}_l$ to $\mathcal{Q}_{g(l)}$  
since $\CN$ acts on the quadrics $U \cap H_l^o$ for each $l$
and $g$ normalizes $\CN$ by Hypothesis \ref{np-h-norm}. 
In addition, $L_g$ is an affine map, since $L_g$ is a projective map sending
a complete affine subspace $H_l^o$ to a complete affine subspace $g(H_l^o)$. 
Finally, by \eqref{np-eqn-kernel}, 
$g$ sends the family of quadrics shifted in the $x_{n+1}$-direction by $\aleph_7(g)$ 
from $l$ to $g(l)$ using the coordinates $(\vec{x}, x_{n+1})$ 
for $\vec{x} \in \bR^{i_0}$. 
That is, 
\begin{equation} \label{np-eqn-qshift}
g: Q_{l, C} \mapsto Q_{g(l), C+ \aleph_7(g)}.
\end{equation}




(B) Now, we give proofs. 
By assumption, $\bGamma_{\tilde E}$ acts on $K = K'' \ast \{k\}$. 
Choose an element $\eta \in \bGamma_{\tilde E, +}$ by Proposition
\ref{prelim-prop-Ben2} such that
$\lambda^{Tr}_{\max}(\eta) > \lambda_2(\eta)$
where $\lambda_{1}(\eta)$ corresponds to a vertex $k$ and $\lambda_{2}(\eta)$ is associated with $K'''$, 
and let $F$ be the fundamental domain in the open convex cone $K^o$ with respect to $\langle \eta \rangle$, which is a bounded domain in $\mathds{A}^{n-i_0-1}$. 
This corresponds to 
a radial subset from $\mbv_{\tilde E}$ 
bounded away at a distance from $K'''$ in $U$.

(i) 
This long proof is divided as follows. 
\begin{description} 
	\item[(i-a)] Forward part: We show that $U$ has a convex hull that is properly convex. 
	\begin{description} 
		\item[(i-a-1)]
		We show that the forward images 
		of the fundamental domain of $U$ under $\langle \eta \rangle$
		is contained in $K^o \times \sigma^{\prime, o}$ for some simplex $\sigma'$. 
		\item[(i-a-2)] Next, we try to show that the backward images the fundamental domain under $\langle \eta \rangle$ is eventually in any neighborhood of $K''' \ast \{\mbv_{\tilde E}\}$.  
		\item[(i-a-3)] Finally, we show that the convex hull of 
		$U$ is in a properly convex domain. 
	\end{description} 
	\item[(i-b)] Converse part: We prove the uniform positive translation condition under the assumption. 
	\end{description}

(i-a) (i-a-1) Let $\lambda_{K'''}(g)$ denote the maximal eigenvalue 
associated with $K'''$ for $g \in \bGamma_{\tilde E }$. 
Choose $x_{0} \in F$. 
Let $\bGamma_{\tilde E, F} := \{g \in \bGamma_{\tilde E}| g(x_0) \in F\}$. For $g \in \bGamma_{\tilde E, F}$, 
\begin{equation} \label{np-eqn-lambdaK''} 
 -C_F < \log \frac{\lambda_{\mbv_{\tilde E}}(g)}{\lambda_{K'''}(g)} < C_F
\end{equation} 
for a uniform $C_F > 0$ a number depending of $F$ only:
Otherwise, we can find 
\begin{itemize} 
\item a sequence $g_i$ with $g_i(x_0)\in F$ such that 
$\left\{\frac{\lambda_k(g_i)}{\lambda_{K'''}(g_i)}\right\} \ra 0$ or 
\item another sequence $g'_i$ with $g'_i(x_0)\in F$ such that 
$\left\{\frac{\lambda_k(g'_i)}{\lambda_{K'''}(g'_i)}\right\} \ra \infty$. 
\end{itemize} 
However, in the first case, let $\tilde x_0 \in U$ be a point mapping to 
$x_0$ under $\Pi_K$. Then 
$\{g_i(\tilde x_0)\}$ accumulates only to points of $K'''$ 
by Proposition \ref{prelim-prop-attract}, which is absurd. 
The second case is also absurd by taking $\{g_i^{-1}(x_0)\}$ instead.

Therefore, given $g \in \bGamma_{\tilde E, F}$, 
we can find a number $i_{0} \in \bZ$ dependent only on $F$ and $g$  
such that $\eta^{i_{0}} g \in \bGamma_{\tilde E, +}$
since $\log \frac{\lambda_k(\eta)}{\lambda_{K'''}(\eta)}$ 
is a constant greater than $1$.
Now, $\aleph_7(\eta^{i_{0}} g)$ is bounded below by some negative number
by the uniform positive eigenvalue condition \eqref{np-eqn-uptc}
and the fact that $\left|\log \frac{\lambda_{\mbv_{\tilde E}}(\eta^{i_{0}}g)}{\lambda_{K'''}(\eta^{i_{0}}g)}\right|$ is also uniformly
bounded.
Since $\aleph_7(\eta^{i_{0}} g) = i_{0}\aleph_7(\eta) + \aleph_7(g)$, 
we obtain 
\begin{equation}\label{np-eqn-a7}
\{\aleph_7(g)| g \in \bGamma_{\tilde E, F}\} > C
\end{equation} 
 for a constant $C$ by \eqref{np-eqn-uptc}. 
 $F$ is covered by $\bigcup_{g\in \bGamma_{\tilde E, F}} g(J)$ for a compact fundamental domain $J$ of $K^{o}$ by $N_{K}$.


We use a fixed affine coordinate system on 
$H^o_{k'}$ parallel under the translations 
preserving $\mathds{A}^{n-i_0-1}$. 
In the affine coordinates for $k'\in F$ of \eqref{np-eqn-affine}, 
the matrix form of \eqref{np-eqn-gformII} shows 
that $g \in \bGamma_{\tilde E}$ 
send paraboloids 
in affine subspaces $H^o_{k'}$ for $k' \in K^o$
to paraboloids in $H^o_{g(k')}$ for $g\in \bGamma_{\tilde E}$. 
(See Section \ref{np-subsub-quadric}.) 
Now, 
\begin{equation} \label{np-eqn-xnlower} 
x_{n}(H^{o}_{k'}\cap U) > C 
\end{equation}
for a uniform constant $C \in \bR$
by \eqref{np-eqn-a7} 
and the fact that $H^o_{k'}$ is the image of 
$H^o_{k''}$ for $k'' \in J$ by an element $g \in 
\bGamma_{\tilde E}$. 
(See Remark \ref{np-rem-alpha7}.)

Since by \eqref{np-eqn-xnlower}, 
\[\bigcup_{k'\in J}\bigcup_{g' \in \bGamma_{\tilde E, F}} g(H^o_{k'} \cap U)\] 
is a lower $x_{n}$-bounded set, its convex hull
$D_{F}$ as a lower-$x_{n}$-bounded subset 
of $K^o\times \mathds{A}^{i_{0}+1} \subset \mathds{A}^{n}$. 
Each region $D_F\cap H_{l'}$
is contained an $(i_0+1)$-dimensional 
simplex $\sigma_0$ with a face in the boundary of $H_{l'}$. 
Since there is a lower $x_{n}$-bound, 
we may use one $\sigma_0$ and translations to contain
every $U \cap H_{l'}$ in $D_F$.


Therefore, 
the convex hull $D_F$ in $\clo(\torb)$ is a properly convex set
contained in a properly convex set $F \times \sigma_0$. 



In $K'''$, the sequence of eigenvalue norms of $\eta^i$ converges to $0$ 
as $i \ra \infty$ and the eigenvalue $\lambda_{\mbv_{\tilde E}}$ at
$\SI_k^{i_0+1}$ goes to $+\infty$. 
Since 
\begin{equation}\label{np-eqn-alpha7gr}  
\aleph_7(\eta^i) = i \aleph_7(\eta) \ra +\infty \hbox{ as } i \ra \infty,
\end{equation} 
we obtain
\[\{\eta^i(D_F)\} \ra \{ \mbv_{\tilde E}\} \hbox{ for } i \ra \infty\]
geometrically, i.e., under the Hausdorff metric $\bdd_{H}$
by \eqref{np-eqn-gplus}. 
Again, for a sufficiently large integer $I$, 
because of the lower bound on the $x_n$-coordinates, we obtain the following: 
\begin{equation} \label{np-eqn-Ksigma}
\bigcup_{i\geq I}\eta^i(D_F) \subset K^o\times \sigma_0, 
\end{equation}
which is our first main result of this proof of the forward part of (i). 


(i-a-2) For each $k' =\llrrparen{\vec{x}_1, 1} \in K^o$, we
can find a point in  $\Pi_K^{-1}(k')$ of 
the form $\llrrparen{\vec{x}_1, 1, \vec{0}, C_{k'}}$ 
in $\Bd U \cap \mathds{A}^n$ for $\vec{0}$ a zero vector in $\bR^{i_0}$
and $C_{k'} \in \bR$. 
Using the $\CN$-action, we can parameterize $\Bd U \cap \mathds{A}^n$ 
starting from the point
$\llrrparen{\vec{x}_1, 1, \vec{0}, C_{k'}}$. 
The $\CN$-orbit of this point is given by
\begin{equation} \label{np-eqn-CNorbit} 
\llrrparen{\vec{x}_1, 1, \vec{v}, \frac{\llrrV{\vec{v}}^2}{2} + C_{k'}}, 
\vec{v} \in \bR^{i_0}.
\end{equation}  

Let 
\[p_i:= \llrrparen{\vec{x}_1(p_i), 1, \vec{v}(p_i), \frac{\llrrV{\vec{v}(p_i)}^2}{2} + C_{k'(p_i)}}.\]
We form a sequence $\{p_i\}$ of points on $\partial U$ for $k'(p_i) = 
\llrrparen{\vec{x}_1(p_i), 1}$. 
Consider $\eta^i$ in the form \eqref{np-eqn-gformII}. 
Since $\bigcup_{i\in \bZ}\eta^i(F)$ covers $K^o$, 
for each $p_i$ there is an integer $j_i$ for which $\eta^{j_i}(F)$ contains $\llrrparen{\vec{x}_1(p_i), 1}$. 
Suppose that $j_i \ra \infty$ as $i \ra \infty$. 
From considerations of $\eta^{j_i}$, 
we deduce that 
\begin{equation} \label{np-eqn-Cvecx1} 
\left\{C_{\llrrparen{\vec{x}_1(p_i), 1}} \right\} \ra +\infty \hbox{ as} 
\llrrV{\vec{x}_1(p_i)} \ra 0 
\end{equation} 
by \eqref{np-eqn-alpha7gr}. 

Similarly, suppose that $j_i \ra -\infty$ as $i \ra \infty$.  
We deduce that 
\begin{equation} \label{np-eqn-Cx1} 
\left\{\frac{C_{\llrrparen{\vec{x}_1(p_i), 1}}}{\llrrV{\vec{x}_1(p_i)}} \right\} \ra 0
	\hbox{ as } \llrrV{\vec{x}_1(p_i)} \ra \infty  
	\end{equation}
since $\aleph_7(\eta^{j_i}) = j_i \aleph_7(\eta) \ra -\infty$ linearly,
and every sequence of $\llrrV{\vec{x}_1(p_i)}$-coordinates of points of $\eta^{j_i}(F)$ grows uniformly exponentially as $j_i \ra -\infty$.



We now try to find all the limit points of $\{p_i\}$: 
\begin{itemize}  
\item Suppose first that  $\llrrV{\vec{x}_1(p_i)} \ra \infty$ and 
$\frac{\llrrV{\vec{x}_1(p_i)}}{\llrrV{\vec{v}(p_i)}^2} \ra \infty$. 
We find that 
\[\left\{\llrrparen{\vec{x}_1(p_i), 1, \vec{v}(p_i), \frac{\llrrV{\vec{v}(p_i)}^2}{2} + C_{k'(p_i)}}\right\},  \vec{v}(p_i) \in \bR^{i_0}\]
has only limit points of form $\llrrparen{\vec{u}, 0, 0, 0}$ 
for a unit vector $\vec{u}$ in the direction of $K''$ by \eqref{np-eqn-Cx1}. 
Hence, the limit is in $K'''$. 
\item Suppose that we have  
$\llrrV{\vec{x}_1(p_i)} \ra +\infty$ with $\llrrV{\vec{v}(p_i)}^2 \ra +\infty$ 
with their ratios bounded between two real numbers. 
Then a limit point is of the form $\llrrparen {\vec{u}, 0, 0, C}$ for some $C> 0$
and a vector $\vec{u}$ in the direction a point of $K'''$ by \eqref{np-eqn-Cx1}. 
In addition, every direction of $K''$ occurs as a direction of $\vec{u}$ for a limit point
by taking a sequence $\{p_i\}$ such that $\{\vec{x}_1(p_i)\}$
converges to $\vec{u}$ in directions. Hence, 
the limits are in $K'''\ast \{\mbv_{\tilde E}\}$. 
\item Suppose that
 $\llrrV{\vec{x}_1(p_i)} \ra +\infty$ with $\llrrV{\vec{v}(p_i)}^2 \ra +\infty$ and
 $\frac{\llrrV{\vec{x}_1(p_i)}}{\llrrV{\vec{v}(p_i)}^2} \ra 0$. 
 Then the only limit point is $\llrrparen{\vec{0}, 0, 0, 1 }$ since the last term dominates the others. 
\item Suppose that $1/C' \leq \llrrV{\vec{x}_1(p_i)} \leq C'$ for a constant $C'$. 
If $\llrrV{\vec{v}(p_i)}$ is uniformly bounded, then
\eqref{np-eqn-xnlower} shows that limit points in $\Bd U$. 
If $\llrrV{\vec{v}(p_i)}\ra \infty$, the only limit point is $\llrrparen{0, 0, 0,  1}$.
\item Suppose that $\llrrV{\vec{x}_1(p_i)} \ra 0$. Then the limit is
$\llrrparen{0, 0, 0,  1}$ by \eqref{np-eqn-Cvecx1}.
\end{itemize} 
These give all the limit points of $U$ in $\Bd \mathds{A}^n$, which we can easily deduce. 
Therefore, 
%
%
\begin{equation}\label{np-eqn-iv} 
\clo(U) \cap \Bd \mathds{A}^n \subset K''' \ast \{\mbv_{\tilde E}\}
\end{equation} 
by Theorem \ref{prelim-thm-Kobayashi2}. 


This also shows that $\eta^i(D_F)$ geometrically converges to
$K''' \ast \{\mbv_{\tilde E}\}$ as $i \ra -\infty$: 
First, the above three items show that 
for every $\eps> 0$, there exists $i_0$ such that 
 $\eta^i(D_F) \subset N_\eps(K''' \ast \{\mbv_{\tilde E}\})$
 for $i > i_0$. 
 Finally,  since $\Pi_K(\eta^i(D_F))$ geometrically converges to $K''$, 
 and we can find a sequence as in the first item that converges to any point of $K'''$, the geometric limit is $K''' \ast \{\mbv_{\tilde E}\}$.  
For every $\eps > 0$, there exists an integer $I$ such that 
\[\bigcup_{i < I} \eta^i(D_F) \subset N_\eps(K''' \ast \mbv_{\tilde E}) \cap \clo(\mathds{A}^n).\]

(i-a-3) By this and \eqref{np-eqn-Ksigma}, we obtain that 
except for finitely many $i$, 
\[\eta^i(D_F) \subset (N_\eps(K''' \ast \mbv_{\tilde E}) \cap \clo(\mathds{A}^n)) 
\cup K \times \sigma_0 \subset \clo(\mathds{A}^n).\]

We use the Fubini-Study metric $\bdd$ on $\SI^n$ where the subspaces spanned 
and $K$ and $\{k\}\times \bR^{i_0}$ are all orthogonal to each other. 

Assume $\eps < \pi/8$.
Let $p$ denote the vertex of the simplex $\{k\}\times \sigma_0$.
We may assume without loss of generality that $p = k\times O$. 
Then we may assume that 
$N_\eps(K''' \ast \mbv_{\tilde E}) \cap \clo(\mathds{A}^n)$ 
is in $\hat K \ast \{p\}$ where $\hat K$ is a  properly compact $(n-1)$-ball 
containing $K''' \ast \{\mbv_{\tilde E}\}$ and is contained
its $2\eps$-$\bdd$-neighborhood. 
Here, we may need to choose sufficiently small $\eps$. 
Now, $\sigma_0 \times \hat K$ is also properly convex for 
a sufficiently small $\eps$. 
For a finite set $L$, 
the convex hull $U_1$ of $\bigcup_{i\in \bZ -L } \eta^i(D_F)$ in $\mathds{A}^n$ 
is properly convex. 
The convex hull of $U_1 \cup U_L$ for
$U_L:= \bigcup_{i\in L } \eta^i(D_F)$ is still properly convex: 
Suppose not. 
Then there exists an antipodal pair in 
\[\clo(\CH(U_1 \cup U_L)) = \CH(\clo(U_1)\cup \clo(U_L)).\] 
The antipodal pair must be in $\Bd \mathds{A}^n$ 
and in $\clo(U_1) \cup \clo(U_L)$
since the interior of 
$\mathds{A}^n$ has no antipodal pair, and 
a point of the convex hull lies on a simplex with vertices in the original set. 
However, since $\clo(U_1) \cap \Bd \mathds{A}^n$ is a properly convex set 
$K'''\ast \{\mbv_{\tilde E}\}$
and $\clo(U_L)\cap \Bd \mathds{A}^n$ is $\{\mbv_{\tilde E}\}$. 
This is a contradiction. 

Let $U'$ denote the convex hull of $U_1\cup U_L$ in $\mathds{A}^n$. 
Hence, we showed that $U'$ is properly convex. 

(i-b) Now we prove the converse part of (i). 
Suppose that $\bGamma_{\tilde E}$ acts 
on a properly convex p-end neighborhood $U'$.

By Lemma \ref{np-lem-i}, we have $\aleph_7(g) \geq 0$ for 
$g \in \bGamma_{\tilde E, +}$. 
Suppose that $\aleph_7(g) =0$ for some $g\in \bGamma_{\tilde E, +}$. 
Then \[\{g^{i}(\clo(U) \cap H_l)\}\ra B \hbox{ as }  i \ra \infty \hbox{ under } \bdd_{H}\]
for a leaf $l$ and a compact domain $B$ at $H_k$ bounded by an ellipsoid.
This contradicts the premise of (i). 
Therefore, 
\begin{equation}\label{np-eqn-mu7ii}  
\mu_7(h) > 0 \hbox{ for every } h \in \bGamma_{\tilde E, +}. 
\end{equation} 

%

Suppose that $\{\mu_7(g_i)\} \ra 0$ for a sequence $g_i \in \bGamma_{\tilde E, +}$. 
We can assume that 
\[ \lambda^{Tr}_{\max}(g_i)/\lambda_{2}(g_i) > 1+\eps 
\hbox{ for a positive constant } \eps > 0\]
since we can take powers of $g_i$ without changing the $\mu_7$-values. 

Since $\{\mu_{7}(g_{i})\} \ra 0$, 
we obtain a nondecreasing sequence $\{n_i\}$, $n_i > 0$, such that 
\[\{\aleph_7(g_i^{n_i})= n_i\aleph_7(g_i)\} \ra 0 \hbox{ and }
 \{\lambda^{Tr}_{\max}(g_i^{n_i})/\lambda_{2}(g_i^{n_i})\} \ra \infty.\]   
 However, from such a sequence, we use  
 \eqref{np-eqn-kernel} to shows that 
 \[\{g_i^{n_i}(\clo(U) \cap H_{l})\} \ra B\]
 to a ball $B$ with a nonempty interior in $H_{k}$.
 Again, the premise contradicts this. 
 Hence, $\mu_7(g) > C$ for all $g \in \bGamma_{\tilde E, +}$ and a uniform constant 
 $C > 0$. 
 This proves the converse part of (i). 

(ii) By the convexity of $U'$, we obtain that it is foliated by radial segments from 
$\mbv_{\tilde E}$.

(iii) From (A), we see that $U' \subset K \times \mathds{A}^{n-i_0+1}$ since $U'$ is the convex hull of $U$.  We have 
$\clo(U') \cap \Bd \mathds{A}^n \subset K''' \ast \{\mbv_{\tilde E}\}$
since $U'$ is a p-end neighborhood and we can apply the same argument as $U$
as in (i).

Consider the projection $K^o \ra (K''')^o$ from the vertex $k$ of $K$. 
Consider \eqref{np-eqn-pik}, and recall that
 $U'$ is an orbifold bundle over $\Sigma_{\tilde E}$ with fibers that are radial rays. 
We see that there is a projection $U' \ra K'''$ that is equivariant with respect to 
the $\bGamma_{\tilde E}$-action since $\bGamma_{\tilde E}$ acts on $K$ fixing the vertex $k$.  
Since $U'$ is open, it contains a generic point projecting to a point of $K^o$. 

Let $N_{L}$ denote the image group of the action of $\bGamma_{\tilde E}$ on $L$
for $L = K$ or $L=K'''$. 
Since $\bGamma_{\tilde E}$ has 
a compact fundamental domain $J$ in $\Sigma_{\tilde E}$. 
Taking the image of $J$ under the projections, we see that there are 
compact subsets $J'$ of $K^o$ and $J'''$ of  $K^{''', o}$ mapping onto $K^o/N_K$ 
and onto $K^{''', o} /N_K$ respectively.  We may assume that $J'''$ is the image of $J'$ under
the projection from $k$. By taking a ray from $k$ with an endpoint at $J'''$,
and taking a sequence of points $p_i$ on it conveying to $p$ in $J'''$. We find $g_i\in N_K$ such that $g_i(p_i)\in J'$. 
This means that $g_i(p) \in J'''$ for all $p$. 
Since $K'''$ is properly convex with a Hilbert metric, $g_i|K'''$ is bounded 
by Proposition \ref{prelim-prop-AutK} 
and we can extract a convergent sequence. Hence, assume that $g_i|K''' \ra g_\infty$. 
We may assume that $g_i|K'''$ is convergent. 
Let $g'_i \in \bGamma_{\tilde E}$ be the one going to $g_i$.
Here, we can see that $\lambda_{\mbv_{\tilde E}}(g'_i) \ra 0$. 
Therefore, an interior point of $K'''$ is in the image of $U'$ under $g_\infty$
and as a limit point of the sequence of the images under $g'_i$ of a point of $U'$. 
By Proposition 3 of \cite{Vey}, the convex hull of any orbit of a point in $K^{''', o}$ equals 
$K^{''', o}$.  We conclude that $\clo(U')\cap \mathds{A}^n =  K''' \ast \{\mbv_{\tilde E}\}$
since $\mbv_{\tilde E}$ is obviously in it and the intersection has to be convex.



(iv)
Suppose that $\aleph_7$ is identically zero for $\bGamma_{\tilde E, +}$. 
Then by \eqref{np-eqn-mu7ii} in the proof of (i), 
  \[\{g^{i}(\clo(U) \cap H_l)\} \ra B \hbox{ as }  i \ra \infty \hbox{ under } \bdd_{H}\]
 for a leaf $l$ and a compact domain $B$ at $H_k$ bounded by an ellipsoid.
 Choose $\hat B$ the maximal domain of the form $B$ arising from the situation. 
Then we may show that $\CH(U) = (K''' \ast \hat B)^o$
by Proposition \ref{prelim-prop-decjoin} since we can find a sequence $g_i$ such that $g_i| K'''$ is bounded and 
$\{\lambda_{K'''}(g)/\lambda_{k}(g)\} \ra \infty$ 
since $(K''\ast k)^o/\bGamma_{\tilde E}$ is compact. 
In addition, $(K''' \ast \hat B)^o$ is properly convex.

In contrast, we have $\aleph_7 \geq 0$ according to Lemma \ref{np-lem-i}.
By premise $\CH(U) = (K''' \ast B)^o$ where $B$ is a convex open ball
in a hemisphere $H_k$ in $\SI^{i_0}_{k}$. Since $\mathcal{N}$ acts on $\CH(U)$, 
$B$ bounded by an ellipsoid. 
If $\aleph_7(g) > 0$ for some $g$, Then $g$ acts on 
$B$ such that $g(B)$ is a translated image of the region $B$
bounded by a paraboloid in the affine subspace $H_{k}^o$.
We obtain $\bigcup_{k=1}^\infty g^{-k}(B) = H_{k}^o$.
This contradicts the proper convexity, and therefore $\aleph_7$ is identically $0$.

\end{proof} 


\begin{definition}\label{np-defn-qjoin} $ $ 
\begin{itemize}
\item Generalizing Example \ref{np-exmp-joined},  
an R-p-end $\tilde E$ satisfying the case (ii) of Proposition \ref{np-prop-qjoin}
is a {\em strictly joined R-p-end} (of {\em a totally geodesic R-end and a horospherical end}\/)
and $\bGamma_{\tilde E}$ now is called a {\em strictly joined end group}. 
In addition, any end finitely covered by a strictly joined R-end is called
a {\em strictly joined R-end}. 
\item An R-p-end $\tilde E$ satisfying the case (i) of Proposition \ref{np-prop-qjoin}, 
 is a {\em quasi-joined} R-p-end (of {\em a totally geodesic R-end and a horospherical end}\/) corresponding to Definition \ref{np-defn-quasijoin}
and $\bGamma_{\tilde E}$ now is a  quasi-joined end holonomy group.
\index{end!quasi-joined|textbf} 
\index{end!holonomy group!quasi-joined|textbf} 
%
\end{itemize}
Also, any p-end $\tilde E$ with $\bGamma_{\tilde E}$ is a finite-index subgroup of $\bGamma_{\tilde E}$ as above is called by
the corresponding names. 
\end{definition}

\subsubsection{The non-existence of strictly joined cases for $\mu\equiv 1$.} 
\label{np-subsub-quasijoin}

\begin{corollary}\label{np-cor-NPNCcase2} 
Let $\Sigma_{\tilde E}$ be the end orbifold of an NPNC R-p-end $\tilde E$ of a
properly convex real projective $n$-orbifold $\orb$ with radial or totally geodesic ends. 
Assume that the holonomy group of $\orb$ is strongly irreducible.
Let $\bGamma_{\tilde E}$ be the p-end holonomy group.
Assume Hypotheses \ref{np-h-norm} only and $\mu_{g} = 1$ for all $g \in \bGamma_{\tilde E}$. 
Then $\tilde E$ is not a strictly joined end. 
\end{corollary}
\begin{proof}
Suppose that $\tilde E$ is a strictly joined end.
By premise, $\mu_g = 1$ for all $g \in \bGamma_{\tilde E}$. 
By Lemma \ref{np-lem-conedecomp1} and Proposition \ref{np-prop-decomposition}, 
every $g \in \bGamma_{\tilde E}$ is of form: 
\renewcommand{\arraystretch}{1.5}
\begin{equation}\label{np-eqn-gform}
\newcommand*{\tempV}{\multicolumn{1}{r|}{}}
\left( \begin{array}{ccccccc} 
S_{g} & \tempV & 0 & \tempV & 0 & \tempV & 0 \\ 
 \cline{1-7}
0 &\tempV & \lambda_g &\tempV & 0 &\tempV & 0 \\ 
 \cline{1-7}
 0 &\tempV & \lambda_g \vec{v}^T_g &\tempV & \lambda_g O_5(g) &\tempV & 0 \\ 
 \cline{1-7}
0 &\tempV & \lambda_g\left(\aleph_7(g) + \frac{\llrrV{\vec{v}_g}^2}{2}\right) &\tempV &
\lambda_g \vec{v}_g O_5(g) &\tempV & \lambda_g 
\end{array} 
\right)
\end{equation}


By Proposition \ref{np-prop-decomposition}, 
we obtain a sequence $\gamma_m$  from step (D) of the proof of the form: 
\renewcommand{\arraystretch}{1.5}
\begin{equation}\label{np-eqn-gammao}
\newcommand*{\tempV}{\multicolumn{1}{r|}{}}
\left( \begin{array}{ccccccc} 
\delta_m S_m & \tempV & 0 & \tempV & 0 & \tempV & 0 \\ 
 \cline{1-7}
0 &\tempV & \lambda_m &\tempV & 0 &\tempV & 0 \\ 
 \cline{1-7}
 0 &\tempV & \lambda_m \vec{v}^T_m &\tempV & \lambda_m O_5(\gamma_m) &\tempV & 0 \\ 
 \cline{1-7}
0 &\tempV & \lambda_m\left(\aleph_7(\gamma_m) + \frac{\llrrV{\vec{v}_m}^2}{2}\right) &\tempV &
\lambda_m \vec{v}_m O_5(\gamma_m) &\tempV & \lambda_m 
\end{array} 
\right)
\end{equation}
as  $C_{1, m} =0$ and $c_{2, m}=0$
where 
\begin{itemize}
\item $\{\lambda_m\} \ra \infty$, 
\item $\{\delta_m\} \ra 0$, 
\item $\{S_m\}$ is in a sequence of bounded matrices in $\SL_{\pm}(n-i_0-1)$, and 
\item $\aleph_7(\gamma_m ) = 0$ by Proposition \ref{np-prop-qjoin} (ii).

\end{itemize}
Moreover, Hypothesis \ref{np-h-qjoin} now holds. 
By Lemma \ref{np-lem-noball}, we obtain a contradiction. 
\end{proof} 

\subsection{The proof for discrete $N_{K}$.} \label{np-sub-discretecase} 


Now, we go to prove Theorem \ref{np-thm-thirdmain} when $N_{K}$ is discrete. 
Taking a finite-index torsion-free subgroup if necessary by Theorem \ref{prelim-thm-vgood}, we may assume that $N_{K}$ acts freely on $K^{o}$. 
We have a corresponding orbifold fibration 
\begin{eqnarray}
l/N & \ra & \tilde \Sigma_{\tilde E}/\bGamma_{\tilde E} \nonumber \\ 
     &       & \, \, \downarrow \nonumber \\ 
   &         & K^o/N_K 
\end{eqnarray}
where the fiber and the quotients are compact orbifolds
since $\Sigma_{\tilde E}$ is compact. 
Here, the fiber equals $l/N$ for generic $l$. 
The action of $N_K$ on $K$ is semisimple by Theorem 3 of Vey \cite{Vey}. 

Since $N$ acts on each leaf $l$ of ${\mathcal F}_{\tilde E}$ in $\tilde \Sigma_{\tilde E}$, 
it also acts on a properly convex domain $\torb$ and $\mbv_{\tilde E}$ in a subspace $\SI^{i_0+1}_l$ in $\SI^n$ 
corresponding to $l$. $l/N \times \bR$ is an open real projective orbifold diffeomorphic to
$(H^{i_0+1}_l \cap \torb)/N$
for an open hemisphere $H^{i_0+1}_l$ corresponding to $l$. 
Since elements of $N$ restricts to $\Idd$ on $K$, 
we obtain 
\[\lambda^{Tr}_{\max}(g) = \lambda^{Tr}_{\min}(g) \hbox{ for all } g \in N:\]
Otherwise, we can easily see that $g$ acts not trivially on $\SI^{n-i_{0}-1}$. 
By Proposition \ref{np-prop-eigSI}, all the norms of eigenvalues are $1$. 
Let $P_l$ denote the smallest subspace containing $\mbv_{\tilde E}$ in the direction of $l$ in 
$\tilde \Sigma_{\tilde E}$. 
\begin{itemize} 
\item Since $l$ is a complete affine subspace, 
Lemma \ref{ce-lem-unithoro} applied to $P_l \cap \torb/N$ shows that
$l$ covers a horospherical end of $(\SI^{i_0+1}_l \cap \torb)/N$.
\item 
By Lemma \ref{ce-lem-unithoro} applied to $P_l \cap \torb/N$, $N$ is virtually unipotent, 
and $N$ is virtually a cocompact subgroup of a unipotent group,
and $N|\SI^{i_0+1}_l$ can be conjugated into a maximal parabolic subgroup of
$\SO(i_0+1, 1)$ in $\Aut(\SI^{i_0+1}_l)$ 
and acting on an ellipsoid of dimension $i_0$ in $H^{i_0+1}_l$. 
\end{itemize} 

We verify Hypothesis \ref{np-h-norm}. 

By the nilpotent Lie group theory of 
Malcev \cite{Malcev49}, the Zariski closure ${\mathcal{Z}}(N)$ of $N$ is a virtually simply connected nilpotent Lie group with finitely many components 
and  ${\mathcal{Z}}(N)/N$ is compact. 
Let $\CN$ denote the identity component of the Zariski closure of $N$
such that $\CN/(\CN \cap N)$ is compact.
$\CN \cap N$ acts on the great sphere $\SI^{i_0+1}_l$ containing $\mbv_{\tilde E}$ and corresponding to $l$.
Since $\CN/(\CN\cap N)$ is compact, 
we can modify $U$ such that $\CN$ acts on $U$ by Lemma \ref{prelim-lem-endnhbd}: 
i.e.,  we take the interior of $\bigcap_{g\in \CN} g(U) = \bigcap_{g\in F} g(U)$
for the fundamental domain $F$ of $\CN$ by $N$.  

Since $\bGamma_{\tilde E}$ normalizes $N$, it also normalizes
the identity component $\mathcal{N}$. 

By the above, $\CN|\SI^{i_0+1}_l$ is conjugate into a parabolic subgroup of
$\SO(i_0+1, 1)$ in $\Aut(\SI^{i_0+1}_l)$, and 
$\CN$ acts on $U \cap \SI^{i_0+1}_l$, which is a horoball for each leaf $l$ of
$\tilde \Sigma_{\tilde E}$. 

Taking a finite-index cover of $U$, we can assume that 
$N \subset \CN$ since  ${\mathcal{Z}}(N)$  is a finite extension of $\CN$. 
We denote the finite-index group by $\bGamma_{\tilde E}$ again. 


Since $\SI^{i_0+1}_l$ corresponds 
to a coordinate $(i_0+2)$-subspace, and 
$\SI^{i_0}_\infty$ and $\{\mbv_{\tilde E}\}$ are 
$\bGamma_{\tilde E}$-invariant, 
we can choose coordinates such that 
\eqref{np-eqn-matstd} and \eqref{np-eqn-nilmatstd}
hold. Hence, Hypothesis \ref{np-h-norm} holds. 

\begin{theorem}\label{np-thm-NPNCcase} 
Let $\Sigma_{\tilde E}$ be the end orbifold of an NPNC R-p-end $\tilde E$ of 
a properly convex real projective $n$-orbifold $\orb$ with radial or totally geodesic ends. 
Assume that the holonomy group of $\pi_{1}(\orb)$ is strongly irreducible.
Let $\bGamma_{\tilde E}$ be the p-end holonomy group
satisfying the transverse weak middle-eigenvalue condition
with respect to R-p-end structure. 
Assume alsuch that $N_{K}$ is discrete, and  $K^{o}/N_{K}$ is compact
and Hausdorff. 
Then
$\tilde E$ is a quasi-joined R-end. 
\end{theorem}
\begin{proof} 
By Lemma \ref{np-lem-similarity}, 
$h(g) \CN(\vec{v}) h(g)^{-1} = \CN( \vec{v} M_g)$ where 
$M_g$ is a scalar multiplied by an element of a copy of an orthogonal group $\Ort(i_0)$. 

The group $\CN$ is isomorphic to $\bR^{i_{0}}$ as a Lie group. 
Since $N \subset \CN$ is a discrete cocompact, $N$ is virtually isomorphic to $\bZ^{i_0}$.
Without loss of generality, we assume that $N$ is a cocompact subgroup of $\CN$. 
By normality of $N$ in $\bGamma_{\tilde E}$, we obtain
$h(g)N h(g)^{-1} = N$ for $g \in \bGamma_{\tilde E}$. 
Since $N$ corresponds to a lattice $L \subset \bR^{i_{0}}$ by the map $\CN$, 
the conjugation by $h(g)$ corresponds 
to an isomorphism $M_g: L \ra L$
by Lemma \ref{np-lem-similarity}.
When we identify $L$ with $\bZ^{i_0}$, 
$M_g: L \ra L$ is represented by an element of 
$\SL_\pm(i_0, \bZ)$ since $(M_g)^{-1} = M_{g^{-1}}$. 
Also, by Lemma \ref{np-lem-similarity}, 
$\{M_g| g \in \bGamma_{\tilde E}\}$ 
is a compact group as their determinants equal $\pm 1$. 
Hence, 
the image of the homomorphism given by
$g \in h(\pi_1(\tilde E)) \mapsto M_g \in \SL_\pm(i_0, \bZ))$ is a finite-order group. 
Moreover, $\mu_{g}= 1$ for every $g\in \bGamma_{\tilde E}$
as we can see from Lemma \ref{np-lem-similarity}. 
Thus, $\bGamma_{\tilde E}$ has a finite-index group $\bGamma_{\tilde E}'$ 
centralizing $\CN$. 

We conclude that Hypothesis \ref{np-h-norm} is satisfied by taking
a finite-index subgroup of $\bGamma_{\tilde E}$ if necessary. 

%

We can now use Proposition \ref{np-prop-decomposition} by letting 
$G$ be $\CN$ since $N_K$ is discrete and $\bar N_K = N_K$. 
We take $\Sigma_{E'}$ as the corresponding cover of $\Sigma_{\tilde E}$. 
By Lemma  \ref{np-lem-conedecomp1} and Proposition 
\ref{np-prop-decomposition}, 
Hypothesis \ref{np-h-qjoin} holds, and 
we have the result needed to apply Proposition \ref{np-prop-qjoin}. 
Finally, Proposition \ref{np-prop-qjoin}(i) and (ii) imply that $\bGamma_{\tilde E}$ virtually is either a join or a quasi-joined group.
Corollary \ref{np-cor-NPNCcase2} shows that a strictly joined end cannot occur.  
\end{proof}



\section{The nondiscrete case} \label{np-sec-nondiscrete}


This is in part a joint work with Y. Carri\`ere. 
Let $\Sigma_{\tilde E}$ be the end orbifold of an NPNC R-end $\tilde E$ of 
a properly convex real projective $n$-orbifold $\orb$ with radial or totally geodesic ends. 
Let $\bGamma_{\tilde E}$ be the p-end holonomy group. 
Let $U$ be a neighborhood of p-end in $\torb$ that corresponds to a p-end vertex $\mbv_{\tilde E}$. 




Recall the exact sequence 
\[ 1 \ra N \ra \pi_1(\tilde E) \stackrel{\Pi^*_K}{\longrightarrow} N_K \ra 1 \] 
where we assume that $N_K \subset \Aut(K)$ is not discrete. 
Since $\tilde \Sigma_{\tilde E}/\bGamma_{\tilde E}$ is compact, 
$N_K$ acts cocompactly on $K^o$. 
However, $N_K$ is not necessarily semisimple. 
(See Section \ref{np-sub-counter}.) 


An element $g \in \bGamma_{\tilde E}$ is of the form: 
\begin{equation} \label{np-eqn-g}
\newcommand*{\tempV}{\multicolumn{1}{r|}{}}
g = \left( \begin{array}{ccc} 
K(g) & \tempV &  0 \\ 
 \cline{1-3}
* &\tempV & U(g)
\end{array} 
\right).
\end{equation}
Here $K(g)$ is an $(n-i_0)\times (n-i_0)$-submatrix and 
$U(g)$ is  an $(i_0+1)\times (i_0+1)$-submatrix acting on $\SI^{i_0}_\infty$. 
We note $\det K(g) \det U(g) = 1$. 

\subsection{Outline of Section \ref{np-sec-nondiscrete}}  



In Section \ref{np-sub-leafclosure}, we take the leaf closure of 
a complete affine $i_0$-dimensional leaf. 
The theory of Molino \cite{Molino88}. 
shows that the space of leaf-closures are
an orbifold. 
In Section \ref{np-subsec-leaves}, we show that 
a leaf-closure is a compact suborbifold in $\Sigma_{\tilde E}$. 
In Section \ref{np-sub-uniN}, 
we show that
the fundamental group of each leaf-closure is virtually solvable.
In Section \ref{np-subsub-syndetic}, we find 
a syndetic closure $S$ according to a theory of Fried-Goldman \cite{FG83}.
From this, we find a subgroup acting on each
complete affine $i_0$-dimensional leaf in Section \ref{np-subsub-US}
which we show to be a cusp group.

In Section \ref{np-sub-nondproof}, we complete the proof of Theorem 
\ref{np-thm-thirdmain} not covered by earlier Theorem \ref{np-thm-NPNCcase}
for discrete $N_K$. 
Proposition \ref{np-prop-NKsemi} shows that $N_K$ is semisimple. 
Proposition \ref{np-prop-mug} shows that $\mu_g =1$ for every $g \in \bGamma_{\tilde E}$. Finally, we prove Theorem \ref{np-thm-thirdmain}. 

\subsection{Taking the leaf closure} \label{np-sub-leafclosure}

\subsubsection{Estimations with $KA \bU$.}

Let $\bU$ denote a maximal nilpotent subgroup of $\SL_\pm(n+1, \bR)_{\SI^{i_0}_\infty, \mbv_{\tilde E}}$  
given 
by lower triangular matrices with diagonal entries equal to $1$.

Recall the foliation on $\tilde \Sigma_{\tilde E}$ given by fibers of $\Pi_K$
has leaves that are $i_0$-dimensional complete affine subspaces.
We denoted it by $\widetilde{\mathcal{F}}_{\tilde E}$. 
Then $K^o$ admits a smooth Riemannian metric $\mu_K$ invariant under $N_K$ by Lemma \ref{np-lem-invmet}. 
We consider the orthogonal frame bundle $\bF K^o$ over $K^o$. 
A metric on each fiber of $\bF K^o$ is induced from $\mu_K$.
Since the action of $N_{K}$ is isometric on $\bF K^o$ with trivial stabilizers, 
$N_{K}$ acts on a smooth orbit submanifold 
of $\bF K^o$ transitively 
with trivial stabilizers. 
(See Lemma 3.4.11 in \cite{Thurston97}.)

There exists a bundle $\bF\tilde \Sigma_{\tilde E}$ from pulling back $\bF K^o$ by the projection map. 
Here, $\bF\tilde \Sigma_{\tilde E}$ covers $\bF\Sigma_{\tilde E}$. 
Since $\bGamma_{\tilde E}$ acts isometrically on $\bF K^o$, 
the quotient space $\bF\tilde \Sigma_{\tilde E}/\bGamma_{\tilde E}$ is a bundle $\bF\Sigma_{\tilde E}$ over 
$\Sigma_{\tilde E}$ 
with compact fibers diffeomorphic to the orthogonal group of dimension $n-i_0$. 
Also, $\bF\tilde \Sigma_{\tilde E}$ is foliated by $i_0$-dimensional affine 
spaces pulled-back from the $i_0$-dimensional leaves on the foliation $\tilde \Sigma_{\tilde E}$. 
We can think of these leaves as the inverse images of points of
$\bF K^{o}$.

\begin{lemma} \label{np-lem-polynomial}
Each leaf $l$ of $\bF\Sigma_{\tilde E}$ is of polynomial growth. That is, each ball $B_R(x)$ in $l$ of radius $R$ for $x\in l$ has an area 
less than equal to $f(R)$ for a polynomial $f$ where we are using an arbitrary Riemannian metric on 
$\bF\tilde \Sigma_{\tilde E}$ induced from one on $\bF\Sigma_{\tilde E}$. 
\end{lemma}
\begin{proof} 
Let us choose a fundamental domain $F$ of $\bF\Sigma_{\tilde E}$. 
Then  for each leaf $l$ 
there exists an index set $I_{l}$ such that 
$l$ is a union of $g_i(D_i)$  $i \in I_l$ for the intersection $D_i$ of a leaf with $F$
and $g_i \in \bGamma_{\tilde E}$.  
We have that $D_i \subset D'_i$ where $D'_i$ is an $\eps$-neighborhood of $D_i$ in the leaf. 
Then \[\{g_i(D'_i)| i \in I_l \}\] cover $l$ in a locally finite manner. 
In fact, the number of elements covering a single point is bounded above by a number 
independent of the point. That is, the cover is so-called uniformly locally finite. 
The subset $G(l):= \{g_i\in \Gamma| i \in I_l\}$ is a discrete subset. 

Choose an arbitrary point $d_i \in D_i$ for every $i \in I_l$. 
The set $\{g_i(d_i)| i \in I_l \}$ and $l$ is quasi-isometric: 
a map from $G(l)$ to $l$ is given by $f_1: g_i \mapsto g_i(d_i)$ 
and the multivalued map $f_2$ from $l$ to $G(l)$ given by sending each point $y \in l$ 
to one of finitely many $g_i$ such that $g_i(D'_i) \ni y$. 
Let $\bGamma_{\tilde E}$ be given the Cayley metric and $\tilde \Sigma_{\tilde E}$ a metric induced 
from $\Sigma_{\tilde E}$. 
Both maps are quasi-isometries since 
these maps are restrictions of quasi-isometries $\bGamma_{\tilde E} \ra \tilde \Sigma_{\tilde E}$
and $\tilde \Sigma_{\tilde E} \ra \bGamma_{\tilde E}$ defined analogously.  

The action of $g_i$ in $K$ is bounded since it sends some points of $\Pi_K(F)$ to those of $\Pi_K(F)$ which is a compact set in $K^o$. 
Thus, $\Pi^*_K(g_i)$ goes to a bounded subset of $\Aut(K)$. 
In the form \eqref{np-eqn-g},
\[K(g_i) = \det(K(g_i))^{1/(n-i_0)} \hat K(g_i) \hbox{ where } \hat K(g_i) \in \SL_{\pm}(n-i_0, \bR)\]
where $\hat K(g_i)$ is uniformly bounded. 
Let $\tilde \lambda_{\max}(g_i)$ and $\tilde \lambda_{\min}(g_i)$ denote the largest norm and the smallest norm of
eigenvalues of $\hat K(g_i)$. 
Since $\Pi^{\ast}_{K}(g_{i})$ are in a bounded set of $\Aut(K)$, 
we obtain
\begin{equation}\label{np-eqn-lambdaB} 
 \frac{1}{C} \leq \tilde \lambda_{\max}(g_i), \tilde \lambda_{\min}(g_i)\leq C
\end{equation}
for $C> 1$ independent of $i$. 
The largest and smallest eigenvalues of $g_i$ are equal.
\[\lambda^{Tr}_{\max}(g_i) = \det(K(g_i))^{1/(n-i_0)}\tilde \lambda_{\max}(g_i) \hbox{ and } 
\lambda^{Tr}_{\min}(g_i) = \det(K(g_i))^{1/(n-i_0)}\tilde \lambda_{\min}(g_i)\]
by Proposition \ref{np-prop-eigSI}.
Denote by $a_j(g_i), j=1, \dots, i_0+1$, the norms of eigenvalues of $g_i$ associated with $\SI^{i_0}_\infty$ 
where $a_1(g_i)\geq \cdots \geq a_{i_0+1}(g_j) > 0$ 
with repetitions allowed. 
Since $\det g_i = 1$, we have 
\[\det(K(g_i)) a_1(g_i) \dots a_{i_0+1}(g_i) = 1.\] 
If $\{|\det(K(g_i))|\} \ra 0$, then 
$\{a_1(g_j)\} \ra \infty$ whereas by \eqref{np-eqn-lambdaB} 
\[\{\det(K(g_i))^{1/(n-i_0)} \tilde \lambda_{\max} (g_i)\} \ra 0\]
contradicting Proposition \ref{np-prop-eigSI}.
If $\{|\det(K(g_i))|\} \ra \infty$, then 
$\{a_{i_0+1}\} \ra 0$ whereas by \eqref{np-eqn-lambdaB} 
\[\{\det(K(g_i))^{1/(n-i_0)} \tilde \lambda_{\min} (g_i)\} \ra \infty\]
contradicting Proposition \ref{np-prop-eigSI}.
Therefore, we obtain 
\[1/C < |\det(K(g_i))| < C\] for a positive constant $C$. 
We deduce that  
the largest norm and the smallest norm of eigenvalues of $g_i$
\[\det(K(g_i))^{1/(n-i_0)}\tilde \lambda_{\max}(g_i) \hbox{ and } 
\det(K(g_i))^{1/(n-i_0)}\tilde \lambda_{\min}(g_i)\]
are bounded above and below by two positive numbers. 
Hence, $\lambda^{Tr}_{\max}(g_i)$ and $\lambda^{Tr}_{\min}(g_i)$ 
are all bounded above and below by a fixed set of positive numbers. 
By Proposition \ref{np-prop-eigSI}, 
the remaining norms of eigenvalues of $g$ 
are bounded above and below by the same fixed set of positive numbers.

By Corollary \ref{prelim-cor-bdunip}, $\{g_i\}$ is of bounded distance from $\bU'$.
Let $N_c(\bU')$ denote a $c$-neighborhood of $\bU'$. 
Then \[G(l)  \subset N_c(\bU') \hbox{ for some } c> 0.\] 

Let $d$ denote the left-invariant metric on $\Aut(\SI^n)$. 
By the discreteness of $\bGamma_{\tilde E}$, the set $G(l)$ is discrete 
and there exists a lower  
bound to \[\{d(g_i, g_j)| g_i, g_j \in G(l), i \ne j\}.\]
Also, given any $g_i \in G(l)$, there exists
an element $g_j \in G(l)$ such that
$d(g_i, g_j) < C$ for a uniform constant $C$.
(We need to choose $g_j \in G(l)$ such that $g_j(F)$ is 
adjacent to $g_i(F)$.)
Let $B_R(\Idd)$ denote the ball in $N_c(\bU')$ of radius $R$ with the center $\Idd$. 
Then $B_R(\Idd) \subset N_c(\bU')$ has a polynomial growth with respect to $R$, and so is $G(l)$ by Lemma \ref{prelim-lem-bU}. 
Since the collection $\{g_i(D'_i)| g_i \in G(l)\}$ 
of uniformly bounded balls covers $l$ in a uniformly locally finite manner, 
and $G(l)$ is discrete with minimal distances between the elements, 
it follows that 
$l$ also has a polynomial growth. 
\end{proof}


\subsubsection{Closures of leaves} \label{np-subsec-leaves}

We consider the orthogonal frame bundle $\bF K^o$ over $K^o$. 
A metric on each fiber of $\bF K^o$ is induced from $\mu_K$.
Since the action of $N_{K}$ is isometric on $\bF K^o$ with trivial stabilizers, 
we find that $N_{K}$ acts on a smooth orbit submanifold 
of $\bF K^o$ transitively 
with trivial stabilizers. 
(See Lemma 3.4.11 in \cite{Thurston97}.)

There exists a bundle $\bF\tilde \Sigma_{\tilde E}$ from pulling back $\bF K^o$ by the projection map. 
Here, $\bF\tilde \Sigma_{\tilde E}$ covers $\bF\Sigma_{\tilde E}$. 
Since $\bGamma_{\tilde E}$ acts isometrically on $\bF K^o$, 
the quotient space $\bF\tilde \Sigma_{\tilde E}/\bGamma_{\tilde E}$ is a bundle $\bF\Sigma_{\tilde E}$ over 
$\Sigma_{\tilde E}$ with a subbundle with compact fibers isomorphic to the orthogonal group of dimension $n-i_0$. 
Also, $\bF\tilde \Sigma_{\tilde E}$ is foliated by $i_0$-dimensional affine 
spaces pulled-back from the $i_0$-dimensional leaves on the foliation $\tilde \Sigma_{\tilde E}$. 
We can think of these leaves as the inverse images of points of
$\bF K^{o}$.



\subsubsection{$\pi_{1}(V_{l}) $ is virtually solvable.} \label{np-sub-uniN}  

Recall the fibration \[\Pi_K:  \tilde \Sigma_{\tilde E} \ra K^o
\hbox{ which induces } \tilde \Pi_K: \bF \tilde \Sigma_{\tilde E} \ra \bF K^o.\]
Since $N_K$ acts as isometries of a Riemannian metric on 
$K^o$, we can obtain a metric on $\Sigma_{\tilde E}$ 
such that the foliation is the Riemannian foliation. 
Let $p_{\Sigma_{\tilde E}}: \bF\tilde \Sigma_{\tilde E} \ra 
\bF \Sigma_{\tilde E}$ be the projective covering map
induced from $\tilde \Sigma_{\tilde E} \ra  \Sigma_{\tilde E}$.
The foliation on $\tilde \Sigma_{\tilde E}$ gives us a foliation of $\bF \tilde \Sigma_{\tilde E}$. 


Let $A_K$ be the identity component of the closure of $N_K$ the image of $\bGamma_{\tilde E}$ in $\Aut(K)$, which is a Lie group of $\dim \geq 1$. 
\begin{proposition} \label{np-prop-Alcentral} 
	$A_K$ is a normal connected nilpotent Lie subgroup of 
	the closure of $N_K$.
	\end{proposition} 
\begin{proof} 
	Since the closure of $N_K$ is 
	normalized by $N_K$, $A_K$ is also normalized by $N_K$. 
	Since $l$ maps to a polynomial growth leaf in $\bF\Sigma_{\tilde E}$ by Lemma \ref{np-lem-polynomial},
The main theorem of	Carri\`ere \cite{Carriere88} shows that $A_K$ is a connected nilpotent Lie group
	in the closure of $N_K$ in $\Aut(K)$ acts on $\bF K^o$ freely.
\end{proof} 
	

Let $l$ be a leaf of $\bF \tilde \Sigma_{\tilde E}$,
and let $p$ be the image of $l$ in $\bF K^o$. 
Moreover, we have 
\begin{alignat}{3} 
\tilde \Pi_K^{-1}(A_K(p))   =: \tilde V_{l} &\, \hookrightarrow \, & \bF \tilde \Sigma_{\tilde E} \nonumber \\
\downarrow \, &                     &\, p_{\Sigma_{\tilde E}} \downarrow \nonumber \\
V_l      & \hookrightarrow & \bF \Sigma_{\tilde E}
\end{alignat} 
for $V_l := \overline{p_{\Sigma_{\tilde E}}(l)}$ in $\bF \Sigma_{\tilde E}$.
Since $\tilde V_l$ is closed and is a component of 
the inverse image of $V_l$ which is a union of copies of 
$\tilde V_l$, the image $V_l$ is a compact submanifold. 
Note that $V_{l}$ has a dimension independent of $l$ since $A_{K}$ acts freely.

Now, $N$ is precisely the subgroup of $\pi_1(V_l)$ 
fixing a leaf $l$ in $\bF K^o$. 
For each closure $V_l$ of a leaf $l$, the manifold $V_l$ is 
a compact submanifold of $\bF \Sigma_{\tilde E}$, and 
we have an exact sequence 
\begin{equation}\label{np-eqn-oiVl}  
 1 \ra N \ra h(\pi_1(V_l)) \stackrel{\Pi^*_K}{\longrightarrow} A'_K \ra 1.
 \end{equation}
Since the leaf $l$ is dense in $V_l$, it follows that $A'_K$ is dense in $A_K$. 
Each leaf $l'$ of $\tilde \Sigma_{\tilde E}$ has a realization as a subset in $\torb$. 
We have the norms of eigenvalues $\lambda_i(g) = 1$ for $g \in N$  
by Proposition \ref{np-prop-eigSI}. 
By Theorem \ref{prelim-thm-orthopotent}, 
$N=N_l$ is virtually unipotent
since the norms of eigenvalues equal $1$ identically and $N_l$ is discrete. 

We take a finite cover of $\Sigma_{\tilde E}$ such that $N$ is nilpotent. 
Hence, $h(\pi_{1}(V_{l}))$ is virtually 
solvable being an extension of a nilpotent group by a nilpotent group. 
%
%
%
%
We summarize in the following. 
\begin{proposition}\label{np-prop-V_l} 
Let $l$ be a generic fiber of $\bF \tilde \Sigma_{\tilde E}$ and $p$ be the corresponding point of $\bF K^o$. 
Then there exists a nilpotent group $A_K$ acting on $\bF K^o$ such that 
$\tilde \Pi_K^{-1}(A_K(p)) =\tilde V_l$ covers a compact suborbifold $V_l$ in $\bF \Sigma_{\tilde E}$, 
a conjugate of the image of the holonomy group of $V_l$ is a dense subgroup of $A_K$, and 
the holonomy group of $V_{l}$ is virtually solvable. Moreover, $\tilde V_{l}$ is homeomorphic to a torus times a cell or a cell.  
\end{proposition} 
\begin{proof} 
We just need to prove the last statement. 
Since $A_K$ is a connected nilpotent group, 
$A_{K}$ is homeomorphic to a torus times a cell or a cell,
and so is the free orbit in $\bF K^o$. 
Since $\tilde \Pi_{K}$ has fibers that are $i_{0}$-dimensional open hemispheres, 
this last statement follows. 
\end{proof}

We note that $A_K$ is nilpotent, but may not be unipotent. 

\begin{remark} \label{np-rem-Molino}
The leaf holonomy acts on $\bF \tilde \Sigma_{\tilde E}/\widetilde{\mathcal{F}}_{\tilde E}$ as a nilpotent killing field group
without any fixed points. 
Hence, each leaf $l$ is in $\tilde V_l$ with a constant dimension. 
Thus, $\mathcal{F}_{\tilde E}$ is a foliation with leaf closures of identical dimensions. 

The leaf closures form another foliation $\overline{\mathcal{F}}_{\tilde E}$ with compact leaves by Lemma 5.2 of Molino \cite{Molino88}. 
We let $\bF \Sigma_{\tilde E}/\overline{\mathcal{F}}_{\tilde E}$ denote the space of closures of leaves. It has an orbifold structure,
where the projection $\bF \Sigma_{\tilde E} \ra \bF \Sigma_{\tilde E}/\overline{\mathcal{F}}_{\tilde E}$ is an orbifold bundle projection
by Proposition 5.2 of \cite{Molino88}.
\end{remark} 


\subsubsection{The holonomy group for a leaf closure is virtually normalized by the end holonomy group.}\label{np-subsub-syndetic}

Note that $\bGamma_l$ is the deck transformation group of $\tilde V_l$ 
over $V_l$. Since $\tilde V_l$ is the inverse image of $A_K(x)$ for $x 
\in \bF K^o$,  
$\bGamma_l$ is the inverse image of $N_K \cap A_K$ under
$\Pi^\ast_K$. 
Since $N_K \cap A_K$ is normal in $N_K$,
$\bGamma_l$ is a normal subgroup of $\bGamma_{\tilde E}$. 

Recall that $\bGamma_{l}$ is virtually solvable, as we showed above. 
We let ${\mathcal{Z}}(\bGamma_{\tilde E})$ and ${\mathcal{Z}}(\bGamma_l)$ denote the Zariski closures in $\Aut(\SI^n)$ of 
$\bGamma_{\tilde E}$ and $\bGamma_l$ respectively. 

By Theorem 1.6 of Fried-Goldman \cite{FG83}, 
there exists a closed virtually solvable Lie group $S_l$ containing $\bGamma_l$ 
with the following four properties: 
\begin{itemize}
\item $S_l$ has finitely many components.
\item $\bGamma_l\backslash S_l$ is compact. 
\item The Zariski closure ${\mathcal{Z}}(S_l)$ is the same as 
${\mathcal{Z}}(\bGamma_l)$. 
\item Finally, we have solvable ranks 
\begin{equation}\label{np-eqn-rankSl}
{\mathrm{rank}}(S_l) \leq {\mathrm{rank}} (\bGamma_l). 
\end{equation}
\end{itemize}
We call this the {\em syndetic hull } of $\bGamma_l$. 
\index{syndetic hull}

We summarize: 
\begin{lemma}\label{np-lem-V} 
$h(\pi_1(V_l))$ is virtually solvable and  is contained in a virtually solvable Lie group $S_l \subset {\mathcal{Z}}(h(\pi_1(V_l))$ 
with finitely many components, and $S_l/h(\pi_1(V_l))$ is compact. 
$S_l$ acts on $\tilde V_l$. 
Furthermore, one can modify a p-end neighborhood $U$ such that 
$S_l$ acts on it. Also the Zariski closure of $h(\pi_1(V_l))$ is the same as that of $S_l$. 
\end{lemma} 
\begin{proof} 
By above, ${\mathcal{Z}}(S_l) = {\mathcal{Z}}(\bGamma_l)$ acts on $\tilde V_l$
and normalizes $\bGamma_l$. 
We need to prove about the p-end neighborhood only. 
Let $F$ be a compact fundamental domain of $S_l$ under the $\Gamma_l$. 
Then we have 
\[ \bigcap_{g\in S_l} g(U) = \bigcap_{g \in F} g(U).\]
By Lemma \ref{prelim-lem-endnhbd}, the latter set contains
a $S_l$-invariant p-end neighborhood. 
\end{proof} 

From now on, we let $S_l$ to denote only the identity component of itself for simplicity
as $S_l$ has finitely many components to begin with. 
This is sufficient for our purposes since we only need a cusp group. 

Since $S_l$ acts on $U$ and hence on $\tilde \Sigma_{\tilde E}$ as shown in Lemma \ref{np-lem-V}, 
we have a homomorphism $S_l \ra \Aut(K)$.  
We define by $S_{l, 0}$ the kernel of this map. 
Then $S_{l, 0}$ acts on each leaf of $\tilde \Sigma_{\tilde E}$.
We have an exact sequence 
\begin{equation}\label{np-eqn-Sl}
1 \ra S_{l, 0} \ra S_l \ra A_K \ra 1.
\end{equation}

\subsubsection{The form of $C_{l, 0}$.} \label{np-subsub-US} 
Recall the parabolic subgroup of the isometry group $\Aut(\bB)$ of 
the hyperbolic space $\bB$ for an $(i_{0}+1)$-dimensional Klein model 
$\bB \subset \SI^{i_{0}+1}$
fixing a point $p$ in the boundary of $\bB$. \index{b@$\bB$}
Such a discrete subgroup of a parabolic subgroup  
is isomorphic to a finite extension of a lattice in $\bR^{i_0}$ 
by the Bieberbach theorem. 
\index{parabolic subgroup}

Let $E$ be an $i_0$-dimensional ellipsoid 
containing the point $\mbv$ in a subspace $P$ of dimension $i_{0}+1$ in 
$\SI^{n}$.  Let $\Aut(P)$ denote the group of projective automorphisms of $P$, 
and let $\SL_{\pm}(n+1, \bR)_{P}$ the subgroup of $\SL_{\pm}(n+1, \bR)$ acting on $P$. 
Let $r_{P}:\SL_\pm(n+1, \bR)_{P} \ra \Aut(P)$ denote the restriction homomorphism $g \ra g| P$. 
An {\em $i_{0}$-dimensional partial cusp subgroup} is 
one mapping  under $r_{P}$ isomorphically
to a cusp subgroup of $\Aut(P)$ acting on $E -\{\mbv\}$, fixing $\mbv$.
\index{cusp group!$i_{0}$-dimensional partial|textbf} 

Suppose now that $\torb \subset \RP^{n}$. 
Let $P'$ denote a subspace of dimension $i_{0}+1$ containing an $i_{0}$-dimensional ellipsoid $E'$ 
containing $\mbv$. 
Let $\PGL(n+1,\bR)_{P'}$ denote the subgroup of $\PGL(n+1, \bR)$ acting on $P'$. 
Let $R_{P'}: \PGL(n+1, \bR)_{P'} \ra \Aut(P')$ denote the restriction $g\mapsto g|P'$.
An {\em $i_{0}$-dimensional partial cusp subgroup} is one mapping  under $R_{P'}$ isomorphically
to a cusp subgroup of $\Aut(P')$ acting on $E'-\{\mbv\}$, fixing $\mbv$. 
When $i_{0} = n-1$, we drop the ``partial'' from the term ``partial cusp group''.

We generalize a bit further. 
 An {\em $i_0$-dimensional cusp group} is a finite extension of a projective conjugate  
of a discrete cocompact subgroup of a group of an $i_{0}$-dimensional partial cusp subgroup acting cocompactly on $E - \{\mbv\}$. 








Let $\SI^{i_{0}+1}_{l}$ denote the $(i_{0}+1)$-dimensional great sphere containing $\SI^{i_{0}}_{\infty}$ 
corresponding to each $i_{0}$-dimensional leaf $l$ of $\widetilde{\mathcal{F}}_{\tilde E}$. 


\begin{proposition}\label{np-prop-ZN} 
Let $l$ be a generic fiber such that $A_K$ acts with trivial stabilizers. 
\begin{enumerate} 
\item[(i)] $S_l$ acts on $\tilde V_l$ cocompactly, acts on $\partial U$ properly,
and acts as isometries on these spaces with respect to some Riemannian metrics. 
\item[(ii)] A closed subgroup $C_{l, 0}$ of unipotent elements
acts transitively on each leaf $l$ with trivial stabilizers, 
and $C_{l, 0}$ acts on an $i_0$-dimensional ellipsoid $\partial U \cap \SI^{i_{0}+1}_{l}$ 
passing $\mbv_{\tilde E}$ with an invariant Euclidean metric. Here, 
we may need to modify $U$ further.
\item[(iii)] $S_{l, 0}$ 
normalizes an $i_0$-dimensional partial cusp group $C_{l, 0}$
where $S_{l, 0}\cap C_{l, 0}$ are cocompact subgroups in
both $S_{l, 0}$ and $C_{l, 0}$. 
\item[(iv)] $C_{l, 0}$  is virtually normalized by $\bGamma_{\tilde E}$
and also by $S_l$. Also, $C_{l, 0}$ acts freely and properly on $m$ for 
each leaf $m$ of $\widetilde{\mathcal{F}}_{\tilde E}$. 

\item[(v)] With setting $\CN := C_{l, 0}$,
Hypothesis \ref{np-h-norm} holds virtually by $\bGamma_{\tilde E}$ for a coordinate system. 
\end{enumerate}
\end{proposition}  
\begin{proof} 

(i) 
By Lemma \ref{ce-lem-transitive},
$S_l$ acts properly on $\tilde V_l$. 
Since $\partial U$ is in one-to-one correspondence with $\tilde \Sigma_{\tilde E}$, 
$S_l$ acts on $\partial U$ properly. Hence, these spaces have 
compact stabilizers with respect to $S_l$. 
The existence of an invariant metric follows from
an argument similar to that in the proof of Lemma \ref{np-lem-invmet}. 
Hence, the action is proper, and the orbit is closed. 

Since $\bGamma_l$ acts on $\tilde V_l$ cocompactly, 
so is the action of $S_l$ on $\tilde V_l$. 


(ii) We may assume that $\bGamma_{\tilde E}$ is torsion-free by
Theorem \ref{prelim-thm-vgood}
taking a finite-index subgroup if necessary. 

Proposition \ref{np-prop-eigSI} implies that  for $g\in \bGamma_{l}$ 
\[\lambda^{Tr}_{\max}(g) \geq \lambda^{\Sio}_{\max}(g) \geq  \lambda^{\Sio}_{\min}(g) \geq \lambda^{Tr}_{\min}(g).\]
Since $S_{l} = F \Gamma_{l}$ for a compact set $F$, 
the inequality 
\begin{gather}
C_{1}\lambda^{Tr}_{\max}(g) \geq \lambda^{\Sio}_{\max}(g) \geq  C_{2}\lambda^{\Sio}_{\min}(g) \geq C_{3}\lambda^{Tr}_{\min}(g), g \in S_{l}, \nonumber \\ 
C_{1}\lambda^{Tr}_{\max}(g) \geq \lambda_1(g) \geq 
\lambda_{n+1}(g) \geq C_2 \lambda^{Tr}_{\min}(g), g \in S_l, 
\label{np-eqn-eigen2} 
\end{gather} 
hold for constants $C_{1}>1, 1> C_{2} > C_{3}> 0$.
Since $S_{l, 0}$ acts trivially on $K^{o}$, we have 
$\lambda^{Tr}_{\max}(g) = \lambda^{Tr}_{\min}(g)$ for $g \in S_{l, 0}$.  
Since the maximal norm $\lambda_1(g)$ of the eigenvalues of $g$ equals 
$\lambda^{Tr}_{\min}(g)$
and the minimal norm of the eigenvalues of $g$ equals $\lambda^{Tr}_{\min}(g)$, 
all the norms of the eigenvalues of $g \in S_{l, 0}$ are bounded above.  
\eqref{np-eqn-eigen2} implies that 
$|\log\lambda^{\Sio}_{\max}(g)|, |\log\lambda_1(g)|, g\in S_{l, 0}$ are both uniformly bounded above.  
Of course, we have
\[|\log\lambda^{\Sio}_{\max}(g^{n})|= |n\log\lambda^{\Sio}_{\max}(g)|,  |\log \lambda_1(g^{n})| = |n\log \lambda_1(g)|, g \in S_{l, 0}.\]
We conclude that
 the norms of eigenvalues of $g \in S_{l, 0}$ are all $1$.

Theorem \ref{prelim-thm-orthopotent} implies that 
$S_{l, 0}$ is a closed orthopotent group and hence a solvable Lie group. 
Lemma \ref{ce-lem-orthnil} gives a unipotent group $C_{l, 0}$ acting 
on $l$ where $C_{l, 0}$ is the Zariski closure of 
the unipotent subgroup $S_{l, 0}^u$ of $S_{l, 0}$. 
We have $C_{l, 0} \cap S_{l, 0} = S_{l, 0}^u$. 
Proposition \ref{ce-prop-orthouni2} shows that 
$C_{l, 0}$ is a cusp group. 
Since $S_l$ normalizes $S_{l, 0}$ and $C_{l, 0} \cap S_{l, 0} = S_{l, 0}^u$ is cocompact 
in $S_{l, 0}$, it follows that $S_l$ normalizes $C_{l, 0}$. 
This also proves (iii).

(iv) Proposition \ref{ce-prop-orthouni2} shows that the action of $C_{l, 0}$ on
any leaf $m$ is a free and proper action. 
Since $C_{m, 0}$ acts on $m$, 
$B_{m}:= H_{m} \cap U$ is again bounded by an ellipsoid. 
Since $B_{m}$ has a hyperbolic metric as a Klein model,
and $C_{l, 0}$ is unipotent acting properly on
horospheres of $B_{m}$ for $\mbv_{\tilde E}$,  
$C_{l, 0}$ must also be a cusp group on $B_{m}$.  
Hence, $C_{l, 0}$ acts as a cusp group on each 
$H_{m} \cap U$. 
%
%
%

Let $g\in \bGamma_{\tilde E}$. 
By using these arguments for $g(l)$ instead of $l$, 
$g C_{l, 0} g^{-1}$ also acts on an ellipsoid $E_{m}$ in the subspace 
corresponding to $m$ from $\mbv_{\tilde E}$ 
as a unipotent Lie group freely and transitively.
Since $E_l$ bounds an $(i_0+1)$-dimensional ball with a hyperbolic metric of the Klein model, such a unipotent group is unique
and hence it follows that $gC_{l, 0} g^{-1}$ and $C_{l, 0}$
restrict to the same group in $H_{m}$.

Let $\hat C$ denote the group generated by $C_{l, 0}$ and its conjugates. 
$\hat C$ is obviously unipotent. 
Also, $\hat C$ acts properly on $\tilde \Sigma_{\tilde E}$ 
since $\bGamma_{\tilde E}$ and 
$C_{l, 0}$ preserve a Riemannian metric. 

Let $g' \in C_{l, 0}$ and $g'' \in gC_{l, 0}g^{-1}$
such that $g'|H_{m} = g''|H_{m}$. 
Then $g^{\prime -1}g''$ fixes every point in $m$. 
Since and the stabilizer of 
the unipotent group acting properly on $\tilde \Sigma_{\tilde E}$ is 
trivial, $g' = g''$. 
Hence, the normality follows.

(v) The first two properties of Hypothesis \ref{np-h-norm} follow from Propositions \ref{prelim-prop-Ben2} and \ref{prelim-prop-sweep}.  
$\bGamma_{\tilde E}$ satisfies the transverse weak middle-eigenvalue condition 
by the premise.  
Since $S_{l, 0}^u$ goes to $\Idd$ under $\Pi_K^\ast$, 
it is in the standard form where $l$ corresponds to a great sphere 
$\SI^{i_0+1}$ containing $\SI^{i_0}_\infty$. 

Since $N$ acts on $\tilde V_l$, $N$ is a subgroup of $\Gamma_l$ virtually.
Hence, $N$ is a subgroup of $S_l$ and therefore of $S_{l, 0}$ virtually. 
Theorem \ref{prelim-thm-orthopotent} tells us that $N$ is virtually unipotent and therefore $N \cap \bGamma_{\tilde E}'$ is in $S_{l, 0}^u$
for a finite-index subgroup $\bGamma_{\tilde E}'$ of $\bGamma_{\tilde E}$.  
(iv) showed that $\CN$ is normalized by $\bGamma_{\tilde E}$.  
(iv) also shows that $\CN$ acts freely and properly on each complete
affine leaf of $\tilde \Sigma_{\tilde E}$. 
\end{proof}



\subsection{The proof for nondiscrete $N_{K}$.} \label{np-sub-nondproof} 

Now, we go to the splitting argument for this case. 
We can parametrize $C_{l, 0}$ by $\CN(\vec{v})$ for $\vec{v} \in \bR^{i_0}$ by Proposition \ref{np-prop-ZN}. 
We showed that Hypothesis \ref{np-h-norm} holds virtually
by Proposition \ref{np-prop-ZN}.  
For convenience, let us assume Hypothesis \ref{np-h-norm} here. 

We outline the proof strategy. 
\begin{itemize} 
	\item $N_K$ is semisimple in Proposition \ref{np-prop-NKsemi}. 
	\item $\mu_g = 1$ for every $g \in \bGamma_{\tilde E}$. 
	\item Hypothesis \ref{np-h-qjoin} holds. Now 
	we use the results in Section \ref{np-subsub-qjoin}. 
	\end{itemize}

Also, $N \cap C_{l, 0}$ has a finite index in $N$ 
since both acts on $l$ and we took many finite-index subgroups in
the processes above. 
Again, using Proposition \ref{prlim-prop-Malcev}, we can consider
the finite covers of the p-end neighborhoods. 
We assume $N \subset C_{l, 0}$. 
Hypothesis \ref{np-h-norm} holds now, as we showed in the previous subsections. 
As in the preceding by Lemmas \ref{np-lem-similarity}, 
we find that the matrices are of form: 
\renewcommand{\arraystretch}{1.5}
\begin{equation} \label{np-eqn-secondm3}
\newcommand*{\tempV}{\multicolumn{1}{r|}{}}
\CN(\vec{v}) = \left( \begin{array}{ccccccc} 
\Idd_{n-i_0-1} & \tempV & 0 & \tempV & 0 & \tempV & 0 \\ 
\cline{1-7}
0 &\tempV & 1 &\tempV & 0 &\tempV & 0 \\ 
\cline{1-7}
C_1(\vec{v}) &\tempV & \vec{v}^T &\tempV & \Idd_{i_0} &\tempV & 0 \\ 
\cline{1-7}
c_2({\vec{v}}) &\tempV & \llrrV{\vec{v}}^2 /2&\tempV & \vec{v} &\tempV & 1 
\end{array} 
\right), 
\end{equation} 
\begin{equation} \label{np-eqn-stdg} 
\newcommand*{\tempV}{\multicolumn{1}{r|}{}}
g = \left( \begin{array}{ccccccc} 
S(g) & \tempV & s_1(g) & \tempV & 0 & \tempV & 0 \\ 
\cline{1-7}
s_2(g) &\tempV & a_1(g) &\tempV & 0 &\tempV & 0 \\ 
\cline{1-7}
C_1(g) &\tempV & a_1(g) \vec{v}^T_g &\tempV & a_5(g) O_5(g) &\tempV & 0 \\ 
\cline{1-7}
c_2(g) &\tempV & a_7(g) &\tempV & a_5(g) \vec{v}^T_g O_5(g) &\tempV & a_9(g) 
\end{array} 
\right)
\end{equation}
where $g \in \bGamma_{\tilde E}$. 
(See \eqref{np-eqn-formgi}.)
Recall 
$\mu_g = a_5(g)/a_1(g) = a_9(g)/a_5(g)$. 
Since $S_l$ is in ${\mathcal{Z}}(\bGamma_l)$ and the orthogonality of normalized $A_5(g)$
is an algebraic condition, the above form also holds for $g \in S_l$.

\begin{proposition}[Semisimple $N_K$] \label{np-prop-NKsemi} 
	Assume hypothesis \ref{np-h-norm}. 
	Suppose that $\pi_1(\tilde E)$ satisfies {\rm (} NS {\rm)} {\rm (} especially if $\dim K = 0, 1$ {\rm)}.
	Then the following hold\/{\rm :}  
	\begin{itemize} 
\item $N_K$ or any of its finite-index subgroup acts semi-simply on $K^o$.
\item There is a finite-index subgroup $N'_K$ of $N_K$ acting on each
$K_i$ irreducibly and has the diagonalizable
commutant $H$ isomorphic to $\bR_+^{\bar l-1}$ for some $\bar l \geq 1$. 
\item $K$ is projectively diffeomorphic to $K_1 \ast \cdots \ast K_{\bar l}$, 
where $H$ acts trivially on each $K_j$ for $j=1, \dots, {\bar l}$. 
\item Let $A_K$ denote the identity component of the closure $\bar N_K$ of $N_K$ in $\Aut(K)$.
Let $A'_K$ denote the image in $A_K$ of $\bGamma_{\tilde E}$. Then 
$A'_K\cap N'_K$ for a subgroup $N'_K$ of $N_K$ 
is free abelian and is a diagonalizable group of 
matrices and is in the virtual center $N_K \cap H$. 
\item $N'_K$ acts on each $K_i$ strongly irreducibly, and 
$N'_K | K_i $ is semisimple and discrete and acts on $K_i^o$ as a dividing action. 
\item $N_K$ contains a free abelian group in $N_K \cap H$ of rank $\bar l -1$.  
\end{itemize} 
\end{proposition}
\begin{proof} 
		It is sufficient to prove for $N_K$ itself since $N'_K$ for a finite-index 
		subgroup $N'_K$ of $N_K$ acts cocompactly. 
	
		If $N_K$ is discrete, then the conclusion follows from Proposition \ref{prelim-prop-Ben2}.

	Suppose $N_K$ is not discrete.  
	For the case where $\dim K=0, 1$, the conclusions are obvious. 
	We prove using induction on $\dim K$.

We now use the notation of Section \ref{np-subsec-leaves}. 
By Theorem \ref{prelim-thm-vgood}, we assume that $\bGamma_{\tilde E}$ is torsion-free.  
	By condition (NS), since $\bGamma_l$ is virtually normal in $\bGamma_{\tilde E}$, 
	$\bGamma_l\cap G$ is central in a finite-index subgroup $G$ of 
	$\bGamma_{\tilde E}$ and is free-abelian. 
	
	Now, we prove this by induction on $\dim K$. 
	We recall the exact sequence. 
	\[1\ra N \ra \bGamma_l \ra A'_K \ra 1. \]
	Here, $A'_K$ is dense in a nilpotent Lie group $A_K$
normalized by $N_K$. 	Let $N'_K$ denote the image of $G$ in $N_K$. 
	Let $\bar N_K'$ denote the closure of $N_K'$ in $\Aut(K)$. 
	
	We give a bit  vague outline of the rest of the proof: 
	\begin{description} 
		\item[(i)] First, we show that there are 
no unipotent actions of some groups on $K^o$
		by length arguments. 
		\item[(ii)] We decompose the space into invariant subspaces on which $N_K$ 
		virtually acts discretely.
		\item[(iii)] Now, we prove that $N_K$ acts semisimply. 
		\item[(iv)] Finally, we show that $N_K$ contains a free abelian 
		group of certain rank. 
\end{description}

(i) We take the unipotent subgroup $\bGamma_{l, u}$,
	$A'_{K, u}$, $A_{K, u}$ of solvable groups $\bGamma_l$, $A'_K$, and $A_K$ 
	respectively. 
	These are normalized by $N_K$. 
	
	
	
	Suppose that $A'_{K, u}\cap N'_K$ is nontrivial. 	
	Choose a nontrivial unipotent element $g_u$ in $A'_{K, u}\cap N'_K$. 
	By Lemma \ref{prelim-lem-unitnolower}, there exists 
	a sequence of elements $x_i\in K^o$ 
	such that $d_{K}(x_i, g_u(x_i)) \ra 0$. 
	Since $N'_K$ is still a sweeping action, 
	let $F$ be a compact set in $K^o$ such that $\bigcup_{g\in N'_K} g(F) = K^o$. 
	Now, we can choose $g_i\in G$ such that $g_i(x_i)\in F$. 
	Then 
	\begin{equation} \label{np-eqn-dOmega}
	\{d_{\Omega}(g_i g_u(x_i), g_i(x_i)) = 
	d_{\Omega}(g_i g_u g_i^{-1}(g_i(x_i)), g_i(x_i))\} \ra 0
	\end{equation} 
	since 
	$g_i$ is a $d_{K}$-isometry. 
	This implies that $\{g_i g_u g_i^{-1}\}$ converges to an element of a stabilizer of 
	a point $f$, $f\in F$ in $\bar N_K$ up to a choice of a subsequence. 
	

Since $\bGamma_l\cap G$ is central in $G$, 
	$g_i g_u g_i^{-1} = g_u$. 
	Since $g_u$ is unipotent, if $g_u$ stabilizes a point of $K^o$, then 
	$g_u$ is the identity element by Lemma \ref{prelim-lem-elliptic}.
	This is a contradiction. 
	Therefore, we conclude that $A'_{K, u}\cap N'_K$ is a trivial group. 

	

(ii)
Since $\bGamma_l$ is central in $G$, 
$A'_K\cap N'_K$ and its closure $A_K\cap \bar N_K'$ 
	are abelian groups. 
	Since $A_K\cap \bar N'_K$ is abelian, 
	we can decompose 
	$\bC^{n-i_0} = V'_1 \oplus \cdots \oplus V'_{l'}$ 
	such that each element $g$ of $A_K \cap \bar N'_K$ acts irreducibly 
	with a single eigenvalue $\lambda_i(g)$ on $V'_i$
	for each $i$ and its conjugate $\bar \lambda_i(g)$
	by the primary decomposition theorem. 
	(See Theorem 12 of Section 6 of \cite{HK71} and Definition \ref{prelim-defn-jordan}.) 
	The map 
	\[g\in A_K \cap \bar N_K' \mapsto 
	(\lambda_1(g), \dots, \lambda_{l'}(g)) \in \bC^{\ast l'}\]
	gives us an isomorphism to the image set
	where we choose a representative eigenvalue $\lambda_i(g)$ for each $V_i$. 
	
	We can define $S_1, \dots, S_{l'}$ 
	to be the set of independent subspaces in $\SI^n$ corresponding to
	a real primary subspace for every $g$
	by the commutativity of elements in $A_K\cap \bar N_K'$. 

	Since $N'_K$ commutes with $A_K$, $N'_K$ also acts on or permutes the 
	corresponding subspaces $S_1, \dots, S_{l'}$ in $\SI^{n-i_0-1}$. 
	We take the finite-index subgroup $N''_K$ of $N'_K$ acting on each $S_i$, 
	and $N''_K$ has a sweeping action on $K^o$. 
		By Proposition \ref{prelim-prop-sweep}, 
	$S_i \cap K \ne \emp$ and $S_i \cap K^o = \emp$. 
	Denote by $K_i:= S_i \cap K$. 
	
	Suppose that $\lambda_i(g)$ is not real for some $i$
	and $g \in A_K\cap \bar N_K'$. 
	Then there is an eigenspace for $\lambda_i(g)$ and one for
	$\bar \lambda_i(g)$ for all $g \in N_K$. 
	Since $K\cap S'$ for a corresponding real subspace $S'$ for the direct sum of 
	two cyclic spaces
	is properly convex, there is a global-fixed point 
	of $g$. The point corresponds to a positive real eigenvalue of $g$. 
	This is a contradiction. The negative eigenvalues of $g$ violate the 
	proper convexity. 
	Therefore, every $\lambda_i(g)> 0$ for $g \in A_K\cap \bar N'_K$.

Let $\bar N''_K$ denote the closure of $N''_K$ in $\bar N_K$. 
	Suppose that $A_K\cap \bar N_K''$ acts on $K_i$ nontrivially. 
	$A_K\cap \bar N_K''|K_i$ can be considered 
as a unipotent action since each element of $A_K\cap \bar N_K'|K_i$ has a single positive eigenvalue affilated with $K_i$. 
	Since $K_i^o/A'_K\cap N''_K$ is compact, we can apply the arguments in 
    the paragraph containing \eqref{np-eqn-dOmega} and the following one, 
    and we obtain a contradiction. 
    Hence, $A_K\cap \bar N_K''|K_i$ is trivial.
	Hence, $A_K\cap \bar N'_K$ is a positive diagonalizable group. 
	

    
    Since $A_K$ is the identity component of $\bar N_K$, 
    and $A_K$ restricts to a trivial group for each $K_i$, 
    $N_{K_i}:= N''_K| K_i$ is discrete. 

(iii)    We are in the case of $N_{K_i}$ being discrete and
$N_{K_i}$ is semisimple 
by Proposition \ref{prelim-prop-Ben2}. 

	Since $N_{K_i}$ is discrete, 
	and $N'_K$ is isomorphic to a subgroup of 
	$N_{K_1}\times \cdots \times N_{K_{l'}} \times \bR_+^{l' -1}$, 
	it follows that the finite extension
	$N_K$ is semisimple. 	
This proves the first to the fifth item of this proposition. 

(iv) Now, we prove the last item in particular. 
	Proposition \ref{prelim-prop-sweep} also shows the existence of 
	a free diagonalizable subgroup $\Lambda$ of $\bar N_K \cap H$ of rank $\bar l-1$. 
	Hence, there must be a lattice in $L \subset \Lambda$ 
	that is Zariski dense in $\bar N_K \cap H$. 
	Choose generators $\eta_1, \dots, \eta_{\bar l}$ of the lattice. 
	For each $\eta_j$, there is a sequence $\{\kappa^j_i\}$ in $N_K$ converging to 
	$\eta_j$ in $\Aut(K)$. Each $\kappa^j_i|K_i$ is in a discrete group $N_{K_i}$. 
	Hence, we may assume that $\kappa^j_i|K_i = \Idd_{K_i}$ for every $i$ since 
	$\eta_j|K_i = \Idd_{K_i}$.  Hence, $\kappa^j_i \in H \cap N_K$ for every $i$. 
	Since $\kappa^j_i$ are sufficiently close to $\eta_j$ for each $j=1, \dots, \bar l$, 
	we can choose a set of generators 
	$\kappa^1_{i'}, \dots, \kappa^{\bar l}_{i'} $ of $H\cap N_K$
for sufficiently large $i'$. 
	This completes the proof.
%
%
%
	%
	%
	\end{proof}

We continue to have Hypothesis \ref{np-h-norm} for $\bGamma_{\tilde E}$
by Proposition \ref{np-prop-ZN}.  
However, we have not shown Hypothesis \ref{np-h-qjoin} yet. 

\begin{proposition} \label{np-prop-mug}
	Suppose that $\orb$ is properly convex. 
We assume that Hypothesis \ref{np-h-norm} and $N_{K}$ is nondiscrete.  
Suppose that $\pi_1(\tilde E)$ satisfies {\rm (}NS\/{\rm )}
{\rm (}in particular if $\dim K = 0, 1${\rm )}.
Then
we have $\mu_g= 1$ for every $g \in \bGamma_{\tilde E}$. 
\end{proposition}
\begin{proof} 
We can take finite-index subgroups for $\bGamma_{\tilde E}$ during the proof
and prove for this group since $\mu$ is a homomorphism to 
the multiplicative group $\bR_+$. 
By Proposition \ref{np-prop-NKsemi}, $N_K$ is semisimple 
and Proposition \ref{np-lem-conedecomp1} and 
Lemma \ref{np-lem-matrix} hold.

Propositions \ref{np-prop-NKsemi} 
and  \ref{prelim-prop-sweep} show that 
$K = K_1 \ast \cdots \ast K_l$ for properly convex sets $K_i$, and 
$N_K$ is virtually isomorphic to a 
subgroup of $N_{K_1}\times \cdots \times N_{K_l}\times \Lambda$ where 
$\Lambda$ is Zariski dense in a diagonalizable group 
$\bR_+^{l-1}$ acting trivially on each $K_i$ and 
$N_{K_i}$ acts semisimply on each $K_i$ for $i=1, \dots, l$. 
We take a finite-index subgroup $N'_K$ such that $N'_K$ acts on
$K_i$ for each $i=1, \dots, l$. We assume that $N_K$ is this $N'_K$ 
without loss of generality. 


We apply Lemma \ref{np-lem-conedecomp1}. 
Then one of $K_i$ is a vertex $k'$.  
Now we can use the coordinates of \eqref{np-eqn-formgi} repeated here. 
\begin{equation} 
\newcommand*{\tempV}{\multicolumn{1}{r|}{}}
\left( \begin{array}{ccccccc} 
S(g) & \tempV & 0 & \tempV & 0 & \tempV & 0 \\ 
 \cline{1-7}
0 &\tempV & a_{1}(g) &\tempV & 0 &\tempV & 0 \\ 
 \cline{1-7}
C_1(g) &\tempV & a_{1}(g) \vec{v}^T_g &\tempV & a_{5}(g) O_5(g) &\tempV & 0 \\ 
 \cline{1-7}
c_2(g) &\tempV & a_7(g) &\tempV &  a_{5}(g) \vec{v}_g O_5(g) &\tempV & a_{9}(g)
\end{array} 
\right)
\end{equation}
defining $\vec{v}_g:= \frac{a_4(g)}{a_1(g)}$.  

Also, Lemma \ref{np-lem-conedecomp1} shows that $C_1(\vec{v}) =0$ for all $\vec{v} \in \bR^{i_0}$ for a coordinate system where $k$ has 
the form \[\llrrparen{0,\dots, 0, 1} \in \SI^{n-i_0-1}.\]
By Proposition \ref{np-prop-decomposition},
we have a coordinate system where 
\begin{gather}
C_1(g)=O, c_2(g)=0
\hbox{ for every } g \in \bGamma_{\tilde E} \hbox{ and } \nonumber \\
C_1(\vec{v}) = O, c_2(\vec{v}) =0
\hbox{ for every }
\mathcal{N}(\vec{v}), \vec{v}\in \bR^{i_0}. 
\label{np-eqn-C1c1}
\end{gather}

Let $\lambda_{S_g}$ denote the maximal norm of the eigenvalues of 
the upper-left block $S_g$ of $g$. 
We define 
\[\bGamma_{\tilde E, +}:= \{g| \lambda_{S_g}(g) < a_1(g) \}.\]
There is always an element like this. 
In particular, we take the inverse image of 
suitable diagonalizable elements of the center  $H \cap N_K$  denoted in Proposition \ref{np-prop-NKsemi}. 
We take the diagonalizable element in $N_K$ with $k$ having a largest norm
eigenvalue. Let $g$ be such an element. Then the transverse weak middle-eigenvalue condition 
shows that $a_1(g)$ is the largest of norms of every eigenvalue
by Proposition \ref{np-prop-eigSI}, and \eqref{np-eqn-a9a5a1} shows   
\[a_{1}(g) \geq a_{9}(g) \hbox{ and equivalently } 
\mu_{g}  \leq 1 \hbox{ for } g \in \bGamma_{\tilde E, +}.\]





By Proposition \ref{np-prop-NKsemi}, 
$N_K \cap H$ contains a free abelian group of rank $\dim H$
which is positive diagonalizable. 
Hence, there exists 
${g_c} \in \bGamma_{\tilde E, +}$ going to 
a center of $N'_K$ with $\mu_{{g_c}} \leq 1$.  


(A) We obtain a nontrivial element of $N$: 
Let us choose $k \in  \mathds{A}^n$ for an p-end neighborhood $U$ of
$\tilde E$. This is a point going to the vertex of $K$ under the extension of 
the projection $\Pi_K$ to $\mathds{A}^n$ which exists. 

We choose a coordinate system such that 
\begin{align} \label{np-eqn-kcoor} 
k& =\llrrparen{0, \dots, 0, x_{n-i_0}, 0, \dots, 0},\,  x_{n-i_0}> 0, \nonumber \\
\SI^{i_0+1}_k & =\{ \llrrparen{0, \dots, 0, x_{n-i_0}, \dots, x_{n+1}} |\, x_i \in \bR  \}. 
\end{align} 
Here, $k$ may be regarded as the origin of $H_k^o$. 

Now, let $g_1$ be any element of $\bGamma_{\tilde E}$. 
We factorize the lower-right $(i_0 +2)\times (i_0 +2)$-submatrix of $g_1$, 
$g_1 \in \bGamma_{\tilde E}$,  
\begin{equation} \label{np-eqn-eta}
\left(
\begin{array}{c|c|c}
a_{1}(g_1) & 0 & 0 \\ 
\hline
a_{1}(g_1) \vec{v}^T_{g_1} & a_{5}({g_1}) O_5({g_1}) & 0 \\ 
\hline
a_7({g_1}) &  a_{5}({g_1}) \vec{v}_{g_1} O_5({g_1}) & a_{9}({g_1})
\end{array} 
\right)  = 
\end{equation} 
\begin{equation}
\left(
\begin{array}{c|c|c}
1 & 0 & 0 \\ 
\hline
0 & \Idd & 0 \\ 
\hline
\frac{a_7({g_1})}{a_{1}({g_1})} -  \frac{\llrrV{\vec{v}_{{g_1}}}^{2}}{2} &  0 & 1
\end{array} 
\right) 
\left(
\begin{array}{c|c|c}
1 & 0 & 0 \\ 
\hline
\vec{v}^T_{g_1} & \Idd & 0 \\ 
\hline
\frac{\llrrV{\vec{v}_{{g_1}}}^{2}}{2} &  \vec{v}_{g_1} & 1
\end{array} 
\right) 
a_{1}({g_1})
\left(
\begin{array}{c|c|c}
1 & 0 & 0 \\ 
\hline
0 & \mu_{{g_1}} O_5({g_1}) & 0 \\ 
\hline
0 &  0 & \mu_{{g_1}}^{2}
\end{array} 
\right).
\end{equation}
using \eqref{np-eqn-a9a5a1}.   
%


Suppose that $\vec{v}_g =0$ for every $\bGamma_{\tilde E}$. 
Then $\bGamma_{\tilde E}$ fixes $k$. 
By Proposition \ref{np-prop-NKsemi} and Lemma \ref{np-lem-conedecomp1},
$K = K'' \ast \{k'\}$ for a compact convex set $K''$ and a point $k'$. 
There is a set $K'''$ in $\Bd \torb$ corresponding to $K''$
by Proposition \ref{np-prop-decomposition}. 
Then the interior $K''' \ast k$ maps to $K^o$ under the extension of the projection
$\Pi_K$ as we can see from the matrix forms. 
By \eqref{np-eqn-C1c1}, $K''' \ast k \ast \mbv_{\tilde E}$ is $\bGamma_{\tilde E}$-invariant. 
Also, under the radial projection to 
$R_{\mbv_{\tilde E}}(\orb) = \tilde \Sigma_{\tilde E}$, 
the interior of $K''' \ast k \ast \mbv_{\tilde E}$ is mapped to the
$\bGamma_{\tilde E}$-invariant subspace
in $\tilde \Sigma_{\tilde E}$ that meets each complete affine leaf
at a point. 
The existence of an invariant manifold of lower dimension
contradicts the cocompactness of 
the action on $\tilde \Sigma_{\tilde E}$. 

Let us take nonidentity $g \in \bGamma_{\tilde E}$ going to $N'_K$. 
with nonzero $\vec{v}_g$. 
Then conjugation $g_c g g_c^{-1}$ gives 
us an element with $\vec{v}_{g_c g g_c^{-1}} = \vec{v}_g \mu_{g_c} O_5(g_c)^{-1}$
by Lemma \ref{np-lem-similarity}. 
This is not equal to $\vec{v}_g$ since $\mu_{g_c}< 1$. 
Hence, a block matrix computation shows that 
$g_c g g_c^{-1} g^{-1}$ is not an identity element in 
$\bGamma_{\tilde E}$ but maps to $\Idd$ in $N_K$. 
We obtain a nontrivial element $n_0$ of $N$. 
By Hypothesis \ref{np-h-norm},
$g_c g g_c^{-1} g^{-1} \in \mathcal{N}$.
Since $n_0:= g_c g g_c^{-1} g^{-1} \ne \Idd$ has the form \eqref{np-eqn-secondm3}, $n_0$ is a unipotent element. 
Since $N \subset \CN$, we may let $n_0= \CN(\vec{v}_0)$ for some 
nonzero vector $v_0$.

(B) Now we show $\mu_g = 1$ for all $g \in \bGamma_{\tilde E}$: 



Suppose that we have an element $g \in \bGamma_{\tilde E, +}$
and $\mu_g < 1$. 
Then we have as above
 \[\vec{v}_{g^k n_0 g^{-k}} =  \vec{v}_{n_0} \mu_{g}^k O_5(g)^n.\] 
Also, $g^k n_0 g^{-k}$ goes to $\Idd$ in $N_K$ since $\Pi_K^\ast(n_0)= \Idd$ in $N_K$. 
Hence, $\{g^k n_0 g^{-k}\} \ra \Idd$ as $k \ra \infty$ 
since $n_0$ is in the forms \eqref{np-eqn-secondm3} 
given by \eqref{np-eqn-C1c1}.
This contradicts the discreteness of $N$. 
Hence, $\mu_g = 1$ for all $g \in \bGamma_{\tilde E,+}$.
 
Since any element of $g \in \bGamma_{\tilde E}$, 
we can take $g'$, $g'\in \bGamma_{\tilde E, +}$ 
such that $g g' \in \bGamma_{\tilde E, +}$ and so
$\mu_{gg'}= \mu_g \mu_{g'}= 1$. Since $\mu_g \leq 1$, 
we obtain $\mu_g = 1$ for all $g\in \bGamma_{\tilde E}$.
\end{proof} 


\begin{proof}[{The proof of Theorem \ref{np-thm-thirdmain}}]
Suppose that $\tilde E$ is an NPNC R-end. 
When $N_{K}$ is discrete, Theorem  \ref{np-thm-NPNCcase}  gives us the result.

When $N_{K}$ is nondiscrete, 
Hypothesis \ref{np-h-norm}  holds by Propositions \ref{np-prop-ZN}. 
Also, $N_K$ is semisimple by Proposition \ref{np-prop-NKsemi}. 

By Proposition \ref{np-prop-mug}, $\mu \equiv 1$ holds.
Lemmas \ref{np-lem-similarity} and \ref{np-lem-conedecomp1},
(iv) and (v) of Proposition \ref{np-prop-ZN} show that the premise of 
Proposition \ref{np-prop-decomposition} holds. 
Proposition \ref{np-prop-decomposition} shows that
Hypothesis \ref{np-h-qjoin} holds. 
Proposition \ref{np-prop-qjoin}
shows that we have a strictly joined or quasi-joined end. 
Corollary \ref{np-cor-NPNCcase2} implies the result. 

Note here that we may prove for finite-index subgroups of 
$\bGamma_{\tilde E}$
by the definition of strictly joined or quasi-joined ends.
\end{proof}



We give a convenient summary. 

\begin{corollary} \label{np-cor-NPNChol} 
Let $\orb$ be a properly convex real projective $n$-orbifold.
Assume that its holonomy group is strongly irreducible. 
Let $\tilde E$ be an NPNC p-end of the universal cover $\torb$ or $\orb$
satisfying the transverse weak middle-eigenvalue condition
for the R-p-end structure of $\tilde E$. 
Suppose that $\pi_1(\tilde E)$ satisfies {\rm (}NS\/{\rm )} {\rm (} in particular if $\dim K = 0, 1$\/{\rm )}.
Then the holonomy group 
$h(\bGamma_{\tilde E})$ is a group whose element under a coordinate system is of form{\,\rm :} 
\begin{equation} \label{np-eqn-finalform} 
\newcommand*{\tempV}{\multicolumn{1}{r|}{}}
g = \left( \begin{array}{ccccccc} 
S(g) & \tempV & 0 & \tempV & 0 & \tempV & 0 \\ 
 \cline{1-7}
0 &\tempV & \lambda(g) &\tempV & 0 &\tempV & 0 \\ 
 \cline{1-7}
0 &\tempV & \lambda(g) \vec{v}(g)^{T}&\tempV & \lambda(g) O_5(g) &\tempV & 0\\ 
 \cline{1-7}
0 &\tempV & \lambda(g)\left(\aleph_7(g)+ \frac{\llrrV{\vec{v}(g)}^{2}}{2}\right) 
&\tempV & \lambda(g) \vec{v}(g) O_5(g) &\tempV & \lambda(g) 
\end{array} 
\right)
\end{equation}
where $\{S(g)|g\in \bGamma_{\tilde E}\}$ acts cocompactly on a properly convex domain in $\Bd \torb$
of dimension $n-i_{0}-1$,  $O_5:\bGamma_{\tilde E } \ra \Ort(i_0)$ is a homomorphism, 
and $\aleph_7(g)$ satisfies the uniform positive translation condition given by \eqref{np-eqn-uptc}. 

And $\bGamma_{\tilde E}$ virtually normalizes the group 
\renewcommand{\arraystretch}{1.5}
\begin{equation} \label{np-eqn-finalCN} 
\newcommand*{\tempV}{\multicolumn{1}{r|}{}}
\Bigg \{\CN(\vec{v}) 
= \left( \begin{array}{ccccccc} 
\Idd_{n-i_0-1} & \tempV & 0 & \tempV & 0 & \tempV & 0 \\ 
 \cline{1-7}
0 &\tempV & 1 &\tempV & 0 &\tempV & 0 \\ 
 \cline{1-7}
0 &\tempV & \vec{v}^T &\tempV & \Idd_{i_0} &\tempV & 0 \\ 
 \cline{1-7}
0 &\tempV & \llrrV{\vec{v}}^2 /2&\tempV & \vec{v} &\tempV & 1 
\end{array} 
\right)\Bigg | \vec{v} \in \bR^{i_{0}} \Bigg\}. 
\end{equation} 
\end{corollary} 
\begin{proof} 
The proof is contained in the proof of Theorem \ref{np-thm-thirdmain}.
\end{proof} 


\index{NPNC-end!dimension of the fiber}


\section{Applications of the NPNC-end theory} \label{np-sec-dualNPNC}

\subsection{The proof of Corollary \ref{np-cor-dualNPNC} } \label{np-sub-dualNPNC}



\begin{proof}[Proof of Corollary \ref{np-cor-dualNPNC}]
We may always take finite-index subgroups for $\bGamma_{\tilde E}$ and 
consider as the end holonomy group. 
By Corollary \ref{np-cor-NPNChol}, we obtain that the dual holonomy group $g^{-1 T} \in \bGamma_{\tilde E}^{\ast}$ 
has form under a coordinate system:  
\begin{equation} \label{np-eqn-dualg}
\newcommand*{\tempV}{\multicolumn{1}{r|}{}}
\Scale[0.85]{g^{-1 T} = \left( \begin{array}{ccccccc} 
S(g)^{-1 T} & \tempV & 0 & \tempV & 0 & \tempV & 0 \\ 
 \cline{1-7}
0 &\tempV & \lambda(g)^{-1} &\tempV & -\lambda(g)^{-1}  O_5(g)^{-1} \vec{v}_g &\tempV & \lambda(g)^{-1}\left(-\aleph_7(g)+ \frac{\llrrV{v(g)}^{2}}{2}\right) \\ 
 \cline{1-7}
0 &\tempV & 0 &\tempV & \lambda(g)^{-1} O_5(g)^{-1} &\tempV & -\lambda(g)^{-1} \vec{v}_g^{T}\\ 
 \cline{1-7}
0 &\tempV & 0 &\tempV & 0 &\tempV & \lambda(g)^{-1} 
\end{array} 
\right).}
\end{equation}

And $\bGamma_{\tilde E}^{\ast}$ virtually normalizes the group 
$\{ \CN(\vec{v})^{-1T}| \vec{v} \in \bR^{i_0} \}$ where 
\renewcommand{\arraystretch}{1.5}
\begin{equation} \label{np-eqn-seconddualm}
\newcommand*{\tempV}{\multicolumn{1}{r|}{}}
\CN(\vec{v})^{-1T} = \left( \begin{array}{ccccccc} 
\Idd_{n-i_0-1} & \tempV &\,\, 0 & \tempV & 0 & \tempV & 0 \\ 
\cline{1-7}
0 &\tempV & \,\, 1 &\tempV &  -\vec{v} &\tempV & \llrrV{\vec{v}}^2 /2 \\ 
\cline{1-7}
0 &\tempV & \, \,0 &\tempV & \Idd_{i_0} &\tempV & -\vec{v}^{T} \\ 
\cline{1-7}
0 &\tempV & \,\, 0 &\tempV & 0 &\tempV & 1 
\end{array} 
\right). 
\end{equation} 
By using coordinate changes by reversing the order of the 
$n-i_{0}+1$-th coordinate to the $n+1$-th coordinate, 
we can make the lower right matrix 
of $\bGamma_{\tilde E}$ and $\CN$
into a lower triangular form. Hence, $\CN^\ast$ is a partial $i_0$-dimensional cusp group.

Suppose that $\langle S(g), g \in \bGamma_{\tilde E}\rangle $  acts on a properly convex set
$K:= K'' \ast \{k\}$ in $\SI^{n-i_{0}-1}$, a strict join,  for a properly convex set 
$K''\ \subset \SI^{n-i_{0}-2} \subset \SI^{n-i_{0}-1}$ 
and $k$ from the proof of 
Proposition \ref{np-prop-qjoin}. 
$\mathcal{N}$ acts on $\SI^{i_{0}+1}$ containing $\SI^{i_{0}}_{\infty}$ and corresponding to a point $k$ 
under the projection $\Pi_{K}:\SI^{n} - \SI^{i_{0}}_{\infty} \ra \SI^{n-i_{0}-1}$. 
Let $K'''$ denote the compact convex set in $\SI^n - \SI^{i_0}_\infty$ mapping 
homeomorphic to $K''$ under $\Pi_{K}$ as we showed in Proposition
\ref{np-prop-decomposition}. 
There is also a subspace $\SI^{n-i_0-2}_{K'''}$ that is the span of $K'''$. 
In addition, $K_4:= K''' \ast \mbv_{\tilde E}$ is in $\SI^{n-i_0-1'}$ the great sphere containing
$\mbv_{\tilde E}$ and projecting to $\SI^{n-i_0-2}$ under the extension of $\Pi_K$. 

Recall Proposition \ref{prelim-prop-NPduality} for the following: 
We have $\bR^{n+1} = V \oplus W$ for subspaces $V$ and $W$ corresponding to $\SI^{n-i_{0}-2}_{K'''}$ and $\SI^{i_{0}+1}$
respectively. 
We may assume that $\mathcal{N}$ acts on both spaces 
and $\bGamma_{\tilde E}$ acts on $K'''$ and both spaces. 
Then $\bR^{n+1\ast} = V^{\dagger}\oplus W^{\dagger}$ for subspaces $V^{\dagger}$ of $1$-forms zero on $W$ and 
$W^{\dagger}$ of $1$-forms zero on $V$. Then $V^{\dagger}$ corresponds to the intrinsic dual subspace $\SI^{n-i_{0}-2\dagger}_{K'''}$ which equals $\SI^{i_0+1\ast}$,  
and $W^{\dagger}$ corresponds to the intrinsic dual $\SI^{i_{0}+1\dagger}$ which equals 
$\SI^{n-i_0-2\ast}_{K'''}$. 
We let $\SI^{n-i_{0}-2\dagger}_{K'''}$ and $\SI^{i_{0}+1\dagger}$ denote the dual subspaces in $\SI^{n\ast}$. 

Hence, $\bGamma_{\tilde E}^\ast$ and $\mathcal{N}^\ast$ act on 
both of these two spaces. 

Since $K''$ is $\Bd \tilde \Sigma_{\tilde E}$ as in Proposition \ref{np-prop-decomposition}, 
we obtain $K_4 \subset \Bd \torb$. 

Let us choose a properly convex p-end neighborhood $U$ where $\CN$ acts on.
$U\cap Q$ for any $i_0+1$-dimensional subspace $Q$ containing $\SI^{i_0}_\infty$ 
is either empty or is an ellipsoid since $\CN$ acts on $U$. 
Any sharply supporting hyperspace $P'$ at $\mbv_{\tilde E}$ to $U$ or $\torb$ must 
contain $\SI^{i_0}_\infty$ since $P' \cap Q$ for any $(i_0+1)$-dimensional 
subspace $Q$ containing $\SI^{i_0}_\infty$ must be disjoint from 
$U \cap Q$ and hence $P'\cap Q \subset \SI^{i_0}_\infty$ and hence 
$P' \supset \SI^{i_0}_\infty$. 

Let $P \subset \SI^{n}$  be an oriented hyperspace sharply supporting $\torb$ at $\mbv_{\tilde E}$
and containing $\SI_{\infty}^{i_0}$ and $\SI^{n-i_0-2}_{K'''}$.
This is unique such one since the hyperspace is the join of the two. 
Hence, $K''' \ast \{\mbv_{\tilde E}\} \subset P$
where  
\[\clo(U)\cap P = \clo(\torb) \cap P = K_4:=  K''' \ast \{\mbv_{\tilde E}\} \] by Proposition \ref{np-prop-qjoin}.

We have the subspaces $\SI^{n-i_0-1'}$ and $\SI^{i_0+1}$ meeting 
at $\mbv_{\tilde E}$ and $\SI^{i_0}$ containing $\mbv_{\tilde E}$ 
meeting with $\SI^{n-i_0-1'}$ at the same point. 


Let $P^\star$ denote the dual space of $P$ under its intrinsic duality. 
Let us denote by $K_{4, P}^\star \subset P^\star$ the dual of $K_4$ with respect to $P$.

Consider a pencil $P_{t, Q}$ of hyperspaces supporting $\torb$ starting from 
$P$ sharing a codimension-two subspace $Q \subset P$. 
Then $P_{t, Q}$ exists for $t$ in a convex interval $I_Q$ in a projective circle. 
$Q$ supports $P \cap \clo(\torb) = {K'''}\ast \{\mbv_{\tilde E}\}$. 
We assume $P_{0, Q}= P$ for all $Q$. 
Thus, the space $S_{P, \Omega}$ of all hyperspaces in one of $P_{t, Q}$ 
is projectively equal to a fibration over $K_4^\star$ in $P^\star$
with fibers a singleton or a compact interval. 

By Proposition \ref{prelim-prop-NPduality}(iii), 
$K_4^\star$ equals $K_4^\dagger \ast \SI^{i_0-1}_P$ 
for a proper subspace dual $K_4^\dagger$ in 
$\SI^{n-i_0-1 \dagger}$ to $K_4$ 
and a great sphere $\SI^{i_0-1}_P$ in $P^\star$. 

Consider $S_{P, \Omega}$ as a subset of $\SI^{n\ast}$ now. 
Then $P_{t, Q}$ under duality goes to a ray in $\torb^\ast$ from $P^\ast$ to 
a boundary point of $\torb^\ast$. 
Hence, the space $R_P(S_{P, \Omega})$ of such open rays is projectively diffeomorphic to 
the interior of $K_4^\dagger \ast \SI^{i_0-1}_P$. 
Let $K^{'''\dagger}$ in $\SI^{n-i_0-2\dagger}$ 
denote the dual of $K'''$ in 
its span $\SI^{n-i_0-2}_{K'''}$. 
Since $K_4^\dagger$ is projectively diffeomorphic to 
$K^{'''\dagger} \ast \{v\}$ for a singleton $v$, 
we have $S_{P, \Omega}$ is projectively diffeomorphic to the interior of
$K^{'''\dagger} \ast H^{i_0}$ for a hemisphere $H^{i_0}$ of dimension $i_0$. 

By the duality argument, this space equals $R_{P^\ast}(\torb^\ast)$ since 
such rays correspond to supporting pencils of $\orb$ and vice versa.

Now, recall that our matrices of $\bGamma_{\tilde E}$ 
in the form of \eqref{np-eqn-finalform}
and matrices in $\CN$ in the form \eqref{np-eqn-finalCN}. 
We can directly show the properness of the action on 
$R_{P^\ast}(\torb^\ast)$: 
 Our dual action is basically that on 
$K^{\dagger, o} \times \bR^{i_0}$, an affine form of the interior of 
$K^{\dagger} \ast \SI^{i_0-1}$, preserving the product structure.
In each factor $K^{\dagger}$, the action is given by
the upper left $(n-i_0)\times (n-i_0)$-matrices 
and on $\bR^{i_0}$ by the lower-right $(i_0+2)\times (i_0+2)$-submatrix of 
\eqref{np-eqn-dualg}.
Let $g_i$ be a sequence of mutually distinct elements of 
$\bGamma_{\tilde E}^\ast$ such that 
$g_i(L)\cap L \ne \emp$ for a compact subset $L$ in
$R_{P^\ast}(\torb^\ast)$. 
Then it must be that $S(g_i)^\ast$ is bounded for our matrices in $\bGamma_{\tilde E}$
since we can project to $K^{\dagger}$ and argue using the invariant Hilbert metric. 
Then $\vec{v}_{g_i}$ blows up: Otherwise, the properness of $\bGamma_{\tilde E}$ does not hold for $\tilde \Sigma_{\tilde E}$. 
This implies that $g_i(L)$ cannot stay in a bounded set, a contradiction. 



Hence, it follows from our matrix form \eqref{np-eqn-dualg} that 
$\bGamma_{\tilde E}^\ast$ acts properly on 
$R_{P^\ast}(\torb^\ast)$ projectively isomorphic to the interior of 
$K_4^\dagger \ast \SI^{i_0-1}$.
By Lemma \ref{intro-lem-actionRend}, $\torb^\ast$ can be considered 
a p-end neighborhood with a radial structure. 
Hence, we can apply our theory of the classification of NPNC-ends. 
The transverse weak uniform middle-eigenvalue condition is satisfied by the form of the matrices. 
Also, the uniform positive translation condition holds by the matrix forms again. 
Proposition \ref{np-prop-qjoin} completes the proof.

Finally, we mention the following: 
Since each $i_0$-dimensional ellipsoid in the fiber in 
$\tilde \Sigma_{\tilde E^{\ast}}$ has a unique fixed point 
that should be common for all ellipsoid fibers, 
the choice of the p-end vertex is uniquely 
determined for $\tilde E^{\ast}$ such that 
$\tilde E^\ast$ is to be quasi-joined.
\end{proof}


\subsection{ A counterexample in a solvable case} \label{np-sub-counter} 

We find some nonsplit NPNC-end where the end holonomy group is solvable.
Our construction is related to the construction of Carri\`ere \cite{Carriere84} and Epstein \cite{Epstein84}. 
A related work is given by Cooper \cite{Cooper17}. 
These are not a quasi-join nor a join and 
do not satisfy our conditions (NS) and neither the end orbifolds
admit properly convex structure nor the fundamental groups are virtually abelian. 
Define
\renewcommand{\arraystretch}{1.5}
\begin{equation} \label{np-eqn-secondm3co}
\newcommand*{\tempV}{\multicolumn{1}{r|}{}}
N(w,v) = \left( \begin{array}{ccccccccc} 
1 & \tempV & w & \tempV & w^2/2 & \tempV & 0 & \tempV & 0 \\ 
\cline{1-9}
0 & \tempV & 1 & \tempV & w & \tempV & 0 & \tempV & 0 \\ 
\cline{1-9}
0 & \tempV & 0 &\tempV & 1 &\tempV & 0 &\tempV & 0 \\ 
\cline{1-9}
0 & \tempV & 0 &\tempV & v &\tempV & 1 &\tempV & 0 \\ 
\cline{1-9}
0 & \tempV &  0 &\tempV & v^2 /2&\tempV & v &\tempV & 1 
\end{array} 
\right), \hbox{ and }
\end{equation} 
\begin{equation} \label{np-eqn-stdgco} 
\newcommand*{\tempV}{\multicolumn{1}{r|}{}}
g_\lambda = \left( \begin{array}{ccccccccc} 
\lambda^2 & \tempV & 0 & \tempV & 0 & \tempV & 0 & \tempV & 0 \\ 
\cline{1-9}
0 & \tempV & \lambda & \tempV & 0 & \tempV & 0 & \tempV & 0 \\ 
\cline{1-9}
0 & \tempV & 0 &\tempV & 1 &\tempV & 0 &\tempV & 0 \\ 
\cline{1-9}
0 & \tempV & 0 &\tempV & 0 &\tempV & 1/\lambda  &\tempV & 0 \\ 
\cline{1-9}
0 & \tempV & 0 &\tempV & 0 &\tempV & 0 &\tempV & 1/\lambda^2 
\end{array} 
\right) 
\hbox{ for } v,  w, \lambda \in \bR, \lambda > 0. 
\end{equation}
We compute 
$g_\lambda N(w, v) g_\lambda^{-1} = N(\lambda w, v/\lambda)$. 
We define the $3$-dimensional group 
\[S:= \langle N(w, v), g_\lambda: v, w \in \bR, \lambda \in \bR_+ \rangle.\]
Consider the affine space $\mathds{A}^4\subset \SI^4$ given by $x_3 > 0$
with coordinates $x_1, x_2, x_4, x_5$ where $S$ acts on. 
Then $\langle N(w, 0), w \in \bR \rangle$, 
acts on an open disk $B_{1,2}^{(C_1)}$ bounded by a quadric $x_1 \geq x_2^2/2+ C_1$ in 
the plane $x_4 =0, x_5 =0$ for each constant $C_1 \in \bR$. (See Section \ref{np-subsub-quadric}.)
$\langle N(0, v), v \in \bR \rangle$, 
acts on an open disk $B_{4, 5}^{(C_2)}$ bounded by a quadric $x_5 \geq x_4^2/2+ C_2$ 
in the plane $x_2=0, x_3 =0$ for each constant $C_2 \in \bR$. 
There is a product of these two types of quadrics in $\mathds{A}^n$. 
Then $g_\lambda$ send to the products of these two types of 
quadrics to such a product. 

Hence, an orbit $S(\llrrparen{1, 0, 1, 0, 1})$ is given by 
the following set as a subset of $\mathds{A}^4$: 
\[ \left\{(x_1, x_2, x_4, x_5)\large| x_1 = \frac{x_2^2}{2} + C^2, x_5= \frac{x_4^2}{2} + \frac{1}{C^2}, \hbox{ for some } C> 0 \right\}.  \]
This is a $3$-cell.
Moreover, 
\begin{multline}
N(w, v)g_{\lambda}\left(\frac{x_2^2}{2} + C^2,x_2,x_4, \frac{x_4^2}{2} + \frac{1}{C^2}\right) 
 \\ 
=\left(\lambda^2 \left(\frac{x_2^2}{2} + C^2\right) + \lambda w x_2 + \frac{w^2}{2}, \lambda x_2 + w, \frac{x_4}{\lambda} + v , \frac{x_4^2}{2\lambda^2} + \frac{1}{C^2\lambda^2} + v \frac{x_4}{\lambda} + \frac{v^2}{2}\right) \\ 
=\left(\frac{(\lambda x_2 + w)^2}{2} +  \lambda^2 C^2, \lambda x_2 + w, \frac{x_4}{\lambda} + v , \frac{\left(\frac{x_4}{\lambda} + v \right)^2}{2} + \frac{1}{C^2\lambda^2}\right) \\
= g_{\lambda} N(w/\lambda, \lambda v)\left(\frac{x_2^2}{2} + C^2,x_2,x_4, \frac{x_4^2}{2} + \frac{1}{C^2}\right).
\end{multline}
Hence, there is an exact sequence 
\[ 1 \ra \{N(w, v)|w, v \in \bR\} \ra S \ra \{g_\lambda|\lambda> 0\} \ra 1,  \]
telling us that $S$ is a solvable Lie group (Thurston's Sol \cite{Thurston97}.).

We find a discrete solvable subgroup. 
We take a lattice $L$ in $\bR^2$ and obtain a free abelian group
$N(L)$ of rank two. 
We can choose $L$ such that the diagonal matrix with diagonal 
$(\lambda_0, 1/\lambda_0)$ for some $\lambda_0$ acts as an automorphism. 
Then the group $S_L$ generated by 
$\langle N(L), g_{\lambda}\rangle$ is a discrete cocompact subgroup of $S$.

We remark that such a group exists by taking a standard lattice in $\bR^2$ and 
choosing an integral Anosov linear map $A$ of determinant $1$
with two eigendirections. 
We choose a new coordinate system such that the eigendirections are parallel to
the $x$-axis and the $y$-axis. Then now $L$ can be read from the new coordinate system,
and $\lambda$ is the eigenvalue of $A$ bigger than $1$. 

The orbit $S(\llrrparen{1,0,1,0,1})$ is a subset of 
$B_{1,2}^{(0)}\times B_{4, 5}^{((0)}$. 

The orbit $S(\llrrparen{1,0,1,0,1})$ is strictly convex:
We work with the affine coordinates. 
We consider this point with affine coordinates $(1,0,0,1)$. 
The tangent hyperspace at this point is given 
by $x_1+ x_5 = 2$. We can show that locally the orbit
meets this hyperspace only at $(1, 0, 0, 1)$ and is otherwise in one side of
the hyperspace. 
Since $Z$ acts transitively, the orbit is strictly convex. 

Also, it is easy to show that the orbit is properly embedded.  
Hence, the orbit is a boundary of a properly convex open domain. 
It is now elementary to show that this is an R-p-end neighborhood
for a choice of p-end vertex 
$\llrrparen{0,0,0, 0, 1}$ or $\llrrparen{1,0,0,0,0}$.

\chapter{Characterization of complete R-ends} \label{ce-sub-horo}



In Chapter \ref{ce-sub-horo}, 
we discuss complete affine ends. 
First, we explain the weak middle-eigenvalue condition. 
We state the main result of the chapter 
Theorem \ref{ce-thm-mainaffine}, which characterizes the 
complete affine ends. We prove it by 
the results in Chapter \ref{ch-np}: 
In Theorem \ref{ce-thm-affinehoro}, 
we show that pre-horospherical ends are horospherical ends. 
In Theorem \ref{ce-thm-comphoro}, we show that 
a complete affine end falls into one of the two classes, 
one of which is a cuspidal and the other one is more 
complicated with two norms of eigenvalues. 
Theorem \ref{ce-thm-mainaffine} shows that the second case
is a quasi-join using the results in Chapter \ref{ch-np}.
In Section \ref{np-sec-misc}, we discuss some miscellaneous results.

The results here overlap with the results of Crampon-Marquis \cite{CM14} and
Cooper-Long-Tillmann \cite{CLT15}. However, the results are of a different direction than theirs 
since they are interested in finite-volume Hilbert metrics, 
and were originally conceived before their papers appeared.  
However, we make use of Crampon-Marquis \cite{CM14}. 

Let $\tilde \orb$ denote a convex domain in $\SI^n$ 
projectively covering a convex real projective $n$-orbifold $\orb$
with a holonomy homomorphism $h: \pi_1(\orb)\ra \SLpm$. 
Let $\tilde E$ be an R-p-end.
A middle-eigenvalue condition
for a p-end holonomy group $h(\pi_{1}(\tilde E))$ with respect to $\mbv_{\tilde E}$ or the R-p-end structure holds if 
for each $g\in h(\pi_{1}(\tilde E)) -\{\Idd\}$, 
the largest norm $\lambda_{1}(g)$ of eigenvalues of $g$ is strictly larger than 
the eigenvalue $\lambda_{\mbv_{\tilde E}}(g)$ associated with the p-end vertex $\mbv_{\tilde E}$. 
\index{middle-eigenvalue condition!weak|textbf} 
\index{lambda@$\lambda_{\mbv_{\tilde E}}(\cdot)$|textbf}


Given an element $g \in h(\pi_{1}(\tilde E))$, 
let $\left(\tilde \lambda_{1}(g), \dots, \tilde \lambda_{n+1}(g)\right)$ be the $(n+1)$-tuple of the eigenvalues
listed with multiplicity given by the characteristic polynomial of $g$
where we repeat each eigenvalue with the multiplicity given by the characteristic polynomial.
The {\em multiplicity} of a norm 
of an eigenvalue of $g$ is 
the number of times the norm occurs among the $(n+1)$-tuples of norms 
\[\left(|\tilde \lambda_{1}(g)|, \dots, |\tilde \lambda_{n+1}(g)|\right).\]

\begin{definition}\label{ce-defn-wmec}
	Let $\tilde E$ be a p-end
	with the holonomy group $h(\pi_{1}(\tilde E))$.  
A {\em weak middle-eigenvalue condition} ({\em wmec}\/)
for an R-p-end $\tilde E$  holds provided for each $g\in h(\pi_{1}(\tilde E))$
the following holds: 
\begin{itemize} 
\item whenever $\lambda_{\mbv_{\tilde E}}(g)$ has the largest norm of all norms of eigenvalues $h(g)$, 
$\lambda_{\mbv_{\tilde E}}(g)$ must have the norm multiplicity $\geq 2$. 
\end{itemize} 
\end{definition}

We note that these definitions depend on the choice of the p-end vertices; 
however, they are well defined once the radial structures are assigned.

\section{Main results}

Our main result classifies CA R-p-ends. We need some facts 
of NPNC-ends that is explained in Section \ref{np-sec-notprop}.


Given a horospherical R-p-end, the p-end holonomy group $\bGamma_{\mbv}$ acts on a p-end neighborhood $U$
and $\bGamma_{\mbv}$ is virtually a subgroup of an $(n-1)$-dimensional cusp group $\mathcal{H}_{\mbv}$. By Lemma \ref{prelim-lem-endnhbd}, 
\[V:=\bigcap_{g\in {\mathcal{H}}_{\mbv }} g(U) \] 
contains a ${\mathcal{H}}_{\mbv}$-invariant p-end neighborhood. 
Hence, $V$ is a horospherical p-end neighborhood of $\tilde E$. 
Thus, a horospherical R-end is pre-horospherical.
(See Definition \ref{intro-defn-R-ends}.)
We clarify by Theorem \ref{ce-thm-affinehoro}:  
\begin{corollary}\label{ce-cor-cusphor} 
Let $\mathcal O$ be a properly convex real projective $n$-orbifold.
Let $E$  be an R-end of its universal cover $\torb$. 
Then $E$ is a pre-horospherical R-end if and only if $\tilde E$ is a horospherical R-end. $\square$
\end{corollary} 


The following classifies the complete affine ends
where we need some results from Chapter \ref{ch-np}.
Since these have virtually abelian holonomy groups
according to Theorem \ref{ce-thm-affinehoro},
they are classified in \cite{BCLp}. 

\begin{theorem}\label{ce-thm-mainaffine} 
Let $\mathcal O$ be a properly convex real projective $n$-orbifold.
Let $\tilde E$  be an R-p-end of its universal cover $\torb \subset \SI^n$ 
{\em (}resp. $\subset \RP^n${\em ).}
Let $\bGamma_{\tilde E}$ denote the end holonomy group. 
Then $\tilde E$ is a complete affine R-p-end if and only if $\tilde E$ is a horospherical R-p-end
or a quasi-joined NPNC-end with fibers of dimension $n-2$ 
with the virtually abelian end-fundamental group
by altering the p-end vertex. 
Furthermore,  if $\tilde E$ is a complete affine R-p-end and 
$\bGamma_{\tilde E}$ satisfies the weak middle-eigenvalue condition
with respect to its R-p-end vertex $\mbv_{\tilde E}$, 
then $\tilde E$ is a horospherical R-p-end. 
\end{theorem}
This is proved at the end of this chapter.






\index{end!pre-horospherical} 
\index{end!horospherical}

\subsection{The Horosphere theorem} 

\begin{theorem}[Horosphere] \label{ce-thm-affinehoro} 
Let $\mathcal O$ be a convex real projective $n$-orbifold, 
$n\geq 2$.  
Let $\tilde E$  be a pre-horospherical R-end of its universal cover $\torb$, $\torb \subset \SI^n$ {\em (}resp. $\subset \RP^n${\em )}, 
and $\bGamma_{\tilde E}$ denote the p-end holonomy group. 
Then the following hold\/{\rm :} 
\begin{enumerate}
\item[{\rm (i)}] The subspace $\tilde \Sigma_{\tilde E} = R_{\mbv_{\tilde E}}(\torb) \subset \SI^{n-1}_{\mbv_{\tilde E}}$
of directions of  
line segments from the endpoint $\mbv_{\tilde E}$ in $\Bd\torb$
into a p-end neighborhood of $\tilde E$ 
forms a complete affine subspace of dimension $n-1$.
\item[{\rm (ii)}]\; The norms of eigenvalues of $g \in \bGamma_{\tilde E}$
are all $1$.
\item[{\rm (iii)}]\;\; $\bGamma_{\tilde E}$ 
virtually is in a conjugate of an abelian parabolic or cusp subgroup of $\SO(n, 1)$
{\rm (}resp. $\PO(n, 1)$\/{\rm )} 
of rank $n-1$ in $\SL_\pm(n+1, \bR)$ {\rm (}resp.  $\PGL(n+1, \bR)${\rm ).}
And hence $\tilde E$ is cusp-shaped.
\item[{\rm (iv)}]\;\; For any compact set $K' \subset \orb$, 
$\orb$ contains a horospherical end neighborhood disjoint from $K'$. 
\end{enumerate}
\end{theorem}
\begin{proof} 
We prove for the case $\torb \subset \SI^n$. 
Let $U$ be a pre-horoball p-end neighborhood with the p-end vertex $\mbv_{\tilde E}$, closed in $\torb$. 
The space of great segments from the p-end vertex passing $U$ forms a convex subset $\tilde \Sigma_{\tilde E}$ 
of a complete affine subspace $\bR^{n-1} \subset \SI^{n-1}_{\mbv_{\tilde E}}$ by Proposition \ref{prelim-prop-projconv}.
The space 
covers an end orbifold $\Sigma_{\tilde E}$ 
with the discrete group $\pi_1(\tilde E)$ acting as a discrete subgroup $\Gamma'_{\tilde E}$ of
the projective automorphisms such that $\tilde \Sigma_{\tilde E}/\Gamma'_{\tilde E}$
 is projectively isomorphic to $\Sigma_{\tilde E}$.

Suppose that 
for at least one $g \in \bGamma_{\tilde E}$,
 the largest norm of eigenvalues and the smallest one of $g$ are both different from $1$. 
By Lemma \ref{prelim-lem-prox}, $g$ is positive bi-semiproximal. 
Therefore, let $\lambda_1(g) >1$ be the greatest norm of the eigenvalues of $g$  and 
$\lambda_2(g)< 1$ be the smallest norm  of the eigenvalues of $g$ as an element of $\SL_\pm(n+1, \bR)$. 
Let $\lambda_{\mbv_{\tilde E}}(g) >0$ 
be the eigenvalue of $g$ associated with $\mbv_{\tilde E}$. 
The possibilities for $g$ are as follows:
\begin{alignat*}{3}
 \lambda_1(g)  && \, = \lambda_{\mbv_{\tilde E}}(g)  && \,> \lambda_2(g), \\ 
  \lambda_1(g)  && \,  > \lambda_{\mbv_{\tilde E}}(g)  &&  \, > \lambda_2(g), \\
\lambda_1(g) && \, > \lambda_2(g)  && \,  = \lambda_{\mbv_{\tilde E}}(g). 
\end{alignat*}
In all cases, at least one of the largest norm or the smallest norm is different from $\lambda_{\mbv_{\tilde E}}(g)$. 
By Lemma \ref{prelim-lem-prox}, this norm is realized by a positive eigenvalue. 
We take $g^{n}(x)$ for a generic point $x \in U$. As $n \ra \infty$ or $n \ra -\infty$, the sequence $\{g^n(x)\}$ limits to a point $x_{\infty}$ 
in $\clo(U)$ distinct from $\mbv_{\tilde E}$.  
In addition, $g$ fixes a point $x_\infty$, and $x_{\infty}$ has a different positive eigenvalue from $\lambda_{\mbv_{\tilde E}}(g)$.
As $x_\infty \not \in U$, it should be $x_\infty =\mbv_{\tilde E}$ by the definition of the pre-horospheres.
This is a contradiction.

Hence, all the norms of eigenvalues of $g, g \in \bGamma_{\tilde E}$ equals $1$. 
The lemma \ref{ce-lem-unithoro} gives us the proof for (i) and (ii). 

Let $U$ denote the domain bounded by the closure of the ellipsoid.
There exist finite elements $g_1, \dots, g_m$ that represent the cosets of
$\bGamma_{\tilde E}/\bGamma_{\tilde E}'$. 
If $g_i(U)$ is a proper subset of $U$, then so is $g_i^n(U)$, and 
hence $g_i^n$ is not in $\bGamma_{\tilde E}'$ for any $n$. 
This is a contradiction. Hence, $\bGamma_{\tilde E}$ also acts on $U$. 
By the same reasoning, $\bGamma_{\tilde E}$ acts on every ellipsoid in
a one-dimensional parameter space containing a unique fixed point, 
and an ellipsoid gives us a horosphere in the interior of a horoball. 
Hence, $\bGamma_{\tilde E}$ is a cusp group. (iii) is proved. 

(iv) We can choose an exiting sequence of p-end horoball neighborhoods $U_i$
where a cusp group acts. We can consider the hyperbolic spaces to understand this. 

\hfill  \SnT {\parfillskip0pt\par}
\end{proof}

\subsection{The forward direction of Theorem \ref{ce-thm-mainaffine}.}








The second case is studied later in Corollary \ref{np-cor-caseiii}.
We show the end $\tilde E$ to be an NPNC-end with the fiber dimension $n-2$ when we choose another point as the new p-end vertex for $\tilde E$.
Clearly, this case is not horospherical. 
(See Crampon-Marquis \cite{CM14} for a similar proof.)

For the following, note that Marquis classified ends for
$2$-manifolds \cite{Marquis12} into cuspidal, hyperbolic, or quasi-hyperbolic ends. Since convex real projective $2$-orbifolds are virtually real convex pojective $2$-manifolds, the following has already been proved for $n=2$
by Marquis. 

\begin{theorem}[Complete affine]\label{ce-thm-comphoro} 
Let $\orb$ be a properly convex real projective $n$-orbifold
for $n \geq 3$. 
Suppose that $\tilde E$ is a complete-affine R-p-end of its universal cover $\torb$ in $\SI^n$ {\rm (} resp. in $\RP^n$\/{\rm ).} Let $\mbv_{\tilde E} \in \SI^n$\/ {\rm (} resp. $\in \RP^n$\/{\rm )} 
be the p-end vertex with the p-end holonomy group 
$\bGamma_{\tilde E}$. Then 
\begin{enumerate} 
\item[{\rm (i)}]  we have following two exclusive alternatives\/{\em :} 
\begin{itemize}
\item $\bGamma_{\tilde E}$ is virtually unipotent where all norms of eigenvalues of elements equal $1$, or 
\item $\bGamma_{\tilde E}$ is virtually abelian where 
\begin{itemize}
\item each $g \in \bGamma_{\tilde E}$ has at most two norms of the eigenvalues, 
\item at least one $g \in \bGamma_{\tilde E}$ has two norms,  and 
\item if $g \in \bGamma_{\tilde E}$ has two distinct norms of the eigenvalues, 
the norm of $\lambda_{\mbv_{\tilde E}}(g)$ has  multiplicity one.
\item $\bGamma_{\tilde E}$ acts as a virtually unipotent group 
on the complete affine space $\tilde \Sigma_{\tilde E}$. 
\end{itemize}
\end{itemize}
\item[{\rm (ii)}] In the first case, $\bGamma_{\tilde E}$ is horospherical, that is, cuspidal. 
\end{enumerate} 
\end{theorem} 

\begin{proof} 
	We first prove for the $\SI^n$-version. 
Using Theorem \ref{prelim-thm-vgood},  we may choose a torsion-free finite-index subgroup. 
We may assume without loss of generality that $\Gamma$ is torsion-free 
since torsion elements have only $1$ as the norms of the eigenvalues, and 
we only need to prove the theorem for a finite-index subgroup. 
Hence, $\Gamma$ does not fix a point in $\tilde \Sigma_{\tilde E}$. 

(i) Since $\tilde E$ is complete affine, $\tilde \Sigma_{\tilde E} \subset \SI^{n-1}_{\mbv_{\tilde E}}$ is identifiable with an affine subspace $\mathds{A}^{n-1}$. 
$\bGamma_{\tilde E}$ induces $\bGamma'_{\tilde E}$ in $\Aff(\mathds{A}^{n-1})$
that are of form $x \mapsto Mx + \vec{b}$ where $M$ is a linear map $\bR^{n-1} \ra \bR^{n-1}$
and $\vec{b}$ is a vector in $\bR^{n-1}$. 
For each $\gamma \in \bGamma_{\tilde E}$, 
\begin{itemize}
\item let $\gamma_{\bR^{n-1}}$ denote this affine transformation, 
\item we denote by $\hat L(\gamma_{\bR^{n-1}})$ the linear part of the affine transformation $\gamma_{\bR^{n-1}}$, and 
\item let $\vec{v}(\gamma_{\bR^{n-1}})$ denote  the translation vector. 

\end{itemize}
A {\em relative eigenvalue} is an eigenvalue of $\hat L(\gamma_{\bR^{n-1}})$.
\index{relative eigenvalue}  

At least one eigenvalue of ${\hat L}(\gamma_{\bR^{n-1}})$ is $1$ 
since $\gamma$ acts without a fixed point on $\bR^{n-1}$.
(See \cite{KS75}.)
Now, ${\hat L}(\gamma_{\bR^{n-1}})$ has a maximal 
invariant vector  subspace $A_\gamma$ of $\bR^{n-1}$ where all norms of the eigenvalues are $1$.


Suppose that $A_\gamma$ is a proper $\gamma$-invariant vector subspace of $\bR^{n-1}$. 
Then $\gamma_{\bR^{n-1}}$ acts on the affine space $\bR^{n-1}/A_\gamma$ 
as an affine transformation with the linear parts without a norm of eigenvalue equal to $1$.
Hence, $\gamma_{\bR^{n-1}}$ has a fixed point in $\bR^{n-1}/A$, and 
$\gamma_{\bR^{n-1}}$ acts on an affine subspace $A'_\gamma$ parallel to $A_\gamma$.

A subspace $H_\gamma$ containing ${\mbv_{\tilde E}}$ corresponds to the direction of $A'_\gamma$ from $\mbv_{\tilde E}$.
The union of open segments with endpoints $\mbv_{\tilde E}, \mbv_{\tilde E-}$ in the directions in $A'_\gamma \subset \SI^{n-1}_{\mbv_{\tilde E}}$
is in an open hemisphere of dimension $< n$. 
Let $H^{+}_\gamma$ denote this space where 
$\partial \clo(H^{+}_\gamma) \ni \mbv_{\tilde E}$ holds.
Since $\bGamma_{\tilde E}$ acts on $A'$, it follows that
$\bGamma_{\tilde E}$ acts on $H^{+}_\gamma$.
Then $\gamma$ has at most two eigenvalues associated with $H^{+}_\gamma$ one of which is $\lambda_{\mbv}(\gamma)$
and the other is to be denoted $\lambda_{+}(\gamma)$.
Since $\gamma$ fixes $\mbv_{\tilde E}$ and there is an eigenvector in the span of 
$H^+_\gamma$ associated with $\lambda_+(\gamma)$,  
$\gamma$ has the matrix form
\[ \gamma = 
\left(
\begin{array}{ccc}
 \lambda_{+}(\gamma) {\hat L}(\gamma_{\bR^{n-1}})  &  \lambda_{+}(\gamma) \vec{v}(\gamma_{\bR^{n-1}}) & 0\\
 0 &  \lambda_{+}(\gamma) &  0\\
 \ast &  \ast &  \lambda_{\mbv_{\tilde E}}(\gamma) 
\end{array}
\right) 
\]
where we have
\begin{equation}\label{ce-eqn-completa} 
\lambda_{+}(\gamma)^{n}\det({\hat L}(\gamma_{\bR^{n-1}})) 
\lambda_{\mbv_{\tilde E}}(\gamma)  = \pm 1.
\end{equation} 
(Note $\lambda_{\mbv_{\tilde E}}(\gamma^{-1}) = \lambda_{\mbv_{\tilde E}}(\gamma)^{-1}$
and $\lambda_{+}(\gamma^{-1}) = \lambda_{+}(\gamma)^{-1}$.)

We show that ${\hat L}(\gamma_{\bR^{n-1}})$ for 
every $\gamma \in \bGamma_{\tilde E}$ is unit-norm-eigenvalued
below. 
As before, $\lambda_1(\gamma)$ denotes the largest norm of the eigevalues of $\gamma$. 
Note that $\lambda_1(\gamma) \geq \lambda_+(\gamma)$ since ${\hat L}(\gamma_{\bR^{n-1}})$ has an eigenvalue equal to $1$. 
There are following possibilities for each $\gamma \in \bGamma_{\tilde E}$:
\begin{itemize}
\item[(a)] $\lambda_{1}(\gamma) > \lambda_{+}(\gamma) \hbox{ and }
\lambda_{1}(\gamma) > \lambda_{\mbv_{\tilde E}}(\gamma)$.
\item[(b)] $\lambda_{1}(\gamma) = \lambda_{+}(\gamma)= \lambda_{\mbv_{\tilde E}}(\gamma)$.
\item[(c)] $\lambda_{1}(\gamma) = \lambda_{+}(\gamma), \lambda_{1}(\gamma) > \lambda_{\mbv_{\tilde E}}(\gamma)$.
\item[(d)]  $\lambda_{1}(\gamma) > \lambda_{+}(\gamma), \lambda_{1}(\gamma) =\lambda_{\mbv_{\tilde E}}(\gamma)$.

\end{itemize}

Suppose that $\gamma$ satisfies (b). The norms of 
the relative eigenvalues of $\gamma$ on $\bR^{n-1}$ are all $\leq 1$. 
Either $\gamma$ is unit-norm-eigenvalued 
or we can take $\gamma^{-1}$ and we are in case (a).  

Suppose that $\gamma$ satisfies (a).
There exists a projective subspace $S$ of dimension $\geq 0$ where 
the points are associated with eigenvalues with the norm $\lambda_{1}(\gamma)$
where $\lambda_{1}(\gamma) >  \lambda_{+}(\gamma), \lambda_{\mbv_{\tilde E}}(\gamma)$. 

Let $S'_\gamma$ be the subspace spanned by $H^+_\gamma$ and $S$. 
Let $U$ be a p-end neighborhood of $\tilde E$. 
Since the space of directions of $U$ is $\bR^{n-1}$, we have 
$U \cap S'_\gamma \ne \emp$. 
We can choose two generic points $y_1$ and $y_2$ of $U \cap S'_\gamma - H^+_\gamma$
such that $\ovl{y_1 y_2}$ meets $H^+_\gamma$ in its interior. 


Then we can choose a subsequence $\{m_i\}$, $\{m_i\} \ra \infty$, such that 
$\{\gamma^{m_i}(y_1)\} \ra f$ and $\{\gamma^{m_i}(y_2)\} \ra f_-$ as $i \ra +\infty$
up to relabeling $y_1$ and $y_2$ for a pair of antipodal points $f, f_- \in S$.
This implies $f, f_- \in \clo(\torb)$, 
and $\torb$ is not properly convex, which is a contradiction. 
Hence, (a) cannot be true. 

We showed that if any $\gamma \in \bGamma_{\tilde E}$ satisfies (a) or (b), then  $\gamma$ is unit-norm-eigenvalued. 

If $\gamma$ satisfies (c), then 
\begin{equation}\label{ce-eqn-c}
\lambda_{1}(\gamma) = \lambda_{+}(\gamma) \geq \lambda_{i}(\gamma) \geq\lambda_{\mbv_{\tilde E}}(\gamma)
\end{equation} 
for all other norms of eigenvalues $\lambda_{i}(\gamma)$: 
Otherwise, $\gamma^{-1}$ satisfies (a), which cannot happen. 

Similarly, if $\gamma$ satisfies (d), then
we have 
\begin{equation}\label{ce-eqn-d}
 \lambda_{1}(\gamma) =\lambda_{\mbv_{\tilde E}}(\gamma)\geq \lambda_{i}(\gamma) \geq \lambda_{+}(\gamma)
 \end{equation}
for all other norms of eigenvalues $\lambda_{i}(\gamma)$. 
We conclude that either $\gamma$ is unit-norm-eigenvalued or satisfies 
\eqref{ce-eqn-c} or \eqref{ce-eqn-d}.





There is a homomorphism 
\[\lambda_{\mbv_{\tilde E}}: \bGamma_{\tilde E} \ra \bR_{+} \hbox{ given by } g \mapsto \lambda_{\mbv_{\tilde E}}(g).\]
This gives us an exact sequence 
\begin{equation} \label{ce-eqn-exxact}
 1 \ra N \ra \bGamma_{\tilde E}  \ra R \ra 1 
\end{equation} 
where $R$ is a finitely generated abelian subgroup of $\bR_{+}$. 
For an element $g \in N$, 
$\lambda_{\mbv_{\tilde E}}(g) = 1$. 
Since the norm of 
the relative eigenvalue corresponding to $\hat L(g_{\bR^{n-1}})|A_g$ is $1$, 
the matrix form shows that $\lambda_{+}(g) = 1$ for $g \in N$. 
\eqref{ce-eqn-c}, \eqref{ce-eqn-d}, 
and the conclusion of the above paragraph 
show that $g$ is unit-norm-eigenvalued. 
%
Thus, $N$ is therefore virtually nilpotent by Theorem \ref{prelim-thm-orthopotent}. 
(See Fried \cite{Fried86}). 
Taking a finite cover again, we may assume that $N$ is nilpotent. 

Since $R$ is a finitely generated abelian group, 
$\bGamma_{\tilde E}$ is solvable by \eqref{ce-eqn-exxact}. 
Since $\tilde \Sigma_{\tilde E}= \bR^{n-1}$ is complete affine, 
Proposition S of Goldman and Hirsch  \cite{GH86} implies 
\[\det(g_{\,\bR^{n-1}}) = 1 \hbox{ for all } g \in \bGamma_{\tilde E}.\]

If $\gamma$ satisfies (c), then all norms of eigenvalues of $\gamma$ 
except for $\lambda_{\mbv_{\tilde E}}(\gamma)$ equal $\lambda_{+}(\gamma)$ 
since otherwise by \eqref{ce-eqn-c}, 
norms of relative eigenvalues 
$\lambda_i(\gamma)/\lambda_+(\gamma)$ are $\leq 1$, and 
the above determinant is less than $1$. 
Similarly, if $\gamma$  satisfies (d),  then similarly all norms of eigenvalues of $\gamma$ except for $\lambda_{\mbv_{\tilde E}}(\gamma)$ equals $\lambda_{+}(\gamma)$.

Therefore, only (b), (c), and (d) hold and $g_{\, \bR^{n-1}}$ is unit-norm-eigenvalued 
for all $g \in \bGamma_{\tilde E}$. 

By Theorem \ref{prelim-thm-orthopotent}, 
$\bGamma_{\tilde E}|\bR^{n-1}$ is an orthopotent group and therefore
is virtually unipotent
by Theorem 3 of Fried \cite{Fried86}. 

Now we return to $\bGamma_{\tilde E}$. 
Suppose that every $\gamma$ is orthopotent. 
Then we have the first case of (i). 
If not, then the second case of (i) holds, and there are only two
norms of eigenvalues for all elements. 
Immediately following 
Collorary \ref{np-cor-caseiii} proves the result.

(ii) This follows from Lemma \ref{ce-lem-unithoro}.
\hfill \SnT {\parfillskip0pt\par}
\end{proof} 



\subsection{Complete affine ends again}  \label{np-sub-compaff}


We now study the second case from the conclusion of Theorem \ref{ce-thm-comphoro}. 


\begin{corollary}[non-cusp complete-affine p-ends] \label{np-cor-caseiii}
Let $\orb$ be a properly convex real projective $n$-orbifold for $n \geq 3$. 
Let $\tilde E$ be a complete affine R-p-end of its universal cover $\torb$ in $\SI^n$. 
Let $\mbv_{\tilde E} \in \SI^n$  be the p-end vertex with the p-end holonomy group 
$\bGamma_{\tilde E}$.  
Suppose that $\tilde E$ is not a cusp p-end. Then we can choose 
a different point as the p-end vertex for $\tilde E$ such that $\tilde E$ is a quasi-joined R-p-end on which a partial cusp group of dimension $n-1$ acts. 
Also, the end fundamental group is virtually abelian. 
\end{corollary} 
\begin{proof} \renewcommand{\qedsymbol}{}  
We use the terminology of the proof of Theorem \ref{ce-thm-comphoro}. 
Theorem \ref{ce-thm-comphoro} 
shows that $\bGamma_{\tilde E}$ is virtually nilpotent and has at most two norms of eigenvalues
for each element. 
By taking a finite-index subgroup, we may assume that $\bGamma_{\tilde E}$ is nilpotent. 
Let $Z$ be the Zariski closure, a nilpotent Lie group. 
We may assume that $Z$ is connected by taking a finite-index subgroup of 
$\bGamma_{\tilde E}$. 
Theorem \ref{ce-thm-comphoro}, says $\bGamma_{\tilde E}$ is 
isomorphic to a virtually unipotent group by restricting to 
the affine space $\tilde \Sigma_{\tilde E}$. Hence, 
$Z$ is simply connected and hence contractible. 
Since $\bGamma_{\tilde E}\cap Z$ is a cocompact lattice in $Z$, 
and $\bGamma_{\tilde E}$ has the virtual cohomological dimension $n-1$, 
it follows that $Z$ is $(n-1)$-dimensional.  

By Lemma \ref{ce-lem-transitive},  
$Z$ acts transitively on the complete affine subspace $\tilde \Sigma_{\tilde E}$ since $\bGamma_{\tilde E}$ acts cocompactly on it. 

The orbit map $Z \ra Z(x)$ for $x \in \tilde \Sigma_{\tilde E}$ is a fiber 
bundle over the contractible space with fiber the stabilizer group. 
Since $\dim Z = n-1$, it follows that it should be discrete. 
Since $\tilde \Sigma_{\tilde E}$ is contractible, the stabilizer is trivial.

Since $Z$ fixes $\mbv_{\tilde E}$, 
we have a homomorphism 
\begin{equation}\label{np-eqn-N}
\lambda_{\mbv_{\tilde E}}: Z \ni g \ra \lambda_{\mbv_{\tilde E}}(g) \in \bR.
\end{equation} 
Let $N$ denote the kernel of the homomorphism, which is a nilpotent Lie group being a subgroup of a nilopotent Lie group. 
Since $Z$ has only two norms of eigenvalues, $N$ has only one norm of eigenvalues.
By Theorem \ref{prelim-thm-orthopotent}, $N$ is orthopotent Lie group
since it has only one norm of the eigenvalues equal to $1$.

Since $\tilde E$ is a complete affine R-p-end, 
$R_{\mbv_{\tilde E}}(U)$ is a complete affine space equal to $\tilde \Sigma_{\tilde E}$. 
Taking a convex hull of finite number of radial rays from 
$\mbv_{\tilde E}$, we may choose 
a properly convex p-end neighborhood $U$ of $\tilde E$.
Also, we may choose such that the closure of $U$ is in another such
p-end neighborhood. Thus, $\Bd U \cap \torb$ is 
in a one-to-one correspondence with $\tilde \Sigma_{\tilde E}$
since a radial segment must escape $U$. 
We modify $U$  to 
\[ \bigcap_{g\in Z} g(U) = \bigcap_{g\in F}g(U).\] 
This contains a nonempty properly convex open set
by Lemma \ref{prelim-lem-endnhbd}.
We may assume that $Z$ acts on a properly convex p-end neighborhood $U$ of $\tilde E$. 
Since $Z$ acts transitively on $\tilde \Sigma_{\tilde E}$, it acts so on 
an embedded convex hypersurface
$\delta U := \Bd U \cap \torb$. This is the set of endpoints of 
maximal segments from $\mbv_{\tilde E}$ in the directions of the complete 
affine space $\tilde \Sigma_{\tilde E}$. 
Since this characterization is independent of $\torb$, 
$\delta U$ is an orbit of $Z$.  
Since $Z$ is a Lie group, $\delta U$ is smooth. 

The smooth convex hypersurface 
$\delta U$ is either strictly convex or has a foliation fibered by totally geodesic submanifolds. Since $\tilde \Sigma_{\tilde E}$ is complete affine, 
these submanifolds must be complete affine subspaces. 
Since $\clo(U)$ contains these and $\clo(U)$ is properly convex,
this is a contradiction. Hence, $\delta U$ is strictly convex. 

Finally, since $\delta U$ is in one-to-one correspondence with 
$\tilde \Sigma_{\tilde E}$, $\delta U/\bGamma_{\tilde E}$ is a compact 
orbifold of codimension $1$. 



Let $A$ be a hyperspace containing $\mbv_{\tilde E}$ in the direction of $\Bd \tilde \Sigma_{\tilde E} = \SI^{n-2}\subset \SI^{n-1}_{\mbv_{\tilde E}}$. 
Then $U_{A}:= A \cap \clo(U)$ is a properly convex compact set on which $Z$ acts.
By Lemma \ref{np-lem-Ucpt}, the $Z$-action on $U_{A}^{o}$ is cocompact. 
By Lemma \ref{np-lem-Ben2}, $U_{A}$ is a properly convex segment. 

We complete the proofs after Lemmas \ref{np-lem-Ben2} and \ref{np-lem-Ucpt}.

\cpr

\end{proof}

\begin{lemma} \label{np-lem-Ucpt} 
The $Z$-action on $U_{A}^{o}$ is cocompact. 
\end{lemma} 
\begin{proof} 
Suppose that $\dim U_{A}^{o} = n-1$. 
The orbit map $Z \ra Z(x)$ for $x\in U_{A}^o$ 
is a fibration over a simply connected domain. 
The stabilizer must be compact since $U_A^o$ has a Hilbert metric.
Since $Z$ being a simply connected nilpotent Lie group is contractible,
the stabilizer has to be trivial.
Since $\dim Z = n-1$,  
$Z$ acts transitively on $U_{A}^{o}$ 
and the $Z$-action on $U_{A}^{o}$ is cocompact by Lemma \ref{ce-lem-transitive2}.


Suppose now that $\dim U_{A}^{o} = j_0 < n-1$. 
Let $L$ be an $(j_0+1)$-dimensional subspace containing $U_{A}$ meeting $A$ transversely. 
Let $l$ be the $j_0$-dimensional affine subspace of the complete affine space $\tilde \Sigma_{\tilde E}$ corresponding to $L$. 
Since \[g(U_{A})=U_{A}, U_{A}\subset L, g(L) \hbox{ and } \dim L = j_0+1, g \in Z\] 
it follows that 
\[ g(L) \cap L = \langle U_{A} \rangle \hbox{ or } g(L) = L, 
\hbox{ which implies }
g(l)= l \hbox{ or } g(l) \cap l = \emp.\] 
Recall that $Z$ acts transitively and freely on the complete affine 
space $\tilde \Sigma_{\tilde E}$ from the beginning of the 
Section \ref{np-sub-compaff}.
Since $\dim l = j_0$, it follows that 
the subgroup $\hat Z_l := \{g \in Z| g(l) = l\}$
has dimension $j_0$. 

Now $\hat Z_l$ acts on $U_{A}^{o}$. As proved above, the stabilizer of 
$\hat Z_l$ of a point of $U_{A}^{o}$ is trivial
since $U_A \cap L $ is properly convex and $\hat Z_l$ is nilpotent 
without a compact subgroup of dimension $> 0$. 
Hence, $\hat Z_l$ acts transitively on $U_{A}^{o}$ 
as in the first paragraph by Lemma \ref{ce-lem-transitive2} 
since $\dim \hat Z_l = \dim U_{A}^{o}$,
and $U_{A}^{o}/\hat Z_l = U_{A}^{o}/Z$ is compact. 
\end{proof}

\begin{lemma}\label{np-lem-Ben2} 
Let a simply connected nilpotent 
Lie group $S$ act cocompactly and effectively
on a properly convex open domain $J$.
Suppose that each element of $S$ has at most two norms of eigenvalues
and fixes a point $p$ in the boundary of $J$. 
Then the dimension of the domain is $0$ or $1$.
\end{lemma}
\begin{proof} 
	Suppose that $\dim J > 1$. 
By Lemma \ref{ce-lem-transitive}, $S$ acts transitively on $J$. 
The action is proper since there is a Hilbert metric on $J$. 
Since $S$ is nilpotent and is simply connected, $S$ is contractible.
Since the stabilizer of $S$ at a point $x\in J$ 
is compact, it is trivial in $S$. 
Hence, $S$ is diffeomorphic to a $(\dim J)$-dimensional cell.

Suppose that $R_p(J)$ is complete affine. Then $N$ acts on 
an affine space consisting of great segments in the directions of $R_p(J)$. 

Suppose that $R_p(J)$ is not complete affine. Then it must be
properly convex or 
is a join of a properly convex domain $K$ and a great sphere of dimension 
$i_0\geq 0$.
In the first case, we take a quotient of $S$ such that it acts effectively. 
Then, by induction, $R_p(J)$ is a segment or a singleton. 
Then $S$ acts on a one- or two-dimensional 
complete affine space containing $J$ since $S$ is connected. 
(That is, it cannot switch the ends of segments).
In the second case, we argue with $K$ and we find that $R_p(J)$ must be a $0$- or $1$-dimensional space
joined with a great sphere
of dimension $i_0$. Again, $S$ acts on a complete affine space
containing the great segments in the direction of $R_p(J)$.

In any case, $S$ acts an affine group on a complete 
affine space containing the segments of directions 
in $R_p(J)$ since $S$ is a connected group. 

Let $\lambda_p(g)$ for $g\in S$ denote the associated 
eigenvalue of $g$ at $p$ for unit-determinant matrix representatives. 
Let $N_S$ denote the kernel of homomorphism
$S \ra \bR_+$ given by $g \mapsto \lambda_p(g)$. 
Then $N_S$ is an orthopotent group of dimension $\dim J -1$
by Theorem \ref{prelim-thm-orthopotent}. 

Then $N_S$ acts on $R_p(J)$ as an orthopotent Lie group. 
We claim that an element $g$ of $N_S$ do not send a point of a line $l$ from $p$ in $J$ to another point on the same line. Otherwise, $g \in N_S$ acts on the line $l'$ that contains $l$. 
Since $g$ has only one norm of the eigenvalue and does not fix a point in $l$, 
we see that $l'$ must be a great segment. Since $l'$ has to be in $J$, this is 
a contradiction. 
Therefore, an orbit of $N_S$ always meets each radial segment from $p$
at most one point. 
Hence, the image of an orbit of $N_S$ to $R_p(J)$ is an open set $O$
by dimension consideration.

 
The stabilizer of a point of 
$R_p(J)$ acts on a segment $s$ from $p$ in $J$. 
The existence of 
two fixed directions of  eigenvalue $1$ implies that 
each point of $s$ is a fixed point, and  
hence the stabilizer is trivial by the above. 
Therefore, the properness and freeness of the action of $N_S$ on 
$R_p(J)$ follow.

We claim that $O = R_p(J)$: 
Since $N_S$ is nilpotent, it has a cocompact discrete subgroup $N'_S$.
Since $O/N'_S$ is a closed affine manifold, $O=R_p(J)$ is a complete affine space by 
Fried \cite{Fried86}. 

By Proposition \ref{ce-prop-orthouni2}, the orbit $N_S(x)$ for $x\in J$ 
is an ellipsoid of dimension $\dim J -1$.  Hence, $S/N_S$ is a $1$-dimensional group. 
The elements of $N_S$ are of form
\begin{equation}\label{np-eqn-form-N_S}
\newcommand*{\tempV}{\multicolumn{1}{r|}{}}
k= 
\left( \begin{array}{ccccc} 
1 &\tempV & 0 &\tempV & 0 \\ 
\cline{1-5}
 \vec{v}_k^T &\tempV & \Idd_{\dim J -1 } &\tempV & 0 \\ 
\cline{1-5}
 \frac{\llrrV{\vec{v}_k}^2}{2} &\tempV & \vec{v}_k &\tempV & 1 
\end{array} 
\right) \hbox{ for } \vec{v}_k \in \bR^{\dim J -1}.
\end{equation} 
Since $N_S$ acts on $J$, this means that $J$ is a standard $(n-1)$-ball $B$
that is tangent to a great sphere $\SI^{n-2}_\infty$of dimension $n-2$ at $p$
up to choice of a coordinate system.

We may write for $g \in S$, 
\begin{equation} \label{np-eqn-form_g}
\newcommand*{\tempV}{\multicolumn{1}{r|}{}}
g= \left( \begin{array}{ccccc} 
a_{1}(g) &\tempV & 0 &\tempV & 0 \\ 
\cline{1-5}
\vec{a}_{4}^T(g)  &\tempV & A_5(g) &\tempV & 0 \\ 
\cline{1-5}
a_7(g) &\tempV &  \vec{a}_{5}(g) &\tempV & a_{9}(g)
\end{array} 
\right).
\end{equation}
This follows since $g$ fixes $p$ and the unique supporting hyperspace at $p$.


Proposition \ref{np-lem-similarity} applies here since we can drop the submatrices $S(g)$, $s_1(g)$, $s_4(g)$ 
by dimension considerations. Also, we do not need the transverse weak middle-eigenvalue 
condition for the arguments to work. 
Lemma \ref{np-lem-conedecomp1} applies since $K$ can be considered as a singleton in 
our case. 
Any element $g \in S$ 
induces an $(\dim J -1)\times(\dim J -1)$-matrix $M_g$ given by
$g \CN(\vec{v}) g^{-1} = \CN(\vec{v}M_g)$ where 
\[M_g = \frac{1}{a_1(g)} (A_5(g))^{-1} = \mu_g O_5(g)^{-1}\]
for $O_5(g)$ in a compact Lie group $G_{\tilde E}$
where $\mu_g=\frac{a_5(g)}{a_1(g)}= \frac{a_9(g)}{a_5(g)}$.

Reasoning as in the proof of Lemma \ref{np-lem-matrix}, 
for every $g \in  N$, we have
\begin{equation} \label{np-eqn-form_g2}
\newcommand*{\tempV}{\multicolumn{1}{r|}{}}
g= \left( \begin{array}{ccccc} 
a_{1}(g) &\tempV & 0 &\tempV & 0 \\ 
\cline{1-5}
a_{1}(g) \vec{v}^T_g &\tempV & a_{5}(g) O_5(g) &\tempV & 0 \\ 
\cline{1-5}
a_7(g) &\tempV &  a_{5}(g) \vec{v}_g O_5(g) &\tempV & a_{9}(g)
\end{array} 
\right), O_5(g)^{-1} = O_5(g)^T,
\end{equation}
and the form of $N_S$ is not changed. 
Also, $a_7(g) = a_1(g)(\alpha_7(g) + \llrrV{v_g}^2/2)$ 
as we can show following the beginning of Section \ref{np-subsub-qjoin}. 
Recall $\alpha_7$ from Section \ref{np-sub-alpha7}. 
Here, $\alpha_7(g) = 0$ since otherwise $J$ cannot be properly convex
since $g$ translates the orbits in the affine space where $p$ 
is the infinity as in Remark \ref{np-rem-alpha7}. 

If there is an element $g$ with $\mu_g \ne 1$, then 
the group $S$ is solvable and not nilpotent as we can see from taking commutators of 
products of elements in the coset $S/N_S$ and $N_S$.
If $\mu_g =1$ for all $g \in S$, then from the matrix form 
we see that $S$ has only $1$ as norms of eigenvalues with $a_1(g)=a_5(g)=a_9(g)$ for 
$g\in S$, and $S$ acts on each ellipsoid orbit of $N_S$. Hence, $J/N = J/N_S$
is not compact. This is a contradiction. 
Therefore, $\dim J = 0, 1$.
%
%
%
\end{proof}

\begin{proof}[Proof of Corollary \ref{np-cor-caseiii} continued] 
If $\dim U_{A} = 0$, then $U$ is a horospherical p-end neighborhood 
where $\bGamma_{\tilde E}$ is unimodular and cuspidal by Theorem \ref{ce-thm-affinehoro}. 

Now assume that $\dim U_{A} = 1$ by Lemma \ref{np-lem-Ben2}. 
Let $q$ denote the other endpoint of the segment $U_{A}$ than $\mbv_{\tilde E}$. 
Since $U$ is convex, $R_q(U)$ is a convex open domain. 
Since an element of $\bGamma_{\tilde E}$ has 
two eigenvalues, 
each radial segment from $q$ maximal in $\torb$ meets the smooth 
strictly convex hypersurface $\delta U$ and transversely since a radial segment cannot be in $\delta U$.  (Recall that $\delta U$ is smooth and strictly convex from the first part of the proof.)


We have $R_q(U) = R_q(\delta U)$:  
Since $\tilde \Sigma_{\mbv_{\tilde E}}$ is complete affine, and $U$ is properly convex, 
each segment from $\mbv_{\tilde E}$ passes $\delta U$ as we lengthen it. 
Hence, $\Bd U = \delta U \cup U_A$ since $U_A$ is precisely 
$\Bd U \cap A$ for the hyperspace $A$ of $\SI^n$ as defined earlier. (Recall that $A$ is the hyperspace containing $\mbv_{\tilde E}$ in the direction of $\Bd \tilde \Sigma_{\tilde E} = \SI^{n-2}\subset \SI^{n-1}_{\mbv_{\tilde E}}$.)
Hence, each segment in 
$U$ from $q$ must end at a point of $\Bd U \cap \torb = \delta U$
since $\delta U$ is strictly convex. 
By the strict convexity of $\delta U$, it is clear that 
each ray in the direction of $R_q(U)$ 
from $q$ meets $\delta U$ transversally. 
By the transversality, a segment from $q$ ending at $\delta U$ must 
have a one-sided neighborhood in $U$ and hence the interior of 
the segment is contained in $U$.

Since the $\bGamma_{\tilde E}$-action on $\delta U$ is proper, 
so is its action on $R_q(U)$. Hence, $U$ can be considered 
as an R-p-end neighborhood with radial lines from $R_q$ foliating 
$U$ by Lemma \ref{intro-lem-actionRend}.  
Lemma \ref{intro-lem-actionRend}
we obtain that $U$ is a R-p-end neighborhood of a p-end vertex $q$, and 
$R_q(U) =R_q(\torb)$ 


There is an embedding from $\delta U$ to $R_q(\delta U) = R_q(U) \subset R_q(\torb)$. 
Since $\delta U/\bGamma_{\tilde E}$ is a compact orbifold, so is 
$R_q(U)/\bGamma_{\tilde E}$. 

By the fourth item in the first item of 
Theorem \ref{ce-thm-comphoro}(i), $R_q(U)$ is not complete-affine: 
The norm $\lambda_{q}(g)$ for some $g \in \bGamma_{\tilde E}$ 
has a multiplicity $n$, $n < n+1$, and $n > 2$ by assumption. 
$\bGamma_{\tilde E}$ cannot act as unimodular group on $R_1(U)$.


Suppose that $R_q(U)$ is properly convex. 
Elements of $\bGamma_{\tilde E}$ have at most two distinct norms of eigenvalues. 
Since $R_q(U)$ is homeomorphic to $\Bd U \cap \torb$ with a compact quotient by 
$\bGamma_{\tilde E}$, 
$R_q(U)$ has a compact quotient by $\bGamma_{\tilde E}$, and 
$\dim R_q(U) = n-1 \geq 2$. 
By Lemma \ref{np-lem-Ben2}, this is a contradiction. 
Hence, $R_q(U)$ is not properly convex. 

Thus, $q$ is the p-end vertex of an NPNC-end for $U$ foliated by 
radial segments from $q$. 
Since the associated upper-left submatrix has only two norms of eigenvalues 
by Lemma \ref{np-lem-Ben2}, and
the properly convex leaf space $K^o$ is $1$-dimensional  
and has a compact quotient, 
the fibers have dimension $n-2=n-1-1$. 
(Also, $K^o$ is a properly convex segment by our definition.) 
Therefore, $R_{q}(U)$ is foliated by $n-2$-dimensional complete affine subspaces. 

%

Recall the Lie group $N$ from \eqref{np-eqn-N}. 
Since $Z$ acts on $R_q(U)$ transitively 
and $\lambda_{\mbv_{\tilde E}}(g) = \lambda_{q}(g), g \in N$, 
the Lie group $N$ acts on each complete affine leaf transitively. 
$N$ is a nilpotent Lie group since it is a subgroup of $Z$. 
Also, $N$ is orthopotent since elements of $N$ 
have only one norm of the eigenvalues by Theorem \ref{prelim-thm-orthopotent}. 
We can apply Proposition \ref{ce-prop-orthouni2}  
to the hyperspace $P$ containing the leaves
with a cocompact subgroup of $N$ acting on it.  
As $U \cap P$ is properly convex, $r_{P}(N)$ is a cusp group. 


Let $x_{1}, \dots, x_{n+1}$ be the coordinates of $\bR^{n+1}$.
Now give coordinates such that $q = \llrrparen{0, 0, \dots, 1}$ and $\mbv_{\tilde E} = \llrrparen{1, 0, \dots, 0}$. 
Since these are fixed points, 
we obtain that elements of $N$ can be put into forms: 
\begin{equation} 
 N(\vec{v}):= \left( \begin{array}{cccc} 
 1 & 0 & 0 & 0 \\ 
 0 &    1 & 0 & 0 \\ 
0 & \vec{v}^T & \Idd & 0 \\ 
0 & \frac{1}{2} \llrrV{\vec{v}}^2 & \vec{v} & 1 
 \end{array} \right) \hbox{ for } \vec{v} \in \bR^{n-2}.
 \end{equation} 


Now, $\bGamma_{\tilde E}$ satisfies the transverse weak middle-eigenvalue condition with respect to $q$
since each element of $\bGamma_{\tilde E}$ has just two 
norms of eigenvalues and $Z$ is generated by $N$ and $g^{t}, t \in \bR,$ for a nonunipotent element $g$ of $\bGamma_{\tilde E}$. This is because 
$g^t, t \ne 0,$ does have multiplicity one eigenvalue at $\mbv_{\tilde E}$.  

Theorem \ref{np-thm-thirdmain} shows that $q$ is a R-p-end vertex for
 a quasi-joined NPNC R-p-end neighborhood $U$ covering the convex 
real projective orbifold $U/\bGamma_{\tilde E}$.   


Finally, $\Gamma_{\tilde E}$ is virtually abelian:
From \eqref{np-eqn-finalform}, 
$S(g)$ is a $1\times 1$-matrix or $0\times 0$-one.
From the matrix form, the Zariski closure $Z$
is an extension of an orthopotent Lie group. 
Since $\tilde \Sigma_{\tilde E}$ equals 
$\mathds{A}^{i_0} \times I$ for an interval or a singleton $I$, 
$i_0= n-2, n-1$, 
$Z$ acts on it. 
$O_5$ extends to a homomorphism
$O_5:Z \ra \Ort(i_0)$. 
Let $Z_K$ denote the kernel. 
Then $Z_K$ also acts properly and cocompactly on 
$\tilde \Sigma_{\tilde E}$
since $Z/Z_K$ is compact. 
It is easy to see $Z_K$ is abelian from the matrix form \eqref{np-eqn-finalform}. 
Also, we can put a $Z$-invariant Euclidean metric on
the complete affine space $\tilde \Sigma_{\tilde E}$ by
the product metric form. 
Then the Bieberbach theorem implies the result. 
\end{proof}

If we require the weak middle-eigenvalue conditions
for a given vertex, the completeness of the end implies 
that the end is cusp. 

\begin{corollary}[cusp and complete affine]\label{np-cor-cusp} 
Let $\orb$ be a properly convex real projective $n$-orbifold. 
Suppose that $\tilde E$ is a complete affine R-p-end of its universal cover $\torb$ in $\SI^n$ {\em (}resp. in $\RP^n${\em ).} 
Let $\mbv_{\tilde E} \in \SI^n$ {\rm (}resp. $\in \RP^n${\rm )} 
be the p-end vertex with the p-end holonomy group 
$\bGamma_{\tilde E}$. Suppose that $\bGamma_{\tilde E}$ satisfies the weak middle-eigenvalue condition with respect to  $\mbv_{\tilde E}$. 
Then $\tilde E$ is a complete affine R-end if and only if $\tilde E$ is a cusp R-end. 
\end{corollary} 
\begin{proof} 
	It is sufficient to prove for the case $\torb \subset \SI^n$. 
Since a horospherical end is a complete affine end 
by Theorem \ref{ce-thm-affinehoro}, 
we need to show the forward direction only.   
In the second possibility of Theorem \ref{ce-thm-comphoro}, the norm of 
$\lambda_{\mbv_{\tilde E}}(\gamma)$ 
has a multiplicity one for a nonunipotent element 
$\gamma$ with $\lambda_{\mbv_{\tilde E}}(\gamma)$ or 
$\lambda_{\mbv_{\tilde E}}(\gamma^{-1})$ equal to the maximal norm.  
This violates the weak middle-eigenvalue condition, and 
 only the first possibility of Theorem \ref{ce-thm-comphoro} holds. 
\hfill  \SnT {\parfillskip0pt\par}
\end{proof}

\begin{proof}[Proof of Theorem \ref{ce-thm-mainaffine}] 
	Again, we assume that $\Omega$ is a domain of $\SI^n$. 
Theorem \ref{ce-thm-comphoro} is the forward direction.
Corollary \ref{np-cor-caseiii} implies 
the second case above. 

Now, we prove the converse using the notation and results of 
Chapter \ref{ch-np}. 
Since a horospherical R-p-end is pre-horospherical, 
Theorem \ref{ce-thm-affinehoro} implies the half of the converse. 
Given an NPNC R-p-end $\tilde E$ with fibers of dimension $n-2$, 
$\tilde \Sigma_{\tilde E}$ is projectively diffeomorphic to an 
affine half-space. Using the notation of Proposition \ref{np-prop-qjoin}, 
$K''$ is zero-dimensional and the end holonomy group
 $\bGamma_{\tilde E}$ acts on $K''\ast \{\mbv\}$ for an end vertex $\mbv$.
 There is a foliation in $\tilde \Sigma_{\tilde E}$ 
 by complete affine spaces of dimension $n-2$ parallel to each other. 
 The space of leaves is projectively diffeomorphic to the interior of $K''\ast \mbv'$ for a point $\mbv'$.
 Let $U$ be the p-end neighborhood for $\tilde E$. 
 Then $\Bd U \cap \torb$ is in one-to-one correspondence with $\tilde \Sigma_{\tilde E}$ 
by radial rays from $\mbv$.
Hence, $\Bd U \cap \torb$ has an induced foliation $\mathcal{F}$ 
with codimension-one leaves.
Each leaf in $\Bd U\cap \torb$ lies in an open hemisphere of dimension $n-1$. 
(See \eqref{np-eqn-pik2} in Section \ref{np-sub-general}.)  
In addition, $\bGamma_{\tilde E}$ acts on an open hemisphere $H^{n-1}_{\mbv'}$ of
dimension $n-1$ with a boundary a great sphere $\SI^{n-2}$ containing $\mbv$ 
in the direction of $\mbv'$. 

 Now we switch the p-end vertex to a singleton $\{k''\}=K''$ from $\mbv$. 
 Then $H^{n-1}_{\mbv'}$ projects to a complete affine space 
 $\mathds{A}^{n-1}_{k''}$ by radial lines from $k''$. 
 Each leaf of $\mathcal{F}$ projects to an ellipsoid in $\mathds{A}^{n-1}_{k''}$ 
 with a common ideal point $\mbv$ corresponding to the direction of 
 $\overline{k''\mbv}$ oriented 
 away from $k''$. The ellipsoids are tangent to a common hyperspace
 in $\SI^{n-1}_{k''}$. Hence, they are parallel paraboloids in 
 an affine subspace $\mathds{A}^{n-1}_{k''}$. 
 The uniform positive translation condition gives 
 us that the union of the parallel paraboloids is $\mathds{A}^{n-1}_{k''}$. 
 Hence, $\tilde E$ is a complete R-end with $k''$ as the vertex. 
The last statement follows from Corollary \ref{np-cor-cusp}. 
\hfill \SnT {\parfillskip0pt\par}
\end{proof}

\section{Some miscellaneous results from the above.} \label{np-sec-misc}

\subsection{Why $\lambda_{\mbv_{\tilde E}}(g) \ne 1$?} 

We need the following to prove Lemma \ref{ex-lem-niceend}.  

\begin{corollary} \label{np-cor-eigenone} 
	Let $\orb$ be a properly convex real projective $n$-orbifold with a strongly
irreducible holonomy group. 
	Suppose that $\tilde E$ is an NPNC R-p-end of its universal cover $\torb$ in $\SI^n$ or 
	{\rm (}resp. in $\RP^n${\rm ).} Let $\mbv_{\tilde E}$  be the p-end vertex with the p-end holonomy group 
	$\bGamma_{\tilde E}$. 
	Suppose that $\pi_1(\tilde E)$ satisfies {\rm (}NS\/{\rm )} for the leaf space $K^o$. 
	Then for some 
	\[ g \in \bGamma_{\tilde E}, 
	\lambda_{\mbv_{\tilde E}}(g) \ne 1.\] 
	\end{corollary} 
\begin{proof}  
	It is sufficient to prove for the case $\torb \subset \SI^n$. 
	Suppose that $\lambda_{\vec{v}_{\tilde E}}(g) = 1$ for all $g \in \bGamma_{\tilde E}$.
		A $\Gamma_{\tilde E}$-invariant $i_{0}$-dimensional subspace $\SI_{\infty}^{i_{0}}$ contains $\mbv_{\tilde E}$ as we discussed in Section \ref{np-sub-general}. 
	
	Suppose that every element of $\bGamma_{\tilde E}$ is unit-norm-eigenvalued. 
	By Theorem \ref{prelim-thm-orthopotent}, $\bGamma_{\tilde E}$ is orthopotent.  
	By Fried \cite{Fried86}, there exists a nonorthopotent element in $\bGamma_{\tilde E}$
	since $\tilde \Sigma_{\tilde E}$ is not complete affine. This is absurd. 
	Hence, there exists an element 
	$g \in \bGamma_{\tilde E}$ that is not unit-norm-eigenvalued. 


	We show that the transverse weak middle-eigenvalue condition for $\tilde E$ holds: 
Suppose not.
We find an element $g$ of 
$\bGamma_{\tilde E}$ not satisfying the condition of the transverse weak 
middle-eigenvalue condition 
with $\lambda_1(g) > 1$. 
	For a real number $\mu$ equal to $\lambda_1(g)$, 
	the subspace $\mathcal{R}_\mu(g)$ projects 
	to $\SI^{i_0}_\infty$ since otherwise the transverse 
	weak middle-eigenvalue condition holds. 
(See \eqref{prelim-eqn-affil}.)  
	There exists a $g$-invariant subspace $\hat P_g \subset \SI^{i_0}_\infty$ 
	that is the projection of 
	$\bigoplus_{\mu = \lambda_1(g)} \mathcal{R}_\mu(g).$
	(See Definition \ref{prelim-defn-jordan}.) 
	This is a proper subspace since $\mbv_{\tilde E}$ has an associated eigenvalue $1$ strictly less than $\lambda_g$. 
	Hence, $\dim \hat P_g \leq i_0 -1$. 
	
	We define $P_{g}$ the projection of 
	\[\bigoplus_{\mu <  \lambda_1(g)} \mathcal{R}_{\mu}(g)\subset \bR^{n+1}.\] 
	Then $P_g$ is complementary to $\hat P_g$. 
	Thus, $\dim P_g \geq n-i_0$, and $P_g \ni \mbv_{\tilde E}$. 

	Under the projection to $\SI^{n-1}(\mbv_{\tilde E})$, 
	$P_{g}$ goes to a subspace $P'_{g}$ of $\dim P_g - 1 \geq n-i_0 -1$. 
	As described in Section \ref{np-sec-notprop},  
	$\tilde S_{\tilde E}$ is foliated by $i_{0}$-dimensional complete affine subspaces. 
	Since $\dim \tilde S_{\tilde E} = n-1$,
	these $i_0$-dimensional leaves must meet $P'_{g}$ of dimension 
	$\geq n-1-i_0$.  
	
	Thus, any p-end neighborhood $U$ of $\tilde E$ meets $P_{g}$. 
	There must be an antipodal pair $\tilde P_g$ in $\hat P_g$ 
	that is the projection of the eigenspace of $g$ with the associated eigenvalue whose norm is $\lambda_1(g) > 1$.
	\[L:=U \cap (P_{g}\ast \tilde P_{g}) \subset P_{g}\ast \tilde P_{g}\]
	is a nonempty open subset of $P_{g}\ast \tilde P_{g}$ meeting $P_{g}$, and 
	$L - P_{g}$ has two components. 
	Let $x, y$ be generic points in distinct components
	in  $L- P_{g}$. 
	Then $\{g^{n}(\{x, y\})\}$ geometrically converges to  
	an antipodal pair of points in $\hat P_{g}$. 
	Since this set is in $\clo(\orb)$, this contradicts the proper convexity of $\orb$. 
	Thus, the transverse weak middle-eigenvalue condition of $\tilde E$ holds. 
	

	The premises of Theorem \ref{np-thm-thirdmain} except for the strong irreducibility of 
	the holonomy group of $\pi_1(\orb)$ 
	are satisfied, and 
	Theorem \ref{np-thm-thirdmain} 
	classifies the ends.
The strong irreducibilty was needed when we used Corollary \ref{np-cor-NPNCcase2}. 
We needed to remove the possibility of (iv) in Proposition \ref{np-prop-qjoin}. 
This was used in Sections \ref{np-sub-discretecase}, and \ref{np-sub-nondproof}. 

We may still obtain strictly joined end without the strong irreducibility. 
Premises and conclusions of Proposition \ref{np-prop-qjoin} still hold. 
	Let $\tilde E$ be a p-end corresponding to one of these.
	Then $\pi_{1}(\tilde E)$ acts on a properly convex domain 
    $K^{\prime \prime o}$ 	disjoint from $\vec{v}_{\tilde E}$. 

We use the commutant $H$ in Proposition \ref{np-prop-decomposition}.
By Propositions \ref{prelim-prop-Ben2} when $N_K$ is discrete and 
\ref{np-prop-NKsemi} when $N_K$ is nondiscrete, 
we may assume $H \subset N_K$. 
	We obtain elements 
 all the norms of whose eigenvalues corresponding to the subspace containing $K''$ 
	are $> 1$ and the rest of the norms of the eigenvalues are $< 1$
	by \eqref{np-eqn-finalform}.  
That is, these are the ones that correspond to the commutant $H$ 
and $N_K$  with these properties and act diagonalizable in $K$ with the smallest norm eigenvalue associated with the vertex of $K$. 
	This is in contradiction to the assumption $\lambda_{\mbv_{\tilde E}}(g) = 1$. 
%

\hfill	\SnT {\parfillskip0pt\par}
\end{proof}

%
%
%

\part{The deformation space of convex real projective structures}

The third part is devoted to understanding 
the deformation spaces of properly convex real projective structures on orbifolds with radial or totally geodesic ends. The end goal is to prove some versions of 
the Ehresmann-Thurston-Weil principle. 

In Chapter \ref{ch-op}, we give the precise definition of the deformation spaces.
We show that the deformation space of 
real projective structures on a strongly tame orbifold with some conditions on the ends is
mapped locally homeomorphically under the holonomy map 
to the character space of the fundamental group of the orbifold
with corresponding conditions. Here, we are not concerned with convexity. 
Thurston's idea of deformation via John Morgan 
as described in Lok \cite{Lok} by charts, works well here as well.  

In Chapter \ref{ch-rh}, we show that a properly convex real projective orbifold is strictly convex with respect to the ends if and only if the fundamental group is hyperbolic with respect to the end fundamental groups. 
Basic tools are from Yaman's work \cite{Yaman04} generalizing the Bowditch's description of 
hyperbolic groups. That is, we look at triples of points in the boundary 
of the universal cover and show that the action is properly discontinuous. 
In addition, we show that the action of the group on the fixed points of end fundamental groups is parabolic in their sense.  
This generalizes the prior work of Cooper-Long-Tillmann \cite{CLT15} and Crampon-Marquis \cite{CM14} for convex real projective 
manifolds with cusp ends. The concept of relative hyperbolic ends depends on 
the types of ends here unfortunately. Our aim was to generalize to orbifolds with our class of ends. 


In Chapter \ref{ch-cl}, we show that the deformation space of properly convex real projective structures on a strongly tame orbifold with some conditions on the ends is identifiable with the union of components of 
the character space of the fundamental group of the orbifold
with corresponding end conditions.

The openness part here continues that  of Chapter \ref{ch-op}. Here, the point is 
to prove the preservation of convexity under small perturbations. 
The proof consists of showing that we can patch the Hessian functions on
the perturbed compact part with the Hessian functions on the end neighborhoods approximating the original Hessian metrics by finding approximating convex domains to the original covering convex domains. 
Cooper-Long-Tillmann \cite{CLT15} 
uses the intrinsic Hessian metrics for ends instead. 

The closedness part generalizes the previous work Choi-Goldman \cite{CG93}, 
Benoist \cite{Benoist05} and Kim \cite{Kimi01}. 
We use the end condition showing us that the sequence of projectively covering convex domains 
can only degenerate to a point or a hemisphere. Then using Benzecri's work \cite{Benzecri60} putting the domains in a fixed ball and containing a fixed smaller ball, we show that the domain has to be actually properly convex. 

Finally, we go to Chapter \ref{ex-sec-nicecase}.
We discuss our nice cases Corollaries \ref{ex-cor-closed2} and \ref{ex-cor-closed3} where the Ehresmann-Thurston-Weil principle holds in
a simple way: the deformation space of the orbifold identifies with a union of 
components of character space of the orbifold fundamental group. 
These include the Coxeter orbifolds admitting complete hyperbolic structures.


\newcommand\DefSO{\Def^s_{\SI^{n}, {\mathcal{E}}, s'_{\mathcal{U}''} }(\mathcal{O})}  
\newcommand\DefEO{\Def^s_{\mathcal{E}, s_{\mathcal{U}}}(\mathcal{O})}

\newcommand\Imm{{\mathsf{Im}}}


\chapter{The openness of deformations   } \label{ch-op}

A real projective structure sometimes admits deformations to parameters of real projective structures. 
We prove the local homeomorphism 
between the deformation space of real projective structures on a strongly tame orbifold with radial or totally geodesic ends 
with various conditions with the $\SLnp$-character space (resp. $\PGL(n+1, \bR)$-character space) of the 
fundamental group with corresponding conditions. 
However, the convexity issues are not studied in this chapter.
In Section \ref{intro-sec-defspace}, 
we state the main
results recalling some definitions such as geometric structures, boundary restrictions, and the deformation spaces. 
In Section \ref{intro-sub-semialg}, we prove the semialgebraic properties of appropriate parts of character varieties. 
In Section \ref{op-sub-endstr}, we introduce a way to compactify our orbifolds and relate these to the local homeomorphism properties. We also define the deformation spaces. 
In Section \ref{op-sec-loch}, we prove the main result of the chapter
Theorem \ref{op-thm-projective}, showing the
openness of the deformation space in the character space. 
We first define the end conditions for real projective structures 
as determined by sections. We describe how to perturb the horospherical ends 
to lens-shaped ones in the affine setting. Then we state the main results. 
In Section \ref{op-sec-relation}, 
we identify the deformation spaces
as defined in our earlier papers  \cite{Choi06} and \cite{CHL12}
 as stated in Chapter \ref{ch-ex} to 
the deformation spaces here. 















\section{Deformation spaces and the spaces of holonomy homomorphisms}  \label{intro-sec-defspace}

Given a real projective $n$-orbifold $\orb$, we add the restriction of the end to be a radial or a totally geodesic type. 
The end is either assigned $\cR$-type or $\cT$-type. 
\index{Rtype@$\cR$-type} \index{Ttype@$\cT$-type}
\begin{itemize} 
\item An $\cR$-type end is required to be radial. 
\item A $\cT$-type end is required to have totally geodesic properly convex ideal boundary 
components or be horospherical. 
\end{itemize} 
\index{end!type-$\cR$|textbf} 
\index{end!type-$\cT$|textbf}
Recall that a strongly tame orbifold always has such an assignment in this monograph, 
and finite-covering maps always respects the types. 
Let $E_1, \dots, E_{e_1}$ be the $\cR$-ends, 
and $E_{e_1+1}, \dots, E_{e_1+e_2}$ be the $\cT$-ends. 



Recall that our strongly tame orbifold $\orb$ comes with 
a compact orbifold $\bar \orb$ with smooth boundary 
whose interior equals $\orb$. (See Hypothesis \ref{intro-hyp-compactification}.)
Each boundary component of $\bar \orb$ is said to be the ideal boundary 
component of $\orb$.

\begin{remark} \label{op-rem-compact}
Recall the compatibility condition of R-end structures and 
T-end structures with the compactification $\bar \orb$ of $\orb$
from Section \ref{intro-sec-ends}.  

For each R-end, 
we require that a vector field tangent to the leaves 
extends to a smooth vector field transverse to the corresponding 
ideal boundary component of $\bar \orb$. 

	The radial foliation gives us a smooth parametrization of 
	an end neighborhood $U$ of 
	a radial end or horospherical end $E$ of $\orb$ by $\Sigma_E \times (0,1)$
	where $x \times (0, 1)$ is the radial line for each $x \in \Sigma_E$
	since we can choose an embedded hypersurface 
	transverse to the radial rays.
	We assume that as $t \ra 1$, the ray escapes to the end.  

	
	
	Let $E$ be a $T$-end. We are given an end neighborhood $U$ diffeomorphic 
	to $S_E \times (0, 1)$ where $S_E$ is the ideal boundary 
	component. We identify $U$ with $S_E \times (0, 1)$ in 
	$S_E\times (0,1]$. 
We requires that $S_E \times \{1\}$ identifies with 
	the ideal boundary component of $\bar \orb$ corresponding to $E$
for this identification.  
	This is the compatibility condition of $\bar \orb$ 
	with a totally geodesic end structure for $E$.

	


	
\end{remark} 
\index{end!compactification} 
\index{end!radial!compactification} 
\index{end!totally geodesic!compactification} 



Also, $\bar \orb$ is a very good orbifold since we can identify
$\bar \orb$ with $\orb -U$ for a union $U$ of open end neighborhoods of product forms as above. 
This is obtained by identifying $\orb$ by $\orb - \clo(U)$ by
an isotopy preserving radial foliations and taking the closures. 
(See Theorem \ref{prelim-thm-vgood}.) 


There is an obvious isomorphism 
$\pi_1(\orb) \ra \pi_1(\bar \orb)$
since we can perturb any $\mathcal{G}$-path in $\bar \orb$ to 
one in $\orb$. 
For the universal cover $\hat \orb$ of $\bar \orb$,
there is an embedding 
$\tilde \orb \ra \hat \orb$ as an inclusion map 
to a dense open subset. 
We always identify $\tilde \orb$ with the dense subset. 
(See \cite{BH99}.)

An {\em isotopy} $i: \orb \ra \orb$ is a self-diffeomorphism such that  \index{isotopy|textbf}
there exists a smooth orbifold map $J: \orb \times [0, 1] \ra \orb$,  
such that 
\[i_{t}:\orb \ra \orb \hbox{ given by } i_{t}(x) = J(x, t)\] 
are self-diffeomorphisms for $t \in [0,1]$
and $i=i_1, i_0 = \Idd_{\orb}$. 
We require $i_t$ to be restrictions of 
isotopies 
\[ \bar i_t: \bar \orb \ra \bar \orb \hbox{ given by } \bar i_{t}(x) = \bar J(x, t)\]
is a self-diffeomorphism for each $t\in [0, 1]$ 
and $\bar J: \bar \orb \ra \bar \orb$ is a smooth orbifold map. 

\begin{definition} \label{op-defn-isolift}
	Let $\iota: \orb \ra \orb$ be an isotopy. 
	We may choose a lift $\tilde \iota: \tilde \orb \ra \tilde \orb$ of 
	$\iota$ such that 
	for the isotopy $F: \orb \times I \ra \orb$ with 
	$F_0 = \Idd_{\orb}$ and $F_t = \iota$ 
	has a lift $\tilde F$ such that $\tilde F_0 = \Idd_{\torb}$ and 
	$\tilde F_1 = \tilde \iota$. 
	We call such a map $\tilde \iota$ an {\em isotopy-lift}.
\end{definition} 
\index{isotopy-lift|textbf} 

Note that the radial structures for each R-end are sent to radial structures and 
the totally geodesic structures for each T-end are preserved 
since we required the radial foliations to extend to $\bar \orb$ smoothly
and the ideal boundary component to be the boundary component of $\bar \orb$ by the compatibility condition above
in Sections \ref{intro-sub-totgeo} and \ref{intro-sub-R-ends}.

%



We define 
{$\Def_{\mathcal E}(\mathcal{O})$} as the deformation space of real projective structures on $\mathcal{O}$ with end structures; more precisely, 
this is the quotient space of the real projective structures on $\mathcal{O}$ satisfying the above conditions for
ends of type $\cR$ and $\cT$  
under the isotopy equivalence relations.
We define the topology more precisely in Section \ref{op-sub-deform}. 
(See \cite{dgorb}, \cite{CEG06} and \cite{Goldman88} for more details. )
\index{deformation space} 
\index{deformation space!topology}
\index{DefE@$\Def_{\mathcal E}(\mathcal{O})$} 

For now, we restrict to compact orbifolds to recall the associated definitions: 
Suppose that $\mathcal{O}$ is compact. 
Let $G$ be a Lie group acting effectively and transitively on a space $X$. 
We define the {\em isotopy-equivalence space} $\widetilde{\Def}_{X, G}(\orb)$ as  \index{isotopy-equivalence space} 
the space of development pairs $(\dev, h)$ quotient by the right action of the isotopy-lifts of the universal cover $\torb$ of $\orb$.
The {\em deformation space} $\Def_{X,G}(\mathcal{O})$ is given by the quotient of $\widetilde{\Def}_{X, G}(\orb')$
by the left action of $G$: $g (\dev, h(\cdot)) = (g \circ \dev, g h(\cdot) g^{-1})$. (See \cite{dgorb} for details.) 
We can also interpret as follows: 
The deformation space $\Def_{X,G}(\mathcal{O})$ of the $(X,G)$-structures is \index{deformation space} 
the space of all $(X,G)$-structures on $\mathcal{O}$ quotient by the isotopy pullback actions.

This space can be thought of as the space of pairs $(\dev, h)$ with the compact open $C^r$-topology for $r \geq 1$ \index{deformation space!topology} 
and the equivalence relation generated by the isotopy relation 
\begin{itemize}
\item $(\dev, h) \sim (\dev', h')$ if $\dev'=\dev\circ \iota$ and $h'=h$ for an isotopy-lift $\iota$ of an isotopy (a right action)
and 
\item $(\dev, h) \sim (\dev', h')$ if $\dev'=k\circ \dev$ and $h(\cdot)=kh(\cdot)k^{-1}$ for $k \in G$.  (a left action)
\end{itemize} 
(See \cite{dgorb} or Chapter 6 of \cite{Cbook}.)

%



 \section{The semi-algebraic properties of 
 	$\rep^s(\pi_1(\mathcal{O}), \PGL(n+1, \bR))$ and related spaces}
 \label{intro-sub-semialg}

Since $\orb$ is strongly tame, the fundamental group $\pi_{1}(\orb)$ is finitely generated. 
Let $\{g_1, \dots, g_m\}$ be a set of generators of $\pi_1(\orb)$. 
As usual $\Hom(\pi_1(\mathcal{O}), G)$ for a Lie group $G$ has an {\em algebraic topology} as a subspace of $G^m$ by the map 
\[  h \in \Hom(\pi_1(\mathcal{O}), G) \mapsto (h(g_1), \dots, h(g_m)) \in  G^m.\]
This topology is given by the notion of {\em algebraic convergence}
\[\{h_i\} \ra h \hbox{ if } \{h_i(g_j)\} \ra h(g_j) \in G \hbox{ for each } j, j=1, \dots, m.\] 
A conjugacy class of a representation is called a {\em character} in this monograph. 

The {\em $\PGL(n+1, \bR)$-character space} ({\em variety}) $\rep(\pi_1(\mathcal{O}), \PGL(n+1,\bR))$ is the quotient space of \index{character space|textbf}
\index{character variety|textbf} 
the homomorphism space \[\Hom(\pi_1(\mathcal{O}), \PGL(n+1,\bR))\] where $\PGL(n+1,\bR)$ acts by conjugation
\[h(\cdot) \mapsto g h(\cdot) g^{-1} \hbox{ for } g \in \PGL(n+1,\bR).\]
Similarly, we define 
\[ \rep(\pi_1(\mathcal{O}), \SLpm):= \Hom(\pi_{1}(\orb), \SLpm)/\SLpm \] as 
the $\SLpm$-character space. This is not really a variety in
the sense of algebraic geometry. We merely define this space 
as the quotient space for now, possibly nonHausdorff one.  


Recall that 
a representation of a group $G$ into $\PGL(n+1, \bR)$ or $\SLpm$ is strongly irreducible 
if the image of every finite-index subgroup of $G$ 
is irreducible. 
Actually, many of the orbifolds have strongly irreducible and stable holonomy homomorphisms by Theorem \ref{intro-thm-sSPC}.

An {\em eigen-$1$-form} of a linear transformation $\gamma$ is a linear functional $\alpha$ in $\bR^{n+1}$ 
such that \index{eigen-$1$-form|textbf}
$\alpha \circ \gamma = \lambda \alpha$ for some $\lambda \in \bR$. 

In the following definitions, $\bG=\PGL(n+1, \bR)$ or $\SLpm$. 


\label{End} 
\begin{itemize} 
\item  \[\Hom_{\mathcal E}(\pi_1(\mathcal{O}), \bG)\] to be the subspace of representations $h$ satisfying  
\index{HomEP@$\Hom_{\mathcal E}(\pi_1(\mathcal{O}), \PGL(n+1,\bR))$|textbf} 
\index{HomES@$\Hom_{\mathcal E}(\pi_1(\mathcal{O}), \SLpm)$|textbf} 
\begin{description}
\item[The vertex condition for $\cR$] $h|\pi_{1}(\tilde E)$ 
has a common eigenvector of positive eigenvalues for 
the images of a  lifting of $h|\pi_1(\tilde E)$ to $\SLpm$ by Theorem \ref{prelim-thm-lifting}
\index{vertex condition}
for each $\cR$-type p-end fundamental group $\pi_{1}(\tilde E)$,
and
\item[The lens-condition for $\cT$] $h|\pi_{1}(\tilde E)$ 
acts on a hyperspace $P$ \index{end!lens condition|textbf} 
for each $\cT$-type p-end fundamental group $\pi_{1}(\tilde E)$
and acts 
discontinuously and cocompactly on a lens $L$, a properly convex domain with $L^o\cap P = L \cap P\ne \emp$
or a horoball tangent to $P$. 
\end{description} 
\item We denote by 
 \[\Hom^{s}(\pi_1(\mathcal{O}), \bG)\]
 the subspace of stable and irreducible representations, and define
  \[\Hom_{\mathcal E}^{s}(\pi_1(\mathcal{O}), \bG)\]
 to be \[\Hom_{\mathcal E}(\pi_1(\mathcal{O}), \bG) \cap \Hom^{s}(\pi_1(\mathcal{O}), \bG).\]
 
 \item We define
 \[\Hom_{\mathcal E, \mathrm{u}}(\pi_1(\mathcal{O}), \bG)\] 
 to be the subspace of 
 representations $h$  where
 \index{HomEu@$\Hom_{\mathcal E, \mathrm{u}}(\pi_1(\mathcal{O}), \bG)$} 
 \begin{itemize}
\item  $h|\pi_{1}(\tilde E)$  has a unique common eigenspace of dimension $1$ in $\bR^{n+1}$ with positive eigenvalues for 
the images of a  lifting of $h|\pi_1(\tilde E)$ to $\SLpm$ by Theorem \ref{prelim-thm-lifting}
 for each p-end holonomy group $h(\pi_{1}(\tilde E))$ 
of an $\mathcal{R}$-type 
\begin{itemize}
\item[(*)] There exists a finitely many elements $g_1, \dots, g_n$ such that 
the intersection $\bigcap_{i=1}^n C_{\lambda_i}(h(g_i))$ is $1$-dimensional 
where $C_{\lambda_i}(g_i)$ is the eigenspace of eigenvalue $\lambda_i$ of $h(g_I)$ 
associated with 
the above common eigenspace. (See Definition \ref{prelim-defn-jordan}.)
\end{itemize} 
and 
\item $h|\pi_{1}(\tilde E)$  has a common null-space $P$ of eigen-$1$-forms 
satisfying the following: 
\begin{itemize}
\item $\pi_{1}(\tilde E)$ acts properly and cocompactly on a lens $L$ and $L \cap P = L^o \cap P$ with a nonempty interior in $P$ with Hausdorff quotients, 
or 
\item $H-\{p\}$ for a horosphere $H$ tangent to $P$ at $p$
\end{itemize} 
and is unique such one 
for each p-end holonomy group $h(\pi_{1}(\tilde E))$ of the p-end of $\mathcal{T}$-type and dual group satisfies above (*) for the dual point of $P$. 
\index{HomEu@$\Hom_{\mathcal E, \mathrm{u}}(\pi_1(\mathcal{O}), \bG)$} 
\label{unique} 
We are putting somewhat strong condition such that some types of stability hold. 
(See Proposition \ref{intro-prop-u_semialg}.)



\end{itemize} 
 For $\mathcal{T}$-ends, the lens condition is satisfied for a hyperspace $P$ and $P$ is unique one satisfying the condition in other words.
We define
\[\Hom^s_{\mathcal E, \mathrm{u}}(\pi_1(\mathcal{O}),\bG)\] 
to be 
\[\Hom^s(\pi_1(\mathcal{O}), \bG \cap 
\Hom_{\mathcal E, \mathrm{u}}(\pi_1(\mathcal{O}), \bG).\] 
\index{HomsEu@$\Hom^s_{\mathcal E, \mathrm{u}}(\pi_1(\mathcal{O}), \bG)$} 
\end{itemize}

\begin{remark} 
The above condition for type $\mathcal{T}$ generalizes the principal boundary condition
for real projective surfaces of Goldman \cite{Goldman90}. 
\end{remark}


Since each $\pi_1(\tilde E)$ is finitely generated and there is only finitely many conjugacy 
classes of $\pi_{1}(\tilde E)$, 
 \[\Hom_{\mathcal E}(\pi_1(\mathcal{O}), \bG)\]
 is a closed semi-algebraic subset. 

\begin{proposition} \label{intro-prop-u_semialg}
	Let $\bG = \PGL(n+1, \bR)$ or $\SLpm$. 
	 \[\Hom_{\mathcal E, \mathrm{u}}(\pi_1(\mathcal{O}), \bG)\] 
	is an open subset of the semi-algebraic subset
	\[\Hom_{\mathcal E}(\pi_1(\mathcal{O}), \bG)\] 
	So is
	 \[\Hom^s_{\mathcal E, \mathrm{u}}(\pi_1(\mathcal{O}), \bG).\] 
	\end{proposition} 
\begin{proof}
%
We prove for the case $\bG = \PGLnp$. The other case is entirely similar. 
For the condition on the uniqueness, we may assume that involved $g_1, \dots, g_n$ are such that 
 $C_{\lambda_i}(h(g_i))$ are transversal and their intersection is $1$-dimensional.
Suppose that there is a sequence of holonomies $h_j:\pi_1(\tilde E)$ converging to $h|\pi_1(\tilde E)$
such that $h_j|\pi_1(\tilde E)$ has more than one fixed points. 
Then the span of the two directions are in the sum of 
$C_{\lambda_{i, j}}(h_j(g_i))$ and $C_{\lambda'_{i, j}}(h_j(g_i))$ for two 
eigenvalues of $\lambda_{i, j}, \lambda'_{i, j}$ of $h_j(g_i)$. 
Clearly, $\lambda_{i, j}, \lambda'_{i, j} \ra \lambda_i$. 
Hence, the limit of the sequence of the sum spaces must converge to the subspace of 
$C_{\lambda_i}(h(g_i))$ by elementary linear algebra. 
However, by the transversality of 
$C_{\lambda_i}(h(g_i)$ for $i=1, \dots, n$, we see that the sum spaces must intersect 
at a $1$-dimensional subspace for $h_i$ sufficiently close to $h$.
This means $h_i$ has only one fixed point, a contradiction. 
Therefore, the uniqueness condition is an open condition. 



Let $\tilde E$ be a T-p-end. 
Let \[h \in \Hom_{\mathcal E}(\pi_1(\mathcal{O}), \PGL(n+1,\bR)),\] 
and let $G:= h(\pi_1(\tilde E))$. 
Assume that $G$ is not a cusp group. 
Let $P$ be a hyperspace where $G$ acts on. 



The dual of the above discussions and 
Proposition \ref{app-prop-koszul} imply that the condition of
the existence of the hyperspace $P$ 
satisfying the lens-property is an open condition in 
$\Hom_{\mathcal E}(\pi_1(\mathcal{O}), \PGL(n+1,\bR))$. 

Suppose that 
there is another hyperspace $P'$ with a lens $L'$ 
satisfying the above properties. Then 
$P\cap P'$ is also $G$-invariant.
Note that $L^o \cap P$ covers a compact end-orbifold. 
By Proposition \ref{prelim-prop-Ben2}, 
we obtain that $\clo(P \cap L)$ is a join $K\ast \{k\}$ for a properly convex domain 
$K$ in $P\cap P'$ and a point $k$ in $P - P'$ since there exists a codimension-one invariant subspace in $P$. 
Similarly, exchanging the role of $P$ and $P'$, we obtain that 
there is a point $k' \in P' - P$ fixed by $G$. 
$G$ acts on the one-dimensional subspace $S_G$ containing $k$ and $k'$. 
There are no other fixed point on $S_G$ since otherwise $S_G$ is the set of 
fixed points and $G$ acts on any hyperspace containing $P\cap P'$ and 
a point on $S_G$. This contradicts our assumption on the first paragraph of the proof. 
Hence, only $k$ and $k'$ are fixed points in $S_G$
and $P\cap P'$ and $\{k, k'\}$ contain all the fixed points of $G$. 

Now, $k'$ is the unique fixed point outside $P$. 
The existence of lens for $P$ tells us that $k'$ must be a fixed point 
outside the closure of the lens. 
By Theorem \ref{pr-thm-equ2}, the existence of a lens for $P$ 
tells us that  every $g \in G$, the maximum norm of eigenvalues of $g$ 
associated with $k$ and $P\cap P'$ is greater than that of $k'$. 

Now, we switch the role of $P$ and $P'$. 
We can take a central element $g'$ with the largest norm of eigenvalue at $k'$ 
by the last item of Proposition \ref{prelim-prop-Ben2}
and the uniform middle-eigenvalue condition from Theorem \ref{pr-thm-equ2}. 
This cannot happen by the above paragraph. 
Hence, $P$ satisfying the lens-condition is unique. 

Suppose that $G$ is a cusp group. Then there exists a unique 
hyperspace $P$ containing the fixed point of $G$ tangent 
to horospheres where $G$ acts on. 
Therefore, \[\Hom_{\mathcal E, \mathrm{u}}(\pi_1(\mathcal{O}), \PGL(n+1,\bR))\] 
is in an open subset of a union of semi-algebraic subsets 
\[\Hom_{\mathcal E}(\pi_1(\mathcal{O}), \PGL(n+1,\bR)). \]
\end{proof} 




We define
\begin{itemize} 
\item  \[\rep_{\mathcal E}(\pi_1(\mathcal{O}), \bG)\] to be the quotient space
 \[\Hom_{\mathcal E}(\pi_1(\mathcal{O}), \bG)/\bG.\]
\item We denote by 
 \[\rep^{s}_{\mathcal E}(\pi_1(\mathcal{O}), \bG)\]
 the subspace of \[\rep_{\mathcal E}(\pi_1(\mathcal{O}), \bG)\]
 of stable and irreducible characters.
 \item 
 We define
 \[\rep_{\mathcal E, \mathrm{u}}(\pi_1(\mathcal{O}), \bG)\] to be
 \[\Hom_{\mathcal E, \mathrm{u}}(\pi_1(\mathcal{O}), \bG)/\bG.\]
 \item We define 
\begin{align} 
&  \rep^s_{{\mathcal E}, \mathrm{u}}(\pi_1(\mathcal{O}), \bG) \nonumber \\ 
 & := \rep^s(\pi_1(\mathcal{O}), \bG) \cap \rep_{{\mathcal E}, \mathrm{u}}(\pi_1(\mathcal{O}), \bG).
  \end{align} 
  \index{repsEu@$\rep^s_{{\mathcal E}, \mathrm{u}}(\pi_1(\mathcal{O}), \bG)$} 
   \index{repEu@$ \rep_{{\mathcal E}, \mathrm{u}}(\pi_1(\mathcal{O}), \bG)$} 
\end{itemize} 

 Let $\rho \in \Hom_{\mathcal E}(\pi_1(E), \bG)$
 where $E$ is an end. 
 Define  
 \[\Hom_{\mathcal E, \mathrm{par}} (\pi_1(E), \bG)\] 
 to be the subspace of representations where 
 $\pi_1(E)$ goes into a cusp group, i.e., 
 a parabolic subgroup in a conjugated copy of $\PO(n, 1)$ or $\Ort(n, 1)$
 By Lemma \ref{intro-lem-parab}, 
 \[\Hom_{\mathcal E, \mathrm{par}} (\pi_1(E), \bG)\] 
 is a closed semi-algebraic set. 
 \index{HomEpar@$\Hom_{\mathcal E, \mathrm{par}} (\pi_1(E), \bG)$}
 
 
 \begin{lemma} \label{intro-lem-parab}
Let $G$ be a finitely presented group. 
 	$\Hom_{\mathcal E, \mathrm{par}} (G, \bG)$ is a semi-algebraic set. 
 \end{lemma}
 \begin{proof}
 	We prove for $\bG = \PGL(n+1, \bR)$. 
 	Let $P$ be a maximal parabolic subgroup of a conjugated copy of $\PO(n+1, \bR)$ that fixes a point $x$. 
 	Then $\Hom(G, P)$ is a closed semi-algebraic set.
 	\[\Hom_{\mathcal E, \mathrm{par}}(G, \PGL(n+1, \bR))\] equals a union 
 	\[\bigcup_{g\in  \PGL(n+1, \bR)} \Hom(G, gPg^{-1})
\]
 	which is a semi-algebraic set
since the map $(g, p) \mapsto g pg^{-1}$ for $g\in \PGL(n, \bR)$ and $p\in P$ is a semi-algebraic map. 

 	When $\bG= \SLpm$, the proof is entirely similar. 
 \end{proof} 
 
 

 

 Let $E$ be an end orbifold of $\orb$. 
 Given \[\rho \in \Hom_{\mathcal E}(\pi_1(E), \bG),\] we define the following sets: 
 \begin{itemize} 
 	\item 
 	Let $E$ be an end of the type $\cR$.
 	Let \[\Hom_{\mathcal E, \mathbb{RL}}(\pi_{1}(E),  \bG)\] denote 
 	the space of  representations $h$ of $\pi_1(E)$
 	where $h(\pi_{1}(E))$ acts on a lens cone $\{p\} \ast L $ for a lens $L$ and $p$ for $p \not\in \clo(L)$ of a p-end $\tilde E$ 
 	corresponding to $E$ and acts properly and cocompactly on the lens $L$ itself.
 	Again, $\{p\}\ast L$ is assumed to be a bounded subset of an affine patch $\mathds{A}^n$.
 	Thus, it is an open subset of the semi-algebraic set
 	$\Hom_{\mathcal E}(\pi_1(E), \bG)$ by Proposition 
\ref{app-thm-qFuch}.
 	\index{HomERL@$\Hom_{\mathcal E, \mathbb{RL}}(\pi_{1}(E),  \PGL(n+1,\bR))$}
 	\index{HomERL@$\Hom_{\mathcal E, \mathbb{RL}}(\pi_{1}(E),  \SLpm)$}
 	\item 
 	Let $E$ denote an end of type $\cT$. 
 	Let \[\Hom_{\mathcal E, \mathbb{TL}}(\pi_{1}(E),  \bG)\] denote 
 	the space of totally geodesic representations $h$ of $\pi_1(E)$ satisfying 
 	the following conditions: 
 	\begin{itemize}
 		\item $h(\pi_{1}(E))$ acts on an lens $L$ and a hyperspace $P$ where
 		\item $L \cap P = L^{o} \cap P \ne \emp$ and 
 		\item $L/h(\pi_{1}(E))$ is a compact orbifold with two strictly convex boundary components.
 	\end{itemize}  
 	\[\Hom_{\mathcal E, \mathbb{TL}}(\pi_{1}(E),  \bG)\] 
 	is again a union of open subsets of the semi-algebraic sets 
 	\[\Hom_{\mathcal{E}}(\pi_{1}(E),  \bG)\] 
 	by Proposition \ref{app-thm-qFuch2}. 
 	\index{HomETLP@$\Hom_{\mathcal E, \mathbb{TL}}(\pi_{1}(E),  \PGL(n+1,\bR)) $} 
 	\index{HomETLS@$\Hom_{\mathcal E, \mathbb{TL}}(\pi_{1}(E),  \SLpm) $} 
 \end{itemize} 
 
 Let 
 \[ R_{E}: \Hom(\pi_1(\mathcal{O}), \bG) \ni h \ra h|\pi_{1}(E) \in \Hom(\pi_1(E), \bG) \] 
 be the restriction map to the p-end holonomy group $h(\pi_1(E))$ corresponding to the end $E$ of $\mathcal{O}$. 
 
 A {\em representative set} of p-ends of $\torb$ is the subset of p-ends where 
 each end of $\orb$ has a corresponding p-end and a unique chosen corresponding p-end.
 Let $\mathcal{R}_{\orb}$ denote the representative set of p-ends of $\torb$ of type $\mathcal{R}$,  
 and let $\mathcal{T}_{\orb}$ denote the representative set of p-ends of $\torb$ of type $\mathcal{T}$. 
 \label{ce} 
 We define a more symmetric space:
 \[\Hom_{\mathcal{E}, \lh}(\pi_1(\mathcal{O}), \bG)\] to be
 {\small
 	\begin{align} \label{intro-eqn-euce2}
 	&  \bigg(\bigcap_{E \in {\mathcal{R}_{\orb}}} R_{E}^{-1}\Big(\Hom_{\mathcal E, \mathrm{par}}(\pi_{1}(E),  \bG) \cup 
 	\Hom_{\mathcal E, \mathbb{RL}}(\pi_{1}(E),  \bG)\Big)\bigg) \cap \nonumber \\
 	&  \bigg(\bigcap_{E \in {\mathcal{T}_{\orb}}} R_{E}^{-1}\Big(\Hom_{\mathcal E, \mathrm{par}}(\pi_{1}(E),  \bG) \cup \Hom_{\mathcal E, \mathbb{TL}}(\pi_{1}(E),  \bG)\Big)\bigg).
 	\nonumber
 	\end{align} 
 }
The quotient space of the space under the conjugation under $\bG$ is denoted by 
\[\rep_{\mathcal{E}, \lh}(\pi_1(\mathcal{O}), \bG).\] 
\index{HomEceP@$\Hom_{\mathcal{E}, \lh}(\pi_1(\mathcal{O}), \PGL(n+1,\bR))$ } 
\index{repEceP@$\rep_{\mathcal{E}, \lh}(\pi_1(\mathcal{O}), \PGL(n+1,\bR))$ } 
\index{HomEceS@$\Hom_{\mathcal{E}, \lh}(\pi_1(\mathcal{O}), \SLpm)$ } 
\index{repEceS@$\rep_{\mathcal{E}, \lh}(\pi_1(\mathcal{O}), \SLpm)$ } 
  We define 
  \[\Hom^s_{\mathcal{E}, \lh}(\pi_1(\mathcal{O}), \bG)\]
  to be 
  	\begin{equation*} 
 	 \Hom^s(\pi_1(\mathcal{O}), \bG)\, \cap   
 	 \Hom_{\mathcal{E}, \lh}(\pi_1(\mathcal{O}), \bG).
 	\end{equation*}
 	Hence, this is an open subset of a semialgebraic subset 
 	\[X:=\Hom^s_{\mathcal E}(\pi_1(\mathcal{O}), \bG).\]
 	These definitions allow for changes between horospherical ends to lens-shaped radial ones and totally geodesic ones.
 	\label{superscripts} 
 	
 	The quotient space of this space 
 	under the conjugation under $\bG$ is denoted by 
 	\[\rep^s_{\mathcal{E}, \lh}(\pi_1(\mathcal{O}), \bG).\] 
 	\index{HomsEceP@$\Hom^s_{\mathcal{E}, \lh}(\pi_1(\mathcal{O}), \PGL(n+1,\bR))$}
 		\index{repsEceP@$\rep^s_{\mathcal{E}, \lh}(\pi_1(\mathcal{O}), \PGL(n+1,\bR))$}
 			\index{HomsEceS@$\Hom^s_{\mathcal{E}, \lh}(\pi_1(\mathcal{O}), \SLpm)$}
 		\index{repsEceS@$\rep^s_{\mathcal{E}, \lh}(\pi_1(\mathcal{O}), \SLpm$}

 Since  \[\rep^s_{{\mathcal E}, \mathrm{u}}(\pi_1(\mathcal{O}), \bG)\]
 is the Hausdorff quotient of the above set with the conjugation $\bG$-action, 
 this is an open subset of a semi-algebraic subset 
 by Proposition \ref{intro-prop-u_semialg} and 
 Proposition 1.1 of \cite{JM87}.

 
 We define \[\Hom^s_{\mathcal E, \mathrm{u}, \lh}(\pi_1(\mathcal{O}), \bG)\] to be the subset 
\[ \Hom^s_{\mathcal E, \mathrm{u}}(\pi_1(\mathcal{O}), \bG) \cap
 \Hom_{\mathcal{E}, \lh}(\pi_1(\mathcal{O}), \bG).\]
\index{HomsEuce@$\Hom^s_{\mathcal{E}, \mathrm{u}, \lh}(\pi_1(\mathcal{O}), \PGL(n+1,\bR))$} 
 
 The above shows:  
 \begin{proposition}\label{intro-prop-uce_semialg} 
Let $\bG=\PGL(n+1, \bR)$ or $\SLpm$. 
 \[\rep^s_{{\mathcal E}, \mathrm{u}, \lh}(\pi_1(\mathcal{O}), \bG)\]
 is an open subset of a semi-algebraic set 
 \[\rep^s_{{\mathcal E}}(\pi_1(\mathcal{O}), \bG).\]
 \end{proposition} 

\index{repsEuce@$\rep^s_{{\mathcal E}, \mathrm{u}, \lh}(\pi_1(\mathcal{O}), 
	\PGL(n+1, \bR))$}

\section{End structures and end compactifications for topological orbifolds}
\label{op-sub-endstr}

Let $\orb$ be a strongly tame smooth $n$-orbifold
with ends $E_1, \dots, E_m, E_{m+1}, \dots, E_{m+l}$. 
For this subsection, we do not consider that $\orb$ is the interior of 
a compact orbifold $\bar \orb$, the associated compactification of 
$\orb$, as we remind from 
Section \ref{intro-sub-endf}.  This is because we wish to consider more than one compactifications. 

An {\em ideal boundary structure} of an end neighborhood $U$ of $E_i$ 
is a pair $(U, f)$ for 
a smooth embedding $f$ of $U$ into a product space $\Sigma \times (0, 1]$ for 
a closed $(n-1)$-orbifold $\Sigma$
where the image is $\Sigma\times (0, 1)$. 
An ideal boundary structure $(U_0, f_0)$ with 
a diffeomorphism $f_0: U_0 \ra \Sigma_0 \times (0, 1)$ for 
a closed $(n-1)$-orbifold $\Sigma_0$ and another one 
$(U_1, f_1)$ 
with a diffeomorphism $f_1: U_1 \ra \Sigma_1\times (0, 1)$ for 
a closed $(n-1)$-orbifold  $\Sigma_1$ 
are {\em compatible} if 
there exists another ideal boundary structure $(U_2, f_2)$ 
such that $U_2 \subset U_0\cap U_1$ with 
a diffeomorphism $f_2: U_2 \ra \Sigma''\times (0, 1)$
such that $f_i\circ f_2^{-1}: \Sigma'' \times (0, 1) \ra \Sigma_i \times (0, 1)$ 
extends to $\Sigma''\times (0, 1]$ as an embedding 
restricting to a diffeomorphism $\Sigma'' \times \{1\}$ 
to $\Sigma_i \times \{0\}$ for $i=0,1 $.



Given an ideal boundary structure on $\orb$ for an end $E_i$, we obtain 
the {\em completion of $\orb$ along $E_i$}. 
We take $U$ an end neighborhood of $E_i$ with an embedding $f:U \ra\Sigma \times (0, 1]$
where the image equals $\Sigma \times (0, 1)$. 
We paste $\orb$ with $\Sigma \times (0, 1]$ by $f$.
The resulting orbifold $\orb'_{(U, f)}$ is said to be the
{\em end compactification of $\orb$ along $E_i$}
using $(U, f)$. 

Let $U'$ and $f': U' \ra \Sigma' \times (0, 1]$ 
be as above with $(U', f')$ compatible to $(U, f)$,
and we obtain an end compactification $\orb'_{(U', f')}$ 
of $\orb$ along $E_i$
using $(U', f')$. 

An {\em isotopy} $\iota$ of $\orb$ with an ideal boundary structure for an end $E_i$ is 
an isotopy of $\orb$ extending to a diffeomorphism 
$\bar \iota: \orb'_{(U, f)} \ra \orb'_{(U', f')}$ for 
at least one compatible pair $(U, f), (U', f')$.

By the following lemma, the definition of an end-structure extendable isotopy is independent of the choice of $(U, f)$ and $(U', f')$. 

\begin{lemma} \label{op-lem-endcompT} 
	Let $U_1$ and $U_2$ be end neighborhoods of a T-end $E$. 
	Let $g$ be an isotopy of $\orb$
	extending to an isotopy 
	$\bar g: \orb'_{(U_1, f_1)} \ra \orb'_{(f_2, U_2)}$
	with diffeomorphisms 
	$f_1: U_1 \ra \Sigma_1 \times (0, 1)$ and $f_2:U_2 \ra \Sigma_2\times (0, 1)$
	for closed $(n-1)$-orbifolds $\Sigma_1$ and $\Sigma_2$. 
Then for any pair 
$(U'_1, f'_1), (U'_2, f'_2)$ for end neighborhoods of $E$ 
compatible to $(U_i, f_i), i=1,2$, 
with diffeomorphisms 
$f'_1: U'_1 \ra \Sigma'_1 \times (0, 1)$ and $f'_2:U'_2 \ra \Sigma'_2\times (0, 1)$
for closed $(n-1)$-orbifolds  $\Sigma'_1$ and $\Sigma'_2$, 
$g$ extends to an isotopy
$\bar g_1: \orb_{(U'_1, f'_1)} \ra \orb_{(U'_2, f'_2)}$. 
\end{lemma} 
\begin{proof} 
	This is straightforward to obtain a diffeomorphism $\bar g_1$
	 since we can take a sufficiently small product neighborhood in each of these end compactifications. 
	To show the isotopy property of $\bar g_1$, we simply take $\orb \times I$ and 
	do the same arguments. 
%
\end{proof}

%
%
%
%
%

A radial structure for ${E_i}$ also gives us an end-compactification 
of $\orb$ along ${E_i}$\/:
Let $U$ be an end neighborhood of ${E_i}$ with a foliation by properly embedded arcs. 
We take a transverse hypersurface $\Sigma_{{E_i}}$ transverse to every leaf,
which is a closed orbifold. Let $U'$ denote a component of 
$U - \Sigma_{{E_i}}$ that is an end neighborhood of ${E_i}$. 
We identify each leaf in $U'$ with a leaf of $\Sigma_{{E_i}}\times (0, 1)$
by a function $f:U' \ra \Sigma_{{E_i}} \times (0, 1)$. 
We call the identification orbifold $\orb'$ of $\orb$ with $\Sigma_{{E_i}}\times (0, 1]$
the {\em end compactification of $\orb$ along ${E_i}$}.   
The suborbifold of $\orb'$ corresponding to $\Sigma_{{E_i}}\times \{1\}$ is 
called the {\em ideal boundary component corresponding to ${E_i}$}. 
This orbifold is independent of the choices up to isotopies of 
the end compactifications
extending isotopies of $\orb$ by Lemma \ref{op-lem-endcompR}.

%
%


\begin{lemma} $ $ \label{op-lem-endcompR} 
	\begin{itemize} 
	\item Let $\orb$ be a strongly tame real projective
	$n$-orbifold with a radial structure at an end $E_i$.
	\item Let $\orb'$ be the end compactification of $\orb$ using $U$ and $\Sigma_{E_i}$ 
	a diffeomorphism $f: U' \ra \Sigma_{E_i}\times (0, 1)$. 
	\item Let $\orb_1$ be the same orbifold with 
	an isotopic radial structure at $E_i$.  
	\item Let $\orb'_1$ be the end compactification of $\orb_1$ 
	for the second radial structure 
	using an end neighborhood $U_1$ and 
	a hypersurface $\Sigma'_{E_i}$ and a diffeomorphism $f': U_1' \ra 
	\Sigma'_{E_i}\times (0, 1)$ for an end neighborhood component $U_1'$ of 
	$U_1 -\Sigma'_{E_i}$. 
\end{itemize} 
	Suppose that an isotopy $\iota$ of $\orb$ sends
	a radial structure of $U$ for end $E_i$ to that of $U_1$ for $\orb_1$ 
	for $E_i$. 
	
	Then an isotopy $\iota': \orb \ra \orb$ extends
	to an isotopy  
	$\hat\iota': \orb' \ra \orb'_1$ sending the ideal boundary component 
	corresponding to $E_i$ in $\orb'$ to one in $\orb'_1$.  
\end{lemma} 
\begin{proof} 
	First, we obtain a diffeomorphism $i'$. 
	We may change $\iota$ such that $\iota(U') \subset U'_1$ 
	by composing with an isotopy supported in $U$ preserving leaves of 
	the radial foliation. 
	
	We consider a diffeomorphism 
	\[f' \circ \iota\circ  f^{-1}|\Sigma_{E_i}\times (0, 1) 
	\ra \Sigma'_{E_i}\times (0, 1)\]
sending each $x \times (0,1)$ to itself. 
	We can find an isotopy $\iota_1: \Sigma'_{E_i}\times (0, 1) 
	\ra \Sigma'_{E_i} \times (0, 1)$ preserving $\{x\} \times (0, 1)$ 
and being the identity on $x \times (0, \eps)$  for 
	each $x \in \Sigma'_{E_i}$ 
	 such that $\iota_1\circ f'\circ \iota \circ f^{-1}$ 
preserving the second parameter of $x \times (1- \eps, \eps)$ for 
each $x\in \Sigma_{\tilde E_i}$ and 
sufficient small $\eps> 0$. 
Hence, $\iota_1\circ f'\circ \iota \circ f^{-1}$ 
	extends to a smooth map at $\Sigma_{E_i}\times \{1\}$. 
	This is fairly simple to show since every self-embedding 
	of $\Sigma_{E_i}\times (0, 1)$ preserving every 
	fiber of form $\{x\}\times (0, 1)$ for $x \in \Sigma'_{E_i}$ are isotopic. 
	Now, we can use $\iota: \orb - U' \ra \orb - U'_1$ 
	and $f^{\prime -1}\circ \iota_1 \circ f'\circ \iota \circ f^{-1}$ on 
	$U'$. Obviously, they extend each other. 
	
\end{proof}

The following shows the well-definedness of the compactification of
a strongly tame orbifold with R- or T-ends.  
\begin{corollary} \label{op-cor-wellend} $ $ 
	\begin{itemize} 
	\item Let $\orb$ be a strongly tame real projective $n$-orbifold with a T-end structure 
	at $E_i$. Then the identity map $\Idd$ of $\orb$ is isotopic 
	to a restriction of a diffeomorphism $\orb'_{(U, f)} \ra \orb'_{(U', f')}$
provided $(U,f)$ and $(U', f')$ are compatible. 
	\item Let $\orb$ be a strongly tame real projective $n$-orbifold with an R-end structure at 
	$E_i$. Then the identity map $\Idd$ of $\orb$ is isotopic to 
	a restriction of 
	a diffeomorphism $\orb'_1 \ra \orb'_2$ for 
	any two end-compactifications $\orb'_1$ and $\orb'_2$ of $\orb$ compatible with 
the R-end structures. 
	\end{itemize} 
	\end{corollary} 
\begin{proof} 
	The first is a corollary of Lemma \ref{op-lem-endcompT}.
	The second one is a corollary of Lemma \ref{op-lem-endcompR}. 
	\end{proof} 

Finally, we say about the compactification $\bar \orb$ 
associated with $\orb$. 

If $\bar \orb$ is the compactification associated with 
$\orb$, the ideal boundary structure 
is given by $(f, N(\Sigma_{E_i})\cap \orb)$ where $f: N(\Sigma_{E_i})\cap \orb \ra \Sigma_{E_i}\times (0, 1]$ is an embedding 
for a tubular neighborhood $N(\Sigma_{E_i})$ of $\Sigma_{E_i}$ in $\bar \orb$. 

Again, if $\bar \orb$ is the associated compactification of $\orb$, 
and the radial structure at ${E_i}$ is compatible with $\bar \orb$ 
as in Section \ref{intro-sub-R-ends}, 
the radial end compactification from $(f, U)$
can be modified to a compatible $(f', N(\Sigma_{E_i})\cap \orb)$ 
for a tubular neighborhood $N(\Sigma_{E_i})$ of $\Sigma_{E_i}$ to 
$\Sigma_{E_i}\times (0, 1]$
where $f'$ extends to a smooth diffeomorphism 
$N(\Sigma_{E_i}) \ra \Sigma_{E_i} \times (0, 1]$. 


\begin{proposition}\label{op-prop-endc} 
We can construct by the above end compactification process
a compact orbifold $\bar \orb$ of which $\orb$
is the interior. 
End compactifications of $\orb$ for all ends 
compatible with the given R-end and T-end structures is always diffeomorphic to 
$\bar \orb$ by a diffeomorphism isotopic to the identity in $\orb$.
\end{proposition} 
\begin{proof} 
Recall from Sections \ref{intro-sub-totgeo} and \ref{intro-sub-R-ends}
the definitions of compatibility. 	
Also, it is straightforward to see that the radial foliation is 
transverse to the added ideal boundary component 
corresponding to $\Sigma_E \times \{1\}$. 
Corollary \ref{op-cor-wellend} completes the proof.
	\end{proof} 


\subsection{Definition of the deformation spaces with end structures} \label{op-sub-deform} 

We extend this notion strongly. 
Two real projective structures $\mu_0$ and $\mu_1$ on $\orb$ with R-ends or T-ends with end structures are {\em isotopic} 
	if there is an isotopy $i$ on $\mathcal{O}$ such that $i^*(\mu_0)=\mu_1$ where $i^*(\mu_0)$ is the induced structure from $\mu_0$ by $i$:
	\begin{itemize}
		\item $i_{\ast}(\mu_{0})$ has a radial end structure for each R-end or 
		horospherical T-end, 
		\item $i$ sends the radial end foliation for $\mu_0$ from an R-end neighborhood or 
		horospherical T-end to the radial end foliation for real projective 
		structure $\mu_1 = i_{\ast}(\mu_{0})$ with corresponding R-end neighborhoods or 
		a horospherical T-end, 
		and
		\item $i$ extends to a diffeomorphism of 
		$\bar \orb$ using the radial foliations 
		and the totally geodesic ideal boundary components for $\mu_0$ and $\mu_1$
		where we use the radial end-compactification for a horospherical T-end. 
  (See Remark  \ref{op-rem-compact}.)
	\end{itemize}
\index{isotopy|textbf}

For noncompact orbifolds with end structures, similar definitions hold except that we have to modify
the notion of isotopies to preserve the end structures. 

\begin{definition}\label{op-defn-isot}
We consider the real projective structures on orbifolds with end structures. 
Let $\orb$ be one of this and $\bar \orb$ be the compactification. 
Let $\hat \orb$ denote the universal cover 
of $\bar \orb$ containing $\tilde \orb$ as a dense open set. 

Let $\dev_\mu$ denote the developing map associated with a convex projective structure $\mu$ with R-end or T-ends. 
The developing map 
 $\dev_\mu:\tilde \orb \ra \rpn$ sometimes extends 
to a map $\overline{\dev}_\mu: \hat \orb \ra \rpn$. 
We need developing maps determined up to isotopies. 
For R-ends, we may choose by a right action by isotopies 
such that $\overline{\dev}_\mu$ is smooth by 
Lemma \ref{op-lem-replacedev}. 
For T-ends, we can always isotopy $\dev_\mu|U$ for a p-end neighborhood $U$ of a p-end 
$\tilde E$ such that it can extend to a smooth map by Lemma \ref{op-lem-Tendisotopy}. 
\end{definition}

\begin{lemma} \label{op-lem-radialisotopy}. 
Let $f_0$ and $f_1$ be two immersions $\tilde\Sigma_{\tilde E}\times (0, 1] \ra \rpn$
{\rm (}resp. $\ra \SI^n${\rm )} 
equivariant with respect to
a holonomy representation $\rho: \pi_1(\Sigma_E) \ra \PGL(n+1, \bR)$ 
{\rm (}resp. $\ra \SL_\pm(n+1, \bR)${\rm )}
fixing a point $p_0$ in $\rpn$ {\rm (}resp. $\SI^n${\rm )}. 
Assume the following\/{\rm :}  
\begin{itemize} 
\item for each $i=0, 1$, 
$f_i|x \times (0, 1] \ra l_x$ is an embedding to a radial segment $l_x$ with endpoint $p_0$ in $\rpn$ 
{\rm (}resp. $\ra \SI^n${\rm )} 
and $f_i(x\times t)$ converges to $p_0$ as $t \ra 0$. Here, $l_x$ is independent of $i=0, 1$. 
\item $f_0| \tilde \Sigma \times [\delta, 1] = f_1| \tilde \Sigma \times [\delta, 1]$ for $\delta > 0$. 
\end{itemize} 
Then $f_0$ and $f_1$ are smoothly isotopic by an isotopy preserving each $x\times (0,1]$
being the identity on $\tilde \Sigma \times [\delta, 1]$. 
\end{lemma} 
\begin{proof} 
	We prove the result for $\rpn$. The proof for the result for $\SI^n$ is entirely similar. 
Let $\mathcal{C}_{x, \delta, f_0}((0, 1], l_x)$ denote the space of embeddings 
$g| (0, 1]$ where $g|(\delta, 1]$ is fixed to be $f_0|x\times [\delta, 1] \ra l_x$, and 
$g(t) \ra p_0$ as $t \ra 0$. This is a contractible space since this is a convex space 
if we identify $l_x$ with a real interval. 
We can build a bundle $\tilde B$ over $\tilde \Sigma_{\tilde E}$ with fiber at $x$ 
equal to $\mathcal{C}_{x, \delta, f_0}((0, 1], l_x)$. 
Then $\pi_1(\tilde E)$ acts on this space where we must act by $\rho$ to
the fibers. Hence, we can find a quotient space $B$ fibering over $\Sigma_{\tilde E}$. 
Then consider the space 
$C_{\rho, x_0, \delta}(\tilde \Sigma_{\tilde E})$ to be the space 
of $\rho$-equivariant sections sending $x \in \tilde \Sigma_{\tilde E}$ to an element of 
$\mathcal{C}_{x, \delta, f_0}((0, 1], l_x)$. 

$f_0$ and $f_1$ give two such sections. 
This induces section $\hat f_0$ and $\hat f_1$ of $B$ over $\Sigma_{\tilde E}$. 
By contractibility of the fibers, we can use the obstruction theory to 
obtain the homotopy between $\hat f_0$ and $\hat f_1$. 
This gives us the $\rho$-equivariant homotopy $f_t, t\in [0, 1]$, 
between them using the contractibility of fibers.  
We use the fact that 
a homotopy $f'_t$ between two homeomorphisms $f'$ and $f''$ of intervals can 
be realized by an isotopy $(f'_t)^{-1}\circ f'$ while each $f'_t$ is a homeomorphism. 
Hence, these homotopies give us the desired isotopies. 
\end{proof}

We now describe the modification of the developing map 
by a process that we call the {\em radial-end projectivization of 
	the developing map with respect to $U$ and $U'$.} 
That is, we prove Lemma \ref{op-lem-replacedev}. 
\index{developing map!radial-end projectivization} 

In the following $\bX = \rpn$ or $\SI^n$ while 
$\bG = \PGL(n+1, \bR)$ or $\SLpm$ respectively
and $\bX^\ast = \RP^{n\ast}$ or $\SI^{n\ast}$ respectively. 
We consider a developing map $\dev: \torb \ra \bX$, and modify it. 
	
Note that $\bR_+ \cup \{\infty\}$ is projectively equivalent to an interval $(0, 1]$ where 
$\infty$ corresponds to $1$. 
	Let $\tilde U$ and $\tilde U'$ denote the closed p-end neighborhoods of 
	$\tilde E$ covering end neighborhoods $U$ and $U'$, $\clo(U) \subset U^{\prime o}$, respectively.
We require $U$ and $U'$ to be compatible product neighborhoods
diffeomorphic to $\Sigma_E \times (0, 1]$. 
	(Recall the compatibility from Section \ref{intro-sec-ends}.)
	Take a maximal radial ray $l_x$ in $\tilde U'$ 
	passing $x \in \Bd \tilde U \cap l$. Then 
	there exists a unique projective diffeomorphism 
$\Pi_x: l_x \ra \bR_+ \cup \{\infty\}$ sending 
	\begin{itemize} 
	\item the endpoint of $l_x$ in $\Bd \tilde U'$ to $\infty$, 
	\item the other end to $0$, and 
	\item $\Bd \tilde U \cap l = \{x\}$ to $1$. 
	\end{itemize} 
We define $\Pi_{\tilde U', \tilde U}: \tilde U' \ra \bR_+$ by 
sending $z \in l_x$ to $\Pi_x(z)$. 
	There is also a unique projective diffeomorphism $P_x: \bR_+\cup \{\infty\} 
\ra \bX$ sending 
	\[ 0 \mapsto \mbv_{\tilde E}, 
	1 \mapsto \dev(l\cap \Bd \tilde U)= \dev(x),  +\infty \mapsto
	\dev(l \cap \Bd \tilde U').  \]
	Define $\vec{v}_x$ 
	to be the vector at $\mbv_{\tilde E}$ 
	of $(\partial P_x(t)/ \partial t)|_{t=0}$. 
	This does depend on $x$ but not on $t$. 
	

Let $\Pi_{\tilde E}: \tilde U' \ra \tilde \Sigma_E$ denote the map sending a point of 
a radial ray in $\tilde U'$ to its equivalence class in $\tilde \Sigma_E$. 
	$\tilde U'$ has coordinate functions 
	\[(\Pi_{\tilde U', \tilde U}, \Pi_{\tilde E}): \tilde U' \ra \bR_+ \times \tilde \Sigma_E,\]
	which is a diffeomorphism. 
	This commutes with the action of $\bGamma_{\tilde E}$ on $\tilde U'$ 
	and the action on $\bR_+\times \tilde \Sigma_E$ is acting on the first factor
	trivially. Also, 
	$\tilde U$ goes to $(0, 1)\times \tilde \Sigma_E$ under 
	the map. 
	
	
	We define a smooth map 
\begin{equation}\label{op-eqn-devN} 
\dev^N: \tilde U' \ra \bX \hbox{ given 	by }  
	\dev^N(y) = P_x \circ \Pi_x(y) \hbox{ for } y \in l_x \subset \tilde U'.
\end{equation}
For each $l_x$, we choose an affine subspace $\mathds{A}^n_x$ containing 
$\dev(l_x)$. These are all related by projective transformations in their intersections 
in $\bX$. 
	Then under the coordinate  of $\tilde U'$ with affine coordinates 
	on $\mathds{A}^n_x$ with a temporary Euclidean norm $\llrrV{\cdot}_x$ containing $\dev(l_x)$ 
	and containing $\dev(\mbv_{\tilde E})$ as the origin, 
	we can write locally
	\begin{equation} \label{np-eqn-fxvx} 
	\dev^N(x, t) = f_x(t) \vec{v}_x, \llrrV{\vec{v}_x}_x=1, 
	\hbox{ for } x \in \Bd \tilde U
	\end{equation} 
	on a neighborhood of $l_{x_0}$ for some $x_0 \in \Bd \tilde U$  
	where $\vec{v}_x$ is a unit vector depending only on $x$ smoothly in the 
	direction of $\overline{\dev(\mbv_{\tilde E}) \dev(x)}$  and 
	\begin{multline}
	f_x: \bR_+ \ra \bR_+, f_x(0)=0, 
	f_x(1) = \llrrV{\dev(x)}_x \\
	f_x(\infty) = \llrrV{\dev(l_x \cap \Bd U')}_x 
	\hbox{ provided } \dev(l_x \cap \Bd U') \in \mathds{A}^n_x
	\end{multline} 
	is a strictly increasing projective function of $t$.
	The coefficients of the $1$-st order 
	rational function $f_x$, as a function of $t$, depend smoothly on $x$
	since $\partial U'$ and $\partial U$ 
	are smooth. 
	It is easy to see that $\dev^N$ extends to 
	$\tilde \Sigma_{E}\times \{0\}$ as a constant map. 
	The expression of $f_x$ shows that 
	$\dev^N$ is a smooth extension also since it has continuous partial derivatives of all orders.

Note that $f_x$ is well-defined for $l_x$ independent of the choice of 
$\mathds{A}_x$ and $\llrrV{\cdot}_x$ by the uniqueness of projective map 
satisfying the properties that $f_x$ does.
	
	Now we change $\dev^N$ on $\tilde U' -\tilde U$ such that 
	it smoothly extends to $\torb - \tilde U'$. 
Actually on $\Bd \tilde U'$ and on $\Bd \tilde U$, $\dev^N$ and $\dev$ agree by our construction. 
	On each $l_x$, the lines $\dev|l_x$ and $\dev^N|l_x$ 
	have the same image to $\dev(l_x)$. 
	Hence, $(\dev^{N}|l_x)^{-1}\circ \dev|l_x$ sends $l_x$ to $l_x$ as a homeomorphism. 
We choose $U''$ to be the open neighborhood of $U$ in $U'$ whose closure is in the interior of $U'$. 
Let $\tilde U''$ denote the component of inverse image containing $\tilde U$.  
	We choose a partition of unity function $\phi$ equal to $1$ on $\tilde U' -\tilde U''$ and 
with support in a small neighborhood of $\tilde U'-\tilde U''$ in $\tilde U' -\clo(\tilde U)$. 
Of course, we need to choose $\phi$ equivariant with respect to $\pi_1(\tilde E)$, 
which can be done as in the proof of Lemma \ref{op-lem-radialisotopy}. 
Also, the directional derivative of $\phi$ in the radial direction of $l_x$ is non-negative. 
We define $g_x: l_x \ra \bR_+$ given by $||\dev(z)||$ for $z \in l_x$, which is a smooth function with
positive radial directional derivatives along $l_x$.
By isotopying for $\dev^N$ towards the p-end vertex equivariantly and taking smaller $U'$ and $U$, 
we may assume without loss of generality that $g_x \geq f_x$ for all $x \in \tilde \Sigma_{\tilde E}$. 
Then we define a new smooth function $h:=(1-\phi)f_x + \phi g_x$ which still has positive directional derivatives in the radial direction of $l_x$. 
Then this agrees  with $g_x$ in $\tilde U' -\tilde U''$ and $f_x$ outside a neighborhood of 
$\tilde U' -\tilde U''$ in $\tilde U' - \clo(\tilde U)$.
We define $\dev': \tilde U' \ra \bX$ by replacing $\Pi_x$ 
 with $h_x$ in \eqref{op-eqn-devN} for every $x \in \tilde \Sigma_{\tilde E}$.
This map is the {\em radial-end projectivization} of $\dev$ with respect to $U'$ and $U$
and agrees with $\dev_N$ outside a neighborhood of $\tilde U'-\tilde U''$ in $\tilde U'$ and 
agrees with $\dev$ on $\tilde U' -\tilde U''$. 
\index{developing map!radial-end projectivization|textbf} 


	
	We define a new developing map $\dev': \torb \ra \bX$ by using $\dev^N$ on 
	$\tilde U$ and letting it equal to 
	$\dev$ on the complement $\torb - \tilde U''$.  
Lemma \ref{op-lem-radialisotopy} implies that $\dev'$ can be obtained from $\dev$
by isotopies. 
	
	
	For other components of form $\gamma(\tilde U)$ for $\gamma \in \pi_1(\orb)$, 
	we do the same constructions. 
	
	

	\begin{lemma}\label{op-lem-replacedev} 
	Let $U$ and $U'$ be a radial end neighborhoods such that 
	$\clo_{\torb}(U) \subset U'$ and compatible to $\bar \orb$. 
	Let $\tilde U$ and $\tilde U'$ denote 
	the p-end neighborhoods of $\tilde E$ covering $U$ and $U'$. 
	We assume that $\Bd U$ and $\Bd U'$ are transverse to radial rays 
	by taking $U$ and $U'$ smaller if necessary. 
	Then we can modify the developing map in $\tilde U$ 
	such that the new developing map $\dev'$ agrees with 
	$\dev$ on $\torb - p_{\orb}(U)$ 
	such that $\dev'$ extends smoothly on the end compactification of $\tilde U$ 
	and $\dev'$ restricts to each radial line segment is a projective map. 
	Finally, $\dev' = \dev \circ \iota $ for an isotopy-lift $\iota$ preserving 
	each radial segment in $\tilde U'$. \hfill {$\square$}
\end{lemma}


Recall the definition of T-end structure from Section \ref{intro-sub-totgeo}. 
\begin{lemma} \label{op-lem-Tendisotopy}
Let $\tilde E$ be a T-p-end. Let $\tilde U$ be a proper T-p-end neighborhood of $\tilde E$. 
Then $\dev_\mu$ is such that $\dev_\mu$ extends to the ideal boundary component 
$\tilde \Sigma_{\tilde E}$ as an immersion. 
\end{lemma} 
\begin{proof} 
For any point $x$ in $\tilde \Sigma_{\tilde E}$, we have a chart $(U, \phi)$ where $x\in U$.
$U - \tilde \Sigma_{\tilde E}$ is inside $\torb$. 
By the definition of a real projective structure as an atlas of compatible charts, 
we obtain $g \in \bG$ 
such that $g\circ \phi$ agrees with $\dev_\mu$ on an open set in $U  - \tilde \Sigma_{\tilde E}$, 
and hence must agree on $U -\tilde \Sigma_{\tilde E}$. 
Hence, $\dev_\mu$ extend to $U$ as well. By continuing with all points of 
$\tilde \Sigma_{\tilde E}$, we obtain the result.  

\end{proof} 


Here, $\overline{\dev}_\mu$ is also equivariant with respect to $h$ if
$\dev_\mu$ was so. We call $(\overline{\dev}_\mu, h)$
an extended developing pair. 

\begin{definition}\label{op-defn-isolift2} 
An {\em extended isotopy} is a diffeomorphism $\bar{\orb} \ra \bar{\orb}$
extending an isotopy of $\orb$. 
An {\em extended isotopy-lift} is an extension of 
an isotopy-lift $\hat{\orb} \ra \hat{\orb}$. 
\end{definition} 
Of course, these are isotopies of $\bar \orb$ and isotopies of 
$\hat \orb$ respectively. 


Then the {\em isotopy-equivalence space } $\widetilde{\Def}_{\mathcal{E}}(\orb)$ is defined 
as the space of extended developing maps $\overline{\dev}_\mu$  
of real projective structures on $\orb$ with ends with radial structures and  
lens-shaped totally geodesic ends with end structures
under the action of the group 
the extended isotopy-lfts where 
an extended isotopy-lift  
$\hat \iota: \hat \orb \ra \hat \orb$ 
acting by 
\[ (\overline{\dev}_\mu, h) \mapsto (\overline{\dev}_\mu \circ \hat\iota, h). \]
\index{isotopy-equivalence space|textbf} 

We explain the topology. 
The space $\mathcal{D}(\hat \orb)$ of maps of form 
$\overline{\dev}_{\mu'}:\hat \orb \ra \rpn$
is given the compact open $C^r$-topology on $\hat \orb$. 

Suppose that we are given two real projective structures 
$\mu$ and $\mu'$ with extendable developing maps 
as above. 
For T-p-ends, this is clear by taking local inverses. 
Suppose we are given 
an isotopy $\iota_{\mu, \mu'}$ such that 
$\iota_{\mu, \mu'}^\ast(\mu') = \mu$.
Then 
the end compactification $\bar \orb$ has an extended isotopy-lift 
$\hat \iota_{\mu, \mu'}: \hat \orb \ra \hat \orb$: 
Suppose we are in an R-p-end neighborhood $U$.
Let $\dev_\mu$ and $\dev_{\mu'}$ denote the corresponding developing maps. 
We give a local spherical coordinate system near $\dev_\mu(\mbv_{\tilde E})$
which we could assume to be $\dev_{\mu'}(\mbv_{\tilde E})$
by $[0, \infty) \times \SI^{n-1}$ where $0$ corresponds to the p-end vertex.
We give coordinates on the closure of $U$ in $\hat \orb$ 
by $[0, \infty) \times \tilde \Sigma_{\tilde E}$
Then $\dev_{\mu}$ can be written locally as $(r(x, t), F(x))$ for 
$t\in [0, \infty), x \in \tilde \Sigma_{\tilde E}$ in the spherical coordinates.
Here, $r$ and $F$ are some functions to $[0, \infty)$.
Also, $\dev_{\mu'}$ can be written locally as $(r'(x, t), F'(x))$ 
for $t\in [0, \infty), x\in \tilde \Sigma_{\tilde E}$. 
And $r', F'$ are again some locally defined functions in $U$ to $\bR$.
Hence, both developing maps on the closures of $U$ locally lifts to 
$\tilde \Sigma_{\tilde E} \times [0, \infty)$ as smooth maps. 
Now, from these expressions, this fact follows. 


Now, $\dev_{\mu'}\circ \iota_{\mu, \mu'}$ is a developing map 
of $\iota_{\mu, \mu'}^\ast(\mu')$ sending the end structures $\torb$ of $\mu'$ 
to radial line. 
Hence, $\overline{\dev}_{\mu'} \circ \hat \iota_{\mu, \mu'}$ is
the unique smooth extension.
Hence, we can reinterpret 
$\mathcal{D}(\hat \orb)$ as the space of extensions of 
developing maps of $\orb$ with a fixed end structure for each end. 

By Lemma \ref{op-lem-replacedev} and \ref{op-lem-Tendisotopy} and 
the above paragraphs, 
we can now give the following definition:  
\begin{definition}\label{op-defn-deforms} 
The quotient space $\mathcal{D(\hat \orb)}/\mathcal{E}(\hat \orb)$ 
of $\mathcal{D}$ under the group of extended isotopy-lifts 
$\mathcal{E}(\hat \orb)$ of form
$\hat \iota_{\mu, \mu'}:\hat \orb \ra \hat \orb$ is
in one-to-one correspondence with  $\widetilde{\Def}_{\mathcal{E}}(\orb)$. 
The topology on $\widetilde{\Def}_{\mathcal{E}}(\orb)$ 
is given as the quotient topology of this space, which is called 
a {\em $C^r$-topology}. 
\index{crtopology@$C^r$-topology|textbf} 

We define
$\Def_{\mathcal{E}}(\orb) := \widetilde{\Def}_{\mathcal{E}}(\orb)/\PGLnp$
by the action 
\[ (\overline{\dev}, h(\cdot)) \mapsto (\phi \circ \overline{\dev}, \phi\circ h(\cdot) \circ \phi^{-1}), \phi \in \PGLnp \]
as in \cite{dgorb} and \cite{Goldman88}. 
The induced quotient topology is called a {\em $C^r$-topology} of
$\Def_{\mathcal{E}}(\orb)$. 

We can define a map 
\[\hol':  \widetilde{\Def}_{\mathcal{E}}(\orb) \ra \Hom_{\mathcal{E}}(\pi_1(\orb), \PGL(n+1, \bR))\] 
by sending the class of $(\dev, h)$ to $h$. This is well-defined since the isotopies do not change $h$. 
There is an induced map: 
\[\hol: {\Def}_{\mathcal{E}}(\orb) \ra \Hom_{\mathcal{E}}(\pi_1(\orb), \PGL(n+1, \bR)/\PGL(n+1, \bR).\] 
\end{definition} 
For these, if $\orb$ is closed, we simply drop the subscripts. 
\index{deformation space|textbf} 
\index{deformation space!topology|textbf} 
\index{DefE@$\Def_{\mathcal E}(\mathcal{O})$|textbf}
 \index{hol@$\hol$} 
\index{Ehresmann-Thurston map} 
\index{Ehresmann-Thurston principle} 
 \index{hol'@$\hol'$} 


It is well-known:
\begin{theorem}[see Choi \cite{dgorb}]  \label{op-thm-dgorb} 
Let $\orb$ be a closed orbifold.  Then 
\[\hol: \widetilde{\Def}(\orb) \ra \Hom(\pi_1(\orb), \PGL(n+1, \bR))\] is a local homeomorphism. 
\end{theorem} 
This map is a so-called Ehresmann-Thurston map, and this theorem 
is often known as Ehresmann-Thurston principle. 
\index{Ehresmann-Thurston map} 
\index{Ehresmann-Thurston principle}

\section{The local homeomorphism theorems}\label{op-sec-loch}

\subsection{The end condition for real projective structures}\label{op-subsec-endreal}
Now, we go over to real projective orbifolds:
We are given a real projective orbifold $\mathcal{O}$ with ends $E_1, \dots, E_{e_1}$ of $\cR$-type and $E_{e_1+1}, \dots, E_{e_1 + e_2}$ 
of $\cT$-type. Let us choose representative p-ends $\tilde E_1, \dots, \tilde E_{e_1}$ and $\tilde E_{e_1+1}, \dots, \tilde E_{e_1 + e_2}$.

In the following $\bX = \rpn$ or $\SI^n$ while 
$\bG = \PGL(n+1, \bR)$ or $\SLpm$ respectively
and $\bX^\ast = \RP^{n\ast}$ or $\SI^{n\ast}$ respectively as well. 

We define a subspace of 
$\Hom_{\mathcal{E}}(\pi_1(\mathcal{O}), \bG)$ to be as in Section
\ref{intro-sub-semialg}. 

Let $\mathcal V$ be an open subset of semi-algebraic subset of
\[\Hom^s_{\mathcal{E}}(\pi_1(\mathcal{O}), \bG)\] invariant under the conjugation action
such that one can choose a continuous section $s_{\mathcal V}^{(1)}: \mathcal V \ra (\bX)^{e_1}$
sending a holonomy homomorphism to a common fixed point of $h(\pi_1(\tilde E_i))$ for $i = 1, \dots, e_1$ and 
 $s_{\mathcal V}^{(1)}$ satisfies 
 \[s_{\mathcal V}^{(1)}(g h(\cdot) g^{-1})  = g \cdot s_{\mathcal V}^{(1)}(h(\cdot)) \hbox{ for } g \in \bG.\] \index{section} 
There might be more than one choice of a section and the domain of definition. 
 $s_{\mathcal V}^{(1)}$ is said to be a {\em fixed-point section}.
 \index{fixed-point section}
 \index{section!fixed-point}
 \index{svone@$s_{\mathcal V}^{(1)}$ } 

Again suppose that one can choose a continuous section $s_{\mathcal V}^{(2)}: \mathcal V \ra (\bX)^{e_2}$
sending a holonomy homomorphism to a common dual fixed point of $\pi_1(\tilde E_i)$ for $i = e_1+1, \dots, e_2$, and 
 $s_{\mathcal V}^{(2)}$ satisfies 
 \[ s_{\mathcal V}^{(2)}(g h(\cdot) g^{-1}) 
 = (g^*)^{-1}\circ s_{\mathcal V}^{(2)}(h(\cdot)) \hbox{ for } g \in \bG).\]
There might be more than one choice of section in certain cases. 
 $s_{\mathcal V}^{(2)}$ is said to be a {\em dual fixed-point section}.
 \index{svtwo@$s_{\mathcal V}^{(2)}$ } 
 \index{section!fixed-point!dual} 
 \index{dual fixed-point section} 
 
 We define $s_{\mathcal V}: \mathcal V \ra (\bX)^{e_1} \times (\bX^\ast)^{e_2}$
 as $ s_{\mathcal V}^{(1)} \times s_{\mathcal V}^{(2)}$ and call it a {\em fixing section}
provided the p-end holonomy group of 
each $\cT$-type p-end $\tilde E_i$ acts on a horosphere tangent to $P$ determined by 
$s_{\mathcal V}^{(2)}$.
\index{fixing section} 
\index{section!fixing} 
\index{sv@$s_{\mathcal V}$ } 

Recall from Section \ref{intro-sub-open}.
We note that the real projective structure with radial and totally geodesic ends with end structures also determines 
a point of $(\bX)^{e_1}\times (\bX^\ast)^{e_2}$. 
Conversely, if the real projective structure with radial and totally geodesic ends 
has the end structure {\em determined by} a section 
$s_{\mathcal{U}}$ if the following hold: 
\begin{itemize} 
\item  $\tilde E_i$ for every $i=1, \dots, e_1$
has a p-end neighborhood with a radial foliation with leaves 
developing into rays ending at the fixed point of the $i$-th 
factor of $s_{\mathcal V}^{(1)}$. 
\item $\tilde E_i$ for every $i= e_1+1, \dots, e_1+ e_2$ 
\begin{itemize} 
	\item has a p-end neighborhood with the ideal boundary component in the hyperspace 
	determined by the $i$-th factor of $s_{\mathcal V}^{(2)}$ provided 
	$\tilde E_i$ is a T-end, or 
	\item has a p-end neighborhood containing a $\bGamma_{\tilde E}$-invariant 
	horosphere tangent to the hyperspace 
	determined by the $i$-th factor of $s_{\mathcal V}^{(2)}$ provided 
	$\tilde E_i$ is a horospherical end. 
	\end{itemize} 
\end{itemize} 

\begin{example}
If $\mathcal{O}$ is a real projective orbifold and has some singularity of dimension one in each end neighborhood of an $\cR$-type end, 
then the  projective universal cover of $\mathcal{O}$ has more than two  lines corresponding 
to singular loci. The developing image of the lines must meet at a point in $\bX$, 
which is a common fixed point of the holonomy group of an end. 
If $\mathcal{O}$ has dimension $3$, this is equivalent to requiring that the end orbifold has corner-reflectors or cone-points. 

Hence, for an open subspace $\mathcal{V}$ of
a semi-algebraic subset of 
 \[\Hom^s_{\mathcal{E}}(\pi_1(\mathcal{O}), \bG)\] 
corresponding to the real projective structures on $\orb$, 
there is a section $s^{(1)}_{\mathcal{V}}$ determined by the common fixed points. 


\end{example} 

\begin{remark}[Cooper]
We do caution the readers that these assumptions are not trivial and exclude some important representations.
For example, these spaces exclude some incomplete hyperbolic structures arising in Thurston's Dehn surgery constructions 
as they have at least two fixed points for the holonomy homomorphism of the fundamental group of a toroidal end as was pointed out by Cooper. 
Hence, the uniqueness condition fails for this class of examples. 
However, if we choose a section on a subset, then we can obtain appropriate results.
Or if we work with particular types of orbifolds, the uniqueness holds.  
See Section \ref{ex-sec-nicecase}.
\end{remark}


\subsection{Perturbing horospherical ends}  
\label{op-sub-perturb}

Theorems \ref{app-thm-qFuch} and \ref{app-thm-qFuch2} study
the perturbation of lens-shaped R-ends and lens-shaped T-ends.  

The following concerns the deformations of $\bGamma_{\tilde E}\ra \bG$ near 
horospherical representations. As long as we restrict to deformed representations
satisfying the lens-condition, there exist  $n$-dimensional properly convex domains
on which the groups act. (This answers a question of Tillmann near 2006.
We also benefited from a discussion with J. Porti in 2011.) 

Let $P$ be an oriented hyperspace of $\SI^n$ with a dual point 
$P^* \in \SI^{n\ast}$ represented by a $1$-form $w_P$ defined on 
$\bR^{n+1}$. 
Let $P^\dagger$ denote the space of oriented hyperspaces in $P$. 
Let $\SI^{n-1 \ast}_{P^*}$ be the space of rays from $P^*$ corresponding to 
hyperspaces in $P$. Then the subspace $P^\dagger$ is dual to $\SI^{n-1 \ast}_{P^*}$: 
each oriented ray in $\SI^{n\ast}$ from $P^*$ defines a hyperspace $S'$ of $P$ as the set of common zeros of the $1$-forms in the ray.
The orientation of $S'$ is given by the open half-space where the $1$-forms near $w_P$ are positive.  
Conversely, an oriented pencil of oriented hyperspaces determined by an 
oriented hyperspace of $P$ is a ray in $\SI^{n-1 \ast}_{P^*}$ from $P^*$. 
(We omit the obvious $\RP^n$-version.)



Let 
\[ \Hom_{\mathcal{E}, \lh, p}(\Gamma', \PGLnp)\:
\left(\hbox{resp. } \Hom_{\mathcal{E}, \lh, p}(\Gamma', \SLnp)\right)\] denote the space of representations $h$ fixing a common fixed point $p$ and 
acts properly and cocompactly on the lens of a lens cone over 
vertex $p$ or is horospherical with a horoball with vertex $p$.

Let \[\Hom_{\mathcal{E}, \lh, P}(\Gamma', \PGLnp)\: \left(
\hbox{resp. } \Hom_{\mathcal{E}, \lh, P}(\Gamma', \SLnp)\right)\] 
denote the space of representations 
where $h(\Gamma')$ for each element $h$ acts on a hyperspace $P$ satisfying the lens-condition.
(See \ref{intro-sub-semialg}.) 


Let a convex cone $B= \partial B \ast \{p\}$ over a point $p$ 
be diffeomorphic to $\partial B \times (0, 1]$. 
Then $B$ with a vertex $p$ has a radial foliation. 
We complete $B$ by identifying with $\partial B \times (0, 1)$
by a diffeomorphism $f$ sending each leaf to $x \times (0, 1)$
and attaching $\partial B \times (0, 1]$ by $f$. 
We denote the partial completion by $\hat B$ diffeomorphic
to $\partial B \times [0, 1]$. 
We call $\hat B$ the {\em p-end completion} of $B$. 
An action of a group $\Gamma$ on $B$ 
extends to $\hat B$ also. 
$\hat B/\Gamma$ is then the end-compactification of $B/\Gamma$.
(See Remark \ref{op-rem-compact}.)

\begin{lemma}[Horospherical-end perturbation] \label{op-lem-horob} $ $
\begin{description}
\item[(A)] Let $B$ be a horoball in $\RP^n$ {\rm (}resp. in $\SI^n${\rm )} and $\Gamma$ be a group of projective automorphisms fixing $p$, 
$p \in \Bd B$ {\rm (}resp. $p \in \SI^n${\rm )},
such that $B/\Gamma$ is a horospherical-end-type orbifold.  
Then there exists a sufficiently small neighborhood $K$ 
of the inclusion homomorphism $h_0$ of $\Gamma$ in 
$\Hom_{\mathcal{E}, \lh, p}(\Gamma, \PGLnp)$ \hfill \break 
{\rm (resp. }$\Hom_{\mathcal{E}, \lh, p}(\Gamma, \SLnp)${\rm )}
where 
\begin{itemize} 
	\item for each $h \in K$, 
	$h(\Gamma')$ acts on a properly convex domain $B_h$ such that $B_h/h(\Gamma')$ 
	is diffeomorphic to $B/\Gamma'$ forming a radial end and fixes  $p$, 
\item $B_h$ forms the  lens-shaped or horospherical R-p-end neighborhood,
\item there is a diffeomorphism 
$f_h: B/\Gamma'  \ra B_h/h(\Gamma'), h\in K$, 
such that 
the lift $\tilde f_h: B \ra B_h$ 
is a continuous family under the $C^r$-topology \index{crtopology@$C^r$-topology}
as a map into $\RP^n$ {\rm (}resp. in $\SI^n${\rm )}
where $\tilde f_{h_0}$ is the identity map. 
\end{itemize} 
Let $\hat B$ and $\hat B_h$ denote the p-end compactifications. Then
$f_h$ extends to the end compactifications 
$\bar f_h: \hat B/\Gamma' \ra \hat B_h/h(\Gamma')$ and 
$\bar f_{h_0}$ is the identity map. 
Furthermore, the lift of this map 
$\hat f_h: \hat B \ra \RP^n$ {\rm (}resp. $\SI^n${\rm )} 
is continuous in the $C^r$-topology
where $\hat f_{h_0}$ is the identity map. 
\item[(B)]  
Let $P$ be a hyperspace in $\RP^n$ {\rm (}resp. in $\SI^n${\rm ).}
Let $\Gamma'$ denote a projective automorphism group 
acting on $P$ and a horoball $B$ tangent to $P$ such that 
$B/\Gamma'$ is a horospherical-end-type orbifold. 
Then there exists a sufficiently small neighborhood $K$ 
of the inclusion homomorphism $h_0$ of $\Gamma'$ in 
$\Hom_{\mathcal{E}, \lh, P}(\Gamma', \PGLnp)$
{\rm (}resp. $\Hom_{\mathcal{E}, \lh, P}(\Gamma', \SLnp)$\/{\rm )} where
\begin{itemize} 
\item  for each $h \in K$,
$h(\Gamma')$ acts on a properly convex domain $B_h$ such that $B_h/h(\Gamma')$ 
is homeomorphic to $B/\Gamma'$,
\item $B_h$ forms a lens-shaped T-p-end or horospherical p-end neighorhood, and 
\item there is a diffeomorphism
$f_h: B/\Gamma' \ra B_h/h(\Gamma')$, $h\in K$,  
such that the lift $\tilde f_h: B \ra B_h$ 
is a continuous family under the $C^r$-topology
where $\tilde f_{h_0}$ is the identity map. 
\end{itemize} 
Let $\hat B$ denote the p-end compactification of $B$ and 
$\hat B_h$ denote 
$B_h$ union with the p-end ideal boundary component of $B_h$ when 
$h$ acts properly and cocompactly on a lens. 
Then $f_h$ extends to the end compactifications 
$\bar f_h: \hat B/\Gamma' \ra \hat B_h/h(\Gamma')$.
Furthermore, the lift of this map 
$\hat f_h: \hat B \ra \RP^n$ {\rm (}resp. $\SI^n${\rm )} 
is a continuous family in the $C^r$-topology
where $\hat f_{h_0}$ is the identity map. 
\end{description} 

\end{lemma} 
\begin{proof} 
We will prove for the $\SI^n$-version. 

(A)  Let us choose a larger horoball $B'$ in $B$ where 
 $B'/\Gamma'$ has a boundary component $S'_{\tilde E}$ so 
$B'/\Gamma$ is diffeomorphic to $S'_{\tilde E} \times [0, 1)$. 
$S'_{\tilde E}$ is strictly convex and transverse to the radial foliation. 
There exists a neighborhood $O_1$ in $\Hom_{\mathcal{E}, \lh, p}(\Gamma', \SLnp)$ corresponding to 
the connection on a fixed compact neighborhood $N$ of $S'_{\tilde E}$ changes only by $\eps$ in the $C^r$-topology, $r \geq 2$, on a compact set containing a compact fundamental domain. 
(See the deformation theorem in \cite{Goldman88} which generalize to the compact orbifolds with boundary.)
 
Let $h\in O_1$.
The universal cover $\tilde S'_{\tilde E}$ is a strictly convex 
codimension-one manifold,
and it deforms to $\tilde S'_{\tilde E, h}$ that is still strictly 
convex for sufficiently small $\eps$. 
Here, $\tilde S'_{\tilde E, h}$ may not be embedded in $\SI^n$ a priori but is a submanifold of the deformed $n$-manifold $N_h$ from $N$
by the change of connections. 
Every ray from $p$ meets $\tilde S'_{\tilde E, h}$ transversely also by the $C^r$-condition.

Let $\orar{v}_{x, h}$ be a vector in the direction of $x$ for $x \in \tilde S'_{\tilde E, h}$
which we choose equivariantly with respect to the action of 
$h(\Gamma')$. We may choose such that 
$(x, h) \mapsto \orar{v}_{x, h}$ is continuous. 
We form a cone 
\[c(\tilde S'_{\tilde E, h}) := \{ [t \orar{v}_p + (1-t) \orar{v}_{x, h}]| t\in [0, 1], x \in \tilde S'_{\tilde E, h} \}. \]
Let $\tilde \Sigma_{\tilde E, h}$ denote the space of rays from $p$ ending at $\tilde S'_{\tilde E, h}$ in $c(\tilde S'_{\tilde E, h})$. 
Here $S'_h := \tilde S'_{\tilde E, h}/h(\Gamma')$ is a compact real projective orbifold of $(n-1)$-dimension.

Since $\Gamma'$ is a cusp group, it is virtually abelian.
Since $h \in \Hom_{\mathcal{E}, \lh}(\Gamma', \SLnp)$, 
Lemma \ref{abelian-lem-brickn} implies that 
$D_h: \tilde S'_{\tilde E, h} \ra \SI^{n-1}_p$ is an embedding to a properly convex  domain or a complete affine domain $\Omega_h$ in $\SI^{n-1}_p$ 
where $h(\Gamma')$ acts properly discontinuously
and cocompactly when $h$ is not the inclusion map. 

There is a one-to-one correspondence from $\tilde S'_{\tilde E, h}$ 
to $\tilde \Sigma_{\tilde E, h}=:\Omega_h$.
By convexity of $\tilde \Sigma_{\tilde E, h}$, the tube domain
$\mathcal{T}_{p}(\Omega_h)$ with vertices $p, -p$ is convex. 
$\tilde S'_{\tilde E, h}$ meets each great segment in the interior of 
the tube domain with vertices at $p, -p$ at a unique point 
transversely since $h$ is in $O_1$ for sufficiently small $\eps$. 
The strict convexity of $\tilde S'_{\tilde E, h}$ implies 
that $B_h$ is convex by Lemma \ref{prelim-lem-locconv}. 
The proper convexity of $B_h$ follows 
since $\tilde S'_{\tilde E, h}$ is strictly convex
and meets each great segment from $p$ in the interior of the tube domain corresponding 
to $\tilde \Sigma_{\tilde E, h}$, 
and hence $\clo(B_h)$ cannot contain a pair of antipodal points. 


By Theorem \ref{ce-thm-affinehoro}, the Zariski closure of 
$h(\Gamma')$ is a cusp group $G_h$ extended by a finite group
and $G_h/h(\Gamma')$ is compact. 
Hence, $\Gamma'$ is virtually abelian by the Bieberbach theorem. 
We take the identity component $\mathcal{N}_h$ of $G_h$, which 
is an abelian group with a uniform lattice 
$h(\Gamma')$. 
The set of orbits of $\mathcal{N}_h$ foliates $B_h$. 
Since $\mathcal{N}_h$ is a normal subgroup of $G_h$, 
$h(\Gamma')$ normalizes $\mathcal{N}_h$. 
Hence, the orbits give us a codimension-one foliation on $B_h/h(\Gamma')$
with compact leaves. The leaves are all diffeomorphic, and 
hence, we obtain a parametrization 
$\partial B/\Gamma' \times (0, 1]$ to $B_h$. 

Now, $h$ induces isomorphism $\hat h_h: \mathcal{N}_0 \ra 
\mathcal{N}_h$ where $\hat h_h \ra \Idd$ as $h \ra h_0:= \Idd$. 

We choose a proper radial path $\alpha_h:I \ra B_h$ 
from a point of $\partial B_h$ and ending at $p$. 
We may assume that $\alpha_h$ is independent of $h$. 
We define a parameterization 
\[\tilde \phi_h: \mathcal{N} \times (0, 1] \ra   B_h, 
(m, t) \mapsto \hat h_h(m)(\alpha_h(t)), t \in (0, 1]. \]
We define $\tilde f_h: B \ra B_h:= \tilde \phi_h \circ \tilde \phi_{h_0}^{-1}$. 
This gives us a map $f_h$. 
(Here, we might be changing $\tilde \Sigma_{\tilde E, h}$.)
Since $\tilde f_h$ sends radial segments to radial segments, 
it extends to a smooth map $\hat f_h: \hat B \ra \hat B_h$. Also, on any compact 
subset $J$ of $\hat B$, a compact foliated set $\hat J$ contains it. 
Let $\hat J_h$ denote the image of $\hat J$ under $\hat f_h$. 
$\hat J_h$ is coordinatized by $\che J\times I$ for a fixed 
compact set $\che J \subset \mathcal{N}$. 
Under these coordinates of $\hat J$ and $\hat J_h$, we can write 
$\hat f_h$ as the identity map. 
Since $\hat h_h \ra \hat h_{h_0}$, we conclude that 
$\hat f_h|\hat J$ uniformly converges to $\Idd$ as $h \ra h_0$. 

We may assume that these are radial-end projectivization
of some developing map by carefully choosing $\alpha_h$. 




Now, we show that $B_h$ is a R-p-end of horospherical or lens type. 
We know from above that $\Omega_h$, $h \in O_1$, is either properly convex or complete affine. 
Suppose that $\Omega_h$ is properly convex. Then we have a tubular action on a tube 
corresponding to $\Omega_h$. 
By Theorem \ref{pr-thm-secondmain}, the lens condition is equivalent to the uniform middle-eigenvalue condition. Hence, we have a distanced action by Proposition \ref{pr-thm-distanced}
and  $\tilde S'_{\tilde E, h}$ must have the boundary in a distanced compact set in the boundary of 
the tube by Proposition \ref{pr-prop-orbit}. By taking the convex hull of $\tilde S'_{\tilde E, h}$, 
we obtain a compact convex set distanced from $p$. Now, looking at from $p_-$, we can 
obtain a smooth lens containing this. Hence, we have a R-p-end of lens type. 

Suppose that $\Omega_h$ is complete affine. Then again $\tilde S'_{\tilde E, h}$ is strictly convex and develops to a complete affine space in $\SI^n_p$. 
Then we have a horospherical end by the premise. 


(B) The second item is the dual of the first one.  If $h(\Gamma')$ acts on 
an open horosphere $B^o$ tangent to $P$ with the vertex in $P$ properly discontinuously, 
then the dual group $h(\Gamma')^*$ acts on a horosphere with a vertex the point $P^*$ dual to $P$. 
By duality Proposition \ref{pr-prop-dualend}, 
$h^\ast$ is in  $\Hom_{\mathcal{E}, \lh, P^\ast}(\Gamma, \SLnp)$. 
We apply the first part, and 
hence, there exists a properly convex domain 
$\Omega_{P, h}$ such that $\Omega_{P, h}/h^\ast(\Gamma')$ is 
an open orbifold for $h \in K$ for some subset $K$ of 
the character space
\[\Hom_{\mathcal{E}, \lh, P^\ast}(\Gamma', \SLnp).\]

By duality Proposition \ref{pr-prop-dualend}, 
$\Omega_{P, h}$ is foliated by radial lines from 
$P^\ast$ and $R_{P^\ast}(\Omega_{P, h}) \subset \SI_{P^\ast}^{n-1}$ is 
a properly convex domain.

Let $B_h^o$ be the properly convex domain dual to $\Omega_{P, h}^o$ and hence is 
a properly convex domain, and $B_h^o/h(\Gamma')$ is a dual orbifold
diffeomorphic to $\Omega_{P, h}/h(\Gamma')^\ast$ by Theorem 
\ref{prelim-thm-dualdiff}.  We have $B_{h_0}= B^o$ since the dual of the dual of
a properly convex open domain is itself.

Let $\Gamma$ act properly on $\partial B$. 
Since $\partial B$ is strictly convex, each point of $\partial B$ 
has a unique hyperspace sharply supporting $B$. 
By Proposition \ref{prelim-prop-duality}, 
there is an open hypersurface $S, S\subset \Bd \Omega_P$ 
dual to $\partial B$. Also, $\Gamma'^\ast$ acts properly 
on $S$ such that $S/\Gamma'^\ast$ is a closed orbifold. 
(The closedness again follows since there is a torsion-free 
subgroup of finite index and hence a finite regular-covering manifold.)
From the first part, there is an open surface 
$S_h, S_h \subset \Bd \Omega_{P, h}$, $h \in K$, meeting 
each radial ray from $P^\ast$ at a unique point.
Also, $S_h/h(\Gamma')^\ast$ is diffeomorphic to $\partial B/\Gamma'^\ast$. 

Again, by Proposition \ref{prelim-prop-duality}, 
we obtain an open surface $S_h^\ast$, 
$S_h^\ast \subset \Bd \Omega_{P, h}^\ast$
 where $\Gamma_h$ acts properly
such that $S_h^\ast/\Gamma_h$ is a closed orbifold. 
We define $B_h:= \Omega_{P, h}^\ast \cup S_h^\ast$. 

Since $\Omega_{P, h}$ is foliated by radial segments from $P^\ast$
with properly convex 
\[R_{P^\ast}(\Omega_{P, h}) \subset \SI_{P^\ast}^{n-1},\]
$D_h:=P\cap \Bd \Omega_{P, h}^\ast$ and 
is a properly convex domain in $P$ by Proposition \ref{pr-prop-dualend}.




Note that $\Gamma'$ is virtually abelian, and when it is not a cusp group, then 
it is lens-type and hence must be virtually diagonalizable. 

Define 
$\hat B_h:=  \Omega_{P, h}^\ast \cup S_h^\ast \cup D_h$, $h \in K$. 
For the second and third items of the second part, 
%
%
we do as above 
but we choose $\alpha_h:I \ra B_h$ to be a single geodesic segment 
starting from 
$x_0\in \partial B_h$ and ending at a point of $D_h$
where as $h \ra h_0$, 
$\alpha_h$ converges as a parameter of functions to a geodesic 
$\alpha_0:I \ra \clo(B) \subset \SI_{\infty}^n$ 
ending at the vertex of the horosphere or a fixed point of $h$ not on $P$ whenever 
$h$ is virtually diagonalizable.
We assume that $\alpha_h$ is a $C^r$-family of geodesics. 
Now, the proof is similar to the above using 
an isomorphism from the identity component $C_h$ of the Zariski closure of 
$\Gamma'$ to that of $h(\Gamma')$, which is an abelian group
since $\Gamma'$ is virtually abelian. 
We denote by $\kappa: C_{h_0} \ra C_h$ the unique homomorphism extending 
$h\circ h_0^{-1}$ on restricted to the abelian subgroup of finite index of $\Gamma'$. 

Here, we need the images of $\alpha_h$ under 
$C_h$ to form a foliation. Since $C_h$ acts on a properly convex set $D_h$, 
it acts as a diagonalizable group on $P$ by Proposition \ref{prelim-prop-Ben2}. 
Being a free abelian group satisfying the uniform middle-eigenvalue condition, 
$C_h$ is a diagonalizable group acting on an $n$-simplex. 
(We dualize the situation and use Theorem \ref{pr-thm-redtot}.) 
We required that the extension of $\alpha_h$ to pass 
a fixed point of $C_h$ not in $\clo(D_h)$. 
The images of $g\circ \alpha_h$, $g \in C_h$, form a foliation of  
$\hat B_h$. 
Using this we define the map $\tilde f_h: B \ra B_h$ sending 
leaves to leaves as given by the function 
\[ g(\alpha_0(t)) \mapsto 
\kappa(g)(\alpha_h(t)) \hbox{ for each } t\in [0, 1], g \in C_{h_0}. \] 
This map extends from $\partial B$ to $\partial B_h$.

Finally notice that our constructions of $f_h$ 
all are smooth from $\bar \orb$. Hence, these are compatible
end neighborhoods. 

We have a lens by reflection about $P$ by the fixed point not on $P$ 
when $\Gamma'$ is virtually diagonalizable. Hence, we have a lens-type end. 
Otherwise, we have a horospherical end as in the end of the proof of the first part.
\hfill \SnT {\parfillskip0pt\par}
\end{proof}

\begin{lemma}[Lens-end perturbation] \label{op-lem-lensp} $ $
\begin{description}
\item[(A)] Let $B$ be a generalized lens cone in $\RP^n$ {\rm (}resp. in $\SI^n${\rm )} and $\Gamma$ be a group of projective automorphisms fixing $p$, 
$p \in \Bd B$ {\rm (}resp. $p \in \SI^n${\rm )},
such that $B/\Gamma$ is a generalized lens-end-type orbifold.  
Then there exists a sufficiently small neighborhood $K$ 
of the inclusion homomorphism $h_0$ of $\Gamma$ in 
$\Hom_{\mathcal{E}, \lh, p}(\Gamma, \PGLnp)$ \hfill \break 
{\rm (}resp. $\Hom_{\mathcal{E}, \lh, p}(\Gamma, \SLnp)${\rm )}
where 
\begin{itemize} 
	\item for each $h \in K$, 
	$h(\Gamma')$ acts on a properly convex domain $B_h$ such that $B_h/h(\Gamma')$ 
	is diffeomorphic to $B/\Gamma'$ forming a radial end and fixes  $p$, 
\item $B_h$ forms the  lens-shaped R-p-end neighborhood, and
\item there is a diffeomorphism 
$f_h: B/\Gamma'  \ra B_h/h(\Gamma'), h\in K$, 
such that 
the lift $\tilde f_h: B \ra B_h$ 
is a continuous family under the $C^r$-topology \index{crtopology@$C^r$-topology}
as a map into $\RP^n$ {\rm (}resp. in $\SI^n${\rm )}
where $\tilde f_{h_0}$ is the identity map. 
\end{itemize} 
Let $\hat B$ and $\hat B_h$ denote the p-end compactifications. Then
$f_h$ extends to the end compactifications 
$\bar f_h: \hat B/\Gamma' \ra \hat B_h/h(\Gamma')$ and 
$\bar f_{h_0}$ is the identity map. 
Furthermore, the lift of this map 
$\hat f_h: \hat B \ra \RP^n$ {\rm (}resp. $\SI^n${\rm )} 
is continuous in the $C^r$-topology
where $\hat f_{h_0}$ is the identity map. 
\item[(B)]  
Let $P$ be a hyperspace in $\RP^n$ {\rm (}resp. in $\SI^n${\rm ).}
Let $\Gamma'$ denote a projective automorphism group 
acting on $P$ and a lens $B$ meeting $P$ in its interior such that 
$B/\Gamma'$ s a lens-end-type orbifold, and a component $B_1$ of $B-P$ is 
a p-end neighborhood of an end. 
Then there exists a sufficiently small neighborhood $K$ 
of the inclusion homomorphism $h_0$ of $\Gamma'$ in 
$\Hom_{\mathcal{E}, \lh, P}(\Gamma', \PGLnp)$ \hfill \break
{\rm (}resp. $\Hom_{\mathcal{E}, \lh, P}(\Gamma', \SLnp)$\/{\rm )} where
\begin{itemize} 
\item  for each $h \in K$,
$h(\Gamma')$ acts on a lens $B_h$ such that $B_h/h(\Gamma')$ 
is homeomorphic to $B/\Gamma'$,
\item A component $B_h -P$ forms lens-shaped T-p-end 
or horospherical p-end neighborhood, and 
\item there is a diffeomorphism
$f_h: B/\Gamma' \ra B_h/h(\Gamma')$, $h\in K$,  
such that the lift $\tilde f_h: B \ra B_h$ 
is a continuous family under the $C^r$-topology
where $\tilde f_{h_0}$ is the identity map. 
\end{itemize} 
Let $\hat B_1$ denote the p-end compactification of $B$ and 
$\hat B_{1, h}$ denote 
$B_h -P$ union with the p-end ideal boundary component of $B_h$ when 
$h$ acts properly and cocompactly on a lens. 
Then $f_h| B_1$ extends to the end compactifications 
$\bar f_h: \hat B_1/\Gamma' \ra \hat B_{1. h}/h(\Gamma')$.
Furthermore, the lift of this map 
$\hat f_h: \hat B_1 \ra \RP^n$ {\rm (}resp. $\SI^n${\rm )} 
is a continuous family in the $C^r$-topology
where $\hat f_{h_0}$ is the identity map. 
\end{description} 

\end{lemma} 
\begin{proof} 
The proofs are very similar to those of Lemma \ref{op-lem-horob} using radial segments and duality.
Here, we need the local homeomorphism property of for the closed orbifolds of 
Choi \cite{dgorb} applied to the end orbifolds.  (See Lok \cite{Lok} also.) 
Hence, for nearby homeomorphisms we can choose developing maps that are very close near 
the hypersurfaces bounding the p-end neighborhoods using the radial-end projectivizations. 

The $\SI^n$-version is obtained by lifting the developing maps to $\SI^n$.
\end{proof} 

We remark that 
we can also reinterpret the parameterization as radial-end projectivization in
Section \ref{op-sub-deform} by taking a second larger end neighborhood
and some modifications of the parameters.

\subsection{Local homeomorphism theorems}

Again $\Def^s_{\mathcal{E}, s_{\mathcal U}}(\mathcal{O})$ is defined to be the subspace of $\Def_{\mathcal{E}}(\orb)$ 
with the stable irreducible holonomy homomorphisms in $\mathcal U$ 
and the end determined by $s_{\mathcal U}$, i.e., 
\begin{itemize} 
\item each $\cR$-type p-end has a p-end neighborhood 
foliated by geodesic leaves that are radial to 
the vector given by $s_{\mathcal U}$ under the developing map, or 
\item each $\cT$-type p-end is totally geodesic of lens-type satisfying the lens-condition or horospherical satisfying the suspended lens condition with respect to the hyperspace determined by $s_{\mathcal U}$. (See Section \ref{intro-sec-ends}.) 
\end{itemize} 

One of the main aim of the chapter is:
\begin{theorem}\label{op-thm-projective} 
Let $\mathcal{O}$ be a noncompact strongly tame real projective $n$-orbifold with lens-shaped radial ends or lens-shaped totally geodesic ends with types assigned. 
In the following $\bX = \rpn$ or $\SI^n$ while 
$\bG = \PGL(n+1, \bR)$ or $\SLpm$ respectively
and $\bX^\ast = \RP^{n\ast}$ or $\SI^{n\ast}$ respectively as well. 
Let
$\mathcal V$ be a conjugation-invariant open subset of the semialgebraic subset \[\Hom^s_{\mathcal{E}}(\pi_1(\mathcal{O}), \bG).\] 
Let $s_{\mathcal V}$ be the fixing section defined on $\mathcal V$ with images in 
$(\bX)^{e_1} \times (\bX^\ast)^{e_2}$. 
Then the map \[\hol:\Def^s_{\mathcal{E}, s_{\mathcal V}}(\mathcal{O}) \ra \rep^s_{\mathcal{E}}(\pi_1(\mathcal{O}), \bG)\]
sending the real projective structures with ends compatible with $s_{\mathcal V}$ to their conjugacy classes of holonomy homomorphisms
is a local homeomorphism  to an open subset of $\mathcal V'$.
\end{theorem}

\begin{corollary} \label{op-cor-A} 
Let $\mathcal{O}$ be a noncompact strongly tame real projective $n$-orbifold with lens-shaped radial ends or lens-shaped totally-geodesic ends with 
end structures and given types $\cR$ or $\cT$. 
Assume $\partial \orb = \emp$. 
Then the following map is a local homeomorphism\,{\rm :}  
\[\hol:\Def^s_{\mathcal{E}, u}(\mathcal{O}) \ra \rep^s_{\mathcal{E}, u}(\pi_1(O), \bG).\]
\end{corollary}  \index{hol@$\hol$} \index{Ehresmann-Thurston map} 
\begin{proof} 
It is clear that the limit of a sequence of unique fixed points of representations 
must be a fixed point of the limit representation. Also, the limit of corresponding cyclic spaces 
must be in the cyclic spaces of the limit representations. 
Hence, if the limit representations have elements with one-dimensional intersections of 
the cyclic spaces, then the corresponding elements of the sequence have one-dimensional intersections of the cyclic spaces. This proves that the map from the representation spaces to the unique fixed points is continuous. Now Theorem \ref{op-thm-projective} implies the conclusion.
\end{proof}

\subsection{The proof of Theorem \ref{op-thm-projective}.} \label{op-sub-proof}


We wish to now prove Theorem \ref{op-thm-projective} following the proof of Theorem 1 in Section 5 of \cite{dgorb}.

\begin{proof}[Proof of Theorem \ref{op-thm-projective}]
Let $\mathcal{O}$ be an affine orbifold with the affine universal covering orbifold $\tilde{\mathcal{O}}$ 
with the affine covering map $p_\orb:\tilde{\mathcal{O}} \ra \mathcal{O}$
and let the fundamental group $\pi_1(\mathcal{O})$ act on it as an automorphism group.

Let $\mathcal U$ and $s_{\mathcal U}$ be as above. 
We now define a map 
\[\hol:\widetilde{\Def}_{{\mathcal{E}}, s_{\mathcal U}}(\mathcal{O}) \ra \Hom_{\mathcal{E}}(\pi_1(\mathcal{O}), \SL_\pm(n+1, \bR))\] 
by sending the projective structure to the pair $(\dev, h)$ and to the conjugacy class of $h$ finally.


We show that $\hol$ is continuous: 
There is a codimension-$0$ compact submanifold $\mathcal{O}'$ of $\mathcal{O}$ 
such that $\pi_1(\mathcal{O}') \ra \pi_1(\mathcal{O})$ is an isomorphism.
The holonomy homomorphism is determined on $\mathcal{O}'$. 
Since the deformation space has the $C^r$-topology, $r \geq 1$,  induced by 
$\dev:\tilde \orb' \ra \bR^{n+1}$, it follows that small changes of $\dev$ on compact domains in $\tilde \orb'$ in the $C^r$-topology imply sufficiently small changes in $h(g'_i)$ for 
generators $g'_i$ of $\pi_1(\mathcal{O}')$ and 
hence sufficiently small change of $h(g_i)$ for generators $g_i$ of $\pi_1(\mathcal{O})$. 
Therefore, $\hol$ is continuous.
(Actually for the continuity, we do not need any condition on ends.)




\label{op-page-thm-affine}
Next, we define the local inverse map from a neighborhood in 
$\mathcal U$ of the image point.
Let $\orb'$ be a compact suborbifold 
of $\orb$ such that $\orb - \orb'$ is a union $U$ of end neighborhoods. 


We show how to change the proof of Theorem 1 of \cite{dgorb}.
Let $h$ be a representation coming from an affine orbifold $\orb$. 
The task is to reassemble $\orb$ with new holonomy homomorphisms as we vary $h$ as in \cite{dgorb} following the approaches of Thurston. 
Suppose that $h'$ is in a neighborhood of $h$ in 
$\Hom_{\mathcal{E}}(\pi_1(\mathcal{O}), \SL_\pm(n+1, \bR))$. 
\begin{itemize} 
\item As in Lok \cite{Lok}, we consider locally finite collections $\mathcal{V}$ 
of open domains that cover $\orb$. 
We find subcollection $\mathcal{V}'$ 
of compact neighborhoods or end neighborhoods in the covering contractible open sets
which covers $\orb$ again. 
Here, each precompact element of $\mathcal{V}$ contains a compact 
domain in $\mathcal{V}'$ forming a cover of $\orb'$. 
The end neighborhoods can meet 
only if they are in the same component of $U$.  
\item We give orders to the open sets covering $\orb$. 
The end neighborhoods have orders higher than all precompact sets. 
\item We regard these as sets in $\SI^{n+1}$ by charts. 
\item We consider the sets that are the intersection \[U_{i_1}\cap \cdots \cap U_{i_k}, 
i_1 > \cdots > i_k, \hbox{ where } U_{i_j}\in \mathcal{V}' \hbox{ for } j=1, \dots, k,\] 
of the largest cardinality of the compact or closed domains in $\mathcal{V}'$
and find the corresponding sets in $\SI^{n+1}$ by charts.
We map it by isotopies to the corresponding intersection of 
deformed collections of domains in $\SI^{n+1}$ 
corresponding to the $h'(\pi_1(\orb))$-action by Lemmas 3 and 4 in \cite{dgorb} and 
using the deformations of $\dev$ by 
Lemmas \ref{op-lem-horob} and \ref{op-lem-lensp}. 
Here, we will follow 
the ordering as above when we deform as in Lok \cite{Lok}. 
That is, we use the isotopy of $U_{i_1}$ restricted to $U_{i_1} \cap \cdots \cap U_{i_k}$ when $U_{i_1}$ has the largest order. 
\item  We extend the isotopies to the sets of intersections of 
smaller number of sets in $\mathcal{V}'$ by Lemma 5 of \cite{dgorb}. 
By induction, we extend it to all the images of compact 
and closed domains in $\mathcal{V}'$. 
\item We patch these open sets to 
build an orbifold $\orb_{h'}$ 
with holonomy $h'$ referring back to $\orb$ by isotopies. 
\item $\orb_{h'}$ is diffeomorphic to $\orb$ by 
the map constructed by the isotopies. 
\end{itemize}

To show that the local inverse is a continuous map 
for the $C^r$-topology of $\hat \orb$, 
we only need to consider compact suborbifolds in $\orb$
since the holonomy representation depends only on any compact submanifold whose complement is a union of proper end neighborhoods. For this, the same argument as in \cite{dgorb} applies.


We now prove the local injectivity of $\hol$.
Given two structures $\mu_0$ and $\mu_1$ in a neighborhood of the deformation space, we show that if their 
holonomy homomorphisms are the same, say $h: \pi_1(\orb) \ra \SL_\pm(n+1, \bR)$, then we can isotopy one in the neighborhood to the other using vector fields as in \cite{dgorb}.

Because of the section $s_{\mathcal U}$ defined on $\mathcal U$, given a holonomy $h$, 
we have a direction of the radial end that is unique for the holonomy homomorphism. 

First assume that $\orb$ has only $\cR$-type ends. 
Recall the compact suborbifold $\orb'$ such that $\orb - \orb'$ is diffeomorphic to 
$E_i \times (0, 1)$ for each end orbifold $E_i$ where each $x \times (0, 1)$ is the image 
of a radial segment. 

We can choose a Riemannian metric on $\bar{\orb}$ such that an 
end neighborhood has a product metric of form $E_i \times (0, 1]$.
Let $\dev_j$ be the developing map of $\mu_j$ for $j=0,1$. 
Then the $C^r$-norm distance of 
extensions $\overline{\dev_0}$ and $\overline{\dev_1}$ to $\bar{\orb}$
is bounded on each compact set $K \subset \bar{\orb}$
by our assumption on the closeness of the two structures. 
Since we chose $\mu_1$ and $\mu_2$ sufficiently close, 
$\overline{\dev_0}$ and $\overline{\dev_1}$ can 
be assumed to be sufficiently close in the $C^r$-topology over $K$. 
The images of $K$ under each of these 
maps can be assumed to lie on a neighborhood of 
the image of a p-end vertex, say $\mbv$. 
Moreover, the radial lines maps to the radial lines ending at the same ideal end vertex. 
We may assume that they have the forms of the radial-end projectivizations. 
We can use the argument in the last part of \cite{dgorb} to show that 
$\overline{\dev}_1$ lifts to an immersion $\hat \orb \ra \hat \orb$ 
equivariant with respect to the deck transformation group. 
Then we use the metrics to equivariantly isotopy it to $\Idd$ 
as in the last section of \cite{dgorb}. 
Hence, $\mu_0$ and $\mu_1$ represent the same point of 
$\widetilde{\Def}_{\mathds{A}, {\mathcal E}, s_{\mathcal U}}^S(\mathcal{O})$. 


Suppose now that $\orb$ has some lens type $\cT$-type ends. 
Suppose that $\mu_0$ and $\mu_1$ have a totally geodesic ideal boundary
component corresponding to an end of $\orb$.
We attach the totally geodesic ideal boundary component for each end, 
and then we can argue as in \cite{dgorb} proving the local injectivity. 

Suppose that $\mu_0$ and $\mu_1$ have horospherical end neighborhoods corresponding to an end of $\orb$.
Then these are radial ends and the same argument as the above one for $\cR$-type ends applies 
to show the local injectivity. 

Finally, we cannot have the situation that $\mu_0$ has the totally geodesic ideal boundary component corresponding 
to an end while $\mu_1$ has a horoball end neighborhood for the same end. 
This follows since the end holonomy group acts on a properly convex domain in a totally geodesic hyperspace
and as such the end holonomy group elements have some norms of eigenvalues $> 1$. (See Proposition 1.1 of \cite{Benoist00} for example.) 
\hfill \SnP {\parfillskip0pt\par}
\end{proof}

\section{Relationship to the deformation spaces in our earlier papers}
\label{op-sec-relation} 


Recall $\mathfrak{D}(\hat P)$ for a Coxeter orbifold $\hat P$ that is 
not necessarily compact in Definition \ref{ex-defn-deformation}.

Let $\orb$ be a strongly tame real projective $n$-orbifolds only radial ends and radial end structure $\mathcal{E}$. 
Recall the definition of $\CDef_{\mathcal{E}}(\orb)$ from Section \ref{op-sub-deform}.

Generalizing this,  let $\mathfrak{D}_{\mathcal{E}}(\orb)$ denote the same set as 
$\CDef_{\mathcal{E}}(\orb)$. We give the topology by $C^r$-topology on 
for the set of all developing maps $\dev: \orb \ra \rpn$
and take the quotient by right actions by the isotopy lifts $\tilde \orb$  
and by the left action by composition with elements of $\PGL(n+1, \bR)$. 
 These were the topology we used before 
as in \cite{Choi06}. 

By following Proposition \ref{op-prop-DPCDef}, we obtain that 
$\mathfrak{D}(\hat P)$ is the same as $\CDef_{{\mathcal E}}(\hat P)$
as topological spaces. 

We say that a diffeomorphism {\em preserves a radial end structure} 
if an end neighborhood with the radial foliation structure contains an end neighborhood
 is mapped into an end neighborhood sending the radial leaves to radial leaves. 

	Let $\tilde \orb$ denote the universal cover of $\orb$. 
	We can form $\bar{\orb}$ by using the end-compactification
	by Remark \ref{op-rem-compact}. 
A radial end structure preserving isotopies can be extended to 
a homeomorphism of end-compactifications. 
We call such diffeomorphism {\em end-structure preserving extended isotopies}. 


\begin{proposition} \label{op-prop-DPCDef}
Let $\orb$ be a stronly tame $n$-orbifold with only radial ends and radial end structure $\mathcal{E}$. 
	Then there is a homeomorphism 
	\[ \mathfrak{D}_{\mathcal{E}}(\orb) \cong \CDef_{{\mathcal E}}(\orb).   \]
	\end{proposition} 
\begin{proof} 
		We denote by $\hat{\orb}$ the universal cover 
		of $\bar{\orb}$ containing $\tilde \orb$ as a dense open set. 
		

	
	Hence, by restricting the structure on $\orb$ only, 
	there is a map 
	\[R'_r: \CDef_{{\mathcal E}}(\orb) \ra  \mathfrak{D}(\orb).\]
	The map is one-to-one and onto since the sets are the same. 
	
Let $D_{\tilde \orb}$ denote the space of all developing maps on $\tilde \orb$ with radial end structures. 
	Let $J^r(f)$ denote the tuples of all jets of $f: \tilde \orb \ra \RP^n$ 
	of order $\leq r$. 
	The topology on $D_{\tilde \orb}$ is given 
	by 	bases of form 
	\[B_{K,\epsilon}(\dev):=\{ f\in D_\orb |\bdd(J^r(\dev)(x), J^r(f)(x)) < \epsilon, x \in K, f \in C^r(\tilde \orb, \RP^n) \}\]
	where $K \subset \tilde \orb$ is a compact set, $\epsilon > 0$,  
	and $\dev \in D_{\tilde \orb}$. 
	The topology on $D_{\hat \orb}$ is 
	given by bases of form 
		\[B'_{K, \epsilon}(\dev):=\{ f\in D_{\hat \orb} |\bdd(J^r(\dev)(x), J^r(f)(x)) < \epsilon, x \in K, f \in C^r(\hat{\orb}, \RP^n) \}\] 
		where  $K \subset \hat \orb$ is a compact set, $\epsilon > 0$, 
		and $\dev \in D_{\hat \orb}$. 
		
	Consider the restriction map 
	\[R_r:C^r(\hat{\orb}, \RP^n) \ra C^r(\tilde{\orb}, \RP^n)\]
	inducing $R'_r$. 
	Since compact subsets of $\tilde \orb$ are compact subsets of $\bar{\orb}$, 
we can find a 
	the inverse image under $R_r$ of a basis element of  $C^r(\tilde{\orb}, \RP^n)$
	is a basis element in $C^r(\bar{\orb}, \RP^n)$. 
	Hence,  the induced map 
\[R'_r:\CDef_{{\mathcal E}}(\orb) =D_{\bar  \orb}/\mathcal{G}_{\bar \orb} \ra  \mathfrak{D}(\orb) = D_{\orb}/\mathcal{G}_\orb\] 
	is continuous. 
	
	Now we find the inverse of $R'_r$: 

Recalling Definition \ref{op-defn-isolift}, we define: 
\begin{itemize} 
\item 	Let $\mathcal{G}_\orb$ denote the group of isotopies of $\orb$ 
	preserving the radial end structures of $\orb$,  and 
\item  let $\mathcal{G}_{\bar \orb}$ denote the group of radial end structure preserving 
	extended isotopies of $\bar \orb$.
\item 	The isotopy-lifts 	$\iota: \tilde \orb \ra \tilde \orb$
form a group which we 
	denote by $\mathcal{G}_{\tilde \orb}$.
\item Denote by $\mathcal{G}_{\hat \orb}$ the group of the extensions of 
isotopy-lifts $\hat \iota: \hat \orb \ra \hat \orb$ of isotopies of $\bar\orb$. 
\end{itemize} 

Fix the union $U$ of mutually disjoint R-end neighborhoods and radial foliations on each component
in $\orb$. 

Also, we choose a union $U'$ of 
such neighborhoods such that $\clo_{\tilde \orb}(U) \subset U'$.    
Let $\tilde U$ denote the inverse image of $U$ and $\tilde U'$ that of $U'$. 
\begin{itemize} 
\item denote by $\mathcal{G}_{\orb, U', U}$ the group of isotopies of form $\iota$ of $\orb$ 
acting on each component of $U$  preserving the radial folations and $U'$. 

\item Let $\mathcal{G}_{\tilde \orb, U', U}$ denote the isotopy-lifts $\tilde \iota$ 
of $\orb$ to $\tilde \orb$
acting on each component of $\tilde U'$ and $\tilde U$ which are lifts of elements of 
$\mathcal{G}_{\orb, U', U}$. 

\item Let 
$\mathcal{G}_{\hat \orb, U', U}$ be the smooth extensions of
isotopy-lifts of $\tilde \orb$ to $\hat \orb$ 
acting on each component of $\tilde U'$ and $\tilde U$.
(These are isotopies of $\hat \orb$. See Definition \ref{op-defn-isolift2}.) 
\end{itemize} 
Clearly there is a natural homomorphism by restriction
$\mathcal{G}_{\hat \orb, U', U} \ra \mathcal{G}_{\orb, U', U}.$ 

Let $D_{\orb, U, U'}$ denote the space of functions of form 
$\dev$ in  a development pair $(\dev, h)$
such that $\dev| \tilde U'$ equals $\dev^N|\tilde U'$ constructed for $U$ and $U'$
by the radial-end projectivization 
in Lemma \ref{op-lem-replacedev}. 
(We may have to construct on a larger union of p-end neighborhood because of
the smoothing process and restrict to $U$ and $U'$.)
By Lemma \ref{op-lem-replacedev}, 
we obtain a natural map $D_{\orb, U, U'}/\mathcal{G}_{\orb, U', U} \ra 
D_{\orb}/\mathcal{G}_{\orb}$ that is a one-to-one onto map. 

Also, denote by 
$D_{\bar \orb, U, U'}$ denote the space of functions of form 
$\overline{\dev}$ which are extended developing maps 
and $\overline{\dev}|\tilde U$ equals $\dev^N| \tilde U$ constructed for $U$ and $U'$ 
by radial-end projectivization
in Lemma \ref{op-lem-replacedev}.  


Recall from the elementary analysis that the $C^r$-topology on a compact set $K$ is  
a metric topology with the metric given by taking the supremum of the distances of jets up to order $\leq r$ at a compact set. 

	We claim that  the canonical map 
	  $F: D_{\orb, U, U'}/\mathcal{G}_{\orb, U',U} \ra 
	  D_{\orb}/\mathcal{G}_\orb$ is a homeomorphism: 
	  Let us take a compact neighborhood $K_F$ of a fundamental domain 
	  $F$ of $\torb - \tilde U$. 
	  Then we define a metric on 
	  $D_{\orb, U, U'}$ by $d_{D_{\orb, U, U'}}$ between 
	  two developing maps $f_1, f_2$ is defined 
	  as $\sup_{x\in K_F} \bdd(J_r(f_1)(x), J_r(f_2)(x))$. 
	  Since the developing map on $\tilde U$ is determined by 
	  the developing map restricted on $K_F$ 
	  as we can see by a radial-end projectivization with respect to $U$ and $U'$, 
	  we obtain a metric $d_{D_{\orb, U, U'}}$ on $D_{\orb, U, U'}$ 
	  giving us the $C^r$-topology on $D_{\orb, U, U'}$. 
	  The map $F$ is continuous since it is induced by 
	  the inclusion map $D_{\orb, U', U} \ra D_\orb$.

	  There is an inverse map 
	  $G: D_\orb/\mathcal{G} \ra D_{\orb, U, U'}/\mathcal{G}_{U', U}$
	  given by taking a developing map $f$ and modifying it 
	  by radial-end projectivization. There are choices 
	  involved, but they are well-defined up to isotopy-lifts. 
	  
	  Any isotopy $\iota$ induces a homeomorphism $\iota^\ast$ in 
	  $D_\orb$. 
	  To show the continuity of $G$, we take 
	  a ball $B_{d_{\orb, U, U'}}(f, \epsilon)$ for $f \in D_{\orb, U, U'}$ 
	  to show that there is a neighborhood of $f\circ \iota$ in $D_\orb$ in 
	  the $C^r$-topology going into it for $\iota \in \mathcal{G}^{\tilde \orb}$. 
	  Since we can take $\iota^{\ast}(B)$ as a neighborhood for $f\circ \iota$, 
	  it is sufficient to find a neighborhood $B$ of $f$ in $D_\orb$ in the $C^r$-topology.
	  Let $K_F$ denote a compact set in $\tilde \orb - \tilde U$.  
	  We take a sufficiently small $\delta$, $\delta > 0$, 
	  the neighborhood $B_{K_F, \delta}(f)$ such that each of 
	  its element $g$ 
	  is still in  
	   $B_{d_{D_{\orb, U, U'}}}(f, \epsilon)$
	   after applying the isotopy of 
	   radial-end projectivization with respect to $U$ and $U'$: 
	   This is because 
	   we only need to worry about the compact set 
	   $K_F \cap \clo_{\tilde \orb}(\tilde U')$ while $g$ after  
	   radial-end projectivization with respect to $U'$ and $U$
	   is determined in this compact set. 
	   If $g$ is sufficiently
	   close to $f$ on $K_F$ in the $C^r$-topology 
	   already in the form of Lemma \ref{op-lem-replacedev}, 
	   then 
	   \begin{itemize} 
	  \item the leaves of the radial foliation of $g$ intersected with $K_F$ 
	   is $C^r$-close to ones for $f$ and 
	  \item the required isotopy $\iota_g$ for radial-end projectivization of 
	  $g$ is also sufficiently $C^r$-close to $\Idd$ on $K_F$ in the uniform $J_r$-topology defined on $U$. 
	  
	  \end{itemize} 
	   Hence, by taking sufficiently small 
	   $\delta$, $g\circ \iota_g$ is in $B_{d_{D_{\orb, U, U'}}}(f, \epsilon)$. 
	   Since we are estimating everything in a compact set $K_F$, 
	   and finding $\iota_g$ depending on $g| K_F$, these methods 
	   are possible. 
	   (By the uniform topology, we mean the topology using 
	   the norm of differences of two functions on $U$ but not just one compact subset of $U$.)

	  
	 We can also show that $D_{\bar \orb, U, U'}/\mathcal{G}_{\bar \orb, U', U} 
	  \ra D_{\bar \orb}/\mathcal{G}_{\bar \orb}$
	  is a homeomorphism similarly.  This is entirely analogous to above. 
	  
It is easy to show that the induced map
$D_{\bar \orb, U, U'}/\mathcal{G}_{\bar \orb, U', U} 
\ra D_{\orb, U, U'}/\mathcal{G}_{\orb, U', U}$ is a homeomorphism. 
This is because we just need to work with the compact orbifold $\orb - U'$
and isotopying with $U$.  
Since there is a commutative diagram 
\begin{alignat}{3} 
D_{\bar \orb, U, U'}/\mathcal{G}_{\bar \orb, U', U} &
\ra  \: & D_{\orb, U, U'}/\mathcal{G}_{\orb, U', U} \nonumber \\
\downarrow \qquad &     & \downarrow \qquad \nonumber \\ 
D_{\bar \orb}/\mathcal{G}_{\bar \orb}  & \: \ra  & D_{\orb}/\mathcal{G}_\orb.
\end{alignat}  	  
Since the downarrows maps are homeomorphisms by above, and 
the upper row map is a homeomorphism, it follows that  
the bottom row is continuous. 

This means that there is a continuous
map $D_{\orb}/\mathcal{G}_\orb \ra D_{\bar \orb}/\mathcal{G}_{\bar \orb}$
giving us the inverse of $R'_r$.  
	\end{proof}

\chapter{Relative hyperbolicity and strict convexity} \label{ch-rh}

We show the equivalence between the relative hyperbolicity of the fundamental groups of the properly convex real projective orbifolds with lens-shaped radial ends, totally geodesic ends, or horospherical ends
with the strict convexity of the orbifolds relative to the ends. 
In Section \ref{rh-sec-constr}, we show how to add some lenses to the T-ends
so that the ends become boundary components and how to remove some open sets to make the R-ends into boundary components. Some constructions preserve strict convexity, and so on. 
In Section \ref{rh-sec-SPC}, we describe the action of the end fundamental group 
on the boundary of the universal cover of a properly convex real projective orbifold. 
In Section \ref{rh-sec-bowditch},  
we prove Theorem \ref{rh-thm-relhyp} that the strict convexity implies the relative hyperbolicity of 
the fundamental group using Yaman's work. We then present the converse of this
Theorem \ref{rh-thm-converse}. 
For this, we use the work of Dru{\c{t}}u and Saphir on tree graded spaces and asymptotic cones.

\section{Some constructions associated with ends} \label{rh-sec-constr} 

We discuss some constructions to begin. 
It is sufficient to prove for the case $\torb \subset \SI^n$ in this chapter
by Theorem \ref{prelim-thm-lifting}
and Proposition \ref{prelim-prop-closureind}. 
Our main results here are Theorem \ref{rh-thm-relhyp}
and Corollary \ref{rh-cor-remhyp}, 
and the proofs fall into \SSn. 

The purpose of this chapter is to prove Corollary \ref{rh-cor-remhyp}, 
the equivalence of the strict convexity of $\orb$ and 
the relative hyperbolicity of $\pi_1(\orb)$ with respect to 
the end fundamental groups. Cooper-Long-Tillmann \cite{CLT15} and
Crampon-Marquis \cite{CM14} proved the same result when we only allowed horospherical ends. 
	Benoist told us at an IMS meeting at the National University of Singapore in 2016 that 
	he has a proof for this theorem for $n=3$ using trees 
	as he had done in closed $3$-dimensional cases in \cite{Benoist062}
	using the Morgan-Shalen's work on trees \cite{MS88}.
For convex cocompact actions, there are some later related works by
Islam and Zimmer \cite{IZp} and \cite{IZp2}, and Weisman \cite{Weisman}  
for relatively hyperbolic groups for convex cocompact orbifold groups. 
These are generalizations of the present work. The approach using Yaman's characterization was first explored by the author (\cite{convsurv},\cite{convMa}). 
However, the present proof simplifies the older versions. 

Recall that properly convex strongly tame real projective orbifolds with generalized lens-shaped or horospherical
ends satisfying (NA) and (IE) 
have strongly irreducible holonomy groups by Theorem \ref{intro-thm-sSPC}.


\subsection{Modifying the T-ends.} \label{rh-subsec-totdual}

From now on in this chapter, we assume $\torb \subset \SI^n$ as we stated the reason above. 

For T-ends, by the lens condition, we only consider the ones that have cocompactly-acted lens neighborhoods in some ambient orbifolds.
First, we discuss the extension to bounded orbifolds. 


\begin{theorem}\label{rh-thm-totgeoext} 
	Suppose that $\mathcal O$ is a strongly tame properly convex real projective $n$-orbifold 
	with generalized lens-shaped or horospherical $\cR$- or $\cT$-ends and satisfies {\em (IE)}. 
	Let $E$ be a lens-shaped T-end, and let $S_E$ be a totally geodesic hypersurface that is the ideal 
	boundary component corresponding to $E$. 
	Let $L$ be a lens-shaped end neighborhood of $S_E$ in an ambient real projective orbifold containing $\orb$. 
	Then 
	\begin{itemize} 
		\item $L \cup \orb$ is a properly convex real projective orbifold 
		and has a strictly convex boundary component 
		corresponding to $E$. 
		\item  Furthermore, if $\orb$ is strictly SPC and $\tilde E$ is a hyperbolic end, then so is $L \cup \orb$, which now has one more boundary component and one less T-ends. 
	\end{itemize} 
\end{theorem} 
\begin{proof} \renewcommand{\qedsymbol}{}
	Let $\torb$ be the universal cover of $\orb$, which we can identify with a properly convex bounded domain in 
	an affine subspace in $\SI^n$. Then $S_E$ corresponds to a T-p-end $\tilde E$ and 
	a totally geodesic hypersurface $S= \tilde S_{\tilde E}$. And $L$ is covered by 
	a lens $\tilde L$ containing $S$. The p-end fundamental group $\pi_1(\tilde E)$ acts on
	$\torb$ and $\tilde L_1$ and $\tilde L_2$ the two components of 
	$\tilde L - \tilde S_{\tilde E}$ in $\torb$ and
	outside $\torb$ respectively. 
\cpr
\end{proof}
	Lemma \ref{rh-lem-commsupp} generalizes Theorem 3.7 of \cite{Goldman90}. 
	\begin{lemma} \label{rh-lem-commsupp} 
		Suppose that $\tilde S_{\tilde E}$ is the totally geodesic ideal boundary component of 
		a lens-shaped T-end $\tilde E$ of a strongly tame real projective orbifold $\orb$.
        Let $\Omega_i$, $i=1, 2$, be a $\pi_1(\tilde E)$-invariant properly convex open domain that satisfies
            $\Bd \Omega_i \cap \SI^{n-1}_\infty = \tilde S_{\tilde E}$. 
Then the following hold\/{\rm :}
		\begin{itemize} 
			\item 
			For each point $p$ of $\Bd \tilde S_{\tilde E}$, and
			any sharply supporting hyperspace $H$ of $\tilde S_{\tilde E}$ at $p$ in $\SI^{n-1}_\infty$, 
			there exists a sharply supporting hyperspace
			to $\Omega_i$ containing $H$.
			\item At each point $p$ of  $\Bd \tilde S_{\tilde E}$, the hyperspace  sharply supporting $\Omega_i$ at $p$ is unique
			if $\pi_1(\tilde E)$ is hyperbolic. 
			\item Suppose that $\Omega_1$ and $\Omega_2$ are on the other sides of 
            $\tilde S_{\tilde E}$. 
			Then $\clo(\Omega_1) \cup \clo(\Omega_2)$ is a convex domain
			with \[\clo(\Omega_1) \cap \clo(\Omega_2) = \clo(\tilde S_{\tilde E})
			\subset \Bd \Omega_1 \cap \Bd \Omega_2, \] 
			and their sharply supporting 
hyperspaces at each point of $\Bd \tilde S_{\tilde E}$ coincide.
	\end{itemize}  
	\end{lemma} 
	\begin{proof} 
		Let $\mathds{A}^n$ denote the affine subspace that is the complement in $\SI^n$ of the hyperspace containing $\tilde S_{\tilde E}$. 
		Because $\pi_1(\tilde E)$ acts properly and cocompactly 
		on a lens-shaped domain, 
		$h(\pi_1(\tilde E))$ satisfies the uniform middle-eigenvalue condition
with respect to $\tilde E$
by Theorem \ref{pr-thm-equ2}. 
		
		The domain $\Omega_1$ has an affine half-space $H(x)$ bounded by 
		an AS-hyperspace for each 
		$x \in \Bd \tilde S_{\tilde E}$ containing $\Omega_1$. 
		Here, $H(x)$ is uniquely determined by $\pi_1(\tilde E)$ and $x$ 
		and $H(x) \cap \SI^{n-1}_\infty$ 
		by Theorems \ref{du-thm-asymnice} and \ref{du-thm-asymniceII}.
		The respective AS-hyperspaces at each point of $\clo(\tilde S_{\tilde E}) - \tilde S_{\tilde E}$ to $\Omega_1$ and $\Omega_2$ 
        must agree by Lemmas \ref{du-lem-inde}
and \ref{du-thm-ASunique} because otherwise, there must be an AS-hyperspace for 
one of $\Omega_i$ that extends to one meeting the interior of $\Omega_j$ for $j \ne i$. 
This proves the first and third items.
		
		
		The second item follows from the third item and Theorem 1.1 of \cite{Benoist04}.
%
%
\end{proof} 

\label{page-thm-totgeoext}
	
	\begin{proof}[Continuation of the proof of Theorem \ref{rh-thm-totgeoext}]
	By Lemma \ref{rh-lem-commsupp}, $\tilde L_2 \cup \tilde S_{\tilde E} \cup \torb$ is a convex domain by the existence of sharply supporting hyperspaces. 
	If $\tilde L_2 \cup \torb$ is not properly convex, then it is a union of 
	two cones over $\tilde S_{\tilde E}$ 
	of $[\pm v_x ] \in \bR^{n+1}, [v_x] = x$. 
	This means that $\torb$ has to be a cone contradicting the strong 
	irreducibility of $h(\pi_1(\orb))$. 
	Hence, it follows that $\tilde L_2 \cup \torb$ is properly convex. 
	
	Suppose that $\orb$ is strictly SPC and $\pi_1(\tilde E)$ is hyperbolic. 
	Then every segment in $\Bd \torb$ or a non-$C^1$-point in $\Bd \torb$
	is in the closure of one of the p-end neighborhoods.
	$\Bd \tilde L_2 - \clo(\tilde S_{\tilde E})$ does not contain any segment in it or not a $C^1$-point. 
	$\Bd \torb - \clo(\tilde S_{\tilde E})$ does not contain any segment or a non-$C^1$-point outside 
	the union of the closures of p-end neighborhoods. 
	$\Bd(\torb \cup \tilde L_2 \cup \tilde S_{\tilde E})$ is $C^1$
	at each point of 
	$\Lambda(\tilde E) := \clo(\tilde S_{\tilde E}) - \tilde S_{\tilde E}$
	by the uniqueness of the sharply supporting hyperspaces of Theorem
\ref{du-thm-ASunique}
	
	Recall that $\tilde S_{\tilde E}$ is strictly convex 
	since $\pi_1(\tilde E)$ is a hyperbolic group. (See Theorem 1.1 of \cite{Benoist04}.)
	Thus, $\Lambda$ does not contain a segment, and hence,
	$\Bd(\torb \cup \tilde L_2 \cup \tilde S_{\tilde E})$ does not contain one. 
	Therefore, $L_2 \cup \orb$ is strictly convex relative to the remaining ends.
	Now, we do this for every copy $g(L_2)$ of $L_2$ for $g\in \pi_1(\orb)$.

	Since $\tilde L_{2}\cup \torb$ has a Hilbert metric by \cite{Kobpaper}, the action is properly discontinuous. 
\end{proof}


\subsection{Shaving the R-ends} 

\begin{corollary} \label{rh-cor-lenssupp}
Suppose that $\mathcal O$ is a noncompact strongly tame properly convex real projective $n$-orbifold with a p-end $\tilde E$, 
and $\pi_1(\tilde E)$ is hyperbolic. 
\begin{enumerate} 
\item[(i)] Let $\tilde E$ be a lens-shaped totally geodesic p-end. 
Let $L$ be a cocompactly-acted lens containing a totally geodesic properly convex hypersurface 
$\tilde E$ such that 
\[\Lambda := \clo(\tilde S_{\tilde E}) - \tilde S_{\tilde E} 
= \Bd L - \partial L.\] Then 
each point of $\Lambda$ has a unique sharply supporting hyperspace of $L$.
and any domain $\Omega_1$ and $\torb$ as in Lemma \ref{rh-lem-commsupp}.  
\item[(ii)] 
Let $\tilde E$ be a lens-shaped radial p-end. 
Let $L$ be a cocompactly-acted lens in the p-end neighborhood.
Define $\Lambda := \Bd L - \partial L$.
Then each point of $\Lambda$ has a unique sharply supporting hyperspace of $L$. 
This holds for the universal cover $\torb$ as well.
\end{enumerate} 
\end{corollary}
\begin{proof} 
(i) is already proved in Lemma \ref{rh-lem-commsupp}
and Theorems \ref{du-thm-asymnice} and \ref{du-thm-asymniceII}. 

(ii) is proved in Proposition \ref{pr-cor-tangency}.
%
\end{proof}

	





We call the following construction {\em shaving the ends}. \index{end!shaving} 

\begin{theorem}\label{rh-thm-remconch} 
	Given an $n$-dimensional 
strongly tame SPC-orbifold $\mathcal{O}$ and its universal cover $\tilde{\mathcal O}$, 
	there exists a collection of mutually disjoint open concave p-end neighborhoods for lens-shaped p-ends.
	We remove a finite union of concave end neighborhoods of some R-ends. 
	Then 
	\begin{itemize}
		\item we obtain a convex domain as the universal 
		cover of a strongly tame orbifold ${\mathcal{O}}_1$ with additional strictly convex smooth boundary components
		that are closed $(n-1)$-dimensional orbifolds.
		\item Furthermore, if $\mathcal{O}$ is strictly SPC with respect 
		to all of its ends, and we remove only some 
		of the concave end neighborhoods of 
		hyperbolic R-ends, then ${\mathcal{O}}_1$ is strictly SPC with respect to the remaining ends. 
	\end{itemize} 
\end{theorem}
\begin{proof} 
We obtained $\mathcal{O}_1$ by removing some 
disjoint collection of concave open p-end neighborhoods. 
	If $\mathcal{O}_1$ is not convex, then there is a triangle $T$ in $\tilde{\mathcal{O}}_1$
	with three segments $s_0, s_1, s_2$ such that $T-s_0^o \subset \tilde{\mathcal{O}}_1$ 
	but $s_0^o - \tilde{\mathcal{O}}_1 \ne \emp$. (See Theorem A.2 of \cite{psconv} for details.)
	Since $\tilde{\mathcal{O}}_1$ is an open manifold, $s_0^o - \tilde{\mathcal{O}}_1$ is a closed subset of
	$s_0^o$. Then a boundary point $x \in s_0^o - \tilde{\mathcal{O}}_1$ 
	is in the boundary of one of the removed concave-open neighborhoods $U$ or is in $\Bd \tilde{\mathcal{O}}$ itself. 
	The second possibility implies that 
	$\mathcal{O}$ is not convex as $\torb_1 \subset \torb$. 
	The first possibility implies that
	there exists an open segment meeting $\Bd U\cap \torb$ at a unique point but 
	disjoint from $U$.
	This is geometrically not possible since $\Bd U \cap \torb$ is strictly convex 
	towards the direction of $U$. 
	These are contradictions. 
	
	Since $\torb$ is properly convex, 
	so is $\torb_1$. 
	Since $\Bd U \cap \torb$ is strictly convex, the new corresponding boundary 
	component of $\torb_1$ is 
	strictly convex. 
	
	
	Now we go to the second part. 
	Suppose that $\mathcal{O}$ is strictly SPC. 
	Let $\mathcal{H}$ denote the set of p-ends with hyperbolic p-end fundamental groups whose concave p-end neighborhoods 
	were removed in the equivariant manner. For each $\tilde E \in \mathcal{H}$, denote by $U_{\tilde E}$ the concave p-end neighborhood that we are removing. 
	
	Any segment in the boundary of the developing image of $\mathcal O$ is in
	the closure of a p-end neighborhood of a p-end vertex. 
	For the p-end-vertex $\mbv_{\tilde E}$ of a p-end $\tilde E$, the domain $R_{\mbv_{\tilde E}}(\torb) \subset \SI^{n-1}_{\mbv_{\tilde E}}$  is strictly convex  by \cite{Benoist04}
	if $\pi_1(\tilde E)$ is hyperbolic. 
	Since $\Bd R_{\mbv_{\tilde E}}(\torb)$ contains no straight segment, 
	only the straight segments in $\clo(U) \cap \Bd \torb$ for the concave p-end neighborhood  $U$ of $\tilde E$
	are in the segments in $\bigcup S({\mbv_{\tilde E}})$. 
	Thus, their interiors are disjoint from $\Bd \tilde{\mathcal{O}}_1$, and hence 
	$\Bd \torb_1$ contains no geodesic segment in $\bigcup_{\tilde E \in {\mathcal{H}}} \clo(U_{\tilde E}) \cap \Bd \torb$. 
	
	Since we removed concave end neighborhoods of the lens-shaped ends with the hyperbolic end fundamental groups, any straight segment in 
	$\Bd \tilde{\mathcal{O}}_1$ lies in the closure of a p-end neighborhood of a remaining p-end vertex. 
	
	A non-$C^1$-point of $\Bd \torb_1$ is neither on the boundary of the concave p-end neighborhood $U_{\tilde E}$ for a hyperbolic 
	p-end $\tilde E$ nor in $\Bd \torb - \bigcup_{\tilde E \in \mathcal{H}} \clo(U_{\tilde E})$.
	We show that points of $\Lambda = L - \partial L$ are $C^1$-points of 
	$\Bd\torb$: 
	$\clo(U_{\tilde E}) \cap \Bd \torb_1$ contains the limit set  $\Lambda = L - \partial L$
	for the cocompactly-acted lens $L$ in a lens-neighborhood. 
	$\torb$ has the same set of sharply supporting hyperspaces as $L$ at points of $\Lambda$
	since they are both $\pi_1(\tilde E)$-invariant convex domains by Corollary \ref{rh-cor-lenssupp}.
	However, the sharply supporting hyperspaces at $\Lambda$ of $L$ are also supporting ones for $\torb_1$ by Corollary \ref{rh-cor-lenssupp}
	since $L \subset \torb_1$ as we removed the outside component $U$ of $\torb - L$. Thus, $\torb_1$ is $C^1$ at points of $\Lambda$. 
	
	Also, points of $\Bd \torb - \bigcup_{\tilde E\in \mathcal{H}} \clo(U_{\tilde E})$ are $C^1$-points of 
	$\Bd \torb$ since $\orb$ is strictly SPC with respect to all of its ends. 
	Let $x$ be a point of this set. 
	Suppose that $x \in s^o$ for a segment $s$ in $\Bd \torb_1$. 
	Then $s\subset \Bd \torb$ and $s$ is not in $\clo(U_{\tilde E})$ for any ${\mbv_{\tilde E}}$ since we removed subsets of $\torb$ to obtain $\torb_1$. 
	Hence, this is not possible. 
	Suppose that $x$ has more than two sharply supporting hyperspaces $P_1, P_2$ 
	to $\torb_1$ at $x$. 
	We may assume that $P_1$ is a sharply supporting hyperspace to $\torb$. 
	Since $P_2$ is not supporting $\torb$, a component $H_2'$ disjoint from 
	$\torb_1$ meets $\torb$. Then $H'_2 \cap \torb$ is a convex domain, which we denote by 
	$\Omega_2$. $\Omega_2 \subset \bigcup_{\tilde E \in \mathcal{H}} \clo(U_{\tilde E})$. 
	Now, it is easy to see that $\clo(U_{\tilde E})$ for at most one p-end $\tilde E$ meets $H'_2$. 
	Since $x \in \torb_2$, $x$ is in the closure of $\clo(U_{\tilde E}) \cap H'_2$. 
	Thus, $x \in \clo(U_{\tilde E}) \cap \clo(\torb_1) \subset \Lambda_{\mbv_{\tilde E}}$ for 
	the limit set $\Lambda_{\mbv_{\tilde E}}$. However, we proved that there is a unique supporting 
	hyperspace at $x$ to $\torb_1$ in the above paragraph. 
	Hence, $\mathcal{O}_1$ is strictly SPC.     
\end{proof}

\section{The strict SPC-structures and the relative hyperbolicity} 
\label{rh-sec-SPC}


\subsection{The Hilbert metric on $\mathcal{O}$.} \label{rh-sub-Hilbert} 

Recall the Hilbert metrics from Section \ref{prelim-sub-Hilbert}. 
A Hilbert metric on an orbifold with an SPC-structure is defined as a distance metric 
given by cross ratios. (We do not assume strictness here.)

Given an SPC-structure on ${\mathcal{O}}$, there is a Hilbert metric which we denote by $d_{\torb}$ 
\index{do@$d_{\torb}$}
on $\tilde{\mathcal{O}}$ and hence on $\tilde {\mathcal{O}}$. 
Actually, we make $\orb$ slightly small by inward perturbations of $\partial \orb$ preserving the strict convexity of $\partial \orb$
by Lemma \ref{prelim-lem-pushing}. 
The Hilbert metric is defined on unchanged $\torb$. (We call this metric the {\em perturbed Hilbert metric}.)
This induces a metric on ${\mathcal{O}}$, including the boundary now. 
We denote the metric by $d_{\orb}$. 
\index{do@$d_{\orb}$}





Given an open properly convex domain $\Omega$, 
we note that given any two points $x, y$ in $\Omega$, there is a geodesic arc \index{Hilbert metric}
$\ovl{xy}$ with endpoints $x, y$ such that its interior is in $\Omega$. 

\begin{figure}[h]
	\centering 
	\includegraphics[height=7cm]{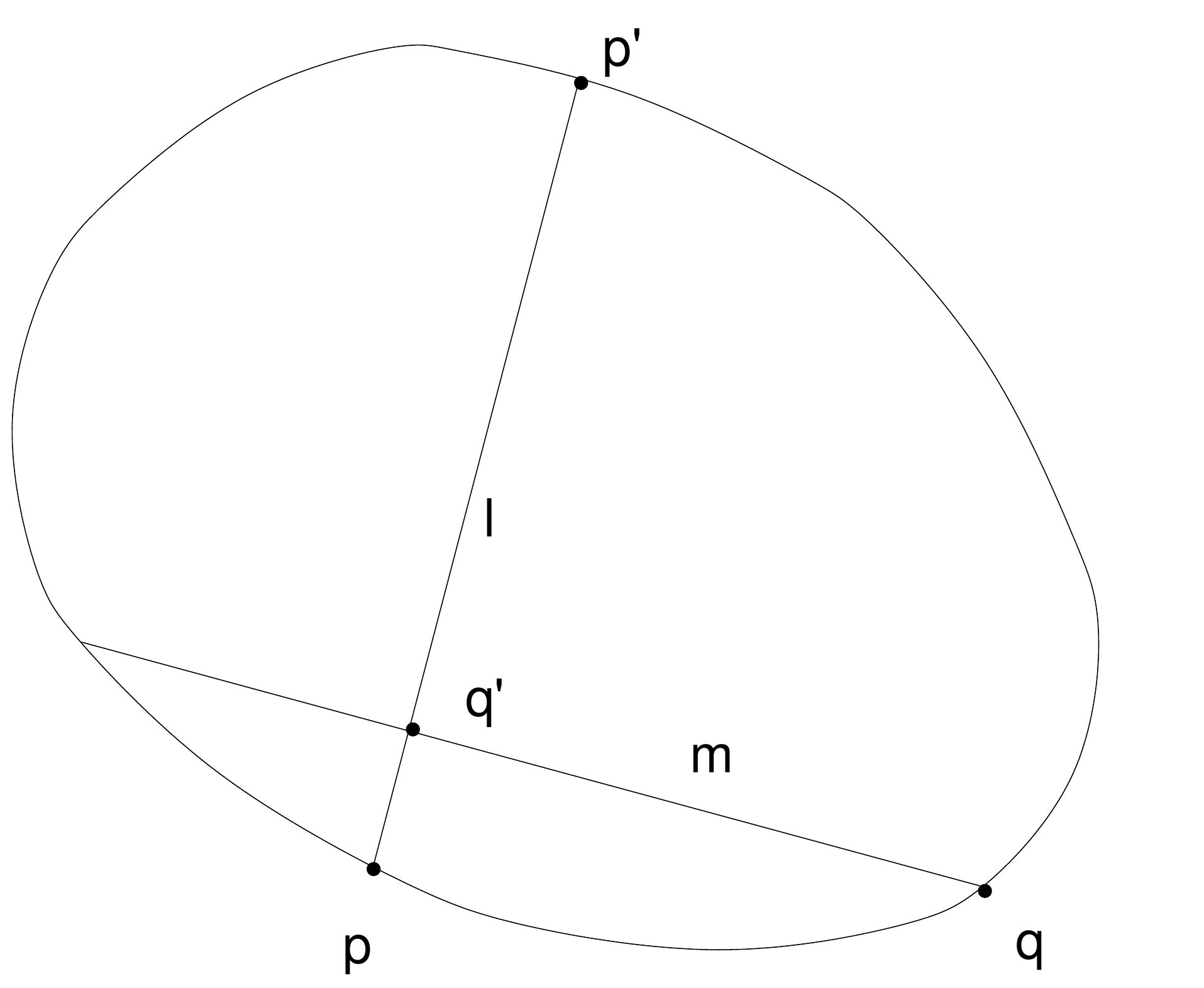}
	
	\caption{The shortest geodesic $m$ to a geodesic $l$.}
	\label{rh-fig:figbow}
\end{figure}

\begin{proposition}\label{rh-prop-shortg} 
Let $\Omega$ be a properly convex open domain. 
Let $P$ be a subspace meeting $\Omega$,
and let $x$ be a point of $\Omega  - P$\,{\rm :}
\begin{description}
\item[(i)] There exists a shortest path $m$ from $x$ to $P \cap \Omega$ that is a line segment. 
\item[(ii)] The set of shortest paths to $P$ from a point $x$ of $\Omega -P$ 
have endpoints in a compact convex subset $K$ of $P\cap \Omega$. 
\item[(iii)] For any line $m'$ containing $m$ and $y \in m'$, the segment in $m'$ from $y$ to the point of $P \cap \Omega$ is 
one of the shortest segments.
\item[(iv)] When $P$ is a complete geodesic in $\Omega$ with $x \in \Omega-P$, 
outside the compact set $K$, $K \subset P$, of endpoints of shortest segments from $x$ to $P$, the distance function from $P-K$ to $x$ is strictly increasing or strictly decreasing. 
\end{description}
\end{proposition}
\begin{proof} 
(i) The distance function $f: P \cap \Omega \ra \bR$ defined by $f(y) = d_\Omega(x,y)$ is 
a proper function where $f(x) \ra \infty$ as $x \ra z$ for any boundary point $z$ of $P \cap \Omega$ in $P$. 
Hence, there exists a shortest segment with an endpoint $x_0$ in $P \cap \Omega$. (iv) is also proved. 

(ii) 
Let $\gamma$ be any geodesic in $P \cap \Omega$ passing $x_0$. 
We need to consider the $2$-dimensional subspace $Q$ containing $\gamma$ and $x$.
The set of endpoints of shortest segments of $\Omega$ in $Q$ is a connected compact subset that contains $x_0$ by 
Proposition 1.4. of \cite{marg}.
Hence, by considering all geodesics in $P \cap \Omega$ passing $x_0$, 
we obtain that the endpoints of the shortest path to $P$ from $x$ 
form a connected compact set. We take two points $z_1, z_2$ on it. 
Then the segment connecting $z_1$ and $z_2$ is also in the set of 
endpoints by Proposition 1.4 of \cite{marg}. Hence, the set is convex. 

(iii) Suppose that there exists $y \in m'$, such that the shortest geodesic $m''$ to $P \cap \Omega$ is not in $m'$. 
Consider the $2$-dimensional subspace $Q$ containing $m'$ and $m''$. Then this is a contradiction by
Corollary 1.5 of \cite{marg}. 

(iv) Again follows from considering a $2$-dimensional subspace containing $P$ and $m$. 
(See Proposition 1.4 of \cite{marg} for details.)
%
\end{proof}

An endpoint in $P$ of a shortest segment is called a {\em foot of the perpendicular} from $x$ to $\gamma$. \index{foot of the perpendicular} 


\subsection{Strict SPC-structures and the group actions}




By Corollary \ref{app-cor-stLens}, strict SPC-orbifolds with generalized 
lens-shaped or horospherical $\cR$- or $\cT$-ends have only lens-shaped or horospherical $\cR$- or $\cT$-ends. 



\begin{lemma} \label{rh-lem-hor} 
Let $\orb$ be an $n$-dimensional strongly tame strictly SPC-orbifold. 
Let $\tilde E$ be a p-end of $\torb$.
\begin{enumerate}
\item[(i)] Suppose that $\tilde E$ is a horospherical p-end.
Let $B$ be a horoball p-end neighborhood 
with a p-end vertex $p$ corresponding to $\tilde E$. 
There exists a homeomorphism $\Phi_{\tilde E}: \Bd B - \{p\} \ra \Bd \torb -\{p\}$ given 
by sending a point $x$ to the endpoint of the maximal convex segment containing $x$ and $p$
in $\clo(\torb)$. 
\item[(ii)] Suppose that $\tilde E$ is a lens-shaped radial p-end. 
Let $U$ be a lens-shaped radial p-end neighborhood with the p-end vertex $p$ corresponding to $\tilde E$. 
There exists a homeomorphism $\Phi_{\tilde E}: \Bd U \cap \torb  \ra \Bd \torb - \clo(U)$ given 
by sending a point $x$ to the other endpoint of the maximal convex segment containing $x$ and $p$
in $\clo(\torb)$.
\end{enumerate}
Moreover, $\Phi_{\tilde E}$ commutes with elements of $h(\pi_1(\tilde E))$. 
\end{lemma} 
\begin{proof} 
(i) 
By Theorem \ref{ce-thm-affinehoro}(i) $\Phi_{\tilde E}$ is well-defined.  
The same proposition implies that  $\Bd B $ is smooth at $p$ and $\Bd \torb$ has a unique sharply supporting hyperspace.
Therefore the map is onto. 

(ii) The second item follows from Theorems \ref{pr-thm-lensclass} and \ref{pr-thm-redtot}
since they imply that the segments in $S(p)$ are maximal ones in $\Bd \torb$ from $p$.
\end{proof} 


We now study the fixed points in $\clo(\torb)$ of elements of $\pi_1(\orb)$. 
Recall that a great segment is a geodesic arc in $\SI^n$ with antipodal p-end vertices. 
\index{great segment} 
It is not properly convex. 

Note that we can replace a generalized lens to a lens for a strongly tame strictly SPC-orbifold by Corollary \ref{app-cor-stLens}.
Also, recall Lemma \ref{prelim-lem-elliptic}.  

\begin{lemma}\label{rh-lem-fix} 
Let $\orb$ be an $n$-dimensional strongly tame strict SPC-orbifold with 
lens-shaped or horospherical $\cR$- or $\cT$-ends. 
Let $g$ be an infinite order element of a p-end fundamental group $\pi_1(\tilde E)$. 
Then every fixed point $x$ of $g$ in $\clo(\torb)$ 
satisfies one of the following\/{\em :} 
\begin{itemize} 
\item $x$ is in the closure of a p-end neighborhood that is a concave end neighborhood of an R-p-end, 
\item $x$ is in the closure of a p-ideal boundary component of a T-p-end or 
\item $x$ is the p-end vertex of a horospherical R-p-end.
\end{itemize} 
\end{lemma}
\begin{proof}
Suppose that the p-end $\tilde E$ is a lens-shaped R-end.
The direction of each segment in the interior of the lens cone 
with an endpoint $\mbv_{\tilde E}$ is fixed by only the finite-order element of 
$\pi_1(\tilde E)$ since $\pi_1(\tilde E)$ acts properly discontinuously on $\tilde S_{\tilde E}$. 
Thus, the fixed points are on the rays in the direction of the boundary of $\tilde E$. 
They are in one of $S(\mbv_{\tilde E})$ for the p-end vertex $\mbv_{\tilde E}$ corresponding to $\tilde E$ by Theorems \ref{pr-thm-lensclass} and \ref{pr-thm-redtot}.
Hence, the fixed points of 
the holonomy homomorphism of $\pi_1(\tilde E)$ is in the closure of the lens cone
with end vertex $\mbv_{\tilde E}$ and nowhere else
in $\clo(\tilde{\mathcal O})$. 

If $\tilde E$ is horospherical, then the p-end vertex $\mbv_{\tilde E}$ is not contained in any segment $s$ in $\Bd \torb$
by Theorem \ref{ce-thm-affinehoro}.  
Hence $\mbv_{\tilde E}$ is the only point of $S \cap \Bd \torb$ for any invariant subset $S$ of $\pi_1(\tilde E)$
by Lemma \ref{rh-lem-hor}. 
Thus, the only fixed point of $\pi_1(E)$ in $\Bd \torb$ is $\mbv_{\tilde E}$.

Suppose that $E$ is a 
lens-shaped T-p-end. 
Since $\tilde E$ is a properly convex real projective orbifold that is closed, 
we obtain an attracting fixed point 
$a$ and a repelling fixed point $r$ of $g|\clo(\tilde S_{\tilde E})$ by \cite{Benoist97}. 
Then $a$ and $r$ respectively 
are attracting and repelling fixed points of $g|\clo(\torb)$ by 
the existence of a cocompactly-acted lens neighborhood of $\tilde S_{\tilde E}$ since 
Theorem \ref{pr-thm-equ2} implies the uniform middle-eigenvalue condition. 

Suppose that we have a fixed point $s \in \Bd \torb$ distinct from $a$ and $r$. 
We claim that $\ovl{as}$ and $\ovl{rs}$ are in $\Bd \torb$. 
The norm of the eigenvalue associated with $s$ is strictly between
those of $r$ and $s$ by the uniform middle-eigenvalue condition. 
Let $P$ denote the two-dimensional subspace containing $r, s, a$. 
Suppose that one of the segment meets $\torb$ at a point $x$. 
We take a convex open-ball neighborhood $B$ of $x$ in $P \cap \torb$.
Suppose that $x \in \ovl{rs}^o$.
Then using the sequence $\{g^n(B)\}$, we obtain a great segment in $\clo(\torb)$ by choosing $n \ra \infty$.
This is a contradiction. If $x \in \ovl{as}^o$, we can use $\{g^{-n}(B)\}$ as $n \ra \infty$, again giving us
a contradiction.  
Hence, $\ovl{as}, \ovl{rs}\subset \Bd \torb$. 

Since $\tilde E$ has a one-sided neighborhood $U$ in a cocompactly-acted lens neighborhood
of $\tilde S_{\tilde E}$ 
by choosing a smaller such neighborhood $U$ if necessary, 
we may assume that 
$\clo(U) \cap \Bd \torb$ is in $\clo(\tilde S_{\tilde E})$.
By the strict convexity of $\torb$, we see that
the nontrivial segments $\ovl{as}$ and $\ovl{rs}$ have to be in 
$\clo(\tilde S_{\tilde E})$. (See Definition \ref{intro-defn-strict}.)
\end{proof}

See Crampon and Marquis \cite{CM14} and Cooper-Long-Tillmann \cite{CLT15} for the work similar to the following.  
We remind the reader that 
generalized lens-shaped R-ends 
are lens-shaped R-ends in the following assumption
by Corollary \ref{app-cor-stLens}, 

\begin{proposition} \label{rh-prop-loxopara}
Suppose that $\mathcal O$ is an $n$-dimensional 
strongly tame strictly SPC-orbifold with lens-shaped ends or horospherical $\cR$- or $\cT$-ends satisfying {\rm (}IE\/{\rm )} and 
{\rm (}NA{\rm ).} 
Then each nonidentity and infinite-order element $g$ of $\pi_1(\mathcal{O})$ has two exclusive possibilities\/{\rm :} 
\begin{itemize}
\item $g|\clo(\torb)$ has exactly two fixed points in $\Bd \tilde{\mathcal O}$ none of which is in the closures of the p-end neighborhoods for distinct ends, and 
$g$ is positive proximal. 
\item $g$ is in a p-end fundamental group $\pi_1(\tilde E)$, and $g|\clo(\torb)$ 
\begin{itemize} 
\item has all fixed points in $\Bd \tilde{O}$ in the closure of a concave p-end neighborhood of a lens-shaped radial p-end $\tilde E$. 
\item has all fixed points  in $\Bd \tilde{O}$ in $\clo(\tilde S_{\tilde E})$ for the ideal boundary component $\tilde S_{\tilde E}$ 
of a lens-shaped totally geodesic p-end $\tilde E$, or
\item has a unique fixed point  in $\Bd \tilde{O}$ at the horospherical p-end vertex. 
\end{itemize} 
\end{itemize}
\end{proposition}
\begin{proof}
Suppose that  $g$ has a fixed point at a horospherical p-end vertex $\mbv$ for a p-end $\tilde E$. 
We can choose the horoball $U$ at $\mbv$ that 
maps into an end neighborhood of $\orb$. 
A horoball p-end neighborhood is either sent to a disjoint one or 
sent to the identical one. 
Since $g(U) \cap U \ne \emp$ by the geometry of a horoball having a smooth boundary at $\mbv$, 
$g$ must act on the horoball, and hence $g$ is in the p-end fundamental group.
The p-end vertex is the unique fixed point of $g$ in $\Bd \torb$ by Lemma \ref{rh-lem-fix}.  

Similarly, suppose that $g \in \pi_1(\orb)$ fixes a point of the closure $U$ of a concave
p-end neighborhood of a p-end vertex $\mbv$ of a lens-shaped end. 
$g(\clo(U))$ and $\clo(U)$ meet at a point. 
By Corollary \ref{app-cor-disjclosure}, $g(\clo(U))$ and $\clo(U)$ share the p-end vertex 
and hence $g(U) = U$ as $g$ is a deck transformation. 
Therefore, $g$ is in the p-end fundamental group of the p-end of $v$. 
Lemma \ref{rh-lem-fix} implies the result. 


Suppose that $g \in \pi_1(\orb)$ fixes a point of $\clo(\tilde S_{\tilde E})$ for a totally geodesic ideal boundary component $\tilde S_{\tilde E}$ 
corresponding to a p-end $\tilde E$. 
Again Corollary \ref{app-cor-disjclosure} and Lemma \ref{rh-lem-fix} imply the result for this case. 




Suppose that an element $g$ of $\pi_1(\mathcal{O})$ is not in 
any p-end fundamental subgroup.
Then by above, $g$ does not fix any of the above types of points. 
We show that $g$ has exactly two fixed points in $\Bd \torb$: 

Suppose that $g  \in \pi_1(\orb)$ fixes a unique point $x$ in the closure of $\Bd \tilde{\mathcal{O}}$ 
and $x$ is not in the closure of p-end neighborhoods by the first part of the proof. 
Then $x$ is a $C^1$-point by the strict convexity. 
(See Definition \ref{intro-defn-strict}.)
Suppose that we have two eigenvalues with the largest norm  
$>1$ and the smallest norm $< 1$ respectively. 
If the largest norm eigenvalue is not positive real, $\clo(\tilde{\mathcal{O}})$ contains a nonproperly convex subset
as we can see by an action of $g^n$ on a generic point of $\torb$. 
Hence, the largest norm eigenvalue is positive and so is the smallest norm
eigenvalue. 
We obtain attracting and repelling subspaces easily with these, and there are at least 
two fixed points. This is a contradiction. 
Therefore, $g$ has only eigenvalues of unit norms. 
However, 
Lemma \ref{prelim-lem-unitnolower} shows that there is a sequence of simple closed 
curves $c_i$ whose the sequence of Hilbert lengths is going to zero.
Hence, $g$ must be freely homotopic to an end neighborhood. 
This is absurd.

We conclude that $g\in \pi_1(\orb)$ not in 
p-end fundamental groups fixes at least two points $a$ and $r$ in $\Bd \tilde{\mathcal{O}}$. 
We choose the two fixed points to have the positive real eigenvalues that are largest and smallest absolute values of 
the eigenvalues of $g$. (As above, the largest and smallest norm eigenvalues must be positive for $\torb$ to be properly convex.) 

No fixed point of $g$ in $\Bd \torb$ is in the closures of p-end neighborhoods
by the first part of the proof. 
By the strict convexity, the interior of $\torb$ contains an open line segment $l$ connecting $a$ and $r$. 



Suppose that there is a third fixed point $t$ in $\clo(\tilde{\mathcal{O}})$, 
which must be a boundary point. 
Let $S$ denote the subspace spanned by $a, r, t$. The point 
$t$ is not in the closures of p-end neighborhoods 
by the above part of the proof
as we assumed that $g$ is not in the p-end fundamental group. 
Then the line segment connecting $t$ to the $a$ or $r$ must be 
in $\Bd \tilde{\mathcal{O}}$: 
Assume without loss of generality $\overline{ta}^o\subset \torb$ 
by taking $g^{-1}$ and switching notation of $a$ and $r$ if necessary. 
Since $g$ acts properly, 
the norms of eigenvalues of $g$ at $t$ or $a$ are distinct. 
We can form a segment $s$ in $\torb \cap S$ transverse to the segment.
Then $\{g^k(s)\}$ geometrically converges to a segment in $\Bd \torb$ containing $t$
with endpoints $r$ and $r_-$ as $k \ra -\infty$.
This is absurd, and $\overline{ta}, \overline{rt} \subset \Bd \torb$. 
Thus, the existence of $t$ contradicts the proper convexity of $\clo(\torb)$. 

Hence, there are exactly two 
fixed points of $g$ in $\Bd \tilde{\mathcal{O}}$ of the  positive real eigenvalues that are largest and smallest absolute values of 
the eigenvalues of $g$.
\end{proof}


\begin{proposition} \label{rh-prop-loxopara2}
Suppose that $\mathcal O$ is an $n$-dimensional 
noncompact strongly tame strictly SPC-orbifold  with lens-shaped ends 
or horospherical $\cR$- or $\cT$-ends. 
Let $\tilde E$ be an end. 
Then for a p-end $\tilde E$, $(\Bd \torb - K)/\pi_1(\tilde E)$ is a compact orbifold where
$K= \bigcup S(\tilde E)$ for a lens-shaped radial p-end $\tilde E$,  
$K= \clo(\tilde S_{\tilde E})$ for totally geodesic p-end $\tilde E$, or
 $K= \{\mbv_{\tilde E}\}$ for horospherical p-end $\tilde E$.
\end{proposition}
\begin{proof} 
Suppose that $\tilde E$ is a lens-shaped R-p-end or a horospherical one.
By Lemma \ref{rh-lem-hor}, the homeomorphism 
$\Phi_{\tilde E}: \tilde S_{\tilde E} \ra \Bd \torb - K$
gives us the result.   

Suppose that $\tilde E$ is a lens-shaped T-p-end. Let $\torb^*$ denote the dual domain.
Then there exists a dual radial p-end $\tilde E^*$ corresponding to $\tilde E$.
Hence, $(\Bd \torb^* - K')/\pi_1(\tilde E^*)$ is compact for $K'$ equal to the closure of p-end neighborhoods of $\tilde E^*$ in the radial case or 
the vertex in the horospherical case.

Recall Section \ref{prelim-sec-duality}.
Let $\Bd^{\Ag} \torb$ be the augmented boundary with the fibration $\Pi_{\Ag}$, and 
let $\Bd^{\Ag} \torb^*$ be the augmented boundary with the fibration map $\Pi^*_{\Ag}$.
Let $K'':= \Pi_{\Ag}^{-1}(K)$ and $K''':= \Pi_{\Ag}^{\ast -1}(K')$. 
The discussion on in the proof of Corollary \ref{pr-cor-duallens} shows that 
there is a duality homeomorphism 
\[ \mathcal{D}_{\torb}: \Bd^{\Ag} \torb - K''\ra \Bd^{\Ag} \torb^* - K'''.\]

Now $(\Bd^{\Ag} \torb^* - K''')/\pi_1(\tilde E^*)$ is compact since $\Bd \torb^*  - K'$ has a compact fundamental domain, 
and the space is the inverse image in $\Bd^{\Ag} \torb^*$
of $\Bd \torb^* - K'$. 
By (iv) of Proposition \ref{prelim-prop-duality}, $(\Bd^{\Ag} \torb - K'')/\pi_1(\tilde E)$ is compact also. 
Since the image of this set under the map induced by a proper map $\Pi_{\Ag}$ is $(\Bd \torb - K)/\pi_1(\tilde E)$. 
Hence, it is is compact.
\end{proof}

\section{Bowditch's method} \label{rh-sec-bowditch} 

\subsection{The strict convexity implies the relative hyperbolicity} 
There are results proved by Cooper, Long, and Tillmann \cite{CLT15} and Crampon and Marquis \cite{CM14} 
similar to below. However, the ends have to be horospherical in their work. 
By Lemma \ref{app-cor-stLens}, for a strongly tame strictly SPC-orbifold, generalized lens-shaped ends are lens-shaped. 
We use Bowditch's result \cite{Bowditch98} to show:

\begin{theorem}\label{rh-thm-relhyp}
Let $\mathcal{O}$ be an $n$-dimensional noncompact strongly tame strictly SPC-orbifold with lens-shaped ends or horospherical $\cR$- or $\cT$-ends 
$E_1, \dots, E_k$ and satisfies  {\rm (}IE\/{\rm )} and {\rm (}NA\/{\rm )}.  
Assume $\partial \orb$ is smooth and strictly convex. 
Let $\tilde {\mathbbm{U}}_i$ be the inverse image ${\mathbbm{U}}_i$ in $\torb$ for a mutually disjoint collection of
 neighborhoods ${\mathbbm{U}}_i$ of the ends $E_i$ for each $i=1, \dots, k$. 
Then the following hold\/{\rm :}
\begin{itemize} 
\item $\pi_1(\mathcal{O})$ is relatively hyperbolic with respect to 
the end fundamental groups \[\pi_1(E_1),\dots, \pi_1(E_k).\] 
Hence $\orb$ is relatively hyperbolic with 
respect to ${\mathbbm{U}}:= {\mathbbm{U}}_1 \cup \cdots \cup  {\mathbbm{U}}_k$
as a metric space.  
\item If $\pi_1(E_{l+1}),\dots, \pi_1(E_k)$ are hyperbolic for some $1 \leq  l \leq k$
{\rm (}possibly some of the hyperbolic ones{\rm )},
then $\pi_1(\mathcal{O})$ is relatively hyperbolic with respect to the end fundamental group 
$\pi_1(E_1), \dots, \pi_1(E_{l})$. 
\end{itemize}
\end{theorem}
\begin{proof}
We show that $\pi_1(\mathcal{O})$ is relatively hyperbolic with respect to 
the end fundamental groups $\pi_1(E_1),\dots, \pi_1(E_k)$. 



Let ${\mathcal C}_B$ denote the following collection:
\begin{itemize}
\item the set of form 
$\clo({\mathbbm{U}}_i) \cap \Bd \torb = \tilde S_{\tilde E}$  for a concave p-end neighborhood ${\mathbbm{U}}_i$ and 
\item the sets of form $\clo(\tilde S_{\tilde E})$ for each lens-shaped totally geodesic end $\tilde E$. 
\item the singletons of vertices of horospherical ends. 
\end{itemize} 
By Corollary \ref{app-cor-disjclosure}, these sets are mutually disjoint balls or singletons in 
$\Bd \torb$. 
 We also denote $C_B:= \bigcup {\mathcal C}_B$. 

We claim that for each closed set $J$ in $\Bd \torb$, the union of $C_J$ of 
elements of ${\mathcal C}_B$ meeting $J$ is also closed: 
Let us choose a sequence $\{x_i \}$ for $x_i \in C_i$, $C_i \cap J \ne \emp$, $C_i \in {\mathcal C}_B$. 
Suppose that $\{x_i\} \ra x$. Let $y_i \in C_i \cap J$. Let $\mbv_i$ be the p-end vertex of $C_i$ if it is from a R-p-end.
Then define $s_i:= \ovl{x_i \mbv_i}\cup \ovl{\mbv_i, y_i} \subset C_i$ if $C_i$ is radial 
or else $s_i:= \ovl{x_iy_i} \subset C_i$. Choose a subsequence such that 
$\{s_i\}$ geometrically converges to a limit containing $x$. 
The limit $s_\infty$ is a singleton, a segment or a union of two segments. 
If $s_\infty$ is a singleton, then $s_\infty$ is a singleton in $J$ clearly. 
Suppose that $s_\infty$ contains a nontrivial segment. 
By the strict convexity of $\torb$, we obtain that 
$s_\infty$ is a subset of an element of ${\mathcal C}_B$
and $s_{\infty}$ meets $J$. 
Thus, $x \in s_{\infty} \subset C_j$ for $C_j \cap J \ne \emp$.
We conclude that $C_J$ is closed. 

We denote this quotient space  $\Bd \tilde{\mathcal{O}}/\sim$ by $B$.
By Proposition \ref{rh-prop-metriz}, $B$ is a metrizable space. 
We denote by $[\cdot]$ the equivalence classes for $\sim$ here. 


We show that $\pi_1(\mathcal{O})$ acts on the metrizable space $B$ 
as a geometrically finite convergence group. 
By Theorem 0.1 of Yaman \cite{Yaman04} following Bowditch \cite{Bowditch98}, this shows that  \index{relative hyperbolicity} 
$\pi_1(\orb)$ is relatively hyperbolic with respect to $\pi_1(E_1), \dots, \pi_1(E_k)$.
The definition of conical limit points and so on are from the article.


Let us choose a basis and the corresponding elements 
$e_1, e_2, \dots, e_{n+1}$ in $\SI^n$. 

(I) We first show that the group acts properly discontinuously on the metrizable space of ordered mutually distinct triples in  the metrizable space $B=\partial \tilde{\mathcal{O}}/\sim$. 
Suppose not. Then 
we can deduce that 
there exists a sequence of mutually distinct triples $\{([p_i], [q_i], [r_i])\}$  
of points $p_i, q_i, r_i$ in $\Bd \torb$  with following properties: 
\begin{itemize} 
\item The sequence converges to a mutually distinct triple $\{([p], [q], [r])\}$  for $p, q, r \in \Bd \torb$ such that 
\[p_i = \gamma_i(p_0^i), q_i = \gamma_i(q_0^i), \hbox{ and } r_i= \gamma_i(r_0^i)\]
\item $\{\gamma_i\}$ is a sequence of mutually distinct elements of $\pi_1(\mathcal{O})$ and 
\item The sequence of 
mutually distinct equivalence classes $[p_0^i], [q_0^i], [r_0^i]$
converge to mutually distinct $[\hat p_0], [\hat q_0], [\hat r_0]$ respectively.
\item We may further assume $p_0^i \ra \hat p_0, q_0^i \ra \hat q_0, r_0^i \ra \hat r_0$
by a choice of subsequences. 
\end{itemize}  
We can deduce this since the space of triples is a metrizable space. 

We find some uniformly bounded sequence of elements $\hat R_i \in \SLpm$
where $\hat R_i(p_0^i), \hat R_i(q_0^i), \hat R_i(r_0^i)$ form a fixed triple
$e_1, e_2, e_3$. 
Also, by multiplying by some uniformly bounded element $R_i$ in $\SLpm$ but not necessarily in $h(\pi_1(\orb))$, we obtain 
that $R_i\circ \gamma_i\circ \hat R_i^{-1}$ 
for each $i$ fixes $e_1, e_2, e_3$ and restricts to a diagonal matrix with 
entries $\lambda_i, \delta_i, \mu_i$ on the span of $e_1, e_2, e_3$. 

Then we can assume that 
\[\lambda_i\delta_i\mu_i =1, \lambda_i \geq \delta_i \geq \mu_i > 0\] 
by restricting to the plane and 
up to choosing subsequences and renaming. 
 
We claim that  $\{\lambda_i\} \ra \infty$ and $\{\mu_i\} \ra 0$: 
Suppose not. Then both of 
these two sequences are bounded. Let $P_i$ denote the $2$-dimension subspace 
spanned by $p_0^i, q_0^i, r_0^i$. 
Then $\gamma_i| P_i$ is a sequence of uniformly bounded 
automorphisms. 

Let $D_i = \clo(\torb\cap P_i)$. Then there exists $\eps>0$ such that 
$D_i$ is not in a $\eps$-$\bdd$-neighborhood of $\Bd \torb$ for all $i$:  
If not, the geometric limit of $D_i$ is a nontrivial 
disk in $\Bd \torb$ containing $\hat p_0, \hat q_0, \hat r_0$. 
Since $\hat p_0, \hat q_0, \hat r_0$ are in mutually distinct equivalence classes, 
this contradicts the strict convexity. 

Hence, for a compact subset $K$ of $\torb$, $K \cap D_i\subset \torb$ 
is not empty for sufficiently large $i$. Choose $s_i \in K \cap D_i$. 
Then $\gamma_i(s_i)$ is in a $d_{\torb}$-bounded neighborhood of $K$ independent of $i$ since otherwise $\gamma_i$ is not uniformly bounded on $\torb \cap P_i$ 
as indicated by the Hilbert metric. To see this, we post-compose by 
a bounded sequence of projective automorphisms such that $P_i$ is fixed. 
Hence, $\gamma_i(K)$ is in 
a $d_{\torb}$-bounded neighborhood of $K$. This contradicts 
the proper discontinuity of the action by $\bGamma$. 

Hence, we showed that  $\{\lambda_i\} \ra \infty$ and $\{\mu_i\} \ra 0$. 

By the strictness of convexity, as we collapsed each of the p-end balls, 
the interiors of the segments $\ovl{p_0^iq_0^i}$, $\ovl{q_0^ir_0^i}$, and $\ovl{r_0^ip_0^i}$ 
are in the interior of $\tilde{\mathcal{O}}$. 
Again, these segments are not contained in some $\eps$-$\bdd$-neighborhoods
$\Bd \torb$ for some $\eps> 0$ since 
$p_0^i \ra \hat p_0, q_0^i \ra \hat q_0, r_0^i \ra \hat r_0$ 
where $[\hat p_0], [\hat q_0], [\hat r_0]$ form a mutually disjoint triple. 

We claim that one of the sequence  $\{\lambda_i/\delta_i\}$ or the sequence $\{\delta_i/\mu_i\}$ are bounded:
Suppose not. 
Then $\{\lambda_i/\delta_i\} \ra \infty$ and $\{\delta_i/\mu_i\} \ra \infty$. 
We choose generic segments $s_0^i$ and $t_0^i$ in $\torb$
with a common endpoint $q_0^i$ and the respective other endpoint $\hat s_0^i$ and $\hat t_0^i$ in different components of 
$P\cap \orb - \ovl{p_0^iq_0^i}$
such that \[\bdd(\hat s_0^i, q_0^i), \bdd(\hat t_0^i, q_0^i) \geq \delta \hbox{ for a uniform } \delta > 0.\] 
We choose $s_0^i$ and $t_0^i$ 
such that their directions from $q_0^i$ differ from 
those of $\ovl{p_0^iq_0^i}$ and $\ovl{q_0^ir_0^i}$
at least by a small constant $\delta' > 0$. 
Then the sequence of lifts to $\SI^n$ of $\{R_i \circ \gamma_i(s_0^i\cup t_0^i)\}$ 
geometrically converges to the segment with endpoint $\hat p_0$ passing $\hat q_0$. 
The segment is a great segment. 
Since $R_i$ is bounded, this implies that there exists such a segment in $\clo(\torb)$. 
This is a contradiction to the proper convexity of $\torb$. 

Suppose now that the sequence $\lambda_i/\delta_i$ is bounded: Now the sequence of segments $\{\ovl{p_iq_i}\}$ 
converges to $\ovl{pq}$ whose interior is in $\tilde{\mathcal{O}}$. 
We can choose a subsegment $s_i \subset \ovl{p_0^iq_0^i}^o$ 
in a compact subset of $\torb$ since $ \ovl{p_0^iq_0^i}^o$ 
converges to a segment $\ovl{\hat p_0\hat q_0}$ whose interior is in $\torb$. 
Then we see that $\ovl{pq}$ must be in the boundary of $\torb$: 
Each point in $\ovl{pq}$ must 
be the limit points of a sequence $\{y_i\}$ for $y_i\in \gamma_i(s_i)$
by the boundedness of the above ratio.
This implies that $\ovl{pq} \subset \Bd \torb$ by  
the proper-discontinuity of the action.
This contradicts the strict convexity as we assumed that 
$p, q,$ and $r$ represent distinct points in $B$. 
If we assume that $\delta_i/\mu_i$ is bounded, then 
we obtain a contradiction similarly.

Hence, we obtain contradictions in all cases. 
This proves the proper discontinuity of the action on the space of distinct triples. 


(II) By Propositions \ref{rh-prop-loxopara} and \ref{rh-prop-loxopara2}, 
each group of form $\Gamma_x$ for a point $x$ of $B=\tilde{\mathcal{O}}/\sim$
corresponding to an end 
is a bounded parabolic subgroup in the sense of Bowditch \cite{Yaman04}. 


Now we take a concave end neighborhood for each radial end of lens-type
and horospherical 
We still choose ones mutually disjoint from themselves and 
nonradial ones in $\mathbbm{U}$.
We denote by $\mathbbm{U}'$ the union of the modified end neighborhoods. 
Let $\mathbbm{U}'_1, \dots, \mathbbm{U}'_k$ denote its components.
Let $\tilde{\mathbbm{U}}'_k$ denote the inverse image of $\mathbbm{U}'_k$ in 
$\torb$ for each $k$. 

(III) Let $p \in \Bd \torb$ be a point that is not in a horospherical endpoint or an equivalence class corresponding to 
a lens-shaped p-end of radial or totally geodesic type of $B$. 
Hence, $[p] = \{p\}$ in $\Bd \torb$. 
That is, there is no segment containing $p$ in one of the collapsed sets. 
We show that $[p]$ is a conical limit point. This completes our proof by Theorem 0.1 of \cite{Yaman04}. 

To show that $[p]$ is a conical limit point, 
we find a sequence of holonomy transformations $\gamma_i$ and 
distinct points $a, b \in B$ 
such that $\{\gamma_i([p])\} \ra a$ and $\{\gamma_i(q)\} \ra b$ locally uniformly for 
$q \in B - \{p\}$:
To do this, we draw a line $l$ in $\torb$ from a point of the fundamental domain to $p$
where as $t \ra \infty$, $l(t) \ra p$ in $\clo(\torb)$.
We may assume that the other endpoint $p'$ of $l$ 
is in a distinct equivalence class from 
$[p]$.  
Since $l(t)$ is not eventually in a p-end neighborhood, there is a sequence $\{t_i\}$ going to $\infty$ 
such that $l(t_i)$ is not in any of the p-end neighborhoods in $\tilde {\mathbbm{U}'}_{1} \cup \cdots \cup  \tilde {\mathbbm{U}'}_{k}$. 
Let $p'$ be the other endpoint of the complete extension of $l(t)$ in $\tilde{\mathcal{O}}$. 
We can assume without loss of generality that 
$p'$ is not in the closure of any p-end neighborhood by choosing the line $l(t)$ differently if necessary. 

Since $(\tilde{\mathcal{O}}- \tilde {\mathbbm{U}'}_1 - \cdots - \tilde {\mathbbm{U}'}_k)/\Gamma$ is compact, 
we have a compact fundamental domain $F$  of 
$\tilde{\mathcal{O}} - \tilde {\mathbbm{U}'}_1 - \cdots - \tilde {\mathbbm{U}'}_k$ with respect to $\Gamma$.
Note that for the minimum distance, 
we have $\bdd(F, \Bd \torb) > C_0$ for some constant $C_0 > 0$. 

We note: 
Given any line $m$ passing $F$, the two endpoints must be in 
distinct equivalence classes because of the convexity of 
each component of $\tilde{\mathbbm{U}'}$.  

We find a sequence of points $z_i \in F$ such that $\gamma_i(l(t_i))= z_i$ for a deck transformation $\gamma_i$.
Then $\{\gamma_i\}$ is an unbounded sequence. 

Using Definition \ref{prelim-defn-convlimit}, we may choose a 
set-convergent subsequence of  
$\{\llrrparen{\gamma_i}\}$ that is convergent in $\SI(M_{n+1}(\bR))$ to 
$\llrrparen{\gamma_\infty}$ for $\gamma_\infty \in M_{n+1}(\bR)$. 
Hence, $A_{\ast}(\{\gamma_i\}) = \SI(\Ima \gamma_{\infty}) \cap \clo(\torb)$. 
Also, on $\clo(\torb ) - N_{\ast}(\{\gamma_i\}) \cap \clo(\torb)$, 
$\{\gamma_i\}$ is convergent to a subset of $A_{\ast}(\{\gamma_i\})$ locally
uniformly as we can easily deduce by linear algebra
and some estimation. 

Since $N_{\ast}(\{\gamma_i\})$ is a convex subset of $\Bd \torb$  
and $A_{\ast}(\{\gamma_i\})$ is a convex subset of $\Bd \torb$ 
by Theorem \ref{prelim-thm-AR},  they are in collapsed sets of 
$\Bd \torb$ by the strictness of the convexity. 


If $p \not\in N_{\ast}(\{\gamma_i\})$, then $\gamma_i(l(t_i))$ is 
also bounded away from $N_{\ast}(\{\gamma_i\})$, and hence 
$\gamma_i(l(t_i))$ accumulates only to $A_{\ast}(\{\gamma_i\})$. 
This is a contradiction since $\gamma_i(l(t_i)) \in F$. 
Thus, $p \in N_{\ast}(\{\gamma_i\})$. 
Since $N_{\ast}(\{\gamma_i\})$ is a convex compact subset of 
$\Bd \torb$, we must have by the strict convexity 
\[N_{\ast}(\{\gamma_i\}) \subset [p].\]
Thus, for all $q \in \Bd \torb - [p]$, we obtain 
a local uniform convergence under $\gamma_i$ to $A_{\ast}(\{\gamma_i\})$.
This shows that $p$ is a conical limit point. 
We let $b$ be the collapsed set containing 
 $A_{\ast}(\{\gamma_i\})$.

Our line $l$ equals the interior of $\overline{pp'}$. 
We choose a subsequence of $\gamma_i$ such that
the corresponding subsequence  $\{\gamma_i(\overline{pp'})\}$ 
geometrically converges to a line passing $F$. 
Since $p$ and $p'$ are in district equivalence classes, 
$[\gamma_i(p')]$ converges to $b$, and 
$\gamma_i(\overline{pp'})$ passes $F$, 
it follows that $\gamma_i(p)$ converges to a point of 
the equivalence class $a$ distinct from $b$ by our note above. 

Finally, we remove concave end neighborhoods for $E_{l+1}, \dots, E_k$ or 
add lens end neighborhoods by Theorems \ref{rh-thm-remconch} and \ref{rh-thm-totgeoext}. 
The resulting orbifold is a strictly SPC-orbifold again and 
we can apply the result (i) to this case and obtain (ii). 
\end{proof}


\subsection{The theorem of Dru{\c{t}}u}

The author obtained a proof of the following theorem from Dru{\c{t}}u. 
See \cite{Drutu09} for more details.  

\begin{theorem}[Dru{\c{t}}u] \label{rh-thm-drutu}
Let $\mathcal{O}$ be a strongly tame properly convex real projective $n$-orbifold with 
generalized lens-shaped ends and horospherical $\cR$- or $\cT$-ends
and satisfies  {\rm (}IE\/{\rm )} and {\rm (}NA\/{\rm )}. 
Let  $\pi_1(E_1),\dots, \pi_1(E_m)$ be end fundamental groups where $\pi_1(E_{l+1}), \dots, \pi_1(E_m)$ for $l\leq m$ 
are hyperbolic groups.
Then $\pi_1(\mathcal{O})$ is a relatively hyperbolic group with respect to 
$\pi_1(E_1),\dots, \pi_1(E_m)$ if and only if 
$\pi_1(\mathcal{O})$ is one with respect to 
 $\pi_1(E_1),\dots, \pi_1(E_l)$.
\end{theorem} 
\begin{proof} 
With the terminology in the paper \cite{Drutu09},   $\pi_1(\mathcal{O})$ is a relatively hyperbolic group with respect to 
the end fundamental groups $\pi_1(E_1),\dots, \pi_1(E_m)$ if and only if 
 $\pi_1(\mathcal{O})$ with a word metric is asymptotically tree graded (ATG) with respect to 
all the left cosets $g\pi_1(E_i)$ for $g \in \pi_1(\mathcal{O})$ and $i=1,\dots, m$. 

We claim that 
 $\pi_1(\mathcal{O})$ with a word metric is asymptotically tree graded (ATG) with respect to 
 all the left cosets $g\pi_1(E_i)$ for $g \in \pi_1(\mathcal{O})$ and $i=1,\dots, m$ if and only if
$\pi_1(\mathcal{O})$ with a word metric is asymptotically tree graded with respect to 
 all the left cosets $g\pi_1(E_i)$ for $g \in \pi_1(\mathcal{O})$ and $i=1,\dots, l$: 
 
 Conditions ($\alpha_1$) and ($\alpha_2$) of Theorem 4.9 in \cite{Drutu09} are satisfied still when we drop end fundamental groups
  $\pi_1(E_{n+1}),\dots, \pi_1(E_m)$ or add  them. 
 (See also Theorem 4.22 in \cite{Drutu09}.)
  
  For the condition $(\alpha_3)$ of Theorem 4.9 of \cite{Drutu09}, it is sufficient to consider only hexagons.
According to Proposition 4.24 of \cite{Drutu09} one can take the fatness constants as large
as one wants, in particular $\theta$ (measuring how fat the hexagon is) much larger than $\chi$ prescribing
how close the fat hexagon is from a left coset.

 If $\theta$ is very large, left cosets containing such hexagons in their
neighborhoods can never be cosets of hyperbolic subgroups 
since hyperbolic groups do not contain fat
hexagons. So the condition $(\alpha_3)$ is satisfied too whether one adds  $\pi_1(E_{n+1}),\dots, \pi_1(E_m)$ or drop them.
\end{proof} 


\subsection{Converse}

We prove the converse to Theorem \ref{rh-thm-relhyp}. We use the theory of tree-graded spaces and asymptotic cones \cite{DS08} and its appendix. 

\begin{itemize} 
\item We shave off every generalized lens-shaped R-ends of $\orb$ by Theorem \ref{rh-thm-remconch}
to obtain $\orb^{(1)}$. 
\item We expand $\orb^{(1)}$ to $\orb^{(2)}$ 
by adding lens neighborhoods to totally geodesic ideal boundary components by Theorem \ref{rh-thm-totgeoext}.
\item We then take out the interior of outside parts of the lens of $\orb^{(2)}$ for every T-ends
to obtain $\orb^{(3)}$. 
\item Next we remove a collection of the mutually disjoint horospherical 
end neighborhoods. 
Let the resulting compact orbifold be denoted $\orb^{(4)}$. 
\end{itemize} 
We can consider their universal covers $\torb^{(i)}$ as subspaces of $\SI^n$. 
Now $\tilde S_{\tilde E}$ for every totally geodesic p-end $\tilde E$ is in $\torb^{(i)}$ for $i=2, 3, 4$. 




\begin{proposition} \label{rh-prop-shaving} 
	Let $\orb$ be a noncompact strongly tame properly convex real projective 
$n$-orbifold with generalized admissible ends. 
Let $\orb^{(4)}$ have the restricted metric $d_{\torb^{(2)}}$.
Then $\tilde \orb^{(4)}$ is quasi-isometric with $\pi_1(\orb)$. 
\index{dtorbM@$d_{\orb^{(2)}}$|textbf} 
\end{proposition} 
\begin{proof} 
Let $\pi_1(\orb)$ have the set of generators $g_1, \dots, g_q$. 
Since we removed all the end neighborhoods of $\orb^e$, our orbifold 
$\orb^M$ is compact. 
Hence,  $\tilde\orb^{(4)}$ has a compact fundamental domain $F$. 
We find a function $\tilde \orb^{(4)} \ra \pi_1(\orb)$ by defining 
$gF^o$ to go to $g$ and defining arbitrarily the faces of $F$ to 
go to $gg_i$. 
and hence there is a function from it to  
$\pi_1({\mathcal{O}})$ decreasing distances up to a positive 
constant. 

Conversely, there is a function from $\pi_1({\mathcal{O}})$ to  $\tilde{\mathcal{O}}^M$
by sending $g$ to $g(x_0)$ for a fixed $x_0 \in F^o$. 
This is also distance decreasing up to a positive constant. 
Hence, this proves the result. 
(This is just Svarc-Milnor Lemma. See \cite{Milnor2} and \cite{Svarc}). 
\end{proof} 




\begin{theorem} \label{rh-thm-converse}
Let $\mathcal{O}$ be a strongly tame properly convex real projective $n$-orbifold with  
generalized lens-shaped or horospherical $\cR$- or $\cT$-ends
and satisfies  {\rm (}IE\/{\rm )} and {\rm (}NA\/{\rm )}. 
Assume $\partial \orb$ is smooth and strictly convex. 
Suppose that $\pi_1(\mathcal{O})$ is a relatively hyperbolic group with respect to 
the end groups $\pi_1(E_1),\dots, \pi_1(E_k)$ where $E_i$ are horospherical 
for $i=1,\dots, m$ and  generalized lens-shaped for $i=m+1,\dots, k$ for $0 \leq m \leq k$.  
Then $\mathcal{O}$ is strictly SPC with respect to the ends $E_1, \dots, E_k$
with lens-type R-ends and T-ends or horospherical ends.
\end{theorem}
\begin{proof} 

Since an $\eps$-mc-p-end neighborhood is always proper by Corollary \ref{app-cor-mcn} for sufficiently small $\epsilon$, 
we choose the end neighborhood of any  generalized lens-shaped 
R-p-end $\tilde E_i$ to be the image of an $\eps$-mc-p-end neighborhood
for some $\epsilon > 0$. 
Assume that all shaved-off parts are inside the union of these. 
We can choose all such neighborhoods and horospherical 
end neighborhoods and lens-shaped end neighborhoods for 
T-ends to be mutually disjoint by Corollaries \ref{app-cor-shrink}
and \ref{app-cor-mcn}.
Let $\tilde{\mathbbm{U}}$ denote the union of the inverse images of the end neighborhoods. 

Suppose that $\mathcal{O}$ is not strictly convex. 
We divide the work into two cases: 
First, we assume that there exists a segment in $\Bd \torb$ not contained in 
the closure of a p-end neighborhood. Second, we assume that there exists a non-$C^1$-point in $\Bd \torb$ not 
contained in the closure of a p-end neighborhood. 

(I) We assume the first case now.
We obtain a triangle with boundary in $\Bd \torb^{(3)}$: 
Let $l$ be a nontrivial maximal segment in $\Bd \torb^{(3)}$ not contained in the closure of a p-end neighborhood intersected with $\Bd \torb^{(3)}$. 
First, $l$ does not meet the closure of a horospherical p-end neighborhood by Theorem \ref{ce-thm-affinehoro}. 
By Theorems \ref{pr-thm-lensclass} and \ref{pr-thm-redtot} if $l^o$ meets the closure of 
a lens-shaped R-p-end neighborhood, then $l^o$ is in the closure. 
Also, suppose that $l^o$ meets $\tilde S_{\tilde E}$ for a totally geodesic p-end $\tilde E$. Then $l^o \cap \partial \clo(\tilde S_{\tilde E}) \ne \emp$. 
$l$ is in the hyperspace $P$ containing $\tilde S_{\tilde E}$ since otherwise 
we have some points of $\tilde S_{\tilde E}$ in the interior of $\clo(\torb)$.
We take a convex hull of $l\cup \tilde S_{\tilde E}$ which 
is a domain $D$ containing $\tilde S_{\tilde E}$ where 
$\pi_1(\tilde E)$ acts on. Then $D$ is still properly convex since so is
$\clo(\torb)$. Since $D^o$ has a Hilbert metric, $\pi_1(\tilde E)$ acts 
properly on $D^o$. By taking a torsion-free subgroup by Theorem \ref{prelim-thm-vgood}, we obtain  that $\tilde S_{\tilde E}/\pi_1(\tilde E) \ra D^o/\pi_1(\tilde E)$ has to be surjective.  Hence, $D^o = \tilde S_{\tilde E}$. 
Therefore, $l^o \subset \clo(\tilde S_{\tilde E})$, a contradiction. 
(See Theorem 4.1 of \cite{CLM18} and \cite{Benoist01}.) 
Therefore, $l$ meets the closures of p-end neighborhoods possibly only at its endpoints. 

%
%


Let $P$ be a $2$-dimensional subspace containing $l$ and meeting the interior of $\torb^{(3)}$ outside $\tilde{\mathbbm{U}}$.
By above, $l^o$ is in the boundary of $P \cap \torb^{(3)}$. 
Draw two segments $s_1$ and $s_2$ in $P \cap \torb^{(3)}$ 
from the two endpoints of $l$ meeting at a vertex $p$ in the interior of $\torb^{(3)}$. 

Let $I$ denote the index set of 
components of $\tilde{\mathbbm{U}}$. Let $\mathbbm{U}_i$ be a component of $\tilde {\mathbbm{U}}$ for $i\in I$. 

Define $A_i$ to be the set of points $x$ of $l^o$ with an open $\bdd$-metric ball-neighborhood in $\clo(\orb) \cap P$ in the closure of 
a single component $\mathbbm{U}_i$. By definition, 
$A_i$ is open in $l^o$. Also, 
$l^o$ is not a subset of single $A_i$ since otherwise $l$ is in the closure
of $U_i$, a contradiction. 
Since $l^o$ is connected, 
$l^o - \bigcup_{i\in I} A_i$ is not empty. 
Choose a point $x$ in it. 
For any open $\bdd$-metric ball-neighborhood $B$ of $x$ in 
$\clo(\torb) \cap P$, we cannot have 
$B \cap \torb \subset \tilde{\mathbbm{U}}$ 
since otherwise $B$ is in a single $\clo(\mathbbm{U}_i)$. 
For each open ball-neighborhood $B$ of $x$, 
$B- \tilde{\mathbbm{U}}$ is not empty. 
We conclude that 
$(\torb - \tilde{\mathbbm{U}}) \cap P$ has 
a sequence of points $\{x_i\}$ converging to a point $x$ of $l^o$. 

Then we claim \[d_{\orb^{(2)}} (x_i, s_1\cup s_2) \ra \infty:\]
Consider any sequence of any maximal straight segment $t_i$ from $x_i$ passing 
a point $y_i$ of $s_1$ or $s_2$. Let us orient it in the direction of $y_i$ from $x_i$.
Then let $\delta_+ t_i$ be the forward endpoint of
$t_i$ and $\delta_- t_i$ the backward one. 
Then the $\bdd$-distance from $y_i$ to $\delta_+ t_i$ goes to zero
by the maximality of $l$, 
which implies the Hilbert metric result by the cross-ratio consideration.

Recall that there is a compact fundamental domain $F$ of $\tilde {\mathcal{O}} -\tilde U$ under the action of $\pi_1(E)$.
Now, we can take $x_i$ to the fundamental domain $F$ by $g_i$. 
We choose $g_i$ to be a sequence of mutually distinct elements of 
$\pi_1(\mathcal{O})$. 
We choose a subsequence such that we assume without loss of generality that 
$\{g_i(T)\}$ geometrically converges to a convex set, 
which could be a point or a segment or a nondegenerate 
triangle. Since $g_i(T) \cap F \ne \emp$, and the sequence $\partial g_i(T)$ exits any compact 
subsets of $\tilde {\mathcal{O}}$ always while 
\[\{d_{\orb^{(2)}}(g_i(x_i), \partial g_i(T))\} \ra \infty\]
and $g_i(T)$ passes $F$, 
we see that a subsequence of $\{g_i(T)\}$ converges to a nondegenerate triangle, say $T_\infty$.

By following Lemma \ref{rh-lem-geogroup}, $T_\infty$ is such that $\partial T_\infty$ is in $\bigcup S(\mbv_{\tilde E})$ 
for a generalized lens-shaped R-p-end $\tilde E$.  


Now, $T_\infty$ is such that $\partial T_\infty \subset \clo(U_1)$ for a p-end neighborhood $U_1$ of a generalized 
lens-shaped end $\tilde E$. Then for sufficiently small $\eps'> 0$, the $\eps'$-$d_{\orb}$-neighborhood of $T_\infty\cap \torb$ 
is a subset of $U_1$ as $U_1$ was chosen to be an $\eps$-mc-p-end neighborhood
(see Lemma \ref{app-lem-mcc}).
However, as $\{g_i(T)\} \ra T_\infty$ geometrically, 
for any compact subset $K$ of $\torb$, $g_i(T) \cap K$ is a subset of $U_1$ for sufficiently large $i$. 
But $g_i(T) \cap F \ne \emp$ for all $i$ and the compact fundamental domain $F$ of $\torb - \tilde U$, 
disjoint from $U_1$. This is a contradiction.


(II) Now we suppose that $\Bd \torb$ has a non-$C^1$-point $x$ outside the closures of p-end neighborhoods. 
Then we go to the dual $\torb^*$ and the dual group
$\Gamma^*$ where $\torb^*/\Gamma^*$ is a strongly tame properly convex real projective $n$-orbifold with horospherical ends, lens-shaped T-ends
or generalized lens-shaped R-ends 
by Corollary \ref{pr-cor-duallens2} and Theorem \ref{prelim-thm-dualdiff}. 
Here the type $\cT$ and $\cR$ are switched for the correspondence between 
the ends of $\orb$ and $\orb^\ast$ by Corollary \ref{pr-cor-duallens2}. 

 Then we have a one-to-one correspondence of the set of p-ends of $\torb$ to the set of p-ends of $\torb^*$, 
and we obtain that $x$ corresponds to a convex subset of $\dim \geq 1$ in $\Bd \torb$ containing a segment $l$ 
not contained in the closure of p-end neighborhoods using the map $\mathcal{D}$ in Proposition \ref{pr-prop-dualend}. 
Thus, the proof reduces to the case (I). 




By Theorem \ref{intro-thm-sSPC}, we obtain that our orbifold is strictly SPC.
By Corollary \ref{app-cor-stLens},
all the generalized lens R-end must also be lens R-ends.
\end{proof}




 Recall that the interior of a triangle has a Hilbert metric called the 
 hex metric by de la Harpe \cite{Harpe93}. The metric space 
 is isometric with a Euclidean space with norms given by regular hexagons. 
 The unit norm sphere of the metric is a regular hexagon ball for this metric. 
 A {\em regular hexagon of side length $l$} is a hexagon in the interior of 
 a triangle $T$ 
 with geodesic edges parallel to the sides of the unit norms  
 and with all edge lengths equal to $l$.  
 The regular hexagon is the boundary of a ball of radius $l$. 
 The {\em center} of a hexagon is the center of the ball. 
 \index{regular hexagon} 
 \index{Hex metric} 
 
\begin{lemma}\label{rh-lem-geogroup} 
Assume the premise of Theorem \ref{rh-thm-converse}.
Let $T$ be a triangle in $\tilde{\mathcal{O}^{(3)}}$ with $T^o \cap \torb^{(3)} \ne \emp$ and $\partial T \subset \Bd \torb^{(3)}$. 
Then $\partial T \subset \bigcup S(\tilde E)$ for an R-p-end $\tilde E$. 
\end{lemma}
\begin{proof}

Let $F$ be the fundamental domain of $\torb^{(3)}$. 



Again, we assume that $\pi_1(\orb)$ is torsion-free by Theorem \ref{prelim-thm-vgood}
since it is sufficient to prove the result for the finite cover of $\orb$. 
Hence, $\pi_1(\orb)$ acts freely on $\torb$.

Let $T'$ be a triangle with $T^{\prime o} \cap \torb^{(3)} \ne \emp$ and 
$\partial T' \subset \Bd \torb^{(3)}$. 
Suppose that $T'$ meets infinitely many horoball p-end neighborhoods in 
$\tilde{\mathbbm{U}}$ of horospherical p-ends, 
and the $d_{\orb^{(2)}}$-diameters of $T'$ intersected with these are not bounded. 
We consider a sequence of such sets $A_i:= T'\cap \tilde{\mathbbm{U}}_{j(i)}$ for some 
$j(i) \in I$ where 
the sequence of $d_{\orb^{(2)}}$-diameter $A_i$ is going to 
$+\infty$, and we choose a deck transformation $g_i$ such that $g_i(\clo(A_i))$ 
intersects
the fundamental domain $F$ of $\torb^{(4)}$. 
We choose a subsequence such that $\{g_i(T')\}$ and $\{g_i(A_i)\}$ geometrically converge
to a triangle $T''$ and a compact set $A_\infty$ respectively. 
Here, $T''$ intersects $F$ and the interior of $T''$ is in $\torb$. 
$g_i(A_i) = g_i(T'') \cap H_i$ for a horoball $H_i$ whose closure meets $F$. 
Since only finitely many closures of the horoball p-end neighborhoods in
$\torb$ meet $F$, 
there are only finitely many such $H_i$, say $H_{i_1}, \dots, H_{i_m}$.  
Now, $T''$ meets one such $H_{i_j}$ such that its vertex is in the boundary 
of $T''$ since the $d_{\orb^{(4)}}$-diameter of $g_i(T'')\cap H_i = g_i(A_i)$ goes to 
$+\infty$. 
This contradicts Theorem \ref{ce-thm-affinehoro}.


Thus, the $d_{\orb^{(2)}}$-diameters of horospherical p-end neighborhoods intersected with $L$ are bounded above uniformly. 
Therefore, by choosing a horospherical end neighborhood sufficiently far inside each horospherical end neighborhood by Corollary \ref{app-cor-shrink}, 
we may assume that $L$ does not 
meet any horospherical p-end neighborhoods. That is we choose a horoball $V'$ inside a one $V$
such that \[d_{\orb^{(2)}}(V', \partial V) >  \frac{1}{2} \sup \{ \hbox{$d_{\orb^{(2)}}$-diam}\{ V \cap T' \}|{V \in \mathcal{V}, T'\in \mathcal{T}}\}\]
where $\mathcal{V}$ is the collection of horoball p-end neighborhoods that we were given in the beginning and $\mathcal{T}$ is the collection of all triangles $T'$ meeting with $\torb^{(3)}$ and with boundary in $\Bd \torb^{(3)}$ --(*).

For $i$ in the index set $I$  of p-ends, we define $L_{1, i}$  to be the following 
subsets of $\torb^{(3)}$:
\begin{itemize} 
	\item $\clo(U(\mbv_{\tilde E})) \cap \torb^{(3)}$ where 
	$U(\mbv_{\tilde E})$ is the open shaved-off 
	concave p-end neighborhood of $\tilde E$
	when $\tilde E$ is a generalized lens-shaped R-p-end, 
	\item $\clo(\tilde S_{\tilde E}) \cap \torb^{(3)}$ if $\tilde E$ is a lens-shaped totally geodesic end, or 
	\item $\clo(U_{\tilde E}) \cap \torb^{(3)}$ for 
	a horoball $U_{\tilde E}$ for a horospherical end $\tilde E$.
\end{itemize}

By Theorem 1.5 of \cite{Drutu09}, $\pi_1(\orb)$ is relatively hyperbolic 
with respect to \[\pi_1(E_1), \dots, \pi_1(E_k)\] if and only if 
every asymptotic cone $\pi_1(\orb)$ is asymptotically tree graded 
with respect to the collection of left cosets of 
\[\mathcal{L} = \{g\pi_1(E_i)| g \in \pi_1(\orb)/\pi_1(E_i), i=1, \dots, k \}.\]
By Theorem 5.1 of \cite{DS08}, $\torb^{(4)}$ with the metric $d_{\torb^{(2)}}$ 
is asymptotically tree-graded with 
respect to $L_{1, i}, i=1,2, \dots,$ since $\pi_1(\orb)$ is quasi-isometric with 
$\torb^{(4)}$ with the cosets of $\pi_1(E_i)$ mapping quasi-isometric into $L_{1, j}$
by Proposition \ref{rh-prop-shaving}.  

Now, we consider $\torb^{(3)}$ and $\Bd \torb^{(3)}$.

For any $\theta > 0,\nu \geq 8$, a regular hexagon in 
$T^{\prime o}$ with side length $l > \nu \theta$ is 
$(\theta, \nu)$-fat according to Definition 5.1 of Dru{\c{t}}u \cite{Drutu09}. 
By Theorem 4.22 of \cite{Drutu09}, 
there is $\chi > 0$ such that a regular hexagon $H_l$ with side length $l > \nu \theta$ is in 
$\chi$-neighborhood $V_1$ of $L_{1, i}$ with respect to $d_{\orb^{(2)}}$ that is either contained in 
a concave p-end neighborhood of an R-p-end, a cocompactly-acted lens of a T-p-end or 
a horoball p-end neighborhood. 

Choose a family of regular hexagons 
\[ \{H_l| H_l \subset T^{'o}, l > \nu \theta\}\]
with a common center in $T^{'o}$.  
Hence,  $\bigcup_{l> \nu \theta} H_l = T^{' o} - K_1$ for a bounded set $K_1, 
K_1\subset T^{'o}$. 
By the above paragraph, $ T^{' o} - K_1 \subset V_1$. 
Now, $\partial T' \subset \Bd \torb^{(3)}$ must be in the closure of 
$L_{1, i}$ by Lemma \ref{rh-lem-chi}.

If $L_{i, 1}$ is from a horospherical p-end, 
$\partial T'$ must be a point. 
This is absurd. 
In the case of a T-p-end, the hyperspace containing $T'$ must coincide 
with one containing $\tilde S_{\tilde E}$. 
This is absurd since $T^{\prime o}$ is a subset of 
$\torb^o$.
In the case of an R-p-end $\tilde E$, it must be 
$\partial T' \subset \bigcup S(\mbv_{\tilde E})$. 
\end{proof}


\begin{lemma} \label{rh-lem-chi} 
	Let $V$ be a $\chi$-neighborhood of $L_{1, i}$ in $\torb^{(2)}$
	under the metric $d_{\orb^{(2)}}$. 
	Then $\clo(V) \cap \Bd \torb^{(2)} = \clo(L_{1, i}) \cap \Bd \torb^{(2)}$. 
	\end{lemma} 
\begin{proof} 
	We recall the metric. We first extend $\torb$ and shave off 
	to $\torb^e$. 
	Then we remove the parts of the lenses outside the ideal end orbifold for T-p-ends  and remove horoballs 
	of ends to obtain $\torb^{(4)}$. 
	
	Suppose that $L_{1, i}$ is from a horospherical p-end. 
	Then the equality is clear since $V$ is contained in a horospherical 
	p-end neighborhood.
	
	Suppose that $L_{1, i}$ is from a T-p-end of lens type. 
	Then there is a cocompactly-acted lens $L$ containing 
	$L_{1, i}$. The closure of $V$ in $\torb^{(3)}$ has a compact fundamental 
	domain $F_V$. Theorem \ref{pr-thm-equ2} and Lemma \ref{du-lem-attracting3} applied to any sequence of images of $F_V$ imply 
	the equality. 
	
	Suppose that $L_{1, i}$ is from an R-p-end of lens type. 
	Then the closure of $L_{1, i}$ in $\torb^{(3)}$ has a compact fundamental
	domain $F_V$. 
	Theorem \ref{pr-thm-equ}, Lemma \ref{pr-lem-attracting2}, and 
	Proposition \ref{pr-prop-orbit}
	again show the equality since the limit sets are independent of 
	the choice of neighborhoods. 
	\end{proof}

We recapitulate the results: 
\begin{corollary}\label{rh-cor-remhyp}
Assume that $\mathcal{O}$ is an $n$-dimensional strongly tame SPC-orbifold with generalized lens-shaped or horospherical $\cR$- or $\cT$-ends 
and satisfies {\rm (}IE\/{\rm )} and {\rm (}NA\/{\rm )}. Let \[E_1, \dots, E_m, E_{m+1}, \dots, E_k\] be the ends of $\orb$
where $E_{m+1}, \dots, E_k$ are some or all of the hyperbolic ends. 
Assume $\partial \orb =\emp$. 
Then $\pi_1(\mathcal{O})$ is a relatively hyperbolic group with respect to
the end groups $\pi_1(E_1),\dots, \pi_1(E_m)$ 
if and only if $\mathcal{O}_1$ as obtained by Theorem \ref{rh-thm-remconch} 
is strictly SPC with respect to the ends $E_1, \dots, E_m$. 
\end{corollary}
\begin{proof} 
If $\pi_1(\mathcal{O})$ is a relatively hyperbolic group with respect to 
the end groups \[\pi_1(E_1),\dots, \pi_1(E_m),\] 
then $\pi_1(\mathcal{O})$ is a relatively hyperbolic group with respect to 
the end groups $\pi_1(E_1),\dots, \pi_1(E_k)$ by Theorem \ref{rh-thm-drutu}. 
By Theorem \ref{rh-thm-converse}, it follows that $\mathcal{O}$ is strictly SPC with 
respect to the ends $E_1,\dots,E_k$. 
Theorem \ref{rh-thm-remconch} shows that $\mathcal O_1$ is strictly SPC with respect to 
$E_1, \dots, E_m$.  

For converse, if $\mathcal O_1$ is strictly SPC with respect to $E_1, \dots, E_m$, 
then $\mathcal O$ is strictly SPC with respect to $E_1, \dots, E_k$. 
By Theorem \ref{rh-thm-relhyp}, $\pi_1(\mathcal{O})$ is a relatively hyperbolic group with respect to 
the end groups $\pi_1(E_1),\dots, \pi_1(E_k)$. 
The conclusion follows from Theorem \ref{rh-thm-drutu}. 
\end{proof} 

\subsection{A topological result} 

\begin{proposition}\label{rh-prop-metriz}
Let $X$ be a compact metrizable space. Let ${\mathcal C}_X$ be a countable collection of mutually disjoint compact connected sets. 
The collection has the property that 
$C_K := \bigcup\{C \in {\mathcal C}_X| C\cap K \ne \emp\}$ is closed for any closed set $K$.
We define the quotient space $X/\sim$ with the equivalence relation 
$x \sim y$ iff $x, y \in C$ for an element $C \in {\mathcal C}_X$. 
Then $X/\sim$ is metrizable. 
\end{proposition} 
\begin{proof} 
We show that $X/\sim$ is Hausdorff, $2$-nd countable, and regular and 
use the Urysohn metrization theorem. 
We define as a basis
 a countable collection $\mathcal{B}$ of open sets of $X$ as follows:
We take an open subset $L$ of $X$ that is an $\eps$-neighborhood for
$\eps \in \bQ, \eps>0$, 
 of an element of ${\mathcal C}_X$ or a point of 
a dense countable set $Y$ in $X - \bigcup {\mathcal C}_X$. 
We form 
\[L - \bigcup\{C \in {\mathcal C}_X| C \cap \Bd L \ne \emp \}\]
for all such $L$ containing an element of ${\mathcal C}_X$
or a point of $Y$. This is an open set. 
The elements of $\mathcal{B}$ 
are neighborhoods of elements of ${\mathcal C}_X$ and $Y$.
Also, each element of ${\mathcal C}_X$ or a point of $Y$ 
is contained in an element of $\mathcal{B}$. 
Furthermore, each element of $\mathcal{B}$ is 
a saturated open set under the quotient map. 
Hence, $X/\sim$ is Hausdorff and $2$-nd countable. 

Now, the proof is reduced to showing that $X/\sim$ is regular. 
For any saturated compact set $K$ in $X$ and a disjoint element $Y$ of ${\mathcal{C}}_X$ or a point of $X$ not in any of ${\mathcal{C}}_X$, 
let $U_K$ and $U_Y$ denote the disjoint neighborhoods of $X$ of $K$ and $Y$ respectively. 
We form 
\[ U:= U_K - \bigcup\{C \in {\mathcal C}_X|C \cap \Bd U_K \ne \emp\},
\hbox{ and }
V := V_K - \bigcup\{C \in {\mathcal C}_X|C \cap \Bd V_K \ne \emp\}.  \]
Then these are disjoint saturated open neighborhoods for the basis $\mathcal{B}$. 
\end{proof}




\chapter{Openness and closedness}\label{ch-cl}

Lastly, we prove the openness and closedness of the properly (resp. strictly) convex real projective structures on the deformation spaces of a class of orbifolds with
generalized lens-shaped or horospherical $\cR$- or $\cT$-ends. 
We need the theory of Crampon and Marquis
and Cooper, Long, and Tillmann on the Margulis lemma for
convex real projective manifolds. 
The theory here partly generalizes that of Benoist on closed real projective orbifolds. 
To begin, we give some definitions and state the main results of the monograph: various Ehresmann-Thurston-Weil principles holding in 
certain circumstances. 
In Section \ref{cl-sec-openness}, we state the openness results to be proved in this chapter. Mainly, we use fixing sections to prove
the results here. 
We show that small deformations preserve convexity. 
Our approach is to work with radiant affine structures of one dimension higher by affine suspension construction and to
prove the results. Then the real projective versions follow easily. 
The idea is to use Hessian functions in
the suspensions of the compact part and to use the approximation of the original domain by
the covering domains of the end neighborhoods. 
In Section \ref{cl-sec-closed}, we show the closedness of the convexity 
under the deformations, first assuming the irreducibility of the holonomy representations. 
In Section \ref{cl-sub-sups}, we show that we do not actually need to assume the irreducibility
a priori. Any sequence of properly convex 
real projective structures converges to one whenever 
 the corresponding sequence of representations converges algebraically. 
In Section \ref{cl-sec-general}, we prove Theorem \ref{CL-THM-D}, the most 
general result of this monograph. Here, we show the natural existence of
the fixing section.

We now state our main results:


\begin{itemize}
\item 
We define $\Def^s_{{\mathcal E}, \lh}(\mathcal{O})$ as the subspace of $\Def_{{\mathcal E}}(\mathcal{O})$ 
consisting of 
real projective structures with generalized lens-shaped or horospherical $\cR$- or $\cT$-ends and stable irreducible holonomy homomorphisms.
\index{DefsEce@ $\Def^s_{{\mathcal E}, \lh}(\mathcal{O})$} 

\item We define 
$\CDef_{{\mathcal E}, \lh}(\mathcal{O})$ as the subspace of 
$\Def_{\mathcal E}(\mathcal{O})$ consisting of SPC-structures with generalized lens-shaped or horospherical $\cR$- or $\cT$-ends.
\index{CDefEce@$\CDef_{{\mathcal E}, \lh}(\mathcal{O})$|textbf} 

\item 
We define 
$\CDef_{{\mathcal E}, \mathrm{u}, \lh}(\mathcal{O})$ as the subspace of 
$\Def_{\mathcal E, \mathrm{u}}(\mathcal{O})$ consisting of SPC-structures with generalized lens-shaped or horospherical $\cR$- or $\cT$-ends.
\index{CDefEuce@$\CDef_{{\mathcal E}, \mathrm{u}, \lh}(\mathcal{O})$|textbf} 


\item We define 
$\SDef_{{\mathcal E}, \lh}(\mathcal{O})$ as the subspace of 
$\Def_{{\mathcal E}, \lh}(\mathcal{O})$ consisting of strict SPC-structures with lens-shaped or horospherical $\cR$- or $\cT$-ends.
\index{SDefEce@$\SDef_{{\mathcal E}, \lh}(\mathcal{O})$} 

\item We define 
$\SDef_{{\mathcal E}, \mathrm{u}, \lh}(\mathcal{O})$ as the subspace of 
$\Def_{{\mathcal E}, \mathrm{u}, \lh}(\mathcal{O})$ consisting of strict SPC-structures with lens-shaped or horospherical $\cR$- or $\cT$-ends.
\index{SDefEuce@$\SDef_{{\mathcal E}, \mathrm{u}, \lh}(\mathcal{O})$}
\end{itemize}

By Theorem \ref{intro-thm-sSPC}, some of these must have stable holonomy representations by topological conditions.
We remark that these strongly tame properly convex real projective $n$-orbifolds are dual to the same type of orbifolds, but we switch the $\mathcal{R}$-end with the $\mathcal{T}$-ends and vice versa.
by Proposition \ref{pr-prop-dualend}.
Also by Corollary \ref{app-cor-stLens}, for strongly tame strictly SPC-orbifolds with generalized 
lens-shaped or horospherical $\cR$- or $\cT$-ends have only lens-shaped or horospherical $\cR$- or $\cT$-ends. 

The following theorems are to be regarded as examples of the
so-called {\em Ehresmann-Thurston-Weil principle}, which identifies the deformation spaces with 
subspaces of the character spaces. 
\index{Ehresmann-Thurston-Weil principle} \index{hol@$\hol$}
\index{Ehresmann-Thurston principle}


\begin{theorem}\label{intro-thm-B} 
Let $\mathcal{O}$ be a noncompact strongly tame 
real projective $n$-orbifold, $n \geq 2$, 
	with generalized lens-shaped or horospherical $\cR$- or $\cT$-ends. 
Assume $\partial \orb =\emp$. 
Suppose that $\mathcal{O}$  satisfies 
{\rm (}IE\/{\rm )} and {\rm (}NA\/{\rm )}.   
Then 
the subspace  
\[\CDef_{{\mathcal E}, \mathrm{u}, \lh}(\mathcal{O}) \subset \Def^s_{\mathcal E, \mathrm{u}, \lh}(\mathcal{O})\]  is open.

Suppose further that every finite-index subgroup of $\pi_1(\mathcal{O})$ contains no nontrivial infinite nilpotent normal subgroup.
Then $\hol$ maps
$\CDef_{{\mathcal E}, \mathrm{u}, \lh}(\mathcal{O})$  homeomorphically to a union of components of 
 \[\rep_{{\mathcal E}, \mathrm{u}, \lh}(\pi_1(\mathcal{O}), \PGL(n+1, \bR)).\]
\end{theorem}


\begin{theorem} \label{intro-thm-C} 
Let $\mathcal{O}$ be a noncompact strongly tame $n$-dimensional strictly SPC-orbifold, $n \geq 2$, with lens-shaped or horospherical $\cR$- or $\cT$-ends 
and satisfies  {\rm (}IE\/{\rm )} and {\rm (}NA\/{\rm )}.  
Assume $\partial \orb =\emp$. 
Then
\begin{itemize}
\item $\pi_1(\mathcal{O})$ is relatively hyperbolic with respect to its end fundamental groups.
\item The subspace  $\SDef_{{\mathcal E}, \mathrm{u}, \lh}(\mathcal{O}) \subset \Def^s_{\mathcal E, \mathrm{u}, \lh}(\mathcal{O})$
of strict SPC-structures with lens-shaped or horospherical $\cR$- or $\cT$-ends is open. 
\end{itemize}
Suppose further that every finite-index subgroup of $\pi_1(\mathcal{O})$ contains no nontrivial infinite nilpotent normal subgroup.
Then $\hol$ maps the deformation space $\SDef_{{\mathcal E}, \mathrm{u}, \lh}(\mathcal{O})$ of
strictly SPC-structures on $\mathcal{O}$ with lens-shaped or horospherical $\cR$- or $\cT$-ends homeomorphically to
a union of components of \[\rep_{{\mathcal E}, \mathrm{u}, \lh}(\pi_1(\mathcal{O}), \PGL(n+1, \bR)).\]
\end{theorem}

Theorems \ref{intro-thm-B} and \ref{intro-thm-C} are proved
 by dividing the result into the openness result in Section \ref{intro-sub-open} 
and the closedness result in Section \ref{intro-sub-closed}.

\subsection{Openness} \label{intro-sub-open}

For the openness of  $\SDef_{{\mathcal E}, \lh}(\mathcal{O})$, 
we use Corollary \ref{rh-cor-remhyp}. 



\begin{theorem}\label{intro-thm-conv} 
Let $\mathcal{O}$ be a noncompact strongly tame real projective $n$-orbifold, $n \geq 2$,  
and satisfies  {\rm (}IE\/{\rm )} and {\rm (}NA\/{\rm )}.  Assume $\partial \orb =\emp$. 
In $\Def^s_{\mathcal E, \mathrm{u}, \lh}(\mathcal{O})$, the subspace  
$\CDef_{{\mathcal E}, \mathrm{u}, \lh}({\mathcal{O}})$ of 
SPC-structures with 
generalized lens-shaped or horospherical $\cR$- or $\cT$-ends is open, and 
so is $\SDef_{{\mathcal E}, \mathrm{u}, \lh}({\mathcal{O}})$.
\end{theorem}
\begin{proof}
	$\Hom_{\mathcal{E}, u}^s(\pi_1(\orb), \PGL(n+1, \bR)$ is 
	an open subset of $\Hom_{\mathcal{E}}(\pi_1(\orb), \PGL(n+1, \bR))$
	by Proposition \ref{intro-prop-u_semialg}. 
	On $\Hom_{\mathcal{E}, u}^s(\pi_1(\orb), \PGL(n+1, \bR))$, one has 
	a uniqueness section defined by Lemma \ref{intro-lem-u_section}. 
	Now, Theorem \ref{intro-thm-conv2} proves the result. 
	\end{proof}


We are given a properly real projective $n$-orbifold $\mathcal{O}$ with ends $E_1, \dots, E_{e_1}$ of $\cR$-type and 
$E_{e_1+1}, \dots, E_{e_1 + e_2}$ 
of $\cT$-type. Let us choose representative p-ends $\tilde E_1, \dots, \tilde E_{e_1}$ and $\tilde E_{e_1+1}, \dots, \tilde E_{e_1 + e_2}$ in their orbit classes under $\pi_(\orb)$.
Again, $e_1$ is the number of $\cR$-type ends,  and 
$e_2$ the number of $\cT$-type ends of $\orb$.

We define a subspace of 
$\Hom_{\mathcal E, \lh}(\pi_1(\mathcal{O}), \PGL(n+1, \bR))$ as in Section \ref{intro-sub-semialg}.

Let $\mathcal V$ be an open subset of a semi-algebraic subset of 
\[\Hom_{\mathcal E}(\pi_1(\mathcal{O}), \PGL(n+1,\bR))\] invariant under the conjugation action of $\PGL(n+1, \bR)$
such that the following hold:
\begin{itemize}
\item one can choose a continuous section $s_{\mathcal V}^{(1)}: \mathcal V \ra (\rpn)^{e_1}$
sending a holonomy homomorphism to a common fixed point of $\bGamma_{\tilde E_i}$ for $i = 1, \dots, e_1$ and 
\item  $s_{\mathcal V}^{(1)}$ satisfies 
 \[s_{\mathcal V}^{(1)}(g h(\cdot) g^{-1})  = g \cdot s_{\mathcal V}^{(1)}(h(\cdot)) \hbox{ for } g \in \PGL(n+1, \bR).\] \index{section|textbf} 
 \end{itemize} 
 $s_{\mathcal V}^{(1)}$ is said to be a {\em fixed-point section}.
 \index{section!fixed-point|textbf} 
 \index{fixed-point section|textbf}
 
 If $\tilde E_i$ for every $i=1, \dots, e_1$
 has a p-end neighborhood that has a radial foliation with leaves 
 developing into rays ending at the point of the $i$-th 
 factor of $s_{\mathcal V}^{(1)}$, 
 we say that the radial end structures are {\em determined by} 
 $s_{\mathcal V}^{(1)}$.
 \index{end!structure!determine} 

Again, we assume that $\mathcal V$ is an open subset of
a semi-algebraic subset.
\[\Hom_{\mathcal E, \lh}(\pi_1(\mathcal{O}), \PGL(n+1,\bR))\]
invariant under the conjugation action by $\PGL(n+1, \bR)$, and 
the following hold: 
\begin{itemize}
\item one can choose a continuous section $s_{\mathcal V}^{(2)}: \mathcal V \ra (\RP^{n\ast})^{e_2}$
sending a holonomy homomorphism to a common dual fixed point of $\bGamma_{\tilde E_i}$ for $i = e_1+1, \dots, e_{1}+e_2$, 
\item $s_{\mathcal V}^{(2)}$ satisfies $s_{\mathcal V}^{(2)}(g h(\cdot) g^{-1}) 
 = (g^*)^{-1}\circ s_{\mathcal V}^{(2)}(h(\cdot))$ for $g \in \PGL(n+1, \bR)$, and
\item letting $P_{\mathcal V}(\tilde E_{i})$ denote the null space of  the $i$-th value of $s_{\mathcal V}^{(2)}$ for 
$i= e_1+1, \dots, e_{1} + e_2$,  
$\bGamma_{\tilde E_i}$ acts on 
the hyperspace $P_{\mathcal V}(\tilde E_{i})$ satisfying the lens-condition for 
$\tilde E_i$. 
 \end{itemize}
 $s_{\mathcal V}^{(2)}$ is said to be a {\em dual fixed-point} section.
 \index{dual fixed-point section|textbf} 
 \index{section!fixed-point!dual|textbf}

If each $\tilde E_i$ for every $i= e_1+1, \dots, e_1+ e_2$ 
\begin{itemize} 
\item has a p-end neighborhood with the ideal boundary component in the hyperspace 
determined by the $i$-th factor of $s_{\mathcal V}^{(2)}$ provided 
$\tilde E_i$ is a T-end, or 
\item has a p-end neighborhood containing a $\bGamma_{\tilde E}$-invariant 
horosphere tangent to the hyperspace 
determined by the $i$-th factor of $s_{\mathcal V}^{(2)}$ provided 
$\tilde E_i$ is a horospherical end, 
\end{itemize} 
 we say that end structures for the totally geodesic ends 
 are {\em determined by} $s_{\mathcal V}^{(2)}$.
 
 We define $s_{\mathcal V}: \mathcal V \ra (\RP^n)^{e_1} \times (\RP^{n\ast})^{e_2}$
 as $ s_{\mathcal V}^{(1)} \times s_{\mathcal V}^{(2)}$ and call it a {\em fixing section}.
 \index{fixing section|textbf} 
 \index{section!fixing|textbf}


\begin{lemma} \label{intro-lem-u_section} 
We can define a section 
\[s_u: \Hom_{\mathcal{E}, \mathrm{u}, \lh}(\pi_1(\orb), \PGLnp)  \ra (\RP^n)^{e_1} \times (\RP^{n \ast})^{e_2} \] by choosing for each holonomy and each p-end 
the unique fixed point and the unique hyperspace as the images. 
\end{lemma} 
\begin{proof} 
$s_u$ is a continuous function since a sequence of fixed points or 
dual fixed points of end holonomy group is a fixed point or 
a dual fixed point of the limit end-holonomy group. 
\end{proof}
We call $s_u$ the {\em uniqueness section}.

Let $\mathcal V$ and $s_{\mathcal V}: {\mathcal{V}} \ra  (\RP^n)^{e_1} \times (\RP^{n \ast})^{e_2}$
 be as above.
\begin{itemize}
\item We define $\Def^s_{{\mathcal E}, s_{\mathcal V}, \lh}(\mathcal{O})$ to be the subspace of 
$\Def_{{\mathcal E}, s_{\mathcal V}}(\mathcal{O})$ 
of real projective structures with generalized lens-shaped or horospherical $\cR$- or $\cT$-end structures determined by $s_{\mathcal V}$,
and stable irreducible holonomy homomorphisms in $\mathcal V$. 

 \label{fixing section} 
\item We define 
$\CDef_{{\mathcal E}, s_{\mathcal V}, \lh}(\mathcal{O})$ to be the subspace consisting of 
SPC-structures with generalized lens-shaped or horospherical $\cR$- or $\cT$-end structures determined by $s_{\mathcal V}$ 
and holonomy homomorphisms in $\mathcal V$
 in $\Def^{s}_{{\mathcal E}, s_{\mathcal V}, \lh}(\mathcal{O})$. 


\index{CDefEsUce@$\CDef_{{\mathcal E}, s_{\mathcal V}, \lh}(\mathcal{O})$|textbf}

\item We define 
$\SDef_{{\mathcal E}, s_{\mathcal V}, \lh}(\mathcal{O})$ to be the subspace of consisting of strictly SPC-structures
with lens-shaped or horospherical $\cR$- or $\cT$-end structures determined by $s_{\mathcal V}$ 
and holonomy homomorphisms in $\mathcal V$
 in $\Def^{s}_{{\mathcal E}, s_{\mathcal V}, \lh}(\mathcal{O})$. 

\end{itemize} 
\index{SDefEsUce@$\SDef_{{\mathcal E}, s_{\mathcal U}, \lh}(\mathcal{O})$}

\begin{theorem}\label{intro-thm-conv2} 
Let $\mathcal{O}$ be a noncompact strongly tame real projective $n$-orbifold, 
$n \geq 2$, 
with generalized lens-shaped or horospherical $\cR$- or $\cT$-ends
and satisfies  {\rm (}IE\/{\rm )} and {\rm (}NA\/{\rm )}.  Assume $\partial \orb =\emp$. 
Choose an open $\PGL(n+1, \bR)$-conjugation invariant subset of a
semialgebraic subset
 \[\mathcal{V} \subset \Hom_{\mathcal E, \lh}(\pi_1(\orb), \PGL(n+1, \bR)),\]
and 
a fixing section
$s_{\mathcal{V}}: {\mathcal{V}} \ra (\RP^n)^{e_1} \times (\RP^{n \ast})^{e_2}$.

Then $\CDef_{{\mathcal E}, s_{\mathcal V}, \lh}(\mathcal{O})$ is open 
 in $\Def^s_{\mathcal E, s_{\mathcal V}, \lh}(\mathcal{O})$,
 and so is $\SDef_{{\mathcal E}, s_{\mathcal V}, \lh}(\mathcal{O})$.
 \end{theorem}

This is proved in Theorem \ref{cl-thm-conv2}.



By Theorems \ref{intro-thm-conv} and \ref{intro-thm-conv2}, we obtain:
\begin{corollary}\label{intro-cor-conv} 
Let $\mathcal{O}$ be a noncompact strongly tame real projective $n$-orbifold, $n \geq 2$, with generalized lens-shaped or horospherical $\cR$- or $\cT$-ends and 
satisfies  {\rm (}IE\/{\rm )} and {\rm (}NA\/{\rm )}.  Assume $\partial \orb =\emp$. 
Then 
\[\hol: \CDef_{{\mathcal E}, \mathrm{u}, \lh}(\mathcal{O}) \ra \rep^{s}_{\mathcal E, \mathrm{u}, \lh}(\pi_1(\mathcal{O}), \PGL(n+1, \bR))\] is a local homeomorphism. 

Furthermore, if $\mathcal{O}$ has a strictly SPC-structure with lens-shaped or horospherical $\cR$- or $\cT$-ends, then so is
\[\hol: \SDef_{{\mathcal E}, \mathrm{u}, \lh}(\mathcal{O}) \ra \rep^{s}_{\mathcal E, \mathrm{u}, \lh}(\pi_1(\mathcal{O}), \PGL(n+1, \bR)).\]
\end{corollary}


%

\subsection{The closedness of convex real projective structures}\label{intro-sub-closed}

The results here are proved in Chapter \ref{ch-cl}. 
We recall  \[\rep^s_{\mathcal E}(\pi_1(\mathcal{O}), \PGL(n+1, \bR))\] that is the subspace of stable irreducible characters
of \[\rep_{\mathcal E}(\pi_1(\mathcal{O}), \PGL(n+1, \bR)).\]
This is shown to be the open subset of a semi-algebraic subset in Section \ref{intro-sub-semialg}, 
and denote by $\rep^s_{{\mathcal E}, \mathrm{u}, \lh}(\pi_1(\mathcal{O}), \PGL(n+1, \bR))$ the subspace of stable irreducible characters
of $\rep_{{\mathcal E}, \mathrm{u}, \lh}(\pi_1(\mathcal{O}), \PGL(n+1, \bR))$,  a union of open subsets of semialgebraic sets. 





\begin{theorem} \label{intro-thm-closed1} 
Let $\mathcal{O}$ be an $n$-dimensional noncompact strongly tame SPC-orbifold, $n \geq 2$, with 
generalized lens-shaped or horospherical $\cR$- or $\cT$-ends and satisfies  {\rm (}IE\/{\rm )} and {\rm (}NA\/{\rm )}.  Assume $\partial \orb =\emp$, and 
that the nilpotent normal subgroups of every finite-index subgroup of $\pi_1(\mathcal{O})$ are trivial.
Then the following hold\,{\rm :} 
\begin{itemize}
\item The deformation space $\CDef_{{\mathcal E}, \mathrm{u}, \lh}(\mathcal{O})$ of SPC-structures on $\mathcal{O}$ with generalized lens-shaped or horospherical $\cR$- or $\cT$-ends maps under $\hol$
homeomorphically to a union of components of $\rep_{{\mathcal E}, \mathrm{u}, \lh}(\pi_1(\mathcal{O}), \PGL(n+1, \bR))$.
\item 
The deformation space $\SDef_{{\mathcal E}, \mathrm{u}, \lh}(\mathcal{O})$ of strictly SPC-structures on $\mathcal{O}$ with lens-shaped or horospherical $\cR$- or $\cT$-ends maps under $\hol$
homeomorphically to the union of components of  $\rep_{{\mathcal E}, \mathrm{u}, \lh}(\pi_1(\mathcal{O}), \PGL(n+1, \bR))$.
\end{itemize}
\end{theorem}
\begin{proof} 
	$\Hom_{\mathcal{E}, u}^s(\pi_1(\orb), \PGL(n+1, \bR)$ is 
	an open subset of $\Hom_{\mathcal{E}}(\pi_1(\orb), \PGL(n+1, \bR)$
	by Proposition \ref{intro-prop-u_semialg}. 
Corollary \ref{cl-cor-closed1} proves this by the existence of 
the uniqueness section
of Lemma \ref{intro-lem-u_section}.  
	\end{proof}

The following is probably the most general result. 


\begin{theorem}\label{CL-THM-D}
	Let $\mathcal{O}$ be an $n$-dimensional noncompact strongly tame SPC-orbifold, $n \geq 2$, with generalized lens-shaped or horospherical $\cR$- or $\cT$-ends 
	and satisfies  {\rm (}IE\/{\rm )} and {\rm (}NA\/{\rm )}.  
	Assume $\partial \orb =\emp$. 
	Then
	\begin{itemize}
		\item  Suppose  that every finite-index subgroup of $\pi_1(\mathcal{O})$ contains no nontrivial infinite nilpotent normal subgroup
		and $\partial \orb =\emp$. 
		Then $\hol$ maps the deformation space $\CDef_{{\mathcal E}, \lh}(\mathcal{O})$ of 
		SPC-structures on $\mathcal{O}$ with generalized lens-shaped or horospherical $\cR$- or $\cT$-ends homeomorphically to 
		a union of components of 
		\[\rep_{{\mathcal E}, \lh}(\pi_1(\mathcal{O}), \PGL(n+1, \bR)).\]
		\item  Suppose  that every finite-index subgroup of $\pi_1(\mathcal{O})$ contains no nontrivial infinite nilpotent normal subgroup
		and $\partial \orb =\emp$. 
		Then $\hol$ maps the deformation space 
		$\SDef_{{\mathcal E}, \lh}(\mathcal{O})$ of 
		strictly SPC-structures on $\mathcal{O}$ with lens-shaped or horospherical $\cR$- or $\cT$-ends homeomorphically to 
		a union of components of \[\rep_{{\mathcal E}, \lh}(\pi_1(\mathcal{O}), \PGL(n+1, \bR)).\]
	\end{itemize}
\end{theorem}

For example, these apply to the projective deformations of
hyperbolic manifolds with the torus boundary
as in \cite{BDL}.

\section{The openness of the convex structures}\label{cl-sec-openness}

In this section also, we only need $\RP^n$ versions. 
Given a strongly tame real projective $n$-orbifold $\orb$ with $e_1$  $\cR$-ends and $e_2$  $\cT$-ends, each end $E_i$, $i=1, \dots, e_1, e_1+1, \dots, e_1+e_2$, has an orbifold structure of dimension $n-1$ and inherits a real projective structure. 

\index{uniqueness section} 

\begin{theorem}\label{cl-thm-conv2} 
Let $\mathcal{O}$ be a strongly tame real projective $n$-orbifold with generalized lens-shaped or horospherical $\cR$- or $\cT$-ends 
and satisfies  {\rm (}IE\/{\rm )} and {\rm (}NA\/{\rm )}.  Assume that $\partial \orb =\emp$.
Let $\mathcal U$ be a $\PGL(n+1, \bR)$-conjugation invariant
open subset of a union of semi-algebraic subsets of
 \[ \Hom^s_{\mathcal{E}, \lh}(\pi_1(\orb), \PGL(n+1, \bR)).\]
We have a $\PGL(n+1, \bR)$-equivariant fixing section 
$s_{\mathcal{U}}: {\mathcal{U}} \ra (\RP^n)^{e_1} \times (\RP^{n \ast})^{e_2}$. 
Then the following are open subspaces satisfying
\begin{gather*} 
\CDef_{{\mathcal E}, s_{\mathcal U}, \lh}(\mathcal{O}) \subset 
 \Def^s_{{\mathcal E}, s_{\mathcal U}, \lh}(\mathcal{O}), \\
\SDef_{{\mathcal E}, s_{\mathcal U}, \lh}(\mathcal{O})\subset 
 \Def^s_{{\mathcal E}, s_{\mathcal U}, \lh}(\mathcal{O}).
\end{gather*} 
 \end{theorem}

For orbifolds such as these, the deformation space of properly convex structures 
may only be a proper subset of the space of the characters.

By Theorem \ref{cl-thm-conv2} and Theorem \ref{op-thm-projective}, 
 we obtain:
\begin{corollary} \label{cl-cor-conv2}
Let $\mathcal{O}$ be a strongly tame real projective $n$-orbifold 
with generalized lens-shaped or horospherical $\cR$- or $\cT$-ends and satisfies  {\rm (}IE\/{\rm )} and {\rm (}NA\/{\rm )}.  Assume that $\partial \orb= \emp$. 
Let $\mathcal{U}$ and $s_{\mathcal{U}}$ be as in Theorem \ref{cl-thm-conv2}.
Suppose that $\mathcal{U}$ has its image ${\mathcal{U}}'$ in 
\[ \rep^s_{{\mathcal E}, \lh}(\pi_1(\mathcal{O}), \PGL(n+1, \bR)).\] 
Then 
\[\hol: \CDef_{{\mathcal E}, s_{\mathcal U}, \lh}(\orb) \ra \mathcal{U'}\] is a local homeomorphism, and so is
\[\hol: \SDef_{{\mathcal E}, s_{\mathcal U}, \lh}(\orb) \ra \mathcal{U'}.\]
\end{corollary}
\begin{proof} 
	Theorem \ref{op-thm-projective} shows that the map 
	\[\hol:\Def^s_{\mathcal{E}, s_{\mathcal V}}(\mathcal{O}) \ra \rep^s_{\mathcal{E}}(\pi_1(\mathcal{O}), \PGL(n+1, \bR))\]
	is an open one. 
Proposition \ref{cl-prop-lenshoro} tells us the openness of the images here. 
Theorem \ref{cl-thm-conv2} completes the proof. 
\end{proof}
Here, in fact, one needs to prove for every possible continuous section. 

Koszul  \cite{Koszul68} proved these facts for closed  affine manifolds, and  Goldman \cite{Goldman90} expanded these
for the closed real projective manifolds. See \cite{dgorb}, \cite{CLM18} and also Benoist \cite{Benoist05}.



\subsection{The proof of the openness}\label{cl-sec-open} 

The major part of showing the preservation of convexity under 
deformation is Proposition \ref{cl-prop-openess} on Hessian function perturbations. 

We mention that our approach to openness is slightly different from
that of Cooper-Long-Tillmann \cite{CLT15} since they 
are using their canonical invariant Hessian metrics for end neighborhoods.
Our Hessian metrics for end neighborhoods are not canonical,
as theirs are. 

Recall that a convex open cone $V$ is a convex cone of $\bR^{n+1}$ containing the origin $O$ in the boundary. 
Recall that a properly convex open cone is a convex cone such that its closure does not contain a pair of $\vec{v}, -\vec{v}$ for a non-zero vector in $\bR^{n+1}$. Equivalently, it does not contain a complete affine line in its interior.

A dual convex cone $V^*$ to a convex open cone is a subset of $\bR^{n+1 *}$ given by condition $\phi \in V^*$ if and only if $\phi(\vec{v})> 0$ for all $\vec{v} \in \clo(V) -\{O\}$.

Recall that $V$ is a properly convex open cone if and only if so is $V^*$ and $(V^*)^* = V$ 
under the identification $(\bR^{n+1 *})^* = \bR^{n+1}$. 
Also, if $V \subset W$ for a properly convex open cone, then $V^* \supset W^*$.

For a properly convex open subset $\Omega$ of $\SI^n$, its dual $\Omega^*$ in $\SI^{n*}$ 
is given by taking a cone $V$ in $\bR^{n+1}$ corresponding to $\Omega$
and taking the dual $V^*$ and projecting it to $\SI^{n *}$. The dual $\Omega^*$ is a properly convex open domain if so was $\Omega$.

Recall the Koszul-Vinberg function for a properly convex cone $V$ and the dual properly convex cone $V^*$ 
\begin{equation}\label{cl-eqn-kv}
f_{V^*}: V \ra \bR_+ \hbox{ defined by }  x \in V \mapsto f_{V^*}(x)= \int_{V^*} e^{-\phi(x)} d\phi
\end{equation}
where the integral is over the Euclidean measure in $\bR^{n+1 *}$. 
\index{Koszul-Vinberg function} 
This function is strictly convex if $V$ is properly convex. $f_{V^*}$ is homogeneous of degree $-(n+1)$. 
Writing $D$ as an affine connection, we write the Hessian $D d\log(f)$. 
The Hessian is positive definite and the norms of the unit vectors are
strictly bounded below in a compact subset $K$ of $V - \{O\}$.
(See Chapter 4 of \cite{goldmanbook} and \cite{Vinberg63}.) 
The metric $D d\log(f)$ is invariant under the group $\Aff(V)$ of affine transformation that acts on $V$. 
In particular, it is invariant under scalar dilatation maps. 
(For an extensive survey, see Shima \cite{Shima07}.)





A {\em Hessian metric} on an open subset $V$ of an affine space 
is a metric of form $\partial^2 f /\partial x_i \partial x_j$ for affine coordinates $x_i$ \index{Hessian metric}
and a function $f: V \ra \bR$ with a positive definite Hessian defined on $V$. 
A Riemannian metric on an affine manifold is a {\em Hessian metric} if the manifold is affinely covered by a cone 
and the metric lifts to a Hessian metric of the cone. 

Let a strongly tame $n$-orbifold 
$\mathcal{O}$ have an SPC-structure $\mu$ with generalized lens-shaped or horospherical $\cR$- or $\cT$-ends.Clearly, $\tilde{\mathcal{O}}$ is a properly convex open domain. 
Then an affine suspension of $\mathcal{O}$ has an affine Hessian metric defined by $D d\phi$ for a function $\phi$ defined on 
the cone in $\bR^{n+1}$ corresponding to $\tilde{\mathcal{O}}$ by above. 

A {\em parameter of real projective structures} $\mu_t, t\in [0, 1]$ on a strongly tame orbifold $\orb$
is a collection such that the restriction $\mu_t|K$ to each compact suborbifold $K$ is a continuous parameter; 
In other words, the associated developing map $\dev_t|\hat K:\hat K \ra \SI^n$   (resp. $\RP^n$)  
for every compact subset $\hat K$ of $\torb$ 
is a continuous family in the $C^r$-topology for the variable $t$. (See
Definition \ref{op-defn-deforms}, Choi \cite{dgorb} and Canary-Epstein-Green \cite{CEG06}.)
\index{parameter of real projective structures} 


\begin{definition} \label{cl-defn-closestr} 
	Let $\orb$ be a strongly tame orbifold with ends. 
	Let $U$ be a union of mutually disjoint end neighborhoods, 
	and let $\mu_0$ and $\mu_1$ be two real projective structures on $\orb$. 
	Let $\overline{\dev}_0, \overline{\dev}_1:\hat \orb \ra \SI^n$ be extended developing maps. 
	We say that $\mu_0$ and $\mu_1$ 
	on $\orb$ are {\em $\delta$-close} in the $C^r$-topology, $r \geq 2$, on 
	the compactification $\bar \orb$ if
	for a compact path-connected domain
	$K$ in $\hat \orb$ mapping onto $\bar \orb$, 
	the associated developing maps $g_0 \circ \overline{\dev}_0| K$ and 
	$g_1\circ \overline{\dev}_1| K$ are $\delta$-close in $C^r$-topology
	for some $g_0$ and $g_1$ in $\PGL(n+1, \bR)$. 
\end{definition} 
Of course, this definition requires us to define $K$ first.

%


Recall the Vinberg metric from Sections 4.4 and 12.3 of  Goldman \cite{goldmanbook}, which is a Hessian metric. 

\begin{proposition}\label{cl-prop-openess}
Let $\mathcal{O}$ be a strongly tame $n$-orbifold with ends and satisfies  {\rm (}IE\/{\rm )} and {\rm (}NA\/{\rm )}. 
Suppose that $\mathcal{O}$ has an SPC structure $\mu_0$ with generalized lens-shaped or horospherical $\cR$- or $\cT$-ends
and the affine suspension of $\mathcal{O}$ with $\mu_0$ has a Hessian metric. 
{\rm (}See Section \ref{prelim-sub-asusp}.{\rm )}
The ends of $\mathcal{O}$ are given $\cR$-type or $\cT$-types. 
Suppose that one of the following holds\/{\rm :} 
\begin{itemize}
\item $\mu_0$ is SPC, and a $C^r$-continuous parameter, 
$r \geq 2$, of real projective structure $\mu_{t}$, $t\in [0, 1]$, radial or totally geodesic ends with end holonomy groups of 
generalized lens-shaped or horospherical $\cR$- or $\cT$-ends 
where the $\cR$-types or $\cT$-types of ends are preserved.
\item We assume that 
$\{\mu_t\} \ra \mu_0$ as $t \ra 0$ in the $C^r$-topology for $r \geq 0$. 
\end{itemize}
Then for sufficiently small $t$, the affine suspension
$C(\torb)$ for $\torb$ with $\mu_t$ also has a Hessian metric invariant under the group of dilatations. 
\end{proposition} 
\begin{proof}\renewcommand{\qedsymbol}{}
We prove for $\SI^n$. It is sufficient since we aim to obtain
the Hessian metric on $C(\torb)$.
We keep $\torb$ and the action of the deck transformation group fixed and only change the structures on it. 
Note that the subsets here remain fixed, and the only changes are on the real projective structures, 
i.e., the atlas of charts to $\SI^n$. 


Let $\tilde{\mathcal{O}}$ in $\SI^n$ denote the universal covering domain corresponding to $\mu_0$. 
Again, $\dev_0$ being an embedding identifies the first with subsets of $\SI^n$
but $\dev_t$ is not known to be so. 



We prove this by steps: 
\begin{enumerate}
	\item[(A)] The first step is to understand the deformations of the end neighborhoods.
	\item[(B)] We change the Hessian function on the cone associated with 
	the universal covers. We need to obtain one for the deformed end neighborhoods
	by Hessian functions from Koszul-Vinberg integrals 
	and another one to the outside 
	of the union of end neighborhoods by isotopies, and 
    we need to patch the two together.  
\end{enumerate}

(A) Let $\tilde E'$ be a p-end of $\torb$, and it also corresponds to a p-end of $\torb'$. Let $E$ be the end of $\orb$ corresponding to $\tilde E$. 
There exists a $C^r$-parameter of real projective structures $\mu_{t}$ 
with generalized lens-shaped or horospherical $\cR$- or $\cT$-ends. 
We can also find a parameter of developing maps $\dev_t$ associated with $\mu_{t}$
where $\dev_t|K$ is continuous with respect to the variable $t$ for each compact 
$K \subset \hat \orb$. 
To begin with, we assume that $\tilde E'$ continues to be a lens-shaped or horospherical p-end. 

Let $h_t$ denote the holonomy homomorphism associated with $\dev_t$ for each $t$. 
Recall that
Theorems \ref{app-thm-qFuch} and \ref{app-thm-qFuch2} study
the perturbation of the lens-shaped R-ends and T-ends.  
Lemma \ref{op-lem-horob} studies the perturbations of 
a horospherical $\cR$- or $\cT$-end to either a R-end or to a T-end. 

In this monograph, 
we do not allow $\cR$-type ends to change to $\cT$-type ends and vice versa,
as this violates the local injectivity property from
the deformation space to
a space of characters. 
(See Theorem \ref{op-thm-projective}.)
Let $\tilde E$ be a p-end. 
We consider the case where the type changes only once since 
the conditions form semialgebraic sets. 
Thus, we need to consider only four cases to prove the openness. 
\begin{description} 
\item[(I)] $h_t(\pi_1(\tilde E))$ changes from 
the holonomy group of a radial p-end to that of a radial p-end in the cases: 
\begin{description} 
\item[(a)] $h_t(\pi_1(\tilde E))$ changes from the holonomy group of a radial p-end of generalized lens-shaped becoming that of a generalized lens-shaped radial p-end.
\item[(b)] $h_t(\pi_1(\tilde E))$ changes from the holonomy group of a horospherical p-end to that of a generalized lens-shaped radial p-end or a horospherical p-end. 
\end{description} 
\item[(II)] $h_t(\pi_1(\tilde E))$ changes from 
the holonomy group of totally geodesic ends of lens type or horospherical $\cR$- or $\cT$-ends changes to that of themselves here. 
\begin{description} 
\item[(a)] $h_t(\pi_1(\tilde E))$ changes from 
the holonomy group of a lens-shaped totally geodesic p-end to that of 
a lens-shaped totally geodesic p-end. 
\item[(b)] $h_t(\pi_1(\tilde E))$ changes from 
the holonomy group of a horospherical p-end to that of a horospherical p-end or to a lens-shaped totally geodesic p-end.

\end{description} 
\end{description}
These hold for the corresponding 
holonomy homomorphisms of the fundamental groups of ends by the premise. 
(The above happens in actuality as well. See \cite{Ballas14},\cite{Ballas15}, and 
\cite{BDL}.) 


We now work on one end at a time: Fix a p-end $\tilde E$ of the R-type of $\tilde{\mathcal{O}}$. 
Let $\mbv$ be the p-end vertex of $\tilde E$ for $\mu_0$ and $\mbv'$ that for $\mu_1$. 
We denote by $\mbv = \mbv_0$ and $\mbv'= \mbv_1$. 
Assume that $\mbv_t$ is the p-end vertex of $\tilde E$ for $\mu_t$. 
Let $\dev_t$ and $h_t$ denote the developing map and the holonomy homomorphism of $\mu_t$. 
Assume first that the corresponding p-end for $\mu$ is of radial or horospherical type. 
By postcomposing the developing map by a transformation near the identity,
we assume that the perturbed vertex $\mbv_t$ of the corresponding p-end $\tilde E$
is mapped to $\mbv_0$, that is, $\mbv = \dev_t(\mbv_t)$.

(I) 
Suppose that $\tilde E$ is a generalized lens-shaped radial p-end or a horospherical p-end for $\mu_0$. Then the holonomy group of $\tilde E$ is that of a generalized lens-shaped radial p-end or a horospherical p-end for $\mu_t$ under (I). 





Let $\Lambda_0$ denote the limit set in the tube of the radial p-end $\tilde E$ for $\torb$ if $\tilde E$ is lens-shaped radial p-end, 
or $\{\mbv_{\tilde E}=\mbv\}$ if $\tilde E$ is a horospherical type for $\mu_0$.  (See Definition \ref{app-defn-limitset}.)

\begin{itemize}
\item Recall that ${R}_{\mbv}(\dev_t(\torb))$ denotes the space of directions of
segments from $\mbv$ in $\dev_t(\torb)$,
\item ${R}_{\mbv}(\dev_t(A_t))$ denotes the space of directions of segments from 
$\mbv$ of $\dev_t(\torb)$ in  ${R}_{\mbv}(\torb)$ passing through the set $\dev(A_t) \subset \dev_t(\torb)$.

\end{itemize}

(i) First, we find domains $\Omega_{s_0}$ with smooth boundary approximating
$\torb$ on which $h(\pi_1(\tilde E))$ acts. 
Here, $s_0$ is a parameter that we use to vary $\Omega_{s_0}$ for
fixed holonomy representations. 
Here, $\Omega_{s_0}$ is a lens cone for $\torb$ with $\mu_0$ 
such that $\partial \Omega_{s_0} = \Bd_{\torb} \Omega_{s_0}$ which 
is the top boundary component of the lens. 


Suppose that $\tilde E$ is a horospherical R-p-end for $\mu_0$. 
Then we obtain $\partial \Omega_{s_0}$ as the boundary of
a convex domain invariant under $h_0(\pi_1(\tilde E))$. 
Again, $\Omega_{s_0} \cup \{\mbv\}$ bounds a properly convex domain
$\Omega_{s_0}$. 


Now, we study the remaining case where $\tilde E$ is a generalized lens-type p-end. 
Choose some small $\epsilon> 0$. 
By Lemma \ref{app-lem-expand},  if $\tilde E$ is a generalized 
lens-shaped R-p-end or a horospherical R-p-end, 
we may choose $\Omega_{s_0} \subset \torb$ such that 
\[\bdd_H(\Omega_{s_0}, \torb) < \eps.\]

Suppose that $\tilde E$ is merely a generalized lens-shaped R-p-end. 
We can show that 
\[\clo(\torb) \cap \mathcal{T}_{\mbv_{\tilde E}}(\tilde \Sigma_{\tilde E})
=\Bd \torb - \bigcup S(\mbv_{\tilde E})
\] 
is a generalized strict lens cone.
Then we can still form a lens cone neighborhood $V$ in 
$\mathcal{T}_{\mbv_{\tilde E}}(\tilde \Sigma_{\tilde E})$ for $\tilde E$. 
We can choose a smooth top boundary component $\Sigma'$ of 
a lens in $V$ that $\eps'$-approximates 
$\Bd \torb - \bigcup S(\mbv_{\tilde E})$ with respect to $d_{V^o}$ inside $V^o$
by Proposition \ref{pr-prop-convhull2} for some sufficiently small $\eps'> 0$.
We now find a new domain $\Omega'_{s_0}$ deformed 
from $\Omega_{s_0}$ as a convex component of $V^o- \Sigma'$. 
We can show that  
for any $\eps > 0$, sufficiently well-chosen $\Omega'_{s_0}$ satisfies
 \[\bdd_H(\Omega'_{s_0}, \torb) < \eps\]
as follows: 
Since $V$ is the interior of a strict lens cone,  
for every $\eps''>0$, we can choose 
a compact set $K$ in $\tilde \Sigma_{\tilde E}$, 
\[\Sigma' - \mathcal{T}_{\mbv_{\tilde E}}(K) \hbox{ and } 
\Bd \torb \cap \mathcal{T}_{\mbv_{\tilde E}}(\tilde \Sigma_{\tilde E}) - \mathcal{T}_{\mbv_{\tilde E}}(K)\]
is $\eps''$-approximate with respect to $\bdd$
by Proposition \ref{pr-prop-orbit}
since $\Lambda_0$ is independent of the choice of lens. 
We choose $\Sigma'$ such that 
$\Sigma' \cap  \mathcal{T}_{\mbv_{\tilde E}}(K)$ $\eps''$-approximates
$\Bd \torb \cap \mathcal{T}_{\mbv_{\tilde E}}(K)$ with respect to $\bdd$
by choosing above $\eps'$ sufficiently small.

Here, 
new $\partial \Omega'_{s_0}$ may not be a subset of $\torb$ unless 
$\tilde E$ is lens-shaped and not just generalized lens-shaped. 
However, this is irrelevant for our purposes. 
From now on, we set $\Omega_{s_0}$ as $\Omega'_{s_0}$. 

For each $\epsilon>0$, we choose the parameter $s_0$ such that 
the above is satisfied. So, $s_0$ is considered to vary for our purposes. 

(ii) Now our purpose is to find a convex domain
$\Omega_{s_0, t}$ where $h_t(\pi_1(E))$ acts on
and approximates $\Omega_{s_0}$, and 
to show that it contains an embedded image of 
a p-end neighborhood. 
We denote by $\tilde \Sigma_{\tilde E, t}$ the universal cover of 
the end orbifold associated with a R-p-end $\tilde E$ of $\torb$ with $\mu_t$. 
By our end holonomy group condition in the premise, 
Corollary \ref{abelian-cor-convhol} 
shows that $\tilde \Sigma_{\tilde E, t}$ is 
again complete affine or properly convex. 


For the $\mathcal{R}$-type end $E$, 
$\orb$ has a concave end neighborhood or a horospherical 
end neighborhood for $E$ bounded by a smooth compact end orbifold $S'_E$ 
transverse to the radial rays. 
Here, $S'_E$ is clearly diffeomorphic to $\Sigma_E$. 

For a sufficiently small $t$ in $\mu_{t}$, we obtain a domain 
$U_t \subset \orb$ 
with $U_0 \subset \Omega_{s_0}/h_0(\pi_1(\tilde E))$ 
bounded by an inverse image 
of a compact orbifold $S'_{E, t}$ diffeomorphic to $S'_E$
still transversal to radial rays 
by Propositions \ref{app-prop-koszul} and  \ref{app-prop-lensP}.
$S'_{E, t}$ is either strictly concave if $S_E$ was a properly convex R-p-end, and
strictly convex if $\Sigma_E$ was horospherical. 
(The strict convexity and the transversality follow since changes in affine connections can be made arbitrarily small, as in the argument of Koszul \cite{Koszul68}.)

Since the change was sufficiently small, we may assume that 
$S'_{E, t}$ still bounds an end neighborhood $U_t$ of a product form
by Lemma \ref{cl-lem-perturbend}. 

The proof continues after Lemma \ref{cl-lem-perturbend}.  \cpr
\end{proof} 

\begin{lemma} \label{cl-lem-perturbend} 
	Suppose that $\tilde U$ is an R-p-end neighborhood of $\torb$ 
	covering an end neighborhood $U$ of an end $E$ in $\orb$. 
	Then for sufficiently small change of real projective structures in 
	$C^r$-sense, $r \geq 2$, in the compact open topology, 
	hypersurface $S_t$ sufficiently close to $S$ in the $C^r$-sense 
	in terms of parameterizing maps
	still bounds an end neighborhood of $E$. 
	Letting $\tilde S_t$ be a component of the inverse image of $S$ 
	on which $\pi_1(\tilde E)$ acts, we still have that 
	$\tilde S_t$ bounds a p-end neighborhood of $\tilde E$. 
	\end{lemma} 
\begin{proof} 
	Straightforward. 
	\end{proof}

\begin{proof}[Proof of Proposition \ref{cl-prop-openess} continued]  \renewcommand{\qedsymbol}{}
Let $\Lambda_t$ denote the limit set in $\bigcup S(\mbv)_t$ for generalized radial p-end cases, and 
$\Lambda_t =\{\mbv\}$ for the horospherical case. 
Let $S(\mbv)_t$ denote the set of maximal segments in the closure of $U_t$ 
from $\mbv$ corresponding to  $\Bd {R}_{\mbv}(\dev_t(\torb))$ of $\mu_t$. 




Suppose that $\tilde E$ is a lens-shaped R-p-end for $\mu_0$. 
We showed above that 
the $C^r$-change $r \geq 2$ of $\mu_{t}$ from $\mu_0$ be sufficiently small such that 
we obtain a region $\Omega_{s_0, t}$ in $B_t^o$ with $\partial \Omega_{s_0, t}$ strictly convex
and transverse to radial rays under $\dev_t$. Here, $\Omega_{s_0, 0} = \Omega_{s_0}$. 

Choose a compact domain $F$ in $\partial \Omega_{s_0}$. 
Let $F_t$ denote the corresponding deformed set in $\partial \Omega_{s_0, t}$.
By Theorem \ref{ce-thm-mainaffine}, $\pi_1(\tilde E)$ is virtually abelian
if $\omega_{s_0}$ was horospherical. 
For sufficiently small $t$, $0 < t <1$, $\dev_t(F_t)$ is a subset of the tube $B_t$ determined by $\dev_t(U_t)$
since $B_t$ and a parameterization of $\dev_t(F_t)$ continuously depends on $t$
by Corollary \ref{abelian-cor-smvar}.
\begin{itemize} 
\item By transversality to the segments mapping to ones from $\mbv$ under $\dev_t$, it follows that 
$\dev_t|\partial \Omega_{s_0, t}$ gives us a smooth immersion to a convex domain $\tilde \Sigma_{\tilde E, t}$ 
that equals the space of maximal segments in $B_t$ with vertices $\mbv$ and $\mbv_-$. 
\item By Corollary \ref{abelian-cor-convhol},  
the immersion
$\dev_t| \partial \Omega_{s_0, t}$ to a properly convex 
domain $\tilde \Sigma_{\tilde E, t}$ is a diffeomorphism 
if $\tilde \Sigma_{\tilde E, t}$ is properly convex. 
It is also so 
if $h_t(\pi_1(\tilde E))$ is horospherical
since this follows from the classical Bieberbach theory
using the Euclidean metric on $\tilde \Sigma_{\tilde E}$ 
where it develops. 

\end{itemize} 
Since $\pi_1(\tilde E)$ is the fundamental group of a generalized lens-shaped or horospherical p-end, it follows that 
\[\dev_t(\clo(\partial \Omega_{s_0, t}) - \partial \Omega_{s_0, t}) \subset \dev_t(\Lambda_t)\]
by Theorems 
\ref{ce-thm-affinehoro},
\ref{pr-thm-lensclass}, and \ref{pr-thm-redtot}. 


Suppose that $\tilde E$ is a generalized lens-shaped or lens-shaped R-p-end for $\mu_t$, $t > 0$. 
Since $\partial \Omega_{s_0, t}$ is convex, each point of $\partial \Omega_{s_0, t} \cup \bigcup S(\mbv)_t$ 
has a neighborhood that is mapped under completion $\widehat{\dev}_t$ to a convex open ball. 
Thus, $\partial \Omega_{s_0, t} \cup \bigcup S(\mbv)_t$ 
bounds a compact ball $\Omega_{s_0, t} \cup \bigcup S(\mbv)_t$ by Lemma \ref{prelim-lem-locconv} 
since the local convexity implies the global convexity
and they are in $B_t$. 

Suppose that $\tilde E$ is a horospherical R-p-end. For $t > 0$, 
$\partial \Omega_{s_0, t} \cup \{\mbv\}$ bounds a convex domain $\Omega_{s_0, t}$
by the local convexity of the boundary set $\partial \Omega_{s_0, t} \cup \{\mbv\}$ and Lemma \ref{prelim-lem-locconv}.
\cpr
\end{proof} 


\begin{proposition} \label{cl-prop-lenshoro} 
	Assume as in Proposition \ref{cl-prop-openess}, and {\rm (}I\/{\rm )} in the proof. 
	Then for $\mu_t$ for sufficiently small $t$, 
	the end corresponding to $E$ is always 
a 	generalized lens-shaped R-end or a horospherical R-end. 
	Also, if $\mu'$ is sufficiently $C^r$-close to 
	$\mu$, then the end of $\orb$ with $\mu'$ 
	is a generalized lens-shaped R-end or a horospherical R-end.  
\end{proposition}
\begin{proof} 
	The above arguments prove this since we are studying arbitrary deformations. 
\end{proof} 

\begin{proof}[Proof of Proposition \ref{cl-prop-openess} continued] 
(iii) We show how these regions deform approximating $\Omega_{s_0}$ 
in the Hausdorff metric sense. 
We define ${\mathbbm{r}}_{\mbv}(K)$ the union of great segments 
with an endpoint $\mbv$ in directions of $K$, $K \subset R_{\mbv}$.
\index{rv@${\mathbbm{r}}_{\mbv}(\cdot)$}
\begin{itemize}
\item Let $K$ be a compact convex subset of  $\Omega_{s_0, 0}$ with smooth boundary,
and $K_t$ the perturbed one in  $\Omega_{s_0, t}$
and $\tilde E$ be the corresponding p-end. We can form a compact set inside 
$\Omega_{s_0, t}$ consisting 
of segments from the p-end vertex to $K$ in the set of radial segments. 
For $\mu_{t}$ from $\mu_0$ changed by a sufficiently small manner in the $C^r$-topology, 
a compact subset $\mathbbm{r}_{\mbv}(K) \subset \mathbbm{r}_{\mbv}(\torb)$ is changed to a compact convex domain 
$\mathbbm{r}_{\mbv}(K_t) \subset \mathbbm{r}_{\mbv}(\tilde \Sigma_{\tilde E, t})$. 

\end{itemize} 




We choose $s_0$ large enough so that $K \subset \Omega_{s_0}^o$. 
For sufficiently small $t$, 
$\Omega_{s_0, t} \cap \mathbbm{r}_{\mbv}(K_t)$ is a convex domain since $\partial \Omega_{s_0, t}$ is strictly convex and 
transverse to great segments from $\mbv$ and therefore 
embeds itself to a convex domain under $\dev_t$. 
We may assume that $\Omega_{s_0, t} \cap \mathbbm{r}_{\mbv}(K_t)$ is sufficiently close to $\Omega_{s_0} \cap \mathbbm{r}_{\mbv}(K)$ 
as we changed the real projective structures sufficiently small in the $C^r$-sense.
See Definition \ref{cl-defn-closestr}. 
%

An {\em $\eps$ thin} space is a space that is an $\eps$ neighborhood of its boundary for small $\eps> 0$.
By Lemma \ref{intro-lem-actionRend}  and Corollary \ref{abelian-cor-smvar}, we may assume that $\clo(\mathbbm{r}_{\mbv}(\torb))$ and $\clo(\mathbbm{r}_{\mbv, t}(\torb))$
are $\eps$-$\bdd$-close convex domains in the Hausdorff sense for sufficiently small $t$.
Thus, 
given an $\eps > 0$, we can choose $K$ and $K'_t$ and a sufficiently small 
deformation of the real projective structures such that 
$\Omega_{s_0}\cap (\mathbbm{r}_{\mbv}(\torb)-\mathbbm{r}_{\mbv}(K))$ 
is an $\eps$-thin space, and so is 
$\Omega_{s_0, t}\cap (\mathbbm{r}_{\mbv}(\torb)-\mathbbm{r}_{\mbv}(K_t))$ 
for sufficiently small changes of $t$.
Moreover, 
\begin{align} \label{cl-eqn-thin}
\clo( \Omega_{s_0}) \cap (\mathbbm{r}_{\mbv}(\torb)-\mathbbm{r}_{\mbv}(K)) \subset N_\eps(\clo( \Omega_{s_0} \cap \mathbbm{r}_{\mbv}(K))) \hbox{ and } \nonumber \\ 
\clo(\Omega_{s_0, t}) \cap (\mathbbm{r}_{\mbv}(\torb)-\mathbbm{r}_{\mbv}(K_t))\subset N_\eps(\clo(\Omega_{s_0, t} \cap \mathbbm{r}_{\mbv}(K_t)));  
\end{align} 
We may assume that the diameters of $\clo(\Omega_{s_0}), \clo(\Omega_{s_0, t})$ are all less than $\pi/2$ by projective automorphisms if necessary.
The reason is that the sharply supporting hyperspaces of $\clo(\Omega_{s_0})$ at points of $\partial \mathbbm{r}_{\mbv}(K) \cap \clo(\Omega_{s_0})$
are in acute angles with the geodesics from ${\mbv}$ and similarly for those of $\clo(\Omega_{s_0, t})$ for sufficiently small $t$
by Corollary \ref{pr-cor-tangency}. 

Therefore, we conclude for (I) that
for any $\eps > 0$, there exists $\delta > 0$ 
\begin{equation} \label{cl-eqn-omegat1} 
\bdd_H(\clo(\Omega_{s_0}), \clo(\Omega_{s_0, t})) < \eps 
\end{equation} 
provided $|t| < \delta$; that is,
 we choose $\mu_{0, t}$ sufficiently close to $\mu_0$:
First, we choose $K$ and the deformation $K_t$ such that
it satisfies \eqref{cl-eqn-thin} for $t < \delta$ for some $\delta > 0$. 
Then we choose $t$ to be sufficiently
small such that $\Omega_{s_0, t} \cap {R}(K_t)$ is
sufficiently close to $\Omega_{s_0}\cap {R}(K)$. 

Also, $\Omega_{s_0, t}$ contains a concave p-end neighborhood of 
$\tilde E$ for $\mu_t$ for sufficiently small $t>0$. 
This can be assured by taking sufficiently large $K$ containing 
a fundamental domain of $R_x(\torb)$ and
sufficiently large $K_t$ containing a fundamental domain $F_t$ of 
$R_t(\dev_t(\orb))$ for $h_t(\pi_1(\tilde E))$ deformed from $F$
by Proposition \ref{cl-prop-lenshoro}.  


(II) Now suppose that $\tilde E$ is a lens-shaped T-p-end or horospherical p-end of type $\mathcal{T}$, 
and we suppose that $h_t(\tilde E)$ is a lens-shaped T-p-end for $\mu_t$ for $t> 0$.
Other cases are similar to (I). 


We take $\Omega_{s_0}$ to be the convex domain obtained as in Lemma \ref{app-lem-expand} with strictly convex boundary 
component $\partial_1 \Omega_{s_0}$ and totally geodesic one 
$\tilde S_{\tilde E, 0}$ in $\Bd \Omega_{s_0}$. 
Now, $\Omega_{s_0}/h(\pi_1(\tilde E))$ has 
strictly convex boundary component $\partial \Omega_{s_0}/h(\pi_1(\tilde E))$
and totally geodesic boundary $\tilde S_{\tilde E, 0}/h(\pi_1(\tilde E))$
when $\tilde E$ is a T-p-end. 
If $\tilde E$ is a horospherical p-end, then
$\partial \Omega_{s_0}/h(\pi_1(\tilde E))$ is still a strictly convex compact $(n-1)$-orbifold. 

Here, we note 
$\Bd \Omega_{s_0} = \partial \Omega_{s_0}\cup \clo(\tilde S_{\tilde E})$. 

Suppose that $\mu_{t}$ is sufficiently close to $\mu_0$. 
Then by Theorem \ref{app-thm-qFuch2}, we deform the lens-shaped 
T-p-end for $\tilde E$: 
\begin{itemize} 
\item we obtain a properly 
convex domain $\Omega_{s_0, t}$ for sufficiently small $t$ 
with a strictly convex boundary 
$\partial_1 \Omega_{s_0, t}$ and $\tilde S_{t}$, 
and 
\item $\partial \Omega_{s_0, t}$ is also
cocompact under the $\pi_1(\tilde E)$-action associated with $\mu_1$ and strictly convex.

\end{itemize}
See Proposition \ref{app-prop-koszul}.
We choose $\Omega_{s_0}$ to be $L_0\cap \torb$ for a lens in an ambient 
orbifold containing $L_0$. 
We also have 
$\Bd \Omega_{s_0, t} = \partial \Omega_{s_0, t}\cup \clo(\tilde S_{\tilde E, t})$
where $\tilde S_{\tilde E, t}$ is the ideal boundary component for 
$\tilde E$ for $\mu_t$. 
By Proposition \ref{app-prop-koszul}, 
$L_0$ deforms to a properly convex domain $L_t$ such that 
$L_t \cap \torb_t$ is $\Omega_{s_0, t}$ where $\torb_t$ is $\torb$ with 
a real projective structure $\mu_t$. 
We have 
\[\clo(\partial \Omega_{s_0, t}) - \partial \Omega_{s_0, t}= 
\partial \clo(\tilde S_{\tilde E, t})\] 
for a totally geodesic ideal boundary 
component $\tilde S_{\tilde E, t}$ by 
Theorem \ref{du-thm-lensn}.
Therefore, the union of $\partial \Omega_{s_0, t}$ and a totally geodesic ideal boundary component
$\clo(\tilde S_{\tilde E, t})$ bounds a properly convex compact $n$-ball in $\SI^n$.
We can find a properly convex lens $L_0$ for $\tilde E$ at $\mu_0$ with $L_0 \cap 
\torb$ in a p-end neighborhood $\Omega_{s_0}$.
Since the change of $\mu_t$ is sufficiently small, 
$L_t \cap \torb_t$ is still in a p-end neighborhood of $\tilde E$
by Lemma \ref{cl-lem-perturbend}. 
We obtain a lens-shaped p-end neighborhood for $\tilde E$ and $\orb$ with $\mu_t$ by 
Theorem \ref{du-thm-lensn}. 


We may assume without loss of generality that 
the hyperspace $V_{\tilde E}$ containing $\tilde S_{\tilde E, t}$ is fixed. 
Moreover, we may assume by choosing a sufficiently large $\Omega_{s_0}$ without loss of generality that
$\bdd_H(\Omega_{s_0}, \torb) < \eps$ for any $\eps> 0$.

By Lemma \ref{pr-lem-attracting2}, a sharply supporting hyperspace at 
a point of $\Bd \tilde S_{\tilde E}$ is 
uniformly bounded away from $V_{\tilde E}$. A sequence of sharply
supporting hyperspaces can converge to a sharply supporting one. 
Let us choose a sufficiently small $\eps > 0$. 
Let $B$ be a compact $\eps$-neighborhood of $\partial \Omega_{s_0}$ such that 
\[ \bdd^H(\partial \Omega_{s_0} - B, 
\partial \clo(\tilde S_{\tilde E})) < \eps. \]
Given a sharply supporting hyperspace $W_x$ of a point $x$ of $\partial \Omega_{s_0}$
of $\Omega_{s_0}$, there exists 
a sharply supporting closed hemisphere $H_x$ bounded by $W_x$. 


We define the {\em shadow} $S$ of $\partial B$ as the set 
\[ \bigcap_{x \in \partial B} H_x \cap V_{\tilde E}. \]
Then we can choose sufficiently small $\eps$ such that 
$\bdd^H(S, \dev_0(\clo(\tilde S_{\tilde E}))) \leq \eps$. 
We can also ensure that $W_x$ meets $V_{\tilde E}$ in the angles in $(\delta, \pi-\delta)$ for some $\delta > 0$
by compactness of $\partial B$ and continuity of map $x \mapsto W_x$. 

Suppose that we change the structure from $\mu_0$ to $\mu_t$ with a small $C^2$-distance. 
Then $B$ changes to $B'_t$ with $W_x$ changing by a small amount. The new shadow $S'_t$ will 
have the property $\bdd^H(S'_t, \clo(\tilde S_{\tilde E, t})) \leq \eps$ for a sufficiently small $C^r$-change, $r\geq 2$, of $\mu_t$ from $\mu_0$. 
Hence, we obtain 
that for each $\epsilon > 0$, there exists $\delta> 0$ such that 
\begin{equation} \label{cl-eqn-omegat2} 
\bdd_H(\partial \Omega_{s_0, t} - B'_t, 
\partial \clo(\tilde S_{\tilde E, t})) < \eps 
\end{equation} 
provided $|t| < \delta$. 
Therefore, by Corollary \ref{abelian-cor-smvar} for
each $\epsilon > 0$, there exists $\delta> 0$ such that 
\[  \bdd_H(\clo(\tilde S_{\tilde E}), \clo(\tilde S_{\tilde E, t})) < \eps.   \]
Recall that 
$\Bd \Omega_{s_0} = \partial \Omega_{s_0}\cup \clo(\tilde S_{\tilde E})$ 
and $\Bd \Omega_{s_0, t} = \partial \Omega_{s_0, t}\cup \clo(\tilde S_{\tilde E, t})$.
Combining with \eqref{cl-eqn-omegat2} and 
the sufficiently small change of $B$ to $B_t$, we obtain that for
each $\epsilon > 0$, there exists $\delta> 0$ such that 
\begin{equation} \label{cl-eqn-omegat3} 
\bdd_H(\clo(\Omega_{s_0}), \clo(\Omega_{s_0, t})) < \eps 
\end{equation} 
provided $|t| < \delta$. 

Suppose that $\tilde E$ was a horospherical p-end of type $\mathcal{T}$. 
Again, the argument is similar. We start with $\Omega_{s_0}$, which is horospherical with a strictly convex boundary and
tangent to a hyperspace $P_0$. 
The deformation gives us a strictly convex hypersurface $\partial \Omega_{s_0, t}$ 
and a hyperspace $P_t$ where $h_t(\pi_1(\tilde E))$ acts. 

We assume without loss of generality that $P_t=P$ for small $t > 0$. 
For sufficiently small $t$, 
we obtain a domain bounded by $\partial \Omega_{s_0, t}$ and 
the closure $\clo(\tilde S_{\tilde E, t})$ 
of a totally geodesic ideal boundary component. 
Taking duals by Corollary \ref{pr-cor-duallens2}, 
we have a lens-shaped or a horospherical R-p-end $\tilde E$.
Proposition \ref{cl-prop-lenshoro} shows that 
we obtain a domain $\Omega_{s_0, t}^\ast$ approximating $\Omega_{s_0}^\ast$
where $\partial \Omega_{s_0, t}^\ast$ is dual to 
$\partial \Omega_{s_0, t}$. 
Then $\Omega_{s_0, t}^{\ast \ast}$ approximates 
$\Omega_{s_0}$ as much as one wishes to
by Lemma \ref{cl-lem-convHaus}.  

We show that we have an embedded image of a p-end neighborhood in 
$\Omega_{s_0, t}$. 
The hypersurface $\partial \Omega_{s_0}$ is embedded in $\torb$ with $\mu_0$. 
Let $F$ be a compact fundamental domain of $\partial \Omega_{s_0}$ by 
$h(\pi_1(\tilde E))$. 
Now, $\partial \Omega_{s_0, t}$ 
contains a compact fundamental domain $F_t$ perturbed 
from $F$. Let $\tilde S_t$ denote a component of the inverse image of
$\partial \Omega_{s_0, t}$ under $\dev_t$ containing the perturbed
fundamental domain $F_t$ deformed from $F$. 
We deduce that  $\pi_1(\tilde E)$ acts on $\tilde S_t$ 
since $h_t(\pi_1(\tilde E))$ acts on $\partial \Omega_{s_0, t}$. 

By the above, $h_t(\pi_1(\tilde E))$ acts on $\partial \Omega_{s_0, t}$ properly
and cocompactly giving us a closed orbifold as a quotient. 
Since $\dev_t|F_t$ is an embedding for sufficiently small $t$, 
the equivariance tells us that $\dev_t: \tilde S_t \ra \partial \Omega_{s_0, t}$ 
is a diffeomorphism. 

Since the change is sufficiently small, $\tilde S_t$ still bounds 
a p-end neighborhood of $\orb$ by Lemma \ref{cl-lem-perturbend}. 

Before going on with part (B) of the proof, we briefly do
a slight generalization. 
\cpr
\end{proof} 

\begin{corollary} \label{cl-cor-dh} 
	We consider all cases {\rm (}I\/{\rm),} {\rm (}II\/{\rm).} 
	For each $\epsilon > 0$, there exists $\delta> 0$
	depending only on $\clo(\Omega_{s_0})$ such that we can choose 
	a convex domain $\Omega_{s_0, 1}$ where $h_1(\pi_1(\tilde E))$ acts on 
	for the holonomy homomorphism $h_1$ 
	such that 
	\begin{itemize} 
		\item $\Omega_{s_0, 1}$ contains a domain that is the embedded
		image of a p-end neighborhood of 
		$\tilde E$ by a developing map $\dev_1$ for $\mu_1$ associated with $h_1$ and 
		\item 
	\begin{equation} \label{cl-eqn-omegat4} 
	\bdd_H(\clo(\Omega_{s_0}), \clo(\Omega_{s_0, 1})) < \eps 
	\end{equation} 
	\end{itemize} 	
	provided $\mu_0$ and $\mu_1$ are $\delta$-close in $C^r$-topology, $r\geq 2$, 
	on the compact set $\orb - U$ for a union of $U$ of end neighborhoods of $\orb$. 
	\end{corollary} 
\begin{proof} 
		Suppose that this is false. There exists a sequence of real projective structures 
	$\mu_{t_i}$ with $\mu_{t_i} \ra \mu_0$ in the $C^2$-topology on $\orb -U$.
	Then letting the associated holonomy of $\mu_{t_i}$ be denoted by $h_{t_i}$,  
	$\{h_{t_i}\}$  converges to $h_0$
	for $\mu_0$. 
	We can apply the argument of cases (I), (II),
	to show that \eqref{cl-eqn-omegat4} holds for every $\eps >0$. 
	\end{proof}


\begin{proof}[Proof of Proposition \ref{cl-prop-openess} continued] 
(B) With $\torb$ with $\mu_t$, we obtain a special affine suspension on $\orb \times \SI^1$ with the affine structure $\hat \mu_t$. 
Let $C(\torb)$ be the cone over $\torb$. Then this covers the special affine suspension. 
Let $\tilde \mu_t$ denote the affine structure on $C(\torb)$ corresponding to $\hat \mu_t$. 
For each $\mu_t$, it has an affine structure $\tilde \mu_t$, different from the induced structure of $\bR^{n+1}$ as for $t = 0$. 
We recall the scalar multiplication 
\[s\cdot \mbv = s\mbv, \mbv \in C(\torb), s \in \bR\] 
for any affine structure $\tilde \mu_t$. 
Also, given a subset $K$ of $\torb$, we denote by $C(K)$ the corresponding set in $C(\torb)$. 
This set is independent of $\tilde \mu_t$ but has different affine structures nearby. 

For $\mu_0$, $\torb$ is a domain in $\SI^n$. 
Recall the Koszul-Vinberg function $f: C(\torb) \ra \bR_+$ homogeneous of degree $-n-1$ as given by \eqref{cl-eqn-kv}. 
(See Lemma \ref{cl-lem-KV}.) \index{Koszul-Vinberg function} 
By our choice above, 
the Hausdorff distance between $\clo(\Omega_{s_0})$ and $\tilde{\mathcal{O}}$ can be made arbitrarily small as
desired for some choice of $0$. 

By the proof of Theorem \ref{op-thm-projective} 
 in Page \pageref{op-page-thm-affine}
constructing the local inverse maps applied to strongly tame orbifolds with boundary, 
there exists a diffeomorphism 
\[F_{t}: \Omega_{s_0}/h_0(\pi_1(E)) \ra \Omega_{s_0, t}/h_t(\pi_1(E)) 
\hbox{ with a lift } \tilde F_t: \Omega_{s_0} \ra \Omega_{s_0, t}\]
such that $\tilde F_t \ra \Idd$ on every compact subset of $C(\torb)$ 
in the $C^r$-topology as $t \ra 0$. 
That is, on every compact subset $K$ of $C(\torb)$, 
\[\{\llrrV{D^j \tilde F_t - D^j \Idd | K}\} \ra 0 \hbox{ for every multi-index } j, 0 \leq |j| \leq r.\]  
We may assume that $\tilde F_t$ commutes with the radial flow $\hat \Psi_s : \torb \ra \torb$ for $s \in \bR$
by restricting $F_t$ to a cross-section of $C(\torb)$ of the radial flow 
and extending radially. 
(See the paragraph after Definition \ref{prelim-defn-affs}.)

By the third item of  Lemma \ref{cl-lem-convHaus} and Lemma \ref{cl-lem-KV},
the Hessian functions $f'_t\circ \tilde F_t$ defined by \eqref{cl-eqn-kv} 
on the inverse image $\tilde F_t^{-1}(C(\Omega_{s_0, t})^o)$
is arbitrarily close to the original Hessian function $f$
in any compact subset of $C(\torb)$ in the $C^r$-topology, $r \geq 2$, 
provided $|t - t_0|$ is sufficiently small. 
By construction, $f'_t$ is homogeneous of degree $-n-1$.

The holonomy groups $h(\pi_1(\mathcal{O}))$ and $h_t(\pi_1(\tilde E))$ 
being in $\SL_{\pm}(n+1, \bR)$ preserve $f$ and $f'_t$ under 
deck transformations respectively.

Now do this for all p-ends.
Let $\mathbbm{U}_t$ be the $\pi_1(\orb)$-invariant mutually disjoint union  of p-end neighborhoods of p-ends of $\torb$.  
 We construct a function $f'_t$ on $C(\mathbbm{U}_t)$
for $\mu_t$ and sufficiently small $|t|$.

Let $\mathbbm{U}$ be the corresponding $\pi_1(\orb)$-invariant union of  proper p-end neighborhoods of $\torb$ for $\mu_0$. 
For each component $\mathbbm{U}_i$ of $\mathbbm{U}$, 
we construct $f'_{i, t}\circ \tilde F_i$ on $C(\mathbbm{U}_i)$ 
using $\Omega_{s_0}$ such that $f'_{i, t}$ satisfies the above properties
and $\tilde F_i$ is constructed as above for $\mathbbm{U}_i$. 
We call $f'_t$ the union of these functions. 

Let $\mathbbm{V}$ be a $\pi_1(\orb)$-invariant compact 
neighborhood of the complement of 
$\mathbbm{U}$ in $\torb$. 
%
\begin{itemize}
\item Let $\partial_s \mathbbm{V}$ be the image of $\partial {\mathbbm{V}}\times \{s\}$ inside the regular neighborhood of $\partial {\mathbbm{V}}$ in $\mathbbm{U}$
parameterized as $\partial {\mathbbm{V}} \times [-1, 1]$ for $s \in [-1, 1]$.
\item We assign $\partial {\mathbbm{V}}= \partial_0 {\mathbbm{V}}$. 
\item Let $\partial_{[s_1, s_2]} {\mathbbm{V}}$ denote the image of $\partial {\mathbbm{V}}_t\times \{[s_1, s_2]\}$ inside 
the regular neighborhood of $\partial {\mathbbm{V}}$ in ${\mathbbm{V}}\cap \mathbbm{U}'$ for a neighborhood $\mathbbm{U}'$ of $\clo(\mathbbm{U}) \cap \torb$.
  
\end{itemize} 
We find a $C^\infty$ map $\phi_t: C(\mathbbm{U}')\cap C({\mathbbm{V}}) \ra \bR_+$ such that $\phi_t(s\orar{v})=\phi_t(\orar{v})$ for every $s >0$
and $f'_t(\orar{v})=\phi_t(\orar{v}) f(\orar{v})$ and $\phi_t$ is very close to the constant value $1$ function.
By making $f'_t/f$ near $1$ and the derivatives of $f'_t/f$ up to two near $0$ as possible, 
we obtain $\phi_t$ that has derivatives up to order two arbitrarily close to $0$ in a compact subset:
This is accomplished by taking a partition of unity functions $p_{1}, p_{2}$ 
invariant under the radial flow such that 
\begin{itemize}
\item $p_{1}=1$ on $C(W)$ for 
\[W:= \partial_{[0, s_1]} {\mathbbm{V}} \cup (\mathbbm{U}'-{\mathbbm{V}})  
\hbox{ for } s_1 < 1,\]
\item $p_{1} = 0$ on $C(\torb -N)$ for a neighborhood $N$ of $W$ 
in $\partial_{(-1, 1)} {\mathbbm{V}} \cup (\mathbbm{U}'-{\mathbbm{V}})$, and 
\item $p_{1}+p_{2}=1$ identically.

\end{itemize}  
We assume that \[1-\eps < f'_t/f < 1+\eps \hbox{ in } C({\mathbbm{U}}' \cap {\mathbbm{V}}),\]
and $f'_t/f$ has derivatives up to order two sufficiently close to $0$
by taking $f'_t$ and $f$ sufficiently close in $C({\mathbbm{U}}') \cap C({\mathbbm{V}})$ by taking sufficiently small $t$. 
We define \[\phi_t = (f'_t/f-(1-\eps))p_{1}+ \eps p_{2} + (1-\eps), 0< t< 1\] 
as $f'_t$ and $f$ are homogeneous of degree $-n-1$. 
Then $1-\eps < \phi_t < 1+ \eps$ and derivatives of $\phi_t$ up to order two are sufficiently close to $0$ by taking sufficiently small $\eps$ 
as we can see easily from computations. 
Thus, using $\phi_t$, we obtain a function $f$ 
obtained from $f'_t$ and $\phi_t f$ on $C(W)$ 
extending them smoothly for sufficiently small $|t|$.

We can check the welded function from $f'_t$ and $\phi_t f$ has the desired Hessian properties for $\mu_t$ for sufficiently small $t$ 
since the derivatives of $\phi_t$ up to order two can be made sufficiently close to zero. 
Now we do this for every p-end of $\tilde{\mathcal{O}}$.


The $-(n+1)$-homogeneity gives us the invariance of the Hessian metric under the scalar dilatations and the affine lifts of the holonomy groups. 
(See Chapter 4 of \cite{goldmanbook}.)  
This completes the proof for Proposition \ref{cl-prop-openess}. 
\hfill \SSn {\parfillskip0pt\par}
\end{proof}


This is a strengthened version of Proposition \ref{cl-prop-openess}. 

\begin{corollary}\label{cl-cor-openness}
	Let $\mathcal{O}$ be a strongly tame $n$-orbifold with ends and satisfies  {\rm (}IE\/{\rm )} and {\rm (}NA\/{\rm )}. 
	Suppose that $\mathcal{O}$ has an SPC-structures $\mu_0$ with generalized lens-shaped or horospherical $\cR$- or $\cT$-ends
	and the suspension of $\mathcal{O}$ with $\mu_0$ has a Hessian metric. 
	Let $U$ be a union of mutually disjoint end neighborhoods of $\orb$. 
	Suppose the following hold\/{\rm :} 
	\begin{itemize}
		\item Let $\mu_1$ be an SPC-structure 
		with generalized lens-shaped or horospherical $\cR$- or $\cT$-ends 
		such that $\mu_0$ and $\mu_1$ are $\eps$-close in 
		$C^r$-topology on $\orb-U$ for $r \geq 2$, and  
		\item the R-end holonomy group of $\mu_1$ is either lens-type or 
		horospherical type and the T-end holonomy group of 
		$\mu_1$ are totally geodesic satisfying the lens conditions.  
	\end{itemize}
	Then for sufficiently small $\eps$, the affine suspension
	$C(\torb)$ for $\torb$ with $\mu_1$ also has a Hessian metric invariant under dilations and the 
	affine suspensions of the holonomy homomorphism for $\mu_{1}$. 
\end{corollary}
\begin{proof} 
We use Corollary \ref{cl-cor-dh} such that \eqref{cl-eqn-omegat4} holds with sufficiently small $\epsilon$. 
Now we use the step (B) of the proof 
of Proposition \ref{cl-prop-openess}.\hfill \SSn  {\parfillskip0pt\par}
	\end{proof} 
\index{crtopology@$C^r$-topology}

%

\begin{lemma}\label{cl-lem-KV} 
Let $V$ be a properly convex cone, 
and let $V^*$ be a dual cone. 
Suppose that a Koszul-Vinberg function 
$f_{V^*}(x)$ is defined on a compact neighborhood $B$ of $x$ contained in a convex cone $V$. 
Let $V_1$ be another properly convex cone
containing the same neighborhood. 
Let $\Omega:= \SI(V^*)$ and $\Omega_1 :=\SI(V_1^*)$
for the dual $V_1^*$ of $V_1$. 
For given any integer $s \geq 1$ and $\eps > 0$, there exists $\delta >0 $ such that 
if  the Hausdorff distance between $\Omega$ and $\Omega_1$ is $\delta$-close, 
then $f_{V^*}(x)$ and $f_{V_1^*}(x)$ are $\eps$-close in $B$ 
in the $C^r$-topology.
\end{lemma}
\begin{proof}
	We prove for $\SI^n$. 
By Lemma \ref{cl-lem-convHaus}, we have 
\begin{align} \label{cl-eqn-OmegaN}
\Omega^* \subset N_\delta(\Omega_1^*), \Omega_1^* \subset N_\delta(\Omega^*), \nonumber \\
(\Omega - N_\delta(\partial \Omega))^* \subset \Omega_1^*, \hbox{ and } \nonumber \\ 
(\Omega_1 - N_\delta(\partial \Omega_1))^* \subset \Omega^*
\end{align} 
provided $\delta$ is sufficiently small. 
We choose sufficiently small $\delta > 0$ such that 
\[ B \subset \Omega - N_\delta(\partial \Omega), \Omega_1 - N_\delta(\partial \Omega_1).\]
Recall the Koszul-Vinberg integral \eqref{cl-eqn-kv}. 
For fixed $\phi \in V^\ast$ or $\in V_1^\ast$,  
the functions $e^{-\phi(x)}$ and the derivatives of $e^{-\phi(x)}$ 
with respect to $x$ in the domains are uniformly bounded on $B$
since $\phi$ and its derivatives are bounded function on $B$. 
The integral is computable from an affine hyperspace meeting $V^*$ and $V_1^*$ in bounded precompact convex sets. Also, the integration is with respect to 
$\phi$. The result follows from \eqref{cl-eqn-OmegaN}.
(See Chapter 4 of \cite{goldmanbook}.) 
\end{proof}

\subsection{The proof of Theorem \ref{cl-thm-conv2}.}

\begin{proof}[The proof of Theorem \ref{cl-thm-conv2}]
Suppose that $\mathcal{O}$ has an SPC-structure $\mu$ with generalized lens-shaped or horospherical $\cR$- or $\cT$-ends. 
Let $\mathbbm{U}$ be the union of end neighborhoods of product form with mutually disjoint 
closures. 
By premises, the end structures are given. 


We assume that $\mu_0$ and $\mu_s$ correspond to elements of 
$\mathcal{U}$ in $\Def^s_{\mathcal{E}, s_{\mathcal{U}}, \lh}(\orb)$. 
We show that a structure $\mu_s$ that has generalized lens-shaped or horospherical $\cR$- or $\cT$-ends sufficiently close $\mu_0$ in $\orb - \mathbbm{U}$
is also SPC. 

Let $h: \pi_1(\orb) \ra \SL_\pm(n+1, \bR)$ be the lift of the holonomy homomorphism corresponding to $\mu_0$ where 
$\torb \subset \SI^n$ is a properly convex domain projectively covering $\orb$. 
Let $h_s: \pi_1(\orb) \ra \SL_\pm(n+1, \bR)$ be the lift of the holonomy homomorphism corresponding to $\mu_s$ sufficiently close to $h$ in 
\[\mathcal{U} \subset \Hom^s_{\mathcal{E}, \lh}(\pi_1(\orb), \SL_\pm(n+1, \bR)).\]  By Theorem \ref{op-thm-projective}, 
it corresponds to a real projective structure $\mu_s$ on $\orb$. 
Since $[\mu_s] \in \Def^s_{{\mathcal E}, \lh}(\orb)$, 
it suffices to show that $\mu_s$ is a properly convex real projective 
structure. 

Let $\mathcal{O}:= C(\torb)/h(\pi_1(\orb))$ with $C(\torb)$ as an universal cover. 
Let $\torb_s$ denote $\torb$ with $\mu_s$. We apply special affine suspension to obtain 
an affine orbifold $\orb \times \SI^1$. (See Section \ref{prelim-sub-asusp}.) 
The universal cover is still $C(\torb)$ and has a corresponding affine structure $\tilde \mu_s$. 
We denote $C(\torb)$ with the lifted affine structure of $\tilde \mu_s$ by $C(\torb)_s$. 
Recall the Kuiper completion $\hat{C}(\torb)_s$ of $C(\torb)_s$. This is a completion of $C(\torb)_s$
with respect to the path metric induced from the pull-back of the standard Riemannian metric on $\SI^{n+1}$ 
by the developing map $\dev_s$ of $\tilde \mu_s$. (Here,the image is in $\bR^{n+1}$ as an affine subspace of $\SI^{n+1}$.) 
The developing maps always extend to ones on $\hat{C}(\torb)_s$ which we denote by $\dev_s$ again. 
(See \cite{cdcr1} and \cite{psconv} for details.) 

By Corollary \ref{cl-cor-openness}, an affine
suspension $\tilde \mu_s$ of $\mu_s$ also has a Hessian function $\phi$
since $\mu_s$ is in a sufficiently small $C^2$-neighborhood of $\mu$
in $\orb - \mathbbm{U}$. 
The Hessian metric $D d\phi$ is invariant under affine automorphism groups of $C(\torb)$ 
by construction. 
We prove that $\tilde \mu_s$ is properly convex, which shows that $\mu_s$ is properly convex:

Suppose that $\tilde \mu_s$ is not convex. Then there exists a triangle $T$ embedded in $\hat{C}(\torb)_s$
where a point in the interior of an edge of $T$ is in the ideal set 
\[\delta_\infty \torb_s:=\hat{C}(\torb)_s - C(\torb)_s\]
while $T^o$ and the union $l'$ of two other edges are in $C(\torb)_s$. 
We can move the triangle $T$ so that the interior of an edge $l$ has a point $x_\infty$ in $\delta_\infty \torb_s$
and $\dev_s(l)$ does not pass the origin.
We form a parameter of geodesics $l_t$, $t \in [0, \eps]$ in $T$ such that 
\[l_0=l \hbox{ and } l_t \subset C(\torb)_s \hbox{ with } \partial l_t \subset l'\]
is close to $l$ in the triangle. 
 (See Theorem A.2 of \cite{psconv} for details.)



Let $p, q$ be the endpoints of $l$. Then the Hessian metric is $D^s d \phi$ for a function $\phi$ defined on 
$C(\torb)_s$. And $d \phi |p$ and $d \phi |q$ are bounded, where $D^s$ is the affine connection of $\mu_s$.
This should be true for 
$p_t$ and $q_t$ for sufficiently small $t$ uniformly. 
Let $u$, $u \in [0, 1]$, be the affine parameter of $l_t$, i.e., $l_t(s)$ is a constant speed line in $\bR^{n+1}$ when developed. 
We assume that $u \in (\eps_t, 1-\eps_t)$ parameterizes $l_t$ for sufficiently small $t$
where $\eps_t \ra 0$ as $t \ra 0$
and $d l_t/ds = \vec{v}$ for a parallel vector $\vec{v}$.
The function $D^s_vd_v \phi( l_t(u))$ is uniformly bounded since its integral 
$d_{\vec{v}} \phi( l_t(u))$ is strictly 
increasing by the strict convexity and converges to certain values as $u \ra \eps_t, 1-\eps_t$. 

Since  \[ \int_{\eps_t}^{1-\eps_t} D^s_vd_v \phi( l_t(u))du = d\phi(p_t)(\vec{v}) - d\phi(q_t)(\vec{v}), \]
the function $\sqrt{D^s_{\vec{v}}d_{\vec{v}} \phi( l_t(u))}$ is also integrable by Jensen's inequality, and 
the length of $l_t$ 
\[ \int_{\eps_t}^{1-\eps_t} \sqrt{D^s_{\vec{v}}d_{\vec{v}} \phi( l_t(u))}du  \]
under the Hessian metric $Dd\phi$ 
has an upper bound $\sqrt{d\phi(p_t)(\vec{v}) - d\phi(q_t)(\vec{v})}$ 
by the same inequality. 
Since 
\[\sqrt{d\phi(p_t)(\vec{v}) - d\phi(q_t)(\vec{v})} \ra 
\sqrt{d\phi(p_0)(\vec{v}) - d\phi(q_0)(\vec{v})} \hbox{ as } t \ra 0,\] 
the length of $l_t$ is uniformly bounded.  (These are the same arguments as in 
\cite{Koszul65}.)

$\mathbbm{U}$ corresponds to an inverse image $\widetilde{\mathbbm{U}}$ in $\torb$ and 
to $C(\tilde{\mathbbm{U}})_s$ the inverse image in $C(\torb)_s$. 
The minimum distance between the components of $\mathbbm{U}$ is bounded below since the metric is invariant under scalar dilatations.
in $C(\torb)_s$. 
Since $(C(\torb_s) - C(\tilde{\mathbbm{U}})_s)/\langle \bR_+\Idd\rangle$ is compact, 
if $l$ meets infinitely many components of $C(\widetilde{\mathbbm{U}})_s$, then the length is infinite.  

As $t \ra 0$, the number is thus bounded, 
$l$ can be divided into finite subsections, each of which meets 
at most one component of $C(\widetilde{\mathbbm{U}})_s$. 

Let $\hat l$ be the subsegment of $l$ in $C(\torb)_s$ containing $x_\infty$ in 
the ideal set of the Kuiper completion of $C(\torb_s)$
with respect to $\dev_s$ 
and meeting only one component $C({\tilde{\mathbbm{U}}}_1)_s$ of $C(\widetilde{\mathbbm{U}})_s$
with $\Bd s \in \Bd C(\mathbbm{U}_1)_s$. 
Let $\hat l_t$ be the subsegment of $l_t$ such that 
the parameter of the endpoints of segments of form $\hat l_t$ 
converges to those of $\hat l$ as $t \ra 0$. 
Let $p'$ and $q'$ be the endpoint of $\hat l$.


Suppose that $C({\tilde{\mathbbm{U}}}_1)_s \subset C(\torb)_s$ corresponds to a lens-shaped or horospherical R-p-end neighborhood $\tilde{\mathbbm{U}}'_1$ 
in $\torb_s$.
and $x_\infty$ is on a line corresponding to the p-end vertex of
$\tilde{\mathbbm{U}}'_1$. 
We project to $\SI^n$ from by the projection $\Pi':\bR^{n+1} -\{O\} \ra \SI^n$.  
Then the proper convexity of $\tilde \Sigma_{\tilde E}$ contradicts this
when $\tilde E$ is an R-p-end of the generalized lens type. 
When $\tilde E$ is a horospherical p-end, 
the whole segment must be in $\Bd \torb$ by convexity.  
Theorem \ref{ce-thm-affinehoro} contradicts this. 

Now suppose that $\Pi'(x_\infty)$ is in the middle of the radial line
from the p-end vertex. 
Then the interior of the triangle 
is transverse to the radial lines. Since our p-end orbifold
$\tilde \Sigma_{\tilde E}$ is convex, there cannot be such a line with a single interior point in the ideal set. 

If $C(\tilde{{\mathbbm{U}}}_1)_s$ is the inverse image in $C(\torb)_s$ of 
a generalized lens-shaped T-p-end neighborhood $\tilde{\mathbbm{U}}_1$ in $\torb_s$, 
then clearly there is no such a segment $l$ containing an ideal 
$x_\infty$ in its interior similarly. 


Now suppose that a subsegment $l_1$ of $l$ contains an ideal point in its interior, but is disjoint from $\tilde{\mathbbm{U}}$. 
There is a connected arc in $l_1 \cap C(\torb - \tilde{\mathbbm{U}})_s$ that ends
at an ideal point $x_\infty$. This is an arc never in a compact subset of 
$C(\torb)_s$. However, we showed above that the Hessian length of $l_t$ is bounded. 
Since for a subarc $l_{1,t}$ of $l_t$, 
the parameter $\{l_{1, t}\}$ converges to $l_1$ as $t \ra 0$. 
Thus, the Hessian length of $l_1$ is also finite. 
Since $C(\torb -\tilde{\mathbbm{U}})$ covers a compact orbifold that is affinely
suspended over $\orb - \mathbbm{U}$, 
the Hessian metric is compatible with any Riemannian metric. 
Since $l_1$ is in a compact orbifold, it cannot have a finite 
Riemannian length. 


This is again a contradiction. Therefore, $\tilde{\mathcal{O}}_s$ is convex. 

Finally,  for sufficiently small deformations, the convex real projective structures are properly convex. 
Suppose not.  
Then there is a sequence $\{\mu_{s_i}\}$ of deformed convex real projective structures that are not properly convex arbitrarily near the original one. 
By Proposition \ref{prelim-prop-classconv}, 
there exists a unique great sphere $\SI^{i_0}$ in the boundary of 
the nonproperly convex set for some $i_0 > 0$ by passing to a subsequence if necessary. 
By uniqueness, the holonomy 
$h_{s_i}$ acts on $\SI^{i_0}$. 

The sequence of structures converges to the beginning $\mu$ 
in $\torb - U$. 
By taking limits, the original holonomy has to be reducible.
This contradicts the premise about $\mathcal{U}$ having only irreducible representations. 

Suppose now that $\mathcal{O}$ with $\mu$ is strictly SPC with lens-shaped or horospherical $\cR$- or $\cT$-ends. 
The relative hyperbolicity of $\tilde{\mathcal{O}}$ with respect to the p-ends is stable under sufficiently small deformations
since it is a metric property invariant under quasi-isometries by Theorem \ref{rh-thm-relhyp}.

The irreducibility and stability follow since these are open conditions
in \[\Hom(\pi_1(\orb), \SL_{\pm}(n+1, \bR)).\] 
Also, the ends are lens-shaped or horospherical. 

By Theorem \ref{intro-thm-sSPC}, the holonomy is not in a parabolic group. 
Projecting back to $\rpn$ and $\PGL(n+1, \bR)$, 
we complete the proof of Theorem \ref{cl-thm-conv2}.
\end{proof}



\section{The closedness of convex real projective structures}\label{cl-sec-closed}

%

%

We recall the subspace $\rep^s_{\mathcal E}(\pi_1(\mathcal{O}), \PGL(n+1, \bR))$ 
of stable irreducible characters
of $\rep_{\mathcal E}(\pi_1(\mathcal{O}), \PGL(n+1, \bR))$
which is shown to be an open subset of a semialgebraic set in Section \ref{intro-sub-semialg}. 
We denote by $\rep^s_{{\mathcal E}, \lh}(\pi_1(\mathcal{O}), \PGL(n+1, \bR))$ the subspace of stable irreducible characters
of $\rep_{{\mathcal E}, \lh}(\pi_1(\mathcal{O}), \PGL(n+1, \bR))$, an open subset of a semialgebraic set. 





\subsection{Preliminary of the section} 


The following shows that the limiting projective orbifolds have nice end decompositions. 
We note that there are a number of similar results by Cooper, Long, Tillmann and so on
on similar settings. Some of these are not yet published. 

 \begin{lemma} \label{cl-lem-endspres} 
	Let $h_i, h \in \Hom_{\mathcal{E}, \lh}(\pi_1(\orb)), \SLpm)$
	$(\hbox{resp. }\in \Hom_{\mathcal{E}, \lh}(\pi_1(\orb)), \PGL(n+1, \bR))).$
	Suppose that the following hold\/{\rm :}
	\begin{itemize}  
		\item Let $\orb$ be a strongly tame 
		real projective orbifold with ends assigned types $\cR$ and $\cT$ 
		satisfying {\rm (}IE\/{\rm )} and {\rm (}NA\/{\rm )}
		with a compatible end compactification. 
	\item Let $\Omega_i$ be a properly convex open domain in $\SI^n$ {\rm (}resp. $\RP^n$\/{\rm ).}
	\item Suppose that $\Omega_i/h_i(\pi_1(\orb))$ is an $n$-dimensional noncompact strongly tame SPC-orbifold with generalized lens-shaped or horospherical $\cR$- or $\cT$-ends.
	\item Assume that each p-end holonomy group $h_i(\pi_1(E_j))$ 
	of $h_i(\pi_1(\orb))$ of 
	type $\mathcal{R}$ has a p-end vertex $\mbv^j_i$ corresponding 
	to the p-end structure where $\{ \mbv^j_i\}$ forms a convergent sequence 
	as $i \ra \infty$. We assume that 
	$h_i \mapsto \mbv^j_i$ extends to an analytic function near 
	$h$.  
	\item Assume that   each p-end holonomy group $h_i(\pi_1(E_j))$ 
	of $h_i(\pi_1(\orb))$ of type $\mathcal{T}$ has a hyperspace $P^j_i$
	containing the p-ideal boundary component 
	where $\{P^j_i\}$ forms 
	a convergent sequence as $i \ra \infty$. 
	We assume that 
	$h_i \mapsto P^j_i$ extends to an analytic function on a neighborhood of 
	$h$.  
	\item Suppose that $\{h_i\} \ra h$ algebraically
	where $h$ is discrete and faithful. 
	\item $\clo(\Omega_i) \ra K$ for a compact properly convex domain 
	$K \subset \SI^n$, $K^o \ne \emp$. 
	\end{itemize} 
	Then the following holds\/{\rm :}
	\begin{itemize}
		\item $\orb_h:= K^o/h(\pi_1(\orb))$ is a strongly tame real projective SPC-orbifold  
		with generalized lens-shaped or horospherical $\cR$- or $\cT$-ends 
		to be denoted $\orb_h$ diffeomorphic to $\orb$. 
		\item For each p-end $\tilde E$ of the universal cover $\torb_h$ of $\orb_h$, 
		$K^o$ has a subgroup $h(\pi_1(\tilde E))$ acting on 
		a $h(\pi_1(\tilde E))$-invariant open set $U_{\tilde E}$
		where $U_{\tilde E}/h(\pi_1(\tilde E))$ is an end neighborhood that 
		is one of the following\/{\rm :}
		\begin{itemize} 
		\item a horospherical or lens-shaped totally geodesic end neighborhood provided 
		$\tilde E$ is a $\cT$-p-end, 
		\item a horospherical or concave end neighborhood provided $\tilde E$ is 
		a $\cR$-p-end.
		\end{itemize} 
		\item Finally, suppose that there is a fixed strongly tame properly convex 
real projective orbifold $\orb'$ with an ideal boundary structure
		and diffeomorphism $f_i: \orb' \ra \Omega_i/h_i(\pi_1(\orb))$ for 
		sufficiently large $i$ extending to
		a diffeomorphism of an end-compactification $\bar \orb'$
		to the end compactification of $\Omega_i/h_i(\pi_1(\orb))$ compatible 
		with R-end and T-end structures given 
		by $\mbv^j_i$ and $P^j_i$. 
		Then $K^o/h(\pi_1(\orb))$ is an orbifold with a diffeomorphism 
		$f$ from $\orb'$ extending to a diffeomorphism 
		from $\bar \orb'$ to an end compactification of 
		$K^o/h(\pi_1(\orb))$ with the above R-end and T-end structures. 
	\end{itemize} 
\end{lemma} 
\begin{proof} 
		Again, it suffices to prove the theorem when $\Omega_i \subset \SI^n$. 
%
The holonomy group $h(\pi_1(\orb))$ acts on $K^o$ with 
a Hilbert metric. 
	Hence, $K^o/h(\pi_1(\orb))$ is an orbifold to be denoted $\orb_h$. 
	(See Lemma 1 of \cite{CG93}.)

We show that $\orb_h$ is strongly tame by constructing the end neighborhoods. 
First, we show the existence of regions projectively
diffeomorphic to p-end neighborhoods of various forms.
	

	
	Since $h \in \Hom_{{\mathcal E}, \lh}(\pi_1(\mathcal{O}), \SLpm)$ holds, 
	each p-end holonomy group $h(\pi_1(\tilde E))$ acts on a horoball $H \subset \SI^n$, 
	a generalized lens cone, or a totally geodesic hypersurface 
	$\tilde S_{\tilde E}$ with a cocompactly-acted lens $L$.
	In the first case, we can choose a sufficiently small horoball $U$ inside $K^o$ and in $H$ since the sharply supporting hyperspaces at the vertex of $H$
	must coincide by the invariance under $h(\pi_1(\tilde E))$
	by a limiting argument. 
Therefore, a component $U$ of $K^o - \Bd U$
	is  projectively diffeomorphic to a generalized lens-shaped p-end neighborhood. 

	Now, we consider the second case. 
	Let $\mbv$ be a limit of the sequence $\{\mbv_{\tilde E, i}\}$ 
	of the fixed p-end vertices of 
	$h_i(\pi_1(\tilde E))$. We obtain $\mbv \in K$.  
	Also, $\mbv \not \in K^o $ since otherwise the elements fixing it has to be of finite order by the proper discontinuity of the Hilbert isometric 
	action of $h(\pi_1(\tilde E))$.
	For each $i$, $h_i(\pi_1(\tilde E))$ acts on a lens cone 
	$L_i\ast \{\mbv_{\tilde E, i}\}$.
		We may assume without loss of generality that $\mbv_{\tilde E, i}$ is 
		constant by changing $\dev_i$ by a convergent sequence $\{g_i\}$ of 
		elements of $\SL_{\pm}(n+1, \bR)$.  
		We may assume that $h_i$ in a continuous parameter 
		converging to $h$ since there are finitely many components in
		the above real algebraic set.
		By Corollary \ref{abelian-cor-smvar} and the condition ``lh", 
		we may assume that 
		\[\{\clo(R_{\mbv_{\tilde E, i}}(\Omega_i))\} \ra K_{\mbv}\] 
		for 
		a properly convex domain $K_{\mbv}$ on which $h(\pi_1(\tilde E))$ acts. 
		Since $\Omega_i$ is a subset of a tube domain 
		for $R_{\mbv_{\tilde E, i}}(\Omega_i)$, we 
		deduce that $\Omega$ is a subset of a tube domain 
		for $K^o_{\mbv}$. Since $K^o_{\mbv}/h(\pi_1(\tilde E))$ is a closed 
		properly convex real projective orbifold, and 
			$R_{\mbv_{\tilde E}}(\Omega) \subset K^o_{\mbv}$, 
			it follows that they are equal by Lemma \ref{prelim-lem-domainIn}. 
			Hence, $R_{\mbv_{\tilde E}}(\Omega)$ is properly convex. 
				By the Hilbert metric on this domain, 
				$h(\pi_1(\tilde E))$ acts properly discontinuously on it.
				Since $h(\pi_1(\tilde E))$ satisfies the uniform 
				middle-eigenvalue condition by the premise, 
				Theorem \ref{pr-thm-equ} shows that 
				the action is distanced in a tubular domain corresponding 
				to $R_{\mbv_{\tilde E}}(K^o)$. 
					Hence, $h(\pi_1(\tilde E))$ acts properly and cocompactly on 
				a generalized lens $L$ in $K$. 
				The group $h(\pi_1(\tilde E))$ acts on 
				an open set $U_L:=L\ast \{\mbv_{\tilde E}\} -L$. 
				We may choose one with sufficiently large $L$ such that 
				the lower boundary component $\partial_- L$ is a subset 
				in $K^o$ since we can make $U_L\cap \mathcal{T}_{\mbv_{\tilde E}}(F)$ be as small as we wish for any compact fundamental domain $F$
				for	$R_{\mbv_{\tilde E}}(K^o)$ and $h(\pi_1(\tilde E))$. 
Therefore, $U_L$  is projectively diffeomorphic to 
a p-end neighborhood of $\tilde E$.  
				
	%

%
	
	In the third case, we can find a cocompactly-acted lens neighborhood of a totally geodesic 
	domain $\tilde S_{\tilde E, i} \subset \clo(\Omega_i)\cap P_i$ 
	in a hypersurface $P_i$  
	on which $h_i(\pi_1(\tilde E))$ acts. 
	We may assume without loss of generality that $P_i$ is constant 
	by changing $\dev_i$ by a convergent sequence $g_i$ in $\SL_{\pm}(n+1, \bR)$. 
	 We may assume by taking  subsequences that 
	 $\{\clo(\tilde S_{\tilde E, i})\} \ra D$ for a properly convex domain $D$ 
	 by Corollary \ref{abelian-cor-smvar}.
	
	By Corollary \ref{prelim-cor-uniqueness}, 
	$D^o/h(\pi_1(\tilde E))$ is a closed orbifold homotopy equivalent to 
	$S_{E, i}$ up to finite manifold covers. 
	By Theorem \ref{pr-thm-equ2}, $h(\pi_1(\tilde E))$ satisfies the uniform 
	middle-eigenvalue condition with respect to the hyperspace containing $D$. 
	By Theorem \ref{du-thm-lensn}, the group $h(\pi_1(\tilde E))$ acts on
	a component $L_1$ of $L - P$ is in $K^o$ for a lens $L$.
	Therefore, $L_1^o$ is diffeomorphic to a p-end neighborhood of $\tilde E$ 

Let us call the constructed objects {\em candidate p-end neighborhoods}. Each one is associated 
with a p-end and so are their images under $h(\pi_1(\orb))$. 
We now show that we can choose the candidate  p-end neighborhoods 
such that their images under $h(\pi_1(\orb)$ are mutually disjoint
in $K^o$.

For each candidate generalized lens-type R-p-end, 
we choose a concave p-end neighborhood in it. 
Then their images are mutually disjoint by Corollary \ref{app-cor-mcn}  where 
we did not require the strong tameness. 

Since the boundary of candidate horospherical p-end neighborhoods and candidate lens-type T-p-end neighborhoods 
cover compact immersed hypersurfaces in $\orb_h$, 
we can choose sufficiently far away concave p-end neighborhoods such that 
their images in $\orb_h$ are disjoint from these. 

We can apply Corollary \ref{app-cor-shrink} for the quotients of these candidate
p-end neighborhoods of a p-end $\tilde E$ by $h(\pi_1(\tilde E))$: 

Suppose that for every choice of shrinking of these p-end neighborhoods, 
there is some pair of 
candidate horospherical p-end neighborhoods or lens T-p-end neighborhoods meeting 
each other.  

Let us choose for representative p-ends $\tilde E_1, \dots, \tilde E_m$ of $\orb$, 
corresponding candidate p-end neighborhoods $U_1, \dots, U_m$ and 
consider their images under $h(\pi_1(\orb))$.  
For each p-end $\tilde E_k$ and another p-end $\tilde E_l$, 
we define the set $C_{k, l}$ of cosets $\gamma h(\pi_1(\tilde E_k))$ 
of $h(\pi_1(\orb))/h(\pi_1(\tilde E_k))$
where $U_k \cap \gamma(U_l) \ne \emp$. 
(Note $C_{k, l} \ne C_{l, k}$. )

Since $h(\tilde E_l)$ acts cocompactly $\Bd U_l \cap K^o$ by our constructions above, 
it follows that $\Bd U_l \cap K^o$ maps to a compact hypersurface $\Sigma_l$ for each $l$. 
We may consider $\Sigma_l$ to be a generic smooth suborbifold by small perturbations.
Their intersection patterns read off these coset classes. 
Since $\Sigma_i$ is compact, their intersections consist 
of finitely many immersed closed suborbifolds 
and above $C_{l, k}$ must be finite for each $l, k$. 

Now, we can apply Corollary \ref{app-cor-shrink} to find mutually disjoint sets of 
candidate p-end neighborhoods for $\tilde E_i$. 
We may suppose that the closures of their images in $\orb_h$ are mutually disjoint as well by a small extension of the above argument. 
Let us call the images of the candidate p-end neighborhoods
the {\em candidate end neighborhoods}.
Hence, we assume that the candidate end neighborhoods can be chosen mutually disjoint. 

Let us call $\hat U$ the union of these candidate end neighborhoods. 
Then $\orb_h -\hat U$ is an orbifold homotopy equivalent to $\bar \orb$
where we can find a homotopy equivalence map restricting to homotopy equivalences
for all boundary components. This means that $\orb_h$ is strongly tame. 

We can also choose a diffeomorphism $f: \orb_h \ra \orb$: 
Since $\partial \hat U$ in $\orb_h$  is a disjoint union of compact hypersurfaces, 
By our constuction, each candidate end neighborhood in $\hat U$ 
bounded by $\Sigma_i$ can be perturbed to
one for $\Sigma_i^h$ for $\orb_i:=\Omega_i/h_i(\pi_1(\orb))$. We call their union $\hat U_i$. 
Then we can find diffeomorphisms $f_i': \hat U \ra \hat U_i$ by deformations
by Lemmas \ref{op-lem-horob} and \ref{op-lem-lensp}, and
$f''_i: \orb_h - \hat U \ra \orb_i - \hat U_i$ by using the argument as in \cite{dgorb} since these 
are compact orbifolds where we know how the boundary components are deforming.
We can patch these to obtain a smooth diffeomorphism $f'''_i: \orb_h \ra \orb_i$
that follows from the existence of collar neighborhoods of the orbifold boundary (see Chapter 4 of \cite{Hirsch} and Chapter 4 of \cite{Cbook}.)
Hence, we obtained a diffeomorphism $f: \orb_h \ra \orb$. 
		
	We apply Theorem \ref{intro-thm-sSPC} to show that the real projective 
	structure is SPC. 
	
	The end compactification $\bar \orb_h$ of $K^o/h(\pi_1(\orb))$ is given 
	by attaching the cover of an ideal boundary component for each T-p-end
	and attaching $\Sigma_E \times [1, 0)$ to 
	an end neighborhood of an R-end $E$ by a diffeomorphism 
	restricting to a proper map as in Section \ref{op-sub-endstr}. 
	
	For the final part, 
	recall from Section \ref{intro-sub-semialg} that 
	$\Hom^s_{\mathcal E, f}(\pi_1(\orb), \SL_\pm(n+1, \bR))$
	is Zariski dense. The homomorphism $h$ is in it 
	by Proposition \ref{cl-prop-section}.
	We choose p-end vertices, and 
	the hyperspace containing the p-ideal boundary components 
	as respective limits of the corresponding sequence of 
    the corresponding objects of $h_i$.
	Based on the premise of analytic extension, we can build a fixed section $s_{\mathcal{U}}$ on a Zariski open subset $\mathcal{U}$ containing $h$ 
    by Proposition \ref{cl-prop-section}.  
	By Theorem \ref{op-thm-projective},  
	we obtain a subspace of the parameter of 
	real projective structures on a strongly tame orbifold $\orb'$ 
	with end structures determined by $s_{\mathcal U}$
	for each point of $\mathcal{U}$. 
	
	By the premise, the p-end vertices of $\Omega_i$ and 
	the hypersurfaces containing the p-ideal boundary components are 
	determined also by $s_{\mathcal{U}}$ for sufficiently large $i$.
%
%
	Now, $\orb':= K^o/h(\pi_1(\orb))$ 
	is realized as a convex real projective orbifold with 
	ends determined by $s_{\mathcal U}$ also. 
	By Theorem \ref{op-thm-projective}, 
	there is a neighborhood $\mathcal{U}' \subset \mathcal{U}$ 
	where every holonomy is realized by 
	a convex real projective structure with end structures 
	determined by $s_{\mathcal U}$. 
	Since $h_i$ may be assumed to be in $\mathcal{U}'$ except for finitely many $i$s,
	a structure $\mu_i$ on $\orb'$ has a holonomy homomorphism $h_i$. 
	By Theorem \ref{closed-thm-holdet}, 
	$\Omega_i/h_i(\pi_1(\orb))$ is projectively 
	diffeomorphic to $\orb'$ with a convex real projective structure $\mu_i$ 
	with identical R-end and T-end structures
	for sufficiently large $i$.
	
%
%


	
We can construct a compactification $\bar \orb_h$ from $\bar h$ using 
the diffeomorphism $f: \orb_h  \ra \orb$. 
Then $f$ extends smoothly as a diffeomorphism from 
	$\bar \orb_h$ to $\bar \orb$. 
%
%

%
%
	%
	%
	%
	%
	%
	%
	Again, the $\RP^n$-version follows from Proposition \ref{prelim-prop-closureind}. 
\hfill	\SnP {\parfillskip0pt\par}
\end{proof}


\begin{theorem}[Uniqueness of domains] \label{closed-thm-holdet} 
	Let $\Gamma$ be a discrete projective automorphism group of a properly 
	convex open domain $\Omega \subset \SI^n$. Suppose that 
	$\Omega/\Gamma$ is an $n$-dimensional strongly tame SPC-orbifold with 
	generalized lens-shaped or horospherical $\cR$- or $\cT$-ends and satisfies  {\rm (}IE\/{\rm )} and {\rm (}NA\/{\rm )}.  Assume $\partial \orb =\emp$. 
	Suppose that for each $\mbv_{\tilde E} \in \Bd \Omega$ for each R-p-end $\tilde E$ is specified up to $\mathcal{A}$ in $\SI^n$,
	and so is each hyperspace for each T-p-end $\tilde E$ meeting 
	$\Bd\Omega$.
	Then $\Omega$ is a unique domain with 
	these properties up to the antipodal map $\mathcal{A}$.
	\end{theorem} 
\begin{proof}
Suppose that $\Omega_1$ and $\Omega_2$ are distinct open domains in $\SI^n$
satisfying the above properties. 
For this, we assume that $\Gamma$ is torsion-free by taking 
a finite-index subgroup by Theorem \ref{prelim-thm-vgood}. 
We claim that $\Omega_1$ and $\Omega_2$ are disjoint:

Suppose that $\Omega':= \Omega_1\cap \Omega_2$ is a nonempty open set. 
Since $\Omega_1, \Omega_2,$ and $\Omega'$ are all $n$-cells, 
the set of p-ends of $\Omega_1$, the set of those of $\Omega_2$, 
in one-to-one correspondences by considering their p-end fundamental groups. 
The types are also preserved by the premise. 

Suppose that $\tilde E_1$ and $\tilde E_2$ are the corresponding R-p-ends of a generalized lens type.   
The p-end vertex $\mbv_{\tilde E_j}$ of a generalized lens-shaped 
R-p-end $\tilde E_j$ of $\Omega_j$, $j=1,2$, is determined up to $\mathcal{A}$
by the premise.

Suppose that $\mbv_{\tilde E_1}$ and $\mbv_{\tilde E_2}$ are antipodal. 
The interiors of $\Omega' \ast \{\mbv_{\tilde E_j}\}$, $j=1,2$, are in
$\Omega_j$ by the convexity of $\clo(\Omega_j)$. 
In addition, $\Omega_1 \cup \Omega_2$ is in a convex tube ${\mathcal{T}}_{\mbv_{\tilde E_1}}(\tilde \Sigma_{\tilde E_1})$ with the
vertices $\mbv_{\tilde E_1}$ and its antipode $\mbv_{\tilde E_2}$ in
the direction of $\tilde \Sigma_{\tilde E_1}$. 
In addition, the convex hull of $\clo(\Omega_1)\cup \clo(\Omega_2)$ is equal to
${\mathcal{T}}_{\mbv_{\tilde E_1}}(\tilde \Sigma_{\tilde E_1})$. 
$h(\pi_1(\orb))$ acts on the unique pair of antipodal points $\{\mbv_{\tilde E_1}, \mbv_{\tilde E_2}\}$.  
Hence, $h(\pi_1(\orb))$ is reducible, contradicting the premise. 

Suppose that $\mbv_{\tilde E_1} = \mbv_{\tilde E_2}$. 
Then  $R_{\mbv_{\tilde E_1}}(\Omega_1)$ and 
$R_{\mbv_{\tilde E_1}}(\Omega_2)$ are not disjoint
since otherwise $\Omega_1$ and $\Omega_2$ are disjoint. 
Lemma \ref{prelim-lem-domainIn} shows that they are equal. 
The generalized lens-cone 
p-end neighborhood $U_1$ in $\Omega_1$ and one $U_2$ in $\Omega_2$ 
must intersect. 
Hence, by Lemma \ref{pr-lem-endnhbd}, 
the intersection of a generalized lens cone p-end neighborhood of 
$\Omega_1$ and that of $\Omega_2$ is one for $\Omega'$:


Suppose that $\tilde E_j$ is a horospherical p-end of $\Omega_j$, $j=1, 2$. 
Then p-end vertices $\mbv_{\tilde E_j}$, $j=1, 2$,  
are either equal or antipodal since there is a unique antipodal pair of fixed points for the cusp group $\Gamma_{\tilde E_j}$, $j=1,2$. 
Since the fixed point in $\mbv_{\tilde E_j}$ is the unique limit point of 
$\{\gamma^n(p)\}$ as $n \ra \infty$ for any $p\in \Omega_j$, 
it follows that $\mbv_{\tilde E_1} = \mbv_{\tilde E_2}$. 
We can verify that $\Omega_1, \Omega_2$, and $\Omega'$ share a horospherical p-end neighborhood 
from this by Lemma \ref{pr-lem-endnhbd}. 

Similarly, consider the ideal boundary component 
$\tilde S_{\tilde E_1}$ for a T-p-end $\tilde E_1$
of $\Omega_1$ and the corresponding $\tilde S_{\tilde E_2}$ for a T-p-end $\tilde E_2$
of $\Omega_2$. 
Since $\bGamma_{\tilde E_1}$ acts on a properly convex domain 
$\Omega'$, Theorem \ref{pr-thm-equ2} and Lemma \ref{du-lem-attracting3}
show that $\clo(\Omega')\cap P$ is a nonempty properly convex set
in $\clo(\tilde S_{\tilde E_1})$. 

We claim that a point $\clo(\tilde S_{\tilde E_2})$ for 
$\Omega_2$ cannot be antipodal to any point of $\clo(\tilde S_{\tilde E_1})$: 
Suppose not. Then $\clo(\tilde S_{\tilde E_2}) = R_2(\clo(\tilde S_{\tilde E_1}))$
for a projective automorphism $R_2$ acting as $\Idd$ or $\mathcal{A}$ on 
a collection of independent subspaces of $P$ by Lemma \ref{prelim-lem-domainIn}. 
There exists a pair of extreme points $p_1 \in \clo(\tilde S_{\tilde E_1})$ and 
$p_2 \in \clo(\tilde S_{\tilde E_2})$, antipodal to each other. 
Here, there exist a point $c \in \tilde S_{\tilde E_j}$ and
a sequence $g^{(j)}_i$ such that $\{g^{(j)}_i(c)\} \ra p_j$ by Lemma 5 of \cite{Vey}. 
By Lemma \ref{du-lem-attracting3}, $\{g^{(j)}_i(d)\} \ra p_j$ for a point 
$d \in \Omega'$. This implies that $\Omega'$ is not properly convex, a contradiction.

Thus, we obtain $\tilde S_{\tilde E} = \tilde S'_{\tilde E}$ again 
by Lemma \ref{prelim-lem-domainIn}. Since $\Omega'$ is a $h(\pi_1(\tilde E))$-invariant 
open set in one side of $P$, it follows that $\Omega'$ contains a one-sided 
lens neighborhood $L_1$ by Lemma \ref{du-thm-asymniceII}. 
By Lemma \ref{pr-lem-endnhbd}, $L_1$ is a p-end neighborhood of $\Omega'$. 

%

We have concave p-end neighborhoods for radial p-ends, lens p-end neighborhoods for totally geodesic p-ends, and horoball p-end 
neighborhoods of p-ends for each of $\Omega_1$, $\Omega_2$, and $\Omega'$. 
%
We verify from the above discussions 
that a p-end neighborhood of $\Omega_1$ exists if and only if a p-end neighborhood of $\Omega_2$ 
exists and their intersection is a p-end neighborhood of $\Omega'$. 
$\Omega'/\Gamma$ is a closed submanifold in $\Omega_1/\Gamma$ and in $\Omega_2/\Gamma$. 
Thus, $\Omega_1/\Gamma, \Omega_2/\Gamma,$ and $\Omega'/\Gamma$  are all homotopy 
equivalent relative to the union of disjoint end neighborhoods. 
The map has to be surjective in order for the map to be 
a homotopy equivalence, as we can show using 
relative homology theories, 
and hence, $\Omega'=\Omega_1=\Omega_2$. 

Suppose that $\Omega_1$ and $\mathcal{A}(\Omega_2)$ meet. Then 
similarly, $\Omega_1= \mathcal{A}(\Omega_2)$. 


Suppose now that $\Omega_1 \cap \Omega_2= \emp$ and 
$\Omega_1 \cap \mathcal{A}(\Omega_2) = \emp$.  
Suppose that $\Omega_1$ has a p-end $\tilde E$ of type $\mathcal{R}$. 
The corresponding pair of the p-end neighborhoods share 
the p-end vertex or have antipodal p-end vertices by the premise. 
Since $\Omega_1$ and $\Omega_2$ are disjoint, it follows that 
$\clo(\Omega_1)\cap \clo(\Omega_2)$ or $\clo(\Omega_1) \cap {\mathcal{A}}(\clo(\Omega_2)$  
is a compact properly convex subset of dimension $< n$.
It is not empty since the p-end vertex or the antipode are in it. 
The minimal hyperspace that contains it is a proper subspace and is invariant under $\Gamma$. 
This contradicts the strong irreducibility of $h(\pi_1(\orb))$ as can be obtained from 
Theorem \ref{intro-thm-sSPC}. 
This also applies to the case when $\tilde E$ is a horospherical end of type $\mathcal{T}$.

Suppose that $\Omega_1$ has a p-end $\tilde E_1$ of type $\mathcal{T}$. 
Then $\Omega_2$ has a p-end $\tilde E_2$ of type $\mathcal{T}$. 
Now, $\Gamma_{\tilde E_1} = \Gamma_{\tilde E_2}$ acts on a hyperspace $P$ containing the 
ideal boundary component
$\tilde S_{\tilde E_i}$ in the boundary of $\Bd \Omega_i$ for $i=1, 2$. 
%
Here, $\tilde S_{\tilde E_1}/\Gamma_{\tilde E_1}$ and 
$\tilde S_{\tilde E_2}/\Gamma_{\tilde E_2}$ are closed $(n-1)$-orbifolds. 
Lemma \ref{prelim-lem-domainIn} shows that 
their closures always meet or are antipodal. 
Hence, up to $\mathcal{A}$, their closures always meet. 
Again, \[\clo(\Omega_1)\cap \clo(\Omega_2)\ne \emp
\hbox{ or }\clo(\Omega_1)\cap \mathcal{A}(\clo(\Omega_2)) \ne \emp\]
while we have $\Omega_1 \cap \Omega_2 =\emp$ and $\Omega_1\cap \mathcal{A}(\Omega_2)=\emp$. 
We obtain a lower-dimensional convex subspace fixed by $\Gamma$. 
This is a contradiction. 
\end{proof} 

\subsection{The main result for the section} 

This generalizes Theorem 4.1 of \cite{CLM18} for closed orbifolds, 
which is actually due to Benoist \cite{Benoist05}.  

\begin{theorem} \label{cl-thm-closed1} 
Let $\mathcal{O}$ be an $n$-dimensional strongly tame SPC-orbifold with 
generalized lens-shaped or horospherical $\cR$- or $\cT$-ends and satisfies  {\rm (}IE\/{\rm )} and {\rm (}NA\/{\rm )}.  Assume $\partial \orb =\emp$. 
We have an open $\PGL(n+1, \bR)$-conjugation invariant set $\mathcal{U}$ 
in a semi-algebraic subset 
\[\Hom^s_{{\mathcal E}, \lh}(\pi_1(\orb), \PGL(n+1, \bR)),\]
and a $\PGL(n+1, \bR)$-equivariant fixing section 
$s_{\mathcal{U}}: {\mathcal{U}} \ra (\RP^n)^{e_1} \times (\RP^{n \ast})^{e_2}$. Let $\mathcal{U}'$ denote the quotient set under 
 $\PGL(n+1, \bR)$. 
Assume that every finite-index subgroup of $\pi_1(\mathcal{O})$ has no nontrivial nilpotent normal subgroup.
Then the following hold\,{\rm :} 
\begin{itemize}
\item The deformation space $\CDef_{{\mathcal E}, s_{\mathcal U}, \lh}(\mathcal{O})$ of SPC-structures on $\mathcal{O}$ with generalized 
lens-shaped or horospherical $\cR$- or $\cT$-ends maps under $\hol$ 
homeomorphically to a union of components of 
\[\mathcal{U}' \subset \rep^s_{{\mathcal E}, \lh}(\pi_1(\mathcal{O}), \PGL(n+1, \bR)).\]
\item 
The deformation space $\SDef_{{\mathcal E}, s_{\mathcal U}, \lh}(\mathcal{O})$ of strictly SPC-structures on $\mathcal{O}$ with lens-shaped or horospherical $\cR$- or $\cT$-ends maps under $\hol$
homeomorphically to the union of components of  
\[\mathcal{U}' \subset \rep^s_{{\mathcal E},  \lh}(\pi_1(\mathcal{O}), \PGL(n+1, \bR)).\]
\end{itemize}
\end{theorem}
\begin{proof}

	Define $\widetilde{\CDef}_{{\mathcal E}, \lh}(\mathcal{O})$ to be the inverse image of $\CDef_{{\mathcal E}, \lh}(\mathcal{O})$ in 
	the isotopy-equivalence space $\widetilde{\Def}_{\mathcal{E}}(\orb)$
	in Definition \ref{op-defn-isot}. 
	Let $\widetilde{\mathcal{U}}$ denote the inverse image  of $\mathcal{U}$ in $\widetilde{\CDef}_{{\mathcal E}, \lh}(\mathcal{O})$
	where the vertices of R-p-ends and hyperspaces of T-p-ends are determined 
	by $s_{\mathcal{U}}$. Then $\widetilde{\mathcal{U}}$ is an open subset by Theorem \ref{op-thm-projective}. 
%

We show that \[\hol: \widetilde{\mathcal{U}} \ra \mathcal{U} \subset 
\Hom^s_{\mathcal{E}, \lh}(\pi_1(\mathcal{O}), \PGL(n+1, \bR))\]
is a homeomorphism onto a union of components. This implies the results. 
%
%
%
Theorem \ref{closed-thm-holdet}  shows that $\hol$ is injective.

Now, $\hol$ is an open map by Theorems \ref{cl-thm-conv2} and  \ref{op-thm-projective}.
 %
 %
 To show that the image is of $\hol$ is closed, we show that 
 the subset of 
$ \mathcal{U} $ corresponding to 
elements in $\widetilde{\mathcal{U} }$ is closed. 
 Let $(\dev_i, h_i)$ be a sequence of development pairs
 such that we have $\{h_i\} \ra h$ algebraically. 
 Let $\Omega_i =\dev_i(\torb)$ denote the corresponding 
  properly convex domains for each $i$. 
 The limit $h$ is a discrete faithful representation 
 by Lemma 1.1 of Goldman-Millson \cite{GM88}.
 Let $\hat \Omega_i$ denote the lift of $\Omega_i$ in 
 $\SI^n$, and let $\hat h_i:\pi_1(\orb)\ra \SLpm$ 
 be the corresponding lift of $h_i$ by Theorem \ref{prelim-thm-lifting}.  
 The sequence $\{\clo(\Omega_i)\}$ also geometrically converges to a compact convex set
 $\hat \Omega$ up to choosing a subsequence 
 by Proposition \ref{prelim-prop-convC} 
 where $\hat h(\pi_1(\mathcal{O}))$ acts on
 as in Lemma 1 of \cite{CG93}. 
 If $\hat \Omega$ has the empty interior, 
 $h$ is reducible, and $h \not\in \mathcal{U}$, contradicting 
 the premise. 
 If $\hat \Omega^o$ is not empty and is not properly convex, 
 then the lift of $\Omega$
 to $\SI^n$ contains a maximal great sphere $\SI^i$, $i \geq 1$, or 
 a unique pair of antipodal points $\{p, p_-\}$ by 
Proposition \ref{prelim-prop-classconv}. 
 In the both cases, $h$ is reducible. 
 Thus, $\hat \Omega^o$ is not empty and is properly convex. 
 Let $\Omega$ denote the image of $\hat \Omega$ under the double covering map.
 As in \cite{CG93}, since $\Omega^o$ has a Hilbert metric, 
 $h(\pi_1(\mathcal{O}))$ acts on $\Omega^o$ properly discontinuously.

 By Lemma \ref{cl-lem-endspres}, 
 the condition of the generalized lens or horospherical 
 condition for $\cR$-ends or lens
 or horospherical condition for $\cT$-ends 
 of the holonomy representation is a closed condition in the 
 \[ \Hom^s_{\mathcal{E}, \lh}(\pi_1(\mathcal{O}),\PGL(n+1,\bR))\]
 as we defined above. 
 
 
 Define $\orb':= \Omega^o/h(\pi_1(\orb))$. 
 We can deform $\orb'$ with holonomy in an open subset of $\widetilde{\mathcal{U}}$
 using the openness of $\hol$ by Theorem \ref{cl-thm-conv2}. 
 We can find a deformed orbifold 
 $\orb''_i$ that has a holonomy $h_i$ for some large $i$.
Now, $\Omega_i/h_i(\orb)$ is diffeomorphic to $\orb$ being 
 in the deformation space. 
 $\orb''_i$ is diffeomorphic to $\orb$ with the corresponding end-compactifications
  since they share the same open domain as the universal cover 
 by the uniqueness for each holonomy group by Theorem \ref{closed-thm-holdet}.  
 By the openness of the map $\hol$ for $\orb'$, $\orb''$ is diffeomorphic to $\orb'$. 
 Hence, $\orb'$ is diffeomorphic to $\orb$.  
 
 

Therefore, we conclude that $\widetilde{\mathcal{U}}$
goes to a closed subset of $\mathcal{U}$. The proof up to here 
imply the first item. 
 
 Now, we go to the second item. 
 Define $\widetilde{\SDef}_{{\mathcal E}, \lh}(\mathcal{O})$ to be the inverse image of $\SDef_{{\mathcal E}, \lh}(\mathcal{O})$ in the isotopy-equivalence space 
$\widetilde{\Def}_{{\mathcal E}}(\mathcal{O})$. 
Let $\widetilde{\mathcal{U}}$ denote the inverse image of $\mathcal{U}$ 
in $\widetilde{\SDef}_{{\mathcal E}, \lh}(\mathcal{O})$. 
We show that \[\hol: \widetilde{\mathcal{U}} \ra \mathcal{U}\]
is a homeomorphism onto a union of components of $\mathcal{U}$. 
 Theorem \ref{cl-thm-conv2} shows that $\hol$ is a local homeomorphism 
 to an open set. The injectivity of $\hol$ follows the same way as in the above item. 
 
 We now show the closedness. 
By Theorem \ref{rh-thm-relhyp}, $\pi_1(\orb)$ is relatively hyperbolic with respect to the end fundamental groups. 
Let $h$ be the limit of a sequence of holonomy representations 
$\{h_i:\pi_1(\orb) \ra \PGL(n+1, \bR)\}$. 
As above, we obtain $\Omega$ as the limit of $\clo(\Omega_i)$ where $\Omega_i$ is the image of 
the developing map associated with $h_i$. 
$\Omega$ is properly convex and $\Omega^o$ is not empty.  
 Since $h$ is irreducible and acts on $\Omega^o$ properly discontinuously, it follows that 
 $\Omega^o/h(\pi_1(\mathcal{O}))$ is a strongly tame properly convex 
$n$-orbifold $\mathcal{O}'$ 
 with generalized lens-shaped or horospherical $\cR$- or $\cT$-ends 
 by the above part of the proof. 
 By Theorem \ref{rh-thm-converse} and Corollary \ref{app-cor-stLens}, $\mathcal{O}'$ is a strictly SPC-orbifold with lens-shaped or horospherical $\cR$- or $\cT$-ends. 
 The rest is the same as above. 
\end{proof}

\begin{remark}[Thurston's example]
We remark that without the end controls we have, there might be counterexamples as we can
see from the examples of geometric limits differing from algebraic limits
for sequences of hyperbolic $3$-manifolds. (See Anderson-Canary \cite{AC96}.) 
\end{remark}


\subsection{Dropping of the superscript $s$.} \label{cl-sub-sups} 

We can drop the superscript $s$ from the above space. 
Hence, the components consist of stable irreducible characters. This is a stronger result. 

%

\begin{corollary} \label{cl-cor-closed1} 
Let $\mathcal{O}$ be an $n$-dimensional noncompact strongly tame SPC-orbifold with lens-shaped or horospherical $\cR$- or $\cT$-ends and satisfies  {\rm (}IE\/{\rm )} and {\rm (}NA\/{\rm )}. 
Assume $\partial \orb =\emp$. 
Assume that no finite-index subgroups $\pi_1(\mathcal{O})$ has a nontrivial nilpotent normal subgroup.
We have a
$\PGL(n+1, \bR)$-conjugation invariant set $\mathcal{U}$ open in
a union of semi-algebraic subsets of 
\[\Hom_{\mathcal{E}, \lh}(\pi_1(\orb), \PGL(n+1, \bR)),\]
and a $\PGL(n+1, \bR)$-equivariant fixing section 
$s_{\mathcal{U}}: {\mathcal{U}} \ra (\RP^n)^{e_1} \times (\RP^{n \ast})^{e_2}$. Let $\mathcal{U}'$ denote the quotient set under 
$\PGL(n+1, \bR)$. 
Suppose that $\mathcal{U}$ contains at least one stable irreducible faithful representation.
Then the following hold\,{\rm :} 
\begin{itemize}
	\item The deformation space $\CDef_{{\mathcal E}, s_{\mathcal U}, \lh}(\mathcal{O})$ of SPC-structures on $\mathcal{O}$ with generalized 
	lens-shaped or horospherical $\cR$- or $\cT$-ends maps under $\hol$ 
	homeomorphically to a union of components of 
	\[\mathcal{U}' \subset \rep_{{\mathcal E}, \lh}(\pi_1(\mathcal{O}), \PGL(n+1, \bR)).\]
	\item 
	The deformation space $\SDef_{{\mathcal E}, s_{\mathcal U}, \lh}(\mathcal{O})$ of SPC-structures on $\mathcal{O}$ with lens-shaped or horospherical $\cR$- or $\cT$-ends maps under $\hol$
	homeomorphically to the union of components of  
	\[\mathcal{U}' \subset \rep_{{\mathcal E}, \lh}(\pi_1(\mathcal{O}), \PGL(n+1, \bR)).\]
\end{itemize}
Furthermore, $\mathcal{U}'$ has to be in 
$\rep^s_{{\mathcal E}, \lh}(\pi_1(\mathcal{O}), \PGL(n+1, \bR))$. 
\end{corollary}
\begin{proof}\let\qed\relax
	We define $\widetilde{\CDef}_{{\mathcal E}, \lh}(\mathcal{O})$ 
	and $\widetilde{\SDef}_{{\mathcal E}, \lh}(\mathcal{O})$ as above. 
	Let $\widetilde{\mathcal{U}}$ be the inverse image of $\mathcal{U}$. 
We show that the image of $\widetilde{\mathcal{U}}$ under $\hol$ in 
$\mathcal{U}$ is closed and consists of stable irreducible characters. 
Now, Theorem \ref{cl-thm-closed1} implies the result. 



We prove this by lifting to $\SI^n$. 
Using Theorem \ref{prelim-thm-lifting}, 
let 
\[\{h_i : \pi_1(\orb) \ra \SL_\pm(n+1, \bR)\}\]
be a sequence of holonomy homomorphisms of real projective structures
corresponding to liftings of elements of $\widetilde{\mathcal{U}}$.
These are stable and strongly irreducible representations by Theorem \ref{intro-thm-sSPC}. 
Let $\Omega_i$ be the sequence of associated properly convex domains in $\SI^n$, and $\Omega_i/h_i(\pi_1(\orb))$ is 
diffeomorphic to $\orb$ and has the structure that lifts an element of 
$\CDef_{{\mathcal E}, \lh}(\mathcal{O})$. 
We assume that $\{h_i\} \ra h$ algebraically; 
that is, for a fixed set of generators $g_1, \dots, g_m$ of $\pi_1(\orb)$, 
$\{h_i(g_j)\} \ra h(g_j) \in \SLpm$ as $i \ra \infty$. 
 The limit $h$ is a discrete representation 
by Lemma 1.1 of Goldman-Millson \cite{GM88}.
We show that $h$ is a lifted holonomy homomorphism of 
an element of $\CDef_{{\mathcal E}, \lh}(\mathcal{O})$, and hence $h$ is stable and strongly irreducible. 

Here, we are using the definition of convexity for $\SI^n$ as given in 
Definition \ref{prelim-defn-convexitySn}. 
Since the Hausdorff metric space $\bdd_H$ of compact subsets of $\SI^n$ is compact, 
we may assume that $\{\clo(\Omega_i)\} \ra K$ for a compact convex set $K$
by taking a subsequence if necessary as in \cite{CG93}.
We take a dual domain $\Omega_i^* \subset \SI^{n \ast}$. Then the sequence $\{\clo(\Omega_i^*)\}$ also geometrically converges  
to a convex compact set $K^*$ by Proposition \ref{prelim-prop-dualHausdorff}. 
(See Section \ref{prelim-sub-Eduality}.) 

	Recall the classification of compact convex 
sets in Proposition \ref{prelim-prop-classconv}. 
For any 
$1$-form $\alpha$ positive on the cone $C_K$, any sufficiently close $1$-form is still positive on $C_K$. 
If $K$ has an empty interior and properly convex, then we can easily show that $K^*$ has a nonempty interior.  
Also, if $K^*$ has an empty interior and properly convex, $K$ has a nonempty interior. 

(I) The first step is to show that at least one of $K$ and $K^*$ has a nonempty interior. We divide the work into four cases (i)-(iv) where the types change for R-ends and T-ends. 

(i) To begin, suppose that 
there exists a radial p-end $\tilde E$ for $\Omega_i$ and $h_i$ and the type does not become horospherical. 
We may assume that $\mbv_{\tilde E, h_i} = \mbv_{\tilde E, h}$ by conjugating 
$h_i$ by a bounded sequence of projective automorphisms. 
For each $i$, $h_i$ acts on a lens cone $\mbv_{\tilde E, h_i} \ast L_i$
in $\Omega_i$ for each p-end $\tilde E$. 
By Theorems \ref{pr-thm-lensclass} and \ref{pr-thm-redtot}, 
$L_i$ can be chosen to be the convex hull of the closure of 
the union of attracting fixed sets of elements of $h_i(\pi_1(\tilde E))$. 
Hence, $h$ is also in it, and  Theorem \ref{pr-thm-equiv} shows that 
there exists a distanced compact convex 
set $L$ distanced away from a point $x$, and the lens cone
$\{\mbv_{\tilde E, h}\} \ast L -\{\mbv_{\tilde E, h} \}$ has a nonempty interior. 

Choose an element $g_1 \in \pi_1(\tilde E)$ such that $h(g_1)$ is positive bi-semiproximal
by Theorem \ref{prelim-thm-semi}.
Then $h_i(g_1)$ is also positive bi-semiproximal for sufficiently large $i$ since 
$\{h_i(g_1)\} \ra h(g_1)$ as a sequence.
We may assume that 
\[\{A_{h_i(g_1)}\} \ra A'_{h(g_1)} \subset A_{h(g_1)} \subset \Bd L\] 
for attracting-fixed-point sets $A_{h_i(g_1)}$ and $A_{h(g_1)}$ and a compact 
subset $A'_{h(g_1)}$ since the limits of sequences of eigenvectors are eigenvectors. 
Since $A'_{h(g_1)}$ is a geometric limit of a sequence compact convex sets, 
it is compact and convex. (See Section \ref{prelim-sub-semi}.)

We claim that the convex hull 
\[ \CH\left(\bigcup_{q\in h(\pi_1(\tilde E))} h(q)A'(h(g_1)) \cup \{\mbv_{\tilde E}\}\right)\]
has a nonempty interior: 
Suppose not. Then it is in a proper subspace where $\pi_1(\tilde E)$ acts. 
This means that $\pi_1(\tilde E)$ is virtually factorizable, and $\tilde E$ is a totally geodesic R-end according to Theorem \ref{pr-thm-redtot}. 
We choose another $g_2$ with $A(h(g_2))$ is not contained in this subspace
by Proposition \ref{prelim-prop-Ben2}. 
Again, we find $A'(h(g_2))\subset A(h(g_2))$. 
Now, 
\[ \CH\left(\bigcup_{q\in h(\pi_1(\tilde E)),g=g_1, g_2} h(q)A'(h(g)) \cup \{\mbv_{\tilde E}\}\right)\]
is in a strictly larger subspace where $h(\pi_1(\tilde E))$ acts on. 
By induction, we stop at a certain point, and we obtain
\[ \CH\left(\bigcup_{q\in h(\pi_1(\tilde E)), g=g_1, \dots, g_m} h(q)A'(h(g)) \cup \{\mbv_{\tilde E}\}\right)\]
that is not contained in a proper subspace. 

It is easy to show $h(q)A'(h(g)) = A'(h(qgq^{-1}))$. 
There exist a finitely many elements $g_1, \dots, g_m$ in $\pi_1(\tilde E)$
such that the attracting fixed set 
\[ \CH\left(\{a_1, \dots, a_m, \mbv_{\tilde E, h} \}\right), a_j \in A'(h(g_j)) \]
has a nonempty interior for some choice of $a_j$. 


We have $A(h_j(g_i)) \subset \clo(L_j)$ for each $j$ and $i$. 
The sequence $\{A(h_j(g_i))  \}$ accumulates only to points of $A'(h(g_j))$
since $\{h_j(g_i)\} \ra h(g_i)$. 
There is a sequence $\{a_{i, j}\}$, $a_{i, j} \in A(h_j(g_i))$, that converges
to $a_j \in A'(h(g_j))$ as $i \ra \infty$. 
Lemma \ref{prelim-lem-fch} implies that  the sequence 
\[\{\CH(\{ a_{1, j}, \dots, a_{m, j},  \mbv_{\tilde E, h_i}  \} )\} \hbox{ in } \Omega_j\] 
converges to $\CH(\{a_1, \dots, a_m, \mbv_{\tilde E}\})$ geometrically. 
Since $\clo(\Omega_j) \ra K$ as $j \ra \infty$, 
\[ \CH(\{a_1, \dots, a_m, \mbv_{\tilde E, h} \}) \subset K,\] 
by Proposition \ref{prelim-prop-BP};  
thus, $K$ has a nonempty interior in the case. (i) is accomplished. 

(ii) Suppose that there is a lens-shaped totally geodesic p-end $\tilde E$ for $\torb$ and the holonomy group $h(\pi_1(\tilde E))$, 
and the type for $\Omega_i$ and $h_i(\tilde E)$ does not
become horospherical. Then the dual 
$\Omega_i^*$ and $K^*$ have a nonempty interior by above arguments 
since $\Omega_i^\ast$ has lens-shaped R-p-ends by Corollary \ref{pr-cor-duallens2}.  
Now, we do the argument for (i) and use the duality at the end
by Proposition \ref{prelim-prop-dualHausdorff}.


(iii) Suppose that there is a lens-shaped R-p-end $\tilde E$ for $\torb$ and the holonomy group $h(\pi_1(\tilde E))$, 
and the type for $\Omega_i$ and $h_i(\tilde E)$ 
becomes horospherical. 

This is sufficient for (I) since when a lens-type T-end changes to a horospherical $\cR$- or $\cT$-end, 
we can use the duality as in (ii).


We have $\tilde \Sigma_{\tilde E, h}\subset \SI_{\mbv_{\tilde E}}^{n-1}$ 
is a complete affine space as $h|\pi_1(\tilde E)$ is parabolic. 
We have a properly convex domain or a complete affine 
space $\tilde \Sigma_{\tilde E, h_i} \subset \SI_{\mbv_{\tilde E_i}}^{n-1}$. 
By multiplying the developing maps by a convergent sequence of elements of 
$\SL_\pm(n+1, \bR)$, we may assume 
$\mbv_{\tilde E} = \mbv_{\tilde E_i}$.  

Recall that the map from the deformation space of real projective structures on
a closed orbifold to the representation space is a local homeomorphism
according to Theorem \ref{op-thm-dgorb}.
(See \cite{dgorb}.) 
Corollary \ref{abelian-cor-smvar} shows that 
$\clo(\tilde \Sigma_{\tilde E, h_i}) \ra \tilde \Sigma_{\tilde E, h}$ 
by passing to subsequences. 


Now $h(\pi_1(\tilde E))$ is the algebraic limit $h_i(\pi_1(\tilde E))$. 
Then $\mathcal{P} \cap h(\pi_1(\tilde E))$ is a lattice in a cusp group
$\mathcal{P}$. We conjugate $\mathcal{P}$ such that it is 
a standard unipotent cusp group in $\SO(n, 1) \subset \SLnp$. 
\index{cusp group} 

We choose any great segment $s_j$ with vertex $\mbv$  in the direction of 
$\tilde \Sigma_{\tilde E, h}$. 
We may assume that these are all inside $\tilde \Sigma_{\tilde E, h_i}$
since this is true for sufficiently large $i$. 
Choose finitely many $s_j$, $j=1, \dots, m$, such that the directions in the convex domain 
$\tilde \Sigma_{\tilde E, h_i}$ has a convex hull with a nonempty interior. 

Suppose that 
\begin{equation} \label{cl-eqn-li0}
 \{\bdd\hbox{-length}(l_{i, j})\} > \eps \hbox{ for } l_{i, j}:= \Omega_i \cap s_j, j=1, \dots, m,
\end{equation}
for a uniform $\eps > 0$. 
Let $l_j$ be the geometric limit of a subsequence $\{l_{i, j}\}$ with a nonzero $\bdd$-length. 
We assume without loss of generality that $\{l_{i,j}\} \ra l_j$ for each $j=1, \dots, m$,
and $h_i \ra h$. 
Then 
$\{h_{i}(g)(l_{j,i})\} \ra h(g)(l_j)$ for $g \in \pi_1(\tilde E)$.
For any finite set $F \subset  \pi_1(\tilde E)$, 
the set \[\{ h_{i}(g)(l_{j,i})| g \in F \} \ra 
\{h(g)(l_j)| g \in F \}\] geometrically. 
By Lemma \ref{prelim-lem-fch}, 
we have the geometric convergence of the sequence of convex hulls 
\[\{\CH( \bigcup_{g \in F} h_{i}(g)(l_{j,i}))\} \ra 
\CH( \bigcup_{g \in F} h(g)(l)).\] 
Since  the latter set has 
a nonempty interior by our assumption on $\bdd$-lengths of $l_{i,j}$, and 
$h(g), g\in F,$ is in a group conjugate to a cusp group, 
it follows that 
the convex hull 
\[\CH(h_{i}(\pi_1(\tilde E))(l_{j,i})) \subset \Omega_{i}\]
contains a fixed open ball $B$ for sufficiently large $i$. 
This means $B \subset K$ showing (I). 

%

Now for the final case, we suppose
\begin{equation} \label{cl-eqn-li02}
\{\bdd\hbox{-length}(l_{j, i})\} \ra 0 \hbox{ for } l_{j, i}:= \Omega_i \cap s_j, j=1, \dots, m
\hbox{ as } i \ra \infty.
\end{equation}
\cpr
\end{proof} 

\begin{lemma} \label{cl-lem-horob2} 
Let $\mbv$, $\mbv =\llrrparen{1, 0, \dots, 0}$ be a fixed point of the standard unipotent cusp group $\mathcal{P}$ and let $L$ be a lattice in $\mathcal{P}$.
Let 
\[ h_0\in \Hom^s_{\mathcal{E}, \lh}(\pi_1(\mathcal{O}),\bG)\]
for $\bG = \PGL(n+1, \bR)$ or $\SL_\pm(n+1, \bR)$
be a representation to a cusp group. 
Let $H$ be a $\mathcal{P}$-invariant hemisphere with $\mbv$ in the boundary, and let $l$ be the maximal  line 
with endpoints $\mbv$ and $\mbv_-$ transversal to $\partial H$ with respect to $\bdd$.  
Then 
there exists a finite subset $F$ of $L$ and a neighborhood $N$ of $h_0$ 
such that the following holds\/{\rm :}
\begin{itemize} 
\item for any point $x \in l$ and a hyperspace at $x$ transversal to $l$ 
bounding a closed hemisphere $H_x$,   
\[I^{h}_x :=\bigcap_{g\in F} h(g)(H_x), h \in N \] 
is a properly convex domain, and 
\item as $x \ra \mbv$, 
$\{I^{h}_x\}$ geometrically converges to $\{\mbv\}$ uniformly for $h\in N$. 
\end{itemize} 
\end{lemma}
\begin{proof} 
We choose the affine coordinate system of $H^o$.
If $F$ is large enough, then $\{h_0(g)(\partial H_x)| g \in F\}$ is in a general position
$x \in l$.   We may also assume the conclusion for $h_0$. 
Now, we choose $N$ so that geometrically $h(g)(\partial H_x)$ is uniformly close to
$h_0(g)(\partial H_x)$ for all $x\in l, h \in N$. 
\end{proof}



\begin{proof}[Proof of Proposition \ref{cl-cor-closed1} continued] \let\qed\relax
We may assume $\mbv_{\tilde E, h_i}=\mbv_{\tilde E, h} = \mbv$ without loss of generality
by changing the developing map by a sequence of bounded automorphisms $g_i$. 
Let $H$ denote the $\mathcal{P}$-invariant hemisphere containing $K$. We assume that $\Omega_i \subset H$.
We choose a finite set $\mathbf{F}$ by Lemma \ref{cl-lem-horob2}. 
Let $l$ be a line of $\bdd$-length $\pi$ from $\mbv$ transverse to $\partial H$.  
Let $x_i$ be the end point of $l \cap \Omega_i$ other than $\mbv$. 
By Lemma \ref{cl-lem-horob2}, 
\[\hat K_i := \bigcap_{g\in F} h_j(g)(H_{x_i}) \cap H\] is properly convex for sufficiently large $j$  since $\{h_j(g)\} \ra h(g), g\in F$. 
This set $\hat K_i$ contains $\clo(\Omega_i)$ since $H_x \supset \clo(\Omega_i)$. 
It follows that $\{\hat K_i\} \ra \{\mbv\}$.
%
Therefore, we conclude that $K$ is a singleton. 


By Proposition \ref{prelim-prop-dualHausdorff}, 
$\{\Omega_i^\ast\}$ geometrically converges to $K^\ast$ dual to $K$
as $i \ra \infty$. 
(See \ref{prelim-sub-Eduality}.) 
If $K$ is a singleton, $K^*$ must be a hemisphere
by Proposition \ref{prelim-prop-NPduality}. 
We now conclude that $K$ or its dual $K^*$ has a nonempty interior. 


Thus, by choosing $h_i^*$ and $h^*$ if necessary, we may assume without loss of generality that $K$ has a nonempty interior. 
We show that $K$ is a properly convex domain and this implies that so is $K^*$. 

(II) The second step is to show that $K$ is properly convex.

Assume that $h(\pi_1(\orb))$ acts on a convex open domain $K^o$.
We may assume that $K^o \subset \mathbb{A}$ for an affine subspace $\mathbb{A}$
and $\Omega_i \subset \mathbb{A}$ by acting by an orthogonal 
$\kappa_i \in \Ort(n+1, \bR)$, where 
$\{\kappa_i\}$ is converging to $\Idd$. 
We can accomplish this by moving $\Omega_i$ into $\mathbb{A}$. 
Since $\{\kappa_i\} \ra \Idd$, we still have $\{\clo(\Omega_i)\} \ra K$ by Lemma \ref{prelim-lem-geoconv}. 
Take a ball $B_{2C}$ of $\bdd$-radius $2C$, $C> 0$, in $K$. 
By Lemma \ref{prelim-lem-bdconv}, 
$\mathbb{A}$ contains a $\bdd$-radius $C$ ball $B_C \subset \Omega_i$
for sufficiently large $i$. 
Without loss of generality, assume $B_C \subset \Omega_i$ for all $i$. 
Choose the $\bdd$-center $x_0$ of $B_C$ as the origin in the affine coordinates. 

Let $g_1, \dots, g_m$ denote the set of generators of $\pi_1(\orb)$. 
Then, by extracting subsequences, we may assume without loss of generality that
$\{h_i(g_j)\} $ converges to $h(g_j)$ for each $j=1, \dots, m$. 
\cpr
\end{proof} 

\begin{lemma} \label{cl-lem-gjC0} 
	For each $g_j$, $j=1, \dots, m$,  
\begin{equation}\label{cl-eqn-gjC0} 
\bdd(h_i(g_j)(x_0), \Bd \Omega_i) \geq C_0 \hbox{ for a uniform constant } C_0.
\end{equation}  
\end{lemma} 
\begin{proof}
Suppose not. Then there is a sequence of $\bdd$-length constant $C$ segments $s_i$, 
$s_i \subset \Omega_i$, with an endpoint $x_0$, and $s_i$
is sent to the segment $h_i(g_j)(s_i)$ in $\Omega_i$ with 
endpoint $h_i(g_j)(x_0)$ and lying on the shortest 
$\bdd$-length segment from $h_i(g_j)(x_0)$ to $\Bd \Omega_i$. 
Thus, the sequence of the $\bdd$-length of $h_i(g_j)(s_i)$ is going to zero. 
This implies that $h_i(g_j)$ is not in a compact subset of $\SL_{\pm}(n+1, \bR)$, a contradiction. 
\end{proof} 

\begin{proof}[Proof of Proposition \ref{cl-cor-closed1} continued]
We may assume without loss of generality that $x_0$ is the origin and $\clo(\Omega_i)$ is 
in a fixed affine space $H^o$ for an $n$-hemisphere $H$. 

By estimation from \eqref{cl-eqn-gjC0}, and the cross-ratio expression of
the Hilbert metric, 
a uniform constant $C$ satisfies
\begin{equation}\label{cl-eqn-bdbel}
d_{\torb_i}(x_0, h_i(g_j)(x_0)) < C.
\end{equation}

Since $\clo(\Omega_i) \subset H^o$ and $\clo(\Omega_i) \ra K'$ for a convex set $K'$ with a non-empty interior in the Hausdorff topology, 
it follows that $\Omega_i \ra K_1 \subset H^o$  for a closed convex domain $K_1$ with a nonempty interior in the Chabauty topology by Proposition \ref{prelim-prop-BP}. 

If $K_1$ is a properly convex domain in $H^o$,  the result follows. 
Hence, we assume that $K_1$ is not bounded. 

We claim that
given a sequence of properly convex domains $\Omega_i \subset H^o$ 
Chabauty converging to $K_1$, 
 there exists a sequence of elements $M_i$ in a common diagonalizable group 
fixing $x_0$ as the largest norm eigenvalues and acting on $H^o$
such that $\{M_i(\Omega_i)\}$ geometrically converges to a properly convex domain with 
$x_0$ in its interior up to choosing a subsequence: 

Give $H^o$ a Euclidean Riemannian metric identical with the Riemannian metric 
for $\bdd$ at $x_0$. We use the product structure. 

We prove this by induction on dimension $n$. 
If $n=1$, it is clear.  Suppose that the conclusion holds for the 
dimension equal to $n-1$.

There exist maximal $r_i$ and minimal $R_i$ 
such that $B_{r_i}(x_0) \subset \Omega_i \subset B_{R_i}(x_0)$. 
We know $R_i \ra \infty$. 
Then $\Bd \Omega_i$ meets $\Bd B_{r_i}(x_0)$. 
We choose a projective automorphism $S_i$ that fixes $x_0$ and each point of $\partial H^0$ with the strongest eigenvalue at $x_0$
such that $S_i(B_{r_i}(x_0)) = B_1(x_0)$.
Up to choosing a subsequence, we may assume that 
the Chaubauty limit of $S_i(\clo(\Omega_i))$ is a closed convex set in $H^o$. 
The Chaubauty limit cannot be complete affine since we have a sharply 
supporting affine hyperspace at a point of $S_i(\Bd(B_{r_i}(x_0)\cap \Bd \Omega_i)$. 

If the limit is properly convex, then the proof is finished.

Suppose not. Then 
the limit equals $K \times A$ for an affine subspace $A$
for a closed properly convex domain $K \subset A'$ where $A$ and $A'$ are complete affine subspaces 
passing $x_0$ orthogonal to each other in $H^o$.

We respectively denote by $\Pi_A:H^o \ra A$ and $\Pi_{A'}: H^o \ra A'$ the orthogonal projections. 


We now apply a projective automorphism 
$L_i$ acting on $H^o$  such that  
as an affine map $L_i$ has a diagonalizable 
linear part with eigenvalues $\lambda_i$ for vectors in $A$ 
and $1$ for vectors in $A'$ only.

Then there exists a maximal ball $B_{r'_i}(x_0) \subset \Pi_A(\Omega_i)$, and
we choose $L_i$ such that $\lambda_i = 1/r'_i$. 
Then $L_i\circ S_i(\Omega_i)$ still contains $B_1(x_0) \cap A'$
and $\Pi_{A}\circ L_i \circ S_i (\Omega_i)$ contains $B_1(x_0) \cap A$. 
Hence, $L_i\circ S_i(\clo(\Omega_i))$ Chaubauty converges to a closed convex set 
containing $B_1(x_0)\cap A'$ and whose projection to $A'$ 
contains $x_0$ in its interior up to a choice of subsequences. 
We choose $\dim A+1$ points in general position in $A$ whose affine convex hull contains $x_0$. 
Then there exists a corresponding point in $L_i \circ S_i(\Omega_i)$. 
The convex hull of these points with $B_1(x_0) \cap A'$ has a nonempty interior and 
is contained in $L_i \circ S_i(\Omega_i)$. 

Since $S_i(\clo(\Omega_i))$ Chaubauty converges to  
$K \times A$ for a properly convex set $K \subset A'$, 
it follows that
$S_i(\clo(\Omega_i))$ is contained in $B \times A$ for a closed properly convex subset $B$ in $A'$
by Proposition E.12 of Benedetti-Petronio \cite{BP92} for sufficiently large $i$. 
Since $L$ acts on every $t\times A$ for $t\in A'$, 
the Chaubauty limit of $L_i \circ S_i(\clo(\Omega_i))$ is a subset of $B\times A$. 

By the induction hypothesis, there is a sequence of projective maps $M'_i$ acting on 
the great sphere that is the span of $A$ such that $M'_i\circ L_i(\Pi_A(\Omega_i))$ 
Chabauty converges to a properly convex domain $D$ in $A$. 
Then we extend $M'_i$  to $M_i$ on $\SI^n$ such that $M_i$ fixes all points of $A'$. 

Note that if $K_1$ or the Chaubauty limit after any steps were already properly convex, then 
we take the automorphism sequences such as $L_i$ and $M_i$ to be identity maps. 

It follows that $M_i\circ L_i \circ S_i(\Omega_i)$ Chaubauty converges to a closed convex 
subset $\hat \Omega$ of $B \times D$. Since $\Pi_A(\hat \Omega) = D$, 
there is at least one point of $\hat \Omega$ for each point of $D$. 
Hence, we may choose $\dim A +1$ points in general position in $D$ whose affine convex hull  contains $\Pi_A(x_0)$ in its interior. 
Then there is a corresponding set of points in $\hat \Omega$. 
The join of this set with $B_1(x_0) \cap A'$ is a subset of $D$ with a nonempty interior.  
  


Hence, up to multiplying $M_i\circ L_i \circ S_i$ by a bounded projective automorphism, 
we obtain a sequence 
$\tau_i \in  \SL_\pm(n+1, \bR)$ acting on $x_0$ 
such that 
\[ B_1(x_0) \subset \tau_i(\Omega_i) \subset B_{R_B}(x_0) \]
up to normalizing.  
Now, $\tau_i h_i(\pi_1(\orb)) \tau_i^{-1}$ acts on $\tau_i(\Omega_i)$.

Let $g_1, \dots, g_m$ denote the generators of $\pi_1(\orb)$.  
Also, we may assume without loss of generality that 
$\tau_i \circ h_i(g_j) \circ \tau_i^{-1}(x_0) \ra y_j$ for each $j=1, \dots, m$
for $y_j$ uniformly bounded away from $\tau_i(\partial \Omega_i)$. 
by \eqref{cl-eqn-bdbel} since $M_i\circ L_i \circ S_i$
is commonly diagonalizable with the largest norm eigenvalue at $x_0$
fixing on $\partial H$.  
By Theorem 7.1 of Cooper-Long-Tillmann \cite{CLT15}, 
we obtain that $\tau_ih_i(g_j)\tau_i^{-1}$ for $j=1, \dots, n$ 
are in a compact subset of $\SL_{\pm}(n+1, \bR)$ independent of $i$. 

Therefore, up to choosing subsequences, we
obtain that 
$\{\tau_i(\Omega_i)\}$ geometrically converges to a properly convex domain $\hat K$ in $B_R$ containing $B_1$, 
and \[\{\tau_i h(\cdot) \tau_i^{-1}:  \pi_1(\orb) \ra \SL_\pm(n+1, \bR) \} \]
algebraically converges to a holonomy homomorphism 
 \[ h': \pi_1(\orb) \ra \SL_\pm(n+1, \bR).\] 
And the image of $h'$ acts on the interior of the properly convex domain $\hat K$. 

Suppose that the sequence $\{\tau_i\}$ is not bounded. Then 
$\tau_i = k_i d_i k'_i$ where $d_i$ is diagonal with respect to a standard 
basis of $\bR^{n+1}$ and $k_i, k'_i \in O(n+1, \bR)$ by the KTK-decomposition of $\SLpm$. 
Then the sequence of the maximum modulus of the eigenvalues of $d_i$ are not bounded above. 
We assume without loss of generality 
\[\{k_i\} \ra k, \{k'_i\} \ra k' \hbox{ in } O(n+1, \bR).\]
Thus, 
$\{k_i' h_i(g_j) k^{ -1}_i\}$ converges to $k' h(g_j) k^{ -1}$ for $k' \in O(n+1, \bR)$. 
Since \[\{k_i d_i k'_i h_i(g_j) k^{\prime -1} d_i^{-1} k^{-1}_i\}\] is convergent
to $h'(g_j)$, we obtain 
\[ \{d_i k'_i h_i(g_j) k^{\prime -1}_i d_i^{-1}\} \ra k^{-1} h'(g_j) k
\hbox{ for each } j. \]
Thus, $\{ d_i k'_i h_i(\pi_1(\orb)) k^{\prime -1}_i d_i^{-1}\}$ algebraically converges to a group $k h'(\pi_1(\orb)) k^{-1}$
acting on $k^{-1}(\hat K)$. 

Since the sequence of 
the norms of $d_i$ is divergent, $k h'(\pi_1(\orb)) k^{-1}$ is reducible: 
We may assume by passing to a subsequence and a change of coordinates that 
the diagonal entries of $d_i$ satisfy 
\[ d_{i, 1}\geq d_{i,2} \geq \cdots   \geq d_{i, n+1}. \]
By passing to a subsequence, 
there is $j$ such that $d_{i, k}/ d_{i, j} \geq 1$ for $k \leq j$ 
and $\{d_{i, k}/ d_{i, j}\} \ra 0$ for $k > j$. 
Then $\{ d_i k'_i h(\pi_1(\orb) k^{\prime -1}_i d_i^{-1}\}$ being a bounded sequence
converges to a matrix with the entries at
$(k+1, \dots, n+1)\times (1, \dots, k)$ being identically zero. 
(Compare to the proof of Lemma 1 of \cite{dgorb}.)

By Lemma \ref{cl-lem-endspres}, 
$k^{-1}(\hat K)^o/ k h'(\pi_1(\orb)) k^{-1}$ is a strongly tame SPC-orbifold with horospherical or generalized lens-shaped ends. 

%
%
%
By Theorem \ref{intro-thm-sSPC}, the algebraic limit
of \[\{ d_i k'_i h(\pi_1(\orb) k^{\prime -1}_i d_i^{-1}\}\]
cannot be reducible. 
Therefore, the sequence of the norms of $d_i$ is uniformly bounded.
This is a contradiction to the unboundedness of $\tau_i$. 


By Lemma \ref{cl-lem-endspres}, we obtain that 
$\orb_h:=\hat K^o/h(\pi_1(\orb))$ is a strongly tame  SPC-orbifold
with generalized lens-shaped or horospherical $\cR$- or $\cT$-ends diffeomorphic to $\orb$. 
This completes the proof for the statement
$\widetilde{\mathcal{U}} \subset \widetilde{\CDef_{{\mathcal E}, \lh}(\mathcal{O})}$.


To prove for $\SDef_{{\mathcal E}, \lh}(\mathcal{O})$, we additionally need
Theorems \ref{rh-thm-relhyp} and \ref{rh-thm-converse} as in 
the last paragraph of the proof of Theorem \ref{cl-thm-closed1}. 
This completes the proof of Corollary \ref{cl-cor-closed1}.  
\end{proof}

\section{General cases without the uniqueness condition: The proof of Theorem \ref{CL-THM-D}.} \label{cl-sec-general} 


We construct a section by the following. 
Let 
\[\Hom_{\mathcal{E}, \lh}(\pi_1(\tilde E), \PGLnp)
\hbox{ (resp. } \Hom_{\mathcal{E}, \lh}(\pi_1(\tilde E), \SLnp))\] 
denote the space of representations $h$ fixing a common fixed point $p$ and acting on a lens $L$ of a lens cone of form $\{p\}\ast L$ with $p \not\in \clo(L)$ or is horospherical with a fixed point $p$. 

Let \[\Hom_{\mathcal{E}, \lh}(\pi_1(\tilde E), \PGLnp) 
\hbox{ (resp. } \Hom_{\mathcal{E}, \lh}(\pi_1(\tilde E), \SLnp))\] 
denote the space of representations 
where $h(\pi_1(\tilde E))$ for each element $h$ acts on $P$ satisfying the lens-condition
or acts on a horosphere tangent to $P$.  
(See Section \ref{intro-sub-semialg}.)

We define the sections
\begin{multline*}
s_R: \Hom_{\mathcal{E}, \lh}(\pi_1(\tilde E), \PGLnp) \ra \RP^n, 
s_T: \Hom_{\mathcal{E}, \lh}(\pi_1(\tilde E), \PGLnp) \ra \RP^{n\ast}
\end{multline*}
\begin{multline*} 
\big(\hbox{resp. } s_R: \Hom_{\mathcal{E}, \lh}(\pi_1(\tilde E), \SLnp) \ra \SI^n, 
s_T: \Hom_{\mathcal{E}, \lh}(\pi_1(\tilde E), \SLnp) \ra \SI^{n\ast}
\big)
\end{multline*} 
by Propositions \ref{cl-prop-section} and \ref{cl-prop-section2}. 

\begin{proposition} \label{cl-prop-uqs}
	The maps $s_R$ and $s_T$ for $\RP^n$ and $\RP^{n\ast}$ 
	{\rm (}resp. $\SI^n$ and $\SI^{n\ast}$\/{\rm )} are both continuous.
	\end{proposition}
\begin{proof} 
	We only prove this for $\RP^n$. The version for $\SI^n$ is similar. 
Consider first the case for $s_R$. 
	Let 
\[h \in \Hom_{\mathcal{E}, \lh}(\pi_1(\tilde E), \PGLnp).\] 
 The  vertex of a lens cone is a common fixed point of all elements of 
 $h(\pi_1(\tilde E))$. Let $F$ be the set of generators of $\pi_1(\tilde E)$
 such that $\{v\} =\{ w| g(w) = w, g \in F \}$.
 Otherwise, we would 
 have a line of fixed points for $\Gamma_E$ and we obtain a contradiction 
 as in the proof of Proposition \ref{cl-prop-section}.
Hence, the holonomies of elements of $F$ determine the vertex.
The continuity follows from a sequence argument. 

For $s_T$,  we take the dual by Proposition \ref{pr-prop-dualend}
and prove continuity.
\end{proof}

\begin{lemma} \label{cl-lem-ceu_section}
	We can construct the uniqueness section of lens-type
	\[s: \Hom_{\mathcal{E}, \lh}(\pi_1(\orb), \PGLnp) \ra (\RP^n)^{e_1} \times 
(\RP^{n \ast})^{e_2}\]
\end{lemma} 
\begin{proof} 
	We can always choose a vertex and the hyperspace by 
	Propositions \ref{cl-prop-section} and \ref{cl-prop-section2}.
	$s$ is continuous by Proposition \ref{cl-prop-uqs}. 
\end{proof} 

\begin{proof}[Proof of Theorem \ref{CL-THM-D}]
	Using the uniqueness section of the lens types, we apply 
    Corollary \ref{cl-cor-closed1}.  
\end{proof}

\chapter{Nice cases}
\label{ex-sec-nicecase}

We now present the cases where the theory presented in
this monograph works best.

In Section \ref{ex-nec-main}, we discuss the results in which the results of this monograph apply.  
In Section \ref{ex-sec-examplesII}, we end with two examples in which the results apply. 

\section{Main results} \label{ex-nec-main}

Let us start with an example: 
\begin{example} \label{ex-exmp-endv2}
Let $M$ be a hyperbolic $3$-orbifold with possibly 
hyperideal ends where each end orbifold has 
a sphere or a disk as the base space. The end fundamental group is 
generated by a finite-order elements.  
By Lemma \ref{ex-lem-niceend}, a properly convex real projective
structure on $M$ has only lens-shaped or horospherical radial ends. 
\end{example}

We need the end-classification results from Chapters \ref{ch-ce}, \ref{ch-pr}, and \ref{ch-np} to prove the following. 
%
Let $g \in \pi_1(\orb)$. 
Using the choice of a representing matrix of $g$ as in Remark \ref{intro-rem-SL},
we let $\lambda_{x}(g)$ denote the eigenvalue of holonomy of $g$ associated with the vector in direction of $x$ if $x$ is a fixed point of $g$. 
\index{lambdax@$\lambda_x(\cdot)$}

Suppose that the holonomy group of $\pi_1(\orb)$ can be lifted to $\SLpm$ such that 
$\lambda_{\mbv_{\tilde E}}(g) = 1$ for the holonomy of every 
$g\in \pi_1(\tilde E)$
where $\mbv_{\tilde E}$ is a p-end vertex of a p-end $\tilde E$ corresponding to $E$. 
Then we say that $E$ or $\tilde E$ satisfies the {\em unit-middle-eigenvalue condition} with respect to $\mbv_{\tilde E}$ or the R-p-end structure. 
\index{unit middle-eigenvalue condition} 

Suppose that $E$ is a $T$-end. If the hyperspace containing 
the ideal boundary component $\tilde S_{\tilde E}$ 
of the p-end $\tilde E$ of $E$
corresponds to $1$ as the eigenvalue of the dual of the holonomy of every 
$g \in \pi_1(\tilde E)$,  then we say 
that $E$ or $\tilde E$ satisfies the {\em unit middle-eigenvalue condition} with respect to $\tilde S_{\tilde E}$ or the T-p-end structure. 
\index{middle-eigenvalue condition!unit|textbf} 

\begin{lemma}\label{ex-lem-niceend} 
	Suppose that $\orb$ is a strongly tame convex real projective orbifold with radial ends. 
	Assume that the end fundamental group $\pi_1(E)$ of an end $E$ satisfies {\rm (}NS\/{\rm)}. 
	Let $E$ be an $R$-end, or is a $T$-end. Suppose that one of the following holds{\rm :}
	\begin{itemize} 
	\item $\pi_1(E)$ is virtually generated by finite order elements or is simple,	or 
	\item the end holonomy group of $E$ satisfies the unit middle-eigenvalue condition. 
	\end{itemize} 
	Then the following hold{\rm :} 
	\begin{itemize} 
	\item the end $E$ is either a properly convex generalized lens-shaped R-end 
	or a lens-shaped T-end, 
	or is horospherical. 
	\item If the end $E$ furthermore has a virtually abelian end holonomy group, then $E$ is a lens-shaped R-end, a lens-shaped T-end, or is a horospherical end. 
	\end{itemize} 
\end{lemma} 
\begin{proof} 
We suppose that $\torb$ is a convex domain in of $\SI^n$. 
First, let $E$ be an R-end. 
The map \[g\in \Gamma_{\tilde E} \mapsto \lambda_{\mbv_{\tilde E}}(g) \in \bR_{+}\]
is a homomorphism. 
Thus, $\lambda_{\mbv_{\tilde E}}(g) = 1$ for $g \in \Gamma_{\tilde E}$
since the end holonomy group is simple or virtually generated by the finite order 
elements. 

Each R-end is either complete, properly convex, or is convex but neither properly convex nor complete affine by Section \ref{ce-sec-ends}.

Suppose that $\tilde E$ is complete. Then Theorem \ref{ce-thm-comphoro} 
shows that either $\tilde E$ is horospherical
or each element $g$, $g\in \pi_{1}(\tilde E)$ has at most two norms of eigenvalues 
where two norms for an element are realized.
Since the multiplication of all eigenvalues equals $1$, 
we obtain $\lambda_{1}^{n+1-r}(g) \lambda_{\mbv_{\tilde E}}(g)^{r} = 1$ for some integer $r, 1 \leq r \leq n$
and the other norm $\lambda_{1}(g)$ of the eigenvalues.
The second case cannot happen. 

Suppose that $\tilde E$ is properly convex. 
Then the uniform middle-eigenvalue condition holds with respect to R-p-end or T-p-end 
$\tilde E$ by Remark
\ref{pr-rem-eigenlem} since $\lambda_{\mbv_{\tilde E}}(g) = 1$ for all $g\in \bGamma_{\tilde E}$. 
(See Definition \ref{pr-defn-umec}.)
By Theorem \ref{pr-thm-secondmain}, $\tilde E$ is of generalized 
lens-type. 



 Finally, Corollary \ref{np-cor-eigenone} rules out the case when 
 $\tilde E$ is convex but neither properly convex nor complete affine. 

Now, let $E$ be a T-end. By dualizing the above, $E^\ast$ is a properly convex radial end, and 
$E^\ast$ satisfies the uniform middle-eigenvalue condition (see Definition \ref{pr-defn-umecD}).  Hence, $E$ satisfies one, and 
Theorem \ref{pr-thm-equ2} implies the result. 
\hfill \SSn {\parfillskip0pt\par}
\end{proof} 

%



\begin{theorem} \label{ex-thm-niceend}
Suppose that $\orb$ is a strongly tame properly convex real projective $n$-orbifold with R-ends or T-ends. 
Suppose that each end fundamental group satisfies property {\rm (}NS\/{\rm)} and 
is virtually generated by finite-order elements, or is simple
or satisfies the unit middle-eigenvalue condition. 
Then the holonomy is in 
\[\Hom^s_{\mathcal E, u, \lh}(\pi_{1}(\orb), \PGL(n+1, \bR)).\] 
\end{theorem} 
\begin{proof} 
Suppose that $E$ is an R-end. Let $\tilde E$ be a p-end corresponding to $E$ 
and $\mbv_{\tilde E}$ be the p-end vertex. 
By Lemma \ref{ex-lem-niceend}, we see that the R end is lens-type or horospherical.
 
We prove the uniqueness of the fixed point under $h(\pi_1(\tilde E))$:
Suppose that $x$ is another fixed point of $h(\pi_1(\tilde E))$. 
Since $\pi_1(\tilde E)$ is as in the premise, the eigenvalue
$\lambda_x(g)$ for every $g \in \pi_1(\tilde E)$ associated with 
$x$ is always $1$.  
In the horospherical case $x = \mbv_{\tilde E}$ since 
the cocompact lattice action on a cusp group fixes a unique point 
in $\RP^n$. 

Now consider the lens case. 
The uniform middle-eigenvalue condition with respect to 
$\mbv_{\tilde E}$ and $x$ holds by Remark
\ref{pr-rem-eigenlem} since $\lambda_{x}(g) = 1$ for all $g$. 
Lemma \ref{ex-lem-niceend} shows that $\pi_1(\tilde E)$ acts on 
a lens cone with vertex at $x$. 
Proposition \ref{cl-prop-section} implies 
the uniqueness of the p-end vertex. 

Suppose that $E$ is a $T$-end. The proof of Proposition \ref{pr-cor-duallens} 
shows that the hyperspace containing $\tilde S_{\tilde E}$ corresponds 
to $\mbv_{\tilde E^\ast}$ for the R-p-end $\tilde E^\ast$ 
corresponding to the dual of the T-p-end $\tilde E$
and vice verse. 
Hence, the result follows from the R-end part of the proof. 
	\end{proof} 

This was proved by Marquis in Theorem A of \cite{Marquis17} when the orbifold is a Coxeter one.

%
%



Theorems \ref{ex-thm-niceend}, \ref{intro-thm-closed1}, and 
\ref{intro-thm-sSPC} 
imply the following: 

\begin{corollary} \label{ex-cor-closed2}
Let $\mathcal{O}$ be an $n$-dimensional 
noncompact strongly tame SPC-orbifold with R-ends and T-ends
and satisfies  {\rm (}IE\/{\rm )} and {\rm (}NA\/{\rm )}.  
Suppose that each end fundamental group is generated by finite-order elements or is simple.  Suppose 
each end fundamental group satisfies {\rm (}NS\/{\rm )}. 
Assume $\partial \orb =\emp$, and 
that the nilpotent normal subgroups of every finite-index subgroup of $\pi_1(\mathcal{O})$ are trivial.
Then 
\[ \CDef_{{\mathcal E}}(\mathcal{O}) = \CDef_{{\mathcal E}, \mathrm{u}, \lh}(\mathcal{O})  \] 
and 
$\hol$
maps the deformation space $\CDef_{{\mathcal E}}(\mathcal{O})$ of SPC-structures
on $\mathcal{O}$ homeomorphic to 
a union of components of 
\[\rep^s_{{\mathcal E}, \mathrm{u}, \lh}(\pi_1(\mathcal{O}), \PGL(n+1, \bR))\]
which is also a union of components of 
\[\rep_{{\mathcal E}, \mathrm{u}, \lh}(\pi_1(\mathcal{O}), \PGL(n+1, \bR))
\hbox{ and }
\rep_{{\mathcal E}}(\pi_1(\mathcal{O}), \PGL(n+1, \bR)).
\]
The same can be said for 
\[\SDef_{{\mathcal E}}(\mathcal{O}) =\SDef_{{\mathcal E}, \mathrm{u}, \lh}(\mathcal{O}) .\]
\end{corollary} 

These types of deformations from structures with cusps to ones with lens-shaped ends are realized in our main examples
as stated in Section \ref{ex-sec-examples}.
We need restrictions on the target space 
since the convexity of $\orb$ is not preserved under the 
hyperbolic Dehn surgery deformations of Thurston, 
as pointed out by Cooper at ICERM in September 2013.

Virtually abelian groups clearly satisfy (NS). 
Since finite-volume hyperbolic $n$-orbifolds satisfy (IE) and (NA) 
(see P.151 of \cite{Marden07} for example),
strongly tame properly convex real projective orbifolds admitting 
compatible complete hyperbolic structures
end fundamental groups generated by finite-order elements 
satisfies the premise. 
Hence,  $2h\underbar{\,\,}1\underbar{\,\,}1$ and the double of the simplex orbifold discussed in Section \ref{ex-sec-examplesII} also do.

Since Coxeter orbifolds satisfy the above properties, we obtain
a simple case: 

\begin{corollary}  \label{ex-cor-closed3}
Let $\mathcal{O}$ be a strongly tame Coxeter $n$-dimensional orbifold, $n \geq 3$, with only $\cR$-ends 
admitting a finite-volume complete hyperbolic structure. 
Then 
\[\SDef_{{\mathcal E}, \mathrm{u}, \lh}(\mathcal{O})\]
is homeomorphic to 
a union of components of 
\[\rep^s_{{\mathcal E}, \mathrm{u}, \lh}(\pi_1(\mathcal{O}), \PGL(n+1, \bR))\]
which is also a union of components of 
\[\rep_{{\mathcal E}}(\pi_1(\mathcal{O}), \PGL(n+1, \bR))\]
Finally, 
 \[\SDef_{{\mathcal E}, \mathrm{u, \lh}}(\mathcal{O}) 
 = \SDef_{{\mathcal E}}(\mathcal{O}).\]
\end{corollary}

\section{Two specific examples}\label{ex-sec-examplesII}

\begin{figure}[t]
	\centering 
	\includegraphics[height=5cm]{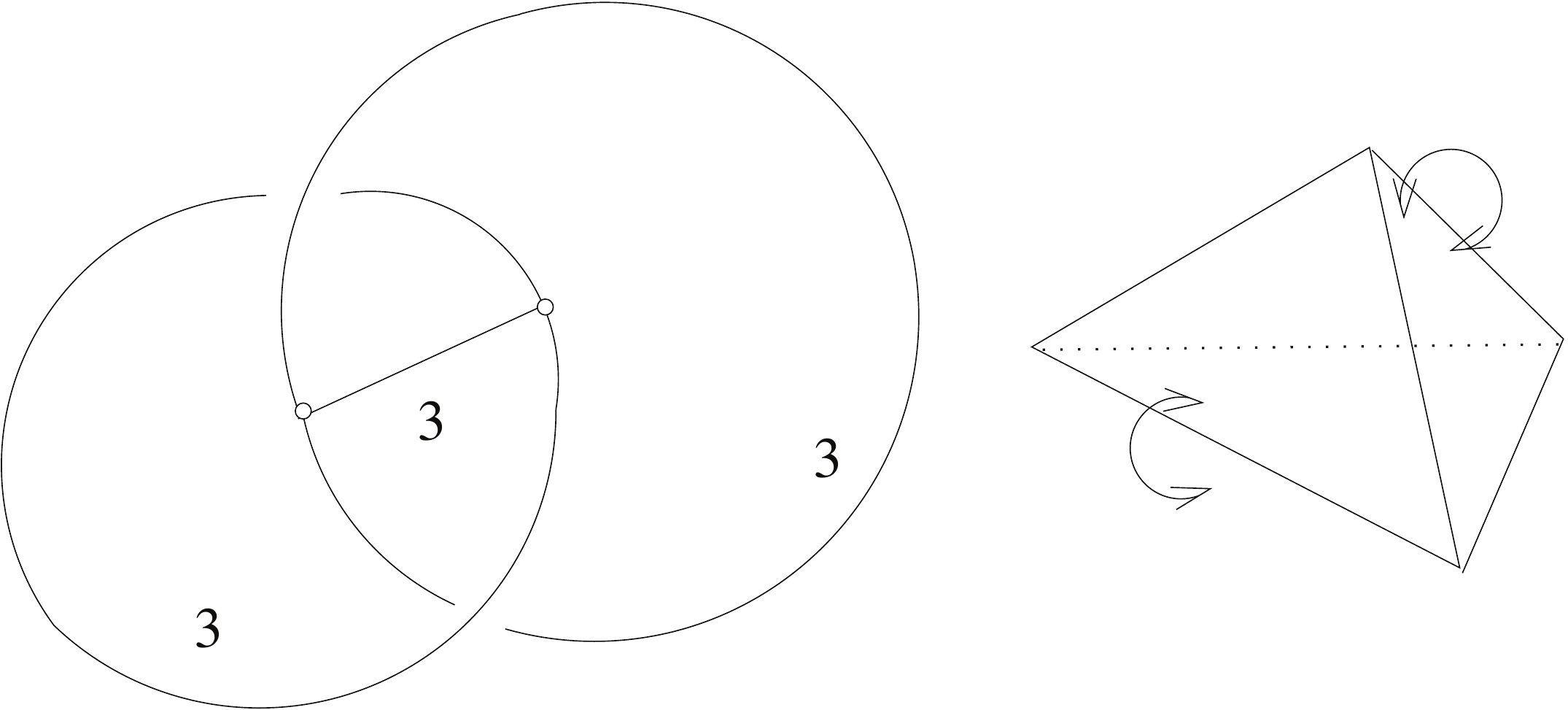}
	
	\caption{Handcuff orbifold.}
	\label{intro-fig-two}
\end{figure}

The example of S. Tillmann is an orbifold on a $3$-sphere with singularity consisting of 
\index{handcuff orbifold} 
two unknotted circles linking each other only once under a projection to a $2$-plane 
and a segment connecting the circles (looking like a linked handcuff) with vertices removed and all arcs 
given as local groups the cyclic groups of order three. (See Figure \ref{intro-fig-two}.)
This is one of the simplest hyperbolic orbifolds in the list of 
Heard, Hodgson, Martelli, and Petronio \cite{HHMP10} labeled $2h\underbar{\,}1\underbar{\,}1$. 
The orbifold admits a complete hyperbolic structure since we can start from a complete hyperbolic 
tetrahedron with four dihedral angles equal to $\pi/6$ and two equal to $2\pi/3$ at a pair of opposite edges $e_1$ and $e_2$.
Then we glue two faces adjacent to $e_i$ by an isometry fixing $e_i$ for $i= 1, 2$. 
The end orbifolds are two $2$-spheres with three cone points of orders equal to $3$ respectively. 
These end orbifolds always have induced convex real projective structures in dimension $2$, 
and real projective structures on them have to be convex. Each of these is either the quotient of a properly convex 
open triangle or a complete affine plane as we saw in Lemma 
\ref{ex-lem-niceend}. 

Porti and Tillmann \cite{PTp} found a two-dimensional solution set from the complete hyperbolic structure by explicit computations. 
Their main questions are the preservation of convexity and realizability as convex real projective structures on the orbifold. Corollary \ref{ex-cor-closed2} answers this since their deformation space identifies with 
$\SDef_{{\mathcal E}, \mathrm{u, \lh}}(\mathcal{O}) 
= \SDef_{{\mathcal E}}(\mathcal{O})$.
(They used different methods due to Cooper-Long-Tilmann \cite{CLT15} for preservation of convexity in \cite{PTp}.)

Another main example can be obtained by doubling a complete hyperbolic Coxeter orbifold based on a convex polytope. 
We take a double $D_T$ of the reflection orbifold based on 
a convex tetrahedron with orders all equal to $3$. This also admits a complete hyperbolic 
structure since we can take the two tetrahedra as the regular complete hyperbolic tetrahedra
and glue them by hyperbolic isometries. The end orbifolds are four $2$-spheres with three singular points of orders $3$. 
Topologically, this is a $3$-sphere with four points removed and six edges connecting them 
all given order $3$ cyclic groups as local groups.

\begin{figure}[t]
	\centering 
	\includegraphics[height=8cm]{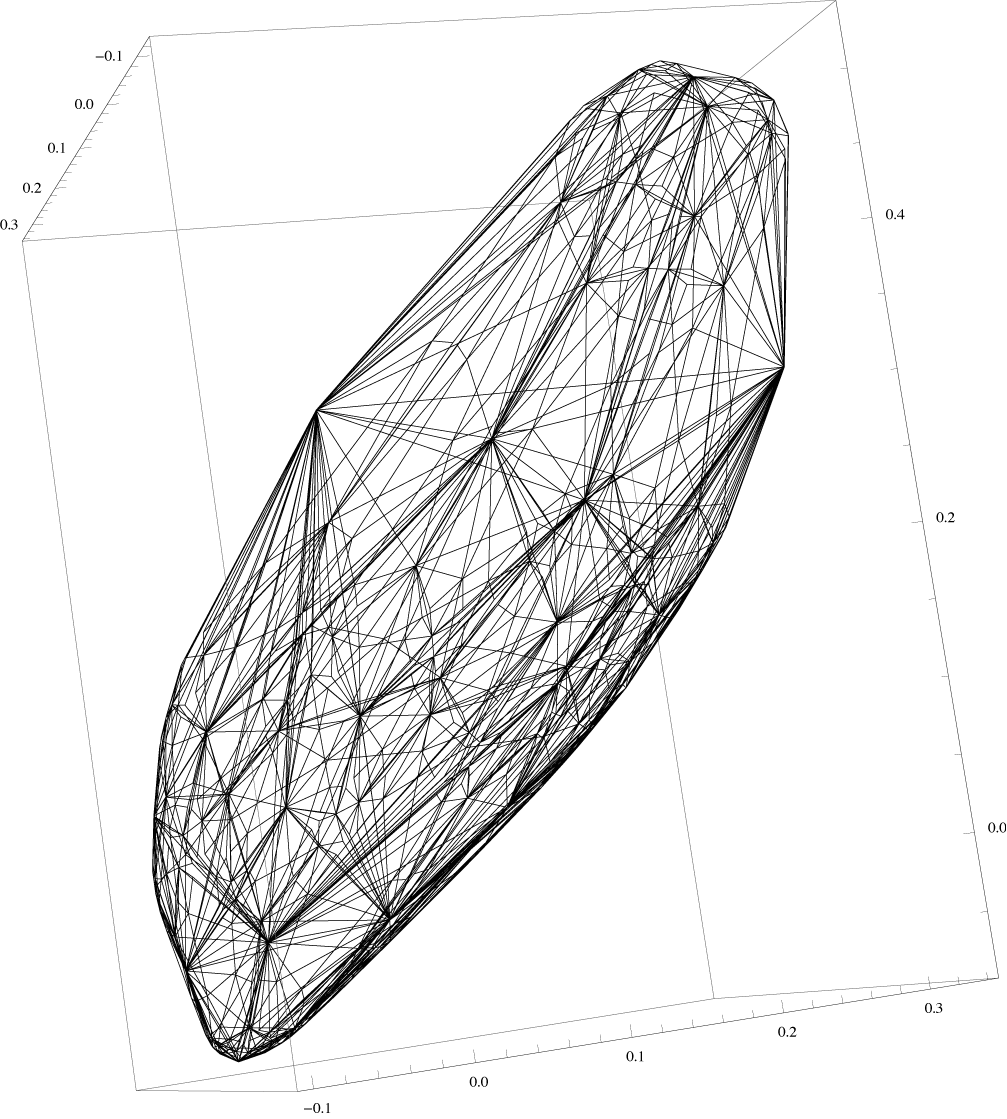}
	
	\caption{A convex developing image example of a tetrahedral orbifold of orders $3,3,3,3,3,3$.}
	\label{ex-fig-good}
\end{figure}

\begin{theorem}\label{ex-thm-tetrbounded} 
	Let $\mathcal{O}$ denote the hyperbolic $3$-orbifold $D_T$.  
	We assign the $\cR$-type to each end. 
	Then  $\SDef_{{\mathcal E}}({\mathcal{O}})$ equals $\SDef_{{\mathcal E}, \romu, \lh}(\orb)$ 
	and $\hol$ maps $\SDef_{{\mathcal E}}({\mathcal{O}})$
	as an onto-map to a component of characters 
\[\rep_{{\mathcal E}}(\pi_1({\mathcal{O}}), \PGL(4, \bR))\] containing 
	a $\PO(3, 1)$-representation which 
	is also a component of 
	\[\rep_{{\mathcal E}, \romu, \lh}(\pi_1({\mathcal{O}}), \PGL(4, \bR)).\] 
	In this case, the component is a cell of dimension $4$. 
\end{theorem}
\begin{proof} 
	A solvable subgroup of $\PO(3, 1)$ fixes a point of the boundary of 
	the Klein ball model $B$.
	Since $\pi_1(\mathcal{O})$ is not elementary, a finite-index subgroup of 
	$\pi_1(\mathcal{O})$ 
	has only trivial normal solvable subgroups. 
	The end orbifolds have zero Euler characteristics,
	and all the singularities are of order $3$. 
	For each end $E$, $\pi_1(E)$ is virtually abelian. Hence, $\pi_1(E)$ satisfies (NS). 
	
Since $\orb$ admits a complete hyperbolic structure with finite volume, 
$\pi_1(\orb)$ is relatively hyperbolic with respect to its  end fundamental groups. 
(This follows from Theorem 0.1 of Yaman \cite{Yaman04} since the group acts on
the sphere of infinity of the hyperbolic $3$-space accordingly. )
By Corollary \ref{rh-cor-remhyp}, any properly convex structures on $\orb$ with 
$\cR$- or $\cT$-ends are strictly convex. 
	By Corollary \ref{ex-cor-closed2}, 
	$\SDef_{{\mathcal E}}({\mathcal{O}})$ equals 
	$\SDef_{{\mathcal E}, u, \lh}(\orb)$.  
	Each of the ends has to be either horospherical, lens-shaped, or totally geodesic radial type. 
	Let $\partial_{\mathcal E} \orb$ denote the union of end orbifolds of $\orb$. 
	
	In \cite{jkms}, we showed that the triangulated real projective structures on the ends determined the real projective structure on $\mathcal{O}$. 
	First, there is a map $\SDef_{{\mathcal E}}({\mathcal{O}}) \ra \CDef(\partial_{\mathcal E} {\mathcal{O}})$ 
	given by sending the real projective structures on $\mathcal{O}$ 
	to the real projective structures of the ends. (Here if $\partial_{\mathcal E} {\mathcal{O}}$ has many components, then
	$\CDef(\partial_{\mathcal E} {\mathcal{O}})$ is the product space of the deformation space of all components.)
	Let $J$ be the image. 
	
	Let $\mu$ be a projective structure corresponding to an element of $\SDef_{{\mathcal E}}({\mathcal{O}})$. The universal cover $\torb$ is identified as a properly convex domain in $\SI^3$. 
	Each singular geodesic arc in $\torb$ connects one of the p-end vertices to the other. The developing image of $\torb$ is a convex open domain,
	and the developing map is a diffeomorphism. 
	The developing images of singular geodesic arcs 
	form geodesics meeting at vertices transversely. 
	There exist two convex tetrahedra $T_1$ and $T_2$ with vertices removed
from which decompose $\torb$. They are adjacent and their images under $\pi_1(\orb)$ tessellate $\torb$. 
	
	
	The end orbifold is such that if given an element of the deformation space, then the geodesic triangulation is uniquely obtained.
	Hence, there is a proper map from $\SDef_{{\mathcal E}}({\mathcal{O}})$ to the space of invariants of the triangulations as in \cite{jkms}, that is, the product space of cross-ratios and Goldman-invariant spaces.
	
	
	Now $\mathcal{O}$ is the orbifold obtained from doubling a tetrahedron with edge orders $3, 3, 3$.
	We consider an element of $\SDef_{\mathcal{E}}(\orb)$. Since it is convex, 
	we triangulate $\mathcal{O}$ into two tetrahedra, and this gives a triangulation for 
	each end orbifold diffeomorphic to $S^{(i)}_{3,3,3}$, $i=1,2,3,4$, corresponding 
	to four ends, 
	each of which gives us triangulations into two triangles.
	We can derive from the result of Goldman \cite{Goldman90} and Choi-Goldman \cite{cg} that given projective invariants 
	$\rho^{(i)}_1, \rho^{(i)}_2, \rho^{(i)}_2, \sigma^{(i)}_1, \sigma^{(i)}_2$ for each of the two triangles satisfying $\rho^{(i)}_1 \rho^{(i)}_2 \rho^{(i)}_3 = \sigma^{(i)}_1 \sigma^{(i)}_2$,
	we can determine the structure on $S^{(i)}_{3, 3, 3}$ for $i=1,2,3,4$
	 completely. 
	 
	For each $S^{(i)}_{3,3,3}$ with a convex real projective structure and divided into two geodesic triangles, 
	we compute respective invariants $\rho^{(i)}_1, \rho^{(i)}_2, \rho^{(i)}_2, \sigma^{(i)}_1, \sigma^{(i)}_2$ for one of the triangles corresponding to 
	the link of $T_1$: 
	\begin{eqnarray} 
	& s_i^{2}+s \tau +1,s_i^{2}+s_i \tau +1,   s_i^{2}+s_i \tau +1,  &\nonumber \\ 
	&   t_i \left(s_i^{2}+s_i \tau +1\right),  \frac{1}{t_i}\left(s_i^{2}+s_i \tau +1\right)  \left(s_i^{2}+s_i \tau +1\right)   & \label{ex-eqn-inv1}
	\end{eqnarray} 
	and for the other triangle corresponding to the link of $T_2$ the respective invariants are 
	\begin{eqnarray}
	& \frac{1}{s_i^{2}}(s_i^{2}+s_i \tau +1),  & \frac{1}{s_i^{2}}(s_i^{2}+s_i \tau+1), \frac{1}{s_i^{2}}(s_i^{2}+s_i \tau +1), \nonumber \\
	& \frac{t_i}{s_i^{3}}\left(s_i^{2}+s_i \tau +1\right), & \frac{1}{s_i^3 t_i}\left(s_i^2+s_i \tau +1\right) \left(s_i^2+s_i \tau +1\right) \label{cl-eqn-inv2}
	\end{eqnarray}
	where $s_i, t_i$, $i=1,2,3,4$, 
	are Goldman parameters and $\tau = 2 \cos \frac{2\pi}{3}$. 
	(See \cite{schoimath1}.)
	\index{projective invariant}
	
	Since $\partial_{\mathcal E} \orb$ is a disjoint union of four spheres with singularities
	$(3, 3, 3)$, $\CDef(\partial_{\mathcal E} {\mathcal{O}})$ is parameterized by 
	$s_i, t_i$ and hence is a cell of dimension $8$. 
	(This can be proved similarly to \cite{cgorb}.)

	The set $J$ is given by projective invariants of the $(3,3,3)$ boundary orbifolds 
	satisfying some equations. 
	Using the method of \cite{jkms} developed by the author, we obtain the equations that $J$ satisfies.  
These are 
\begin{align*}
s_i^2 + s_i \tau +1 &= s_j^2 + s_j \tau +1, i, j = 1,2,3,4 \\
\frac{1}{s_i^{2}}(s_i^{2}+s_i \tau +1) &= 
\frac{1}{s_j^{2}}(s_j^{2}+s_i \tau +1),  i, j = 1,2,3,4 \\
t_1 t_2 t_3 t_4 \prod_{i=1}^4 \left(s_i^{2}+s_i \tau +1\right) 
&= \frac{1}{t_1 t_2 t_3 t_4}\prod_{i=1}^4\left(s_i^{2}+s_i \tau +1\right)^2 \\ 
\prod_{i=1}^4 \frac{t_i}{s_i^{3}}\left(s_i^{2}+s_i \tau +1\right) &= 
\prod_{i=1}^4 \frac{1}{s_i^3 t_i}\left(s_i^2+s_i \tau +1\right)^2.
\end{align*}
The first and second lines of equations are from matching the cross ratios 
$\rho_l^{(i)}$ with $\rho_l^{(j)}$ for any pair $i, j$ 
corresponding to an edge connecting the $i$-th vertex to the $j$-vertex in $D_T$ and four faces containing the edge. 
(See (5) of \cite{jkms}). 
The third and fourth lines of the equations are from the equations that match
the products $\prod_{i=1}^4\sigma^{(i)}_1  = \prod_{i=1}^4 \sigma^{(i)}_2$ 
for $T_1$ and $T_2$, respectively.

	The equation is solvable: 
	\[s_1 = s_2 = s_3 = s_4=s, t_1 t_2 t_3 t_4 = C(s) \hbox{ for a constant } C(s)> 0 \hbox{ depending on } s.\]
	Thus, $J$ is contained in the solution subspace $C$, a $4$-dimensional cell in $\CDef(\partial_{\mathcal E} {\mathcal{O}})$.

	Conversely, given an element of $C$, we can assign invariants at each edge of the tetrahedron and Goldman $\sigma$-invariants at the vertices
if the invariants satisfy the equations.
	This is given by starting from the first convex tetrahedron and gluing one by one 
	using the projective invariants (see \cite{jkms} and \cite{schoimel}):
	Let the first one always be  the standard tetrahedron with vertices 
	\[[1,0,0,0],[0,1,0,0],[0,0,1,0], \hbox{ and } [0,0,0,1]\] and we let $T_2$ a fixed adjacent tetrahedron with vertices 
	\[[1,0,0,0],[0,1,0,0],[0,0,1,0] \hbox{ and } [2,2,2,-1].\] Then projective invariants 
	determine all other tetrahedron triangulating $\tilde{\mathcal{O}}$. 
	Given any deck transformation $\gamma$, $T_1$ and $\gamma(T_1)$ is connected by a sequence of adjacency-related tetrahedrons, and their pasting maps are wholly
determined by the projective invariants, where cross-ratios do not equal $0$. 
	Therefore, 
	as long as the projective invariants are bounded, the holonomy transformations of the generators are bounded. 
	Corollary \ref{ex-cor-closed3} shows that these correspond to elements of $\SDef_{\mathcal{E}, \mathrm{u}, \lh}({\mathcal{O}})
	= \SDef_{\mathcal{E}}({\mathcal{O}})$. 
	(This method was discussed in our talk in Melbourne, May 18, 2009 \cite{schoimel}.)
	Hence, we showed that $\SDef_{{\mathcal E}, \mathrm{u}, \lh}({\mathcal{O}})$
is parameterized by the set of solutions $C$. 
	Thus, $J= C$ since each element of $C$ gives us an element of
    $ \SDef_{\mathcal{E}}({\mathcal{O}})$. 
\hfill	\SSn  {\parfillskip0pt\par}
\end{proof}

The dimension is one higher than 
that of the deformation space of the reflection $3$-orbifold based on the tetrahedron. Thus 
we have examples not arising from reflection ones here as well.
See the Mathematica files \cite{schoimath} for a different explicit method of 
solutions. Also, see \cite{schoimath2} to see how to draw 
Figure \ref{ex-fig-good}.

We remark that the above theorem can be generalized to orders $\geq 3$ with hyperideal ends
with similar computations. 
See \cite{schoimath} for examples to modify orders and so on.





\appendix

\chapter{Projective abelian group actions on convex domains} 
\label{ch-abelian}


	We explore some theories of projective abelian group actions on convex domains, which is not necessarily properly convex. 
In Section \ref{abelian-sub-abelian}, we show that 
a free abelian group action decomposes the space into joins. 
In Section \ref{abelian-sub-convex}, we discuss 
convex projective orbifolds with free abelian holonomy groups. 
In Lemma \ref{abelian-lem-realeign},  
we show a decomposition similar 
to the Benoist decomposition 
for the dividing projective actions on 
properly convex domains. 
We also show that a parameter of 
orbits of free abelian groups geometrically 
converges to an orbit in Lemma \ref{abelian-lem-gconv}.  
In Section \ref{abelian-sub-deform}, we show that such actions always immediately deform
to the abelian group actions on properly convex domains.
In Section \ref{abelian-sub-geoconv}, we prove geometric convergences 
of convex real projective orbifolds slightly more general than those explored by Benoist. 
In Section \ref{abelian-sec-wmec}, we give some justification of 
why we are using the weak middle-eigenvalue conditions. 


%
%
%

\section{Convex real projective orbifolds} 

We explore a class of convex real projective 
orbifolds a little bit more general than the properly convex ones. 
Also see Leitner \cite{Leitner161}, \cite{Leitner162}, 
and \cite{Leitner163} for a similar work, where she explores 
representations of abelian groups; however, these 
do not act cocompactly on convex domains. 


Recall the Cartan decomposition $\SL_\pm(n, \bR) = KTK$ where 
$K = \Ort(n, \bR)$ and $T$ is the group of positive diagonal matrices.  
Note that the endomorphisms in $M_n(\bR)$ may have null spaces. 
(See Theorem 7.39 of \cite{Knapp} and \cite{Helgason}.)

A Cartan decomposition $g = k_{1, g} A_g k_{2, g}$ for $k_{1, g}, k_{2, g}\in \Ort(n, \bR)$ and a
diagonal matrix $A_g$ with a nonnegative nonincreasing set of diagonals
for each element $g$ of $M_n(\bR)$ exist since each element is
a limit of elements of $\GL(n, \bR)$ admitting a Cartan decomposition.

Each induced projective 
endomorphism $g$ for $\SI^{n-1}$ may have a nonempty subspace $V_g$ where it is not defined.
We call the projectivization of the null space the {\em undefined subspace} of $g$.
It could be an empty set. 

	Let $N_g$ be the projectivization of the null space of $A_g$. 
	Then $V_g:= k_{2, g}^{-1}N_g$ is the undefined subspace of $f$.  

Let $M_g$ be the matrix of $g$ written as $k_{1, g} \hat A_g k_{2, g}$ where $\hat A_g$ is the diagonal matrix with the maximum entry being $1$.  We call this the normalization of $g$. 

Use the Riemannian metric of $\SI^{n-1}$ to compute the norms of 
differentials. We will call these $\bdd$-norms.




\subsection{A connected free abelian group with positive eigenvalues only} 
\label{abelian-sub-abelian}

Recall that 
for a matrix $A$, we denote by $|A|$ the maximum of the norms of entries of $A$. 

We can deform the unipotent abelian representation to diagonalizable ones that are arbitrarily close to the original one. 
\begin{lemma}\label{abelian-lem-dense} 
	Let $h:\bZ^l \ra \SL_\pm(n, \bR)$ be a representation to unipotent elements. 
	Let $g_1, \dots, g_l$ denote the generators. 
	Then given $\eps> 0$
	there exists a positive diagonalizable representation $h': \bZ^l \ra \SL_\pm(n, \bR)$ with matrices 
	satisfying  $|h'(g_i) - h(g_i)| < \eps, i=1, \dots, l$.
	Furthermore, we may choose a continuous parameter of $h'_t$ 
	such that $h'_0=h$ and $h'_t$ is positive diagonalizable for $t> 0$. 
\end{lemma} 
\begin{proof} 
	First, assume that every $h(g_i)$, $i=1, \dots, l$, has matrices that are upper triangular matrices with diagonal elements equal to $1$
since the Zariski closure is in a nilpotent Lie group
and Theorem 3.5.4 of \cite{Varadarajan84}. 
	
	Let $\eps > 0$ be given. 
	We will inductively prove 
	that we can find $h'$ as above with eigenvalues of $h'(g_1)$ are all positive and mutually distinct. 
	For $n=2$, we can simply change the diagonal elements to positive numbers not equal to $1$. 
	Then the group embeds in $\Aff(\bR^1)$. We choose a positive constant $a_i$ such that $|a_i - 1| < \eps$ for each $i$. 
	Let $g_i$ be given as $x \mapsto a_i x + b_i$. 
	The commutativity reduces to equations $a_j b_i = a_i b_j$ for all $i, j$. 
	Then the solutions are given by $b_i = a_1^{-1} a_i b_1$ for any given $b_1$. 
We can construct the diagonalizable representations. 
	
	Suppose that the conclusion
	is true for dimension $k-1$. We will now consider a unipotent homomorphism $h: \bZ^l \ra \SL_\pm(k, \bR)$. We conjugate such that every $h(g_i)$ is upper-triangular. 
	Since $h(g_i)$ is upper triangular, let $h_1(g_i)$ denote the upper-left $(k-1)\times (k-1)$-matrix. 
	By the induction hypothesis, 
	we find a positive diagonalizable representation $h'_1: \bZ^l \ra \SL_\pm(k-1, \bR)$.
	Also, assume $|h'_1(g_i) - h_1(g_i)| < \eps/2$ for $i=1, \dots, l$.  
	We change 
	\[h(g_i) =  \left( 
	\begin{array}{cc}
	h_1(g_i)  & b(g_i)   \\
	0   & 1
	\end{array}
	\right) 
	\hbox{ to } 
	h'(g_i) = \left( 
	\begin{array}{cc}
	\frac{1}{\lambda'(g_i)^{\frac{1}{k-1}}} h'_1(g_i)  & b'(g_i)   \\
	0   & \lambda'(g_i)   
	\end{array}
	\right)
	\]
	for some choice of $h_1'(g), b'(g), \lambda'(g_i) > 0$ for $i, j=1, \dots, l$. 
	For commutativity, we need to solve for $b'(g_i)$
	for $i, j=1, \dots, l$, 
	\[\left(\frac{1}{\lambda'(g_i)^{\frac{1}{k-1}}}  h'_1(g_i) - \lambda'(g_j)\Idd \right) b'(g_j) = \left(\frac{1}{\lambda'(g_j)^{\frac{1}{k-1}}} h'_1(g_j) - \lambda'(g_i)\Idd \right) b'(g_i). \]
	We denote by 
 \[M_i := \left(\frac{1}{\lambda'(g_i)^{\frac{1}{k-1}}}  h'_1(g_i) - \lambda'(g_j)\Idd \right).\]
 Note $M_i M_j = M_j M_i$. 
 By generic choices  of $\lambda'(g_i)$s, we may assume that $M_i$ are invertible. 
	The solution is given by 
	\[ b'(g_i) = M(g_1)^{-1}M(g_i) b'(g_1)\]
	for an arbitrary choice of $b'(g_1)$. 
	We choose $b'(g_1)$ arbitrarily near $b(g_1)$. 
	Here, $\lambda(g_1)$ has to be chosen generically to make all the eigenvalues distinct and sufficiently near $1$ 
	such that $|h'(g_i) - h(g_i)|< \eps$, $i=1, \dots, l$.  
	We can check the solution by the commutativity. 
	Hence, we complete the induction steps. 
	
	To find a parameter denoted $h'_t$, we simply repeat the induction process	building a parameter of $h'_t$.
\end{proof}

\begin{lemma} \label{abelian-lem-Abconv}
	Let $L$ be a connected projective abelian group acting on a properly convex domain $K$
	cocompactly and effectively. Then $L$ is positive diagonalizable and the domain is a simplex.
	\end{lemma}
\begin{proof}
$L$ contains a cocompact lattice $L'$. By the Hilbert metric of $K^o$, 
$L'$ acts properly discontinuously on $K^o$. Proposition \ref{prelim-prop-Ben2} applies now. 
Since $L$ is the Zariski closure of positive diagonalizable $L'$, the proof is finished.
\end{proof}

Fix $g \in A, g \ne \Idd$. 
The minimal polynomial of $g$ is of the form 
$\prod_{i=1}^m (x-\lambda_i) $ where each nonreal $\lambda_i$ 
pairs with exactly one other $\lambda_j $ equal to $\bar \lambda_i$. 
We can write it as 
\begin{equation}\label{abelian-eqn-minpoly}
(x- \lambda_1)^{r_1} \cdots (x- \lambda_s)^{r_s}
\prod_{t=0}^{(m-s)/2-1 }(x^2 - 2\Re\lambda_{s+ 2t + 1} x - |\lambda_{s+ 2t + 1}|^2)^{r_{s+2t+1}}, 
\end{equation} 
where $\lambda_1, \dots, \lambda_s$ are all the real eigenvalues and 
$\Re z$ denotes the real part of any complex number $z$.
We define
a primary decomposition subspace $C_i$ as the kernel of 
\[M_i(g):= (g - \lambda_i \Idd )^{r_i} \hbox{ for } i=1, \dots, s\]
and $C_t$ to be the kernel of 
\[M_{s+2t+1}(g):= (g^2 - 2 \Re\lambda_{s+2t+1} x - |\lambda_{s+2t+1}|^2 \Idd)^{r_{s+2t+1}} \hbox{ for } t =0, \dots, (m-s)/2 -1. \]
(See \cite{HK71}.)


\begin{lemma}\label{abelian-lem-cyclic} 
	Let $g \in \SLpm$ be a nonidentity element. 
	\begin{itemize} 
	\item Given a primary decomposition space $C_i$ of $g$, 
	we have $h C_i = C_i$ for any $h \in \SLpm$ commuting with $g$.
	\item Given a primary decomposition subspace $C$ of $g$ and 
	$D$ of $h$ for $g, h \ne \Idd$, 
	$C\cap D$ are both $h$ and $g$ invariant
	provided $h$ and $g$ commute with each other.
	\item Given a free abelian group $A$ of finite rank, there exists a maximal collection of invariant subspaces 
	\[C_1, \dots, C_m \hbox{ satisfying } \bR^{n+1} = C_1\oplus \cdots \oplus C_m \]
	where each $C_j$ is $g$-invariant in a primary decomposition space of every $g$, $g \in A$. 
	\end{itemize} 
\end{lemma} 
\begin{proof} 
Corollary to Theorem 12 of Section 6.8 in \cite{HK71} implies the first statement. 
Let $g_1, \dots, g_k$ denote the generators of $A$. 
We obtain $C_1, \dots, C_m$ by taking a primary decomposition
space $C_{i, j}$ for $g_j$ and taking intersections of the 
arbitrary collections of $C_{i, j}$ for all $i, j$.
\end{proof} 
	


A {\em scalar group} is a group acting by $s\Idd$ for $s \in \bR$ and $s> 0$. 
A {\em scalar unipotent group} is a subgroup of the product of a scalar group with a unipotent group. 
Hence, on each $A| C_i$ is a scalar unipotent group
for each $i$. 

\begin{lemma} \label{abelian-lem-A} 
	Let $A$ be a connected free abelian group acting on $\bR^n$
	with positive eigenvalues only. 
	Then there exists a decomposition 
	$\bR^n = V_0 \oplus V_1 \oplus \cdots \oplus V_m$ 
	where $A$ acts as a positive diagonalizable group on $V_0$ 
	and acts as a positive scalar unipotent group on each $V_i$ for $1 \leq i \leq m$.  
\end{lemma} 
\begin{proof} 
We obtain $C_1, \dots, C_m$ by Lemma \ref{abelian-lem-cyclic}. 
	On $C_i$, $A$ acts as a scalar group acting 
	on a one-dimensional space or a scalar unipotent group
	since the corresponding factor of the minimal 
	polynomial is $(x - \lambda_i \Idd)^{r_i}$. 
\end{proof} 

\begin{proposition}\label{abelian-prop-cocompact}
Suppose that $\Gamma$ is a discrete free abelian group whose 
Zariski closure is $A$.
Suppose that elements of $\Gamma$ have only positive eigenvalues. 
Then $A/\Gamma$ is compact. 
\end{proposition} 
\begin{proof} 
	By Lemma \ref{abelian-lem-A}, 
	there is a decomposition 
	$\bR^n = V_0 \oplus V_1 \oplus \cdots \oplus V_m$ 
	where $\Gamma$ acts as a positive diagonalizable group on $V_0$ 
	and acts as a positive scalar unipotent group on each $V_i$ for $1 \leq i \leq m$.

	Let $\Gamma_i$ denote the restriction $\Gamma| \SI(V_i)$. 
	The restriction $A_i$ of $A|\SI(V_i)$ is a unipotent group or a diagonalizable group, and $A_i/\Gamma_i$ is compact. 
	Hence, $\Gamma$ and $A$ are subgroups of 
	$A_0\times A_1 \times \cdots \times A_m \times Q$ where 
	$Q$ is a maximal diagonalizable group acting 
	trivially on $\SI(V_i)$ for $i \geq 1$ 
	and acting diagonally on $\SI(V_0)$.	

Since the commutativity is also an algebraic condition, $A_i$ are abelian for $i=0, 1, \dots, m$ and 
must be free by the unipotence, and hence $A$ is a free abelian group.
	
	Let $q =\dim V_0$. 
    Let $x_1(g), \dots, x_q(g)$ denote 
    eigenvalues of $g$, $g\in A$, for $V_0$ and 
	\[x_{q+1}(g), \dots, x_{q+m}(g)\] denote 
	respective ones for $V_1,\dots, V_m$.
	Let $D$ be a positive diagonalizable group 
	acting as a scalar group on each $V_1, \dots, V_m$ 
	and positive diagonalizable group on $V_0$
	defined as a subgroup of $\bR_+^{q+ m}$
	given by the equation $x_1(g)\cdots x_{q+m}(g) =1$

	There is a homomorphism $c:A \ra  D$ 
	given by sending $g$ to the tuples of eigenvalues respectively 
	associated with $V_0$ and $V_i$ for $i=1,\dots, m$. 
	Then $\Gamma$ has a cocompact image in $c(A)$ under $c$
	since $c(A)$ is a connected diagonalizable group 
	that is a Zariski closure of $c(\Gamma)$ and $c$ is continuous. 
	
	Let $K$ be the kernel of $c:A \ra D$ which is an algebraic group. 
	This is a unipotent subgroup of $A$ and contains
	$K \cap \Gamma$. 
	Let $K_1$ be the Zariski closure of $K\cap \Gamma$ in $K$. 
	Since $K\cap \Gamma$ is normal in $\Gamma$, 
	and $K_1$ is the minimal algebraic group containing $K\cap \Gamma$, 
$K_1$ is normalized by $\Gamma$ and hence by $A$. 
	
	If $K_1$ is a proper subgroup of $K$, then
	there is a proper algebraic subgroup of $A$ containing 
	$\Gamma$ since $A$ is a product of $K$ and a group isomorphic to 
	$c(A)$. This is a contradiction. 
Hence, $K$ is the Zariski closure of $K \cap \Gamma$. 
	
	Since $K$ is unipotent, $K \cap \Gamma$ is cocompact in $K$ by Malcev. 
	Hence, $\Gamma$ is cocompact in $A$. 
	\end{proof}

\subsection{Convex real projective structures} \label{abelian-sub-convex}

One can think of the following lemma as a classification of convex real projective orbifolds
with abelian fundamental groups. Benoist \cite{Benoist94}, \cite{Benoist99} investigated these in a more general setting. 

\begin{lemma} \label{abelian-lem-realeign} 
	Let $\Gamma$ be a finitely generated free abelian group acting on a convex domain $\Omega$ of \,  $\SI^{n-1}$ \rlp resp. 
	$\RP^{n-1}$ \rrp\, properly and cocompactly. 
	Then the following hold:
	\begin{enumerate} 
	\item[(i)] the Zariski closure $L$ of a finite-index subgroup $\Gamma'$ of $\Gamma$ is 
	such that $L/\Gamma'$ is compact, and $L$ has only positive eigenvalues
	\rlp resp. a lift of $L$ to $\SL_{\pm}(n, \bR)$ has\rrp. 
	\item[(ii)]  $\Omega$ is an orbit of the abelian Lie group $L$ of dimension $n-1$ 
	acting properly and freely on it. 
	\item[(iii)] 
	$\Omega = (A_1\ast \dots \ast A_p \ast \{p_1\} \ast \dots \ast \{p_q\})^o$ for a complete affine subspace $A_i$, $i=1, \dots, p$, and points $p_j$, $j=1, \dots, q$. 
	Here, $\langle A_1 \rangle, \dots, \langle A_p \rangle, p_1, \dots, p_q$ are independent. 
	\item[(iv)] $L$ contains a central Lie subgroup $Q$ of rank 
	$p+q-1$ acting trivially on 
	$A_j$ and $p_k$ for $j =1, \dots, p, k=1, \dots, q$. 
	\end{enumerate} 
\end{lemma}
\begin{proof}
		We will prove for the case $\Omega \subset \SI^{n-1}$. The other case is implied by this. 
	If $\Omega$ is properly convex, then Proposition \ref{prelim-prop-Ben2} gives us a diagonal matrix group $L$ 
	acting on a simplex. (i) to (iv) follow in this case.

	(i) 
	%
	Assume that $\Omega$ is not properly convex. 
	
	Now, $\Gamma$ has no invariant lower-dimensional subspace $P$ meeting $\Omega$:
	otherwise, $\Gamma$ acts on $P\cap \Omega$ properly such that 
	$P\cap \Omega/\Gamma$ is virtually homeomorphic to 
	a lower-dimensional manifold homotopy
	equivalent to a cover of $\Omega/\Gamma$ by Theorem \ref{prelim-thm-vgood}. 
	This is a contradiction. 
	
	The positivity of the eigenvalues will be proved:
	Let $C_1, \dots, C_{p+q}$ denote the subspaces of $\Gamma$ obtained by 
	Lemma \ref{abelian-lem-cyclic}
	where $\dim C_1,..., \dim C_p \geq 2, \dim C_{p+1}= \dots = \dim C_{p+q} = 1$. 
	We also denote $C_{j} = \{p_{j-p}\}$ for $j=p+1,\dots, p+q$. 
	Let $\lambda_1(g), \dots, \lambda_{p+q}(g)$ denote the norms of 
	eigenvalues of each element $g$ of $L$ restricted to $C_1, \dots, C_{p+q}$.
	
	Define \[ \hat \SI_j:= 
	\SI(C_1)\ast \cdots \ast \SI(C_{j}) \ast \SI(C_{j+1}) \ast \cdots \ast \SI(C_{p+q}).\] 
	By the second paragraph after (i), we have
	\[ \Omega \cap \hat \SI_j  =\emp.\] 
	Let 
	\[\Pi_j: \SI^{n-1}  - \hat \SI_j \ra \SI(C_j), j=1, \dots, p+q \]
	denote the projection. Consider $\Omega_j$ denote the image of $\Omega$ under $\Pi_i$. 
	The image is a convex subset of $\SI(C_j)$
	since convex segments go to convex segments or a point under $\Pi$.
	The image is open since otherwise the dimension of $\Omega < n$. 
	Now $\Gamma$ acts cocompactly on $\Omega_j$ since $\Gamma$ acts cocompactly  
	on $\Omega$ which maps to $\Omega_j$ under $\Pi_j$. 
	Then $\Omega \subset \Omega_1\ast \cdots \ast \Omega_{p+q}$. 
	

We may prove by induction that 
	$\Omega_j$ is contained in a hemisphere 
	and is a convex open domain for each $j$: 
This is clear for $p=1$ by Proposition \ref{prelim-prop-classconv}.  
Now assume that this holds for $p=m$. 
Suppose that any $\Omega_j$ contains an antipodal pair of points. Then 
there is a pair $r, s$ in $\Omega$ that must be on a great segment passing $\hat \SI_j$ by 
our definition of the projection $\Pi_i: \hat \SI_j \ast \SI(C_j) \ra \SI(C_j)$. 
Then there is a properly convex segment 
connecting two points of $\Omega$ in the great segment meeting $\hat \SI_j$. 
Since $\Omega$ does not meet $\hat \SI_j$, this is not possible.
Hence, $\Omega_j$ is contained in a hemisphere by Proposition \ref{prelim-prop-classconv}. 
	
	%
	%
	

	We can consider the action of $\Gamma$ on $C_j$ has the norm of 
	the eigenvalue equal to $1$ only by 
multiplying by a representation $\Gamma \ra \bR_+$ and we are working on 
$\SI^{n-1}$. 
	The action of $\Gamma$ on $C_j$ is orthopotent 
	by Theorem \ref{prelim-thm-orthopotent}. 
	By Conze-Guivarc'h \cite{CG74} 
	or Moore \cite{Moore68}, there is an orthopotent flag in $\SI_j$
	and hence a proper $\Gamma$-invariant subspace. 	
	Let $\Gamma_j$ denote the image of $\Gamma$ by the restriction homomorphism to 
	$\SI(C_j)$. 
	Since $\Gamma_j$ is abelian, $\Gamma_j$ contains a uniform lattice $L_j'$
	in the Zariski closure of $\Gamma_j$. 
Since $L_j'$ is discrete, 
Theorem \ref{prelim-thm-orthopotent} shows that $L'_j$ has a virtually unipotent
action on $\SI(C_j)$  and 
so is its Zariski closure. Hence
	$\Gamma_j$ is virtually unipotent. (See Theorem 3 of Fried \cite{Fried86}.)
Let $\Gamma'_j$ be the unipotent subgroup of $\Gamma_j$ of finite-index. 
We can take $\Gamma':= \bigcap_{j=1}^{p+q} \Pi_j^{\ast -1}(\Gamma'_j) \cap \Gamma$. 
	The finite-index subgroup $\Gamma'$ of 
	$\Gamma$ has only positive eigenvalue at $\SI(C_j)$ for each $j$, $j=1, \dots, p+q$. 
	Also, the Zariski closure $Z_j$ of $\Gamma'_j$ is isomorphic to 
	$\bR^{n_j}$ for some $n_j$. 

		We assume that $\Gamma'$ is torsion-free by taking a finite-index subgroup 
by Selberg's lemma, i.e, Theorem \ref{prelim-thm-vgood}. 
		The Zariski closure $L'$ of $\Gamma$ is in $Z_1 \times \cdots \times Z_p \times Q$
for the positive diagonalizable group $Q$ acting trivially on $\SI(C_j)$ for each $j=p+1, \dots, p+q$, 
		and hence $L'$ is free abelian. 
%
$L'/\Gamma$ is a closed manifold by Proposition \ref{abelian-prop-cocompact}. 
We take a connected component $L$ of $L'$, and we redefine 
$\Gamma' = L \cap \Gamma$. Now, 
$L/\Gamma'$ is a manifold, and $\Omega/\Gamma'$ is a closed manifold. 
Since they are both $K(\Gamma', 1)$-spaces, it follows that 
$\dim L = \dim \Omega = n-1$. 
	This proves (i). 
	
	
	(ii) Now, we will let $\Gamma$ be $\Gamma'$ above without loss of generality.  
	Suppose that $p=1$ and $q=0$. 
	Since $\Omega$ is a convex domain in 
	an affine subspace in $\SI^{n-1}$, $\Omega$ is in a complete affine subspace.
	We can change $\Gamma$ to be unipotent by changing scalars.
	A unipotent group acts on a half-space in $\bR^{n}$ since
	its dual must fix a point in $\bR^{n\ast}$ being solvable.  
	Thus $\Gamma$ acts on an affine subspace $\mathds{A}^{n-1}$ in $\SI^{n-1}$, 
	and $\Gamma$ acts as an affine transformation group of $\mathds{A}^{n-1}$.  
	 Proposition T of \cite{GH86} proves our result. 
	
If $p> 1$, then $\Omega$ cannot be complete affine since 
some of $\SI(C_j)$ meets $\Omega$ otherwise. 
Hence, it must be that  $\Omega$ is not complete affine and not properly convex. 	

	There exists a great sphere 
	$\SI^{i-1}$ in the boundary of $\Omega$ where $L$ acts on 
	and is the common boundary of $i$-dimensional affine spaces foliating $\Omega$
	by Proposition \ref{prelim-prop-classconv}
	as in Section \ref{np-sub-general}.
	There is a projective projection 
	\[\Pi_{\SI^{i-1}}: \SI^{n-1} - \SI^{i-1} \ra \SI^{n-i-1}.\]  
	Then the image $\Omega_1$ of $\Omega$ is properly convex,
	and $\Omega$ is the inverse image $\Pi_{\SI^{i-1}}^{-1}(\Omega_1)$. 
	Since $L$ acts on $\Omega_1$, it follows that $L$ acts on $\Omega$. 
	Since $\dim L =\dim \Omega$, and $\Gamma$ acts on $\Omega$ properly 
    discontinuously with a compact fundamental domain, 
	$L$ acts properly and cocompactly on $\Omega$. (See Section 3.5 of \cite{Thurston97}.)
	

	Let $N$ denote the kernel of $L$ going to a
	connected Lie group $L_1$ acting on $\Omega_1$ properly and cocompactly. 
	\[ 1 \ra N \ra L \ra L_1 \ra 1.\]
	By Lemma \ref{abelian-lem-Abconv},
	$\Omega_1$ is a simplex. Hence, $L_1$ is a positive diagonalizable group. 
	Since $\Omega_1/L_1$ is compact, 
	$L_1$ acts simply transitively on $\Omega_1$ 
	by  Lemma 2.5 of \cite{Benoist03}.
	$\dim L_1 = n-i-1$. Thus, $\dim N = i$ and the abelian group $N$ acts on each complete affine $i$-dimensional 
	affine space $A_l$ that is a leaf. Since the action of $N$ is proper, $N$ acts on $A_l$ transitively by the proof of Lemma 2.5 of \cite{Benoist03}. 
	The action is simple since $\dim A_l \leq  \dim l=i$. 
	Hence, $L$ acts transitively on $\Omega$.
	Since the action of $L$ is proper, $L$ acts simply transitively 
	by the dimension count.

(iv)	By Lemma \ref{abelian-lem-A} and (i), 
	we decompose $\bR^n$ into subspaces 
	$\bR^{n}= V_1 \oplus \cdots \oplus V_p \oplus V_0$ where
	$V_j$ corresponds to $C_j$ for $j=1, \dots, p$, 
	$V_0$ corresponds to $C_{p+1}\ast \cdots \ast C_{p+q}$, 
	 $L$ acts on $V_{0}$ as a positive diagonalizable linear group and 
	$L$ acts on $V_j$ as elements of an abelian positive 
	scalar unipotent group for $j=1, \dots, p$. 
	(Here $V_{0}$ can be $0$ and $V_i$ for $i \geq 1$ equals $C_j$ for some $j$. )
	
	Since $L$ acts on $V_0$ as a positive diagonalizable group, 
	it fixes points $p_1, \dots, p_q$ in general position in 
	$\SI^{q-1}:= \SI(V)$ with $q = \dim V_0$. 	
	We claim that $L$ contains an abelian Lie subgroup $Q$ of rank $p+q-1$
	acting trivially on each $\Omega_j$, $j=1, \dots, p$, 
	and fixing $p_i$, $i=1, \dots, q$. 
	Suppose that $\Omega$ is properly convex.
	Then $\Omega$ is the interior of a simplex. 
	The cocompactness of a lattice of $L$ shows that $L$ contains 
	a discrete free abelian central group of rank $p+q-1$ by the last 
	part of Proposition \ref{prelim-prop-Ben2}. 
	The central group is contained in $Q$. Since $L$ is the Zariski closure 
	of the lattice, $Q$ is a subgroup of $L$. 
	
	Suppose that $\Omega$ is not properly convex. 
	We deform a lattice of $L$ to a diagonalizable one by a generalization of 
	Lemma \ref{abelian-lem-dense} to the direct sum of scalar unipotent 
	representations, 
	and use the limit argument. 
	
	Actually, $Q$ is the maximal diagonalizable group 
	acting on the lines in directions of 
	$p_1, \dots, p_q$ and eigenspaces $V_j$, $j=1, \dots, p$.
	(iv) again follows from the limit argument.

	
	
	
	(iii)	We choose a generic point $\llrrparen{\vec{x}}$, $\llrrparen{\vec{x}}\in \Omega$, 
		in the complement of $\SI(V_0)$, $\SI(V_j)$ for 
	$j=1, \dots, p$. Since these are independent spaces, 
	$\vec{x} = \vec{x}_0 + \sum_{i=1}^{p} \vec{x}_j$ where $\vec{x}_0\in V_0$, 
	$\vec{x}_j \in V_j$. 
	We choose a parameter of element $\eta_t$ of $Q$ fixing $V_0$ or $V_j$ for some $j$ 
	with largest norm eigenvalues and $\{\eta_t\}$ converging to $0$-maps on other subspaces
	as $t \ra \infty$. By Theorem \ref{prelim-thm-converg}, 
	we obtain a projection to $V_0$ or $V_j$ for each $j$ as a limit 
	in $\SI(M_n(\bR^n))$, and 
	$\llrrparen{\vec{x}_0}, \llrrparen{\vec{x}_j}$ are in 
	the closure of $L(\llrrparen{x})$. 
	
%

Since $\Pi_j(L(x)) = L(\Pi_j(x))$, we obtain
\begin{equation} \label{eqn:Lx}
L(\llrrparen{\vec{x}}) \subset L(\llrrparen{\vec{x}_0}) \ast L(\llrrparen{\vec{x}_1}) \ast \cdots \ast L(\llrrparen{\vec{x}_p}). 
\end{equation}
From the above paragraph, we can show that
$L(\llrrparen{\vec{x}_j})$ is contained in $\clo(L(\llrrparen{\vec{x}}))$.
Hence, 
\begin{equation} \label{eqn:Lx2}
\clo(L(\llrrparen{x})) \supset L(\llrrparen{\vec{x}_0}) \ast L(\llrrparen{\vec{x}_1}) \ast \cdots \ast L(\llrrparen{\vec{x}_p}). 
\end{equation}

Since $L$ acts transitively on $\Omega$, 
$L$ acts so on the projection $\Omega_j$ under $\Pi_j$.
Hence, $L(\llrrparen{\vec{x}}) = \Omega$ and $L(\llrrparen{\vec{x}_j}) = \Omega_j$. 
By convexity of the domain $\clo(\Omega)$, 
$(\clo(\Omega))^o = \Omega$. 
We obtain 
\[\Omega= (\{p_1\} \ast \cdots \ast \{p_q\} \ast \Omega_1 \ast \cdots \ast \Omega_p)^o. \]

Recall from above that $L|\Omega_j$ is a unipotent abelian group and hence has a distal action. The proof of Theorem 2 of \cite{Fried86} applies here since $L|\Omega_j$ contains a cocompact lattice, 
and it follows that $L(\vec{x}_j)$ is a complete affine space in $\SI(V_j)$.
This proves (iii).
\hfill \SnT  {\parfillskip0pt\par}
\end{proof} 




\begin{lemma} \label{abelian-lem-gconv} 
	Let $t_0 \in I$ for an interval $I$. 
	Suppose that we have a parameter of compact convex domains $\tri_t \subset \SI^{n-1}$ for $t < t_0$, $t \in I$, 
	and a transitive group action 
	$\Phi_t:L \times \tri^o_t \ra \tri^o_t, t \in I$ by 
	a connected free abelian group $L$ of rank $n-1$ for each $t \in I$. 
	Suppose that $\Phi_t$ depends continuously on $t$ and $\Phi_t$ is given by a continuous parameter of 
	homomorphisms $h_t: L \ra \SL_{\pm}(n, \bR)$. 
	Then $\{\tri_t\} \ra \tri_{t_0}$ geometrically where 
	$\triangle_{t_0}$ is a convex domain, and $L$ 
	acts transitively on $\triangle_{t_0}^o$. 
\end{lemma} 
\begin{proof} 
	Let $L \cong \bR^{n-1}$ have coordinates $(z_1, \dots, z_{n-1})$. 
	$\Phi_t(g, \cdot): \SI^{n-1} \ra \SI^{n-1}$ is represented 
	as a matrix 
	\begin{equation}\label{abelian-eqn-exp} 
	h_t(g) = \exp ( H_t( \sum_{i=1}^{n-1} z_i(g) e_i))
	\end{equation}
	where $\{H_t: \bR^{n-1} \ra \mathfrak{sl}(n, \bR)\}$ is 
	a uniformly bounded collection of linear maps.
	
Since $\Phi_{t_0}$ is an isomorphism, 
we may assume that some open set is contained in an orbit of $\Phi_{t_0}(L)$ for a point $x_0$. 
Since $\Phi_t \ra \Phi_{t_0}$ algebraically, 
	we may assume without loss of generality that $\bigcap_{t\in I} \tri^o_t \ne \emp$ and 
it contains an open neighborhood of $x_0$ by taking a smaller $I$ containing $t_0$.

	Let $\Delta_{t_0}$ denote the interior of a geometric limit of 
	$\clo(\Delta_{t_i})$ for some subsequence $t_i \ra t_0$ as we can see from 
 Section \ref{prelim-sub-Haus}. 
By Proposition \ref{prelim-prop-convC}, 
$\Delta_{t_0}$ is a convex open set. 

Obviously, we hace
	\[\{\Phi_t(g, x_0)\} \in \tri^o_t \ra \Phi_{t_0}(g, x_0) \hbox{ as } t \ra t_0.\]
Hence,  $\Phi_{t_0}(L)(x_0)$ is contained in $\Delta_{t_0}$ 
by \eqref{prelim-eqn-characterization}.

Now we show that $\Delta_{t_0}^{o}$ is a unique open orbit of $L$ under $\Phi_{t_0}$. 


	First suppose that $\Delta_t$ is properly convex for $t \in I - \{t_0\}$. 
	Suppose that $\Delta'_{t_0}$ contains more than two open orbits	under $\Phi_{t_0}$. 
	Then there exists a point $y$ in the interior of $\Delta'_{t_0}$ and in 
	an orbit of dimension $< n$. 
By (iv) of Lemma \ref{abelian-lem-realeign}, there exists a one parameter group 
of element of form $g^s\in L, s\in \bR$ fixing each point of 
a hyperspace $P$ passing $y$ and having a largest norm eigenvalue at another point $x$
of multiplicity one outside $P$. (To see, we just need to consider the diagonalizable group
acting trivially on each of the subspaces and find the kernel acting trivially on the join producing 
$P$.) 

Let $B(y)$ denote a compact ball which is contained in $\Delta'_{t_0}$ and in $\Delta_t$ for sufficiently close $t$ to $t_0$. 
For sufficiently close $t$ to $t_0$, $g$ still has a largest norm eigenvalue of multiplicity 
one and a $g$-invariant subspace $P_t$ meeting the interior of $B(y)$. 
Then acting by $g^s, s\in \bR$, we see that $\Delta_t$ cannot be properly convex.
This is a contradiction. 
Hence, $\Delta_{t_0}^o$ is an orbit.

	Suppose now that $\Delta_t$ is not properly convex for $t$ in $J -\{t_0\}$ 
	for some interval $J \subset I$. 
	Then $\clo(\Delta_t) = \SI^{i-1}_t\ast K_t$ for a properly convex domain $K_t$
	by Proposition \ref{prelim-prop-classconv}. 
	$L$ acts on a great sphere 
	$\SI^{i-1}_t$ in the boundary of $\Delta_t$. Now, 
	$\SI^{i-1}_t$ is the common boundary of $i$-dimensional affine spaces foliating $\Delta_t$ by Proposition \ref{prelim-prop-classconv}
	as in Section \ref{np-sub-general}.
	We may assume without loss of generality that $\SI^{i-1}_t$ is a fixed sphere
	$\SI^{i-1}$ acting by an element $g_t$ where $\{g_t\}$ converging to $\Idd$
	as $t \ra t_0$. 
		There is a projection 
		\[\Pi_{\SI^{i-1}}: \SI^{n-1} - \SI^{i-1} \ra \SI^{n-i-1}.\] 
		Then we consider $\Pi_{\SI^{i-1}}(K_t) \subset \SI^{n-i-1}$.
		Now, the discussion reduces to the above by taking a subgroup $L' \subset L$ 
		acting transitively on the interior of $\Pi_{\SI^{i-1}}(K_t)$ 
		for $t \in J -\{t_0\}$.

For each $i$-hemisphere in $\Delta_{t}^o$ for $t$ near $t_0$ with boundary $\SI^{i-1}$, 
$L$ acts transitively. We may assume that this is true for $t= t_0$ by the limit argument.
Hence, we show that $\Delta_{t_0}$ is an orbit. 
\end{proof}

\subsection{Deforming convex real projective structures}
\label{abelian-sub-deform} 

\begin{lemma} \label{abelian-lem-injectivity} 
	Let $\mu$ be a real projective structure on a closed $(n-1)$-orbifold $M$ 
	with a developing map $\dev: \tilde M \ra \SI^{n-1}$ {\rm (}resp. $ \RP^{n-1}$\/{\rm )}
	is not injective. Then for any structure $\mu'$ sufficiently close to 
	$\mu$, its developing map $\dev'$ is not injective. 
	\end{lemma} 
\begin{proof} 
	We prove for $\SI^{n-1}$. 
	We take two open sets $U_1$ and $U_2$ 
	respectively containing two points $x, y\in\tilde M$ 
	with $\dev(x) = \dev(y)$ where $\dev| U_i$ is an embedding for 
	each $i=1,2.$
	Then for any developing map $\dev'$ for $\mu'$ perturbed from $\dev$
	under the $C^r$-topology, 
	$\dev'(U_1)\cap \dev(U_2)\ne \emp$. Hence, $\dev'$ is not injective.
	\hfill \SnP {\parfillskip0pt\par}
	\end{proof} 

\begin{lemma} \label{abelian-lem-oneorbit} 
Suppose that $\mu_0$ is a convex real projective structure on a closed orbifold $M$ with 
a developing map $\dev_0: \tilde M \ra \SI^{n-1}$ {\rm (}resp. $ \RP^{n-1}$\/{\rm )}
mapping into an orbit of the identity component of 
the Zariski closure of its free abelian holonomy group given by 
the holonomy homomorphism $h_0: \pi_1(M) \ra \SLpm$. 
Let $\mu_i$ be a sequence of real projective structures on $M$ converging to $\mu_0$, let
$\dev_i$ be one of the corresponding 
developing maps converging to $\dev_0$, and let $h_i$ be one of 
be the corresponding holonomy homorphisms converging to $h_0$. 
Suppose that a sequence of the orbits  of the identity component 
of Zariski closure of the holonomy group of $\dev_i$ geometrically converges to that of $h_0$. 
Then $\dev_i$ must map into 
an orbit of the identity component of 
the Zariski closure of the corresponding holonomy group of $\dev_i$. 
\end{lemma}
\begin{proof} 
We prove the result for $\SI^{n-1}$. 
Now $\dev_t(F)$ geometrically converges to $\dev_0(F)$ by Proposition \ref{prelim-prop-BP}. 
Also, the orbit of the identity component of the Zariski closure geometrically converges.  
For the fundamental group $F$ of $\tilde M$, 
$\dev_0(F)$ and hence $\dev_i(F)$ must map into the orbit for sufficiently large $i$. 
Hence, $\dev_i$ must develop into an orbit. 
\hfill \SnP {\parfillskip0pt\par}
\end{proof}

A convex real projective structure $\mu_0$ on an orbifold $\Sigma$ is {\em virtually immediately deformable to} a properly convex real structure
if there exists a parameter $\mu_t$ of real projective structures on a finite cover $\hat \Sigma$ 
of $\Sigma$ such that $\hat \Sigma$ with induced structures $\hat \mu_t$ is 
properly convex for $t > 0$. 
\index{virtually immediately deformable} 

\begin{proposition}\label{abelian-prop-vip} 
	A convex real projective structure on a closed $(n-1)$-orbifold $M$ 
	with virtually free abelian holonomy subgroup of  a finite index is always virtually  immediately deformable to	a properly convex real projective structure. 
\end{proposition} 
\begin{proof} 
	Again, we prove the result for $\SI^{n-1}$. 
	Let $\bZ^l$ denote the fundamental group of a finite cover $M'$ of $M$. 
	Let $h \in \Hom(\bZ^l, \SL_\pm(n, \bR))$ be the restriction of the holonomy homomorphism to $\bZ^l$. 
	Every $h$ has a neighborhood containing a positive
	diagonalizable holonomy $h': \bZ^l \ra \SL_{\pm}(n, \bR)$ by Lemmas \ref{abelian-lem-dense} and \ref{abelian-lem-realeign}.
We may assume that $h'$ has a continuous path to $h$. 
	By the deformation theory of \cite{dgorb}, $h'$ is realized as a holonomy of 
	a real projective manifold $M'$ diffeomorphic to $M'$
provided $h'$ is sufficiently close to $h$.  

In addition, the universal cover of $M'$ is a union of orbits of an abelian Lie group $L$
by Benoist \cite{Benoist94}.  
	Here, $h'(\bZ^l)$ is a lattice in $L$ by Lemma \ref{abelian-lem-realeign}. 
	
	By Lemma \ref{abelian-lem-brickn}, 
	$M'$ is a properly convex real projective orbifold.
%
\hfill \SSn {\parfillskip0pt\par}
\end{proof} 

Let us recall a work of Benoist: 
Let $M$ be a real projective $(n-1)$-orbifold with nilpotent holonomy. 
(Here we are not working with orbifolds.)
Let $N$ be the nilpotent identity component of 
the Zariski closure of the holonomy group. 
Benoist \cite{Benoist94}, \cite{Benoist99} showed that $\tilde M$ decomposes into 
a union of connected open submanifolds $D_i$, $i\in I$ for an index set $I$,  such that 
$\dev:D_i \ra \dev(D_i)$ is a diffeomorphism to an orbit of 
a nilpotent Lie group $N$. 
The brick number of $M$ is the number of 
the $(n-1)$-dimensional open orbits that map to mutually distinct connected open strata in 
$M$.  



\begin{lemma} \label{abelian-lem-brickn} 
	Let $M$ be a closed $(n-1)$-orbifold. 
	Let $(M, \mu)$ be a real convex projective orbifold with a virtually abelian fundamental group. Suppose $\mu$ is deformed into a continuous parameter $\mu_t$ of real projective structures such that $\mu_0 = \mu$. 
	Let $h_t: \pi_1(M) \ra \SL_{\pm}(n, \bR)$ {\em (}resp. $\PGL(n, \bR)$\/{\rm )}, 
	$t\in [0, 1]=I$, 
	be a continuous family of the associated holonomy homomorphism
	with $\dev_t: \tilde M \ra \SI^{n-1}$  {\em (}resp. $\RP^n$\/{\rm )}. 
	Suppose $\mu_t$ has the holonomy group $h_t(\pi_1(M))$
	with following properties for $t> 0$\/{\rm :} 
	\begin{description} 
		\item[{\rm (}A{\rm )}] it is virtually diagonalizable for 
		or, more generally, it acts on some properly convex domain $D_t$, or 
		\item[{\rm (}B{\rm )}] it acts properly on a complete affine space $D_t$
	\end{description} 
where $D_t$ has no proper $h_t(\pi_1(M))$-invariant open domain.

	Then $(M, \mu_t)$ is properly convex or is complete affine for $t > 0$. 
\end{lemma} 
\begin{proof} 
	Again, we prove for $\SI^{n-1}$. 
	For $t =0$, $\dev_0$ is a diffeomorphism into a complete affine space 
	or a properly convex domain by our conditions. 
	
	
	
	Define the following sets: 	
	\begin{itemize} 
		\item $A$ is the subset of $t$ satisfying (A)  and $\mu_t$ is a properly convex structure, and 
		\item $B$ is the subset of $t$ for (B). 
		
	\end{itemize} 
	We will decompose $[0, 1] = A \cup B$. 
	
	By Lemma \ref{abelian-lem-gconv}, the set $\hat A$ of points of $I$ satisfying (A) is open. 
	By Koszul \cite{Koszul70}, the set $A$ is open and $A \subset \hat A$. 
	We have $\hat A \cup B =[0, 1]$. 
	
	(A) Suppose that we have an open connected 
	subset $U$ with $U \subset \hat A$ and $t_0 \in \clo(U)$ with properly convex or complete affine $\mu_{t_0}$. 
	Then we claim $U \subset A$. 
	
	Let $t \in U$. 
	Since the holonomy group $h_t(\pi_1(M))$ 
	is virtually abelian, a finite cover $M'$ of $M$ is octantizable by Proposition 2 of \cite{Benoist94}. 
	$\dev_t$ maps to a union of orbits of 
	a connected abelian Zariski closure $\Delta_t$ of a finite-index abelian subgroup $H_t$ 
	of $h(\pi_1(M'))$ in $\SI^{n-1}$ by \cite{Benoist99}. 
	Benoist also shows that $\Delta_t$ acts on $\tilde M$.
	
	By our assumption for $h_t$, $h_t(\pi_1(M'))$ 
	acts on a properly convex domain in $\SI^{n-1}$.
	$\Delta_t$ is positive diagonalizable by Lemma \ref{abelian-lem-Abconv}. 
	Now orbits of $\Delta_t$ are convex simplexes in $\SI^{n-1}$
	by Section 3.1 of Benoist \cite{Benoist99}
	where he explains the classification of such structures
	by Smillie \cite{Smillie} and \cite{SmillieTh}.
	
Now, $\dev_t$ cannot map to a more than one orbit of $\Delta_t$ for sufficiently small 
$|t - t_0|$ by Lemma \ref{abelian-lem-oneorbit}. 
By Theorem \ref{prelim-thm-Kobayashi}, $\dev_t$ is a diffeomorphism to an orbit of $\Delta_t$. 
	(In other words, a sequence of real projective structures with more than one bricks 
	cannot converge to a properly convex structure or 
	a complete affine  structure.) 
		Since the orbits of $\Delta_t$ on $\SI^{n-1}$ are properly convex, 
		we conclude that $\mu_t \in A$ for a sufficiently small 
		$|t-t_0|$ by Theorem \ref{prelim-thm-Kobayashi}. 
		Hence, $U \cap A$ is a nonempty open set.

	Also, we claim that $U$ is in $A$: for any sequence $\{t_i\}$ converging to $t_0'$ 
	in $U$ in $\hat A$, choose a point $x_0$ 
	in the developing image of $\dev_{t_i}$ for sufficiently large $i$
	since $\mu_{t_i}$ are sufficiently $C^r$-close. 
	Then $\Delta_{t_i}(x_0) \ra \Delta_{t_0'}(x_0)$ by Lemma \ref{abelian-lem-gconv}. 
	Since $\tilde M$ with $\mu_{t_i}$ develops into $\Delta_{t_i}(x_0)$,
the fundamental domain of 
	$\tilde M$ with $\mu_{t_0'}$ develops into $\Delta_{t_0'}(x_0)$,
and hence $\dev_{t_0'}$ develops into $\Delta_{t_0'}(x_0)$.  
	The developing map is a diffeomorphism to $\Delta_{t_0'}(x_0)$ 
	by Theorem \ref{prelim-thm-Kobayashi}. 
	Hence, $t_0' \in A$ also. 
	Thus, $U\cap A$ is open and closed and hence $U\subset A$. 
	
	Hence, for connected open $U$, $U \subset \hat A$, with $t_0 \in \clo(U)$ for 
	properly convex or complete affine $\mu_{t_0}$, we have 
	$U \subset A$. 

	
(B)	Let $I'$ denote the subset of $[0, 1]$ consisting 
	of $t$ with $\mu_t$ satisfying conclusions of the lemma. 
		Since $\mu_0$ is properly convex or complete affine, 
	$I'$ is not empty. 
	Also, $I'$ contains all components of $\hat A$ meeting it
	by the above argument. Also, $A' = I' \cap \hat A$ is open. 
	
		We will show that $I'$ is open and closed in $[0, 1]$ and hence 
		$I'=[0,1]$: 
	
	First, we show that $I'$ is open. 
	Also, for $t_0 \in I'\cap B$, we claim 
	that a neighborhood is in $I'$:
	Otherwise, there is a sequence $\{t_i\}$ for $t_i \not\in I'$ converging to $t_0$
	with $t_i \in I$ and $\dev_{t_i}$ is not a diffeomorphism 
	to a complete affine space where $h_{t_i}$ acts on. 
Without the loss of generality, we may assume the following: 
\begin{itemize} 
\item 	If $t_i \in \hat A$ for infinitely many $i$, 
	then $t_i \in A\subset I'$ for sufficiently large $i$ 
by Lemma \ref{abelian-lem-oneorbit} 
and Theorem \ref{prelim-thm-Kobayashi}. 
	This is a contradiction. 
\item 	Assume now $t_i \in B$ for sufficiently large $i$. 
Then again $\dev_{t_i}$ maps into an orbit of the Zariski closures of the image of
 $h_{t_i}$. Since $\dev_{t_i}$ are diffeomorphisms to a open hemisphere that is an orbit,
it follows that so is $\dev_{t_0}$: The injectivity follows from Lemma \ref{abelian-lem-injectivity}
and the surjectivity follows from the condition that there is no holonomy 
invariant open subset.

	
	\end{itemize} 
	Hence, we showed that $I'$ is an open subset of $I$.

	Now, we show that $I'$ is closed: 
	For any boundary point $t_0$ of $I'$ in $\hat A$, 
	we have $t_0\in I'$ since $I'$ contains a component of $\hat A$ it meets.  
	
	For any boundary point $t_0$ of $I'$ that is in $A$ or $B$, 
	we have $\dev_{t_0}$ must be injective by Lemma \ref{abelian-lem-injectivity} since $\dev_t$ is injective for $t< t_0$ or $t> t_0$ 
	for $t$ in an open interval with boundary point $t_0$. 
	Since the image of $\dev_{t}$ for $t\in I'$ is in an open hemisphere or 
	a compact properly convex domain $D_i$, 
$\dev_t$ must be a diffeomorphism to the interior since there are no holonomy invariant 
open subset. 
	The image of $\dev_{t_0}$ is in an open hemisphere or a properly convex 
	domain that is the geometric limit of $D_i$ 
	up to a choice of subsequences. 
	Hence, $t_0 \in I'$. 
	

	Thus, $I'$ is open and closed. This completes the proof. 
	\hfill \SnT {\parfillskip0pt\par}
\end{proof}

\begin{lemma} \label{abelian-lem-Ben3} 
	Let $\Sigma$ be a closed $(n-1)$-orbifold with a properly convex real projective 
structure $\mu_t, t\in [0, 1]$. 
	Let $\dev_t$ be a continuous parameter of developing maps for $\mu_t$
to $\SI^{n-1}$. 
	Then for $t$ in a open subset of $[0, 1]$, 
	\[\clo(\dev_t(\tilde \Sigma)) = S_{1, t} \ast \cdots \ast S_{m_t, t}\]
	for properly convex domains $S_{1, t}, \dots, S_{m_t,t}$ 
	where each $S_{j, t}$ span a subspace $P_{j, t}$. 
	The finite-index subgroup of $h_t(\pi_1(\Sigma))$ acting on $P_{j, t}$ acts strongly irreducibly for each $t$. 
	Furthermore, $m_t, \dim S_{j, t}$ are constant, and $P_{j, t}, j=1, \dots, m_t$ 
	are always independent, and $P_{j, t}$ forms a continuous parameter in
	the Grassmannian spaces $G(n, \dim P_{j, t})$ up to reordering. 
	\end{lemma} 
\begin{proof} 
For an open susbet $O$ of $[0, 1]$, $\dev_t$, $t\in O$, is a diffeomorphism by Proposition 
\ref{app-prop-koszul}.  
The decomposition follows from Proposition \ref{prelim-prop-Ben2}. 
There is a virtual center $Z$, a free abelian group of rank $m_t-1$, 
mapping to a positive diagonalizable group $Z_t$ acting trivially on each $S_{j, t}$.  	
	By Theorem 1.1 of \cite{Benoist03}, an infinite-order virtually central element cannot have 
	nontrivial action on $S_{j, t}$ since otherwise the Zariski closure 
	of the subgroup of $h_t(\pi_1(\Sigma))$ acting on $P_{j, t}$ cannot be simple as claimed immediately after 
	that theorem. Hence, any infinite-order virtually central elements are in 
	a maximal free abelian group of rank $m_t$. 
	Hence, for each $t$, there are subspaces $S_{j, t}$ for $j=1, \dots, m_t$ such that 
	the above decomposition holds. Now, we need to show that the dimensions are constant.
	
	We decompose $I$ into mutually disjoint subsets  where 
	\[
	\dim S_{j, t} = n_{j}^{(t)} \hbox{ where } n_{1}^{(t)}+ \cdots + n_{m_t}^{(t)} + m_t^{(t)} = n+1
\]
	by reordering the indices. 
	Then each of these sets is closed as we can see from a sequence argument since
	the above rank argument shows that there cannot be further decomposition in
	the limit increasing the rank of the virtual center.
	Since $I$ is connected, there is only one such set equal to $I$. 
	Now, the conclusion follows up to reordering. 
	\end{proof}

We generalize Proposition \ref{abelian-lem-brickn}: 
\begin{corollary} \label{abelian-cor-convhol} 
	Suppose that a real projective orbifold $\Sigma$  is a closed $(n-1)$-orbifold 
	with the structure $\mu$.
	 Let $\mu_t$, $t\in [0, 1]$, be a parameter of projective structures on 
	$\Sigma$ such that $\mu_0$ is properly convex or complete affine 
	and $\mu_1 = \mu$. 
	Let $h_t$ denote the associated holonomy homomorphisms. 
	\begin{itemize} 
	\item Suppose that the holonomy group 
	$h_t(\pi_1(\Sigma))$ in $\PGL(n, \bR)$ {\rm (}resp. $\SL_{\pm}(n, \bR)$\,{\rm )}
	acts on a properly convex domain or a complete affine subspace $D_t$.
		\item Suppose $D_t$ is a minimal $h_t(\pi_1(\Sigma))$-invariant convex 
		open domain. 
	\item We require $\pi_1(\Sigma)$ to be virtually abelian if
	$D_{t_0}$ is complete affine for at least one $t_0$. 
	\end{itemize} 
	Then $\Sigma$ is also properly convex or complete affine
	where the following hold\/{\em :}
	\begin{itemize} 
	\item a developing map 
	$\dev_t$  is a diffeomorphism to $R_t(D_t)$ for every $t\in [0, 1]$
	where $R_t$ is a uniformly bounded projective automorphism for each $t$ 
	and is a composition of reflections commuting with one another. 
	\item Furthermore, 
	$\Sigma$ is properly convex if so is $D_t$.
	\end{itemize} 
\end{corollary}
\begin{proof} 
	Again, we prove in $\SI^{n-1}$. 
	If $\pi_1(\Sigma)$ is virtually abelian, then it follows from Proposition
	\ref{abelian-lem-brickn}  
	
	Now, suppose that $\pi_1(\Sigma)$ is not virtually abelian. 
	Then $D_{t}$ is properly convex for every $t$, $t \in I$, by the premise.

	 The set $A$ where $\mu_t$ is properly convex is open by 
	 Koszul \cite{Koszul68}. Let $t_0$ be the supremum of 
	 the connected component $A'$ of $A$ containing $0$.
	 There is a developing map $\dev_t: \tilde \Sigma \ra \SI^{n-1}$ for $t \in A$
	 is a diffeomorphism to a properly convex domain $D'_t$
	 where $D'_t = R_t(D_t)$ for a projective automorphism $R_t$ by
	 Lemma \ref{prelim-lem-domainIn}.
	 

	 Since the associated developing map $\dev_t$ maps into $D'_t$ for $t < t_0$, 
	 $\dev_t| K$ for a compact fundamental domain $F\subset \tilde \Sigma$ maps into 
	 a compact subset of $D_{t}^{\prime o}$ for $t < t_0$. 
	 Since $\dev_t:\tilde \Sigma$ is injective, 
	 $\dev_{t_0}: \tilde \Sigma \ra \SI^{n-1}$ 
	 is also injective by Lemma \ref{abelian-lem-injectivity}.
	  By the injectivity and the invariance of domain, 
	  $\dev_{t_0}$ is a diffeomorphism to an open domain $\Omega$. 
	  Since every point of the image of $\dev_{t_0}$ is 
	  approximated by the points in the image $D_t^{\prime o}$
	  of $\dev_t$ for $t < t_0$.  
      Hence, $\clo(\Omega)$ is contained 
	  in the geometric limit of a convergent subsequence of any sequence 
	  $\clo(D'_{t_i})$ by Proposition \ref{prelim-prop-BP}. 
	  By Lemma \ref{prelim-lem-domainIn}, 
	  $\Omega$ is a properly convex domain since 
	  the holonomy group acts on a properly convex domain $D_{t_0}$
	  and $\Omega = R_{t_0}(D_{t_0})$ for a projective automorphism $R_{t_0}$
	  that is a composition of reflections commuting with one another.

	 Hence, $A$ is also closed, and $A =[0, 1]$. 
%
%
	 By Lemma \ref{abelian-lem-Ben3}  the uniform boundedness of $R_t$ follows
	 since the subspaces $P_{j, t}$ are continuous and $R_t$ are either $\Idd$ or 
	 $\mathcal{A}$ on it. \hfill \SnT {\parfillskip0pt\par}
	\end{proof}

\subsection{Geometric convergence of convex real projective orbifolds} 
\label{abelian-sub-geoconv}


Note that the third item of the premise below
is automatically true by Theorem \ref{ce-thm-mainaffine} 
if $\Sigma$ is an end-orbifold of a properly convex affine $n$-orbifold
for any $t$. 

\begin{corollary}\label{abelian-cor-smvar} 
	Suppose that $\Sigma$ is a closed $(n-1)$-orbifold.
	We are given a path $\mu_t$, $t \in [0, 1]$, 
	of $\RP^{n-1}$-structures on $\Sigma$ equipped with the $C^r$-topology, 
	$r\geq 2$. \index{crtopology@$C^r$-topology} 
	Suppose that $\mu_0$ is properly convex or complete affine.  
	\begin{itemize} 
	\item Suppose that the holonomy group 
	$h_t(\pi_1(\Sigma))$ in $\PGL(n, \bR)$ {\rm (}resp. $\SL_{\pm}(n, \bR)$\,{\rm )}
	acts on a properly convex domain or a complete affine subspace $D_t$. 
	\item Suppose that $D_t$ is the minimal holonomy invariant domain. 
	\item We require that if $\mu_t$ is complete affine for at least one $t$, then 
	the holonomy group is virtually abelian.
	\end{itemize} 
	Then the following holds\/{\rm :} 
	\begin{itemize} 
		\item We can find a  family of developing maps $\dev_t$ to $\RP^{n-1}$ \/  
		{\rm (}resp. in $\SI^{n-1}$\/{\rm )} continuous in the $C^r$-topology
		and a continuous family of holonomy homomorphisms $h_t: \Gamma \ra \Gamma_t$ 
		such that $K_t := \clo(\dev_t(\tilde \Sigma))$ is a continuous family of convex domains in $\RP^{n-1}$ {\rm (}resp. in $\SI^{n-1}$\/{\rm )} under 
		the Hausdorff metric topology of the space of closed subsets of $\RP^{n-1}$ {\rm (}resp. $\SI^{n-1}$\/{\rm ).}
		\item In other words, 
		given $0< \eps < 1/2$ and $t_0, t_1 \in [0, 1]$, we can find $\delta > 0$ such that if $|t_0 -t _1| < \delta$, then 
		$K_{t_1} \subset N_\eps(K_{t_0})$ and $K_{t_0} \subset N_\eps(K_{t_1})$.
		\item Also, given $0 < \eps < 1/2$ and $t_0, t_1 \in [0,1]$, we can find $\delta >0$
		such that if $|t_0 -t _1| < \delta$, then 
		$\partial K_{t_1} \subset N_\eps(\partial K_{t_0})$ and $\partial K_{t_0} \subset N_\eps(\partial K_{t_1})$ where we assume  $\dev_t$ maps to $\SI^{n-1}$ exclusively here. 
		\item Finally, $\mu_t$ is always virtually immediately deformable to a properly convex structure. 
	\end{itemize} 
\end{corollary}   
\begin{proof} 
	We will prove for the result for 
    $\SI^{n-1}$. 
	Suppose first that $\pi_1(\Sigma)$ is not virtually abelian. 
	Then $D_t$ is never a complete affine space and hence 
	is always properly convex by the premise. 
By Corollary \ref{abelian-cor-convhol},  
$\dev_t$ is an embedding to the interior of properly convex $K_t$ for each $t$. 
	
%
%
%
%
%
%
	
	First, for any sequence $\{t_i\}$ converging to $t_0$, 
	we can choose a subsequence $\{t_{i_j}\}$ so 
	that $\{K_{t_{i_j}}\}$ converges to a compact convex set $K_\infty$ in the Hausdorff metric.
	$h_{t_0}(\pi_1(\Sigma))$
	 acts on $K_\infty$ by Corollary \ref{prelim-cor-uniqueness}. 
%
	 
	By Lemma \ref{abelian-lem-Ben3}, we define 
	$K_{l, t} := P_{l, t}\cap K_t$ where $h_t(\pi_1(\Sigma))$ acts 
	strongly irreducibly on $P_{l, t}$, and 
	$K_t = K_{1, t}\ast \cdots \ast K_{m_1}$.  
	We obtain $K_{l, t_{i_j}} \ra K_{l}$ as $j \ra \infty$ for a compact convex 
	set $K_l$ in a subspace $P_{l, t_0}$ where a finite-index subgroup of 
	$h_{t_0}(\pi_1(\Sigma))$ acts on
	by Lemma \ref{abelian-lem-Ben3} and Proposition \ref{prelim-prop-BP}. 
	Proposition \ref{prelim-prop-Ben2} shows that 
	the action on $P_{l, t_0}$ by $\pi_1(\Sigma')$ is irreducible. 
	Hence, $K_l$ must be properly convex by Proposition \ref{prelim-prop-projconv}.
	We have $K_\infty \subset K_1 \ast \cdots \ast K_m$ by Proposition \ref{prelim-prop-BP} since 
	$K_{t_{i_j}} \subset K_{1, t_{i_j}}\ast \cdots \ast K_{m, t_{i_j}}$ for each $j$.
	Hence, $K_\infty$ is properly convex.  
	 
	Now, for any sequence $\{t'_i\}$ converging to $t_0$, 
	suppose that a convergent subsequence $\{K_{t'_{i_j}}\}$ 
	converges to $K'_\infty$. Then we claim that 
	$K_\infty = K'_\infty$: 
	Now, $K'_\infty$ is properly convex also. 
	Choose a torsion-free finite-index subgroup $\Gamma'$ of $h_{t_0}(\pi_1(\Sigma))$
	by Theorem \ref{prelim-thm-vgood}. 
	$K_\infty^o/\Gamma' $ and $K_\infty^{\prime o}/\Gamma'$
	are homotopy equivalent. 
	Since $\dev_{t_i}$ and $\dev_{t'_i}$ are close, 
	we may assume that $K_\infty^o \cap K_\infty^{\prime o} \ne \emp$. 
	Lemma \ref{prelim-lem-domainIn} shows that  
	$K_\infty^o = K_\infty^{\prime o}$.
	This implies the first item 
	for this case. 
	
	
	
	

	Suppose now that $\Gamma$ is virtually abelian. 
	Then $\Omega_t$ is determined by the generators of the free abelian subgroup $\Gamma'$ of 
	a finite index with only positive eigenvalues by Lemma \ref{abelian-lem-realeign}. 
	$\Gamma'$ determines the connected abelian Lie group $\Delta_t$ containing $h_t(\Gamma')$ 
	and $\Omega_t$ is an orbit of $\Delta_t$ by Lemma \ref{abelian-lem-realeign}.  Now Lemma \ref{abelian-lem-gconv} implies the first item. 

	The second item follows from the first one. The third one can be deduced by  Theorem \ref{prelim-thm-partialKi}. 
	The fourth item follows from Proposition \ref{abelian-prop-vip}. 
%
\hfill	\SnT {\parfillskip0pt\par}
\end{proof}

\begin{remark} \label{ablian-rem-finalquestion} 
	Of course, we wish to generalize Lemma \ref{abelian-lem-brickn} and 
	Corollaries \ref{abelian-cor-convhol} and \ref{abelian-cor-smvar} 
	 for fully general cases without the restriction on 
	 the domains of actions starting from any properly convex 
	projective orbifold and show the similar results. 
Then we can allow NPNC-ends into the discussions. 
	We leave this as a question of whether 
	one can achieve such results.
	\end{remark} 






\section{ The justification for weak middle-eigenvalue conditions} 
\label{abelian-sec-wmec}

\begin{theorem} \label{abelian-thm-justify}
	Let $\orb$ be a properly convex real projective $n$-orbifold
	with ends. 
	Let $\Sigma_{\tilde E}$ be an end orbifold of an R-p-end $\tilde E$ of $\orb$
	with the virtually abelian end-fundamental group $\pi_1(\tilde E)$. 
	\begin{itemize} 
	\item Suppose that $\mu_i$ be a sequence of properly convex 
	structure on an R-p-end neighborhood $U_{\tilde E}$ of $\tilde E$ 
	corresponding to a generalized lens-shaped
	R-p-ends satisfying 
	the uniform middle-eigenvalue conditions. 
	\item Suppose that $\mu_i$ limits to $\mu_\infty$ 
	in the $C^{r}$-topology, $r \geq 2$. 
	\end{itemize} 
	Then $\mu_\infty$ satisfies the weak middle-eigenvalue 
	condition for $\tilde E$. 
	Furthermore, the holonomy group for $\mu_\infty$ virtually 
	satisfies the transverse weak middle-eigenvalue condition for $\tilde E$
	if it is NPNC and $\pi_1(\tilde E)$ is virtually abelian. 
	\end{theorem} 
\begin{proof} 
	Assume first that $\torb \subset \SI^n$. 
	We may assume that the p-end vertex $\mbv_{\tilde E}$ is independent of $\mu_i$
	by conjugation of the holonomy homomorphisms. 
	Let $h_i: \pi_1(\Sigma_{\tilde E}) \ra \SLpm_{\mbv_{\tilde E}}$.
	Since these satisfy the uniform middle-eigenvalue conditions, 
	we have 
	\[ C^{-1}\cwl(g) < \log \left( \frac{\lambda_1(h_i(g))}{\lambda_{\mbv_{\tilde E}}(h_i(g))} \right)  < C \cwl(g), C > 1, g \in \pi_1(\Sigma_{\tilde E}) \]
	where $C$ is a constant which may depend on $h_i$. 
	Let $h_\infty$ denote the holonomy homomorphism for $\mu_\infty$,
	which is an algebraic limit of $h_i$. 
	By taking limits, we obtain that $h_\infty$ satisfies the weak middle
	eigenvalue condition. 

	Suppose now that $\mu_\infty$ is NPNC.	
	For convenience, we may assume without loss of generality 
	that $\pi_1(\tilde E)$ is free abelian. 
	Let $\Gamma$ denote $h_\infty(\pi_1(\tilde E))$. 
	

	By Lemma \ref{abelian-lem-realeign}, 
	$\tilde \Sigma_E$ is the interior of a strict join of 
	hemispheres and a properly convex domain
	\[H_1 \ast \cdots \ast H_m \ast K_0 \subset \SI^{n-1}_{\mbv_{\tilde E}}\]
	where 
	\begin{itemize} 
	\item $\Gamma|H_j$, $j=1, \dots, m$, has the Zariski closure 
	a unipotent Lie group for a finite-index subgroup $\Gamma$ of 
	$h_\infty(\pi_1(\tilde E))$, 
	\item $\Gamma|K_0$ is a diagonalizable group acting so, 
	and 
	\item $\Gamma$ has a center $Q$ of rank $m+ \dim K_0 - 1$
	acting trivially on each $H_i$, $i=1, \dots, m$ 
	and fixing the vertices of $K_0$. 
	
	\end{itemize} 
	
	Given any $i$-dimensional hemisphere $V$ of $\SI^{n-1}_{\mbv_{\tilde E}}$ 
	for $0 \leq i \leq n-1$, 
	there exists unique $(i+1)$-dimensional hemisphere $\hat V$ in $\SI^n$ 
	in the direction of $V$ from $\mbv_{\tilde E}$ and 
	containing $\mbv_{\tilde E}$ in $\partial \hat V$ 


We denote by $\hat H_i$ the hemispheres in $\SI^n$ corresponding to 
the directions of $H_i$ for 
$i=1, \dots, m$ and $\hat p_i$ the great segments in $\SI^n$ corresponding 
to the directions of vertices $p_1, \dots, p_k$ of $K_0$. 
Let $g \in h(\pi_1(\tilde E))$.
Since a cocompactly-acted lens for $h_i(\pi_1(\tilde E))$ 
contains the points affiliated with the largest norm of eigenvalues for 
$h_j(g)$ for each $g \in \pi_1(\tilde E)$, a limiting argument shows that 
points in $\hat H_i^o$ or $\hat p_i$ must be affiliated with 
the largest norm $\lambda_1(h(g))$ of the eigenvalues. 
(Of course, these are not all such points necessarily)

By Proposition \ref{prelim-prop-classconv}, 
the maximal dimensional great sphere $\SI^{i_0-1}_\infty$ in 
$\Bd \tilde \Sigma_{\tilde E}\subset \SI^{n-1}_{\mbv_{\tilde E}}$ 
corresponding the boundary of complete affine leaves 
in $\tilde \Sigma_{\tilde E}$ equals 
 $\partial \hat H_1 \ast \cdots \ast \partial \hat H_m$. 
Since these points are not in the directions of
 $\partial \hat H_1 \ast \cdots \ast \partial \hat H_m$, 
	the desired inequality
	$\lambda^{Tr}_{\mbv_{\tilde E}}(g) \geq \lambda_{\mbv_{\tilde E}}(g)$
	holds. 
%
%
%
	%
\hfill	\SnT  {\parfillskip0pt\par}
	\end{proof}


%
%




%
%


\chapter*{Index of Notations}

These are not the definitions. Please see the pages 
to find precise definitions. 

\begin{description} 
\item[$\bR$] The real number field.

\item[$\bR_+$ ] The set of positive real numbers.

\item[$\bC$ ] The complex number field.

\item[$\mathcal{A}$] The antipodal map $\SI^n \ra \SI^n$. 

\item[$\bdd$] The Fubini-Study metric on $\SI^n$. p.\pageref{bdd} 

\item[$\Hom(G, H)$] The set of homomorphisms $G \ra H$ for two groups $G, H$. 

\item[$\rep(G, H)$] The set of conjugacy classes of homomorphisms $G \ra H$ for two groups $G, H$. 

\item[$| \cdot |$] The maximal norm of the entries of a matrix. 
p.\pageref{twoverts}

\item[$\llrrV{\cdot}$] The Euclidean metric on a vector space over $\bR$, 
(also we use $\llrrV{\cdot}_E$ for emphasis).   
p.\pageref{twovert2}

\item[$\llrrV{\cdot}_{\mathrm{fiber}}$] A fiberwise metric on a vector bundle over 
an orbifold.  p.\pageref{vertfiber} 

\item[$\pi_1(\cdot)$] The orbifold fundamental group of an orbifold. 

\item[${\mathcal{Z}(\cdot)}$] The Zariski closure of a group. p.\pageref{zariski} 
	
\item[$\bZ(\cdot)$] The center of a group. p.\pageref{center} 

\item[$\Aut(K)$] The group of projective automorphisms of a set $K$ in $\rpn$ or 
$\SI^n$. p.\pageref{Aut} 

\item[$\RP^n$] The $n$-dimensional projective space. p.\pageref{Pparanthesis}

\item[$\RP^{n\ast}$] The dual $n$-dimensional projective space. p.\pageref{dualrpn}
 
\item[$\mathds{A}^n$] The $n$-dimensional affine space. p.\pageref{affine space}

\item[$\SI^n$] The sphere.  p.\pageref{sphere} 

\item[$p_{\SI^n}$] The double covering map 
$\SI^n \ra \RP^n$.  p. \pageref{sphere} 

\item[$\SI^{n\ast}$] The dual sphere. p.\pageref{Snast} 

\item[$\Pi$] $\bR^{n+1}-\{O\} \ra \rp^n$ projection. p.\pageref{pi} 

\item[$\Pi'$] $\bR^{n+1} - \{O\} \ra \SI^n$ projection. p.\pageref{piprime}

\item[$\Bd^\Ag \Omega$] 
$
\{ (x, H)| x \in \Bd \Omega, x \in H,  
H \hbox{ is an oriented sharply supporting hyperspace of } \Omega \}  
\subset \SI^{n}\times \SI^{n\ast}.
\hbox{ p. \pageref{PiAg}} 
$

\item[$\Pi_\Ag$] projection $\Pi_\Ag: \Bd^\Ag \Omega \ra \Bd \Omega$ given 
by $(x, H) \mapsto x$. p. \pageref{PiAg}

\item[$\partial M$] manifold or orbifold boundary of a manifold or orbifold $M$. p.\pageref{manifoldboundary}

\item[$\Bd X$] topological boundary of $X$ in an ambient space.
p.\pageref{manifoldboundary}

\item[$\Bd_X Y$] topological boundary of $Y$ in an ambient space $X$.
p.\pageref{manifoldboundary}

\item[$M^o$] the manifold or orbifold interior of a manifold or orbifold $M$ or 
the relative interior of a convex domain in a spanning projective or affine subspace.
 p.\pageref{manifoldboundary}


\item[$\bP(V)$] the projectivization of a vector space $V$. p. \pageref{Pparanthesis}

\item[$\SI(V)$] the sphericalization of a vector space $V$. p. \pageref{Sparanthesis} 

\item[$\overline{pq}$] the geodesic segment in $\rpn$ or $\SI^n$ 
connecting $p$ and $q$ not antipodal to $p$.
p.\pageref{overline} 

\item[$\overline{pzq}$] the geodesic segment in $\SI^n$ connecting $p$ and $q$ antipodal to $p$ and passing $z$. p.\pageref{overline} 

\item[$()^\ast$] Duals of vector spaces or convex sets or linear groups. 
p.\pageref{dualrpn}

\item[$()^\dagger$] The intrinsic subspace dual of a properly convex domain 
in a subspace. 
p.\pageref{dagger}

\end{description}


\backmatter

\bibliographystyle{plain} 
\bibliography{choibib}

\begin{thebibliography}{100}

\bibitem{Ablondi25}
A.~Ablondi.
\newblock Affine deformations of divisible convex cones and affine spacetimes.
\newblock arXiv:2506.12172.

\bibitem{ALR07}
A.~Adem, J.~Leida, and Y.~Ruan.
\newblock {\em Orbifolds and stringy topology}, volume 171 of {\em Cambridge
  Tracts in Mathematics}.
\newblock Cambridge University Press, Cambridge, 2007.

\bibitem{AC96}
J.~Anderson and R.~Canary.
\newblock Algebraic limits of {K}leinian groups which rearrange the pages of a
  book.
\newblock {\em Invent. Math.}, 126(2):205--214, 1996.

\bibitem{Andreev70}
E.~Andreev.
\newblock Convex polyhedra in {L}oba\v cevski\u\i\ spaces.
\newblock {\em Mat. Sb. (N.S.)}, 81 (123):445--478, 1970.

\bibitem{Ballas14}
S.~Ballas.
\newblock Deformations of noncompact projective manifolds.
\newblock {\em Algebr. Geom. Topol.}, 14(5):2595--2625, 2014.

\bibitem{Ballas15}
S.~Ballas.
\newblock Finite volume properly convex deformations of the figure-eight knot.
\newblock {\em Geom. Dedicata}, 178:49--73, 2015.

\bibitem{BCLp}
S.~Ballas, D.~Cooper, and A.~Leitner.
\newblock Generalized cusps in real projective manifolds: classification.
\newblock {\em J. Topol.}, 13(4):1455--1496, 2020.

\bibitem{BCL22}
S.~Ballas, D.~Cooper, and A.~Leitner.
\newblock The moduli space of marked generalized cusps in real projective
  manifolds.
\newblock {\em Conform. Geom. Dyn.}, 26:111--164, 2022.

\bibitem{BDL}
S.~Ballas, J.~Danciger, and G.-S. Lee.
\newblock Convex projective structures on nonhyperbolic three-manifolds.
\newblock {\em Geom. Topol.}, 22(3):1593--1646, 2018.

\bibitem{BM20}
S.~Ballas and L.~Marquis.
\newblock Properly convex bending of hyperbolic manifolds.
\newblock {\em Groups Geom. Dyn.}, 14(2):653--688, 2020.
\newblock arXiv:1609.03046.

\bibitem{BaoBonahon}
X.~Bao and F.~Bonahon.
\newblock Hyperideal polyhedra in hyperbolic 3-space.
\newblock {\em Bull. Soc. Math. France}, 130(3):457--491, 2002.

\bibitem{Barbot00}
T.~Barbot.
\newblock Vari\'et\'es affines radiales de dimension 3.
\newblock {\em Bull. Soc. Math. France}, 128(3):347--389, 2000.

\bibitem{Barbot2005}
T.~Barbot.
\newblock Globally hyperbolic flat space-times.
\newblock {\em J. Geom. Phys.}, 53(2):123--165, 2005.

\bibitem{Barden}
D.~Barden.
\newblock {\em The structure of manifolds}.
\newblock PhD thesis, Cambridge University, Cambridge, England, 1963.

\bibitem{Batyrev}
V.V. Batryrev.
\newblock Toric varieties and smooth convex approximations of a polytope.
\newblock {\em RIMS Kokyuroku}, 776:20--20, 1992.

\bibitem{Ben-Tal}
A.~Ben-Tal and M.~Teboulle.
\newblock A smoothing technique for nondifferentiable optimization problems.
\newblock In {\em Optimization ({V}aretz, 1988)}, volume 1405 of {\em Lecture
  Notes in Math.}, pages 1--11. Springer, Berlin, 1989.

\bibitem{BP92}
R.~Benedetti and C.~Petronio.
\newblock {\em Lectures on hyperbolic geometry}.
\newblock Universitext. Springer-Verlag, Berlin, 1992.

\bibitem{Benoist94}
Y.~Benoist.
\newblock Nilvari\'et\'es projectives.
\newblock {\em Comment. Math. Helv.}, 69(3):447--473, 1994.

\bibitem{Benoist97}
Y.~Benoist.
\newblock Propri\'et\'es asymptotiques des groupes lin\'eaires.
\newblock {\em Geom. Funct. Anal.}, 7(1):1--47, 1997.

\bibitem{Benoist00}
Y.~Benoist.
\newblock Automorphismes des c\^ones convexes.
\newblock {\em Invent. Math.}, 141(1):149--193, 2000.

\bibitem{Benoist99}
Y.~Benoist.
\newblock Tores affines.
\newblock In {\em Crystallographic groups and their generalizations
  ({K}ortrijk, 1999)}, volume 262 of {\em Contemp. Math.}, pages 1--37. Amer.
  Math. Soc., Providence, RI, 2000.

\bibitem{Benoist01}
Y.~Benoist.
\newblock Convexes divisibles.
\newblock {\em C. R. Acad. Sci. Paris S\'er. I Math.}, 332(5):387--390, 2001.

\bibitem{Benoist03}
Y.~Benoist.
\newblock Convexes divisibles. {II}.
\newblock {\em Duke Math. J.}, 120(1):97--120, 2003.

\bibitem{Benoist04}
Y.~Benoist.
\newblock Convexes divisibles. {I}.
\newblock In {\em Algebraic groups and arithmetic}, pages 339--374. Tata Inst.
  Fund. Res., Mumbai, 2004.

\bibitem{Benoist05}
Y.~Benoist.
\newblock Convexes divisibles. {III}.
\newblock {\em Ann. Sci. \'Ecole Norm. Sup. (4)}, 38(5):793--832, 2005.

\bibitem{Benoist062}
Y.~Benoist.
\newblock Convexes divisibles. {IV}. {S}tructure du bord en dimension 3.
\newblock {\em Invent. Math.}, 164(2):249--278, 2006.

\bibitem{Benzecri60}
J.-P. Benz{\'e}cri.
\newblock Sur les vari\'et\'es localement affines et localement projectives.
\newblock {\em Bull. Soc. Math. France}, 88:229--332, 1960.

\bibitem{Berger}
M.~Berger.
\newblock {\em Geometry {I}}.
\newblock Universitext. Springer-Verlag, Berlin, 2009.
\newblock Translated from the 1977 French original by M. Cole and S. Levy,
  Fourth printing of the 1987 English translation [MR0882541].

\bibitem{Bobb}
M.~Bobb.
\newblock Convex projective manifolds with a cusp of any non‐diagonalizable
  type.
\newblock {\em J. London Math. Soc.}, 100(1):183--202, 2019.

\bibitem{BFp2024}
M.~Bobb and J.~Farre.
\newblock Affine laminations and coaffine representations.
\newblock arXiv. 2404.14284.

\bibitem{Bonahon}
F.~Bonahon.
\newblock The geometry of {T}eichm\"{u}ller space via geodesic currents.
\newblock {\em Invent. Math.}, 92(1):139--162, 1988.

\bibitem{Borzellino12}
Joseph~E. Borzellino and Victor Brunsden.
\newblock Elementary orbifold differential topology.
\newblock {\em Topology Appl.}, 159(17):3583--3589, 2012.

\bibitem{Borzellino15}
Joseph~E. Borzellino and Victor Brunsden.
\newblock On the notions of suborbifold and orbifold embedding.
\newblock {\em Algebr. Geom. Topol.}, 15(5):2789--2803, 2015.

\bibitem{Bowditch98}
B.~Bowditch.
\newblock A topological characterisation of hyperbolic groups.
\newblock {\em J. Amer. Math. Soc.}, 11(3):643--667, 1998.

\bibitem{BH99}
M.~Bridson and A.~Haefliger.
\newblock {\em Metric spaces of non-positive curvature}, volume 319 of {\em
  Grundlehren der Mathematischen Wissenschaften [Fundamental Principles of
  Mathematical Sciences]}.
\newblock Springer-Verlag, Berlin, 1999.

\bibitem{Bryant}
R.~Bryant.
\newblock An answer given at ``{L}eft invariant metrics on {$\mathrm{SL}(n,
  \mathbb{R})$}''.
\newblock
  \url{https://mathoverflow.net/questions/108280/left-invariant-metric-on-rm-sl-n-mathbbr},
  [Accessed April 6, 2024].

\bibitem{CEG06}
R.~Canary, D.~Epstein, and P.~Green.
\newblock Notes on notes of {T}hurston.
\newblock In {\em Fundamentals of hyperbolic geometry: selected expositions},
  volume 328 of {\em London Math. Soc. Lecture Note Ser.}, pages 1--115.
  Cambridge Univ. Press, Cambridge, 2006.
\newblock With a new foreword by Canary.

\bibitem{Carriere84}
Y.~Carri\`ere.
\newblock Flots {R}iemanniens.
\newblock {\em Ast\'{e}risque}, (116):31--52, 1984.
\newblock Transversal structure of foliations (Toulouse, 1982).

\bibitem{Carriere88}
Y.~Carri{\`e}re.
\newblock Feuilletages {R}iemanniens \`a croissance polyn\^omiale.
\newblock {\em Comment. Math. Helv.}, 63(1):1--20, 1988.

\bibitem{ChCh93}
Y.~Chae, S.~Choi, and C.~Park.
\newblock Real projective manifolds developing into an affine space.
\newblock {\em Internat. J. Math.}, 4(2):179--191, 1993.

\bibitem{ChengYau}
S.~Cheng and S.-T. Yau.
\newblock Complete affine hypersurfaces. {I}. {T}he completeness of affine
  metrics.
\newblock {\em Comm. Pure Appl. Math.}, 39(6):839--866, 1986.

\bibitem{End2}
S.~Choi.
\newblock The classification of radial or totally geodesic ends of real
  projective orbifolds. {II}: properly convex ends and totally geodesic ends.
\newblock arXiv:1501.00352.

\bibitem{convMa}
S.~Choi.
\newblock The convex real projective orbifolds with radial or totally geodesic
  ends: The closedness and openness of deformations.
\newblock arXiv:1011.1060.v3.

\bibitem{schoimath1}
S.~Choi.
\newblock The mathematica file to compute real projective structures on
  ${S}^2(3,3,3)$.
\newblock \url{http://mathsci.kaist.ac.kr/~schoi/triangle4reversal.nb}.

\bibitem{schoimath}
S.~Choi.
\newblock The mathematica files to compute tetrahedral real projective
  orbifolds.
\newblock \url{http://mathsci.kaist.ac.kr/~schoi/exampletet3.nb},
  \url{http://mathsci.kaist.ac.kr/~schoi/equations.nb}.

\bibitem{quasijoin}
S.~Choi.
\newblock The mathematica files to drawing the end neighorhood of the
  quasi-joined end.
\newblock \url{http://mathsci.kaist.ac.kr/~schoi/quasijoin2.nb}.

\bibitem{schoimath2}
S.~Choi.
\newblock The mathematica files to drawing the tetrahedral real projective
  orbifolds.
\newblock \url{http://mathsci.kaist.ac.kr/~schoi/coefficientsol8-inpaper.nb}.

\bibitem{schoimel}
S.~Choi.
\newblock Real projective structures on 3-orbifolds and projective invariants.
\newblock the presentation for Melbourne Talk, May 18, 2009.
  \url{http://mathsci.kaist.ac.kr/~schoi/melb2009handout.pdf}.

\bibitem{cdcr1}
S.~Choi.
\newblock Convex decompositions of real projective surfaces. {I}.
  {$\pi$}-annuli and convexity.
\newblock {\em J. Differential Geom.}, 40(1):165--208, 1994.

\bibitem{cdcr2}
S.~Choi.
\newblock Convex decompositions of real projective surfaces. {II}. {A}dmissible
  decompositions.
\newblock {\em J. Differential Geom.}, 40(2):239--283, 1994.

\bibitem{marg}
S.~Choi.
\newblock The {M}argulis lemma and the thick and thin decomposition for convex
  real projective surfaces.
\newblock {\em Adv. Math.}, 122(1):150--191, 1996.

\bibitem{psconv}
S.~Choi.
\newblock {\em The convex and concave decomposition of manifolds with real
  projective structures}, volume~78 of {\em M\'em. Soc. Math. Fr. {\rm
  (}N.S.{\rm )}}.
\newblock Soc. Math. France, Paris, 1999.

\bibitem{rdsv}
S.~Choi.
\newblock The decomposition and classification of radiant affine 3-manifolds.
\newblock {\em Mem. Amer. Math. Soc.}, 154(730):viii+122, 2001.

\bibitem{jkms}
S.~Choi.
\newblock Geometric structures on low-dimensional manifolds.
\newblock {\em J. Korean Math. Soc.}, 40(2):319--340, 2003.

\bibitem{dgorb}
S.~Choi.
\newblock Geometric structures on orbifolds and holonomy representations.
\newblock {\em Geom. Dedicata}, 104:161--199, 2004.

\bibitem{Choi06}
S.~Choi.
\newblock The deformation spaces of projective structures on 3-dimensional
  {C}oxeter orbifolds.
\newblock {\em Geom. Dedicata}, 119:69--90, 2006.

\bibitem{Cbook}
S.~Choi.
\newblock {\em Geometric structures on 2-orbifolds\/{\em :} exploration of
  discrete symmetry}, volume~27 of {\em MSJ Memoirs}.
\newblock Mathematical Society of Japan, Tokyo, 2012.

\bibitem{End1}
S.~Choi.
\newblock A classification of radial or totally geodesic ends of real
  projective orbifolds. {I}: a survey of results.
\newblock In {\em Hyperbolic Geometry and Geometric Group Theory}, Advanced
  Studies in Pure Mathematics, pages 69--134, Tokyo, 2017. Mathematical Society
  of Japan.

\bibitem{convsurv}
S.~Choi.
\newblock The convex real projective orbifolds with radial or totally geodesic
  ends: a survey of some partial results.
\newblock In A.~Reid A.~Basmanjian, Y.~Minsky, editor, {\em In the Tradition of
  Ahlfors–Bers, VII}, volume 696 of {\em Contemp. Math.}, pages 51-- 86,
  Providence, RI., 2017. Amer. Math. Soc.

\bibitem{CG93}
S.~Choi and W.~Goldman.
\newblock Convex real projective structures on closed surfaces are closed.
\newblock {\em Proc. Amer. Math. Soc.}, 118(2):657--661, 1993.

\bibitem{cg}
S.~Choi and W.~Goldman.
\newblock The classification of real projective structures on compact surfaces.
\newblock {\em Bull. Amer. Math. Soc. (N.S.)}, 34(2):161--171, 1997.

\bibitem{cgorb}
S.~Choi and W.~Goldman.
\newblock The deformation spaces of convex {$\mathbb{RP}^2$}-structures on
  2-orbifolds.
\newblock {\em Amer. J. Math.}, 127(5):1019--1102, 2005.

\bibitem{CGLMp}
S.~Choi, R.~Green, G.~Lee, and L.~Marquis.
\newblock Projective deformations of hyperbolic coxeter 3-orbifolds of finite
  volume.
\newblock in preparation.

\bibitem{CHL12}
S.~Choi, C.~Hodgson, and G.~Lee.
\newblock Projective deformations of hyperbolic {C}oxeter 3-orbifolds.
\newblock {\em Geom. Dedicata}, 159:125--167, 2012.

\bibitem{Choi12}
S.~Choi, C.~Hodgson, and G.-S. Lee.
\newblock Projective deformations of hyperbolic {C}oxeter 3-orbifolds.
\newblock {\em Geom. Dedicata}, 159:125--167, 2012.

\bibitem{CL15}
S.~Choi and G.~Lee.
\newblock Projective deformations of weakly orderable hyperbolic {C}oxeter
  orbifolds.
\newblock {\em Geom. Topol.}, 19(4):1777--1828, 2015.

\bibitem{CLM18}
S.~Choi, G.~Lee, and L.~Marquis.
\newblock Convex real projective structures on manifolds and orbifolds.
\newblock In S.-T.~Yau L.~Ji, A.~Papadopoulos, editor, {\em The Handbook of
  Group Actions III}, chapter~10, pages 263--310. Higher Education Press and
  International Press, Boston, 2018.
\newblock arXiv:1605.02548.

\bibitem{CG74}
J.-P. Conze and Y.~Guivarc'h.
\newblock Remarques sur la distalit\'e dans les espaces vectoriels.
\newblock {\em C. R. Acad. Sci. Paris S\'er. A}, 278:1083--1086, 1974.

\bibitem{Cooper17}
D.~Cooper.
\newblock The {H}eisenberg group acts on a strictly convex domain.
\newblock {\em Conform. Geom. Dyn.}, 21:101--104, 2017.

\bibitem{CLT06}
D.~Cooper, D.~Long, and M.~Thistlethwaite.
\newblock Computing varieties of representations of hyperbolic 3-manifolds into
  {${\rm SL}(4,\mathbb{R})$}.
\newblock {\em Experiment. Math.}, 15(3):291--305, 2006.

\bibitem{CLT07}
D.~Cooper, D.~Long, and M.~Thistlethwaite.
\newblock Flexing closed hyperbolic manifolds.
\newblock {\em Geom. Topol.}, 11:2413--2440, 2007.

\bibitem{CLT18}
D.~Cooper, D.~Long, and S.~Tillman.
\newblock Deforming convex projective manifolds.
\newblock {\em Geom. Topol.}, 22:1349--1404, 2018.
\newblock arXiv:1511.06206.

\bibitem{CLT15}
D.~Cooper, D.~Long, and S.~Tillmann.
\newblock On convex projective manifolds and cusps.
\newblock {\em Adv. Math.}, 277:181--251, 2015.

\bibitem{CM14}
M.~Crampon and L.~Marquis.
\newblock Finitude g\'eom\'etrique en g\'eom\'etrie de {H}ilbert.
\newblock {\em Ann. Inst. Fourier (Grenoble)}, 64(6):2299--2377, 2014.

\bibitem{DGK21}
J.~Danciger, F.~Gu\'eritaud, and F.~Kassel.
\newblock Convex cocompact acts on real projective geometry.
\newblock {\em Ann. Sci. \'Ec. Norm. Sup\'er.}, 57:1753--1843, 2024.

\bibitem{DGKLM}
J.~Danciger, F.~Gu\'eritaud, F.~Kassel, G.~Lee, and L.~Marquis.
\newblock Convex cocompactness for {C}oxeter groups.
\newblock {\em J. Eur. Math. Soc.}, 27(1):119--181, 2025.

\bibitem{GV58}
J.~de~Groot and H.~de~Vries.
\newblock Convex sets in projective space.
\newblock {\em Compositio Math.}, 13:113--118, 1958.

\bibitem{Harpe93}
P.~de~la Harpe.
\newblock On {H}ilbert's metric for simplices.
\newblock In {\em Geometric group theory, {V}ol.\ 1 ({S}ussex, 1991)}, volume
  181 of {\em London Math. Soc. Lecture Note Ser.}, pages 97--119. Cambridge
  Univ. Press, Cambridge, 1993.

\bibitem{deRham}
G.~de~Rham.
\newblock {\em Differentiable manifolds}, volume 266 of {\em Grundlehren der
  Mathematischen Wissenschaften [Fundamental Principles of Mathematical
  Sciences]}.
\newblock Springer-Verlag, Berlin, 1984.
\newblock Forms, currents, harmonic forms, Translated from the French by F. R.
  Smith, With an introduction by S. S. Chern.

\bibitem{DK18}
C.~Dru\c{t}u and M.~Kapovich.
\newblock {\em Geometric group theory}, volume~63 of {\em American Mathematical
  Society Colloquium Publications}.
\newblock American Mathematical Society, Providence, RI, 2018.
\newblock With an appendix by Bogdan Nica.

\bibitem{Drutu09}
C.~Dru{\c{t}}u.
\newblock Relatively hyperbolic groups: geometry and quasi-isometric
  invariance.
\newblock {\em Comment. Math. Helv.}, 84(3):503--546, 2009.

\bibitem{DS08}
C.~Dru{\c{t}}u and M.~Sapir.
\newblock Tree-graded spaces and asymptotic cones of groups.
\newblock {\em Topology}, 44(5):959--1058, 2005.
\newblock With an appendix by Denis Osin and Mark Sapir.

\bibitem{Efremovic}
V.~Efremovi\v{c}.
\newblock On the proximity geometry of {R}iemannian manifolds.
\newblock {\em Uspekhi Math Nauk}, 8(189), 1953.

\bibitem{Epstein84}
D.~B.~A. Epstein.
\newblock Transversely hyperbolic {$1$}-dimensional foliations.
\newblock {\em Ast\'{e}risque}, (116):53--69, 1984.
\newblock Transversal structure of foliations (Toulouse, 1982).

\bibitem{FM}
B.~Farb and D.~Margalit.
\newblock {\em A primer on mapping class groups}, volume~49 of {\em Princeton
  Mathematical Series}.
\newblock Princeton University Press, Princeton, NJ, 2012.

\bibitem{Fried82}
D.~Fried.
\newblock The geometry of cross sections to flows.
\newblock {\em Topology}, 21(4):353--371, 1982.

\bibitem{Fried86}
D.~Fried.
\newblock Distality, completeness, and affine structures.
\newblock {\em J. Differential Geom.}, 24(3):265--273, 1986.

\bibitem{Friedlens}
D.~Fried.
\newblock Choi's lenses.
\newblock informal note, 2012.

\bibitem{FG83}
D.~Fried and W.~Goldman.
\newblock Three-dimensional affine crystallographic groups.
\newblock {\em Adv. in Math.}, 47(1):1--49, 1983.

\bibitem{GR70}
H.~Garland and M.~S. Raghunathan.
\newblock Fundamental domains for lattices in ({R}-)rank {$1$} semisimple {L}ie
  groups.
\newblock {\em Ann. of Math. (2)}, 92:279--326, 1970.

\bibitem{Ghomi}
M.~Ghomi.
\newblock Optimal smoothing for convex polytopes.
\newblock {\em Bull. London Math. Soc.}, 36(4):483--492, 2004.

\bibitem{Gthesis}
W.~Goldman.
\newblock Affine manifolds and projective geometry on surfaces.
\newblock senior thesis, Princeton University, 1977.

\bibitem{Goldman88}
W.~Goldman.
\newblock Geometric structures on manifolds and varieties of representations.
\newblock In {\em Geometry of group representations ({B}oulder, {CO}, 1987)},
  volume~74 of {\em Contemp. Math.}, pages 169--198. Amer. Math. Soc.,
  Providence, RI, 1988.

\bibitem{Goldman90}
W.~Goldman.
\newblock Convex real projective structures on compact surfaces.
\newblock {\em J. Differential Geom.}, 31(3):791--845, 1990.

\bibitem{goldmanbook}
W.~Goldman.
\newblock {\em Geometric structures on manifolds}, volume 227 of {\em Graduate
  Studies in Mathematics}.
\newblock American Mathematical Society, 2022.

\bibitem{GH86}
W.~Goldman and M.~Hirsch.
\newblock Affine manifolds and orbits of algebraic groups.
\newblock {\em Trans. Amer. Math. Soc.}, 295(1):175--198, 1986.

\bibitem{GLM09}
W.~Goldman, F.~Labourie, and G.~Margulis.
\newblock Proper affine actions and geodesic flows of hyperbolic surfaces.
\newblock {\em Ann. of Math. (2)}, 170(3):1051--1083, 2009.

\bibitem{GM88}
W.~Goldman and J.~Millson.
\newblock The deformation theory of representations of fundamental groups of
  compact {K}\"ahler manifolds.
\newblock {\em Inst. Hautes \'Etudes Sci. Publ. Math.}, (67):43--96, 1988.

\bibitem{greene}
R.~Greene.
\newblock {\em The deformation theory of discrete reflection groups and
  projective structures}.
\newblock PhD thesis, Ohio State University, 2013.

\bibitem{Gromov81}
M.~Gromov.
\newblock Groups of polynomial growth and expanding maps.
\newblock {\em Inst. Hautes \'Etudes Sci. Publ. Math.}, (53):53--73, 1981.

\bibitem{Gromov86}
M.~Gromov.
\newblock {\em Partial diﬀerential relations}, volume~9 of {\em Ergebnisse
  der Mathematik und ihrer Grenzgebiete (3)}.
\newblock Springer-Verlag, Berlin, 1986.

\bibitem{Guichard05}
O.~Guichard.
\newblock Sur la r\'egularit\'e {H}\"older des convexes divisibles.
\newblock {\em Ergodic Theory Dynam. Systems}, 25(6):1857--1880, 2005.

\bibitem{GW12}
O.~Guichard and A.~Wienhard.
\newblock Anosov representations: domains of discontinuity and applications.
\newblock {\em Invent. Math.}, 190(2):357--438, 2012.

\bibitem{HHMP10}
D.~Heard, C.~Hodgson, B.~Martelli, and C.~Petronio.
\newblock Hyperbolic graphs of small complexity.
\newblock {\em Experiment. Math.}, 19(2):211--236, 2010.

\bibitem{HK78}
E.~Heintze and H.~Karcher.
\newblock A general comparison theorem with applications to volume estimates
  for submanifolds.
\newblock {\em Ann. Sci. \'Ec. Norm. Sup\'er.}, 11(4):451--470, 1978.

\bibitem{Helgason}
S.~Helgason.
\newblock {\em Differential geometry, {L}ie groups, and symmetric spaces},
  volume~34 of {\em Graduate Studies in Mathematics}.
\newblock American Mathematical Society, Providence, RI, 2001.
\newblock Corrected reprint of the 1978 original.

\bibitem{HP11}
M.~Heusener and J.~Porti.
\newblock Infinitesimal projective rigidity under {D}ehn filling.
\newblock {\em Geom. Topol.}, 15(4):2017--2071, 2011.

\bibitem{Hirsch59}
M.~W. Hirsch.
\newblock Immersions of manifolds.
\newblock {\em Trans. Amer. Math. Soc.}, 93:242--276, 1959.

\bibitem{Hirsch}
M.~W. Hirsch.
\newblock {\em Differential topology}, volume~33 of {\em Graduate Texts in
  Mathematics}.
\newblock Springer-Verlag, New York, 1994.
\newblock Corrected reprint of the 1976 original.

\bibitem{HK71}
K.~Hoffman and R.~Kunze.
\newblock {\em Linear algebra}.
\newblock Second edition. Prentice-Hall, Inc., Englewood Cliffs, N.J., 1971.

\bibitem{Humphreys}
J.~Humphreys.
\newblock {\em Linear algebraic groups}.
\newblock Springer-Verlag, New York-Heidelberg, 1975.
\newblock Graduate Texts in Mathematics, No. 21.

\bibitem{IZp}
M.~Islam and A.~Zimmer.
\newblock A flat torus theorem for convex co-compact actions of projective
  linear groups.
\newblock {\em J. London. Math. Soc.}, 103:470--489, 2021.

\bibitem{IZp2}
M.~Islam and A.~Zimmer.
\newblock Convex co-compact actions of relatively hyperbolic groups.
\newblock {\em Geometry and Topology}, 27:417--511, 2023.

\bibitem{Jenkins}
J.~Jenkins.
\newblock Growth of connected locally compact groups.
\newblock {\em J. Functional Analysis}, 12:113--127, 1973.

\bibitem{JM87}
D.~Johnson and J.~Millson.
\newblock Deformation spaces associated to compact hyperbolic manifolds.
\newblock In {\em Discrete groups in geometry and analysis ({N}ew {H}aven,
  {C}onn., 1984)}, volume~67 of {\em Progr. Math.}, pages 48--106. Birkh\"auser
  Boston, Boston, MA, 1987.

\bibitem{Vinberg67}
V.~Kac and {\`E}.~Vinberg.
\newblock Quasi-homogeneous cones.
\newblock {\em Mat. Zametki}, 1:347--354, 1967.

\bibitem{Kapovich09}
M.~Kapovich.
\newblock {\em Hyperbolic manifolds and discrete groups}.
\newblock Modern Birkh\"auser Classics. Birkh\"auser Boston, Inc., Boston, MA,
  2009.
\newblock Reprint of the 2001 edition.

\bibitem{KP2002}
F.~Kassel and R.~Potrie.
\newblock Eigenvalue gaps for hyperbolic groups and semigroups.
\newblock {\em J. Mod. Dyn.}, 18:161--208, 2022.
\newblock arXiv:2002.07015.

\bibitem{KH95}
A.~Katok and B.~Hasselblatt.
\newblock {\em Introduction to the modern theory of dynamical systems},
  volume~54 of {\em Encyclopedia of Mathematics and its Applications}.
\newblock Cambridge University Press, Cambridge, 1995.
\newblock With a supplementary chapter by Katok and Leonardo Mendoza.

\bibitem{Kimi01}
I.~Kim.
\newblock Rigidity and deformation spaces of strictly convex real projective
  structures on compact manifolds.
\newblock {\em J. Differential Geom.}, 58(2):189--218, 2001.

\bibitem{Knapp}
A.~Knapp.
\newblock {\em Lie groups beyond an introduction}, volume 140 of {\em Progress
  in Mathematics}.
\newblock Birkh\"auser Boston, Inc., Boston, MA, second edition, 2002.

\bibitem{Kobpaper}
S.~Kobayashi.
\newblock Projectively invariant distances for affine and projective
  structures.
\newblock In {\em Differential geometry ({W}arsaw, 1979)}, volume~12 of {\em
  Banach Center Publ.}, pages 127--152. PWN, Warsaw, 1984.

\bibitem{Kostant73}
B.~Kostant.
\newblock On convexity, the {W}eyl group and the {I}wasawa decomposition.
\newblock {\em Ann. Sci. \'Ecole Norm. Sup. (4)}, 6:413--455 (1974), 1973.

\bibitem{KS75}
B.~Kostant and D.~Sullivan.
\newblock The {E}uler characteristic of an affine space form is zero.
\newblock {\em Bull. Amer. Math. Soc.}, 81(5):937--938, 1975.

\bibitem{Koszul65}
J.-L. Koszul.
\newblock Vari\'et\'es localement plates et convexit\'e.
\newblock {\em Osaka J. Math.}, 2:285--290, 1965.

\bibitem{Koszul68}
J.-L. Koszul.
\newblock D\'eformations de connexions localement plates.
\newblock {\em Ann. Inst. Fourier (Grenoble)}, 18(fasc. 1):103--114, 1968.

\bibitem{Koszul70}
J.-L. Koszul.
\newblock Trajectoires convexes de groupes affines unimodulaires.
\newblock In {\em Essays on {T}opology and {R}elated {T}opics ({M}\'emoires
  d\'edi\'es \`a {G}eorges de {R}ham)}, pages 105--110. Springer, New York,
  1970.

\bibitem{Kuiper53}
N.~H. Kuiper.
\newblock On convex locally-projective spaces.
\newblock In {\em Convegno {I}nternazionale di {G}eometria {D}ifferenziale,
  {I}talia, 1953}, pages 200--213. Edizioni Cremonese, Roma, 1954.

\bibitem{Labourie07}
F.~Labourie.
\newblock Flat projective structures on surfaces and cubic holomorphic
  differentials.
\newblock {\em Pure Appl. Math. Q.}, 3(4, part 1):1057--1099, 2007.

\bibitem{Leitner162}
A.~Leitner.
\newblock A classification of subgroups of {$SL(4,\mathbb R)$} isomorphic to
  {$\mathbb R^3$} and generalized cusps in projective 3 manifolds.
\newblock {\em Topology Appl.}, 206:241--254, 2016.

\bibitem{Leitner161}
A.~Leitner.
\newblock Conjugacy limits of the diagonal {C}artan subgroup in
  {$SL_3(\mathbb{R})$}.
\newblock {\em Geom. Dedicata}, 180:135--149, 2016.

\bibitem{Leitner163}
A.~Leitner.
\newblock Limits under conjugacy of the diagonal subgroup in
  {$SL_n(\mathbb{R})$}.
\newblock {\em Proc. Amer. Math. Soc.}, 144(8):3243--3254, 2016.

\bibitem{Lok}
W.~Lok.
\newblock {\em Deformations of locally homogeneous spaces and Kleinian groups}.
\newblock PhD thesis, Columbia University, New York, New York, 1984.
\newblock pp. 1--178.

\bibitem{Malcev}
A.~Malcev.
\newblock On isomorphic matrix representations of infinite groups.
\newblock {\em Rec. Math. [Mat. Sbornik] N.S.}, 8 (50):405--422, 1940.

\bibitem{Malcev49}
A.~I. Malcev.
\newblock On a class of homogeneous spaces.
\newblock {\em Izvestiya Akad. Nauk. SSSR. Ser. Mat.}, 13:9--32, 1949.

\bibitem{Marden07}
A.~Marden.
\newblock {\em Outer circles}.
\newblock Cambridge University Press, Cambridge, 2007.

\bibitem{Marquis10}
L.~Marquis.
\newblock Espace des modules de certains poly\`edres projectifs miroirs.
\newblock {\em Geom. Dedicata}, 147:47--86, 2010.

\bibitem{Marquis12}
L.~Marquis.
\newblock Surface projective convexe de volume fini.
\newblock {\em Ann. Inst. Fourier (Grenoble)}, 62(1):325--392, 2012.

\bibitem{Marquis17}
L.~Marquis.
\newblock Coxeter group in {H}ilbert geometry.
\newblock {\em Groups Geom. Dyn.}, 11(3):819--877, 2017.

\bibitem{Mazur}
B.~Mazur.
\newblock Relative neighborhoods and the theorems of {S}male.
\newblock {\em Ann. of Math. (2)}, 77:232--249, 1963.

\bibitem{Mess07}
G.~Mess.
\newblock Lorentz spacetimes of constant curvature.
\newblock {\em Geom. Dedicata}, 126:3--45, 2007.

\bibitem{Milnor}
J.~Milnor.
\newblock Whitehead torsion.
\newblock {\em Bull. Amer. Math. Soc.}, 72:358--426, 1966.

\bibitem{Milnor2}
J.~Milnor.
\newblock A note on curvature and fundamental group.
\newblock {\em J. Differential Geometry}, 2:1--7, 1968.

\bibitem{Molino88}
P.~Molino.
\newblock {\em Riemannian foliations}, volume~73 of {\em Progress in
  Mathematics}.
\newblock Birkh\"auser Boston, Inc., Boston, MA, 1988.
\newblock Translated from the French by Grant Cairns, With appendices by
  Cairns, Y. Carri{\`e}re, {\'E}. Ghys, E. Salem and V. Sergiescu.

\bibitem{Moore68}
C.~Moore.
\newblock Distal affine transformation groups.
\newblock {\em Amer. J. Math.}, 90:733--751, 1968.

\bibitem{Morgan79}
J.~Morgan.
\newblock On {T}hurston's uniformization theorem for three-dimensional
  manifolds.
\newblock In {\em The {S}mith conjecture ({N}ew {Y}ork, 1979)}, volume 112 of
  {\em Pure Appl. Math.}, pages 37--125. Academic Press, Orlando, FL, 1984.

\bibitem{MS88}
J.~Morgan and P.~Shalen.
\newblock Degenerations of hyperbolic structures. {III}. {A}ctions of
  {$3$}-manifold groups on trees and {T}hurston's compactness theorem.
\newblock {\em Ann. of Math. (2)}, 127(3):457--519, 1988.

\bibitem{Morris15}
D.~Morris.
\newblock {\em Introduction to arithmetic groups}.
\newblock Deductive Press, 2015.

\bibitem{Munkres75}
J.~Munkres.
\newblock {\em Topology: a first course}.
\newblock Prentice-Hall, Inc., Englewood Cliffs, N.J., 1975.

\bibitem{NY}
T.~Nagano and K.~Yagi.
\newblock The affine structures on the real two-torus. {I}.
\newblock {\em Osaka J. Math.}, 11:181--210, 1974.

\bibitem{NieSeppi23}
X.~Nie and A.~Seppi.
\newblock Affine deformations of quasi-divisible convex cones.
\newblock {\em Proc. Lond. Math. Soc. (3)}, 127(1):35--83, 2023.

\bibitem{Porti20}
J.~Porti.
\newblock Dimension of representation and character varieties for two- and
  three-orbifolds.
\newblock {\em Algebraic \& Geometric Topology}, 22:1905--1967, 2020.

\bibitem{PTp}
J.~Porti and S.~Tillman.
\newblock Projective structures on a hyperbolic 3--orbifold.
\newblock {\em Acta Mathematica Vietnamica}, 46(2):347--355, 2021.

\bibitem{Roeder07}
R.~Roeder, J.~Hubbard, and W.~Dunbar.
\newblock Andreev's theorem on hyperbolic polyhedra.
\newblock {\em Ann. Inst. Fourier (Grenoble)}, 57(3):825--882, 2007.

\bibitem{Rousset}
M.~Rousset.
\newblock Sur la rigidit\'e{} de poly\`edres hyperboliques en dimension 3: cas
  de volume fini, cas hyperid\'eal cas fuchsien.
\newblock {\em Bull. Soc. Math. France}, 132(2):233--261, 2004.

\bibitem{RS}
D.~Ruelle and D.~Sullivan.
\newblock Currents, flows and diffeomorphisms.
\newblock {\em Topology}, 14(4):319--327, 1975.

\bibitem{Schlenker}
J.-M. Schlenker.
\newblock Des immersions isom\'etriques de surfaces aux vari\'et\'es
  hyperboliques \`a{} bord convexe.
\newblock In {\em S\'eminaire de {T}h\'eorie {S}pectrale et {G}\'eom\'etrie.
  {V}ol. 21. {A}nn\'ee 2002--2003}, volume~21 of {\em S\'emin. Th\'eor. Spectr.
  G\'eom.}, pages 165--216. Univ. Grenoble I, Saint-Martin-d'H\`eres, 2003.

\bibitem{Selberg}
A.~Selberg.
\newblock On discontinuous groups in higher-dimensional symmetric spaces.
\newblock In {\em Contributions to function theory ({I}nternat. {C}olloq.
  {F}unction {T}heory, {B}ombay, 1960)}, pages 147--164. Tata Institute of
  Fundamental Research, Bombay, 1960.

\bibitem{Shima07}
H.~Shima.
\newblock {\em The geometry of {H}essian structures}.
\newblock World Scientific Publishing Co. Pte. Ltd., Hackensack, NJ, 2007.

\bibitem{Smale59}
S.~Smale.
\newblock The classification of immersions of spheres in euclidean spaces.
\newblock {\em Ann. Math.}, 69, 1959.

\bibitem{Smillie}
J.~Smillie.
\newblock Affine manifolds with diagonal holonmy.
\newblock preprint, 1977.

\bibitem{SmillieTh}
J.~Smillie.
\newblock {\em Affinely flat manifolds}.
\newblock PhD thesis, University of Chicago, Chicago, USA, 1977.

\bibitem{Sullivan82}
D.~Sullivan.
\newblock Discrete conformal groups and measurable dynamics.
\newblock {\em Bull. Amer. Math. Soc. (N.S.)}, 6(1):57--73, 1982.

\bibitem{ST83}
D.~Sullivan and W.~Thurston.
\newblock Manifolds with canonical coordinate charts: some examples.
\newblock {\em Enseign. Math. (2)}, 29(1-2):15--25, 1983.

\bibitem{Thurston80}
W.~Thurston.
\newblock {\em The geometry and topology of three-manifolds}.
\newblock Lecture notes. Princeton University, 1980.
\newblock \, Available at https://library.slmath.org/nonmsri/gt3m/.

\bibitem{Thurston97}
W.~Thurston.
\newblock {\em Three-dimensional geometry and topology. {V}ol. 1}, volume~35 of
  {\em Princeton Mathematical Series}.
\newblock Princeton University Press, Princeton, NJ, 1997.
\newblock Edited by Silvio Levy.

\bibitem{Varadarajan84}
V.~S. Varadarajan.
\newblock {\em Lie groups, {L}ie algebras, and their representations}, volume
  102 of {\em Graduate Texts in Mathematics}.
\newblock Springer-Verlag, New York, 1984.
\newblock Reprint of the 1974 edition.

\bibitem{Vey}
J.~Vey.
\newblock Sur les automorphismes affines des ouverts convexes saillants.
\newblock {\em Ann. Scuola Norm. Sup. Pisa (3)}, 24:641--665, 1970.

\bibitem{Vinberg63}
\`E. Vinberg.
\newblock The theory of homogeneous convex cones.
\newblock {\em Trudy Moskov. Mat. Ob\v s\v c.}, 12:303--358, 1963.

\bibitem{Vinberg71}
{\`E}.~Vinberg.
\newblock Discrete linear groups that are generated by reflections.
\newblock {\em Izv. Akad. Nauk SSSR Ser. Mat.}, 35:1072--1112, 1971.

\bibitem{Svarc}
A.~\v{S}varc.
\newblock A volume invariant of coverings.
\newblock {\em Dokl. Akad. Nauk SSSR (N.S.)}, 105:32--34, 1955.

\bibitem{Weil62}
A.~Weil.
\newblock On discrete subgroups of {L}ie groups. {II}.
\newblock {\em Ann. of Math. (2)}, 75:578--602, 1962.

\bibitem{Weisman}
T.~Weisman.
\newblock Dynamical properties of convex cocompact actions in projective space.
\newblock {\em J. Topology}, 16, 2023.

\bibitem{Witte95}
D.~Witte.
\newblock Superrigidity of lattices in solvable {L}ie groups.
\newblock {\em Invent. Math.}, 122(1):147--193, 1995.

\bibitem{Yaman04}
A.~Yaman.
\newblock A topological characterisation of relatively hyperbolic groups.
\newblock {\em J. Reine Angew. Math.}, 566:41--89, 2004.

\bibitem{Zassenhaus}
H.~Zassenhaus.
\newblock Beweis eines satzes \"uber diskrete gruppen.
\newblock {\em Abh. Math. Sem. Univ. Hamburg}, 12(1):289--312, 1937.

\bibitem{Zimmer23}
A.~Zimmer.
\newblock A higher-rank rigidity theorem for convex real projective manifolds.
\newblock {\em Geom. Topol.}, 23(7):2899--2936, 2023.

\end{thebibliography}


\printindex



\end{document}